\documentclass[a4paper, reqno, 12pt, oneside]{book}

 \usepackage{color}
 \definecolor{darkblue}{rgb}{0,0,0.75}
 \definecolor{darkgreen}{rgb}{0,0.75,0}
 \usepackage[
  colorlinks, linkcolor=darkblue,
    citecolor=darkgreen, urlcolor=blue, bookmarks=false, breaklinks]{hyperref}

\usepackage{cmap} 

\usepackage[T2A]{fontenc}
\usepackage[cp1251]{inputenc} 
\usepackage[russian]{babel}

\usepackage{lscape} 

\usepackage{indentfirst}
\usepackage[left=2.5cm, right=1cm, top=2cm, bottom=2cm]{geometry}
\usepackage{fancyhdr}
\fancypagestyle{plain}{
\fancyhf{}
\chead{\thepage}
}
\usepackage[ddmmyyyy]{datetime}

\usepackage{titlesec}
\titleformat{\chapter}[display]{\normalfont\huge\bfseries\centering}{\chaptertitlename\ \thechapter}{20pt}{\Huge}
\titleformat{\section}[block]{\bfseries\large\centering}{\thesection}{1ex}{}
\titleformat{\subsection}[block]{\bfseries\normalsize\centering}{\thesubsection}{1ex}{}
\usepackage{amssymb, amsmath, amsthm, longtable, bbm}
\usepackage[matrix, arrow, curve]{xy}

\usepackage[geometry]{ifsym}
\renewcommand{\square}{\text{\SmallSquare}}

\newtheoremstyle{mytheoremstyle} 
    {\topsep}                    
    {\topsep}                    
    {\itshape}                   
    {5ex}                           
    {\bfseries}                   
    {.}                          
    {.5em}                       
    {}  

\newtheoremstyle{myremarkstyle} 
    {\topsep}                    
    {\topsep}                    
    {}                   
    {5ex}                           
    {\bfseries}                   
    {.}                          
    {.5em}                       
    {}  

\theoremstyle{mytheoremstyle}
\newtheorem{theorem}{Теорема}[chapter]
\newtheorem{proposition}[theorem]{Предложение}
\newtheorem{corollary}[theorem]{Следствие}
\newtheorem{lemma}[theorem]{Лемма}
\newtheorem*{conjecture}{Гипотеза}

\theoremstyle{myremarkstyle}
\newtheorem{remark}[theorem]{Замечание}

\newtheorem{example}[theorem]{Пример}

\newtheorem*{question}{Вопрос}
\newtheorem{definition}[theorem]{Определение}
\newtheorem{problem}[theorem]{Задача}

\numberwithin{equation}{chapter}

\makeatletter
\renewenvironment{proof}[1][\proofname]{\par\indent {\bfseries #1\@addpunct{.} }}{\qed}
\makeatother

\DeclareMathOperator{\Id}{Id}
\DeclareMathOperator{\Der}{Der}
\DeclareMathOperator{\id}{id}
\DeclareMathOperator{\chr}{char}
\DeclareMathOperator{\im}{im}
\DeclareMathOperator{\ad}{ad}
\DeclareMathOperator{\tr}{tr}
\DeclareMathOperator{\GL}{GL}
\DeclareMathOperator{\UT}{UT}

\DeclareMathOperator{\diff}{diff}
\DeclareMathOperator{\supp}{supp}
\DeclareMathOperator{\cosupp}{cosupp}
\DeclareMathOperator{\PGL}{PGL}

\DeclareMathOperator{\height}{ht}
\DeclareMathOperator{\End}{End}
\DeclareMathOperator{\Aut}{Aut}
\DeclareMathOperator{\diag}{diag}

\DeclareMathOperator{\sign}{sign}
\DeclareMathOperator{\Alt}{Alt}
\DeclareMathOperator{\Ann}{Ann}
\DeclareMathOperator{\Hom}{Hom}
\DeclareMathOperator{\op}{op}
\DeclareMathOperator{\PIexp}{PIexp}
\newcommand{\hatotimes}{\mathbin{\widehat{\otimes}}}

\DeclareRobustCommand{\No}{\ifmmode{\nfss@text{\textnumero}}\else\textnumero\fi}

\begin{document}

{\fontsize{14pt}{18} \selectfont

\thispagestyle{empty}
\begin{center}

\vfill\vfill \ \\ {
 МОСКОВСКИЙ ГОСУДАРСТВЕННЫЙ УНИВЕРСИТЕТ \\
имени М.В.~ЛОМОНОСОВА

МЕХАНИКО-МАТЕМАТИЧЕСКИЙ ФАКУЛЬТЕТ

}

\vfill \hfill \textit{На правах рукописи}


\bigskip

\bigskip

\vfill {\Large
 \textbf{Гордиенко Алексей Сергеевич}}

\bigskip

\bigskip

\bigskip

{\LARGE
\textbf{(Ко)модульные алгебры 
 и их обобщения} 
}

%
%
%
%
%
%
%


\vfill\vfill\vfill\vfill Москва~"--- 2020
\end{center}
}

\newpage

\pagestyle{plain}

\tableofcontents

\newpage

\chapter*{Предисловие}
\addcontentsline {toc}{chapter}{Предисловие}

Данная рукопись является расширенной версией диссертации на соискание учёной степени доктора физико-математических наук, подготовленной автором в 2019--2020 годах, что и объясняет подбор материала на основе результатов автора. Автор надеется, что рукопись окажется полезной всем, кто решил познакомиться с данной темой на русском языке. От окончательного варианта диссертации эта работа отличается, в основном, наличием глав~\ref{ChapterGradEquiv}, \ref{ChapterOmegaAlg} и~\ref{ChapterHequiv}, включающих в себя результаты статей~\cite{ASGordienko20ALAgoreJVercruysse, ASGordienko21ALAgoreJVercruysse, ASGordienko18Schnabel, ASGordienko19Schnabel}, посвящённых эквивалентностям градуировок и (ко)модульных структур на алгебрах и связанным с ними категорным вопросам.

\medskip

  Автор благодарит своего научного консультанта,
  доктора физико-математических наук,
  профессора
  Михаила Владимировича Зайцева за поддержку и внимание к результатам.
 
  Автор выражает признательность
  Ю.\,А.~Бахтурину, Эрику Йесперсу, М.\,В.~Кочетову, Ане Агоре, Йоосту Веркрёйссе, Йорку Зоммерхойзеру, А.\,А.~Клячко, В.\,К.~Харченко и Офиру Шнабелю  за плодотворное обсуждение результатов данной работы и конструктивные замечания.
  Автор благодарен участникам семинара <<Избранные вопросы алгебры>> кафедры высшей алгебры механико-математического факультета МГУ за
   обсуждение результатов диссертации.
  Автор признателен сотрудникам математических кафедр университетов Memorial University of Newfoundland (Сент-Джонс, Канада) и Vrije Universiteit Brussel (Брюссель, Бельгия), на которых автору посчастливилось заниматься исследованиями в 2010--2018 годах в рамках грантов AARMS и Fonds Wetenschappelijk Onderzoek~"--- Vlaanderen, за
   обсуждение результатов диссертации и творческую атмосферу,
   которая способствовала научной работе. Наконец, автор выражает благодарность коллективу кафедры высшей математики МГТУ ГА, на которой он работает с февраля 2019 года, за дружескую атмосферу и возможность завершить работу над диссертацией.
   
   Автор признателен своей семье за помощь и поддержку при написании данной работы.
   
   \bigskip

   Зеленоград, декабрь 2020 года. \hfill Гордиенко А.\,С.

\newpage

\chapter*{Введение}
\addcontentsline {toc}{chapter}{Введение}

Во многих областях математики и физики (см., например, \cite{ArnoldBook, ModernGeometry, MurphyBook, HaagKastler}) находят своё применение \textit{алгебры}, то есть векторные пространства над некоторым полем $\mathbbm{k}$ (например, $\mathbbm{k}$ может быть полем $\mathbb{R}$ вещественных или $\mathbb{C}$ комплексных чисел), в которых задана бинарная операция внутреннего умноже­ния, линейная по каждому аргументу. 

Часто алгебры, встречающиеся в приложениях, наделены некоторой дополнительной структурой или (обобщёнными) симметриями: действием (полу)группы эндоморфизмами и антиэндоморфизмами, (полу)групповой градуировкой или действием алгеб­ры Ли дифференцированиями (см., например, \cite{PolyakovBook, HaagBook, KakuBook, MajidBook}). Для работы с такими дополнительными структурами оказывается полезным понятие модульной и комодульной алгебры над алгеброй Хопфа и даже более общее понятие алгебры с обобщённым $H$-действием. В частности, понятие (ко)модульной алгебры позволяет изучать различные виды дополнительных структур на алгебрах одновременно.

Кроме того, (ко)модульные алгебры естественным образом возникают в геометрии: если некоторая аффинная алгебраическая группа действует морфизмами на аффинном алгебраи­ческом многообразии (например, рассматриваются симметрии некоторой поверхности, задан­ной алгебраическими уравнениями), то алгебра регулярных функций (множество функций, которые можно определить при помощи многочленов от координат точки,
с операциями сложения, умножения между собой и на скаляры) будет модульной алгеброй над групповой алгеброй этой группы и над универсальной обёртывающей алгебры Ли этой группы и комодульной алгеброй над алгеброй регулярных функций на аффинной алгебраической группе~\cite{Abe}.

Обратим внимание на то, что в классическом случае алгебры регулярных функций на многообразиях коммутативны, так как для умножения функций выполнен перестановочный закон. Однако в новом направлении, которое получило название некоммутативной геометрии, рассмат­риваются <<некоммутативные пространства>>, т. е. такие пространства, алгебры регулярных функций которых некоммутативны. Поэтому изучение (ко)действий необязательно коммутатив­ных алгебр Хопфа на необязательно коммутативных алгебрах можно интерпретировать как изучение кван­товых симметрий некоммутативных пространств. Последние находят своё применение в теорети­ческой физике (см., например, \cite{ConnesMarcolli, Donatsos}).

Одним из важнейших примеров комодульных алгебр являются алгебры, градуированные группами и полугруппами.
В 2002--2008 гг. Ю.\,А. Бахтуриным, М.\,В. Зайцевым и С.\,К. Сегалом были классифицированы групповые градуировки на матричных алгебрах и градуированно простые алгебры, градуированные группами~\cite{BahturinZaicevSegal, BahZaicAllGradings}.
В 2013 году вышла монография А.~Эльдуке и М.\,В.~Кочетова, подводящая итог исследованиям
групповых градуировок на простых алгебрах Ли~\cite{ElduqueKochetov}.
Многими авторами также исследовались кольца, градуированные полугруппами.
В частности, Э.~Йесперсом, А.\,В.~Келаревым и М.~Клэйзом~\cite{Clase,KelarevPI, KelarevBook}
изучались градуированные полугруппами PI-кольца, а
 Ю.\,А.~Бахтуриным, М.\,В.~Зайцевым, Э.~Йесперсом, П.~Нистедтом, Й.~Ойнертом и П.~Уаутерсом~\cite{BahturinZaicevSem, oinert3, jespers3}
 были получены результаты, касающиеся простых и градуированно простых колец и алгебр, градуированных полугруппами.
 Как показано в \S\ref{SectionGeneralReductionSemigroup}
данной работы, при изучении градуированно простых ассоциативных колец или алгебр, градуированных
полугруппами, достаточно ограничиться градуировками $0$-простыми полугруппами,
т.е. полугруппами с $0$, не содержащими нетривиальных идеалов. Для конечных 
$0$-простых полугрупп $T$ справедлива структурная теорема Риса~\cite[теорема 3.5]{clipre},
которая утверждает, что полугруппа $T$ изоморфна 
полугруппе Риса матричного типа
над группой с нулём $G^0$
для некоторой максимальной подгруппы $G$ полугруппы $T$.
При этом двумя диаметрально противоположными важнейшими частными случаями
оказываются случай, когда $T=G^0$, т.е. кольцо или алгебра градуированы конечной группой
(этот случай как раз и был исследован Ю.\,А. Бахтуриным, М.\,В. Зайцевым и С.\,К. Сегалом в~\cite{BahturinZaicevSegal}), и случай, когда все максимальные подгруппы
градуирующей полугруппы тривиальны. Второму случаю посвящены \S\ref{SectionGeneralReductionSemigroup}--\ref{SectionTGradedReesExistence} настоящей работы.

Оказывается, что для многих приложений, например, при изучении градуированных подпространств и идеалов, градуированных тождеств и т.д., является несущественным, какая именно группа градуирует алгебру, т.е. реализует заданную градуировку, и компоненты градуировки можно переиндексировать элементами другой группы. Получившаяся градуировка называется \textit{эквивалентной}\footnote{или \textit{слабо эквивалентной}~\cite{ASGordienko18Schnabel}} \cite{ElduqueKochetov} исходной градуировке. 
Здесь следует обратить внимание на то, что группы, реализующие одну и ту же градуировку, могут очень сильно отличаться друг от друга: быть конечными и бесконечными, абелевыми и неабелевыми и т.д.
В силу того, что для ряда задач (см., например, \cite{GiaLa,ASGordienko3}) возможность заменить градуирующую группу на конечную или, например, на абелеву, может существенно упростить решение, представляет интерес вопрос о том, когда такая замена возможна. Первый пример конечномерной алгебры, градуированной бесконечной группой, которая не может быть заменена на конечную, был построен М.~Клэйзом, Э.~Йесперсом и А.~Дель~Рио~\cite{ClaseJespersDelRio} в 1996 году (см. также~\cite{DasNasDelRioVanOyst}). Алгебра в этом примере имела нетривиальный радикал Джекобсона, поэтому вопрос о том, существует ли конечномерная полупростая алгебра, градуированная бесконечной группой, которая не может быть заменена на конечную, оставался открытым.

Оказывается, что среди среди всех групп, реализующих заданную градуировку, существует группа, являющаяся по отношению к остальным группам универсальной. Впервые понятие универсальной группы градуировки было введено
И.~Патерой и Г.~Цассенхаузом~\cite{PZ89} в 1989 году. В силу того, что универсальная группа может быть достаточно прозрачно задана порождающими и определяющими соотношениями, это позволяет переформулировать проблемы, связанные с переградуированием градуировок, на теоретико-групповом языке.
В главе~\ref{ChapterGradEquiv} даётся критерий возможности реализации заданной градуировки, а также её огрублений градуировкой конечной группой в терминах универсальной группы градуировки. В теореме~\ref{TheoremGivenFinPresGroupExistence} показано, что любая конечно представленная группа может быть реализована в качестве универсальной группы некоторой элементарной градуировки на матричной алгебре. Как следствие, строится такая элементарная градуировка на матричной алгебре, которую невозможно реализовать конечной группой.  Кроме того, в главе~\ref{ChapterGradEquiv} понятие эквивалентности градуировок переформулируется
в терминах линейных операторов, что позволяет ввести понятие эквивалентности действий групп.
Оставшаяся часть главы~\ref{ChapterGradEquiv} посвящена изучению категорий и функторов, которые естественным образом возникают при исследовании эквивалентных градуировок. 

Главы~\ref{ChapterH(co)modAssoc} и~\ref{ChapterH(co)modLie} посвящены структурным вопросам
теории (ко)модульных ассоциативных алгебр и алгебр Ли над алгебрами Хопфа. Центральную роль здесь играет вопрос о (ко)инвариантности радикалов и существовании (ко)инвариантных разложений Веддербёрна~"--- Артина, Веддербёрна~"--- Мальцева, Леви и разложения полупростой $H$-(ко)модульной алгебры Ли в прямую сумму (ко)инвариантных идеалов, являющихся $H$-простыми алгебрами.  Достаточные условия инвариантности радикала Джекобсона в ассоциативных алгебрах были получены В.\,В.~Линченко, С.~Монтгомери и Л.~Смоллом~\cite{LinchenkoJH, LinchenkoMontgomery, LinchenkoMontgomerySmall}.
 Проблема градуированности радикала Джекобсона
 в алгебрах, градуированных группами и полугруппами, изучалась, в частности, Дж.~Бергманом, Э.~Йесперсом,  А.\,В.~Келаревым, М.~Коэн, С.~Монтгомери, Я.~Окнинским, Э.~Пучиловским~\cite{Bergman, CohenMontgomery, JespersPuczylowski, KelarevOkninski} и другими авторами.
 Достаточные условия существования инвариантных разложений
 Веддербёрна~"--- Мальцева и Леви в ассоциативных алгебрах и алгебрах Ли
 были получены Э.~Тафтом~\cite{Taft}.
Достаточные условия существования $H$-(ко)инвариантного разложения Веддербёрна~"--- Мальцева
для конечномерных ассоциативных $H$-(ко)модульных алгебр были получены
А.\,В.~Сидоровым~\cite{SidorovSplitting} и
Д.~Штефаном и Ф.~Ван Ойстайеном~\cite{SteVanOyst}.
$H$-инвариантный аналог теоремы Веддербёрна~"--- Артина был доказан Ф.~Ван Ойстайеном и С.\,М.~Скрябиным~\cite{SkryabinHdecomp, SkryabinVanOystaeyen}.
В случае, когда некоторая алгебра $A$ является $H$-комодульной для некоторой
бесконечномерной алгебры Хопфа,  на алгебре $A$ действует алгебра $H^*$, двойственная к коалгебре $H$,
 однако это действие не является действием алгебры Хопфа,
поскольку на $H^*$ не удаётся определить коумножение. В данной работе предлагается
метод, позволяющий устранить это затруднение в случае конечномерных алгебр $A$.
При помощи этого метода, в частности, доказываются достаточные условия
$H$-коинвариантности радикалов в $H$-комодульных ассоциативных алгебрах и алгебрах Ли.
Кроме того, в главе~\ref{ChapterH(co)modLie} доказываются
достаточные условия $H$-инвариантности радикалов в $H$-модульных алгебрах Ли,
существование $H$-(ко)инвариантного разложения Леви и (ко)инвариантные аналоги теоремы Вейля.
Наконец, в главах~\ref{ChapterH(co)modAssoc} и~\ref{ChapterH(co)modLie} классифицируются
конечномерные ассоциативные алгебры и алгебры Ли, простые по отношению к действию алгебр Тафта $H_{m^2}(\zeta)$, что является частью более обширного
исследования алгебр с действием расширений Оре.
Конечномерные ассоциативные $H_{m^2}(\zeta)$-модульные алгебры, которые не содержат ненулевых
нильпотентных элементов были ранее классифицированы в теореме 2.5 из работы~\cite{MontgomerySchneider}.
Точные модульные категории над категорией $\mathrm{Rep}(H_{m^2}(\zeta))$ рассматривались в теореме~4.10 из работы~\cite{EtingofOstrik}, а $H_{m^2}(\zeta)$-действия
на алгебрах путей колчанов исследовались в~\cite{KinserWalton}.
Действия точечных алгебр Хопфа на матричных алгебрах также изучались в~\cite{BahturinMontgomery}.

Оказывается, что если ослабить требования к (ко)действующей алгебре и разрешить ей быть просто биалгеброй, а не алгеброй Хопфа, то многие структурные результаты оказываются неверны. Так, например, в главе~\ref{ChapterGradSemigroup} рассматриваются ассоциативные алгебры, градуированные полугруппами из двух элементов (что эквивалентно кодействию соответствующей полугрупповой биалгебры),
радикал Джекобсона которых не является градуированным идеалом, и алгебры, для которых не существует градуированных разложений Веддербёрна~"--- Артина и Веддербёрна~"--- Мальцева. Кроме того, получаются некоторые достаточные условия существования таких градуированных разложений.

(Ко)модульные алгебры~"--- это область, бурно развивающаяся в настоящее время. Как уже было отмечено выше, с одной стороны, доказываются результаты о структуре таких алгебр. С другой стороны, в ряде случаев показано, что действие одних алгебр Хопфа сводится к действию других алгебр Хопфа~\cite{CuadraEtingofWalton, EtingofWalton}. Как и в случае градуировок, для многих приложений (структурная теория алгебр, полиномиальные тождества, ...) бывает несущественно, какая конкретно алгебра Хопфа (ко)действует на данной алгебре, и зачастую эту алгебру Хопфа удаётся заменить другой (иногда более просто устроенной), (ко)действие которой эквивалентно (ко)действию первоначальной алгебры Хопфа. Понятие эквивалентности (ко)действий  является естественным обобщением понятия эквивалентности градуировок. Как и в случае групповых градуировок на алгебрах, где среди всех групп, реализующих данную градуировку, существует универсальная, для заданного (ко)действия существует алгебра Хопфа, являющаяся универсальной среди всех остальных алгебр Хопфа, действие которых эквивалентно заданному. Существование такой алгебры Хопфа доказывается в \S\ref{SectionActions} и \S\ref{SectionCoactions}. 
При этом результаты П.\,И.~Этингофа, Х.~Куадры и Ч.~Уолтон~\cite{CuadraEtingofWalton, EtingofWalton} могут проинтерпретированы следующим образом: при определённых условиях на модульную структуру её универсальная алгебра Хопфа изоморфна групповой алгебре некоторой группы.

В теории (ко)модульных алгебр известны и другие конструкции универсальных (ко)действующих биалгебр и алгебр Хопфа. В своей книге, опубликованной в 1969 году, М.Е.~Свидлер строит универсальную измеряющую коалгебру
из алгебры $A$ в алгебру $B$, которая при совпадении алгебр $A$ и $B$ становится универсальной действующей биалгеброй.  
В своих работах Ю.И.~Манин~\cite{Manin} и Д.~Тамбара~\cite{Tambara} также рассматривали универсальные кодействующие биалгебры и алгебры Хопфа. Биалгебры и алгебры Хопфа Свидлера, Манина и Тамбары являются обобщением групп автоморфизмов и полугрупп эндоморфизмов и могут быть названы квантовыми (полу)группами квантовых симметрий. 
Если рассматривать только классические симметрии алгебры $A$ над полем $\mathbbm{k}$, т.е. действия групп автоморфизмами, то все такие действия пропускаются
через подгруппу автоморфизмов $\Aut(A)$, которая является подмножеством в множестве всех линейных операторов
$\End_\mathbbm{k}(A)$ на алгебре $A$. Если же теперь принять во внимание и квантовые симметрии, т.е. (ко)модульные структуры на $A$, то в алгебре $A$ не окажется внутренней квантовой группы таких симметрий, через которую любое такое (ко)действие пропускается. Поэтому универсальные алгебры Хопфа (когда они существуют) являются теми внешними объектами, который призваны стать квантовыми аналогами групп автоморфизмов.
Однако они универсальны среди всех кодействующих биалгебр и алгебр Хопфа, а не только среди тех, которые реализуют кодействие эквивалентное заданному, и их кодействие не обязано быть эквивалентным кодействию первоначальной алгебры Хопфа. 
Более того, если не накладывать на кодействие дополнительных ограничений вроде сохранения определённых
градуировок, как это делается, например, в~\cite{Manin}, существование универсальных кодействующих биалгебр и алгебр Хопфа доказано лишь
для конечномерных алгебр. В~\S\ref{SectionUnivCoactHopfNonExist} приводятся примеры бесконечномерных алгебр, для которых биалгебры Тамбары и алгебры Хопфа Манина не существуют. Напротив, универсальная алгебра Хопфа  заданной комодульной структуры всегда существует (см. \S\ref{SectionCoactions}).

В главе~\ref{ChapterOmegaAlg} строится единая теория, которая объединяет все вышеуказанные виды универсальных (ко)действующих биалгебр и алгебр Хопфа. В силу того, что в теории алгебр Хопфа находят своё применение не только (ко)модульные алгебры, но и (ко)модульные коалгебры~\cite{RadfordHopf},
представляет интерес изучение эквивалентности (ко)действий и универсальных (ко)действующих алгебр Хопфа
на коалгебрах. Для того, чтобы одновременно изучать (ко)действия на алгебрах и коалгебрах, а также включить случай заплетённых\footnote{braided vector spaces (англ.)} векторных пространств,
в работе рассматриваются так называемые \textit{$\Omega$-алгебры}, т.е. векторные пространства $A$ с фиксированным множеством $\Omega$ линейных отображений между тензорными степенями (см.~\S\ref{SectionOmegaAlgebras}).
Для того, чтобы включить случай универсальной алгебры Хопфа заданной (ко)модульной структуры на $\Omega$-алгебре $A$, вводится понятие $V$-универсальной алгебры Хопфа, где $V\subseteq \End_\mathbbm{k}(A)$.
$V$-универсальная действующая алгебра Хопфа всегда существует, а для $V$-универсальной кодействующей алгебры Хопфа в работе выявляются условия, гарантирующие её существование, см. теорему~\ref{TheoremBHsquareVExistence}.
Оказывается, что при этих условиях между $V$-универсальной действующей и кодействующей алгебрами Хопфа существует двойственность, являющаяся обобщением и уточнением двойственности между универсальной действующей и кодействующей биалгебрами для конечномерных алгебр, полученной Д.~Тамбарой~\cite{Tambara}.
Одно из возможных применений этой двойственности~"--- вычисление универсальных алгебр Хопфа в случаях, когда это не удаётся сделать другими методами.

В главе~\ref{ChapterHequiv} понятие эквивалентности рассматривается для важнейших классов (ко)модульных структур: действий кокоммутативных алгебр, рациональных действий связных аффинных алгебраических групп, групповых градуировок и расширений Хопфа--Галуа. Оказывается, например, что любая комодульная структура, эквивалентная некоторой групповой градуировке также сводится к некоторой градуировке, а для расширений Хопфа--Галуа универсальные алгебры Хопфа совпадают с исходными (ко)действующими алгебрами.

Одной из важных сторон исследования  алгебраических
систем  является  изучение тех тождеств, которые выполняются
в этих алгебраических системах.
<<Хотя тождества представляют собой простейшие
замкнутые высказывания логического языка, язык
тождеств все же достаточно богатый, чтобы на нём
можно было выражать многие тонкие свойства систем и их классов.>>
(А.И.~Мальцев~\cite[с.~337]{MalcevAIAlgSyst})
При исследовании тождеств в алгебрах естественным образом
возникают их числовые характеристики.
Одними из важнейших числовых характеристик являются
коразмерности.
Коразмерности оказываются  полезным инструментом
при решении различных задач, например, при доказательстве наличия 
или отсутствия нетривиальных тождеств~\cite{RegTensor,
RegExplicit}. Кроме того, коразмерности естественным образом возникают 
при вычислении базиса тождеств в алгебре над полем характеристики $0$.
 Коразмерности
служат своеобразной оценкой количества тождеств,
 которым удовлетворяет конкретная алгебра.
 Существует также интерпретация коразмерностей, напрямую не связанная с полиномиальными
 тождествами: $n$-я коразмерность~"--- это размерность пространства $n$-линейных
 функций на алгебре, представимых в виде (вообще говоря, некоммутативного или даже неассоциативного)
 многочлена от своих аргументов.

 Асимптотическое поведение
коразмерностей вызывает дополнительный интерес
 в связи c тем, что это поведение тесно связано со структурой
изучаемой алгебры~\cite{ZaiExp, ZaiGia}, на чём мы подробно остановимся в главах~\ref{ChapterGenHAssocCodim}--\ref{ChapterHModLieCodim}.

В 1984~году А.~Регев показал~\cite{RegMatrix},
 что коразмерности~$c_n(M_k(\mathbbm{k}))$
 полиномиальных тождеств алгебры~$M_k(\mathbbm{k})$ всех матриц~$k\times k$
 над произвольным полем~$\mathbbm{k}$ характеристики~$0$
 имеют следующую асимптотику (здесь и далее~$f \sim g$,
  если~$\lim \frac{f}{g}=1$):
 \begin{equation*} 
 c_n(M_k(\mathbbm{k})) \sim \alpha_k
   n^{-\frac{k^2-1}{2}}
   k^{2n}
 \text{ при } n \to \infty,
 \end{equation*}
 где $\alpha_k = \left(\frac{1}{\sqrt{2\pi}}\right)^{k-1}
 \left(\frac{1}{2}\right)^{(k^2-1)/2} \cdot 1! \cdot 2! \cdot\ldots
 \cdot (k-1)! \cdot k^{(k^2+4)/2}$, $k \in \mathbb N$ фиксировано.

 Основываясь на этом результате, Ш.~Амицур
 выдвинул следующую гипотезу:

\begin{conjecture}[Ш.~Амицур {\cite[гипотеза~6.1.3]{ZaiGia}}]
Пусть $A$~"--- ассоциативная PI-алгебра над полем характеристики~$0$,
а
$c_n(A)$~"--- последовательность коразмерностей ее полиномиальных
тождеств. Тогда существует \textit{PI-экспонента}
$\PIexp(A):=\lim\limits_{n \to \infty} \sqrt[n]{c_n(A)} \in \mathbb Z_+$.
\end{conjecture}

 Гипотеза Ш.~Амицура была
доказана М.\,В.~Зайцевым и А.~Джамбруно~\cite{ZaiExp}
 в~1999 году для всех ассоциативных алгебр.
 Наверное, самым замечательным в этом результате оказалось то,
 что в случае, когда алгебра $A$ конечномерна, а основное поле алгебраически замкнуто,
 для PI-экспоненты была получена явная формула, в которой участвовал радикал Джекобсона $J(A)$
 и простые компоненты разложения факторалгебры $A/J(A)$.
Кроме того, в~2002 году М.\,В.~Зайцев~\cite{ZaicevLie}
доказал аналог гипотезы Амицура для коразмерностей полиномиальных тождеств
 конечномерных алгебр Ли. В случае алгебр Ли также была получена формула для PI-экспоненты,
 однако эта формула оказалась сложнее, чем в ассоциативном случае, и включала в себя аннуляторы неприводимых
 факторов присоединённого представления алгебры Ли.
   В 2011 году А.~Джамбруно, М.\,В.~Зайцев и И.\,П.~Шестаков
доказали аналог гипотезы Амицура для конечномерных йордановых и альтернативных алгебр~\cite{GiaSheZai}.
В 2012 году А.~Джамбруно и М.\,В.~Зайцев
доказали существование  PI-экспоненты
для любой конечномерной простой необязательно ассоциативной алгебры~\cite[теорема 3]{ZaiGiaFinDimSuperAlgebras}.

В общем же случае аналог гипотезы Амицура для произвольных неассоциативных
и даже для бесконечномерных алгебр Ли может быть неверен.
Во-первых, рост коразмерностей может оказаться сверхэкспоненциальным~\cite{Volichenko}.
Во-вторых, экспонента этого роста может оказаться дробным числом~\cite{VerZaiMishch, GiaMishZaiCodimGrFun, ZaiMishchFracPI}. В-третьих, в 2014 году М.\,В.~Зайцев
построил пример бесконечномерной неассоциативной алгебры $A$,
для которой
$\mathop{\underline\lim}_{n\to\infty}\sqrt[n]{c_n(A)}=1$,
а $\mathop{\overline\lim}_{n\to\infty}\sqrt[n]{c_n(A)} > 1$~\cite{ZaicevPIexpDoesNotExist}.
Вопрос о том, существует ли такая алгебра Ли $L$, что $$\mathop{\underline\lim}_{n\to\infty}\sqrt[n]{c_n(L)}\ne\mathop{\overline\lim}_{n\to\infty}\sqrt[n]{c_n(L)},$$
по-прежнему остаётся открытым.
 
 В связи с результатами, перечисленными выше, представляет интерес дальнейшее изучение связи между
 строением алгебр и асимптотическим поведением их полиномиальных тождеств, особенно в случае алгебр, наделённых какой-то дополнительной структурой.
 При изучении алгебр с дополнительной структурой естественно ввести эту дополнительную структуру в сигнатуру полиномиальных тождеств и рассматривать градуированные, $G$-, дифференциальные и $H$-тождества~\cite{Kharchenko, BahtGiaZai, BahtDrensky}. Такой подход оказывается весьма плодотворным, например, при изучении градуированных алгебр:
 в 2012 году Э.~Альхадефф и О.~Давид~\cite{AljaDavidOnRegGGrad}, используя градуированные тождества,
 показали, что для всякой алгебры $A$ над алгебраически замкнутым полем, градуированной конечной
 группой, при условии, что градуировка минимальна и регулярна, порядок градуирующей группы является фиксированным числом и совпадает с экспонентой роста обычных тождеств алгебры $A$.
 Кроме того, в \S\ref{SectionAssocExamples} и \S\ref{SectionLieExamples} данной работы
 формулируются критерии градуированной, $G$- и $H$-простоты ассоциативных алгебр и алгебр Ли в терминах
 соответствующих коразмерностей.
 
 
Для того, чтобы одновременно изучать разные типы дополнительных структур на алгебрах,
оказывается удобным рассматривать так называемые обобщённые $H$-действия,
определение которых даётся в главе~\ref{ChapterGenHActions}.
Скорее всего, первым, кто начал рассматривать такие действия и соответствующие
полиномиальные $H$-тождества, был Аллан Берел~\cite[замечание после теоремы~15]{BereleHopf} в 1996 году.
Обобщённые $H$-действия особенно полезны тем, что к ним сводятся многие структуры, которые, вообще говоря, не являются модульными структурами над алгебрами Хопфа, например, действия групп не только автоморфизмами, но и антиавтоморфизмами, градуировки бесконечными (полу)группами на конечномерных алгебрах,
градуированные действия групп, суперинволюции и псевдоинволюции (см. главу~\ref{ChapterGenHActions}).

Обозначим через $\bigl(c_n^H(A)\bigr)_{n=1}^{\infty}$
последовательность коразмерностей полиномиальных $H$-тождеств в алгебре $A$ с обобщённым $H$-действием.
Следующий вопрос возникает естественным образом:

\begin{question}
При каких условиях на $H$-действие справедлив аналог гипотезы Амицура,
т.е. существует предел $$\PIexp^H(A):=\lim_{n\to\infty}\sqrt[n]{c_n(A)},$$
который является целым числом? Какие условия необходимо для этого наложить на структуру алгебры $A$ ?
\end{question}

Здесь нужно отметить, что при такой постановке вопроса сам аналог гипотезы Амицура становится
инструментом исследования и понимания дополнительных структур на алгебрах.

Ответу на вопрос, поставленный выше, посвящены главы~\ref{ChapterFreePICodim}--\ref{ChapterSGGrAssocCodim}, причём в указанных главах интенсивно применяется структурная теория, развитая
в главах~\ref{ChapterGradEquiv}--\ref{ChapterGenHActions}.

Как показано в примере~\ref{ExampleInfCodim},
в случае, когда обе алгебры $A$ и $H$ бесконечномерны,
сами коразмерности $c_n^H(A)$ могут оказаться бесконечными.
Поэтому при изучении $H$-коразмерностей имеет смысл ограничиться
случаем, когда одна из двух алгебр $A$ и $H$ конечномерна.

Аналог гипотезы Амицура для $*$-коразмерностей конечномерных ассоциативных 
алгебр с инволюцией был доказан в 1999 году М.\,В.~Зайцевым и А.~Джамбруно~\cite{ZaiGiaInvolution}.
Данный результат был обобщён в 2017 году А.~Джамбруно, С.~Полсино Милисом и А.~Валенти
на случай бесконечномерных PI-алгебр с инволюцией~\cite{GiaPolMilVal}.

Аналог гипотезы Амицура для $\mathbb Z/2\mathbb Z$-градуированных коразмерностей
конечномерных ассоциативных супералгебр был доказан в 2003 году 
Ф.~Бенанти, А.~Джамбруно и М.\,Пипитоне~\cite{BenantiGiambrunoPipitone}.
В случае ассоциативных алгебр, градуированных произвольной конечной группой,
аналог гипотезы Амицура для градуированных полиномиальных тождеств
был доказан в 2010--2011 гг. Э.~Альхадеффом,  А.~Джамбруно и Д.~Ла~Маттиной~\cite{AljaGia, AljaGiaLa, GiaLa}.
Отсюда сразу же получалась справедливость и аналога гипотезы Амицура для полиномиальных $G$-тождеств
в случая действия конечной абелевой группы $G$ автоморфизмами.
Доказательство оценки снизу было основано на классификации конечномерных
градуированно простых алгебр.
Здесь нужно отметить, что в~\cite{ASGordienko18Schnabel} (см. также~\cite{ClaseJespersDelRio}) было показано, что не всякая градуировка даже на матричной алгебре эквивалентна градуировке конечной группой.

Также М.\,В.~Зайцевым
и Д.\,Реповшем было доказано существование градуированной PI-экспоненты
в конечномерных градуированно простых алгебрах, градуированных
коммутативными полугруппами~\cite[теорема 2]{ZaicevRepovs}. Полиномиальная оценка сверху для кодлин полиномиальных $H$-тождеств ассоциативных PI-алгебр  с обобщённым действием конечномерной ассоциативной алгебры $H$ с единицей была получена А.~Берелом~\cite{BereleHopf}.

Наконец, в последнее время
получили популярность исследования асимптотического поведения коразмерностей соответствующих типов тождеств
в конечномерных ассоциативных алгебрах с супер- и псевдоинволюциями~\cite{dosSantos,
GiaIopLaMa1, GiaIopLaMa2, Ioppolo, IopMar}. Как было впервые отмечено Р.\,Б. дос Сантосом~\cite{dosSantos}, всякая алгебра с суперинволюцией является
алгеброй с обобщённым $H$-действием. В теореме~\ref{TheoremGradGenActionReplace} ниже будет показано, что к обобщённым $H$-действиям
сводятся не только супер- и псевдоинволюции, но и для любые действия на алгебре, согласованные с некоторой градуировкой, имеющей конечный носитель.

В главах~\ref{ChapterGenHAssocCodim} и~\ref{ChapterHModLieCodim} аналог гипотезы Амицура для градуированных полиномиальных тождеств доказывается для всех конечномерных ассоциативных алгебр и алгебр Ли, градуированных
произвольными группами.
Для доказательства аналога гипотезы Амицура в случае ассоциативных алгебр автором диссертации был разработан метод, который
заключается в замене градуировки на обобщённое $H$-действие для некоторой ассоциативной алгебры $H$ с единицей. Этот метод, в частности, позволяет отказаться от использования классификации градуированно простых
алгебр и доказывать оценку снизу для коразмерностей с использованием теоремы плотности и центральных многочленов. Случай конечномерных алгебр Ли, градуированных произвольной группой, сводится
к случаю алгебр Ли, градуированных конечнопорождённой абелевой группой. Градуировка
конечнопорождённой абелевой группой $G$ затем заменяется на действие алгебры Хопфа $\mathbbm{k}\hat G$,
где $\mathbbm{k}$~"--- основное поле, а $\hat G := \Hom(G, \mathbbm{k}^\times)$. В главе~\ref{ChapterGenHAssocCodim} аналог гипотезы Амицура
также доказывается для коразмерностей полиномиальных $G$-тождеств конечномерных
ассоциативных алгебр с действием произвольной группы $G$ автоморфизмами и антиавтоморфизмами,
а в главе~\ref{ChapterHModLieCodim} аналогичный результат доказывается для алгебр Ли при условии, что
либо группа $G$ является редуктивной аффинной алгебраической группой, действующей рационально, либо
разрешимый радикал алгебры Ли является нильпотентным идеалом.
 Кроме того, в главе~\ref{ChapterGenHAssocCodim} аналог гипотезы Амицура
доказывается для коразмерностей дифференциальных полиномиальных тождеств конечномерных
ассоциативных алгебр с действием произвольной алгебры Ли дифференцированиями,
а в главе~\ref{ChapterHModLieCodim} аналогичный результат доказывается для алгебр Ли при условии, что
либо алгебра Ли, действующая дифференцированиями, конечномерная полупростая, либо
разрешимый радикал алгебры Ли, на которой действуют дифференцированиями, является нильпотентным идеалом.
Наконец, в главе~\ref{ChapterGenHAssocCodim} аналог гипотезы Амицура доказывается для
 ассоциативных алгебр c градуированным действием произвольной группы (см. определение~\ref{DefGradedAction}).

 В главе~\ref{ChapterSGGrAssocCodim} изучается
рост коразмерностей в конечномерных ассоциативных алгебрах, градуированных полугруппами, и строится 
семейство таких алгебр, имеющих дробную градуированную PI-экспоненту. Эти примеры оказались первыми примерами ассоциативных алгебр с дополнительной структурой, имеющих дробную экспоненту
роста коразмерностей соответствующих тождеств. В \S\ref{SectionHPIexpExistHSimple}
доказывается существование (необязательно целой) градуированной PI-экспоненты
у любой конечномерной градуированно простой алгебры, градуированной произвольным множеством.
Здесь также существенно используется метод замены градуировки на обобщённое $H$-действие для некоторой ассоциативной алгебры $H$ с единицей.

Для полиномиальных $H$-тождеств $H$-модульных ассоциативных алгебр аналог гипотезы Амицура
может быть сформулирован в следующей форме, которая принадлежит Ю.\,А.~Бахтурину:

\begin{conjecture} Пусть $A$~"--- конечномерная
ассоциативная $H$-модульная алгебра, где $H$~"--- алгебра Хопфа на полем характеристики $0$.
Тогда существует предел $\PIexp^H(A):=\lim\limits_{n\to\infty}
 \sqrt[n]{c^H_n(A)}$, который является целым числом.
\end{conjecture}

В главе~\ref{ChapterGenHAssocCodim} гипотеза Амицура~"--- Бахтурина
доказывается для конечномерных ассоциативных $H$-модульных
алгебр $A$, где $H$~"--- произвольная алгебра Хопфа, в случае, когда радикал Джекобсона
$J(A)$ является $H$-подмодулем. Кроме того, гипотеза Амицура~"--- Бахтурина
доказывается для конечномерных ассоциативных $H$-модульных
алгебр в случае, когда алгебра Хопфа $H$ либо полупроста и конечномерна, либо получена при помощи (возможно, многократного) расширения Оре конечномерной полупростой алгебры Хопфа косопримитивными элементами.
Класс $\mathcal C$ таких алгебр Хопфа достаточно широк. С одной стороны, он включает алгебры Тафта $H_{m^2}(\zeta)$, также как алгебры Хопфа $H(C,n,c,c^*,a,b)$ (см.~\cite[определение 5.6.15]{Danara}),
которые были использованы для опровержения гипотезы Капланского о конечном числе
попарно неизоморфных алгебр Хопфа заданной размерности~\cite{AndrusSchneider, BDG1, BDG2}. 
С другой стороны, класс $\mathcal C$ содержит и неточечные алгебры Хопфа,
например, две алгебры Хопфа размерности $16$ (см.~\cite[теорема~5.1]{CalDascMasMen}).

Здесь нужно отметить, что Я.~Карасик, используя результаты автора данной работы, доказал в 2015 году
гипотезу Амицура~"--- Бахтурина для необязательно конечномерных
ассоциативных $H$-модульных PI-алгебр в случае, когда $H$~"--- конечномерная полупростая алгебра Хопфа~\cite{Karasik}.

Все аналоги гипотезы Амицура для ассоциативных алгебр с различной дополнительной структурой
получаются в данной работе как следствия единственного утверждения (теоремы~\ref{TheoremHmodHRadAmitsurPIexpHBdimB}), которое выводит требуемое асимптотическое поведение
коразмерностей из существования для простых алгебр соответствующей сигнатуры полилинейного многочлена, кососимметричного по достаточному числу наборов переменных, а также существования инвариантного аналога разложения Веддербёрна~"--- Артина. Выполнение данных условий напрямую зависит от того, насколько хорошей
является дополнительная структура на алгебре.

В главе~\ref{ChapterHModLieCodim} аналог гипотезы Амицура доказывается
для конечномерных $H$-модульных алгебр Ли $L$ в следующих трёх случаях:
\begin{enumerate}
 \item когда  $H$~"--- конечномерная полупростая алгебра Хопфа;
 \item когда $H$~"--- произвольная алгебра Хопфа, но 
разрешимый радикал алгебры Ли $L$ является нильпотентным $H$-инвариантным идеалом;
\item когда $H$ является алгеброй Тафта $H_{m^2}(\zeta)$, а $L$~"--- $H_{m^2}(\zeta)$-простая алгебра Ли.
(Это единственный известный случай справедливости аналога гипотезы Амицура в $H$-модульной
алгебре Ли с не $H$-инвариантным разрешимым радикалом.)
\end{enumerate}

В целом схема доказательства аналогов гипотезы Амицура следует плану, разработанному А.~Джамбруно и М.\,В.~Зайцевым для обычных тождеств, однако для того, чтобы учесть дополнительную структуру
на алгебрах, требуются новые методы, которые и разрабатываются в работе.

\medskip

Данную работу, если не учитывать вспомогательные разделы, можно условно разделить на две части, тесно связанные между собой. В первой части, состоящей из глав~\ref{ChapterGradEquiv}--\ref{ChapterGenHActions}, рассматриваются структурные вопросы теории (ко)модульных алгебр, а также понятие эквивалентности (ко)модульных структур.
Во второй части, состоящей из глав~\ref{ChapterFreePICodim}--\ref{ChapterSGGrAssocCodim}, исследуется рост коразмерностей полиномиальных тождеств в алгебрах с дополнительной структурой.
Кроме того, во вводной главе~\ref{VvodGlav} приводится список обозначений и перечисляются основные используемые понятия и факты.

\newpage

\chapter{Основные понятия}
\label{VvodGlav}

\section{Список обозначений}

В работе используются следующие обозначения:

\medskip

\begin{longtable}{p{2.75cm} p{13cm}}
$1_M$ & единица моноида (группы, алгебры, \ldots) $M$ \\ \\
$\delta^i_j$,
$\delta_{ij}$   & символы Кронекера,
$\delta^i_j = \delta_{ij} = \left\lbrace
\begin{array}{ll}
1, & \text{ при } i=j, \\
0, & \text{ при } i\ne j
\end{array}
\right.$  \\ \\
$|M|$ & число элементов во множестве $M$ \\ \\
$A^2$                    & по определению,
                  $A^2=\langle ab \mid a,b \in A \rangle_\mathbbm{k}$ \\ \\
 $A \hookrightarrow B$ & вложение (инъективное отображение) множества $A$
во множество $B$ \\ \\
$A \mathrel{\widetilde\to} B$ & изоморфизм алгебраических структур $A$ и $B$ \\ \\
$A \times B$ & произведение объектов $A$ и $B$ некоторой категории, с.~\pageref{DefProd};
например, прямое произведение алгебраических структур $A$ и $B$
или декартово произведение множеств $A$ и $B$ \\ \\
$A \sqcup B$ & копроизведение объектов $A$ и $B$ некоторой категории, с.~\pageref{DefCoprod};
например, объединение непересекающихся множеств $A$ и $B$ \\ \\
$ f \sim g$              & запись означает,
                  что~$\lim \frac{f}{g} = 1$ \\ \\
$V \otimes W$ & тензорное произведение векторных пространств $V$ и $W$ над основным полем $\mathbbm{k}$ \\ \\
$ \lambda \vdash n$      & разбиение числа~$n$,
                      с.~\pageref{DefPartition}  \\ \\
$\pi_{\Gamma_1 \to \Gamma_2}$      & гомоморфизм универсальных групп градуировок $G_{\Gamma_1} \twoheadrightarrow G_{\Gamma_2}$, порождённый огрублением, с.~\pageref{NotationPiCoarser} \\ \\
$\binom{n}{n_1, n_2, \ldots, n_N}$ & полиномиальный
                   коэффициент,
                    равный~$\frac{n!}{n_1! n_2! \ldots n_N!}$ \\ \\
$\binom{n}{k}_\zeta$  & квантовый биномиальный коэффициент, с.~\pageref{DefQuantumBinom} \\ \\                  
$V^*$ & векторное пространство, двойственное к векторному пространству $V$, т.е. пространство линейных функций на $V$ со значениями в основном поле $\mathbbm{k}$ \\ \\ 
$C_{\mathrm{coc}}$ & сумма всех кокоммутативных подкоалгебр коалгебры $C$\\ \\
$A^\mathrm{op}$ & алгебра $A$, в которой в умножении меняются местами аргументы: $\mu^\mathrm{op}(a\otimes b):= \mu(b\otimes a)$ \\ \\                   
$C^\mathrm{cop}$ & коалгебра $C$, в которой в коумножении меняются местами тензорные множители в разложении результата: $\Delta^\mathrm{op}(c):= c_{(2)}\otimes c_{(1)}$\\ \\                   
$A^\circ$ & коалгебра, конечная двойственная к алгебре $A$, с.~\pageref{DefACirc} \\ \\
 $H^\circ$ & алгебра Хопфа, конечная двойственная к алгебре Хопфа $H$, с.~\pageref{DefHCirc} \\ \\
 $\langle Q \rangle_\mathbbm{k}$ & $\mathbbm{k}$-линейная оболочка подмножества $Q$ некоторого векторного
 пространства над полем $\mathbbm{k}$ \\ \\
$\tilde\varphi$ & изоморфизм $\tilde\varphi \colon \End_{\mathbbm{k}}(A_1) \mathrel{\widetilde\to} \End_{\mathbbm{k}}(A_2)$,
индуцированный изоморфизмом $\varphi \colon A_1 \mathrel{\widetilde\to} A_2$, 
задаётся при помощи равенства $\tilde\varphi(\psi)(a)=\varphi\Bigl(\psi\bigl(\varphi^{-1}(a)\bigr)\Bigr)$ для $\psi\in \End_{\mathbbm{k}}(A_1)$ и $a\in A_2$
\\ \\
$\hat G$ & если $G$~"--- группа, то $\hat G := \Hom(G, \mathbbm{k}^\times)$~"--- группа гомоморфизмов
из $G$ в мультипликативную группу $\mathbbm{k}^\times$ поля $\mathbbm{k}$, т.е. группа всех линейных характеров группы $G$ \\ \\
$\hat\rho$ & если $\rho \colon  A \to B \otimes Q$~"--- линейное
отображение, то линейное отображение $\hat\rho \colon Q^* \otimes A \to B$
определяется по формуле $\hat\rho(q^*\otimes a):= q^*(a_{(1)})a_{(0)}$, где $q^*\in Q^*$ и $a\in A$, с.~\pageref{DefHatRho} \\ \\
$\Phi(x_1, \ldots, x_q)$ & 
$\Phi(x_1, \ldots, x_q):=\frac{1}{x_1^{x_1} \ldots x_s^{x_q}}$ для $x_1, \ldots, x_q \geqslant 0$,
где $0^0 := 1$ \\
$\chi^H_n(A)$ & $n$-й
  кохарактер полиномиальных $H$-тождеств алгебры $A$, с.~\pageref{DefHCochar} \\ \\
$\chi^{T\text{-}\mathrm{gr}}_n(A)$  & $n$-й
  кохарактер $T$-градуированных полиномиальных тождеств алгебры $A$, с.~\pageref{DefCocharGr} \\ \\
$\chi(\lambda)$ &  характер неприводимого представления группы $S_n$, отвечающий разбиению $\lambda \vdash n$ \\ \\
$[x_1, \ldots, x_n]$ & длинный коммутатор
$[[\ldots [[x_1, x_2], x_3], \ldots], x_n]$ \\ \\
$\mathcal{A}^\square(A,B,V)$ & $V$-универсальная коизмеряющая алгебра, с.~\pageref{NotationAsquareABV} \\ \\
$\ad$ & присоединённое представление алгебры Ли: $(\ad a)b:=[a,b]$ \\ \\
$\mathbf{Alg}_\mathbbm{k}$ & категория ассоциативных алгебр с единицей над полем~$\mathbbm{k}$ \\ \\
$\Ann_R(M)$ & аннулятор модуля $M$ над кольцом $R$\\ \\
$\Aut(A)$ & группа автоморфизмов алгебраической структуры $A$ \\ \\
$\Aut^*(A)$ & группа автоморфизмов и антиавтоморфизмов алгебраической структуры $A$ \\ \\
$\Aut_{\mathbf{Alg}}(A)$ & группа линейных биекций $A \mathrel{\widetilde{\to}} A$,
сохраняющих умножение \\ \\
$\mathcal{B}^\square(A)$ & универсальная кодействующая биалгебра, с.~\pageref{NotationBsquareA} \\ \\
${}_\square \mathcal{B}(A)$ & универсальная действующая биалгебра, с.~\pageref{NotationsquareBA} \\ \\
$\mathcal{B}^\square(A,V)$ & $V$-универсальная кодействующая биалгебра, с.~\pageref{NotationBsquareAV} \\ \\
${}_\square \mathcal{B}(A,V)$ & $V$-универсальная действующая
 биалгебра, с.~\pageref{NotationsquareBAV} \\ \\
$\mathbf{Bialg}_\mathbbm{k}$ & категория биалгебр над полем~$\mathbbm{k}$ \\ \\
${}_\square \mathcal{C}(A,B,V)$ & $V$-универсальная измеряющая коалгебра,
  с.~\pageref{NotationsquareCABV} \\ \\
$\mathop\mathrm{char} \mathbbm{k}$ & характеристика поля~$\mathbbm{k}$ \\ \\
$\mathbf{Coalg}_\mathbbm{k}$ & категория коалгебр над полем~$\mathbbm{k}$ \\ \\
$\cosupp \rho$ & коноситель отображения $\rho \colon A \to B \otimes Q$, с.~\pageref{DefCosupportABQ} \\ \\
$\cosupp \psi$ & коноситель отображения $\psi \colon P \otimes A \to B$, с.~\pageref{DefCosupportPAB} \\ \\
$C_n$ & циклическая группа порядка $n$ в мультипликативной записи \\ \\
$c_n(A)$ & $n$-я коразмерность обычных полиномиальных тождеств алгебры $A$,
с.~\pageref{RemarkOrdinaryCodim}  \\ \\
$c_n^H(A)$ & $n$-я коразмерность полиномиальных $H$-тождеств алгебры $A$ с обобщённым $H$-действием,
с.~\pageref{DefHCodim} \\ \\
$c_n^{T\text{-}\mathrm{gr}}(A)$ & $n$-я коразмерность $T$-градуированных полиномиальных тождеств алгебры $A$,
с.~\pageref{DefCodimGr} \\ \\
$C_{T_\lambda}$          & подгруппа симметрической группы,
                     переводящая в себя множество чисел
                      каждого столбца таблицы Юнга~$T_\lambda$ \\ \\
$D_\lambda$    & диаграмма Юнга, отвечающая разбиению~$\lambda$,
                                с.~\pageref{DefDlambda} \\ \\
$\deg P$       & степень многочлена~$P$ \\ \\
$\Der(A)$       & алгебра Ли дифференцирований алгебры~$A$ \\ \\
$\diff U$      & если $U$~"--- подмножество группы $G$, то \par $\diff U := \lbrace uv^{-1} \mid u,v\in U,\ u\ne v\rbrace$, с. \pageref{NotationDiff} \\ \\
$e_{ij}$       & матричная единица, т.е. матрица~$n\times n$,
в которой на пересечении $i$-й строки и $j$-го столбца стоит~$1$,
 а все остальные клетки заполнены нулями \\ \\
 $\End_\mathbbm{k}(V)$ & алгебра $\mathbbm{k}$-линейных операторов из векторного пространства $V$ над полем $\mathbbm{k}$ в себя \\ \\
$e_{T_\lambda}$ & симметризатор, построенный
                 по таблице Юнга~$T_\lambda$,
                                с.~\pageref{DefeTlambda} \\ \\
$\mathbbm{k}G$ & групповая алгебра группы $G$ над полем $\mathbbm{k}$, с.~\pageref{DefGroupAlg}\\ \\
$\mathbbm{k}T$ & полугрупповая алгебра полугруппы $T$ над полем $\mathbbm{k}$, с.~\pageref{DefSemigroupAlg} \\ \\                                
$\mathbbm{k}^T$ & алгебра функций из множества $T$ в поле $\mathbbm{k}$ с поточечными операциями \\ \\
$\mathbbm{k}\langle X \rangle$  & свободная ассоциативная алгебра без единицы на
         множестве~$X$ 
             над полем~$\mathbbm{k}$,
                                с.~\pageref{DefFX} \\ \\
$\mathbbm{k}\lbrace X \rbrace$  & (абсолютно) свободная неассоциативная алгебра на
         множестве~$X$ 
             над полем~$\mathbbm{k}$,
                                с.~\pageref{DefFXbrace} \\ \\
$\mathbbm{k}\langle X|H \rangle$ & свободная ассоциативная алгебра без единицы на множестве $X$ с символами операторов из алгебры $H$ над полем~$\mathbbm{k}$, с.~\pageref{DefFXHangle} \\ \\
$\mathbbm{k}\lbrace X|H \rbrace$ & свободная неассоциативная алгебра на множестве $X$ с символами операторов из алгебры $H$ над полем~$\mathbbm{k}$, с.~\pageref{DefFXHbrace} \\ \\
$\mathcal F(X)$ & свободная группа со множеством $X$
 свободных порождающих \\ \\
 $G(H)$ & группа группоподобных элементов алгебры Хопфа $H$, с.~\pageref{DefGroupLike} \\ \\
 $\GL_n(\mathbbm{k})$ & группа всех невырожденных квадратных матриц $n\times n$ над полем $\mathbbm{k}$ \\ \\
 $\mathfrak{gl}_n(\mathbbm{k})$ & алгебра Ли всех квадратных матриц размера $n\times n$ \\ \\
$\mathbf{Grp}$ & категория групп \\ \\
 $H_4$ & алгебра Свидлера, с.~\pageref{ExampleTaftAlgebra} \\ \\
$\mathcal{H}^\square(A)$ & универсальная кодействующая алгебра Хопфа, с.~\pageref{NotationHsquareA} \\ \\
${}_\square \mathcal{H}(A)$ & универсальная действующая алгебра Хопфа, с.~\pageref{NotationsquareHA} \\ \\
$\mathcal{H}^\square(A,V)$ & $V$-универсальная кодействующая алгебра Хопфа, с.~\pageref{NotationHsquareAV} \\ \\
${}_\square \mathcal{H}(A,V)$ & $V$-универсальная действующая алгебра Хопфа, с.~\pageref{NotationsquareHAV} \\ \\
${}_\square \mathcal{H}(A,V)_\mathrm{coc}$ & $V$-универсальная кокоммутативная действующая  алгебра Хопфа, с.~\pageref{NotationsquareHAVcoc} \\ \\
 $H_{m^2}(\zeta)$ & алгебра Тафта, с.~\pageref{ExampleTaftAlgebra} \\ \\
 $\mathbf{Hopf}_\mathbbm{k}$ & категория алгебр Хопфа над полем~$\mathbbm{k}$ \\ \\
$\Hom(A,B)$ & гомоморфизмы из алгебраической структуры $A$ в алгебраическую структуру $B$ \\ \\
$\Hom_\mathbbm{k}(V,W)$ & пространство $\mathbbm{k}$-линейных операторов (или, что то же, $\mathbbm{k}$-линейных отображений) из векторного пространства $V$ в векторное
пространство $W$ над полем $\mathbbm{k}$ \\ \\
$q_t$ & элемент алгебры $\mathbbm{k}^T$, заданный равенством
$q_t(s):=\left\lbrace\begin{smallmatrix} 1 & \text{при} & s=t,\\ 0 & \text{при} & s\ne t,\end{smallmatrix} \right.$ с.~\pageref{DefHT} \\ \\
$\id_M$ & 
тождественный морфизм объекта $M$ некоторой категории в себя; например, тождественное отображение множества $M$ в себя \\ \\
$\Id(A)$ & идеал обычных полиномиальных тождеств алгебры $A$, с.~\pageref{RemarkOrdinaryCodim} \\ \\
$\Id^H(A)$ & идеал полиномиальных $H$-тождеств алгебры $A$ с обобщённым $H$-действием, с.~\pageref{DefIdHA} \\ \\
$\Id^{T\text{-}\mathrm{gr}}(A)$ & $T$-градуированных полиномиальных тождеств алгебры $A$, с.~\pageref{DefIdTgrA} \\ \\
$\im \varphi$ & образ гомоморфизма или линейного отображения $\varphi$ \\ \\
$J(A)$ & радикал Джекобсона ассоциативной алгебры $A$ \\ \\
$J^H(A)$ & $H$-радикал конечномерной ассоциативной $H$-модульной алгебры $A$,
 с.~\pageref{DefJHA} \\ \\
$\ker \varphi$ & ядро гомоморфизма или линейного отображения $\varphi$ \\ \\
$\ell_n(A)$ & $n$-я кодлина  обычных полиномиальных тождеств алгебры $A$, с.~\pageref{RemarkOrdinaryCodim} \\ \\
$\ell_n^H(A)$ & $n$-я кодлина полиномиальных $H$-тождеств алгебры $A$, с.~\pageref{DefHColength}\\ \\
$\ell_n^{T\text{-}\mathrm{gr}}(A)$ & $n$-я кодлина $T$-градуированных полиномиальных тождеств алгебры $A$, с.~\pageref{DefColengthGr}\\ \\
$L(X)$ & свободная алгебра Ли на множестве $X$, с.~\pageref{DefLX} \\ \\
$L(X|H)$ & свободная $H$-модульная алгебра Ли на множестве $X$, с.~\pageref{DefLXH} \\ \\
$m(A, \lambda)$ & кратность неприводимого характера $\chi(\lambda)$
  в разложении кохарактера $\chi_n(A)$ \\ \\
$m(A, H, \lambda)$ & кратность неприводимого характера $\chi(\lambda)$
  в разложении кохарактера $\chi^H_n(A)$ \\ \\
$m(A, T\text{-}\mathrm{gr}, \lambda)$ & кратность неприводимого характера $\chi(\lambda)$
  в разложении кохарактера $\chi^{T\text{-}\mathrm{gr}}_n(A)$ \\ \\
$M_{k\times \ell}(\mathbbm{k})$ & пространство всех матриц размера $k\times \ell$ над полем $\mathbbm{k}$ \\ \\
$M_{s,\ell}(\mathbbm{k})$ & алгебра $M_{s+\ell}(\mathbbm{k})$
с $\mathbb Z/2\mathbb Z$-градуировкой, заданной равенствами
$M_{s,\ell}^{(0)}(\mathbbm{k}) := \left( \begin{smallmatrix}
 M_s(\mathbbm{k}) & 0 \\
 0 & M_\ell(\mathbbm{k})\\
 \end{smallmatrix}\right)$ и
 $M_{s,\ell}^{(1)}(\mathbbm{k}) := \left( \begin{smallmatrix}
 0 & M_{s\times \ell}(\mathbbm{k}) \\
 M_{\ell\times s}(\mathbbm{k}) & 0\\
 \end{smallmatrix}\right)$
 \\ \\
$M_n(\mathbbm{k})$         &  алгебра всех матриц размера~$n\times n$ над полем~$\mathbbm{k}$ \\ \\
$M(\lambda)$      &  $\mathbbm{k}S_n$-модуль, изоморфный минимальному
                      левому идеалу~$\mathbbm{k}S_n e_{T_\lambda}$ \\ \\
$\mathbb N$      & множество~$\lbrace 1, 2, 3, \ldots \rbrace$ натуральных
                  чисел \\ \\
$\mathbf{nuAlg}_\mathbbm{k}$ & категория ассоциативных алгебр над полем~$\mathbbm{k}$ (необязательно с единицей) \\ \\                  
$\mathcal O(G)$ & алгебра Хопфа регулярных функций на аффинной алгебраической
группе $G$, с.~\pageref{ExampleOAffAlgGrp} \\ \\
$\mathcal O(V)$  & алгебра регулярных функций на аффинном алгебраическом
многообразии $V$ \\ \\
$\PGL_n(\mathbbm{k})$ & проективная линейная группа, т.е. факторгруппа группы $\GL_n(\mathbbm{k})$ по её центру \\ \\
$P_n$ & пространство полилинейных ассоциативных многочленов от переменных $x_1,\ldots, x_n$\\ \\
$P_n^H$ & пространство полилинейных ассоциативных $H$-многочленов от переменных $x_1,\ldots, x_n$, с.~\pageref{DefPnH} \\ \\
$P^{T\text{-}\mathrm{gr}}_n$ & пространство полилинейных ассоциативных $T$-градуированных многочленов $n$-й степени \\ \\
$\PIexp^H(A)$ & $H$-PI-экспонента алгебры $A$, с.~\pageref{DefHPIexp},\\ & $\PIexp^H(A):=\lim\limits_{n\rightarrow\infty} \sqrt[n]{c^H_n(A)}$ \\ \\
$\PIexp^{T\text{-}\mathrm{gr}}(A)$ & $T$-градуированная PI-экспонента алгебры $A$, с.~\pageref{DefTgrPIexp},\\ & $\PIexp^{T\text{-}\mathrm{gr}}(A):=\lim\limits_{n\rightarrow\infty} \sqrt[n]{c^{T\text{-}\mathrm{gr}}_n(A)}$ \\ \\
$R_{T_\lambda}$          & подгруппа симметрической группы,
                      переводящая в себя множество чисел
                      каждой строки таблицы Юнга~$T_\lambda$ \\ \\
$S_n$             & группа подстановок (симметрическая группа)
                     на множестве~$\lbrace 1,2, \ldots, n \rbrace$ \\ \\
$\mathbf{Sets}$ & категория множеств   \\ \\
$\mathfrak{sl}_n(\mathbbm{k})$ & алгебра Ли всех квадратных матриц размера $n\times n$ с нулевым следом \\ \\
$\supp \Gamma$ & носитель градуировки $\Gamma$, с.~\pageref{DefSupport} \\ \\
$\supp \rho$ & носитель отображения $\rho \colon A \to B \otimes Q$, с.~\pageref{DefSupportABQ} \\ \\
$\mathfrak{t}_n(\mathbbm{k})$ & алгебра Ли верхнетреугольных матриц размера $n\times n$ \\ \\
$T_\lambda$    & таблица Юнга, отвечающая разбиению~$\lambda$,
         с.~\pageref{DefTlambda} \\ \\
$\mathbf{Top}$ & категория, в которой объектами являются всевозможные топологические пространства,
а морфизмами~"--- непрерывные отображения \\ \\
$\tr \mathcal A$ & след линейного оператора $\mathcal A$ \\ \\
$\lambda^T$              & разбиение, транспонированное по отношению
                            к разбиению~$\lambda$,
                                с.~\pageref{DeflambdaT} \\ \\
$\mathcal U(M)$  & группа обратимых элементов моноида $M$ \\ \\
$U(L)$  & универсальная обёртывающая алгебра алгебры Ли $L$ \\ \\
$\UT_n(\mathbbm{k})$ & алгебра всех верхнетреугольных матриц размера $n\times n$ над полем $\mathbbm{k}$ \\ \\
 $\mathbf{Vect}_\mathbbm{k}$ & категория векторных пространств над полем $\mathbbm{k}$   \\ \\
$V_n$ & пространство полилинейных лиевских многочленов от переменных $x_1,\ldots, x_n$ \\ \\
$V_n^H$ & пространство полилинейных лиевских $H$-многочленов от переменных $x_1,\ldots, x_n$, с.~\pageref{DefVnH} \\ \\
$V^{T\text{-}\mathrm{gr}}_n$ & пространство полилинейных лиевских $T$-градуированных многочленов $n$-й степени \\ \\
$W_n$ & пространство полилинейных неассоциативных многочленов от переменных $x_1,\ldots, x_n$, с.~\pageref{RemarkOrdinaryCodim} \\ \\
$W_n^H$ & пространство полилинейных неассоциативных $H$-многочленов от переменных $x_1,\ldots, x_n$, с.~\pageref{DefWnH} \\ \\
$W^{T\text{-}\mathrm{gr}}_n$ & пространство полилинейных неассоциативных $T$-градуированных многочленов $n$-й степени, с.~\pageref{DefWTgrn} \\ \\
$X^{T\text{-}\mathrm{gr}}$ &  объединение $X^{T\text{-}\mathrm{gr}}:=\bigsqcup_{t \in T}X^{(t)}$
непересекающихся множеств $X^{(t)} = \{ x^{(t)}_1,
x^{(t)}_2, \ldots \}$, с.~\pageref{DefXTgr} \\ \\
$\mathbb Z_+$      & множество~$\lbrace 0, 1, 2, 3, \ldots \rbrace$
                   целых неотрицательных  чисел \\ \\
                   $Z(A)$ & \textit{центр} алгебраической
                   структуры $A$: если $A$~"--- группа или ассоциативное кольцо,
                   то $$Z(A):= \lbrace b\in A \mid ab = ba \text{ для всех } a\in A \rbrace;$$
      если же $A$~"--- алгебра Ли,
                   то $$Z(A):= \lbrace b\in A \mid [a,b]=0 \text{ для всех } a\in A \rbrace$$            
 \end{longtable}

Алгебры, рассматриваемые в работе, необязательно ассоциативны и необязательно содержат единицу.
Все коалгебры, которые рассматриваются в работе, коассоциативные с коединицей.

 Если $A$ и $B$~"--- объекты категории $\mathcal A$, то через $\mathcal A(A,B)$
обозначается множество морфизмов $A \to B$ в категории $\mathcal A$. Например, $\mathbf{Vect}_\mathbbm{k}(V,W)$~"--- это то же самое, что и $\Hom_\mathbbm{k}(V,W)$, т.е. множество всех линейных отображений $V \to W$, где $V$ и $W$~"--- векторные пространства над полем $\mathbbm{k}$.
  
  Забывающие функторы, например $\mathbf{Alg}_\mathbbm{k} \to \mathbf{Vect}_\mathbbm{k}$, и функторы вложения, 
  например, $\mathbf{Hopf}_\mathbbm{k} \subseteq \mathbf{Bialg}_\mathbbm{k}$, в тексте работы обычно никак не обозначаются. 
  Множества $\mathbf{Vect}_\mathbbm{k}(V,W)$ считаются по умолчанию наделёнными структурой векторного пространства над полем $\mathbbm{k}$.

%

%

\section{(Ко)алгебры, алгебры Ли, биалгебры и алгебры Хопфа}

В данном параграфе мы напомним основные понятия, связанные с (ко)алгебрами, алгебрами Ли, биалгебрами и алгебрами Хопфа.
Подробнее с этими понятиями можно познакомиться в монографиях~\cite{HumphreysLieAlg, Danara, Montgomery, Sweedler}.

\subsection{Алгебры и коалгебры}

Прежде всего напомним понятие алгебры над полем в удобной для нас форме, а именно, на языке линейных отображений и коммутативных диаграмм. 

\begin{definition} 
 \textit{Алгеброй над полем $\mathbbm{k}$} называется пара $(A, \mu)$, состоящая из
 векторного пространства $A$ над $\mathbbm{k}$ и линейного отображения
$\mu \colon A \otimes A \to A$.
\end{definition}

При помощи отображения $\mu$ на векторном пространстве $A$ задаётся операция внутреннего умножения, линейная по каждому аргументу:
$ab := \mu(a\otimes b)$ для всех $a,b\in A$.

Алгебра $(A, \mu)$ называется \textit{ассоциативной}, если
следующая диаграмма коммутативна: $$\xymatrix{ A \otimes A \otimes A \ar[d]_{\mu \otimes \id_A} \ar[rr]^{\id_A \otimes \mu} && A\otimes A
\ar[d]_\mu \\
A\otimes A \ar[rr]^\mu && A  } $$

В случае алгебр Ли существует традиция обозначать операцию внутреннего умножения при помощи коммутатора:

\begin{definition} 
 \textit{Алгеброй Ли над полем $\mathbbm{k}$} называется векторное пространство $L$ над полем $\mathbbm{k}$ с заданным на нём билинейным отображением $(a,b)\mapsto [a,b]$, удовлетворяющим следующим аксиомам:
 \begin{enumerate}
\item (антикоммутативность)\qquad 
$[a,a]=0$;
\item (тождество Якоби)\qquad $[a,[b,c]]+[[b,c],a]+[[c,a],b]=0
\text{ для всех }a,b,c \in L$.
\end{enumerate}
\end{definition}
\begin{remark}
В случае, когда $\chr \mathbbm{k} \ne 2$, условие 1 эквивалентно условию $[a,b]=-[b,a]$ для всех $a,b\in L$.
\end{remark}

Дадим теперь определение алгебры с единицей:

\begin{definition} 
 Набор $(A, \mu, u)$ называется \textit{алгеброй с единицей над полем $\mathbbm{k}$}, если
 $(A,\mu)$~"--- алгебра над $\mathbbm{k}$, а
 $u \colon \mathbbm{k} \to A$~"--- линейное отображение и, кроме того, следующие диаграммы коммутативны:
 $$\xymatrix{  A \otimes \mathbbm{k} \ar[rr]^{\id_A \otimes u} \ar[rd]^{\sim} 
  & & A\otimes A \ar[ld]^{\mu}\\
 & A  & } \text{\qquad и \qquad}\xymatrix{  \mathbbm{k} \otimes A \ar[rr]^{u \otimes \id_A} \ar[rd]^{\sim} 
  & & A\otimes A \ar[ld]^{\mu}\\
 & A  & }$$
(Здесь мы использовали естественные отождествления
 $A \otimes \mathbbm{k} \cong \mathbbm{k} \otimes A \cong A$.)
\end{definition}
Если ввести обозначение $1_A := u(1_\mathbbm{k})$, то элемент $1_A\in A$ будет удовлетворять классическому определению единицы в алгебре.

\medskip

Определение коассоциативной коалгебры с коединицей получается из определения ассоциативной алгебры с единицей формальным обращением стрелок:

\begin{definition} Набор $(C, \Delta, \varepsilon)$
называется \textit{коассоциативной коалгеброй с коединицей},
если $C$~"--- векторное пространство на полем $\mathbbm{k}$,
$\Delta \colon C \to C \otimes C$ и $\varepsilon \colon C \to \mathbbm{k}$
являются линейными отображениями и, кроме того, следующие диаграммы коммутативны:
 \begin{enumerate}
\item (коассоциативность) $$\xymatrix{ C \ar[d]^\Delta \ar[rr]^\Delta && C\otimes C
\ar[d]^{\id_C \otimes \Delta} \\
C\otimes C \ar[rr]_{\Delta \otimes \id_C} && C \otimes C \otimes C  } $$
\item (наличие коединицы) $$\xymatrix{  & C \ar[ld]_{\Delta} \ar[rd]^{\sim} &  \\
C \otimes C \ar[rr]^{\id_C \otimes \varepsilon} && C \otimes \mathbbm{k}
}\text{\qquad и \qquad}\xymatrix{  & C \ar[ld]_{\Delta} \ar[rd]^{\sim} &  \\
C \otimes C \ar[rr]^{\varepsilon \otimes \id_C} && \mathbbm{k} \otimes C
}$$

(Здесь мы использовали естественные отождествления $C \otimes \mathbbm{k} \cong \mathbbm{k} \otimes C \cong C$.)
\end{enumerate}
\end{definition}

\medskip

В работе под \textit{коалгебрами} понимаются именно коассоциативные коалгебры с коединицей.

Отображение $\Delta$ называется \textit{коумножением}, а отображение $\varepsilon$ называется \textit{коединицей} коалгебры $(C, \Delta, \varepsilon)$.

Подпространство $I$ коалгебры $(C,\Delta,\varepsilon)$ называется \textit{коидеалом},
если $\varepsilon(I)=0$ и $\Delta(I)\subseteq I \otimes C + C \otimes I$.
На факторпространстве $C/I$ можно естественным образом ввести структуру коалгебры.
Коидеалы~"--- это в точности ядра гомоморфизмов коалгебр.

Подпространство $D$ коалгебры $(C,\Delta,\varepsilon)$ называется \textit{подкоалгеброй},
если $\Delta(D)\subseteq D \otimes D$.

В работе используются обозначения М.~Свидлера $\Delta c = c_{(1)}\otimes c_{(2)}$, где $c\in C$
и опущен знак суммы. 
Если оператор $\Delta$ применяется к элементу несколько раз, то в силу коассоциативности
неважно, к какому именно множителю тензорного произведения применяется каждое последующее
$\Delta$. В частности, мы можем определить $$c_{(1)}\otimes \ldots \otimes c_{(n)} := (\underbrace{\id_C \otimes \ldots \otimes \id_C}_{n-k-1} \otimes \Delta \otimes
\underbrace{\id_C \otimes \ldots \otimes \id_C}_{k})\, c_{(1)}\otimes \ldots \otimes c_{(n-1)}$$ для $c\in C$, где результат не зависит от $0\leqslant k \leqslant n-1$.

Говорят, что коалгебра $(C,\Delta,\varepsilon)$ \textit{кокоммутативна},
если $c_{(1)}\otimes c_{(2)}=c_{(2)}\otimes c_{(1)}$ для всех $c\in C$.

Если $(C,\Delta,\varepsilon)$~"--- коалгебра, то алгебра $(C^*, \Delta^*, \varepsilon^*)$, \textit{двойственная} к $(C,\Delta,\varepsilon)$, определяется при помощи ограничения
отображения $ \Delta^* \colon (C\otimes C)^* \to C^*$,  двойственного к $\Delta$,
на подпространство $C^*\otimes C^* \subseteq (C\otimes C)^*$ и композиции отображения
$\varepsilon^* \colon C^* \to \mathbbm{k}^*$ с естественным отождествлением $\mathbbm{k}^* \cong \mathbbm{k}$.
(Для удобства отображения $\Delta^*$ и $\varepsilon^*$, изменённые таким образом, обозначаются
по-прежнему через $\Delta^*$ и $\varepsilon^*$.)
Очевидно, что коалгебра $(C, \Delta, \varepsilon)$ кокоммутативна, если и только если алгебра $(C^*, \Delta^*, \varepsilon^*)$ коммутативна.

Если $(A,\mu,u)$~"--- конечномерная ассоциативная алгебра с единицей, то,
используя двойственные отображения $\mu^* \colon A^* \to (A\otimes A)^*$
и $u^* \colon A^* \to \mathbbm{k}^*$ и естественные отождествления $(A\otimes A)^* \cong A^*\otimes A^*$ и $\mathbbm{k}^* \cong \mathbbm{k}$, можно определить коалгебру $(A^*,\mu^*,u^*)$, \textit{двойственную} к $(A,\mu,u)$.
Если алгебра $A$ бесконечномерна, то определить двойственную коалгебру таким образом уже не получается.
В этом случае рассматривают помножество $A^\circ \subseteq A^*$ таких линейных функций на алгебре $A$, которые содержат в своём ядре некоторый идеал конечной коразмерности.
Тогда $\mu^*(A^\circ)\subseteq A^\circ \otimes A^\circ$. Обозначив через
$\mu^\circ$ и $u^\circ$, соответственно, ограничения отображений $\mu^*$ и $u^*$
на $A^\circ$, получаем коалгебру $(A^\circ,\mu^\circ,u^\circ)$\label{DefACirc}, \textit{конечную двойственную} к $(A,\mu,u)$.

В дальнейшем коалгебры $(C, \Delta, \varepsilon)$ для краткости обозначаются просто через $C$. Аналогичное соглашение действует и для остальных алгебраических структур.

Коалгебра $C$ называется \textit{простой},
если её единственными подкоалгебрами являются $0$ и $C$.
Например, основное поле $\mathbbm{k}$ с тривиальными коумножением $\Delta$ и коединицей $\varepsilon$,
заданными равенствами $\Delta(1)=1\otimes 1$ и $\varepsilon(1)=1$,
является простой коалгеброй.

\begin{proposition}\label{PropositionCocommAlgClosedSimple} Если $C$~"--- простая кокоммутативная коалгебра над алгебраически замкнутым полем $\mathbbm{k}$,
то $C \cong \mathbbm{k}$.
\end{proposition}
\begin{proof} Известно (см., например, теорему 1.4.7 из~\cite{Danara}), что все простые коалгебры конечномерны. В силу соответствия между идеалами в $C^*$ и подкоалгебрами в $C$ коммутативная алгебра $C^*$ проста. Поскольку поле $\mathbbm{k}$ алгебраически замкнуто, алгебра  $C^*$, а значит, и коалгебра $C$ одномерна.
\end{proof}

\subsection{Биалгебры и алгебры Хопфа}

Биалгеброй называется алгебра, на которой, кроме этого, ещё задана структура коалгебры, и эти две структуры согласованы между собой:

\begin{definition} Набор $(B, \mu, u, \Delta, \varepsilon)$ называется \textit{биалгеброй}
над полем $\mathbbm{k}$, если $(B, \mu, u)$~"--- ассоциативная алгебра с единицей, а
$(B, \Delta, \varepsilon)$~"--- (коассоциативная) коалгебра (с коединицей),
причём $\Delta \colon B \to B\otimes B$
и $\varepsilon \colon B \to \mathbbm{k}$ являются гомоморфизмами алгебр с единицей.
\end{definition}

Подпространство $I$ биалгебры $B$ называется \textit{биидеалом}, если $B$~"--- одновременно и идеал, и коидеал. Факторпространство $B/I$ естественным образом наследует с $B$ структуру биалгебры.

Если от моноида потребовать, чтобы все его элементы были обратимы, то мы получим определение группы.
Алгебры Хопфа выделяются среди класса всех биалгебр схожим способом\footnote{Многочисленные сходства между группами и алгебрами Хопфа объясняются тем, что оба этих понятия являются частными случаями одного и того же понятия для разных категорий: категории $\mathbf{Sets}$ множеств и категории $\mathbf{Vect}_\mathbbm{k}$ векторных пространств над некоторым полем $\mathbbm{k}$.}:
\begin{definition} Набор $(H, \mu, u, \Delta, \varepsilon, S)$ называется \textit{алгеброй Хопфа}
над полем $\mathbbm{k}$, если $(H, \mu, u, \Delta, \varepsilon)$~"--- биалгебра,
а $S \colon H \to H$~"--- такое линейное отображение, что
 $$(Sh_{(1)})h_{(2)}=h_{(1)}(Sh_{(2)})=\varepsilon(h)1_H\text{ для всех }h\in H.$$
\end{definition}

Отображение $S$ называется \textit{антиподом} алгебры $H$. Можно показать, что $S$ является антигомоморфизмом алгебры $H$.
Более того, $(Sh)_{(1)}\otimes (Sh)_{(2)}=Sh_{(2)}\otimes Sh_{(1)}$ для всех $h\in H$ \cite [предложение~4.2.6]{Danara}.

\begin{example}\label{ExampleFTBialgebra} Пусть $T$~"--- полугруппа, а $\mathbbm{k}$"--- поле.
Рассмотрим алгебру $\mathbbm{k}T$ с базисом, состоящим из элементов полугруппы $T$ и умножением,
продолженным по линейности с умножения в $T$. Тогда $\mathbbm{k}T$ называется \textit{полугрупповой алгеброй}\label{DefSemigroupAlg}
полугруппы $T$. Определим отображения $\Delta \colon \mathbbm{k}T \to \mathbbm{k}T \otimes \mathbbm{k}T$
и $\varepsilon \colon \mathbbm{k}T \to \mathbbm{k}$ при помощи формул
$\Delta(t)=t\otimes t$, $\varepsilon(t)=1$ при $t\in T$, а затем продолжим их по линейности на всё $\mathbbm{k}T$.
Тогда $\mathbbm{k}T$~"--- биалгебра, если и только если
$T$~"--- \textit{моноид}, т.е. полугруппа с единицей.
Более того, $\mathbbm{k}T$ является алгеброй Хопфа,
если и только если $T$~"--- группа, причём в последнем случае справедливо равенство $S(t)=t^{-1}$ для всех $t\in T$.
Полугрупповая алгебра $\mathbbm{k}G$ группы $G$ называется \textit{групповой алгеброй}
группы $G$.\label{DefGroupAlg}
\end{example}

\begin{example} Если $L$~"--- алгебра Ли, то её универсальная обёртывающая алгебра $U(L)$ является алгеброй Хопфа, где $\Delta(v)=1\otimes v + v\otimes 1$, $\varepsilon(v)=0$,
$S(v)=-v$ для $v \in L$. На всю алгебру $U(L)$ отображения $\Delta, \varepsilon, S$
продолжаются таким образом, чтобы $\Delta$ и $\varepsilon$ были гомоморфизмами алгебр с $1$, а $S$~"--- антигомоморфизмом алгебр с $1$.
\end{example}

\begin{example}\label{ExampleOAffAlgGrp} Пусть $G$~"--- аффинная алгебраическая группа над алгебраически замкнутым полем $\mathbbm{k}$.
Тогда алгебра регулярных функций $\mathcal O(G)$ на $G$ является алгеброй Хопфа, где
коумножение $\Delta \colon \mathcal O(G) \to  \mathcal O(G) \otimes  \mathcal O(G)$~"--- это гомоморфизм алгебр, индуцированный умножением $G \to G\times G$ в группе, коединица $\varepsilon \colon \mathcal O(G) \to \mathbbm{k}$ определяется равенством $\varepsilon(f)=f(1_G)$ для всех $f\in \mathcal O(G)$,
а антипод $S \colon \mathcal O(G) \to  \mathcal O(G)$~"--- это гомоморфизм алгебр, индуцированный морфизмом аффинных алгебраических многообразий $G \to G$, $g\mapsto g^{-1}$.
\end{example}

\begin{example}\label{ExampleTaftAlgebra}
Пусть $m \geqslant 2$~"--- натуральное число, а $\zeta$~"--- примитивный корень $m$-й степени из единицы
в поле $\mathbbm{k}$. (Такой корень может существовать в поле $\mathbbm{k}$, только если $\chr \mathbbm{k} \nmid m$.)
Рассмотрим алгебру $H_{m^2}(\zeta)$ с $1$, порождённую
элементами $c$ и $v$, которые удовлетворяют соотношениям $c^m=1$, $v^m=0$, $vc=\zeta cv$.
Тогда элементы $(c^i v^k)_{0 \leqslant i, k \leqslant m-1}$ образуют базис алгебры $H_{m^2}(\zeta)$.
Введём на алгебре $H_{m^2}(\zeta)$
структуру коалгебры при помощи равенств
  $\Delta(c)=c\otimes c$,
$\Delta(v) = c\otimes v + v\otimes 1$, $\varepsilon(c)=1$, $\varepsilon(v)=0$.
Тогда $H_{m^2}(\zeta)$~"--- алгебра Хопфа с антиподом~$S$, где $S(c)=c^{-1}$
и $S(v)=-c^{-1}v$. Алгебра Хопфа $H_{m^2}(\zeta)$ называется \textit{алгеброй Тафта}.
Алгебра $H_4(-1)$ называется \textit{алгеброй Свидлера}. Так как $(-1)$~"--- единственный примитивный корень $2$-й степени из единицы, мы будем обозначать алгебру Свидлера просто через $H_4$.
\end{example}

\begin{definition}
Элемент $h\ne 0$ алгебры Хопфа $H$ называется \textit{группоподобным},\label{DefGroupLike}
если $\Delta h = h\otimes h$. Группоподобные элементы алгебры Хопфа образуют группу $G(H)$ относительно операции умножения.
\end{definition}
\begin{definition}
Элемент $h$ алгебры Хопфа $H$ называется
\textit{примитивным},
если $\Delta h = h \otimes 1 + 1 \otimes h$.
Примитивные элементы алгебры Хопфа образуют алгебру Ли относительно коммутатора $[h_1,h_2]=h_1 h_2 - h_2 h_1$.
\end{definition}

Если $H$~"--- биалгебра или алгебра Хопфа, то можно рассмотреть коалгебру $H^\circ$, конечную двойственную к
алгебре $H$. При этом $H^\circ$ также оказывается подалгеброй в алгебре $H^*$. Биалгебра $H^\circ$
называется \textit{конечной двойственной} к биалгебре $H$.\label{DefHCirc} Если $H$~"--- алгебра Хопфа, то $H^\circ$ также является алгеброй Хопфа. Если $H$ конечномерна, то $H^\circ=H^*$, и $H^*$ называется,
соответственно, биалгеброй или алгеброй Хопфа, \textit{двойственной} к $H$.

Подпространство $I$ алгебры Хопфа $H$ называется \textit{идеалом Хопфа}, если $I$~"--- биидеал и, кроме того, $SI \subseteq I$. Факторпространство $H/I$ естественным образом наследует с $H$ структуру алгебры Хопфа.

\subsection{ad-инвариантные левые интегралы на алгебрах Хопфа}

Пусть $H$~"--- алгебра Хопфа. Напомним, что линейная функция $t \in H^*$ называется
 \textit{левым интегралом на $H$}, если $t(h_{(2)})h_{(1)}= t(h)1_H$ для всех $h \in H$.
 Говорят, что левый интеграл $t$ является \textit{$\ad$-инвариантным},
 если $t(a_{(1)}\ b\ S(a_{(2)}))=\varepsilon(a)t(b)$ для всех $a,b \in H$.
 
  В теореме~\ref{TheoremHcoLevi}
будет доказано существование $H$-коинвариантного
разложения Леви в предположении, что для алгебры Хопфа $H$
существует $\ad$-инвариантный левый интеграл
 $t \in H^*$, такой, что $t(1)=1$.
 Приведём три главных примера алгебр Хопфа $H$, обладающих таким свойством~\cite{SteVanOyst}.

Прежде всего заметим, что существование такого интеграла $t \in H^*$, что $t(1)=1$,
эквивалентно кополупростоте алгебры Хопфа $H$ (см., например, \cite[пример~5.5.9]{Danara}).

\begin{example}\label{ExampleIntegralHSS}
Пусть $H$~"--- конечномерная (ко)полупростая
алгебра Хопфа над полем характеристики $0$.
Тогда существует такой $\ad$-инвариантный левый интеграл $t \in H^*$,
что $t(1)=1$.
\end{example}
\begin{proof} В силу теоремы Ларсона~"--- Рэдфорда (см., например, \cite[теорема 7.4.6]{Danara}),
алгебра Хопфа $H$ полупроста, если и только если она кополупроста.
Третьим эквивалентным условием является
 $S^2 = \id_H$. Как мы уже отметили выше, если $H$ кополупроста, существует левый интеграл $t \in H^*$, такой, что $t(1)=1$. Всякий интеграл на конечномерной полупростой алгебре Хопфа кокоммутативен (см., например, \cite[упражнение~7.4.7]{Danara}), т.е. $t(ab)=t(ba)$ для всех $a,b \in H$.  Следовательно, $$t(a_{(1)}\ b\ S(a_{(2)}))=t(b (Sa_{(2)})a_{(1)})=
t\Bigl(b\ S\bigl((Sa_{(1)})a_{(2)}\bigr)\Bigr)=\varepsilon(a)t(b)\text{ для всех }a,b \in H$$ и интеграл $t$ является $\ad$-инвариантным.\end{proof}

\begin{example}\label{ExampleIntegralFG}
Пусть $G$~"--- группа. Определим $t \in (\mathbbm{k}G)^*$
при помощи равенств
$$t(g)=\left\lbrace\begin{array}{ccc} 0 & \text{ при } & g\ne 1, \\ 
1 & \text{ при } & g = 1. \end{array}\right.$$
Тогда $t$ является левым $\ad$-инвариантным интегралом для алгебры Хопфа $\mathbbm{k}G$, причём $t(1)=1$.
\end{example}

\begin{example}\label{ExampleIntegralAffAlgGr}
Пусть $G$~"--- редуктивная аффинная алгебраическая группа над алгебраически замкнутым полем $\mathbbm{k}$ характеристики $0$.
Тогда существует такой $\ad$-инвариантный левый интеграл $t \in \mathcal O(G)^*$, что $t(1)=1$.
\end{example}
\begin{proof}
Все рациональные представления
группы $G$ вполне приводимы (см. \cite{Nagata}). 
Отсюда в силу теоремы 4.2.12 из~\cite{Abe} алгебра Хопфа $\mathcal O(G)$ кополупроста.
Следовательно, существует такой левый интеграл $t \in \mathcal O(G)^*$, что $t(1)=1$. 
Его $\ad$-инвариантность следует из коммутативности алгебры Хопфа $\mathcal O(G)$.
\end{proof}

Приведём теперь пример алгебры Хопфа, у которой не существует ненулевых интегралов.

\begin{example}\label{ExampleIntegralUnivEnv}
Пусть $L$~"--- алгебра Ли над полем $\mathbbm{k}$ характеристики $0$ и $t \in U(L)^*$ является
левым интегралом, тогда $t=0$.
\end{example}
\begin{proof}
В силу теоремы Пуанкаре~"--- Биркгофа~"--- Витта
достаточно показать, что $t(v_1^{m_1}\ldots v_k^{m_k})=0$
для всех линейно независимых
элементов $v_1, \ldots, v_k \in L$ и всех $m_1, \ldots, m_k \geqslant 0$,
$k \geqslant 0$. Фиксируем элементы $v_1, \ldots, v_k \in L$ и 
введём лексикографическое упорядочение $\prec$ на  наборах $(m_1, \ldots, m_k)$.
Будем доказывать утверждение индукцией по этому упорядочению.

Сперва заметим, что $t(v)1 = t(1)v+t(v)1$ для всех $v \in L$. Следовательно, $t(1)=0$
и база индукции доказана. Пусть $m_1, \ldots, m_k \geqslant 0$.
Введём обозначения $$\Lambda := \lbrace (\ell_1, \ldots, \ell_k) \mid 0 \leqslant \ell_i \leqslant m_i,
\ 1 \leqslant i \leqslant k-1;\ 0 \leqslant \ell_k \leqslant m_k+1 \rbrace$$
и $$\Lambda_1 := \Lambda\backslash \lbrace (m_1, \ldots, m_k+1),  (m_1, \ldots, m_k) \rbrace
= \lbrace  (\ell_1, \ldots, \ell_k) \in \Lambda \mid (\ell_1, \ldots, \ell_k) \prec (m_1, \ldots, m_k)\rbrace.$$
Тогда \begin{equation*}\begin{split}
t(v_1^{m_1}\ldots v_k^{m_k+1})1=t\left((v_1^{m_1}\ldots v_k^{m_k+1})_{(2)}\right)
(v_1^{m_1}\ldots v_k^{m_k+1})_{(1)}=\\=\sum_{(\ell_1, \ldots, \ell_k)\in\Lambda} \tbinom{m_1}{\ell_1}\tbinom{m_2}{\ell_2}\ldots \tbinom{m_{k-1}}{\ell_{k-1}}
\tbinom{m_k+1}{\ell_k}\ t(v_1^{\ell_1}\ldots v_k^{\ell_k})
\ v_1^{m_1-\ell_1}\ldots v_k^{(m_k+1)-\ell_k}
= \\
=t(v_1^{m_1}\ldots v_k^{m_k+1})1+(m_k+1)t(v_1^{m_1}\ldots v_k^{m_k})v_k+
\\
+\sum_{(\ell_1, \ldots, \ell_k)\in\Lambda_1} \tbinom{m_1}{\ell_1}\tbinom{m_2}{\ell_2}\ldots \tbinom{m_{k-1}}{\ell_{k-1}}
\tbinom{m_k+1}{\ell_k}\ t(v_1^{\ell_1}\ldots v_k^{\ell_k})
\ v_1^{m_1-\ell_1}\ldots v_k^{(m_k+1)-\ell_k} = \\
=t(v_1^{m_1}\ldots v_k^{m_k+1})1+(m_k+1)t(v_1^{m_1}\ldots v_k^{m_k})v_k+0 \end{split}
\end{equation*}
в силу предположения индукции.
Следовательно, $t(v_1^{m_1}\ldots v_k^{m_k})=0$.
\end{proof}

\section{Градуировки, их эквивалентность и универсальные группы}\label{SectionGradEquivUnivGroup}
         
В данном параграфе мы следуем монографиям~\cite{ElduqueKochetov, KelarevBook}.
         
         Пусть $\Gamma \colon A=\bigoplus_{t\in T} A^{(t)}$~"--- разложение (необязательно ассоциативной) алгебры $A$
над некоторым полем $\mathbbm{k}$ в прямую сумму подпространств, проиндексированных элементами некоторого множества $T$. Говорят, что $\Gamma$~"--- \textit{градуировка}, если для любых элементов $t_1, t_2 \in T$ существует такой элемент $t \in T$, что $A^{(t_1)}A^{(t_2)}\subseteq  A^{(t)}$. Если на алгебре $A$ задана градуировка множеством $T$, 
алгебра $A$ называется \textit{$T$-градуированной}. Подпространства $A^{(t)}$, где $t\in T$, называются
\textit{однородными компонентами} или \textit{компонентами градуировки}.

Если $T$ наделена структурой (полу)группы и для всех $t_1, t_2 \in T$ выполнено
$A^{(t_1)}A^{(t_2)}\subseteq A^{(t_1 t_2)}$, то градуировка называется \textit{(полу)групповой}.

Пусть $V\subseteq A$~"--- некоторое подпространство градуированной алгебры $A$. Говорят, что $V$~"--- \textit{градуированное} или \textit{однородное} относительно градуировки $\Gamma$ подпространство, если 
$V=\bigoplus_{t\in T} V \cap A^{(t)}$. Аналогично вводятся понятия \textit{однородных (градуированных) идеалов} и \textit{подалгебр}.

Элементы подмножества $\bigcup_{t\in T} A^{(t)}$ называются \textit{градуированными} или \textit{однородными}. Если $a\in A^{(t)}$, то говорят, что $a$~"--- \textit{элемент $T$-степени $t$}.
 
Пусть алгебры $A_1$ и $A_2$ градуированы множествами $T_1$ и $T_2$, соответственно.
Говорят, что гомоморфизм алгебр $\varphi \colon A_1\to A_2$ \textit{(нестрого) градуированный},
если он переводит однородные элементы в однородные, т.е. для всякого $t_1\in T_1$
существует такое $t_2 \in T_2$, что $\varphi\left(A_1^{(t_1)}\right) \subseteq A_2^{(t_2)}$.
Говорят, что $\varphi$ \textit{строго градуированный},
если $T_1=T_2$ и $\varphi\left(A_1^{(t)}\right) \subseteq A_2^{(t)}$ для всех $t\in T_1$.

Важнейшим примером групповой градуировки является элементарная градуировка на матричной алгебре:
\begin{definition}
Пусть $\mathbbm{k}$~"--- поле, $G$~"--- группа, $n\in\mathbb N$, а $(g_1, \ldots, g_n)$ 
некоторый набор из $n$ элементов группы $G$.
Определим градуировку на алгебре  $M_n(\mathbbm{k})$ матриц размера $n\times n$,
потребовав, чтобы для всех $1\leqslant i,j \leqslant n$ матричная единица $e_{ij}$ принадлежала
 однородной компоненте, отвечающей элементу $g_i g_j^{-1}$.
Такая градуировка называется \textit{элементарной $G$-градуировкой, заданной набором $(g_1, \ldots, g_n)$}.
\end{definition}
\begin{remark}\label{RemarkElementary}
Заметим, что элементарная градуировка однозначно определяется $G$-степенями элементов $e_{i,i+1}$, $1\leqslant i \leqslant n-1$, так как эти матричные единицы вместе с матричными единицами $e_{i+1,i}$ порождают алгебру $M_n(\mathbbm{k})$. Обратно, если $G$~"--- произвольная группа, а $(h_1, \ldots, h_{n-1})$~"--- произвольный набор из $(n-1)$ её элемента, тогда элементарная градуировка, обладающая свойством $e_{i,i+1} \in M_n(\mathbbm{k})^{(h_i)}$ может быть задана набором $(g_1, \ldots, g_n)$, где $g_i = \prod_{j=i}^{n-1} h_j$.
\end{remark}

         При работе с градуировками важно определиться, когда мы считаем две градуировки схожими. Самое узкое определение, которое здесь можно дать~"--- это определение изоморфизма градуировок.
         
        Пусть
\begin{equation}\label{EqTwoSemiGroupGradings}\Gamma_1 \colon A_1=\bigoplus_{t_1 \in T_1} A_1^{(t_1)},\qquad \Gamma_2 \colon A_2=\bigoplus_{t_2 \in T_2} A_2^{(t_2)}\end{equation}
~"--- две градуировки, где $A_1$ и $A_2$~"--- алгебры, а $T_1$ и $T_2$~"--- два множества.

         \begin{definition} 
Пусть $\varphi \colon A_1 \mathrel{\widetilde\to} A_2$~"--- изоморфизм алгебр.         
         Говорят, что $\varphi$~"--- \textit{изоморфизм градуировок}, если
         $T_1 = T_2$ и $\varphi\left(A_1^{(t)}\right)\subseteq A_2^{(t)}$ для всех $t\in T_1$.
         В этом случае говорят, что алгебры $A_1$ и $A_2$ \textit{градуированно изоморфны}.
         \end{definition}
         
         Для многих приложений~\cite{AljaGia, AljaGiaLa,
BahtZaiGradedExp, GiaLa, ASGordienko9} бывает несущественным, элементами какого конкретного множества
или (полу)группы алгебра градуирована. (Здесь следует, однако отметить, что то, является ли $T$ группой или полугруппой, оказывает существенное влияние на структуру градуированной алгебры, см. главу~\ref{ChapterGradSemigroup}.) Замена градуирующего множества оставляет градуированные подпространства и идеалы градуированными. В случае градуированных тождеств (см. главу~\ref{ChapterGenHAssocCodim}) эта замена просто приводит к замене меток у переменных, а все числовые характеристики остаются прежними. Таким образом, оказывается полезным следующее более широкое определение:
       \begin{definition}[\cite{ElduqueKochetov}]
Пусть заданы две градуировки~\eqref{EqTwoSemiGroupGradings}  и изоморфизм алгебр $\varphi \colon A_1 \mathrel{\widetilde\to} A_2$.         
         Говорят, что $\varphi$~"--- \textit{эквивалентность градуировок}, если
         для любого $t_1 \in T_1$ с $A_1^{(t_1)}\ne 0$ существует такое $t_2\in T_2$,
         что $\varphi\left(A_1^{(t_1)}\right)=A_2^{(t_2)}$.
         В этом случае говорят, что градуировки $\Gamma_1$ и $\Gamma_2$ \textit{эквивалентны при помощи изоморфизма $\varphi$}.
         \end{definition}
          \begin{definition}
         Если $A_1=A_2$ и $\varphi=\id_{A_1}$, то говорят, что множества $T_1$ и $T_2$ \textit{реализуют на $A_1$ одну и ту же градуировку}, подразумевая здесь под градуировкой лишь разложение алгебры $A_1$ в прямую сумму градуированных компонент.
  \end{definition}
          \begin{definition}\label{DefSupport}
         \textit{Носителем} градуировки $\Gamma \colon A=\bigoplus_{t\in T} A^{(t)}$ называется
         множество $\supp \Gamma := \lbrace t\in T \mid  A^{(t)}\ne 0\rbrace$.
  \end{definition}
  Очевидно что, всякая эквивалентность $\varphi$ градуировок $\Gamma_1$ и $\Gamma_2$ определяет
  биекцию $\psi \colon \supp \Gamma_1 \mathrel{\widetilde\to} \supp \Gamma_2$,
  где $\varphi\left(A_1^{(t)}\right)=A_2^{\left(\psi(t)\right)}$ для всех $t\in \supp \Gamma_1$.
  
  Если $\varphi$~"--- эквивалентность градуировок $\Gamma_1$ и $\Gamma_2$, то можно
определить градуировку $\Gamma_3 \colon A_2=\bigoplus_{t_1 \in T_1} \varphi\left(A_1^{(t_1)}\right)$. 
При этом $\varphi$ будет являться изоморфизмом градуировок $\Gamma_1$ и $\Gamma_3$, а $\Gamma_2$
и $\Gamma_3$ будут эквивалентны при помощи тождественного изоморфизма алгебры $A_2$, т.е. $\Gamma_2$
будет получаться из $\Gamma_3$ переименованием однородных компонент при помощи элементов множества $T_2$.
 Таким образом, всякая эквивалентность градуировок сводится к изоморфизму градуировок и переименованию однородных компонент и, с точностью до градуированного изоморфизма, всегда можно считать, что все градуировки, эквивалентные данной, заданы на одной и той же алгебре, а эквивалентности являются тождественными изоморфизмами.

Рассмотрим подробнее групповые градуировки. Среди всех групп, реализующих заданную градуировку, существует особая группа, называемая универсальной группой градуировки~\cite[определение~1.17]{ElduqueKochetov}, \cite{PZ89}.

\begin{definition}
Пусть $\Gamma$~"--- групповая градуировка, заданная на алгебре $A$, причём $\Gamma$ допускает реализацию
в качестве $G_\Gamma$-градуировки для некоторой группы $G_\Gamma$.
Обозначим через $\varkappa_\Gamma$ соответствующее вложение носителя $\supp \Gamma \hookrightarrow G_\Gamma$.
Говорят, что пара $(G_\Gamma,\varkappa_\Gamma)$~"--- \textit{универсальная группа градуировки},
если для любой реализации $\Gamma$ в качестве $G$-градуировки для некоторой группы $G$
с вложением носителя $\psi \colon \supp \Gamma \hookrightarrow G$ существует
единственный гомоморфизм групп $\varphi \colon G_\Gamma \to G$, такой, что
следующая диаграмма коммутативна:
$$\xymatrix{ \supp \Gamma \ar[r]^(0.6){\varkappa_\Gamma} \ar[rd]^\psi & G_\Gamma \ar@{-->}[d]^\varphi \\
& G
}
$$
 \end{definition}
 \begin{remark} Универсальная группа градуировки является начальным объектом категории
 $\mathcal C_\Gamma$, заданной следующим образом: \begin{itemize}
 \item объектами в $\mathcal C_\Gamma$ являются всевозможные пары $(G,\psi)$, где $G$~"--- группа, реализующая градуировку $\Gamma$, а $\psi \colon \supp \Gamma \hookrightarrow G$~"--- соответствующее вложение носителя; 
 \item множество морфизмов из объекта $(G_1,\psi_1)$ в объект $(G_2,\psi_2)$
 состоит из всех гомоморфизмов групп $f \colon G_1 \to G_2$, таких, что диаграмма 
 $$\xymatrix{ \supp \Gamma \ar[r]^(0.6){\psi_1} \ar[rd]^{\psi_2} & G_1 \ar[d]^f \\
& G_2
}
$$
коммутативна.
\end{itemize}
 \end{remark}

Пусть $X$~"--- множество. Обозначим через $\mathcal F(X)$ свободную группу с множеством
 свободных порождающих $X$. Легко видеть, что если $\Gamma \colon A = \bigoplus_{g\in G} A^{(g)}$~"--- градуировка группой $G$, то $$G_\Gamma \cong \mathcal F([\supp \Gamma])/N,$$ где $[\supp \Gamma]:=\lbrace [g] \mid g\in \supp \Gamma \rbrace$, а $N$~"--- нормальная подгруппа, порождённая словами $[g][h][gh]^{-1}$ для всех $g,h \in \supp \Gamma$, таких, что $A^{(g)}A^{(h)}\ne 0$.

По определению, универсальная группа градуировки $\Gamma$~"--- это пара $(G_\Gamma,\varkappa_\Gamma)$.
Из теоремы~\ref{TheoremGivenFinPresGroupExistence}, которую мы докажем в главе~\ref{ChapterGradEquiv},
следует, что первой компонентой этой пары может выступать любая \textit{конечно представленная} группа, т.е. группа, которую можно задать конечным числом порождающих и определяющих соотношений.
Более того, для такой группы можно выбрать в качестве $\Gamma$ элементарную градуировку на матричной
алгебре.

Дадим ещё несколько полезных определений:

\begin{definition}
Пусть $\Gamma_1 \colon A=\bigoplus_{g \in G} A^{(g)}$
и $\Gamma_2 \colon A=\bigoplus_{h \in H} A^{(h)}$~"--- некоторые градуировки,
где $G$ и $H$~"--- группы, а $A$~"--- алгебра.
Будем говорить, что $\Gamma_2$ \textit{грубее}, чем $\Gamma_1$,
а $\Gamma_1$ \textit{тоньше}, чем $\Gamma_2$,
 если для всех $g\in G$, таких, что $A^{(g)}\ne 0$, существует $h\in H$,
 такое, что $A^{(g)}\subseteq A^{(h)}$. В этом случае $\Gamma_2$ называется \textit{огрублением}
 градуировки $\Gamma_1$, а $\Gamma_1$ называется \textit{утончением} градуировки $\Gamma_2$.
Будем обозначать через $\pi_{\Gamma_1 \to \Gamma_2} \colon G_{\Gamma_1} \twoheadrightarrow G_{\Gamma_2}$
гомоморфизм групп, заданный условием $\pi_{\Gamma_1 \to \Gamma_2}(\varkappa_{\Gamma_1}(g))=\varkappa_{\Gamma_2}(h)$ для $g\in\supp\Gamma_1$ и $h\in\supp \Gamma_2$, таких, что $A^{(g)} \subseteq A^{(h)}$.\label{NotationPiCoarser}
\end{definition}

С лингвистической точки зрения правильней было бы говорить <<грубее или эквивалентна>> вместо просто <<грубее>> и  <<тоньше или эквивалентна>> вместо просто <<тоньше>>. Однако для краткости слова <<или эквивалентна>> будут опускаться. Легко видеть, что отношение <<тоньше>> задаёт на классе всех градуировок \textit{предпорядок}, т.е. рефлексивное и транзитивное отношение, а сопоставление градуировке $\Gamma$ её универсальной группы
$G_\Gamma$ является функтором из этого предпорядка в категорию групп.
Этот функтор сопоставляет всякой паре градуировок
$\Gamma_1$ и $\Gamma_2$, таких, что
$\Gamma_1$ тоньше, чем $\Gamma_2$,
 гомоморфизм
$\pi_{\Gamma_1 \to \Gamma_2}$.

\section{(Ко)модули и (ко)модульные алгебры}\label{Section(Co)modules}

\subsection{Модули и комодули}

Подобно тому, как определение ассоциативной алгебры с единицей было дано на языке линейных отображений, дадим определение модуля над алгеброй.

\begin{definition}
Пусть $M$~"--- векторное пространство, $(A, \mu)$~"--- ассоциативная алгебра,
а $\psi \colon A \otimes M \to M$~"--- линейное отображение.
Говорят, что $(M,\psi)$~"--- \textit{(левый) модуль} над алгеброй $A$,
если коммутативна следующая диаграмма:
$$\xymatrix{ A \otimes A \otimes M \ar[d]_{\mu \otimes \id_M} \ar[rr]^{\id_A \otimes \psi} && A\otimes M
\ar[d]_\psi \\
A\otimes M \ar[rr]^\psi && M  } $$ 
 \end{definition}

Действие алгебры $A$ на $M$ обозначается через $am := \psi(a\otimes m)$ для всех $a\in A$, $m\in M$.
 
\begin{definition} Пусть $(A, \mu, u)$~"--- ассоциативная алгебра  с единицей над полем $\mathbbm{k}$.
Говорят, что $(M,\psi)$~"--- \textit{унитальный модуль} над $(A, \mu, u)$,
если $(M,\psi)$~"--- модуль над  $(A, \mu)$ и, кроме того,
коммутативна диаграмма:
$$\xymatrix{  \mathbbm{k} \otimes M \ar[rr]^{u \otimes \id_M} \ar[rd]^{\sim} 
  & & A\otimes M \ar[ld]^{\mu}\\
 & M  & }$$
(Здесь мы использовали естественное отождествление
 $ \mathbbm{k} \otimes M \cong M$.)
  \end{definition}

 В дальнейшем под модулями над ассоциативными алгебрами с единицей будут пониматься именно унитальные модули.
 
 В случае модулей над алгебрами Ли согласованность действия алгебры Ли и её операции коммутатора выглядит следующим образом:
 
 \begin{definition}
Пусть $V$~"--- векторное пространство, а $L$~"--- алгебра Ли над полем $\mathbbm{k}$.
Говорят, что $V$~"---  \textit{(левый) модуль} над $L$,
если задано билинейное отображение $L\times V \to V$, $(a,v)\mapsto av$, такое, что
$[a,b]v=a(bv)-b(av)$ для всех $a,b\in L$, $v\in V$. 
 \end{definition}
 
Эквивалентным языком является является язык линейных представлений алгебр Ли.
\textit{Линейным представлением} $\psi$ алгебры Ли $L$ на векторном пространстве
$V$ над полем $\mathbbm{k}$ называется гомоморфизм алгебр Ли $L \to \mathfrak{gl}(V)$.
(Напомним, что алгебра Ли $\mathfrak{gl}(V)$ совпадает как множество с множеством $\End_\mathbbm{k}(V)$ линейных операторов $V\to V$, а коммутатор на $\mathfrak{gl}(V)$ определяется по формуле $[\mathcal A, \mathcal B]:= \mathcal A \mathcal B - \mathcal B \mathcal A$ для всех $\mathcal A, \mathcal B \in \End_\mathbbm{k}(V)$.)
Ясно, что любое представление $\psi \colon L \to \mathfrak{gl}(V)$ задаёт на $V$ структуру
$L$-модуля через равенство $av := \psi(a)v$, где $a\in L$, $v\in V$.

Важнейшим примером представления всякой алгебры Ли $L$ является \textit{присоединённое представление}
$\ad \colon L \to \mathfrak{gl}(L)$, где $(\ad a)(b):=[a,b]$ для всех $a,b\in L$.

Теперь дадим определение комодуля над коалгеброй:

\begin{definition} Пусть $(C, \Delta, \varepsilon)$~"--- коалгебра на полем $\mathbbm{k}$.
Пара $(M,\rho)$, где $M$~"--- векторное пространство над $\mathbbm{k}$,
$\rho \colon M \to M \otimes C$~"--- линейное отображение,
называется \textit{(правым) комодулем} над коалгеброй $(C, \Delta, \varepsilon)$,
если следующие диаграммы коммутативны:
 $$\xymatrix{ M \ar[d]^\rho \ar[rr]^\rho && M\otimes C
\ar[d]^{\id_M \otimes \Delta} \\
M\otimes C \ar[rr]_{\rho \otimes \id_C} && M \otimes C \otimes C  } 
\text{\qquad и \qquad}\xymatrix{  & M \ar[ld]_{\rho} \ar[rd]^{\sim} &  \\
M \otimes C \ar[rr]^{\id_M \otimes \varepsilon} && M \otimes \mathbbm{k}
}$$
(Здесь мы использовали естественное отождествление $M \otimes \mathbbm{k} \cong M$.)
\end{definition}

При работе с комодулями в работе используются обозначения Свидлера: $\rho(m) = m_{(0)}\otimes m_{(1)}$, где $m\in M$
и опущен знак суммы. 

\begin{example}\label{ExampleRegularAffActAutComodule}  Пусть $G$~"--- аффинная алгебраическая группа над алгебраически замкнутым полем $\mathbbm{k}$
и пусть $\mathcal O(G)$~"--- алгебра регулярных
функций на $G$, которая, как было отмечено в примере~\ref{ExampleOAffAlgGrp}, является алгеброй Хопфа.
Пусть на некотором конечномерном векторном пространстве $V$ над $\mathbbm{k}$ задано \textit{рациональное} представление
группы $G$, т.е. 
 для некоторого базиса $v_1, \ldots, v_n$ пространства $V$ существуют элементы $\omega_{ij} \in \mathcal O(G)$, где $1\leqslant i,j \leqslant n$, такие, что $g v_j = \sum_{i=1}^n \omega_{ij}(g) v_i$
для всех $1\leqslant j \leqslant n$ и $g\in G$.
Тогда $V$ является $\mathcal O(G)$-комодулем, где $\rho(v_j) := \sum_{i=1}^n 
 v_i \otimes \omega_{ij}$ при $1\leqslant j \leqslant n$.
 \end{example}

Если $M$~"--- правый $C$-комодуль, то $M$ является левым $C^*$-модулем,
где $$c^*m:=c^*(m_{(1)})m_{(0)}\text{ для всех }c^*\in C^*\text{ и }m\in M.$$

Действие алгебры $C^*$ позволяет описать подкомодули, порождённые заданными подпространствами:

  \begin{lemma}\label{LemmaGenSubcomodules} Пусть $V$~"--- подпространство $C$-комодуля $M$, где $C$~"--- коалгебра над полем $\mathbbm{k}$.
  Тогда подпространство $C^*V := \langle c^* v \mid c^*\in C^*,\ v\in V\rangle_\mathbbm{k}$ является $C$-подкомодулем.
  \end{lemma}
 \begin{proof} Достаточно доказать, что для любых $c^*\in C^*$ и $v\in V$
 выполнено $\rho(c^*v) \in C^*V \otimes C$. Выберем базис $(c_\alpha)_{\alpha}$ в коалгебре $C$
 и фиксируем произвольные элементы $c^*\in C^*$ и $v\in V$.
 Тогда существуют такие элементы $m_{\alpha\beta} \in M$ (из которых только конечное число ненулевые), что
 $$(\rho \otimes \id_C)\rho(v)=\sum_{\alpha,\beta} m_{\alpha\beta} \otimes c_\alpha \otimes c_\beta.$$
 
 Определим линейные функции $c^\alpha \in C^*$ при помощи формул 
 $$c^\alpha(c_\beta):=\left\lbrace \begin{array}{ccc} 
0 & \text{ при } & \alpha \ne \beta, \\
1 & \text{ при } & \alpha = \beta.
 \end{array}\right.$$
Тогда для всех $\alpha$ и $\beta$
имеем $$m_{\alpha\beta}=(\id_M \otimes c_\alpha \otimes c_\beta)(\rho \otimes \id_C)\rho(v)=
(\id_M \otimes c_\alpha \otimes c_\beta)(\id_M \otimes \Delta)\rho(v)=
(c_\alpha c_\beta)v \in C^*V.$$ 
Отсюда $\rho(c^*v)= \sum_{\alpha,\beta} c^*(c_\beta)m_{\alpha\beta} \otimes c_\alpha \in C^*V\otimes C$. 
 \end{proof}

Если $M$~"--- левый модуль над конечномерной алгеброй $A$, то $M$~"--- правый $A^*$-комодуль,
где соответствующее отображение $\rho \colon M \to M \otimes A^*$ определяется равенством
 $$am = m_{(1)}(a)\, m_{(0)}\text{ для всех }a\in A\text{ и }m\in M.$$

\subsection{Модульные и комодульные алгебры}\label{SubsectionModComodAlg}

Наличие в биалгебре коумножения и коединицы позволяет формулировать дополнительные условия
на действие такой биалгебры на некоторой алгебре:

\begin{definition}
Говорят, что (необязательно ассоциативная) алгебра $A$ над полем $\mathbbm{k}$
является \textit{(левой) $H$-модульной алгеброй} для биалгебры $H$,
если $A$~"--- это левый $H$-комодуль и 
\begin{equation}\label{EqModCompat} h(ab)=(h_{(1)}a)(h_{(2)}b) \text{
для всех }a,b\in A,\ h\in H.\end{equation}
\end{definition}

Обозначим через $\psi$ линейное отображение $H \otimes A \to A$,
заданный равенством $\psi(h\otimes a)=ha$ для всех $h\in H$ и $a\in A$.
Отображение $\psi$ называется \textit{$H$-модульной структурой}
или \textit{$H$-действием} на $A$.

Мы будем говорить, что $A$~"--- \textit{$H$-модульная алгебра с единицей},
если $A$~"--- алгебра с единицей $1_A$
и $h 1_A = \varepsilon(h) 1_A$ для всех $h\in H$.

\begin{example}\label{ExampleFTmodule} Если $T$~"--- моноид, то $\mathbbm{k}T$-модульная алгебра~"--- это в точности алгебра,
на которой моноид $T$ действует эндоморфизмами. Если $G$~"--- группа, то $\mathbbm{k}G$-модульная алгебра~"--- это в точности алгебра, на которой группа $G$ действует автоморфизмами.
\end{example}

\begin{example}\label{ExampleUgModule} Если $L$~"--- алгебра Ли, то $U(L)$-модульная алгебра~"--- это в точности такая алгебра $A$, на которой алгебра Ли $L$ действует дифференцированиями,
т.е. $v(ab)=(va)b+a(vb)$ для всех $v\in L$ и $a,b\in A$.
\end{example}

\begin{example}\label{ExampleHModEnd} Пусть $M$~"--- $H$-модуль, где $H$~"--- алгебра Хопфа над полем $\mathbbm{k}$. Тогда алгебра
$\End_\mathbbm{k}(M)$ является ассоциативной $H$-модульной алгеброй с единицей, где \begin{equation}\label{EqHActionOnEnd}(h\varphi)(m):= h_{(1)}\varphi\bigl((Sh_{(2)})m\bigr)\text{ для всех }m\in M\text{ и }h\in H.\end{equation}
\end{example}

\begin{definition}
Говорят, что (необязательно ассоциативная) алгебра $A$ над полем $\mathbbm{k}$
является \textit{(правой) $H$-комодульной алгеброй} для биалгебры $H$,
если $A$~"--- это правый $H$-комодуль
и
 $(ab)_{(0)}\otimes (ab)_{(1)} = a_{(0)}b_{(0)}\otimes a_{(1)}b_{(1)}$
для всех $a,b \in A$.
\end{definition}

Отображение $\rho \colon A \to A \otimes
H$, задающее на $A$ структуру $H$-комодуля, называется \textit{$H$-комодульной структурой}
или \textit{$H$-кодействием} на $A$.

Мы будем говорить, что $A$~"--- \textit{$H$-комодульная алгебра с единицей},
если $A$~"--- алгебра с единицей $1_A$
и $\rho(1_A) = 1_A \otimes 1_H$.

\begin{example}\label{ExampleComoduleGraded} Если $T$~"--- моноид, то $\mathbbm{k}T$-комодульная алгебра~"--- это в точности алгебра,
градуированная моноидом $T$. При этом комодульная структура $\rho \colon A \to A \otimes \mathbbm{k}T$
на алгебре $A=\bigoplus_{t\in T} A^{(t)}$
задаётся следующим образом: $\rho(a):=a\otimes t$ для всех $a\in A^{(t)}$, где $t\in T$.
\end{example}

Для того, чтобы сформулировать второй пример, напомним, что 
 если $V$ и $W$~"--- векторные пространства над полем $\mathbbm{k}$, причём $\dim V < +\infty$,
то существует естественная линейная биекция
$$\Hom_\mathbbm{k}(V, V \otimes W) \cong V^* \otimes V \otimes W \cong \End_\mathbbm{k}(V)\otimes W.$$

\begin{example}\label{ExampleHComodEnd}
Пусть $M$~"--- конечномерный $H$-комодуль, где $H$~"--- алгебра Хопфа над полем $\mathbbm{k}$. Тогда алгебра
$\End_\mathbbm{k}(M)$ является ассоциативной $H$-комодульной алгеброй с единицей, где
комодульная структура задаётся равенством
 $$\varphi_{(0)}m \otimes \varphi_{(1)}
= \varphi(m_{(0)})_{(0)} \otimes \varphi(m_{(0)})_{(1)}(Sm_{(1)})\text{
для всех }\varphi \in \End_\mathbbm{k}(M)\text{ и }m\in M.$$
(Здесь мы пользуемся линейной биекцией $\Hom_\mathbbm{k}(M, M \otimes H) \cong \End_\mathbbm{k}(M)\otimes H$.)
\end{example}

Соответствие между модулями и комодулями переносится и на модульные и комодульные алгебры:
\begin{itemize}
\item если $A$~"--- правая $H$-комодульная алгебра для некоторой биалгебры $H$, то $A$ является левой $H^\circ$-модульной алгеброй, где $h^\circ a:=h^\circ(a_{(1)})a_{(0)}$ для всех $h^\circ\in H^\circ$ и $a\in A$;
\item если $A$~"--- левая модульная алгебра над конечномерной биалгеброй $H$, то $A$~"--- правая $H^*$-комодульная алгебра,
где $H^*$-комодульная струкутра $\rho \colon A \to A \otimes H^*$ определяется равенством
 $ha = a_{(1)}(h)\, a_{(0)}$ для всех $h\in H$ и $a\in A$.
\end{itemize}

\begin{example}\label{ExampleRegularAffActAutAlgebra}  Пусть $G$~"--- аффинная алгебраическая группа над алгебраически замкнутым полем $\mathbbm{k}$
и пусть $\mathcal O(G)$~"--- алгебра Хопфа регулярных
функций на $G$.
Предположим, что $G$ действует рационально автоморфизмами на конечномерной алгебре $A$, т.е. для некоторого базиса $a_1, \ldots, a_n$ алгебры $A$ существуют элементы $\omega_{ij} \in \mathcal O(G)$, где $1\leqslant i,j \leqslant n$, такие, что $g a_j = \sum_{i=1}^n \omega_{ij}(g) a_i$
для всех $1\leqslant j \leqslant n$ и $g\in G$.
Тогда $A$ является $\mathcal O(G)$-комодульной алгеброй, где $\rho(a_j) := \sum_{i=1}^n 
 a_i \otimes \omega_{ij}$ при $1\leqslant j \leqslant n$.
  В то же время $A$~"--- $\mathcal O(G)^\circ$-модульная алгебра:
$f^* a_j = \sum_{i=1}^n f^*(\omega_{ij}) a_i$ для всех $1\leqslant j \leqslant n$ и $f^* \in \mathcal O(G)^\circ$.  
  Алгебра Ли $\mathfrak g$ группы $G$~"--- это подпространство алгебры Хопфа $\mathcal O(G)^\circ$,
  состоящее из всех её примитивных элементов, а $\mathfrak g$-действие на $A$
  дифференцированиями~"--- это просто ограничение 
  $\mathcal O(G)^\circ$-действия. В то же время, сама группа $G$
может быть отождествлена с множеством всех группоподобных элементов алгебры Хопфа.
    Таким образом, на $A$ действуют три алгебры Хопфа:
  $\mathcal O(G)^\circ$, $\mathbbm{k}G$ и $U(\mathfrak g)$.
\end{example}

Пусть $A$~"--- ассоциативная $H$-модульная алгебра для некоторой алгебры Хопфа $H$
над полем $\mathbbm{k}$. 
\textit{Смэш-произведением} $A\mathbin{\#}H$ называется векторное пространство
$A\otimes H$, для элементов которого мы используем обозначение $a\mathbin{\#} h := a\otimes h$, со структурой ассоциативной алгебры, заданной умножением $(a\mathbin{\#} h)(b\mathbin{\#} k):=a(h_{(1)})b \mathbin{\#} h_{(2)}k$
для всех $a,b\in A$, $h,k\in H$.

Пусть $G$~"--- группа, которая действует автоморфизмами на алгебре Ли $L$.
Тогда $U(L)$~"--- $\mathbbm{k}G$-модульная алгебра, а смэш-произведение $U(L)\mathbin{\#}\mathbbm{k}G$
является алгеброй Хопфа, где $U(L)\mathbin{\#}\mathbbm{k}G \cong U(L)\otimes \mathbbm{k}G$ как коалгебра, т.е. коумножение и коединица задаются
при помощи тензорного произведения соответствующих отображений коалгебр $U(L)$ и $\mathbbm{k}G$.
Антипод $S$ определяется равенством $S(w\mathbin{\#} g):= (1 \mathbin{\#} g^{-1})(Sw\mathbin{\#} 1)$
при $g\in G$ и $w\in U(L)$. 

Для кокоммутативных алгебр Хопфа над алгебраически замкнутым полем характеристики $0$
справедлива следующая структурная теорема:

\begin{theorem}[Габриэль--Картье--Костант--Милнор--Мур]\label{TheoremCartierGabrielKostant}
Пусть $H$~"--- кокоммутативная алгебра Хопфа над алгебраически замкнутым полем $\mathbbm{k}$ характеристики $0$,
причём $G$~"--- группа её обратимых элементов, а $L$~"--- алгебра Ли её примитивных элементов,
на которой группа $G$ действует сопряжениями.
Тогда $H \cong U(L)\mathbin{\#} \mathbbm{k}G$.
\end{theorem}
\begin{proof} С учётом предложения~\ref{PropositionCocommAlgClosedSimple} данной главы
утверждение теоремы следует, например,  из следствия~5.6.4 и теоремы~5.6.5 монографии \cite{Montgomery}.
\end{proof}

Нам также потребуется обратное утверждение:

\begin{proposition}\label{PropositionPrimitiveGrouplikeSmash}
Пусть $G$~"--- группа, которая действует автоморфизмами на алгебре Ли $L$ над полем $\mathbbm{k}$ характеристики $0$.
Тогда группа группоподобных элементов алгебры Хопфа $U(L)\mathbin{\#} \mathbbm{k}G$ совпадает с подмножеством $1 \otimes G$,
а алгебра Ли примитивных элементов~~"--- с подмножеством $L \otimes 1$. 
\end{proposition}
\begin{proof}
Пусть $h\in U(L)\mathbin{\#} \mathbbm{k}G$. Тогда существуют однозначно определённые элементы $w_g \in U(L)$,
только конечное число из которых ненулевые, такие, что $h = \sum_{g\in G} w_g \otimes g$.
Тогда $\Delta h = \sum_{g\in G} (w_g)_{(1)} \otimes g \otimes (w_g)_{(2)} \otimes g$.

Если $w_{g_1}\ne 0$ и $w_{g_2}\ne 0$ для каких-то $g_1 \ne g_2$,
то разложении элемента $h\otimes h$ обязательно будет присутствовать элемент $w_{g_1}\otimes g_1 \otimes w_{g_2}\otimes g_2$, которого не будет в разложении для $\Delta h$.
Отсюда $\Delta h = h \otimes h$, если и только если $w_g = 0$ для всех $g\in G$, кроме, быть может, одного  элемента $g_0 \in G$, и для этого элемента справедливо равенство $\Delta w_{g_0} = w_{g_0}\otimes w_{g_0}$.

Аналогично, равенство $\Delta h = h \otimes 1 + 1 \otimes h$ справедливо, если и только если
$h=w\otimes 1$, где $w$~"--- примитивный элемент алгебры Хопфа $U(L)$.
Теперь доказываемое утверждение следует из предложения~5.5.3 монографии~\cite{Montgomery}.
\end{proof}

\section{Сопряжённые функторы, единицы и коединицы сопряжения}\label{SectionAdjoint}

Мы предполагаем знакомство читателя с базовыми понятиями теории категорий, такими как
\textit{категория}, \textit{функтор} и \textit{естественное преобразование функторов}.
(См., например, \cite{MacLaneCatWork}.)

Пусть $F \colon \mathcal X \to \mathcal A$ и $G \colon \mathcal X \to \mathcal A$~"---
функторы. Говорят, что функтор $F$~"--- \textit{левый сопряжённый} для функтора  $G$,
а функтор $G$~"--- \textit{правый сопряжённый} для функтора  $F$,
если существует биекция \begin{equation}\label{EqAdjunction}\mathcal A(FX,A)\mathrel{\widetilde\to}\mathcal X(X,GA),\end{equation}
где $X$ и $A$, соответственно, объекты категорий $\mathcal X$ и $\mathcal A$,
причём эта биекция является естественным преобразованием функторов, если рассматривать
$\mathcal A(F-,-)$ и $\mathcal X(-,G-)$ как функторы
$\mathcal X\times \mathcal A \to \mathbf{Sets}$.
Используется обозначение $F \dashv G$.

Для всякого объекта $X$ категории $\mathcal X$ обозначим через $i_X \colon X \to GFX$ элемент множества $\mathcal X(X,GFX)$,
соответствующий при этой биекции элементу $\id_{FX}\in \mathcal A(FX,FX)$.
Тогда естественное преобразование $i \colon \id_{\mathcal X} \Rightarrow GF$ называется
\textit{единицей сопряжения} $F \dashv G$.

Для всякого объекта $A$ категории $\mathcal A$ обозначим через $\nu_A \colon FGA \to A$ элемент множества $\mathcal A(FGA,A)$,
соответствующий элементу $\id_{GA}\in \mathcal X(GA,GA)$.
Естественное преобразование $\mu \colon FG \Rightarrow \id_{\mathcal A}$ называется
\textit{коединицей сопряжения} $F \dashv G$.

Единица и коединица сопряжения удовлетворяют следующим универсальным свойствам:
\begin{enumerate}\item для любого морфизма  $f \colon FX \to A$ существует единственный морфизм
$g \colon X \to GA$, такой, что  коммутативна диаграмма
 $$\xymatrix{FX \ar[rd]_f \ar@{-->}[r]^{Fg} &  FGA \ar[d]_(0.4){\nu_A} \\
                    & A }$$
                    \item для любого морфизма $g \colon X \to GA$ существует единственный морфизм
 $f \colon FX \to A$, такой, что коммутативна диаграмма
 $$\xymatrix{X \ar[rd]_g \ar[r]^{i_X} &  GFX \ar@{-->}[d]_(0.4){Gf} \\
                    & GA }$$
                    \end{enumerate}
 При этом $f\in \mathcal A(FX,A)$ и
 $g\in\mathcal X(X,GA)$ соответствуют друг другу при биекции~\eqref{EqAdjunction}.
 
 Обратно, наличие естественного преобразования, удовлетворяющего одному из универсальных свойств, приведённых выше, говорит о наличии сопряжения между функторами.

Приведём важнейшие примеры пар сопряженных функторов, встречающиеся в работе.

\begin{example} Пусть $\mathbf{nuAlg}_\mathbbm{k}$~"--- категория ассоциативных алгебр над полем $\mathbbm{k}$, необязательно
содержащих единицу. Обозначим через $U$ \textit{забывающий} функтор $\mathbf{nuAlg}_\mathbbm{k}\to \mathbf{Sets}$,
т.е. функтор, который ставит в соответствие каждой алгебре множество её элементов.
Для всякого множества $X$ определим \textit{свободную ассоциативную алгебру $\mathbbm{k}\langle X \rangle$ без единицы} на множестве $X$ как алгебру многочленов от некоммутирующих переменных из множества $X$ без свободного члена.
\label{DefFX}
Всякое отображение $X\to A$, где $A$~"--- ассоциативная алгебра над полем $\mathbbm{k}$, однозначно продолжается до гомоморфизма алгебр $\mathbbm{k}\langle X \rangle \to A$.
Отсюда имеем естественную биекцию $\mathbf{nuAlg}_\mathbbm{k}(\mathbbm{k}\langle X \rangle, A)\mathrel{\widetilde\to} \mathbf{Sets}(X,UA)$, откуда функтор $\mathbbm{k}\langle - \rangle$ является левым сопряжённым к функтору $U$.
Единицей этого сопряжения является вложение $X \hookrightarrow \mathbbm{k}\langle X \rangle$.
\end{example}

Для забывающего функтора  $\mathbf{Alg}_\mathbbm{k}\to \mathbf{Sets}$ левым сопряжённым является функтор,
сопоставляющий каждому множеству его \textit{свободную ассоциативную алгебру с $1$}, т.е. алгебру многочленов со свободным членом от некоммутирующих переменных из этого множества.

\begin{example}\label{ExampleTensorAlg} Левым сопряжённым к забывающему функтору $\mathbf{Alg}_\mathbbm{k}\to \mathbf{Vect}_\mathbbm{k}$ является функтор $T\colon\mathbf{Vect}_\mathbbm{k}\to \mathbf{Alg}_\mathbbm{k}$,
который ставит в соответствие всякому векторному пространству
\textit{тензорную алгебру} $T(V):= \bigoplus_{n=0}^\infty V^{\otimes n}$, где $V^{\otimes 0}:= \mathbbm{k}$,
а умножение определяется при помощи равенства $uv := u\otimes w$
для всех $u\in V^{\otimes k}$, $w\in V^{\otimes m}$, $m,n\in\mathbb Z_+$. Алгебра $T(V)$ изоморфна свободной ассоциативной алгебре с $1$, построенной по базису пространства $V$.\label{DefTV}
Единицей сопряжения является вложение $V \hookrightarrow T(V)$.
\end{example}

\begin{example}\label{ExampleCofreeCoalg} Забывающий функтор $\mathbf{Coalg}_\mathbbm{k}\to \mathbf{Vect}_\mathbbm{k}$ имеет правый сопряжённый $K\colon\mathbf{Vect}_\mathbbm{k}\to \mathbf{Alg}_\mathbbm{k}$ (см., например,
\cite[\S 1.6]{Danara}).
Коалгебра $K(V)$ называется \textit{косвободной коалгеброй} пространства $V$
и удовлетворяет следующему универсальному свойству:
для любого линейного отображения $f \colon C \to V$
некоторой коалгебры $C$ в векторное пространство $V$ существует единственный гомоморфизм
 $g \colon C \to K(V)$ коалгебр, такой, что коммутативна диаграмма
 $$\xymatrix{C \ar[rd]_f \ar@{-->}[r]^(0.4){g} &  K(V) \ar[d]^{{\nu\,}_V} \\
                    & V }$$
(Здесь ${\nu\,}_V \colon K(V) \to V$~"--- коединица сопряжения.)
\end{example}

\section{Сопряжённые функторы, связанные с градуировками}\label{SectionGrAdjoint}

В главе~\ref{ChapterGradEquiv} будут рассматриваться категории и функторы, связанные с понятием 
эквивалентности градуировок. Для того, чтобы дать читателю возможность сравнить конструкции главы~\ref{ChapterGradEquiv} с уже известными конструкциями, приведём пример двух пар сопряжённых функторов, связанных с градуировками на алгебрах. Ниже для определённости рассматриваются ассоциативные алгебры необязательно с единицей, хотя, конечно, аналоги таких сопряжений существуют и для категорий ассоциативных алгебр с единицей, а также для категорий неассоциативных алгебр.

Обозначим через $\mathbf{nuAlg}_\mathbbm{k}^{G\text{-}\mathrm{gr}}$ 
категорию ассоциативных алгебр на полем $\mathbbm{k}$, необязательно содержащих единицу, градуированных
группой $G$.
Морфизмами категории $\mathbf{nuAlg}_\mathbbm{k}^{G\text{-}\mathrm{gr}}$
являются всевозможные гомоморфизмы алгебр
$$\psi \colon A=\bigoplus_{g\in G}A^{(g)}\ \longrightarrow\ B=\bigoplus_{g\in G}B^{(g)},$$
такие, что $\psi(A^{(g)})\subseteq B^{(g)}$ для всех $g\in G$.
Для произвольного гомоморфизма групп $\varphi \colon G \to H$
обозначим через $U_\varphi$ функтор $ \mathbf{nuAlg}_\mathbbm{k}^{G\text{-}\mathrm{gr}} \to \mathbf{nuAlg}_\mathbbm{k}^{H\text{-}\mathrm{gr}}$,
который сопоставляет всякой $G$-градуировке $A=\bigoplus_{g\in G}A^{(g)}$ на алгебре $A$ 
 градуировку $A = \bigoplus_{h\in H} A^{(h)}$ группой $H$, где $A^{(h)} := \bigoplus\limits_{\substack{g\in G,\\ \varphi(g)=h}} A^{(g)}$,
 и не затрагивает гомоморфизмы.

Функтор $U_\varphi  \colon \mathbf{nuAlg}_\mathbbm{k}^{G\text{-}\mathrm{gr}} \to \mathbf{nuAlg}_\mathbbm{k}^{H\text{-}\mathrm{gr}}$
обладает правым сопряжённым функтором $K_\varphi \colon \mathbf{nuAlg}_\mathbbm{k}^{H\text{-}\mathrm{gr}} \to \mathbf{nuAlg}_\mathbbm{k}^{G\text{-}\mathrm{gr}}$,
который определяется следующим образом:
для $H$-градуировки $B=\bigoplus_{h\in H} B^{(h)}$ 
положим
 $K_\varphi (B) := \bigoplus_{g\in G} \bigl(K_\varphi(B)\bigr)^{(g)}$,
 где $\bigl(K_\varphi(B)\bigr)^{(g)} := \lbrace (g, b) \mid b\in B^{\left(\varphi(g)\right)} \rbrace$.
Структура векторного пространства на $\bigl(K_\varphi(B)\bigr)^{(g)}$ переносится с векторного пространства $B^{\left(\varphi(g)\right)}$, а произведение
на $K_\varphi(B)$ задаётся равенством $(g_1,a)(g_2,b):=(g_1 g_2, ab)$ при $g_1, g_2\in G$, $a\in B^{\left(\varphi(g_1)\right)}$, $b\in B^{\left(\varphi(g_2)\right)}$. Для $\psi \in \mathbf{nuAlg}_\mathbbm{k}^{H\text{-}\mathrm{gr}}(B_1,B_2)$
морфизм $K_\varphi (\psi)\in \mathbf{nuAlg}_\mathbbm{k}^{G\text{-}\mathrm{gr}}(K_\varphi (B_1), K_\varphi (B_2))$ задаётся равенством $K_\varphi (\psi)(g,b):=(g,\psi(b))$
при $g\in G$, $b\in B_1^{\left(\varphi(g)\right)}$.
Существует естественная биекция
$$\mathbf{nuAlg}_\mathbbm{k}^{H\text{-}\mathrm{gr}} (U_\varphi (A), B)\to \mathbf{nuAlg}_\mathbbm{k}^{G\text{-}\mathrm{gr}} (A, K_\varphi(B)),$$
где $A$~"--- объект категории $\mathbf{nuAlg}_\mathbbm{k}^{G\text{-}\mathrm{gr}}$, а $B$~"--- объект категории $\mathbf{nuAlg}_\mathbbm{k}^{H\text{-}\mathrm{gr}}$ (см., например, \cite[\S~1.2]{NastasescuVanOyst}).

Другой пример сопряжения~"--- свободно-забывающее.

Пусть $G$~"--- группа, а $\mathbf{Sets}^{G\text{-}\mathrm{gr}}_*$~"--- категория,
объектами которой являются множества
$X$, содержащие выделенный элемент $0$, с фиксированным разложением $X=\lbrace 0 \rbrace \sqcup \bigsqcup_{g\in G} X^{(g)}$ в объединение непересекающихся множеств.
Морфизмами между множествами $X=\lbrace 0 \rbrace \sqcup \bigsqcup_{g\in G} X^{(g)}$ и
 $Y=\lbrace 0 \rbrace \sqcup \bigsqcup_{g\in G} Y^{(g)}$  
 являются такие отображения $\varphi \colon X \to Y$,
 что $\varphi(0)=0$ и $\varphi(X^{(g)})\subseteq \lbrace 0 \rbrace \sqcup Y^{(g)}$
 для всех $g\in G$.
Существует очевидный забывающий функтор $U \colon \mathbf{nuAlg}_\mathbbm{k}^{G\text{-}\mathrm{gr}} \to \mathbf{Sets}^{G\text{-}\mathrm{gr}}_*$, который сопоставляет
каждой градуированной алгебре $A$ объект $U(A)=\lbrace 0 \rbrace \sqcup \bigsqcup_{g\in G} \left(A^{(g)} \backslash \lbrace 0 \rbrace \right)$.
Функтор $U$ обладает левым сопряжённым функтором $\mathbbm{k}\langle (-)\backslash \lbrace 0 \rbrace  \rangle$,
который сопоставляет множеству $X=\lbrace 0 \rbrace \sqcup \bigsqcup_{g\in G} X^{(g)} \in \mathbf{Sets}^{G\text{-}\mathrm{gr}}_*$
свободную ассоциативную алгебру $\mathbbm{k}\langle X\backslash \lbrace 0 \rbrace  \rangle$ 
без единицы на множестве $X\backslash \lbrace 0 \rbrace$,
причём алгебра $\mathbbm{k}\langle X\backslash \lbrace 0 \rbrace  \rangle$
снабжена градуировкой, заданной условием $x_1 \ldots x_n \in
\mathbbm{k}\langle X\backslash \lbrace 0 \rbrace  \rangle^{(g_1 \ldots g_n)}$
при $x_i \in X^{(g_i)}$, $1\leqslant i \leqslant n$.
Существует естественная биекция
$$\mathbf{nuAlg}_\mathbbm{k}^{G\text{-}\mathrm{gr}} (\mathbbm{k}\langle X\backslash \lbrace 0 \rbrace\rangle, A)\to \mathbf{Sets}^{G\text{-}\mathrm{gr}}_* (X, UA),$$
где $A$ является объектом категории $\mathbf{nuAlg}_\mathbbm{k}^{G\text{-}\mathrm{gr}}$,
а $X$ является объектом категории $\mathbf{Sets}^{G\text{-}\mathrm{gr}}_*$.

Заметим, что для всех объектов любой категории, которая рассматривается в данном параграфе, фиксирована одна и та же градуирующая группа.

\section{(Ко)произведения в категориях $\mathbf{Alg}_\mathbbm{k}$, $\mathbf{Coalg}_\mathbbm{k}$, $\mathbf{Bialg}_\mathbbm{k}$}\label{SectionCategoryConstructions}

В этом и следующем параграфах мы вкратце изложим теоретические сведения, необходимые для построения
левого $H_l$ и правого $H_r$ сопряжённых функторов для функтора вложения категории $\textbf{Hopf}_\mathbbm{k}$
в категорию $\textbf{Bialg}_\mathbbm{k}$, которые нам потребуется в главе~\ref{ChapterOmegaAlg}. В данных параграфах
мы следуем работам~\cite{AgorePhDThesis, Chirvasitu, Pareigis, Takeuchi}.

Напомним, что если $A_\alpha$, где $\alpha \in \Lambda$,~"--- объекты некоторой категории $\mathcal C$,
то их \textit{произведением}\label{DefProd} называется такой объект $\prod_{\alpha \in \Lambda} A_\alpha$
вместе с морфизмами $\pi_\alpha \colon \prod_{\alpha \in \Lambda} A_\alpha \to A_\alpha$,
что для любого объекта $A$ и морфизмов $\varphi_\alpha \colon A \to A_\alpha$ существует единственный
морфизм $\varphi \colon A \to \prod_{\alpha \in \Lambda} A_\alpha$,
делающий для любого $\alpha \in \Lambda$ коммутативной диаграмму
 $$\xymatrix{A \ar[rd]_{\varphi_\alpha} \ar@{-->}[r]^(0.3){\varphi} &  \prod_{\alpha \in \Lambda} A_\alpha \ar[d]^{\pi_\alpha} \\
                    & A_\alpha }$$
                    
Определение копроизведения получается обращением стрелок.
Если $A_\alpha$, где $\alpha \in \Lambda$,~"--- объекты некоторой категории $\mathcal C$,
то их \textit{копроизведением}\label{DefCoprod} называется такой объект $\coprod_{\alpha \in \Lambda} A_\alpha$
вместе с морфизмами $i_\alpha \colon A_\alpha \to \coprod_{\alpha \in \Lambda} A_\alpha$,
что для любого объекта $A$ и морфизмов $\varphi_\alpha \colon A_\alpha \to A$ существует единственный
морфизм $\varphi \colon \coprod_{\alpha \in \Lambda} A_\alpha \to  A$,
делающий для любого $\alpha \in \Lambda$ коммутативной диаграмму
 $$\xymatrix{A_\alpha \ar[rd]_{\varphi_\alpha} \ar[r]^(0.4){i_\alpha} &  \coprod_{\alpha \in \Lambda} A_\alpha \ar@{-->}[d]^{\varphi} \\
                    & A }$$

В категории $\mathbf{Alg}_\mathbbm{k}$ произведением алгебр является их декартово произведение как множеств
с покомпонентными операциями, причём результатом применения отображения $\pi_\alpha$ является компонента с индексом $\alpha$ соответствующего элемента произведения.

Копроизведением алгебр $A_\alpha$ является алгебра $T\left(\bigoplus_{\alpha \in \Lambda} A_\alpha\right)/I$,
где $I$~"--- пересечение всех идеалов $Q$ алгебры $T\left(\bigoplus_{\alpha \in \Lambda} A_\alpha\right)$,
таких, что композиции всех вложений $A_\alpha \hookrightarrow T\left(\bigoplus_{\alpha \in \Lambda} A_\alpha\right)$
с естественным сюръективным гомоморфизмом $T\left(\bigoplus_{\alpha \in \Lambda} A_\alpha\right) \twoheadrightarrow T\left(\bigoplus_{\alpha \in \Lambda} A_\alpha\right)/Q$
являются гомоморфизмами алгебр.

Рассмотрим категорию $\mathbf{Coalg}_\mathbbm{k}$. Копроизведением коалгебр $C_\alpha$
является их прямая сумма $\bigoplus_{\alpha \in \Lambda} C_\alpha$ как векторных пространств,
причём коумножение и коединица продолжаются с соответствующих отображений на $C_\alpha$.

Произведением коалгебр $C_\alpha$
является сумма всех подкоалгебр $Q$ коалгебры $K\left(\prod_{\alpha \in \Lambda} C_\alpha\right)$
(здесь $\prod_{\alpha \in \Lambda} C_\alpha$~"--- декартово произведение векторных пространств),
таких, что для любого $\alpha \in \Lambda$ композиция вложения $Q \hookrightarrow
K\left(\prod_{\alpha \in \Lambda} C_\alpha\right)$ с коединицей сопряжения $${\nu\,}_{\prod_{\alpha \in \Lambda} C_\alpha} \colon K\left(\prod_{\alpha \in \Lambda} C_\alpha\right) \to \prod_{\alpha \in \Lambda} C_\alpha \quad\text{(см. пример~\ref{ExampleCofreeCoalg})}$$ и проекцией $\prod_{\alpha \in \Lambda} C_\alpha \twoheadrightarrow C_\alpha$ является гомоморфизмом коалгебр.

Рассмотрим теперь произведение и копроизведение в категории $\mathbf{Bialg}_\mathbbm{k}$.
Пусть $B_\alpha$, $\alpha \in \Lambda$,~"--- биалгебры. Произведением биалгебр $B_\alpha$
является их произведение $B$ как коалгебр,
причём умножение $\mu_B \colon B \otimes B \to B$ и единица $u_B \colon \mathbbm{k} \to B$
определяются как единственные отображения, делающие для любого $\alpha \in \Lambda$ коммутативными
диаграммы 
 $$\xymatrix{ B \otimes B \ar@{-->}[r]^(0.6){\mu_B} \ar[d]_{\pi_\alpha \otimes \pi_\alpha} & B \ar[d]_{\pi_\alpha}  \\
 B_\alpha \otimes B_\alpha \ar[r]_(0.6){\mu_{B_\alpha}} & B_\alpha
 } \text{\qquad и \qquad} \xymatrix{ \mathbbm{k} \ar@{-->}[r]^{u_B} \ar[rd]_{u_{B_\alpha}} & B \ar[d]^{\pi_\alpha}  \\
  & B_\alpha  }$$
  
  Копроизведением биалгебр $B_\alpha$, где $\alpha \in \Lambda$,
является их копроизведение $B$ как алгебр,
причём коумножение $\Delta_B \colon B  \to B\otimes B$ и коединица $\varepsilon_B \colon  B \to \mathbbm{k}$
определяются как единственные отображения, делающие для любого $\alpha \in \Lambda$ коммутативными
диаграммы 
 $$\xymatrix{ B_\alpha \ar[r]^(0.4){\Delta_{B_\alpha}}\ar[d]^{i_\alpha}  & B_\alpha \otimes B_\alpha \ar[d]^{i_\alpha \otimes i_\alpha} \\
 B  \ar@{-->}[r]_(0.4){\Delta_B}  & B\otimes B 
 } \text{\qquad и \qquad} \xymatrix{ B_\alpha \ar[r]^{i_\alpha}\ar[rd]_{\varepsilon_{B_\alpha}} &  B  \ar@{-->}[d]^{\varepsilon_B} \\
 & \mathbbm{k} }$$

\section{Функторы $H_l$ и $H_r$}\label{SectionFunctorsHlHr}

Напомним, что через $A^\mathrm{op}$ обозначается алгебра $A$, в которой в умножении меняются местами аргументы, а через $C^\mathrm{cop}$~"--- коалгебра, в которой в коумножении меняются местами тензорные множители в разложении результата.

Конструкция функтора  $H_l \colon \mathbf{Bialg}_\mathbbm{k} \to \mathbf{Hopf}_\mathbbm{k}$ восходит к работе М.~Такеучи~\cite{Takeuchi}. В явном виде она было изложена Б.~Парейгисом в~\cite{Pareigis}.

Для каждой биалгебры $B$ определим биалгебры $B_n$, где $n\in\mathbb Z_+$,
следующим образом: $B_{2k} \cong B$, а $B_{2k+1} \cong B^{\mathrm{op,cop}}$
для всех $k\in\mathbb Z_+$. 

На копроизведении биалгебр $\coprod_{n=0}^\infty  B_n$ определим гомоморфизм 
$$S \colon \coprod_{n=0}^\infty  B_n \to \left(\coprod_{n=0}^\infty  B_n\right)^{\mathrm{op,cop}}$$
как единственный гомоморфизм биалгебр, делающий коммутативными
диаграммы
 $$\xymatrix{  B_n \ar[r]^{\sim}\ar[d]^{i_n} &  B_{n+1}^{\mathrm{op,cop}}\ar[d]^{i_{n+1}} \\
   \coprod_{n=0}^\infty  B_n \ar@{-->}[r]^(0.4)S & \left(\coprod_{n=0}^\infty  B_n\right)^{\mathrm{op,cop}}}$$
 Пусть $I$~"--- идеал биалгебры $\coprod_{n=0}^\infty  B_n$,
 порождённый элементами $S(h_{(1)})h_{(2)}-\varepsilon(h)1$
 и $h_{(1)}S(h_{(2)})-\varepsilon(h)1$ для всех $h\in \coprod_{n=0}^\infty  B_n$. Тогда $I$~"--- биидеал, инвариантный относительно
 отображения $S$, и $H_l(B):= \left(\coprod_{n=0}^\infty  B_n\right)/I$~"--- алгебра Хопфа. Обозначим через $i_B \colon B \to H_l(B)$ отображение, индуцированное отображением
 $i_0 \colon B=B_0 \to \coprod_{n=0}^\infty  B_n$.
 Тогда $H_l \colon \mathbf{Bialg}_\mathbbm{k} \to \mathbf{Hopf}_\mathbbm{k}$ является функтором, а $i$ является естественным преобразованием тождественного функтора категории $\mathbf{Bialg}_\mathbbm{k}$
 в композицию функтора $H_l$ и функтора вложения $\mathbf{Hopf}_\mathbbm{k} \to \mathbf{Bialg}_\mathbbm{k}$.
 Естественное преобразование $i$ удовлетворяет следующему универсальному свойству:
 для любой алгебры Хопфа $H$ и гомоморфизма биалгебр $\varphi \colon B \to H$
 существует единственный гомоморфизм алгебр Хопфа $\varphi_0 \colon H_l(B) \to H$,
 делающий коммутативными диаграммы
 
 $$\xymatrix{ B \ar[r]^(0.4){i_B}\ar[rd]_{\varphi} & H_l(B) \ar@{-->}[d]^{\varphi_0} \\
 & H }$$
 
 Иными словами, функтор $H_l$~"--- левый сопряжённый к функтору вложения $\mathbf{Hopf}_\mathbbm{k} \to \mathbf{Bialg}_\mathbbm{k}$.
 
 Функтор  $H_r \colon \mathbf{Bialg}_\mathbbm{k} \to \mathbf{Hopf}_\mathbbm{k}$, который строится ниже двойственным образом,
впервые рассматривался в работах А.~Агоре и А.~Чирваситу~\cite{AgoreCatConstr, Chirvasitu}.

 На произведении биалгебр $\prod_{n=0}^\infty  B_n$ определим гомоморфизм 
$$S \colon \prod_{n=0}^\infty  B_n \to \left(\prod_{n=0}^\infty  B_n\right)^{\mathrm{op,cop}}$$
как единственный гомоморфизм биалгебр, делающий коммутативной
диаграмму
 $$\xymatrix{  \prod_{n=0}^\infty  B_n \ar[d]^{\pi_{n+1}}\ar@{-->}[r]^(0.4)S & \left(\prod_{n=0}^\infty  B_n\right)^{\mathrm{op,cop}\ar[d]^{\pi_n}} \\
  B_{n+1} \ar[r]^{\sim} &  B_n^{\mathrm{op,cop}} }$$

 Обозначим через $H_r(B)$
сумму всех таких подкоалгебр биалгебры $\prod_{n=0}^\infty  B_n$, что
для всех их элементов $h$ справедливы равенства $S(h_{(1)})h_{(2)}=\varepsilon(h)1$
и $h_{(1)}S(h_{(2)})=\varepsilon(h)1$.
Тогда $H_r(B)$~"--- подалгебра Хопфа. Обозначим через ${\nu\,}_B \colon H_r(B) \to B$ отображение, индуцированное отображением
 $\pi_0 \colon \coprod_{n=0}^\infty  B_n \to B_0=B$.
 Тогда $H_r \colon \mathbf{Bialg}_\mathbbm{k} \to \mathbf{Hopf}_\mathbbm{k}$ является функтором, а $\nu$ является естественным преобразованием 
 композиции функтора $H_r$ и функтора вложения $ \mathbf{Hopf}_\mathbbm{k} \to \mathbf{Bialg}_\mathbbm{k}$ в тождественный функтор категории $\mathbf{Bialg}_\mathbbm{k}$.
 Естественное преобразование $\nu$ удовлетворяет следующему универсальному свойству:
 для любой алгебры Хопфа $H$ и гомоморфизма биалгебр $\varphi \colon H \to B$
 существует единственный гомоморфизм алгебр Хопфа $\varphi_0 \colon H \to H_r(B)$,
 делающий коммутативной диаграмму
 
 $$\xymatrix{ H \ar@{-->}[r]^(0.4){\varphi_0}\ar[rd]_{\varphi} & H_r(B) \ar[d]^{{\nu\,}_B} \\
 & B }$$
 
 Иными словами, функтор $H_r$~"--- правый сопряжённый к функтору вложения $\mathbf{Hopf}_\mathbbm{k} \to \mathbf{Bialg}_\mathbbm{k}$.
 
\section{Критерий нильпотентности линейного оператора}

В данном параграфе формулируется и доказывается обобщение
 известного критерия нильпотентности линейного оператора в терминах следа на случай поля произвольной характеристики. Хотя это утверждение и не является новым, но его доказательство в требуемой нам форме сложно найти в общедоступной литературе.

\begin{theorem}\label{TheoremTraceCriterionForNilpotence}
Пусть $\mathcal A \colon V \to V$~"--- линейный оператор на конечномерном векторном пространстве $V$ над полем $\mathbbm{k}$, причём либо $\chr \mathbbm{k} = 0$, либо $\chr \mathbbm{k} > \dim V$.
Тогда оператор $\mathcal A$ нильпотентен, если и только если $\tr\left(\mathcal A^k\right) = 0$
для всех $1 \leqslant k \leqslant \dim V$.
\end{theorem}
\begin{proof}
Расширим основное поле $\mathbbm{k}$ до его алгебраического замыкания $\overline{\mathbbm{k}}$ и рассмотрим
оператор $\bar{\mathcal A} \colon \overline{V} \to \overline{V}$, где $\overline{V} := \overline{\mathbbm{k}} \mathbin{\otimes_\mathbbm{k}} V$, заданный при помощи равенства $\bar{\mathcal A} (\alpha \otimes v) := \alpha \otimes {\mathcal A} v$ для всех $\alpha \in \overline{\mathbbm{k}}$ и $v\in V$.

Для любого базиса $v_1, \ldots, v_n$ пространства $V$ элементы $1 \otimes v_1, \ldots, 1\otimes v_n$
образуют базис пространства $\overline V$. Более того, оператор $\bar{\mathcal A}$ в базисе 
$1 \otimes v_1, \ldots, 1\otimes v_n$ имеет ту же матрицу, что и оператор $\mathcal A$ в базисе $v_1, \ldots, v_n$. Следовательно, $\tr\left(\mathcal A^k\right)=\tr\left(\bar{\mathcal A}^k\right)$. Более того, оператор $\mathcal A$ нильпотентен, если и только если нильпотентен оператор $\bar{\mathcal A}$.
Отсюда без ограничения общности можно считать, что основное поле $\mathbbm{k}$ алгебраически замкнуто.

Приведём оператор $\mathcal A$ к жордановой нормальной форме, сгруппировав клетки, отвечающие одинаковым собственным значениям. Пусть в матрице $A$ оператора $\mathcal A$ в соответствующем базисе на главной диагонали стоит $n_1$ элементов $\lambda_1$, $n_2$ элементов $\lambda_2$, \ldots, $n_s$ элементов $\lambda_k$, причём элементы $\lambda_1, \ldots, \lambda_s \in \mathbbm{k}$ попарно различны. Поскольку матрица $A$ верхнетреугольная, при возведении $\mathcal A$ в степень в эту же степень возводятся и элементы $\lambda_i$. Поэтому, если оператор $\mathcal A$ нильпотентен, все элементы $\lambda_i$ равны нулю, откуда
все  $\tr\left(\mathcal A^k\right) = \sum_{i=1}^s n_i \lambda_i^k$ также равны $0$.

Обратно, предположим, что $\tr\left(\mathcal A^k\right) = \sum_{i=1}^s n_i \lambda_i^k = 0$ для всех 
$1 \leqslant k \leqslant \dim V=: n$.

Тогда $$ \begin{pmatrix} 1 & 1 & \ldots  & 1\\
\lambda_1 & \lambda_2 & \ldots  & \lambda_s\\
\lambda_1^2 & \lambda_2^2 & \ldots  & \lambda_s^2\\
\hdotsfor{4} \\
\lambda_1^{n-1} & \lambda_2^{n-1} & \ldots  & \lambda_s^{n-1}
 \end{pmatrix}
\begin{pmatrix}n_1\lambda_1 \\ n_2\lambda_2 \\ \vdots \\ n_s\lambda_s\end{pmatrix}
 = \begin{pmatrix}0 \\ 0 \\ \vdots \\ 0\end{pmatrix}.
 $$
 Отсюда, учитывая невырожденность матрицы Вандермонда, получаем, что
 $$\begin{pmatrix}n_1\lambda_1 \\ n_2\lambda_2 \\ \vdots \\ n_s\lambda_s\end{pmatrix}
 = \begin{pmatrix}0 \\ 0 \\ \vdots \\ 0\end{pmatrix}.
 $$

Если $\chr \mathbbm{k} = 0$, то $\lambda_1 = \ldots = \lambda_s  = 0$. 
 
 В случае $\chr \mathbbm{k} > 0$ справедливы неравенства $n_i \leqslant \dim A < \chr \mathbbm{k}$, откуда также следует, что
 $\lambda_1 = \ldots = \lambda_s  = 0$.
 
 Таким образом, оператор $\mathcal A$ имеет верхнетреугольную матрицу с нулями на главной диагонали и, следовательно, нильпотентен.
\end{proof}
\begin{remark}\label{RemarkTraceCriterionForNilpotence} В теореме~\ref{TheoremTraceCriterionForNilpotence} недостаточно требовать только того, чтобы число $\dim V$ не делилось на
характеристику поля. Действительно, рассмотрим для простого числа $p = \chr \mathbbm{k}$ оператор $\mathcal A $,
который в некотором базисе имеет матрицу $$\begin{pmatrix}
1 &  0 &  \ldots    &  0  &  0  \\
0 &  1 &  \ldots    &  0  &  0  \\
\hdotsfor{5} \\
 0 & 0 &  \ldots &  1 & 0   \\      
 0 & 0 & \ldots &   0 &  0 
\end{pmatrix},$$ где число единиц на главной диагонали равно $p$. Тогда $\tr\left(\mathcal A^k\right) = p = 0$
для всех $k > 0$, однако оператор $\mathcal A$ не является нильпотентным, несмотря на то,
что $\dim V = p+1$ не делится на $p$.
\end{remark}

\section{Представления симметрической группы}
\label{Represent}

Приведем необходимые определения и результаты из теории
представлений симметрической группы~\cite{Bahturin, Fulton, FulHar, ZaiGia}.
В данном параграфе через~$\mathbbm{k}$ обозначается произвольное
 поле характеристики~$0$.

{\itshape Разбиением} числа~$n \in \mathbb N$
называется набор~$\lambda=(\lambda_1,
\lambda_2, \ldots, \lambda_s)$ из целых
 чисел~$\lambda_1 \geqslant \lambda_2 \geqslant \ldots
 \geqslant \lambda_s > 0$, $s\in\mathbb N$, такой,
  что~$\sum\limits_{k=1}^s \lambda_k = n$.
  Будем обозначать это
  через~$\lambda \vdash n$.\label{DefPartition}
Удобно считать, что $\lambda_j=0$ при $j > s$ или $j < 0$.

По каждому разбиению~$\lambda$ можно построить
{\itshape диаграмму Юнга} $$D_\lambda =
\lbrace (i,j) \in \mathbb Z \times \mathbb Z \mid
1 \leqslant i \leqslant s, 1 \leqslant j \leqslant \lambda_i \rbrace.
\label{DefDlambda}$$

Диаграмму Юнга можно изобразить в виде набора клеток
(числа в скобках указывают количество клеток
в соответствующей строчке):
$$ D_\lambda =
\begin{array}{c|c|c|c|c|c|c|c|}
\cline{2-8}
(\lambda_1) &  \multicolumn{3}{c|}{\ldots\ldots} & ~  & ~ & ~ & ~ \\
\cline{2-8}
(\lambda_2) &  \multicolumn{3}{c|}{\ldots\ldots} & ~ \\
\cline{2-5}
\ldots  & \multicolumn{3}{c|}{\ldots\ldots} \\
\cline{2-4}
(\lambda_{s-1}) & ~ & ~ \\
\cline{2-3}
(\lambda_s) & ~ \\
\cline{2-2}
\end{array}.$$

Будем придерживаться соглашения, при котором паре~$(i,j)$
соответствует клетка, стоящая в~$i$-й строке и~$j$-м столбце.

 Диаграмма Юнга, заполненная в некотором порядке
 числами~$1$, $2$, \ldots, $n$, называется
 {\itshape таблицей Юнга}:\label{DefTlambda}
 $$T_\lambda =
\begin{array}{|c|c|c|c|c|}
\hline
2 &  1 & 3  & 9 & 7 \\
\hline
5 &  6 \\
\cline{1-2}
8 \\
\cline{1-1}
4 \\
\cline{1-1}
\end{array}.$$

В данном примере~$\lambda = (5,2,1,1)$.

 В обозначении~$T_\lambda$
  подразумевается, что таблица~$T_\lambda$
  получена заполнением диаграммы~$D_\lambda$,
 отвечающей разбиению~$\lambda \vdash n$.
Будем говорить также, что таблица~$T_\lambda$
 имеет форму~$\lambda$.

Пусть $S_n$~"--- симметрическая группа, т.е.
группа биекций множества~$\lbrace 1,2,\ldots, n \rbrace$.
Обозначим через~$C_{T_\lambda}$ подгруппу в~$S_n$,
переводящую в себя множество чисел каждого столбца,
а через~$R_{T_\lambda}$~"--- подгруппу, переводящую
в себя множество чисел каждой строки
 диаграммы~$T_\lambda$. Пусть~$\mathbbm{k}S_n$~"--- групповая алгебра
 группы~$S_n$ над полем~$\mathbbm{k}$.
Для каждой таблицы Юнга~$T_\lambda$ определим {\itshape
симметризатор Юнга} $$e_{T_\lambda} =
 \sum\limits_{\substack{\rho \in R_{T_\lambda}, \\
 \sigma \in C_{T_\lambda} }} (\sign \sigma) \rho \sigma.
 \label{DefeTlambda}$$

 \begin{theorem}[{\cite[лемма~4.26, с.~54]{FulHar}},
  {\cite[следствие из леммы~3.2.5,
   с.~110]{Bahturin}}]\label{SymmYoungMul}
 Пусть~$T_\lambda$ и~$T_\mu$~"--- таблицы Юнга,
 отвечающие некоторым разбиениям~$\lambda, \mu \vdash n$
 соответственно. Тогда~$e_{T_\lambda} (\mathbbm{k}S_n) e_{T_\mu} = 0$
  при~$\lambda \ne \mu$
 и~$e_{T_\lambda}^2 = \beta e_{T_\lambda}$ для
  некоторого~$\beta \in \mathbb N$.
 \end{theorem}

  В силу последнего свойства
 симметризаторы Юнга называются квазиидемпотентами.

 \begin{theorem}[{\cite[теорема~4.3, с.~46]{FulHar}},
  {\cite[теорема~3.2.6,
   с.~111]{Bahturin}}]\label{IrrSymm}
 Для всякого разбиения~$\lambda \vdash n$, $n \in \mathbb N$,
  и таблицы Юнга~$T_\lambda$
 идеал~$\mathbbm{k}S_n e_{T_\lambda}$ является минимальным левым идеалом
 алгебры~$\mathbbm{k}S_n$. Всякий неприводимый~$\mathbbm{k}S_n$-модуль изоморфен
 одному из модулей~$\mathbbm{k}S_n e_{T_\lambda}$ для
  некоторого~$\lambda\vdash n$. При
   этом~$\mathbbm{k}S_n$-модули~$\mathbbm{k}S_n e_{T'_\lambda}$
 и~$\mathbbm{k}S_n e_{T''_\mu}$, где $\lambda,\mu \vdash n$, изоморфны
 тогда и только тогда, когда~$\lambda=\mu$ (таблицы~$T'_\lambda$
 и~$T''_\mu$ могут быть произвольными).
 \end{theorem}

 Для произвольного разбиения~$\lambda=(\lambda_1,\lambda_2,
 \ldots, \lambda_s) \vdash n$
 определим~{\itshape транспонированное разбиение} $\lambda^T=
 (\lambda^T_1,\lambda^T_2,
 \ldots, \lambda^T_r) \vdash n$.\label{DeflambdaT}
  Здесь~$\lambda^T_j = \max \lbrace k \mid \lambda_k \geqslant j \rbrace$,
  $r = \max \lbrace j \mid \lambda^T_j > 0 \rbrace = \lambda_1$.
  Иными словами,~$\lambda^T_j$~"--- число клеток
  в~$j$-м столбце диаграммы~$D_\lambda$, а
  диаграмма~$D_{\lambda^T}$ получается из диаграммы~$D_\lambda$
  отражением относительно главной диагонали.
  Обозначим через~$h_{ij}$ длину крюка диаграммы~$D_\lambda$
   с углом в клетке~$(i,j)$:
$$\begin{array}{|c|c|c|c|c|c|c|c|}
\hline
\qquad & \qquad &     ~ &       ~ & \qquad &      ~      & \qquad & \qquad \\
\hline
\qquad & \qquad & (i,j) & (i,j+1) & \ldots & (i,\lambda_i) \\
\cline{1-6}
\qquad & \qquad & (i+1,j) & ~   \\
\cline{1-4}
\qquad & \qquad & \ldots \\
\cline{1-3}
\qquad & \qquad & \phantom{\Bigl|} (\lambda^T_j,j)\\
\cline{1-3}
\qquad \\
\cline{1-1}
\end{array}.$$

Понятно, что~$h_{ij} = (\lambda_i - i)+(\lambda^T_j - j)+1$.

  \begin{theorem}[формула крюков~{\cite[предложение~4.12, с.~50]{FulHar}},
   {\cite[\S 4.3, с.~68]{Fulton}}]\label{HookFormula}
 Пусть~$M_\lambda$~"---  $\mathbbm{k}S_n$-модуль, изоморфный минимальному
 левому идеалу~$\mathbbm{k}S_n e_{T_\lambda}$.
Тогда
$$\dim M_\lambda = \frac{n!}{\strut\prod
\limits_{(i,j)\in D_\lambda} h_{ij}}.$$

 \end{theorem}

  \begin{lemma}[{\cite[лемма~4.21, с.~53]{FulHar}},
   {\cite[лемма~3.2.5, с.~110]{Bahturin}}]
 \label{vonNeumann} Пусть~$T_\lambda$~"---
 таблица Юнга, отвечающая
 разбиению~$\lambda \vdash n$,
   и для некоторого элемента~$a \in \mathbbm{k}S_n$
   при всех~$\rho \in R_{T_\lambda}$ и~$\sigma \in C_{T_\lambda}$
   выполняется равенство~$\rho a \sigma = (\sign \sigma)a$.
Тогда~$a = \alpha e_{T_\lambda}$ для некоторого~$\alpha \in \mathbbm{k}$.
Обратно, для всякого элемента~$a = \alpha e_{T_\lambda}$
 справедливо указанное равенство.
 \end{lemma}

 Следующая теорема позволяет вычислять кратности вхождения
 неприводимых слагаемых.

\begin{theorem}\label{ThKratnost}
Пусть
$T_{\lambda}$~"--- таблица Юнга,
 отвечающая фиксированному разбиению~$\lambda \vdash n$,
   $M=\bigoplus\limits_{i=1}^m M_i  \oplus
   \bigoplus\limits_{j=1}^k L_j $~"---
   разложение $\mathbbm{k}S_n$-модуля~$M$ в прямую сумму неприводимых
   подмодулей,
   причем~$M_i \cong M_{\lambda}$ при всех~$1 \leqslant i \leqslant m$,
   а $L_j \ncong M_{\lambda}$ при всех~$1 \leqslant j \leqslant k$.
 Тогда~$m = \dim \left( e_{T_\lambda} M \right)$.
В частности, если~$e_{T_\lambda} M = 0$, то модули~$M_{\lambda}$
не входят в разложение модуля~$M$.
\end{theorem}
\begin{proof}
 Для всякого~$j$ по теореме~\ref{IrrSymm} существует
 такое~$\mu\vdash n$, $\mu \ne \lambda$, что~$L_j \cong \mathbbm{k}S_n e_{T_\mu}$.
 Поэтому~$e_{T_\lambda} L_j = e_{T_\lambda} (\mathbbm{k}S_n) e_{T_\mu} = 0$
в силу теоремы~\ref{SymmYoungMul}. Отсюда для завершения доказательства
достаточно показать, что~$\dim e_{T_\lambda} M_\lambda = 1$.
Это следует из леммы~\ref{vonNeumann}
 и того, что,
в силу теоремы~\ref{SymmYoungMul},
 $e_{T_\lambda} M_\lambda =
 e_{T_\lambda} (\mathbbm{k}S_n) e_{T_\lambda}
 \ni e_{T_\lambda}^2 \ne 0$.
\end{proof}

\newpage

\chapter{Эквивалентность градуировок и действий групп}
\label{ChapterGradEquiv}

Результаты этой главы были опубликованы в работах~\cite{ASGordienko18Schnabel, ASGordienko19Schnabel, ASGordienko20ALAgoreJVercruysse}.

         \section{Реализация элементарных градуировок на матричных алгебрах конечной группой}
\label{SectionRegradingMatrixFinite}

         \subsection{Различные формулировки задачи}
         
         Как будет показано ниже в следствии~\ref{CorollaryFiniteRegradingImpossible349}, при больших $n$ на алгебре $M_n(\mathbbm{k})$ существует элементарная градуировка (см. \S\ref{SectionGradEquivUnivGroup}) бесконечной группой, не эквивалентная никакой градуировке конечной группой. В связи с этим представляет интерес следующая задача:
         
         \begin{problem}\label{ProblemMinNElAGradNonFin} 
Определить множество $\Omega$ таких натуральных чисел $n\in\mathbb N$,
что для любая элементарная градуировка на $M_n(\mathbbm{k})$ реализуется градуировкой конечной группой.
\end{problem}

В подпараграфе~\ref{SubsectionProofElementaryNonRegrad} мы докажем следующий результат:

\begin{theorem}\label{TheoremRegradeElementaryOmega}
Множество $\Omega$, введённое в задаче~\ref{ProblemMinNElAGradNonFin},
имеет следующую структуру: $$\Omega=\lbrace n \in \mathbb N \mid 1 \leqslant n \leqslant n_0\rbrace$$
для некоторого $3\leqslant n_0 \leqslant 348$.
В частности, для любого $n \geqslant 349$ 
существует элементарная градуировка на $M_n(\mathbbm{k})$, которую нельзя реализовать градуировкой конечной группой.
\end{theorem}

В этом подпараграфе мы покажем, что задача~\ref{ProblemMinNElAGradNonFin} может быть переформулирована в чисто теоретико-групповых терминах.

Для произвольного подмножества $U$ группы $G$ введём обозначение $\diff U := \lbrace uv^{-1} \mid u,v\in U,\ u\ne v\rbrace$. \label{NotationDiff}

Напомним также, что если градуировка $\Gamma_2$ является огрублением градуировки $\Gamma_2$,
то через $\pi_{\Gamma_1 \to \Gamma_2} \colon G_{\Gamma_1} \twoheadrightarrow G_{\Gamma_2}$
мы обозначаем гомоморфизм групп, заданный условием $\pi_{\Gamma_1 \to \Gamma_2}(\varkappa_{\Gamma_1}(g))=\varkappa_{\Gamma_2}(h)$ для $g\in\supp\Gamma_1$ и $h\in\supp \Gamma_2$, таких, что $A^{(g)} \subseteq A^{(h)}$.

\begin{lemma}\label{LemmaCoarseningSubset} 
Пусть $\Gamma_2 \colon A=\bigoplus_{h \in H} A^{(h)}$~"--- огрубление градуировки $\Gamma_1 \colon A=\bigoplus_{g \in G} A^{(g)}$.
Обозначим через $Q \triangleleft  G_{\Gamma_1}$ нормальную подгруппу, порождённую множеством \begin{equation*}\begin{split}W := \lbrace \varkappa_{\Gamma_1}(g_1) \varkappa_{\Gamma_1}(g_2)^{-1} \mid
g_1,g_2 \in \supp \Gamma_1,\text{ такие, что } \\ A^{(g_1)}\oplus A^{(g_2)} \subseteq A^{(h)}\text{ для некоторых }h\in \supp \Gamma_2\rbrace.\end{split}\end{equation*}
 Тогда $\ker \pi_{\Gamma_1 \to \Gamma_2} = Q$. Кроме того, $Q \cap \diff\varkappa_{\Gamma_1}(\supp \Gamma_1)= W$.
\end{lemma}
\begin{proof}
Очевидно, что $Q\subseteq \ker \pi_{\Gamma_1 \to \Gamma_2}$.

Пусть $\pi_{\Gamma_1} \colon \mathcal F([\supp \Gamma_1]) \twoheadrightarrow G_{\Gamma_1}$
и $\pi_{\Gamma_2} \colon \mathcal F([\supp \Gamma_2]) \twoheadrightarrow G_{\Gamma_2}$~"---
естественные сюръективные гомоморфизмы. Обозначим через $\varphi \colon \mathcal F([\supp \Gamma_1]) \twoheadrightarrow \mathcal F([\supp \Gamma_2])$ сюръективный гомоморфизм, заданный условием $\varphi([g]) = [h]$ для $g\in\supp\Gamma_1$ и $h\in\supp \Gamma_2$, таких, что $A^{(g)} \subseteq A^{(h)}$.
Тогда следующая диаграмма коммутативна:
\begin{equation}\label{EqCommDiagGamma1Gamma2}\xymatrix{ \mathcal F([\supp \Gamma_1]) \ar[r]^(0.65){\pi_{\Gamma_1}} \ar[d]^{\varphi} & G_{\Gamma_1} \ar[d]^{\pi_{\Gamma_1 \to \Gamma_2}}\\ \mathcal F([\supp \Gamma_2])
\ar[r]^(0.65){\pi_{\Gamma_2}}  & G_{\Gamma_2}}\end{equation}

Заметим, что $\ker \varphi$ совпадает с нормальной подгруппой группы $\mathcal F([\supp \Gamma_1])$,
порождённой
 всевозможными элементами $[g_1][g_2]^{-1}$, где $g_1,g_2 \in \supp \Gamma_1$ и $A^{(g_1)}\oplus A^{(g_2)} \subseteq A^{(h)}$ для некоторого $h\in \supp \Gamma_2$. Следовательно, $\pi_{\Gamma_1}(\ker\varphi) = Q$.

Предположим, что $\pi_{\Gamma_1 \to \Gamma_2} \pi_{\Gamma_1}(w) = 1_{G_{\Gamma_2}}$ для
некоторого $w\in \mathcal F([\supp \Gamma_1])$. Тогда из~(\ref{EqCommDiagGamma1Gamma2}) следует, что $\varphi(w)\in \ker \pi_{\Gamma_2}$.
Следовательно, элемент $\varphi(w)$ 
принадлежит нормальной подгруппе, порождённой словами
$[h_1][h_2] [h_1 h_2]^{-1}$ для всех $h_1, h_2 \in \supp \Gamma_2$, таких, что $A^{(h_1)}A^{(h_2)}\ne 0$.
Однако последнее включение справедливо, если и только если $A^{(g_1)}A^{(g_2)}\ne 0$
для некоторых $g_1, g_2 \in \supp \Gamma_1$, таких, что $A^{(g_1)}\subseteq A^{(h_1)}$, $ A^{(g_2)}\subseteq A^{(h_2)}$.
Следовательно, $w=w_0 w_1$, где $w_0 \in \ker\varphi$ и $w_1 \in \ker \pi_{\Gamma_1}$.
В частности, $\pi_{\Gamma_1}(w) = \pi_{\Gamma_1}(w_0)\in Q$. В силу того, что $\pi_{\Gamma_1}$
сюръективно,  $\ker \pi_{\Gamma_1 \to \Gamma_2} = Q$.
Равенство
  $Q \cap \diff\varkappa_{\Gamma_1}(\supp \Gamma_1)= W$ теперь следует из очевидного равенства
 $\ker \pi_{\Gamma_1 \to \Gamma_2} \cap \diff\varkappa_{\Gamma_1}(\supp \Gamma_1)= W$.
\end{proof}

\begin{lemma}\label{LemmaSubsetCoarseningExists} Пусть $\Gamma_1 \colon A=\bigoplus_{g \in G} A^{(g)}$ ~"--- градуировка группой $G$.
Тогда для всякого подмножества $W \subseteq \diff\varkappa_{\Gamma_1}(\supp \Gamma_1)$
существует огрубление $\Gamma_2 \colon A=\bigoplus_{h \in H} A^{(h)}$  градуировки $\Gamma_1$, такое,
что $\ker\pi_{\Gamma_1 \to \Gamma_2} = Q$,
где $Q$~"--- нормальная подгруппа группы $G_{\Gamma_1}$, порождённая множеством $W$.
\end{lemma}
\begin{proof}
Пусть $\pi_{\Gamma_1} \colon \mathcal F([\supp \Gamma_1]) \twoheadrightarrow G_{\Gamma_1}$ и $\pi \colon G_{\Gamma_1} \twoheadrightarrow G_{\Gamma_1}/Q$~"--- естественные сюръективные
гомоморфизмы. Рассмотрим градуировку $\Gamma_2 \colon A=\bigoplus_{u \in G_{\Gamma_1}/Q} A^{(u)}$,
где $A^{(u)} := \bigoplus_{\substack{g\in \supp \Gamma_1,\\ \pi(\varkappa_{\Gamma_1}(g))=u}} A^{(g)}$. 
Докажем, что $G_{\Gamma_1}/Q$~"--- универсальная группа градуировки $\Gamma_2$.

Предположим, что $\Gamma_2$
может быть реализована как градуировка группой $H$ 
и $\psi \colon \supp\Gamma_2\hookrightarrow H$~"--- соответствующее вложение носителя.
Тогда
существует единственный гомоморфизм $\varphi \colon \mathcal F([\supp \Gamma_1]) \to H$,
такой, что $\varphi([g])=\psi(\pi(\varkappa_{\Gamma_1}(g)))$
для всех $g\in \supp \Gamma_1$.
Заметим, что $\ker (\pi\pi_{\Gamma_1})$~"--- нормальная подгруппа
группы $\mathcal F([\supp \Gamma_1])$, порождённая следующими словами:
\begin{enumerate}
\item $[g][h][gh]^{-1}$ для всех  $g,h \in \supp \Gamma_1$,
таких, что $A^{(g)}A^{(h)}\ne 0$;
\item $[g_1] [g_2]^{-1}$ для всех $\varkappa_{\Gamma_1}(g_1) \varkappa_{\Gamma_1}(g_2)^{-1} \in W$.
\end{enumerate}
Следовательно, $\ker (\pi\pi_{\Gamma_1}) \subseteq \ker \varphi$ и существует
гомоморфизм $\bar\varphi \colon G_{\Gamma_1}/Q \to H$, такой, что
$\bar \varphi\pi\pi_{\Gamma_1} = \varphi$. В частности, $\bar\varphi(u)=\psi(u)$ для всех $u\in \supp \Gamma_2$.
В силу того, что группа $G_{\Gamma_1}/Q$ порождена множеством $\supp \Gamma_2$, гомоморфизм $G_{\Gamma_1}/Q \to H$, обладающий этим свойством,
единственный. Поэтому $G_{\Gamma_1}/Q$ является универсальной группой градуировки $\Gamma_2$,
а $\pi_{\Gamma_1 \to \Gamma_2}$ может быть отождествлён с $\pi$.
\end{proof}
\begin{remark}
Заметим, что включение $W \subseteq Q \cap \diff\varkappa_{\Gamma_1}(\supp \Gamma_1)$
в лемме~\ref{LemmaSubsetCoarseningExists}
может быть строгим.
\end{remark}

Напомним, что группа $G$ называется \textit{финитно аппроксимируемой} или \textit{остаточно конечной},
если для любого $g\in G$, где $g \ne 1_G$, существует такая нормальная подгруппа $N \triangleleft G$ конечного индекса, что образ $g$ в $G/N$ отличен от единицы группы $G/N$, т.е. $g\notin N$. В силу того, что пересечение подгрупп конечного индекса является подгруппой конечного индекса, для любого конечного подмножества $W$ элементов группы $G$, такого, что $1_G \notin W$, существует такая нормальная подгруппа $N \triangleleft G$ конечного индекса, что $W \cap N = \varnothing$.

Нам потребуются более слабые ограничения на группу $G$, при которых множество $W$ не произвольно, а фиксировано. Для удобства формулировок разрешим, чтобы множество $W$ содержало единицу $1_G$.

\begin{definition}
Пусть $G$~"--- группа, а $W$~"--- некоторое её подмножество. 
Будем говорить, что группа $G$ \textit{финитно аппроксимируема по отношению к $W$},
если  существует такая нормальная подгруппа $N \triangleleft G$ конечного индекса,
что $W \cap N = W \cap \lbrace 1_G \rbrace$.
\end{definition}

Если потребовать, чтобы не только сама группа $G$ была финитно аппроксимируема по отношению к $W$, но и её гомоморфные образы $\overline G$ былы финитно аппроксимируемы по отношению к $\overline W$,
где $\pi \colon G \twoheadrightarrow \overline G$~"--- соответствующий сюръективный гомоморфизм, $\overline W=\pi(W)$, а группа $\ker\pi$ порождена как нормальная группа некоторым подмножеством в $W$, то мы придём к понятию группы, наследственно финитно аппроксимируемой по отношению к $W$.

\begin{definition}
Будем говорить, что группа $G$ \textit{наследственно финитно аппроксимируема по отношению к $W$}, 
если для нормального замыкания $N_1 \triangleleft G$ любого (возможно, пустого) подмножества
множества $W$ существует такая нормальная подгруппа $N \triangleleft G$ конечного индекса,
что $N_1 \subseteq N$ и $W \cap N = W \cap N_1$.
\end{definition}

В теореме~\ref{TheoremFinGradingsToGroupsTransition}, которая доказывается ниже,
 задача о том, может ли градуировка и её всевозможные огрубления быть
реализованы конечными группами, переформулируется на чисто теоретико-групповом языке.

\begin{theorem}\label{TheoremFinGradingsToGroupsTransition} Пусть $G$~"--- группа, $A$~"--- алгебра, а
$\Gamma \colon A = \bigoplus_{g\in G} A^{(g)}$~"--- $G$-градуировка на $A$.
Тогда \begin{enumerate}
\item градуировка $\Gamma$ может быть реализована конечной группой, если и только если
$G_\Gamma$ финитно аппроксимируема по отношению к $\diff\varkappa_\Gamma(\supp \Gamma)$;
\item градуировка $\Gamma$ и все её огрубления могут быть реализованы конечными группами, если и только если
$G_\Gamma$ наследственно финитно аппроксимируема по отношению к $\diff\varkappa_\Gamma(\supp \Gamma)$.
\end{enumerate}
\end{theorem}
\begin{proof}
Градуировка $\Gamma$ 
может быть реализована любой такой факторгруппой группы $G_\Gamma$,
что элементы носителя градуировки имеют в этой факторгруппе различные образы.
Следовательно, если группа $G_\Gamma$ финитно аппроксимируема по отношению к $\diff\varkappa_\Gamma(\supp \Gamma)$,
то такая конечная факторгруппа существует. 
Обратно, если градуировка $\Gamma$
реализуется конечной группой $G$, 
то подгруппа группы $G$,
порождённая носителем градуировки является конечной факторгруппой группы $G_\Gamma$,
в которой элементы носителя градуировки имеют различные образы.
Следовательно, $G_\Gamma$ финитно аппроксимируема по отношению к $\diff\varkappa_\Gamma(\supp \Gamma)$,
и первая часть теоремы доказана.

Допустим, градуировка $G_\Gamma$ 
наследственно финитно аппроксимируема по отношению ко множеству 
 $\diff\varkappa_\Gamma(\supp \Gamma)$.
В силу леммы~\ref{LemmaCoarseningSubset} для любого огрубления $\Gamma_1$ градуировки $\Gamma$
существует нормальная подгруппа $Q$, которая является нормальным замыканием множества
\begin{equation*}\begin{split}Q \cap
\diff\varkappa_\Gamma(\supp \Gamma) = \lbrace \varkappa_{\Gamma}(g_1) \varkappa_{\Gamma}(g_2)^{-1} \mid
g_1,g_2 \in \supp \Gamma,\text{ такие, что } \\ A^{(g_1)}\oplus A^{(g_2)} \subseteq A^{(h)}\text{ для некоторых }h\in \supp \Gamma_1\rbrace,\end{split}\end{equation*}
такая, что $G_{\Gamma_1} \cong G_{\Gamma}/Q$. В силу того, что $G_\Gamma$
наследственно финитно аппроксимируема по отношению ко множеству
 $\diff\varkappa_\Gamma(\supp \Gamma)$, существует нормальная подгруппа $N \triangleleft G_{\Gamma}$ 
конечного индекса, такая, что $Q \subseteq N$ и
\begin{equation}\label{EqNQDiffKappa} N \cap \diff\varkappa_\Gamma(\supp \Gamma) = Q \cap \diff\varkappa_\Gamma(\supp \Gamma).\end{equation}
Введём обозначение $\bar N=\pi_{\Gamma \to \Gamma_1}(N)$. Предположим, что $\varkappa_{\Gamma_1}(h_1)\varkappa_{\Gamma_1}(h_2)^{-1} \in \bar N$ для некоторых $h_1, h_2 \in \supp \Gamma_1$. 
В силу изоморфизма $G_\Gamma/N \cong G_{\Gamma_1}/\bar N$  для всех $g_1, g_2 \in \supp \Gamma$, таких, что $A^{(g_1)} \subseteq A^{(h_1)}$ и $A^{(g_2)} \subseteq A^{(h_2)}$ 
выполнено $\varkappa_{\Gamma}(g_1)\varkappa_{\Gamma}(g_2)^{-1} \in N$.
Поэтому из~(\ref{EqNQDiffKappa})
следует, что $\varkappa_{\Gamma_1}(g_1)\varkappa_{\Gamma_1}(g_2)^{-1} \in Q$,
$\varkappa_{\Gamma_1}(h_1)=\varkappa_{\Gamma_1}(h_2)$ и $h_1 = h_2$. Следовательно, 
в факторгруппе группы $G_{\Gamma_1}$ по
нормальной подгруппе $\bar N$ различные элементы носителя градуировки $\Gamma_1$ имеют различные образы.
 В силу того, что $|G_{\Gamma_1}/\bar N| = |G_\Gamma/N|< +\infty$,
градуировка $\Gamma_1$ реализуется при помощи конечной группы $G_{\Gamma_1}/\bar N$.

Предположим, что всякое огрубление $\Gamma_1$ градуировки $\Gamma$ 
реализуется конечной группой.
Докажем, что $G_\Gamma$ 
наследственно финитно аппроксимируема по отношению ко множеству $\diff\varkappa_\Gamma(\supp \Gamma)$.
Действительно, пусть $N_1 \triangleleft G_\Gamma$~"---
нормальное замыкание некоторого подмножества множества $\diff\varkappa_\Gamma(\supp \Gamma)$. В силу леммы~\ref{LemmaSubsetCoarseningExists} существует такая градуировка $\Gamma_1$, что $G_{\Gamma_1} \cong G_{\Gamma}/N_1$.
В силу того, что градуировку $\Gamma_1$ 
можно реализовать конечной группой, существует нормальная подгруппа $Q \triangleleft G_{\Gamma_1}$
конечного индекса, такая, что $\varkappa_{\Gamma_1}(h_1) \varkappa_{\Gamma_1}(h_2)^{-1} \notin Q$
для всех $h_1, h_2 \in \supp \Gamma_1$, $h_1 \ne h_2$. 
Пусть $N := \pi_{\Gamma \to \Gamma_1}^{-1}(Q)$.
Тогда $|G_\Gamma/N| = |G_{\Gamma_1}/ Q| < +\infty$, $N \supseteq \ker\pi_{\Gamma \to \Gamma_1}=N_1$, и
$$N \cap \diff\varkappa_\Gamma(\supp \Gamma) \supseteq N_1 \cap \diff\varkappa_\Gamma(\supp \Gamma).$$
Предположим, что $\varkappa_{\Gamma}(g_1) \varkappa_{\Gamma}(g_2)^{-1} \in N$ для некоторых элементов $g_1, g_2 \in \supp \Gamma$.
Выберем $h_1, h_2 \in \supp \Gamma_1$, такие, что $A^{(g_1)} \subseteq A^{(h_1)}$ и $A^{(g_2)} \subseteq A^{(h_2)}$. Тогда $\varkappa_{\Gamma_1}(h_1) \varkappa_{\Gamma_1}(h_2)^{-1} \in Q$
и $h_1 = h_2$. Следовательно, $\varkappa_{\Gamma}(g_1) \varkappa_{\Gamma}(g_2)^{-1} \in N_1$
и $$N \cap \diff\varkappa_\Gamma(\supp \Gamma) = N_1 \cap \diff\varkappa_\Gamma(\supp \Gamma).$$
Отсюда $G_\Gamma$ наследственно финитно аппроксимируема по отношению ко множеству
$\diff\varkappa_\Gamma(\supp \Gamma)$.
\end{proof}

Предложение~\ref{prop:abelian}, которое приводится ниже в силу его важности, не является новым. Хотя это предложение можно было бы вывести из теоремы~\ref{TheoremFinGradingsToGroupsTransition}, докажем его непосредственно.

\begin{proposition}\label{prop:abelian}
Пусть $\Gamma \colon A=\bigoplus_{g \in G} A^{(g)}$
~"--- градуировка конечномерной алгебры $A$ абелевой группой $G$. 
Тогда $\Gamma$ эквивалентна градуировке конечной группой.
\end{proposition}
\begin{proof}
Группу $G$ можно заменить её подгруппой, порождённой носителем $\supp \Gamma$.
Поэтому в силу конечности $\supp \Gamma$ 
можно без ограничения общности считать $G$
конечно порождённой абелевой группой. Тогда $G$
является прямым произведением свободных и примарных циклических групп.
Заменяя свободные циклические группы конечными циклическими группами
достаточно большого порядка (см. пример~\ref{ExampleAbelianRegrading} ниже),
получаем градуировку конечной группой.
\end{proof}

\begin{example}\label{ExampleAbelianRegrading}
Пусть $n\in \mathbb N$, а $\Gamma$~"--- элементарная $\mathbb Z$-градуировка на $M_n(\mathbbm{k})$,
заданная набором $(1,2,\ldots,n)$, т.е. $e_{ij}\in M_n(\mathbbm{k})^{(i-j)}$.
Тогда $$\supp \Gamma = \lbrace -(n-1), -(n-2), \ldots, -1, 0, 1, 2, \ldots, n-1 \rbrace$$
и $\Gamma$ эквивалентна элементарной $\mathbb Z/(2n\mathbb Z)$-градуировке,
заданной набором $(\bar 1,\bar 2,\ldots,\bar n)$, т.е. $e_{ij}\in M_n(\mathbbm{k})^{(\overline{i-j})}$.
\end{example}

Рассмотрим такую градуировку $\Gamma_0$ на алгебре $M_n(\mathbbm{k})$ свободной
группой $\mathcal F(x_1, \ldots, x_{n-1})$,
что  $e_{r,r+1} \in M_n(\mathbbm{k})^{(x_r)}$ для всех $1\leqslant r \leqslant n-1$.
Иными словами, эта градуировка задана набором $(x_1 x_2 \ldots x_{n-1},\ x_2 \ldots x_{n-1},\ \ldots, \ x_{n-1},\ 1)$. Заметим, что однородная компонента градуировки $\Gamma_0$, отвечающая единице группы, 
является линейной оболочкой матричных единиц $e_{ii}$, $1\leqslant i \leqslant n$,
и имеет размерность $n$, а все остальные однородные компоненты одномерны.
В силу того, что для любой элементарной градуировки диагональные матричные единицы $e_{ii}$ принадлежат
однородной компоненте, отвечающей единице группы, и все матричные единицы являются однородными элементами,
любая элементарная градуировка на  $M_n(\mathbbm{k})$ является огрублением градуировки $\Gamma_0$. В силу того, что группа $\mathcal F(x_1, \ldots, x_{n-1})$ свободна и все её свободные 
порождающие принадлежат множеству $\supp \Gamma_0$,
существует изоморфизм $G_{\Gamma_0} \cong \mathcal F(x_1, \ldots, x_{n-1})$.
Заметим также, что \begin{equation*}\begin{split}\supp \Gamma_0 = \lbrace x_i x_{i+1}\ldots x_j  \mid 1\leqslant i \leqslant j \leqslant n-1 \rbrace\cup{ } \\  { }\cup
\lbrace x_j^{-1} x_{j-1}^{-1}\ldots x_i^{-1}  \mid 1\leqslant i \leqslant j \leqslant n-1 \rbrace \cup \lbrace 1\rbrace.\end{split}\end{equation*}
В силу теоремы~\ref{TheoremFinGradingsToGroupsTransition} задача~\ref{ProblemMinNElAGradNonFin}
эквивалентна задаче~\ref{ProblemMinFnHRFRx1x2x3} ниже:
\begin{problem}\label{ProblemMinFnHRFRx1x2x3} 
Определить множество $\Omega$ таких натуральных чисел $n\in\mathbb N$,
что группа $\mathcal F(x_1, \ldots, x_{n-1})$ наследственно финитно аппроксимируема по отношению
ко множеству $\diff W_n$,
где \begin{equation}\begin{split}\label{EqWnFreeGroup}W_n=\lbrace x_i x_{i+1}\ldots x_j  \mid 1\leqslant i \leqslant j \leqslant n-1 \rbrace\cup{ } \\  { }\cup
\lbrace x_j^{-1} x_{j-1}^{-1}\ldots x_i^{-1}  \mid 1\leqslant i \leqslant j \leqslant n-1 \rbrace \cup \lbrace 1\rbrace.\end{split}\end{equation}
\end{problem}

\subsection{Получение оценок}\label{SubsectionProofElementaryNonRegrad}

В этом подпараграфе доказывается теорема~\ref{TheoremRegradeElementaryOmega},
в которой даются оценки на множество $\Omega$ из задач~\ref{ProblemMinNElAGradNonFin} и~\ref{ProblemMinFnHRFRx1x2x3}.

Напомним, что универсальная группа градуировки $\Gamma$~"--- это пара $(G_\Gamma,\varkappa_\Gamma)$,
где $G_\Gamma$~"--- это группа, а $\varkappa_\Gamma\colon\supp \Gamma \hookrightarrow G_\Gamma$~"--- соответствующее вложение носителя.
Ниже мы доказываем теорему~\ref{TheoremGivenFinPresGroupExistence}, из которой 
следует, что первой компонентой этой пары может выступать любая конечно представленная группа.
Более того, для такой группы можно выбрать в качестве $\Gamma$ элементарную градуировку на матричной
алгебре. Возможность включить в носитель произвольное подмножество $V$ группы $G$ будет использована позже,
в доказательстве теоремы~\ref{TheoremFiniteRegradingImpossible}.

\begin{theorem}\label{TheoremGivenFinPresGroupExistence}
Пусть $\mathbbm{k}$~"--- поле, $G$~"--- конечно представленная группа, а $V \subseteq G$~"--- конечное (возможно, пустое) подмножество.
Тогда для некоторого $n\in\mathbb N$, которое зависит лишь от длины записи определяющих соотношений
группы $G$ и элементов множества $V$ через элементы, порождающие группу $G$,
существует элементарная градуировка $\Gamma$ на $M_n(\mathbbm{k})$, такая, что $G_\Gamma \cong G$ и 
 $V \subseteq \supp \Gamma$.
\end{theorem}
\begin{proof}
Пусть $G \cong \mathcal F(X)/N$, где $X=\lbrace x_1, x_2, \ldots, x_\ell\rbrace$~"--- конечное множество порождающих, а $N$~"--- нормальная подгруппа,
порождённая словами $w_1, \ldots, w_m$.
Пусть $V=\lbrace g_1, \ldots, g_s\rbrace \subseteq G$.
Выберем такие слова $w_{m+1},\ldots, w_{m+s}\in \mathcal F(X)$,
что $g_i=\bar w_{m+i}$, $1\leqslant i \leqslant s$, где через $\bar u$ обозначается образ слова $u\in \mathcal F(X)$ в $G$. Пусть $w_i = y_{i1}\ldots y_{ik_i}$, где $y_{ij}\in X \cup X^{-1}$, $1\leqslant i \leqslant m+s$.

Меняя, если нужно, нумерацию элементов $x_1, \ldots, x_\ell$, можно считать, что
первые $\ell_0$ порождающих
$x_1, \ldots, x_{\ell_0}$, где $0 \leqslant \ell_0 \leqslant \ell$, не встречаются
среди элементов $y_{ij}$ и $y_{ij}^{-1}$,
а для любого $\ell_0 < k \leqslant \ell$ существуют такие $i,j$, что либо $y_{ij}=x_k$, либо $y_{ij}=x_k^{-1}$.
 Пусть $n:=\ell_0+1+\sum_{i=1}^{m+s} k_i$.
 
 Обозначим через $$\Gamma \colon M_n(\mathbbm{k})=\bigoplus_{g\in G} M_n(\mathbbm{k})^{(g)}$$
элементарную $G$-градуировку, определённую следующим образом:
 $$e_{r,r+1} \in M_n(\mathbbm{k})^{(\bar y_{ij})} \text{ при } r=k_1+\ldots+k_{i-1}+j,
\ 1\leqslant j \leqslant k_i,\ 1\leqslant i \leqslant m+s$$ и $$e_{r,r+1} \in M_n(\mathbbm{k})^{(\bar x_j)}\text{ при }r=\sum_{i=1}^{m+s} k_i+j,\ 1\leqslant j \leqslant \ell_0$$
(соответствующая элементарная градуировка существует в силу замечания~\ref{RemarkElementary}).

Заметим, что произведение \begin{equation*}\begin{split}e_{\sum_{i=1}^{m+r-1} k_i+1, \sum_{i=1}^{m+r} k_i+1}
= \prod_{j=1}^{k_{m+r}} e_{\sum_{i=1}^{m+r-1} k_i + j,\sum_{i=1}^{m+r-1} k_i + j+1}  \in
\\
 \in M_n(\mathbbm{k})^{\left(\prod_{j=1}^{k_{m+r}}\bar y_{m+r,j}\right)}=M_n(\mathbbm{k})^{(g_r)}\end{split}
 \end{equation*}
 ненулевое. Следовательно $g_r \in \supp \Gamma$ для всех $1\leqslant r \leqslant s$,
 и $V\subseteq \supp \Gamma$. 

Докажем, что $G_\Gamma \cong G$.

Предположим, что градуировка $\Gamma$ реализована в качестве градуировки некоторой группой $H$.
Тогда существует вложение $\psi \colon \supp \Gamma \hookrightarrow H$,
заданное условием $M_n(\mathbbm{k})^{(g)}= M_n(\mathbbm{k})^{(\psi(g))}$.
При этом \begin{equation}\label{EqPsiPartialHom} \psi(g_1 g_2)=\psi(g_1)\psi(g_2) \end{equation}
 для всех $g_1, g_2 \in G$, таких, что $M_n(\mathbbm{k})^{(g_1)}M_n(\mathbbm{k})^{(g_2)}\ne 0$.
В силу того, что $M_n(\mathbbm{k})$~"--- алгебра с единицей, справедливо равенство $\psi(1_G)=1_H$.

Из того, как строилась градуировка $\Gamma$, следует, что для всех $x\in X$ выполнено
$\bar x, \bar x^{-1}\in\supp\Gamma$. Следовательно, 
в силу~\eqref{EqPsiPartialHom} $\psi(\bar x^{-1})=\psi(\bar x)^{-1}$
и элементы $\psi(\bar x)$ и $\psi(\bar x^{-1})$ определены для всех $x\in X$.
По индукции $$\psi(\bar y_{i1})\ldots \psi(\bar y_{ik_i})=\psi(\bar y_{i1}\ldots \bar y_{ik_i})=\psi(\bar w_i)=1_H$$
для всех $1\leqslant i \leqslant m$.
Следовательно, элементы $\psi(\bar x)$ удовлетворяют определяющим соотношениям группы $G$.
В силу этого существует гомоморфизм $\varphi \colon G \to H$, такой, что $\varphi(\bar x) = \psi(\bar x)$ для каждого $x\in X$. Из того, что $\lbrace \bar x \mid x \in X\rbrace$ порождает $G$, следует, что этот гомоморфизм единственный.

Осталось доказать, что $\varphi \bigr|_{\supp \Gamma} = \psi$.
Каждый элемент $g$ множества $\supp \Gamma$ 
соответствует матричной единице
$e_{ij}$, где либо $i > j$ и $e_{ij}=e_{i, i-1} e_{i-1,i-2}\ldots e_{j+1,j}$, либо $i < j$ и $e_{ij}=e_{i, i+1} e_{i+1,i+2}\ldots e_{j-1,j}$,
либо $i=j$ и $g=1_G$. Так как для всех $1\leqslant i,j \leqslant n$, таких, что
$|i-j|=1$, выполнено $e_{ij} \in M_n(\mathbbm{k})^{(\bar x)}$ для некоторого $x \in X \cup X^{-1}$ и $\varphi(x)=\psi(x)$ для $x \in X \cup X^{-1}$,
индукция по $|i-j|$ с использованием~(\ref{EqPsiPartialHom}) показывает, что $\varphi(g)=\psi(g)$.
 Следовательно, $G \cong G_\Gamma$.
\end{proof}

В лемме~\ref{LemmaSemigroupMatrixElementaryGroup}, которая доказывается ниже, используется идея из~\cite[\S 4]{ClaseJespersDelRio}.
Из леммы~\ref{LemmaSemigroupMatrixElementaryGroup}, в частности, следует, что если некоторая элементарная
градуировка на алгебре $M_n(\mathbbm{k})$ реализуется некоторой полугруппой, то эта полугрупповая градуировка
всё равно сводится к некоторой групповой градуировке.

\begin{lemma}\label{LemmaSemigroupMatrixElementaryGroup}
Пусть $\mathbbm{k}$~"--- поле, $n\in\mathbb N$.
Если $\Gamma \colon M_n(\mathbbm{k})=\bigoplus_{t\in T} M_n(\mathbbm{k})^{(t)}$~"--- градуировка
некоторой полугруппой $T$, причём все матричные единицы $e_{ij}$
являются элементами, однородными в этой градуировке, и существует элемент $e \in T$,
такой, что все $e_{ii}\in M_n(\mathbbm{k})^{(e)}$, тогда $e^2=e$
и $\supp \Gamma \subseteq \mathcal U(eTe)$, где
$\mathcal U(eTe)$~"--- группа обратимых элементов моноида $eTe$.
\end{lemma}
\begin{proof} Из равенства
$e_{11}^2=e_{11}$ следует, что $e^2=e$. 
Поскольку по условию единичная матрица принадлежит компоненте
$M_n(\mathbbm{k})^{(e)}$, справедливо включение $\supp \Gamma \subseteq eTe$.
При этом из равенства $e_{ij}e_{ji}=e_{ii}$ следует, что $\supp \Gamma \subseteq 
\mathcal U(eTe)$.
\end{proof}

Покажем теперь, что если для некоторого $m\in\mathbb N$
справедливо $m\notin \Omega$ (определение множества $\Omega$ см. в задаче~\ref{ProblemMinNElAGradNonFin}), то $n\notin \Omega$ для всех $n\geqslant m$.

\begin{lemma}\label{LemmaGradedSubalgebra}
Пусть $\Gamma \colon A = \bigoplus_{t\in T} A^{(t)}$~"---
градуировка (полу)группой $T$ на некоторой алгебре $A$ над полем $\mathbbm{k}$,
а $B\subseteq A$~"--- градуированная подалгебра.
Тогда если индуцированная на $B$ градуировка не может быть переградуирована
никакой	конечной группой, то и исходная градуировка $\Gamma$
не может быть переградуирована никакой конечной (полу)группой.
\end{lemma}
\begin{proof}
Всякая замена градуирующей группы градуировки $\Gamma$
приводит к соответствующей замене градуирующей группы на $B$.
Отсюда, если можно было бы реализовать $\Gamma$
конечной (полу)группой, это можно было бы сделать и для градуировки на $B$.
Однако последнее по условию невозможно.
\end{proof}

\begin{corollary}\label{CorollaryMatrixElemGreaterSize}
Если для некоторого $m\in\mathbb N$ и группы $G$
существует элементарная $G$-градуировка $\Gamma_0$ на алгебре $M_m(\mathbbm{k})$,
не эквивалентная никакой градуировке конечной (полу)группой,
то $G$-градуировки с таким свойством существуют
на алгебрах $M_n(\mathbbm{k})$ для всех $n\geqslant m$.
\end{corollary}
\begin{proof}
Предположим, что $\Gamma_0$ определяется
набором элементов
$(g_1, g_2, \ldots, g_m) \in G^m$.
Рассмотрим элементарную $G$-градуировку $\Gamma$ на алгебре $M_n(\mathbbm{k})$,
заданную набором $(g_1, g_2, \ldots, g_m,\underbrace{g_m, \ldots, g_m}_{n-m})$.
Тогда алгебра $M_m(\mathbbm{k})$ оказывается градуированной подалгеброй
алгебры $M_n(\mathbbm{k})$ (с другой единицей).
Следовательно, если $\Gamma$
была бы эквивалентна некоторой градуировке конечной (полу)группой,
то можно было бы переиндексировать элементами этой (полу)группы
и компоненты градуировки $\Gamma_0$, что по условию невозможно.
Отсюда градуировка $\Gamma$
неэквивалентна никакой градуировке конечной группой.
\end{proof}

Поскольку задачи~\ref{ProblemMinNElAGradNonFin} и~\ref{ProblemMinFnHRFRx1x2x3}
эквивалентны, отсюда немедленно получается следующее утверждение:
\begin{corollary}
Если для некоторого $m\in\mathbb N$ свободная группа $\mathcal F(x_1, \ldots, x_{m-1})$ 
не является наследственно финитно аппроксимируемой по отношению
ко множеству $\diff W_m$, то и для любого $n \geqslant m$ группа
$\mathcal F(x_1, \ldots, x_{n-1})$ не является наследственно финитно аппроксимируемой по отношению
ко множеству $\diff W_n$. (Определение множества $W_n$ см. в~(\ref{EqWnFreeGroup}).)
\end{corollary}

Примеры элементарных градуировок, которые не реализуются никакими конечными (полу)группами,
можно строить, используя финитно неаппроксимируемые группы.

\begin{theorem}\label{TheoremFiniteRegradingImpossible} Пусть $\mathbbm{k}$~"---
поле, а $G$~"--- конечно представленная группа,
не являющаяся финитно аппроксимируемой. (Например, $G$~"--- конечно представленная бесконечная простая
группа, см., например, \cite{Higman}.)
Тогда для некоторого $n\in\mathbb N$ на алгебре $M_n(\mathbbm{k})$ существует $G$-градуировка,
которая не эквивалента $H$-градуировке ни для какой конечной (полу)группы~$H$.
\end{theorem}
\begin{proof} Пусть $g_0\ne 1_G$~"--- произвольный элемент
из пересечения всех нормальных подгрупп группы $G$ конечного индекса. (Если группа $G$
проста, можно взять произвольный элемент $g_0\ne 1_G$.)

В силу теоремы~\ref{TheoremGivenFinPresGroupExistence}
для некоторого $n\in\mathbb N$, которое зависит только от длины записи определяющих соотношений
группы $G$
и выражения для $g_0$ через порождающие,
существует такая элементарная $G$-градуировка $\Gamma$ на алгебре $M_n(\mathbbm{k})$, что $G_\Gamma \cong G$ и $g_0 \in \supp \Gamma$.

Предположим, что $\Gamma$
эквивалентна градуировке конечной полугруппой $H$ и $\psi \colon \supp \Gamma \hookrightarrow H$~"--- соответствующее вложение носителя. Тогда существует такой элемент $e\in H$, что все $e_{ii} \in M_n(\mathbbm{k})^{(e)}$. Поскольку все матричные единицы являются в градуировке $\Gamma$ однородными элементами,
в силу леммы~\ref{LemmaSemigroupMatrixElementaryGroup} 
можно считать, что $H$~"--- группа. Из того, что $G_\Gamma \cong G$, следует существование единственного
такого гомоморфизма $\varphi \colon G \to H$, что $\varphi|_{\supp \Gamma} = \psi$.
Однако группа $H$ конечна, поэтому нормальная подгруппа $\ker \varphi$ имеет конечный индекс, откуда $\varphi(g_0)=1_H$. В силу того, что $g_0 \ne 1_G$, а $g_0, 1_G \in \supp \Gamma$, получаем противоречие.
Следовательно, никакая $H$-градуировка не может быть эквивалентна градуировке $\Gamma$.
\end{proof}

Из теоремы~\ref{TheoremFiniteRegradingImpossible} следует верхняя оценка из формулировки теоремы~\ref{TheoremRegradeElementaryOmega}:

\begin{corollary}\label{CorollaryFiniteRegradingImpossible349} Если $n \geqslant 349$, то на полной матричной алгебре $M_n(\mathbbm{k})$ существует градуировка бесконечной группой, не эквивалентная никакой градуировке конечной (полу)группой.
\end{corollary}
\begin{proof}
Рассмотрим конечно представленную бесконечную простую группу $G_{2,1}$ Томпсона.
Воспользуемся представлением группы $G_{2,1}$ из~\cite[\S~8]{Higman}. В качестве $g_0$ можно взять любой из порождающих группы $G_{2,1}$, поскольку все порождающие
оказываются в носителе градуировки, которая строится в теореме~\ref{TheoremGivenFinPresGroupExistence} для $G=G_{2,1}$.
Следовательно, достаточно применить теорему~\ref{TheoremGivenFinPresGroupExistence} для $V = \varnothing$.
Суммируя длины порождающих, получаем, что число $n$ из доказательства
теоремы~\ref{TheoremGivenFinPresGroupExistence}
равно $349$. Теперь достаточно применить следствие~\ref{CorollaryMatrixElemGreaterSize}.
\end{proof}


Теперь укажем класс градуировок, которые эквивалентны градуировкам конечными группами.
Рассмотрим элементарные градуировки, для которых различные недиагональные матричные единицы $e_{ij}$
принадлежат различным однородным компонентам, отвечающим элементам $g\ne 1$.

\begin{theorem}\label{TheoremFiniteRegradingAllDifferent} Пусть $G$~"--- группа, $\mathbbm{k}$~"--- поле,
$n\in\mathbb N$, а $(g_1, \ldots, g_n)$~"--- набор элементов группы $G$, такой, что $g_i g_j^{-1} = g_k g_\ell^{-1}$, если и только если
$\left\lbrace \begin{smallmatrix} i=k, \\ j=\ell \end{smallmatrix}\right.$ или
$\left\lbrace \begin{smallmatrix} i=j, \\ k=\ell. \end{smallmatrix}\right.$
Тогда элементарная градуировка на $M_n(\mathbbm{k})$,
заданная набором $(g_1, \ldots, g_n)$ 
эквивалентна  элементарной $S_{n+1}$-градуировке, заданной набором $(\gamma_1,\ldots,\gamma_n)$,
где $S_{n+1}$~"--- симметрическая группа, действующая на множестве $\lbrace 1,2,\ldots,n+1\rbrace$,
а $\gamma_i := (1, i+1) \in S_{n+1}$~"--- транспозиция, меняющая местами элементы $1$ и $i+1$.
Эта же градуировка на $M_n(\mathbbm{k})$ эквивалентна $\mathbb Z/(2^{n+1}\mathbb Z)$-градуировке, заданной набором $(\bar 2, \bar 2^2, \bar 2^3, \ldots, \bar 2^n)$.
\end{theorem}
\begin{proof} 
Для того, чтобы доказать первую часть теоремы, достаточно показать, что $\gamma_i \gamma_j^{-1} = \gamma_k \gamma_\ell^{-1}$, если только если $\left\lbrace \begin{smallmatrix} i=k, \\ j=\ell \end{smallmatrix}\right.$ или
$\left\lbrace \begin{smallmatrix} i=j, \\ k=\ell. \end{smallmatrix}\right.$ Однако  при $i=j$ элемент $\gamma_i \gamma_j^{-1}$ является тождественным преобразованием, а при $i\ne j$ элемент $\gamma_i \gamma_j^{-1}$ является $3$-циклом $(1, j+1, i+1)$.

Для того, чтобы доказать вторую часть теоремы, достаточно заметить, что
 $\bar 2^i-\bar 2^j=\bar 2^k-\bar 2^\ell$ в $\mathbb Z/(2^{n+1}\mathbb Z)$, если и только если $2^i-2^j=2^k-2^\ell$, что, в свою очередь, выполнено, если и только если или $\left\lbrace \begin{smallmatrix} i=k, \\ j=\ell \end{smallmatrix}\right.$, или 
$\left\lbrace \begin{smallmatrix} i=j, \\ k=\ell. \end{smallmatrix}\right.$
Действительно, разделив равенство $2^i-2^j=2^k-2^\ell$ на $2^{\min(i,j,k,\ell)}$, 
получим, что, по крайней мере, два числа из $i,j,k,\ell$ должны совпадать, что и доказывает требуемое.
\end{proof}

В теореме~\ref{TheoremFiniteRegradingImpossible}
была построена элементарная градуировка на $M_n(\mathbbm{k})$,
которая неэквивалентна никакой градуировке конечной группой.
Однако эта градуировка является огрублением элементарной градуировки, которая соответствует набору элементов $(1,z_1,\ldots, z_{n-1})$ свободной группы $\mathcal F(z_1, \ldots, z_{n-1})$
и которая в силу теоремы~\ref{TheoremFiniteRegradingAllDifferent}
эквивалентна градуировке конечной группой. Другими словами, существуют такие градуировки,
которые сами могут быть переградуированы конечными группами, а некоторые из их огрублений~"--- нет.

\begin{theorem}\label{TheoremFiniteRegradingSmallMatrix}
Пусть $\Gamma$~"--- элементарная $G$-градуировка
на полной матричной алгебре $M_n(\mathbbm{k})$,
где $n \leqslant 3$, $\mathbbm{k}$~"--- поле, а $G$~"--- группа.
Тогда $\Gamma$
эквивалентна градуировке конечной группой.
\end{theorem}
\begin{proof}
Если $n=1$, то $M_n(\mathbbm{k})=\mathbbm{k}$, и эта алгебра может быть переградуирована 
тривиальной группой.
Поэтому можно считать, что $n=2,3$.
Напомним, что $e_{ij} \in M_n(\mathbbm{k})^{(g_i g_j^{-1})}$,
где $(g_1, \ldots, g_n)$~"--- набор элементов группы $G$, задающий
градуировку $\Gamma$.
Обозначим через $A$ матрицу размера
$n\times n$, в клетке $(i,j)$ которой находится элемент $g_i g_j ^{-1}$.
Очевидно, что такая матрица полностью определяет всякую градуировку на алгебре $M_n(\mathbbm{k})$,
матричные единицы которой являются однородными элементами.
Более того, если $\Gamma_1$ и $\Gamma_2$~"--- элементарные градуировки 
на $M_n(\mathbbm{k})$, которым соответствуют матрицы $A_1=(g_{1ij})_{i,j}$ и $A_2=(g_{2ij})_{i,j}$,
то $\Gamma_1$ и $\Gamma_2$ эквивалентны при помощи изоморфизма $\id_{M_n(\mathbbm{k})}$,
если и только если для всех $1\leqslant i_1,j_1,i_2,j_2 \leqslant n$
выполнено условие $g_{1i_1j_1}=g_{1i_2j_2} \Leftrightarrow g_{2i_1j_1}=g_{2i_2j_2}$.

В случае $n=2$ справедливо равенство $A=\left(\begin{smallmatrix}
1 & g \\
g^{-1} & 1
\end{smallmatrix}\right)$, где $1=1_G$. 
Здесь есть две возможности: $g=g^{-1}$
и $g\ne g^{-1}$. В случае $g=g^{-1}$ алгебра $M_2(\mathbbm{k})$
может быть переградуирована группой $\mathbb Z/2\mathbb Z$ (соответствующая элементарная градуировка задаётся набором $(\bar 0,\bar 1)$). В случае $g\ne g^{-1}$ алгебра $M_2(\mathbbm{k})$
может быть переградуирована группой $\mathbb Z/3\mathbb Z$ (соответствующая элементарная градуировка задаётся набором $(\bar 0,\bar 1)$).

Рассмотрим теперь случай $n=3$.
Матрица, составленная из элементов $g_ig_j^{-1}$
имеет вид
$A=\left(\begin{smallmatrix}
1      & g      & gh \\
g^{-1} & 1      & h \\
h^{-1}g^{-1} & h^{-1} & 1
\end{smallmatrix}\right)$,
причём эта же градуировка может быть задана набором $(gh,h,1)$.
Если $g\in \lbrace 1, h, h^{-1}, gh, h^{-1}g^{-1} \rbrace$
или $h\in \lbrace 1, g, g^{-1}, gh, h^{-1}g^{-1} \rbrace$, тогда подгруппа группы $G$,
порождённая элементами $g$ и $h$ является абелевой. В силу предложения~\ref{prop:abelian}
в этом случае алгебра $M_3(\mathbbm{k})$ 
может быть переградуирована конечной группой.
Поэтому ниже мы будем считать, что элементы $1, g, h, gh$
попарно различны и условие $\lbrace g,h,gh\rbrace \cap \lbrace g^{-1},h^{-1},h^{-1}g^{-1}\rbrace
\ne \varnothing$ выполнено, если и только если верно одно из равенств $g=g^{-1}$, $h=h^{-1}$, $gh=h^{-1}g^{-1}$, другими словами, только если по крайней мере один из элементов $g$, $h$, $gh$ имеет порядок $2$.

Для трёх возможных равенств имеется $2^3=8$ случая в зависимости от того, выполняются каждое из этих равенств или нет. 

\begin{enumerate}
\item $\lbrace g,h,gh\rbrace \cap \lbrace g^{-1},h^{-1},h^{-1}g^{-1}\rbrace
= \varnothing$.
В этом случае можно воспользоваться теоремой~\ref{TheoremFiniteRegradingAllDifferent},
поскольку во всех клетках матрицы $A$, кроме диагональных, стоят различные элементы группы $G$.

\item $gh = h^{-1}g^{-1}$, но $\lbrace g,h\rbrace \cap \lbrace g^{-1},h^{-1}\rbrace
= \varnothing$.
В этом случае градуировка $\Gamma$ эквивалентна
$\mathbb Z/6\mathbb Z$-градуировке, заданной условиями $g\mapsto \bar 1$ и $h\mapsto \bar 2$,
т.е. при помощи набора $(\bar 3, \bar 2, \bar 0)$.

\item $gh = h^{-1}g^{-1}$, $g=g^{-1}$, но $h\ne h^{-1}$.
В этом случае $\Gamma$
эквивалентна элементарной $S_3$-градуировке, заданной условиями  $g\mapsto (12)$ и $h\mapsto (123)$.

\item $gh = h^{-1}g^{-1}$, $g\ne g^{-1}$, $h = h^{-1}$. В этом случае $\Gamma$
эквивалентна элементарной $S_3$-градуировке, заданной условиями  $g\mapsto (123)$ и $h\mapsto (12)$.

\item $gh = h^{-1}g^{-1}$, $g=g^{-1}$, $h=h^{-1}$. В данном случае градуировка $\Gamma$
эквивалентна элементарной $(\mathbb Z/2\mathbb Z) \times (\mathbb Z/2\mathbb Z)$-градуировке,
заданной условиями $g\mapsto(\bar 0, \bar 1)$ и $h\mapsto(\bar 1, \bar 0)$.
(Поскольку в этом случае $g$ и $h$ коммутируют, мы могли бы вместо этого воспользоваться предложением~\ref{prop:abelian}.)

\item $gh \ne h^{-1}g^{-1}$, $g=g^{-1}$, $h\ne h^{-1}$.
В этом случае градуировка $\Gamma$ эквивалентна
$\mathbb Z/6\mathbb Z$-градуировке, заданной условиями $g\mapsto \bar 3$ и $h\mapsto \bar 1$.

\item $gh \ne h^{-1}g^{-1}$, $g\ne g^{-1}$, $h=h^{-1}$.
Этот случай рассматривается аналогичным образом.

\item $gh \ne h^{-1}g^{-1}$, $g=g^{-1}$, $h=h^{-1}$. В этом случае $\Gamma$
эквивалентна элементарной $S_3$-градуировке, заданной условиями $g\mapsto (12)$, $h \mapsto (13)$.
\end{enumerate}

Все возможные случаи разобраны, откуда градуировка $\Gamma$ действительно эквивалентна
элементарной градуировке конечной группой.
\end{proof}

\begin{remark}
Элементарная $S_3$-градуировка на алгебре $M_3(\mathbbm{k})$, определённая выше при помощи соответствия  $g\mapsto (12)$, $h \mapsto (13)$ неэквивалентна никакой градуировке абелевой группой, поскольку $(12)(13)\ne (13)(12)$.
\end{remark}

\begin{proof}[Доказательство теоремы~\ref{TheoremRegradeElementaryOmega}]
Достаточно воспользоваться следствиями~\ref{CorollaryMatrixElemGreaterSize}, \ref{CorollaryFiniteRegradingImpossible349} и теоремой~\ref{TheoremFiniteRegradingSmallMatrix}.
\end{proof}

         \section{Носитель градуировки как op-слабый 2-функтор
         из категорий $\mathbf{GrnuAlg}_\mathbbm{k}$ и $\mathbf{GrAlg}_\mathbbm{k}$} \label{SectionOplaxSupport}
         
         В \S\ref{SectionGrAdjoint} были приведены два известных примера пар сопряжённых функторов,
         связанных с градуировками. Здесь прежде всего нужно отметить, что в обоих примерах
         используются категории
         $\mathbf{nuAlg}_\mathbbm{k}^{G\text{-}\mathrm{gr}}$, объекты которых градуированы одними и теми же группами $G$, т.е. понятие изоморфизма в этих категориях совпадает с понятием изоморфизма градуировок.
         Кроме того, хотя функтор $U_\varphi$ и меняет градуирующую группу, градуировка на алгебре $U_\varphi(A)$ не обязана быть эквивалентна исходной градуировке на алгебре $A$.

Для работы с эквивалентностями градуировок введём категорию $\mathbf{GrnuAlg}_\mathbbm{k}$ 
всех ассоциативных алгебр необязательно с единицей
над полем $\mathbbm{k}$, градуированных произвольными группами. При этом понятие изоморфизма в $\mathbf{GrnuAlg}_\mathbbm{k}$
будет совпадать с понятием эквивалентности градуировок.

Напомним, что гомоморфизм $\psi \colon A \to B$
между алгебрами $A=\bigoplus_{g \in G} A^{(g)}$
и $B=\bigoplus_{h \in H} B^{(h)}$ называется \textit{(нестрого) градуированным},
если для любого $g\in G$ существует такое $h\in H$, что $\psi(A^{(g)})\subseteq B^{(h)}$.

Любой (нестрого) градуированный гомоморфизм алгебр индуцирует отображение между носителями градуировок.
Такие отображения порождают ряд функторов, которые исследуются ниже.

Итак, пусть $\mathbf{GrnuAlg}_\mathbbm{k}$~"--- категория,
объектами которой являются все групповые градуировки на ассоциативных алгебрах необязательно с единицей
над полем $\mathbbm{k}$, а морфизмами~"--- градуированные гомоморфизмы между соответствующими алгебрами.

Через $\mathbf{GrAlg}_\mathbbm{k}$ обозначим категорию,
объектами которой являются все групповые градуировки на ассоциативных алгебрах с единицей
над полем $\mathbbm{k}$, а морфизмами~"--- градуированные гомоморфизмы алгебр с единицей.

Кроме этого, рассмотрим категорию $\mathcal C$, в которой
\begin{itemize}
\item объектами являются тройки $(G, S, P)$, где $G$~"--- группа, а $S \subseteq G$
 и $P \subseteq S\times S$~"--- такие подмножества, что $\lbrace gh \mid (g,h)\in P \rbrace \subseteq S$;
\item морфизмами  $(G_1, S_1, P_1)\to (G_2, S_2, P_2)$ являются всевозможные тройки $(\psi, R, Q)$,
где $R \subseteq S_1$, $Q\subseteq P_1 \cap (R\times R)$, а $\psi \colon R \to S_2$~"--- такое отображение,
что $\lbrace gh \mid (g,h)\in Q\rbrace
\subseteq R$, $\psi(g)\psi(h)=\psi(gh)$ для всех $(g,h)\in Q$ и $\lbrace(\psi(g),\psi(h))
\mid (g,h)\in Q \rbrace \subseteq P_2$;
\item тождественный морфизм объекта $(G, S, P)$~"--- это тройка $(\id_{S}, S, P)$;
\item если $(\psi_1, R_1, Q_1) \colon (G_1, S_1, P_1)\to (G_2, S_2, P_2)$
и $(\psi_2, R_2, Q_2) \colon (G_2, S_2, P_2)\to (G_3, S_3, P_3)$,
то $$(\psi,R,Q)=(\psi_2, R_2, Q_2)(\psi_1, R_1, Q_1) \colon (G_1, S_1, P_1)\to (G_3, S_3, P_3)$$
задаётся формулами $R=\lbrace g\in R_1 \mid \psi_1(g)\in R_2 \rbrace$, $Q = \lbrace (g,h)\in Q_1 \mid (\psi_1(g),\psi_1(h))\in Q_2 \rbrace$ и $\psi(g)=\psi_2(\psi_1(g))$ при $g\in R$.
\end{itemize}

Существует очевидное функтороподобное\footnote{Под \textit{функтороподобным отображением} $\mathcal F$ из категории $\mathcal A$
в категорию $\mathcal B$ будем понимать такое соответствие, при котором объекту $A$ категории $\mathcal A$ соответствует объект $\mathcal F(A)$ категории $\mathcal B$,
а морфизму $f \colon A_1 \to A_2$ в $\mathcal A$~"--- морфизм $\mathcal F (f) \colon \mathcal F (A_1) \to \mathcal F (A_2)$ в $\mathcal B$. Никаких других ограничений на $\mathcal F$, вообще говоря, не накладывается.} отображение $L \colon \mathbf{GrnuAlg}_\mathbbm{k} \to \mathcal C$,
где $$L(\Gamma):=(G, \supp \Gamma, \lbrace
(g_1,g_2)\in G \times G \mid A^{(g_1)}A^{(g_2)}\ne 0\rbrace)$$
для $\Gamma \colon A=\bigoplus_{g\in G} A^{(g)}$,
а для градуированного морфизма $\varphi \colon \Gamma \to
\Gamma_1$,
где $\Gamma_1 \colon B=\bigoplus_{h\in H} B^{(h)}$,
тройка $L(\varphi)=(\psi, R, Q)$ определяется следующим образом: $R=\lbrace g\in G \mid \varphi\left(A^{(g)}\right)\ne 0 \rbrace$, $$Q=\lbrace(g_1,g_2)\in R\times R\mid \varphi\left(A^{(g_1)}\right)\varphi\left(A^{(g_2)}\right)\ne 0\rbrace,$$
отображение $\psi$ задано при помощи включения $\varphi\left(A^{(g)}\right) \subseteq B^{(\psi(g))}$
 при всех $g\in R$.
 
 Как будет показано ниже  в примере~\ref{ExampleLnotAFunctor},
 отображение $L$ не является обычным функтором: $L(\varphi_1)L(\varphi_2)$
 не всегда совпадает с $L(\varphi_1 \varphi_2)$.
 Для преодоления этого затруднения наделим множество морфизмов
$\mathcal C((G_1, S_1, P_1),(G_2, S_2, P_2))$ частичным порядком $\preccurlyeq$:
будем говорить, что $(\psi_1, R_1, Q_1) \preccurlyeq (\psi_2, R_2, Q_2)$,
если $R_1 \subseteq R_2$, $Q_1 \subseteq Q_2$ и $\psi_1 = \psi_2 \bigl|_{R_1}$. Для любой пары $\varphi_1, \varphi_2$ перемножаемых морфизмов из категории $\mathcal C$ выполняется неравенство \begin{equation}\label{EqLSucc}L(\varphi_1)L(\varphi_2)
\succcurlyeq L(\varphi_1 \varphi_2),\end{equation}
которое, как показывает пример ниже, может быть строгим:

\begin{example}\label{ExampleLnotAFunctor}
Пусть $A=A^{(\bar 0)} \oplus A^{(\bar 1)}$~"--- $\mathbb Z/2\mathbb Z$-градуированная алгебра,
где $$A^{(\bar 0)}=\mathbbm{k}1_A,\quad A^{(\bar 1)}=\mathbbm{k}a\oplus \mathbbm{k}b,\quad a^2=ab=ba=b^2=0.$$
Пусть $\varphi \colon A \to A$~"--- градуированный морфизм, заданный равенствами $\varphi(1_A)=1_A$, $\varphi(a)=b$, $\varphi(b)=0$.
Тогда $$L(A)=(\mathbb Z/2\mathbb Z, \mathbb Z/2\mathbb Z, (\mathbb Z/2\mathbb Z)^2 \backslash \lbrace (\bar 1, \bar 1)\rbrace),$$ $$L(\varphi)=(\id_{\mathbb Z/2\mathbb Z}, \mathbb Z/2\mathbb Z,
(\mathbb Z/2\mathbb Z)^2 \backslash \lbrace (\bar 1, \bar 1)\rbrace)$$
и $L(\varphi)^2=L(\varphi)$, однако $$L(\varphi^2)=(\id_{\lbrace \bar 0 \rbrace}, \lbrace \bar 0 \rbrace, \lbrace (\bar 0, \bar 0) \rbrace)\prec L(\varphi)^2.$$
\end{example}

Подберём для $L$
соответствующее категорное понятие, используя теорию обогащённых категорий.
Категория $\mathcal A$ называется
\textit{обогащённой} над категорией $\mathcal B$, если 
 hom-объекты $\mathcal A(A,B)$ категории $\mathcal A$ являются объектами категории $\mathcal B$,
 а композиция морфизмов и назначение тождественного морфизма являются
 морфизмами в категории $\mathcal B$. (См. точное определение в~\cite[\S 1.2]{KellyEnriched}.)

Напомним, что всякое частично упорядоченное множество $(M, \preccurlyeq)$
является категорией, в которой объектами являются элементы $m\in M$ и если $m \preccurlyeq n$, то существует единственный морфизм $m \to n$. Если же $m \not\preccurlyeq n$, то морфизмов $m \to n$ нет.
Обозначим через $\mathbf{Cat}$ категорию малых категорий. Поскольку понятие $\mathbf{Cat}$-обогащённой категории совпадает с понятием $2$-категории, всякая категория обогащённая над
категорией частично упорядоченных множеств является $2$-категорией.
Частичное упорядочение $\preccurlyeq$ превращает $\mathcal C$ 
в категорию, обогащённую над категорией частично упорядоченных множеств и, следовательно,
в $2$-категорию, где
$0$-клетками являются тройки $(G, S, P)$, $1$-клетками являются тройки $(\psi, R, Q)$
и между $1$-клеткой $(\psi_1, R_1, Q_1)$ и $1$-клеткой $(\psi_2, R_2, Q_2)$ существует $2$-клетка,
если и только если $(\psi_1, R_1, Q_1) \preccurlyeq (\psi_2, R_2, Q_2)$.

 Между $2$-категориями можно рассматривать \textit{строгие} $2$-функторы,
 от которых требуется, чтобы они сохраняли операцию умножения и тождественные морфизмы $1$- и $2$-клеток, 
а  можно рассматривать \textit{(op-)слабые} $2$-функторы $F$, для которых
требование $F(\varphi_1)F(\varphi_2)=F(\varphi_1 \varphi_2)$
 для всех $1$-клеток $\varphi_1,\varphi_2$ заменяется требованием существования $2$-клетки
 между $1$-клетками
 $F(\varphi_1)F(\varphi_2)$ и $F(\varphi_1 \varphi_2)$.
Разница между слабым $2$-функтором и op-слабым $2$-функтором
заключается в направлении, в котором идёт эта $2$-клетка.
 (См. детали и точное определение (op-)слабого 2-функтора в \cite[с.~83]{KellyTwoCat} и \cite[с.~29]{BenabouBi}.)

В этой терминологии неравенство~(\ref{EqLSucc})
означает, что существует $2$-клетка, которая идёт из $1$-клетки $L(\varphi_1)L(\varphi_2)$
в $1$-клетку $L(\varphi_1 \varphi_2)$. Это превращает $L$ в op-слабый $2$-функтор из $\mathbf{GrnuAlg}_\mathbbm{k}$
в $\mathcal C$, если мы будем рассматривать $\mathbf{GrnuAlg}_\mathbbm{k}$ в качестве $2$-категории
с дискретными hom-категориями (т.е. единственными $2$-клетками в $\mathbf{GrnuAlg}_\mathbbm{k}$ являются тождественные $2$-клетки между морфизмами).
Все требуемые равенства из определения op-слабого $2$-функтора выполняются, 
поскольку в категориях, являющихся частично упорядоченными множествами все диаграммы являются коммутативными.

Ясно, что можно ограничить область определения op-слабого $2$-функтора $L$ до категории $\mathbf{GrAlg}_\mathbbm{k}$
и рассматривать то же явление в случае градуированных алгебр с единицей.

\section{Универсальная группа градуировки как функтор}\label{SectionUniversalGradingGroupFunctors}

В предыдущем параграфе было введено функтороподобное отображение $L$, которое ставит в соответствие
каждой градуировке её носитель, и было показано, что $L$~"--- это не обычный, а op-слабый 2-функтор.
В данном параграфе мы так изменим $L$, чтобы получить обычный функтор. Во-первых, мы ограничим множество допустимых морфизмов. Во-вторых, вместо носителя градуировки каждой градуировке мы будем ставить в соответствие группу, порождённую носителем, в которой элементы носителя удовлетворяют соотношениям, накладываемым градуировкой, т.е. универсальную группу градуировки.

Будем называть градуированный гомоморфизм \textit{градуированно инъективным},
если его ограничение на каждую однородную компоненту является инъективным отображением.
Например,  любой градуированный гомоморфизм алгебр с единицей групповой алгебры в любую градуированную алгебру градуированно инъективный, поскольку все групповые элементы отображаются в обратимые элементы.

При этом отображение аугментации $\varepsilon \colon \mathbbm{k}G \to \mathbbm{k}$, где $\varepsilon\left(\sum_g \alpha_g g \right)
\mapsto \sum_g \alpha_g$, будучи всегда градуированно инъективным гомоморфизмом, является инъективным гомоморфизмом только тогда, когда $G$~"--- тривиальная группа.

Введём категорию $\widetilde{\mathbf{GrnuAlg}_\mathbbm{k}}$,
объектами которой являются всевозможные групповые градуировки на ассоциативных алгебрах над полем $\mathbbm{k}$ необязательно с единицей, а морфизмами являются градуированно инъективные гомоморфизмы между соответствующими алгебрами. Тогда соответствие, описанное в начале параграфа, задаёт функтор $R \colon \widetilde{\mathbf{GrnuAlg}_\mathbbm{k}} \to \mathbf{Grp}$
где $R(\Gamma):=G_{\Gamma}$ для всякого $\Gamma \colon A=\bigoplus_{g\in G} A^{(g)}$,
а для всякого градуированно инъективного гомоморфизма $\varphi \colon
\Gamma \to \Gamma_1$, где $\Gamma_1 \colon B=\bigoplus_{h\in H} B^{(h)}$,
гомоморфизм групп $R(\varphi) \colon G_{\Gamma} \to G_{\Gamma_1}$ определяется через $R(\varphi)(g):=h$,
где $\varphi\left(A^{(g)}\right) \subseteq B^{(h)}$, $g\in G$.

Можно ограничится лишь алгебрами с единицей.
Обозначим через $\widetilde{\mathbf{GrAlg}_\mathbbm{k}}$
категорию, объектами которой являются всевозможные групповые градуировки на ассоциативных алгебрах
с единицей над полем $\mathbbm{k}$, а морфизмами являются градуированно инъективные гомоморфизмы алгебр с единицей.
Обозначим через $R_1$ 
ограничение функтора $R$ на категорию $\widetilde{\mathbf{GrAlg}_\mathbbm{k}}$. Будем называть $R$ и $R_1$
\textit{функторами универсальной группы градуировки}.

В случаях, когда ясно, какая именно градуировка задана на алгебре $A$,
мы будем отождествлять $A$ с градуировкой, заданной на $A$, и считать $A$ объектом соответствующей категории.

В следующих параграфах изучаются категории $\mathbf{GrnuAlg}_\mathbbm{k}$,
$\mathbf{GrAlg}_\mathbbm{k}$, $\widetilde{\mathbf{GrnuAlg}_\mathbbm{k}}$
и $\widetilde{\mathbf{GrAlg}_\mathbbm{k}}$, а также функторы $R$ и $R_1$, введённые выше.

\section{Категорные конструкции в категориях $\mathbf{GrnuAlg}_\mathbbm{k}$ и $\mathbf{GrAlg}_\mathbbm{k}$}\label{SectionLimColimGrAlg}

В этом разделе исследуются (ко)произведения и (ко)уравнители
в категориях
$\mathbf{GrnuAlg}_\mathbbm{k}$ и $\mathbf{GrAlg}_\mathbbm{k}$.
Особый интерес именно к этим примерам (ко)пределов
связан с тем хорошо известным фактом, что если в какой-то
категории существуют все (ко)пределы и (ко)уравнители,
то в ней существуют и все (ко)пределы. (См., например, теорему 2 из \S 2 главы V \cite{MacLaneCatWork}.)
В предложении~\ref{PropositionCriterionForMonomorphismGrAlg} ниже
описываются мономорфизмы в категориях
$\mathbf{GrnuAlg}_\mathbbm{k}$ и $\mathbf{GrAlg}_\mathbbm{k}$ и сравниваются с инъективными и градуированно инъективными гомоморфизмами.

Легко видеть, что произведение градуированной алгебры $A$ и нулевой алгебры изоморфно $A$.
Кроме того, если две алгебры $A$ и $B$ наделены тривиальной градуировкой, то их произведение как неградуированных алгебр (снова с тривиальной градуировкой) будет их произведением и в $\mathbf{GrnuAlg}_\mathbbm{k}$,
а если $A$ и $B$~"--- алгебры с единицей, то и в $\mathbf{GrAlg}_\mathbbm{k}$.

Оказывается, в любом другом случае произведение двух объектов в $\mathbf{GrnuAlg}_\mathbbm{k}$ и $\mathbf{GrAlg}_\mathbbm{k}$
не существует:

\begin{theorem}\label{TheoremAbsenceOfProductsGrAlg}
Пусть $\Gamma_1 \colon A = \bigoplus_{g\in G} A^{(g)}$ и $\Gamma_2 \colon B = \bigoplus_{h\in H} B^{(h)}$~"---
групповые градуировки на алгебрах над полем $\mathbbm{k}$, такие, что 
$\supp \Gamma_1$ состоит по крайней мере из двух различных элементов, а $B\ne 0$.
Тогда произведение $A$ и $B$ не существует ни в $\mathbf{GrnuAlg}_\mathbbm{k}$, ни в $\mathbf{GrAlg}_\mathbbm{k}$.
\end{theorem}
\begin{proof}
Будем доказывать утверждение для обоих случаев $\mathbf{GrnuAlg}_\mathbbm{k}$ и
$\mathbf{GrAlg}_\mathbbm{k}$ одновременно. В последнем случае будем считать, что алгебры $A$ и $B$ с единицей.

Выберем $g,t\in G$, $h\in H$, $a^{(g)} \in A^{(g)}$, $a^{(t)} \in A^{(t)}$,
$b^{(h)} \in B^{(h)}$, такие, что $g\ne t$, $a^{(g)}\ne 0$, $a^{(t)}\ne 0$,
$b^{(h)}\ne 0$.

Предположим, что произведение $A \times B$ существует. Пусть $\pi \colon A\times B \to A$
и $\rho \colon A\times B \to B$~"--- соответствующие проекции.

Обозначим через $D$ 
свободную ассоциативную некоммутативную алгебру (с единицей, если рассматривается категория $\mathbf{GrAlg}_\mathbbm{k}$
и без единицы, если рассматривается категория $\mathbf{GrnuAlg}_\mathbbm{k}$) над полем $\mathbbm{k}$ от переменных $x, y$
с $\mathbb Z$-градуировкой по общей степени одночленов.
Из определения произведения следует, что для любых морфизмов $\alpha \colon D \to A$ и $\beta \colon D \to B$ существует единственный морфизм $\psi \colon D \to A\times B$,
делающий диаграмму ниже коммутативной:

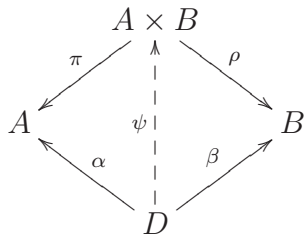
$$\xymatrix{ & A \times B
\ar[ld]_{\pi} \ar[rd]^{\rho} & \\
A & & B \\
& D \ar[lu]_{\alpha} \ar[ru]^{\beta} \ar@{-->}[uu]^{\psi} &
}
$$

Определим градуированные гомоморфизмы $\alpha_1, \alpha_2 \colon D \to A$
и $\beta \colon D \to B$ при помощи формул $\alpha_1(x)=\alpha_2(x)=0$, $\alpha_1(y)=a^{(g)}$,
 $\alpha_2(y)=a^{(t)}$, $\beta(x)=\beta(y)=b^{(h)}$.
 Пусть теперь $\psi_1, \psi_2 \colon D \to A\times B$~"--- единственные такие градуированные
 гомоморфизмы, что $\pi \psi_i = \alpha_i$ и $\rho\psi_i = \beta$.

Обозначим через $C$ 
алгебру многочленов от переменной $x$ с коэффициентами из $\mathbbm{k}$ без свободного члена, если рассматривается категория $\mathbf{GrnuAlg}_\mathbbm{k}$, и со свободным членом, если рассматривается категория $\mathbf{GrAlg}_\mathbbm{k}$. 
Введём на $C$ градуировку группой $\mathbb Z$ по степени одночленов. 
Тогда всякий морфизм $\alpha \colon C \to A$ и $\beta \colon C \to B$
однозначно определяется выбором однородного элемента $\alpha(x)$ и $\beta(x)$.
Пусть $\tau \colon C \hookrightarrow D$~"--- вложение, заданное формулой $\tau(x)=x$.
Поскольку $\pi\psi_1\tau(x)=\alpha_1(x)=\alpha_2(x)=\pi\psi_2\tau(x)$
и $\rho\psi_1\tau(x)=\beta(x)=\rho\psi_2\tau(x)$,
из универсального свойства произведения следует, что $\psi_1\tau = \psi_2\tau$ and $\psi_1(x)=\psi_2(x)$.
В то же время из $\rho(\psi_1(x))=b^{(h)}\ne 0$ получаем $\psi_1(x)=\psi_2(x)\ne 0$.
Используя аналогичные рассуждения, заключаем, что $\psi_1(y)\ne 0$ и $\psi_2(y)\ne 0$.
Поскольку $x$ и $y$ принадлежат одной и той же однородной компоненте алгебры $D$,
$\psi_1(x)$ и $\psi_1(y)$ принадлежат одной и той же однородной компоненте алгебры $A\times B$.
Аналогично, $\psi_2(x)$ и $\psi_2(y)$ принадлежат одной и той же однородной компоненте алгебры $A\times B$.
Поскольку все вышеуказанные элементы ненулевые и $\psi_1(x)=\psi_2(x)$, получаем, что $\psi_1(y)$ и $\psi_2(y)$ принадлежат одной и той же однородной компоненте алгебры $A\times B$. В частности, $a^{(g)}=\pi(\psi_1(y))$
и $a^{(t)}=\pi(\psi_2(y))$ принадлежат одной и той же однородной компоненте алгебры $A$.
Тогда $g=t$, и мы приходим к противоречию.
Следовательно, $A\times B$ не существует ни в
$\mathbf{GrnuAlg}_\mathbbm{k}$, ни в $\mathbf{GrAlg}_\mathbbm{k}$.
\end{proof}

Приведём теперь пример градуированных алгебр, у которых не существует копроизведения.

\begin{proposition}
Пусть $A$~"--- ассоциативная алгебра над полем $\mathbbm{k}$ с базисом $1_A, a, b, c, d, cd$,
причём все попарные произведения элементов
$a,b,c,d$, за исключением $cd=dc\ne 0$, нулевые.
Пусть $B=\langle 1_B, v\rangle_\mathbbm{k}$, $v^2=0$. Обозначим через $\xi$ порождающий циклической группы $C_3$
порядка $3$. Зададим $C_3$-градуировку на $A$ при помощи равенств $A^{(1)}=\langle 1_A, cd\rangle_\mathbbm{k}$, $A^{(\xi)}=\langle a,c\rangle_\mathbbm{k}$, $A^{\left(\xi^2\right)}=\langle b,d\rangle_\mathbbm{k}$
и тривиальную градуировку на $B$.
Тогда копроизведение градуированных алгебр $A$ и $B$ не существует
ни в $\mathbf{GrnuAlg}_\mathbbm{k}$, ни в $\mathbf{GrAlg}_\mathbbm{k}$.
\end{proposition}

\begin{proof}
Предположим, что копроизведение $A\sqcup B$ существует.
Обозначим через $i\colon A \to A\sqcup B$ и $j\colon B \to A\sqcup B$
соответствующие морфизмы.

Пусть $f\colon B \to A$~"--- гомоморфизм алгебр с единицей, заданный равенством $f(v)=0$.
Тогда существует единственный такой морфизм $h\colon A \sqcup B \to A$, что диаграмма
ниже коммутативна:

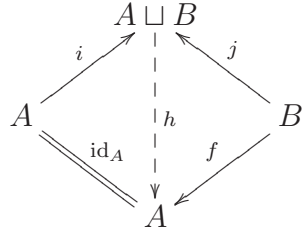
$$\xymatrix{ & A \sqcup B
\ar@{-->}[dd]^h & \\
A \ar[ru]^{i}\ar@{=}[rd]^{\id_A} & & B \ar[lu]_{j}\ar[ld]_{f}
 \\
& A   &
}
$$

Поскольку $cd=h(i(cd))\ne 0$, справедливо равенство \begin{equation}\label{EqIcIdNe0}
i(c)i(d)=i(d)i(c)\ne 0.
\end{equation}

Рассмотрим ассоциативную алгебру $C$ с базисом $\lbrace 1_C, x, y, z, xy, yz, zy, yx, xyz, zyx\rbrace$,
где все прочие произведения элементов $x,y,z$ нулевые.
Наделим алгебру $C$ градуировкой свободной группой $\mathcal F (X,Z)$ при помощи равенств $x\in C^{(X)}$; $1_C,y\in C^{(1)}$; $z \in C^{(Z)}$. 
Определим градуированные гомоморфизмы алгебр с единицей $\alpha \colon A \to C$
и $\beta \colon B \to C$ при помощи равенств
$$\alpha(a)=x,\text{ } \alpha(b)=z,\text{ } \alpha(c)=0, \text{ } \alpha(d)=0, \text{ }\beta(v)=y.$$

Тогда существует единственный такой морфизм $\psi \colon A\sqcup B \to C$,
что диаграмма ниже коммутативна:

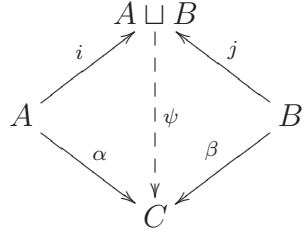
$$\xymatrix{ & A \sqcup B
\ar@{-->}[dd]^\psi & \\
A \ar[ru]^{i}\ar[rd]^{\alpha} & & B \ar[lu]_{j}\ar[ld]_{\beta}
 \\
& C    &
}
$$

Обозначим через $G$
 градуирующую группу для алгебры $A \sqcup B$. Поскольку гомоморфизмы $i$ и $j$
отображают однородные элементы в однородные элементы, существуют такие $g,h, t \in G$,
что
$$i(a)\in (A \sqcup B)^{(g)};\text{ } i(b)\in (A \sqcup B)^{(h)}; \text{ } j(1_B), j(v)\in (A \sqcup B)^{(t)}.$$ Из того, что
$a,c \in A^{(\xi)}$ и $b,d \in A^{\left(\xi^2\right)}$,
следует, что $i(c)\in (A \sqcup B)^{(g)}$ и
$i(d)\in (A \sqcup B)^{(h)}$.
Используя~(\ref{EqIcIdNe0}), получаем, что $gh=hg$.
В силу того, что $\psi(j(v))=\beta(v)=y\ne 0$,
справедливо неравенство $j(v) \ne 0$.
В то же время $j(1_B)j(v)=j(v)$. Следовательно, $t^2=t$ и $t=1_G$. Отсюда оба элемента $i(a)j(v)i(b)$
и $i(b)j(v)i(a)$ принадлежат $(A \sqcup B)^{(s)}$, где $s=g t h = gh = hg = htg$.
Однако элементы $\psi(i(a)j(v)i(b))=xyz$ и $\psi(i(b)j(v)i(a))=zyx$
принадлежат разным однородным компонентам алгебры $C$.
Получаем противоречие, откуда копроизведение $A \sqcup B$ не существует.
\end{proof}

Обратимся теперь к описанию уравнителей.

\begin{proposition}\label{PropositionEqualizerGrAlg}
Пусть $\alpha, \beta \colon A \to B$~"--- гомоморфизмы градуированных
алгебр.
Обозначим через $C$
линейную оболочку всех таких однородных элементов $a\in A$,
что $\alpha(a)=\beta(a)$. Пусть $i \colon C \to A$~"--- соответствующее вложение. 
Тогда $i$ является уравнителем морфизмов $\alpha$ и $\beta$
как в $\mathbf{GrnuAlg}_\mathbbm{k}$, так и в $\mathbf{GrAlg}_\mathbbm{k}$.
(Как обычно, в случае категории $\mathbf{GrAlg}_\mathbbm{k}$ гомоморфизмы
 $\alpha$ и $\beta$ должны быть гомоморфизмами алгебр с единицей.)
\end{proposition}
\begin{proof}
Очевидно, что $C$~"--- градуированная подалгебра алгебры $A$.
Требуется доказать, что для любого гомоморфизма $\gamma \colon D \to A$
градуированных алгебр, такого, что $\alpha\gamma=\beta\gamma$,
существует единственный такой градуированный гомоморфизм
$\varphi \colon D \to C$, что $\gamma = i\varphi$:

$\xymatrix{ C \ar[r]^i  & A \ar@<-0.5ex>[r]_\beta \ar@<0.5ex>[r]^\alpha & B \\
                       & D \ar[u]_\gamma \ar@{-->}[lu]^\varphi}$

Однако это утверждение следует из того факта, что для любого однородного элемента $d\in D$
справедливо равенство $\alpha(\gamma(d))=\beta(\gamma(d))$,
откуда $\gamma(d)\in C$. Поскольку $D$ является линейной оболочкой
своих однородных компонент, $\varphi(D)\subseteq C$
и $\varphi$ получается из $\gamma$ ограничением области значений.
\end{proof}

Приведём теперь пример параллельной\footnote{Пара морфизмов $\alpha_1 \colon A_1 \to B_1$ и $\alpha_2 \colon A_2 \to B_2$ называется \textit{параллельной}, если $A_1=A_2$ и $B_1=B_2$.} пары градуированных гомоморфизмов,
у которой не существует коуравнителя.
Напомним, что в наших обозначениях $\mathcal F (x,y)$~"--- это свободная группа со свободными порождающими $x,y$.

\begin{proposition}
Пусть $G=\mathcal F (x,y)$. Обозначим через $B$ алгебру
над полем $\mathbbm{k}$ с базисом $1, a,b,c,d$ и единицей $1$,
умножение в которой задаётся равенством $\langle a,b,c,d \rangle_\mathbbm{k}^2=0$.
Определим $G$-градуировку на $B$ при помощи равенств
$$B^{(1)}=\mathbbm{k}1_B,\text{ } B^{(x)}=\langle a,c \rangle_\mathbbm{k}, \text{ }B^{(y)}=\langle b,d \rangle_\mathbbm{k}.$$
Пусть $A:=\langle 1, a, b\rangle_\mathbbm{k} \subset B$.
Обозначим через $\alpha \colon A \to B$ естественное вложение,
 а через $\beta \colon A \to B$
гомоморфизм алгебр с единицей, заданный равенствами $\beta(a)=b$
и $\beta(b)=a$.
Тогда коуравнитель
 гомоморфизмов $\alpha$ и $\beta$ не существует ни в $\mathbf{GrnuAlg}_\mathbbm{k}$, ни в $\mathbf{GrAlg}_\mathbbm{k}$.
\end{proposition}
\begin{proof}
Предположим, что $\gamma \colon B \to C$~"--- коуравнитель морфизмов $\alpha$ и $\beta$.
Тогда \begin{equation}\label{EqNoCoeqInGrAlg}
\gamma(a)=\gamma(\alpha(a))=\gamma(\beta(a))=\gamma(b).\end{equation}

 Обозначим через $\varphi \colon B \to D$, где $D:=\langle 1,a\rangle_\mathbbm{k} \subset B$,
гомоморфизм алгебр с единицей, заданный при помощи равенств
$$\varphi(a)=\varphi(b)=\varphi(c)=\varphi(d)=a.$$
Тогда существует единственный такой морфизм $\psi \colon C \to D$, что $\psi\gamma = \varphi$:

$\xymatrix{ A \ar@<-0.5ex>[r]_\beta \ar@<0.5ex>[r]^\alpha & B  \ar[r]^\gamma \ar[d]^\varphi & C  \ar@{-->}[ld]^\psi \\
           &  D  }$

Отсюда $\psi(\gamma(a))=\varphi(a)=a\ne 0$ и $\gamma(a)\ne 0$.

Пусть $\theta \colon B \to A$~"--- гомоморфизм алгебр с единицей, заданный при помощи равенств $\theta(a)=\theta(b)=0$,
$\theta(c)=a$, $\theta(d)=b$. Тогда, в силу наших предположений должен существовать единственный такой морфизм $\mu \colon C \to A$,
что $\mu\gamma = \theta$.
Из того, что элементы $\mu(\gamma(c))=a \ne 0$ и $\mu(\gamma(d))=b\ne 0$ 
принадлежат различным компонентам алгебры $A$, следует, что
элементы
$\gamma(c)$ и $\gamma(d)$ принадлежат различным компонентам алгебры $C$. Отсюда
элементы
$\gamma(a)$ и $\gamma(b)$ принадлежат различным компонентам алгебры $C$,
что противоречит~(\ref{EqNoCoeqInGrAlg}), поскольку выше было показано, что $\gamma(a)\ne 0$.
\end{proof}

Опишем теперь мономорфизмы и сравним их с инъективными и градуированно инъективными гомоморфизмами.

\begin{proposition}\label{PropositionCriterionForMonomorphismGrAlg}
Пусть $f\colon A=\bigoplus_{g\in G} A^{(g)} \to B$~"--- морфизм в $\mathbf{GrnuAlg}_\mathbbm{k}$ или в $\mathbf{GrAlg}_\mathbbm{k}$. Тогда $f$~"--- мономорфизм, если и только если из неравенства $a \ne b$ для некоторых $a,b\in \bigcup_{g\in G}
A^{(g)}$ всегда следует $f(a)\ne f(b)$. В частности, всякий мономорфизм в $\mathbf{GrnuAlg}_\mathbbm{k}$ и в $\mathbf{GrAlg}_\mathbbm{k}$ градуированно инъективный.
\end{proposition}
\begin{proof}
Пусть $C$~"--- алгебра всех многочленов
от переменной $x$ с коэффициентами из поля $\mathbbm{k}$  без свободного члена, если рассматривается категория $\mathbf{GrnuAlg}_\mathbbm{k}$, и со свободным членом, если рассматривается категория $\mathbf{GrAlg}_\mathbbm{k}$,
наделённая $\mathbb Z$-градуировкой по степеням.
Тогда для любого однородного $a\in A$ 
существует единственный градуированный гомоморфизм $\lambda \colon C \to A$,
такой, что $\lambda(x)=a$. Следовательно, если бы существовали такие $a,b\in \bigcup_{g\in G}
A^{(g)}$, $a \ne b$, что
$f(a)=f(b)$, можно было бы построить 
градуированные гомоморфизмы $\lambda, \mu \colon
C \to A$, такие, что $\lambda(x)=a$, $\mu(x)=b$, т.е. $\lambda\ne \mu$, но $f\lambda=f\mu$.
Таким образом, необходимость доказана.

Предположим теперь, что $f\colon A \to B$~"--- такой градуированный гомоморфизм,
что для любых однородных $a,b\in A$, $a \ne b$, справедливо неравенство $f(a)\ne f(b)$.
Пусть теперь $\alpha,\beta \colon D \to A$~"--- два градуированных
гомоморфизма, таких, что
$f\alpha = f\beta$. Тогда для любого однородного элемента $d\in D$
имеем $\alpha(d),\beta(d)\in \bigcup_{g\in G} A^{(g)}$
и $f(\alpha(d))=f(\beta(d))$. Следовательно, $\alpha(d)=\beta(d)$. Поскольку $D$ является линейной оболочкой своих градуированных компонент, справедливо равенство $\alpha = \beta$ и $f$ действительно является мономорфизмом.
\end{proof}

Отметим, что не всякий градуированно инъективный гомоморфизм является мономорфизм
в $\mathbf{GrnuAlg}_\mathbbm{k}$ или в $\mathbf{GrAlg}_\mathbbm{k}$.  Например, отображение аугментации $\varepsilon \colon \mathbbm{k}G \to \mathbbm{k}$, $\varepsilon\left(\sum_{g\in G} \alpha_g g\right):=\sum_{g\in G} \alpha_g$,
градуированно инъективно, но не является мономорфизмом в случае нетривиальной группы $G$.

Дадим теперь пример мономорфизма, который не является инъективным отображением.

\begin{example}
Пусть $A=\mathbbm{k}1_A\oplus \mathbbm{k}a_1\oplus \mathbbm{k}a_2 \oplus \mathbbm{k}a_3$, где $\langle a_1, a_2, a_3 \rangle_\mathbbm{k}^2=0$.
Определим $\mathbb Z/4\mathbb Z$-градуировку на $A$ при помощи условий $a_i \in A^{(\bar i)}$
при $i=1,2,3$. Пусть алгебра $B = A/(a_1+a_2+a_3)$ наделена тривиальной градуировкой.
Тогда естественная сюрьекция $\pi \colon A \twoheadrightarrow B$
является мономорфизмом как в $\mathbf{GrnuAlg}_\mathbbm{k}$, так и в $\mathbf{GrAlg}_\mathbbm{k}$.
\end{example}
\begin{proof}
Пусть $a,b \in \bigcup_{\bar i \in \mathbb Z/4\mathbb Z} A^{(\bar i)}$,
такие, что $\pi(a)=\pi(b)$. Тогда $a-b = \alpha(a_1+a_2+a_3)$
для некоторого $\alpha \in \mathbbm{k}$. Теперь заметим, что коэффициент $\alpha$
обязан быть нулевым, поскольку иначе элемент
$a-b$ имел бы ненулевые компоненты в каждом из подпространств $A^{(\bar i)}$, $\bar i \ne \bar 0$.
Отсюда $a=b$ и из предложения~\ref{PropositionCriterionForMonomorphismGrAlg}
следует, что $\pi$ является мономорфизмом как в $\mathbf{GrnuAlg}_\mathbbm{k}$, так и в $\mathbf{GrAlg}_\mathbbm{k}$.
\end{proof}

\section{Категорные конструкции в категориях $\widetilde{\mathbf{GrnuAlg}_\mathbbm{k}}$
и $\widetilde{\mathbf{GrAlg}_\mathbbm{k}}$}\label{SectionLimColimTildeGrAlg}

Начнём с примера градуированных алгебр, произведение которых существует как в $\widetilde{\mathbf{GrnuAlg}_\mathbbm{k}}$, так и в $\widetilde{\mathbf{GrAlg}_\mathbbm{k}}$,
но сильно отличается от их произведения как неградуированных алгебр.

\begin{proposition}\label{PropositionProductGroupAlgebrasTildeGrAlg} Пусть $G$ и $H$~"---
группы, а $\mathbbm{k}$~"--- поле.
Тогда групповая алгебра $\mathbbm{k}(\mathbbm{k}^{\times} \times G \times H)$, наделённая стандартной градуировкой,
является произведением групповых алгебр $\mathbbm{k}G$ и $\mathbbm{k}H$
как в $\widetilde{\mathbf{GrnuAlg}_\mathbbm{k}}$, так и в $\widetilde{\mathbf{GrAlg}_\mathbbm{k}}$.
\end{proposition}
\begin{proof}
Пусть $\pi_1 \colon \mathbbm{k}(\mathbbm{k}^{\times} \times G \times H) \to \mathbbm{k}G$
и $\pi_2 \colon \mathbbm{k}(\mathbbm{k}^{\times} \times G \times H) \to \mathbbm{k}H$~"--- гомоморфизмы,
заданные равенствами $\pi_1\bigl((\alpha,g,h)\bigr)=\alpha g$
и $\pi_2\bigl((\alpha,g,h)\bigr)=h$
для $\alpha \in \mathbbm{k}^\times$, $g\in G$, $h\in H$.
Ясно, что они являются градуированно инъективными.

Предположим, что существует градуированная алгебра $A$ и градуированно инъективные
гомоморфизмы
$\varphi_1 \colon A \to \mathbbm{k}G$ и $\varphi_2 \colon A \to \mathbbm{k}H$.
Докажем, что существует единственный градуированно инъективный гомоморфизм
$\varphi \colon A \to \mathbbm{k}(\mathbbm{k}^{\times} \times G \times H)$,
такой, что диаграмма ниже коммутативна:

$$\xymatrix{ & \mathbbm{k}(\mathbbm{k}^{\times} \times G \times H)
\ar[ld]_{\pi_1} \ar[rd]^{\pi_2} & \\
\mathbbm{k}G & & \mathbbm{k}H \\
& A \ar[lu]_{\varphi_1} \ar[ru]^{\varphi_2} \ar@{-->}[uu]^\varphi &
}
$$
Вначале заметим, что поскольку всякая однородная компонента алгебры $\mathbbm{k}G$
имеет размерность $1$ и гомоморфизм $\varphi_1$ градуированно инъективен,
каждая однородная компонента алгебры $A$ 
должна иметь размерность не более $1$.
Предположим теперь, что градуированно инъективный гомоморфизм $\varphi$
действительно существует. Пусть $a \in A$~"--- ненулевой
однородный элемент.
Тогда $\varphi(a)$ также должен быть ненулевым однородным элементом,
т.е. $\varphi(a)=\alpha {(\beta, g,h)}$
для некоторых элементов поля $\alpha, \beta \in \mathbbm{k}^\times$ и элементов групп $g\in G$, $h\in H$.
При этом $\varphi_1(a)=\pi_1\varphi(a)= \alpha\beta g$
и $\varphi_2(a)=\pi_2\varphi(a)= \alpha h$, т.е. элемент $\varphi(a)$
однозначно определён элементами $\varphi_1(a)$ и $\varphi_2(a)$.

Отсюда при заданных $\varphi_1$ и $\varphi_2$ гомоморфизм $\varphi$
определяется следующим образом:
если $\varphi_1(a)=\lambda g$ и $\varphi_2(a)=\mu h$, то
$\varphi(a)=\mu {(\lambda/\mu, g,h)}$.
\end{proof}
\begin{corollary}
Если поле $\mathbbm{k}$ состоит более чем из $2$ элементов, то функторы $R$ и $R_1$ не имеют левых сопряжённых.
\end{corollary}
\begin{proof} Всякий функтор, обладающий левым сопряжённым, сохраняет пределы (см., например, теорему 1 из \S 5 главы V в~\cite{MacLaneCatWork})
и, в частности, произведения. Однако, взяв произвольные конечные группы $G$ и $H$, получаем $$R(\mathbbm{k}(\mathbbm{k}^{\times} \times G \times H))=
R_1(\mathbbm{k}(\mathbbm{k}^{\times} \times G \times H))=\mathbbm{k}^{\times} \times G \times H
\ncong R(\mathbbm{k}G) \times R(\mathbbm{k}H)=G\times H.$$
\end{proof}

В следующем параграфе в теореме~\ref{TheoremAbsenceOfLeftAdjointUGGF}
будет доказано, что ограничение на мощность поля $\mathbbm{k}$ является излишним,
т.е. функторы $R$ и $R_1$ не имеют левых сопряжённых в случае любого поля $\mathbbm{k}$.

Приведём пример алгебр, произведение которых не существует ни в
$\widetilde{\mathbf{GrnuAlg}_\mathbbm{k}}$, ни в $\widetilde{\mathbf{GrAlg}_\mathbbm{k}}$.

\begin{theorem}
Пусть $A=\mathbbm{k}[a_1,a_2]$ и $B=\mathbbm{k}[b_1,b_2,b_3]$~"--- алгебры многочленов от коммутирующих переменных с коэффициентами из поля $\mathbbm{k}$, снабжённые $\mathbb Z$-градуировкой по степеням одночленов.
Тогда произведение объектов $A$ и $B$ не существует ни в $\widetilde{\mathbf{GrnuAlg}_\mathbbm{k}}$, ни в $\widetilde{\mathbf{GrAlg}_\mathbbm{k}}$.
\end{theorem}
\begin{proof}
Пусть $A\times B$~"--- произведение объектов $A$ и $B$, а
$\pi \colon A\times B \to A$ и $\rho \colon A\times B \to B$~"--- соответствующие проекции.
Снова обозначим через $C$ 
алгебру многочленов от переменной $x$ с коэффициентами из $\mathbbm{k}$ без свободного члена, если рассматривается категория $\mathbf{GrnuAlg}_\mathbbm{k}$, и со свободным членом, если рассматривается категория $\mathbf{GrAlg}_\mathbbm{k}$,
и введём на $C$ градуировку группой $\mathbb Z$ по степени одночленов. 

Будем использовать тот же трюк, что и в доказательстве теоремы~\ref{TheoremAbsenceOfProductsGrAlg}.
Для любых ненулевых (и, следовательно, ненильпотентных)
однородных элементов $a \in A$ и $b\in B$
существуют единственные такие гомоморфизмы $\alpha \colon C \to A$
и $\beta \colon C \to B$, что $\alpha(x)= a$ и $\beta(x)=b$.

Из определения произведения следует, что существует единственный градуированно инъективный гомоморфизм
 $\psi \colon C \to A\times B$, такой, что $\pi(\psi(x))=a$
и $\rho(\psi(x))=b$:
$$\xymatrix{ & A\times B
\ar[ld]_{\pi} \ar[rd]^{\rho} & \\
A & & B \\
&  C \ar[lu]^\alpha \ar[ru]_\beta \ar@{-->}[uu]^\psi &
}
$$
 Поскольку гомоморфизм $\psi$  однозначно определяется образом элемента $x$,
 существует единственный ненильпотентный однородный элемент $c:= \psi(x)\in A\times B$,
 такой, что $\pi(c)=a$, $\rho(c)=b$.

Отображая $x$ в элементы $a_i$, $b_k$, $a_i + a_j$, $b_k+b_\ell$,
получим единственные ненильпотентные однородные элементы $c_{ik}, d_{ijk\ell} \in A \times B$,
где $1\leqslant i,j \leqslant 2$, $1\leqslant k,\ell \leqslant 3$, $i\ne j$, $k\ne \ell$,
такие, что
$$\pi(c_{ij})=a_i,\text{ } \pi(d_{ijk\ell})=a_i+a_j,\text{ } \rho(c_{ij})=b_j, \text{ }\rho(d_{ijk\ell})=b_k+b_\ell.$$
Из того, что
$$\pi(c_{ik}+c_{j\ell})=\pi(d_{ijk\ell})=a_i+a_j,\text{ }
\rho(c_{ik}+c_{j\ell})=\rho(d_{ijk\ell})=b_k+b_\ell$$
и элементы $a_i+a_j$ и $b_k+b_\ell$ ненильпотентны,
следует, что элемент
$c_{ik}+c_{j\ell}$ также ненильпотентен и справедливо равенство $c_{ik}+c_{j\ell} = d_{ijk\ell}$.
Поскольку все элементы $c_{ik},c_{j\ell},d_{ijk\ell}$ ненулевые,
для любого набора индексов $(i,j,k,\ell)$ элементы $c_{ik},c_{j\ell},d_{ijk\ell}$
принадлежат одной и той же однородной компоненте алгебры $A\times B$. Меняя $i,j,k,\ell$,
получим, что элементы $c_{ik},d_{ijk\ell}$ для всех значений $(i,j,k,\ell)$
принадлежат одной и той же однородной компоненте алгебры $A\times B$.

В силу того, что
 $\pi(c_{11})=\pi(c_{12})=a_1$ и гомоморфизм
 $\pi$ градуированно инъективен, справедливо равенство $c_{11}=c_{12}$.
 Однако $\rho(c_{11})=b_1\ne b_2=\rho(c_{12})$, откуда получаем противоречие
 и произведение $A\times B$ не существует.
\end{proof}

Покажем, что копроизведения в $\widetilde{\mathbf{GrnuAlg}_\mathbbm{k}}$
и $\widetilde{\mathbf{GrAlg}_\mathbbm{k}}$ также не всегда существуют.

\begin{proposition} Пусть $G$ и $H$~"--- группы, а $\mathbbm{k}$~"--- поле.
Тогда копроизведение в категории $\widetilde{\mathbf{GrnuAlg}_\mathbbm{k}}$ групповых алгебр $\mathbbm{k}G$ и $\mathbbm{k}H$ (наделённых стандартными градуировками) не существует.
\end{proposition}
\begin{proof}
Пусть $A$~"--- копроизведение объектов $\mathbbm{k}G$ и $\mathbbm{k}H$ в $\widetilde{\mathbf{GrnuAlg}_\mathbbm{k}}$, а $i_1 \colon \mathbbm{k}G \to A$
и $i_2 \colon \mathbbm{k}H \to A$~"--- соответствующие морфизмы.
Пусть $\varphi_1 \colon \mathbbm{k}G \to \mathbbm{k}(G\times H)$ и
$\varphi_2 \colon \mathbbm{k}H \to \mathbbm{k}(G\times H)$~"--- естественные вложения.
Тогда существует единственный градуированно инъективный гомоморфизм $\varphi \colon A \to \mathbbm{k}(G\times H)$,
такой, что диаграмма ниже коммутирует:
$$\xymatrix{ & A
\ar@{-->}[dd]^\varphi & \\
\mathbbm{k}G \ar[ru]^{i_1}\ar[rd]^{\varphi_1} & & \mathbbm{k}H \ar[lu]_{i_2}\ar[ld]_{\varphi_2} \\
& \mathbbm{k}(G \times H)    &
}
$$

В частности, $$\varphi(i_1(g)i_2(h))={(g,1_H)}{(1_G,h)}= {(g,h)}\ne 0,$$
откуда $i_1(g)i_2(h)\ne 0$.

Пусть $\mathbbm{k}G \times \mathbbm{k}H := \lbrace (a, b) \mid a \in \mathbbm{k}G,\ b\in \mathbbm{k}H\rbrace$~"---
алгебра с покомпонентными операциями,
на которой задана $G\times H$-градуировка при помощи условий $(g, 0)\in (\mathbbm{k}G \times \mathbbm{k}H)^{\bigl((g, 1_H )\bigr)}$ и $(0, h)\in (\mathbbm{k}G \times \mathbbm{k}H)^{\bigl((1_G, h)\bigr)}$ для всех $g\in G$ и $h\in H$.

Пусть теперь $\psi_1 \colon \mathbbm{k}G \to \mathbbm{k}G \times \mathbbm{k}H$ и
$\psi_2 \colon \mathbbm{k}H \to \mathbbm{k}G\times \mathbbm{k}H$~"--- естественные вложения.
Тогда существует единственный градуированно инъективный гомоморфизм $\psi \colon A \to \mathbbm{k}G \times \mathbbm{k}H$, такой,
что диаграмма ниже коммутирует:
$$\xymatrix{ & A
\ar@{-->}[dd]^\psi & \\
\mathbbm{k}G \ar[ru]^{i_1}\ar[rd]^{\psi_1} & & \mathbbm{k}H \ar[lu]_{i_2}\ar[ld]_{\psi_2} \\
& \mathbbm{k}G \times \mathbbm{k}H    &
}
$$
В частности, $$\psi(i_1(g)i_2(h))=(g, 0)(0, h)= 0$$ и мы получаем противоречие, поскольку элемент $i_1(g)i_2(h)\ne 0$ однороден как произведение однородных элементов, а гомоморфизм $\psi$ градуированно инъективен.
\end{proof}

В то же время копроизведение алгебр $\mathbbm{k}G$ и $\mathbbm{k}H$ в категории $\widetilde{\mathbf{GrAlg}_\mathbbm{k}}$ равно $\mathbbm{k}(G*H)$,
групповой алгебре копроизведения групп $G$ и $H$ в категории групп.
В предложении~\ref{PropositionCoproductDoesntExistTilde}, которое доказывается ниже,
строится пример градуированных алгебр, у которых не существует копроизведения ни в $\widetilde{\mathbf{GrnuAlg}_\mathbbm{k}}$, ни в $\widetilde{\mathbf{GrAlg}_\mathbbm{k}}$.

\begin{proposition}\label{PropositionCoproductDoesntExistTilde} Пусть $\mathbbm{k}$~"--- поле,
а $A_i=\langle 1, a_i\rangle_\mathbbm{k}$, где $a_i^2=0$, $i=1,2$,~"--- две двумерные алгебры,
на которых заданы $\mathbb Z/2\mathbb Z$-градуировки при помощи условий $a_i \in A_i^{(\bar 1)}$.
Тогда копроизведение объектов $A_1$ и $A_2$ не существует ни в категории
$\widetilde{\mathbf{GrnuAlg}_\mathbbm{k}}$, ни в категории $\widetilde{\mathbf{GrAlg}_\mathbbm{k}}$.
\end{proposition}
\begin{proof} Предположим противное:
пусть существует копроизведение $A$ объектов $A_1$ и $A_2$, а $i_j \colon A_j \to A$, где $j=1,2$,~"--- соответствующие морфизмы.
Обозначим через $A_0=\langle 1, a_1,a_2\rangle_\mathbbm{k}$
трёхмерную градуированную алгебру с единицей $1$, порождающие которой удовлетворяют соотношениям $a_1^2=a_2^2=a_1a_2=a_2a_1 = 0$, а $\mathbb Z/3\mathbb Z$-градуировка задаётся при помощи условий $a_j \in A_0^{(\bar j)}$, $j=1,2$.
Обозначим через $\varphi_j \colon A_j \to A_0$ естественные вложения.

В силу определения копроизведения существуют единственные градуированно инъективные гомоморфизмы $\varphi \colon A \to A_0$, такие, что диаграмма ниже коммутирует:
$$\xymatrix{ & A
\ar@{-->}[dd]^\varphi & \\
A_1 \ar[ru]^{i_1}\ar[rd]^{\varphi_1} & & A_2 \ar[lu]_{i_2}\ar[ld]_{\varphi_2}
 \\
& A_0    &
}
$$

В частности, $\varphi(i_1(a_1)i_2(a_2))=a_1 a_2=0$.
Поскольку и $i_1(a_1)$, и $i_2(a_2)$
являются однородными элементами, а гомоморфизм $\varphi$ градуированно инъективен, получаем
\begin{equation}\label{eq:equal0}
i_1(a_1)i_2(a_2) = 0.
\end{equation}

Пусть теперь $B=\langle 1, b_1,b_2, b_1b_2\rangle_\mathbbm{k}$~"---
четырёхмерная градуированная алгебра с единицей $1$,
порождающие которой удовлетворяют соотношениям $b_1^2=b_2^2=b_2b_1 = 0$,
а $\mathbb Z/2\mathbb Z\times \mathbb Z/2\mathbb Z$-градуировка задаётся условиями $b_1 \in B^{(\bar 1,\bar 0)}$, $b_2 \in B^{(\bar 0,\bar 1)}$.
Пусть $\psi_j \colon A_j \to B$, где $j=1,2$,~"--- такие вложения, что $a_j \mapsto b_j$, $1\mapsto 1$.

В силу определения копроизведения существуют единственный градуированно инъективный гомоморфизм
$\psi \colon A \to B$, делающий диаграмму ниже коммутативной:
$$\xymatrix{ & A
\ar@{-->}[dd]^\psi & \\
A_1\ar[ru]^{i_1}\ar[rd]^{\psi_1} & & A_2 \ar[lu]_{i_2}\ar[ld]_{\psi_2} \\
& B & }$$
В частности, $\psi(i_1(a_1)i_2(a_2))=b_1 b_2 \ne 0$ и мы получаем противоречие с~\eqref{eq:equal0}.
Отсюда копроизведение объектов $A_1$ и $A_2$ не существует.
\end{proof}

Заметим, что отображение $i$ из предложения~\ref{PropositionEqualizerGrAlg}
инъективно и, следовательно, градуированно инъективно.
Отсюда вложения $\widetilde{\mathbf{GrnuAlg}_\mathbbm{k}} \subset \mathbf{GrnuAlg}_\mathbbm{k}$ и $\widetilde{\mathbf{GrAlg}_\mathbbm{k}} \subset \mathbf{GrAlg}_\mathbbm{k}$ сохраняют уравнители.

\begin{proposition}
Пусть $f,g \colon A \to B$~"--- различные морфизмы в категории $\widetilde{\mathbf{GrnuAlg}_\mathbbm{k}}$ или $\widetilde{\mathbf{GrAlg}_\mathbbm{k}}$, а алгебра $B$ наделена тривиальной градуировкой.
Тогда коуравнитель морфизмов $f$ и $g$ не существует.
\end{proposition}
\begin{proof}
Предположим, что $h\colon B \to C$~"--- коуравнитель морфизмов $f$ и $g$. Тогда $hf = hg$. 
В силу того, что гомоморфизм $h$ градуированно инъективен, а градуировка на $B$ тривиальна, гомоморфизм $h$
инъективен, откуда $f=g$. Получаем противоречие.
\end{proof}

Покажем теперь, что мономорфизмы в категориях $\widetilde{\mathbf{GrnuAlg}_\mathbbm{k}}$ и $\widetilde{\mathbf{GrAlg}_\mathbbm{k}}$
допускают то же описание, что и в категориях $\mathbf{GrnuAlg}_\mathbbm{k}$ and $\mathbf{GrAlg}_\mathbbm{k}$.

\begin{proposition}\label{PropositionCriterionForMonomorphismTildeGrAlg}
Пусть $f\colon A=\bigoplus_{g\in G} A^{(g)} \to B$~"--- морфизм
в категории $\widetilde{\mathbf{GrnuAlg}_\mathbbm{k}}$ или $\widetilde{\mathbf{GrAlg}_\mathbbm{k}}$.
Тогда $f$ является мономорфизмом, если и только если $f(a)\ne f(b)$ для всех $a,b\in \bigcup_{g\in G}
 A^{(g)}$, таких, что $a \ne b$.
\end{proposition}
\begin{proof}
Предположим, $f$~"--- мономорфизм.
Выберем такие $a,b\in \bigcup_{g\in G} A^{(g)}$, что $f(a)= f(b)$.

Рассмотрим сперва случай, когда $f(a)^k \ne 0$ для всех $k\in\mathbb N$.
Тогда $a^k\ne 0$, $b^k \ne 0$ для всех $k\in\mathbb N$.
Снова обозначим через $C$ 
алгебру многочленов от переменной $x$ с коэффициентами из $\mathbbm{k}$ без свободного члена, если рассматривается категория $\widetilde{\mathbf{GrnuAlg}_\mathbbm{k}}$, и со свободным членом, если рассматривается категория $\widetilde{\mathbf{GrAlg}_\mathbbm{k}}$,
и введём на $C$ градуировку группой $\mathbb Z$ по степени одночленов.
Существуют единственные градуированно инъективные гомоморфизмы
$\alpha,\beta \colon C \to A$,
такие, что $\alpha(x)=a$, $\beta(x)=b$. При этом $f\alpha=f\beta$,
так как $f(\alpha(x)) = f(a)=f(b)=f(\beta(x))$. Поскольку $f$~"--- мономорфизм,
отсюда следует, что $\alpha = \beta$ и $a=\alpha(x)=\beta(x)= b$.

Рассмотрим теперь случай, когда
элемент $f(a)$ нильпотентен.
Определим число $k\in \mathbb N$ при помощи условий $f(a)^k=0$, $f(a)^{k-1}\ne 0$. (Если $f(a)=0$, положим $k:=1$.)
Поскольку $f$ градуированно инъективен, а все элементы $a^i$ и $b^j$ однородны,
справедливо равенство $a^k=b^k=0$ и неравенства $a^{k-1}\ne 0$ и $b^{k-1}\ne 0$.
Если $k=1$, то $a=b=0$ и необходимость доказана.
Предположим, что $k > 1$. Пусть $\bar C = C/(x^k)$. Обозначим через $\bar x$ образ элемента $x$ в алгебре $\bar C$.
Тогда существуют единственные градуированно инъективные гомоморфизмы $\alpha,\beta \colon \bar C \to A$,
такие, что $\alpha(\bar x)=a$, $\beta(\bar x)=b$. Следовательно, $f\alpha=f\beta$,
так как $f(\alpha(\bar x)) = f(a)=f(b)=f(\beta(\bar x))$. Поскольку $f$~"--- мономорфизм,
справедливо равенство $\alpha = \beta$ и $a=\alpha(\bar x)=\beta(\bar x)= b$. Теперь необходимость доказана полностью.

Достаточность следует из предложения~\ref{PropositionCriterionForMonomorphismGrAlg},
поскольку в случае выполнения условия из формулировки предложения~\ref{PropositionCriterionForMonomorphismTildeGrAlg} отображение $f$ 
является мономорфизмом в б\'ольшей категории $\mathbf{GrnuAlg}_\mathbbm{k}$ (соответственно, $\mathbf{GrAlg}_\mathbbm{k}$).
\end{proof}

\section{Отсутствие сопряжённых у функторов универсальной группы градуировки}\label{SectionAbsenceGrAdjunctions}

В параграфе~\S\ref{SectionGrAdjoint} были введены две пары сопряжённых функторов, связанных с градуировками. 
В этом параграфе мы покажем, что функторы
$R \colon \widetilde{\mathbf{GrnuAlg}_\mathbbm{k}} \to \mathbf{Grp}$ и $R_1\colon \widetilde{\mathbf{GrAlg}_\mathbbm{k}}\to \mathbf{Grp}$, которые были введены в~\S\ref{SectionUniversalGradingGroupFunctors}
и ставят в соответствие каждой градуировке её универсальную группу,
не имеют ни левых, ни правых сопряжённых.

 \begin{theorem}\label{TheoremAbsenceOfLeftAdjointUGGF}
Функторы универсальной группы градуировки
$$R\colon \widetilde{\mathbf{GrnuAlg}_\mathbbm{k}} \rightarrow \mathbf{Grp} \text{ и }
R_1\colon \widetilde{\mathbf{GrAlg}_\mathbbm{k}}\rightarrow \mathbf{Grp}$$
не имеют ни левых, ни правых сопряжённых.
\end{theorem}

Разобьём доказательство теоремы~\ref{TheoremAbsenceOfLeftAdjointUGGF} на три предложения.
\begin{proposition}\label{PropositionAbsenceOfLeftAdjointUGGF}
Функтор $R$ не имеет левого сопряжённого.
\end{proposition}
\begin{proof} Предположим, что $K$~"--- левый сопряжённый для функтора $R$.
Тогда существует естественная биекция $$\widetilde{\mathbf{GrnuAlg}_\mathbbm{k}}(K(H), \Gamma) \to \mathbf{Grp}(H, R(\Gamma)).$$

Докажем, что для любой группы $H$ градуировка $K(H)$
является градуировкой на нулевой алгебре.
Пусть $K(H) \colon A=\bigoplus_{g\in G} A^{(g)}$ и $A\ne 0$.
Выберем произвольное множество индексов $\Lambda$, такое, что $|\Lambda| > |\Hom(H,G_{K(H)})|$.
Рассмотрим прямую сумму $\bigoplus_{\lambda \in \Lambda} A$
копий алгебры $A$, где каждая копия наделена $G$-градуировкой $K(H)$.
Обозначим через $\Xi$ полученную таким образом $G$-градуировку на алгебре $\bigoplus_{\lambda \in \Lambda} A$. Тогда $G_{\Xi} \cong G_{K(H)}$, однако $$|\widetilde{\mathbf{GrnuAlg}_\mathbbm{k}}(K(H), \Xi)|
\geqslant |\Lambda| > |\Hom(H,G_{K(H)})|=|\mathbf{Grp}(H, R(\Xi))|,$$
что противоречит существованию естественной биекции.
Следовательно, для любой группы
 $H$ градуировка $K(H)$ 
является градуировкой на нулевой алгебре.
В частности, любое множество $\widetilde{\mathbf{GrnuAlg}_\mathbbm{k}}(K(H), \Gamma)$
содержит в точности один элемент.

Если $H$~"--- нетривиальная группа, а $\mathbbm{k}H$~"--- её групповая алгебра, наделённая стандартной
градуировкой $\Gamma$, то множество $\mathbf{Grp}(H, R(\Gamma))=\Hom(H,H)$
содержит больше одного элемента (по крайней мере, тождественное отображение
и гомоморфизм, отображающий всю группу в $1_H$), и мы снова приходим к противоречию.
Отсюда левый сопряжённый для функтора 
$R$ не существует.
\end{proof}

Приём с бесконечной суммой работает только в случае, когда разрешено, чтобы алгебра, которая получается в результате, была без единицы. С другой стороны, для алгебры $A$ с единицей $1_A$ можно использовать существование гомоморфного вложения $\mathbbm{k} \cdot 1_A \to A$.

\begin{proposition}\label{PropositionAbsenceOfLeftAdjointUGGF1}
Функтор $R_1$ не имеет левого сопряжённого.
\end{proposition}
\begin{proof} Предположим, что для функтора $R_1$ существует левый сопряжённый~$K$.
Тогда существует естественная биекция $$\widetilde{\mathbf{GrAlg}_\mathbbm{k}}(K(H), \Gamma) \to \mathbf{Grp}(H, R_1(\Gamma)).$$

Докажем, что для любой группы $H$ градуировка $K(H)$
является градуировкой на алгебре, изоморфной основному полю $\mathbbm{k}$.
Пусть $H$~"--- группа, а $K(H) \colon A=\bigoplus_{g\in G} A^{(g)}$.
Обозначим через $\Upsilon$ градуировку на $\mathbbm{k}$ тривиальной группой. Тогда множество $\mathbf{Grp}(H, R_1(\Upsilon))$, а следовательно, и $\widetilde{\mathbf{GrAlg}_\mathbbm{k}}(K(H), \Upsilon)$ состоит из одного элемента.
В частности, существует градуированно инъективный гомоморфизм
  $\varphi \colon A \to \mathbbm{k}$ алгебр с единицей.
Следовательно, существует идеал $\ker \varphi \subsetneqq A$ коразмерности $1$.
Обозначим через $\Xi$ градуировку на алгебре $A$ тривиальной группой.
  Если $\ker \varphi\ne 0$, то множество $\widetilde{\mathbf{GrAlg}_\mathbbm{k}}(K(H), \Xi)$
  состоит по крайней мере из двух различных элементов:
  \begin{enumerate}
    \item тождественного отображения $A\to A$;
    \item композиции градуированно инъективного гомоморфизма $\varphi$
    и вложения $\mathbbm{k} \cdot 1_A \to A$.
  \end{enumerate}
   Поскольку группа $R_1(\Xi)$ тривиальна и множество $\mathbf{Grp}(H, R_1(\Xi))$
   состоит из единственного элемента, множество $\widetilde{\mathbf{GrAlg}_\mathbbm{k}}(K(H), \Xi)$
также состоит из единственного элемента,
    $\ker \varphi=0$ и $A \cong \mathbbm{k}$.
   
   В частности, $\widetilde{\mathbf{GrAlg}_\mathbbm{k}}(K(H), \Gamma)$ 
содержит в точности один элемент для всех $H$ и $\Gamma$.
  Однако если $\Gamma$~"--- стандартная градуировка
  на групповой алгебре $\mathbbm{k}H$
  произвольной нетривиальной группы $H$, то множество $\mathbf{Grp}(H, R_1(\Gamma))=\Hom(H,H)$ 
содержит больше одного элемента. Получаем противоречие, откуда
левого сопряжённого функтора  для функтора $R_1$  не существует.
\end{proof}

Доказательства отсутствия правого сопряжённого функтора для функторов $R$ и $R_1$ совпадают:

\begin{proposition}\label{PropositionAbsenceOfRightAdjointsUGGF}
У функторов $R$ и $R_1$ отсутствуют правые сопряжённые.
\end{proposition}
\begin{proof}
Для определённости рассмотрим случай функтора $R$.

Пусть $K$~"--- правый сопряжённый для функтора $R$. Тогда
существует естественная биекция $$\mathbf{Grp}(R(\Gamma), H) \to \widetilde{\mathbf{GrnuAlg}_\mathbbm{k}}(\Gamma, K(H)).$$

Заметим, что множество $\mathbf{Grp}(R(\Gamma), H)$ никогда не пусто,
поскольку, по крайней мере, содержит гомоморфизм, который отображает все элементы группы $H$ в $1_H$. 
Фиксируем группу~$H$. Пусть $K(H) \colon A=\bigoplus_{g\in G} A^{(g)}$.
Возьмём $\mathbbm{k}$-алгебру $B$ мощности $|B|$ большей, чем
мощность $|A|$ алгебры $A$. Например,
можно взять $B = \End_\mathbbm{k}(V)$, где $V$~"--- векторное пространство
с базисом, мощность которого больше $|A|$. Обозначим через $\Gamma$ градуировку на $B$
тривиальной группой. Тогда ни для какого $g\in G$ не существует инъективных
отображений $B \to A^{(g)}$, откуда множество $\widetilde{\mathbf{GrnuAlg}_\mathbbm{k}}(\Gamma, K(H))$
пусто. Получаем противоречие. Следовательно, правого сопряжённого функтора  для функтора $R$  не существует.
\end{proof}

Теорема~\ref{TheoremAbsenceOfLeftAdjointUGGF} является следствием 
предложений~\ref{PropositionAbsenceOfLeftAdjointUGGF}--\ref{PropositionAbsenceOfRightAdjointsUGGF}.
Из доказательств видно, что для того, чтобы получить пару сопряжённых функторов, необходимо ограничить категорию алгебр до таких алгебр, которые полностью определяются своими градуирующими группами.

\section{Сопряжение в случае групповых алгебр}\label{SectionFunctorGroupAlgebra}

Пусть $\mathbbm{k}$~"--- некоторое поле.  Обозначим через $\mathbf{Grp}'_\mathbbm{k}$
категорию, объектами которой являются все группы $G$,
которые не имеют нетривиальных одномерных представлений (другими словами, $H^1(G, \mathbbm{k}^{\times}) = 0$),
а морфизмами являются всевозможные гомоморфизмы групп.
Пусть $\mathbf{GrpAlg}'_\mathbbm{k}$~"--- категория,
объектами которой являются групповые алгебры $\mathbbm{k}G$ групп $G$ из категории $\mathbf{Grp}'_\mathbbm{k}$
со стандартными градуировками, а морфизмами являются все ненулевые градуированные
гомоморфизмы.
Обозначим через $U$ функтор $\mathbf{GrpAlg}'_\mathbbm{k} \to \mathbf{Grp}'_\mathbbm{k}$
заданный равенством $U(\mathbbm{k}G)=G$ и условием $\varphi\left(\mathbbm{k}G_1^{(g)}\right)\subseteq \mathbbm{k}G_2^{\left(U(\varphi)(g) \right)}$
 для всех $\varphi \colon \mathbbm{k}G_1 \to \mathbbm{k}G_2$.
Обозначим через $\mathbbm{k}-$ функтор, который ставит в соответствие каждой группе её групповую алгебру над $\mathbbm{k}$.

\begin{proposition}\label{PropositionEquivGroupAlgGroups}
Для любого объекта $G$ категории $\mathbf{Grp}'_\mathbbm{k}$
и любого объекта $A$ категории $\mathbf{GrpAlg}'_\mathbbm{k}$ существует биекция $\theta_{G, A} \colon \mathbf{GrpAlg}'_\mathbbm{k}(\mathbbm{k}G, A)
\to \mathbf{Grp}'_\mathbbm{k}(G, U(A))$, естественная по $A$ и $G$. Более того, $\mathbbm{k}U(-)= 1_{\mathbf{GrpAlg}'_\mathbbm{k}}$
и $U(\mathbbm{k}-)= 1_{\mathbf{Grp}'_\mathbbm{k}}$, т.е. категории
$\mathbf{Grp}'_\mathbbm{k}$ и $\mathbf{GrpAlg}'_\mathbbm{k}$ изоморфны.
\end{proposition}
\begin{proof}
Пусть $\varphi\colon \mathbbm{k}G \to \mathbbm{k}H$~"--- ненулевой градуированный гомоморфизм.
Тогда из того, что $\varphi({g_0}) \ne 0$ для некоторого $g_0\in G$,
следует, что
$\varphi({g_0})=\varphi({g_0} g^{-1}) \varphi(g)\ne 0$
и $\varphi(g)\ne 0$ для всех $g\in G$. Отсюда гомоморфизм $\varphi$ градуированно инъективен.
Следовательно, $\varphi$ однозначно определяется групповыми гомоморфизмами $\psi \colon G \to H$ и $\alpha \colon G \to \mathbbm{k}^\times$, такими, что $\varphi(g) = \alpha(g){\psi(g)}$ для всех $g\in G$.
Поскольку у $G$ нет нетривиальных одномерных представлений,
гомоморфизм $\alpha$ тривиален и существует естественная биекция $\theta_{G, A}$. Равенства $\mathbbm{k}U(-)= 1_{\mathbf{GrpAlg}'_\mathbbm{k}}$
и $U(\mathbbm{k}-)= 1_{\mathbf{Grp}'_\mathbbm{k}}$ проверяются непосредственно.
\end{proof}

         \section{Критерий эквивалентности в терминах линейных операторов}
         
         В теоремах~\ref{TheoremGradEquivCriterion} и~\ref{TheoremGradFinerCoarserCriterion}, которые доказываются ниже, мы формулируем критерии эквивалентности градуировок и наличия отношения <<тоньше-грубее>> между ними в терминах линейных операторов. Эти критерии затем  будут использованы в \S\ref{SectionGroupActions} 
и главе~\ref{ChapterOmegaAlg}        
         для того, чтобы ввести понятие эквивалентности 
и отношение <<тоньше-грубее>> для действий групп и (ко)модульных структур.
         
Пусть $G$~"--- группа, а $\mathbbm{k}$~"--- поле. Тогда алгебра $(\mathbbm{k}G)^*$, двойственная к коалгебре $\mathbbm{k}G$, изоморфна алгебре всех функций $G \to \mathbbm{k}$ с поточечными операциями.
Если алгебра $A = \bigoplus_{g\in G} A^{(g)}$ над полем $\mathbbm{k}$ градуирована группой $G$, тогда
алгебра $(\mathbbm{k}G)^*$ действует на алгебре $A$ следующим образом: $ha=h(g)a$ для всех $h\in (\mathbbm{k}G)^*$, $a\in A^{(g)}$, $g\in G$. 

Обе теоремы данного параграфа являются следствиями следующего утверждения:

\begin{lemma}\label{LemmaGradEquivCriterion}
Пусть $\Gamma_1 \colon A_1 = \bigoplus_{g\in G_1} A_1^{(g)}$ и 
$\Gamma_2 \colon A_2 = \bigoplus_{g\in G_2} A_2^{(g)}$~"--- групповые градуировки,
а $\varphi \colon A_1 \mathrel{\widetilde\to} A_2$~"--- изоморфизм алгебр.
Обозначим через $\zeta_i \colon (\mathbbm{k}G_i)^* \to \End_{\mathbbm{k}}(A_i)$, где $i=1,2$, гомоморфизм из алгебры $(\mathbbm{k}G_i)^*$ в алгебру $\End_{\mathbbm{k}}(A_i)$  линейных операторов на $A_i$, соответствующий $(\mathbbm{k}G_i)^*$-действию, определённому выше, а через $\tilde\varphi$~"---  изоморфизм алгебр  $\End_{\mathbbm{k}}(A_1) \mathrel{\widetilde\to} \End_{\mathbbm{k}}(A_2)$,
заданный при помощи равенства $\tilde\varphi(\psi)(a)=\varphi\Bigl(\psi\bigl(\varphi^{-1}(a)\bigr)\Bigr)$ для всех $\psi\in \End_{\mathbbm{k}}(A_1)$ и $a\in A_2$.
Тогда включение
\begin{equation}\label{EqEmbeddingOfImagesOfFG}
\tilde\varphi\Bigl(\zeta_1\bigl((\mathbbm{k}G_1)^*\bigr)\Bigr)\supseteq \zeta_2\bigl((\mathbbm{k}G_2)^*\bigr),
\end{equation}
справедливо, если и только если
для любого $g_1 \in G_1$ существует такое $g_2\in G_2$,
что $\varphi\left(A_1^{(g_1)}\right)\subseteq A_2^{(g_2)}$.
\end{lemma}
\begin{proof}
Если
для любого $g_1 \in G_1$ существует такое $g_2\in G_2$,
что $\varphi\left(A_1^{(g_1)}\right)\subseteq A_2^{(g_2)}$,
то всякое подпространство
$A_2^{(g_2)}$ является прямой суммой некоторых из подпространств
$\varphi\left(A_1^{(g_1)}\right)$, поскольку
$$\bigoplus_{g_1 \in G_1} \varphi\left(A_1^{(g_1)}\right)=\varphi(A_1)=A_2=\bigoplus_{g_2 \in G_2} A_2^{(g_2)}.$$
Заметим, что множество $\zeta_i\bigl((\mathbbm{k}G_i)^*\bigr)$ 
состоит из всех линейных операторов, которые действуют скалярно на каждой
однородной компоненте $A_i^{(g_i)}$, где $g_i \in G_i$.
 Следовательно, множество 
$\tilde\varphi\Bigl(\zeta_1\bigl((\mathbbm{k}G_1)^*\bigr)\Bigr)$ состоит из всех линейных операторов
на $A_2$, которые действуют скалярно на каждом
подпространстве $\varphi\Bigl(\zeta_1\left(A_1^{(g_1)}\right)\Bigr)$.
Поскольку  всякое подпространство
$A_2^{(g_2)}$ является прямой суммой некоторых из подпространств
$\varphi\left(A_1^{(g_1)}\right)$, все операторы из множества $\zeta_2 \bigl((\mathbbm{k}G_2)^*\bigr)$
также действуют скалярно на каждом
подпространстве $\varphi\left(A_1^{(g_1)}\right)$,
откуда и следует включение~(\ref{EqEmbeddingOfImagesOfFG}).

Обратно, пусть справедливо включение~(\ref{EqEmbeddingOfImagesOfFG}).
Обозначим через $p_{g_2}$, где $g_2 \in \supp \Gamma_2$, проекцию на подпространство $A_2^{(g_2)}$
вдоль подпространства $\bigoplus\limits_{\substack{h\in \supp\Gamma_2, \\ h \ne g_2}} A_2^{(h)}$,
т.е. $p_{g_2} a := a$ для всех $a\in A_2^{(g_2)}$
и $p_{g_2} a := 0$ для всех $a \in \bigoplus\limits_{\substack{h\in \supp\Gamma_2, \\ h \ne g_2}} A_2^{(h)}$. Тогда из~(\ref{EqEmbeddingOfImagesOfFG}) следует, что $p_{g_2}\in \tilde\varphi\Bigl(\zeta_1\bigl((\mathbbm{k}G_1)^*\bigr)\Bigr)$
для всех $g_2 \in \supp \Gamma_2$. В частности, оператор $p_{g_2}$ действует скалярно на всех компонентах $\varphi\left(A_1^{(g)}\right)$,
где $g\in \supp\Gamma_1$. Поскольку $p_{g_2}^2 = p_{g_2}$,
то для всякого $g\in \supp\Gamma_1$ либо $p_{g_2} \varphi\left(A_1^{(g)}\right) = 0$,
либо $p_{g_2} a = a$
для всех $a\in \varphi\left(A_1^{(g)}\right)$.
Следовательно, подпространство $A_2^{(g_2)}=\im(p_{g_2})$ является прямой суммой некоторых
компонент $\varphi\left(A_1^{(g)}\right)$ для некоторых
$g\in\supp \Gamma_1$. В силу того, что $$
\bigoplus\limits_{g_1\in \supp\Gamma_1} \varphi\left(A_1^{(g_1)}\right) =
\bigoplus\limits_{g_2\in \supp\Gamma_2} A_2^{(g_2)},$$
для любого $g_1 \in G_1$ существует такое $g_2\in G_2$,
что $\varphi\left(A_1^{(g_1)}\right)\subseteq A_2^{(g_2)}$.
 \end{proof}

Из леммы~\ref{LemmaGradEquivCriterion} 
 вытекают следующие утверждения:

\begin{theorem}\label{TheoremGradEquivCriterion}
Пусть $\Gamma_1 \colon A_1 = \bigoplus_{g\in G_1} A_1^{(g)}$ и 
$\Gamma_2 \colon A_2 = \bigoplus_{g\in G_2} A_2^{(g)}$~"--- групповые градуировки,
а $\varphi \colon A_1 \mathrel{\widetilde\to} A_2$~"--- изоморфизм алгебр.
Тогда $\varphi$~"--- эквивалентность градуировок, если и только если
$$
\tilde\varphi\Bigl(\zeta_1\bigl((\mathbbm{k}G_1)^*\bigr)\Bigr)=\zeta_2\bigl((\mathbbm{k}G_2)^*\bigr),
$$
где $\zeta_i \colon (\mathbbm{k}G_i)^* \to \End_{\mathbbm{k}}(A_i)$~"--- гомоморфизм из алгебры $(\mathbbm{k}G_i)^*$ в алгебру $\End_{\mathbbm{k}}(A_i)$  линейных операторов на $A_i$, соответствующий $(\mathbbm{k}G_i)^*$-действию, определённому выше, $i=1,2$,
а изоморфизм алгебр $\tilde\varphi \colon \End_{\mathbbm{k}}(A_1) \mathrel{\widetilde\to} \End_{\mathbbm{k}}(A_2)$
задан при помощи равенства $\tilde\varphi(\psi)(a)=\varphi\Bigl(\psi\bigl(\varphi^{-1}(a)\bigr)\Bigr)$ для $\psi\in \End_{\mathbbm{k}}(A_1)$ и $a\in A_2$.
\end{theorem}
\begin{proof}
Достаточно применить лемму~\ref{LemmaGradEquivCriterion} к $\varphi$ и $\varphi^{-1}$.
\end{proof}
\begin{theorem}\label{TheoremGradFinerCoarserCriterion}
Пусть $\Gamma_1 \colon A = \bigoplus_{g\in G_1} A^{(g)}$ и 
$\Gamma_2 \colon A = \bigoplus_{g\in G_2} A^{(g)}$~"--- градуировки на одной и той же алгебре~$A$
и пусть $\zeta_i \colon (\mathbbm{k}G_i)^* \to \End_{\mathbbm{k}}(A)$~"--- соответствующие этим градуировкам гомоморфизмы.
Тогда $\Gamma_1$ тоньше, чем $\Gamma_2$, если и только если
$\zeta_1\bigl((\mathbbm{k}G_1)^*\bigr) \supseteq \zeta_2\bigl((\mathbbm{k}G_2)^*\bigr)$.
\end{theorem}
\begin{proof}
Достаточно применить лемму~\ref{LemmaGradEquivCriterion} к $A_1=A_2=A$ и $\varphi = \id_A$.
\end{proof}

 \section{Эквивалентность действий групп}\label{SectionGroupActions}

Предположим, что мы находимся в двойственной ситуации: 
группы $G_i$, где $i=1,2$, действуют на алгебрах $A_i$ автоморфизмами.
Пусть $\zeta_i \colon G_i \to \Aut(A_i)$~"--- соответствующие гомоморфизмы групп.
Отталкиваясь от теоремы~\ref{TheoremGradEquivCriterion}, дадим следующее определение:

\begin{definition}\label{DefGroupActionEquivalence}
Пусть $\varphi \colon A_1 \mathrel{\widetilde\to}A_2$~"--- изоморфизм алгебр.
Будем говорить, что действия $\zeta_1$ и $\zeta_2$  \textit{эквивалентны при помощи изоморфизма $\varphi$},
если \begin{equation*}
\tilde\varphi\Bigl(\langle \zeta_1\bigl(G_1 \bigr) \rangle_\mathbbm{k}\Bigr)=\langle \zeta_2\bigl(G_2 \bigr) \rangle_\mathbbm{k},
\end{equation*}
где $\mathbbm{k}$-линейная оболочка $\langle \cdot \rangle_\mathbbm{k}$ берётся в соответствующей алгебре $\End(A_i)$ 
и изоморфизм $\tilde\varphi \colon \End_{\mathbbm{k}}(A_1) \mathrel{\widetilde\to} \End_{\mathbbm{k}}(A_2)$ задаётся равенством $\tilde\varphi(\psi)(a)=\varphi\Bigl(\psi\bigl(\varphi^{-1}(a)\bigr)\Bigr)$ для $\psi\in \End_{\mathbbm{k}}(A_1)$ и $a\in A_2$.
\end{definition}

Если $\Gamma \colon A=\bigoplus_{g\in G} A^{(g)}$~"--- градуировка алгебры $A$ группой $G$, 
существует стандартное $\Hom(G,\mathbbm{k}^\times)$-действие на алгебре $A$ автоморфизмами: \begin{equation}\label{EqGroupActionChi}\chi a := \chi(g)a\text{ для всех }\chi \in \Hom(G,\mathbbm{k}^\times),\ a\in A^{(g)},\ g\in G.\end{equation}
Каждое отображение $\chi \in \Hom(G,\mathbbm{k}^\times)$ продолжается по линейности до отображения
$\mathbbm{k}G \to \mathbbm{k}$, что задаёт вложение $\Hom(G,\mathbbm{k}^\times) \hookrightarrow (\mathbbm{k}G)^*$. Отсюда $\Hom(G,\mathbbm{k}^\times)$-действие~\eqref{EqGroupActionChi}
является ограничением отображения $\zeta \colon (\mathbbm{k}G)^* \to \End_\mathbbm{k}(A)$, соответствующего градуировке~$\Gamma$.

В случае, когда поле $\mathbbm{k}$ алгебраически замкнуто характеристики $0$,
а $G$~"--- конечная абелева группа, группа $\Hom(G,\mathbbm{k}^\times)$ 
обычно обозначается  через $\widehat G$
и называется \textit{группой линейных характеров группы $G$}.
Из классической теоремы о структуре конечно порождённых абелевых
групп следует, что $\widehat G \cong G$. Более того, равенство (\ref{EqGroupActionChi}) 
задаёт взаимно однозначное соответствие между $G$-градуировками и $\widehat G$-действиями. (См., например, \cite[\S 3.2]{ZaiGia}.)
 Докажем, что эквивалентные $G$-градуировки соответствуют эквивалентным $\widehat G$-действиям.

\begin{proposition}
Пусть $\Gamma_1 \colon A = \bigoplus_{g\in G_1} A^{(g)}$ и 
$\Gamma_2 \colon A = \bigoplus_{g\in G_2} A^{(g)}$~"--- градуировки
конечными абелевыми группами $G_1$ и $G_2$,
а $\zeta_i$~"--- соответствующие 
$(\mathbbm{k}G_i)^*$-действия, $i=1,2$.
 Предположим, что основное поле $\mathbbm{k}$ алгебраически замкнуто характеристики $0$. 
Тогда $\Gamma_1$ и $\Gamma_2$ эквивалентны как групповые градуировки, если и только если $\zeta_1 \bigr|_{\widehat G_1}$
и $\zeta_2 \bigr|_{\widehat G_2}$ эквивалентны как действия групп.
\end{proposition}
\begin{proof} В силу соотношений ортогональности для характеров,
 элементы группы $\widehat G_i$ образуют базис в алгебре $(\mathbbm{k}G_i)^*$. Отсюда $\zeta_i((\mathbbm{k}G_i)^*)= \langle\zeta_i(\widehat G_i )\rangle$ при $i=1,2$,
и предложение является следствием теоремы~\ref{TheoremGradEquivCriterion}.
\end{proof}

Вернёмся к случаю произвольных групп $G_1$ и $G_2$ и произвольного поля $\mathbbm{k}$.
Как и в случае градуировок, можно отождествить $A_1$ и $A_2$ при помощи изоморфизма $\varphi$. 
Тогда эквивалентность действий $\zeta_1$ и $\zeta_2$ означает, что
образы групп $G_1$ и $G_2$ порождают одну и ту же подалгебру в алгебре $\mathbbm{k}$-линейных операторов
на $A_1=A_2$.

\begin{definition}\label{UniversalGroupOfTheAction}
Пусть $\zeta \colon G \to \Aut(A)$~"--- действие группы $G$
на алгебре $A$ и пусть $\varkappa_{\zeta} \colon G_{\zeta} \to \Aut(A)$~"--- действие группы $G_\zeta$, эквивалентное $\zeta$ при помощи тождественного изоморфизма $\id_A$.  Будем говорить, что пара $(G_{\zeta}, \varkappa_{\zeta})$~"--- \textit{универсальная группа действия $\zeta$},
если для любого другого действия $\zeta_1 \colon G_1 \to \Aut(A)$, эквивалентного действию $\zeta$ 
при помощи $\id_A$, существует единственный гомоморфизм групп $\varphi \colon G_1 \to G_{\zeta}$,
такой, что следующая диаграмма коммутативна:
$$\xymatrix{ \Aut(A)  & G_{\zeta} \ar[l]_(0.3){\varkappa_{\zeta}}  \\
& G_1 \ar[lu]^{\zeta_1} \ar@{-->}[u]_\varphi
}
$$
\end{definition}

\begin{remark}\label{RemarkUniversalGroupOfAnAction}
Рассмотрим группу $\mathcal U(\langle \zeta(G) \rangle_\mathbbm{k})
\cap \Aut(A)$, где $\mathcal U(\langle \zeta(G) \rangle_\mathbbm{k})$~"---
группа обратимых элементов алгебры $\langle \zeta(G) \rangle_\mathbbm{k} \subseteq \End_\mathbbm{k}(A)$.
В силу того, что для любой группы $G_1$ из определения~\ref{UniversalGroupOfTheAction} образ $G_1$ в $\Aut(A)$ 
является подмножеством  группы $\mathcal U(\langle \zeta_1(G_1) \rangle_\mathbbm{k}) \cap \Aut(A) = \mathcal U(\langle \zeta(G) \rangle_\mathbbm{k})\cap \Aut(A)$,
 универсальная группа действия $\zeta$~"--- это
 (с точностью до изоморфизма) пара $(G_{\zeta}, \varkappa_{\zeta})$, где $$G_{\zeta}:= \mathcal U(\langle \zeta(G) \rangle_\mathbbm{k})
\cap \Aut(A)$$
и $\varkappa_{\zeta}$~"--- естественное вложение $G_{\zeta} \subseteq \Aut(A)$.
\end{remark}

Пусть $\zeta_i \colon G_i \to \Aut(A)$, $i=1,2$, действия групп на алгебре $A$ автоморфизмами.
Будем говорить, что $\zeta_1$ \textit{тоньше}, чем $\zeta_2$, а $\zeta_2$ \textit{грубее}, чем $\zeta_1$,
если $\langle \zeta_2(G) \rangle_\mathbbm{k} \subseteq \langle \zeta_1(G) \rangle_\mathbbm{k}$.
Как и в случае градуировок, легко видеть, что это отношение является препорядком,
а $\zeta_1$ одновременно и тоньше, и грубее, чем  $\zeta_2$, если и только
если $\id_A$~"--- эквивалентность действий $\zeta_1$ и $\zeta_2$. 
Более того, универсальная группа действия является функтором из этого предпорядка в категорию
групп:
если $\zeta_1$ грубее, чем $\zeta_2$, функтор ставит
в соответствие  этой стрелке вложение
$\mathcal U\bigl(\langle \zeta(G_2) \rangle_\mathbbm{k}\bigr)
\cap \Aut(A) \subseteq \mathcal U\bigl(\langle \zeta(G_1) \rangle_\mathbbm{k}\bigr)
\cap \Aut(A)$.

\begin{proposition}
Пусть $\Gamma_1 \colon A = \bigoplus_{g\in G_1} A^{(g)}$ и 
$\Gamma_2 \colon A = \bigoplus_{g\in G_2} A^{(g)}$"--- градуировки на алгебре $A$
конечными абелевыми группами $G_1$ и $G_2$,
а $\zeta_i$~"--- соответствующие 
$(\mathbbm{k}G_i)^*$-действия, $i=1,2$.
 Предположим, что основное поле $\mathbbm{k}$ алгебраически замкнуто характеристики $0$. 
Тогда $\Gamma_1$ тоньше, чем $\Gamma_2$, если и только если $\zeta_1 \bigr|_{\widehat G_1}$
тоньше, чем $\zeta_2 \bigr|_{\widehat G_2}$.
\end{proposition}
\begin{proof}
Опять воспользуемся тем фактом, что, в силу соотношений ортогональности для характеров,
 элементы группы $\widehat G_i$ образуют базис в алгебре $(\mathbbm{k}G_i)^*$. Отсюда $\zeta_i((\mathbbm{k}G_i)^*)= \langle\zeta_i(\widehat G_i )\rangle$ при $i=1,2$,
и предложение является следствием теоремы~\ref{TheoremGradFinerCoarserCriterion}.
\end{proof}

\newpage

\chapter{Ассоциативные (ко)модульные алгебры}\label{ChapterH(co)modAssoc}

В данной главе рассматриваются структурные вопросы теории ассоциативных (ко)модульных алгебр.
(Отдельные результаты остаются справедливыми и для необязательно ассоциативных алгебр.)
Доказанные утверждения будут затем использованы в главе~\ref{ChapterGenHAssocCodim} при изучении полиномиальных $H$-тождеств. Для работы с конечномерными $H$-комодульными алгебрами разрабатывается специальная техника, которая позволяет компенсировать отсутствие коумножения в алгебре $H^*$ в случае бесконечномерных алгебр Хопфа~$H$.

Результаты главы были опубликованы в работах~\cite{ASGordienko6Kochetov, ASGordienko11,ASGordienko12,ASGordienko9, ASGordienko15}.

\section{(Ко)инвариантность радикала Джекобсона}

Пусть $A$~"--- $H$-модульная алгебра для некоторой биалгебры $H$.
Подпространство $V \subseteq A$ называется  \textit{инвариантным
относительно $H$-действия}, если $HV=V$, т.е. если $V$ является $H$-подмодулем.

Аналогично, если $A$~"--- $H$-комодульная алгебра для некоторой биалгебры $H$, то
подпространство $V \subseteq A$ называется  \textit{коинвариантным
относительно $H$-кодействия}, если $\rho(V)\subseteq V \otimes H$, т.е. если $V$ является $H$-подкомодулем.

Если $A^2\ne 0$ и алгебра $A$ не содержит нетривиальных $H$-(ко)инвариантных двусторонних идеалов, то $A$
называется \textit{$H$-простой} алгеброй. В частности, алгебра $A$, градуированная группой $G$, называется \textit{$G$-градуированно простой}, если $A^2 \ne 0$ и $A$ не содержит нетривиальных градуированных идеалов. Аналогично, алгебра $A$ с действием группы $G$ автоморфизмами называется \textit{$G$-простой}, если $A^2 \ne 0$ и $A$ не содержит нетривиальных $G$-инвариантных идеалов.

В.\,В.~Линченко было доказано достаточное условие $H$-инвариантности радикала в $H$-модульной ассоциативной алгебре. Для удобства читателя приведём не только уточнённую\footnote{В работе В.\,В.~Линченко имеется неточность связанная с условием на характеристику основного поля, см. замечание~\ref{RemarkTraceCriterionForNilpotence}.} формулировку этого результата, но и его доказательство.

\begin{theorem}[{В.\,В.~Линченко~\cite[теорема 2.1]{LinchenkoJH}}]\label{TheoremRadicalHSubMod}
Пусть $A$~"--- конечномерная
ассоциативная $H$-модульная алгебра над полем $\mathbbm{k}$
для некоторой алгебры Хопфа $H$ с антиподом $S$, таким, что $S^2=\id_H$.
Предположим, что либо $\chr \mathbbm{k} = 0$, либо $\chr \mathbbm{k}>\dim A$. Тогда радикал Джекобсона $J:=J(A)$
алгебры $A$ является её $H$-подмодулем.
\end{theorem}

В доказательстве используется следующее утверждение, представляющее и самостоятельный интерес:

\begin{proposition}[{\cite[лемма 1.2]{LinchenkoJH}}]\label{PropositionHIHSubMod}
Пусть $A$~"--- (необязательно ассоциативная) $H$-модульная алгебра над полем $\mathbbm{k}$
для некоторой алгебры Хопфа $H$ с антиподом $S$, являющимся обратимым оператором.
Тогда для всякого двухстороннего идеала $I \subseteq A$ подмодуль $HI := \langle ha \mid h\in H, a\in I \rangle_\mathbbm{k}$ также является двухсторонним идеалом.
\end{proposition}
\begin{proof} Пусть $a\in I$, $b\in A$ и $h\in H$. Тогда \begin{equation}\label{Eqhab}(ha)b=(h_{(1)}a)(\varepsilon(h_{(2)})b)=
(h_{(1)}a)(h_{(2)}(Sh_{(3)})b)=h_{(1)}(a((Sh_{(2)})b)) \in HI,\end{equation}
поскольку $a(gb) \in I$ для любого $g \in H$.

Применив к равенствам $$(Sh_{(1)}) h_{(2)} =  h_{(1)}(Sh_{(2)})=\varepsilon(h)1_H$$
антигомоморфизм $S^{-1}$, для всякого $h\in H$
получаем
\begin{equation}\label{EqInverseAntipode}h_{(2)} (S^{-1}h_{(1)}) =  (S^{-1}h_{(2)})h_{(1)}=\varepsilon(h)1_H.
\end{equation}

Отсюда получаем \begin{equation}\label{Eqahb}b(ha)=(\varepsilon(h_{(1)})b)(h_{(2)}a)=
(h_{(2)}(S^{-1}h_{(1)})b)(h_{(3)}a)=h_{(2)}(((S^{-1}h_{(1)})b)a) \in HI.\end{equation}
\end{proof}
\begin{proof}[Доказательство теоремы~\ref{TheoremRadicalHSubMod}.]
Поскольку всякий ниль-идеал содержится в радикале Джекобсона (см., например, \cite[лемма 1.2.2]{Herstein}), достаточно доказать, что 
$H$-подмодуль $HJ$, который в силу предложения~\ref{PropositionHIHSubMod} является идеалом,
состоит из нильпотентных элементов.

Обозначим через $\Phi \colon A \to \End_\mathbbm{k}(A)$ \textit{левое регулярное представление} алгебры $A$,
то есть гомоморфизм,
заданный равенством $\Phi(a)b := ab$ для всех $a\in A$ и $b\in A$.
Для всякого элемента $h \in H$ обозначим через $\rho \colon H \to \End_\mathbbm{k}(A)$ гомоморфизм,
заданный равенством $\rho(h)a := ha$ для всех $a\in A$.

Докажем, что для всех $v \in HJ$ выполнено $\tr\left(\Phi(v)\right)=0$.
Действительно, если $h\in H$ и $a\in J$, то
$$\tr\left(\Phi(ha)\right)\stackrel{\eqref{Eqhab}}{=}\tr\left(\rho(h_{(1)})\Phi(a)\rho(Sh_{(2)})\right)=
\tr\left(\rho(Sh_{(2)})\rho(h_{(1)})\Phi(a)\right)\stackrel{\eqref{EqInverseAntipode}}{=}
\varepsilon(h)\tr\left(\Phi(a)\right)=0,$$
поскольку элемент $a$, а следовательно, и оператор $\Phi(a)$ нильпотентны.

Из $v^k \in HJ$ для всех $v\in HJ$ и $k\in \mathbb N$ получаем $\tr\left(\Phi(v^k)\right)=0$. Следовательно, в силу теоремы~\ref{TheoremTraceCriterionForNilpotence} оператор $\Phi(v)$ нильпотентен, т.е.
$\Phi(v)^m=0$ для некоторого $m\in \mathbb N$. Отсюда $v^{m+1}=\Phi(v)^m v = 0$
и $HJ$ является ниль-идеалом. Поэтому $HJ\subseteq J$.
\end{proof}

\begin{corollary}\label{CorollaryRadicalHSSSubMod}
Пусть $A$~"--- конечномерная
ассоциативная $H$-модульная алгебра над полем характеристики $0$
для некоторой конечномерной полупростой (или кополупростой) алгебры Хопфа $H$. Тогда радикал Джекобсона
алгебры $A$ является её $H$-подмодулем.
\end{corollary}
\begin{proof}
В силу теоремы Ларсона~"--- Рэдфорда всякая конечномерная алгебра Хопфа $H$ над полем характеристики $0$ полупроста, если и только если она кополупроста, что, в свою очередь, справедливо,
если и только если $S^2=\id_H$ (см., например, \cite[теорема 7.4.6]{Danara}). Отсюда для доказательства утверждения достаточно применить теорему~\ref{TheoremRadicalHSubMod}.
\end{proof}

Оказывается, аналогичный результат справедлив и для комодульных алгебр:
     
     \begin{theorem}\label{TheoremRadicalHSubComod}
Пусть $A$~"--- конечномерная
ассоциативная $H$-комодульная алгебра над полем $\mathbbm{k}$
для некоторой алгебры Хопфа $H$ с антиподом $S$, таким, что $S^2=\id_H$.
Предположим, что либо $\chr \mathbbm{k} = 0$, либо $\chr \mathbbm{k}>\dim A$. Тогда радикал Джекобсона $J:=J(A)$
алгебры $A$ является её $H$-подкомодулем.
\end{theorem}

Отсюда мы получаем новое доказательство\footnote{Другое доказательство в случае поля характеристики $0$ см., например, в~\cite[теорема 2.1]{KelarevOkninski}.} градуированности радикала Джекобсона в конечномерных градуированных алгебрах:

\begin{corollary}\label{CorollaryRadicalGraded} Пусть $A$~"---
конечномерная ассоциативная алгебра над полем $\mathbbm{k}$, градуированная произвольной
группой $G$, причём либо $\chr \mathbbm{k} = 0$, либо $\chr \mathbbm{k} > \dim A$. 
Тогда радикал Джекобсона алгебры $A$ является её градуированным идеалом.
\end{corollary}

Прежде чем доказывать теорему~\ref{TheoremRadicalHSubComod}, напомним
что всякая $H$-комодульная алгебра $A$, где $H$~"--- некоторая биалгебра, является $H^*$-модулем, где $H^*$~"--- алгебра, двойственная к коалгебре~$H$, и $h^*a := h^*(a_{(1)})a_{(0)}$ при $a\in A$ и $h^* \in H^*$.

Если биалгебра $H$ бесконечномерна, далеко не всегда можно определить на $H^*$ структуру коалгебры, двойственной к алгебре $H$. 
Однако на $H^*$ можно тем не менее определить некий аналог коумножения:

\begin{lemma}\label{LemmaSubstituteComult}
Пусть $H_1$~"--- конечномерное подпространство биалгебры $H$ над полем $\mathbbm{k}$.
Тогда для любого $h^*\in H^*$ существуют такие $s\in \mathbb N$, ${h_i^*}', {h_i^*}''\in H^*$, $1\leqslant i\leqslant s$, что $$h^*(hq) = \sum_{i=1}^s{h_i^*}'(h){h_i^*}''(q) \text{
для всех } h,q \in H_1.$$
В частности, если $A$~"--- (необязательно ассоциативная) $H$-комодульная алгебра, а 
$a,b \in A$~"--- такие элементы, что $\rho(a),\rho(b)\in A \otimes H_1$,
то \begin{equation}\label{EqSubstituteComult}h^*(ab)=\sum_{i=1}^s ({h_i^*}'a)({h_i^*}''b).\end{equation}
\end{lemma}
\begin{proof}
Рассмотрим линейное отображение $\Xi \colon H^* \to (H_1 \otimes H_1)^*$,
заданное при помощи равенства $\Xi(h^*)(h \otimes q) = h^*(hq)$ для всех $h^* \in H^*$,
$h,q \in H$. Используя естественное отождествление $(H_1 \otimes H_1)^*=H_1^*\otimes H_1^*$,
получаем $\Xi(h^*) = \sum_{i=1}^s {h_i^*}' \otimes {h_i^*}''$
для некоторых $s\in \mathbb N$, ${h_i^*}', {h_i^*}''\in H_1^*$, $1\leqslant i\leqslant s$.
Продолжая линейные функции с $H_1$ на $H$, можно считать, что
${h_i^*}', {h_i^*}''\in H^*$.
Тогда $h^*(hq) = \sum_{i=1}^s{h_i^*}'(h){h_i^*}''(q)$
для всех $h,q \in H_1$.
Первая часть утверждения доказана.

Если $A$~"--- $H$-комодульная алгебра и $\rho(A)\subseteq A \otimes H_1$,
то
$$h^*(ab)=h^*(a_{(1)}b_{(1)})a_{(0)}b_{(0)} = 
\sum_{i=1}^s {h_i^*}'(a_{(1)}){h_i^*}''(b_{(1)}) a_{(0)}b_{(0)}
=\sum_{i=1}^s ({h_i^*}'a)({h_i^*}''b)$$ для всех $a,b \in A$.
\end{proof}


Докажем теперь, что подкомодуль (см. лемму~\ref{LemmaGenSubcomodules}), порождённый идеалом, также является идеалом.
 
  \begin{lemma}\label{LemmaHComoduleIdeal}
 Пусть $A$~"--- (необязательно ассоциативная) $H$-комодульная алгебра над полем $\mathbbm{k}$
для некоторой алгебры Хопфа $H$ с антиподом $S$, являющимся обратимым оператором.
Тогда для всякого двухстороннего идеала $I \subseteq A$ подпространство $H^*I := \langle h^*a \mid h^*\in H^*, a\in I \rangle_\mathbbm{k}$ также является двухсторонним идеалом.
  \end{lemma}
  \begin{proof} Пусть $a \in I$, $b\in A$, $h^*\in H^*$.
Выберем такое конечномерное подпространство $H_2 \subseteq H$, что $\rho(a),\rho(b)\in A \otimes H_2$
и $(\id_A \otimes \Delta)\rho(b) \in A \otimes H_2 \otimes H_2$.
Теперь выберем для $h^*$  элементы ${h^*_i}'$ и ${h^*_i}''$
в соответствии с леммой~\ref{LemmaSubstituteComult}, где $H_1 := H_2 + H_2^2+SH_2+S^{-1}H_2$.
  
   Тогда
  $$(h^*a)b=h^*(a_{(1)})a_{(0)}b=h^*(\varepsilon(b_{(1)})a_{(1)})a_{(0)}b_{(0)}
  =h^*(a_{(1)}b_{(1)}Sb_{(2)})a_{(0)}b_{(0)}
  =$$ \begin{equation}\label{EqMoveHFirst}
  =\sum_i {h^*_i}'(a_{(1)}b_{(1)}){h^*_i}''(Sb_{(2)})a_{(0)}b_{(0)}
  =\sum_i {h^*_i}'(a (S^*{h^*_i}'') b) \in H^*I.
  \end{equation}
  Аналогично,
  $$b(h^*a)=h^*(a_{(1)})ba_{(0)}=h^*(\varepsilon(b_{(1)})a_{(1)})b_{(0)}a_{(0)}
  =h^*((S^{-1}b_{(2)})b_{(1)} a_{(1)})b_{(0)}a_{(0)}
  =$$ \begin{equation}\label{EqMoveHSecond}
  =\sum_i {h^*_i}'(S^{-1}b_{(2)}){h^*_i}''(b_{(1)}a_{(1)})b_{(0)}a_{(0)}
  =\sum_i {h^*_i}''((((S^{-1})^*{h^*_i}') b) a) \in H^*I.\end{equation}
  \end{proof}
  
  Теперь докажем одно важное свойство следа регулярного представления:
  
 \begin{lemma}\label{LemmaTraceHRad}
 Пусть $A$~"--- конечномерная ассоциативная $H$-комодульная алгебра над полем $\mathbbm{k}$
для некоторой алгебры Хопфа $H$ с антиподом $S$, таким, что $S^2=\id_H$.
 Тогда для левого регулярного представления $\Phi \colon A \to \End_\mathbbm{k}(A)$, $\Phi(a)b := ab$,
 справедливо равенство
 $$\tr(\Phi(h^* a))=h^*(1) \tr(\Phi(a))\text { для всех } h^{*}\in H^{*} \text{ и } a \in A.$$
 \end{lemma}
 \begin{proof} 
Зададим на алгебре 
$\End_\mathbbm{k}(A)$ структуру $H$-комодульной алгебры в соответствии с примером~\ref{ExampleHComodEnd}.
Кроме того, выберем такое подпространство $H_2 \subseteq H$, что $\rho(A)\subseteq A \otimes H_2$.
Положим $H_1 := H_2 + H_2^2+SH_2$ и в соответствии с леммой~\ref{LemmaSubstituteComult}
для каждого $h^*\in H^*$ фиксируем такие ${h^*_i}', {h^*_i}'' \in H^*$, что для всех $a,b\in A$ справедливо равенство~\eqref{EqSubstituteComult}.
 Тогда для всех $h^*\in H^*$, $\psi \in \End_\mathbbm{k}(A)$ и $a\in A$
 \begin{equation}\label{EqMoveHstarRight}(h^*\psi)a = h^*(\psi(a_{(0)})_{(1)}(Sa_{(1)})) \psi(a_{(0)})_{(0)}
 = \sum_i {h^*_i}'\psi((S^*{h^*_i}'') a).\end{equation}
 Следовательно, $h^*\psi = \sum_i \zeta({h^*_i}') \psi \zeta(S^*{h^*_i}'')$
 для всех $h^* \in H^*$ и $\psi \in \End_\mathbbm{k}(A)$,
 где $\zeta \colon H^* \to \End_\mathbbm{k}(A)$~"--- отображение, задающее  на $A$
 структуру $H^*$-модуля. 
  
 Заметим, что отображение $\Phi$ является гомоморфизмом $H$-комодулей.
 Отсюда $\Phi$ также является гомоморфизмом $H^*$-модулей.
 Из~\eqref{EqMoveHstarRight} следует, что
 $$\tr(\Phi(h^* a))= \sum_i \tr\zeta({h^*_i}') \Phi(a) \zeta(S^*{h^*_i}'')
 = \sum_i \tr(\zeta((S^*{h^*_i}''){h^*_i}') \Phi(a)).$$
 Учитывая, что $$\sum_i \zeta((S^*{h^*_i}''){h^*_i}') b = 
\sum_i(S^*{h^*_i}'')({h^*_i}'(b_{(1)})b_{(0)})=$$ $$=
\sum_i({h^*_i}'')(Sb_{(1)}){h^*_i}'(b_{(2)})b_{(0)}
=h^*(b_{(2)} Sb_{(1)}) b_{(0)} =h^*(1)b$$ для всех $h^*\in H^*$
и $b\in A$, получаем
 $\tr(\Phi(h^* a)) = h^*(1) \tr(\Phi(a))$.
  \end{proof}

\begin{proof}[Доказательство теоремы~\ref{TheoremRadicalHSubComod}.]
В силу лемм~\ref{LemmaGenSubcomodules} и~\ref{LemmaHComoduleIdeal} подпространство $H^* J \supseteq J$ является $H$-инвариантным двустронним идеалом. Поэтому для доказательства теоремы достаточно показать, что $H^* J$ является ниль-идеалом.
 
 Напомним, что в конечномерных алгебрах радикал Джекобсона является нильпотентным идеалом,
 откуда $\tr(\Phi(a))=0$ для всех $a \in J$. Применяя лемму~\ref{LemmaTraceHRad}, получаем,
 что $\tr(\Phi(v))=0$ для всех $v \in H^* J$ и, в частности, $\tr(\Phi(w)^k)=0$ для всех $w\in H^* J$ и $k\in\mathbb N$. В силу теоремы~\ref{TheoremTraceCriterionForNilpotence} оператор $\Phi(a)$ является нильпотентным, и $H^* J$~"--- ниль-идеал. Отсюда $H^* J = J$.
\end{proof}

   \section{$U(\mathfrak g)$-простые и $G$-простые алгебры}\label{SectionUgGSimple}
   
   В данном параграфе доказывается, что для конечномерных ассоциативных алгебр над полем характеристики нуль
   с действием некоторой алгебры Ли $\mathfrak g$ дифференцированиями понятия $U(\mathfrak g)$-простой
   и простой алгебры совпадают. Аналогичный результат верен и для алгебр с рациональным действием связной аффинной алгебраической группы.
   
 \begin{lemma}\label{LemmaDerSimpleSum}
 Пусть $B=B_1 \oplus \ldots \oplus B_q$ (прямая сумма идеалов)~"--- необязательно ассоциативная
 алгебра над произвольным полем $\mathbbm{k}$, где $B_i$~"--- простые алгебры,
 а $\delta$~"--- некоторое дифференцирование алгебры $B$.
  Тогда все алгебры $B_i$ инвариантны относительно дифференцирования $\delta$.
 \end{lemma}
 \begin{proof}
 Пусть $a \in B_i$ для некоторого $1 \leqslant i \leqslant q$. Тогда $\delta(a)=\sum_{i=1}^q b_i$,
 где $b_j \in B_j$ при $1 \leqslant j \leqslant q$.
 Для всех $b \in B_j$, где $j\ne i$, справедливо равенство $$0=\delta(ab)=\delta(a)b+a\delta(b)=b_jb+a\delta(b),$$
 откуда $b_jb = - a\delta(b) \in B_i \cap B_j$, т.е. $b_j b = 0$. Аналогично доказывается, что $b b_j = 0$ для всех $b\in B_j$. Поскольку множество $\lbrace v \in B_j \mid vB_j = B_j v = 0 \rbrace $
 образует двусторонний идеал в простой алгебре $B_j$, получаем, что  $b_j = 0$ для всех $j\ne i$ и $\delta(a)\in B_i$.
  \end{proof}
  \begin{lemma}\label{LemmaDerTrivIdealsAssoc}
  Пусть $A$~"--- конечномерная ассоциативная алгебра над полем $\mathbbm{k}$ характеристики $0$,
  на которой действует дифференцированиями некоторая алгебра Ли $\mathfrak g$.
   Предположим, что $A$ и $0$~"--- единственные $\mathfrak g$-инвариантные
   идеалы алгебры $A$. Тогда либо алгебра $A$ полупроста, либо $A^2=0$.
  \end{lemma}
  \begin{proof} В силу теоремы~\ref{TheoremRadicalHSubMod}
  радикал Джекобсона является $U(\mathfrak g)$-инвариантным идеалом,
  откуда либо $J(A)=0$ и лемма доказана, либо $A=J(A)$ является нильпотентной
  алгеброй. Во втором случае $A^2 \ne A$~"--- $\mathfrak g$-инвариантный идеал, и, как следствие, $A^2=0$.
   \end{proof}
  \begin{theorem}\label{TheoremDerSimpleAssoc}
  Если $B$~"--- конечномерная $\mathfrak g$-простая ассоциативная алгебра над полем $\mathbbm{k}$ характеристики $0$,
  на которой действует дифференцированиями некоторая алгебра Ли~$\mathfrak g$, то алгебра $B$ проста.
  \end{theorem}
  \begin{proof} Согласно лемме~\ref{LemmaDerTrivIdealsAssoc} алгебра $B$ полупроста, т.е. в силу теоремы Веддербёрна~"--- Артина $B=B_1 \oplus \ldots \oplus B_q$ (прямая сумма идеалов) для некоторых простых алгебр $B_i$. В силу леммы~\ref{LemmaDerSimpleSum} всякий идеал $B_i$ является $\mathfrak g$-инвариантным, откуда $q=1$ и $B=B_1$.
   \end{proof}
   
   Формально алгебры, на которых некоторая группа действовать не только автоморфизмами, но и антиавтоморфизмами, не являются модульными алгебрами над алгебрами Хопфа.
   Однако в случае действия связной аффинной алгебраической группы всё сводится к случаю лишь автоморфизмов:
  \begin{proposition}\label{PropositionConnectedAffAlgAutNoAnti}
Пусть $A$~"--- конечномерная (необязательно ассоциативная) алгебра над алгебраически замкнутым полем $\mathbbm{k}$,
на которой рационально действует  автоморфизмами и антиавтоморфизмами некоторая
связная аффинная алгебраическая группа $G$. Тогда
группа $G$ действует на $A$ автоморфизмами.
  \end{proposition}
  \begin{proof}  Если $\Aut(A)\ne \Aut^*(A)$, то $G$-действие индуцирует гомоморфизм групп $G \to \Aut^*(A)/\Aut(A) \cong \mathbb Z/2\mathbb Z$. При этом прообразы элементов группы
  $\mathbb Z/2\mathbb Z$ являются непересекающимися подмножествами в $G$, замкнутыми в топологии Зарисского.
  В силу связности группы $G$ она целиком отображается в нейтральный элемент группы $\Aut^*(A)/\Aut(A)$,
  откуда группа $G$ действует на $A$ автоморфизмами.
\end{proof}

Теперь докажем теорему об алгебрах, простых по отношению к действию таких групп:
  \begin{theorem}\label{TheoremGSimpleAssoc}
Пусть $A$~"--- конечномерная $G$-простая ассоциативная алгебра над алгебраически замкнутым полем $\mathbbm{k}$,
на которой рационально действует  автоморфизмами некоторая
связная аффинная алгебраическая группа $G$. Тогда алгебра $A$ проста в обычном смысле.
  \end{theorem}
  \begin{proof} Предположим теперь, что алгебра $A$ является $G$-простой.
  Поскольку радикал Джекобсона инвариантен относительно автоморфизмов, всякая $G$-простая алгебра
  полупроста. Отсюда  в силу теоремы Веддербёрна~"--- Артина $A=B_1 \oplus \ldots \oplus B_q$ (прямая сумма идеалов) для некоторых простых алгебр $B_i$. Заметим, что $B_i$ являются единственными простыми
  идеалами алгебры $A$ (достаточно рассмотреть для каждого простого идеала $I$
  алгебры $A$ произведения $IB_j, B_j I \subseteq I \cap B_j$).
  Отсюда $G$-действие индуцирует гомоморфизм $\varphi \colon G \to S_n$, где $g B_i = B_{\varphi(g)(i)}$ для всех $1 \leqslant i \leqslant q$ и $g\in G$, а $S_n$~"--- $n$-я группа подстановок.
  Следовательно, группа $G$ раскладывается в объединение попарно непересекающихся замкнутых
  подмножеств, отвечающих разным $\varphi(g) \in S_n$.
  Поскольку группа $G$ связна, гомоморфизм $\varphi$ отображает всю группу $G$
  в тождественную подстановку и  $g B_i = B_i$ для всех $1 \leqslant i \leqslant q$ и $g\in G$.
  Отсюда $q=1$, $B=B_1$, и алгебра $B$ простая.
    \end{proof}
    
    Аналоги теорем~\ref{TheoremDerSimpleAssoc} и~\ref{TheoremGSimpleAssoc} справедливы
    и для алгебр Ли. (См. теоремы~\ref{TheoremDerSimpleLie} и~\ref{TheoremGSimpleLie}).

   \section{$H$-(ко)инвариантные аналоги теорем Веддербёрна~"--- Артина и Веддербёрна~"--- Мальцева}\label{SectionH(co)invWedArtMalcev}
   
   Сперва докажем достаточные условия для (ко)инвариантности аннуляторов под(ко)модулей в (ко)модульных алгебрах:
      
  \begin{lemma}\label{LemmaAnnIdealHmod}
  Пусть $M$~"--- $H$-подмодуль в ассоциативной $H$-модульной алгебре $A$ для некоторой алгебры Хопфа $H$ с биективным антиподом $S$ над полем $\mathbbm{k}$,
  тогда подпространство $\Ann_{\mathrm{lr}}(M) := \lbrace a \in A \mid ab=ba=0 \text{ для всех } b \in M \rbrace$
  также
является $H$-подмодулем.
  \end{lemma}
  \begin{proof}
  Пусть $a \in \Ann_{\mathrm{lr}}(M)$, $b\in M$, $h\in H$. Тогда из~(\ref{Eqhab})
  следует, что $$(ha)b=h_{(1)}(a((Sh_{(2)})b))=0.$$
В то же время   
  из~(\ref{Eqahb})
  следует, что $$b(ha)=h_{(2)}(((S^{-1}h_{(1)})b)a)=0.$$
  
  Таким образом, $ha \in \Ann_{\mathrm{lr}}(M)$.
\end{proof}

  \begin{lemma}\label{LemmaAnnIdealHcomod}
  Пусть $M$~"--- $H$-подкомодуль в ассоциативной $H$-комодульной алгебре $A$ для некоторой алгебры Хопфа $H$ с биективным антиподом $S$ над полем $\mathbbm{k}$,
  тогда подпространство $\Ann_{\mathrm{lr}}(M)$
также является $H$-подкомодулем.
  \end{lemma}
  \begin{proof}
  Пусть $a \in \Ann_{\mathrm{lr}}(M)$, $b\in M$, $h^*\in H^*$. Тогда
\begin{equation*}\begin{split} (h^*a)b=h^*(a_{(1)})a_{(0)}b = h^*\bigl(a_{(1)}\varepsilon(b_{(1)})1_H\bigr)a_{(0)}b_{(0)}= \\ = h^*(a_{(1)}b_{(1)}Sb_{(2)})a_{(0)}b_{(0)} = h^*\bigl((ab_{(0)})_{(1)}Sb_{(1)}\bigr)(ab_{(0)})_{(0)}=0,\end{split}\end{equation*}
поскольку для каждого слагаемого $ab_{(0)}=0$.
Аналогично,
\begin{equation*}\begin{split} b(h^*a)=h^*(a_{(1)})ba_{(0)} = h^*\bigl(\varepsilon(b_{(1)})a_{(1)}\bigr)b_{(0)}a_{(0)}= \\ = h^*\bigl((S^{-1}b_{(2)}) b_{(1)}a_{(1)}\bigr)b_{(0)}a_{(0)} = h^*\bigl((S^{-1}b_{(1)})(b_{(0)}a)_{(1)}\bigr)(b_{(0)}a)_{(0)}=0,\end{split}\end{equation*}
поскольку для каждого слагаемого $b_{(0)}a=0$.
    Таким образом, $h^*a \in \Ann_{\mathrm{lr}}(M)$, откуда $\Ann_{\mathrm{lr}}(M)$
является $H$-подкомодулем.
\end{proof}

Рассмотрим ассоциативные алгебры, \textit{полупростые} в обычном смысле, т.е. такие алгебры,
радикал Джекобсона которых равен $0$. В доказательстве теоремы ниже мы используем идею из предложения 3.1 из работы~\cite{PaReZai}:
  
  \begin{theorem}\label{TheoremWedderburnHmod}
  Пусть $B$~"--- конечномерная полупростая ассоциативная $H$-модульная алгебра над полем $\mathbbm{k}$
для некоторой алгебры Хопфа $H$ с биективным антиподом.
Тогда $B = B_1 \oplus B_2 \oplus \ldots \oplus B_s$
  (прямая сумма идеалов, являющихся $H$-подмодулями) для некоторых $H$-простых
  $H$-модульных алгебр $B_i$.
  \end{theorem}
  \begin{proof} Докажем теорему индукцией по $\dim B$. Если алгебра $B$ является $H$-простой, доказывать
  нечего. Предположим, что $B$ содержит нетривиальные $H$-инвариантные идеалы.
В силу обычной теоремы Веддербёрна~--- Артина справедливо равенство $B=A_1 \oplus \ldots \oplus A_s$ (прямая сумма идеалов), где $A_i$~"--- простые алгебры, которые необязательно являются $H$-подмодулями.
Пусть $B_1$~"--- минимальный идеал алгебры $B$, являющийся $H$-подмодулем.
 Тогда $B_1 = A_{i_1} \oplus \ldots \oplus A_{i_k}$ для некоторых $i_1, i_2, \ldots, i_k \in \lbrace 1, 2, \ldots, s\rbrace$.  (Для доказательства этого равенства достаточно рассмотреть идеалы
 $B_1 A_j$ для всех $j$, которые либо равны $0$, либо совпадают с $A_j$ в силу минимальности последних.)
Рассмотрим $\tilde B_1 = \lbrace a \in A \mid ab=ba=0 \text{ для всех } b \in B_1 \rbrace$.
Тогда $\tilde B_1$ равна сумме всех идеалов $A_j$, где $j \notin \lbrace 
i_1, i_2, \ldots, i_k\rbrace$, и $A = B_1 \oplus \tilde B_1$.
 В силу леммы~\ref{LemmaAnnIdealHmod} идеал $\tilde B_1$ является $H$-подмодулем. 
Теперь достаточно применить предположение индукции к алгебре $\tilde B_1$.
  \end{proof}

Используя вместо леммы~\ref{LemmaAnnIdealHmod} лемму~\ref{LemmaAnnIdealHcomod}
получаем аналогичное утверждение и для комодульных алгебр:
  
  \begin{theorem}\label{TheoremWedderburnHcomod}
  Пусть $B$~"--- конечномерная полупростая ассоциативная $H$-комодульная алгебра над полем $\mathbbm{k}$
для некоторой алгебры Хопфа $H$ с биективным антиподом.
Тогда $B = B_1 \oplus B_2 \oplus \ldots \oplus B_s$
  (прямая сумма идеалов, являющихся $H$-подкомодулями) для некоторых $H$-простых $H$-комодульных алгебр $B_i$.
  \end{theorem}
  
  Для конечномерных $H$-модульных алгебр $A$, радикал Джекобсона которых не является $H$-подмодулем,
  естественно ввести понятие \textit{$H$-радикала} $J^H(A)$\label{DefJHA}, т.е. максимального $H$-инвариантного нильпотентного идеала. Такой идеал существует, поскольку сумма всех $H$-инвариантных нильпотентных идеалов
  лежит в $J(A)$ и снова является $H$-инвариантным нильпотентным идеалом.
  
  В главе~\ref{ChapterGenHAssocCodim} нам потребуется разложение и для таких алгебр, для которых $J^H(A)=0$,
  хотя $J(A)\ne 0$:
\begin{theorem}[С.\,М.~Скрябин~"--- Ф.~Ван Ойстайен]\label{TheoremSkryabinVanOystaeyen}
Пусть $A$~"--- левая $H$-модульная правоартинова ассоциативная алгебра
с единицей для некоторой алгебры Хопфа $H$ над полем $\mathbbm{k}$ такая, что $J^H(A)=0$.
Тогда $A=B_1 \oplus \ldots \oplus B_q$ (прямая сумма $H$-инвариантных идеалов) для некоторого $q\in \mathbb Z_+$ и некоторых $H$-простых $H$-модульных алгебр $B_i$.
\end{theorem}
\begin{proof}
См. \cite[теорема 1.1]{SkryabinHdecomp}, \cite[лемма 4.2]{SkryabinVanOystaeyen}.
\end{proof}

Для того, чтобы применять теорему~\ref{TheoremSkryabinVanOystaeyen},
не предполагая наличие единицы, 
сделаем следующее наблюдение, которое представляет и самостоятельный интерес:

\begin{lemma}\label{LemmaHSemiSimpleIsUnital}
Пусть $A$~"--- левая $H$-модульная правоартинова ассоциативная алгебра
с единицей для некоторой алгебры Хопфа $H$ над полем $\mathbbm{k}$ такая, что $J^H(A)=0$.
Тогда $A$~"---  $H$-модульная алгебра с единицей.
\end{lemma}
\begin{proof}
Пусть $A^+ := A \oplus \mathbbm{k}1_{A^+}$, где $1_{A^+}$~"--- присоединённая единица алгебры~$A$.
Определим действие $h 1_{A^+} := \varepsilon(h)1_{A^+}$ для всех $h\in H$.
Тогда $A^+$~"--- $H$-модульная алгебра с единицей, а $A$~"--- её двусторонний $H$-инвариантный идеал.
Поскольку $A^+/A$ и $A$ являются правоартиновыми $A^+$-модулями, $A^+$~"---
правоартинова алгебра, и мы можем применить теорему~\ref{TheoremSkryabinVanOystaeyen}.
Тогда $A^+=B_1 \oplus \ldots \oplus B_q$ (прямая сумма $H$-инвариантных идеалов) для некоторых $q\in \mathbb Z_+$ некоторых $H$-простых $H$-модульных алгебр $B_i$.
Поскольку $A$~"---  двусторонний $H$-инвариантный идеал алгебры $A^+$
и каждый идеал $AB_i$ равен либо $0$, либо $B_i$, 
мы получим, что алгебра $A$
является прямой суммой всех идеалов $B_i$, за исключением одного. Теперь утверждение леммы следует из того, что алгебры $B_i$ с единицей и $h 1_{B_i} = \varepsilon(h)
1_{B_i}$ для всех $h\in H$.
\end{proof}

Как следствие, если $A$~"--- конечномерная $H$-модульная 
ассоциативная алгебра для некоторой алгебры Хопфа $H$,
то $A/J^H(A) = B_1 \oplus B_2 \oplus \ldots
\oplus B_q$ (прямая сумма $H$-инвариантных идеалов) для некоторых $H$-простых $H$-модульных алгебр $B_i$.

Нам также понадобится следующее обобщение теоремы Веддербёрна~"--- Мальцева:

\begin{theorem}[Д.~Штефан~"--- Ф.~Ван Ойстайен]
Пусть $H$~"--- алгебра Хопфа над произвольным полем, для которой существует ad-инвариантный левый интеграл $t\in H^*$,
такой, что $t(1)=1$, а $A$~"--- конечномерная $H$-комодульная ассоциативная алгебра, причём $J(A)$ является
$H$-подкомодулем, а $A/J(A)$~"--- сепарабельная алгебра.
Тогда существует такая $H$-коинвариантная максимальная полупростая подалгебра $B$, что $A=B\oplus J(A)$ (прямая сумма $H$-подкомодулей).
\end{theorem}
\begin{proof}
Доказательство для случая, когда алгебра $A$ с единицей см. в~\cite[теорема~2.4]{SteVanOyst}.

В случае, когда алгебра $A$ не содержит единицы, как и в доказательстве леммы~\ref{LemmaHSemiSimpleIsUnital},
рассмотрим алгебру $A^+ := A \oplus \mathbbm{k}1_{A^+}$, где $1_{A^+}$~"--- присоединённая единица алгебры~$A$,
и определим кодействие $\rho(1_{A^+}) := 1_{A^+} \otimes 1_H$.
Тогда $A^+$~"--- $H$-комодульная алгебра с единицей, а $A$~"--- её двусторонний $H$-коинвариантный идеал,
причём в силу своей нильпотентности $J(A^+)\subseteq A$, т.е. $J(A^+)=J(A)$.
Поскольку алгебра $A^+$ содержит единицу, а $A^+/J(A)\cong \mathbbm{k} \oplus A/J(A)$~"--- сепарабельная алгебра (см., например, следствие из \S 10.6 в~\cite{PierceAssoc}),
в силу сделанных замечаний существует разложение $A^+ = B_0 \oplus J(A)$ (прямая сумма $H$-подкомодулей), где
$B_0$~"--- некоторая $H$-коинвариантная максимальная подалгебра алгебры $A^+$.
Отметим при этом, что единица алгебры $B_0$ совпадает с единицей алгебры $A^+$
(например, это следует из нильпотентности идеала $J(A)$).
Отсюда $B_0 =  \mathbbm{k}1_{A^+} \oplus B$
и $A=B\oplus J(A)$ (прямая сумма $H$-подкомодулей), где $B:= A \cap B_0$.
\end{proof}

\begin{corollary}
Пусть $H$~"--- алгебра Хопфа над полем характеристики $0$, для которой существует ad-инвариантный левый интеграл $t\in H^*$,
такой, что $t(1)=1$, а $A$~"--- конечномерная $H$-комодульная ассоциативная алгебра, причём $J(A)$ является
$H$-подкомодулем.
Тогда существует такая $H$-коинвариантная максимальная полупростая подалгебра $B$, что $A=B\oplus J(A)$ (прямая сумма $H$-подкомодулей).
\end{corollary}

С учётом примеров~\ref{ExampleIntegralHSS}--\ref{ExampleIntegralAffAlgGr}, \ref{ExampleComoduleGraded}, \ref{ExampleRegularAffActAutAlgebra},
  следствий~\ref{CorollaryRadicalHSSSubMod}, \ref{CorollaryRadicalGraded}
 и двойственности  между $H$-действиями и $H^*$-кодействиями при $\dim H < +\infty$ получаются также такие следствия:

\begin{corollary}\label{CorollaryGradedWedderburnMalcev}
Пусть $A$~"--- конечномерная ассоциативная алгебра над полем характеристики $0$, 
градуированная произвольной группой.
Тогда существует такая градуированная максимальная полупростая подалгебра $B$, что $A=B\oplus J(A)$ (прямая сумма градуированных
подпространств).
\end{corollary}

\begin{corollary}\label{CorollaryGReductWedderburnMalcev}
Пусть $A$~"--- конечномерная ассоциативная алгебра над алгебраически замкнутым полем характеристики $0$, 
на которой рационально действует автоморфизмами редуктивная аффинная алгебраическая группа $G$.
Тогда существует такая $G$-инвариантная максимальная полупростая подалгебра $B$, что $A=B\oplus J(A)$ (прямая сумма $G$-инвариантных подпространств).
\end{corollary}

\begin{corollary}\label{CorollaryHSSWedderburnMalcev}
Пусть $A$~"--- конечномерная ассоциативная $H$-(ко)модульная алгебра над полем характеристики $0$, 
где $H$~"--- конечномерная (ко)полупростая алгебра Хопфа.
Тогда существует такая $H$-(ко)инвариантная максимальная полупростая подалгебра $B$, что $A=B\oplus J(A)$ (прямая сумма $H$-(ко)подмодулей).
\end{corollary}

Приведём теперь примеры $H$-модульных ассоциативных алгебр, для которых $H$-инвариантное
разложение Веддербёрна~"--- Мальцева не существует.

\begin{example}[Ю.\,А.~Бахтурин]\label{ExampleGnoninvWedderburnMalcev}
Пусть $$A = \left\lbrace\left(\begin{array}{cc} C & D \\
0 & 0
  \end{array}\right) \mathrel{\biggl|} C, D\in M_m(\mathbbm{k})\right\rbrace
  \subseteq M_{2m}(\mathbbm{k}),$$ $m \geqslant 2$.
  Тогда  \begin{equation}\label{EqGnoninvWedderburnMalcev}
  J(A)=\left\lbrace\left(\begin{array}{cc} 0 & D \\
0 & 0
  \end{array}\right) \mathrel{\biggl|} D\in M_m(\mathbbm{k})\right\rbrace.
  \end{equation}  
      Определим $\varphi \in \Aut(A)$ по формуле
  $$\varphi\left(\begin{array}{cc} C & D \\
0 & 0
  \end{array}\right)=\left(\begin{array}{cc} C & C+D \\
0 & 0
  \end{array}\right).$$
  Тогда группа $G=\langle \varphi \rangle
  \cong \mathbb Z$ действует на алгебре $A$ автоморфизмами, т.е. $A$ является $\mathbbm{k}G$-модульной
  алгеброй.
  Однако в $A$ не существует такой $\mathbbm{k}G$-инвариантной полупростой подалгебры
   $B$, что $A=B\oplus J(A)$ (прямая сумма $\mathbbm{k}G$-подмодулей).
\end{example}
\begin{proof} 
Пусть $a \in A$. Тогда $\varphi(a)-a \in J(A)$. Предположим, что $B$~"--- такое $G$-инвариантное 
подпространство, что $B \cap J(A) = \lbrace 0 \rbrace$. Тогда $\varphi(b)-b=0$ для всех $b \in B$
и $B \subseteq J(A)$. Следовательно, $B=0$ и $\mathbbm{k}G$-инвариантного разложения Веддербёрна~"--- Мальцева не существует.
\end{proof}

\section{Связь между дифференцированиями и автоморфизмами}

В данном параграфе доказывается, что можно заменить всякое действие конечномерной полупростой
алгебры Ли дифференцированиями на рациональное действие связной редуктивной аффинной алгебраической группы
автоморфизмами, что позволяет перенести ряд результатов на случай алгебр с действием алгебр Ли.

Рассмотрим сперва классическую связь между автоморфизмами и дифференцированиями:

\begin{proposition}
\label{PropositionDerAutConnection}
Пусть $A$~"--- необязательно ассоциативная конечномерная алгебра
над алгебраически замкнутым полем $\mathbbm{k}$ характеристики $0$, на которой
задано рациональное представление некоторой 
связной аффинной алгебраической группы $G$.
Обозначим через~$\mathfrak g$ алгебру Ли группы $G$
и определим представление алгебры Ли $\mathfrak g$ на $A$
как дифференциал представления группы $G$.
Тогда \begin{enumerate}\item  алгебра Ли $\mathfrak g$ действует на $A$ дифференцированиями,
если и только если группа $G$ действует на $A$ автоморфизмами;
\item все $\mathfrak g$-подмодули в $A$
являются $G$-инвариантными подпространствами и наоборот.\end{enumerate}
\end{proposition}
\begin{proof} Умножение $\mu \colon A \otimes A \to A$
можно отождествить с элементом $\mu = \sum_{i} \mu_{1i} \otimes \mu_{2i} \otimes \mu_{3i} \in A^*\otimes A^* \otimes A$. Группа $G$ и алгебра Ли $\mathfrak g$
действуют на $A^*\otimes A^* \otimes A$ следующим образом: $$g(u(\cdot) \otimes v(\cdot) \otimes w)=u(g^{-1}(\cdot))
\otimes v(g^{-1}(\cdot)) \otimes (gw),$$ $$\delta (u(\cdot) \otimes v(\cdot) \otimes w)
=u(\cdot) \otimes v(\cdot) \otimes \delta w - u(\delta (\cdot)) \otimes v(\cdot) \otimes w
- u(\cdot) \otimes v(\delta (\cdot)) \otimes w,$$
где $u, v \in A^*$, $w \in A$, $\delta \in \mathfrak g$, $g\in G$.
Алгебра Ли $\mathfrak g$ действует на $A$ дифференцированиями, если и только если $\delta (bc)=(\delta b)c+b(\delta c)$ для всех $b,c \in A$, $\delta \in \mathfrak g$,
что в свою очередь эквивалентно равенствам
$\sum_{i} \mu_{1i}(b)  \mu_{2i}(c) (\delta \mu_{3i})
= \sum_{i} (\mu_{1i}(\delta b)  \mu_{2i}(c) \mu_{3i}+\mu_{1i}(b)  \mu_{2i}(\delta c) \mu_{3i})$,
 $\delta  \mu = 0$ для всех $\delta  \in \mathfrak g$ и $\mathfrak g \mu = 0$. В силу теоремы~13.2 из~\cite{HumphreysAlgGr} последнее равенство эквивалентно
  равенству $G\mu = \mu$, что в свою очередь равносильно равенству
 $g(bc)=(gb)(gc)$, которое означает, что группа $G$ действует на $A$ автоморфизмами. Используя теорему~13.2 из~\cite{HumphreysAlgGr} ещё раз, получаем,
 что $G$ и $\mathfrak g$ имеют в $A$ одни и те же инвариантные подпространства. 
\end{proof}

Теперь докажем требуемый результат о замене действия:

\begin{theorem}
\label{TheoremLieDiffActionReplacement}
Пусть $A$~"--- необязательно ассоциативная конечномерная алгебра
над алгебраически замкнутым полем $\mathbbm{k}$ характеристики $0$, на которой
действует дифференцированиями конечномерная полупростая алгебра Ли $\mathfrak g$.
Тогда существует такое рациональное действие односвязной
полупростой аффинной алгебраической группы $G$ на $A$ автоморфизмами,
что выполняются следующие условия:
 \begin{enumerate}\item алгебра Ли группы $G$ совпадает с $\mathfrak g$;
\item $\mathfrak g$-действие на $A$ является дифференциалом $G$-действия на $A$;
\item все $\mathfrak g$-подмодули в $A$
являются $G$-инвариантными подпространствами и наоборот.\end{enumerate}
\end{theorem}
\begin{proof} В силу теоремы~5.1 из главы XVIII монографии~\cite{Hochschild}
существует односвязная аффинная алгебраическая группа~$G$,
алгебра Ли которой изоморфна алгебре Ли~$\mathfrak g$.

 В силу теоремы 13.5 из~\cite{HumphreysAlgGr} группа $G$ полупроста, а значит и редуктивна.

 Алгебра $A$ как $\mathfrak g$-модуль  является прямой суммой неприводимых $\mathfrak g$-подмодулей $M_i$,
 которые отвечают некоторым доминантным весам $\lambda_i$ алгебры Ли $\mathfrak g$.
Определим на $M_i$ неприводимые рациональные представления группы $G$
со старшими весами~$\lambda_i$. Теперь утверждение теоремы следует из предложения~\ref{PropositionDerAutConnection}.
\end{proof}
\begin{corollary}\label{CorollaryDiffWedderburn} Пусть $A$~"--- конечномерная ассоциативная алгебра над алгебраически замкнутым полем характеристики $0$, 
на которой
действует дифференцированиями конечномерная полупростая алгебра Ли $\mathfrak g$.
Тогда существует такая $\mathfrak g$-инвариантная максимальная полупростая подалгебра $B$, что $A=B\oplus J(A)$ (прямая сумма $\mathfrak g$-инвариантных подпространств).
\end{corollary}
\begin{proof} Достаточно воспользоваться теоремой~\ref{TheoremLieDiffActionReplacement} и следствием~\ref{CorollaryGReductWedderburnMalcev}.
\end{proof}

Как показывает пример ниже, в случае неполупростых алгебр Ли $\mathfrak g$ подобный результат может быть неверен.

\begin{example}\label{ExampleU(L)noninvWedderburnMalcev}
Пусть $A$~"---  та же ассоциативная алгебра, что и в примере~\ref{ExampleGnoninvWedderburnMalcev}.
Определим на векторном пространстве $A$ структуру алгебры Ли при помощи коммутатора  $[x,y]=xy-yx$
и обозначим соответствующую алгебру Ли через $\mathfrak g$.
Рассмотрим стандартное представление алгебры Ли $\mathfrak g$ на $A$ дифференцированиями.
Тогда $A$ оказывается $U(\mathfrak g)$-модульной ассоциативной алгеброй
(см. пример~\ref{ExampleUgModule}), но в $A$ не существует такой $U(\mathfrak g)$-инвариантной полупростой подалгебры $B$, что $A=B\oplus J(A)$ (прямая сумма $U(\mathfrak g)$-подмодулей).
\end{example}
\begin{proof}
Предположим, что $A=B\oplus J(A)$, где $B$~"--- некоторый $U(\mathfrak g)$-подмодуль. Тогда $B$ является идеалом алгебры Ли $\mathfrak g$, а $J(A)$~"--- абелевым идеалом алгебры Ли~$\mathfrak g$. Отсюда $J(A)$
содержится в центре алгебры Ли $\mathfrak g$, что не соответствует действительности. 
Получаем противоречие. Следовательно, $U(\mathfrak g)$-инвариантного разложения
Веддербёрна~"--- Мальцева для алгебры $A$ не существует.
\end{proof}

     \section{Алгебры с действием расширений Оре}
     \label{SectionHSimpleOreExtStructure}
     
В \S\ref{SectionHSimpleOreExtMul}--\ref{SectionTaftSimpleSemisimple}
будут классифицированы конечномерные $H_{m^2}(\zeta)$-простые $H_{m^2}(\zeta)$-модульные ассоциативные алгебры,
где $H_{m^2}(\zeta)$~"--- алгебра Тафта.
Однако определённые структурные результаты могут быть получены и в более общем случае
алгебр, простых по отношению к действию алгебр Хопфа, полученных при помощи расширений Оре.
Эти результаты будут затем использованы в главе~\ref{ChapterGenHAssocCodim} для
доказательства гипотезы Амицура~"--- Бахтурина для $H$-модульных алгебр в случае, когда алгебра Хопфа $H$
получена при помощи (возможно, многократного)
расширения Оре конечномерной полупростой алгебры Хопфа косопримитивными элементами (теорема~\ref{TheoremHOreAmitsur}).

В данном параграфе доказывается, что если $A$~"---
конечномерная $H$-простая $H$-модульная
алгебра, где $H$~"--- алгебра Хопфа,
полученная при помощи расширения Оре некоторой алгебры Хопфа 
 $\tilde H$  косопримитивным элементом $v$,
то алгебра $A/J^{\tilde H}(A)$ является $\tilde H$-простой.
Если при этом  $A/J^{\tilde H}(A) \ne 0$, то,
применяя оператор $v$ достаточное число раз,
можно перевести всякий ненулевой элемент
в элемент, не принадлежащий $J^{\tilde H}(A)$.
Это свойство будет затем использовано в доказательстве
теоремы~\ref{TheoremHOrePropertyStar}.

Напомним сперва определение расширения Оре.

\begin{definition}
Пусть заданы ассоциативная алгебра $C$ над полем $\mathbbm{k}$, 
 автоморфизм $\varphi \colon C\mathrel{\widetilde{\to}}C$ и
\textit{$\varphi$-косое дифференцирование} $\delta \colon C \to C$, т.е.
такое отображение $\delta$, что
$\delta(ab)=\varphi(a)\delta(b)+\delta(a)b$ для всех $a,b\in C$.
Тогда \textit{расширение Оре} $C[v,\varphi,\delta]$
состоит из всевозможных формальных многочленов $\sum_{i=0}^n a_i v^n$, где $a_i\in C$, $n\in\mathbb Z_+$,
умножение которых индуцировано умножением в алгебре $C$ с учётом соотношения
$va - \varphi(a) v = \delta(a)$ для всех $a\in C$.
\end{definition}

Заметим, что если $B$~"--- алгебра, порождённая элементом $v$
 и подалгеброй $C \subseteq B$,
а $\varphi \colon C\mathrel{\widetilde{\to}}C$~"--- такой автоморфизм алгебры, что $v a - \varphi(a)v \in C$ для всех $a\in C$,
то
отображение $\delta$, определённое при помощи равенства $$\delta(a)=v a - \varphi(a)v \text{ для } a\in C,$$ всегда является $\varphi$-косым дифференцированием на $C$.
В этом случае $B \cong C[v,\varphi,\delta]/I$, где $I$~"--- некоторый такой идеал алгебры $B$,
что $C \cap I = 0$.

\begin{definition}\label{DefinitionOreConstructed}
Мы будем говорить, что алгебра $B$ \textit{получена при помощи расширения Оре алгебры $C$ элементом $v$}.
\end{definition}

Элемент $v\in H$ называется \textit{косопримитивным},
если $\Delta v = g_1 \otimes v + v \otimes g_2$
для некоторых $g_1, g_2 \in G(H)$.

\begin{example}
Алгебра Тафта $H_{m^2}(\zeta)$ получена при помощи расширения Оре групповой алгебры $\mathbbm{k}C_m$, где $C_m=\langle c\rangle_m$~"--- циклическая группа порядка $m$, косопримитивным элементом $v$. Как алгебра
$H_{m^2}(\zeta) \cong \mathbbm{k}C_n[v,\varphi,0]/I$, где автоморфизм $\varphi$
задан равенством $\varphi(c)=\zeta c$, а идеал $I$ порождён элементом $v^m$. 
\end{example}

С другими примерами алгебр Хопфа, полученных при помощи расширений Оре можно познакомиться, например, в~\cite[\S 5.6]{Danara}.

\begin{lemma}\label{LemmaHOreRadical}
Пусть алгебра Хопфа $H$ над полем $\mathbbm{k}$
порождена как алгебра своей подалгеброй Хопфа $\tilde H$
и косопримитивным элементом $v\in H$, где $\Delta v = g_1 \otimes v + v \otimes g_2$, $g_1, g_2 \in G(\tilde H)$. Пусть $A$~"--- конечномерная ассоциативная $H$-простая $H$-модульная алгебра.
Тогда либо $A$ является $\tilde H$-простой, либо $J^{\tilde H}(A)\ne 0$.
\end{lemma}
\begin{proof}
Предположим, что $J^{\tilde H}(A)= 0$.
В силу теоремы~\ref{TheoremSkryabinVanOystaeyen}
и леммы~\ref{LemmaHSemiSimpleIsUnital}
существует разложение
$A=N_1 \oplus \ldots \oplus N_s$ (прямая сумма $\tilde H$-инвариантных идеалов)
для некоторых $\tilde H$-простых алгебр $N_i$ с единицей. 
Докажем, что все $N_i$ являются $H$-подмодулями.
Достаточно показать, что
$vN_i\subseteq N_i$ для всех $1\leqslant i \leqslant s$.
Пусть $a\in N_i$ и $b\in N_j$, $i\ne j$.
Тогда $$(va)b=v(a (g_2^{-1}b))-(g_1a)(vg_2^{-1}b)=-(g_1 a)(vg_2^{-1}b)
\in N_i.$$ В то же время $b\in N_j$ и $(va)b \in N_j$.
Следовательно, $(va)b=0$. Аналогично получаем, что $b(va)=0$. Поскольку элемент $b\in N_j$, где $j\ne i$, 
является произвольным, получаем $va \in N_i$. Следовательно, все $N_i$ являются $H$-инвариантными,
$s=1$, и алгебра $A$ является $\tilde H$-простой.
\end{proof}

Для всякой алгебры Хопфа $H$ будем обозначать через 
$\Aut_{\mathbf{Alg}}(H)$ группу всевозможных
автоморфизмов $H$ как алгебры, т.е. таких линейных биекций,
от которых требуется лишь сохранение умножения.

\begin{theorem}\label{TheoremHSimpleOreExtStructure}
Пусть алгебра Хопфа $H$ над полем $\mathbbm{k}$
порождена как алгебра своей подалгеброй Хопфа $\tilde H$
и косопримитивным элементом $v\in H$, где $\Delta v = g_1 \otimes v + v \otimes g_2$, $g_1, g_2 \in G(\tilde H)$. Предположим также, что существует такой автоморфизм
$\varphi\in\Aut_{\mathbf{Alg}}(\tilde H)$, что
что $vh-\varphi(h)v\in \tilde H$ для всех $h\in\tilde H$.
Пусть $A$~"--- конечномерная ассоциативная $H$-простая $H$-модульная алгебра с $J^{\tilde H}(A)\ne 0$. 
Обозначим через $q$ такое натуральное число, что 
$J^{\tilde H}(A)^q=0$, а $J^{\tilde H}(A)^{q-1}\ne 0$.
Пусть $\tilde J \subseteq J^{\tilde H}(A)^{q-1}$~"--- минимальный двухсторонний
 $\tilde H$-инвариантный идеал.
Тогда существует такое $m\in\mathbb N$,
что $$I_k := \bigoplus_{i=0}^{k-1} v^i \tilde J$$
при всех $1 \leqslant k \leqslant m$ являются $\tilde H$-инвариантными
идеалами,
$I_m = A$, $I_{m-1} = J^{\tilde H}(A)$.
Также для всех $0 \leqslant k \leqslant m-1$ отображение $v^k \tilde J \twoheadrightarrow v^{k+1}\tilde J$,
заданное при помощи формулы $a \mapsto va$, является $\mathbbm{k}$-линейной
биекцией. Более того, $A/J^{\tilde H}(A)$~"--- $\tilde H$-простая алгебра.
\end{theorem}

Перед тем, как перейти к доказательству
теоремы~\ref{TheoremHSimpleOreExtStructure},
дадим несколько определений.

Пусть $A$~"--- (левая) $H$-модульная алгебра. Будем говорить, что $M$~"--- \textit{$(H,(A,A))$-бимодуль},
если $M$~"--- (левый) $H$-модуль и $(A,A)$-бимодуль и при этом $h(am)=(h_{(1)}a)(h_{(2)}m)$,
 $h(ma)=(h_{(1)}m)(h_{(2)}a)$ для всех $a\in A$ и $m\in M$.

Пусть $\psi \colon M_1 \to M_2$~"--- $\mathbbm{k}$-линейное отображение, где $M_1$ и $M_2$
являются $(H,(A,A))$-бимодулями. Пусть $g_1, g_2 \in G(H)$
и пусть $\varphi\in\Aut_{\mathbf{Alg}}(H)$. Будем говорить, что $\psi$~"--- \textit{$(\varphi, g_1, g_2)$-гомоморфизм}, если $\psi(am)=(g_1 a)\psi(m)$, $\psi(ma)=\psi(m) (g_2 a)$ и
 $\psi(hm)=\varphi(h)\psi(m)$.
 Если такое отображение $\psi$ биективно, будем говорить, что $\psi$~"---
\textit{$(\varphi, g_1, g_2)$-изоморфизм}. Ясно, что в этом случае $\psi^{-1}$
является $(\varphi^{-1}, g_1^{-1}, g_2^{-1})$-изоморфизмом.

\begin{proof}[Доказательство теоремы~\ref{TheoremHSimpleOreExtStructure}.]
Если идеал $\tilde J$ инвариантен относительно действия $v$, то
 $\tilde J$ является $H$-инвариантным идеалом алгебры $A$, откуда $A=\tilde J$ в силу $H$-простоты алгебры $A$. Тогда из $\tilde J\subseteq J(A)$ следует, что $A^2=\tilde J^2 = 0$, что противоречит определению
  $H$-простой алгебры.
Отсюда $v \tilde J \subsetneqq \tilde J$.

Введём обозначение $$I_k := \sum\limits_{i=0}^{k-1} v^i \tilde J\text{ при }k\geqslant 1.$$
Поскольку $hv a \in v\varphi^{-1}(h)a + \tilde H a$
для всех $h\in \tilde H$ и $a\in A$,
подпространства $I_k$ являются $\tilde H$-подмодулями.

Более того, поскольку элемент $v$ косопримитивен,
справедливы равенства $$b(va) = v ((g_1^{-1} b)a) - (v g_1^{-1} b)(g_2 a)$$
и $$(va)b = v(a(g_2^{-1}b)) - (g_1a)(vg_2^{-1}b)$$ для всех $a,b \in A$.
Рассматривая произвольные $a\in I_k$, по индукции получаем, что все $I_k$ являются двусторонними
идеалами.

В силу того, что алгебра $A$ конечномерна, существует такое число  $m\in\mathbb N$,
что $I_k \subsetneqq I_{k+1}$ для всех $1\leqslant k \leqslant m-1$,
а $I_{m+1}=I_m$. Тогда $\tilde H$-модуль $I_m$ инвариантен относительно действия $v$,
откуда $A=I_m$ в силу $H$-простоты алгебры $A$.

Пусть $I_0:=0$. Определим сюрьективные $\mathbbm{k}$-линейные отображения $$\theta_k \colon I_k/I_{k-1} \twoheadrightarrow I_{k+1}/I_k,\text{ где }1 \leqslant k \leqslant m-1,$$ при помощи равенств $$\theta_k(a + I_{k-1}):= va+ I_k,\quad a\in I_k.$$

Заметим, что $$\theta_k((a+I_{k-1})b)=v(ab)+ I_k=(g_1 a)(vb)+(v a)(g_2b)+ I_k =\theta_k(a+I_{k-1})(g_2b)$$
и $$\theta_k(b(a+I_{k-1}))= v(ba)+ I_k
=(g_1 b)(va)+(v b) (g_2a)+ I_k =
(g_1 b)\theta_k(a+I_{k-1})$$
для всех $a\in I_k$, $b\in A$.

Более того, $$\theta_k(h(a+I_{k-1}))=vha + I_k=
\varphi(h) va +I_k
= \varphi(h)\theta_k(a+I_{k-1}).
$$
Следовательно, $\theta_k$ является $(\varphi,g_1,g_2)$-гомоморфизмом.

Поскольку $I_1=\tilde J$ является неприводимым $(\tilde H,(A,A))$-бимодулем,
либо отображение $\theta_1$ биективно, либо $I_2/I_1=0$, т.е. $m=1$.
В первом случае $I_2/I_1$
также является неприводимым $(\tilde H,(A,A))$-бимодулем.
Продолжая эту процедуру, мы получаем, что $I_k/I_{k-1}$, $1\leqslant k \leqslant m$,
являются неприводимыми $(\tilde H,(A,A))$-бимодулями, $(\varphi^{1-k},g_1^{1-k},g_2^{1-k})$-изоморфными~$\tilde J$. 
Сравнивая их размерности, получаем, что отображение $v^k \tilde J \twoheadrightarrow v^{k+1}\tilde J$,
определённое при помощи формулы $a \mapsto va$ 
должно быть $\mathbbm{k}$-линейной биекцией
для всех $0 \leqslant k \leqslant m-1$ и сумма в определении идеалов $I_k$ прямая.

Из $(\tilde H,(A,A))$-бимодульной версии теоремы Жордана~"--- Гёльдера
следует, что в любом композиционном ряде $(\tilde H,(A,A))$-бимодулей
в $A$ всякий фактор $(\tilde \varphi,\tilde g_1, \tilde g_2)$-изоморфен
$(\tilde H,(A,A))$-бимодулю $\tilde J$ для подходящих $\tilde\varphi \in \Aut_{\mathbf{Alg}}(\tilde H)$, $\tilde g_1, \tilde g_2\in G(\tilde H)$. 
В силу теоремы~\ref{TheoremSkryabinVanOystaeyen}
существует разложение
$A/J^{\tilde H}(A)=N_1 \oplus \ldots \oplus N_s$ (прямая сумма $\tilde H$-инвариантных идеалов)
для некоторых $\tilde H$-простых алгебры $N_i$. 
Предположим, что $s\geqslant 2$. Поскольку $N_1$ и $N_2$
являются неприводимыми факторами в композиционном
ряде $(\tilde H,(A,A))$-бимодулей алгебры $A$,
существует $(\tilde\varphi,\tilde g_1,\tilde g_2)$-изоморфизм $\psi \colon N_1 \mathrel{\widetilde\to} N_2$
для некоторых $\tilde g_1, \tilde g_2 \in G(\tilde H)$
и $\tilde\varphi\in\Aut_{\mathbf{Alg}}(H)$.

Обозначим через $\bar a$ образ элемента $a\in A$ в $A/J^{\tilde H}(A)$.
Тогда для всех $\bar a, \bar b\in N_2$ 
имеем $$\bar a \bar b = \bar a b=\psi(\psi^{-1}(\bar a)) b
=\psi(\psi^{-1}(\bar a)(\tilde g_2^{-1} b))
= \psi(\psi^{-1}(\bar a)(\tilde g_2^{-1} \bar b))
= 0,$$ поскольку $\psi^{-1}(\bar a)\in N_1$ и $\tilde g_2^{-1} \bar b \in N_2$.
Отсюда $N_2^2=0$ и мы получаем
противоречие с тем, что идеал $N_2$ является $\tilde H$-простой алгеброй. Следовательно, $s=1$
и $A/J^{\tilde H}(A)$
является $\tilde H$-простой алгеброй. В частности, $J^{\tilde H}(A)$ является
максимальным $\tilde H$-инвариантным идеалом.

Докажем, что $J^{\tilde H}(A)$ является единственным максимальным $\tilde H$-инвариантным
идеалом. Действительно, если $I \subseteq A$ был бы другим $\tilde H$-инвариантным
идеалом, таким, что $I \subsetneqq J^{\tilde H}(A)$, тогда было бы справедливо равенство $A=I+J^{\tilde H}(A)$
и существовал бы изоморфизм $$A/(I\cap J^{\tilde H}(A)) \cong I/(I\cap J^{\tilde H}(A)) \oplus
J^{\tilde H}(A) /(I\cap J^{\tilde H}(A))$$ (прямая сумма $\tilde H$-инвариантных идеалов).
Поскольку в силу леммы~\ref{LemmaHSemiSimpleIsUnital} алгебра $A$ с единицей, алгебра $A/(I\cap J^{\tilde H}(A))$ и, следовательно,
алгебра $J^{\tilde H}(A)/(I\cap J^{\tilde H}(A))$ также
должны были бы быть с единицей, что противоречит нильпотентности
идеала 
$J^{\tilde H}(A)$. Следовательно, $J^{\tilde H}(A)$ действительно
является единственным максимальным $\tilde H$-инвариантым идеалом алгебры $A$.

 Поскольку $A/I_{m-1}$ как $(\tilde H,(A,A))$-бимодуль  
  $(\varphi^{1-m},g_1^{1-m},g_2^{1-m})$-изоморфен $\tilde J$,
  а $\tilde J$ является неприводимым $(\tilde H,(A,A))$-бимодулем, $A/I_{m-1}$~"--- $\tilde H$-простая алгебра.
 Следовательно,
$I_{m-1}$ также является максимальным $\tilde H$-инвариантным идеалом, и $I_{m-1}=J^{\tilde H}(A)$.
\end{proof}

При помощи теоремы~\ref{TheoremHSimpleOreExtStructure}
можно определить $\mathbbm{k}$-линейное отображение $\psi 
\colon A \to A$ при помощи формул $\psi(v^k a) := v^{k-1} a$
для $1\leqslant k \leqslant m-1$ и $a\in \tilde J$, и $\psi(\tilde J) :=0$.
Заметим, что $\psi^m=0$, а $\psi(A)= J^{\tilde H}(A)$.

Определим теперь $\mathbbm{k}$-линейные отображения $\psi_k \colon A/J^{\tilde H}(A)\to A$, где $0\leqslant k \leqslant m-1$, следующим образом.
Для всякого $\bar a \in \bar A := A/J^{\tilde H}(A)$
существует единственный такой элемент $a\in v^{m-1} \tilde J$,
что $\bar a = a + J^{\tilde H}(A)$. Определим 
$\psi_k(\bar a) := \psi^k(a)$, $0\leqslant k \leqslant m-1$.
Тогда $A=\bigoplus_{k=0}^{m-1}\psi_k(\bar A)$.
Другими словами, алгебра $A$ составлена из образов алгебры $\bar A$.
Произведения элементов $\psi_k(a)$ будут вычислены в следующем разделе.

Исследуем, как ведёт себя отображение $\psi_{m-1}$ по отношению к $\tilde H$-действию.

\begin{proposition}\label{PropositionHOreActionOnm1}
Предположим, что выполнены предположения теоремы~\ref{TheoremHSimpleOreExtStructure}.
Тогда
$h \psi_{m-1}(\bar a)= \psi_{m-1}(\varphi^{m-1}(h)\bar a)$ для всех $\bar a\in \bar A$ и $\tilde h \in \tilde H$. 
\end{proposition}
\begin{proof} Поскольку $\psi_{m-1}(\bar a) \in \tilde J$, а $\tilde J$~"--- $\tilde H$-инвариантный идеал,
$h \psi_{m-1}(\bar a)= b $ для некоторого $b\in \tilde J$.
Пусть $\delta(h):=vh-\varphi(h)v$, где $h\in\tilde H$.
 Тогда 
\begin{equation*}\begin{split}v^{m-1} b = v^{m-1} h \psi_{m-1}(\bar a)=v^{m-2} \varphi(h) v \psi_{m-1}(\bar a)+v^{m-2} \delta(h)
\psi_{m-1}(\bar a) =\\ = \varphi^{m-1}(h) v^{m-1} \psi_{m-1}(\bar a)+ \sum_{i=0}^{m-2}  v^{m-i-2} \delta\left(\varphi^{i}(h)\right) v^i \psi_{m-1}(\bar a)=\\=
\varphi^{m-1}(h) \psi_0(a) + \sum_{i=0}^{m-2}  v^{m-i-2} \delta\left(\varphi^{i}(h)\right) v^i \psi_{m-1}(\bar a).\end{split}\end{equation*}
В силу теоремы~\ref{TheoremHSimpleOreExtStructure} второе слагаемое принадлежит $J^{\tilde H}(A)$.
Следовательно, $\pi(v^{m-1} b)= \varphi^{m-1}(h) \bar a$, где $\pi \colon A \twoheadrightarrow A/J^{\tilde H}(A)$~"--- естественный сюръективный гомоморфизм. Поскольку $v^{m-1}b \in v^{m-1} \tilde J$, получаем $\psi_0(\varphi^{m-1}(h) \bar a)=v^{m-1}b$ и $$h \psi_{m-1}(\bar a)=b= \psi_{m-1}(\varphi^{m-1}(h)\bar a).$$
\end{proof}

\section{Умножение в алгебрах, простых по отношению к действию расширений Оре}\label{SectionHSimpleOreExtMul}

В этом параграфе мы покажем, что при некоторых дополнительных предположениях относительно $H$ умножение
в $H$-простой алгебре $A$ определяется умножением в $\tilde H$-простой алгебре
$A/J^{\tilde H}(A)$. Примерами алгебр Хопфа $H$, для которых справедливы результаты
параграфа по-прежнему служат алгебра Тафта $H_{m^2}(\zeta)$ и алгебры Хопфа из~\cite[\S 5.6]{Danara}.
В частности, удаётся получить полную классификацию конечномерных
неполупростых $H_{m^2}(\zeta)$-простых алгебр (теоремы~\ref{TheoremTaftSimpleNonSemiSimplePresent} и~\ref{TheoremTaftSimpleNonSemiSimpleClassify}).

В теореме~\ref{TheoremHSimpleOreExtMul}, которая доказывается ниже,
используются \textit{квантовые биномиальные коэффициенты}:\label{DefQuantumBinom}
 $$\binom{n}{k}_\zeta := \frac{n!_\zeta}{(n-k)!_\zeta\ k!_\zeta},$$ где $n!_\zeta := n_\zeta (n-1)_\zeta \cdot \ldots \cdot 1_\zeta$ и
$n_\zeta := 1 + \zeta + \zeta^2 + \ldots + \zeta^{n-1}$, $n\in\mathbb N$, $0_\zeta :=1$.
При этом даже при $n!_\zeta = 0$ мы считаем, что $\binom{n}{n}_\zeta:=\binom{n}{0}_\zeta:=1$.

Пусть $H$~"--- алгебра Хопфа над полем $\mathbbm{k}$,
порождённая как алгебра подалгеброй Хопфа $\tilde H$
и косопримитивным элементом $v\in H$, где $\Delta v = g \otimes v + v \otimes 1$, $g \in G(\tilde H)$.
Предположим также, что существует автоморфизм $\varphi \in \Aut_{\mathbf{Alg}}(\tilde H)$,
такой, что $vh-\varphi(h)v\in \tilde H$ для всех $h\in\tilde H$,
причём $vg=\varphi(g)v$ и $\varphi(g)=\zeta g$ для примитивного корня $\zeta$ степени $t$ из единицы,
где $v^t \tilde J \subseteq \tilde J$, например, $v^t \in \tilde H$.
Пусть $A$~"--- конечномерная ассоциативная $H$-простая $H$-модульная алгебра с $J^{\tilde H}(A)\ne 0$. 
Обозначим через $\tilde J \subseteq J^{\tilde H}(A)^{q-1}$
некоторый минимальный двухсторонний $\tilde H$-инвариантный идеал,
где натуральное число $q$ задано условиями $J^{\tilde H}(A)^q=0$ и $J^{\tilde H}(A)^{q-1}\ne 0$.
Введём обозначение $v^{-1} \tilde J:= 
\lbrace a\in A \mid va\in  \tilde J\rbrace$.

\begin{theorem}\label{TheoremHSimpleOreExtMul}
Число $t$ равно числу $m$ из теоремы~\ref{TheoremHSimpleOreExtStructure}, подпространство $B := v^{-1} \tilde J$ является $g$-инвариантной подалгеброй, которая совпадает с $v^{m-1} \tilde J$,
имеет место разложение $A=B\oplus J^{\tilde H}(A)$
(прямая сумма подпространств), отображение $\psi_0 \colon A / J^{\tilde H}(A) \mathrel{\widetilde\to}
B$ является изоморфизмом алгебр
и для любых $0\leqslant k,\ell < m$ и $\bar a, \bar b \in A / J^{\tilde H}(A)$
справедлива формула \begin{equation}\label{EqQuantumCoefMulHSimpleOreExt}
\psi_k(\bar a)\psi_\ell(\bar b)=\tbinom{k+\ell}{k}_\zeta \psi_{k+\ell}
((g^\ell \bar a) \bar b).
\end{equation} (Здесь $\psi_k$~"--- отображения, определённые в
конце \S\ref{SectionHSimpleOreExtStructure}.)
\end{theorem}

Для того, чтобы доказать теорему~\ref{TheoremHSimpleOreExtMul}, нам потребуется несколько лемм.

Прежде всего заметим, что
$\psi(ga)=\zeta^{-1}g\psi(a)$ для всех $a\in A$, где $\psi$~"---
отображение, определённое в предыдущем параграфе.

\begin{lemma}\label{LemmaHSimpleOreExtMul0}
Для всех $a \in A$ и $b \in v^{-1} \tilde J$ справедливы равенства $\psi(ab)=\psi(a)b$ и $\psi(ba)=(g^{-1}b)\psi(a)$.
\end{lemma}
\begin{proof} 
В случае $a\in \tilde J$ утверждение леммы очевидно, так как $\psi(\tilde J)=0$.
Поскольку в силу теоремы~\ref{TheoremHSimpleOreExtStructure}
справедливо равенство $A=vJ^{\tilde H}(A) + \tilde J$,
достаточно доказать равенства из формулировки
в случае, когда $a = vu$ для некоторого $u\in J^{\tilde H}(A)$.
Заметим, что по определению $\psi(vu)=u$ для всех $u\in J^{\tilde H}(A)$.
Более того,
$J^{\tilde H}(A)$ является $\tilde H$-инвариантным
двухсторонним идеалом и
$(gu)(vb)\in J^{\tilde H}(A) \tilde J = 0$.
Следовательно, $$\psi(ab)=\psi((vu)b)=\psi(v(ub))- \psi((gu)(vb)) = ub=\psi(vu)b=\psi(a)b.$$
Аналогично, \begin{equation*}\begin{split}\psi(ba)=\psi(b(vu))=\psi(v((g^{-1}b)u))- \psi((vg^{-1}b)u) = \\ = (g^{-1}b)u-\zeta^{-1}\psi((g^{-1}vb)u)=(g^{-1}b)\psi(vu)-0=(g^{-1}b)\psi(a).\end{split}\end{equation*}
\end{proof}

Теперь мы можем доказать формулу, в которою входят степени отображения $\psi$. 

\begin{lemma}\label{LemmaHSimpleOreExtMul}
Пусть $a,b \in v^{-1} \tilde J$, $0\leqslant k,\ell < t$.
Тогда \begin{equation}\label{EqQuantumCoefMulHSimpleOreExt0}
\psi^k(a)\psi^\ell(b)=\tbinom{k+\ell}{k}_\zeta \psi^{k+\ell}
((g^\ell a) b).
\end{equation}
\end{lemma}
\begin{proof}
Докажем лемму индукцией по $k+\ell$. Если $k=0$ и $\ell = 0$,
то (\ref{EqQuantumCoefMulHSimpleOreExt0})
является следствием леммы~\ref{LemmaHSimpleOreExtMul0}.
Предположим, что $k,\ell \geqslant 1$.
Тогда $\psi^k(a), \psi^\ell(b) \in J^{\tilde H}(A)$.
Для всех $u\in A$ выполнено $v\psi(u)-u \in \tilde J$.
Поскольку $ J^{\tilde H}(A)\tilde J = \tilde J J^{\tilde H}(A) = 0$,
выполнено $(g\psi^k(a))v\psi^\ell(b) = (g \psi^k(a))\psi^{\ell-1}(b)$
и $(v\psi^k(a))\psi^\ell(b) = \psi^{k-1}(a)\psi^\ell(b)$.
Кроме того, для всех $u \in J^{\tilde H}(A)$ справедливо равенство $\psi(vu)=u$.
Следовательно,
 \begin{equation*}\begin{split}\psi^k(a)\psi^\ell(b)
= \psi(v(\psi^k(a)\psi^\ell(b)))=\\=\psi((g\psi^k(a))v\psi^\ell(b)
+(v\psi^k(a))\psi^\ell(b))=\psi(\zeta^k(\psi^k(g a))\psi^{\ell-1}(b)
+ \psi^{k-1}(a)\psi^\ell(b))=\\=\psi\left(
 \zeta^k \tbinom{k+\ell-1}k_\zeta \psi^{k+\ell-1}
((g^\ell a) b)
+ \tbinom{k+\ell-1}{k-1}_\zeta \psi^{k+\ell-1}
((g^\ell a) b)\right)=\tbinom{k+\ell}{k}_\zeta \psi^{k+\ell}
((g^\ell a) b).\end{split}\end{equation*}
\end{proof}
\begin{corollary}Пусть $m$~"--- число из теоремы~\ref{TheoremHSimpleOreExtStructure}. Тогда $t=m$
и $\zeta$~"--- примитивный корень степени $m$ из единицы.
\end{corollary}
\begin{proof} Сперва заметим, что из
$v^t \tilde J \subseteq \tilde J$ следует, что $I_t$ является $H$-инвариантным
идеалом алгебры $A$ и $m \leqslant t$.
Кроме того, $m \geqslant 2$, так как $I_{m-1}=J^{\tilde H}(A) \ne 0$.
В силу леммы~\ref{LemmaHSemiSimpleIsUnital}
алгебра $A$ содержит единицу $1_A$,
причём $1_A \notin I_{m-1}$,
поскольку $I_{m-1}$ является нетривиальным идеалом.
 Отсюда $\psi^{m-1}(1_A) \in J^{\tilde H}(A) \backslash \lbrace 0 \rbrace$.
Так как $v\psi(a)-a \in \tilde J$ для всех $a\in A$,
получаем $v\psi^{m-1}(1_A)=\psi^{m-2}(1_A)+j_1$
и $v\psi(1_A)=1_A+j_2$ для некоторых $j_1, j_2 \in \tilde J$.
Заметим, что $\psi^{m-1}(1_A) \psi(1_A)=\binom{m}{m-1}_\zeta \psi^m(1_A) = 0$.
Однако \begin{equation*}\begin{split}0=v(\psi^{m-1} (1_A)\psi(1_A))
=(v\psi^{m-1}(1_A))\psi(1_A)+(g\psi^{m-1}(1_A))v\psi(1_A)
=\\ =
(\psi^{m-2}(1_A)+j_1)\psi(1_A)+\zeta^{m-1}\psi^{m-1}(1_A)(1_A+j_2)
=\psi^{m-2}(1_A)\psi(1_A)+\zeta^{m-1}\psi^{m-2}(1_A)1_A
=\\ = \left(\binom {m-1}{m-2}_\zeta+\zeta^{m-1}\right)\psi^{m-1}(1_A)
=m_\zeta\ \psi^{m-1}(1_A),\end{split}\end{equation*}
так как $ J^{\tilde H}(A)\tilde J = \tilde J J^{\tilde H}(A) = 0$.
Следовательно, $m_\zeta = 0$ и $t=m$.
\end{proof}
\begin{proof}[Доказательство теоремы~\ref{TheoremHSimpleOreExtMul}.] 
По условию $v^{m-1} \tilde J \subseteq v^{-1} \tilde J$.
Заметим, что $$vJ^{\tilde H}(A)=v I_{m-1}=\bigoplus_{i=1}^{m-1} v^i\tilde J$$
и $v^{-1} \tilde J \cap J^{\tilde H}(A)=0$.
Следовательно, $v^{-1} \tilde J = v^{m-1} \tilde J$.
Так как $\tilde J$ является $\tilde H$-инвариантным идеалом,
$$v(ab)=(ga)(vb)+(va)b \in \tilde J\text{ для всех }a,b \in v^{-1} \tilde J,$$
и
 $v^{-1} \tilde J$~"--- подалгебра. Теперь из $$A=v^{m-1} \tilde J \oplus J^{\tilde H}(A)\qquad\text{(прямая сумма подпространств)}$$
следует, что $A=B\oplus J^{\tilde H}(A)$, а из~(\ref{EqQuantumCoefMulHSimpleOreExt0}) следует~(\ref{EqQuantumCoefMulHSimpleOreExt}).
\end{proof}

Теперь применим полученные результаты в случае, когда
$H$ является алгеброй Тафта $H_{m^2}(\zeta)$, т.е. $\tilde H = \mathbbm{k}C_m$ и $vc=\zeta cv$.
При этом будем пользоваться соответствием между $\mathbbm{k}C_m$-модульными алгебрами
и $\mathbb Z/m\mathbb Z$-градуированными алгебрами, которое основано на изоморфизмах $(\mathbbm{k}C_m)^*\cong \mathbbm{k}C_m$ (см., например, \cite[\S 3.2]{ZaiGia})
и $(\mathbb Z/m\mathbb Z,+)\cong (C_m, \cdot)$. Для того, чтобы не перегружать обозначения,
будем обозначать компоненту $\mathbb Z/m\mathbb Z$-градуированной алгебры $A$, отвечающую элементу
$\bar k \in \mathbb Z/m\mathbb Z$, через $A^{(k)}$, где $0\leqslant k \leqslant m-1$, а не через $A^{(\bar k)}$.

Начнём с того, что приведём примеры ассоциативных $H_{m^2}(\zeta)$-простых алгебр, а затем докажем, что любая из конечномерных неполупростых ассоциативных $H_{m^2}(\zeta)$-простых алгебр изоморфна одной из $H_{m^2}(\zeta)$-модульных алгебр из теоремы ниже:

\begin{theorem}\label{TheoremTaftSimpleNonSemiSimplePresent} Пусть $B$~"--- некоторая $\mathbb Z/m\mathbb Z$-градуированно простая ассоциативная алгебра над полем $\mathbbm{k}$, содержащим
примитивный корень $\zeta$ степени $m$ из единицы.
Положим $W_0 := B$ и определим векторные пространства
$W_i$, где $1\leqslant i \leqslant m-1$, как копии векторного пространства~$B$.
Обозначим одной и той же буквой $\psi$ соответствующие линейные биекции  $W_{i-1} \to W_i$
и зададим на каждом $W_i$ градуировку группой $\mathbb Z/m\mathbb Z$
по формуле $W_i^{(\ell+1)} := \psi\left(W_{i-1}^{(\ell)}\right)$.
Положим $\psi(W_{m-1}):=0$.
Рассмотрим $H_{m^2}(\zeta)$-модуль $A=\bigoplus_{i=0}^{m-1} W_i$ (прямая сумма подпространств),
где $v\psi(a):=a$ для всех $a \in W_i$,
$0\leqslant i \leqslant m-2$, $vB:=0$,
а $c a^{(i)}:=\zeta^i a^{(i)}$ для всех  $a^{(i)} \in A^{(i)}$, $A^{(i)} := \bigoplus_{k=0}^{m-1} W_k^{(i)}$ (прямая сумма подпространств).
  Определим на $A$ умножение при помощи формулы $$\psi^k(a)\psi^\ell(b):=\binom{k+\ell}{k}_\zeta\ \psi^{k+\ell}((c^\ell a)b) \text{ для всех }a, b\in B \text{ и } 0 \leqslant k,\ell < m.$$ Тогда $A$~"---
  ассоциативная $H_{m^2}(\zeta)$-простая $H_{m^2}(\zeta)$-модульная алгебра.
\end{theorem}
\begin{proof} То, что формулы выше действительно
задают на $A$ структуру $H_{m^2}(\zeta)$-модульной алгебры, проверяется непосредственно.

Пусть $I$~"--- $H_{m^2}(\zeta)$-инвариантный двухсторонний идеал алгебры $A$.
 В силу того, что $v^m=0$, справедливо равенство $v^m I = 0$. Обозначим через $t \in \mathbb Z_+$ такое число, что $v^t I \ne 0$, однако $v^{t+1} I = 0$. Тогда $0 \ne v^t I \subseteq I \cap \ker v$. Однако $\ker v = B$ является градуированно простой алгеброй. Следовательно, $\ker v \subseteq I$. Поскольку $1_A=1_B \in I$,
 получаем $I = A$. Отсюда алгебра $A$ является $H_{m^2}(\zeta)$-простой.
\end{proof}

Теперь докажем, что мы действительно получили все неполупростые ассоциативные
 $H_{m^2}(\zeta)$-простые алгебры.

\begin{theorem}\label{TheoremTaftSimpleNonSemiSimpleClassify}
Пусть $A$~"--- конечномерная ассоциативная $H_{m^2}(\zeta)$-простая $H_{m^2}(\zeta)$-модульная
алгебра над полем $\mathbbm{k}$, содержащим примитивный корень $\zeta$ степени $m$ и $J(A)\ne 0$.
Тогда $A$ изоморфна алгебре из теоремы~\ref{TheoremTaftSimpleNonSemiSimplePresent}.
\end{theorem}
\begin{proof}
Из теоремы~\ref{TheoremHSimpleOreExtMul}
следует, что $v^{-1}\tilde J = v^{m-1} \tilde J$. Из равенства $v^m=0$, которое выполняется в алгебре Тафта,
получаем, что $$v^{m-1} \tilde J \subseteq \ker v  \subseteq v^{-1}\tilde J$$
и $$\ker v = v^{-1}\tilde J = v^{m-1} \tilde J \cong A/J^{\mathbbm{k}C_m}(A)$$~"--- градуированно простая алгебра. 
Теперь достаточно положить $W_i := v^{m-i+1} \tilde J$ и снова применить теорему~\ref{TheoremHSimpleOreExtMul}.
\end{proof}

\begin{remark}\label{RemarkTaftNonSSAut}
Поскольку максимальная полупростая подалгебра $\ker v$ определена однозначно,
любые две такие $H_{m^2}(\zeta)$-простые алгебры $A$
изоморфны как $H_{m^2}(\zeta)$-модульные алгебры,
если и только если их подалгебры $\ker v$ изоморфны как $\mathbb Z/m\mathbb Z$-градуированные алгебры.
Более того, все автоморфизмы алгебры $A$ как $H_{m^2}(\zeta)$-модульной
алгебры порождаются автоморфизмами алгебры $\ker v$ как $\mathbb Z/m\mathbb Z$-градуированной алгебры.
Действительно, пусть $\theta \colon A \to A$~"--- автоморфизм алгебры $A$ как $H_{m^2}(\zeta)$-модульной
алгебры. Тогда из $$\tilde J = \left(\psi(1_A)\right)^{m-1}(\ker v) = J(A)^{m-1}$$ 
следует, что $\theta(\tilde J) = \tilde J$ и $$v^{m-1}\theta(\psi^{m-1}(a))=\theta(a)
\text{ для всех }a\in \ker v,$$ откуда $$\theta(\psi^{m-1}(a))=\psi^{m-1}(\theta(a)).$$
Применяя $v^{m-k-1}$, получаем $\theta(\psi^k(a))=\psi^k(\theta(a))$
для всех $a\in\ker v$ и $0 \leqslant k < m$.
\end{remark}

     \section{Полупростые $H_{m^2}(\zeta)$-простые алгебры}
\label{SectionTaftSimpleSemisimple}

Для завершения классификации конечномерных ассоциативных $H_{m^2}(\zeta)$-простых алгебр
осталось рассмотреть случай, когда такие алгебры полупросты. Оказывается,
все полупростые $H_{m^2}(\zeta)$-простые алгебры являются $\mathbbm{k}C_m$-простыми, т.е. не имеют нетривиальных идеалов, инвариантных
относительно действия оператора $c$:

\begin{theorem}\label{TheoremTaftSimpleSemisimple}
Пусть $A$~"---  конечномерная полупростая ассоциативная $H_{m^2}(\zeta)$-простая $H_{m^2}(\zeta)$-модульная алгебра
над алгебраически замкнутым полем~$\mathbbm{k}$.
Тогда
 $$A \cong \underbrace{M_k(\mathbbm{k}) \oplus M_k(\mathbbm{k}) \oplus \ldots \oplus M_k(\mathbbm{k})}_t\qquad \text{(прямая сумма идеалов)}$$ для некоторых $k,t\in\mathbb N$, $t \mid m$, и
 существуют такие матрицы $P \in M_k(\mathbbm{k})$ и $Q \in \GL_k(\mathbbm{k})$,  где  $Q^{\frac{m}{t}} = E_k$, 
 $E_k$~"--- единичная матрица размера $k\times k$, $Q P Q^{-1}=\zeta^{-t} P$, $P^m = \alpha E_k$ для некоторого $\alpha \in \mathbbm{k}$, что
\begin{equation}\label{EqSSTaftSimple1} c\, (a_1, a_2, \ldots, a_t) = (Q a_t Q^{-1}, a_1, \ldots, a_{t-1}),
\end{equation}\begin{equation}\label{EqSSTaftSimple2} v\,(a_1, a_2, \ldots, a_t)=(Pa_1 - (Q a_t Q^{-1}) P, \zeta (P a_2 - a_1 P), \ldots, \zeta^{t-1}(Pa_t - a_{t-1} P))\end{equation}
для всех $a_1, a_2, \ldots, a_t \in M_k(\mathbbm{k})$.
\end{theorem}

\begin{remark}\label{RemarkTaftSimpleSemisimpleDiagonalizeQ}
Диагонализируя $Q$, можно считать, что $$Q = \diag\Bigl\lbrace\underbrace{1,\ldots,1}_{k_1},
\underbrace{\zeta^t,\ldots,\zeta^t}_{k_2}, \ldots, \underbrace{\zeta^{t\left(\frac{m}{t}-1\right)},\ldots,\zeta^{t\left(\frac{m}{t}-1\right)}}_{k_{\frac{m}{t}}}\Bigr\rbrace$$ 
для некоторых $k_1,\ldots,k_{\frac{m}{t}}\in \mathbb Z_+$, $k_1+\ldots+k_{\frac{m}{t}}=k$.
Теперь из $QPQ^{-1} = \zeta^{-t} P$ следует, что $P=(P_{ij})$
является блочной матрицей, где $P_{ij}$~"--- матрица размера $k_{i-1} \times k_{j-1}$,
причём $P_{ij} = 0$ для всех $j \ne i+1$ и $(i,j)\ne \left(\frac{m}{t}, 1\right)$.
\end{remark}

Прежде, чем доказывать теорему~\ref{TheoremTaftSimpleSemisimple},
докажем три вспомогательные леммы.
В первых двух мы докажем все утверждения теоремы~\ref{TheoremTaftSimpleSemisimple},
за исключением равенства $P^m=\alpha E_k$. В лемме~\ref{LemmaTaftSimpleMatrix}
рассматривается случай $t=1$, т.е. когда алгебра $A$ изоморфна алгебре всех квадратных матриц.

\begin{lemma}\label{LemmaTaftSimpleMatrix}
Пусть $A$~"--- $H_{m^2}(\zeta)$-модульная
алгебра над алгебраически замкнутым полем $\mathbbm{k}$, изоморфная как алгебра алгебре $M_k(\mathbbm{k})$ для некоторого $k\in \mathbb N$. Тогда существуют такие матрицы $P\in M_k(\mathbbm{k})$, $Q\in \GL_k(\mathbbm{k})$, $Q^m=E_k$, что $Q P Q^{-1}=\zeta^{-1} P$ и
$A$ изоморфна как $H_{m^2}(\zeta)$-модульная алгебра алгебре $M_k(\mathbbm{k})$, наделённой
следующим
$H_{m^2}(\zeta)$-действием: $ca=QaQ^{-1}$ и $va=Pa-(QaQ^{-1})P$ для всех $a\in M_k(\mathbbm{k})$.
\end{lemma}
\begin{proof}
Все автоморфизмы полных матричных алгебр внутренние (см., например, \cite[\S 4.3]{Herstein}).
Отсюда $ca=QaQ^{-1}$ для некоторой матрицы $Q\in \GL_k(\mathbbm{k})$. Поскольку $c^m = 1$, 
матрица $Q^m$ скалярная. Умножая $Q$ на корень $m$-й степени
из соответствующего скаляра, мы можем считать, что $Q^m = E_k$.

 Напомним, что элемент $v$ действует на $A$ как косое дифференцирование.
Докажем\footnote{Данный результат является <<фольклорным>>. 
Автор благодарен В.\,К. Харченко, который сообщил ему простое доказательство данного результата.},
что это косое дифференцирование является \textit{внутренним},  т.е. существует такая матрица $P \in A$,
что $va = Pa-(ca)P$ для всех $a\in A$. Действительно,
для всех $a,b \in A$ справедливы равенства
 $$Q^{-1}(v(ab)) = Q^{-1}((ca)(vb)+(va)b)
= Q^{-1}((QaQ^{-1})(vb)+(va)b)=a(Q^{-1}(vb))+(Q^{-1}(va))b.$$
 Следовательно, $Q^{-1}(v(\cdot))$~"--- обычное дифференцирование. Поскольку все дифференцирования полных матричных алгебр внутренние (см., например, \cite[\S 4.3]{Herstein}),
существует такая матрица  $P_0 \in A$, что
 $Q^{-1}(va)=P_0a-aP_0$
 для всех $a\in A$.
 При этом $$va = QP_0a- Q a P_0= QP_0 a - Q aQ^{-1}QP_0= Pa-(QaQ^{-1})P
 \text{ для всех }a\in A,$$ где $P=QP_0$, т.е. $v$
действует как внутреннее косое дифференцирование.
 
Теперь заметим, что из $vc=\zeta cv$ следует, что $c^{-1}v=\zeta v c^{-1}$,
$$Q^{-1}(Pa-(QaQ^{-1})P)Q = \zeta \bigl(P(Q^{-1}aQ)-aP\bigr),$$ 
$$Q^{-1}PaQ-aQ^{-1}PQ = \zeta PQ^{-1}aQ-\zeta aP,$$ 
$$Q^{-1}Pa-aQ^{-1}P = \zeta PQ^{-1}a-\zeta aPQ^{-1},$$ 
$$Q^{-1}Pa-\zeta PQ^{-1}a = aQ^{-1}P-\zeta aPQ^{-1},$$ 
$$(Q^{-1}P-\zeta PQ^{-1})a = a(Q^{-1}P-\zeta PQ^{-1}) \text{ для всех } a\in A.$$
Следовательно, $Q^{-1}P-\zeta PQ^{-1} = \alpha E_k$ для некоторого $\alpha \in \mathbbm{k}$. Теперь
заменим $P$ на $(P-\frac{\alpha}{1-\zeta} Q)$.
Тогда равенство $va = Pa-(ca)P$ будет по-прежнему справедливо, но, кроме этого,
будут также справедливы равенства $Q^{-1}P-\zeta PQ^{-1} = 0$ и $Q PQ^{-1} = \zeta^{-1} P$.
\end{proof}

Рассмотрим теперь общий случай.

\begin{lemma}\label{LemmaTaftSimpleSemisimpleFirst}
Пусть $A$~"--- конечномерная полупростая ассоциативная $H_{m^2}(\zeta)$-простая $H_{m^2}(\zeta)$-модульная
алгебра над алгебраически замкнутым полем $\mathbbm{k}$.
Тогда
 $A \cong \underbrace{M_k(\mathbbm{k}) \oplus M_k(\mathbbm{k}) \oplus \ldots \oplus M_k(\mathbbm{k})}_t$ (прямая сумма идеалов) для некоторых $k, t \in\mathbb N$, где $t \mid m$,
 и существуют такие матрицы $P \in M_k(\mathbbm{k})$ и $Q \in \GL_k(\mathbbm{k})$, $Q^{\frac{m}{t}} = E_k$,
 где $Q P Q^{-1}=\zeta^{-t} P$, что для всех $a_1, a_2, \ldots, a_t \in M_k(\mathbbm{k})$
 справедливы равенства~(\ref{EqSSTaftSimple1}) и~(\ref{EqSSTaftSimple2}).
\end{lemma}
\begin{proof}
Поскольку алгебра $A$ полупроста, она является прямой суммой идеалов, каждый из которых, в свою очередь,
является $\mathbb Z/m\mathbb Z$-градуированно простой
подалгеброй. (Это следует, например, из теоремы~\ref{TheoremWedderburnHcomod}.)
Пусть $B$~"--- один из таких идеалов.
Тогда $vb= v(1_B b)=(c1_B)(vb)+(v1_B)b \in B$
для всех $b\in B$. Отсюда $B$~"--- $H_{m^2}(\zeta)$-подмодуль, $A=B$,
и $A$~"--- $\mathbb Z/m\mathbb Z$-градуированно простая или, что то же, $\mathbbm{k}C_m$-простая алгебра.
В силу классической теоремы Веддербёрна~"--- Артина алгебра $A$ является прямой суммой идеалов, каждый из которых изоморфен алгебре $M_k(\mathbbm{k})$ для некоторого $k$, при этом оператор $c$, действуя как автоморфизм, переводит всякий минимальный идеал $I$, где $I\cong M_k(\mathbbm{k})$, в другой идеал $cI$.
Отсюда $I+ cI + \ldots + c^{m-1} I$ является $\mathbbm{k}C_m$-инвариантным идеалом,
т.е., c учётом $\mathbbm{k}C_m$-простоты алгебры $A$, должно выполняться равенство $$A=I+ cI + \ldots + c^{m-1} I.$$ 
Выбирая такое минимальное $t\in\mathbb N$, что $c^t I = I$, получаем, что $$A \cong \underbrace{M_k(\mathbbm{k}) \oplus M_k(\mathbbm{k}) \oplus \ldots \oplus M_k(\mathbbm{k})}_t\text{ (прямая сумма
идеалов)},$$
причём элемент $c$ отображает $i$-ю компоненту в $(i+1)$-ю.
Из равенства $c^m I = I$ следует, что $t \mid m$.

Для случая $t=1$
утверждение данной леммы доказано в лемме~\ref{LemmaTaftSimpleMatrix}.
Рассмотрим случай $t\geqslant 2$.
Заметим, что $c^t$ отображает всякое прямое слагаемое в себя.
Поскольку всякий автоморфизм полной матричной алгебры внутренний,
существует такая матрица $Q$, что $c^t(a, 0, \ldots, 0) = (QaQ^{-1}, 0, \ldots, 0)$
для всех $a \in M_k(\mathbbm{k})$. Теперь из $c^m = \id_A$ следует, что $Q^{\frac{m}{t}}$~"---
скалярная матрица. Без ограничения общности можно считать, что $Q^{\frac{m}{t}} = E_k$,
поскольку поле $\mathbbm{k}$ алгебраически замкнуто и мы можем домножить
матрицу $Q$ на корень $m$-й степени из соответствующего скаляра.
Отсюда мы можем считать, что справедливо равенство~(\ref{EqSSTaftSimple1}).

Обозначим через $\pi_i \colon A \to M_k(\mathbbm{k})$
проекцию на $i$-ю компоненту. Рассмотрим отображения
 $\rho_{ij} \in \End_\mathbbm{k}(M_k(\mathbbm{k}))$, где $1 \leqslant i,j\leqslant t$,
 заданные формулами $\rho_{ij}(a):= \pi_i\bigl(v\,(\underbrace{0, \ldots, 0}_{j-1}, a,0, \ldots, 0)\bigr)$ для $a \in M_k(\mathbbm{k})$.
Тогда \begin{equation*}\begin{split}\rho_{ij}(ab)=\pi_i\bigl(v\,(0,\ldots, 0, ab,0, \ldots, 0)\bigr)=\\=
\pi_i\Bigl(v\bigl((0,\ldots, 0, a,0, \ldots, 0)(0,\ldots, 0, b,0, \ldots, 0)\bigr)\Bigr)=\\=
\pi_i\Bigl(\bigl(c(0,\ldots, 0, a,0, \ldots, 0)\bigr)v(0,\ldots, 0, b,0, \ldots, 0)\Bigr)+ \\+
\pi_i\Bigl(\bigl(v(0,\ldots, 0, a,0, \ldots, 0)\bigr)(0,\ldots, 0, b,0, \ldots, 0)\Bigr)= \\=
\delta_{ij}\,\rho_{ii}(a)b+\delta_{j, i-1}\,a\rho_{i,i-1}(b)+\delta_{i1}\delta_{jt}\,QaQ^{-1}\rho_{1t}(b)
\end{split}\end{equation*}
для всех $a,b \in M_k(\mathbbm{k})$, где $\delta_{ij}$~"--- символ Кронекера.

Пусть $\rho_{ii}(E_k)=P_i$, $\rho_{i, i-1}(E_k)=Q_i$, $\rho_{1t}(E_k)=Q_1$, где
$P_i, Q_i \in M_k(\mathbbm{k})$.
Тогда из равенства $a_i E_k = a_i $ для всех $a_i \in M_k(\mathbbm{k})$
следует, что
 \begin{equation*}\begin{split}v(a_1, \ldots, a_t)=\Bigl(\pi_1(v\,(a_1,\ldots, a_t)\bigr), \ldots, \pi_t\bigl(v\,(a_1,\ldots, a_t)\bigr)\Bigr)= \\ =
\bigl(P_1 a_1 + (Qa_tQ^{-1}) Q_1, P_2 a_2 + a_1 Q_2, \ldots, P_t a_t + a_{t-1} Q_t\bigr).
\end{split}\end{equation*}
Заметим теперь, что \begin{equation*}\begin{split}0=v((\underbrace{0, \ldots, 0}_{i-1}, E_k,0, \ldots, 0)(\underbrace{0, \ldots, 0}_i, E_k,0, \ldots, 0))=\\= (c(\underbrace{0, \ldots, 0}_{i-1}, E_k,0, \ldots, 0))(v(\underbrace{0, \ldots, 0}_i, E_k,0, \ldots, 0))+\\+ (v(\underbrace{0, \ldots, 0}_{i-1}, E_k,0, \ldots, 0))(\underbrace{0, \ldots, 0}_i, E_k,0, \ldots, 0)=\\= (\underbrace{0, \ldots, 0}_i, P_{i+1}+Q_{i+1}, 0, \ldots, 0).\end{split}\end{equation*}
Следовательно, $P_{i+1}+Q_{i+1}=0$, откуда $Q_j=-P_j$ для всех $1\leqslant j \leqslant t$.

Теперь заметим, что
\begin{equation*}\begin{split}(-\zeta QP_tQ^{-1}, 0,  \ldots, 0, \zeta P_{t-1})=\zeta cv(0, \ldots, 0, E_k, 0)=vc(0, \ldots, 0, E_k, 0)=\\= v(0, \ldots, 0, E_k)=(-P_1,0, \ldots, 0, P_t),\end{split}\end{equation*}
\begin{equation*}\begin{split}(\zeta QP_tQ^{-1}, -\zeta P_1,0, \ldots, 0)=\zeta cv(0, \ldots, 0, E_k)=vc(0, \ldots, 0, E_k)=\\= v(E_k,0, \ldots, 0)=(P_1, -P_2,0, \ldots, 0),\end{split}\end{equation*}
и $P_1 = \zeta Q P_t Q^{-1}$, $P_2 = \zeta P_1$, $P_t = \zeta P_{t-1}$.

Кроме того, если $t > 2$,  то \begin{equation*}\begin{split}(\underbrace{0, \ldots, 0}_i, \zeta P_i, -\zeta P_{i+1},0, \ldots, 0)=\zeta cv(\underbrace{0, \ldots, 0}_{i-1}, E_k,0, \ldots, 0)=\\= vc(\underbrace{0, \ldots, 0}_{i-1}, E_k,0, \ldots, 0)=v(\underbrace{0, \ldots, 0}_i, E_k,0, \ldots, 0)=(\underbrace{0, \ldots, 0}_i, P_{i+1}, -P_{i+2},0, \ldots, 0)\end{split}\end{equation*}
для всех $1\leqslant i \leqslant t-2$.
Следовательно, $P_{i+1}=\zeta P_i$ для всех $1\leqslant i \leqslant t-1$.
Пусть $P:= P_1$. Тогда $P_i = \zeta^{i-1} P$, $\zeta^t QP Q^{-1} = P$,  справедливо равенство (\ref{EqSSTaftSimple2}), и лемма доказана.
\end{proof}

Докажем теперь следующий результат:

\begin{lemma}\label{LemmaTaftSimpleSemisimpleFormula}
Пусть $v$~"--- оператор на алгебре $M_k(\mathbbm{k})^t$, заданный при помощи формулы~(\ref{EqSSTaftSimple2}),
где $QPQ^{-1}=\zeta^{-t} P$.
Тогда $$v^\ell (a_1, a_2, \ldots, a_t) =(b_1, b_2, \ldots, b_t),$$
где \begin{equation}\label{EqSSTaftSimple3} b_k = \zeta^{\ell(k-1)} \sum_{j=0}^\ell  (-1)^j \zeta^{-\frac{j(j-1)}{2}} \binom{\ell}{j}_{\zeta^{-1}}
P^{\ell-j} a_{k-j} P^j \end{equation}
и $a_{-j} := Q a_{t-j} Q^{-1}$, $j \geqslant 0$, $a_i \in M_k(\mathbbm{k})$, $1\leqslant \ell \leqslant m$.
\end{lemma}
\begin{proof}
Докажем утверждение индукцией по $\ell$. База $\ell=1$ очевидна в силу~(\ref{EqSSTaftSimple2}).
Предположим, что (\ref{EqSSTaftSimple3}) справедливо для $\ell$. Тогда
$$v^{\ell+1} (a_1, a_2, \ldots, a_t) =(\tilde b_1, \tilde b_2, \ldots, \tilde b_t),$$
где $\tilde b_k = \zeta^{k-1}(P b_k - b_{k-1} P)$, $1\leqslant k \leqslant t$ и $b_0 := Q b_t Q^{-1}$,
причём в силу условия $QPQ^{-1}=\zeta^{-t} P$ результат вычисления $b_k$
по формуле~(\ref{EqSSTaftSimple3}) при $k=0$ также совпадает с $b_0$.
Следовательно,
\begin{equation*}\begin{split}
\tilde b_k =  \zeta^{k-1}\left( \zeta^{\ell(k-1)} \sum_{j=0}^\ell (-1)^j \zeta^{-\frac{j(j-1)}{2}} \binom{\ell}{j}_{\zeta^{-1}}
P^{\ell-j+1} a_{k-j} P^j  -\right. \\ \left.- \zeta^{\ell(k-2)} \sum_{j=0}^\ell (-1)^j \zeta^{-\frac{j(j-1)}{2}} \binom{\ell}{j}_{\zeta^{-1}}
P^{\ell-j} a_{k-j-1} P^{j+1} \right)= \\= \zeta^{k-1}\left( \zeta^{\ell(k-1)} \sum_{j=0}^\ell (-1)^j \zeta^{-\frac{j(j-1)}{2}} \binom{\ell}{j}_{\zeta^{-1}}
P^{\ell-j+1} a_{k-j} P^j  -\right. \\ \left.- \zeta^{\ell(k-2)} \sum_{j=1}^{\ell+1} (-1)^{j-1} \zeta^{-\frac{(j-2)(j-1)}{2}} \binom{\ell}{j-1}_{\zeta^{-1}}
P^{\ell-j+1} a_{k-j} P^j \right)= \end{split}\end{equation*} \begin{equation*}\begin{split}
= \zeta^{(\ell+1)(k-1)}\left( \sum_{j=0}^\ell (-1)^j \zeta^{-\frac{j(j-1)}{2}} \binom{\ell}{j}_{\zeta^{-1}}
P^{\ell-j+1} a_{k-j} P^j  +\right. \\ \left.+  \sum_{j=1}^{\ell+1} (-1)^j \zeta^{-\frac{j(j-1)}{2}} \zeta^{j-\ell-1} \binom{\ell}{j-1}_{\zeta^{-1}}
P^{\ell-j+1} a_{k-j} P^j \right)=\\ =
\zeta^{(\ell+1)(k-1)} \sum_{j=0}^{\ell+1} (-1)^j \zeta^{-\frac{j(j-1)}{2}} \binom{\ell+1}{j}_{\zeta^{-1}}
P^{\ell-j+1} a_{k-j} P^j.\end{split}\end{equation*}  Отсюда (\ref{EqSSTaftSimple3}) справедливо для всех $1\leqslant \ell \leqslant m$.
\end{proof}
\begin{proof}[Доказательство теоремы~\ref{TheoremTaftSimpleSemisimple}.]
Напомним, что $v^m=0$ и $\binom{m}{j}_{\zeta^{-1}}=0$ для всех $1\leqslant j \leqslant m-1$.
Если $m$ нечётно, то $$(-1)^m \zeta^{-\frac{m(m-1)}{2}}=(-1)\cdot 1^{-\frac{m-1}2}=-1.$$
Если $m$ чётно, то $\zeta^{\frac{m}{2}}=-1$, т.е.  по-прежнему
$$(-1)^m \zeta^{-\frac{m(m-1)}{2}}=1 \cdot (-1)^{-(m-1)}=-1.$$
Следовательно, из лемм~\ref{LemmaTaftSimpleSemisimpleFirst}
и~\ref{LemmaTaftSimpleSemisimpleFormula} следует, что
\begin{equation*}\begin{split}v^m(a_1, \ldots, a_t)=(P^m a_1 - a_{1-m}P^m, P^m a_2 - a_{2-m}P^m, \ldots, P^m a_t - a_{t-m}P^m)=\\ =
([P^m, a_1], [P^m, a_2], \ldots, [P^m, a_t]) = 0\end{split}\end{equation*}
для всех $a_i \in M_k(\mathbbm{k})$, поскольку $Q^{\frac{m}{t}}=E_k$. Отсюда $P^m = \alpha E_k$ для некоторого $\alpha\in \mathbbm{k}$ и теорема доказана.
\end{proof}

\begin{remark}
Обратно, для всех $k, t \in\mathbb N$, где $t \mid m$, и таких матриц $P \in M_k(\mathbbm{k})$ и $Q \in \GL_k(\mathbbm{k})$,
что $Q^{\frac{m}{t}} = E_k$, $Q P Q^{-1}=\zeta^{-t} P$, $P^m = \alpha E_k$ для некоторого $\alpha \in \mathbbm{k}$,
можно определить структуру $H_{m^2}(\zeta)$-простой алгебры на алгебре 
$A = \underbrace{M_k(\mathbbm{k}) \oplus M_k(\mathbbm{k}) \oplus \ldots \oplus M_k(\mathbbm{k})}_t$ (прямая сумма идеалов)
при помощи формул~(\ref{EqSSTaftSimple1}) и~(\ref{EqSSTaftSimple2}).
Более того, алгебра $A$ будет $\mathbbm{k}C_m$-простой, а значит, и $\mathbb Z/m\mathbb Z$-градуированно простой.
\end{remark}

Исследуем теперь изоморфизмы между такими алгебрами. Начнём с того, что в
алгебре $\underbrace{M_k(\mathbbm{k}) \oplus M_k(\mathbbm{k}) \oplus \ldots \oplus M_k(\mathbbm{k})}_t$ (прямая сумма идеалов)
 существует ровно $t$
минимальных идеалов, причём каждый из них изоморфен алгебре $M_k(\mathbbm{k})$ (см., например, \cite[лемма~1.4.4]{Herstein}). Отсюда, если две какие-то $H_{m^2}(\zeta)$-простые алгебры изморфны, то числа $t$ у них
должны совпадать. Поэтому достаточно ограничиться изучением изоморфизмов между алгебрами
с одинаковыми $t$:

\begin{theorem}\label{TheoremTaftSimpleSSIso}
Пусть $A_1 \cong \underbrace{M_k(\mathbbm{k}) \oplus M_k(\mathbbm{k}) \oplus \ldots \oplus M_k(\mathbbm{k})}_t$ (прямая сумма идеалов)~"---
полупростая $H_{m^2}(\zeta)$-простая $H_{m^2}(\zeta)$-модульная алгебра над произвольным полем $\mathbbm{k}$,
заданная матрицами $P_1 \in M_k(\mathbbm{k})$ и $Q_1 \in \GL_k(\mathbbm{k})$
при помощи формул~(\ref{EqSSTaftSimple1}) и~(\ref{EqSSTaftSimple2}),
и пусть $A_2$~"--- другая такая алгебра, заданная матрицами $P_2 \in M_k(\mathbbm{k})$ и $Q_2 \in \GL_k(\mathbbm{k})$.
Тогда $A_1 \cong A_2$ как алгебры и $H_{m^2}(\zeta)$-модули,
если и только если $P_2 =\zeta^s\, T P_1 T^{-1}$ и $Q_2 = \beta T Q_1 T^{-1}$ для некоторых $s\in\mathbb Z$, $\beta \in \mathbbm{k}$ и $T \in \GL_k(\mathbbm{k})$.
\end{theorem}
\begin{proof} Рассмотрим произвольный изоморфизм 
$H_{m^2}(\zeta)$-модульных алгебр $\varphi \colon A_1 \to A_2$.
Как уже было отмечено перед теоремой, единственными минимальными идеалами обеих алгебр являются компоненты соответствующих прямых сумм. Отсюда каждая копия алгебры $M_k(\mathbbm{k})$ из алгебры $A_1$
отображается на одну из копий алгебры $M_k(\mathbbm{k})$ в $A_2$.
Поскольку каждый автоморфизм алгебры $M_k(\mathbbm{k})$ внутренний,
существует такая матрица $T \in \GL_k(\mathbbm{k})$ и число $1 \leqslant r \leqslant t$,
что $$\varphi(a, 0, \ldots, 0) = (\underbrace{0,\ldots, 0}_{r-1}, TaT^{-1}, \underbrace{0,\ldots, 0}_{t-r})
\text{ для всех }a\in M_k(\mathbbm{k}).$$
По условию $\varphi$ коммутирует с оператором $c$, который передвигает компоненты вправо. Отсюда для всех $1\leqslant j \leqslant t-r+1$
и $a\in M_k(\mathbbm{k})$ справедливо равенство
$$\varphi(\underbrace{0,\ldots, 0}_{j-1}, a, 0, \ldots, 0)
= \varphi\bigl(c^{j-1}(a, 0, \ldots, 0) \bigr)
=  c^{j-1}\varphi(a, 0, \ldots, 0)
 = \varphi(\underbrace{0,\ldots, 0}_{r+j-2}, TaT^{-1}, \underbrace{0,\ldots, 0}_{t-r-j+1}).$$
Для всех $t-r+2 \leqslant j \leqslant t$
это равенство выглядит следующим образом: $$\varphi(\underbrace{0,\ldots, 0}_{j-1}, a, 0, \ldots, 0)
= c^{j-1}\varphi(a, 0, \ldots, 0)
 = \varphi(\underbrace{0,\ldots, 0}_{r+j-t-2}, Q_2 TaT^{-1} Q_2^{-1}, \underbrace{0,\ldots, 0}_{2t-r-j+1}).$$
Следовательно,
\begin{equation}\begin{split}\label{EqSSTaftSimple4}\varphi(a_1, \ldots, a_t)=(Q_2 T a_{t-r+2} T^{-1} Q_2^{-1},
Q_2 T a_{t-r+3} T^{-1} Q_2^{-1}, \ldots,
Q_2 T a_t T^{-1} Q_2^{-1}, \\ T a_1 T^{-1}, T a_2 T^{-1}, T a_3 T^{-1},
\ldots, T a_{t-r+1} T^{-1})
\end{split}\end{equation}
для всех $a_i \in M_k(\mathbbm{k})$.
При этом \begin{equation*}\begin{split}c\varphi(a_1, \ldots, a_t)=(
Q_2 T a_{t-r+1} T^{-1}Q_2^{-1},
Q_2 T a_{t-r+2} T^{-1} Q_2^{-1},
 \ldots,
Q_2 T a_{t-1} T^{-1} Q_2^{-1}, \\ Q_2 T a_t T^{-1} Q_2^{-1}, T a_1 T^{-1}, T a_2 T^{-1},
\ldots, T a_{t-r} T^{-1}),
\end{split}\end{equation*}
а \begin{equation*}\begin{split}\varphi\bigl(c(a_1, \ldots, a_t)\bigr)=(Q_2 T a_{t-r+1} T^{-1} Q_2^{-1},
Q_2 T a_{t-r+2} T^{-1} Q_2^{-1}, \ldots,
Q_2 T a_{t-1} T^{-1} Q_2^{-1}, \\ T Q_1 a_t Q_1^{-1}  T^{-1}, T a_1  T^{-1}, T a_2 T^{-1},
\ldots, T a_{t-r} T^{-1}).
\end{split}\end{equation*}
Учитывая, что $c\varphi=\varphi c$, получаем, что
$Q_2 T a_t T^{-1} Q_2^{-1} =  T Q_1 a_t Q_1^{-1} T^{-1}$ для всех $a_t\in M_k(\mathbbm{k})$,
т.е. матрица $Q_1^{-1} T^{-1}Q_2 T $ коммутирует со всеми квадратными матрицами
и, следовательно, является скалярной.
Отсюда $Q_2 = \beta T Q_1 T^{-1}$ для некоторого $\beta \in \mathbbm{k}$. Отметим, что $\beta \ne 0$
в силу того, что $Q_1, Q_2, T \in \GL_k(\mathbbm{k})$.

Используя равенство $v\varphi = \varphi v$, получаем
\begin{equation}\begin{split}\label{EqSSTaftSimple5}
v\varphi(a_1, \ldots, a_t) = \bigl(P_2 Q_2 T a_{t-r+2} T^{-1} Q_2^{-1} - Q_2 T a_{t-r+1} T^{-1} Q_2^{-1} P_2, \\
\zeta(P_2 Q_2 T a_{t-r+3} T^{-1} Q_2^{-1} - Q_2 T a_{t-r+2} T^{-1} Q_2^{-1} P_2), \ldots,
\\ \zeta^{r-2}(P_2 Q_2 T a_t T^{-1} Q_2^{-1} - Q_2 T a_{t-1} T^{-1} Q_2^{-1} P_2), 
\zeta^{r-1}(P_2 T a_1 T^{-1} - Q_2 T a_t T^{-1} Q_2^{-1} P_2), \\ \zeta^{r}(P_2 T a_2 T^{-1} - T a_1 T^{-1} P_2), 
\ldots,   \zeta^{t-1}(P_2 T a_{t-r+1} T^{-1}-T a_{t-r} T^{-1} P_2)\bigr) = \\ =
\varphi\bigl(v(a_1, \ldots, a_t)\bigr) = \\ =
\varphi\Bigl(\bigl(P_1a_1 - (Q_1 a_t Q_1^{-1}) P_1, \zeta (P_1 a_2 - a_1 P_1), \ldots, \zeta^{t-1}(P_1a_t - a_{t-1} P_1) \bigr)\Bigr)= \\ = \bigl(
\zeta^{t-r+1} Q_2 T(P_1 a_{t-r+2} - a_{t-r+1} P_1)T^{-1} Q_2^{-1},  
\ldots, \zeta^{t-1} Q_2 T(P_1 a_t - a_{t-1} P_1)T^{-1} Q_2^{-1},\\
 T(P_1 a_1 - Q_1 a_t Q_1^{-1} P_1)T^{-1}, 
 \zeta T(P_1 a_2 - a_1 P_1)T^{-1},
 \ldots, \zeta^{t-r} T(P_1 a_{t-r+1} - a_{t-r} P_1)T^{-1}\bigr)
\end{split}\end{equation}
для всех $a_i \in M_k(\mathbbm{k})$.

Если $t \geqslant 2$, то, подставляя $a_1=E_k$ и
$a_j = 0$ при $j\ne 1$,
получаем $P_2 = \zeta^{1-r}\, T P_1 T^{-1}$.

В случае $t=1$ равенство \eqref{EqSSTaftSimple5}
принимает вид $$P_2 T a_1 T^{-1} - Q_2 T a_1 T^{-1} Q_2^{-1} P_2 = 
T(P_1 a_1 - Q_1 a_1 Q_1^{-1} P_1)T^{-1}.$$
Умножая слева на $T^{-1} Q_2^{-1}$,
справа на $T$ и используя
равенство $Q_2 = \beta T Q_1 T^{-1}$, получаем
$$T^{-1} Q_2^{-1} P_2 T a_1  - a_1 T^{-1} Q_2^{-1} P_2 T -
\beta^{-1} Q_1^{-1} P_1 a_1 + \beta^{-1} a_1 Q_1^{-1} P_1 = 0,$$
$$\beta^{-1}Q_1^{-1} T^{-1} P_2 T a_1  - \beta^{-1}a_1 Q_1^{-1} T^{-1} P_2 T -
\beta^{-1} Q_1^{-1} P_1 a_1 + \beta^{-1} a_1 Q_1^{-1} P_1 = 0,$$
$$(Q_1^{-1} T^{-1} P_2 T-  Q_1^{-1} P_1) a_1  - a_1 (Q_1^{-1} T^{-1} P_2 T -  Q_1^{-1} P_1) = 0.$$
Поскольку последнее равенство справедливо для всех $a_1 \in M_k(\mathbbm{k})$,
существует такое $\gamma \in \mathbbm{k}$,
что $$Q_1^{-1} T^{-1} P_2 T-  Q_1^{-1} P_1 = \gamma E_k,$$
т.е. \begin{equation}\label{EqSSTaftSimple6} P_2 = T P_1 T^{-1} + \gamma T Q_1 T^{-1} = T P_1 T^{-1} + \gamma\beta^{-1} Q_2.\end{equation}
С другой стороны, умножая \eqref{EqSSTaftSimple6} слева на $Q_2$, а справа на $Q_2^{-1}$,
получаем
$$ Q_2 P_2 Q_2^{-1}  = T Q_1 P_1 Q_1^{-1} T^{-1} + \gamma\beta^{-1} Q_2.$$
Учитывая равенства $Q_i P_i Q_i^{-1} = \zeta^{-1} P_i$, отсюда следует, что
$$ \zeta^{-1} P_2  = \zeta^{-1} T  P_1  T^{-1} + \gamma\beta^{-1} Q_2,$$
\begin{equation}\label{EqSSTaftSimple7} P_2  = T  P_1  T^{-1} + \zeta \gamma \beta^{-1} Q_2. \end{equation}
Сравнивая \eqref{EqSSTaftSimple6} и \eqref{EqSSTaftSimple7},
получаем, что $\gamma = 0$
и $P_2  = T  P_1  T^{-1}$. Поскольку в случае $t=1$
справедливо равенство $r=1$, это означает, что 
$P_2 = \zeta^{1-r}\, T P_1 T^{-1}$ и при $t=1$.
 Таким образом, необходимость условий из формулировки теоремы доказана.

Обратно, пусть $P_2 =\zeta^s\, T_0 P_1 T_0^{-1}$ и
$Q_2 = \beta T_0 Q_1 T_0^{-1}$
 для некоторых $s \in \mathbb Z$, $\beta \in \mathbbm{k}$
и $T_0 \in \GL_n(\mathbbm{k})$.
Представим $(1-s)$ в виде $n t + r$, где $n\in\mathbb Z$ и $1 \leqslant r \leqslant t$.
Положив $T:=T_0 Q_1^n$, получим $P_2 =\zeta^{1-r}\, T P_1 T^{-1}$
и $Q_2 = \beta T Q_1 T^{-1}$.
Теперь для завершения доказательства достаточно определить искомый изоморфизм $\varphi \colon A_1 \mathrel{\widetilde\to} A_2$ по формуле \eqref{EqSSTaftSimple4}.
\end{proof}

\begin{remark}\label{RemarkTaftSSAut}
Пусть полупростая ассоциативная $H_{m^2}(\zeta)$-простая $H_{m^2}(\zeta)$-модульная алгебра $A$ задаётся числом $t \in \mathbb N$ и матрицами $P \in M_k(\mathbbm{k})$ и $Q \in \GL_k(\mathbbm{k})$.
Из доказательства теоремы~\ref{TheoremTaftSimpleSSIso}
следует, что всякий автоморфизм $\varphi \colon A \mathrel{\widetilde\to} A$
задаётся формулой
\begin{equation*}\begin{split}\varphi(a_1, \ldots, a_t)=(Q T a_{t-r+2} T^{-1} Q^{-1},
Q T a_{t-r+3} T^{-1} Q^{-1}, \ldots,
Q T a_t T^{-1} Q^{-1}, \\ T a_1 T^{-1}, T a_2 T^{-1}, 
\ldots, T a_{t-r+1} T^{-1})
\end{split}\end{equation*}
для некоторых $1\leqslant r \leqslant t$ и $T\in\GL_k(\mathbbm{k})$.
Для того, чтобы описать умножение в группе $\Aut(A)$ таких автоморфизмов,
положим $$q:=\left\lbrace\begin{array}{ccc} 0 & \text{ при } & r = 1, \\
 t-r+1 & \text{ при } & r > 1 \end{array}\right.
  \text{ и } R :=\left\lbrace\begin{array}{ccc} T & \text{ при } & r=1, \\
Q T  & \text{ при } & r > 1. \end{array}\right.$$
Тогда
\begin{equation*}\begin{split}
\varphi(a_1, \ldots, a_t)=(R a_{q+1} R^{-1},
R a_{q+2} R^{-1}, \ldots,
R a_t R^{-1}, \\ Q^{-1} R a_1 R^{-1} Q, Q^{-1} R a_2 R^{-1} Q,
\ldots, Q^{-1} R a_q R^{-1} Q)
\end{split}\end{equation*}
для всех $a_i \in M_k(\mathbbm{k})$.  
Очевидно, что число $q$ определено однозначно, а матрица $R$~"--- с точностью до скалярного множителя.
Отсюда множество $\Aut(A)$ можно отождествить с множеством
пар $(\overline R, q)$, где $R\in \GL_k(\mathbbm{k})$, а $0\leqslant q < t$,
причём $P =\zeta^q\, R P R^{-1}$ и $QRQ^{-1}R^{-1} = \beta E_k$ для некоторого $\beta \in \mathbbm{k}$.
(Здесь при помощи $\bar R$ обозначен класс матрицы $R \in \GL_k(\mathbbm{k})$ в группе $\PGL_k(\mathbbm{k})$.)
Если перенести операцию умножения из группы автоморфизмов
на множество таких пар, получается, что
\begin{equation}\label{EqAutTaftSS}(\overline R_1,q_1)(\overline R_2,q_2)=\left\lbrace\begin{array}{llcc} \Bigl(\overline{R_1 R_2}, & q_1+q_2\Bigr) & \text{ при }& q_1+q_2 < t, \\
\Bigl(\overline{R_1 R_2 Q^{-1}}, & q_1+q_2-t\Bigr) & \text{ при }& q_1+q_2 \geqslant t. \end{array}\right.\end{equation}

Отсюда существует точная последовательность

$$\xymatrix{0 \ar[r] & G_0 \ar[r] & \Aut(A) \ar[r]^\rho & \mathbb Z/t\mathbb Z,}$$
где $\rho(\varphi):= \bar q$ (класс числа $q$ в группе $\mathbb Z/t\mathbb Z$),
а $$G_0 := \bigl\lbrace \bar R \in \PGL_n(\mathbbm{k}) \mathrel{\bigl|}
P = R P R^{-1},\ QRQ^{-1}R^{-1} = \beta E_k
\text{ для некоторого }\beta \in \mathbbm{k}
\bigr\rbrace.$$

При $P=0$ это описание становится проще. Группа $G_0$ является централизатором элемента $\bar Q$
в  $\PGL_k(\mathbbm{k})$, а гомоморфизм $\rho$ сюрьективен:
$$\xymatrix{0 \ar[r] & G_0 \ar[r] & \Aut(A) \ar[r]^\rho & \mathbb Z/t\mathbb Z \ar[r] & 0.}$$
Действительно, для всякого $0\leqslant q < t$ пара $(\bar E_k, q)$ задаёт автоморфизм
алгебры $A$, причём $\rho(\bar E_k, q) = \bar q \in \mathbb Z/t\mathbb Z$.

При $P\ne 0$ можно дать следующее описание группы $\Aut(A)$, заметив, что в этом случае в паре $(\bar R, q)$
число $q$ однозначно определено классом $\bar R \in \PGL_k(\mathbbm{k})$
при помощи равенства $P =\zeta^q\, R P R^{-1}$. Обозначим через $G$ подгруппу
в $\PGL_k(\mathbbm{k})$ состоящую из всех таких классов $\bar R$, что $P =\zeta^q\, R P R^{-1}$
и $QRQ^{-1}R^{-1} = \beta E_k$ для некоторых $q\in\mathbb Z$ и $\beta \in\mathbb Z$.
Тогда циклическая группа $\langle \bar Q \rangle$ является нормальной подгруппой группы $G$,
причём $\Aut(A) \cong G/\langle \bar Q \rangle$.
Действительно, поскольку в силу теоремы~\ref{TheoremTaftSimpleSemisimple}
справедливы равенства~$QPQ^{-1}=\zeta^{-t} P$ и $Q^{\frac{m}{t}}=E_k$, в каждом классе~"--- элементе группы $G/\langle \bar Q \rangle$~"--- существует ровно один представитель $\bar R$, для которого $0\leqslant q < t$, причём умножение таких представителей осуществляется по формуле~\eqref{EqAutTaftSS}.
\end{remark}

    \section{Алгебры, простые по отношению к действию алгебры Свидлера}
\label{SectionSweedlerSimple}

В данном параграфе мы рассмотрим случай $m=2$ подробнее,
получив соответствующие следствия из теорем~\ref{TheoremTaftSimpleNonSemiSimpleClassify}, \ref{TheoremTaftSimpleSemisimple} и~\ref{TheoremTaftSimpleSSIso} и замечаний~\ref{RemarkTaftNonSSAut} и~\ref{RemarkTaftSSAut}. Напомним, что в случае, когда мы рассматриваем алгебру Свидлера $H_4$,
мы всегда предполагаем, что $(-1)$ является примитивным корнем степени $2$ из единицы, т.е. характеристика
основного поля $\mathbbm{k}$ отлична от $2$.

Будем через $M_{k\times \ell}(\mathbbm{k})$ обозначать пространство всех матриц размера $k\times \ell$ над полем $\mathbbm{k}$.

\begin{theorem}\label{TheoremSweedlerSimpleMatrix}
Пусть $A$~"--- $H_4$-модульная алгебра над некоторым полем $\mathbbm{k}$, 
изоморфная как алгебра алгебре $M_k(\mathbbm{k})$ для некоторого $k\in \mathbb N$.
Тогда 
\begin{enumerate}
\item
либо $A$ изоморфна как алгебра и $H_4$-модуль
 алгебре
$$M_{s,\ell}(\mathbbm{k}):=M_{s,\ell}^{(0)}(\mathbbm{k})\oplus M_{s,\ell}^{(1)}(\mathbbm{k})\text{ (прямая сумма подпространств)},$$
где $c a = a$ для всех $a \in M_{s,\ell}^{(0)}(\mathbbm{k})$,
$c a = -a$ для всех $a \in M_{s,\ell}^{(1)}(\mathbbm{k})$,
$$M_{s,\ell}^{(0)}(\mathbbm{k}) := \left( \begin{array}{cc}
 M_s(\mathbbm{k}) & 0 \\
 0 & M_\ell(\mathbbm{k})\\
 \end{array}\right),$$
 $$M_{s,\ell}^{(1)}(\mathbbm{k}) := \left( \begin{array}{cc}
 0 & M_{s\times \ell}(\mathbbm{k}) \\
 M_{\ell\times s}(\mathbbm{k}) & 0\\
 \end{array}\right)$$ и $va=Pa-(ca)P$ для всех $a\in A$, где $P=\left( \begin{array}{cc}
 0 & P_1 \\
 P_2 & 0\\
 \end{array}\right)$, $P_1 P_2 = \alpha E_s$, $P_2 P_1 = \alpha E_\ell$
  для некоторых $s,\ell > 0$, $s+\ell = k$, $s\geqslant \ell$, $P_1 \in M_{s\times \ell}(\mathbbm{k})$,
  $P_2 \in M_{\ell\times s}(\mathbbm{k})$, $\alpha \in \mathbbm{k}$,
  \item либо $ca=a$, $va=0$ для всех $a\in A$.
  \end{enumerate}
\end{theorem}

\begin{remark}
Обратно, для любых таких $P$ и $s \geqslant \ell$,
существует $H_4$-модульная алгебра, изоморфная как алгебра алгебре $M_k(\mathbbm{k})$.
\end{remark}

\begin{remark}
Если $s > \ell$, то матрица $P_1 P_2$ вырождена и $P^2=0$.
\end{remark}
\begin{proof}[Доказательство теоремы~\ref{TheoremSweedlerSimpleMatrix}.]
Достаточно применить теорему~\ref{TheoremTaftSimpleSemisimple} и замечание~\ref{RemarkTaftSimpleSemisimpleDiagonalizeQ}. Условие $\ell > 0$
следует из того, что если $\ell = 0$, то $ca=a$ для всех $a\in A$ и
$va=vca=-cva=-va$, т.е. $va=0$.
\end{proof}

\begin{theorem} \label{TheoremSweedlerSimpleMatrixIso}
Пусть $A_1$ и $A_2$~"--- $H_4$-модульные алгебры над некоторым полем $\mathbbm{k}$,
$A_i=A_i^{(0)}\oplus A_i^{(1)}$ 
 (прямая сумма подпространств),
 где $c a = a$ для всех $a \in A_i^{(0)}$,
$c a = -a$ для всех $a \in A_i^{(1)}$,
$$A_i^{(0)} = \left( \begin{array}{cc}
 M_{s_i}(\mathbbm{k}) & 0 \\
 0 & M_{\ell_i}(\mathbbm{k})\\
 \end{array}\right),$$
 $$A_i^{(1)} = \left( \begin{array}{cc}
 0 & M_{s_i\times \ell_i}(\mathbbm{k}) \\
 M_{\ell_i\times s_i}(\mathbbm{k}) & 0\\
 \end{array}\right),$$
и $va=P_ia-(ca)P_i$ для всех $a\in A_i$, где $P_i=\left( \begin{array}{cc}
 0 & P_{i1} \\
 P_{i2} & 0\\
 \end{array}\right)$,  $P_{i1} P_{i2} = \alpha_i E_{s_i}$, $P_{i2} P_{i1} = \alpha_i E_{\ell_i}$
  для некоторых $s_i,\ell_i > 0$, $s_i\geqslant \ell_i$, $P_{i1} \in M_{s_i\times \ell_i}(\mathbbm{k})$,
  $P_{i2} \in M_{\ell_i\times s_i}(\mathbbm{k})$, $\alpha_i\in \mathbbm{k}$, $i=1,2$.
Тогда $A_1 \cong A_2$ как алгебры и $H_4$-модули, если и только если $s_1 = s_2$, $\ell_1 = \ell_2$
и существуют такие $T_1 \in \GL_{s_1}(\mathbbm{k})$, $T_2 \in \GL_{\ell_1}(\mathbbm{k})$,
что
\begin{enumerate}\item
 либо $P_{21}=T_1P_{11}{T_2}^{-1}$ и $P_{22}=T_2 P_{12}{T_1}^{-1}$,
\item 
либо $s_1=\ell_1$,
$P_{21}=T_1P_{12}{T_2}^{-1}$ и $P_{22}=T_2 P_{11}{T_1}^{-1}$.
\end{enumerate}
\end{theorem}
\begin{proof}
В силу теоремы~\ref{TheoremTaftSimpleSemisimple} алгебрам $A_i$
соответствуют матрицы $Q_i = \left( \begin{array}{cc}
 E_{s_i} & 0 \\
 0 & E_{\ell_i}\end{array}\right)$, причём в силу теоремы~\ref{TheoremTaftSimpleSSIso} изоморфизм
  $A_1 \mathrel{\widetilde\to} A_2$ существует, если и только если существует
 такая матрица $T\in\GL_k(\mathbbm{k})$, что $Q_2 = \beta TQ_1 T^{-1}$
 и $P_2 = \pm TP_1 T^{-1}$ для некоторого $\beta \in \mathbbm{k}$, причём, меняя $T$ на $TQ_1$, можно считать, что
 $P_2 = TP_1 T^{-1}$.
 
 Поскольку размерности собственных подпространств, отвечающих одному и тому же собственному
 значению, у линейных операторов с подобными матрицами совпадают и $s_i \geqslant \ell_i$,
 в случае существования изоморфизма $A_1 \mathrel{\widetilde\to} A_2$ выполнены условия $\beta = \pm 1$ и $s_1 = s_2$, $\ell_1 = \ell_2$,
 причём равенство $\beta = -1$ возможно лишь при $s_1 = \ell_1= s_2 = \ell_2$.
 
 Если $\beta = 1$, то
  $T=\left( \begin{array}{cc}
 T_1 & 0 \\
 0 & T_2 \end{array}\right)$,
 а если $\beta = -1$, то
 $T=\left( \begin{array}{cc}
 0 & T_1 \\
 T_2 & 0 \end{array}\right)$
  для некоторых $T_1 \in \GL_{s_1}(\mathbbm{k})$, $T_2 \in \GL_{\ell_1}(\mathbbm{k})$.
  Теперь утверждение теоремы следует из равенства $P_2 = TP_1 T^{-1}$. Обратное утверждение проверяется непосредственно.
\end{proof}
 
 \begin{example}
 В случае матриц размера $2\times 2$ возможны следующие варианты:
 \begin{enumerate}
 \item $A = A^{(0)}= M_2(\mathbbm{k})$, $A^{(1)}= 0$, $ca=a$, $va=0$ для всех $a\in A$;
 \item $A = A^{(0)} \oplus A^{(1)}$, где $$A^{(0)}=\left\lbrace\left(\begin{array}{cc}
\alpha & 0 \\
0      & \beta 
 \end{array}\right) \mathbin{\biggl|} \alpha,\beta \in \mathbbm{k}\right\rbrace$$
 и $$A^{(1)}=\left\lbrace\left(\begin{array}{cc}
0 & \alpha \\
\beta      & 0 
 \end{array}\right) \mathbin{\biggl|} \alpha,\beta \in \mathbbm{k}\right\rbrace,$$
 $ca=(-1)^{i}a$, $va=0$ для всех $a\in A^{(i)}$;
  \item $A = A^{(0)} \oplus A^{(1)}$, где $$A^{(0)}=\left\lbrace\left(\begin{array}{cc}
\alpha & 0 \\
0      & \beta 
 \end{array}\right) \mathbin{\biggl|} \alpha,\beta \in \mathbbm{k}\right\rbrace$$
 и $$A^{(1)}=\left\lbrace\left(\begin{array}{cc}
0 & \alpha \\
\beta      & 0 
 \end{array}\right) \mathbin{\biggl|} \alpha,\beta \in \mathbbm{k}\right\rbrace,$$
 $ca=(-1)^{i}a$, $va = Pa-(ca)P$ для всех $a\in A^{(i)}$, где $P=\left(\begin{array}{cc}
0 & 1 \\
\gamma      & 0 
 \end{array}\right)$ и $\gamma \in \mathbbm{k}$~"--- фиксированный элемент поля.
  \end{enumerate} 
 Действительно, случай $P=0$ относится к п.1. Если же $P\ne 0$, то, сопрягая, если нужно,
 матрицу $P$ матрицей $T$ вида $\left(\begin{array}{cc}
 \lambda & 0\\
 0 & \mu 
 \end{array}\right)$  или $\left(\begin{array}{cc}
0 & \lambda \\
\mu  & 0 
 \end{array}\right)$, где $\lambda,\mu \in \mathbbm{k}^\times$, можно считать, что элемент в правом верхнем углу равен $1$. 
В силу теоремы~\ref{TheoremSweedlerSimpleMatrixIso} различным элементам $\gamma$ соответствуют неизоморфные
алгебры $A$.
 \end{example}

Непосредственным следствием теоремы~\ref{TheoremTaftSimpleSemisimple} является теорема ниже:

\begin{theorem}\label{TheoremSweedlerSimpleSemisimple}
Пусть $A$~"--- конечномерная полупростая ассоциативная $H_4$-простая $H_4$-модульная алгебра над алгебраически замкнутым полем $\mathbbm{k}$.
Тогда либо
\begin{enumerate}
\item
 $A$ изоморфна алгебре $M_k(\mathbbm{k})$
для некоторого $k \geqslant 1$, либо
\item
 $A \cong M_k(\mathbbm{k}) \oplus M_k(\mathbbm{k})$ (прямая сумма идеалов) для некоторого $k \geqslant 1$
и существует такая матрица $P \in M_k(\mathbbm{k})$, что $P^2=\alpha E_k$ для некоторого $\alpha \in \mathbbm{k}$ и
\begin{equation}\label{EqSSSweedlerSimple} c\, (a, b) = (b,a),\qquad v\,(a,b)=(Pa-bP,aP-Pb)\end{equation}
для всех $a,b \in M_k(\mathbbm{k})$.
\end{enumerate}
\end{theorem}
\begin{remark}
Обратно, для любой матрицы $P \in M_k(\mathbbm{k})$, такой, что
$P^2=\alpha E_k$ для некоторого $\alpha \in \mathbbm{k}$,
на алгебре $M_k(\mathbbm{k}) \oplus M_k(\mathbbm{k})$ при помощи~(\ref{EqSSSweedlerSimple})
можно определить структуру $H_4$-простой алгебры,
которая, более того, является $\mathbb Z/2\mathbb Z$-градуированно простой.
\end{remark}
Из теоремы~\ref{TheoremTaftSimpleSSIso} получаем:
\begin{theorem}\label{TheoremSweedlerSimpleSSIso}
Пусть $A_1 = M_k(\mathbbm{k}) \oplus M_k(\mathbbm{k})$~"--- полупростая $H_4$-простая алгебра над полем $\mathbbm{k}$,
заданная матрицей $P_1 \in M_k(\mathbbm{k})$, где
$P_1^2=\alpha_1 E_k$ для некоторого $\alpha_1 \in \mathbbm{k}$, при помощи формул~(\ref{EqSSSweedlerSimple}),
а $A_2$~"--- другая такая алгебра, заданная матрицей $P_2 \in M_k(\mathbbm{k})$.
Тогда $A_1 \cong A_2$ как алгебры и $H_4$-модули,
если и только если $P_2 = \pm T P_1 T^{-1}$
для некоторого $T \in \GL_k(\mathbbm{k})$.
\end{theorem}

\begin{remark}
Из теорем~\ref{TheoremSweedlerSimpleSemisimple} и~\ref{TheoremSweedlerSimpleSSIso}
следует, что любая полупростая ассоциативная $H_4$-простая алгебра $A$
над алгебраически замкнутым полем $\mathbbm{k}$,
которая не является простой как обычная алгебра,
изоморфна
$M_k(\mathbbm{k}) \oplus M_k(\mathbbm{k})$ (прямая сумма идеалов) для некоторого $k \geqslant 1$,
где $$ c\, (a, b) = (b,a),\qquad v\,(a,b)=(Pa-bP,aP-Pb)$$
для всех $a,b \in M_k(\mathbbm{k})$
и
\begin{enumerate}
\item
либо $P=(\underbrace{\alpha,\alpha,\ldots, \alpha}_s, 
\underbrace{-\alpha,-\alpha,\ldots, -\alpha}_\ell)$ для некоторых $\alpha \in \mathbbm{k}$
и $s \geqslant \ell$, $s+\ell=k$,
\item либо $P$~"--- блочно диагональная матрица с несколькими блоками $\left(\begin{smallmatrix} 
 0 & 1 \\
 0 & 0\\
 \end{smallmatrix}\right)$
по главной диагонали (остальные клетки заполнены нулями),
 \end{enumerate}
 причём эти алгебры не изоморфны для различных $P$.
\end{remark}

Теоремы~\ref{TheoremSweedlerNonSemiSimple}, \ref{TheoremSweedlerNonSemiSimpleExistence}
и замечание~\ref{RemarkSweedlerNonSemiSimpleIso} являются
следствиями, соответственно, теорем~\ref{TheoremTaftSimpleNonSemiSimpleClassify},
\ref{TheoremTaftSimpleNonSemiSimplePresent} и замечания~\ref{RemarkTaftNonSSAut}:
\begin{theorem}\label{TheoremSweedlerNonSemiSimple}
Пусть $A$~"--- конечномерная ассоциативная $H_4$-модульная
алгебра над некоторым полем $\mathbbm{k}$.
Предположим, что $A$~"--- $H_4$-проста, но не полупроста.
Тогда существует такой $\mathbb Z/2\mathbb Z$-градуированный идеал $J \subset A$,
что $A = vJ \oplus J$ (прямая сумма $\mathbb Z/2\mathbb Z$-градуированных подпространств), $J^2 = 0$,
$vJ$~"--- $\mathbb Z/2\mathbb Z$-градуированно простая алгебра. Более того, существует такая линейная биекция $\psi \colon vJ \to J$,
что $v\psi(b)=b$, $a \psi(b) = \psi((ca)b)$,
$\psi(a)b = \psi(ab)$ для всех $a,b\in vJ$
 и $\psi((vJ)^{(0)})=J^{(1)}$, $\psi((vJ)^{(1)})=J^{(0)}$.
\end{theorem}

\begin{theorem}\label{TheoremSweedlerNonSemiSimpleExistence}
Пусть $B$~"--- $\mathbb Z/2\mathbb Z$-градуированно простая алгебра над некоторым полем $\mathbbm{k}$, $\chr \mathbbm{k} \ne 2$, а $\psi \colon B \to J$ линейная биекция из $B$
на некоторое векторное пространство $J$.
Зададим  на $J$ градуировку группой $\mathbb Z/2\mathbb Z$, положив $J^{(0)}:=\psi(B^{(1)})$
и $J^{(1)}:=\psi(B^{(0)})$, и определим на пространстве $A := B \oplus J$  (прямая сумма $\mathbb Z/2\mathbb Z$-градуированных подпространств)
действие оператора $c$ при помощи формул $c a := a$ для всех $a\in A^{(0)}:=B^{(0)} \oplus J^{(0)}$
и $c a := -a$ для всех $a\in A^{(1)}:=B^{(1)} \oplus J^{(1)}$
и действие оператора $v$ при помощи формулы $v(a+\psi(b))=b$ для всех $a,b\in B$.
Кроме того, определим на
$A$ операцию умножения, положив
 $a \psi(b) := \psi((ca)b)$,
$\psi(a)b := \psi(ab)$ для всех $a,b\in vJ$ и $J^2:=0$.
Тогда $A$ является $H_4$-простой ассоциативной $H_4$-модульной алгеброй.
\end{theorem}

\begin{remark}\label{RemarkSweedlerNonSemiSimpleIso}
Две такие алгебры $A$ изоморфны как $H_4$-модульные алгебры,
если и только если их максимальные $\mathbb Z/2\mathbb Z$-градуированно простые
подалгебры $B$ изоморфны как $\mathbb Z/2\mathbb Z$-градуированные алгебры.
\end{remark}

Рассмотрим теперь автоморфизмы конечномерных ассоциативных $H_4$-простых
алгебр. В первых пяти случаях из шести, которые приводятся ниже, группы автоморфизмов вычисляются явным образом.
В последнем случае предлагается метод описания группы автоморфизмов.

\begin{enumerate}
\item \label{ItemSweedlerSimpleMkl} Если $A=M_{s,\ell}(\mathbbm{k})$ для некоторых $s > \ell$,  $c a = a$ при $a \in M_{s,\ell}^{(0)}(\mathbbm{k})$,
$c a = -a$ при $a \in M_{s,\ell}^{(1)}(\mathbbm{k})$, $va=0$ для всех $a\in A$, тогда
$\Aut(A) \cong (\GL_s(\mathbbm{k})\times \GL_\ell(\mathbbm{k}))/\mathbbm{k}^{\times}E_{s+\ell}$,
где $\mathbbm{k}^{\times}E_n$~"--- группа невырожденных скалярных матриц размера $n\times n$
(это следует из замечания~\ref{RemarkTaftSSAut} и теоремы~\ref{TheoremSweedlerSimpleMatrixIso}).

\item \label{ItemSweedlerSimpleMkk} Если $A=M_{s,s}(\mathbbm{k})$ для некоторого $s\in \mathbb N$,  $c a = a$ для всех $a \in M_{s,s}^{(0)}(\mathbbm{k})$,
$c a = -a$ для всех $a \in M_{s,s}^{(1)}(\mathbbm{k})$, $va=0$ для всех $a\in A$, тогда
$$\Aut(A) \cong \left(\bigl(\GL_s(\mathbbm{k})\times \GL_s(\mathbbm{k})\bigr) \leftthreetimes \left\langle\left(\begin{smallmatrix}0 & E_s \\ E_s & 0 \end{smallmatrix}\right)\right\rangle_2 \right)/\mathbbm{k}^{\times}E_{2s}$$ 
(это также следует из замечания~\ref{RemarkTaftSSAut} и теоремы~\ref{TheoremSweedlerSimpleMatrixIso}).

\item \label{ItemSweedlerSimpleDoubleMk} Пусть $A \cong M_k(\mathbbm{k}) \oplus M_k(\mathbbm{k})$ (прямая сумма идеалов)~"---$H_4$-простая алгебра  для некоторого $k \geqslant 1$ и
$$ c\, (a, b) = (b,a),\qquad v\,(a,b)=(0,0)$$
для всех $a,b \in M_k(\mathbbm{k})$.
Тогда $\Aut(A) \cong \PGL_k(\mathbbm{k}) \times C_2$ 
(это следует из замечания~\ref{RemarkTaftSSAut}).

\item Если $A$~"--- конечномерная неполупростая $H_4$-простая алгебра, то, как было отмечено в замечании~\ref{RemarkTaftNonSSAut}, $\Aut(A)=\Aut(vJ)$, где
$J:=J(A)$, а $\Aut(vJ)$~"--- группа
автоморфизмов $vJ$ как $\mathbb Z/2\mathbb Z$-градуированной алгебры.
Группа $\Aut(vJ)$ описана в пп.~\ref{ItemSweedlerSimpleMkl}--\ref{ItemSweedlerSimpleDoubleMk}.

\item  \label{ItemSweedlerSimpleSSP}
Пусть $A \cong M_k(\mathbbm{k}) \oplus M_k(\mathbbm{k})$ (прямая сумма идеалов), где $k \geqslant 1$,~"--- полупростая $H_4$-простая алгебра, где 
$$ c\, (a, b) = (b,a),\qquad v\,(a,b)=(Pa-bP,aP-Pb)$$
для всех $a,b \in M_k(\mathbbm{k})$, а $P \in M_k(\mathbbm{k})$~"--- такая ненулевая матрица, что $P^2=\alpha E_k$ для некоторого $\alpha \in \mathbbm{k}$.  Тогда из замечания~\ref{RemarkTaftSSAut}
следует, что $\Aut(A)$ изоморфна подгруппе группы $\PGL_{k}(\mathbbm{k})$,
состоящей из образов всех таких невырожденных матриц $T$ размера $k\times k$,
что $T P T^{-1} = \pm P$. Без ограничения общности можно считать, что матрица $P$ приведена к жордановой
нормальной форме.
Заметим, что в силу теоремы о ранге матрицы пространство $M_{s,\ell}(\mathbbm{k})^{(1)}$
содержит невырожденные матрицы только при $s=\ell$.
  Следовательно,
\begin{enumerate}
\item[а)]  
  если $P=(\underbrace{\alpha,\alpha,\ldots, \alpha}_s, 
\underbrace{-\alpha,-\alpha,\ldots, -\alpha}_\ell)$ для некоторого $\alpha \in \mathbbm{k}^\times$
и $s \geqslant \ell$, $s+\ell=k$, то группа $\Aut(A)$ изоморфна
подгруппе группы $\PGL_{k}(\mathbbm{k})$,
которая состоит из образов всех обратимых матриц
из множества $M_{s,\ell}(\mathbbm{k})^{(0)}$ при $s > \ell$ и из множества $M_{s,s}(\mathbbm{k})^{(0)} \cup M_{s,s}(\mathbbm{k})^{(1)}$ при $s=\ell$.
\item[б)] Если $$P=\left(
\begin{array}{ccccccc}
\begin{array}{|cc|}
\hline
0 & 1 \\
0 & 0 \\
\hline
\end{array} & \begin{array}{cc}
0 & 0 \\
0 & 0 \\
\end{array} & \ldots & \begin{array}{cc}
0 & 0 \\
0 & 0 \\
\end{array} & \begin{array}{c}
0  \\
0  
\end{array}  & \ldots & \begin{array}{c}
0  \\
0  
\end{array}\\
       \begin{array}{cc}
0 & 0 \\
0 & 0 \\
\end{array}    & \begin{array}{|cc|}
                  \hline
                   0 & 1 \\
                    0 & 0 \\
                \hline
                \end{array} &  \ldots & \begin{array}{cc}
0 & 0 \\
0 & 0 \\
\end{array} &  \begin{array}{c}
0  \\
0  
\end{array} & \ldots & \begin{array}{c}
0  \\
0  
\end{array}\\
\multicolumn{7}{c}{\dotfill}\\
\begin{array}{cc}
0 & 0 \\
0 & 0 \\
\end{array} & \begin{array}{cc}
0 & 0 \\
0 & 0 \\
\end{array} & \ldots & \begin{array}{|cc|}
\hline
0 & 1 \\
0 & 0 \\
\hline
\end{array} & \begin{array}{c}
0  \\
0  
\end{array}  & \ldots & \begin{array}{c}
0  \\
0  
\end{array}\\
\begin{array}{cc}
0  & 0  
\end{array} & \begin{array}{cc}
0  & 0  
\end{array} & \ldots & \begin{array}{cc}
0  & 0  
\end{array} & 0 & \ldots & 0 \\ 
\multicolumn{7}{c}{\dotfill}\\
\begin{array}{cc}
0  & 0  
\end{array} & \begin{array}{cc}
0  & 0  
\end{array} & \ldots & \begin{array}{cc}
0  & 0  
\end{array} & 0 & \ldots & 0 \\ 
\end{array}
\right),$$
где число клеток $\left(\begin{smallmatrix} 
 0 & 1 \\
 0 & 0\\
 \end{smallmatrix}\right)$ равно $\ell\in\mathbb N$,
 то из~\cite[\S 69]{MalcevAIFoundations} следует,
 что группа $\Aut(A)$ изоморфна
 подгруппе $\PGL_{k}(\mathbbm{k})$, состоящей из образов всех обратимых матриц
 вида
$$\left(
\begin{array}{cc|cc|c|cc|ccc}
 \alpha_{11} & \beta_{11} & \alpha_{12} & \beta_{12} &   \ldots & \alpha_{1\ell} & \beta_{1\ell}
 & \gamma_{1,\ell+1} & \ldots &  \gamma_{1,{k-2\ell}}\\
 0           &  \alpha_{11} &       0   & \alpha_{12} &  \ldots &  0 & \alpha_{1\ell}
 & 0 & \ldots &  0 \\
\hline
 \alpha_{21} & \beta_{21} & \alpha_{22} & \beta_{22} &   \ldots & \alpha_{2\ell} & \beta_{2\ell}
 & \gamma_{2,\ell+1} & \ldots &  \gamma_{2,{k-2\ell}}\\
 0           &  \alpha_{21} &       0   & \alpha_{22} &  \ldots &  0 & \alpha_{2\ell}
 & 0 & \ldots &  0 \\
\hline
\multicolumn{2}{c|}{\dotfill} & \multicolumn{2}{c|}{\dotfill} & \ldots & \multicolumn{2}{c|}{\dotfill}
& \multicolumn{3}{c}{\dotfill} \\
\hline
 \alpha_{\ell1} & \beta_{\ell1} & \alpha_{\ell2} & \beta_{\ell2} &   \ldots & \alpha_{\ell\ell} & \beta_{\ell\ell}
 & \gamma_{\ell,\ell+1} & \ldots &  \gamma_{\ell,{k-2\ell}}\\
 0           &  \alpha_{\ell 1} &       0   & \alpha_{\ell 2} &  \ldots &  0 & \alpha_{\ell \ell }
 & 0 & \ldots &  0 \\
\hline
0 & \gamma_{\ell +1, 1} & 0 & \gamma_{\ell +1, 2} & \ldots & 0 &  \gamma_{\ell +1, \ell }
&  \gamma_{\ell +1, \ell +1} & \ldots & \gamma_{\ell +1, k-2\ell } \\
\multicolumn{2}{c|}{\dotfill} & \multicolumn{2}{c|}{\dotfill} & \ldots & \multicolumn{2}{c|}{\dotfill}
& \multicolumn{3}{c}{\dotfill} \\
0 & \gamma_{k-2\ell , 1} & 0 & \gamma_{k-2\ell , 2} & \ldots & 0 &  \gamma_{k-2\ell , \ell }
&  \gamma_{k-2\ell , \ell +1} & \ldots & \gamma_{k-2\ell , k-2\ell } \\
\end{array}
\right)$$
и 
$$\left(
\begin{array}{cc|cc|c|cc|ccc}
 \alpha_{11} & \beta_{11} & \alpha_{12} & \beta_{12} &   \ldots & \alpha_{1\ell } & \beta_{1\ell }
 & \gamma_{1,\ell +1} & \ldots &  \gamma_{1,{k-2\ell }}\\
 0           &  -\alpha_{11} &       0   & -\alpha_{12} &  \ldots &  0 & -\alpha_{1\ell }
 & 0 & \ldots &  0 \\
\hline
 \alpha_{21} & \beta_{21} & \alpha_{22} & \beta_{22} &   \ldots & \alpha_{2\ell } & \beta_{2\ell }
 & \gamma_{2,\ell +1} & \ldots &  \gamma_{2,{k-2\ell }}\\
 0           &  -\alpha_{21} &       0   & -\alpha_{22} &  \ldots &  0 & -\alpha_{2\ell }
 & 0 & \ldots &  0 \\
\hline
\multicolumn{2}{c|}{\dotfill} & \multicolumn{2}{c|}{\dotfill} & \ldots & \multicolumn{2}{c|}{\dotfill}
& \multicolumn{3}{c}{\dotfill} \\
\hline
 \alpha_{\ell 1} & \beta_{\ell 1} & \alpha_{\ell 2} & \beta_{\ell 2} &   \ldots & \alpha_{\ell \ell } & \beta_{\ell \ell }
 & \gamma_{\ell ,\ell +1} & \ldots &  \gamma_{\ell ,{k-2\ell }}\\
 0           &  -\alpha_{\ell 1} &       0   & -\alpha_{\ell 2} &  \ldots &  0 & -\alpha_{\ell \ell }
 & 0 & \ldots &  0 \\
\hline
0 & \gamma_{\ell +1, 1} & 0 & \gamma_{\ell +1, 2} & \ldots & 0 &  \gamma_{\ell +1, \ell }
&  \gamma_{\ell +1, \ell +1} & \ldots & \gamma_{\ell +1, k-2\ell } \\
\multicolumn{2}{c|}{\dotfill} & \multicolumn{2}{c|}{\dotfill} & \ldots & \multicolumn{2}{c|}{\dotfill}
& \multicolumn{3}{c}{\dotfill} \\
0 & \gamma_{k-2\ell , 1} & 0 & \gamma_{k-2\ell , 2} & \ldots & 0 &  \gamma_{k-2\ell , \ell }
&  \gamma_{k-2\ell , \ell +1} & \ldots & \gamma_{k-2\ell , k-2\ell } \\
\end{array}
\right).$$
\end{enumerate}

\item \label{ItemSweedlerSimpleMatrixQ} Пусть $A=M_{s,\ell}(\mathbbm{k})$~"---$H_4$-модульная алгебра для некоторых $s,\ell\in \mathbb N$, $s\geqslant \ell$,
где $c a = a$ для всех $a \in M_{s,\ell}^{(0)}(\mathbbm{k})$,
$c a = -a$ для всех $a \in M_{s,\ell}^{(1)}(\mathbbm{k})$,
и $va=Pa-(ca)P$ для всех $a\in A$, где $P=\left( \begin{array}{cc}
 0 & P_1 \\
 P_2 & 0\\
 \end{array}\right)$, $P_1 P_2 = \alpha E_s$, $P_2 P_1 = \alpha E_\ell$, $P_1 \in M_{s\times \ell}(\mathbbm{k})$,
  $P_2 \in M_{\ell\times s}(\mathbbm{k})$, $\alpha \in \mathbbm{k}$.
  Тогда  из замечания~\ref{RemarkTaftSSAut} и теоремы~\ref{TheoremSweedlerSimpleMatrixIso}
  следует, что 
\begin{enumerate}
\item  если $s\ne \ell$, то группа
  $\Aut(A)$ изоморфна подгруппе группы $\PGL_{s+\ell}(\mathbbm{k})$,
  состоящей из образов всех матриц $R=\left(\begin{array}{cc}
  R_1 & 0 \\  0 & R_2 \end{array} \right)$, где $R_1 \in \GL_s(\mathbbm{k})$, $R_2\in \GL_\ell(\mathbbm{k})$,
которые коммутируют с $P$;
  \item  
  если $s=\ell$, то
  группа $\Aut(A)$ изоморфна подгруппе группы $\PGL_{2k}(\mathbbm{k})$
  состоящей из образов всех матриц $R=\left(\begin{array}{cc}
  R_1 & 0 \\  0 & R_2 \end{array} \right)$ и $R=\left(\begin{array}{cc}
 0 & R_1 \\  R_2 & 0 \end{array} \right)$, где $R_1, R_2\in \GL_s(\mathbbm{k})$, которые коммутируют с $P$.
\end{enumerate}
  
  Как и в предыдущем случае, для любой конкретной матрицы $P$ матрицы $R$ определяются с использованием жордановой нормальной формы матрицы $P$
 (см., например, \cite[\S 69]{MalcevAIFoundations}).
\end{enumerate}

Рассмотрим некоторые частные случаи пп.~\ref{ItemSweedlerSimpleSSP}
и~\ref{ItemSweedlerSimpleMatrixQ}:

\begin{example}
Пусть $A=M_{1,1}(\mathbbm{k})$, причём $c a = a$ для всех $a \in M_{1,1}^{(0)}(\mathbbm{k})=\left\lbrace\left(
\begin{smallmatrix} 
\alpha & 0 \\
0 & \beta
\end{smallmatrix}
 \right)\right\rbrace$,
$c a = -a$ для всех $a \in M_{1,1}^{(1)}(\mathbbm{k})=\left\lbrace\left(
\begin{smallmatrix} 
0 & \alpha \\
\beta & 0
\end{smallmatrix}
 \right)\right\rbrace$
и $va=Pa-(ca)P$ для всех $a\in A$, где $P=\left(\begin{smallmatrix} 
 0 & 1 \\
 1 & 0\\
 \end{smallmatrix}\right)$.
  Заметим, что все невырожденные матрицы из $M_{1,1}^{(0)}(\mathbbm{k})$ и $M_{1,1}^{(1)}(\mathbbm{k})$,
  коммутирующие с $P$, имеют вид $\left(\begin{smallmatrix} 
 \alpha & 0 \\
 0 & \alpha\\
 \end{smallmatrix}\right)$ и $\left(\begin{smallmatrix} 
 0 & \alpha \\
 \alpha & 0\\
 \end{smallmatrix}\right)$, соответственно, где $\alpha \in \mathbbm{k}^\times$.
 Факторизуя по подгруппе скалярных матриц, получаем $\Aut(A)\cong C_2$.
\end{example}

\begin{example}
Пусть $A \cong M_2(\mathbbm{k}) \oplus M_2(\mathbbm{k})$ (прямая сумма идеалов) "---$H_4$-простая алгебра над полем $\mathbbm{k}$, где
 $$ c\, (a, b) = (b,a),\qquad v\,(a,b)=(Pa-bP,aP-Pb)$$
для всех $a,b \in M_k(\mathbbm{k})$, где $P=\left(\begin{smallmatrix} 
 0 & 1 \\
 0 & 0 \\
 \end{smallmatrix}\right)$, $\lambda \in \mathbbm{k}$.
 Все невырожденные матрицы $T$, коммутирующие с $P$, имеют вид $T=\left(\begin{smallmatrix} 
 \mu & \gamma \\
 0 & \mu \\
 \end{smallmatrix}\right)$, $\gamma \in \mathbbm{k}$, $\mu\in \mathbbm{k}^\times$.
  Все невырожденные матрицы $T$, такие, что $TPT^{-1}=-P$, имеют вид $T=\left(\begin{smallmatrix} 
 \mu & \gamma \\
 0 & -\mu \\
 \end{smallmatrix}\right)$, $\gamma \in \mathbbm{k}$, $\mu\in \mathbbm{k}^\times$.
 Факторизуя по подгруппе скалярных матриц, получаем $$\Aut(A)\cong 
 \left\lbrace \left(\begin{smallmatrix} 
 1 & \mu \\
 0 & 1 
 \end{smallmatrix} \right) \mathbin{\bigr|} \mu \in \mathbbm{k}  \right\rbrace
 \leftthreetimes \left\langle \left(\begin{smallmatrix} 
 1 & 0 \\
 0 & -1 
 \end{smallmatrix} \right) \right\rangle_2 \cong
(\mathbbm{k}, +) \leftthreetimes C_2,$$
где результат сопряжения элемента $\alpha \in (\mathbbm{k}, +)$ порождающим группы $C_2$ равен  $(-\alpha)$.
\end{example}

\begin{example}
Пусть $A \cong M_2(\mathbbm{k}) \oplus M_2(\mathbbm{k})$ (прямая сумма идеалов)~"---$H_4$-простая алгебра над полем $\mathbbm{k}$,
 $$ c\, (a, b) = (b,a),\qquad v\,(a,b)=(Pa-bP,aP-Pb)$$
для всех $a,b \in M_2(\mathbbm{k})$, где $P=\left(\begin{smallmatrix} 
 1 & 0 \\
 0 & -1 \\
 \end{smallmatrix}\right)$. Все невырожденные матрицы, коммутирующие с $P$, диагональны.
 Все невырожденные матрицы $T$, такие, что $TPT^{-1}=-P$, имеют вид $T=\left(\begin{smallmatrix} 
 0 & \gamma \\
 \mu & 0 \\
 \end{smallmatrix}\right)$, $\gamma, \mu\in \mathbbm{k}^\times$.
 Факторизуя по подгруппе скалярных матриц, получаем $$\Aut(A) \cong \left\lbrace \left(\begin{smallmatrix} 
 1 & 0 \\
 0 & \alpha 
 \end{smallmatrix} \right) \mathbin{\bigr|} \alpha \in \mathbbm{k}^\times  \right\rbrace
 \leftthreetimes \left\langle \left(\begin{smallmatrix} 
 0 & 1 \\
 1 & 0 
 \end{smallmatrix} \right) \right\rangle_2\cong
 \mathbbm{k}^\times \leftthreetimes C_2,$$ где результат сопряжения элемента $\alpha \in \mathbbm{k}^\times$
  порождающим группы $C_2$ равен $\alpha^{-1}$.
\end{example}

\newpage

\chapter{(Ко)модульные алгебры Ли}\label{ChapterH(co)modLie}

В данной главе рассматриваются структурные вопросы теории (ко)модульных алгебр Ли.
Доказанные утверждения будут затем использованы в главе~\ref{ChapterHModLieCodim} при изучении полиномиальных $H$-тождеств.

Результаты главы были опубликованы в работах~\cite{ASGordienko2, ASGordienko4, ASGordienko5, ASGordienko9, ASGordienko17}.

\section{(Ко)инвариантность радикалов}\label{SectionLieStability}

В данном параграфе доказываются достаточные условия (ко)инвариантности радикалов в (ко)модульных алгебрах Ли.

Для начала докажем следующее утверждение:

\begin{lemma}\label{LemmaNJadL} Пусть $L$~"--- алгебра Ли над некоторым полем $\mathbbm{k}$, а $N\subseteq L$~"---
нильпотентный идеал. Обозначим через $A$ ассоциативную подалгебру алгебры $\End_\mathbbm{k}(L)$,
порождённую подпространством $(\ad L)$. Тогда $(\ad N) \subseteq J(A)$.
(Напомним, что через $\ad \colon L \to \mathfrak{gl}(L)$ мы обозначаем присоединённое 
представление алгебры $L$ и $(\ad a)b := [a,b]$.)
\end{lemma}
\begin{proof}
Пусть $N^m=0$, $m \in \mathbb N$.
Тогда \begin{equation}\label{EqBp0} b_1 \ldots b_m = 0 \text{ для всех } b_i \in \ad N.
\end{equation}

Обозначим через $I$ двухсторонний идеал алгебры $A$, порождённый подпространством $(\ad N)$.
Тогда подпространство $I^n$, $n\in \mathbb N$, состоит из линейных комбинаций 
элементов
$$(a_{i_{01}} \ldots a_{i_{0s_0}}) b_1 (a_{i_{11}} \ldots a_{i_{1s_1}})
b_2 (a_{i_{21}} \ldots a_{i_{2s_2}})  \ldots b_{n-1}(a_{i_{n-1,1}} \ldots a_{i_{{n-1},s_{n-1}}}) b_n  (a_{i_{n1}} \ldots a_{i_{ns_n}}),$$
где $b_i \in \ad N$, $a_{ij} \in \ad L$.
Используя равенство $[\ad a, \ad b] = \ad [a, b]$, можно перенести все элементы $a_{ij}$
вправо и увидеть, что $I^n$ состоит из линейных комбинаций
элементов $a(b_1 b_2 \ldots b_n)c$, где $a,c \in A$, $b_i \in \ad N$.
Тогда из~(\ref{EqBp0}) следует, что $I^m = 0$ и $(\ad N) \subseteq I \subseteq J(A)$.
\end{proof}

\subsection{Модульные алгебры Ли}

Напомним, что алгебра Ли $L$ называется \textit{$H$-модульной} для некоторой алгебры Хопфа $H$,
если алгебра Ли $L$ является $H$-модулем и
 \begin{equation}\label{EqHmoduleLieAlgebra}
h\,[a,b]=[h_{(1)}a, h_{(2)}b]
\text{ для всех }h \in H,\ a,b \in L.
\end{equation}

 Докажем, что $H$-подмодуль, порождённый идеалом, также является идеалом:

\begin{lemma}\label{LemmaHISubModLie}
Пусть $I$~"--- идеал в $H$-модульной алгебре Ли $L$, где $H$~"--- алгебра Хопфа
над некоторым полем $\mathbbm{k}$.
Тогда $HI$ является $H$-инвариантным идеалом алгебры Ли $L$.
\end{lemma}
\begin{proof} Пусть $a \in I$, $h\in H$, $b \in L$. Тогда
\begin{equation}\label{EqPerebrosLie}
[ha, b] = [h_{(1)}a, \varepsilon(h_{(2)})b]
= [h_{(1)}a, h_{(2)} (Sh_{(3)}) b] = h_{(1)}[a, (Sh_{(2)})b].
\end{equation}
Отсюда $[ha, b]\in HI$.
\end{proof}

Теперь мы можем доказать достаточное условие инвариантности радикалов в алгебрах Ли:

\begin{theorem}\label{TheoremHModRadicalsLie}
Пусть $L$~"--- конечномерная $H$-модульная алгебра Ли над полем $\mathbbm{k}$
характеристики $0$, а $H$~"--- алгебра Хопфа с антиподом $S$, таким, что $S^2=\id_H$. 
Тогда разрешимый $R$ и нильпотентный $N$ радикалы алгебры Ли $L$ являются её $H$-подмодулями.
\end{theorem} 
\begin{proof} Определим действие алгебры Хопфа $H$ на алгебре $\End_\mathbbm{k}(L)$
при помощи равенства
$(h\psi)(a) := h_{(1)} \psi((Sh_{(2)})a)$,
где $a \in L$, $h \in H$, $\psi \in \End_\mathbbm{k}(L)$. (См. пример~\ref{ExampleHModEnd}.)
Тогда отображение $\ad$ оказывается гомоморфизмом $H$-модулей.
Отсюда ассоциативная подалгебра $A$ алгебры $\End_\mathbbm{k}(L)$, порождённая подпространством $(\ad L)$, является $H$-подмодулем.
В силу леммы~\ref{LemmaHISubModLie} $H$-подмодули $HN$ и $HR$ являются идеалами алгебры Ли $L$.
Из леммы~\ref{LemmaNJadL} следует, что $(\ad (HN)) \subseteq HJ(A)$.
 В то же время, согласно теореме~\ref{TheoremRadicalHSubMod} справедливо равенство $HJ(A)=J(A)$. Отсюда идеал $HN$ нильпотентен и $HN = N$.

В силу предложения 2.1.7 из \cite{GotoGrosshans} справедливо включение
 $[L, R] \subseteq N$. Применяя~(\ref{EqPerebrosLie}),
 получаем $$[HR, HR] \subseteq [HR, L] \subseteq H [R, HL] \subseteq H[R,L] \subseteq HN = N.$$
 Следовательно, идеал $HR$ разрешим и $HR = R$.
\end{proof}
\begin{corollary}\label{CorollaryHModRadicalsLie}
Пусть $L$~"--- конечномерная $H$-модульная алгебра Ли над полем $\mathbbm{k}$
характеристики $0$, а $H$~"--- конечномерная
(ко)полупростая алгебра Хопфа. 
Тогда разрешимый $R$ и нильпотентный $N$ радикалы алгебры Ли $L$ являются $H$-подмодулями.
\end{corollary}
\begin{proof}
В силу теоремы Ларсона~"--- Рэдфорда всякая конечномерная алгебра Хопфа $H$ над полем характеристики $0$ полупроста, если и только если она кополупроста, что в свою очередь справедливо,
если и только если $S^2=\id_H$ (см., например, \cite[теорема 7.4.6]{Danara}). Отсюда для доказательства утверждения достаточно применить теорему~\ref{TheoremHModRadicalsLie}.
\end{proof}

Также получаем новое доказательство известного результата (см., например, \cite[глава~III, \S6, теорема~7]{JacobsonLie}) о инвариантности радикалов относительно дифференцирований:

\begin{corollary}\label{CorollaryLieRadicalsDiffInvariant}
Пусть $L$ и $\mathfrak g$~"--- алгебры Ли над полем характеристики $0$,
причём~$\mathfrak g$ действует на~$L$ дифференцированиями и $\dim L < + \infty$.
Тогда $R$ и $N$ являются
$\mathfrak g$-подмодулями.
\end{corollary}
\begin{proof}
Алгебра Ли $L$ является $U(\mathfrak g)$-модульной алгеброй Ли (см. пример~\ref{ExampleUgModule}),
причём $S^2=\id_{U(\mathfrak g)}$. Отсюда в силу теоремы~\ref{TheoremHModRadicalsLie}
идеалы $R$ и $N$ являются $U(\mathfrak g)$-подмодулями, а следовательно, и $\mathfrak g$-подмодулями.
\end{proof}

Используя следствие~\ref{CorollaryLieRadicalsDiffInvariant} и разложение
конечномерных полупростых алгебр Ли в прямую сумму простых алгебр Ли (см., например, \cite[теорема 2.1.4]{GotoGrosshans}),
получаем по аналогии с теоремами~\ref{TheoremDerSimpleAssoc} и~\ref{TheoremGSimpleAssoc}
 следующие утверждения:
  \begin{theorem}\label{TheoremDerSimpleLie}
  Если $B$~"--- конечномерная $\mathfrak g$-простая алгебра Ли над полем $\mathbbm{k}$ характеристики $0$,
  на которой действует дифференцированиями некоторая алгебра Ли~$\mathfrak g$, то алгебра Ли $B$ проста.
  \end{theorem}
  \begin{theorem}\label{TheoremGSimpleLie}
Пусть $L$~"--- конечномерная $G$-простая алгебра Ли над алгебраически замкнутым полем $\mathbbm{k}$ характеристики $0$, на которой рационально действует автоморфизмами некоторая
связная аффинная алгебраическая группа $G$. Тогда алгебра Ли $L$ проста в обычном смысле.
  \end{theorem}

Приведём теперь пример $H$-модульной алгебры Ли,
радикалы которой не являются её $H$-подмодулями:

\begin{example}\label{ExampleLieSweedlerNonStableRadical}
Пусть $H_4$~"--- алгебра Свидлера (см. пример~\ref{ExampleTaftAlgebra}) над полем $\mathbbm{k}$.
Заметим, что $S^2 v = -S(c^{-1}v)=(c^{-1}v)c=-v \ne v$.
Пусть $W$~"--- некоторое трёхмерное векторное пространство.
Фиксируем линейную биекцию $\varphi \colon \mathfrak{sl}_2(\mathbbm{k}) \to W$.
Рассмотрим алгебру Ли $L = \mathfrak{sl}_2(\mathbbm{k})\oplus W$ (прямая сумма подпространств),
в которой коммутатор задаётся формулой $$[a+\varphi(b), u + \varphi(w)]=[a,u]+\varphi([a,w]+[b,u])
\text{ для всех } a,b,u,w \in \mathfrak{sl}_2(\mathbbm{k}).$$
Таким образом, $W$~"--- абелев идеал алгебры Ли $L$,
который совпадает с нильпотентным и разрешимым идеалами алгебры Ли $L$.
Определим $H_4$-действие через $c(a+\varphi(b))=a-\varphi(b)$ и $v(a+\varphi(b))=b$
для всех $a,b \in \mathfrak{sl}_2(\mathbbm{k})$.
Тогда $L$~"--- $H_4$-модульная алгебра Ли, однако её радикал $W$ не является $H_4$-подмодулем.
\end{example}

\subsection{Комодульные алгебры Ли}
Ниже мы получим достаточные условия коинвариантности радикалов и для комодульных алгебр Ли.

Напомним, что \textit{$H$-комодульной алгеброй Ли} для алгебры Хопфа $H$ 
называется алгебра Ли $L$, которая является $H$-комодулем,
причём $$\rho([a,b])=[a_{(0)},b_{(0)}] \otimes a_{(1)}b_{(1)}
\text{ для всех } a,b \in L.$$ Здесь $\rho \colon L \to L \otimes H$~"--- отображение,
задающее на $L$ структуру $H$-комодуля, $\rho(a)= a_{(0)}\otimes a_{(1)}$.

\begin{theorem}
\label{TheoremLieRadicalHSubComod}
Пусть $L$~"--- конечномерная $H$-комодульная алгебра Ли над полем $\mathbbm{k}$
характеристики $0$, а $H$~"--- алгебра Хопфа с антиподом $S$, таким, что $S^2=\id_H$. 
Тогда разрешимый $R$ и нильпотентный $N$ радикалы алгебры Ли $L$ являются её $H$-подкомодулями.
\end{theorem}
\begin{proof}
Рассмотрим присоединённое представление $\ad \colon L \to \mathfrak{gl}(L)$,
 $(\ad a)b:=[a,b]$.
Введём на ассоциативной алгебре $\End_\mathbbm{k}(L)$ структуру $H$-комодульной алгебры,
основываясь на примере~\ref{ExampleHComodEnd}.
Тогда отображение $\ad$ является гомоморфизмом $H$-комодулей и $H^*$-модулей.

Обозначим через $A$ ассоциативную подалгебру алгебры $\End_\mathbbm{k}(L)$, порождённую подпространством $(\ad L)$.
Поскольку $(\ad L)$ является $H$-подкомодулем, $A$ также является $H$-подкомодулем. 
В силу леммы~\ref{LemmaHComoduleIdeal} подпространства $H^*N$ и $H^*R$ являются идеалами алгебры Ли $L$.
Из теоремы~\ref{TheoremRadicalHSubComod} и леммы~\ref{LemmaNJadL}
следует, что $(\ad (H^*N)) \subseteq H^*J(A)=J(A)$. Следовательно, идеал $H^*N$ нильпотентен и $H^*N = N$.

В силу предложения 2.1.7 из \cite{GotoGrosshans} справедливо включение
 $[L, R] \subseteq N$. Применяя тот же приём, что и в~(\ref{EqMoveHFirst}),
 получаем $$[H^*R, H^*R] \subseteq [H^*R, L] \subseteq H^* [R, H^*L] \subseteq H^*[R,L] \subseteq H^*N = N.$$
 Отсюда  идеал $H^*R$ разрешим и $H^*R = R$.
\end{proof}

Из теоремы~\ref{TheoremLieRadicalHSubComod}
получается следующее обобщение результата
М.\,В.~Зайцева, Д.~Пагона и Д.~Реповша~\cite[предложение 3.3]{PaReZai}:
\begin{corollary}\label{CorollaryLieRadicalsGraded}
Пусть $L$~"--- конечномерная алгебра Ли над полем
характеристики $0$, градуированная произвольной группой. 
Тогда разрешимый $R$ и нильпотентный $N$ радикалы алгебры Ли $L$ являются градуированными идеалами.
\end{corollary}

\section{$(H,L)$-модули над $H$-модульными алгебрами Ли $L$}

Как будет видно из дальнейшего, при изучении $H$-(ко)модульных алгебр Ли $L$
важную роль играют $L$-модули, наделённые дополнительно структурой $H$-(ко)модуля,
в которых структуры $L$-модуля и $H$-(ко)модуля согласованы. Мы будем называть
такие модули $(H,L)$-модулями.

Пусть $L$~"--- $H$-модульная алгебра Ли, $V$~"--- $H$-модуль для некоторой алгебры Хопфа $H$,
а $\psi \colon L \to \mathfrak{gl}(V)$~"--- гомоморфизм, задающий на 
$V$ структуру $L$-модуля.
Будем говорить, что $(V, \psi)$ является
 \textit{$(H,L)$-модулем}, если \begin{equation*}h(\psi(a)v)=\psi(h_{(1)}a)(h_{(2)}v)
\text{ для всех }a\in L,\ h\in H,\ v\in V.\end{equation*}
Если для всех $a\in L$, $v\in V$ положить $av:= \psi(a)v$,
последнее равенство переписывается в виде
\begin{equation}\label{EqHLModule}h(av)=(h_{(1)}a)(h_{(2)}v)
\text{ для всех }a\in L,\ h\in H,\ v\in V.\end{equation}
   Будем называть $(H,L)$-модуль $(V, \psi)$ \textit{неприводимым},
   если он не содержит нетривиальных $H$-инвариантных $L$-подмодулей.
   Будем говорить, что $(H,L)$-модуль $(V, \psi)$ \textit{вполне приводим},
   если он является прямой суммой неприводимых $(H,L)$-подмодулей.
   Будем называть $(H,L)$-модуль $(V, \psi)$  \textit{точным}, если $(V, \psi)$ точен как $L$-модуль, т.е. если $\ker \psi = 0$.
Часто для краткости мы будем опускать вторую компоненту в обозначении
$(V, \psi)$ и писать просто $V$.

В силу~(\ref{EqHmoduleLieAlgebra}) для всякой $H$-модульной алгебры Ли $L$ пара $(L, \ad)$
является $(H,L)$-модулем, причём $(H,L)$-подмодули $(H,L)$-модуля $(L, \ad)$~"--- это
в точности $H$-инвариантные идеалы алгебры Ли $L$.
 
Если $V$~"--- модуль над некоторой алгеброй Хопфа $H$, то
алгебра Ли $\mathfrak{gl}(V)$ наследует структуру $H$-модуля от ассоциативной $H$-модульной
алгебры $\End_\mathbbm{k}(V)$ (см. пример~\ref{ExampleHModEnd}).
Однако если $H$ некокоммутативна, нельзя утверждать, что в $\mathfrak{gl}(V)$ выполяется условие~(\ref{EqHmoduleLieAlgebra}) и что $\mathfrak{gl}(V)$ является $H$-модульной алгеброй.

Пусть $(V,\psi)$ является $(H,L)$-модулем для некоторой $H$-модульной алгебры Ли $L$.
Тогда $\psi$~"--- гомоморфизм $H$-модулей.
   Более того, если через $\zeta \colon H \to \End_\mathbbm{k}(V)$ мы обозначим гомоморфизм, задающий на $V$ структуру $H$-модуля, то условие~(\ref{EqHLModule}) станет эквивалентым условию
\begin{equation}\label{EqHLModule2}\zeta(h)\psi(a) = \psi(h_{(1)}a) \zeta(h_{(2)})
\text{ для всех }
h\in H,\ a\in L.
\end{equation}

Докажем теперь некоторые свойства $(H,L)$-модулей,
которые будут затем использованы в \S\ref{SectionHSS}, \S\ref{SectionClassTaftSLieNSS} и главе~\ref{ChapterHModLieCodim}.

Как и в случае ассоциативных алгебр (см. предложение~\ref{PropositionHIHSubMod}), начнём с того, что докажем $H$-инвариантность аннуляторов $H$-подмодулей, в данном случае $H$-подмодулей $(H,L)$-модулей. Благодаря тождеству антикоммутативности в случае алгебр Ли этот результат справедлив без дополнительных предположений относительно антипода алгебры Хопфа.

\begin{lemma}\label{LemmaAnnHLmoduleHMod}
Пусть $V$~"--- $(H,L)$-модуль для некоторой алгебры Хопфа $H$ и $H$-модульной алгебры Ли $L$
над произвольным полем $\mathbbm{k}$, а $M \subseteq V$~"--- его $H$-подмодуль. Тогда аннулятор $$\Ann_L(M):=\lbrace a\in L \mid aM=0\rbrace$$ $H$-подмодуля $M$ является $H$-подмодулем алгебры Ли $L$.
Если, кроме этого, $M$ является ещё и $L$-подмодулем, то $\Ann_L(M)$ является $H$-инвариантным идеалом алгебры Ли $L$.
\end{lemma}
\begin{proof}
Заметим, что \begin{equation*}\begin{split}(ha)v=(h_{(1)}\varepsilon(h_{(2)})a)v
= (h_{(1)}a)((\varepsilon(h_{(2)})1)v)=(h_{(1)}a)(h_{(2)} (Sh_{(3)})v)=
\\ = (h_{(1)(1)}a)(h_{(1)(2)} (Sh_{(2)})v)
=h_{(1)}(a(Sh_{(2)})v)=0\end{split}\end{equation*}  для всех $v \in M$, $a \in \Ann_L(M)$,  $h\in H$,
поскольку для каждого слагаемого $(Sh_{(2)})v \in M$.
Отсюда $\Ann_L(M)$ действительно является $H$-подмодулем.
Если $M$~"--- $L$-подмодуль, то $$[a,b]v=a(bv)-b(av)=0\text{ для всех }v \in M,\ a \in \Ann_L(M),\ b \in L,$$
т.е. $\Ann_L(M)$ является идеалом.
\end{proof}

В случае, когда $V=L$, аннулятор $\Ann_L(M)$ подмножества $M\subseteq L$ называется \textit{централизатором} подмножества $M$.

В леммах~\ref{LemmaLR}--
\ref{LemmaRedIrr},
которые доказываются ниже,
 $V$ является конечномерным
$(H,L)$-модулем, где $H$~"--- алгебра Хопфа
над алгебраически замкнутым полем $\mathbbm{k}$ характеристики $0$, а $L$~"--- $H$-модульная
алгебра Ли с $H$-инвариантным радикалом $R$.
Через $\zeta \colon H \to \End_\mathbbm{k}(V)$ и
$\psi \colon L \to \mathfrak{gl}(V)$
будем обозначать гомоморфизмы, отвечающие $(H,L)$-модульной структуре.
Обозначим через $A$ ассоциативную подалгебру алгебры $\End_\mathbbm{k}(V)$,
порождённую операторами из $\psi(L)$
и $\zeta(H)$.

\begin{lemma}\label{LemmaLR} Справедливо включение
$\psi([L, R]) \subseteq J(A)$, где $J(A)$~"--- радикал Джекобсона алгебры $A$.
\end{lemma}
\begin{proof}
Пусть $$V = W_0 \supseteq W_1 \supseteq W_2 \supseteq \ldots \supseteq W_t = \left\lbrace 0 \right\rbrace$$
"--- композиционный ряд в $V$ из (необязательно $H$-инвариантных) $L$-подмодулей.
Тогда все факторы $W_i/W_{i+1}$ являются неприводимыми $L$-модулями.
Обозначим через $\psi_i \colon L \to
\mathfrak{gl}(W_i/W_{i+1})$ соответствующие гомоморфизмы.
Тогда в силу теоремы
 Э.~Картана~\cite[предложение~1.4.11]{GotoGrosshans},
каждая из алгебр Ли $\psi_i (L)$
либо полупроста, либо
является прямой суммой полупростого идеала и центра алгебры Ли $\mathfrak{gl}(W_i/W_{i+1})$.
Следовательно, алгебра Ли $\psi_i ([L, L])$ полупроста
и $\psi_i ([ L, L ] \cap R) = 0$.
Поскольку $[L, R] \subseteq [L, L] \cap R$,
получаем $\psi_i([L, R])=0$ и $[L,R]W_i \subseteq W_{i+1}$.

Обозначим через $I$ ассоциативный идеал алгебры $A$,
порождённый всеми элементами пространства $\psi([L, R])$.
Тогда $I^t$~"--- ассоциативный идеал, порождённый элементами вида
\begin{equation*}\begin{split}
a_1 \bigl(\zeta(h_{10})b_{11}\zeta(h_{11})b_{12} \ldots \zeta(h_{1,s_1-1})b_{1,s_1}
\zeta(h_{1,s_1})\bigr) a_2 \bigl(\zeta(h_{20})b_{21}\zeta(h_{21})b_{22} \ldots \zeta(h_{2,s_2-1})b_{2,s_2}
\zeta(h_{2,s_2})\bigr)\cdot 
\\ 
\cdot \ldots \cdot
a_{t-1}\bigl(\zeta(h_{t-1,0})b_{t-1,1}\zeta(h_{t-1,1})b_{t-1,2} \ldots
 \zeta(h_{t-1,s_{t-1}-1})b_{t-1,s_{t-1}}\zeta(h_{t-1,s_{t-1}})\bigr) a_t,\end{split}\end{equation*}
 где $a_i \in \psi([L, R])$, $b_{ij}\in\psi(L)$, $h_{ij}\in H$.
Используя~(\ref{EqHLModule2}), переместим все $\zeta(h_{ij})$
вправо и получим, что $I^t$ порождён элементами
  \begin{equation*}\begin{split}b =
a_1 \bigl((h'_{11}b_{11}) \ldots (h'_{1,s_1} b_{1,s_1})
\bigr)\ (h_2a_2)\left((h'_{21}b_{21}) \ldots (h'_{2,s_2} b_{2,s_2})\right)\cdot 
\\ 
\cdot \ldots \cdot
(h_{t-1}a_{t-1})
\bigl((h'_{t-1,1}b_{t-1,1}) \ldots
 (h'_{t-1,s_{t-1}}b_{t-1,s_{t-1}})\bigr) (h_t a_t) \zeta(h_{t+1}),
\end{split}\end{equation*} где $h_i, h'_{ij} \in H$.
Однако все $h_i a_i \in \psi([L, R])$, так как согласно нашим предположениям разрешимый радикал $R$ является $H$-подмодулем.
Следовательно, $(h_i a_i) W_{k-1} \subseteq W_k$ для всех $1 \leqslant k \leqslant t$.
Отсюда $b=0$, $I^t=0$ и $\psi([L, R]) \subseteq J(A)$.
\end{proof}

\begin{lemma}\label{LemmaHSolvableElementsJacobsonRadical}
Предположим, что радикал Джекобсона $J(A_1)$ любой $H$-инвариантной
ассоциативной подалгебры $A_1 \subseteq \End_\mathbbm{k}(V)$
является $H$-подмодулем.
Обозначим через $A_2$ ассоциативную подалгебру в $End_\mathbbm{k}(V)$,
порождённую операторами из $\psi(R)$.
Тогда для любой ассоциативной $H$-инвариантной подалгебры
$A_1 \subseteq A_2$ справедливо включение $J(A_1) \subseteq J(A)$.
\end{lemma}
\begin{proof} Заметим, что в силу $H$-инвариантности разрешимого радикала $R$ 
подалгебра $A_2$ является $H$-подмодулем. 
Рассмотрим $H$-инвариантную подалгебру Ли $\psi(R)+J(A)\subseteq \mathfrak{gl}(V)$.
Заметим, что алгебра Ли $\psi(R)+J(A)$ разрешима, так как радикал Джекобсона $J(A)$
нильпотентен, а $(\psi(R)+J(A))/J(A)
\cong \psi(R)/(\psi(R) \cap J(A))$ является гомоморфным образом разрешимого радикала $R$.
В силу теоремы Ли в пространстве $V$
существует базис, в котором все операторы из $\psi(R)+J(A)$ 
имеют верхнетреугольные матрицы. Обозначим соответствующее вложение
 $A  \hookrightarrow M_s(\mathbbm{k})$, где $s := \dim V$, через $\theta$.
 При этом $\theta(A_2+J(A)) \subseteq \UT_s(\mathbbm{k})$,
 где $\UT_s(\mathbbm{k})$~"--- ассоциативная алгебра верхнетреугольных матриц $s\times s$.

Докажем, что $H$-инвариантный идеал $I$ алгебры $A$, порождённый подпространством $J(A_1)+J(A)$,
нильпотентен. Отсюда будет следовать, что $J(A_1) \subseteq J(A)$.

Прежде всего заметим, что $\theta(J(A_1))$ и $\theta(J(A))$ содержатся в $\UT_s(\mathbbm{k})$
и состоят из нильпотентных элементов. 
Следовательно, в соответствующих матрицах
на главной диагонали стоят нулевые элементы и
$\theta(J(A_1)), \theta(J(A)) \subseteq \tilde N$,
где $$\tilde N := \langle e_{ij} \mid 1 \leqslant i < j \leqslant s \rangle_\mathbbm{k}.$$

Введём обозначение $$\tilde N_k := \langle e_{ij} \mid i+k \leqslant j \rangle_\mathbbm{k} \subseteq \tilde N.$$
Тогда $$\tilde N = \tilde N_1 \supsetneqq \tilde N_2 \supsetneqq \ldots \supsetneqq \tilde N_{m-1} \supsetneqq \tilde N_s = \lbrace 0\rbrace.$$
Пусть $\height_{\tilde N} a := k$, если $\theta(a) \in \tilde N_k$, $\theta(a) \notin \tilde N_{k+1}$.

Поскольку идеал $J(A)$ нильпотентен, $(J(A))^p =0$ для некоторого $p \in\mathbb N$.
Докажем, что $I^{s+p} = 0$.
Пользуясь~(\ref{EqHLModule2}), переместим все $\zeta(h)$, где $h \in H$, вправо
и получим, что пространство $I^{s+p}$
является линейной оболочкой элементов
 $b_1 j_1 b_2 j_2 \ldots j_{s+p} b_{s+p+1} \zeta(h)$, где
 $j_k \in J(A_1) \cup J(A)$, $b_k \in A_3 \cup \lbrace 1\rbrace$,
 $h \in H$. Здесь через $A_3$ обозначена подалгебра алгебры $\End_\mathbbm{k}(V)$,
 порождённая  элементами пространства $\psi(L)$.
 Если по крайней мере $p$ элементов $j_k$ принадлежат $J(A)$,
 такое произведение равно $0$.
 Следовательно, мы можем считать, что по крайней мере $s$
элементов $j_k$ принадлежат $J(A_1)$.

 Пусть $j_i \in J(A_1)$, $b_i \in A_3 \cup \lbrace 1\rbrace$.
Воспользуемся индукцией по $\ell$, чтобы доказать, что элемент
$j_1 b_1 j_2 b_2 \ldots b_{\ell-1} j_{\ell}$
может быть представлен в виде суммы $\tilde j_1 \tilde j_2 \ldots \tilde j_\alpha j'_1 j'_2\ldots j'_\beta
a$, где $\tilde j_i \in J(A_1)$, $j'_i \in J(A)$,
$a \in A_3 \cup \lbrace 1\rbrace$,
 а $\alpha+\sum_{i=1}^\beta \height_{\tilde N} j'_i \geqslant \ell$.
 Предположим, что элемент
  $j_1 b_1 j_2 b_2 \ldots b_{\ell-2} j_{\ell-1}$
 может быть представлен в виде суммы $\tilde j_1 \tilde j_2 \ldots \tilde j_\gamma j'_1 j'_2\ldots j'_\varkappa
a$,
где $\tilde j_i \in J(A_1)$, $j'_i \in J(A)$,
$a \in A_3 \cup \lbrace 1\rbrace$,
 а $\gamma+\sum_{i=1}^\varkappa \height_{\tilde N} j'_i \geqslant \ell-1$.
 Тогда
 $j_1 b_1 j_2 b_2 \ldots j_{\ell-1} b_{\ell-1}j_{\ell}$
 является суммой элементов
 $$\tilde j_1 \tilde j_2 \ldots \tilde j_\gamma j'_1 j'_2\ldots j'_\varkappa
a b_{\ell-1}j_{\ell} =
\tilde j_1 \tilde j_2 \ldots \tilde j_\gamma j'_1 j'_2\ldots j'_\varkappa
[ab_{\ell-1}, j_{\ell}] + \tilde j_1 \tilde j_2 \ldots \tilde j_\gamma j'_1 j'_2\ldots j'_\varkappa
j_{\ell} (a b_{\ell-1}).$$
Заметим, что в силу леммы~\ref{LemmaLR} и тождества Якоби $[ab_{\ell-1}, j_{\ell}] \in
J(A)$. Следовательно, достаточно рассматривать только второе слагаемое.
Однако
\begin{equation*}\begin{split}\tilde j_1 \tilde j_2 \ldots \tilde j_\gamma j'_1 j'_2\ldots j'_\varkappa
j_{\ell} (a b_{\ell-1})
= \tilde j_1 \tilde j_2 \ldots \tilde j_\gamma j_{\ell} j'_1 j'_2\ldots j'_\varkappa
 (a b_{\ell-1})+\\+\sum_{i=1}^{\varkappa}
\tilde j_1 \tilde j_2 \ldots \tilde j_\gamma  j'_1 j'_2\ldots j'_{i-1}[j'_{i}, j_\ell]
j'_{i+1}\ldots j'_\varkappa (a b_{\ell-1}).\end{split}\end{equation*}
Поскольку $[j'_{i}, j_\ell] \in J(A)$ и $\height_{\tilde N} [j'_{i}, j_\ell] \geqslant 1+ \height_{\tilde N} j'_i$, все слагаемые имеют требуемый вид.
 Следовательно, $$j_1 b_1 j_2 b_2 \ldots j_{s-1} b_{s-1}j_{s}
 \in \theta^{-1}(\tilde N_s) = \lbrace 0 \rbrace,$$ $I^{s+p}=0$ и $$
 J(A) \subseteq J(A_1)+J(A) \subseteq I \subseteq J(A).$$
\end{proof}

Нам также потребуется следующий вариант разложения Жордана:

\begin{lemma}\label{LemmaGeneralizedJordan} 
Предположим, что для любой $H$-инвариантной
ассоциативной подалгебры $A_1 \subseteq \End_\mathbbm{k}(V)$
существует $H$-инвариантное разложение
Веддербёрна~"--- Мальцева.
Пусть $W=\langle a_1, \ldots, a_t \rangle_\mathbbm{k}$~"--- некоторое подпространство разрешимого радикала $R$,
такое, что $\psi(W)$ является $H$-подмодулем $H$-модуля $\psi(R)$.
 Тогда для всех $1 \leqslant i \leqslant t$ существует разложение
   $\psi(a_i) = c_i+d_i$, где
$c_i$ и $d_i$ являются ассоциативными многочленами от $\psi(a_j)$, $1 \leqslant j
\leqslant t$, без свободного члена, $c_i$~"--- коммутирующие
диагонализуемые операторы на $V$,
а $d_i \in J(A)$. Более того, $\langle c_1, \ldots, c_t \rangle_\mathbbm{k}$
и $\langle d_1, \ldots, d_t \rangle_\mathbbm{k}$ являются $H$-подмодулями алгебры $\End_\mathbbm{k}(V)$.
\end{lemma}
\begin{proof}
Как и в доказательстве леммы~\ref{LemmaHSolvableElementsJacobsonRadical},
выберем такое вложение $\theta \colon A  \hookrightarrow M_s(\mathbbm{k})$,
что $\theta(\psi(R)+J(A)) \subseteq \UT_s(\mathbbm{k})$.
Обозначим через $A_1$ ассоциативную подалгебру алгебры $\End_\mathbbm{k}(V)$, порождённую элементами $\psi(a_i)$,
где
$1 \leqslant i \leqslant t$.
Подалгебра $A_1$ является $H$-инвариантной, так как $\psi(W)$~"--- $H$-подмодуль в $\psi(R)$.
Следовательно, $A_1$~"--- ассоциативная $H$-модульная алгебра.
Согласно нашим предположениям радикал Джекобсона $J(A_1)$ является $H$-инвариантным идеалом
и существует $H$-инвариантное разложение Веддербёрна~"--- Мальцева
   $A_1 = \tilde A_1 \oplus J(A_1)$
 (прямая сумма $H$-подмодулей),
 где $\tilde A_1$~"--- $H$-инвариантная полупростая подалгебра алгебры $A_1$.
Поскольку $\theta(\psi(R)) \subseteq \mathfrak{t}_s(\mathbbm{k})$, где $\mathfrak{t}_s(\mathbbm{k})$~"--- алгебра Ли
верхнетреугольных матриц $s\times s$, имеет место включение
$\theta(A_1) \subseteq \UT_s(\mathbbm{k})$. 
 Воспользуемся разложением $$\UT_s(\mathbbm{k}) = \mathbbm{k}e_{11}\oplus \mathbbm{k}e_{22}\oplus
 \ldots\oplus \mathbbm{k}e_{ss}\oplus \tilde N,$$
 где $$\tilde N := \langle e_{ij} \mid 1 \leqslant i < j \leqslant s \rangle_\mathbbm{k}$$
 "--- нильпотеный идеал. Отсюда алгебра $A_1$
 не содержит подалгебр, изоморфных $M_2(\mathbbm{k})$, и
  $\tilde A_1=\mathbbm{k}e_1 \oplus \ldots \oplus \mathbbm{k}e_q$
    для некоторых идемпотентов $e_i \in A_1$.
Для всякого $a_j$ обозначим через $d_j$ его компоненту в $J(A_1)$,
а через $c_j$ его компоненту $\mathbbm{k}e_1 \oplus \ldots \oplus \mathbbm{k}e_q$.
Заметим, что $e_i$~"--- коммутирующие диагонализуемые операторы.
Следовательно, для них существует базис в~$V$, состоящий из общих собственных векторов,
и $c_i$ также являются коммутирующими диагонализуемыми операторами.
Более того,
$$hc_j+hd_j=h\psi(a_j) \in
 \langle \psi(a_i)
\mid 1 \leqslant i \leqslant t \rangle_\mathbbm{k}
\subseteq \langle c_i
\mid 1 \leqslant i \leqslant t \rangle_\mathbbm{k}
\oplus \langle d_i
\mid 1 \leqslant i \leqslant t \rangle_\mathbbm{k}
\subseteq \tilde A_1 \oplus J(A_1)$$
для всех $h \in H$. Однако $\tilde A_1$ и $J(A_1)$ являются $H$-подмодулями, откуда $hc_j \in \tilde A_1$, а $hd_j \in J(A_1)$.
Следовательно, $hc_j\in\langle c_1, \ldots, c_t \rangle_\mathbbm{k}$
и $hd_j\in\langle d_1, \ldots, d_t \rangle_\mathbbm{k}$. Отсюда $\langle c_1, \ldots, c_t \rangle_\mathbbm{k}$
и $\langle d_1, \ldots, d_t \rangle_\mathbbm{k}$ являются $H$-подмодулями алгебры $\End_\mathbbm{k}(V)$.
В силу леммы~\ref{LemmaHSolvableElementsJacobsonRadical}
справедливы включения $J(A_1) \subseteq J(A)$ и 
$\langle d_1, \ldots, d_t \rangle_\mathbbm{k} \subseteq J(A)$.
\end{proof}

\begin{lemma} \label{LemmaRedIrr}
Пусть $V$~"--- конечномерный неприводимый $(H,L)$-модуль.
Предположим, что для любой $H$-инвариантной
ассоциативной подалгебры $A_1 \subseteq \End_\mathbbm{k}(V)$
существует $H$-инвариантное разложение
Веддербёрна~"--- Мальцева.
Тогда \begin{enumerate}
\item \label{RedSum} $V=V_1 \oplus \ldots \oplus V_q$ для некоторых
$L$-подмодулей $V_i$;
\item \label{RedScalar}  на каждом $V_i$ элементы разрешимого радикала $R$ действуют
скалярными операторами.
\end{enumerate}
\end{lemma}
\begin{proof}
Пусть $\psi(r_1), \ldots, \psi(r_t)$~"--- базис в $\psi(R)$. В силу леммы~\ref{LemmaGeneralizedJordan} существует разложение
 $\psi(r_i)=r'_i+r''_i$,
где  $r'_i$~"--- коммутирующие диагонализуемые операторы на $V$, а
 $r''_i \in J(A)$.
 Заметим, что в силу теоремы плотности $A = \End_\mathbbm{k}(V)$. Следовательно, $J(A)=0$ и $\psi(r_i)=r'_i$.
Отсюда для операторов $\psi(r_i)$
существует общий базис из собственных векторов, и
мы можем
выбрать подпространства
 $V_i$, где $1 \leqslant i \leqslant q$, $q\in \mathbb N$,
такие, что $$V=V_1 \oplus \ldots \oplus V_q,$$ и всякое $V_i$
является пересечением собственных подпространств операторов $\psi(r_i)$.
Заметим, что в силу леммы~\ref{LemmaLR} $$[\psi(r_i),\psi(a)]\in J(a)=0\text{ для всех }a\in L.$$
Следовательно, $V_i$ являются $L$-подмодулями, и лемма доказана.
\end{proof}

Рассмотрим теперь случай, когда $H=\mathbbm{k}G$ для некоторой группы $G$.

Пусть $L$~"--- алгебра Ли, на которой
действует автоморфизмами
некоторая группа $G$, а $V$~"--- $\mathbbm{k}G$-модуль.
Предположим, что гомоморфизм $\psi \colon L \to \mathfrak{gl}(V)$
задаёт на $V$ структуру $L$-модуля.
Будем говорить, что $(V, \psi)$ является \textit{$(G,L)$-модулем}, если $g(\psi(a)v)=\psi(ga)(gv)$
для всех $a\in L$, $g\in G$ и $v\in V$.
   Будем называть $(G,L)$-модуль $(V, \psi)$ \textit{неприводимым},
   если он не содержит нетривиальных
$G$-инвариантных $L$-подмодулей.

В случае действий групп лемма~\ref{LemmaRedIrr}
допускает следующее уточнение:

\begin{lemma} \label{LemmaRedIrrG}
Пусть $V$~"--- конечномерный неприводимый $(G,L)$-модуль,
а $L$~"--- алгебра Ли над алгебраически замкнутым полем $\mathbbm{k}$ характеристики $0$, на которой
действует автоморфизмами
некоторая группа $G$.
Предположим, что для любой $G$-инвариантной
ассоциативной подалгебры $A_1 \subseteq \End_\mathbbm{k}(V)$
существует $G$-инвариантное разложение
Веддербёрна~"--- Мальцева. (Это всегда выполнено, если, например, группа $G$ конечна.)
Тогда \begin{enumerate}
\item \label{RedSumG} $V=V_1 \oplus \ldots \oplus V_q$ для некоторых
$L$-подмодулей $V_i$;
\item \label{RedScalarG} на каждом $V_i$ элементы разрешимого радикала $R$ действуют
скалярными операторами;
\item \label{RedTransitiveG}
 для любого $g\in G$ существует такое $1 \leqslant j \leqslant q$, что $g V_i = V_j$,
 и это действие группы $G$ на множестве $\lbrace V_1, \ldots, V_q \rbrace$ транзитивно. 
\end{enumerate}
\end{lemma}
\begin{proof}
Предложения~\ref{RedSumG} и~\ref{RedScalarG} являются следствиями предложений~\ref{RedSum} и~\ref{RedScalar} 
леммы~\ref{LemmaRedIrr}.

Докажем предложение~\ref{RedTransitiveG}. Для всякого $V_i$
определим линейную функцию $\alpha_i \colon R \to \mathbbm{k}$, 
такую, что $\psi(r)v=\alpha_i(r)v$ для всех $r \in R$ и $v \in V$.
Тогда $V_i = \bigcap_{r\in R} \ker(\psi(r)-\alpha_i(r)\id_V)$,
а $$g V_i = \bigcap_{r\in R} \ker(\psi(g r)-\alpha_i(r)\id_V)
= \bigcap_{\tilde r\in R}
 \ker(\psi(\tilde r)-\alpha_i({g^{-1}}\tilde r) \id_V),$$
 где $\tilde r=gr$.
Следовательно, подмодуль $g V_i$ должен совпадать с одним из подмодулей $V_j$ для некоторого $1 \leqslant j \leqslant q$. Теперь транзитивность $G$-действия на множестве $\lbrace V_1, \ldots, V_q \rbrace$ следует из того, что $(G,L)$-модуль $V$ неприводим.
\end{proof}

\section{$(H,L)$-модули над $H$-комодульными алгебрами Ли $L$}

Пусть $L$~"---  $H$-комодульная алгебра Ли для некоторой алгебры Хопфа $H$, а
$\psi \colon L \to \mathfrak{gl}(V)$~"--- её представление.
Будем говорить, что $(V, \psi)$ является \textit{$(H,L)$-модулем}, если $V$ является $H$-комодулем
и $$\rho_V(\psi(a)v)=\psi(a_{(0)}) v_{(0)}\otimes a_{(1)}v_{(1)} \text{ для всех } a \in L,\ v \in V,$$
где $\rho_V \colon V \to V \otimes H$~"--- отображение, задающее на $V$ структуру комодуля.
Назовём $(V, \psi)$ \textit{симметрическим $(H,L)$-модулем},
если $$\rho_V(\psi(a)v)=\psi(a_{(0)}) v_{(0)}\otimes a_{(1)}v_{(1)} =
 \psi(a_{(0)}) v_{(0)}\otimes v_{(1)}a_{(1)} \text{ для всех } a \in L,\ v \in V.$$
\begin{example}
Если $L$~"--- $H$-комодульная алгебра Ли, 
то присоединённое представление
$\ad \colon L \to \mathfrak{gl}(L)$ определяет на $L$
структуру симметрического $(H,L)$-модуля,
поскольку
 $$\rho((\ad a)b)=\rho([a,b])= -\rho([b,a]) =-[b_{(0)}, a_{(0)}]\otimes b_{(1)}a_{(1)} =
 (\ad a_{(0)}) b_{(0)}\otimes b_{(1)}a_{(1)}$$ для всех $a,b \in L$.
\end{example}

Будем говорить, что $(H,L)$-модуль $(V, \psi)$ \textit{неприводим},
если он не содержит нетривиальных $H$-коинвариантных $L$-подмодулей.

 Докажем результат, двойственный к лемме~\ref{LemmaAnnHLmoduleHMod}.
\begin{lemma}\label{LemmaAnnHLmoduleHComod}
Пусть $V$~"--- $(H,L)$-модуль для некоторой алгебры Хопфа $H$ и $H$-комодульной алгебры Ли $L$
над произвольным полем $\mathbbm{k}$, а $M \subseteq V$~"--- его $H$-подкомодуль. Тогда $\Ann_L(M)$  является $H$-подкомодулем алгебры Ли $L$.
Если, кроме этого, $M$ является ещё и $L$-подмодулем, то $\Ann_L(M)$ является $H$-коинвариантным идеалом алгебры Ли $L$.
\end{lemma}
\begin{proof}
Достаточно доказать, что для всех $h^* \in H^*$, $a \in \Ann_L(M)$,
 $m \in M$,
 справедливо равенство $h^*(a_{(1)})a_{(0)}m=0$. Заметим, что
  \begin{equation*}\begin{split}h^*(a_{(1)})a_{(0)} m=h^*(a_{(1)})a_{(0)}m=h^*(a_{(1)})a_{(0)} \varepsilon(m_{(1)}) m_{(0)}=
  h^*(a_{(1)}\varepsilon(m_{(1)})1_H) a_{(0)} m_{(0)}= \\ =
  h^*(a_{(1)}m_{(1)} S(m_{(2)}))a_{(0)} m_{(0)}=
  h^*([a,m_{(0)}]_{(1)} S(m_{(1)}))(am_{(0)})_{(0)}
   =0,\end{split}\end{equation*} поскольку для каждого слагаемого
 $am_{(0)}=0$. Следовательно, $\Ann_L(M)$ является $H$-подкомодулем.
 Для доказательства второй части леммы достаточно повторить соответствующие рассуждения
 из доказательства леммы~\ref{LemmaAnnHLmoduleHComod}.
 \end{proof}

\section{$H$-(ко)инвариантное разложение полупростых алгебр}\label{SectionHSS}

Докажем теперь аналоги известного разложения для полупростых (в обычном смысле, т.е. таких, что их разрешимый радикал равен нулю) $H$-модульных и $H$-комодульных алгебр Ли.

\begin{theorem}\label{TheoremHLieSemiSimple}
Пусть $B$~"--- конечномерная полупростая $H$-модульная алгебра Ли,
где $H$~"--- произвольная алгебра Хопфа
над полем характеристики $0$.
Тогда $B=B_1 \oplus B_2 \oplus \ldots \oplus B_s$ (прямая сумма $H$-инвариантных идеалов)
для некоторых $H$-простых $H$-модульных алгебр Ли $B_i$.
\end{theorem} 
\begin{proof} Докажем теорему индукцией по $\dim B$. Если алгебра $B$ является $H$-простой, доказывать
  нечего. Предположим, что $B$ содержит нетривиальные $H$-инвариантные идеалы.
  В силу классической теоремы о разложении полупростой алгебры Ли
  (см., например, \cite[теорема 2.1.4]{GotoGrosshans})
   $$B=\tilde B_1 \oplus \tilde B_2 \oplus \ldots \oplus \tilde B_q
\text{ (прямая сумма идеалов)}$$
для некоторых простых алгебр Ли $\tilde B_j$.
Фиксируем произвольный минимальный $H$-инвариантный идеал 
$B_1 \subset B$.
Тогда $B_1 = \tilde B_{i_1} \oplus \tilde B_{i_2}
 \oplus \ldots \oplus \tilde B_{i_\ell}$ для некоторых $i_k$,
 а централизатор $B_0$ идеала $\tilde B$ в $B$ 
состоит из прямой суммы остальных простых алгебр Ли $\tilde B_j$.
  В силу леммы~\ref{LemmaAnnHLmoduleHMod} идеал $B_0$ является $H$-подмодулем.
  Теперь достаточно применить к $B_0$ предположение индукции
  и воспользоваться равенством
   $B=B_1 \oplus B_0$.
\end{proof}

Двойственный результат выглядит следующим образом:
\begin{theorem}\label{TheoremCoHLieSemiSimple}
Пусть $B$~"--- конечномерная полупростая $H$-комодульная алгебра Ли,
где $H$~"--- произвольная алгебра Хопфа
над полем характеристики $0$.
Тогда $B=B_1 \oplus B_2 \oplus \ldots \oplus B_s$ (прямая сумма $H$-коинвариантных идеалов)
для некоторых $H$-простых $H$-комодульных алгебр Ли $B_i$.
\end{theorem}
\begin{proof}
Повторим дословно доказательство теоремы~\ref{TheoremHLieSemiSimple}, использовав
вместо леммы~\ref{LemmaAnnHLmoduleHMod} лемму~\ref{LemmaAnnHLmoduleHComod}.
\end{proof}

 \section{Когомологии алгебр Ли и (ко)инвариантное разложение Леви}
 
Сперва напомним основные понятия когомологий алгебр Ли (см.~\cite{GotoGrosshans, Postnikov}).

Пусть $\psi \colon L \to \mathfrak{gl}(V)$~"--- представление
алгебры Ли $L$ на некотором векторном пространстве $V$ над полем $\mathbbm{k}$.
Обозначим через $C^k(L; V) \subseteq \Hom_\mathbbm{k}(L^{{}\otimes k}; V)$, где $k \in \mathbb N$, 
подпространство, состоящее из всех кососимметрических отображений. При этом считаем, что $C^0(L; V) := V$.

Напомним, что элементы пространства $C^k(L; V)$ называются \textit{$k$-коцепями}
с коэффициентами в $V$.   \textit{Кограничные операторы} $d \colon C^k(L;V) \to C^{k+1}(L; V)$
определяются на этих пространствах таким образом, чтобы было справедливо равенство $d^2 = 0$.
Элементы подпространства $$Z^k(L;\psi):= \ker(d \colon C^k(L;V) \to C^{k+1}(L; V))\subseteq C^k(L;V)$$ называются \textit{$k$-коциклами},
а элементы подпространства $$B^k(L;\psi):= d(C^{k-1}(L;V)) \subseteq C^k(L;V)$$
называются \textit{$k$-кограницами}. Пространство $H^k(L;\psi) := Z^k(L;\psi)/B^k(L;\psi)$
называется \textit{$k$-й группой когомологий}.

Для доказательства колинейного аналога теоремы Леви нам 
потребуются колинейные коцепи с коэффициентами в $(H,L)$-модулях.

Пусть $L$~"---  $H$-комодульная алгебра Ли для некоторой алгебры Хопфа $H$,
а $(V, \psi)$~"--- $(H,L)$-модуль.
Обозначим через $\tilde C^k(L; V)$ подпространство \textit{$H$-колинейных
коцепей}, т.е. таких отображений
 $\omega \in C^k(L; V)$, что $$\rho_V(\omega(a_1, a_2, \ldots, a_k))
 =\omega({a_1}_{(0)}, {a_2}_{(0)}, \ldots, {a_k}_{(0)})\otimes
 {a_1}_{(1)} {a_2}_{(1)} \ldots {a_k}_{(1)} \text{ для всех }a_i \in L.$$
  
 Если $(V,\psi)$~"--- $(H,L)$-модуль и алгебра Хопфа $H$ коммутативна,
то, очевидно, кограница всякой $H$-колинейной коцепи снова является $H$-колинейной коцепью.
Однако для $1$-коцепей и симметрического
 $(H,L)$-модуля $(V,\psi)$ это утверждение справедливо даже для некоммутативных
 алгебр Хопфа $H$:

\begin{lemma}\label{LemmaHCoboundary}
Если $(V,\psi)$~"--- симметрический $(H,L)$-модуль, тогда $$d(\tilde C^1(L;V)) 
\subseteq \tilde C^2(L;V).$$
\end{lemma}
\begin{proof}
Пусть $\omega \in \tilde C^1(L;V)$.
Тогда $$(d\omega)(x,y)
:=\psi(x)\omega(y)-\psi(y)\omega(x) - \omega([x,y])$$
и \begin{equation*}\begin{split}\rho_V((d\omega)(x,y))= \\ = \psi(x_{(0)})\omega(y)_{(0)}
\otimes x_{(1)}\omega(y)_{(1)} -\psi(y_{(0)})\omega(x)_{(0)}
\otimes {\omega(x)_{(1)}} y_{(1)} - \omega([x,y]_{(0)}) \otimes [x,y]_{(1)}
=\\= \psi(x_{(0)})\omega(y_{(0)})
\otimes x_{(1)}y_{(1)} -\psi(y_{(0)})\omega(x_{(0)})
\otimes x_{(1)}y_{(1)} - \omega([x_{(0)},y_{(0)}]) \otimes x_{(1)}y_{(1)}
=\\= (d\omega)(x_{(0)},y_{(0)}) \otimes x_{(1)}y_{(1)}.\end{split}\end{equation*}
\end{proof}

Пусть $\tilde Z^2(L; \psi) := Z^2(L; \psi) \cap \tilde C^2(L;V)$
и $\tilde B^2(L; \psi) := d(\tilde C^1(L;V))$.
Лемма~\ref{LemmaHCoboundary} позволяет определить
\textit{вторую группу $H$-колинейных когомологий} $\tilde H^2(L; \psi) :=
 \tilde Z^2(L; \psi)/\tilde B^2(L; \psi)$.
 
 Пусть $V$ и $W$~"--- $H$-комодули  для некоторой алгебры Хопфа $H$.
 Будем говорить, что $\mathbbm{k}$-линейное отображение $\varphi \colon V \to W$
является \textit{$H$-колинейным}, если $$\varphi(v)_{(0)}\otimes \varphi(v)_{(1)}=\varphi(v_{(0)})
\otimes v_{(1)}\text{ для всех }v \in V.$$

В~\cite{Taft} Э.\,Дж.~Тафт использовал оригинальный приём Машке
для того, чтобы превратить неинвариантное отображение в инвариантное.
В~\cite{SteVanOyst} Д. Штефан и Ф. Ван Ойстайен 
использовали приём Машке, адаптированный для алгебр Хопфа с левым интегралом.
Мы также воспользуемся последним приёмом:
 
 \begin{lemma}\label{LemmaHcolinear}
 Пусть $r \colon V \to W$~"--- $\mathbbm{k}$-линейное отображение, где $V$ и $W$~"--- $H$-комодули,
 а $H$~"--- алгебра Хопфа над полем $\mathbbm{k}$. Пусть алгебра Хопфа $H$ обладает левым интегралом $t \in H^*$. Тогда
 отображение $\tilde r \colon V \to W$,
 заданное равенством
 $$\tilde r(x) = t\bigl(r(x_{(0)})_{(1)}S(x_{(1)})\bigr)r(x_{(0)})_{(0)} \text{ при } x \in V$$ является $H$-колинейным. Если, кроме этого, $\pi  r = \id_V$ для некоторого
 $H$-колинейного отображения $\pi \colon W \to V$ и $t(1)=1$, то $\pi  \tilde r = \id_V$.
 \end{lemma}
 \begin{example} Если $G$~"--- группа, а $H=\mathbbm{k}G$, тогда $V=\bigoplus_{g \in G}V^{(g)}$ и $W=\bigoplus_{g \in G}W^{(g)}$ являются градуированными пространствами.
 Предположим, что левый интеграл $t$ взят из примера~\ref{ExampleIntegralFG}.
 Тогда отображение $\tilde r$, заданное равенством $\tilde r(x) = \sum_{g \in G} p_{W,g}\ r(p_{V,g} x)$ при $x \in V$,
 является градуированным.
 Здесь $p_{V,g}$~"--- проектор пространства $V$ на $V^{(g)}$
 с ядром $\bigoplus_{\substack{h \in G, \\ h \ne g}} V^{(h)}$,
 а $p_{W,g}$~"--- проектор пространства $W$ на $W^{(g)}$
 с ядром $\bigoplus_{\substack{h \in G, \\ h \ne g}} W^{(h)}$.
 \end{example}
 \begin{proof}[Доказательство леммы~\ref{LemmaHcolinear}]
 Заметим, что \begin{equation*}\begin{split}
 \tilde r(x_{(0)}) \otimes x_{(1)}=
 t\bigl(r(x_{(0)})_{(1)}S(x_{(1)})\bigr)\ r(x_{(0)})_{(0)} \otimes x_{(2)}
 = \\ =
 r(x_{(0)})_{(0)} \otimes \Bigl(t\bigl(r(x_{(0)})_{(1)}S(x_{(1)})\bigr)1\Bigr)x_{(2)}=
 \\ = r(x_{(0)})_{(0)} \otimes t\bigl(r(x_{(0)})_{(1)(2)}(Sx_{(1)})_{(2)}\bigr)\ r(x_{(0)})_{(1)(1)}(Sx_{(1)})_{(1)}
  x_{(2)} = \\ =
  r(x_{(0)})_{(0)} \otimes t\bigl(r(x_{(0)})_{(2)}S(x_{(1)})\bigr)\ r(x_{(0)})_{(1)}(Sx_{(2)})x_{(3)}
  = \\ =
 t\bigl(r(x_{(0)})_{(2)}S(x_{(1)})\bigr)\ r(x_{(0)})_{(0)} \otimes r(x_{(0)})_{(1)}
  = \rho_W(\tilde r(x)).\end{split}\end{equation*}
 Тогда отображение $\tilde r$ является $H$-колинейным, и первая часть
 леммы доказана.
 
 Предположим, что $\pi  r = \id_V$ для некоторого
 $H$-колинейного отображения $\pi \colon W \to V$. Пусть $x \in V$.
 Тогда \begin{equation*}\begin{split}(\pi  \tilde r)(x)=
 t\bigl(r(x_{(0)})_{(1)}S(x_{(1)})\bigr)\ \pi\bigl(r(x_{(0)})_{(0)}\bigr)
 = t\bigl((\pi r)(x_{(0)})_{(1)}S(x_{(1)})\bigr)\ (\pi r)(x_{(0)})_{(0)}
= \\ =
  t\bigl(x_{(0)(1)}S(x_{(1)})\bigr)\ x_{(0)(0)} =
 t\bigl(x_{(1)}S(x_{(2)})\bigr)\ x_{(0)} = t(1) x = x.\end{split}\end{equation*}
 \end{proof}

\begin{lemma}\label{LemmaH2zero}
Пусть $(V,\psi)$~"--- конечномерный симметрический $(H,L)$-модуль,
где $L$~"--- конечномерная $H$-комодульная полупростая алгебра Ли над полем $\mathbbm{k}$ характеристики $0$,
а $H$~"--- алгебра Хопфа с таким $\ad$-инвариантным левым интегралом $t \in H^*$, что $t(1)=1$.
 Тогда $\tilde H^2(L; \psi) = 0$.
\end{lemma}
\begin{proof} Напомним, что в силу второй леммы Уайтхеда
 (см., например, \cite[упражнение~3.5]{GotoGrosshans}
 или \cite[лекция 19, следствие 2]{Postnikov}), $H^2(L; \psi) = 0$.
Следовательно, если $\omega \in \tilde Z^2(L; \psi)$, то существует такая $1$-коцепь $\nu \in C^1(L; V)$,
что $\omega = d \nu$. Пусть $\tilde \nu$~"--- отображение,
полученное из $\nu$ согласно лемме~\ref{LemmaHcolinear}. 
Тогда $\tilde \nu \in \tilde C^1(L; V)$. Докажем, что $d\tilde\nu = \omega$.

Пусть $a,b \in L$. Тогда
\begin{equation*}\begin{split}(d\tilde\nu)(a, b) = \psi(a)\tilde\nu(b)
- \psi(b)\tilde\nu(a)-\tilde\nu([a,b])= t(\nu(b_{(0)})_{(1)} S(b_{(1)}))\ \psi(a) \nu(b_{(0)})_{(0)}
- \\ - t(\nu(a_{(0)})_{(1)} S(a_{(1)}))\ \psi(b) \nu(a_{(0)})_{(0)}
-t(\nu([a,b]_{(0)})_{(1)} S([a,b]_{(1)}))\ \nu([a,b]_{(0)})_{(0)}= \\ =
  t(\nu(b_{(0)})_{(1)} S(b_{(1)}))\ \psi(\varepsilon(a_{(1)})a_{(0)}) \nu(b_{(0)})_{(0)}
- \\ - t(\nu(a_{(0)})_{(1)} S(a_{(1)}))\ \psi(\varepsilon(b_{(1)})b_{(0)}) \nu(a_{(0)})_{(0)}
-t(\nu([a_{(0)},b_{(0)}])_{(1)} S(a_{(1)}b_{(1)}))\ \nu([a_{(0)},b_{(0)}])_{(0)}= \\ =
  \varepsilon(a_{(1)})t\Bigl( \nu(b_{(0)})_{(1)}S(b_{(1)})\Bigr)\ \psi(a_{(0)}) \nu(b_{(0)})_{(0)}
- \\ - t\Bigl(\nu(a_{(0)})_{(1)}\varepsilon(b_{(1)})S(a_{(1)})\Bigr)\ \psi(b_{(0)}) \nu(a_{(0)})_{(0)}
-\\ - t\Bigl(\nu([a_{(0)},b_{(0)}])_{(1)}S(a_{(1)}b_{(1)})\Bigr)\ \nu([a_{(0)},b_{(0)}])_{(0)}.\end{split}\end{equation*}

Поскольку $t$ является $\ad$-инвариантным, справедливо равенство
\begin{equation*}\begin{split}(d\tilde\nu)(a, b) =  t\Bigl(a_{(1)} \nu(b_{(0)})_{(1)}(Sb_{(1)})S(a_{(2)})\Bigr)\ \psi(a_{(0)}) \nu(b_{(0)})_{(0)}
- \\ - t\Bigl(\nu(a_{(0)})_{(1)}b_{(1)}(Sb_{(2)})S(a_{(1)})\Bigr)\ \psi(b_{(0)}) \nu(a_{(0)})_{(0)}
-\\ - t\Bigl(\nu([a_{(0)},b_{(0)}])_{(1)}S(a_{(1)}b_{(1)})\Bigr)\ \nu([a_{(0)},b_{(0)}])_{(0)}
= \\ = t\Bigl(\bigl(a_{(0)(1)} \nu(b_{(0)})_{(1)}\bigr)(Sb_{(1)})S(a_{(1)})\Bigr)\ \psi(a_{(0)(0)}) \nu(b_{(0)})_{(0)}
- \\ -t\Bigl(\bigl(\nu(a_{(0)})_{(1)}b_{(0)(1)}\bigr)(Sb_{(1)})S(a_{(1)})\Bigr)\ \psi(b_{(0)(0)}) \nu(a_{(0)})_{(0)}
-\\ - t\Bigl(\nu([a_{(0)},b_{(0)}])_{(1)}S(a_{(1)}b_{(1)})\Bigr)\ \nu([a_{(0)},b_{(0)}])_{(0)}
.\end{split}\end{equation*}
Поскольку $(H,L)$-модуль $(V, \psi)$ симметрический, имеем
\begin{equation*}\begin{split}(d\tilde\nu)(a, b)  =
  t\Bigl(\bigl(\psi(a_{(0)}) \nu(b_{(0)})\bigr)_{(1)}(Sb_{(1)})S(a_{(1)})\Bigr)\ \bigl(\psi(a_{(0)}) \nu(b_{(0)})\bigr)_{(0)}
- \\ - t\Bigl(\bigl(\psi(b_{(0)}) \nu(a_{(0)})\bigr)_{(1)}(Sb_{(1)})S(a_{(1)})\Bigr)\ \bigl(\psi(b_{(0)}) \nu(a_{(0)})\bigr)_{(0)}
-\\ - t\Bigl(\nu([a_{(0)},b_{(0)}])_{(1)}S(a_{(1)}b_{(1)})\Bigr)\ \nu([a_{(0)},b_{(0)}])_{(0)}= \\ =
  t\Bigl(\bigl(\psi(a_{(0)}) \nu(b_{(0)})\bigr)_{(1)}S(a_{(1)}b_{(1)})\Bigr)\ \bigl(\psi(a_{(0)}) \nu(b_{(0)})\bigr)_{(0)}
- \\ -t\Bigl(\bigl(\psi(b_{(0)}) \nu(a_{(0)})\bigr)_{(1)}S(a_{(1)}b_{(1)})\Bigr)\ \bigl(\psi(b_{(0)}) \nu(a_{(0)})\bigr)_{(0)}
-\\ - t\Bigl(\nu([a_{(0)},b_{(0)}])_{(1)} S(a_{(1)}b_{(1)})\Bigr)\ \nu([a_{(0)},b_{(0)}])_{(0)}= \\ =
 t\bigl(\omega(a_{(0)}, b_{(0)})_{(1)}S(a_{(1)}b_{(1)})\bigr)\ \omega(a_{(0)}, b_{(0)})_{(0)}
=t\bigl(a_{(1)}b_{(1)} S(a_{(2)}b_{(2)})\bigr)\ \omega(a_{(0)}, b_{(0)})=\omega(a,b),\end{split}\end{equation*} поскольку $\omega \in \tilde Z^2(L; \psi)$ и $t(1)=1$. Следовательно, 
$\tilde Z^2(L; \psi) = \tilde B^2(L; \psi)$ и $\tilde H^2(L; \psi)=0$.
\end{proof}

Теорема~\ref{TheoremHcoLevi} является $H$-комодульной версией теоремы Леви:

\begin{theorem}\label{TheoremHcoLevi}
Пусть $L$~"--- конечномерная $H$-комодульная алгебра Ли над полем $\mathbbm{k}$ характеристики $0$,
где $H$~"--- алгебра Хопфа.  Предположим, что разрешимый радикал $R$
алгебры Ли $L$ является $H$-подкомодулем и 
существует такой $\ad$-инвариантный левый интеграл $t \in H^*$, что $t(1)=1$.
Тогда существует такая  $H$-коинвариантная максимальная полупростая подалгебра $B\subseteq L$,
что
$L=B\oplus R$ (прямая сумма $H$-подкомодулей).
\end{theorem}
\begin{proof}
Проведём доказательство по той же схеме,
по которой доказывается обычная теорема Леви.

Сперва рассмотрим случай, когда разрешимый радикал $R$ является абелевым идеалом. 
Пусть $\pi \colon L \to L/R$~"--- естественный сюръективный гомоморфизм.
Заметим, что $L/R$~"--- полупростая $H$-комодульная алгебра Ли, а
 $\pi$~"--- $H$-колинейное отображение, поскольку $R$ является $H$-подкомодулем.  Рассмотрим произвольное 
 $\mathbbm{k}$-линейное отображение $r \colon L/R \to L$, такое, что $\pi  r = \id_{L/R}$.
В силу леммы~\ref{LemmaHcolinear} можно считать, что отображение $r$ является $H$-колинейным.

Введём обозначение $$\Phi(x,y):=[r(x),r(y)]-r([x,y]).$$ Заметим, что $\pi(\Phi(a,b))=[a,b]-[a,b]=0$
для всех $a,b \in L/R$. Следовательно, $\Phi(a,b) \in R$.
Пусть $\psi \colon L/R \to \mathfrak{gl}(R)$~"--- линейное отображение,
заданное равенством $\psi(a)(v)=[r(a), v]$, где $a\in L/R$, $v\in R$.
При этом
 \begin{equation*}\begin{split}\psi(a)\psi(b)(v)-\psi(b)\psi(a)(v)-\psi([a,b])(v)=[r(a), [r(b), v]]-[r(b), [r(a), v]]
 -[r([a, b]), v]=\\=[[r(a),r(b)],v]-[r([a, b]), v]=[\Phi(a,b), v]=0
 \text{ для всех } a,b \in L/R,\end{split}\end{equation*} поскольку $\Phi(a,b) \in R$ и $[R,R]=0$.
 Следовательно, отображение $\psi$ является представлением алгебры Ли $L/R$.
Более того, $(R, \psi)$~"--- симметрический $(H, L/R)$-модуль.
Заметим, что  $\Phi \in \tilde Z^2(L/R; \psi)$,  поскольку $d \Phi = 0$.
Следовательно, в силу леммы~\ref{LemmaH2zero},
 $\Phi=d \omega$ для некоторой $H$-колинейной $1$-коцепи $\omega \in \tilde C^1(L/R; R)$.
 Отсюда \begin{equation*}\begin{split}[(r-\omega)(a),(r-\omega)(b)]-(r-\omega)([a,b])=\\=
 ([r(a),r(b)]-r([a,b]))-([r(a),\omega(b)]-[r(b),\omega(a)]-\omega([a,b]))+[\omega(a),\omega(b)]=\\=
 \Phi(a,b)-(d\omega)(a,b)+0=0\end{split}\end{equation*}  для всех $a,b \in L/R$ и $\pi  (r-\omega) = \pi  r = \id_{L/R}$. Следовательно, $(r-\omega)$ является $H$-колинейным гомоморфным вложением
алгебры Ли $L/R$ в $L$, и $L=B\oplus R$ (прямая сумма $H$-подкомодулей), где $B=(r-\omega)(L/R)$, т.е.
в этом случае $H$-коинвариантная теорема Леви доказана.
 
 Докажем теперь общий случай индукцией по $\dim R$.
Теорема уже доказана в случае $[R,R]=0$. Предположим, что $[R,R]\ne 0$. Заметим, что $[R,R]\ne R$, поскольку
идеал $R$ разрешим.
 Более того, $[R,R]$ является $H$-подкомодулем. Рассмотрим алгебру Ли $L/[R,R]$.
 Поскольку алгебра Ли $(L/[R,R])/(R/[R,R])\cong L/R$ полупроста,
 идеал $R/[R,R]$ является разрешимым радикалом агебры Ли $L/[R,R]$.
Применим теперь предположение индукции и получим, что $L/[R,R] = L_1/[R,R] \oplus R/[R,R]$
 (прямая сумма $H$-подкомодулей)
 для некоторой $H$-коинвариантной подалгебры Ли $L_1 \subset L$,
 где $L_1/[R,R] \cong L/R$.
Применим теперь предположение индкуции к
 $L_1$ и получим, что $L_1 = B \oplus [R,R]$ (прямая сумма $H$-подкомодулей),
 где $B \cong L_1/[R,R] \cong L/R$~"--- полупростая подалгебра.
 Следовательно, $L = B \oplus R$ (прямая сумма $H$-подкомодулей), и теорема доказана.
 \end{proof}
 
 Получим теперь некоторые важные следствия из теоремы~\ref{TheoremHcoLevi}:
 
\begin{theorem}\label{TheoremHLevi}
Пусть $L$~"--- конечномерная $H$-(ко)модульная алгебра Ли над полем $\mathbbm{k}$ характеристики $0$,
где $H$~"--- конечномерная (ко)полупростая алгебра Хопфа.
Тогда существует такая $H$-(ко)инвариантная максимальная полупростая подалгебра $B\subseteq L$,
что $L=B\oplus R$ (прямая сумма $H$-под(ко)модулей).
\end{theorem}
\begin{proof}
Достаточно использовать двойственность между $H$-действиями и $H^*$-кодействиями,
применить пример~\ref{ExampleIntegralHSS}, следствие~\ref{CorollaryHModRadicalsLie} и теорему \ref{TheoremHcoLevi}.
\end{proof}

Д.~Пагон, Д.~Реповш и М.\,В.~Зайцев~\cite{PaReZai}
доказали градуированную версию теоремы Леви
для конечной группы.
Используя теорему \ref{TheoremHcoLevi} мы можем доказать это утверждение для любой группы:

\begin{theorem}\label{TheoremGradLevi}
Пусть $L$~"--- конечномерная алгебра Ли над полем $\mathbbm{k}$ характеристики $0$,
градуированная произвольной группой $G$. 
Тогда существует такая градуированная максимальная полупростая подалгебра $B\subseteq L$,
что $L=B\oplus R$ (прямая сумма градуированных подпространств).
\end{theorem}
\begin{proof}
Воспользуемся примерами~\ref{ExampleIntegralFG}, \ref{ExampleComoduleGraded}, следствием~\ref{CorollaryLieRadicalsGraded} и теоремой~\ref{TheoremHcoLevi}.
\end{proof}

Применим теперь теорему~\ref{TheoremHcoLevi}
к алгебрам Ли с рациональным действием редуктивной
аффинной алгебраической группы:

\begin{theorem}\label{TheoremAffAlgGrLevi}
Пусть $L$~"--- конечномерная 
алгебра Ли над алгебраически замкнутым полем $\mathbbm{k}$ характеристики $0$,
а $G$~"--- редуктивная аффинная
алгебраическая группа над $\mathbbm{k}$.
Предположим, что группа $G$ рационально действует на $L$
автоморфизмами. Тогда существует такая $G$-инвариантная подалгебра $B\subseteq L$,
что $L=B\oplus R$ (прямая сумма $G$-инвариантных подпространств).
\end{theorem}
\begin{proof}
Нужно заметить, что разрешимый радикал $R$
инвариантен относительно всех автоморфизмов
и использовать примеры~\ref{ExampleIntegralAffAlgGr}, \ref{ExampleRegularAffActAutAlgebra} и теорему~\ref{TheoremHcoLevi}.
\end{proof}

\begin{theorem}\label{TheoremSSLieActionLevi}
Пусть $L$~"--- конечномерная 
алгебра Ли над алгебраически замкнутым полем $\mathbbm{k}$ характеристики $0$, на которой
действует дифференцированиями конечномерная полупростая алгебра Ли $\mathfrak g$. Тогда существует такая $\mathfrak g$-инвариантная подалгебра $B\subseteq L$,
что $L=B\oplus R$ (прямая сумма $\mathfrak g$-инвариантных подпространств).
\end{theorem}
\begin{proof} Достаточно воспользоваться теоремами~\ref{TheoremLieDiffActionReplacement} и~\ref{TheoremAffAlgGrLevi}.
\end{proof}

Приведём теперь примеры $H$-модульных алгебр Ли, для которых $H$-инвариантное
разложение Леви не существует.

\begin{example}[Ю.\,А.~Бахтурин]\label{ExampleGnoninvLevi}
Пусть $$L = \left\lbrace\left(\begin{array}{cc} C & D \\
0 & 0
  \end{array}\right) \mathrel{\biggl|} C \in \mathfrak{sl}_m(\mathbbm{k}), D\in M_m(\mathbbm{k})\right\rbrace
  \subseteq \mathfrak{sl}_{2m}(\mathbbm{k}),\ m \geqslant 2.$$
  Тогда идеал $$R=\left\lbrace\left(\begin{array}{cc} 0 & D \\
0 & 0
  \end{array}\right) \mathrel{\biggl|} D\in M_m(\mathbbm{k})\right\rbrace
  $$
  является разрешимым (и нильпотентным) идеалом алгебры Ли $L$.
    Определим $\varphi \in \Aut(L)$ по формуле
  $$\varphi\left(\begin{array}{cc} C & D \\
0 & 0
  \end{array}\right)=\left(\begin{array}{cc} C & C+D \\
0 & 0
  \end{array}\right).$$
  Тогда группа $G=\langle \varphi \rangle
  \cong \mathbb Z$ действует на $L$ автоморфизмами, т.е. $L$ является $\mathbbm{k}G$-модульной
  алгеброй.
  Однако в $L$ не существует такой $\mathbbm{k}G$-инвариантной полупростой подалгебры
   $B$, что $L=B\oplus R$ (прямая сумма $\mathbbm{k}G$-подмодулей).
\end{example}
\begin{proof} 
Пусть $a \in L$. Тогда $\varphi(a)-a \in R$. Предположим, что $B$~"--- такое $G$-инвариантное 
подпространство, что $B \cap R = \lbrace 0 \rbrace$. Тогда $\varphi(b)-b=0$ для всех $b \in B$
и $B \subseteq R$. Следовательно, $B=0$ и $\mathbbm{k}G$-инвариантного разложения Леви не существует.
\end{proof}
\begin{example}\label{ExampleDiffnoninvLevi}
Пусть $L$~"---  та же алгебра Ли, что и в примере~\ref{ExampleGnoninvLevi}.
Рассмотрим присоединённое представление алгебры Ли $L$ на себе самой дифференцированиями.
Тогда $L$ оказывается $U(L)$-модульной алгеброй Ли
(см. пример~\ref{ExampleUgModule}), и
все $U(L)$-подмодули в $L$~"--- это в точности идеалы алгебры Ли $L$.
Однако в $L$ не существует такой $U(L)$-инвариантной полупростой подалгебры $B$,
что $L=B\oplus R$ (прямая сумма $U(L)$-подмодулей).
\end{example}
\begin{proof}
Предположим, что $L=B\oplus R$, где $B$~"--- некоторый $U(L)$-подмодуль. Тогда $B$ является идеалом алгебры Ли $L$, а $R$~"--- центром алгебры Ли $L$, поскольку $[R,R]=0$. 
Получаем противоречие. Следовательно, $U(L)$-инвариантного разложения
Леви для алгебры Ли $L$ не существует.
\end{proof}

\section{$H$-(ко)инвариантный аналог теоремы Вейля}\label{SectionHWeyl}

Прежде всего нам потребуется следующее дополнение к лемме~\ref{LemmaHcolinear}.
\begin{lemma}\label{LemmaHLlinear}
Пусть $\pi \colon V \to W$~"--- гомоморфизм
 $L$-модулей, где $(V, \varphi)$ и $(W, \psi)$ являются $(H,L)$-модулями
 для $H$-комодульной алгебры Ли $L$ и алгебры Хопфа $H$. Пусть $t \in H^*$~"--- $\ad$-инвариантный
 левый интеграл на $H$. Тогда $\tilde \pi \colon V \to W$, где
 $$\tilde \pi(x) = t\bigl(\pi(x_{(0)})_{(1)}S(x_{(1)})\bigr)\pi(x_{(0)})_{(0)} \text{ при } x \in V,$$ 
 является
 $H$-колинейным гомоморфизмом $L$-модулей. Более того, если $t(1)=1$, $W \subseteq V$, 
 а $\pi$ является проектором $V$ на $W$, то $\tilde\pi$ также является проектором $(H,L)$-модуля $V$ на $W$.
 
\end{lemma}
\begin{proof} Отображение $\tilde \pi$ является $H$-колинейным согласно лемме~\ref{LemmaHcolinear}. 
Пусть $a \in L$, $x\in V$. Тогда \begin{equation*}\begin{split}\tilde \pi(\varphi(a)x)= 
t\bigl(\pi((\varphi(a)x)_{(0)})_{(1)}S((\varphi(a)x)_{(1)})\bigr)\pi((\varphi(a)x)_{(0)})_{(0)}=\\=
t\bigl(\pi(\varphi(a_{(0)})x_{(0)})_{(1)}S(a_{(1)}x_{(1)})\bigr)\pi(\varphi(a_{(0)})x_{(0)})_{(0)}=\\=
t\bigl((\psi(a_{(0)})\pi(x_{(0)}))_{(1)}S(a_{(1)}x_{(1)})\bigr)(\psi(a_{(0)})\pi(x_{(0)}))_{(0)}=\\=
t\bigl(a_{(1)}\pi(x_{(0)})_{(1)}S(a_{(2)}x_{(1)})\bigr)\psi(a_{(0)})\pi(x_{(0)})_{(0)}=\\=
t\bigl(a_{(1)}\pi(x_{(0)})_{(1)}(Sx_{(1)})S(a_{(2)})\bigr)\psi(a_{(0)})\pi(x_{(0)})_{(0)}=\\=
t\bigl(\pi(x_{(0)})_{(1)}S(x_{(1)})\bigr)\psi(a)\pi(x_{(0)})_{(0)}
=\psi(a)\tilde \pi(x),\end{split}\end{equation*}
поскольку $t$ является $\ad$-инвариантным левым интегралом. Следовательно, $\tilde\pi$~"---
 $H$-колинейный гомоморфизм $L$-модулей.
 
 Пусть $t(1)=1$, $W \subseteq V$ и $\pi$ является проектором $(H,L)$-модуля $V$ на $W$.
Рассмотрим произвольный элемент $x \in W$.
Поскольку $W$ является $H$-подкомодулем, $x_{(0)}\otimes x_{(1)} \in W \otimes H$.
 Следовательно, \begin{equation*}\begin{split}\tilde \pi(x) = t\bigl(\pi(x_{(0)})_{(1)}S(x_{(1)})\bigr)\pi(x_{(0)})_{(0)}
 =\\=t\bigl(x_{(0)(1)}S(x_{(1)})\bigr)x_{(0)(0)}
 =t\bigl(x_{(1)}S(x_{(2)})\bigr)x_{(0)}
 =t(1)x=x\end{split}\end{equation*} и $\tilde\pi$ также является проектором $(H,L)$-модуля $V$ на $W$.
\end{proof}

\begin{theorem}\label{TheoremHcoWeyl}
Пусть $L$~"--- $H$-комодульная алгебра Ли над полем характеристики $0$,
а $H$~"--- алгебра Хопфа, обладающая таким $\ad$-инвариантным интегралом $t \in H^*$, что $t(1)=1$.
Пусть $(V, \psi)$~"--- конечномерный $(H,L)$-модуль,
который вполне приводим как обычный $L$-модуль (без учёта $H$-кодействия). Тогда $V=V_1 \oplus V_2 \oplus \ldots \oplus V_s$ для некоторых неприводимых $(H,L)$-подмодулей $V_i$. 
\end{theorem}
\begin{proof}
Снова воспользуемся приёмом Машке.
Достаточно показать, что для любого $H$-коинвариантного
$L$-подмодуля $W \subseteq V$ существует проектор $\tilde \pi \colon V \to W$, который
является $H$-колинейным гомоморфизмом $L$-модулей. Тогда $W = V \oplus \ker \tilde\pi$,
и можно воспользоваться индукцией по $\dim V$.

Поскольку $V$ является вполне приводимым $L$-модулем, 
существует проектор $\pi \colon V \to W$, являющийся гомоморфизмом $L$-модулей.
Теперь достаточно определить проектор $\tilde \pi$ в соответствии с
леммой~\ref{LemmaHLlinear}. Тогда $\tilde \pi$
является $H$-колинейным
гомоморфизмом $L$-модулей.
\end{proof}

Теперь докажем аналог теоремы Вейля (см., например, \cite[теорема 6.3]{HumphreysLieAlg}) для $H$-комодульных
алгебр Ли.
\begin{corollary}
Пусть $L$~"--- конечномерная полупростая $H$-комодульная алгебра Ли над полем характеристики $0$,
а $H$~"--- алгебра Хопфа, обладающая таким $\ad$-инвариантным интегралом $t \in H^*$, что $t(1)=1$.
Пусть $(V, \psi)$~"--- конечномерный $(H,L)$-модуль. Тогда $V=V_1 \oplus V_2 \oplus \ldots \oplus V_s$ для некоторых неприводимых $(H,L)$-подмодулей $V_i$. 
\end{corollary}
\begin{proof}
Нужно сперва заметить, что в силу обычной теоремы Вейля модуль $V$
вполне приводим как $L$-модуль, а затем применить теорему~\ref{TheoremHcoWeyl}.
\end{proof}

Точно так же, как и у теорем, доказанных ранее, у теоремы~\ref{TheoremHcoWeyl} имеются варианты для случая $H=\mathbbm{k}G$ и для случая конечномерной алгебры Хопфа $H$.

Пусть $G$~"--- группа, $L=\bigoplus_{g\in G} L^{(g)}$~"--- $G$-градуированная алгебра Ли,
а $V=\bigoplus_{g\in G} V^{(g)}$~"--- $G$-градуированное векторное пространство.
Пусть $\psi \colon L \to \mathfrak{gl}(V)$~"--- гомоморфизм, задающий на 
$V$ структуру $L$-модуля.
 Будем говорить, что $(V, \psi)$~"--- \textit{градуированный $L$-модуль},
 если $\psi(a^{(g)})v^{(h)} \in V^{(gh)}$
 для всех $g,h \in G$, $a^{(g)} \in L^{(g)}$, $v^{(h)} \in V^{(h)}$.
 Будем называть градуированный $L$-модуль $(V, \psi)$ \textit{неприводимым},
 если он не содержит нетривиальных градуированных $L$-подмодулей.

\begin{theorem}\label{TheoremGradWeyl}
Пусть $L$~"--- алгебра Ли над полем характеристики $0$,
градуированная произвольной группой, а $(V, \psi)$~"---
конечномерный градуированный $L$-модуль, вполне приводимый как обычный $L$-модуль
без учёта градуировки.
 Тогда $$V=V_1 \oplus V_2 \oplus \ldots \oplus V_s$$
для некоторых неприводимых градуированных $L$-подмодулей $V_i$. 
\end{theorem} 
\begin{proof}
Воспользуемся примерами~\ref{ExampleIntegralFG}, \ref{ExampleComoduleGraded}  и теоремой~\ref{TheoremHcoWeyl}.
\end{proof}
\begin{corollary}
Пусть $L$~"--- конечномерная полупростая алгебра Ли над полем характеристики $0$,
градуированная произвольной группой, а $(V, \psi)$~"---
конечномерный градуированный $L$-модуль.
 Тогда $$V=V_1 \oplus V_2 \oplus \ldots \oplus V_s$$
для некоторых неприводимых градуированных $L$-подмодулей $V_i$. 
\end{corollary}

\begin{theorem}\label{TheoremHWeyl}
Пусть $L$~"--- $H$-(ко)модульная алгебра Ли над полем характеристики $0$,
$H$~"--- конечномерная (ко)полупростая алгебра Хопфа,
а $(V, \psi)$~"--- конечномерный $(H,L)$-модуль, вполне приводимый как обычный $L$-модуль
(без учёта $H$-(ко)действия).
 Тогда $$V=V_1 \oplus V_2 \oplus \ldots \oplus V_s$$
 для некоторых неприводимых $(H,L)$-подмодулей $V_i$. 
\end{theorem}
\begin{proof}
Достаточно использовать двойственность между $H$-действиями и $H^*$-кодействиями,
применить пример~\ref{ExampleIntegralHSS} и теорему~\ref{TheoremHcoWeyl}.
\end{proof}
\begin{corollary}\label{CorollaryHWeylLSS}
Пусть $L$~"--- конечномерная полупростая $H$-(ко)модульная
 алгебра Ли над полем характеристики $0$,
$H$~"--- конечномерная (ко)полупростая алгебра Хопфа,
а $(V, \psi)$~"--- конечномерный $(H,L)$-модуль.
 Тогда $$V=V_1 \oplus V_2 \oplus \ldots \oplus V_s$$
 для некоторых неприводимых $(H,L)$-подмодулей $V_i$. 
\end{corollary}

\begin{theorem}\label{TheoremAffAlgGrWeyl}
Пусть $L$~"--- конечномерная алгебра Ли над алгебраически замкнутым полем $\mathbbm{k}$
характеристики $0$, на которой рационально действует автоморфизмами
редуктивная аффинная алгебраическая группа $G$.
Пусть $(V, \psi)$~"--- конечномерный $(G,L)$-модуль
с рациональным $G$-действием,
вполне приводимый как обычный $L$-модуль
без учёта $G$-действия.
 Тогда $$V=V_1 \oplus V_2 \oplus \ldots \oplus V_s$$
для некоторых неприводимых $(G,L)$-подмодулей $V_i$. 
\end{theorem}
\begin{proof}
Воспользуемся примерами~\ref{ExampleIntegralAffAlgGr}, \ref{ExampleRegularAffActAutAlgebra}
и теоремой~\ref{TheoremHcoWeyl}.
\end{proof}
\begin{corollary}
Пусть $L$~"--- конечномерная полупростая алгебра Ли над алгебраически замкнутым полем $\mathbbm{k}$
характеристики $0$, на которой рационально действует автоморфизмами
редуктивная аффинная алгебраическая группа $G$.
Пусть $(V, \psi)$~"--- конечномерный $(G,L)$-модуль
с рациональным $G$-действием.
 Тогда $$V=V_1 \oplus V_2 \oplus \ldots \oplus V_s$$
для некоторых неприводимых $(G,L)$-подмодулей $V_i$. 
\end{corollary}

\section{$H$-(ко)инвариантное разложение разрешимого радикала}\label{SectionHLBQN}

Сперва рассмотрим случай, когда алгебра Хопфа $H$ необязательно конечномерна.

\begin{theorem}\label{TheoremHcoLBQN}
Пусть $L$~"--- конечномерная $H$-комодульная алгебра Ли над полем $\mathbbm{k}$ характеристики $0$,
где $H$~"--- алгебра Хопфа с $\ad$-инвариантным левым интегралом $t \in H^*$, где $t(1)=1$.
Предположим, что нильпотентный радикал $N$ и разрешимый радикал $R$ являются $H$-подкомодулями.
Обозначим через $B$ такую $H$-коинвариантную максимальную полупростую подалгебру, что $L=B \oplus R$.
Тогда существует такой $H$-подкомодуль $Q$, что $R = Q \oplus N$ и $[B,Q]=0$. В частности, $L=B\oplus Q\oplus N$.
\end{theorem}
\begin{proof}
Рассмотрим присоединённое представление подалгебры $B$ на $L$. Тогда $L$ является $(H,B)$-модулем.
В силу следствия~\ref{CorollaryHWeylLSS}
этот $(H,B)$-модуль вполне приводим.
Более того, идеалы $N$ и $R$ являются $(H,B)$-подмодулями $(H,B)$-модуля $L$.
Отсюда существует такой $(H,B)$-подомодуль $Q$, что $R = Q \oplus N$.
В частности, $[B, Q] \subseteq Q$. 
Однако в силу предложения 2.1.7 из \cite{GotoGrosshans}
справедливо включение
$[B, Q] \subseteq [L,R]\subseteq N$. Отсюда $[B,Q]=0$.
\end{proof}

Используя двойственность между $H$-действиями и $H^*$-кодействиями
и применяя пример~\ref{ExampleIntegralHSS}, следствие~\ref{CorollaryHModRadicalsLie}
и теорему~\ref{TheoremHcoLBQN}, получаем слудующий результат:

\begin{theorem}\label{TheoremHLBQN}
Пусть $L$~"--- конечномерная полупростая $H$-(ко)модульная
 алгебра Ли над полем характеристики $0$, а
$H$~"--- конечномерная (ко)полупростая алгебра Хопфа.
Обозначим через $N$ нильпотентный, а через $R$~"--- разрешимый радикалы алгебры Ли $L$.
Тогда существуюет такой $H$-под(ко)модуль $Q$, что $R = Q \oplus N$,
$L=B\oplus Q\oplus N$ (прямая сумма $H$-под(ко)модулей) и $[B,Q]=0$,
где $B$~"--- $H$-(ко)инвариантная максимальная полупростая подалгебра.
\end{theorem}
 
Используя примеры~\ref{ExampleIntegralFG}, \ref{ExampleComoduleGraded},
следствие~\ref{CorollaryLieRadicalsGraded} и теорему~\ref{TheoremHcoLBQN},
получаем:

\begin{theorem}\label{TheoremGradLBQN}
Пусть $L$~"--- конечномерная $G$-градуированная алгебра Ли над полем характеристики $0$. Обозначим через $N$ нильпотентный, а через $R$~"--- разрешимый радикалы алгебры Ли $L$.
Тогда существует такое градуированное подпространство $Q$, что $R = Q \oplus N$,
$L=B\oplus Q\oplus N$ (прямая сумма градуированных подпространств) и $[B,Q]=0$,
где $B$~"--- градуированная максимальная полупростая подалгебра.
\end{theorem}

Также получаем следующую теорему:

\begin{theorem}\label{TheoremAffAlgGrLBQN}
Пусть редуктивная аффинная алгебраическая группа $G$ рационально действует 
автоморфизмами на конечномерной алгебре Ли $L$ над алгебраически замкнутым полем характеристики $0$.
Обозначим через $N$ нильпотентный, а через $R$~"--- разрешимый радикалы алгебры Ли $L$.
Тогда существует такое $G$-инвариантное подпространство $Q$, что $R = Q \oplus N$,
$L=B\oplus Q\oplus N$ (прямая сумма $G$-инвариантных подпространств) и $[B,Q]=0$,
где $B$~"--- $G$-инвариантная максимальная полупростая подалгебра.
\end{theorem}
\begin{proof}
Достаточно заметить, что $R$ и $N$ инвариантны
относительно всех автоморфизмов, и использовать примеры~\ref{ExampleIntegralAffAlgGr}, \ref{ExampleRegularAffActAutAlgebra} и теорему~\ref{TheoremHcoLBQN}.
\end{proof}

     \section{Полупростые $H_{m^2}(\zeta)$-простые алгебры Ли}\label{SectionClassTaftSLieSS}
     
     Оставшаяся часть данной главы посвящена изучению конечномерных $H_{m^2}(\zeta)$-простых
     алгебр Ли, где $H_{m^2}(\zeta)$~"--- алгебра Тафта. При этом тождество антикоммутативности, 
     которое выполняется во всех алгебрах Ли, а также тот факт, что любая простая алгебра Ли неприводима как модуль над собой, накладывают  по сравнению со случаем ассоциативных алгебр дополнительные ограничения.
Если ассоциативные полупростые $H_{m^2}(\zeta)$-простые алгебры
      параметризуются матрицами, а число минимальных идеалов в разложении может быть, вообще говоря, произвольным делителем числа $m$, то полупростые $H_{m^2}(\zeta)$-простые алгебры Ли параметризуются элементами основного поля,
      а число минимальных идеалов в разложении должно совпадать с $m$, если только, конечно, $v\in H_{m^2}(\zeta)$ не действует на алгебре Ли как нулевой оператор.     
      Если ассоциативные неполупростые $H_{m^2}(\zeta)$-простые алгебры строятся на основе $\mathbb Z/m\mathbb Z$-градуированно простых алгебр $B$, то
      неполупростые $H_{m^2}(\zeta)$-простые алгебры Ли строятся на основе алгебр Ли $B$, простых в
      обычном смысле. 

В этом параграфе классифицируются полупростые $H_{m^2}(\zeta)$-простые алгебры Ли,
которые не являются простыми в обычном смысле.

Напомним, что если на некоторой алгебре $A$ действует алгебра Тафта~$H_{m^2}(\zeta)$, то циклическая группа $C_m :=\langle c \rangle_m$, где $c \in H_{m^2}(\zeta)$, действует на $A$ автоморфизмами.
Существование $C_m$-действия эквивалентно наличию $\mathbb Z/m\mathbb Z$-градуировки
$A=\bigoplus_{k=0}^{m-1} A^{(k)}$, причём действие и градуировка связаны соотношением
$ca=\zeta^k a$ для всех $a\in A^{(k)}$.

Пусть $B$~"--- простая алгебра Ли над полем $\mathbbm{k}$.
Предположим, что $\mathbbm{k}$
содержит примитивный корень $\zeta$ 
степени $m$ из единицы. 
Пусть $\alpha \in \mathbbm{k}$. Введём обозначение
$$L_\alpha(B) := \underbrace{B\oplus \ldots \oplus B}_{m} \text{ (прямая сумма идеалов)}$$
и определим действие элементов $c,v \in H_{m^2}(\zeta)$ на алгебре Ли $L_\alpha(B)$
при помощи формул
\begin{equation}\label{EqLieTaftSimpleSSDefc}c(a_1, a_2, \ldots, a_{m-1}, a_m) = (a_m, a_1, a_2, \ldots, a_{m-1})\end{equation}
и \begin{equation}\label{EqLieTaftSimpleSSDefv}v(a_1, \ldots, a_m)=\alpha\left(a_1 - a_m, \zeta(a_2 - a_1), \ldots, \zeta^{m-1}(a_m - a_{m-1})\right)\end{equation} для всех $a_1, \ldots, a_m \in B$.

Аналогично тому, как это было сделано
в лемме~\ref{LemmaTaftSimpleSemisimpleFormula},
для произвольных $a_1, \ldots, a_m \in B$
получаем $$v^\ell (a_1, a_2, \ldots, a_m) =(b_1, b_2, \ldots, b_m),$$
где \begin{equation*} b_k = \alpha^\ell \zeta^{\ell(k-1)} \sum_{j=0}^\ell  (-1)^j \zeta^{-\frac{j(j-1)}{2}} \binom{\ell}{j}_{\zeta^{-1}}
a_{k-j} \end{equation*}
и $a_{-j} := a_{m-j} $ при $j \geqslant 0$.
В доказательстве теоремы~\ref{TheoremTaftSimpleSemisimple}
было получено равенство $(-1)^m \zeta^{-\frac{m(m-1)}{2}}=-1$, из которого теперь следует, что $v^m (a_1, a_2, \ldots, a_m) = 0$ для всех $a_i \in B$.
Непосредственная проверка показывает, что формулы~(\ref{EqLieTaftSimpleSSDefc}) и~(\ref{EqLieTaftSimpleSSDefv})
задают  $H_{m^2}(\zeta)$-действие на $L_\alpha(B)$ корректно. Поскольку коммутатор любого идеала
алгебры Ли $L_\alpha(B)$ со любой копией алгебры Ли $B$ либо равен $0$, либо совпадает
с этой копией, алгебра Ли $L_\alpha(B)$ является $C_m$-простой, а следовательно, $\mathbb Z/m\mathbb Z$-градуированно простой  и $H_{m^2}(\zeta)$-простой алгеброй Ли.

Другое описание алгебр Ли $L_\alpha(B)$ при $\alpha \ne 0$
будет дано ниже в теореме~\ref{TheoremTaftSimpleLieEquivDef}.

\begin{theorem}\label{TheoremTaftSimpleSemiSimpleLieClassify}
Пусть $L$~"--- конечномерная полупростая $H_{m^2}(\zeta)$-простая $H_{m^2}(\zeta)$-модульная  алгебра Ли
над алгебраически замкнутым полем $\mathbbm{k}$ характеристики $0$.
Предположим, что $L$ полупростая, но не простая алгебра Ли. Тогда
$L$ является $\mathbb Z/m\mathbb Z$-градуированно простой алгеброй Ли.
Если $vL \ne 0$, то $L \cong L_\alpha(B)$ для некоторой простой алгеброй Ли $B$ и элемента $\alpha \in \mathbbm{k}$.
\end{theorem}
\begin{proof}
В силу полупростоты алгебры $L$ существует разложение $L=B_1 \oplus \ldots \oplus B_t$ (прямая сумма идеалов), где $B_i$~"--- простые алгебры Ли. Тогда для любого $1\leqslant i \leqslant t$
существует такое число $1\leqslant j(i) \leqslant t$, что
$cB_i = B_{j(i)}$.
Кроме того, $v[a,b]=[ca, vb]+[va, b] \in B_i \oplus B_{j(i)}$
для всех $a,b \in B_i$. Поскольку $[B_i, B_i] = B_i$, получаем $vB_i \subseteq  B_i \oplus B_{j(i)}$.
В частности, идеал $\sum_{k=0}^{m-1} c^k B_1$ инвариантен относительно действия обоих элементов $c$ и $v$.
В силу того, что $L$ является $H_{m^2}(\zeta)$-простой алгеброй Ли, $\sum_{k=0}^{m-1} c^k B_1 = L$
и, очевидно, алгебра Ли $L$ является $C_m$-простой и $\mathbb Z/m\mathbb Z$-градуированно простой.
Без ограничения общности можно считать, что $B_i = c^{i-1} B_1$.
Тогда в силу единственности минимальных идеалов $c^t B_1 = B_1$.

Обозначим через $\pi_i \colon B \twoheadrightarrow B_i$ естественные сюръективные гомоморфизмы.
Определим $\rho_i \colon B_i \to B_i$ и $\theta_i \colon B_i \to cB_i$
при помощи формул $\rho_i(a)=\pi_i(va)$ и $\theta_i(a)=\pi_{i+1}(va)$
для всех $a\in B_i$. (В доказательстве 
теоремы будем считать, что все нижние индексы определены по модулю $t$, например, $\pi_{t+1} := \pi_1$.)
Тогда $$\rho_i[a,b]=\pi_i(v[a,b])=\pi_i([ca,vb]+[va,b])=[\rho_i(a),b].$$
Аналогично, $$\theta_i[a,b]=\pi_{i+1}(v[a,b])=\pi_{i+1}([ca,vb]+[va,b])=[ca, \theta_i(b)]$$
для всех $a,b\in B_i$.
В частности, как $\rho_i$, так и $\theta_i$
являются гомоморфизмами $B_i$-модулей.

Поскольку $B_i$ являются простыми алгебрами Ли, $B_i$ неприводимы как $B_i$-модули.
Отсюда в силу леммы Шура справедливы равенства
$\rho_i = \alpha_i \id_{B_i}$
и $\theta_i = \beta_i \left(c\bigr|_{B_i}\right)$ для некоторых $\alpha_i, \beta_i \in \mathbbm{k}$.
Поскольку $vc=\zeta cv$, получаем $$\alpha_{i+1} ca =\rho_{i+1}(ca)=\zeta\pi_{i+1}(c(v(a)))=
\zeta c(\rho_i(a))=\zeta \alpha_i ca$$
и $$\beta_{i+1} c^2 a =\theta_{i+1}(ca)=\zeta\pi_{i+2}(c(v(a)))=
\zeta c(\theta_i(a))=\zeta \beta_i c^2 a$$ для всех $1\leqslant i \leqslant t$
и $a\in B_i$.
 Следовательно, $\alpha_i = \zeta^{i-1} \alpha_1$
и $\beta_i = \zeta^{i-1} \beta_1$ для всех $1\leqslant i \leqslant t$.
Более того, если хотя бы один из элементов $\alpha_1$ и $\beta_1$ ненулевой, то $\zeta^t = 1$
и $t=m$.

Заметим, что для всех $1\leqslant i \leqslant t$, $a \in B_i$ и $b \in B_{i+1}$
справедливо равенство $$0=v[a,b]=[ca,vb]+[va, b]=[ca,\rho_{i+1}(b)]+[\theta_i(a),b]=(\alpha_{i+1}+\beta_i) [ca, b].$$
Поскольку $[B_{i+1}, B_{i+1}]=B_{i+1}$,
получаем $\beta_i = -\alpha_{i+1}$ для всех $1\leqslant i \leqslant t$.

Если $\alpha_1 = 0$, то $vL=0$ и теорема доказана.
Предположим, что $\alpha_1 \ne 0$. Тогда $t=m$.
Поскольку $B_i = c^{i-1} B_1$ и $va=\rho_i(a)+\theta_i(a)$ для всех $a\in B_i$, 
можно отождествить все $B_i$ между собой  и считать, что $L=\underbrace{B\oplus \ldots \oplus B}_{m}$ (прямая сумма идеалов)
для простой алгебры Ли $B:=B_1$,
причём для $\alpha := \alpha_1$
справедливы равенства (\ref{EqLieTaftSimpleSSDefc}) и~(\ref{EqLieTaftSimpleSSDefv}).
\end{proof}
\begin{remark}
Если $vL=0$, то из доказательства теоремы~\ref{TheoremTaftSimpleSemiSimpleLieClassify} 
следует, что существуют такое число $t\in\mathbb N$, где $t \mid m$,
и простая алгебра Ли $B$ с действием циклической группы 
порядка $\frac{m}{t}$ с порождающим $d$, что
$$L \cong \underbrace{B\oplus \ldots \oplus B}_{t} \text{ (прямая сумма идеалов)},$$
\begin{equation*}c(a_1, a_2, \ldots, a_{m-1}, a_m) = (d a_m, a_1, a_2, \ldots, a_{m-1})\end{equation*}
и \begin{equation*}v(a_1, \ldots, a_m)=0\end{equation*} для всех $a_1, \ldots, a_m \in B$.
\end{remark}

\medskip

В теореме~\ref{TheoremTaftSimpleSemiSimpleLieIso}, которую мы доказываем
ниже, приводятся необходимые и достаточные условия для
существования изоморфизма $H_{m^2}(\zeta)$-модульных
алгебр Ли $L_{\alpha_1}(B_1) \cong L_{\alpha_2}(B_2)$.

\begin{theorem}\label{TheoremTaftSimpleSemiSimpleLieIso}
Пусть $B_1, B_2$~"--- простые алгебры Ли над полем $\mathbbm{k}$, $\alpha_1, \alpha_2 \in \mathbbm{k}$,
а $\zeta$~"--- примитивный корень из единицы степени $m$.
Пусть $\theta \colon L_{\alpha_1}(B_1) \mathrel{\widetilde{\to}} L_{\alpha_2}(B_2)$~"---
изоморфизм алгебр Ли и $H_{m^2}(\zeta)$-модулей.
Тогда существует такое число $0\leqslant k \leqslant m-1$ и изоморфизм алгебр Ли
$\varphi \colon B_1 \mathrel{\widetilde{\to}} B_2$, что
\begin{equation}\label{EqLieTaftSimpleSSDefTheta}\theta(b_1, \ldots, b_m)=(\varphi(b_{k+1}), \ldots, \varphi(b_m), \varphi(b_1), \ldots, \varphi(b_k))\end{equation}
для всех $b_i \in B_1$.
Кроме того, $\alpha_2 = \zeta^k\alpha_1$.
Обратно, если $B_1 \cong B_2$ как обычные алгебры Ли и $\alpha_2 = \zeta^k\alpha_1$ для
некоторого $k\in\mathbb Z$, то $L_{\alpha_1}(B_1) \cong L_{\alpha_2}(B_2)$ как $H_{m^2}(\zeta)$-модульные
алгебры Ли.
\end{theorem}
\begin{proof}
Заметим, что любой минимальный идеал алгебры $L_{\alpha_2}(B_2)$ совпадает
с одной из копий алгебры $B_2$. Отсюда либо существует такое $1\leqslant k \leqslant m-1$,
что $$\theta(B_1, 0, \ldots, 0)=
(\underbrace{0,\ldots, 0}_{m-k}, B_2, 0, \ldots, 0),$$
либо $$\theta(B_1, 0, \ldots, 0)=
(B_2, 0, \ldots, 0).$$ В последнем случае полагаем $k:= 0$.

Обозначим индуцированный изоморфизм $B_1 \mathrel{\widetilde{\to}} B_2$ через $\varphi$. Тогда
$$\theta(b, 0, \ldots, 0)=
(\underbrace{0,\ldots, 0}_{m-k}, \varphi(b), 0, \ldots, 0)$$
для всех $b\in B$.
Теперь~(\ref{EqLieTaftSimpleSSDefTheta}) следует из~(\ref{EqLieTaftSimpleSSDefc})
и равенства $\theta(ca)=c\theta(a)$, которое выполняется для всех $a\in L_{\alpha_1}(B_1)$.
Используя~(\ref{EqLieTaftSimpleSSDefv}) 
и то, что $\theta(va)=v\theta(a)$ для всех $a\in L_{\alpha_1}(B_1)$,
получаем $\alpha_2 = \zeta^k\alpha_1$.
Таким образом, прямое утверждение доказано.
Обратное утверждение очевидно.
\end{proof}
\begin{remark}
В частности, при $\alpha \ne 0$ все автоморфизмы $H_{m^2}(\zeta)$-модульной
алгебры Ли $L_{\alpha}(B)$ индуцированы автоморфизмами обычной алгебры Ли $B$
и соответствующие группы автоморфизмов $\Aut(L_{\alpha}(B))$ и $\Aut(B)$
могут быть отождествлены друг с другом.
Если $\alpha = 0$, то $\Aut(L_{\alpha}(B)) \cong \Aut(B) \times \mathbb Z/m\mathbb Z$.
\end{remark}

\section{Алгебры Ли $L(B, \gamma)$ и $H_{m^2}(\zeta)$-действия на простых алгебрах Ли}
\label{SectionClassTaftSLieSimple}

Следующим шагом в классификации конечномерных $H_{m^2}(\zeta)$-простых алгебр Ли
является изучение $H_{m^2}(\zeta)$-действий на простых алгебрах Ли.
Как мы покажем ниже в теореме~\ref{TheoremTaftSimpleSimpleLieClassify},
во всех конечномерных простых алгебрах Ли, наделённых $H_{m^2}(\zeta)$-действием,
справедливо равенство $va=0$ для всех $a\in L$.
Для того, чтобы доказать этот результат, определим $H_{m^2}(\zeta)$-простые алгебры Ли $L(B, \gamma)$.

\begin{theorem}\label{TheoremTaftSimpleLiePresent} Пусть $B$~"---
простая алгебра Ли над полем $\mathbbm{k}$,
а $\gamma\in \mathbbm{k}$~"--- некоторый элемент поля.
Предположим, что $\mathbbm{k}$
содержит примитивный корень $\zeta$ степени $m$
из единицы. Рассмотрим векторные пространства $L^{(i)}$, где $1\leqslant i \leqslant m-1$,
которые являются копиями векторного пространства $L^{(0)} := B$. Обозначим
одной и той же буквой $\psi$ соответствующие $\mathbbm{k}$-линейные биекции $L^{(i-1)} \mathrel{\widetilde\to} L^{(i)}$, где $1 \leqslant i \leqslant m-1$.
Пусть $\psi(L^{(m-1)}):=0$.
Рассмотрим $H_{m^2}(\zeta)$-модуль $L(B,\gamma) :=\bigoplus_{i=0}^{m-1} L^{(i)}$ (прямая сумма подпространств), где $v\psi(a):=a$ для всех $a \in L^{(i)}$ и
$0\leqslant i \leqslant m-2$, $vB:=0$ и $c a^{(i)}:=\zeta^i a^{(i)}$,  $a^{(i)} \in L^{(i)}$.
  Определим лиевский коммутатор на $L(B,\gamma)$ по формуле
  \begin{equation}\label{EqMultTaftSimpleLiePresent}
  [\psi^k(a),\psi^\ell(b)]:=\left\lbrace
\begin{array}{rrr}
  \binom{k+\ell}{k}_\zeta\ \psi^{k+\ell}[a,b]  & \text{при} & k+\ell < m,\\
  \gamma\frac{(k+\ell-m)!_\zeta}{k!_\zeta \ell!_\zeta}\ \psi^{k+\ell-m}[a,b] 
  & \text{при} & k+\ell \geqslant m
  \end{array}\right.
  \end{equation}
  для всех $a, b\in B$ и $0 \leqslant k,\ell < m$.
    Тогда $L(B,\gamma)$~"--- $H_{m^2}(\zeta)$-простая алгебра Ли.
\end{theorem}
\begin{proof}
То, что приведённые выше формулы действительно задают на $L(B,\gamma)$
структуру $H_{m^2}(\zeta)$-модульной алгебры Ли, доказывается непосредственно.
Проверим только, что $v[u,w]=[v_{(1)}u, v_{(2)}w]$
при $u,w\in L(B,\gamma)$.

Пусть $0\leqslant k,\ell < m$, а $a,b\in B$.
Если $k=\ell=0$, то $v[a,b]=0=[ca,vb]+[va,b]$.

Если $k=0$, $\ell > 0$, то
  $$v[a, \psi^\ell(b)]=v \psi^\ell[a,b]= \psi^{\ell-1}[a,b]=
 [a,\psi^{\ell-1}(b)]
 =  [ca, v\psi^\ell(b)]+[va, \psi^\ell(b)].$$

Если $k > 0$, $\ell=0$, 
 то $$v[\psi^k(a), b]=v \psi^k[a,b]= \psi^{k-1}[a,b]=
 [\psi^{k-1}(a), b]
 =  [c\psi^k(a), vb]+[v\psi^k(a), b].$$

Если $k,\ell > 0$, $k+\ell < m$, то
 \begin{equation*}\begin{split}v[\psi^k(a),\psi^\ell(b)] = \binom{k+\ell}{k}_\zeta \psi^{k+\ell-1}[a,b]
 = \left(\zeta^k\binom{k+\ell-1}{k}_\zeta+\binom{k+\ell-1}{k-1}_\zeta\right) \psi^{k+\ell-1}[a,b]=\\= [c\psi^k(a),\psi^{\ell-1}(b)] + [\psi^{k-1}(a),\psi^\ell(b)]= [c\psi^k(a),v\psi^\ell(b)] + [v\psi^k(a),\psi^\ell(b)],\end{split}\end{equation*}
поскольку \begin{equation}\begin{split}\label{EquationMainQuantumBinomial}
\zeta^k \binom{k+\ell-1}{k}_\zeta + \binom{k+\ell-1}{k-1}_\zeta
= \frac{\zeta^k (k+\ell-1)!_\zeta}{k!_\zeta (\ell-1)!_\zeta}+
\frac{(k+\ell-1)!_\zeta}{(k-1)!_\zeta \ell!_\zeta}=\\=(\zeta^k \ell_\zeta + k_\zeta)\frac{(k+\ell-1)!_\zeta}{k!_\zeta \ell!_\zeta}=(k+\ell)_\zeta \frac{(k+\ell-1)!_\zeta}{k!_\zeta \ell!_\zeta}
= \frac{(k+\ell)!_\zeta}{k!_\zeta \ell!_\zeta}=\binom{k+\ell}{k}_\zeta.\end{split}\end{equation}

Пусть $k+\ell \geqslant m$. Тогда $k,\ell > 0$.
Если $k+\ell =m$, то $(k+\ell)_\zeta = m_\zeta = 0$
и
 \begin{equation*}\begin{split}v[\psi^k(a),\psi^\ell(b)] = \frac{\gamma}{k!_\zeta \ell!_\zeta}\ v[a,b]= 0 = \binom{k+\ell}{k}_\zeta \psi^{k+\ell-1}[a,b]
 =\\= \left(\zeta^k\binom{k+\ell-1}{k}_\zeta+\binom{k+\ell-1}{k-1}_\zeta\right) \psi^{k+\ell-1}[a,b]= [c\psi^k(a),\psi^{\ell-1}(b)] + [\psi^{k-1}(a),\psi^\ell(b)]=\\= [c\psi^k(a),v\psi^\ell(b)] + [v\psi^k(a),\psi^\ell(b)].\end{split}\end{equation*}

Если $k+\ell >m$, то $(k+\ell)_\zeta = (k+\ell-m)_\zeta + \zeta^{k+\ell-m} m_\zeta = (k+\ell-m)_\zeta$
и
 \begin{equation*}\begin{split}v[\psi^k(a),\psi^\ell(b)] = \gamma\frac{(k+\ell-m)!_\zeta}{k!_\zeta \ell!_\zeta}\ \psi^{k+\ell-m-1}[a,b] 
= \\=
\gamma(k+\ell)_\zeta\frac{(k+\ell-m-1)!_\zeta}{k!_\zeta \ell!_\zeta}\ \psi^{k+\ell-m-1}[a,b] 
= \\= 
  \gamma\left(\zeta^k \frac{(k+\ell-m-1)!_\zeta}{k!_\zeta (\ell-1)!_\zeta}+\frac{(k+\ell-m-1)!_\zeta}{(k-1)!_\zeta \ell!_\zeta}\right) \psi^{k+\ell-m-1}[a,b]=\\= [c\psi^k(a),\psi^{\ell-1}(b)] + [\psi^{k-1}(a),\psi^\ell(b)]=\\= [c\psi^k(a),v\psi^\ell(b)] + [v\psi^k(a),\psi^\ell(b)].\end{split}\end{equation*}

Таким образом, были рассмотрены все возможные варианты  при $0\leqslant k,\ell < m$.
Отсюда $v[u,w]=[v_{(1)}u, v_{(2)}w]$
для всех $u,w\in L(B,\gamma)$.

Предположим, что $I$~"--- $H_{m^2}(\zeta)$-инвариантный идеал алгебры Ли $L$. Тогда $v^m I = 0$.
Обозначим через $t \in \mathbb Z_+$ такое число, что $v^t I \ne 0$, а $v^{t+1} I = 0$. Тогда $0 \ne v^t I \subseteq I \cap \ker v$. Однако $\ker v = B$~"--- простая алгебра Ли. Следовательно, $I \cap \ker v =\ker v$ и $\ker v \subseteq I$. Поскольку $$[\ker v, L^{(i)}]=[B,\psi^i(B)]=\psi^i[B,B]=\psi^i(B)=L^{(i)}\text{ для всех }1\leqslant i \leqslant m-1,$$
получаем, что $I = L(B,\gamma)$. Следовательно, $L(B,\gamma)$ является
  $H_{m^2}(\zeta)$-простой алгеброй Ли.
\end{proof}
\begin{remark}
Алгебра Ли $L(B,0)$ не является полупростой. Разрешимый радикал алгебры Ли $L(B,0)$
совпадает с нильпотентным радикалом
и равен $\bigoplus_{i=1}^{m-1} L^{(i)}$.
\end{remark}

Из теоремы~\ref{TheoremTaftSimpleLieEquivDef},
которая доказывается ниже, следует, что если поле $\mathbbm{k}$ алгебраически замкнуто и $\gamma \ne 0$,
то $L(B, \gamma)$ как $H_{m^2}(\zeta)$-модульная алгебра Ли
изоморфна одной из непростых $\mathbb Z/m\mathbb Z$-градуированно простых $H_{m^2}(\zeta)$-модульных алгебр Ли $L_\alpha(B)$, определённых в \S\ref{SectionClassTaftSLieSS}.

\begin{theorem}\label{TheoremTaftSimpleLieEquivDef}
Пусть $B$~"--- простая алгебра Ли над полем $\mathbbm{k}$.
Предположим, что поле $\mathbbm{k}$ содержит примитивный корень $\zeta$ степени $m$
из единицы. Пусть $\alpha \in \mathbbm{k}$, $\alpha \ne 0$.
Тогда $L(B, \frac{1}{\alpha^m (1-\zeta)^m}) \cong L_\alpha(B)$ 
как $H_{m^2}(\zeta)$-модульные алгебра Ли.
 \end{theorem}
\begin{proof}
Заметим, что $$L_\alpha(B)^{(k)}=\left\lbrace\left(b, \zeta^{-k} b, \zeta^{-2k} b, \ldots, \zeta^{-(m-1)k}b\right)
\mathrel{\bigl|} b \in B\right\rbrace$$ для всех $0\leqslant k \leqslant m-1$. В частности, $L_\alpha(B)^{(0)}\cong B$.
Положим $$\psi(b, \zeta^{-k} b, \zeta^{-2k} b, \ldots, \zeta^{-(m-1)k}b)
:= \frac{1}{\alpha(1-\zeta^{k+1})}\left(b, \zeta^{-(k+1)} b, \zeta^{-2(k+1)} b, \ldots, \zeta^{-(m-1)(k+1)}b\right)$$
для всех $b\in B$ и $0\leqslant k < m-1$.
Тогда  \begin{equation*}\begin{split}\psi^k(b,b,\ldots, b) = \frac{1}{\alpha^k(1-\zeta)(1-\zeta^2)\ldots (1-\zeta^k)}
(b, \zeta^{-k} b, \zeta^{-2k} b, \ldots, \zeta^{-(m-1)k}b)=\\=\frac{1}{\alpha^k(1-\zeta)^k k!_\zeta}
(b, \zeta^{-k} b, \zeta^{-2k} b, \ldots, \zeta^{-(m-1)k}b).\end{split}\end{equation*}
Заметим, что
$L_\alpha(B)^{(k)}=\psi^k(L_\alpha(B)^{(0)})$ и $v\psi(a)=a$ для всех $a\in L_\alpha(B)^{(k)}$, где $0\leqslant k < m-1$. 
Более того, элемент $[\psi^k(a),\psi^\ell(b)]$ может быть вычислен с использованием~(\ref{EqMultTaftSimpleLiePresent})
при $\gamma=\frac{1}{\alpha^m (1-\zeta)^m}$ для всех $a,b \in L_\alpha(B)^{(0)}$ и $0\leqslant k,\ell < m$.
Отсюда $L_\alpha(B) \cong L(B, \frac{1}{\alpha^m (1-\zeta)^m})$ как $H_{m^2}(\zeta)$-модульные алгебры Ли.
\end{proof}

Докажем теперь несколько лемм об $H_{m^2}(\zeta)$-модульных
алгебрах Ли.

\begin{lemma}\label{LemmaLieTaftSimpleSGradVGrad}
Пусть $L$~"--- $H_{m^2}(\zeta)$-модульная алгебра Ли над полем $\mathbbm{k}$. 
Тогда
\begin{equation}\label{EqLieTaftSimpleSGradMultLR}(\zeta^k-1)[a^{(k)},v b^{(\ell)}] = (\zeta^\ell-1) [v a^{(k)}, b^{(\ell)}],\end{equation}
\begin{equation}\label{EqLieTaftSimpleSGradMultProd}(\zeta^\ell-1)v[a^{(k)},b^{(\ell)}]=(\zeta^{k+\ell}-1) [a^{(k)},v b^{(\ell)}]\end{equation}
для всех $a^{(k)}\in L^{(k)}$, $b^{(\ell)}\in L^{(\ell)}$,
где $L=\bigoplus_{k=0}^{m-1} L^{(k)}$~"--- $\mathbb Z/m\mathbb Z$-градуировка,
индуцированная действием группы $C_m=\langle c \rangle_m$.
Более того, если, будучи наделённой этой градуировкой, алгебра Ли $L$ является $\mathbb Z/m\mathbb Z$-градуированно простой алгеброй Ли, то $vL^{(0)}=0$.
\end{lemma}
\begin{proof}
Заметим, что  $$v[a^{(k)},b^{(\ell)}]=[ca^{(k)}, vb^{(\ell)}]+[va^{(k)}, b^{(\ell)}]=\zeta^k [a^{(k)},v b^{(\ell)}]
+[va^{(k)}, b^{(\ell)}]$$ для всех $a^{(k)}\in L^{(k)}$
и $b^{(\ell)} \in L^{(\ell)}$.
В то же время 
\begin{equation*}\begin{split}v[a^{(k)},b^{(\ell)}]=-v[b^{(\ell)},a^{(k)}]=-[cb^{(\ell)}, va^{(k)}]-[vb^{(\ell)}, a^{(k)}]=\\ = -\zeta^\ell [b^{(\ell)},v a^{(k)}]
-[vb^{(\ell)}, a^{(k)}]=[a^{(k)}, vb^{(\ell)}]+\zeta^\ell [v a^{(k)}, b^{(\ell)}].\end{split}\end{equation*}
Отсюда получаем~(\ref{EqLieTaftSimpleSGradMultLR})
и~(\ref{EqLieTaftSimpleSGradMultProd}).
Если $\ell \ne 0$, то мы также получаем
$$ v[a^{(k)},b^{(\ell)}]=\frac{\zeta^{k+\ell}-1}{\zeta^\ell-1} [a^{(k)},v b^{(\ell)}].
$$

В частности, $v[L^{(k)}, L^{(m-k)}]=0$ для всех $1\leqslant k \leqslant m-1$.

Предположим теперь, что $L$~"--- $\mathbb Z/m\mathbb Z$-градуированно простая алгебра Ли.
Если $L= L^{(0)}$, то элемент $c$ действует как тождественный оператор и из $vc=\zeta cv$
следует $vL=0$.
Отсюда без ограничения общности можно считать, что $L\ne L^{(0)}$.
Пусть $a \in L^{(k)}$, $k\ne 0$, $a\ne 0$.
Поскольку алгебра Ли $L$ является $\mathbb Z/m\mathbb Z$-градуированно простой, а
элемент $a$ однороден, идеал, порождённый элементом $a$, совпадает с $L$.
Следовательно, $L^{(0)}$ является $\mathbbm{k}$-линейной оболочкой элементов
$[a, a_1, \ldots, a_n]$, где $a_i \in L^{(k_i)}$, $0\leqslant k_i \leqslant m-1$,
$m \mid (k+k_1+\ldots+k_n)$, $n\in \mathbb N$.
(Здесь мы используем длинные коммутаторы $[x_1, \ldots, x_n]
:= [[\ldots [[x_1, x_2], x_3], \ldots], x_n]$.)
Если $k_n \ne 0$,
то из $[a, a_1, \ldots, a_n]=[[a, a_1, \ldots, a_{n-1}], a_n]$
следует, что $[a, a_1, \ldots, a_n]\in [L^{(m-k_n)}, L^{(k_n)}]$
и $v[a, a_1, \ldots, a_n]=0$.
В случае $k_n = 0$
применим тождество Якоби
и перепишем $[a, a_1, \ldots, a_n]$ как сумму элементов
$[[a,a_n], a_1, \ldots, a_{n-1}]$, $[a, [a_1,a_n], \ldots, a_{n-1}]$ и
 $[a, a_1, \ldots, [a_{n-1},a_n]]$. Если $a_{n-1} \in L^{(0)}$, то мы продолжим эту процедуру.
В конце концов мы придём к ситуации, когда элементы,
стоящие в длинных коммутаторах на последних местах,
принадлежат $L^{(k_i)}$, $k_i \ne 0$.
Применяя тот же приём, что и выше,
получим $v[a, a_1, \ldots, a_n]=0$. 
 Таким образом, $v L^{(0)}=0$.
\end{proof}

\begin{lemma}\label{LemmaLieTaftSimpleSGradTrois}
Пусть $L$~"--- $H_{m^2}(\zeta)$-модульная алгебра Ли над полем $\mathbbm{k}$ характеристики $\chr \mathbbm{k} \nmid m$, $\chr \mathbbm{k} \ne 2$. Пусть $a^{(\ell)} \in L^{(\ell)}$, $b^{(k)}\in L^{(k)}$, $u^{(m-k)}\in L^{(m-k)}$ для некоторых $1\leqslant k, \ell \leqslant m-1$.
Тогда $$v[a^{(\ell)},[b^{(k)}, u^{(m-k)}]] = \frac{\zeta^\ell-1}{\zeta^{m-k}-1}[a^{(\ell)},[b^{(k)}, vu^{(m-k)}]].$$
\end{lemma}
\begin{proof} В силу~(\ref{EqLieTaftSimpleSGradMultLR}), (\ref{EqLieTaftSimpleSGradMultProd}) и тождества Якоби
\begin{equation*}\begin{split}[a^{(\ell)},[b^{(k)}, vu^{(m-k)}]]=-\frac{\zeta^{m-k}-1}{\zeta^{k}-1}[a^{(\ell)},[u^{(m-k)}, vb^{(k)}]]
=\\ =-\frac{\zeta^{m-k}-1}{\zeta^{k}-1}([[a^{(\ell)},u^{(m-k)}], vb^{(k)}]+ [u^{(m-k)}, [a^{(\ell)},vb^{(k)}]]) =\\ =-\frac{\zeta^{m-k}-1}{\zeta^{\ell}-1}v[a^{(\ell)},u^{(m-k)}, b^{(k)}]
+\frac{\zeta^{m-k}-1}{\zeta^{\ell}-1}[u^{(m-k)}, [b^{(k)},va^{(\ell)}]]
=\\=-\frac{\zeta^{m-k}-1}{\zeta^{\ell}-1}v[a^{(\ell)},u^{(m-k)}, b^{(k)}]
+\frac{\zeta^{m-k}-1}{\zeta^{\ell}-1}([[u^{(m-k)}, b^{(k)}],va^{(\ell)}]
+ [b^{(k)},[u^{(m-k)}, va^{(\ell)}]])=\\=\frac{\zeta^{m-k}-1}{\zeta^{\ell}-1}v\left(-[a^{(\ell)},u^{(m-k)}, b^{(k)}]+[u^{(m-k)}, b^{(k)},a^{(\ell)}]\right)-[b^{(k)},[a^{(\ell)}, vu^{(m-k)}]]=\\=
\frac{\zeta^{m-k}-1}{\zeta^{\ell}-1}v\left(-[a^{(\ell)},u^{(m-k)}, b^{(k)}]+[u^{(m-k)}, b^{(k)},a^{(\ell)}]\right)-\\-[[b^{(k)},a^{(\ell)}], vu^{(m-k)}]-[a^{(\ell)}, [b^{(k)},vu^{(m-k)}]]
=\\=\frac{\zeta^{m-k}-1}{\zeta^{\ell}-1}v\left(-[a^{(\ell)},u^{(m-k)}, b^{(k)}]+[u^{(m-k)}, b^{(k)},a^{(\ell)}]-[b^{(k)},a^{(\ell)}, u^{(m-k)}]\right)-\\-[a^{(\ell)}, [b^{(k)},vu^{(m-k)}]]
=\\= 2\frac{\zeta^{m-k}-1}{\zeta^{\ell}-1} v[a^{(\ell)}, [b^{(k)}, u^{(m-k)}]]-[a^{(\ell)}, [b^{(k)},vu^{(m-k)}]].
\end{split}\end{equation*} Таким образом, лемма доказана.
\end{proof}

\begin{lemma}\label{LemmaLieTaftSimpleSGradVMany}
Пусть $L$~"--- $H_{m^2}(\zeta)$-модульная алгебра Ли над полем $\mathbbm{k}$ характеристики $\chr \mathbbm{k} \nmid m$, $\chr \mathbbm{k} \ne 2$. 
Пусть $s \geqslant 2$, $0\leqslant k_i \leqslant m-1$ при $1\leqslant i \leqslant s$,
$k_s > 0$, а $a_i^{(k_i)} \in L_i^{(k_i)}$ при $1\leqslant i \leqslant s$.
Тогда 
\begin{equation}\label{EqLieTaftSimpleSGradVMany}
v[a_1^{(k_1)}, [a_2^{(k_2)}, \ldots, [a_{s-1}^{(k_{s-1})}, a_s^{(k_s)}]\ldots]
 = \frac{\zeta^{\sum_{i=1}^s k_i}-1}{\zeta^{k_s}-1}[a_1^{(k_1)}, [a_2^{(k_2)}, \ldots, [a_{s-1}^{(k_{s-1})}, v a_s^{(k_s)}]\ldots].\end{equation}
\end{lemma}
\begin{proof}
Проведём доказательство индукцией по $s$. База $s=2$ является следствием равенства~(\ref{EqLieTaftSimpleSGradMultProd}).
Предположим, что $s > 2$.

 Если $m \nmid \sum_{i=2}^s k_i$, то
 в силу~(\ref{EqLieTaftSimpleSGradMultProd}) и предположения индукции \begin{equation*}\begin{split}v[a_1^{(k_1)}, [a_2^{(k_2)}, \ldots, [a_{s-1}^{(k_{s-1})}, a_s^{(k_s)}]\ldots]
 =  \frac{\zeta^{\sum_{i=1}^s k_i}-1}{\zeta^{\sum_{i=2}^s k_i}-1}[a_1^{(k_1)}, v[a_2^{(k_2)}, \ldots, [a_{s-1}^{(k_{s-1})}, a_s^{(k_s)}]\ldots]
  = \\ =\frac{\zeta^{\sum_{i=1}^s k_i}-1}{\zeta^{k_s}-1}[a_1^{(k_1)}, [a_2^{(k_2)}, \ldots, [a_{s-1}^{(k_{s-1})}, v a_s^{(k_s)}]\ldots]\end{split}\end{equation*}
и в этом случае лемма доказана.

Предположим, что $k_i = 0$ для некоторого $1\leqslant i < s$.
Тогда в силу тождества Якоби (символ $\widehat{a_i^{(k_i)}}$
обозначает пропуск элемента $a_i^{(k_i)}$ в выражении), \begin{equation*}\begin{split}a:=v[a_1^{(k_1)}, [a_2^{(k_2)}, \ldots, [a_{s-1}^{(k_{s-1})}, a_s^{(k_s)}]\ldots]  = \\
 = \sum_{j=i+1}^{s-1} v[a_1^{(k_1)}, [a_2^{(k_2)}, \ldots, [\widehat{a_i^{(k_i)}}, \ldots, [a_{j-1}^{(k_{j-1})}, [[a_i^{(k_i)}, a_j^{(k_j)}],
 [a_{j+1}^{(k_{j+1})}, \ldots [a_{s-1}^{(k_{s-1})}, a_s^{(k_s)}]\ldots]
 +\\+ v[a_1^{(k_1)}, [a_2^{(k_2)}, \ldots, [\widehat{a_i^{(k_i)}}, \ldots, [a_{s-1}^{(k_{s-1})}, [a_i^{(k_i)}, a_s^{(k_s)}]]\ldots].\end{split}\end{equation*}

Поскольку $k_i=0$, всякий элемент $[a_i^{(k_i)}, a_j^{(k_j)}]$ снова имеет степень $k_j$ в $\mathbb Z/m\mathbb Z$-градуировке.
Будем рассматривать такой элемент $[a_i^{(k_i)}, a_j^{(k_j)}]$ как единое целое.
Применяя предположение индукции для $s-1$, получаем
\begin{equation*}\begin{split}a 
 =  \frac{\zeta^{\sum_{i=1}^s k_i}-1}{\zeta^{k_s}-1} \cdot \\ \cdot\left(\sum_{j=i+1}^{s-1} [a_1^{(k_1)}, [a_2^{(k_2)}, \ldots, [\widehat{a_i^{(k_i)}}, \ldots, [a_{j-1}^{(k_{j-1})}, [[a_i^{(k_i)}, a_j^{(k_j)}],
 [a_{j+1}^{(k_{j+1})}, \ldots [a_{s-1}^{(k_{s-1})}, v a_s^{(k_s)}]\ldots]
 +\right. \\ \left.+ [a_1^{(k_1)}, [a_2^{(k_2)}, \ldots, [\widehat{a_i^{(k_i)}}, \ldots, [a_{s-1}^{(k_{s-1})}, v[a_i^{(k_i)}, a_s^{(k_s)}]]\ldots]\right).\end{split}\end{equation*}
 Из~(\ref{EqLieTaftSimpleSGradMultProd}) и тождества Якоби следует, что
  \begin{equation*}\begin{split}a 
 = \frac{\zeta^{\sum_{i=1}^s k_i}-1}{\zeta^{k_s}-1} \cdot \\ \cdot\left(\sum_{j=i+1}^{s-1} [a_1^{(k_1)}, [a_2^{(k_2)}, \ldots, [\widehat{a_i^{(k_i)}}, \ldots, [a_{j-1}^{(k_{j-1})}, [[a_i^{(k_i)}, a_j^{(k_j)}],
 [a_{j+1}^{(k_{j+1})}, \ldots [a_{s-1}^{(k_{s-1})}, v a_s^{(k_s)}]\ldots]
  +\right. \\ \left.+[a_1^{(k_1)}, [a_2^{(k_2)}, \ldots, [\widehat{a_i^{(k_i)}}, \ldots, [a_{s-1}^{(k_{s-1})}, [a_i^{(k_i)}, v a_s^{(k_s)}]]\ldots]\right)
  =\\= \frac{\zeta^{\sum_{i=1}^s k_i}-1}{\zeta^{k_s}-1}[a_1^{(k_1)}, [a_2^{(k_2)}, \ldots, [a_{s-1}^{(k_{s-1})}, v a_s^{(k_s)}]\ldots],\end{split}\end{equation*}
 т.е. лемма доказана и в этом случае.

Единственный случай, который ещё остался неразобранным,~"--- это когда $1\leqslant k_i \leqslant m-1$
для всех $i$, а $m \mid \sum_{i=2}^s k_i$. Однако в этом случае $m \mid \left(k_2 + \sum_{i=3}^s k_i\right)$,
и мы можем применить лемму~\ref{LemmaLieTaftSimpleSGradTrois}:
\begin{equation*}\begin{split} v[a_1^{(k_1)}, [a_2^{(k_2)}, [a_3^{(k_3)} \ldots, [a_{s-1}^{(k_{s-1})}, a_s^{(k_s)}]\ldots] = \\ = \frac{\zeta^{\sum_{i=1}^s k_i}-1}{\zeta^{\sum_{i=3}^s k_i}-1}[a_1^{(k_1)}, [a_2^{(k_2)}, v\bigl[a_3^{(k_3)} \ldots, [a_{s-1}^{(k_{s-1})}, a_s^{(k_s)}]\ldots\bigr]]]. 
 \end{split}\end{equation*}
 Теперь утверждение леммы следует из предположения индукции для $s-2$.\end{proof}

\begin{lemma}\label{LemmaLieTaftSimpleSGradVKer}
Пусть $L$~"--- $H_{m^2}(\zeta)$-модульная алгебра Ли над полем $\mathbbm{k}$ характеристики $\chr \mathbbm{k} \nmid m$, $\chr \mathbbm{k} \ne 2$. Предположим, что $L$~"--- $\mathbb Z/m\mathbb Z$-градуированно простая алгебра Ли, $vL\ne 0$.
Тогда $\ker v = L^{(0)}$.
\end{lemma}
\begin{proof}
В силу леммы~\ref{LemmaLieTaftSimpleSGradVGrad} справедливо равенство $vL^{(0)}=0$. Поскольку $vc=\zeta cv$, подпространство $\ker v$ является $\mathbb Z/m\mathbb Z$-градуированным. Предположим, что $\ker v \supsetneqq L^{(0)}$. Тогда существует такой элемент $a^{(k)} \subseteq L^{(k)}$, где
$1\leqslant k \leqslant m-1$, что $a^{(k)}\ne 0$, но $va^{(k)}=0$.
Поскольку алгебра Ли $L$ является $\mathbb Z/m\mathbb Z$-градуированно простой, а элемент $a^{(k)}$
однороден, справедливо равенство $L=\sum_{n \geqslant 0}[\underbrace{L, [L,\ldots [L,}_{n} \mathbbm{k} a^{(k)}]\ldots]$. Теперь из леммы~\ref{LemmaLieTaftSimpleSGradVMany}
следует, что $vL=0$.
\end{proof}

\begin{lemma}\label{LemmaLieTaftSimpleSGradVIm}
Пусть $L$~"--- $H_{m^2}(\zeta)$-модульная алгебра Ли над полем $\mathbbm{k}$ характеристики $\chr \mathbbm{k} \nmid m$, $\chr \mathbbm{k} \ne 2$. Предположим, что $L$~"--- $\mathbb Z/m\mathbb Z$-градуированно простая алгебра Ли, $vL\ne 0$.
Тогда  $v L^{(k)}=L^{(k-1)}$ для всех $1\leqslant k \leqslant m-1$.
\end{lemma}
\begin{proof}
Во-первых, докажем, что $L = vL + \sum_{k= 1}^{m-1} [L^{(k)}, vL^{(m-k)}]$.
Заметим, что $$I=vL + \sum_{k= 1}^{m-1} [L^{(k)}, vL^{(m-k)}]$$~"--- $\mathbb Z/m\mathbb Z$-градуированное подпространство. Утверждается, что $I$ является идеалом. 

В силу~(\ref{EqLieTaftSimpleSGradMultProd}) получаем
$[L^{(\ell)}, vL^{(k)}] \subseteq vL$ для всех $0\leqslant k,\ell < m$,
таких, что $m \nmid (k+\ell)$.
Отсюда $[L, vL] \subseteq I$.

Теперь докажем, что $\left[L, \sum_{k= 1}^{m-1} [L^{(k)}, vL^{(m-k)}]\right] \subseteq I$.

Во-первых, в силу тождества Якоби \begin{equation}\label{EqLieTaftSimpleSGradJacobi}[a^{(\ell)}, [b^{(k)}, vu^{(m-k)}]]
= [[a^{(\ell)}, b^{(k)}], vu^{(m-k)}]+ [ b^{(k)}, [a^{(\ell)},vu^{(m-k)}]]\end{equation}
для всех $0\leqslant k,\ell < m$ и $a^{(\ell)} \in L^{(\ell)}$, $b^{(k)} \in L^{(k)}$, $u^{(m-k)}\in L^{(m-k)}$.

Если $\ell = 0$, то из~(\ref{EqLieTaftSimpleSGradMultProd})
и~(\ref{EqLieTaftSimpleSGradJacobi}) получаем
$$[a^{(\ell)}, [b^{(k)}, vu^{(m-k)}]]
= [[a^{(\ell)}, b^{(k)}], vu^{(m-k)}]+ [ b^{(k)}, v[a^{(\ell)},u^{(m-k)}]] \in [L^{(k)}, vL^{(m-k)}].$$

Если $\ell\ne 0$ и $k \ne \ell$, то из~(\ref{EqLieTaftSimpleSGradMultProd})
и~(\ref{EqLieTaftSimpleSGradJacobi}) получаем
\begin{equation*}\begin{split}[a^{(\ell)}, [b^{(k)}, vu^{(m-k)}]]=
 [[a^{(\ell)}, b^{(k)}], vu^{(m-k)}]+[ b^{(k)}, [a^{(\ell)},vu^{(m-k)}]] = \\ =
 \frac{\zeta^{m-k}-1}{\zeta^{\ell}-1}v\left([[a^{(\ell)}, b^{(k)}], u^{(m-k)}]+[ b^{(k)}, [a^{(\ell)},u^{(m-k)}]]\right) \in vL.\end{split}\end{equation*}

Предположим, что $\ell \ne 0$ и $k=\ell$.
В этом случае нужно доказать, что $$[L^{(k)}, [L^{(k)},vL^{(m-k)}]] \subseteq vL.$$

Если $m \ne 2k$, то $m-k \ne k$. Из~(\ref{EqLieTaftSimpleSGradMultLR}), (\ref{EqLieTaftSimpleSGradMultProd}) и тождества Якоби следует, что
\begin{equation*}\begin{split}[L^{(k)}, [L^{(k)},vL^{(m-k)}]] \subseteq [L^{(k)}, [L^{(m-k)}, vL^{(k)}]]
\subseteq \\ \subseteq [[L^{(k)}, L^{(m-k)}], vL^{(k)}]
+ [L^{(m-k)}, [L^{(k)}, vL^{(k)}]] \subseteq vL. \end{split}\end{equation*}

Если $m=2k$, то $\chr \mathbbm{k} \ne 2$ и включение $[L^{(k)}, [L^{(k)},vL^{(k)}]] \subseteq vL$
является следствием леммы~\ref{LemmaLieTaftSimpleSGradTrois}.
 
 Отсюда $I$ действительно является $\mathbb Z/m\mathbb Z$-градуированным идеалом и $$L = vL + \sum_{k\ne 0} [L^{(k)}, vL^{(m-k)}],$$ поскольку алгебра Ли
$L$  градуированно проста.

Поскольку $vc = \zeta cv$, справедливо включение $vL^{(k)} \subseteq L^{(k-1)}$, откуда $$\sum_{k\ne 0} [L^{(k)}, vL^{(m-k)}]\subseteq L^{(m-1)}.$$ Следовательно, $\bigoplus_{k=0}^{m-2} L^{(k)} \subseteq vL$.
В силу леммы~\ref{LemmaLieTaftSimpleSGradVGrad} справедливо равенство $vL^{(0)}=0$,
откуда
$vL \cap L^{(m-1)}=0$ и $vL = \bigoplus_{k=0}^{m-2} L^{(k)}$. В частности, $v L^{(k)}=L^{(k-1)}$ для всех $1\leqslant k \leqslant m-1$.
\end{proof}

\begin{lemma}\label{LemmaLieTaftSimpleSGradPhiMult}
Пусть $L$~"--- $H_{m^2}(\zeta)$-модульная алгебра Ли над полем $\mathbbm{k}$ характеристики $\chr \mathbbm{k} \nmid m$, $\chr \mathbbm{k} \ne 2$. Предположим, что $L$~"--- $\mathbb Z/m\mathbb Z$-градуированно простая алгебра Ли, $vL\ne 0$.
Определим отображения $\psi \colon L^{(k)} \to L^{(k+1)}$ (мы будем обозначать их одной и той же буквой)
 при помощи равенства $\psi(va)=a$, где $a \in L^{(k+1)}$, $0\leqslant k \leqslant m-2$.
Положим $\lbrace a, b \rbrace := (m-1)!_\zeta [\psi(a), \psi^{m-1}(b)]$.
Тогда 
  \begin{equation}\label{EqLieTaftSimpleSGradPhiMult}
  [\psi^k(a),\psi^\ell(b)]:=\left\lbrace
\begin{array}{rrr}
  \binom{k+\ell}{k}_\zeta\ \psi^{k+\ell}[a,b]  & \text{при} & k+\ell < m,\\
  \frac{(k+\ell-m)!_\zeta}{k!_\zeta \ell!_\zeta}\ \psi^{k+\ell-m}\lbrace a,b \rbrace 
  & \text{при} & k+\ell \geqslant m
  \end{array}\right.
  \end{equation}
  для всех $a, b\in L^{(0)}$ и $0 \leqslant k,\ell < m$.
\end{lemma}
\begin{proof}
В силу лемм~\ref{LemmaLieTaftSimpleSGradVKer} и~\ref{LemmaLieTaftSimpleSGradVIm}
отображение
$\psi$ определено корректно. Более того,
 $v\psi(a)=a$ для всех $a\in L^{(i)}$, $0\leqslant i \leqslant m-2$.

Если $k=\ell=0$, то справедливость~(\ref{EqLieTaftSimpleSGradPhiMult}) очевидна.
В силу леммы~\ref{LemmaLieTaftSimpleSGradVGrad}
для всех $1\leqslant k \leqslant m-1$ и $a,b\in L^{(0)}$ справедливо равенство $v[\psi^k(a),b]=[\psi^{k-1}(a),b]$.
Следовательно, $[\psi^k(a),b]=\psi^k[a,b]$ и
в случае, когда одно из чисел $k,\ell$ равно нулю, равенство~(\ref{EqLieTaftSimpleSGradPhiMult}) также доказано.

В случае произвольных $k,\ell > 0$, где $k+\ell < m$, равенство~(\ref{EqLieTaftSimpleSGradPhiMult})
доказывается по индукции с использованием~(\ref{EquationMainQuantumBinomial}):
\begin{equation*}\begin{split}
[\psi^k(a),\psi^\ell(b)]= \psi \left(v[\psi^k(a),\psi^\ell(b)]\right)
=\psi([c\psi^k(a), \psi^{\ell-1}(b)]+[\psi^{k-1}(a),\psi^\ell(b)])= \\
=\psi\left([\zeta^k\psi^k(a), \psi^{\ell-1}(b)]+[\psi^{k-1}(a),\psi^\ell(b)]\right)= \\
=\psi\left(\zeta^k\binom{k+\ell-1}{k}_\zeta\ \psi^{k+\ell-1}[a,b]+\binom{k+\ell-1}{k-1}_\zeta\ \psi^{k+\ell-1}[a,b] \right)= \\
=\left(\zeta^k \binom{k+\ell-1}{k}_\zeta + \binom{k+\ell-1}{k-1}_\zeta\right) \psi^{k+\ell}[a,b]
= \\ =\binom{k+\ell}{k}_\zeta\ \psi^{k+\ell}[a,b].
\end{split}\end{equation*}

Предположим, что $k+\ell = m$.
Докажем равенство~(\ref{EqLieTaftSimpleSGradPhiMult}) индукцией по $k$.
 Если $k=1$ и $\ell = m-1$, то~(\ref{EqLieTaftSimpleSGradPhiMult})
следует из определения операции $\lbrace , \rbrace$.
Если $k>1$,
то $\ell < m-1$. В силу~(\ref{EqLieTaftSimpleSGradMultLR})
и предположения индукции при $k-1$ получаем
\begin{equation*}\begin{split}[\psi^k(a),\psi^\ell(b)]=
[\psi^k(a),v\psi^{\ell+1}(b)]
=\frac{\zeta^{\ell+1}-1}{\zeta^k-1}
[v\psi^k(a), \psi^{\ell+1}(b)]=\frac{\zeta^{\ell+1}-1}{\zeta^k-1}[\psi^{k-1}(a), \psi^{\ell+1}(b)]
=\\=\frac{\zeta^{\ell+1}-1}{(\zeta^k-1)(k-1)!_\zeta (\ell+1)!_\zeta } \lbrace a,b\rbrace =
\frac{(\ell+1)_\zeta}{k_\zeta(k-1)!_\zeta (\ell+1)!_\zeta } \lbrace a,b\rbrace=\frac{1}{k!_\zeta \ell!_\zeta } \lbrace a,b\rbrace.\end{split}\end{equation*}

Если $k+\ell > m$, то мы воспользуемся индукцией по $(k+\ell)$:
\begin{equation*}\begin{split}
[\psi^k(a),\psi^\ell(b)]= \psi \left(v[\psi^k(a),\psi^\ell(b)]\right)
=\psi([c\psi^k(a), \psi^{\ell-1}(b)]+[\psi^{k-1}(a),\psi^\ell(b)])=\\
=\psi\left([\zeta^k\psi^k(a), \psi^{\ell-1}(b)]+[\psi^{k-1}(a),\psi^\ell(b)]\right)=\\
=\psi\left(\zeta^k\frac{(k+\ell-m-1)!_\zeta}{k!_\zeta
(\ell-1)!_\zeta}\ \psi^{k+\ell-m-1}\lbrace a,b\rbrace+\frac{(k+\ell-m-1)!_\zeta}{(k-1)!_\zeta
\ell!_\zeta}\ \psi^{k+\ell-m-1}\lbrace a,b\rbrace \right)=\\=
\left(\zeta^k \frac{(k+\ell-m-1)!_\zeta}{k!_\zeta
(\ell-1)!_\zeta} + \frac{(k+\ell-m-1)!_\zeta}{(k-1)!_\zeta
\ell!_\zeta}\right) \psi^{k+\ell-m}\lbrace a,b\rbrace
= \\ =\frac{(k+\ell)_\zeta (k+\ell-m-1)!_\zeta}{k!_\zeta
\ell!_\zeta}\ \psi^{k+\ell-m}\lbrace a,b\rbrace=\frac{(k+\ell-m)!_\zeta}{k!_\zeta
\ell!_\zeta}\ \psi^{k+\ell-m}\lbrace a,b\rbrace,
\end{split}\end{equation*}
поскольку $(k+\ell)_\zeta=(k+\ell-m)_\zeta$ при $m<k+\ell < 2m$.
\end{proof}

\begin{lemma}\label{LemmaLieTaftSimpleSGradIsoToLBGamma}
Пусть $L$~"--- конечномерная $H_{m^2}(\zeta)$-модульная алгебра Ли над алгебраически замкнутым полем $\mathbbm{k}$ характеристики $0$. Предположим, что $L$~"--- $\mathbb Z/m\mathbb Z$-градуированно простая алгебра Ли, $vL\ne 0$.
Тогда $L^{(0)}$~"--- простая алгебра Ли и существует такое $\gamma \in \mathbbm{k}$, что $\gamma \ne 0$
и $\lbrace a, b \rbrace = \gamma [a,b]$ для всех $a,b \in L^{(0)}$.
Другими словами, $L \cong L(L^{(0)}, \gamma)$ как $H_{m^2}(\zeta)$-модульные алгебры Ли.
\end{lemma}
\begin{proof}
Заметим, что в силу тождества Якоби $$[\psi^{m-1}(a), \psi(b), u]
+ [\psi(b), u, \psi^{m-1}(a)]
+ [u, \psi^{m-1}(a), \psi(b)]=0$$
для всех $a,b,u \in L^{(0)}$.
Теперь из леммы~\ref{LemmaLieTaftSimpleSGradPhiMult}
следует, что 
\begin{equation}\label{EqLieTaftSimpleSGradIsoToLBGamma1}
[\lbrace a, b \rbrace, u]+\lbrace [b,u], a \rbrace
+ \lbrace [u,a], b \rbrace = 0.
\end{equation}

Снова используя тождество Якоби, получаем $$[\psi^{m-1}(a), \psi^{m-1}(b), \psi(u)]
+ [\psi^{m-1}(b), \psi(u), \psi^{m-1}(a)]
+ [\psi(u), \psi^{m-1}(a), \psi^{m-1}(b)]=0$$
для всех $a,b,u \in L^{(0)}$.
В силу леммы~\ref{LemmaLieTaftSimpleSGradPhiMult}
отсюда следует, что
\begin{equation}\label{EqLieTaftSimpleSGradIsoToLBGamma2}
[\lbrace a, b \rbrace, u]+[\lbrace b,u \rbrace, a]
+ [\lbrace u, a \rbrace, b] = 0.
\end{equation}

Меняя в~(\ref{EqLieTaftSimpleSGradIsoToLBGamma1}) значения $a,b,u$ местами,
получаем
$$[\lbrace b, u \rbrace, a]+\lbrace [u,a], b \rbrace
+ \lbrace [a,b], u \rbrace = 0, $$
$$[\lbrace u, a \rbrace, b]+\lbrace [a,b], u \rbrace
+ \lbrace [b,u], a \rbrace = 0.$$
Складывая эти равенства с~(\ref{EqLieTaftSimpleSGradIsoToLBGamma1})
и используя равенство~(\ref{EqLieTaftSimpleSGradIsoToLBGamma2}), а также условие $\chr \mathbbm{k} \ne 2$,
получаем
$$\lbrace [a,b], u \rbrace +\lbrace [b,u], a \rbrace
+ \lbrace [u,a], b \rbrace =0$$ для всех $a,b,u \in L^{(0)}$.
В силу~(\ref{EqLieTaftSimpleSGradIsoToLBGamma1})
отсюда следует, что \begin{equation}\label{EqLieTaftSimpleSGradIsoToLBGamma3}\lbrace [a,b], u \rbrace = [\lbrace a,b\rbrace, u].\end{equation}
В силу леммы~\ref{LemmaLieTaftSimpleSGradPhiMult}
справедливо равенство
$\lbrace a, b \rbrace = -\lbrace b,a \rbrace$.
Используя~(\ref{EqLieTaftSimpleSGradIsoToLBGamma2}) и~(\ref{EqLieTaftSimpleSGradIsoToLBGamma3}),
получаем отсюда, что
$$\lbrace [a,b], u \rbrace = [\lbrace a, u\rbrace, b] + [a, \lbrace b, u \rbrace] $$ для всех $a,b,u \in L^{(0)}$.
Другими словами, $\lbrace \cdot, u \rbrace$ является дифференцированием
для всех $u \in L^{(0)}$.

Докажем теперь, что $L^{(0)}$~"--- простая алгебра Ли.
Предположим сперва, что $I \ne 0$~"--- идеал алгебры Ли $L^{(0)}$,
такой, что $\lbrace I, u \rbrace \subseteq I$ для всех $u \in L^{(0)}$.
Тогда в силу леммы~\ref{LemmaLieTaftSimpleSGradPhiMult} подпространство
$\bigoplus_{k=0}^{m-1} \psi^k(I)$ является ненулевым $\mathbb Z/m\mathbb Z$-градуированным идеалом
алгебры Ли $L$.
Отсюда $L = \bigoplus_{k=0}^{m-1} \psi^k(I)$ и $I=L^{(0)}$.
В силу теоремы 7 из~\cite[глава~III, \S6]{JacobsonLie} (это также следует из
теоремы~\ref{TheoremHModRadicalsLie}) разрешимый и нильпотентный радикалы
алгебры~$L^{(0)}$ 
инвариантны относительно всех дифференцирований.
Отсюда алгебра Ли $L^{(0)}$  либо полупроста, либо
нильпотентна. 

Предположим, что алгебра Ли $L^{(0)}$ нильпотентна.
Пусть \begin{equation}\label{EqLieTaftSimpleSGradIsoToLBGammaNilp1}[a_1, \ldots, a_t] = 0\end{equation} для всех
$a_1, \ldots, a_t \in L^{(0)}$.
Докажем, что \begin{equation}\label{EqLieTaftSimpleSGradIsoToLBGammaNilp2}[\psi^{s_1}(a_1), \ldots, \psi^{s_{m(t+1)}}(a_{m(t+1)})] = 0\end{equation}
для всех $a_i \in L^{(0)}$ и $0\leqslant s_i \leqslant m-1$, где $1\leqslant i \leqslant m(t+1)$.
Действительно, согласно лемме~\ref{LemmaLieTaftSimpleSGradPhiMult}
$$[\psi^{s_1}(a_1), \ldots, \psi^{s_{m(t+1)}}(a_{m(t+1)})]=\alpha \psi^{s_1+\ldots+s_{m(t+1)}-\ell m}(b),$$
где $\alpha \in \mathbbm{k}$, число $\ell \in \mathbb Z_+$ подобрано
таким образом,  чтобы обеспечить справедливость неравенства $$0 \leqslant s_1+\ldots+s_{m(t+1)}-\ell m \leqslant m-1,$$
а $b$ является смешанным коммутатором длины $m(t+1)$, являющимся композицией
$m(t+1)-\ell-1$ операций $[\cdot,\cdot]$ и $\ell$ операций $\lbrace \cdot,\cdot \rbrace$.
При этом $$m(t+1)-\ell-1 \geqslant m(t+1)-\frac{s_1+\ldots+s_{m(t+1)}}{m}-1 \geqslant m(t+1)-\frac{m(t+1)(m-1)}{m}- 1 \geqslant t.$$
Учитывая, что в силу~\eqref{EqLieTaftSimpleSGradIsoToLBGamma3}
мы можем менять операции $[\cdot,\cdot]$ и $\lbrace \cdot,\cdot \rbrace$
местами, без ограничения общности 
можно считать, что в выражении $b$ подряд идёт не менее $t$ коммутаторов $[\cdot,\cdot]$, откуда $b=0$ 
в силу~\eqref{EqLieTaftSimpleSGradIsoToLBGammaNilp1}.
Следовательно, справедливо равенство~\eqref{EqLieTaftSimpleSGradIsoToLBGammaNilp2} 
и алгебра Ли $L$ нильпотентна. Получаем противоречие с $\mathbb Z / m \mathbb Z$-простотой
алгебры Ли $L$.

Итак, алгебра  Ли $L^{(0)}$ полупроста.
Если алгебра Ли $L^{(0)}$ не является простой, то 
$L^{(0)}=I_1 \oplus I_2$ для некоторых ненулевых идеалов $I_1$ и $I_2$.
Пусть $\delta$~"--- дифференцирование алгебры Ли $L^{(0)}$.
Тогда $\delta(I_i)=\delta[I_i, I_i]\subseteq[\delta(I_i), I_i]+[I_i, \delta(I_i)]\subseteq I_i$.
Другими словами, идеалы $I_i$, где $i=1,2$, инвариантны 
относительно всех дифференцирований, откуда $L^{(0)}=I_1 = I_2$,
т.е. мы получаем противоречие. Отсюда $L^{(0)}$~"--- простая алгебра Ли.

Поскольку алгебра Ли $L^{(0)}$ проста, все её дифференцирования внутренние (см., например, \cite[теорема 5.3]{HumphreysLieAlg}).
Следовательно, существует такое $\mathbbm{k}$-линейное отображение $\theta \colon L^{(0)} \to L^{(0)}$,
что $\lbrace a,b \rbrace = [a, \theta(b)]$
для всех $a,b\in L^{(0)}$. Теперь из~(\ref{EqLieTaftSimpleSGradIsoToLBGamma3})
следует, что $$[\theta[a,b], u] = \lbrace [a,b], u \rbrace
= [\lbrace a,b\rbrace, u] = [[a,\theta(b)], u]$$
для всех $a,b,u \in L^{(0)}.$
Поскольку центр алгебры Ли $L^{(0)}$ нулевой, справедливо равенство $\theta[a,b] = [a,\theta(b)]$
для всех $a,b \in L^{(0)}$. Другими словами,
$\theta \colon L^{(0)} \to L^{(0)}$~"--- гомоморфизм
 $L^{(0)}$-модулей. Поскольку $L^{(0)}$~"--- неприводимый
 $L^{(0)}$-модуль, в силу леммы Шура $\theta$ является скалярным отображением, т.е. $\lbrace a,b \rbrace = \gamma[a, b]$
 для некоторого $\gamma \in \mathbbm{k}$. Из леммы~\ref{LemmaLieTaftSimpleSGradPhiMult}
 теперь следует, что $L$ как $H_{m^2}(\zeta)$-модульная алгебра Ли
 изоморфна алгебре Ли
 $L(L^{(0)}, \gamma)$.
 Поскольку $L$ является $\mathbb Z/m\mathbb Z$-градуированно простой, а следовательно, и полупростой,
 выполнено условие $\gamma \ne 0$.
\end{proof}

Ниже в теореме~\ref{TheoremTaftSimpleSimpleLieClassify}
показывается, что всякая конечномерная $H_{m^2}(\zeta)$-модульная
алгебра Ли, простая в обычном смысле, является
$\mathbb Z/m\mathbb Z$-градуированной алгеброй Ли
с тривиальным $v$-действием.

\begin{theorem}\label{TheoremTaftSimpleSimpleLieClassify}
Пусть $L$~"--- конечномерная $H_{m^2}(\zeta)$-модульная алгебра Ли над алгебраически замкнутым полем $\mathbbm{k}$ характеристики $0$. Предположим, что $L$~"--- простая (в обычном смысле) алгебра Ли. Тогда $vL = 0$.
\end{theorem}
\begin{proof} 
Предположим, что $vL\ne 0$.
Тогда в силу леммы~\ref{LemmaLieTaftSimpleSGradIsoToLBGamma}
существует изоморфизм $H_{m^2}(\zeta)$-модульных алгебр Ли $L \cong L(L^{(0)}, \gamma)$ 
для некоторого $\gamma \ne 0$.
Согласно теореме~\ref{TheoremTaftSimpleLieEquivDef}
алгебра Ли
$L(L^{(0)}, \gamma)$ изоморфна алгебре Ли $L_\alpha(L^{(0)})$,
где $\alpha =  \frac{(1-\zeta)^{-1}}{\sqrt[m]{\gamma}}$. Однако $L_\alpha(L^{(0)})$
не проста как обычная алгебра Ли. Отсюда получаем противоречие
и $vL=0$. \end{proof}

\section{Неполупростые $H_{m^2}(\zeta)$-простые алгебры Ли}\label{SectionClassTaftSLieNSS}

В данном параграфе доказывается, что все конечномерные неполупростые
 $H_{m^2}(\zeta)$-простые алгебры Ли
 изоморфны алгебрам Ли из теоремы~\ref{TheoremTaftSimpleLiePresent} при $\gamma = 0$.

\begin{theorem}\label{TheoremTaftSimpleNonSemiSimpleLieClassify}
Пусть $L$~"--- конечномерная $H_{m^2}(\zeta)$-простая $H_{m^2}(\zeta)$-модульная алгебра Ли над алгебраически замкнутым полем $\mathbbm{k}$ характеристики $0$, причём разрешимый радикал алгебры Ли $L$ не равен $0$. Тогда $L$ изоморфна как $H_{m^2}(\zeta)$-модульная алгебра Ли алгебре Ли $L(B, 0)$
для некоторой конечномерной простой алгебры Ли $B$.
\end{theorem}

Для того, чтобы доказать теорему~\ref{TheoremTaftSimpleNonSemiSimpleLieClassify},
нам потребуется несколько лемм.

Пусть $M_1,M_2$~"---  два $\mathbb Z/m\mathbb Z$-градуированных 
модуля над $\mathbb Z/m\mathbb Z$-градуированной алгеброй Ли $L$. 
Будем говорить, что $\mathbbm{k}$-линейная биекция $\psi \colon M_1 \mathrel{\widetilde\to} M_2$ является \textit{$c$-изоморфизмом}
из $M_1$ в $M_2$, если существует такое число $r\in\mathbb Z$, что $c\psi(b) = \zeta^{-r} \psi(cb)$ и $\psi(ab)=a\psi(b)$ для всех $b\in M_1$, $a\in L$.

Напомним, что для любой конечномерной алгебры Ли $L$
над полем характеристики $0$ справедливо включение
$[L,R]\subseteq N$ (см., например, \cite[предложение 2.1.7]{GotoGrosshans}),
где $R$, $N$~"--- соответственно,  разрешимый и нильпотентный радикалы. 
Отсюда, если $N=0$, то $R \subseteq Z(L)\subseteq N=0$,
где $Z(L)$~"--- центр алгебры Ли $L$. Напомним также,
что если $L$~"--- $H_{m^2}(\zeta)$-модульная алгебра Ли, то $R$ и $N$ являются $\mathbb Z/m\mathbb Z$-градуированными идеалами. Например, это следует из того, что $R$ и $N$ инвариантны относительно всех автоморфизмов алгебры Ли $L$ и, в частности, относительно $c$-действия. 

\begin{lemma}\label{LemmaTaftSimpleNonSemiSimpleClassifyLieSumDirect}
Пусть $L$~"--- конечномерная $H_{m^2}(\zeta)$-простая $H_{m^2}(\zeta)$-модульная алгебра Ли над полем $\mathbbm{k}$, причём разрешимый радикал $R$ алгебры Ли $L$ не равен $0$.
Определим число $\ell \in\mathbb N$ при помощи условий
 $N^\ell = 0$, $N^{\ell-1} \ne 0$. Выберем минимальный $\mathbb Z/m\mathbb Z$-градуированный идеал $\tilde N \subseteq N^{\ell-1}$ алгебры Ли $L$.
Тогда для любого $k$ подпространство $N_k := \sum_{i=0}^{i=k} v^i \tilde N$ является $\mathbb Z/m\mathbb Z$-градуированным идеалом алгебры Ли $L$, причём $L = \bigoplus_{i=0}^t v^i \tilde N$ (прямая сумма $\mathbb Z/m\mathbb Z$-градуированных подпространств) для некоторого $1 \leqslant t \leqslant m-1$.
Более того, для всех $0 \leqslant k \leqslant t$ пространства $N_k/N_{k-1}$ являются неприводимыми $\mathbb Z/m\mathbb Z$-градуированными $L$-модулями, $c$-изоморфными друг другу. (Здесь $N_{-1} := 0$.)
\end{lemma}
\begin{proof}
Поскольку для любых $a \in \tilde N$ и $b\in L$ элемент $$[v^k a, b]=v[v^{k-1} a, b]-[cv^{k-1} a, vb]$$
мождет быть представлен в виде $\mathbbm{k}$-линейной комбинации
элементов $v^i[c^{k-i} a, v^{k-i}b]$, всякое подпространство $N_k := \sum_{i=0}^{i=k} v^i \tilde N$
является $\mathbb Z/m\mathbb Z$-градуированным идеалом алгебры Ли $L$.

 Напомним, что $v^m =0$. Следовательно, $N_{m-1}$ является $H_{m^2}(\zeta)$-инвариантным идеалом алгебры Ли $L$, откуда $L=N_{m-1}$.

Пусть $\theta_k \colon N_k/N_{k-1} \twoheadrightarrow N_{k+1}/N_k$, где $0 \leqslant k \leqslant m-2$,~"---
отображения, заданные равенствами  $$\theta_k (b + N_{k-1}) := vb + N_k.$$ 
Введём обозначение $\bar b:= b + N_{k-1}$.
Тогда $c\theta_k (\bar b) = \zeta^{-1}\theta_k (c\bar b)$,
\begin{equation*}\begin{split}
 \theta_k(a \bar b) = v[a,b]+N_k = -v[b,a] + N_k= -[cb,va]-[vb,a]+N_k = \\=
-[vb,a]+N_k=[a,vb]+N_k = a\theta_k (\bar b)
\text{ для всех }a\in L,\ b \in N_k.\end{split}
\end{equation*}
Заметим, что $\tilde N = N_0/N_{-1}$~"--- неприводимый $\mathbb Z/m\mathbb Z$-градуированный $L$-модуль.
Следовательно, для всех $0 \leqslant k < m-1$ пространства $N_{k+1}/N_k = \theta_k(N_k/N_{k-1})$
либо также являются неприводимыми $\mathbb Z/m\mathbb Z$-градуированными $L$-модулями, либо равны $0$.
Следовательно, если $L = N_t$, $L \ne N_{t-1}$,
то $\dim N_t = (t+1)\dim \tilde N$ и $L = \bigoplus_{i=0}^t v^i \tilde N$ (прямая сумма $\mathbb Z/m\mathbb Z$-градуированных подпространств).
\end{proof}

\begin{lemma}\label{LemmaTaftSimpleNonSemiSimpleLieClassifyUnity}
Предположим, что выполнены условия леммы~\ref{LemmaTaftSimpleNonSemiSimpleClassifyLieSumDirect}.
Кроме того, предположим, что основное поле $\mathbbm{k}$
алгебраически замкнуто характеристики $0$.
Тогда $R=N=N_{t-1}$, $L^{(i)} = v^{t-i} \tilde N$ для всех $0\leqslant i \leqslant t$ и $L^{(0)}\cong L/R$ является простой алгеброй Ли. 
Более того, $\dim (N_k/N_{k-1})=\dim (L/R)$ для всех $0 \leqslant k\leqslant t$.
Также $\ker v = L^{(0)}$.
\end{lemma}
\begin{proof}
Сперва заметим, что $[L,L]$~"--- $H_{m^2}(\zeta)$-инвариантный идеал. Следовательно, $L=[L,L]$ и
$L\ne R$.

В силу, например, теоремы~\ref{TheoremGradLevi}, существует максимальная $\mathbb Z/m\mathbb Z$-градуированная полупростая подалгебра $B\subseteq L$,
такая, что $L = B \oplus R$ (прямая сумма $\mathbb Z/m\mathbb Z$-градуированных подпространств), а $B \cong L/R$. Заметим, что $N$ аннулирует все неприводимые $\mathbb Z/m\mathbb Z$-градуированные $L$-модули,
которые являются факторами присоединённого представления алгебры Ли $L$.
Кроме того, $[L, R] \subseteq N$ (см., например, \cite[предложение 2.1.7]{GotoGrosshans}).
Отсюда $L/N$~"--- редуктивная $\mathbb Z/m\mathbb Z$-градуированная
алгебра Ли, а $\tilde N$~"--- неприводимый $\mathbb Z/m\mathbb Z$-градуированный
$L/N$-модуль.

 В силу леммы~\ref{LemmaRedIrrG}
 существует такой $L$-подмодуль $M\subseteq \tilde N$, что $\tilde N = \bigoplus_{i=0}^s c^i M$ для некоторого
 $s \in\mathbb N$ и разрешимый радикал $R$ действует на $M$ скалярными операторами.
 Поскольку $\mathbb Z/m\mathbb Z$-градуированный идеал $\tilde N$ минимален,
 $L$-подмодуль $M$ можно считать неприводимым.
Все $B$-подмодули в $M$ являются $L$-подмодулями,
так как $R$ действует на $M$ скалярными операторами.
Следовательно, $M$~"--- неприводимый $B$-модуль.
 
Поскольку $\tilde N=N_0/N_{-1}$, все $N_k/N_{k-1}$ $c$-изоморфны друг другу и алгебра Ли $B$ полупроста,
алгебра Ли $L$ является прямой суммой неприводимых $B$-подмодулей,
изоморфных $c^j M$, где $j \in\mathbb Z$.
Заметим, что $B$-действие на $M$ и, следовательно, на всяком $c^j M$
ненулевое, поскольку алгебра Ли $B$ сама является $B$-подмодулем алгебры Ли $L$
с ненулевым $B$-действием.
С другой стороны, существует такой $B$-подмодуль $Q \subseteq R$,
что $R=N\oplus Q$. В силу того, что $[L,R]\subseteq N$, справедливо включение $[B,Q]\subseteq N \cap Q = 0$,
т.е. $Q\subseteq L$ является подмодулем с нулевым $B$-действием.
Следовательно, $Q=0$ и $R=N$. Учитывая, что $[N,N_k]\subseteq N_{k-1}$, все
$N_k/N_{k-1}$ являются неприводимыми $\mathbb Z/m\mathbb Z$-градуированными
$B$-модулями,
$c$-изоморфными друг другу.
Однако $B\subseteq L$~"--- $\mathbb Z/m\mathbb Z$-градуированный
 $B$-подмодуль.
 Если алгебра Ли $B$ не была бы $\mathbb Z/m\mathbb Z$-градуированно простой алгеброй Ли,
  то $B$ была бы прямой суммой своих идеалов,
  являющихся $\mathbb Z/m\mathbb Z$-градуированно простыми
  алгебрами Ли (это следует, например, из теоремы~\ref{TheoremCoHLieSemiSimple}),
которые не $c$-изоморфны как $B$-модули.
Следовательно, $B$~"--- $\mathbb Z/m\mathbb Z$-градуированно простая алгебра Ли
и все $N_k/N_{k-1}$ $c$-изоморфны модулю $B$ как $\mathbb Z/m\mathbb Z$-градуированные
$B$-модули.
Определим число $q \in\mathbb Z_+$ при помощи условий $B \subseteq N_q$, $B \subsetneqq N_{q-1}$.
Если $q < t$, то $[B, N_t] \subseteq N_{t-1}$
и $B \left( N_t/N_{t-1} \right) = 0$. Получаем противоречие
с тем, что все $N_k/N_{k-1}$ $c$-изоморфны друг другу. Отсюда $B \cap N_{t-1} = 0$, $L=B \oplus N_{t-1}$ (прямая сумма подпространств) и $B \cong L/N_{t-1}$.

Докажем, что $ca=a$ для всех $a\in B$ и, следовательно, $B$ проста как обычная алгебра Ли.
В лемме~\ref{LemmaTaftSimpleNonSemiSimpleClassifyLieSumDirect}
было доказано, что $\theta_k(a\bar b) = a\theta_k(\bar b)$ 
для всех $0 \leqslant k \leqslant m-2$, $a \in L$ и $b\in N_k$.
Аналогично доказывается, что $$\theta_k(a\bar b) = v[a,b]+N_k = [ca, vb] + [va, b]+ N_k = 
[ca, vb] + N_k = (ca)\theta_k(\bar b).$$
Другими словами, $((ca)-a)$ действует как нулевой оператор $0$ на всех пространствах $N_k/N_{k-1}$
для любого $a\in B$.
 В частности, элемент $((ca)-a)$ принадлежит центру алгебры Ли $B$.
Поскольку алгебра Ли $B$ полупроста, получаем $ca=a$ для всех $a\in B$,
откуда алгебра Ли $B \subseteq L^{(0)}$ имеет тривиальную градуировку и
проста как обычная алгебра Ли.

Заметим, что $B \cong L/N_{t-1} \cong v^t \tilde N$ как $\mathbb Z/m\mathbb Z$-градуированные пространства.
 Следовательно, $v^t \tilde N \subseteq L^{(0)}$.  Используя
$vc=\zeta cv$, получаем $v^i \tilde N \subseteq L^{(t-i)}$.
Поскольку $L=\bigoplus_{i=0}^{m-1} L^{(i)} = \bigoplus_{i=0}^t v^i \tilde N$, при $0\leqslant i \leqslant t$ справедливо равенство $L^{(i)} = v^{t-i} \tilde N$,
а при $t+1\leqslant i \leqslant m-1$ подпространство $L^{(i)}$ нулевое.
В частности, $B=L^{(0)}=v^t \tilde N$.

Напомним, что всякое подпространство $N_j = \bigoplus_{i=0}^j v^i \tilde N=\bigoplus_{i=t-j}^t L^{(i)}$ 
является идеалом.
Следовательно, при всех
$0 \leqslant j \leqslant t$ и $0 \leqslant i \leqslant m-1$
справедливо неравенство $0 \leqslant t-j+i < t-j+m$ и включения
$$[L^{(i)}, L^{(t-j)}] \subseteq N_j \cap L^{(t-j+i)} \subseteq N_{j-i}$$
и $[L^{(i)}, N_k] \subseteq N_{k-i}$. (Считаем, что $N_k := 0$ при $k<0$.)
В частности, идеал $N_{t-1}$ нильпотентен и $R=N=N_{t-1}$.

Если $t=m-1$, то $vL^{(0)}=v^m \tilde N = 0$. Поскольку $N_{t-1} \cap (\ker v) = 0$, получаем $\ker v =L^{(0)}$.
Если $t < m-1$, то $vL^{(0)} \subseteq L^{(m-1)} = 0$. Снова получаем, что $\ker v =L^{(0)}$.
\end{proof}

\begin{lemma}\label{LemmaTaftSimpleNonSemiSimpleClassifyLieFormula}
Предположим, что выполнены условия леммы~\ref{LemmaTaftSimpleNonSemiSimpleLieClassifyUnity}.
Определим $\mathbbm{k}$-линейное отображение $\psi \colon L \to L$ при помощи равенства $\psi(v^k a) = v^{k-1} a$
для всех $a \in\tilde N$, $1 \leqslant k \leqslant t$, $\psi(\tilde N) = 0$.
Тогда \begin{equation}\label{EqQuantumBinomPhiLie}
[\psi^k(a),\psi^\ell(b)]=\binom{k+\ell}{k}_\zeta\ \psi^{k+\ell}[a,b] \text{ для всех }a, b\in L^{(0)} \text{ и } 0 \leqslant k,\ell \leqslant t.
\end{equation}
\end{lemma}
\begin{proof}
Заметим, что $\psi(va)=a$ для всех $a\in N_{t-1}$.
Следовательно, для всех $b\in L^{(0)}$ и $a=vu$, где $u\in N_{t-1}$,
справедливы равенства $$c\psi(a)=c\psi(vu)=cu=\psi(vcu)=\zeta\psi(cvu)=\zeta\psi(ca),$$
$$\psi[a,b]=\psi[vu,b]=\psi(v[u,b]-[cu, vb])=[u,b]=[\psi(a),b].$$
Поскольку  $\psi(\tilde N) = 0$,
 $L = v N_{t-1} \oplus \tilde N$ (прямая сумма $\mathbb Z/m\mathbb Z$-градуированных подпространств),
а $\tilde N$ является идеалом,
получаем равенства $c\psi(a) = \zeta \psi(ca)$ и  
$\psi[a,b]=[\psi(a),b]$  для всех $a\in L$, $b\in L^{(0)}$.
Этим~(\ref{EqQuantumBinomPhiLie}) доказано в случае, когда либо $k=0$, либо $\ell=0$.

В случае $k,\ell \geqslant 1$ формула~(\ref{EqQuantumBinomPhiLie})
доказывается по индукции с использованием равенства
\begin{equation*}\begin{split}
[\psi^k(a),\psi^\ell(b)]= \psi \left(v[\psi^k(a),\psi^\ell(b)]\right)
=\psi([c\psi^k(a), \psi^{\ell-1}(b)]+[\psi^{k-1}(a),\psi^\ell(b)])\end{split}\end{equation*}
точно так же, как это было сделано выше в лемме~\ref{LemmaLieTaftSimpleSGradPhiMult}.
\end{proof}

\begin{proof}[Доказательство теоремы~\ref{TheoremTaftSimpleNonSemiSimpleLieClassify}.]
В силу леммы~\ref{LemmaTaftSimpleNonSemiSimpleLieClassifyUnity}
алгебра Ли
$L^{(0)}$ является простой.
Возьмём произвольные элементы $a,b \in L^{(0)}$, такие, что $[a,b]\ne 0$.
Тогда $\psi^t [a,b]\ne 0$.
В то же время $[\psi^t(a), \psi(b)]=\binom{t+1}{t}_\zeta \psi^{t+1}[a,b] = 0$.
Однако \begin{equation*}\begin{split}0=v[\psi^t(a), \psi(b)]
=[v\psi^t(a), \psi(b)]+[c\psi^t(a), v\psi(b)]
=\\= \left(\binom t{t-1}_\zeta+\zeta^t\right)\psi^t[a,b]
=(t+1)_\zeta\ \psi^t[a,b].\end{split}\end{equation*} Следовательно, $(t+1)_\zeta=0$ и $m=t+1$.
Теперь утверждение теоремы следует из~(\ref{EqQuantumBinomPhiLie}) и леммы~\ref{LemmaTaftSimpleNonSemiSimpleLieClassifyUnity}.
\end{proof}

\begin{remark}
Поскольку максимальная полупростая подалгебра $\ker v$ определена однозначно,
любые две такие $H_{m^2}(\zeta)$-простые алгебры Ли $L$ изоморфны как $H_{m^2}(\zeta)$-модульные
алгебры Ли, если и только если их подалгебры Ли $\ker v$
изоморфны как обычные алгебры Ли.
Более того, все автоморфизмы алгебры Ли $L$ как $H_{m^2}(\zeta)$-модульной
алгебры Ли индуцируются автоморфизмами алгебры Ли $\ker v$ как обычной алгебры Ли.
Действительно, пусть $\theta \colon L \mathrel{\widetilde\to} L$~"---
произвольный автоморфизм алгебры Ли $L$ как $H_{m^2}(\zeta)$-модульной алгебры Ли.
Поскольку $\tilde N = N^{m-1}$, получаем $\theta(\tilde N) = \tilde N$
и $\theta(v^k \tilde N)=v^k \tilde N$ для всех $0\leqslant k < m$. Теперь
из $v^k \theta(\psi^k(a))=\theta(a)$
для всех $a\in \ker v$ следует, что $\theta(\psi^k(a))=\psi^k(\theta(a))$ и $\theta$
однозначно определено своим ограничением на $\ker v$.
\end{remark}

\newpage

\chapter{Ассоциативные алгебры, градуированные полугруппами}
\label{ChapterGradSemigroup}

В данной главе рассматриваются структурные вопросы теории ассоциативных алгебр, градуированных полугруппами.

Результаты \S\ref{SectionSemigroupTwoElements}--\ref{SectionSemigroupWedderburn}
касаются градуированности радикала Джекобсона и градуированных вариантов теорем Веддербёрна~"--- Артина и  Веддербёрна~"--- Мальцева, при этом основное внимание уделено алгебрам, градуированным полугруппами из двух элементов и лентами левых и правых нулей. Эти результаты были опубликованы
в работе~\cite{ASGordienko13}.

Результаты~\S\ref{SectionGeneralReductionSemigroup}--\ref{SectionTGradedReesExistence}
были получены автором совместно с Э.~Йесперсом и Дж.~Янссенсом
и опубликованы в~\cite{ASGordienko14JanssensJespers}.
Они посвящены классификации градуированно простых алгебр, градуированных конечными полугруппами
с тривиальными максимальными подгруппами.

Доказанные утверждения будут затем использованы в главе~\ref{ChapterSGGrAssocCodim} при изучении градуированных полиномиальных тождеств.

\section{Полугруппы, состоящие из двух элементов}\label{SectionSemigroupTwoElements}

Как будет показано в следующих параграфах, для многих интересных примеров достаточно полугрупп, состоящих из двух элементов. Докажем, что таких попарно неизоморфных полугрупп существует ровно пять.

Обозначим через $Q_1 =(\lbrace 0,1\rbrace, \cdot)$ мультипликативную полугруппу поля $\mathbb Z/2\mathbb Z$.

Пусть $Q_2 =(\lbrace 0,v\rbrace, \cdot)$~"--- полугруппа, заданная соотношениями $v^2 =0^2= 0\cdot v=v\cdot 0=0$.

Полугруппа $T$ называется \textit{лентой левых нулей},
если $t_1 t_2 = t_1$ для всех $t_1, t_2 \in T$
и \textit{лентой правых нулей},
если $t_1 t_2 = t_2$ для всех $t_1, t_2 \in T$.

Обозначим через $Q_3$ ленту правых нулей, состоящую из двух элементов.

\begin{proposition}\label{PropositionSemigroupTwoElements}
Пусть $T$~"--- полугруппа, состоящая из двух элементов. Тогда $T$ 
изоморфна
одной из полугрупп из списка $\lbrace Q_1, Q_2, Q_3, Q_3^{\,\op}, (\mathbb Z/2\mathbb Z,+)\rbrace$
и любые две полугруппы из этого списка неизоморфны. (Здесь через  $Q_3^{\,\op}$
обозначена полугруппа, антиизоморфная
полугруппе $Q_3$.)
\end{proposition}
\begin{proof} Сперва рассмотрим случай, когда $T = \lbrace a, a^2 \rbrace$ для некоторого $a\in T$.
Тогда если $a^3=a$, то $a^4=a^2$ и элемент $a^2$ является единицей полугруппы $T$, т.е.
$T \cong (\mathbb Z/2\mathbb Z,+)$. Если $a^3=a^2$, то $a^4=a^3=a^2$ и элемент $a^2$ является нулём полугруппы $T$,
т.е. $T \cong Q_2$.

Теперь рассмотрим случай, когда $T \ne \lbrace a, a^2 \rbrace$ ни для какого элемента $a\in T$.
Тогда $T = \lbrace a, b \rbrace$, $a^2=a$, $b^2=b$. Если $ab=ba$, то $T \cong Q_1$.
Если $ab\ne ba$, то $T\cong Q_3$ при $ab = b$, $ba=a$ и $T\cong Q_3^{\,\mathrm{op}}$
при $ab = a$, $ba=b$.
\end{proof}

\section{Градуированность радикала Джекобсона}\label{SectionSemigroupJacobson}

Алгебры, градуированные полугруппой $Q_1$,~"--- это в точности алгебры
с фиксированным разложением в прямую сумму двустороннего идеала и подалгебры.
Если $T$~"--- лента правых нулей, то
$T$-градуированные алгебры~"--- это в точности алгебры, разложенные в прямую сумму левых идеалов,
проиндексированных элементами полугруппы $T$.

Известно (см., например, \cite[пример 4.2]{KelarevBook}),
что в алгебрах, градуированных полугруппами,
радикал Джекобсона необязательно является градуированным  идеалом.
(См. обзор достаточных условий градуированности в~\cite[\S 4.4]{KelarevBook}.)
В данном параграфе мы приведём примеры и докажем результаты, касающиеся полугрупп из двух 
элементов и лент правых и левых нулей.

В алгебре $M_k(\mathbbm{k})$ всех квадратных матриц размера $k\times k$ зафиксируем
базис, состоящий из матричных единиц $e_{i\ell}$, где $1\leqslant i,\ell \leqslant k$.

\begin{example}\label{ExampleQ1} 
Пусть $A = M_k(\mathbbm{k})\oplus \UT_k(\mathbbm{k})$ (прямая сумма идеалов), где $\mathbbm{k}$~"--- некоторое поле, а $k\geqslant 2$. Определим $Q_1$-градуировку на $A$  при помощи формул $A^{(0)}=(M_k(\mathbbm{k}),0)$, $A^{(1)}=\lbrace (\varphi(a), a) \mid a \in \UT_k(\mathbbm{k}) \rbrace$, где $\varphi \colon \UT_k(\mathbbm{k}) \hookrightarrow M_k(\mathbbm{k})$~"---
естественное вложение. Тогда $$J(A) = \lbrace (0,e_{ij}) \mid 1 \leqslant i < j \leqslant k\rbrace \subset (0,\UT_k(\mathbbm{k})),$$ $J(A) \cap A^{(0)} = J(A) \cap A^{(1)}=0$
и $J(A)$ не является градуированным идеалом.
\end{example}

\begin{example}\label{ExampleQ2}
Пусть $A = M_k(\mathbbm{k})\oplus V$ (прямая сумма идеалов), где $V\cong M_k(\mathbbm{k})$ как векторные пространства, $k\in \mathbb N$, $V^2 = 0$, а $\mathbbm{k}$~"--- некоторое поле. Обозначим через $\varphi \colon V \mathrel{\widetilde{\rightarrow}} M_k(\mathbbm{k})$
соответствующую биекцию.
 Определим $Q_2$-градуировку на $A$ при помощи формул $A^{(0)}=(M_k(\mathbbm{k}),0)$, $A^{(v)}=\lbrace (\varphi(a), a) \mid a \in V \rbrace$. Тогда $$J(A) = (0,V),\ J(A) \cap A^{(0)} = J(A) \cap A^{(v)}=0$$
и $J(A)$ не является градуированным идеалом.
\end{example}

\begin{example} \label{ExampleQ3}
Пусть $A = M_k(\mathbbm{k})\oplus V$ (прямая сумма левых идеалов), где $V\cong M_k(\mathbbm{k})$ как левые
$M_k(\mathbbm{k})$-модули, $k\in \mathbb N$, $V^2=V M_k(\mathbbm{k}) = 0$, а $\mathbbm{k}$~"--- поле. Обозначим через $\varphi \colon V \mathrel{\widetilde{\rightarrow}} M_k(\mathbbm{k})$
соответствующий изоморфизм. Определим $Q_3$-градуировку на $A$ при помощи формул $A^{(e_1)}=(M_k(\mathbbm{k}),0)$, $$A^{(e_2)}=\lbrace (\varphi(a), a) \mid a \in V \rbrace.$$ Тогда $$J(A) = (0,V),\ J(A) \cap A^{(e_1)} = J(A) \cap A^{(e_2)}=0$$
и $J(A)$ не является градуированным идеалом.
\end{example}

\begin{remark}
Рассматривая вместо $A$ алгебру $A^{\,\mathrm{op}}$, получим пример $Q_3^{\,\mathrm{op}}$-градуированной алгебры, радикал Джекобсона которой не является градуированным.
\end{remark}

Однако если $Q_3$- или $Q_3^{\,\mathrm{op}}$-градуированная алгебра содержит единицу,
то её радикал Джекобсона градуирован. В действительности справедлив более общий результат:

\begin{proposition}\label{PropositionTIdemGradedIdeals}
Пусть $A$~"--- $T$-градуированная ассоциативная алгебра с $1$
над полем $\mathbbm{k}$ для некоторой ленты $T$ левых или правых нулей.
Тогда любой двусторонний идеал алгебры $A$ градуированный.
\end{proposition}
\begin{proof}
Рассмотрим случай, когда $t_1 t_2 = t_2$ для всех $t_1, t_2 \in T$. (Случай ленты левых нулей рассматривается аналогично.)
Тогда все $A^{(t)}$, где $t\in T$, являются левыми идеалами алгебры $A$ и $1 = \sum_{t\in T}e_t$ для некоторых $e_t \in A^{(t)}$.
Пусть $I$~"--- произвольный идеал алгебры $A$. Тогда для любого $a\in I$ имеет место разложение $a=\sum_{t\in T}a e_t$,
где $ae_t\in I \cap A^{(t)}$ для любого $t\in T$. Следовательно, идеал $I=\bigoplus_{t\in T} I \cap A^{(t)}$
градуированный. 
\end{proof}

\section{Градуированные аналоги теорем Веддербёрна и $T$-градуированная простота}\label{SectionSemigroupWedderburn}

Теперь исследуем справедливость градуированных аналогов теорем Веддербёрна
в $T$-градуированных алгебрах, где $T$~"--- полугруппа. 

\begin{example}\label{ExampleQ1Wedderburn}
Пусть $B=M_k(\mathbbm{k})\oplus M_k(\mathbbm{k})$
(прямая сумма идеалов), $k\in \mathbb N$, а $\mathbbm{k}$~"--- поле.
Определим $Q_1$-градуировку на $A$ при помощи равенств $$A^{(0)}=(M_k(\mathbbm{k}),0),\qquad A^{(1)}=\lbrace (a, a) \mid a \in M_k(\mathbbm{k}) \rbrace.$$
Тогда $B$ нельзя представить в виде прямой суммы $Q_1$-градуированных идеалов,
которые являлись бы $Q_1$-градуированно простыми алгебрами,
т.е.  $Q_1$-градуированный аналог теоремы Веддербёрна~"--- Артина неверен.
\end{example}
\begin{proof}  Достаточно заметить, что у полупростой алгебры $B$ только четыре
идеала: $0$, $B$, $(0,M_k(\mathbbm{k}))$ и $(M_k(\mathbbm{k}),0)$. Три из них, а именно $0$, $B$ и $(M_k(\mathbbm{k}),0)$,  $Q_1$-градуированы, и только $(M_k(\mathbbm{k}),0)$ является $Q_1$-градуированно простой алгеброй.
\end{proof}
\begin{remark}
Поскольку компонента $B^{(0)}$ всегда является градуированным идеалом,
всякая $Q_1$-градуированно простая алгебра $B$ имеет тривиальную градуировку, т.е. $B=B^{(0)}$.
Следовательно,  всякая $Q_1$-градуированно простая алгебра
является простой.
\end{remark}

\begin{proposition}\label{PropositionQ2Wedderburn}
Пусть $B$~"--- конечномерная
полупростая ассоциативная $Q_2$-градуированная алгебра над полем $\mathbbm{k}$.
Тогда $B=B^{(0)}$ и в силу обычной теоремы Веддербёрна~"--- Артина
алгебра $B$ является прямой суммой $Q_2$-градуированных идеалов,
являющихся простыми алгебрами (с тривиальной градуировкой).
Более того, всякая $Q_2$-градуированно простая алгебра проста.
В частности, справедлив $Q_2$-градуированный 
аналог теоремы Веддербёрна~"--- Артина.
\end{proposition}
\begin{proof}
Предположим, что $B \ne B^{(0)}$.
Учитывая, что $B^{(0)}$~"--- идеал, из обычной
теоремы Веддербёрна~"--- Артина
следует, что $B=B^{(0)} \oplus I$ некоторого полупростого идеала $I$ алгебры $B$. Однако $(B/B^{(0)})^2=0$,
поскольку $(B^{(1)})^2 \subseteq B^{(0)}$ и алгебра $I \cong B/B^{(0)}$ не может быть полупростой.
Отсюда $B=B^{(0)}$, $Q_2$-градуировка на алгебре $B$ является тривиальной и к $B$
можно применить обычную теорему Веддербёрна~"--- Артина.
\end{proof}

\begin{proposition}\label{PropositionTIdemWedderburn}
Пусть $B$~"--- конечномерная полупростая ассоциативная $T$-градуированная
алгебра над полем $\mathbbm{k}$ для некоторой ленты $T$
левых или правых нулей.
Тогда $B$ расскладывается в прямую сумму $T$-градуированных идеалов,
являющихся простыми алгебрами.
В частности, справедлив $T$-градуированный аналог теоремы
Веддербёрна~"--- Артина, и любая конечномерная полупростая $T$-градуированно простая
алгебра является простой.
\end{proposition}
\begin{proof} 
Рассмотрим случай, когда $T$ является лентой правых нулей.
(Случай ленты левых нулей разбирается аналогично.)
 
В силу обычной теоремы Веддербёрна~"--- Артина  $$B=B_1\oplus B_2 \oplus \ldots \oplus B_s \text{ (прямая сумма идеалов)}$$ для некоторых простых алгебр $B_i$.
Теперь утверждение теоремы следует из того, что в силу
теоремы~\ref{PropositionTIdemGradedIdeals}
 каждый из идеалов
$B_i$ является $T$-градуированным.
\end{proof}

Здесь нужно отметить, что существуют неполупростые градуированно простые
алгебры над лентами левых и правых нулей. (См. теорему~\ref{TheoremImagesOfGradedComponentsReesExistence} ниже.)

В силу предложения~\ref{PropositionTIdemGradedIdeals} радикал Джекобсона
(точно так же, как и все остальные идеалы)
любой $T$-градуированной алгебры с единицей, где $T$~"--- лента левых или правых нулей,
является градуированным. 
Поэтому естественным образом возникает вопрос о справедливости $T$-градуированного
аналога теоремы Веддербёрна~"--- Мальцева.
Оказывается, $T$-градуированный аналог этой теоремы действительно справедлив.

\begin{theorem}\label{TheoremTIdemGradedWeddMalcev}
Пусть $A$~"--- конечномерная ассоциативная
 $T$-градуированная алгебра с единицей над полем $\mathbbm{k}$,
 где $T$~"--- лента левых или правых нулей, а алгебра $A/J(A)$ сепарабельна. (Например, поле $\mathbbm{k}$
совершенно.)
Тогда существует такая максимальная полупростая подалгебра $B$, что $A = B \oplus J$
 (прямая сумма градуированных подпространств), где $J := J(A)$.
\end{theorem}

\begin{proof} Без ограничения общности можно считать, что $T$~"--- лента правых нулей.
Рассмотрим сперва случай $J^2=0$.

Заметим, что $1_A = \sum_{t\in T} e_t$ для некоторых $e_t\in A^{(t)}$. 
Более того, поскольку $ e_r 1_A=
\sum_{t\in T} e_r e_t $ и $e_r e_t \in A^{(t)}$ для всех $t\in T$, для всех $r\ne t$
выполняются равенства $e_r e_t = 0$ и
  $e_r^2 = e_r$.

Используя обычную теорему Веддербёрна~"--- Мальцева, выберем максимальную полупростую
подалгебру $B$, такую, что $A = B \oplus J$ (прямая сумма подпространств).
Пусть $\pi \colon A \twoheadrightarrow A/J$~"--- естественный
сюръективный гомоморфизм, который является градуированным, поскольку
идеал $J$ градуированный. Пусть $\varphi \colon A/J \hookrightarrow A$~"--- такое гомоморфное
вложение, что $\varphi(A/J)=B$ и $\pi\varphi = \id_{A/J}$.
Заметим, что $\pi(1_A)=\sum_{t\in T} \pi(e_t)$~"--- единица алгебры $A/J$.
Кроме того, $1_A = 1_B = \varphi\pi(1_A)$.

В силу того, что алгебра конечномерна, носитель градуировки также конечен.
Поскольку любое подмножество ленты правых нулей является подполугруппой,
можно считать, что полугруппа $T$ конечна.
Пусть $T=\lbrace t_1, \ldots, t_s\rbrace$.
Если $\varphi\pi(e_{t_i})=e_{t_i}$ для всех $1\leqslant i\leqslant s$,
то $B{e_{t_i}} \subseteq B$ для всех $1\leqslant i\leqslant s$, $B=\bigoplus_{t\in T} B{e_t}$
является градуированной подалгеброй и теорема доказана.

Предположим, что $\varphi\pi(e_{t_i}) \ne e_{t_i}$ по крайней мере для одного $1\leqslant i\leqslant s$.
Выберем $0\leqslant k \leqslant s-1$, такое, что $\varphi\pi(e_{t_i})=e_{t_i}$
для всех $1\leqslant i \leqslant k$
и $\varphi\pi(e_{t_{k+1}}) \ne e_{t_{k+1}}$.
 Заметим, что $\pi(\varphi\pi(e_{t_{k+1}})-e_{t_{k+1}})=0$, т.е. $\varphi\pi(e_{t_{k+1}}) = e_{t_{k+1}}+j$ для некоторого $j\in J$. Кроме того, $$j e_{t_i} = (\varphi\pi(e_{t_{k+1}})-e_{t_{k+1}})e_{t_i}
 = \varphi\pi(e_{t_{k+1}} e_{t_i})-e_{t_{k+1}} e_{t_i} = 0\text{ для всех }1\leqslant i\leqslant k.$$
 Аналогично, $e_{t_i}j=0$ для всех $1\leqslant i\leqslant k$.
 Более того, из  $(e_{t_{k+1}} +j)^2=e_{t_{k+1}} +j$ следует, что
 $j=e_{t_{k+1}}j + j e_{t_{k+1}}$ и $e_{t_{k+1}} j e_{t_{k+1}} = 0$.
 
 Пусть $\tilde \varphi \colon A/J \hookrightarrow A$~"---
 гомоморфное вложение, заданное  равенством \begin{equation*}\begin{split}\tilde\varphi(a) = (1_A + e_{t_{k+1}}j - je_{t_{k+1}})\varphi(a)(1_A + e_{t_{k+1}}j - je_{t_{k+1}})^{-1}
 =\\=(1_A + e_{t_{k+1}}j - je_{t_{k+1}})\varphi(a)(1_A - e_{t_{k+1}}j + je_{t_{k+1}}).\end{split}\end{equation*}
 Заметим, что $\pi\tilde\varphi=\id_{A/J}$ и $$\tilde \varphi \pi(e_{t_i})=(1_A + e_{t_{k+1}}j - je_{t_{k+1}})e_{t_i}(1_A - e_{t_{k+1}}j + je_{t_{k+1}})=e_{t_i} \text{ для всех } 1\leqslant i \leqslant k.$$
 Более того, \begin{equation*}\begin{split}\tilde \varphi \pi(e_{t_{k+1}})=
 (1_A + e_{t_{k+1}}j - je_{t_{k+1}})\varphi \pi(e_{t_{k+1}})(1_A - e_{t_{k+1}}j + je_{t_{k+1}})
 =\\= 
 (1_A  + e_{t_{k+1}}j - je_{t_{k+1}})(e_{t_{k+1}}+j)(1_A  - e_{t_{k+1}}j + je_{t_{k+1}})=\\=
 (e_{t_{k+1}}+j  - je_{t_{k+1}})(1_A  - e_{t_{k+1}}j + je_{t_{k+1}})=\\=
e_{t_{k+1}}+j  - je_{t_{k+1}} - e_{t_{k+1}}j =
 e_{t_{k+1}}.\end{split}\end{equation*}
 Следовательно, $\tilde B = \tilde\varphi(A/J)$ 
является максимальной  полупростой подалгеброй,
такой, что $A = \tilde B \oplus J$ (прямая сумма подпространств)
 и $\tilde \varphi \pi(e_{t_i}) = e_{t_i}$ для всех $1\leqslant i \leqslant k+1$.
 
 Рассуждая по индукции, можно считать, что 
 $e_t=\tilde \varphi \pi(e_t) \in \tilde B$ для всех $t\in T$.
Следовательно, $\tilde B=\bigoplus_{t\in T} \tilde B{e_t}$~"--- градуированная подалгебра
алгебры $A$.

В случае $J^2=0$ теорема доказана.
Общий случай доказывается индукцией по $\dim A$.
Предположим, что $J^2\ne 0$. Тогда $A/J^2 = B_0 \oplus J/J^2$ (прямая сумма градуированных подпространств)
для некоторой градуированной максимальной полупростой подалгебры $B_0$ алгебры $A/J^2$.  Заметим, что $1_{A/J^2} \in B_0$. 
Рассмотрим прообраз $B_1$ алгебры $B_0$ в $A$ относительно естественного сюръективного гомоморфизма $\pi_1 \colon A \twoheadrightarrow A/J^2$. Тогда $1_A \in B_1$.
Поскольку алгебра $B_0 \cong A/J$ полупроста, справедливо равенство $J(B_1)=J^2$.
Более того, $\dim B_1 < \dim A$ и в силу предположения индукции $B_1 = B \oplus J^2$
(прямая сумма градуированных подпространств) для некоторой градуированной максимальной полупростой подалгебры $B$ алгебры $A$.
Следовательно, $A = B \oplus J$ (прямая сумма градуированных подпространств) и теорема доказана.
\end{proof}

Полугруппа $T$ называется \textit{полугруппой с сокращениями}, если для всех $a,b,c \in T$
из любого из равенств $ac = bc$ или $ca=cb$ следует, что $a=b$.

\begin{proposition}\label{PropositionTCancelWedderburn}
Пусть $B$~"--- конечномерная полупростая ассоциативная $T$-градуированная 
алгебра над полем $\mathbbm{k}$
для некоторой полугруппы $T$ с сокращениями.
Тогда $B$ является прямой суммой $T$-градуированных идеалов,
являющихся $T$-градуированно простыми алгебрами.
Иными словами, справедлив $T$-градуированный аналог  теоремы Веддербёрна~"--- Артина.
\end{proposition}
\begin{proof}
В силу обычной теоремы Веддербёрна~"--- Артина  $B = B_1 \oplus \ldots \oplus B_s$ (прямая сумма идеалов)
для некоторых простых алгебр $B_i$.
Пусть $I$~"--- минимальный $T$-градуированный идеал алгебры $B$.
Тогда $I=B_{i_1}\oplus\ldots \oplus B_{i_k}$ для некоторых $i_1, \ldots, i_k$.
Введём обозначение $N:=\bigoplus_{i \in \lbrace 1,\ldots, s\rbrace \backslash \lbrace i_1, \ldots, i_k\rbrace} B_i$.
Поскольку все идеалы $B_i$ полупросты, справедливо равенство $N=\lbrace b\in B \mid ba =0
 \text{ для всех } a\in I\rbrace$.
В силу того, что полугруппа $T$ с сокращениями, а идеал $I$ градуированный, идеал $N$
также является градуированным.
Следовательно, $B = I \oplus N$ (прямая сумма градуированных идеалов),
где $I$~"--- $T$-градуированно простая алгебра. Для завершения доказательства достаточно применить
к $N$ индукционные рассуждения.
\end{proof}

\section{Кольца, градуированные конечными полугруппами}\label{SectionGeneralReductionSemigroup}

Оставшаяся часть главы посвящена изучению градуированно простых колец и алгебр.

До этого момента в работе использовалось только понятие градуированной алгебры над полем. Аналогичным образом вводится понятие градуированного кольца.

Пусть $T$~"--- полугруппа, а $R$~"--- ассоциативное кольцо.
Разложение $\Gamma \colon R=\bigoplus_{t\in T} R^{(t)}$ кольца $R$ в прямую сумму
аддитивных подгрупп $R^{(t)}$, где $t\in T$, называется \textit{$T$-градуировкой},
если $R^{(s)} R^{(t)} \subseteq R^{(st)}$ для всех $s,t\in T$.
Подгруппа $P$ аддитивной группы кольца $R$ называется \textit{однородной} или \textit{градуированной}, если $P=\bigoplus_{t\in T} (R^{(t)}\cap P)$.
Говорят, что кольцо $R$ является  \textit{$T$-градуированно простым}, если  $R^{2}\neq 0$ и $0$ и $R$
являются единственными градуированными идеалами кольца $R$.
Полугруппу $T$ всегда можно заменить на подполугруппу, порождённую
\textit{носителем} $\supp \Gamma=\{ t\in T \mid R^{(t)}\neq  0\} $ градуировки $\Gamma$.

Заметим, что если полугруппа $T$ содержит нулевой элемент $0$, то компонента $R^{(0)}$ является градуированным идеалом кольца $R$, т.е. если $R^{(0)}\neq 0$  и кольцо $R$
градуированно простое, то  $R^{(0)}=R$ и $\supp\Gamma=\{ 0 \}$, т.е. градуировка $\Gamma$ тривиальна.
Отсюда при изучении градуированно простых колец можно без ограничения общности считать, что если
 $0 \in T$, то $R^{(0)}=0$. 
Заменяя, если нужно, $T$ на $T^0:=T \cup \lbrace 0 \rbrace$, мы можем считать, что $T$ всегда содержит нулевой элемент, $R^{(0)}=0$ и $T=\langle\supp\Gamma\rangle \cup \{ 0 \}$.

Предположим теперь, что кольцо $R$ является $T$-градуированно простым и $T=\langle \supp\Gamma\rangle  \cup \{ 0 \}$. Заметим, что если $I$~"--- \textit{идеал} полугруппы $T$, т.е. $st,ts \in I$ для всех $t\in T$, $s\in I$, то подмножество $R_{I}:=\bigoplus_{t\in I} R^{(t)}$ является идеалом кольца $R$.
Следовательно, в силу градуированной простоты кольца $R$
либо $R_{I}=0$, либо $R_{I}=R$. 
Из последнего следует, что $I=T$.
Отсюда если $I$~"--- нетривиальный идеал полугруппы $T$, то  $I\cap \supp\Gamma = \varnothing$. 
В этом случае полугруппу $T$ можно заменить на \textit{факторполугруппу Риса} $T/I$, т.е.
на полугруппу $T$, в которой подмножество элементов $I$ заменено единственным элементом, являющимся
нулём по умножению. Таким образом, без ограничения общности можно считать, что
полугруппа $T$ не содержит идеалов, отличных от $0$ и $T$ и $T^2 \ne 0$.
Такая полугруппа $T$ называется \textit{$0$-простой}.
Отсюда если $T$-градуированное кольцо $R$ является $T$-градуированно простым, тогда
без ограничения общности мы можем считать, что полугруппа $T$ является $0$-простой.
Разумеется, те же допущения можно сделать и при изучении градуированных полугруппами градуированно простых алгебр над полями.

Пусть $G$~"--- группа, $I,J$~"--- множества, а
$P=(p_{ji})$~"--- матрица размера $J\times I$
с элементами из множества $G^{0}:=G\cup \{ 0 \}$.
Обозначим через $\mathcal{M}(G^{0};I,J;P)$ множество $\{ (g,i,j) \mid i\in I,\, j\in J,\, g\in G^{0}\}$,
в котором все элементы $(0 ,i,j)$ отождествлены с нулевым элементом $0$. 
Определим на $\mathcal{M}(G^{0};I,J;P)$ ассоциативное умножение
следующим образом: $(g,i,j)(h,k,\ell):=(gp_{jk}h,i,\ell)$.
Полугруппа $\mathcal{M}(G^{0};I,J;P)$ называется \textit{полугруппой Риса матричного типа
над группой с нулём $G^0$ с сэндвич-матрицей $P$}.

Всякая конечная $0$-простая полугруппа $T$ является \textit{вполне $0$-простой}, т.е.
содержит ненулевой примитивный идемпотент (см. \cite[\S 2.7]{clipre}).
Отсюда в силу теоремы Риса~\cite[теорема 3.5]{clipre}
полугруппа $T$ изоморфна 
$\mathcal{M}(G^{0};I,J;P)$
для некоторой максимальной подгруппы $G$ полугруппы $T$, множеств $I,J$
и матрицы $P$, в каждой строчке и каждом столбце которой находится как минимум один ненулевой
элемент.

Таким образом, если градуированно простое кольцо $R$ градуировано
конечной полугруппой $T$, без ограничения общности можно считать,
что $$T=\mathcal{M}(G^{0},m,n;P)=\{ (g,i,j) \mid 1\leqslant i\leqslant m,\; 1\leqslant j \leqslant n,\;  g\in G^{0}\},$$ где $G$~"--- группа,
$m,n\in\mathbb N$, $P$~"--- матрица размера $n\times m$ с элементами из множества $G^0$.
Заметим, что кольцо $R$
градуировано также полугруппой $T'=\mathcal{M}(\{ e\}^{0},m,n;P')$,
где $e^2=e$, а матрица $P'$ задана следующим условием: элемент, стоящий на месте $(i,j)$, равен $e$, если $p_{ij}\neq 0$, и нулю в противном случае.

\begin{example}\label{ExampleRightZeroBandZeroSimple}
Если $T$~"--- лента правых нулей, $|T|=n$, то $T^0 \cong \mathcal{M}(\{ e\}^{0},1,n;P)$,
где $p_{i1}=e$ для всех $1\leqslant i \leqslant n$.
\end{example}

В данной главе мы классифицируем все $T$-градуированно простые кольца и алгебры над полями, градуированные
 конечными $0$-простыми полугруппами с тривиальными максимальными подгруппами. Напомним, что Ю.\,А.~Бахтурин, М.\,В.~Зайцев и С.\,К.~Сегал в своих работах \cite{BahturinZaicevSegal, BahZaicAllGradings}
 классифицировали градуированно простые конечномерные алгебры над алгебраически замкнутыми полями, градуированные группами. Таким образом, комбинирование результатов данной главы с результатами Ю.\,А.~Бахтурина, М.\,В.~Зайцева и С.\,К.~Сегала может привести к решению общей проблемы классификации градуированно простых колец и алгебр, градуированных конечными полугруппами с необязательно тривиальными подгруппами.

Покажем теперь, что разница между понятиями простого и $T$-градуированно простого кольца для таких полугрупп $T$ заключается
в нетривиальности левых и правых аннуляторов.

Будем называть кольцо \textit{точным}, если его левый и правый аннуляторы тривиальны,
т.е. из $aR = 0$ или $Ra=0$ для некоторого $a\in R$ следует, что $a=0$.

\begin{proposition}\label{simple iff}
Пусть  $R$~"--- кольцо, градуированное  конечной $0$-простой полугруппой $T$ с тривиальными максимальными подгруппами. Тогда кольцо $R$ простое, если и только если $R$ является точным $T$-градуированно простым кольцом.
\end{proposition}
\begin{proof}
Необходимость условий предложения очевидна, так как левый и правый аннуляторы являются двухсторонними аннуляторами.

Предположим теперь, что кольцо $R$ является точным и $T$-градуированно простым.
Поскольку $T$~"--- конечная $0$-простая полугруппа с тривиальными максимальными подгруппами, 
существует изоморфизм
$T \cong \mathcal{M}(\{ e\}^{0};I,J;P)$. Если  $I$~"--- ненулевой идеал кольца $R$, то $RIR$ является $T$-градуированным идеалом. В силу точности кольца $R$ идеал $I$ ненулевой, откуда $RIR=R$ и $RIR \subseteq I=R$, т.е. кольцо $R$ действительно является простым.
\end{proof}

Следующий несложный пример показывает, что условие точности кольца не может быть опущено.

Пусть $T=\lbrace e,f\rbrace$~"--- лента правых нулей, состоящая из двух элементов.
Нетрудно видеть, что
конечная полугруппа $T^0$ является $0$-простой, а
полугрупповая алгебра $\mathbbm{k}T$, где $\mathbbm{k}$~"--- произвольное поле, градуированно проста.
Однако эта алгебра не является простой, поскольку она содержит собственный
двухсторонний идеал $\mathbbm{k}(e-f)$. Заметим, что $(e-f)\mathbbm{k}T=0$, т.е. алгебра $\mathbbm{k}T$ не точна.

Покажем теперь, что
если $R \ne J(R)$ и кольцо $R$ является $T$-градуированно простым для некоторой полугруппы $T$,
то радикал Джекобсона $J(R)$ не содержит никакой особой информации не только о $T$-градуировке,
но и о структуре кольца $R$.

Начнём с того, что докажем, что никакой нетривиальный идеал 
не может содержать ненулевые однородные элементы.

\begin{lemma}\label{LemmaRadicalSemigroupGradedSimple}
Пусть $I \ne R$~"--- двухсторонний идеал $T$-градуированно простого кольца $R=\bigoplus_{t\in T} R^{(t)}$
для некоторой полугруппы $T$. Тогда $R^{(t)} \cap I = 0$
для всех $t\in T$.
\end{lemma}
\begin{proof}
Предположим, что $r \in R^{(t)} \cap I$ для некоторого $t\in T$.
Тогда наименьший двухсторонний идеал $I_0$, содержащий $r$, градуирован.
Поскольку $I_0 \subseteq I \subsetneqq R$, получаем $I_0=0$ и $r=0$.
\end{proof}

Напомним, что гомоморфизм  $T$-градуированных колец $\varphi \colon R_1 \rightarrow R_2$
называется \textit{строго градуированным}, если $\varphi\left(R_1^{(t)} \right) \subseteq R_2^{(t)}$
 для всех $t\in T$. Два $T$-градуированных кольца $R_1$ и $R_2$ называются \textit{градуированно изоморфными},
 если существует строго градуированный изоморфизм $R_1 \mathrel{\widetilde\rightarrow} R_2$.

Теорема~\ref{TheoremEquivalenceSemigroupGradedSimple}, которая доказывается ниже, оказывается особенно полезной в случае, когда кольцо $R/J(R)$ просто.

\begin{theorem}\label{TheoremEquivalenceSemigroupGradedSimple}
Пусть $T$~"---  полугруппа, а $R_i=\bigoplus_{t\in T} R_i^{(t)}$, где $i=1,2$, "--- $T$-градуированно простые кольца, $R_i \ne J(R_i)$.
Если существует изоморфизм колец $\bar\varphi \colon R_1/J(R_1) \mathrel{\widetilde\rightarrow}  R_2/J(R_2)$, такой,
что $\bar\varphi\left(\pi_1\left(R_1^{(t)}\right)\right)=\pi_2\left(R_2^{(t)}\right)$
для всех $t\in T$, где $\pi_i \colon R_i \twoheadrightarrow R_i/J(R_i)$, $i=1,2$, "--- естественные
сюръективные гомоморфизмы,
тогда существует изоморфизм $\varphi \colon R_1 \mathrel{\widetilde\rightarrow} R_2$ градуированных колец, такой,
что $\pi_2\varphi=\bar\varphi\pi_1$.

Обратно, если $\varphi \colon R_1 \mathrel{\widetilde\rightarrow} R_2$~"--- изоморфизм градуированных колец,
мы всегда можем определить изоморфизм колец $\bar\varphi \colon R_1/J(R_1) \mathrel{\widetilde\rightarrow} R_2/J(R_2)$
при помощи равенств $\bar\varphi(\pi_1(a)) = \pi_2\varphi(a)$ для всех $a \in R_1$.
При этом $\bar\varphi\left(\pi_1\left(R_1^{(t)}\right)\right)=\pi_2\left(R_2^{(t)}\right)$
для любого $t\in T$.
\end{theorem}
\begin{proof} Предположим, что существует
изоморфизм $\bar\varphi$.
Из леммы~\ref{LemmaRadicalSemigroupGradedSimple}
следует, что $$\pi_i\bigr|_{R^{(t)}_i} \colon R^{(t)}_i \rightarrow \pi_i\left(R^{(t)}_i\right)$$
является изоморфизмом аддитивных групп для любого $t\in T$ и $i=1,2$.
Определим $\varphi \colon R_1 \rightarrow R_2$ по формуле  $$\varphi\left(r\right):=\left(\pi_2\bigr|_{R^{(t)}_2}\right)^{-1}\bar\varphi\pi_1(r)\text{
при }r\in R_1^{(t)}\text{ и }t\in T,$$
продолжив этот гомоморфизм по аддитивности. Ясно, что $\varphi \left(R^{(t)}_1\right)=R^{(t)}_2$ и $\varphi$~"--- градуированный сюръективный гомоморфизм аддитивных групп. Более того, справедливо равенство $\pi_2\varphi=\bar\varphi\pi_1$.

Предположим, что $\varphi\left( \sum_{t\in T} r^{(t)}\right)=0$
для некоторых $r^{(t)}\in R_1^{(t)}$, где $t\in T$. Поскольку гомоморфизм $\varphi$ градуированный,  $\varphi\left(r^{(t)}\right)=0$ для всех $t\in T$. Отсюда $\pi_1\left(r^{(t)}\right)=0$
и $r^{(t)}=0$, так как согласно лемме~\ref{LemmaRadicalSemigroupGradedSimple} для каждого $t\in T$
справедливо равенство $R_1^{(t)}\cap J(R_1)=0$. Следовательно, $\varphi$~"--- биекция.

Докажем теперь, что $\varphi$~"--- изоморфизм колец. Для этого достаточно показать, что 
$\varphi$ сохраняет умножение.
Пусть $r^{(s)} \in R_1^{(s)}$ и $r^{(t)} \in R_1^{(t)}$.
Тогда $$\pi_2\varphi\left(r^{(s)}r^{(t)}\right)=\bar\varphi\pi_1\left(r^{(s)}r^{(t)}\right)
=\bar\varphi\pi_1\left(r^{(s)}\right)\bar\varphi\pi_1\left(r^{(t)}\right)=\pi_2\left(
\varphi\left(r^{(s)}\right)\varphi\left(r^{(t)}\right)\right).$$
Поскольку и $\varphi\left(r^{(s)}r^{(t)}\right)$, и $\varphi\left(r^{(s)}\right)\varphi\left(r^{(t)}\right)$ принадлежат одной и той же компоненте $ R_2^{(st)}$, а
 $\pi_2\bigr|_{R^{(st)}_2}$~"--- изоморфизм, получаем отсюда, что $$\varphi\left(r^{(s)}r^{(t)}\right)=\varphi\left(r^{(s)}\right)\varphi\left(r^{(t)}\right)
 \text{ для всех }r^{(s)} \in R_1^{(s)}\text{ и }r^{(t)} \in R_1^{(t)},$$
 и первое утверждение доказано.
 
 Второе утверждение очевидно, поскольку радикал Джекобсона под действием изоморфизма переходит в радикал Джекобсона.
\end{proof}

\begin{remark} Очевидно, что теорема~\ref{TheoremEquivalenceSemigroupGradedSimple} 
справедлива не только для радикала Джекобсона, но и для любого радикала.
(См. общее определение, например, в~\cite[глава IV, \S 6]{Skornyakov} или в~\cite[\S 4.5]{KelarevBook}.)
\end{remark}

\section{Односторонние идеалы матричных алгебр}\label{SectionLeftIdealsOfMatrixAlgebras}

В данном параграфе доказываются утверждения, которые будут затем использованы при работе
с конечномерными градуированно простыми алгебрами.

\begin{lemma}\label{LemmaKerMnFActDuality} Пусть $\mathbbm{k}$~"--- поле, а $k\in\mathbb N$. 
Рассмотрим естественное $M_k(\mathbbm{k})$-действие
на арифметическом векторном пространстве $\mathbbm{k}^k$ линейными операторами.
Тогда существует взаимно однозначное соответствие между левыми идеалами $I$ алгебры $M_k(\mathbbm{k})$
и подпространствами
 $W \subseteq \mathbbm{k}^k$, такими, что \begin{equation}\label{EquationIWonetoone}
I=\Ann W :=\lbrace a\in M_k(\mathbbm{k}) \mid aW = 0\rbrace,\qquad W=\bigcap_{a\in I} \ker a\end{equation}
и $\dim I = k(k-\dim W)$. Кроме того, если $I_1=\Ann W_1$ и $I_2=\Ann W_2$,
то $I_1+I_2= \Ann (W_1 \cap W_2)$ и $I_1 \cap I_2= \Ann(W_1+W_2)$.
\end{lemma}
\begin{proof} 
Пусть $W$~"--- подпространство пространства $\mathbbm{k}^k$. Выберем некоторый базис $w_{k+1-\dim W}, \ldots, w_k$ в $W$, а также такие элементы
$w_1, \ldots, w_{k-\dim W} \in \mathbbm{k}^k$, что $w_1, \ldots, w_k$ является базисом в $\mathbbm{k}^k$.
Тогда $\Ann W$ состоит из всех таких операторов $a\in M_k(\mathbbm{k})$, последние $\dim W$ столбцов матриц которых в базисе $w_1, \ldots, w_k$ состоят из нулей.
Заметим, что $\bigcap_{a\in \Ann W} \ker a = W$ и $\dim \Ann W=k(k-\dim W)$.

Пусть $I\subseteq M_k(\mathbbm{k})$~"--- левый идеал.
Поскольку $I$~"--- левый идеал в полупростой артиновой алгебре $M_k(\mathbbm{k})$, в силу теоремы~1.4.2
из~\cite{Herstein} существует такой идемпотент $e\in I$, что $I=M_k(\mathbbm{k})e$. Отсюда $I(\ker e)=0$. 
Заметим, что $e$ действует на пространстве $\mathbbm{k}^k$ как проектор. Следовательно, $\mathbbm{k}^k = \ker e \oplus \im e$. Выберем в $\mathbbm{k}^k$ базис, который является объединением базисов
пространств $\im e$ и $\ker e$. Тогда оператор $e$ имеет в этом базисе матрицу $\left(\begin{smallmatrix}
E & 0 \\
0 & 0 
\end{smallmatrix}\right)$, и левый идеал $I$
состоит из всех операторов с нулевыми последними $\dim\ker e$ столбцами.
Следовательно, $\bigcap_{a\in I} \ker a = \ker e$ и $\Ann\ker e = I$.
Теперь из первого абзаца доказательства получаем, что~(\ref{EquationIWonetoone}) действительно
является взаимно однозначным соответствием.

Предположим, что $I_1=\Ann W_1$ и $I_2=\Ann W_2$.
Тогда $$I_1 \cap I_2 = \Ann W_1 
\cap \Ann W_2 = \Ann(W_1+W_2).$$
Более того, $(I_1 + I_2)(W_1\cap W_2)=0$ и $I_1 + I_2 \subseteq 
\Ann (W_1\cap W_2)$.
Теперь утверждение леммы следует из равенств \begin{equation*}\begin{split}\dim(I_1 + I_2)=\dim I_1 + \dim I_2 - \dim(I_1\cap I_2)
=\\= k(2k-\dim W_1 - \dim W_2)-k(k-\dim(W_1+W_2))=\\= k(k-(\dim W_1 + \dim W_2-\dim(W_1+W_2)))
= \dim \Ann(W_1\cap W_2).\end{split}\end{equation*}
\end{proof}

\begin{theorem}\label{TheoremSumLeftIdealsMatrix} Пусть
$k,s \in \mathbb N$, $\mathbbm{k}$~"--- поле и
заданы такие левые идеалы $I_i$ алгебры $M_k(\mathbbm{k})$,
что $M_k(\mathbbm{k})=\bigoplus_{i=1}^s I_i$, причём
$\dim I_i = n_i k$, где $n_i\in\mathbb Z_+$.
Тогда существует такая матрица $P\in \GL_k(\mathbbm{k})$, что левый идеал $P^{-1} I_i P$
состоит из всех матриц $(\alpha_{k\ell})$, у которых 
$\alpha_{k\ell}=0$ при $\ell \leqslant \sum_{j=1}^{i-1} n_j$
и при $\ell > \sum_{j=1}^i n_j$.
\end{theorem}
\begin{proof}
Рассмотрим стандартное действие алгебры
$M_k(\mathbbm{k})$
на арифметическом векторном пространстве $\mathbbm{k}^k$.
 В силу леммы~\ref{LemmaKerMnFActDuality} для некоторого $V_i\subseteq \mathbbm{k}^k$
 справедливо равенство $I_i=\Ann V_i$.
 Применяя двойственность из леммы~\ref{LemmaKerMnFActDuality} к~$M_k(\mathbbm{k})=\bigoplus_{i=1}^s I_i$, получаем,
 что $\bigcap_{i=1}^s V_i = 0$ и $$V_i + \bigcap_{\substack{j=1,\\ j\ne i}}^s V_j = 
\mathbbm{k}^k\text{ для всех }1\leqslant i\leqslant s.$$ 
 Введём обозначение $W_i := \bigcap\limits_{\substack{j=1,\\ j\ne i}}^s V_j$.
 Тогда \begin{equation}\label{EquationFkViWi}\mathbbm{k}^k=V_i\oplus W_i.\end{equation} 
 Заметим, что $$\Ann W_i = \bigoplus_{\substack{j=1,\\ j\ne i}}^s I_j.$$
 Поскольку $\bigcap\limits_{\substack{j=1,\\ j\ne i}}^s \Ann W_j = 
 I_i$,
 получаем, что $V_i = \sum\limits_{\substack{j=1,\\ j\ne i}}^s W_j$
 и $\mathbbm{k}^k = \sum\limits_{j=1}^s W_j$.
 В силу~\eqref{EquationFkViWi} сумма $\mathbbm{k}^k = \bigoplus\limits_{j=1}^s W_j$ прямая.

 Теперь выберем в $\mathbbm{k}^k$ базис, который является объединением базисов в $W_i$. 
 Обозначим через $P\in \GL_k(\mathbbm{k})$ матрицу перехода от стандартного базиса в $\mathbbm{k}^k$
 к данному базису. Тогда каждый идеал $P^{-1} I_i P$
 состоит из всех матриц, у которых все столбцы нулевые, за исключением, быть может, тех столбцов,
 которые соответствуют подпространству $W_i$.  
  \end{proof}
  
  \begin{lemma}\label{LemmaLeftIdealMatrix}
  Пусть $I$~"--- минимальный левый идеал алгебры $M_k(\mathbbm{k})$, где $k\in\mathbb N$, а $\mathbbm{k}$~"--- некоторое поле.
  Тогда существуют такие фиксированные элементы $\mu_j\in \mathbbm{k}$,  где $1\leqslant j \leqslant k$, что
  $I = \left\langle \sum\limits_{j=1}^k \mu_j e_{ij} \mathrel{\biggl|} 1\leqslant i \leqslant k  \right\rangle_\mathbbm{k}$.
  \end{lemma}
  \begin{proof} Пусть $a=\sum\limits_{i,j=1}^k \mu_{ij} e_{ij} \in I \backslash \lbrace 0 \rbrace$.
В силу того, что $\sum_{\ell=1}^k e_{\ell\ell} a = a$, справедливо равенство $e_{\ell\ell}a \ne 0$ для некоторого $1\leqslant \ell \leqslant k$. Обозначим $\mu_j := \mu_{\ell j}$ для всех $1\leqslant j \leqslant k$. Тогда $$\left\langle e_{i\ell}a \mid 1\leqslant i \leqslant k  \right\rangle_\mathbbm{k}=\left\langle \sum\limits_{j=1}^k \mu_j e_{ij} \mathrel{\biggl|} 1\leqslant i \leqslant k  \right\rangle_\mathbbm{k}
$$ является ненулевым левым идеалом, содержащимся в $I$. Теперь утверждение леммы следует из минимальности левого идеала $I$.
  \end{proof}

\begin{lemma}\label{LemmaMatrixIdealsProduct}
Пусть $D$~"--- конечномерная алгебра с делением над полем $\mathbbm{k}$, а $k\in \mathbb N$.
Пусть $I$ и $V$~"--- соответственно, левый и правый идеалы алгебры $M_k(D)$.
Тогда $\dim_\mathbbm{k} (VI)=\frac{\dim_\mathbbm{k} V \dim_\mathbbm{k} I}{k^2 \dim D}$.
\end{lemma}
\begin{proof} 
Заметим, что $$I \cong \underbrace{M_k(D)e_{11}\oplus \ldots \oplus M_k(D)e_{11}}_{{\dim_\mathbbm{k} I}/(k\dim_\mathbbm{k} D)}\text{\quad и\quad}V \cong \underbrace{e_{11}M_k(D)\oplus \ldots \oplus e_{11}M_k(D)}_{{\dim_\mathbbm{k} V}/(k\dim_\mathbbm{k} D)}$$
как, соответственно, левый и правый $M_k(D)$-модули.
Следовательно \begin{equation*}\begin{split}\dim_\mathbbm{k} (VI)= \frac{\dim_\mathbbm{k} I}{k\dim_\mathbbm{k} D} \dim_\mathbbm{k}(V M_k(D)e_{11})
= \\ = \frac{\dim_\mathbbm{k} V \dim_\mathbbm{k} I}{k^2(\dim_\mathbbm{k} D)^2} \dim_\mathbbm{k}(e_{11} M_k(D)e_{11})
=  \frac{\dim_\mathbbm{k} V \dim_\mathbbm{k} I}{k^2 \dim D}.\end{split}\end{equation*}
\end{proof}

\section{Структура градуированно простых алгебр}\label{SectionReesSemigroupSimpleDescription}

В данном параграфе через $A$ обозначается некоторая $T$-градуированная алгебра над
полем $\mathbbm{k}$, где $T=\mathcal{M}(\{ e\}^{0},m,n;P)$~"--- произвольная конечная $0$-простая полугруппа с тривиальными максимальными подгруппами. 


Обозначим однородную компоненту алгебры $A$, отвечающую 
элементу $(e,i,j)$, через $A_{ij}$.
Тогда
    $$A=\bigoplus_{\substack{1\leqslant i \leqslant m,\\ 1\leqslant j \leqslant n}} A_{ij}$$
и 
  $$A_{ij}A_{k\ell}\subseteq A_{i\ell}.$$
  В частности, любая однородная компонента $A_{ij}$ является подалгеброй алгебры $A$.
  Пусть $ P = (p_{jk})_{j,k}$. Тогда  $A_{ij}A_{k\ell} = 0$
  для всех таких $1\leqslant i,k \leqslant m$ и $1\leqslant j,\ell \leqslant n$, что $p_{jk}=0$.

Заметим, что $A$~"--- $\mathcal{M}(\{ e\}^{0},m,n;P)$-градуированно простая алгебра для некоторой матрицы $P$, если и только если $A$~"--- $\mathcal{M}(\{ e\}^{0},m,n;P')$--градуированно простая алгебра для матрицы
$P'$, все элементы которой равны $e$.

Сделаем несколько предварительных замечаний:

\begin{lemma} \label{LemmaReesGradedFirstProperties} 
В произвольной $T$-градуированно простой алгебре $A$
выполнены следующие свойства:
\begin{enumerate}
\item $A_{ij} \cap J(A)= 0$ для всех $i,j$;
\item если $I\subseteq A$~"--- произвольное подмножество, то подмножество $AIA$ является градуированным идеалом, откуда $AIA$ равно либо $0$, либо $A$;
\item $AJ(A)A=0$.
\end{enumerate}
\end{lemma}
\begin{proof}
Предложение 1 является непосредственным следствием
леммы~\ref{LemmaRadicalSemigroupGradedSimple}.
Предложение 2 очевидно, а предложение 3 немедленно следует из предложения 2.
\end{proof}

Определим левые идеалы $L_j:=\bigoplus_{k=1}^m A_{kj}$ и правые идеалы
$R_i:=\bigoplus_{k=1}^n A_{ik}$. Тогда $L_j \cap R_i = A_{ij}$.

Следующая теорема является градуированной версией теоремы Веддербёрна~"--- Мальцева.
В ней доказывается существование ортогональных столбцовых и строчных идемпотентов,
которые задают полупростое дополнение к радикалу.
Этот результат~"--- первый шаг в классификации $T$-градуированно простых
алгебр.

\begin{theorem}\label{TheoremBGradedReesSemigroupGrSimple}
Пусть $A=\bigoplus_{i,j} A_{ij}$~"--- конечномерная $T$-градуированная алгебра над полем $\mathbbm{k}$, такая, что $AJ(A)A=0$. Тогда существуют такие ортогональные идемпотенты
$f_{1}, \ldots, f_n$ и такие ортогональные идемпотенты $f_1', \ldots, f_m'$ (некоторые из них могут быть нулевыми), что
$$B=\bigoplus_{i,j} f_i' A f_j =\bigoplus_{i,j} (B\cap A_{ij})$$~"--- $T$-градуированная максимальная полупростая подалгебра алгебры $A$, $f_i' \in B \cap R_i$ при $1\leqslant i \leqslant m$, $f_j \in B \cap L_j$ при $1\leqslant j \leqslant n$, $\sum_{i=1}^m f'_i=\sum_{j=1}^n f_j=1_B$ и 
 $A=B\oplus J(A)$ (прямая сумма подпространств).
\end{theorem}

\begin{proof} Через $\bar X$ мы будем обозначать образ подмножества $X$ алгебры $A$ в факторалгебре $A/J(A)$ под действием естественного сюръективного гомоморфизма
$A \rightarrow A/J(A)$.

Заметим, что $\bar A = \sum_{j=1}^n \bar L_j$.
Поскольку алгебра $\bar A = A/J(A)$ полупроста и вполне приводима как левый $A/J(A)$-модуль,
существуют левые идеалы $\tilde L_i \subseteq \bar  L_i$, которые дополняют $\bar{L_i} \bigcap \sum\limits_{j=1}^{i-1} \bar{L_j}$ до $\bar{L_i}$. Ясно, что $\bar A= \bigoplus\limits_{i=1}^{m} \tilde L_i$.
 Компоненты $\bar{\omega_i} \in \tilde{L_i}$ разложения $1_{\bar{A}} = \sum\limits_{i=1}^{m} \bar{\omega_i}$ единицы алгебры $\bar A$ являются ортогональными идемпотентами. Идемпотенты $\bar \omega_i$ 
являются образами некоторых идемпотентов $\omega_i \in L_i$ алгебры $A$
под действием естественных сюръективных гомоморфизмов $\pi\bigl|_{L_i} \colon L_i \rightarrow L_i / L_i \cap J(A)$, поскольку идеал $J(A)$ нильпотентен. Идемпотенты $\omega_1,\ldots, \omega_m$
также ортогональны, так как $\omega_i \omega_j = \omega_i(\omega_i \omega_j)\omega_j \in AJ(A)A= 0$.
 
 Аналогичным образом получаем такие ортогональные идемпотенты $\omega_1', \ldots, \omega_m' \in A$, где $\omega_i' \in R_i$ при $1\leqslant i \leqslant m$, что $\sum_{i=1}^m \bar f'_i = 1_{\bar A}$.
Пусть теперь $B := \bigoplus\limits_{\substack{1\leqslant i \leqslant m, 
\\ 1 \leqslant j \leqslant n}} \omega_i^{'} A \omega_j$.
Заметим, что $\omega_i^{'} A \omega_j \subseteq L_j \cap R_i= A_{ij}$, т.е. $B$ является $T$-градуированной подалгеброй алгебры $A$.
Пусть $a = \sum\limits_{\substack{1\leqslant i \leqslant m, 
\\ 1 \leqslant j \leqslant n}} \omega_i^{'} a_{ij} \omega_j \in J(A)$
для некоторых $a_{ij} \in A$. Тогда из $AJ(A)A=0$ следует, что $\omega_i^{'} a_{ij} \omega_j = \omega_i^{'} a \omega_j = 0$
для всех $i,j$. Отсюда $a=0$ и $J(A) \cap B = 0$. Более того, $\bar B = 1_{\bar A}  \bar A1_{\bar A}=\bar A$. Следовательно, $B$~"--- $T$-градуированная максимальная 
полупростая подалгебра алгебры $A$ и
$A = B \oplus J(A)$ (прямая сумма подпространств).

 Раскладывая $1_B$ по левым идеалам $\bigoplus_{i=1}^n \omega_i' A \omega_j$, где $1\leqslant j \leqslant n$, и по правым идеалам $\bigoplus_{j=1}^m \omega_i' A \omega_j$, где $1\leqslant i \leqslant m$, получаем
такие ортогональные идемпотенты $f_i \in B \cap L_i$, где $1\leqslant i \leqslant m$,
 и ортогональные идемпотенты $f_j' \in B \cap R_j$, где $1\leqslant j \leqslant n$,
 что $$\sum_{i=1}^m f'_i = \sum_{j=1}^n f_j = 1_B.$$
 При этом $$B = 1_B B 1_B = \bigoplus\limits_{\substack{1\leqslant i \leqslant m, 
\\ 1 \leqslant j \leqslant n}} f_i' B f_j = \bigoplus\limits_{\substack{1\leqslant i \leqslant m, 
\\ 1 \leqslant j \leqslant n}} f_i' A f_j,$$ поскольку $A=B \oplus J(A)$ и $AJ(A)A=0$.
\end{proof}

Пример~\ref{ExampleTwoBGradingsNonIso} показывает, что градуировки на разных подалгебрах $B$ 
из теоремы~\ref{TheoremBGradedReesSemigroupGrSimple} могут быть неэквивалентны.

\begin{example}\label{ExampleTwoBGradingsNonIso}
Пусть $\mathbbm{k}$~"--- поле, $I$~"--- левый $M_2(\mathbbm{k})$-модуль, изоморфный $\langle e_{12}, e_{22}\rangle_\mathbbm{k}$,
а $\varphi \colon I \mathrel{\widetilde{\rightarrow}} \langle e_{12}, e_{22}\rangle_\mathbbm{k}$~"--- соответствующий
изоморфизм.
Пусть $A=M_2(\mathbbm{k})\oplus I$ (прямая сумма $M_2(\mathbbm{k})$-модулей), где $IM_2(\mathbbm{k})=I^2=0$.
Зададим на $A$ следующую $Q_3$-градуировку:
$A^{(e_1)}= (M_2(\mathbbm{k}),0)$ и $A^{(e_2)}=\lbrace (\varphi(a),a) \mid a\in I\rbrace$.
Тогда алгебра $A$ является $Q_3$-градуированно простой и обе алгебры $B_1 = A^{(e_1)}$ и $B_2=\langle (e_{11},0),(e_{21},0)\rangle_\mathbbm{k}\oplus
A^{(e_2)}$
являются градуированными максимальными полупростыми подалгебрами алгебры $A$. Однако градуировки на $B_1$
и $B_2$ неэквивалентны. Тем более $B_1$
и $B_2$ неизоморфны как градуированные алгебры.
\end{example}

Построим теперь пример конечномерной $T$-градуированной алгебры,
которая не является градуированно простой и у которой
не существует $T$-градуированных максимальных полупростых
подалгебр, которые дополняют радикал до всей алгебры.
Из этого примера будет видно, что условие $AJ(A)A=0$ в теореме~\ref{TheoremBGradedReesSemigroupGrSimple} является существенным.

\begin{example}
Пусть
$R = \mathbbm{k}[X] / (X^2)$, а $A= M_{2}(R)$. Положим
$$v_1 = \left(\begin{smallmatrix}
1 \\
0 
\end{smallmatrix}\right),
v_2 = \left(\begin{smallmatrix}
0 \\
1 
\end{smallmatrix}\right),
w_1 = \left(\begin{smallmatrix}
1 & X 
\end{smallmatrix}\right), 
w_2 = \left(\begin{smallmatrix}
0 & 1
\end{smallmatrix}\right).$$
Рассмотрим следующие $\mathbbm{k}$-подпространства в $A$:
 $$\begin{array}{ll}
A_ {11} = R v_1 w_1 = R \left(\begin{smallmatrix}
1 & X\\
0 & 0
\end{smallmatrix}\right),
&
A_{12} = R v_1 w_2 = R \left(\begin{smallmatrix}
0 & 1\\
0 & 0
\end{smallmatrix}\right), \\
 &\\
A_{21} = R v_2 w_1 = R \left(\begin{smallmatrix}
0 & 0 \\
1 & X
\end{smallmatrix}\right),
&
A_{22} = R v_2 w_2 = R \left(\begin{smallmatrix}
0 & 0\\
0 & 1
\end{smallmatrix}\right).
\end{array}$$ Тогда разложение $$A = A_{11} \oplus A_{12} \oplus A_{21} \oplus A_{22}$$ является $T$-градуировкой,
где $T=\mathcal{M}\left(\{e \}^0, 2,2, \left( \begin{smallmatrix}
1 &0 \\
0 &1
\end{smallmatrix} \right)\right)$. 
При этом в $A$ не существует $T$-градуированной максимальной полупростой
подалгебры $B$, такой, что $A= B \oplus J(A)$.
\end{example}

\begin{proof}
Сперва заметим, что разложение  $A = A_{11} \oplus A_{12} \oplus A_{21} \oplus A_{22}$
действительно является $T$-градуировкой, поскольку $(a\, v_i w_j) (b\, v_\ell w_k) = ab (w_j v_\ell)\, v_i w_k$ для всех $a,b \in R$ и $1 \leqslant i,j,k,\ell \leqslant 2$, а $w_j v_\ell$~"--- матрица размера $1\times 1$, которая может быть отождествлена с соответствующим элементом поля $\mathbbm{k}$.

 Ясно, что $J(A) = \left(\begin{smallmatrix}
(X) & (X) \\
(X) & (X)
\end{smallmatrix}\right)$
и $A/J(A) \cong M_2(\mathbbm{k})$.
При этом $\left(\begin{smallmatrix}
X & 0 \\
0 & 0
\end{smallmatrix}\right) = X \left(\begin{smallmatrix}
1 & X \\
0 & 0
\end{smallmatrix}\right) \in A_{11}$
и $\left(\begin{smallmatrix}
0 & 0 \\
X & 0
\end{smallmatrix}\right)=X \left(\begin{smallmatrix}
0 & 0 \\
1 & X
\end{smallmatrix}\right) \in A_{21}$,
т.е. радикал Джекобсона $J(A)$ является градуированным идеалом,
а $A/J(A) \cong M_2(\mathbbm{k})$~"--- градуированной алгеброй.

Предположим, что существует такая $T$-градуированно простая максимальная подалгебра $B=\bigoplus_{1\leqslant i,j \leqslant 2} B_{ij}$, что $A = B \oplus J(A)$
и $B_{ij} \subseteq A_{ij}$.
Тогда ограничение на $B$ естественного сюръективного
гомоморфизма $\pi \colon A \twoheadrightarrow A/J(A)$
является градуированным изоморфизмом алгебр $B$ и $A/J(A)$.
В частности, $B_{ii}=\langle b_{ii} \rangle_\mathbbm{k}$, где $\pi(b_{ii})=e_{ii}$ при $i=1,2$.
Следовательно, $b_{11}=(1+\alpha X)\left(\begin{smallmatrix}
1 & X \\
0 & 0
\end{smallmatrix}\right)$ для некоторого $\alpha \in \mathbbm{k}$,
а $b_{22}=(1+\beta X)\left(\begin{smallmatrix}
0 & 0 \\
0 & 1
\end{smallmatrix}\right)$
для некоторого $\beta \in \mathbbm{k}$.
Тогда $$b_{11} b_{22} = (1+(\alpha+\beta)X)\left(\begin{smallmatrix}
0 & X \\
0 & 0
\end{smallmatrix}\right) = \left(\begin{smallmatrix}
0 & X \\
0 & 0
\end{smallmatrix}\right) \in J(A),$$
и мы получаем противоречие.
\end{proof}

Теорема~\ref{TheoremBGradedReesSemigroupGrSimple}
описывает полупростую часть $T$-градуированно простой алгебры.
Ниже изучается строение радикала, что завершает исследование всей $T$-градуированно простой алгебры.
В \S\ref{SectionTGradedReesExistence} 
будет показано, что полученное описание действительно является полной классификацией $T$-градуированно простых алгебр.

Для любого $r\in A$ будем обозначать элемент $x-xr$ (соответственно, $x-rx$) через $x(1-r)$ (соответственно,
через $(1-r)x$). Если алгебра $A$ без единицы, под символом $1$ можно понимать присоединённую единицу
алгебры $A^+:= \mathbbm{k}1\oplus A$.

\begin{lemma} \label{LemmaReesGrSimpleOtherProperties} 
Пусть $A$~"--- конечномерная $T$-градуированно простая алгебра, а $A=B\oplus J(A)$ (прямая сумма подпространств)~"---
разложение из теоремы~\ref{TheoremBGradedReesSemigroupGrSimple}. 
Тогда выполняются следующие свойства:
\begin{enumerate}
\item $J(A)^2 A = A J(A)^2=0$;
\item $B$~"--- простая подалгебра;
\item $A=A1_B A$;
\item $J(A)=(1-1_B)A 1_B \oplus 1_B A(1-1_B)\oplus J(A)^2$ (прямая сумма подпространств);
\item $J(A)^2= (1-1_B)A 1_B A (1-1_B) = (1-1_B)A (1-1_B)$.
\end{enumerate}
\end{lemma}
\begin{proof} Предложение 3 является немедленно следует из предложения 2 леммы~\ref{LemmaReesGradedFirstProperties}.

Пусть $f$~"--- примитивный центральный идемпотент алгебры $B$. 
Тогда $A=AfA$ и в силу предложения 3  леммы~\ref{LemmaReesGradedFirstProperties} $$B= 1_B (B\oplus J(A)) 1_B =  1_B A 1_B= 1_B AfA 1_B = 1_B A 1_B f 1_B A 1_B = BfB=Bf=fB.$$ Следовательно, $f=1_B$.
Отсюда алгебра
 $B$ проста, и мы получаем предложение 2. 

 Заметим, что $A(1-1_B),\, (1-1_B)A \subseteq J(A)$, так как
оба подпространства переходят в нуль при естественном сюръективном гомоморфизме $A \twoheadrightarrow A/J(A)$.
Используя равенство $B=1_B A 1_B$ и 
разложение Пирса по отношению к идемпотенту $1_B$, получаем
 $$J(A)=(1-1_B)A 1_B \oplus 1_B A (1-1_B) \oplus (1-1_B)A (1-1_B) \text{ (прямая сумма подпространств)}.$$
 Из предложения 3 теперь следует, что $$(1-1_B)A (1-1_B) = (1-1_B)A 1_B A (1-1_B) \subseteq J(A)^2.$$
Отсюда $$(1-1_B)A (1-1_B) J(A),\quad J(A)(1-1_B)A (1-1_B) \subseteq J(A)^3=0.$$  
  Поскольку $$A (1-1_B)A = A 1_B A (1-1_B)A \subseteq AJ(A)A = 0,$$
  получаем, что $J(A)^2 \subseteq (1-1_B)A (1-1_B)$,
  откуда и следуют предложения 1, 4 и 5.
\end{proof}

Заметим, что если матрица в определении полугруппы $T$ 
состоит всего лишь из одной строки, все $A_{1j}$
являются левыми идеалами, откуда
 $a - 1_B a \subseteq J(A)\cap A_{1j}=0$ для всех $a \in A_{1j}$ и $1 \leqslant j \leqslant n$. В этом случае элемент $1_B$
действует как левая единица на радикале $J(A)$.

Докажем теперь, что справедливость условия 1 леммы~\ref{LemmaReesGradedFirstProperties}
и условия 2 леммы~\ref{LemmaReesGrSimpleOtherProperties}
вместе с равенством $A^2=A$ эквивалентны градуированной $T$-простоте.

\begin{proposition}\label{PropositionReesGradedSimplicityCriterion}
Предположим, что основное поле $\mathbbm{k}$ совершенно, конечномерная алгебра $A/J(A)$ проста, $A^2=A$ и $A_{ij} \cap J(A) = 0$ для всех $1\leqslant i \leqslant m$ и $1\leqslant j \leqslant n$. Тогда алгебра $A$ является $T$-градуированно простой.
\end{proposition}
\begin{proof}
Пусть $I$~"--- ненулевой двухсторонний градуированный идеал алгебры $A$. Обозначим через $\pi \colon A \twoheadrightarrow A/J(A)$ естественный сюръективный гомоморфизм. Тогда $\pi(I)\ne 0$. Поскольку алгебра $A/J(A)$ проста, справедливы $\pi(I)=A/J(A)$
и $A=I+J(A)$. В силу обычной теоремы Веддербёрна~"--- Мальцева
существует максимальная полупростая подалгебра $B \subseteq I$, такая, что $I=B\oplus J(I)$ (прямая сумма подпространств).
Напомним, что $J(I)=J(A)\cap I$. Отсюда $A = B \oplus J(A)$ и  $\pi(A (1-1_B) A)=0$. Следовательно, $A(1-1_B) A \subseteq J(A)$.
Поскольку идеал $A (1-1_B) A$ градуированный, получаем $A(1-1_B)A=0$
и $ab = a 1_B b \in I$ для всех $a,b \in A$. Отсюда $A=A^2\subseteq I$
и $I=A$.
\end{proof}

Из леммы~\ref{LemmaReesGrSimpleOtherProperties} следует, что 
$$J(A)=1_B A (1-1_B) \oplus  (1-1_B)A 1_B \oplus J(A)^2=
\sum_{j=1}^n 1_B L_j (1-1_B) \oplus \sum_{i=1}^m   (1-1_B) R_i 1_B
\oplus J(A)^2.$$
Для $1\leqslant i\leqslant m$ и $1\leqslant j \leqslant n$ положим
   $$J^{10}_{ij}:=f'_i L_j (1-1_B) \quad \text{ и } \quad  J^{01}_{ij}:=(1-1_B)R_i f_j.$$
Кроме того, введём обозначение
   $$J^{10}_{*j} :=\bigoplus_{1\leqslant i \leqslant m} J^{10}_{ij} = 1_B L_j (1-1_B) \quad \text{ и } \quad  
     J^{01}_{i*}  : = \bigoplus_{1\leqslant j \leqslant n} J^{01}_{ij}  = (1-1_B) R_i 1_B.$$
     (Суммы являются прямыми в силу ортогональности соответствующих идемпотентов.)
    
     Ясно, что каждое подпространство $J^{10}_{*j}$ является левым $B$-подмодулем
радикала $J(A)$, а каждое каждое подпространство $J^{01}_{i*}$ является правым $B$-подмодулем
радикала $J(A)$.

 Покажем, что подмодули $J^{01}_{i*}$ и $J^{10}_{*j}$ являются компонентами, из которых строится
     радикал Джекобсона $J(A)$. Для начала определим отображение $\varphi$ из этих подмодулей
     в $B$, которое, как мы увидим ниже в теореме~\ref{TheoremReesGradedRadicalStructure},
     будет биективным отображением каждого из подмодулей на соответствующий односторонний идеал алгебры $B$.

Для всякого элемента $a \in J^{01}_{ij}$ или $a \in J^{10}_{ij}$
обозначим через $\varphi(a)$ проекцию на $B$ с ядром $J(A)$ однородной $A_{ij}$-компоненты элемента $a$.
Иными словами, $$a = (\varphi(a) + v) + w,\quad\text{где }\varphi(a) \in B,\ v\in J(A),\ \varphi(a)+v\in A_{ij},\ w\in \bigoplus\limits_{\substack{r\ne i \\ \text{или} \\ \ell \ne j}} A_{r\ell}.$$

\begin{theorem} \label{TheoremReesGradedRadicalStructure}
Пусть $A$~"--- конечномерная $T$-градуированно простая алгебра над полем $\mathbbm{k}$, а
$B$ и 
$f_{1}, \ldots, f_m,\ f_1', \ldots, f_n'$~"--- соответственно, градуированная 
подалгебра и системы ортогональных идемпотентов из теоремы~\ref{TheoremBGradedReesSemigroupGrSimple}.  
Тогда $$J(A)=\bigoplus_{i=1}^m J^{01}_{i*} \oplus \bigoplus_{j=1}^n J^{10}_{*j}
\oplus J(A)^2 \quad \text{ и } \quad
J(A)^2=\bigoplus_{i=1}^m \bigoplus_{j=1}^n  J^{01}_{i*}J^{10}_{*j}$$
(прямые суммы подпространств), а
если продолжить отображение $\varphi$ по $\mathbbm{k}$-линейности до отображения 
$\bigoplus_{i=1}^m J^{01}_{i*} \oplus \bigoplus_{j=1}^n J^{10}_{*j} \rightarrow B$,
то  оказывается, что $\varphi\bigl|_{\bigoplus_{j=1}^n J^{10}_{*j}}$~"--- гомоморфизм левых $B$-модулей,
  \begin{equation}\label{EquationReesSimpleJ10*jfj0}
     J^{10}_{*j}\cap \ker \varphi = 0,\qquad \varphi (J^{10}_{*j}) \cap B f_j =\varphi (J^{10}_{*j}) f_j= 0
        \quad    \text{ для всех } 1\leqslant j \leqslant n,
  \end{equation}
$\varphi\bigl|_{\bigoplus_{i=1}^m J^{01}_{i*}}$~"--- гомоморфизм правых $B$-модулей,
\begin{equation}\label{EquationReesSimplefi'J01i*0}
       J^{01}_{i*} \cap \ker \varphi = 0,\qquad \varphi (J^{01}_{i*}) \cap f_i' B =
        f_i'\varphi (J^{01}_{i*}) = 0 \quad \text{ для всех }1\leqslant i \leqslant m.
   \end{equation}
Более того, 
   \begin{equation}\begin{split}\label{EquationDecompAijReesSimple}
      A_{ij}=f_i' B f_j\oplus \left\{ \varphi(v)+v \mid  
          v \in J^{10}_{ij}\oplus J^{01}_{ij} \right\} \oplus \\
          \oplus \left\langle \varphi(v)\varphi(w)+v\varphi(w) +  \varphi(v) w + vw\mid  
            v \in J^{01}_{i*},\ w\in J^{10}_{*j} \right\rangle_\mathbbm{k}
            \end{split}\end{equation}
(прямая сумма подпространств) для всех  $1\leqslant i \leqslant m$, $1\leqslant j \leqslant n$.

Если $s\in\mathbb N$, $v_\ell \in J^{01}_{i*}$ и $w_\ell \in J^{10}_{*j}$
для всех $1\leqslant \ell \leqslant s$,
то $\sum_{\ell=1}^s v_\ell w_\ell = 0$, если только если $\sum_{\ell=1}^s \varphi(v_\ell) \varphi(w_\ell) = 0$.

Наконец, $B \cong M_k(D)$ для некоторого $k \in \mathbb N$ и алгебры с делением $D$,
причём
 \begin{eqnarray}\label{EquationDimJ01}
   \dim_{\mathbbm{k}} \bigoplus_{i=1}^m J^{01}_{i*} &\leqslant & (m - 1) \dim_{\mathbbm{k}} B = (m-1)k^2\dim_\mathbbm{k} D,\\
 \label{EquationDimJ10}
     \dim_{\mathbbm{k}} \bigoplus_{j=1}^n J^{10}_{*j}  &\leqslant &  (n - 1) \dim_{\mathbbm{k}} B = (n-1)k^2\dim_\mathbbm{k} D,\\
  \label{EquationDimJAReesSemiGr}
    \dim_{\mathbbm{k}} J(A) &\leqslant & (nm - 1) \dim_{\mathbbm{k}} B =(|T|-1) \dim_{\mathbbm{k}} B = (|T|-1)k^2\dim_\mathbbm{k} D.    
  \end{eqnarray}
\end{theorem}

\begin{proof}
В силу леммы~\ref{LemmaReesGrSimpleOtherProperties}
$$J(A)=1_B A (1-1_B) \oplus  (1-1_B)A 1_B \oplus J(A)^2=
\sum_{j=1}^n 1_B L_j (1-1_B) \oplus \sum_{i=1}^m   (1-1_B) R_i 1_B
\oplus J(A)^2.$$
Заметим, что если $\sum_{j=1}^n 1_B a_j (1-1_B) = 0$ для некоторых
$a_j \in L_j$, то $\sum_{j=1}^n 1_B a_j 1_B = \sum_{j=1}^n 1_B a_j$.
Так как $1_B A 1_B = B$ является градуированной подалгеброй, получаем, что $1_B a_j \in B \cap L_j$,
$1_B a_j 1_B = 1_B a_j$ и  все $1_B a_j (1-1_B) = 0$.
Следовательно, сумма $\bigoplus_{j=1}^n 1_B L_j (1-1_B)$ является прямой.
Аналогично, прямой является и сумма $\bigoplus_{i=1}^m (1-1_B) R_i 1_B$, причём
\begin{equation*}\begin{split}J(A)=\bigoplus_{j=1}^n 1_B L_j (1-1_B) \oplus \bigoplus_{i=1}^m (1-1_B) R_i 1_B
\oplus J(A)^2  = \bigoplus_{j=1}^n J^{10}_{*j}\oplus\bigoplus_{i=1}^m J^{01}_{i*} 
\oplus J(A)^2.\end{split}\end{equation*}
Используя предложение 5 из леммы~\ref{LemmaReesGrSimpleOtherProperties}, получаем \begin{equation}\label{EquationJASquareSum} J(A)^2 = \sum_{i=1}^m\sum_{j=1}^n
(1-1_B)R_i 1_B 1_B L_j (1-1_B) = \sum_{i=1}^m\sum_{j=1}^n J^{01}_{i*}J^{10}_{*j}.\end{equation}

Докажем теперь формулы \begin{equation}\label{EqVarphiLj(1-1B)}\varphi(a(1-1_B))=a(1_B-f_j) \text{ для всех } a \in 1_B L_j\end{equation}
и \begin{equation}\label{EqVarphi(1-1B)Ri}\varphi((1-1_B)a)=(1_B-f'_i)a \text{ для всех } a \in R_i1_B. \end{equation}
Из формул~\eqref{EqVarphiLj(1-1B)} и~\eqref{EqVarphi(1-1B)Ri} будет немедленно следовать, что $\varphi\bigl|_{\bigoplus_{j=1}^n J^{10}_{*j}}$~"--- гомоморфизм левых $B$-модулей,
а $\varphi\bigl|_{\bigoplus_{i=1}^m J^{01}_{i*}}$~"--- гомоморфизм правых $B$-модулей.

В силу линейности по $a$ достаточно доказать формулу~\eqref{EqVarphiLj(1-1B)}
для любого $a \in f_i' L_j$, а формулу~\eqref{EqVarphi(1-1B)Ri}~"---
для любого $a \in R_i f_j$.

Итак, пусть для определённости $a \in f_i' L_j$. Тогда $a-a f_j \in A_{ij}$, $a (1_B-f_j) \in B$, $a(1-1_B) \in J(A)$,  $\sum\limits_{\ell \ne j} a f_\ell \in \bigoplus\limits_{\ell \ne j} A_{i\ell}$. Таким образом,
$$ a(1-1_B)=(a-a f_j)-\sum_{\ell \ne j} a f_\ell = \bigl(a (1_B-f_j)+a(1-1_B)\bigr)-\sum_{\ell \ne j} a f_\ell$$
и формула~\eqref{EqVarphiLj(1-1B)} доказана. Формула~\eqref{EqVarphi(1-1B)Ri} доказывается аналогично. 

Предположим, что $\varphi(1_B a (1-1_B))=0$ для некоторого $a\in L_j$.
В силу~\eqref{EqVarphiLj(1-1B)} это означает, что $1_B a(1_B - f_j) = 0$, т.е. $1_B a f_j = 1_B a 1_B$
и $f_i' a f_j = f_i' a 1_B$ для всех $1\leqslant i \leqslant m$. Отсюда $$f_i' a (1-1_B)=f_i'a - f_i' a f_j \in A_{ij} \cap J(A)=0.$$
Следовательно, $1_B a (1-1_B) = 0$ и $$1_B L_j (1-1_B)\cap \ker \varphi = J^{10}_{*j} \cap \ker \varphi = 0.$$
Аналогично, $$(1-1_B) R_i 1_B \cap \ker \varphi = J^{01}_{i*} \cap \ker \varphi= 0.$$

Более того, $$\varphi(J^{10}_{*j}) \cap B f_j =\varphi (1_B L_j (1-1_B))
\cap B f_j \subseteq \varphi (1_B L_j (1-1_B))f_j = 1_B L_j (1_B-f_j) f_j = 0$$ для всех $1\leqslant j \leqslant n$,
и равенство~(\ref{EquationReesSimpleJ10*jfj0}) доказано.
Аналогично, $$\varphi(J^{01}_{i*}) \cap f_i' B = \varphi ((1-1_B) R_i 1_B)
\cap f_i' B  \subseteq f_i'\varphi ((1-1_B) R_i 1_B)
  =  f_i'(1_B-f_i') R_i 1_B = 0$$ для всех $1\leqslant i \leqslant m$, и равенство~(\ref{EquationReesSimplefi'J01i*0}) также доказано.
  
В силу предложения 2 леммы~\ref{LemmaReesGrSimpleOtherProperties}
существует изоморфизм $B \cong M_k(D)$ для некоторого $k \in \mathbb N$ и алгебры с делением $D$.
  Отсюда  \begin{equation*}\begin{split}\dim_{\mathbbm{k}} \bigoplus_{i=1}^m J^{01}_{i*}
= \sum_{i=1}^m \dim_\mathbbm{k} \varphi(J^{01}_{i*})
\leqslant \sum_{i=1}^m (\dim_\mathbbm{k} B  - \dim_\mathbbm{k} f_i' B)
    =\\ =(m - 1) \dim_{\mathbbm{k}} B = (m-1)k^2\dim_\mathbbm{k} D\end{split}\end{equation*}
   и равенство~(\ref{EquationDimJ01}) доказано.
   Равенство~(\ref{EquationDimJ10}) доказывается аналогично.
  
  Докажем теперь равенство~(\ref{EquationDecompAijReesSimple}).
  Пусть
  \begin{equation*}
  \begin{split} \tilde A_{ij}=\left(f_i' B f_j\oplus \left\lbrace \varphi(v)+v \mid  
v \in J^{10}_{ij}\oplus J^{01}_{ij} \right\rbrace\right) +\\+
    \left\langle \varphi(v)\varphi(w)+v\varphi(w) +  \varphi(v) w + vw\mid  
       v \in J^{01}_{i*},\ w\in J^{10}_{*j} \right\rangle_\mathbbm{k} .
  \end{split}\end{equation*}
  
Сперва покажем, что $\tilde A_{ij}=A_{ij}$. Для этого заметим, что если $a \in f'_i L_j$, то $$\varphi(a(1-1_B))+a(1-1_B) = a(1_B-f_j)+a(1-1_B)
= a- a f_j \in A_{ij}.$$ Следовательно, $\varphi(w)+w \in A_{ij}$
для всех $w\in J^{10}_{ij}$.
Аналогично, $\varphi(v)+v \in A_{ij}$
для всех $v\in J^{01}_{ij}$. Отсюда $\tilde A_{ij} \subseteq A_{ij}$.

Очевидно, что $$1_B \left(\bigoplus_{i=1}^m\bigoplus_{j=1}^n \tilde A_{ij} \right)  \subseteq \bigoplus_{i=1}^m\bigoplus_{j=1}^n \tilde A_{ij}$$
и $$ \left(\bigoplus_{i=1}^m\bigoplus_{j=1}^n \tilde A_{ij} \right) 1_B \subseteq \bigoplus_{i=1}^m\bigoplus_{j=1}^n \tilde A_{ij}.$$
Отсюда в силу леммы~\ref{LemmaReesGrSimpleOtherProperties} справедливы равенства $$1_B \left(\bigoplus_{i=1}^m\bigoplus_{j=1}^n \tilde A_{ij} \right) (1-1_B) = \bigoplus_{j=1}^n J^{10}_{*j}$$
и $$(1-1_B) \left(\bigoplus_{i=1}^m\bigoplus_{j=1}^n \tilde A_{ij} \right) 1_B = \bigoplus_{i=1}^m J^{01}_{i*}.$$
Следовательно,
$$\bigoplus_{i=1}^m J^{01}_{i*} \oplus \bigoplus_{j=1}^n J^{10}_{*j}
\subseteq \bigoplus_{i=1}^m\bigoplus_{j=1}^n \tilde A_{ij}.$$
Кроме того, из (\ref{EquationJASquareSum}) следует, что $$(1-1_B) \left(\bigoplus_{i=1}^m\bigoplus_{j=1}^n \tilde A_{ij} \right) (1-1_B)= J^2(A)$$ и
$J(A)^2 \subseteq \bigoplus_{i=1}^m\bigoplus_{j=1}^n \tilde A_{ij}$.
Отсюда $\bigoplus_{i=1}^m\bigoplus_{j=1}^n \tilde A_{ij} = A$ и
$\tilde A_{ij}=A_{ij}$. Теперь для завершения доказательства справедливости равенства~(\ref{EquationDecompAijReesSimple})
достаточно показать, что сумма в определении пространства $\tilde A_{ij}$
является прямой.

Пусть $v_\ell \in J^{01}_{i*}$ и $w_\ell \in J^{10}_{*j}$, где $1\leqslant \ell \leqslant s$, $s\in\mathbb N$.
Предположим, что $\sum_{\ell=1}^s \varphi(v_\ell) \varphi(w_\ell) = 0$.
Поскольку $\varphi\bigl|_{\bigoplus_{j=1}^n J^{01}_{i*}}$~"--- гомоморфизм правых $B$-модулей,
справедливо равенство 
$\varphi\left(\sum_{\ell=1}^s v_\ell \varphi(w_\ell)\right)=0$, откуда в силу предложения 3
получаем, что
$\sum_{\ell=1}^s v_\ell \varphi(w_\ell) = 0$.
Аналогично, $\sum_{\ell=1}^s \varphi(v_\ell) w_\ell = 0$.
Следовательно, $$\sum_{\ell=1}^s v_\ell w_\ell = \sum_{\ell=1}^s (\varphi(v_\ell)+v_\ell) (\varphi(w_\ell)+w_\ell)\in A_{ij}\cap J(A)=0.$$

Обратно, предположим, что $\sum_{\ell=1}^s v_\ell w_\ell = 0$
для некоторых $v_\ell \in J^{01}_{i*}$ и $w_\ell \in J^{10}_{*j}$, где $1\leqslant \ell \leqslant s$, $s\in\mathbb N$. Пусть $a=\sum_{\ell=1}^s (\varphi(v_\ell)+v_\ell) (\varphi(w_\ell)+w_\ell)$, $$b=\sum_{\ell=1}^s  (\varphi(\varphi(v_\ell)w_\ell)+\varphi(v_\ell)w_\ell)
 + \sum_{\ell=1}^s  (\varphi(v_\ell\varphi(w_\ell))+v_\ell\varphi(w_\ell)) - \sum_{\ell=1}^s \varphi(v_\ell)\varphi(w_\ell).$$
Тогда $a - b=\sum_{\ell=1}^s v_\ell w_\ell = 
 0$. Отсюда $b=a\in A_{ij}$.
 Однако
 \begin{equation*}
 \begin{split}
   b = \sum_{q=1}^m \sum_{\ell=1}^s  (f_q'\varphi(\varphi(v_\ell)w_\ell)+f_q'\varphi(v_\ell)w_\ell)
   +\\
   + \sum_{r=1}^n\sum_{\ell=1}^s  (\varphi(v_\ell\varphi(w_\ell))f_r+v_\ell\varphi(w_\ell)f_r) -
 \sum_{q=1}^m\sum_{r=1}^n \sum_{\ell=1}^s f'_q \varphi(v_\ell)\varphi(w_\ell) f_r.\end{split}
 \end{equation*}
 Рассматривая однородную компоненту элемента $b$, принадлежащую подпространству $A_{ij}$, т.е. 
слагаемое с $q=i$ и $r=j$, получаем, что $a=b=0$, поскольку в силу~(\ref{EquationReesSimpleJ10*jfj0}) и~(\ref{EquationReesSimplefi'J01i*0}) справедливы равенства
 $$f_i'\varphi(v_\ell)=\varphi(w_\ell)f_j = 0.$$
 Проектируя элемент  $a$ на $B$ с ядром $J(A)$, получаем $\sum_{\ell=1}^s \varphi(v_\ell) \varphi(w_\ell)=0$.
 Следовательно, $$\sum_{\ell=1}^s \varphi(v_\ell) \varphi(w_\ell) = \sum_{\ell=1}^s v_\ell \varphi(w_\ell)
= \sum_{\ell=1}^s \varphi(v_\ell) w_\ell = 0.$$

Теперь всё готово для того, чтобы доказать, что
сумма в определении подпространства $\tilde A_{ij}$
прямая.

Предположим, что
\begin{equation*}\begin{split}\sum_{\ell=1}^s \left(\varphi(v_\ell)\varphi(w_\ell)+v_\ell \varphi(w_\ell)
+  \varphi(v_\ell) w_\ell + v_\ell w_\ell \right) \in \\ \in f_i' B f_j\oplus \left\lbrace \varphi(v)+v \mid  
v \in J^{10}_{ij}\oplus J^{01}_{ij} \right\rbrace \subseteq B \oplus \bigoplus_{j=1}^n J^{10}_{*j}\oplus\bigoplus_{i=1}^m J^{01}_{i*}\end{split}\end{equation*} для некоторых $  
v_\ell \in J^{01}_{i*}$ и $w_\ell \in J^{10}_{*j}$.
В силу того, что $\sum_{\ell=1}^s v_\ell w_\ell \in J(A)^2$,
а $$\sum_{\ell=1}^s \left(\varphi(v_\ell)\varphi(w_\ell)+v_\ell \varphi(w_\ell)
+  \varphi(v_\ell) w_\ell\right) \in B \oplus \bigoplus_{j=1}^n J^{10}_{*j}\oplus\bigoplus_{i=1}^m J^{01}_{i*},$$ справедливо равенство $\sum_{\ell=1}^s  v_\ell w_\ell=0$.
В силу замечаний, сделанных выше, $$\sum_{\ell=1}^s \varphi(v_\ell)\varphi(w_\ell)
=\sum_{\ell=1}^s v_\ell\varphi(w_\ell) = \sum_{\ell=1}^s \varphi(v_\ell)w_\ell = 0,$$
а значит, и $\sum_{\ell=1}^s \left(\varphi(v_\ell)\varphi(w_\ell)+v_\ell \varphi(w_\ell)
+  \varphi(v_\ell) w_\ell + v_\ell w_\ell \right) = 0$.

В частности, сумма в определении подпространства $\tilde A_{ij}$
прямая, и равенство~(\ref{EquationDecompAijReesSimple}) доказано.

Докажем теперь, что сумма $J(A)^2=\sum_{i=1}^m \sum_{j=1}^n  J^{01}_{i*}J^{10}_{*j}$
также является прямой. Действительно, предположим, что $\sum_{i=1}^m \sum_{j=1}^n u_{ij} = 0$
для некоторых $u_{ij}\in J^{01}_{i*}J^{10}_{*j}$.
В силу~(\ref{EquationDecompAijReesSimple}) существует разложение $u_{ij}= a_{ij}-v_{ij}$, где $a_{ij}\in A_{ij}$, а $v_{ij}$
является линейной комбинацией однородных элементов из $B$
 и однородных элементов из  
$\left\lbrace\varphi(v)+ v\mid  
v \in J^{01}_{i*} \oplus J^{10}_{*j} \right\rbrace$.
Сгруппировав в $\sum_{i=1}^m \sum_{j=1}^n (a_{ij}-v_{ij}) = 0$
элементы по однородным компонентам, получаем, что всякий элемент $a_{ij}$
является линейной комбинацией элементов
из $f_i' B f_j$ и элементов из
$\left\lbrace\varphi(v)+ v\mid  
v \in J^{01}_{ij} \oplus J^{10}_{ij} \right\rbrace$.
Отсюда $$u_{ij}=(1-1_B)u_{ij}(1-1_B)=(1-1_B)(a_{ij}-v_{ij})(1-1_B) = 0,$$
и сумма $J(A)^2=\bigoplus_{i=1}^m \bigoplus_{j=1}^n  J^{01}_{i*}J^{10}_{*j}$
действительно является прямой.

Осталось доказать только неравенство~(\ref{EquationDimJAReesSemiGr}).
Однако это неравенство следует из леммы~\ref{LemmaMatrixIdealsProduct}
и неравенств~(\ref{EquationDimJ01}) и~(\ref{EquationDimJ10}):
    \begin{equation*}\begin{split}
      \dim_\mathbbm{k} J(A)^2 =  \sum_{i=1}^m \sum_{j=1}^n \dim_\mathbbm{k}  J^{01}_{i*}J^{10}_{*j}
                                =\sum_{i=1}^m \sum_{j=1}^n \dim_\mathbbm{k}  \varphi(J^{01}_{i*})\varphi(J^{10}_{*j})=\\
                               =\sum_{i=1}^m \sum_{j=1}^n \frac{\dim_\mathbbm{k} \varphi( J^{01}_{i*})\dim_\mathbbm{k} 
                               \varphi(J^{10}_{*j})}{k^2 \dim_\mathbbm{k} D} \leqslant \\
                                \leqslant \sum_{i=1}^m \sum_{j=1}^n \frac{(\dim_\mathbbm{k} B - \dim_\mathbbm{k} f_i' B )
                                    (\dim_\mathbbm{k} B - \dim_\mathbbm{k} B f_j)}{k^2 \dim_\mathbbm{k} D}=\\
                                =\frac{(\dim_\mathbbm{k} B)^2(n-1)(m-1)}{k^2 \dim_\mathbbm{k} D}=\\
                                =(n-1)(m-1) \dim_\mathbbm{k} B.   
   \end{split}\end{equation*}
\end{proof}

\section{Теоремы существования для градуированно простых алгебр}\label{SectionTGradedReesExistence}

В теоремах~\ref{TheoremBGradedReesSemigroupGrSimple} и~\ref{TheoremReesGradedRadicalStructure}
было получено описание $T$-градуированно простых алгебр $A$.
Было показано, что для $A$ существует разложение Веддербёрна~"--- Мальцева $B\oplus J(A)$,
в котором подалгебра $B$ градуирована, и что $J(A)$ является, грубо говоря,
прямой суммой левых и правых $B$-модулей, изоморфных некоторым левым и правым идеалам алгебры $B$,
которые, кроме этого, удовлетворяют некоторым другим ограничениям.
Для того, чтобы завершить классификацию, докажем теперь, что
любой такой набор левых и правых идеалов конечномерной
простой алгебры  $B$ задаёт $T$-градуированно простую алгебру $A$
с градуированной максимальной подалгеброй $B$.

Пусть $k,m,n\in\mathbb N$, а $D$~"--- тело.
Предположим, что $B \cong M_k(D)$, а
 $f_1, \ldots, f_n \in B$ и  $f'_1, \ldots, f'_n \in B$~"--- два набора идемпотентов
(некоторые из которых могут быть равны), таких, что $\sum_{i=1}^m f_i' = \sum_{j=1}^n f_j = 1_B$
и в каждом наборе идемпотенты попарно ортогональны.

Пусть $J^{10}_{*1}, \ldots, J^{10}_{*m}$ и 
$J^{01}_{1*}, \ldots, J^{01}_{n*}$~"--- соответственно, левые и правые $B$-модули, такие,
что существуют вложения 
   $$\varphi \colon J^{10}_{*j} \hookrightarrow B \quad \text{ и } \quad 
   \varphi \colon J^{01}_{i*} \hookrightarrow B,$$
   которые являются гомоморфизмами, соответственно,
   левых и правых $B$-модулей, причём 
   $$\varphi(J^{10}_{*j}) f_j = 0 \quad \text{ и  } \quad  f_i' \varphi (J^{01}_{i*})= 0.$$
(Для удобства будем обозначать оба отображения одной и той же буквой $\varphi$,
а также считать, что $\mathbb Z$-линейное отображение $\varphi$ определено на
аддитивной группе $\bigoplus_{i=1}^m J^{01}_{i*} \oplus \bigoplus_{j=1}^n J^{10}_{*j}$.)
Определим аддитивные группы $J_{ij}$ как изоморфные копии групп $\varphi(J^{01}_{i*})\varphi(J^{10}_{*j}) \subseteq B$ при $1\leqslant i \leqslant m$, $1\leqslant j \leqslant n$.
Пусть
    $$\Theta_{ij} \colon \varphi(J^{01}_{i*})\varphi(J^{10}_{*j}) 
    \mathrel{\widetilde\rightarrow}
J_{ij}$$ "--- соответствующие изоморфизмы, а 
  $$\mu \colon J^{01}_{i*} \times J^{10}_{*j}
\rightarrow J_{ij}$$ "--- $\mathbb Z$-билинейное отображение,
заданное равенством
$$\mu(v,w):=\Theta_{ij}(\varphi(v)\varphi(w))$$
при $v\in J^{01}_{i*}$ и $w\in J^{10}_{*j}$.
Продолжим $\mu$ по $\mathbb Z$-линейности до отображения 
   $$\mu \colon \bigoplus_{i=1}^m J^{01}_{i*} \times
\bigoplus_{j=1}^n J^{10}_{*j}
\rightarrow \bigoplus_{i=1}^m \bigoplus_{j=1}^n J_{ij} .$$
Обозначим через 
$Q$ матрицу размера $n\times m$, в каждой клетке которой стоит
формальная единица $e$.

\begin{theorem}\label{TheoremTGradedReesExistence}
Зададим на аддитивной группе $$A=B \oplus \bigoplus_{i=1}^m J^{01}_{i*} \oplus \bigoplus_{j=1}^n J^{10}_{*j} 
\oplus \bigoplus_{i=1}^m \bigoplus_{j=1}^n J_{ij}$$
умножение $$(b_1,v_1,w_1,u_1)(b_2,v_2,w_2,u_2) =(b_1 b_2, v_1 b_2, b_1 w_2, \mu(v_1, w_2))$$ 
для всех $b_1,b_2 \in B$, $v_1,v_2\in \bigoplus_{i=1}^m J^{01}_{i*}$, $w_1,w_2 \in \bigoplus_{j=1}^n J^{10}_{*j}$, $u_1,u_2 \in \bigoplus_{i=1}^m \bigoplus_{j=1}^n J_{ij}$ и
$\mathcal{M}(\{ e\}^{0},m,n;Q)$-градуировку   
  \begin{equation*}\begin{split}
    A_{ij}=(f_i' B f_j,0,0,0) \oplus \left\lbrace (\varphi(v),v,0, 0) \mid v \in 
J^{01}_{i*}f_j \right\rbrace 
 \oplus \left\lbrace (\varphi(w),0,w, 0) \mid w \in 
 f_i'J^{10}_{*j} \right\rbrace \oplus \\
  \oplus \left\langle 
(\varphi(v)\varphi(w), v\varphi(w), \varphi(v)w, \mu(v,w)) \mid 
v \in 
J^{01}_{i*},\ w \in J^{10}_{*j} \right\rangle_{\mathbb Z}.
  \end{split}\end{equation*}
Тогда $A$ является $\mathcal{M}(\{ e\}^{0},m,n;Q)$-градуированно простым
кольцом.
   \end{theorem}

\begin{proof} 
Непосредственная проверка показывает, что 
$A = \bigoplus_{i=1}^m \bigoplus_{j=1}^n A_{ij}$ и что $A$ действительно является
$\mathcal{M}(\{ e\}^{0},m,n;Q)$-градуированным
кольцом.

Ясно, что $J(A)=(0,\bigoplus_{i=1}^m J^{01}_{i*}, \bigoplus_{j=1}^n J^{10}_{*j}, 
 \bigoplus_{i=1}^m \bigoplus_{j=1}^n J_{ij})$, поскольку третья степень правой части равенства нулевая.
 
  Заметим, что для всех $1\leqslant i \leqslant m$ и $1\leqslant j \leqslant n$
  $$1_B (f_i' B f_j,0,0,0) 1_B = \left(f_i' B f_j,0,0,0\right),$$
 $$1_B \left\lbrace (\varphi(v),v,0, 0) \mid v \in 
J^{01}_{i*}f_j \right\rbrace 1_B \subseteq \left(\bigoplus_{\ell \ne i} f_\ell' B f_j,0,0,0\right),$$
$$1_B \left\lbrace (\varphi(w),0,w, 0) \mid w \in 
 f_i'J^{10}_{*j} \right\rbrace 1_B \subseteq \left(\bigoplus_{r \ne j} f_i' B f_r,0,0,0\right),$$
 $$1_B \left\langle 
(\varphi(v)\varphi(w), v\varphi(w), \varphi(v)w, \mu(v,w)) \mid 
v \in 
J^{01}_{i*},\ w \in J^{10}_{*j} \right\rangle_{\mathbb Z} 
  1_B \subseteq \left(\bigoplus_{\substack{\ell\ne i, \\ r \ne j}} f_\ell' B f_r,0,0,0\right).$$
  Таким образом, образы разных слагаемых прямой суммы из определения пространства $A_{ij}$
  лежат в разных градуированных компонентах алгебры $B$.
  Учитывая, что ограничения $\varphi$
  на $J^{10}_{*j}$ и 
$J^{01}_{i*}$ являются инъекциями, получаем, что если $1_B a 1_B  = 0$ для некоторого $a\in A_{ij}$, то $a = 0$. Из того, что $1_B J(A) 1_B = 0$, теперь следует, что $A_{ij} \cap J(A)=0$.
  
 Предположим, что $I$~"--- градуированный двухсторонний идеал алгебры $A$. Пусть $a\in I$, где $a\ne 0$,~"--- однородный элемент. В силу вышесказанного  $a=(b,u,v,w)$, где  $b\ne 0$. Отсюда $(1_B,0,0,0) a (1_B,0,0,0) = (b,0,0,0)\in I$.
Поскольку кольцо $B$ простое, справедливо включение $(B,0,0,0)\subseteq I$.
Тогда $(1_B,0,0,0) A \subseteq I$ и $A(1_B,0,0,0)\subseteq I$. Поскольку $$A = (B,0,0,0)+ (1_B,0,0,0) A
+ A(1_B,0,0,0) + (1_B,0,0,0) A^2 (1_B,0,0,0),$$
получаем, что $I=R$, т.е кольцо $A$ градуированно простое.
\end{proof}

В случае, когда тело $D$ является алгеброй над полем $\mathbbm{k}$, все указанные модули являются модулями над $\mathbbm{k}$-алгеброй $B$, а участвующие в определении отображения $\mathbbm{k}$-линейны, конструкция из теоремы~\ref{TheoremTGradedReesExistence} приводит к $T$-градуированно простой алгебре.
Если алгебра $B$ конечномерна, вложение левых $B$-модулей 
$\varphi  \colon J^{10}_{*j} \rightarrow \bigoplus\limits_{\substack{ r=1,\\r \ne j}}^n B f_r$
существует, если и только если $\dim_\mathbbm{k} J^{10}_{*j} \leqslant \dim_\mathbbm{k} B - \dim_\mathbbm{k} ( B f_j).$
Вложение правых $B$-модулей $$\varphi \colon J^{01}_{*i} \rightarrow \bigoplus_{\ell \ne i} f'_\ell B$$ существует, если и только если $$\dim_\mathbbm{k} J^{01}_{i*} \leqslant \dim_\mathbbm{k} B - \dim_\mathbbm{k} (f'_i B).$$

Из теоремы~\ref{TheoremEquivalenceSemigroupGradedSimple} 
следует, что градуировка на алгебре $A$ полностью определена
образами её компонент в алгебре $A/J(A)$.
Покажем теперь, что всякое такой набор образов задаёт некоторую $T$-градуировку.

\begin{theorem}\label{TheoremImagesOfGradedComponentsReesExistence}
Пусть $D$~"--- алгебра с делением над полем $\mathbbm{k}$, а $B \cong M_k(D)$.
Предположим, что $B = \sum_{i=1}^m \sum_{j=1}^n B_{ij}$,
где $B_{ij}$~"--- некоторые  такие подпространства алгебры $B$, что
 $B_{ij} B_{\ell r} \subseteq B_{ir}$
для всех $1\leqslant i,\ell \leqslant m$ и $1\leqslant j,r \leqslant n$.
Пусть $P=(p_{ij})_{i,j}$~"--- такая матрица размера $n\times m$, где $p_{ij}\in \{ 0,e\}$, что $B_{ij} B_{\ell r}=0$ для всех $(j,\ell)$ с $p_{j\ell} = 0$.
Тогда существует такая $\mathcal{M}(\{ e\}^{0},m,n;P)$-градуированно простая алгебра $A=\bigoplus_{i=1}^m \bigoplus_{j=1}^n A_{ij}$ и такой сюръективный гомоморфизм алгебр $\psi \colon A \rightarrow  B$, что $\ker \psi = J(A)$ и $\psi(A_{ij})=B_{ij}$ для всех 
 $1\leqslant i \leqslant m$ и $1\leqslant j \leqslant n$.
\end{theorem}
\begin{proof} Пусть $\bar L_j := \bigoplus_{i=1}^m B_{ij}$ для всех $1\leqslant j \leqslant n$
и $\bar R_i := \bigoplus_{j=1}^n B_{ij}$ для всех $1\leqslant i \leqslant m$. Тогда
$\bar L_1, \ldots, \bar L_n$~"--- левые идеалы, а $\bar R_1, \ldots, \bar R_m$~"--- правые идеалы.
Более того, поскольку алгебра $B \cong M_k(D)$ полупроста,
алгебра $B$ вполне приводима как левый и правый $B$-модуль.
Определим $\tilde L_j$ как левый $B$-подмодуль, дополняющий $\bar L_j\cap \sum_{\ell=1}^{j-1} \bar L_\ell$ до $\bar L_j$,  где $1\leqslant j \leqslant n$. Аналогично, определим $\tilde R_i$ 
как правый $B$-подмодуль, дополняющий $\bar R_i\cap \sum_{\ell=1}^{i-1} \bar R_\ell$ до $\bar R_i$, где $1\leqslant i \leqslant m$.
Тогда $\bigoplus_{\ell=1}^i \tilde R_\ell = \sum_{\ell=1}^i \bar R_\ell$ для всех $1\leqslant i \leqslant m$
и $\bigoplus_{\ell=1}^j \tilde L_\ell = \sum_{\ell=1}^j \bar L_\ell$ для всех $1\leqslant j \leqslant n$.
В частности, $B = \bigoplus_{i=1}^m \tilde R_i = \bigoplus_{j=1}^n \tilde L_j$.
Раскладывая $1_B$ в сумму элементов, соответственно, подпространств $\tilde L_1, \ldots, \tilde L_m$ и $\tilde R_1, \ldots, \tilde R_n$,
получаем два множества ортогональных идемпотентов $f_1, \ldots, f_m$ и $f_1', \ldots, f_n'$,
таких, что $\tilde L_j = B f_j$, $\tilde R_i = B f_i'$, $\sum_{i=1}^m f_i' = \sum_{j=1}^n f_j = 1_B$.

Пусть $W^{10}_{*j}:=\bar L_j (1_B-f_j) \subseteq \bar L_j$ для $1\leqslant j \leqslant n$ 
и $W^{01}_{i*}:=(1_B-f_i')\bar R_i \subseteq \bar R_i$ для $1\leqslant i \leqslant m$.
Тогда $\bar L_j = \bar L_j f_j \oplus W^{10}_{*j}$ и $\bar R_i = f_i' \bar R_i \oplus W^{01}_{i*}$ (прямые суммы, соответственно, левых и правых идеалов). Обозначим через $J^{10}_{*j}$
изоморфную копию модуля $W^{10}_{*j}$, а через
$J^{01}_{i*}$~"--- изоморфную копию модуля $W^{01}_{i*}$.
Обозначим соответствующие изоморфизмы $\varphi \colon J^{10}_{*j} \mathrel{\widetilde\rightarrow} W^{10}_{*j}$
левых $B$-модулей и $\varphi \colon J^{01}_{i*} \mathrel{\widetilde\rightarrow} W^{01}_{i*}$
правых $B$-модулей одной и той же буквой $\varphi$.
Теперь продолжим $\varphi$ до $\mathbbm{k}$-линейного отображения
$$\varphi \colon 
\bigoplus_{i=1}^m J^{01}_{i*} \oplus \bigoplus_{j=1}^n J^{10}_{*j} \to
\sum_{i=1}^m W^{01}_{i*} + \sum_{j=1}^n W^{10}_{*j} \subseteq B.$$

Пусть $A=\bigoplus_{i=1}^m\bigoplus_{j=1}^n A_{ij}$~"--- кольцо,
построенное в теореме~\ref{TheoremTGradedReesExistence}, которое в силу замечания
после теоремы является алгеброй над $\mathbbm{k}$.
Докажем, что алгебра $A$ удовлетворяет всем условиям теоремы~\ref{TheoremImagesOfGradedComponentsReesExistence}.
Действительно, определим $\psi$ как проекцию алгебры $$A=B \oplus \bigoplus_{i=1}^m J^{01}_{i*} \oplus \bigoplus_{j=1}^n J^{10}_{*j} 
\oplus \bigoplus_{i=1}^m \bigoplus_{j=1}^n J_{ij}$$ на пространство $B$ с ядром $$J(A)=\bigoplus_{i=1}^m J^{01}_{i*} \oplus \bigoplus_{j=1}^n J^{10}_{*j} 
\oplus \bigoplus_{i=1}^m \bigoplus_{j=1}^n J_{ij}.$$
Тогда 
   \begin{equation}\label{EqSemigrBijOne}\begin{split}
      \psi(A_{ij}) = f_i' B f_j' + \varphi(J^{01}_{i*}f_j) + \varphi(f_i'J^{10}_{*j})
+\varphi(J^{01}_{i*}) \varphi(J^{10}_{*j}) 
    =\\ = f_i' B f_j + (1-f_i')\bar R_i f_j + f_i' \bar L_j (1-f_j) +(1-f_i')\bar R_i \bar L_j (1-f_j)=\\
       =f_i' \bar R_i f_j + (1-f_i')\bar R_i f_j +f_i' \bar R_i \bar L_j (1-f_j) +(1-f_i')\bar R_i \bar L_j (1-f_j)
   =\\= \bar R_i f_j + \bar R_i \bar L_j (1-f_j) = \bar R_i \bar L_j f_j + \bar R_i \bar L_j (1-f_j) = \bar R_i \bar L_j.\end{split}
   \end{equation}
Поскольку алгебра $B$ полупроста, существуют идемпотенты $e_i'$ и $e_j$,
такие, что $\bar R_i = e'_i B$ и $\bar L_j = B e_j$.   
Следовательно,
   \begin{equation}\label{EqSemigrBijTwo}\begin{split}
B_{ij} \subseteq \bar R_i \cap \bar L_j = e_i' B e_j = \bar R_i \bar L_j \subseteq B_{ij},\end{split}
   \end{equation}
   откуда
   \begin{equation}\label{EqSemigrBijThree}
   \bar R_i \cap \bar L_j = \bar R_i \bar L_j = B_{ij}.
   \end{equation}
 Из~\eqref{EqSemigrBijOne} и~\eqref{EqSemigrBijThree} теперь следует,
 что $\psi(A_{ij}) = \bar R_i \bar L_j = B_{ij}$.
   Более того, при $p_{j\ell}=0$ справедливы равенства $B_{ij}B_{\ell r} = 0$, $\psi(A_{ij}A_{\ell r})=0$ и $A_{ij}A_{\ell r} = 0$. Последнее равенство следует из того, что $(\ker \psi)  \cap A_{ir}=J(A)\cap A_{ir}=0$.
\end{proof}

В силу равенств~\eqref{EqSemigrBijThree},
 всякое разложение алгебры $B=M_{k}(D)$ в сумму компонент $B_{ij}$, удовлетворяющих условиям теоремы~\ref{TheoremImagesOfGradedComponentsReesExistence},
однозначно задаётся таким набором левых идеалов
$\bar L_1, \ldots, \bar L_n$ и таким набором правых идеалов $\bar R_1,\ldots, \bar R_m$ алгебры
$B$, что  $B=\sum_{i=1}^m \bar R_i = \sum_{j=1}^n \bar L_j$. 
Обратно, любые такие наборы левых и правых идеалов задают подпространства $B_{ij}:=\overline{R}_{i} \cap \overline{L}_{j}$, удовлетворяющие условиям теоремы~\ref{TheoremImagesOfGradedComponentsReesExistence}
для некоторой матрицы $P$. В частности,
 получается следующее утверждение, дополняющее критерий изоморфности
 градуированно простых алгебр, доказанный в теореме~\ref{TheoremEquivalenceSemigroupGradedSimple}:
\begin{theorem}\label{TheoremIsoReesMatrixGradedSimple}
Пусть $A=\bigoplus_{i=1}^m \bigoplus_{j=1}^n A_{ij}$
и $A'=\bigoplus_{i=1}^m \bigoplus_{j=1}^n A'_{ij}$~"--- конечномерные $\mathcal{M}(\{ e\}^{0},m,n;P)$-градуированно простые алгебры над полем $\mathbbm{k}$ для некоторой матрицы $P$.
Введём обозначения $$L_j:=\bigoplus_{k=1}^m A_{kj},\ L'_j:=\bigoplus_{k=1}^m A'_{kj},\ 
R_i:=\bigoplus_{k=1}^n A_{ik}\text{ и }R'_i:=\bigoplus_{k=1}^n A'_{ik}.$$
Тогда если существует такой изоморфизм алгебр $\bar\varphi \colon A/J(A) \mathrel{\widetilde\rightarrow}  A'/J(A')$, что $\bar\varphi\left(\pi\left(R_i\right)\right)=\pi'\left(R'_i\right)$
и $\bar\varphi\left(\pi\left(L_j\right)\right)=\pi'\left(L'_j\right)$
для всех $1\leqslant i \leqslant m$ и $1\leqslant j \leqslant n$, где $\pi \colon A \twoheadrightarrow A/J(A)$ и $\pi' \colon A' \twoheadrightarrow A'/J(A')$ "--- естественные
сюръективные гомоморфизмы,
то существует изоморфизм $\varphi \colon A \mathrel{\widetilde\rightarrow} A'$ градуированных алгебр, такой,
что $\pi'\varphi=\bar\varphi\pi$.

Обратно, если $\varphi \colon A \mathrel{\widetilde\rightarrow} A'$~"--- изоморфизм градуированных алгебр,
мы всегда можем определить изоморфизм алгебр $\bar\varphi \colon A/J(A) \mathrel{\widetilde\rightarrow}  A'/J(A')$
при помощи равенств $\bar\varphi(\pi(a)) = \pi'\varphi(a)$ для всех $a \in A$.
При этом $\bar\varphi\left(\pi\left(R_i\right)\right)=\pi'\left(R'_i\right)$
и $\bar\varphi\left(\pi\left(L_j\right)\right)=\pi'\left(L'_j\right)$
для всех $1\leqslant i \leqslant m$ и $1\leqslant j \leqslant n$.
\end{theorem}
\begin{proof} Достаточно воспользоваться теоремой~\ref{TheoremEquivalenceSemigroupGradedSimple}, леммой~\ref{LemmaReesGrSimpleOtherProperties} 
и равенствами $\pi(A_{ij})= \pi(R_i) \cap \pi(L_j)$ и $\pi'(A'_{ij})= \pi'(R'_i) \cap \pi'(L'_j)$, которые следуют из~\eqref{EqSemigrBijThree}. 
\end{proof}

\newpage

\chapter{$\Omega$-алгебры, эквивалентность (ко)действий и $V$-универсальные (ко)действующие алгебры Хопфа} \label{ChapterOmegaAlg}

В данной главе мы, отталкиваясь от теоремы~\ref{TheoremGradEquivCriterion}, вводим понятие эквивалентности (ко)действий. Это понятие будет затем использовано в параграфе~\ref{SectionDoubleNumbers} для классификации модульных структур и в параграфе~\ref{SectionEquivApplToPolyIden} для доказательства гипотезы Ш.~Амицура~"--- Ю.\,А.~Бахтурина для всевозможных модульных структур на конкретной алгебре.
  Строится общая теория, которая позволяет изучать с единых позиций как универсальные биалгебры и алгебры Хопфа М.~Свидлера--Ю.\,И.~Манина--Д.~Тамбары, так и универсальные алгебры Хопфа, которые возникают при изучении эквивалентных (ко)модульных структур. При этом свойство двойственности, которое было известно в случае универсальной измеряющей коалгебры на алгебре и универсальной коизмеряющей алгебры~\cite{Tambara},
  продолжается на случай универсальных действующих и кодействующих биалгебр и алгебр Хопфа, в том числе для 
  эквивалентных (ко)модульных структур. Для того, чтобы одновременно изучать (ко)действия как на алгебрах, так и на коалгебрах, вводится понятие $\Omega$-алгебры.
  
  Результаты главы вошли в статью~\cite{ASGordienko21ALAgoreJVercruysse}, написанную совместно автором, Аной Агоре и Йоостом Веркрёйссе. Ане Агоре принадлежит идея рассмотреть аналоги биалгебр и алгебр Хопфа М.~Свидлера--Ю.\,И.~Манина--Д.~Тамбары, (ко)действующие на коалгебрах, а также первый вариант доказательства двойственности для универсальных действующих и кодействующих биалгебр и алгебр Хопфа. Автором работы результаты Аны Агоре были обобщены на случай $\Omega$-алгебр и было введено понятие $V$-универсальных (ко)действующих биалгебр и алгебр Хопфа, позволяющее включить в рассмотрение универсальные алгебры Хопфа эквивалентных (ко)модульных структур. Доказательства результатов на языке коммутативных диаграмм, приводимые ниже, также принадлежат автору данной работы. Йоосту Веркрёйссе принадлежит идея рассмотреть двойственность между универсальными (ко)действующими биалгебрами на $\Omega$- и $\Omega^*$-алгебрах, а также идея получить универсальные кодействующие биалгебры и алгебры Хопфа Ю.\,И.~Манина
  как $V$-универсальные кодействующие биалгебры и алгебры Хопфа для соответствующих алгебр~$V$.

\section{Линейные отображения, их (ко)носители и конечная топология}\label{SectionLinearMaps}
         
В этом параграфе вводятся понятия, которые затем используются в случае (ко)действий
и (ко)измерений.         
         
\subsection{Линейные отображения $P \otimes A \to B$}\label{SubsectionLinearMapsMod}

Пусть $\psi \colon P \otimes A \to B$~"--- линейное отображение для некоторых векторных пространств $P,A,B$ над полем $\mathbbm{k}$.
Назовём подпространство $\cosupp \psi := \psi(P \otimes (-)) \subseteq \mathbf{Vect}_\mathbbm{k}(A,B)$
 \textit{коносителем} отображения $\psi$.\label{DefCosupportPAB} 
 Здесь $\psi(p \otimes (-))$
рассматривается как линейное отображение $ A \to B$ для фиксированного элемента $p\in P$,
а $\psi(P \otimes (-)) := \lbrace \psi(p \otimes (-)) \mid p \in P \rbrace$.
Другими словами, коноситель $\psi$~"--- это пространство
всех линейных операторов
 $ A \to B$, отвечающих отображению $\psi$.

\begin{definition}
Пусть $\psi_i \colon P_i \otimes A_i \to B_i$~"--- линейные отображения для некоторых векторных пространств $P_i,A_i,B_i$, где $i=1,2$, и пусть $\varphi \colon A_1 \mathrel{\widetilde{\to}} A_2$
и $\xi \colon B_1 \mathrel{\widetilde{\to}} B_2$~"--- линейные биективные отображения.
Будем говорить, что пара $(\varphi, \xi)$ является \textit{эквивалентностью} отображений
$\psi_1$ и $\psi_2$, если $$\xi(\cosupp\psi_1)\varphi^{-1} = \cosupp\psi_2.$$
В этом случае будем говорить, что $\psi_1$ и $\psi_2$ \textit{эквивалентны} при помощи $(\varphi, \xi)$.
\end{definition}

Рассмотрим теперь линейные такие отображения $\psi_i \colon P_i \otimes A \to B$, $i=1,2$, где пространства $A$ и $B$ общие для обоих отображений.
Будем говорить, что $\psi_1$ \textit{грубее}, чем $\psi_2$, а 
$\psi_2$ \textit{тоньше}, чем $\psi_1$, и использовать обозначение $\psi_1 \preccurlyeq \psi_2$,
если $\cosupp \psi_1 \subseteq \cosupp \psi_2$.

\begin{remark}
Отображения $\psi_i \colon P_i \otimes A \to B$, $i=1,2$, эквивалентны при помощи $(\id_{A},\id_{B})$, если
и только если $\psi_1 \preccurlyeq \psi_2$
и $\psi_1 \succcurlyeq \psi_2$.
\end{remark}
\begin{remark}
Пусть группы $G_1$ и $G_2$ действуют автоморфизмами на алгебрах $A_1$ и $A_2$, соответственно.
Тогда эти действия индуцируют линейные отображения $\psi_i \colon \mathbbm{k}G_i \otimes A_i \to A_i$.
Предположим, что алгебры $A_1$ и $A_2$ изоморфны. Обозначим через
$\varphi \colon A_1 \to A_2$ один из таких изоморфизмов.
Тогда отображения $\psi_1$ и $\psi_2$ эквивалентны при помощи $(\varphi, \varphi)$,
если и только если эквивалентны соответствующие действия (см. \S\ref{SectionGroupActions}). 
Если $A_1=A_2$, то $\psi_1 \preccurlyeq \psi_2$, если и только если $G_1$-действие грубее,
чем $G_2$-действие.
\end{remark}

\subsection{Линейные отображения $A \to B \otimes Q$}\label{SubsectionLinearMapsComod}

Пусть $\rho \colon  A \to B \otimes Q$~"--- линейное
отображение для некоторых векторных пространств $A,B,Q$ над полем $\mathbbm{k}$.
Как и в случае комодульных структур, 
для таких отображений мы используем обозначения М.~Свидлера $\rho(a)=a_{(0)}\otimes a_{(1)}$,
в которых знак суммы для краткости не пишется.

Пусть $(a_\alpha)_\alpha$~"--- базис в алгебре $A$, а $(b_\beta)_\beta$~"--- базис в алгебре $B$.
Определим элементы $q_{\beta\alpha} \in Q$ при помощи равенств $\rho(a_\alpha)=\sum_{\beta} b_\beta \otimes q_{\beta\alpha}$. (Для фиксированного $\alpha$ только конечное число элементов $q_{\beta\alpha}$ ненулевое.)
Назовём линейную оболочку элементов  $q_{\beta\alpha}$ \textit{носителем отображения} $\rho$
и будем обозначать её через $\supp \rho$. \label{DefSupportABQ}
Вычисляя значения линейных функций, двойственных к $b_\beta$, 
легко видеть, что подпространство $\supp \rho$ 
не зависит от выбора базисов в пространствах $A$ и $B$, поскольку $\supp \rho$
совпадает с наименьшим подпространством $Q_1$ пространства $Q$, таким, что $\rho(A)\subseteq B \otimes Q_1$.
Для произвольного отображения $\rho$ определим линейное отображение $\hat\rho \colon Q^* \otimes A \to B$
по формуле $\hat\rho(q^*\otimes a):= q^*(a_{(1)})a_{(0)}$, где $q^*\in Q^*$ и $a\in A$. \label{DefHatRho}
Назовём подпространство $$\cosupp \rho := \cosupp \hat\rho = \hat\rho(Q^* \otimes (-)) \subseteq \mathbf{Vect}_\mathbbm{k}(A,B)$$
\textit{коносителем} отображения $\rho$.\label{DefCosupportABQ}

\begin{definition}
Пусть $\rho_i \colon  A_i \to B_i \otimes Q_i$~"--- линейные отображения для некоторых векторных пространств $A_i,B_i,Q_i$, где $i=1,2$, и пусть $\varphi \colon A_1 \mathrel{\widetilde{\to}} A_2$
и $\xi \colon B_1 \mathrel{\widetilde{\to}} B_2$~"--- линейные биективные отображения.
Будем говорить, что пара $(\varphi, \xi)$ является \textit{эквивалентностью} отображений
$\rho_1$ и $\rho_2$, если $$\xi(\cosupp\rho_1)\varphi^{-1} = \cosupp\rho_2.$$
В этом случае будем говорить, что $\rho_1$ и $\rho_2$ \textit{эквивалентны} при помощи $(\varphi, \xi)$.
\end{definition}

Рассмотрим теперь линейные такие отображения $\rho_i \colon  A \to B \otimes Q_i$, $i=1,2$, где пространства $A$ и $B$ общие для обоих отображений.
Будем говорить, что $\rho_1$ \textit{грубее}, чем $\rho_2$, а 
$\rho_2$ \textit{тоньше}, чем $\rho_1$, и использовать обозначение $\rho_1 \preccurlyeq \rho_2$,
если $\cosupp \rho_1 \subseteq \cosupp \rho_2$.

\begin{remark}
Отображения $\rho_i \colon  A \to B \otimes Q_i$, $i=1,2$, эквивалентны при помощи $(\id_{A},\id_{B})$, если
и только если $\rho_1 \preccurlyeq \rho_2$
и $\rho_1 \succcurlyeq \rho_2$.
\end{remark}
\begin{remark}
Пусть  $\Gamma_i \colon A_i=\bigoplus_{g\in G_i} A_i^{(g)}$~"--- групповые градуировки.
Обозначим через $\rho_i \colon A_i \to A_i \otimes \mathbbm{k}G_i$ соответствующие комодульные структуры.
(См. пример~\ref{ExampleComoduleGraded}.)
Предположим, что алгебры $A_1$ и $A_2$ изоморфны. Обозначим через
$\varphi \colon A_1 \to A_2$ один из таких изоморфизмов.
Тогда в силу теоремы~\ref{TheoremGradFinerCoarserCriterion} отображения $\rho_1$ и $\rho_2$ эквивалентны при помощи $(\varphi, \varphi)$,
если и только если эквивалентны градуировки $\Gamma_1$ и $\Gamma_2$. 
Если $A_1=A_2$, то $\rho_1 \preccurlyeq \rho_2$, если и только если градуировка $\Gamma_1$ грубее,
чем градуировка $\Gamma_2$.
\end{remark}

Следующее предложение и два следствия из него позволяют переформулировать отношение <<тоньше-грубее>> и понятие эквивалентности для отображений на языке, близком к тому языку, на котором эти понятия формулировались для градуировок, т.е. на языке носителей.

\begin{proposition}\label{PropositionEquivPrecLinearMapComodCriterion}
Пусть $\rho_i \colon  A_i \to B_i \otimes Q_i$~"--- линейные отображения для некоторых векторных пространств $A_i,B_i,Q_i$, где $i=1,2$, и пусть $\varphi \colon A_1 \mathrel{\widetilde{\to}} A_2$
и $\xi \colon B_1 \mathrel{\widetilde{\to}} B_2$~"--- линейные биективные отображения.
Тогда $$\xi(\cosupp\rho_1)\varphi^{-1} \supseteq \cosupp\rho_2,$$
если и только если существует линейное отображение
$\tau \colon \supp \rho_1 \to \supp \rho_2$,
такое, что следующая диаграмма коммутативна (отображения  $A_i \to B_i \otimes \supp \rho_i$
получены из $\rho_i$ ограничением области значений):
\begin{equation}\label{EqLinearMapComodABTauQ} \xymatrix{ A_1 \ar[r] \ar[d]^\varphi
 & B_1 \otimes \supp \rho_1 \ar@{-->}[d]^{\xi \otimes \tau} \\
A_2 \ar[r] & B_2 \otimes \supp \rho_2} \end{equation}
Такое отображение $\tau$, если оно существует, единственно и сюръективно.
\end{proposition}
\begin{proof}
Выберем базис $(a_\alpha)_\alpha$ в пространстве  $A_1$ и базис $(b_\beta)_\beta$ в пространстве $B_1$.
Тогда совокупности $(\varphi(a_\alpha))_\alpha$ и $(\xi(b_\beta))_\beta$ являются базисами
в пространствах $A_2$ и $B_2$, соответственно.

Определим элементы $p_{\beta\alpha} \in Q_1$ и $q_{\beta\alpha} \in Q_2$ при
помощи равенств $\rho_1(a_\alpha)=\sum_{\beta} b_\beta \otimes p_{\beta\alpha}$ и $\rho_2(\varphi(a_\alpha))=\sum_{\beta} \xi(b_\beta) \otimes q_{\beta\alpha}$.
Если линейное отображение $\tau\colon
\supp \rho_1 \to \supp \rho_2$, делающее диаграмму~\eqref{EqLinearMapComodABTauQ} коммутативной,
действительно существует, тогда
$\tau(p_{\beta\alpha})=q_{\beta\alpha}$ для всех $\alpha$ и $\beta$. В частности, отображение $\tau$ единственное и сюръективное.
Для того, чтобы показать, что $\tau$ действительно существует,
достаточно доказать, что из
$\sum_{\alpha,\beta}\lambda_{\beta\alpha} p_{\beta\alpha}=0$
(где только конечное число элементов $\lambda_{\beta\alpha} \in \mathbbm{k}$ не равно нулю)
всегда следует, что $\sum_{\alpha,\beta}\lambda_{\beta\alpha} q_{\beta\alpha}=0$.

Предположим, что $\sum_{\alpha,\beta}\lambda_{\beta\alpha} p_{\beta\alpha}=0$.
Обозначим через $b^\gamma \in B_1^*$ линейные функции, заданные условиями $b^{\gamma}(b_\beta)=\delta_{\beta\gamma}$, где
$\delta_{\beta\gamma}$~"--- символ Кронекера.
Так как для всех $\alpha$ и $\beta$ выполнено $p_{\beta\alpha}=(b^\beta \otimes \id_{Q_1})(\rho_1(a_\alpha))$, то
 $$\sum_{\alpha,\beta}\lambda_{\beta\alpha} (b^\beta \otimes \id_{Q_1})(\rho_1(a_\alpha))=0,$$
$$\sum_{\alpha,\beta}\lambda_{\beta\alpha} (b^\beta \otimes p^*)(\rho_1(a_\alpha))=0\text{ для всех }p^*\in Q_1^*,$$
$$\sum_{\alpha,\beta}\lambda_{\beta\alpha} b^\beta (\hat\rho_1(p^*\otimes a_\alpha))=0 \text{ для всех }p^*\in Q_1^*,$$
$$\sum_{\alpha,\beta}\lambda_{\beta\alpha} b^\beta (T_1 a_\alpha)=0 \text{ для всех }T_1\in \cosupp \rho_1,$$
Если $\xi(\cosupp\rho_1)\varphi^{-1} \supseteq \cosupp\rho_2$, то
$$\sum_{\alpha,\beta}\lambda_{\beta\alpha} b^\beta\bigl(\xi^{-1}T_2 \varphi(a_\alpha)\bigr)=0 \text{ для всех }T_2\in \cosupp \rho_2,$$
$$\sum_{\alpha,\beta}\lambda_{\beta\alpha} (b^\beta\xi^{-1} \otimes q^*)(\rho_2(\varphi(a_\alpha)))=0\text{ для всех }q^*\in Q_2^*,$$
$$\sum_{\alpha,\beta}\lambda_{\beta\alpha} (b^\beta\xi^{-1} \otimes \id_{Q_2})(\rho_2(\varphi(a_\alpha)))=0,$$
и $$\sum_{\alpha,\beta}\lambda_{\beta\alpha} q_{\beta\alpha}=0.$$
Следовательно, можно корректно определить  линейное отображение $\tau\colon
\supp \rho_1 \to \supp \rho_2$,
положив $\tau(p_{\beta\alpha})=q_{\beta\alpha}$ для всех $\alpha$ и $\beta$,
и диаграмма~\eqref{EqLinearMapComodABTauQ} станет коммутативной.

Обратно, предположим, что существует такое линейное отображение $\tau\colon
\supp \rho_1 \to \supp \rho_2$, что диаграмма~\eqref{EqLinearMapComodABTauQ} коммутативна. 
Для любого линейного отображения $T_2 \in \cosupp \rho_2$ существует такое $q^* \in Q_2^*$, что $T_2a_2 = (\id_{B_2} \otimes q^*) \rho_2(a_2)$
для всех $a_2\in A_2$. Поэтому для любого $a_1 \in A_1$ справедливо равенство
$$\xi^{-1} T_2 \varphi(a_1) = \xi^{-1}(\id_{B_2} \otimes q^*) \rho_2(\varphi(a_1))\stackrel{\text{\eqref{EqLinearMapComodABTauQ}}}{=} (\id_B \otimes q^*\tau)\rho_1(a_1)=
(\id_{B_1} \otimes \tau^*(q^*))\rho_1(a_1).$$
Следовательно, $\xi^{-1} T_2 \varphi \in \cosupp \rho_1$ и $\xi(\cosupp\rho_1)\varphi^{-1} \supseteq \cosupp\rho_2$.
\end{proof}

\begin{corollary}\label{CorollaryEquivLinearMapComodCriterion}
Пусть $\rho_i \colon  A_i \to B_i \otimes Q_i$~"--- линейные отображения для некоторых векторных пространств $A_i,B_i,Q_i$, где $i=1,2$, и пусть $\varphi \colon A_1 \mathrel{\widetilde{\to}} A_2$
и $\xi \colon B_1 \mathrel{\widetilde{\to}} B_2$~"--- линейные биективные отображения.
Тогда пара $(\varphi, \xi)$ является эквивалентностью отображений
$\rho_1$ и $\rho_2$, если и только если существует линейное биективное отображение
$\tau \colon \supp \rho_1 \mathrel{\widetilde{\to}} \supp \rho_2$,
такое, что следующая диаграмма коммутативна:
\begin{equation*} \xymatrix{ A_1 \ar[r] \ar[d]^\varphi
 & B_1 \otimes \supp \rho_1 \ar@{-->}[d]^{\xi \otimes \tau} \\
A_2 \ar[r] & B_2 \otimes \supp \rho_2} \end{equation*}
Такое отображение $\tau$, если оно существует, единственно.
\end{corollary}
\begin{proof} Если такое линейное биективное отображение
$\tau \colon \supp \rho_1 \mathrel{\widetilde{\to}} \supp \rho_2$ действительно существует, то достаточно применить предложение~\ref{PropositionEquivPrecLinearMapComodCriterion} к $\tau$, $\varphi$, $\xi$
и к $\tau^{-1}$, $\varphi^{-1}$, $\xi^{-1}$.

Обратно, если $(\varphi, \xi)$~"--- эквивалентность отображений $\rho_1$ и $\rho_2$, то, опять применяя
предложение~\ref{PropositionEquivPrecLinearMapComodCriterion} к $(\varphi, \xi)$ и к $(\varphi^{-1}, \xi^{-1})$, получаем, что существует как линейное отображение $\tau$, такое, что диаграмма с $\tau$, $\varphi$ и $\xi$ становится коммутативной, так и линейное отображение $\tau_0 \colon \supp \rho_2 \mathrel{\widetilde{\to}} \supp \rho_1$, такое, что становится коммутативной диаграмма с $\tau_0$, $\varphi^{-1}$ и $\xi^{-1}$.
Отсюда $(\id_{A_1} \otimes \tau_0\tau)\rho_1 = \rho_1$ и $(\id_{A_2} \otimes \tau\tau_0)\rho_2 = \rho_2$.
Учитывая, что, кроме этого, $(\id_{A_1} \otimes \id_{\supp \rho_1})\rho_1 = \rho_1$ и $(\id_{A_2} \otimes \id_{\supp \rho_2})\rho_2 = \rho_2$, и применяя условие единственности в предложении~\ref{PropositionEquivPrecLinearMapComodCriterion} для $\id_{\supp \rho_1}$, $\id_{A_1}$, $\id_{B_1}$ и $\id_{\supp \rho_2}$, $\id_{A_2}$, $\id_{B_2}$, получаем,
что $\tau_0\tau = \id_{\supp \rho_1}$ и $\tau\tau_0 = \id_{\supp \rho_2}$, т.е. $\tau$ действительно 
линейное биективное отображение.
\end{proof}

\begin{corollary}\label{CorollaryPrecLinearMapComodCriterion}
Пусть $\rho_i \colon  A \to B \otimes Q_i$~"--- линейные отображения для некоторых векторных пространств $A,B,Q_i$, где $i=1,2$. Тогда 
 $\rho_1 \succcurlyeq \rho_2$, если и только если существует линейное отображение $\tau\colon
\supp \rho_1 \to \supp \rho_2$, такое, что следующая диаграмма коммутативна:
\begin{equation*} \xymatrix{ A \ar[r] 
\ar[rd]
 & B \otimes \supp \rho_1 \ar@{-->}[d]^{\id_B \otimes \tau} \\
& B \otimes \supp \rho_2} \end{equation*}
Отображение $\tau$, если оно существует, единственно и сюръективно.
\end{corollary}
\begin{proof} Достаточно в предложении~\ref{PropositionEquivPrecLinearMapComodCriterion}
положить $A_1=A_2=A$, $B_1=B_2=B$, $\varphi = \id_A$ и $\xi = \id_B$.
\end{proof}

\subsection{Замкнутые поточечно конечномерные подпространства в $\mathbf{Vect}_\mathbbm{k}(A,B)$}\label{SubsectionClLocFinSubspaces}

В этом подпараграфе мы получим необходимые и достаточные условия (см. теорему~\ref{TheoremVectDualization})
для того, чтобы подпространство $V \subseteq \mathbf{Vect}_\mathbbm{k}(A,B)$ было дуализируемо,
т.е. совпадало с $\cosupp \rho$ для некоторого линейного отображения $\rho \colon A \to B \otimes Q$.

Напомним, что если $X$ и $Y$~"--- топологические пространства,
 тогда на множестве $\mathbf{Top}(X,Y)$ непрерывных отображений $X\to Y$ 
можно задать \textit{компактно-открытую топологию},
предбазой которой являются всевозможные множества $$W(K,U)=\lbrace f \in \mathbf{Top}(X,Y) \mid
f(K) \subseteq U \rbrace,$$ где $K$~"--- компактное подмножество топологического пространства $X$,
а $U$~"--- открытое подмножество топологического пространства $Y$.

Если $X$ и $Y$ снабжены \textit{дискретной топологией}, т.е. 
все подмножества множеств $X$ и $Y$ одновременно и открыты, и замкнуты,
компактно-открытая топология на множестве $\mathbf{Top}(X,Y)$ (которое теперь совпадает со множеством $\mathbf{Sets}(X,Y)$ всех отображений $X\to Y$) называется \textit{конечной топологией}.

Поскольку подмножество дискретного топологического пространства компактно,
если и только если оно конечно,
база конечной топологии состоит из всех множеств $$W_{x_1,\ldots, x_n; y_1,\ldots,y_n}=\lbrace f \in \mathbf{Sets}(X,Y) \mid
f(x_1)=y_1,\ldots, f(x_n)=y_n \rbrace,$$ где $n\in\mathbb N$, $x_1,\ldots,x_n \in X$,
$y_1,\ldots,y_n \in Y$. Легко видеть, что множества $W_{x_1,\ldots, x_n; y_1,\ldots,y_n}$
одновременно и открыты, и замкнуты.

Ниже мы считаем множества $\mathbf{Vect}_\mathbbm{k}(A,B)$, где $A$ и $B$~"--- векторные пространства
над полем~$\mathbbm{k}$,
наделёнными конечной топологией.

\begin{lemma}\label{LemmaDualizationClosed}
Пусть $\rho \colon  A \to B \otimes Q$~"--- линейное
отображение для некоторых векторных пространств $A,B,Q$ над полем $\mathbbm{k}$.
Тогда
подмножество $\hat\rho(Q^* \otimes (-))\subseteq \mathbf{Vect}_\mathbbm{k}(A,B)$ замкнуто в конечной топологии.
\end{lemma}
\begin{proof}
Выберем базис $(a_\alpha)_\alpha$ в пространстве  $A$ и базис $(b_\beta)_\beta$ в пространстве $B$.
Как и прежде, определим элементы $q_{\beta\alpha} \in Q$ при
помощи равенств $\rho(a_\alpha)=\sum_{\beta} b_\beta \otimes q_{\beta\alpha}$.
(Напомним, что для фиксированного $\alpha$ только конечное число элементов $q_{\beta\alpha}$ ненулевое.)

Пусть $f \in \mathbf{Vect}_\mathbbm{k}(A,B)$. Определим элементы $f_{\beta\alpha} \in \mathbbm{k}$ при
помощи равенств
$f(a_\alpha)=\sum_{\beta} f_{\beta\alpha} b_\beta$. (Как и в случае с $q_{\beta\alpha}$,
для фиксированного $\alpha$ только конечное число элементов $f_{\beta\alpha}$ ненулевое.)
Тогда $f=\hat\rho(q^* \otimes (-))$ для некоторого $q^* \in Q^*$,
если и только если \begin{equation}\label{EqSystemFqalphabeta} q^*(q_{\beta\alpha})= f_{\beta\alpha}\text{ для всех }\alpha\text{ и }\beta.\end{equation}

Если существуют $n\in\mathbb N$ и $\alpha_1, \ldots, \alpha_n$, $\beta_1, \ldots, \beta_n$,
такие, что для любого $q^* \in Q^*$ хотя бы одно из равенств $q^*(q_{\beta_i\alpha_i})= f_{\beta_i\alpha_i}$,
где $1\leqslant i \leqslant n$, нарушается, то открытая окрестность $$\lbrace
g \in \mathbf{Vect}_\mathbbm{k}(A,B) \mid g(a_{\alpha_1})=f(a_{\alpha_1}),
\ldots, g(a_{\alpha_n})=f(a_{\alpha_n}) \rbrace$$ отображения $f$ не пересекается с
 $\hat\rho(Q^* \otimes (-))$. Поэтому достаточно доказать, что либо для некоторого $q^* \in Q^*$
 справедливы все равенства \eqref{EqSystemFqalphabeta}, либо среди них можно выбрать
 конечное число таких равенств, что при любом $q^* \in Q^*$ среди них хотя бы одно всё равно нарушается.

Предположим, что для любых конечных подмножеств $\alpha_1, \ldots, \alpha_n$, $\beta_1, \ldots, \beta_n$
индексов система линейных уравнений \begin{equation*}
q^*(q_{\beta_i\alpha_i})= f_{\beta_i\alpha_i},\ 1\leqslant i \leqslant n,\end{equation*} имеет решение $q^* \in Q^*$. Тогда для всех $\lambda_i \in \mathbbm{k}$, таких, что $\sum_{i=1}^n \lambda_i q_{\beta_i\alpha_i} = 0$,
выполнено $\sum_{i=1}^n \lambda_i f_{\beta_i\alpha_i} = 0$.
Другими словами, мы можем исключить из системы те уравнения, в которых левые части линейно зависимы с левыми частями оставшихся уравнений.
Обозначим через $(q_\gamma)_\gamma$ максимальную линейно независимую подсистему системы элементов $(q_{\beta\alpha})_{\alpha,\beta}$. (Последняя существует в силу леммы Цорна.) Обозначим через $f_\gamma$ элемент $f_{\beta\alpha}$, соответствующий
элементу $q_\gamma = q_{\beta\alpha}$.
Тогда система линейных уравнений 
\begin{equation*} q^*(q_\gamma)= f_\gamma\text{ для всех }\gamma\end{equation*}
эквивалентна системе \eqref{EqSystemFqalphabeta}.
Поскольку элементы $q_\gamma$ линейно независимы, мы можем включить $(q_\gamma)_\gamma$
в базис пространства $Q$, и система~\eqref{EqSystemFqalphabeta} будет иметь решение.
В этом случае $f \in \hat\rho(Q^* \otimes (-))$.

Отсюда для $f \in \mathbf{Vect}_\mathbbm{k}(A,B)$ есть только две возможности: либо $f \in \hat\rho(Q^* \otimes (-))$, либо у элемента $f$ существует такая окрестность, которая не пересекается с $\hat\rho(Q^* \otimes (-))$. Следовательно, подмножество $\hat\rho(Q^* \otimes (-))$ действительно замкнуто в конечной топологии.
\end{proof}

Дадим ещё одно важное определение.

Будем говорить, что подпространство $V\subseteq \mathbf{Vect}_\mathbbm{k}(A,B)$ \textit{поточечно конечномерно},
если для любого $a\in A$
подпространство $Va := \lbrace va \mid v\in V \rbrace \subseteq B$ конечномерно.

\begin{lemma}\label{LemmaDualizationPointwiseFinDim}
Пусть $\rho \colon  A \to B \otimes Q$~"--- линейное
отображение для некоторых векторных пространств $A,B,Q$ над полем $\mathbbm{k}$.
Тогда
подпространство $\hat\rho(Q^* \otimes (-))\subseteq \mathbf{Vect}_\mathbbm{k}(A,B)$ поточечно конечномерно.
\end{lemma}
\begin{proof} Для любого $a\in A$ существуют $n\in\mathbb N$, $a_i \in A$ и $q_i \in Q$, такие,
что
$$\rho(a)=a_1\otimes q_1 + \ldots + a_n \otimes q_n.$$ Следовательно $\hat\rho(q^*\otimes a)=
\sum_{k=1}^n q^*(q_j)a_j \in \langle a_1, \ldots, a_n \rangle_\mathbbm{k}$ для всех $q^* \in Q^*$
и $\dim_\mathbbm{k}(Va)\leqslant n$.
\end{proof}

Хотя следующий результат не является новым, мы приводим его доказательство для удобства читателя:
\begin{lemma}\label{LemmaVVstarstarOnFinite} Пусть $V$~"--- векторное пространство. 
Тогда для любого конечномерного подпространства $T \subseteq V^*$ и любого элемента $v^{**} \in V^{**}$ существует элемент $v \in V$, такой, что $v\bigr|_{T}$ (рассматриваемый как элемент пространства $V^{**}$)
совпадает с $v^{**}\bigr|_{T}$.
\end{lemma}
\begin{proof}
Выбирая базис в $V$, можно представить любую линейную функцию $v^{*} \in V^{*}$ как (возможно, бесконечную)
строку значений функции $v^{*}$ на элементах базиса.
Выберем базис $v^*_1, \ldots, v^*_n$ в  $T$. 
Применяя метод Жордана~"--- Гаусса к строчкам, соответствующим
линейным функциям  $v^*_1, \ldots, v^*_n$, мы можем перейти к другому такому базису
$\tilde v^*_1, \ldots, \tilde v^*_n$ пространства $V^{*}$,
что соответствующие строки содержат единичную подматрицу размера $n\times n$.
Выберем базисные элементы $v_1, \ldots, v_n$, отвечающие столбцам этой матрицы.
Тогда $\tilde v^*_i(v_j)=\delta_{ij}$, и
нам достаточно положить $v:= \sum_{i=1}^n v^{**}(\tilde v^*_i)v_i$.
\end{proof}

Теперь мы готовы доказать основную теорему, характеризующую дуализируемые подпространства
$V \subseteq \mathbf{Vect}_\mathbbm{k}(A,B)$

\begin{theorem}\label{TheoremVectDualization}
Пусть $V \subseteq \mathbf{Vect}_\mathbbm{k}(A,B)$~"--- подпространство
для некоторых векторных пространств $A$ и $B$
над полем $\mathbbm{k}$. Тогда $V=\hat\rho(Q^* \otimes (-))$
для некоторого векторного пространства $Q$ и линейного отображения $\rho \colon A \to B \otimes Q$,
если и только если $V$ поточечно конечномерно и замкнуто в конечной топологии.
\end{theorem}
\begin{proof}
Необходимость была доказана в леммах~\ref{LemmaDualizationClosed} и~\ref{LemmaDualizationPointwiseFinDim}.

Докажем достаточность. Предположим, что пространство $V$
поточечно конечномерно и замкнуто в конечной топологии.
Снова выберем базис $(a_\alpha)_\alpha$ в $A$ и базис $(b_\beta)_\beta$ в $B$.
Определим отображение $\rho  \colon A \to B \otimes V^*$
следующим образом. В силу поточечной конечномерности
пространства $V$ для любого $\alpha$
существуют элементы $n\in\mathbb N$ и $\beta_1, \ldots, \beta_n$,
такие, что $Va_\alpha \subseteq \langle b_{\beta_1}, \ldots, b_{\beta_n} \rangle_\mathbbm{k}$.
Следовательно, существуют $v_i^* \in V^*$, такие, что
$f(a_\alpha)=\sum_{i=1}^n v_i^*(f)b_{\beta_i}$ для всех $f\in V$.
Пусть $\rho(a_\alpha) := \sum_{i=1}^n b_{\beta_i} \otimes v_i^*$.

Очевидно, что $V \subseteq \hat\rho(V^{**} \otimes (-))$.
Докажем обратное включение. Пусть $v^{**} \in V^{**}$.
Мы утверждаем, что тогда $\hat\rho(v^{**} \otimes (-)) \in V$.

Определим элементы $v^*_{\beta\alpha} \in V^*$ при помощи равенства $\rho(a_\alpha) = \sum_{\beta} b_{\beta} \otimes v^*_{\beta\alpha}$.
(Для каждого $\alpha$ только конечное число элементов $v^*_{\beta\alpha}$ ненулевое.)
Предположим, что $\hat\rho(v^{**} \otimes (-)) \notin V$. 
Поскольку $V$ замкнуто в конечной топологии,
существуют такие
 $n\in\mathbb N$ и $\alpha_1, \ldots, \alpha_n$,
 что для любого $f \in V$ существует такой номер $1\leqslant i \leqslant n$,
 что $$\sum_{\beta} v^*_{\beta\alpha_i}(f) b_{\beta}=f(a_{\alpha_i}) \ne \hat\rho(v^{**} \otimes a_{\alpha_i})= \sum_{\beta} v^{**}(v^*_{\beta\alpha_i}) b_{\beta}.$$
Другими словами, система линейных уравнений
\begin{equation}\label{EqSystemFVstaralphabeta} v^*_{\beta\alpha_i}(f)=v^{**}(v^*_{\beta\alpha_i})\text{ для всех }\beta\text{ и всех }1\leqslant i \leqslant n\end{equation}
не имеет решения $f \in V$. 
Поскольку для фиксированного $i$ только конечное число элементов $v^*_{\beta\alpha_i}$ ненулевое,
система~\eqref{EqSystemFVstaralphabeta} состоит из конечного числа уравнений.
В силу леммы~\ref{LemmaVVstarstarOnFinite} мы получаем противоречие, и $\hat\rho(v^{**} \otimes (-)) \in V$.
Следовательно, $V = \hat\rho(V^{**} \otimes (-))$.
\end{proof}

\begin{remark}
Если оба пространства $A$ и $B$ конечномерны, то конечная топология на $\mathbf{Vect}_\mathbbm{k}(A,B)$
дискретна и все подпространства в $\mathbf{Vect}_\mathbbm{k}(A,B)$ замкнуты и поточечно конечномерны.
\end{remark}

Леммы, которые доказываются ниже, будут использованы в \S\ref{SectionDuality(Co)measurings} и~\S\ref{SectionCoactions}.
\begin{lemma}\label{LemmaClosurePwFD}
Пусть $V \subseteq \mathbf{Vect}_\mathbbm{k}(A,B)$~"--- подпространство
для некоторых векторных пространств $A$ и $B$
над полем $\mathbbm{k}$. Если $V$ поточечно конечномерно,
то его замыкание $\overline V$ в конечной топологии также является поточечно конечномерным.
\end{lemma}
\begin{proof}
Пусть $f\in \overline V$. Тогда для всех $n\in\mathbb N$ и $a_1, \ldots, a_n \in A$ существует $g\in V$,
такое, что $f(a_i)=g(a_i)$, $1\leqslant i \leqslant n$.
Следовательно, для любого $a\in A$ справедливо равенство $\overline V a = Va$, и $\overline V$
действительно поточечно конечномерно.
\end{proof}

\begin{lemma}\label{LemmaLinearMap(Co)modClosure}
Пусть $A,B,P$~"--- векторные пространства над полем $\mathbbm{k}$.
Пусть $\rho \colon A \to B \otimes P^*$~"--- линейное отображение, $\rho(a)=a_{(0)}\otimes a_{(1)}$
при $a\in A$. Определим линейное отображение $\rho^\vee \colon P \otimes A \to B$ по формуле $\psi(p\otimes a):=a_{(1)}(p)a_{(0)}$
для всех $a\in A$, $p\in P$.
Тогда $\cosupp \rho = \overline{\vphantom{|}\cosupp\rho^\vee}$ (замыкание берётся в конечной топологии).
\end{lemma}
\begin{proof}
Очевидно, $\cosupp \rho^\vee \subseteq \cosupp \rho$. В силу леммы~\ref{LemmaDualizationClosed}
подпространство $\cosupp \rho$ замкнуто. Следовательно, достаточно доказать, что
для любого $f \in \cosupp \rho$, $n\in\mathbb N$, $a_1, \ldots, a_n \in A$ существует
такое $g \in \cosupp \psi$, что $f(a_i)=g(a_i)$ для всех $1\leqslant i \leqslant n$. 

Снова выберем базис $(a_\alpha)_\alpha$ в $A$ и базис $(b_\beta)_\beta$ в $B$
и определим  $q_{\beta\alpha} \in P^*$ при помощи равенств $\rho(a_\alpha)=\sum_{\beta} b_\beta \otimes q_{\beta\alpha}$.
В силу леммы~\ref{LemmaVVstarstarOnFinite}
для любого $n\in\mathbb N$, любых индексов $\alpha_1, \ldots, \alpha_n$ и любого $p^{**}\in P^{**}$
существует
  $p\in P$, такое, что
  $q_{\beta\alpha_i}(p)=p^{**}(q_{\beta\alpha_i})$
для всех $1\leqslant i \leqslant n$ и всех $\beta$  
  (здесь снова используется тот факт, что для любого $i$ только конечное число элементов $q_{\beta\alpha}$
 ненулевое).
 Следовательно, значения линейных отображений $A \to B$,
 соответствующих элементам $p$ и $p^{**}$ совпадают на элементах $a_{\alpha_1}, \ldots, a_{\alpha_n}$.
 Поскольку всякий элемент пространства $A$ 
может быть записан как конечная линейная комбинация некоторых из элементов $a_\alpha$,
 утверждение леммы доказано.
\end{proof}

         \section{$\Omega$-алгебры}\label{SectionOmegaAlgebras}

Пусть $\Omega$~"--- это такое множество, что 
для любого его элемента $\omega \in \Omega$ заданы два целых неотрицательных числа $s(\omega)$
и $t(\omega)$. 
Будем называть $s(\omega)$ \textit{арностью} символа $\omega$,
а $t(\omega)$~"--- \textit{коарностью} символа $\omega$.
 Назовём \textit{$\Omega$-алгеброй над полем $\mathbbm{k}$} векторное пространство $A$ над полем $\mathbbm{k}$,
 на котором для любого $\omega \in \Omega$ заданы линейные отображения $\omega_A \colon A^{\otimes s(\omega)} \to
A^{\otimes t(\omega)}$. (Когда это не будет приводить к недоразумениям, мы будем опускать
индекс $A$ и обозначать соответствующее отображение также через $\omega$.) Здесь мы полагаем $A^{\otimes 0} := \mathbbm{k}$.
 
 \begin{example} $\varnothing$-алгебры~"--- это просто векторные пространства.
 \end{example}
 \begin{example}\label{ExampleAlgebra} $\lbrace \mu \rbrace$-алгебры, где $s(\mu)=2$, $t(\mu)=1$,"---
 это в точности необязательно ассоциативные алгебры $A$
 с умножением $\mu \colon A  \otimes A \to A$.
 \end{example}
  \begin{example}\label{ExampleUnitalAlgebra} Алгебры $A$
 с умножением $\mu \colon A  \otimes A \to A$
 и единицей $1_A$
 являются примером $\lbrace \mu, u \rbrace$-алгебр, где $s(\mu)=2$, $t(\mu)=1$, $s(u)=0$, $t(u)=1$,
 а линейное отображение $u \colon \mathbbm{k} \to A$ задаётся условием $u(\alpha)=\alpha 1_A$ для всех $\alpha \in \mathbbm{k}$.
 \end{example}
   \begin{example}\label{ExampleCoalgebra} Коалгебры $C$
 с коумножением $\Delta \colon C \to C  \otimes C$
 и коединицей $\varepsilon \colon C \to \mathbbm{k}$
 являются примерами $\lbrace \Delta, \varepsilon \rbrace$-алгебр, где $s(\Delta)=1$, $t(\Delta)=2$, $s(\varepsilon)=1$, $t(\varepsilon)=0$.
 \end{example}
    \begin{example} Заплетённые векторные пространства $W$
 с заплетением $$\tau \colon W \otimes W \to W \otimes W$$
  являются примерами $\lbrace \tau \rbrace$-алгебр, где $s(\tau)=2$, $t(\tau)=2$.
 \end{example}
   \begin{example}
   Биалгебры и алгебры Хопфа также являются $\Omega$-алгебрами для подходящих $\Omega$.
    \end{example}
    
    Для множества $\Omega$ определим множество $\Omega^* := \lbrace \omega^* \mid \omega \in \Omega \rbrace$, где $s(\omega^*):= t(\omega)$ и $t(\omega^*):= s(\omega)$. Очевидно, что если $A$~"--- конечномерная $\Omega$-алгебра, то $A^*$~"--- это $\Omega^*$-алгебра,
    где $\omega^*_{A^*}:= \left(\omega_A\right)^*$.

Если задана моноидальная категория\footnote{Определение см. в~\cite{EGNObook} или~\cite{MacLaneCatWork}.} $\mathcal C$ с моноидальным произведением $\otimes$
и нейтральным объектом $\mathbbm 1$,
то можно ввести понятие $\Omega$-алгебры в категории.
 Назовём \textit{$\Omega$-алгеброй} в  категории $\mathcal C$ объект $A$
 с морфизмами $\omega_A \colon A^{\otimes s(\omega)} \to
A^{\otimes t(\omega)}$, заданными для всех $\omega \in \Omega$. Здесь мы полагаем $A^{\otimes 0} := \mathbbm 1$.
 Очевидно, что $\Omega$-алгебра над полем $\mathbbm{k}$~"--- это в точности $\Omega$-алгебра в категории $\mathbf{Vect}_\mathbbm{k}$. 

         \section{Измерения}\label{SectionMeasurings}
         
         Наша конечная цель~"--- действия и кодействия биалгебр и алгебр Хопфа на $\Omega$-алгебрах.
         Однако, следуя М.~Свидлеру~\cite[глава VII]{Sweedler}, мы сперва введём более общее понятие
         измерения, которое позволяет сделать наши построения яснее с теоретико-категорной точки зрения. 

Пусть $A$ и $B$~"--- $\Omega$-алгебры над полем $\mathbbm{k}$ и пусть $P$~"---
$\mathbbm{k}$-коалгебра (в обычном смысле) с коумножением $\Delta_P$ и коединицей $\varepsilon_P$.
Если $\psi \colon P \otimes A \to B$~"--- линейное отображение,
то для любого $n\in\mathbb Z_+$ можно определить следующие линейные отображения $\psi_n \colon P \otimes A^{\otimes n} \to B^{\otimes n}$: 
$$\psi_0(p)=\varepsilon_P(p)$$ и $$\psi_n(p\otimes a_1 \otimes \ldots \otimes a_n):= \psi(p_{(1)} \otimes a_1) \otimes \ldots \otimes \psi(p_{(n)}\otimes a_n)\text{
для всех }a_i \in A,\ p\in P.$$
Здесь мы используем обозначения $$p_{(1)} \otimes p_{(2)}\otimes \ldots \otimes p_{(n)} := \Delta_n (p)\quad \text{(опущен знак суммы),}$$ 
где $\Delta_1 := \id_P$, $\Delta_{k+1}(p):=(\Delta_k  \otimes \id_P)\Delta_P(p)$
при $1\leqslant k \leqslant n$, $n\geqslant 1$.

 Будем говорить, что $\psi$~"--- \textit{измерение} из $A$ в $B$,
 если это отображение сохраняет операции из множества $\Omega$, т.е. для любого $\omega \in \Omega$
 справедливо равенство
$$\psi_{t(\omega)}(p \otimes \omega_A(a_1 \otimes \ldots \otimes a_{s(\omega)}))
=\omega_B(\psi_{s(\omega)}(p \otimes a_1 \otimes \ldots \otimes a_{s(\omega)}))$$
для всех $a_i \in A$, $p\in P$.


\begin{example} Если $\Omega = \lbrace \mu, u \rbrace$
и $A$, $B$~"--- алгебры с единицей, мы получаем классическое определение~\cite[глава VII]{Sweedler}
измерения из алгебры $A$ в алгебру $B$:

\hspace{1cm}$\psi(p \otimes ab)= \psi(p_{(1)} \otimes a)\psi(p_{(2)}\otimes b)$
и $\psi(p\otimes 1_A)= \varepsilon(p)1_B$ для всех $a,b \in A$ и $p\in P$.
\end{example}

Зафиксируем $\Omega$-алгебры $A$ и $B$ над полем $\mathbbm{k}$ и подпространство $V\subseteq \mathbf{Vect}_\mathbbm{k}(A,B)$. Рассмотрим категорию $\mathbf{Meas}(A,B,V)$, в которой объектами являются измерения $\psi \colon P \otimes A \to B$ для всевозможных $\mathbbm{k}$-коалгебр $P$, такие, что $\cosupp \psi \subseteq V$, а морфизмами из $\psi_1 \colon P_1 \otimes A \to B$
в $\psi_2 \colon P_2 \otimes A \to B$ являются гомоморфизмы коалгебр $\varphi \colon P_1 \to P_2$,
для которых коммутативна
диаграмма

$$\xymatrix{ P_1 \otimes A \ar[r]^(0.6){\psi_1} \ar[d]_{\varphi \otimes \id_A} & B \\
P_2 \otimes A  \ar[ru]_{\psi_2}  } $$

Будем называть коалгебру ${}_\square \mathcal{C}(A,B,V)$,
отвечающую  терминальному объекту категории $\mathbf{Meas}(A,B,V)$ (в теореме~\ref{TheoremUnivMeasExistence}  мы доказываем его существование), \textit{$V$-универсальной измеряющей коалгеброй} из $A$ в $B$.\label{NotationsquareCABV}

\begin{theorem}\label{TheoremUnivMeasExistence} В категории $\mathbf{Meas}(A,B,V)$ существует терминальный объект.
\end{theorem}
\begin{proof} Напомним, что через $K$ мы обозначаем правый сопряжённый функтор к забывающему функтору $\mathbf{Coalg}_\mathbbm{k} \to \mathbf{Vect}_\mathbbm{k}$ (см. пример~\ref{ExampleCofreeCoalg}).
Обозначим через ${}_\square \mathcal{C}(A,B,V)$ сумму всех подкоалгебр $C$ коалгебры $K(V)$, 
таких, что композиция вложения $C \subseteq K(V)$, коединицы $K(V)\to V$ сопряжения (см. \S\ref{SectionAdjoint})
и вложения $V \subseteq \mathbf{Vect}_\mathbbm{k}(A,B)$
задаёт измерение. Тогда, в силу универсального свойства коалгебры $K(V)$, линейное отображение $\psi_{A,B,V} \colon {}_\square \mathcal{C}(A,B,V) \otimes A \to B$,
отвечающее композиции вложения $ {}_\square\mathcal{C}(A,B,V) \subseteq K(V)$, коединицы $K(V)\to V$ сопряжения и вложения $V \subseteq \mathbf{Vect}_\mathbbm{k}(A,B)$, является терминальным объектом категории
$\mathbf{Meas}(A,B,V)$, и ${}_\square \mathcal{C}(A,B,V)$~"--- $V$-универсальная измеряющая коалгебра из $A$ в $B$. 
\end{proof}

\begin{remark} Легко видеть, что коноситель отображения $\psi_{A,B,V}$
является подпространством в $V$. Коноситель $\cosupp \psi_{A,B,V}$ совпадает с $V$,
если и только если существует измерение из $A$ в $B$ с коносителем, равным $V$. 
\end{remark}

\begin{example}
Когда $A$ и $B$~"--- ассоциативные алгебры с единицей над полем $\mathbbm{k}$, коалгебра $ {}_\square\mathcal{C}(A,B,\mathbf{Vect}_\mathbbm{k}(A,B))$~"--- это в точности универсальная измеряющая коалгебра $M(A,B)$ М.~Свидлера~\cite[глава VII]{Sweedler}.
\end{example}

         \section{Коизмерения}\label{SectionComeasurings}

Коизмерения, которые играют для кодействий ту же роль, что измерения играют для действий,
вводятся двойственным образом.

Пусть $A$ и $B$~"--- $\Omega$-алгебры и пусть $Q$~"--- ассоциативная алгебра с единицей над полем $\mathbbm{k}$.
Для всякого линейного отображения $\rho \colon A \to B \otimes Q$ и
числа $n\in\mathbb Z_+$ определим линейные отображения $\rho_n \colon A^{\otimes n} \to B^{\otimes n} \otimes Q$
следующим образом: $\rho_0(\alpha)=\alpha 1_Q$ для всех $\alpha \in \mathbbm{k}$ и
 $$\rho_n(a_1 \otimes \ldots \otimes a_n):= (a_1)_{(0)} \otimes \ldots \otimes (a_{n})_{(0)}
\otimes (a_1)_{(1)}  \ldots  (a_{n})_{(1)}$$
для всех $a_i \in A$.

Назовём линейное отображение $\rho \colon A \to B \otimes Q$ \textit{коизмерением} из $A$ в $B$,
если оно сохраняет операции из множества $\Omega$, т.е. если для всех  $\omega \in \Omega$ и
$a_i \in A$ справедливо равенство
$$\rho_{t(\omega)} \left( \omega_A(a_1 \otimes \ldots \otimes a_{s(\omega)}) \right)
= (\omega_B \otimes \id_Q)\rho_{s(\omega)}(a_1 \otimes \ldots \otimes a_{s(\omega)}).$$

 \begin{example} В случае, когда $\Omega = \lbrace \mu, u \rbrace$ и $A$, $B$~"---
 алгебры с единицей, определение коизмерения из $A$ в $B$ выглядит следующим образом:
 
\hspace{1.7cm}$\rho(1_A)= 1_B \otimes 1_Q$ и $\rho(ab)= a_{(0)}b_{(0)}\otimes a_{(1)}b_{(1)}$  для всех $a,b \in A$.
\end{example}

\begin{remark} Ограничение отображения $\hat\rho$
на подпространство $Q^\circ \otimes A$ является измерением.
В частности, если алгебра $Q$ конечномерна, то отображение $\hat\rho$ само является измерением, причём $\cosupp \rho = \cosupp \hat\rho$.
\end{remark}

Зафиксируем $\Omega$-алгебры $A$ и $B$ над полем $\mathbbm{k}$ и подпространство $V\subseteq \mathbf{Vect}_\mathbbm{k}(A,B)$. Рассмотрим категорию $\mathbf{Comeas}(A,B,V)$, объектами которой являются такие коизмерения $\rho \colon A \to B \otimes Q$ для произвольных ассоциативных $\mathbbm{k}$-алгебр $Q$ с единицей, что выполнено условие $\cosupp \rho \subseteq V$, а морфизмами из $\rho_1 \colon A \to B \otimes Q_1$
в $\rho_2 \colon A \to B \otimes Q_2$  являются такие гомоморфизмы $\varphi \colon Q_1 \to Q_2$
алгебр с единицей, что коммутативна диаграмма

$$\xymatrix{ A \ar[r]^(0.4){\rho_1} 
\ar[rd]_{\rho_2}
 & B \otimes Q_1 \ar[d]^{\id_B \otimes \varphi} \\
& B \otimes Q_2} $$

Назовём алгебру $\mathcal{A}^\square(A,B,V)$,
соответствующую начальному объекту в категории $\mathbf{Comeas}(A,B,V)$ (если он существует), \textit{$V$-универсальной коизмеряющей алгеброй} из $A$ в $B$.\label{NotationAsquareABV}

\begin{theorem}\label{TheoremUnivComeasExistence}  Если $V$ 
поточечно конечномерно и замкнуто в конечной топологии, то в категории $\mathbf{Comeas}(A,B,V)$
существует начальный объект.
\end{theorem}
\begin{proof} В силу теоремы~\ref{TheoremVectDualization}
существуют векторное пространство~$W$
 и такое линейное отображение $\rho_0 \colon A \to B \otimes W$, что $\hat\rho_0(W^* \otimes (-)) = V$.

Пусть $(a_\alpha)_\alpha$~"--- базис $\Omega$-алгебры $A$ и пусть $(b_\beta)_\beta$~"--- базис $\Omega$-алгебры $B$.
Определим элементы $p_{\beta\alpha} \in W$ при помощи равенств $$\rho_0(a_\alpha)=\sum_{\beta} b_\beta \otimes p_{\beta\alpha}.$$
Пусть $W_0 := \supp \rho_0=\langle p_{\beta\alpha}  \mid \alpha, \beta \rangle_\mathbbm{k}$.

Для каждого $\omega \in \Omega$ определим элементы $a^{\beta_1,\ldots,\beta_{t(\omega)}}_{\omega; \alpha_1,\ldots,\alpha_{s(\omega)}}, b^{\beta_1,\ldots,\beta_{t(\omega)}}_{\omega; \alpha_1,\ldots,\alpha_{s(\omega)}}\in \mathbbm{k}$
при помощи равенств
$$\omega_A(a_{\alpha_1} \otimes \ldots \otimes a_{\alpha_{s(\omega)}})
= \sum_{\beta_1,\ldots,\beta_{t(\omega)}} a^{\beta_1,\ldots,\beta_{t(\omega)}}_{\omega; \alpha_1,\ldots,\alpha_{s(\omega)}}a_{\beta_1} \otimes \ldots \otimes a_{\beta_{t(\omega)}}$$
и
$$\omega_B(b_{\alpha_1} \otimes \ldots \otimes b_{\alpha_{s(\omega)}})
= \sum_{\beta_1,\ldots,\beta_{t(\omega)}} b^{\beta_1,\ldots,\beta_{t(\omega)}}_{\omega; \alpha_1,\ldots,\alpha_{s(\omega)}}b_{\beta_1} \otimes \ldots \otimes b_{\beta_{t(\omega)}}.$$

Другими словами, $a^{\beta_1,\ldots,\beta_{t(\omega)}}_{\omega; \alpha_1,\ldots,\alpha_{s(\omega)}}$
и $b^{\beta_1,\ldots,\beta_{t(\omega)}}_{\omega; \alpha_1,\ldots,\alpha_{s(\omega)}}$
являются структурными константами, соответственно, $\Omega$-алгебр $A$ и $B$.

Напомним, что через $T$ мы обозначаем левый сопряжённый функтор к забывающему функтору $\mathbf{Alg}_\mathbbm{k} \to \mathbf{Vect}_\mathbbm{k}$ (см. пример~\ref{ExampleTensorAlg}).
Пусть $i_{W_0} \colon W_0 \to T(W_0)$~"--- единица этого сопряжения (см. \S\ref{SectionAdjoint}).
 Обозначим через $I$ идеал алгебры $T(W_0)$, порождённый элементами
\begin{equation*}\begin{split} \sum_{\gamma_1,\ldots,\gamma_{t(\omega)}} a^{\gamma_1,\ldots,\gamma_{t(\omega)}}_{\omega; \alpha_1,\ldots,\alpha_{s(\omega)}}
i_{W_0}(p_{\beta_1 \gamma_1})i_{W_0}(p_{\beta_2 \gamma_2})
\ldots i_{W_0}(p_{\beta_{t(\omega)} \gamma_{t(\omega)}}) - \\
- \sum_{\mu_1,\ldots,\mu_{s(\omega)}}
b^{\beta_1,\ldots,\beta_{t(\omega)}}_{\omega; \mu_1,\ldots,\mu_{s(\omega)}}
i_{W_0}(p_{\mu_1 \alpha_1})i_{W_0}(p_{\mu_2 \alpha_2})\ldots i_{W_0}(p_{\mu_{s(\omega)} \alpha_{s(\omega)}})\end{split}\end{equation*}
 для всех $\omega \in \Omega$ и $\alpha_1,\ldots,\alpha_{s(\omega)}$, $\beta_1,\ldots,\beta_{t(\omega)}$.
 (Произведение нулевого числа множителей $i_{W_0}(\ldots)$ считаем
 по определению равным $1_{T(W_0)}$.)

Пусть $ \mathcal{A}^\square(A,B,V):=T(W_0)/I$. Определим
$\rho_{A,B,V} \colon A \to B
\otimes \mathcal{A}^\square(A,B,V)$
как композицию отображения $\rho_0$, область значений которого ограничена до $B \otimes W_0$ и
отображений $i_{W_0} \colon W_0 \to T(W_0)$ и $T(W_0)\twoheadrightarrow T(W_0)/I$,
тензорно умноженных на $\id_B$.
В силу выбора идеала $I$ линейное отображение $\rho_{A,B,V}$ является коизмерением.
Пусть $\rho \colon A \to B \otimes Q$~"--- другое такое коизмерение, что $\cosupp \rho \subseteq V$,
а $Q$~"--- некоторая ассоциативная алгебра с единицей. 
Тогда в силу следствия~\ref{CorollaryPrecLinearMapComodCriterion} 
существует такое линейное отображение $\tau \colon W_0 \to Q$,
что коммутативна диаграмма
\begin{equation*} \xymatrix{ A \ar[r] 
\ar[rd]_\rho
 & B \otimes W_0 \ar@{-->}[d]^{\id_B \otimes \tau} \\
& B \otimes Q} \end{equation*}

Поскольку $Q$~"--- ассоциативная алгебра с единицей,
существует такой гомоморфизм алгебр с единицей
$\tau_1 \colon T(W_0) \to Q$,
что коммутативна диаграмма
\begin{equation*} \xymatrix{ A \ar[r] 
\ar[rd]_\rho
 & B \otimes T(W_0) \ar@{-->}[d]^{\id_B \otimes \tau_1} \\
& B \otimes Q} \end{equation*}

Поскольку $\rho$~"--- коизмерение, справедливо включение $I \subseteq \ker(\tau_1)$, 
и существует такой гомоморфизм алгебр
$\bar\tau_1 \colon \mathcal{A}^\square(A,B,V) \to Q$,
что коммутативна диаграмма
\begin{equation*} \xymatrix{ A \ar[rr]^(0.3){\rho_{A,B,V}}
\ar[rrd]_\rho
 & & B \otimes \mathcal{A}^\square(A,B,V) \ar@{-->}[d]^{\id_B \otimes \bar\tau_1} \\
& & B \otimes Q} \end{equation*}

В силу того, что алгебра $\mathcal{A}^\square(A,B,V)$ порождена образами элементов $p_{\beta\alpha}$, гомоморфизм $\bar \tau_1$~"--- единственный, и $\rho_{A,B,V}$~"--- начальный объект категории $\mathbf{Comeas}(A,B,V)$.
\end{proof}

\begin{remark} Легко видеть, что коноситель отображения $$\rho_{A,B,V} \colon A \to B
\otimes \mathcal{A}^\square(A,B,V)$$
является подпространством пространства $V$ (коноситель может стать собственным подпространством пространства $V$ после факторизации алгебры $T(W_0)$ по идеалу $I$) и совпадает с $V$, если и только если существует
коизмерение из $A$ в $B$ с коносителем, равным $V$. 
\end{remark}

\begin{example}
Если $A$ и $B$~"--- ассоциативные алгебры с единицей, причём алгебра $B$ конечномерна,
то всё пространство $\mathbf{Vect}_\mathbbm{k}(A,B)$, очевидно, поточечно конечномерно и замкнуто в конечной топологии. Тогда $\mathcal{A}^\square(A,B,\mathbf{Vect}_\mathbbm{k}(A,B))$ 
совпадает с алгеброй $\alpha(B,A)$ Д.~Тамбары~\cite[\S 1]{Tambara}. 
\end{example}

         \section{Двойственность между измерениями и коизмерениями} \label{SectionDuality(Co)measurings}

В этом параграфе мы устанавливаем двойственность между коалгеброй
${}_\square \mathcal{C}(A,B,V)$ и алгеброй
$\mathcal{A}^\square(A,B,V)$. 

Сперва докажем следующую лемму:

\begin{lemma}\label{LemmaMeasComeasBijection}
Пусть $A$ и $B$~"--- $\Omega$-алгебры над полем $\mathbbm{k}$, $P$~"--- коалгебра, а $V \subseteq \mathbf{Vect}_\mathbbm{k}(A,B)$~"---
поточечно конечномерное подпространство, замкнутое в конечной
топологии.
Обозначим через $\mathrm{Meas}(P,A,B,V)$
множество всех измерений $\psi \colon P \otimes A \to B$, для которых $\cosupp \psi \subseteq V$,
а через $\mathrm{Comeas}(P^*,A,B,V)$
множество всех коизмерений $\rho \colon A \to B \otimes P^*$,  для которых $\cosupp \rho \subseteq V$.
Тогда существует биекция $\mathrm{Meas}(P,A,B,V) 
\cong \mathrm{Comeas}(P^*,A,B,V)$,
естественная по коалгебре $P$, если рассматривать $\mathrm{Meas}(-,A,B,V)$ и $\mathrm{Comeas}((-)^*,A,B,V)$
как функторы $\mathbf{Coalg}_\mathbbm{k} \to \mathbf{Sets}$.
\end{lemma}
\begin{proof} Любое коизмерение $\rho \in \mathrm{Comeas}(P^*,A,B,V)$ порождает измерение $\rho^\vee \colon P \otimes A \to B$, где $\rho^\vee(p\otimes a):=a_{(1)}(p)a_{(0)}$ для всех $a\in A$ и $p\in P$.
Тогда $\cosupp \rho^\vee \subseteq \cosupp \rho \subseteq V$ и $\rho^\vee \in \mathrm{Meas}(P,A,B,V)$.

Обратно, если $\psi \in \mathrm{Meas}(P,A,B,V)$, то $\cosupp \psi \subseteq V$,
т.е. $\cosupp \psi$ является поточечно конечномерным пространством.
Отсюда для каждого $a \in A$ существует конечное число элементов
$b_1, \ldots, b_n \in B$ и линейных функций $p_1^*, \ldots, p_n^* \in P^*$, таких, что
для любого $p\in P$ справедливо равенство $\psi(p\otimes a)=\sum_{i=1}^n p_i^*(p)b_i$,
причём при фиксированных линейно независимых элементах $b_1, \ldots, b_n$ элементы $p_1^*, \ldots, p_n^* \in P^*$ определены однозначно. Поэтому существует единственное коизмерение $\rho \colon A \to B \otimes P^*$, такое, что $\psi = \rho^\vee$. Для любого $a\in A$ справедливо равенство $\rho(a)=\sum_{i=1}^n b_i \otimes p_i^*$. 
Поскольку $V$ замкнуто, в силу леммы~\ref{LemmaLinearMap(Co)modClosure}
справедливо включение $\cosupp \rho = \overline{\cosupp \psi} \subseteq V$, и $\rho \in \mathrm{Comeas}(P^*,A,B,V)$.

\end{proof}

Напомним, что для алгебры $A$ коалгебра, конечная двойственная к ней (см. с.~\pageref{DefACirc}), обозначается через $A^\circ$.

\begin{theorem}\label{TheoremUniv(Co)measDuality} Пусть $A$ и $B$~"--- $\Omega$-алгебры над полем $\mathbbm{k}$, а $V \subseteq \mathbf{Vect}_\mathbbm{k}(A,B)$~"---
поточечно конечномерное подпространство, замкнутое в конечной
топологии. Обозначим ограничение отображения $\hat \rho_{A,B,V} \colon 
\mathcal{A}^\square(A,B,V)^* \otimes A \to B$ на подпространство $\mathcal{A}^\square(A,B,V)^\circ \otimes A$ также через $\hat \rho_{A,B,V}$, а через $\theta$ обозначим  единственный такой гомоморфизм коалгебр,
что диаграмма
$$\xymatrix{ {}_\square \mathcal{C}(A,B,V) \otimes A \ar[dr]^(0.6){\psi_{A,B,V}}& \\
                   &   B \\
\mathcal{A}^\square(A,B,V)^\circ \otimes A \ar[ru]_(0.6){\hat\rho_{A,B,V}} \ar@{-->}[uu]^{\theta \otimes \id_A} } $$
коммутативна. Тогда $\theta$~"--- изоморфизм коалгебр.
\end{theorem}
\begin{proof}
Отталкиваясь от замечания~1.3 из~\cite{Tambara},
мы можем добавить
биекцию из леммы~\ref{LemmaMeasComeasBijection}
к следующим биекциям, естественным по коалгебре~$P$:
\begin{equation}\label{EqTambaraBijections}\begin{split}\mathbf{Coalg}_\mathbbm{k}(P, {}_\square \mathcal{C}(A,B,V))\cong \mathrm{Meas}(P,A,B,V) \\ \cong \mathrm{Comeas}(P^*,A,B,V)\cong
\mathbf{Alg}_\mathbbm{k}(\mathcal{A}^\square(A,B,V), P^*)
\cong \mathbf{Coalg}_\mathbbm{k}(P, \mathcal{A}^\square(A,B,V)^\circ).\end{split}
\end{equation}

Заменим теперь $P$ коалгеброй $\mathcal{A}^\square(A,B,V)^\circ$.
Тогда гомоморфизм коалгебр $\theta \colon \mathcal{A}^\square(A,B,V)^\circ \to {}_\square \mathcal{C}(A,B,V)$ будет соответствовать отображению $\id_{\mathcal{A}^\square(A,B,V)^\circ}$.
Если мы заменим $P$ коалгеброй ${}_\square \mathcal{C}(A,B,V)$,
то в силу естественности биекций, гомоморфизм $\id_{{}_\square \mathcal{C}(A,B,V)}$ 
будет соответствовать отображению $\theta^{-1}$. Следовательно, $\theta$~"--- изоморфизм коалгебр.
\end{proof}

\begin{remark} Можно сформулировать и вариант двойственности в виде сопряжения. Будем называть измерение $\psi \colon P \otimes A \to B$ \textit{рациональным},
если $\psi = \rho_0^\vee$ для некоторого линейного отображения $\rho_0 \colon A \to B \otimes P^*$. (См. лемму~\ref{LemmaMeasComeasBijection} выше.)
Введём обозначение $\tilde\psi := \rho_0$.
Отметим, что отображение $\tilde\psi$ определено однозначно и является измерением.
Обозначим через $\mathbf{Meas}^{\mathrm{rat}}(A,B,V)$
полную подкатегорию категории $\mathbf{Meas}(A,B,V)$,
объектами которой являются все рациональные измерения $\psi \colon P \otimes A \to B$.
Тогда для любого подпространства $V \subseteq \mathbf{Vect}_\mathbbm{k}(A,B)$, замкнутого
в конечной топологии (замкнутость требуется для того, чтобы можно было применить лемму~\ref{LemmaLinearMap(Co)modClosure}), любого рационального измерения $\psi \colon P \otimes A \to B$ 
и любого коизмерения $\rho \colon A \to B \otimes Q$, таких, что $\cosupp \psi, \cosupp\rho \subseteq V$,
существует естественная биекция $$\mathbf{Meas}^{\mathrm{rat}}(A,B,V)(\psi, \hat\rho\bigr|_{Q^\circ \otimes A}) \cong \mathbf{Comeas}(A,B,V)(\rho, \tilde\psi)$$
(ограничение биекции $\mathbf{Coalg}_\mathbbm{k}(P,Q^\circ)\cong \mathbf{Alg}_\mathbbm{k}(Q,P^*)$),
что означает наличие сопряжения между соответствующими контравариантными функторами.
Из этого сопряжения получается, в частности, другое доказательство теоремы~\ref{TheoremUniv(Co)measDuality}.
(Достаточно использовать контравариантную версию того факта, что правый сопряжённый ковариантный функтор сохраняет терминальный объект.)
\end{remark}

          \section{Действия}\label{SectionActions}

Пусть $B$~"---  биалгебра, а $A$~"--- $\Omega$-алгебра над полем $\mathbbm{k}$.
Будем говорить, что линейное отображение $\psi \colon B \otimes A \to A$
является \textit{$B$-действием} на $A$ или, что то же самое, $\psi$ задаёт на $A$ структуру \textit{(левой) $B$-модульной $\Omega$-алгебры}, если выполняются два следующих условия: \begin{enumerate}
\item $\psi$ определяет на $A$ структуру левого (унитального) $B$-модуля;
\item $\psi$ является измерением из $A$ в $A$.
\end{enumerate}

Иными словами, левая $B$-модульная $\Omega$-алгебра~"--- это не что иное, как $\Omega$-алгебра
в моноидальной категории левых $B$-модулей.

Если $\Omega=\lbrace\mu, u\rbrace$ и $A$~"--- алгебра с единицей, то мы получаем обычное определение модульной алгебры с единицей (см.~\S\ref{Section(Co)modules}).
Если $\Omega=\lbrace\mu\rbrace$, то алгебра $A$ является модульной алгеброй, необязательно с единицей.
Если, наконец, $\Omega=\lbrace\Delta, \varepsilon\rbrace$ и $A$~"--- коалгебра,
то мы получаем определение модульной коалгебры.

Для действий мы будем использовать понятия коносителя и эквивалентности, а также предпорядок,
введённые в~\S\ref{SectionLinearMaps}.

Заметим, что для действия $\psi \colon B \otimes A \to A$ 
его коноситель $\cosupp\psi$~"--- это в точности образ соответствующего гомоморфизма алгебр $B \to \End_\mathbbm{k}(A)$. 
В частности, $\cosupp\psi$~"--- это подалгебра в $\End_\mathbbm{k}(A)$, содержащая тождественный оператор.

\begin{theorem}\label{TheoremsquareBBialgebra}
Пусть $A$~"--- это $\Omega$-алгебра над полем $\mathbbm{k}$, а $V\subseteq \End_\mathbbm{k}(A)$~"--- подалгебра, содержащая тождественный оператор.
Тогда на $V$-универсальной измеряющей коалгебре ${}_\square \mathcal{B}(A,V) := {}_\square \mathcal{C}(A,A,V)$ 
можно задать структуру биалгебры таким образом, что
для любой биалгебры $B$ и любого действия
$\psi \colon B \otimes A \to A$, такого, что $\cosupp \psi \subseteq V$, единственный гомоморфизм коалгебр
$\varphi$, такой, что диаграмма ниже является коммутативной,
является ещё и гомоморфизмом биалгебр:
\begin{equation}\label{EqTheoremsquareBBialgebra}\xymatrix{ B \otimes A \ar[r]^(0.6){\psi} \ar@{-->}[d]_{\varphi \otimes \id_A} & A \\
{}_\square \mathcal{B}(A,V) \otimes A  \ar[ru]_{\psi_{A,V}}  } \end{equation}
(Здесь $\psi_{A,V} := \psi_{A,A,V}$.)
\end{theorem}

Основываясь на теореме~\ref{TheoremsquareBBialgebra}, будем называть биалгебру ${}_\square \mathcal{B}(A,V)$  \textit{$V$-универсальной действующей биалгеброй}.\label{NotationsquareBAV}

\begin{proof}[Доказательство теоремы~\ref{TheoremsquareBBialgebra}]
Заметим, что $$\psi_{A,V}(\id_{{}_\square \mathcal{B}(A,V)} \otimes \psi_{A,V})
\colon {}_\square \mathcal{B}(A,V)\otimes {}_\square \mathcal{B}(A,V) \otimes A \to A$$
также является измерением, причём $\cosupp(\psi_{A,V}(\id_{{}_\square \mathcal{B}(A,V)} \otimes \psi_{A,V}))
\subseteq V$, поскольку $V$~"--- подалгебра.
Отсюда существует единственный гомоморфизм коалгебр $\mu \colon 
{}_\square \mathcal{B}(A,V)\otimes {}_\square \mathcal{B}(A,V) \to {}_\square \mathcal{B}(A,V)$,
такой, что следующая диаграмма коммутативна:
$$\xymatrix{ {}_\square \mathcal{B}(A,V)\otimes {}_\square \mathcal{B}(A,V) \otimes A \ar[rrr]^(0.6){\id_{{}_\square \mathcal{B}(A,V)} \otimes \psi_{A,V}} \ar@{-->}[d]_{\mu \otimes \id_A} && & {}_\square \mathcal{B}(A,V) \otimes A  \ar[d]^{\psi_{A,V}} \\
{}_\square \mathcal{B}(A,V) \otimes A  \ar[rrr]_{\psi_{A,V}} &&  & A  } $$
Определим теперь умножение в ${}_\square \mathcal{B}(A,V)$ при помощи отображения  $\mu$.

Для того, чтобы задать единичный элемент, рассмотрим тривиальное отображение $$\mathbbm{k} \otimes A
\mathrel{\widetilde\to}  A.$$
В силу того, что это отображение является измерением,
существует единственный гомоморфизм коалгебр $u \colon \mathbbm{k} \to {}_\square \mathcal{B}(A,V)$,
такой, что коммутативна диаграмма
$$\xymatrix{ \mathbbm{k} \otimes A \ar[r]^(0.6){\sim} \ar@{-->}[d]_{u \otimes \id_A} & A \\
{}_\square \mathcal{B}(A,V) \otimes A  \ar[ru]_{\psi_{A,V}}  } $$
Определим $1_{{}_\square \mathcal{B}(A,V)} := u(1_\mathbbm{k})$.

Теперь нужно доказать, что умножение $\mu$
ассоциативно и $1_{{}_\square \mathcal{B}(A,V)}$ действительно является
единицей алгебры ${}_\square \mathcal{B}(A,V)$.

Рассмотрим диаграмму

{\scriptsize
$$\xymatrix{  {}_\square \mathcal{B}(A,V)\otimes {}_\square \mathcal{B}(A,V) \otimes {}_\square\mathcal{B}(A,V) \otimes A 
\ar[rr]^{ \id_{{}_\square \mathcal{B}(A,V)} \otimes \mu\otimes \id_A} \ar[dd]^(0.7){\id_{{}_\square \mathcal{B}(A,V)} \otimes\id_{{}_\square \mathcal{B}(A,V)} \otimes\psi_{A,V}} \ar[rd]^(0.6){\qquad\mu\otimes \id_{{}_\square \mathcal{B}(A,V)} \otimes \id_A} &&  {}_\square \mathcal{B}(A,V) \otimes {}_\square\mathcal{B}(A,V) \otimes A \ar'[d][dd]^{\id_{{}_\square \mathcal{B}(A,V)} \otimes\psi_{A,V}} \ar[rd]^{\mu\otimes \id_A} & \\
& {}_\square \mathcal{B}(A,V) \otimes {}_\square\mathcal{B}(A,V) \otimes A \ar[rr]^(0.4){\mu\otimes \id_A} \ar[dd]^(0.7){\id_{{}_\square \mathcal{B}(A,V)} \otimes\psi_{A,V}}
&& {}_\square\mathcal{B}(A,V) \otimes A \ar[dd]^{\psi_{A,V}} \\
{}_\square \mathcal{B}(A,V) \otimes {}_\square\mathcal{B}(A,V) \otimes A \ar'[r][rr]^(0.4){\id_{{}_\square \mathcal{B}(A,V)} \otimes\psi_{A,V}} \ar[rd]^{\mu\otimes \id_A} && {}_\square\mathcal{B}(A,V) \otimes A \ar[rd]^{\psi_{A,V}} & \\
& {}_\square\mathcal{B}(A,V) \otimes A \ar[rr]_{\psi_{A,V}} && A } $$ }

Левая грань коммутативна, поскольку обе композиции равны 
$\mu \otimes \psi_{A,V}$.
Коммутативность остальных граней, за исключением верхней, следует из
определения отображения $\mu$.
Отсюда композиции, отвечающие верхней грани, также становятся равными друг другу
после их композиции с $\psi_{A,V}$.
Теперь из универсального свойства коалгебры ${}_\square \mathcal{B}(A,V)$
следует коммутативность диаграммы
 $$\xymatrix{ {}_\square \mathcal{B}(A,V)\otimes {}_\square \mathcal{B}(A,V) \otimes {}_\square\mathcal{B}(A,V) \ar[rrr]^(0.6){\id_{{}_\square \mathcal{B}(A,V)} \otimes \mu} \ar[d]_{\mu \otimes \id_{{}_\square \mathcal{B}(A,V)}} && & {}_\square \mathcal{B}(A,V) \otimes {}_\square\mathcal{B}(A,V)  \ar[d]^{\mu} \\
{}_\square \mathcal{B}(A,V) \otimes {}_\square \mathcal{B}(A,V)  \ar[rrr]_{\mu} &&  & {}_\square\mathcal{B}(A,V)  },$$
откуда умножение в ${}_\square \mathcal{B}(A,V)$ действительно ассоциативно.

Рассмотрим диаграмму

$$\xymatrix{ \mathbbm{k} \otimes {}_\square \mathcal{B}(A,V)\otimes A \ar[rr]^{u \otimes \id_{{}_\square \mathcal{B}(A,V)}\otimes\id_A} \ar[rd]^{\sim} \ar[dd]_(0.7){\id_\mathbbm{k}  \otimes\psi_{A,V}}
  & & {}_\square \mathcal{B}(A,V)\otimes {}_\square \mathcal{B}(A,V)\otimes A  \ar[dd]^{\id_{{}_\square \mathcal{B}(A,V)} \otimes \psi_{A,V}}  \ar[ld]^{\mu\otimes\id_A}\\
 & {}_\square \mathcal{B}(A,V)\otimes A \ar[dd]^(0.3){\psi_{A,V}} & \\
\mathbbm{k} \otimes A \ar'[r][rr]^(0.3){u\otimes\id_A} \ar[rd]^{\sim}
  & & {}_\square \mathcal{B}(A,V)\otimes  A  \ar[ld]^{\psi_{A,V}}\\
 &  A  & 
 } $$

Коммутативность задней и левой граней очевидна.
Коммутативность нижней грани следует из определения отображения $u$, а коммутативность правой грани следует из определения отображения $\mu$.
Следовательно, гомоморфизмы коалгебр, образующие верхнюю грань, становятся равными после их композиции с
 $\psi_{A,V}$. 
 Применяя универсальное свойство коалгебры ${}_\square \mathcal{B}(A,V)$,
 получаем отсюда
коммутативность диаграммы
 $$\xymatrix{ \mathbbm{k} \otimes {}_\square \mathcal{B}(A,V) \ar[rr]^{u \otimes \id_{{}_\square \mathcal{B}(A,V)}} \ar[rd]^{\sim} 
  & & {}_\square \mathcal{B}(A,V)\otimes {}_\square \mathcal{B}(A,V) \ar[ld]^{\mu}\\
 & {}_\square \mathcal{B}(A,V)  & }$$
 Отсюда элемент $1_{{}_\square \mathcal{B}(A,V)}$ является левой единицей алгебры ${}_\square \mathcal{B}(A,V)$.
 
 Коммутативность диаграммы
 $$\xymatrix{  {}_\square \mathcal{B}(A,V) \otimes \mathbbm{k} \ar[rr]^{\id_{{}_\square \mathcal{B}(A,V)} \otimes u} \ar[rd]^{\sim} 
  & & {}_\square \mathcal{B}(A,V)\otimes {}_\square \mathcal{B}(A,V) \ar[ld]^{\mu}\\
 & {}_\square \mathcal{B}(A,V)  & }$$
 получается при помощи аналогичных рассуждений, применённых к диаграмме
  $$\xymatrix{  {}_\square \mathcal{B}(A,V)\otimes \mathbbm{k} \otimes A \ar[rr]^{ \id_{{}_\square \mathcal{B}(A,V)}\otimes u \otimes\id_A} \ar[rd]^{\sim} \ar[dd]^{\sim}
  & & {}_\square \mathcal{B}(A,V)\otimes {}_\square \mathcal{B}(A,V)\otimes A  \ar[dd]^{\id_{{}_\square \mathcal{B}(A,V)} \otimes \psi_{A,V}}  \ar[ld]^{\mu\otimes\id_A}\\
 & {}_\square \mathcal{B}(A,V)\otimes A \ar[dd]^(0.3){\psi_{A,V}} & \\
{}_\square \mathcal{B}(A,V) \otimes A \ar@{=}'[r][rr] \ar[rd]_{\psi_{A,V}}
  & & {}_\square \mathcal{B}(A,V)\otimes  A  \ar[ld]^{\psi_{A,V}}\\
 &  A  & 
 } $$
 
 Следовательно, $1_{{}_\square \mathcal{B}(A,V)}$ является не только левой, но и правой единицей
 алгебры ${}_\square \mathcal{B}(A,V)$, и ${}_\square \mathcal{B}(A,V)$~"--- действительно биалгебра.
 
 Предположим, что $B$~"--- биалгебра, а $\psi \colon B \otimes A \to A$~"--- такое
 действие, что $\cosupp \psi \subseteq V$.
 Обозначим через $\varphi \colon B \to {}_\square \mathcal{B}(A,V)$ 
единственный гомоморфизм коалгебр, делающий диаграмму~\eqref{EqTheoremsquareBBialgebra}
коммутативной. Такой гомоморфизм существует в силу теоремы~\ref{TheoremUnivMeasExistence}. Докажем, что в этом случае $\varphi$ является гомоморфизмом биалгебр.

Рассмотрим диаграмму
$$\xymatrix{  B \otimes B \otimes A 
\ar[rr]^(0.4){\varphi\otimes\varphi\otimes\id_A} \ar[dd]_{\id_B \otimes \psi}  \ar[rd]^{\mu_B\otimes\id_A} &&  {}_\square \mathcal{B}(A,V) \otimes {}_\square\mathcal{B}(A,V) \otimes A \ar'[d][dd]^{\id_{{}_\square \mathcal{B}(A,V)}\otimes \psi_{A,V}} \ar[rd]^{\mu\otimes\id_A} & \\
& B \otimes A \ar[rr]^(0.3){\varphi\otimes\id_A} \ar[dd]^(0.3){\psi} && {}_\square\mathcal{B}(A,V) \otimes A \ar[dd]^{\psi_{A,V}} \\
B \otimes A \ar'[r][rr]^(0.3){\varphi\otimes\id_A} \ar[rd]^{\psi} && {}_\square\mathcal{B}(A,V) \otimes A  \ar[rd]^{\psi_{A,V}} \\
& A \ar@{=}[rr] && A } $$ 
 (Здесь $\mu_B \colon B\otimes B \to B$~"--- умножение в $B$.)
 
 Нижняя и передняя грани коммутативны по определению гомоморфизма $\varphi$. 
Задняя грань коммутативна, так как она совпадает с нижней (и передней) гранью, тензорно умноженной
слева на $\varphi$.
Левая и правая грани коммутативны, поскольку и $\psi$, и $\psi_{A,V}$ определяют на $A$
структуру левого модуля. 
  Отсюда после композиции с $\psi_{A,V}$
верхняя грань также становится коммутативной, и из
универсального свойства коалгебры ${}_\square\mathcal{B}(A,V)$
следует коммутативность диаграммы
 $$\xymatrix{ B \otimes B   \ar[r]^(0.3){\varphi\otimes\varphi}\ar[d]_{\mu_B} & {}_\square \mathcal{B}(A,V) \otimes {}_\square\mathcal{B}(A,V)  \ar[d]^{\mu} \\ 
 B \ar[r]_{\varphi} & {}_\square\mathcal{B}(A,V) }$$ Следовательно, отображение $\varphi$ 
сохраняет операцию умножения.
 
Рассмотрим диаграмму $$\xymatrix{ & \mathbbm{k} \otimes A \ar[lddd]_{u_B \otimes \id_A} 
 \ar[dd]^(0.7){u \otimes \id_A}  \ar[rddd]^{\sim} & \\ & & \\
 & {}_\square\mathcal{B}(A,V) \otimes A  
 \ar[rd]_{\psi_{A,V}} & \\
 B \otimes A \ar[rr]_{\psi} \ar[ru]_{\varphi \otimes \id_A}& & A }$$ (Здесь отображение $u_B \colon \mathbbm{k}  \to B$ задаётся равенством $u(\alpha)=\alpha 1_B$, где $\alpha\in \mathbbm{k}$.)
 
 Большой треугольник и правый треугольник коммутативны в силу того, что единицы биалгебр $B$ и ${}_\square\mathcal{B}(A,V)$ действуют на $A$ как тождественный оператор.
Нижний треугольник коммутативен по определению отображения $\varphi$.
Следовательно, левый треугольник также становится коммутативным после композиции его отображений
с отображением $\psi_{A,V}$.
Следовательно, из универсального свойства коалгебры
 ${}_\square\mathcal{B}(A,V)$ 
следует коммутативность диаграммы
 $$\xymatrix{  & \mathbbm{k} \ar[ld]_{u_B} \ar[rd]^(0.4){u} & \\
 B \ar[rr]_(0.4){\varphi} & & {}_\square\mathcal{B}(A,V) }$$
 Отсюда отображение $\varphi$ сохраняет единицу и поэтому действительно является
 гомоморфизмом биалгебр.
\end{proof}

Будем называть биалгебру ${}_\square \mathcal{B}(A,\End_\mathbbm{k}(A))$ \textit{универсальной действующей биалгеброй} на $\Omega$-алгебре $A$ и обозначать её через ${}_\square \mathcal{B}(A)$.\label{NotationsquareBA}

\begin{remark}
Если $A$~"--- ассоциативная алгебра с единицей над полем $\mathbbm{k}$, то биалгебра ${}_\square \mathcal{B}(A)$~"--- это в точности универсальная измеряющая биалгебра $M(A,A)$ М.~Свидлера~\cite[глава~7]{Sweedler}.
\end{remark}

Обратимся теперь к действиям алгебр Хопфа.

 Напомним, что для функтора вложения $\mathbf{Hopf}_\mathbbm{k} \to \mathbf{Bialg}_\mathbbm{k}$
 существует правый сопряжённый функтор
$H_r \colon \mathbf{Bialg}_\mathbbm{k} \to \mathbf{Hopf}_\mathbbm{k}$.
(См. \S\ref{SectionFunctorsHlHr}.) Пусть ${}_\square \mathcal{H}(A,V) := H_r({}_\square \mathcal{B}(A,V))$.
Определим действие $\psi^{\mathbf{Hopf}}_{A,V} \colon {}_\square \mathcal{H}(A,V) \otimes A \to A$
как композицию коединицы ${}_\square \mathcal{H}(A,V)
\to {}_\square \mathcal{B}(A,V)$ сопряжения, тензорно умноженной на $\id_A$,
и  действия $\psi_{A,V}$.
Тогда для любой алгебры Хопфа $H$ и действия $\psi \colon H \otimes A \to A$, где $\cosupp \psi \subseteq V$, 
существует единственный такой гомоморфизм алгебр Хопфа $\varphi$,
что следующая диаграмма коммутативна:
$$\xymatrix{ H \otimes A \ar[r]^(0.6){\psi} \ar@{-->}[d]_{\varphi \otimes \id_A} & A \\
{}_\square \mathcal{H}(A,V) \otimes A  \ar[ru]_(0.6){\psi^{\mathbf{Hopf}}_{A,V}}  } $$ 
 Будем называть ${}_\square \mathcal{H}(A,V)$ \textit{$V$-универсальной
 действующей алгеброй Хопфа}.\label{NotationsquareHAV}

 Алгебру Хопфа ${}_\square \mathcal{H}(A,\End_\mathbbm{k}(A))$ будем называть \textit{универсальной действующей алгеброй Хопфа} на $\Omega$-алгебре $A$ и обозначать её через ${}_\square \mathcal{H}(A)$.\label{NotationsquareHA}
 
\begin{definition}\label{DefinitionUnivHopfMod}
Если $\psi \colon H \otimes A \to A$~"--- действие некоторой алгебры Хопфа $H$ на $\Omega$-алгебре $A$,
 то алгебра Хопфа ${}_\square \mathcal{H}(A,\cosupp \psi)$
называется \textit{универсальной алгеброй Хопфа} модульной структуры $\psi$.
\end{definition}

Это понятие было введено автором в работе~\cite{ASGordienko20ALAgoreJVercruysse}.
Из определения действия  $\psi^{\mathbf{Hopf}}_{A,\cosupp \psi}$
следует, что оно
является универсальным среди всех действий алгебр Хопфа на $A$, эквивалентных $\psi$.
В частности,
$\cosupp \psi \subseteq \cosupp\left(\psi^{\mathbf{Hopf}}_{A,\cosupp \psi}\right)$.
Из определения ${}_\square \mathcal{C}(A,A,\cosupp \psi)$
 следует обратное включение, откуда 
$\cosupp \psi^{\mathbf{Hopf}}_{A,\cosupp \psi} = \cosupp \psi$,
т.е. действие $\psi^{\mathbf{Hopf}}_{A,\cosupp \psi}$ 
эквивалентно $\cosupp \psi$.

Если $A$~"--- алгебра с единицей, то её можно рассматривать и как $\lbrace \mu \rbrace$-алгебру,
и как $\lbrace \mu, u \rbrace$-алгебру. 
Это приводит к двум различным определениям алгебры Хопфа ${}_\square \mathcal{H}(A,V)$.
Ниже в теореме~\ref{TheoremHopfActionAlgebraUnital}
мы доказываем, что в случае, когда $V$ является коносителем отображения $\psi 
\colon H \otimes A \to A$, задающего на $A$ структуру модульной алгебры с единицей,
т.е. $\psi(h \otimes 1_A)=\varepsilon(h) 1_A$ для всех $h\in H$, эти две $V$-универсальные
действующие алгебры Хопфа изоморфны.

Сперва докажем следующее предложение:

\begin{proposition}\label{PropositionHopfUnivModUnitality} Пусть $A$~"---
$H$-модульная алгебра над некоторой алгеброй Хопфа $H$ над полем $\mathbbm{k}$. Предположим, что
существует единичный элемент $1_A \in A$, который является общим собственным вектором
для всех операторов из $H$.
Тогда $A$~"--- $H$-модульная алгебра с единицей.
В частности, структура модульной алгебры с единицей над алгеброй Хопфа может быть
эквивалентна только модульной структуре с единицей.
\end{proposition}
\begin{proof}
Обозначим через $\lambda \in H^*$ такую линейную функцию, что $h1_A = \lambda(h)1_A$. 
Поскольку $A$~"--- $H$-модульная алгебра, 
для всех $h_1,h_2 \in H$ справедливо равенство $(h_1 h_2) 1_A=h_1 (h_2 1_A)$,
откуда $\lambda(h_1 h_2)=\lambda(h_1)\lambda(h_2)$.
Из равенства $h(1_A 1_A)=h 1_A$
следует, что $\lambda(h_{(1)})\lambda(h_{(2)}) = \lambda(h)$ для всех $h\in H$.
Равенство $\lambda(1_H)=1_\mathbbm{k}$ следует из равенства $1_H 1_A = 1_A$.
 Следовательно, \begin{equation*}\begin{split}
\lambda(h) =\lambda(h_{(1)})\varepsilon(h_{(2)})\lambda(1_H)=\lambda(h_{(1)})\lambda(\varepsilon(h_{(2)})1_H)=
  \lambda(h_{(1)})\lambda(h_{(2)} Sh_{(3)})=\\=
  \lambda(h_{(1)})\lambda(h_{(2)})\lambda(Sh_{(3)}) =
 \lambda(h_{(1)})\lambda(Sh_{(2)})  =\lambda(h_{(1)}(Sh_{(2)}))
 =\varepsilon(h)\lambda(1_H)=\varepsilon(h)\end{split}\end{equation*}
 и $A$~"--- $H$-модульная алгебра с единицей.
\end{proof}

\begin{theorem}\label{TheoremHopfActionAlgebraUnital}
Пусть $A$~"--- алгебра с единицей над полем $\mathbbm{k}$
и пусть $V\subseteq \End_\mathbbm{k}(A)$~"--- подалгебра, содержащая тождественный оператор.
Предположим, что $V 1_A = \mathbbm{k} 1_A$. (Например, $V=\cosupp \psi$ для некоторого действия $\psi$, задающего на $A$ структуру модульной алгебры с единицей.)
Тогда алгебра Хопфа ${}_\square \mathcal{H}(A,V)$ и её действие $\psi^{\mathbf{Hopf}}_{A,V}$
не зависят от того, рассматриваем ли мы $A$ как $\lbrace \mu \rbrace$-алгебру или как $\lbrace \mu, u \rbrace$-алгебру.
\end{theorem}
\begin{proof}
Достаточно доказать, что отображение $\psi^\mathbf{Hopf}_{A,V} \colon {}_\square \mathcal{H}(A,V) \otimes A
\to A$ задаёт на $A$ структуру модульной алгебры с единицей, даже если мы рассматриваем
  $A$ как $\lbrace \mu \rbrace$-алгебру.
Равенство $V 1_A = \mathbbm{k} 1_A$ означает, что $1_A$
является общим собственным вектором для всех операторов из ${}_\square \mathcal{H}(A,V)$.
В силу предложения~\ref{PropositionHopfUnivModUnitality} действие $\psi^\mathbf{Hopf}_{A,V} \colon {}_\square \mathcal{H}(A,V) \otimes A
\to A$ задаёт на $A$ структуру модульной алгебры с единицей и, следовательно,
 ${}_\square \mathcal{H}(A,V)$~"--- $V$-универсальная действующая алгебра Хопфа для алгебры $A$,
 рассматриваемой и как $\lbrace \mu, u \rbrace$-алгебра.
 \end{proof}

\section{Кодействия}
\label{SectionCoactions}

Пусть $B$~"--- биалгебра, а $A$~"--- $\Omega$-алгебра над полем $\mathbbm{k}$.
Будем говорить, что линейное отображение $\rho \colon A \to A \otimes B$
является \textit{$B$-кодействием} на алгебре $A$ или, что то же самое, $\rho$ задаёт на $A$
структуру \textit{(правой) $B$-комодульной $\Omega$-алгебры},
если выполнены следующие два условия: \begin{enumerate}
\item $\rho$ задаёт на $A$ структуру (коунитального) правого $B$-комодуля;
\item $\rho$ является коизмерением из $A$ в $A$.
\end{enumerate}

Иными словами, правая $B$-комодульная $\Omega$-алгебра~"--- это не что иное, как $\Omega$-алгебра
в моноидальной категории правых $B$-комодулей.

Если $\Omega=\lbrace\mu, u\rbrace$, а $A$~"--- алгебра с единицей, то
мы получаем обычное определение комодульной алгебры с единицей.
Если $\Omega=\lbrace\mu\rbrace$, то алгебра $A$ является модульной алгеброй необязательно с единицей.
Наконец, если $\Omega=\lbrace\Delta, \varepsilon\rbrace$, а $A$~"--- коалгебра,
мы получаем определение комодульной коалгебры.

Для кодействий мы будем использовать понятия коносителя, эквивалентности и предпорядок,
введённые в \S~\ref{SectionLinearMaps}.

Заметим, что для кодействия $\rho \colon A \to A \otimes B$ его коноситель $\cosupp\rho$~"--- это
в точности образ соответствующего гомоморфизма алгебр $B^* \to \End_\mathbbm{k}(A)$. 
В частности, $\cosupp\rho$~"--- это подалгебра в $\End_\mathbbm{k}(A)$, содержащая тождественный оператор.
Из того, что $A$~"--- это $B$-комодуль, следует, что $\supp\rho$~"--- это подкоалгебра в $B$.
Учитывая, что ядро гомоморфизма $B^* \to \End_\mathbbm{k}(A)$ совпадает с ядром 
сюръективного гомоморфизма $B^* \to (\supp\rho)^*$, полученного
ограничением линейных функций, заданных на векторном пространстве
$B$, на подпространство $\supp\rho$, получаем изоморфизм алгебр $\cosupp\rho \cong (\supp\rho)^*$. 

\begin{theorem}\label{TheoremBsquareBialgebra}
Пусть $A$~"--- $\Omega$-алгебра над полем $\mathbbm{k}$, а $V\subseteq \End_\mathbbm{k}(A)$~"---
замкнутая в конечной топологии
подалгебра, содержащая единицу, такая, что существует $V$-универсальная коизмеряющая алгебра $\mathcal{B}^\square(A,V) := \mathcal{A}^\square(A,A,V)$.
Тогда на $\mathcal{B}^\square(A,V)$
можно так задать структуру биалгебры, что для любой биалгебры $B$
и любого кодействия $\rho \colon  A \to A \otimes B$, такого, что $\cosupp \rho \subseteq V$,  единственный гомоморфизм алгебр $\varphi$,
делающий диаграмму ниже коммутативной, на самом деле является гомоморфизмом биалгебр:
\begin{equation}\label{EqTheoremBsquareBialgebra}\xymatrix{ A \ar[r]^(0.3){\rho_{A,V}} \ar[rd]_\rho & A \otimes \mathcal{B}^\square(A,V) \ar@{-->}[d]^{\id_A\otimes
\varphi} \\
& A\otimes B} \end{equation}
(Здесь $\rho_{A,V} := \rho_{A,A,V}$.)
\end{theorem}

Основываясь на теореме~\ref{TheoremBsquareBialgebra}, будем называть биалгебру $\mathcal{B}^\square(A,V)$  \textit{$V$-универсальной кодействующей биалгеброй}.\label{NotationBsquareAV}

\begin{proof}[Доказательство теоремы~\ref{TheoremBsquareBialgebra}]
Заметим, что $$(\rho_{A,V}\otimes \id_{\mathcal{B}^\square(A,V)})\rho_{A,V}
\colon A \to A\otimes \mathcal{B}^\square(A,V) \otimes \mathcal{B}^\square(A,V)$$
также является коизмерением и, поскольку $V$~"--- подалгебра, замкнутая в конечной топологии, справедливо включение $\cosupp((\rho_{A,V}\otimes \id_{\mathcal{B}^\square(A,V)})\rho_{A,V})\subseteq V$.
Следовательно, существует единственный гомоморфизм алгебр
$\Delta \colon \mathcal{B}^\square(A,V) \to \mathcal{B}^\square(A,V) \otimes \mathcal{B}^\square(A,V)$,
делающий диаграмму ниже коммутативной:
$$\xymatrix{ A \ar[rr]^{\rho_{A,V}} \ar[d]_{\rho_{A,V}} && A\otimes \mathcal{B}^\square(A,V)
\ar@{-->}[d]^{\id_A\otimes \Delta} \\
A\otimes \mathcal{B}^\square(A,V) \ar[rr]_(0.4){\rho_{A,V}\otimes \id_{\mathcal{B}^\square(A,V)}} &&  A\otimes \mathcal{B}^\square(A,V) \otimes \mathcal{B}^\square(A,V)}$$
Определим теперь коумножение в $\mathcal{B}^\square(A,V)$
при помощи отображения~$\Delta$.

Для того, чтобы задать коединицу, рассмотрим тривиальное отображение $$A
\mathrel{\widetilde\to} A  \otimes \mathbbm{k},$$ которое является коизмерением.
Существует единственный гомоморфизм алгебр $\varepsilon \colon   \mathcal{B}^\square(A,V) \to \mathbbm{k}$,
делающий диаграмму ниже коммутативной:
$$\xymatrix{ A \ar[r]^(0.3){\rho_{A,V}} \ar[rd]^{\sim} & A \otimes \mathcal{B}^\square(A,V)
\ar[d]^{\id_A \otimes \varepsilon}\\
& A \otimes \mathbbm{k}}$$
Определим теперь коединицу в $\mathcal{B}^\square(A,V)$ при помощи отображения~$\varepsilon$.

Теперь нужно доказать, что коумножение $\Delta$ коассоциативно и что $\varepsilon$
и $\Delta$ действительно удовлетворяют аксиомам коединицы.

Рассмотрим диаграмму

{\scriptsize
$$\xymatrix{ A \ar[rr]^{\rho_{A,V}} \ar[dd]_{\rho_{A,V}} \ar[rd]^{\rho_{A,V}} && A \otimes \mathcal{B}^\square(A,V) \ar[rd]^{\qquad \rho_{A,V} \otimes \id_{\mathcal{B}^\square(A,V)}} \ar'[d][dd]_{\rho_{A,V} \otimes \id_{\mathcal{B}^\square(A,V)}} & \\
& A \otimes \mathcal{B}^\square(A,V) \ar[rr]^(0.3){\id_A \otimes \Delta} \ar[dd]_(0.3){\rho_{A,V} \otimes \id_{\mathcal{B}^\square(A,V)}} && A \otimes \mathcal{B}^\square(A,V) \otimes \mathcal{B}^\square(A,V)\ar[dd]_(0.3){\rho_{A,V} \otimes \id_{\mathcal{B}^\square(A,V) \otimes \mathcal{B}^\square(A,V)}} \\
A \otimes \mathcal{B}^\square(A,V)  \ar'[r][rr]_(0.3){\id_A \otimes \Delta}
 \ar[rd]^{\id_A \otimes \Delta} && A \otimes \mathcal{B}^\square(A,V) \otimes \mathcal{B}^\square(A,V) \ar[rd]_(0.35){\id_A \otimes \Delta \otimes \id_{\mathcal{B}^\square(A,V)}\qquad} &\\
 & A \otimes \mathcal{B}^\square(A,V) \otimes \mathcal{B}^\square(A,V) 
 \ar[rr]_(0.4){\id_A  \otimes \id_{\mathcal{B}^\square(A,V)}\otimes \Delta}  && A \otimes \mathcal{B}^\square(A,V) \otimes \mathcal{B}^\square(A,V) \otimes \mathcal{B}^\square(A,V) } $$}

Правая грань коммутативна, поскольку эта диаграмма является диаграммой, определяющей
отображение $\Delta$, тензорно умноженной на $\id_{\mathcal{B}^\square(A,V)}$.
Коммутативность левой, верхней и задней граней также следует из определения
отображения $\Delta$.
Передняя грань коммутативна, так как  обе композиции равны $\rho_{A,V} \otimes \Delta$.
Следовательно, композиции, отвечающие нижней грани, после их композиции с $\rho_{A,V}$
также оказываются равными.

Отсюда из универсального свойства алгебры $\mathcal{B}^\square(A,V)$ следует коммутативность диаграммы
 $$\xymatrix{ \mathcal{B}^\square(A,V) \ar[rrr]_(0.4)\Delta\ar[d]^\Delta &&& \mathcal{B}^\square(A,V) \otimes \mathcal{B}^\square(A,V) \ar[d]^{\Delta\otimes \id_{\mathcal{B}^\square(A,V)}}\\ 
 \mathcal{B}^\square(A,V) \otimes \mathcal{B}^\square(A,V) \ar[rrr]_(0.42){ \id_{\mathcal{B}^\square(A,V)}\otimes\Delta} &&& \mathcal{B}^\square(A,V) \otimes \mathcal{B}^\square(A,V) \otimes \mathcal{B}^\square(A,V) }$$
 Это в точности означает, что коумножение $\Delta$ коассоциативно.

Рассмотрим диаграмму
$$\xymatrix{ A \ar[rr]^{\rho_{A,V}} \ar[rd]^{\rho_{A,V}} \ar[dd]_{\rho_{A,V}} && A \otimes \mathcal{B}^\square(A,V) \ar@{=}[ld]\ar[dd]^{\rho_{A,V} \otimes \id_{\mathcal{B}^\square(A,V)}} \\
& A \otimes \mathcal{B}^\square(A,V) \ar[dd]^(0.3){\sim} & \\
A \otimes \mathcal{B}^\square(A,V)\ar[rd]^{\sim} \ar'[r][rr]^(0.3){\id_A\otimes \Delta} && A \otimes \mathcal{B}^\square(A,V) \otimes \mathcal{B}^\square(A,V) \ar[ld]^{\qquad\id_A\otimes\varepsilon\otimes\id_{\mathcal{B}^\square(A,V)}} \\
& A \otimes \mathbbm{k} \otimes \mathcal{B}^\square(A,V) & } $$

Коммутативность верхней и левой граней очевидна.
Коммутативность правой грани следует из определения
отображения $\varepsilon$,
а коммутативность задней грани 
следует из определения отображения $\Delta$. 
Следовательно, гомоморфизмы алгебр, образующие нижнюю грань, становятся равными после их композиции
с $\rho_{A,V}$. 
 Отсюда из универсального свойства алгебры $\mathcal{B}^\square(A,V)$
следует коммутативность диаграммы
 $$\xymatrix{ \mathcal{B}^\square(A,V) \ar[rr]^{\Delta} \ar[rd]^{\sim} 
  & & \mathcal{B}^\square(A,V)\otimes \mathcal{B}^\square(A,V) \ar[ld]^{\qquad\varepsilon\otimes \id_{\mathcal{B}^\square(A,V)}}\\
 & \mathbbm{k} \otimes \mathcal{B}^\square(A,V)  & }$$
 Поэтому отображение $\varepsilon$ 
удовлетворяет аксиоме левой коединицы.
 
 Аналогично, рассматривая диаграмму
  $$\xymatrix{ A \ar[rr]^{\rho_{A,V}}\ar[dd]_{\rho_{A,V}} \ar[rd]^{\sim} && A\otimes \mathcal{B}^\square(A,V)
  \ar[ld]^{\qquad\id_A\otimes \varepsilon} \ar[dd]^{\rho_{A,V} \otimes \id_{\mathcal{B}^\square(A,V)}} \\
  & A\otimes \mathbbm{k} \ar[dd]^(0.3){\rho_{A,V} \otimes \id_\mathbbm{k}} & \\
 A \otimes \mathcal{B}^\square(A,V)\ar[rd]^{\sim} \ar'[r][rr]_(0.3){\id_A\otimes \Delta} && A \otimes \mathcal{B}^\square(A,V) \otimes \mathcal{B}^\square(A,V) \ar[ld]^{\qquad\id_A\otimes\id_{\mathcal{B}^\square(A,V)}\otimes\varepsilon} \\
& A \otimes \mathcal{B}^\square(A,V) \otimes \mathbbm{k} & } $$
  получаем коммутативность диаграммы
 $$\xymatrix{ \mathcal{B}^\square(A,V) \ar[rr]^{\Delta} \ar[rd]^{\sim} 
  & & \mathcal{B}^\square(A,V)\otimes \mathcal{B}^\square(A,V) \ar[ld]^{\qquad \id_{\mathcal{B}^\square(A,V)}\otimes\varepsilon}\\
 & \mathcal{B}^\square(A,V) \otimes  \mathbbm{k} & }$$
 Следовательно, отображение $\varepsilon$
удовлетворяет аксиоме левой единицы, и $\mathcal{B}^\square(A,V)$ действительно биалгебра.
 
 Пусть $B$~"---  биалгебра, а $\rho \colon  A \to A\otimes B$ такое кодействие,
 что $\cosupp \rho \subseteq V$.
 Обозначим через $\varphi \colon \mathcal{B}^\square(A,V) \to B$ 
единственный гомоморфизм алгебр, делающий диаграмму~\eqref{EqTheoremBsquareBialgebra} коммутативной.
Такой гомоморфизм существует в силу теоремы~\ref{TheoremUnivComeasExistence}.
Докажем, что в этом случае $\varphi$ является гомоморфизмом биалгебр.

Рассмотрим диаграмму
$$\xymatrix{  A \ar@{=}[rr]\ar[dd]_{\rho_{A,V}}\ar[rd]^{\rho_{A,V}} && A \ar'[d][dd]_{\rho}  \ar[rd]^{\rho}& \\
& A \otimes \mathcal{B}^\square(A,V) \ar[rr]^(0.4){\id_A\otimes\varphi}\ar[dd]_(0.3){\rho_{A,V} \otimes \id_{\mathcal{B}^\square(A,V)}} && A\otimes B \ar[dd]^{\rho \otimes \id_{B}} \\
A \otimes \mathcal{B}^\square(A,V) \ar'[r][rr]^(0.3){\id_A\otimes\varphi}\ar[rd]^{\id_A \otimes \Delta} && A\otimes B \ar[rd]^{\id_A \otimes \Delta_B}& \\
& A \otimes \mathcal{B}^\square(A,V)\otimes \mathcal{B}^\square(A,V) \ar[rr]^{\id_A\otimes\varphi\otimes\varphi} & & A\otimes B\otimes B } $$
 (Здесь $\Delta_B \colon B \to B\otimes B$~"--- коумножение в $B$.)
 Верхняя и задняя грани коммутативны по определению гомоморфизма $\varphi$.
Передняя грань коммутативна, поскольку она совпадает с верхней (а также с задней) гранью,
тензорно умноженной на $\varphi$ справа.
Левая и правая грани коммутативны, поскольку и $\rho$, и $\rho_{A,V}$ задают на $A$
структуру правого комодуля. 
 
 Следовательно, после композиции с $\rho_{A,V}$
нижняя грань также становится коммутативной, а из универсального свойства
алгебры $\mathcal{B}^\square(A,V)$ 
следует коммутативность диаграммы
 $$\xymatrix{\mathcal{B}^\square(A,V)\ar[d]_{\Delta}  \ar[r]^{\varphi} & B \ar[d]^{\Delta_B}\\
   \mathcal{B}^\square(A,V) \otimes \mathcal{B}^\square(A,V)  \ar[r]_(0.7){\varphi\otimes\varphi} & B \otimes B   \\ 
  }$$ Отсюда отображение $\varphi$ сохраняет коумножение.
 
 Рассмотрим диаграмму $$\xymatrix{ & A \ar[lddd]_{\rho} 
 \ar[dd]^(0.7){\rho_{A,V}}  \ar[rddd]^{\sim} & \\ & & \\
 &  A  \otimes \mathcal{B}^\square(A,V) 
 \ar[rd]_{\id_A\otimes\varepsilon\ } \ar[ld]^{\ \id_A \otimes \varphi}& \\
 A \otimes B \ar[rr]_{\id_A\otimes\varepsilon_B} & & A\otimes \mathbbm{k} }$$ (Здесь $\varepsilon_B$~"---
 коединица в $B$.)
 
 Большой и правый треугольники коммутативны в силу коунитальности комодульных структур на $A$.
 Левый треугольник коммутативен по определению отображения $\varphi$.
Следовательно, нижний треугольник становится коммутативным после композиции с отображением $\rho_{A,V}$.
Отсюда из универсального свойства алгебры $\mathcal{B}^\square(A,V)$ 
следует коммутативность диаграммы
 $$\xymatrix{ \mathcal{B}^\square(A,V)\ar[rr]^(0.6){\varphi}\ar[rd]_{\varepsilon} & & B \ar[ld]^{\varepsilon_B}  \\
  & \mathbbm{k}   &}$$
 Это означает, что отображение $\varphi$
сохраняет коединицу и действительно является гомоморфизмом биалгебр.
\end{proof}

Биалгебру $\mathcal{B}^\square(A,\End_\mathbbm{k}(A))$ будем называть \textit{универсальной кодействующей биалгеброй} на $\Omega$-алгебре $A$ и обозначать её через $\mathcal{B}^\square(A)$.\label{NotationBsquareA} 

\begin{remark}
Если $A$~"--- ассоциативная алгебра с единицей над полем~$\mathbbm{k}$,
то биалгебра $\mathcal{B}^\square(A)$~"--- это в точности
универсальная кодействующая биалгебра $\alpha(A,A)$ Д.~Тамбары~\cite{Tambara}.
Заметим, что если алгебра $A$ конечномерна, то, очевидно, справедливы условия теоремы~\ref{TheoremUnivComeasExistence}, обеспечивающие существование алгебры $\mathcal{B}^\square(A)$. 
\end{remark}
\begin{remark}\label{RemarkManinEnd}
Пусть $A$~"--- ассоциативная $\mathbb Z$-градуированная алгебра с единицей над полем~$\mathbbm{k}$,
причём $\dim A^{(n)} < +\infty$ для всех $n\in\mathbb Z$.
Обозначим через $V\subseteq \End_\mathbbm{k}(A)$ подалгебру всех линейных операторов, сохраняющих градуировку.
Очевидно, $V$ поточечно конечномерна. Кроме того, если какой-то оператор $f_0\in\End_\mathbbm{k}(A)$ нарушает
градуировку, то $f_0(a)=b$ для некоторых $n\in\mathbb Z$, $a\in A^{(n)}$ и $b\notin A^{(n)}$.
При этом окрестность $\lbrace f\in \End_\mathbbm{k}(A) \mid f(a)=b \rbrace$ точки $f_0$ в конечной топологии
не пересекается с $V$. Отсюда множество $V$ замкнуто в конечной топологии
и в силу теоремы~\ref{TheoremUnivComeasExistence} существует биалгебра $\mathcal{B}^\square(A, V)$, которая является
универсальной кодействующей биалгеброй $\underline{\mathrm{end}}(A)$ Ю.\,И.~Манина~\cite{Manin}.
\end{remark}

Рассмотрим теперь кодействия алгебр Хопфа.

 Напомним, что для функтора вложения $\mathbf{Hopf}_\mathbbm{k} \to \mathbf{Bialg}_\mathbbm{k}$
 существует левый сопряжённый функтор $H_l \colon \mathbf{Bialg}_\mathbbm{k} \to \mathbf{Hopf}_\mathbbm{k}$.
(См. \S\ref{SectionFunctorsHlHr}.) Обозначим $\mathcal{H}^\square(A,V) := H_l(\mathcal{B}^\square(A,V))$.
Определим кодействие $\rho^{\mathbf{Hopf}}_{A,V} \colon A \to A \otimes \mathcal{H}^\square(A,V)$
как композицию отображения $\rho_{A,V}$ и единицы $\mathcal{B}^\square(A,V) \to \mathcal{H}^\square(A,V)$
 сопряжения, тензорно умноженной на $\id_A$.
 Тогда для любой алгебры Хопфа  $H$ и любого кодействия
$\rho \colon A \to A \otimes H$, такого, что $\cosupp \rho \subseteq V$,
существует единственный гомоморфизм алгебр Хопфа $\varphi$, делающий коммутативной
диаграмму
$$\xymatrix{ A \ar[rd]_\rho \ar[r]^(0.3){\rho^{\mathbf{Hopf}}_{A,V}} & A \otimes \mathcal{H}^\square(A,V)
\ar@{-->}[d]^{\id_A \otimes \varphi} \\
& A \otimes H }$$
 Назовём алгебру Хопфа $\mathcal{H}^\square(A,V)$ \textit{$V$-универсальной кодействующей алгеброй Хопфа}.
 \label{NotationHsquareAV} 
 
 Алгебру Хопфа $\mathcal{H}^\square(A,\End_\mathbbm{k}(A))$ будем называть \textit{универсальной кодействующей  алгеброй Хопфа} на $\Omega$-алгебре $A$ и обозначать её через $\mathcal{H}^\square(A)$.\label{NotationHsquareA} 

\begin{remark}\label{RemarkManinAut}
Пусть $A$~"--- ассоциативная $\mathbb Z$-градуированная алгебра с единицей над полем~$\mathbbm{k}$,
причём $\dim A^{(n)} < +\infty$ для всех $n\in\mathbb Z$. Обозначим через
$V\subseteq \End_\mathbbm{k}(A)$ подалгебру всех линейных операторов, сохраняющих градуировку.
Как было показано в замечании~\ref{RemarkManinEnd},
подалгебра $V$ поточечно конечномерна и замкнута в конечной топологии.
Поэтому существует алгебра Хопфа $\mathcal{H}^\square(A, V)$, которая является
универсальной кодействующей алгеброй Хопфа $\underline{\mathrm{aut}}(A)$ Ю.\,И.~Манина~\cite{Manin}.
\end{remark}

 \begin{definition}\label{DefinitionUnivHopfComod}
Если $\rho \colon  A \to A \otimes H$~"--- кодействие некоторой алгебры Хопфа $H$ на $\Omega$-алгебре $A$,
 то алгебра Хопфа $\mathcal{H}^\square(A,\cosupp \rho)$
называется \textit{универсальной алгеброй Хопфа} комодульной структуры $\rho$.
\end{definition}

Это понятие было введено автором в работе~\cite{ASGordienko20ALAgoreJVercruysse}.
Из определения кодействия  $\rho^{\mathbf{Hopf}}_{A,\cosupp \rho}$
следует, что оно
является универсальным среди всех кодействий алгебр Хопфа на $A$, эквивалентных $\rho$.
В частности,
$\cosupp \rho \subseteq \cosupp\left(\rho^{\mathbf{Hopf}}_{A,\cosupp \rho}\right)$.
Из определения $\mathcal{A}^\square(A,A,\cosupp \rho)$
 следует обратное включение, откуда 
$\cosupp \rho^{\mathbf{Hopf}}_{A,\cosupp \rho} = \cosupp \rho$,
т.е. кодействие $\rho^{\mathbf{Hopf}}_{A,\cosupp \rho}$ 
эквивалентно $\cosupp \rho$.
 
 \begin{theorem}\label{TheoremBHsquareVExistence}
 Пусть $A$~"--- $\Omega$-алгебра над полем $\mathbbm{k}$, а $V\subseteq \End_\mathbbm{k}(A)$~"--- подалгебра, содержащая тождественный оператор и замкнутая в конечной топологии.
 Тогда биалгебра $\mathcal{B}^\square(A,V)$ 
существует, если и только если подалгебра, порождённая
подалгебрами $\cosupp\rho$ для всевозможных биалгебр $B$ и кодействий $\rho \colon  A \to A\otimes B$,
таких, что $\cosupp \rho \subseteq V$, поточечно конечномерна.
 Аналогично, алгебра Хопфа  $\mathcal{H}^\square(A,V)$ существует, если и только если подалгебра, порождённая
подалгебрами $\cosupp\rho$ для всевозможных алгебр Хопфа $H$ и кодействий $\rho \colon  A \to A\otimes H$,
таких, что $\cosupp \rho \subseteq V$, поточечно конечномерна.
 \end{theorem}
\begin{proof} Доказательство в обоих случаях одинаково. В связи с этим
приведём доказательство лишь для случая биалгебр.

Необходимость очевидна: если биалгебра $\mathcal{B}^\square(A,V)$ действительно существует, то $\cosupp \rho \subseteq \cosupp \rho_{A,V}$ для всех биалгебр $B$ и кодействий $\rho \colon  A \to A\otimes B$, таких, что
 $\cosupp \rho \subseteq V$. При этом пространство $\cosupp \rho_{A,V}$ поточечно конечномерно.

Обратно, пусть $W \subseteq \End_\mathbbm{k}(A)$~"--- подалгебра, порождённая всеми такими $\cosupp\rho\subseteq V$.
Если подалгебра $W$ поточечно конечномерна, тогда в силу леммы~\ref{LemmaClosurePwFD}
замыкание $\overline W$ алгебры $W$ в конечной топологии также поточечно конечномерно.
Отсюда согласно теореме~\ref{TheoremUnivComeasExistence} существует биалгебра $\mathcal{B}^\square(A,\overline W)$. Поскольку $\overline W \subseteq V$, получаем $\mathcal{B}^\square(A,V)=\mathcal{B}^\square(A,\overline W)$.
\end{proof} 
\begin{corollary}\label{CorollaryManinExistence} Универсальная кодействующая алгебра Хопфа $\mathcal{H}^\square(A)$ 
существует для $\Omega$-алгебры $A$, если и только если подалгебра, порождённая
подалгебрами $\cosupp\rho$ для всевозможных алгебр Хопфа $H$ и кодействий $\rho \colon  A \to A\otimes H$, поточечно конечномерна.
\end{corollary}

Предложение~\ref{PropositionHopfCoactionAlgebraUnital},
которое доказывается ниже, является двойственным к предложению~\ref{PropositionHopfUnivModUnitality}.

\begin{proposition}\label{PropositionHopfCoactionAlgebraUnital}
Пусть $A$~"--- алгебра с единицей, а $\rho \colon A \to A \otimes H$~"--- кодействие
алгебры Хопфа $H$ над полем $\mathbbm{k}$.
Предположим, что $\rho(1_A)=1_A \otimes h$ для некоторого $h\in H$.
Тогда $h=1_H$, т.е. кодействие $\rho$ задаёт на $A$ структуру $H$-комодульной алгебры с единицей.
\end{proposition}
\begin{proof} Поскольку $\rho$~"--- кодействие, из определения комодуля следует, что $\varepsilon(h)=1_\mathbbm{k}$ и
$\Delta(h) = h \otimes h$. Рассматривая равенство $1_A^2=1_A$, получаем
$h^2=h$. Из равенства $(Sh)h=h(Sh)=1_H$ следует, что элемент $h$ обратим.
Отсюда, используя равенство $h^2=h$, заключаем, что $h=1_H$.
\end{proof}

Из теоремы~\ref{TheoremHopfCoactionAlgebraUnital}, которая доказывается ниже, следует, что если пространство $V$ является коносителем некоторой структуры комодульной алгебры с единицей, то нет разницы, включаем ли мы символ $u$ отображения, отвечающего единице, в $\Omega$ или нет.

\begin{theorem}\label{TheoremHopfCoactionAlgebraUnital}
Пусть $A$~"--- алгебра с единицей, а $V\subseteq \End_\mathbbm{k}(A)$~"--- подалгебра, содержащая тождественный оператор, такая, что существует $\mathcal{H}^\square(A,V)$.
Предположим, что $V 1_A = \mathbbm{k} 1_A$. (Например, $V=\cosupp \rho$, где $\rho$~"--- некоторое кодействие, задающее на $A$ структуру комодульной алгебры с единицей.)
Тогда алгебра Хопфа $\mathcal{H}^\square(A,V)$ 
не зависит от того, рассматриваем ли мы алгебру $A$ как $\lbrace \mu \rbrace$-алгебру или как $\lbrace \mu, u \rbrace$-алгебру.
\end{theorem}
\begin{proof} Сперва рассмотрим $A$ как $\lbrace \mu \rbrace$-алгебру.
Из равенства $V 1_A = \mathbbm{k} 1_A$ следует, что
$\rho^\mathbf{Hopf}_{A,V}(1_A)=1_A \otimes h$ для некоторого $h\in \mathcal{H}^\square(A,V)$.
В силу предложения~\ref{PropositionHopfCoactionAlgebraUnital} кодействие $\rho^\mathbf{Hopf}_{A,V}$ задаёт на $A$ структуру комодульной алгебры с единицей. Отсюда $\mathcal{H}^\square(A,V)$~"--- $V$-универсальная кодействующая алгебра Хопфа на $A$, даже если рассматривать $A$ как $\lbrace \mu, u \rbrace$-алгебру.
\end{proof}

         \section{Двойственность между действиями и кодействиями}\label{SectionDualityActionsCoactions}

Оказывается, что в случае (ко)действий изоморфизм $\theta$
из теоремы~\ref{TheoremUniv(Co)measDuality} в действительности является изоморфизмом биалгебр.

\begin{theorem}\label{TheoremUniv(Co)actDuality}
Пусть $A$~"--- $\Omega$-алгебра над полем $\mathbbm{k}$, а $V\subseteq \End_\mathbbm{k}(A)$~"--- поточечно конечномерная
подалгебра, содержащая тождественный оператор и замкнутая в конечной топологии.
Обозначим через $\theta$ единственный гомоморфизм биалгебр,
делающий диаграмму ниже коммутативной:
$$\xymatrix{ {}_\square \mathcal{B}(A,V) \otimes A \ar[dr]^(0.6){\psi_{A,V}} & \\
                   &   A \\
\mathcal{B}^\square(A,V)^\circ \otimes A \ar[ru]_(0.6){\widehat{\rho_{A,V}}}\ar@{-->}[uu]^{\theta \otimes \id_A} } $$
Тогда $\theta$~"--- изоморфизм биалгебр.
В частности, для любой конечномерной $\Omega$-алгебры $A$
существует изоморфизм биалгебр $
{}_{\square}\mathcal{B}(A) \cong \mathcal{B}^{\square}(A)^{\circ}$.
\end{theorem}
\begin{proof}
Существование гомоморфизма биалгебр $\theta$ следует из теоремы~\ref{TheoremsquareBBialgebra}. 
В силу теоремы~\ref{TheoremUniv(Co)measDuality} гомоморфизм $\theta$ биективен.
\end{proof}

То же верно для $V$-универсальных алгебр Хопфа.

\begin{theorem}\label{TheoremUnivHopf(Co)actDuality}
Пусть $A$~"--- $\Omega$-алгебра над полем $\mathbbm{k}$, а $V\subseteq \End_\mathbbm{k}(A)$~"--- поточечно конечномерная
подалгебра, содержащая тождественный оператор и замкнутая в конечной топологии.
Тогда единственный гомоморфизм $\theta^\mathbf{Hopf} \colon \mathcal{H}^\square(A,V)^\circ \to {}_\square \mathcal{H}(A,V) $ алгебр Хопфа, делающий коммутативной  диаграмму $$\xymatrix{ {}_\square \mathcal{H}(A,V) \otimes A \ar[dr]^(0.6){\psi^\mathbf{Hopf}_{A,V}} & \\
                   &   A \\
\mathcal{H}^\square(A,V)^\circ \otimes A \ar[ru]_(0.6){\widehat{\rho^\mathbf{Hopf}_{A,V}}}\ar@{-->}[uu]^{\theta^\mathbf{Hopf} \otimes \id_A} } $$ в действительности является изоморфизмом.
В частности, для любой конечномерной $\Omega$-алгебры $A$
существует изоморфизм алгебр Хопфа
$
{}_{\square}\mathcal{H}(A) \cong \mathcal{H}^{\square}(A)^{\circ}$.
\end{theorem}
\begin{proof}Пусть $H$~"--- произвольная алгебра Хопфа над полем $\mathbbm{k}$.
Рассмотрим естественные биекции
\begin{equation*}
\begin{split}
\mathbf{Hopf}_\mathbbm{k}(H, {}_\square \mathcal{H}(A,V)) = 
\mathbf{Hopf}_\mathbbm{k}(H, H_r({}_\square \mathcal{B}(A,V)))\cong
\mathbf{Bialg}_\mathbbm{k}(H, {}_\square \mathcal{B}(A,V)) \cong \\ 
\cong \mathbf{Bialg}_\mathbbm{k}(H, \mathcal{B}^\square(A,V)^\circ)
\cong \mathbf{Bialg}_\mathbbm{k}(\mathcal{B}^\square(A,V), H^\circ)
\cong \mathbf{Hopf}_\mathbbm{k}(H_l(\mathcal{B}^\square(A,V)), H^\circ) = \\
= \mathbf{Hopf}_\mathbbm{k}(\mathcal{H}^\square(A,V), H^\circ)
\cong \mathbf{Hopf}_\mathbbm{k}(H, \mathcal{H}^\square(A,V)^\circ).\end{split}\end{equation*}
Следовательно, $\mathcal{H}^\square(A,V)^\circ \cong {}_\square \mathcal{H}(A,V)$,
где изоморфизм соответствует элементу $$\id_{\mathcal{H}^\square(A,V)^\circ}
\in \mathbf{Hopf}_\mathbbm{k}(\mathcal{H}^\square(A,V)^\circ, \mathcal{H}^\square(A,V)^\circ),$$
если взять $H=\mathcal{H}^\square(A,V)^\circ$.
Однако соответствующий элемент множества $\mathbf{Hopf}_\mathbbm{k}(\mathcal{H}^\square(A,V)^\circ, {}_\square \mathcal{H}(A,V))$~"--- это в точности гомоморфизм $\theta^\mathbf{Hopf}$, заданный
универсальными свойствами биалгебр ${}_\square \mathcal{B}(A,V)$ и алгебр Хопфа ${}_\square \mathcal{H}(A,V)$:
$$\xymatrix{ {}_\square \mathcal{H}(A,V) \otimes A \ar[r] \ar@/^5pc/[rrd]^(0.6){\psi^\mathbf{Hopf}_{A,V}} &
{}_\square \mathcal{B}(A,V) \otimes A \ar[rd]^{\psi_{A,V}} & \\
                 &  &   A \\
\mathcal{H}^\square(A,V)^\circ \otimes A \ar[r] \ar@/_5pc/[rru]_(0.6){\widehat{\rho^\mathbf{Hopf}_{A,V}}}\ar@{-->}[uu]^{\theta^\mathbf{Hopf} \otimes \id_A} &  \mathcal{B}^\square(A,V)^\circ \otimes A
\ar[uu]^{\theta \otimes \id_A}
 \ar[ru]_(0.6){\widehat{\rho_{A,V}}} & }$$
Единственность гомоморфизма $\theta^\mathbf{Hopf}$
также следует из универсального свойства алгебры Хопфа ${}_\square \mathcal{H}(A,V)$.
\end{proof}

         \section{Двойственность между конечномерными $\Omega$- и $\Omega^*$-алгебрами}\label{SectionDualityAlgCoalg}
         
В этом параграфе мы изучим, что произойдёт с универсальными (ко)действующими биалгебрами,
если заменить $\Omega$-алгебру $A$ на $\Omega^*$-алгебру $A^*$ (см. \S\ref{SectionOmegaAlgebras}).
Напомним, что через $B^\mathrm{op}$ обозначается биалгебра $B$, в которой в умножении меняются местами аргументы (коумножение остаётся прежним), а через $B^\mathrm{cop}$~"--- биалгебра, в которой в коумножении меняются местами тензорные множители в разложении результата (умножение остаётся прежним).

Для произвольной $\Omega$-алгебры $A$ и подалгебры $V\subseteq \End_\mathbbm{k}(A)$, содержащей тождественный оператор, определим категорию $\mathbf{Actions}(A,V)$,
в которой объектами являются всевозможные действия $\psi \colon B \otimes A \to A$ для всех $\mathbbm{k}$-биалгебр $B$, где $\cosupp \psi \subseteq V$, а морфизмами из $\psi_1 \colon B_1 \otimes A \to A$
в $\psi_2 \colon B_2 \otimes A \to A$  являются всевозможные гомоморфизмы биалгебр $\varphi \colon B_1 \to B_2$, делающие диаграмму ниже коммутативной:

$$\xymatrix{ B_1 \otimes A \ar[r]^(0.6){\psi_1} \ar[d]_{\varphi \otimes \id_A} & A \\
B_2 \otimes A  \ar[ru]_{\psi_2}  } $$

Тогда терминальным объектом категории $\mathbf{Actions}(A,V)$ является действие биалгебры ${}_\square \mathcal{B}(A,V)$ на $A$.

Аналогично, определим категорию $\mathbf{Coactions}(A,V)$,
в которой объектами являются всевозможные кодействия $\rho \colon  A \to A \otimes B$ для всех $\mathbbm{k}$-биалгебр $B$, где $\cosupp \psi \subseteq V$, а морфизмами из $\rho_1 \colon A \to A \otimes B_1$
в $\rho_2 \colon A \to A \otimes B_2$ являются всевозможные гомоморфизмы биалгебр $\varphi \colon B_1 \to B_2$, делающие диаграмму ниже коммутативной:

$$\xymatrix{ A \ar[r]^(0.4){\rho_1} 
\ar[rd]_{\rho_2}
 & A \otimes B_1 \ar[d]^{\id_A \otimes \varphi} \\
& A \otimes B_2} $$

Тогда  начальным объектом категории $\mathbf{Coactions}(A,V)$ является действие биалгебры $\mathcal{B}^\square(A,V)$ на $A$ (если $\mathcal{B}^\square(A,V)$ существует).

\begin{theorem}\label{TheoremDualityActions} Пусть $A$~"--- конечномерная $\Omega$-алгебра над полем $\mathbbm{k}$, а $V\subseteq \End_\mathbbm{k}(A^*)$~"--- подалгебра, содержащая тождественный оператор.
Тогда $${}_\square \mathcal{B}(A^*,V)
\cong {}_\square \mathcal{B}(A,V^{\#})^{\mathrm{op}}\text{\quad(изоморфизм биалгебр),}$$
где $V^{\#} := \lbrace f^* \mid f \in V \rbrace$.
\end{theorem}
\begin{proof} 
Обозначим через $\widetilde{\ }$ функтор $\mathbf{Actions}(A,V^{\#}) \to \mathbf{Actions}(A^*,V)$,
заданный следующим образом: если $\psi \colon B \otimes A \to A$~"--- объект
категории $\mathbf{Actions}(A,V^{\#})$, то
 $\tilde\psi \colon B^\mathrm{op} \otimes A^* \to A^*$ определяется через равенство
$$\tilde\psi(b\otimes a^*)(a):=a^*\left(\psi(b\otimes a)\right)\text{ для всех }b\in B,\ a\in A,\ a^*\in A^*.$$ Другими словами, $\tilde\psi(b\otimes(-)) = (\psi(b\otimes(-)))^*$ и $\cosupp \tilde\psi = (\cosupp\psi)^{\#}$. Алгебра $A^*$ становится правым $B$-модулем и, следовательно, левым $B^\mathrm{op}$-модулем.

Непосредственно проверяется, что $\widetilde{\ }$ 
является изоморфизмом категорий. Следовательно, 
$\widetilde{\ }$ отображает терминальный объект в терминальный объект, откуда и следует утверждение теоремы.
\end{proof}

Аналогичный результат справедлив и в случае кодействий:

\begin{theorem}\label{TheoremDualityCoactions}
Пусть $A$~"--- конечномерная $\Omega$-алгебра над полем $\mathbbm{k}$, а $V\subseteq \End_\mathbbm{k}(A^*)$~"--- подалгебра, содержащая тождественный оператор.
Тогда $$\mathcal{B}^\square(A^*,V)
\cong \mathcal{B}^\square(A,V^{\#})^{\mathrm{cop}}\text{\quad(изоморфизм биалгебр).}$$
где $V^{\#} := \lbrace f^* \mid f \in V \rbrace$.
\end{theorem}
\begin{proof}
Обозначим через $\widetilde{\ }$ функтор $\mathbf{Coactions}(A,V^{\#}) \to \mathbf{Coactions}(A^*,V)$,
заданный следующим образом: если $\rho \colon  A \to A \otimes B$~"--- объект
категории $\mathbf{Coactions}(A,V^{\#})$,
то $\tilde\rho \colon  A^* \to A^* \otimes B^\mathrm{cop}$ определяется через равенство
$$\tilde\rho(a^*)(a):=(a^*\otimes \id_B)\rho(a)\text{ для всех }a\in A,\ a^*\in A^*.$$
Заметим, что $\hat{\tilde\rho}(b^*\otimes(-)) = \left(\hat\rho(b^*\otimes(-))\right)^*$
для всех $b^*\in B^*$, где $\hat\rho(b^*\otimes a):=b^*(a_{(1)})a_{(0)}$.
Следовательно, $\cosupp \tilde\rho = (\cosupp\rho)^{\#}$.
Алгебра $A$ становится левым $B$-комодулем и, следовательно, правым $B^\mathrm{cop}$-комодулем.

Непосредственно проверяется, что $\widetilde{\ }$ 
является изоморфизмом категорий. Следовательно, 
$\widetilde{\ }$ отображает начальный объект в начальный объект, откуда и следует утверждение теоремы.
\end{proof}

         \section{Примеры несуществования универсальной \\ кодействующей алгебры Хопфа}
         \label{SectionUnivCoactHopfNonExist}

Как было впервые показано Д.~Тамбарой и следует, например, из следствия~\ref{CorollaryManinExistence}
и примеров~\ref{ExampleAlgebra}--\ref{ExampleCoalgebra}, 
универсальная кодействующая алгебра Хопфа на конечномерной (ко)алгебре всегда 
существует.

Однако для произвольных (ко)алгебр это не всегда так, как показывают примеры, которые строятся ниже.

\begin{example}
Пусть $A$~"--- алгебра над алгебраически замкунутым полем $\mathbbm{k}$ характеристики $0$
со счётным базисом $1_A, v_1, v_2, v_3, \ldots$, где $v_i v_j = 0$ для всех $i,j \in\mathbb N$.
Тогда не существует ни $\mathcal{B}^{\square}(A)$, ни $\mathcal{H}^{\square}(A)$.
\end{example}
\begin{proof}
Обозначим через $C_n=\langle c_n \rangle_n$ циклическую группу порядка $n$, $n\in\mathbb N$.
Рассмотрим следующее $C_n$-действие на алгебре $A$ автоморфизмами: $$c_n v_i = \left\lbrace \begin{array}{lll}
 v_{i+1} & \text{ при } &  i<n, \\
 v_1 & \text{ при } &  i=n,\\
 v_i & \text{ при } &  i>n.
  \end{array} \right.$$
  Пусть $\zeta_n$~"--- примитивный корень $n$-й степени из единицы.
  Легко видеть, что  элементы $$1_A,$$ $$w_{nj}=\sum_{i=1}^n \zeta_n^{(i-1)(j-1)} v_i, \text{ где }j=1,\ldots, n,$$
  $$v_{n+1}, v_{n+2},\ldots$$ образуют базис из собственных векторов для этого $C_n$-действия.
   
  Пусть $C^*_n = \Hom(C_n, \mathbbm{k}^\times) = \langle \chi_n \rangle_n$~"--- группа,
  двойственная к $C_n$,
  $\chi_n(c_n):= \zeta_n$. 
 Тогда  $c_n w_{nj} = \zeta_n^{1-j} w_{nj} = \chi_n^{1-j}(c_n) w_{nj}$,
  и на  
  $A$ можно задать $C^*_n$-градуировку $A=\bigoplus_{j=0}^n A^{\left(\chi_n^j \right)}$,
   причём $w_{nj} \in A^{\left(\chi_n^{1-j} \right)}$
   и $1_A, v_{n+1}, v_{n+2},\ldots \in A^{\left(1\right)}$.
   
  Заданные выше $C_n$-действие и $C^*_n$-градуировка соответствуют $\mathbbm{k}C^*_n$-комодульной структуре $\rho_n \colon A \to A \otimes \mathbbm{k}C^*_n$,
  где  $\rho_n(1_A)=1_A\otimes 1$, $\rho_n(w_{nj}) := w_{nj} \otimes \chi_n^{1-j}$ при $j=1,\ldots, n$ и $\rho_n(v_i)=v_i \otimes 1$
   при $i > n$. Заметим, что
   $v_1 = \frac{1}{n}\sum_{j=1}^n w_{nj}$, а  
      $$\rho_n(v_1)= \frac{1}{n} \sum_{j=1}^n w_{nj} \otimes \chi_n^{1-j}
      = \frac{1}{n} \sum_{j=1}^n \sum_{i=1}^n \zeta_n^{(i-1)(j-1)} v_i \otimes \chi_n^{1-j}
      = \sum_{i=1}^n v_i \otimes \chi_{ni},$$
      где $\chi_{ni} := \frac{1}{n} \sum_{j=1}^n \zeta_n^{(i-1)(j-1)} \chi_n^{1-j}$, $1\leqslant i \leqslant n$.
      Предположим теперь, что существует $\mathcal{B}^{\square}(A)$. Обозначим через $\rho \colon A \to A \otimes
      \mathcal{B}^{\square}(A)$ соответствующее кодействие.
      Тогда $\rho(v_1) = 1_A \otimes h_0 + \sum_{i=1}^m v_i \otimes h_i$
      для некоторого $m\in\mathbb N$ и некоторых $h_i \in \mathcal{B}^{\square}(A)$.
      Поскольку при фиксированном $n$ элементы $\chi_{ni}$, где $i=1,\ldots,n$, 
линейно независимы (это следует из рассуждения с использованием определителя Вандермонда), 
при $n > m$ не существует гомомоморфизма биалгебр $\varphi \colon \mathcal{B}^{\square}(A) \to
      \mathbbm{k}C^*_n$, делающего коммутативной диаграмму
      $$\xymatrix{ A \ar[r]^(0.3){\rho} \ar[rd]_{\rho_n} &  A \otimes \mathcal{B}^{\square}(A) \ar@{-->}[d]^{\id_A \otimes \varphi} \\ 
              &  A \otimes \mathbbm{k}C^*_n}$$
              Те же рассуждения можно провести и для  $\mathcal{H}^{\square}(A)$.
      \end{proof}
      
     \begin{remark}
     Аналогичные рассуждения показывают, что для $A=\mathbbm{k}[v_1, v_2, v_3, \ldots]$
     также не существует ни универсальной кодействующей алгебры Хопфа, ни  универсальной кодействующей биалгебры.
     \end{remark}
      
\begin{example}
Пусть $(C,\Delta,\varepsilon)$~"--- коалгебра над алгебраически замкунутым полем $\mathbbm{k}$ характеристики $0$
со счётным базисом $v_0, v_1, v_2, v_3, \ldots$, где
$\varepsilon(v_0)=1$, $\Delta(v_0)=v_0 \otimes v_0$,
$\Delta(v_i)=v_0 \otimes v_i + v_i\otimes v_0$, $\varepsilon(v_i) = 0$  для всех $i \in\mathbb N$.
Тогда не существует ни $\mathcal{B}^{\square}(C)$, ни $\mathcal{H}^{\square}(C)$.
\end{example}
\begin{proof} Нужно использовать $C_n$-действие
$$c_n v_i = \left\lbrace \begin{array}{lll}
 v_{i+1} & \text{ при } &  i<n, \\
 v_1 & \text{ при } &  i=n,\\
 v_i & \text{ при } &  i>n \text{ или } i=0,
  \end{array} \right.$$
  а затем дословно повторить доказательство из предыдущего примера.
\end{proof}

\newpage

\chapter{Эквивалентность (ко)действий для некоторых классов (ко)модульных алгебр} \label{ChapterHequiv}

В данной главе понятие эквивалентности (см. определение в \S\ref{SubsectionLinearMapsMod} и \S\ref{SubsectionLinearMapsComod}) рассматривается для важнейших классов (ко)модульных структур:
действий кокоммутативных алгебр, рациональных действий связных аффинных алгебраических групп, групповых градуировок и расширений Хопфа--Галуа. Оказывается, например, что любая комодульная структура, эквивалентная некоторой групповой градуировке также сводится к некоторой градуировке, а для расширений Хопфа--Галуа универсальные алгебры Хопфа совпадают с исходными (ко)действующими алгебрами.

Результаты главы вошли в статью~\cite{ASGordienko21ALAgoreJVercruysse}.

    \section{Универсальные алгебры Хопфа комодульных структур, индуцированных градуировками}
    
    В теореме~\ref{TheoremUniversalHopfOfAGradingIsJustUniversalGroup}
    мы показываем, что универсальная алгебра Хопфа комодульной структуры (см. определение~\ref{DefinitionUnivHopfComod}), индуцированной градуировкой~"--- это в точности групповая алгебра универсальной группы градуировки (см. \S\ref{SectionGradEquivUnivGroup}).
 
\begin{theorem}\label{TheoremUniversalHopfOfAGradingIsJustUniversalGroup}
Пусть $\Gamma \colon A = \bigoplus_{g\in G} A^{(g)}$~"--- градуировка на алгебре $A$ над полем $\mathbbm{k}$ группой $G$. Обозначим через $\rho \colon A \to A \otimes \mathbbm{k}G$ отображение, задающее на $A$ соответствующую структуру $\mathbbm{k}G$-комодуля.
Пусть $G_\Gamma$~"--- универсальная группа градуировки $\Gamma$, а $\rho_\Gamma \colon A \to A \otimes \mathbbm{k}G_\Gamma$~"--- отображение, задающее на $A$  соответствующую структуру $\mathbbm{k}G_\Gamma$-комодуля.
Тогда $\mathbbm{k}G_\Gamma$~"--- это универсальная алгебра Хопфа комодульной структуры $\rho$,
причём соответствующее кодействие алгебры Хопфа $\mathbbm{k}G_\Gamma$ задаётся отображением $\rho_\Gamma$.
\end{theorem}
\begin{proof} Обратимся к конструкции универсальной коизмеряющей алгебры
$\mathcal{A}^\square(A,A,\cosupp \rho)$ из теоремы~\ref{TheoremUnivComeasExistence}.
В качестве $W=W_0$ можно взять векторное пространство $\supp\rho$, которое является линейной оболочкой множества $\supp\Gamma$ в алгебре $\mathbbm{k}G$.
Выберем базис $(a_\alpha)_\alpha$ алгебры $A$ однородным относительно градуировки.
В этом случае $p_{\alpha\beta}=0$ при $\alpha\ne \beta$
и всякий элемент $p_{\alpha\alpha}$ является групповым элементом, отвечающим однородной
компоненте, содержащей $a_\alpha$.
Отсюда порождающие идеала $I$ алгебры $T(W_0)$ имеют вид
\begin{equation*} \sum_{u} a^{u}_{\alpha\beta}i_{W_0} (p_{\gamma u}) - \sum_{r,q} a^\gamma_{rq} i_{W_0}(p_{r\alpha})i_{W_0}(p_{q\beta})=  a^{\gamma}_{\alpha\beta}
\bigl(i_{W_0} (p_{\gamma \gamma}) - i_{W_0}(p_{\alpha\alpha})i_{W_0}(p_{\beta\beta})
\bigr)
\end{equation*}  для всевозможных значений индексов $\alpha,\beta,\gamma$.
(Напомним, что $a^{\gamma}_{\alpha\beta}\in \mathbbm{k}$~"--- структурные константы алгебры $A$.)
Отсюда $\mathcal{A}^\square(A,A,\cosupp \rho)=T(W_0)/I$~"--- полугрупповая алгебра моноида, порождённого элементами
носителя $\supp\Gamma$ градуировки и теми же соотношениями, которым удовлетворяют
порождающие группы $G_\Gamma$. Из конструкции, использованной в теореме~\ref{TheoremBsquareBialgebra},
следует, что $\mathcal{B}^\square(A,\cosupp \rho) := \mathcal{A}^\square(A,A,\cosupp \rho)$
как биалгебра изоморфна полугрупповой биалгебре моноида, заданного порождающими и соотношениями
группы $G_\Gamma$. В свою очередь, из определения функтора $H_l$ (см. \S\ref{SectionFunctorsHlHr})
следует, что $\mathcal{H}^\square(A,\cosupp \rho)=H_l(\mathcal{B}^\square(A,\cosupp \rho))
\cong \mathbbm{k}G_\Gamma$.
\end{proof}

\begin{theorem}\label{TheoremGradingCanBeEquivalentToAGradingOnly}
Пусть $\rho \colon A \to A \otimes H$~"--- $H$-комодульная структура на алгебре $A$,
где $H$~"--- алгебра Хопфа над полем $\mathbbm{k}$, причём $\rho$ эквивалентна некоторой групповой градуировке.
Тогда существует подалгебра Хопфа $H_1 \subseteq H$, изоморфная некоторой групповой алгебре, такая
что $\rho(A)\subseteq A \otimes H_1$.
\end{theorem}
\begin{proof} Для того, чтобы вывести эту теорему из теоремы~\ref{TheoremUniversalHopfOfAGradingIsJustUniversalGroup},
достаточно рассмотреть гомоморфный образ $H_1$ групповой алгебры $\mathcal{H}^\square(A,\cosupp \rho)$
в $H$ и использовать тот факт, что любой гомоморфизм алгебр Хопфа переводит
группоподобные элементы в группоподобные.
\end{proof}

\begin{remark}
Результаты этого параграфа естественным образом обобщаются на $\Omega$-алгебры,
если ввести понятие $G$-градуированной $\Omega$-алгебры как $\mathbbm{k}G$-комодульной алгебры для группы $G$.
Кроме того, рассуждения, приведённые выше, показывают, что если рассматривать градуировки моноидами и по аналогии с универсальной группой градуировки ввести понятие универсального моноида градуировки, то
универсальная биалгебра комодульной структуры, отвечающей моноидной градуировке,~"--- это в точности полугрупповая биалгебра универсального моноида.
\end{remark}

    \section{Универсальные алгебры Хопфа расширений Хопфа--Галуа}
  \subsection{Комодульные расширения Хопфа--Галуа}
      
    Пусть $H$~"--- алгебра Хопфа, $A$~"--- ненулевая $H$-комодульная алгебра с 1,
    а $\rho \colon A \to A\otimes H$~"--- отображение, задающее на $A$ структуру такой алгебры. 
Обозначим через $A^{\mathrm{co}H}$ подалгебру \textit{коинвариантов}: 
$$A^{\mathrm{co}H} := \lbrace a\in A \mid \rho(a)=a\otimes 1_H \rbrace.$$

Напомним, что алгебра $A$ называется \textit{расширением Хопфа--Галуа} алгебры $A^{\mathrm{co}H}$,
если линейное отображение $\mathsf{can} \colon A \mathbin{\otimes_{A^{\mathrm{co}H}}} A
\to A \mathbin{\otimes} H$, определённое через
$$\mathsf{can}(a\otimes b) := 
ab_{(0)}\otimes b_{(1)},$$
 биективно.

В теореме ниже вычисляется универсальная алгебра Хопфа расширения Хопфа--Галуа. 

\begin{theorem}\label{TheoremHComodHopfGaloisUnivHopfAlg}
Пусть $A/A^{\mathrm{co}H}$~"--- расширение Хопфа--Галуа. Тогда $H$~"---
универсальная алгебра Хопфа для $\rho$, причём кодействие алгебры Хопфа $H$ как универсальной алгебры Хопфа совпадает с $\rho$.
\end{theorem}
\begin{proof}
Из сюрьективности отображения $\mathsf{can}$ следует, что $C(\rho)=H$. Следовательно,
для всякой $H_1$-комодульной структуры $\rho_1$ на $A$,
которая эквивалентна $\rho$, где $H_1$~"--- некоторая алгебра Хопфа, существует единственный гомоморфизм коалгебр $\tau \colon H \to H_1$,
такой, что \begin{equation*}\xymatrix{ A  \ar[r]^(0.4){\rho} \ar[rd]_{\rho_1}& A \otimes H \ar[d]^{\id_A \otimes \tau} \\
& A\otimes H_1}
\end{equation*}

Единственное, что нужно теперь доказать,~"--- это то,  что $\tau$~"--- гомоморфизм алгебр Хопфа.

Во-первых, в силу теоремы~\ref{TheoremHopfCoactionAlgebraUnital}
структуры комодульных алгебр с единицей могут быть эквивалентны
только комодульным алгебрам с единицей, откуда $\rho_1(1_A)=1_A \otimes 1_{H_1}$.
В то же время $\rho_1(1_A)=(\id_A \otimes \tau)\rho(1_A)=1_A \otimes \tau(1_H)$.
Следовательно, $\tau(1_H)=1_{H_1}$.

Во-вторых, для любых $a,b\in A$ имеем \begin{equation}\label{EqHopfGaloisPart0}a_{(0)}b_{(0)}\otimes \tau(a_{(1)}b_{(1)}) = \rho_1(ab)=
\rho_1(a)\rho_1(b)=a_{(0)}b_{(0)}\otimes \tau(a_{(1)})\tau(b_{(1)}).\end{equation}
Докажем, что
\begin{equation}\label{EqHopfGaloisPart1} a b_{(0)}\otimes \tau(h b_{(1)}) = a b_{(0)}\otimes \tau(h) \tau(b_{(1)})
\text{ для всех }a, b\in A,\ h\in H.\end{equation}
Выберем базис $(h_\alpha)_\alpha$ в $H$ и рассмотрим произвольный элемент $h_\beta$ этого базиса.
В силу сюръективности отображения $\mathsf{can}$ существуют
$a_i, b_i \in A$, такие, что $a \otimes h_\beta = \sum_i a_i b_{i(0)} \otimes b_{i(1)}$.
При этом
$\rho(b_i)= \sum_\alpha b_{i\alpha} \otimes h_\alpha$ для некоторых $b_{i\alpha} \in A$,
где для любого $i$ только конечное число элементов $b_{i\alpha}$ ненулевые.
В этих обозначениях $a \otimes h_\beta = \sum_{i,\alpha} a_i b_{i\alpha} \otimes h_\alpha$
и \begin{equation}\label{EqHopfGaloisSokr}\sum_i a_i b_{i\alpha} = \left\lbrace\begin{array}{rrr} a & \text{ при } &  \alpha = \beta, \\ 
0 & \text{ при } &  \alpha \ne \beta.  \end{array}\right.\end{equation}
Следовательно,
\begin{equation*}\begin{split}  a b_{(0)}\otimes \tau(h_\beta) \tau(b_{(1)})
= \sum_{i,\alpha} a_i b_{i\alpha} b_{(0)} \otimes \tau(h_\alpha) \tau(b_{(1)})
\stackrel{(\ref{EqHopfGaloisPart0})}{=} \\ =
\sum_{i,\alpha} a_i b_{i\alpha} b_{(0)} \otimes \tau(h_\alpha b_{(1)})
\stackrel{(\ref{EqHopfGaloisSokr})}{=} 
a b_{(0)} \otimes \tau(h_\beta b_{(1)}).
\end{split}\end{equation*}
Таким образом, равенство~(\ref{EqHopfGaloisPart1}) справедливо в случае $h=h_\beta$.
В силу того, что $\beta$ было произвольным индексом и обе части равенства~(\ref{EqHopfGaloisPart1})
линейны по $h$, равенство (\ref{EqHopfGaloisPart1}) справедливо и для произвольного $h$.

Рассмотрим произвольный элемент $h\in H$ и некоторый элемент базиса $h_\gamma$.
Снова используя сюрьективность $\mathsf{can}$, получаем, что существуют
$c_i, d_i \in A$, такие, что $1_A \otimes h_\gamma = \sum_i c_i d_{i(0)} \otimes d_{i(1)}$.
При этом $\rho(d_i)=\sum_{\alpha} d_{i\alpha} \otimes h_\alpha$
для некоторых $d_{i\alpha} \in A$. В этих обозначениях $1_A \otimes h_\gamma = \sum_{i,\alpha} c_i d_{i\alpha} \otimes h_\alpha$
и \begin{equation}\label{EqHopfGaloisSokr2}\sum_i c_i d_{i\alpha} = \left\lbrace\begin{array}{rrr} 1_A & \text{ при } &  \alpha = \gamma, \\ 
0 & \text{ при } &  \alpha \ne \gamma.  \end{array}\right.\end{equation}
Поэтому \begin{equation*}\begin{split} 1_A \otimes \tau(h)\tau(h_\gamma)= \sum_{i,\alpha} c_i d_{i\alpha} \otimes \tau(h)\tau(h_\alpha)
\stackrel{(\ref{EqHopfGaloisPart1})}{=} \\ = \sum_{i,\alpha} c_i d_{i\alpha} \otimes \tau(h h_\alpha)
\stackrel{(\ref{EqHopfGaloisSokr2})}{=} 1_A \otimes \tau(h h_\gamma).
\end{split}\end{equation*}
Следовательно, $\tau(h h_\gamma)=\tau(h)\tau(h_\gamma)$. В силу того, что индекс $\gamma$ 
был произвольным, $\tau$~"--- гомоморфизм биалгебр, откуда $\tau$ является и гомоморфизмом алгебр Хопфа (см., например, \cite[предложение 4.2.5]{Danara}). Как следствие, $H$~"--- универсальная алгебра Хопфа комодульной структуры $\rho$.
\end{proof}

\begin{remark}
В доказательстве теоремы~\ref{TheoremHComodHopfGaloisUnivHopfAlg}
мы использовали только сюръективность $\mathsf{can}$.
\end{remark}

Если мы рассмотрим стандартную $G$-градуировку на групповой алгебре $\mathbbm{k}G$ группы $G$,
тогда универсальная группа этой градуировки изоморфна $G$. В случае комодульных структур справедлив аналогичный результат:

\begin{corollary}\label{CorollaryUniversalHopfOfCoactionOnItselfIsItself}
Пусть $H$~"--- алгебра Хопфа.
Тогда универсальная алгебра Хопфа $H$-комодульной структуры на алгебре $H$,
заданной коумножением $\Delta \colon H \to H \otimes 
H$ совпадает с $H$.
\end{corollary}
\begin{proof}
Заметим, что $H^{\mathrm{co}H}  \cong \mathbbm{k}$, а $H/\mathbbm{k}$~"--- расширение Хопфа--Галуа в силу \cite[примеры 6.4.8, 1)]{Danara}. Поэтому требуемое утверждение сразу же вытекает из теоремы~\ref{TheoremHComodHopfGaloisUnivHopfAlg}.
\end{proof}

\subsection{Модульные расширения Хопфа--Галуа}

В данном пункте доказывается результат, аналогичный теореме~\ref{TheoremHComodHopfGaloisUnivHopfAlg}, но для
$H$-модульных алгебр.
Пусть $A$~"--- ненулевая $H$-модульная алгебра с $1$ для некоторой алгебры Хопфа $H$.
Обозначим через $\psi \colon H \otimes A \to A$ соответствующее отображение. 
Пусть $A^H$~"--- подалгебра \textit{инвариантов}, т.е. 
$A^H := \lbrace a\in A \mid ha=\varepsilon(h)a \text{ для всех }
h\in H \rbrace$. Говорят, что алгебра $A$~"--- \textit{расширение Хопфа--Галуа} алгебры $A^H$,
если линейное отображение $\mathsf{can} \colon A \mathbin{\otimes_{A^H}} A
\to \Hom_\mathbbm{k}(H, A)$, заданное при помощи равенства $$\mathsf{can}(a\otimes b)(h) := 
a(hb),$$ инъективно и имеет плотный образ в конечной топологии.

\begin{theorem}\label{TheoremHModHopfGaloisUnivHopfAlg}
Пусть $A/A^H$~"--- \textit{расширение Хопфа--Галуа}. Тогда $H$~"--- универсальная алгебра Хопфа модульной структуры $\psi$, причём действие алгебры Хопфа $H$ как универсальной алгебры Хопфа совпадает с $\psi$.
\end{theorem}
\begin{proof}
Обозначим через $\zeta \colon H \to \End_\mathbbm{k}(A)$ гомоморфизм, заданный равенством
$\zeta(h)(a)=\psi(h\otimes a)$ для всех $h\in H$, $a\in A$. 

Из плотности образа отображения $\mathsf{can}$ следует, что для всякого $h\ne 0$
существуют элементы $a,b\in A$, такие, что $a(hb)\ne 0$. В частности, $\ker\zeta=0$.
Поэтому для любой алгебры Хопфа $H_1$ и любой $H_1$-модульной структуры $\psi_1$ на $A$
эквивалентной $\psi$ при помощи $\id_A$
существует единственный гомоморфизм алгебр $\tau \colon H_1 \to H$, такой,
что диаграмма
\begin{equation*}\xymatrix{ H  \ar[r]^(0.35){\zeta} & 
\End_\mathbbm{k}(A)  \\
H_1 \ar[ru]_(0.4){\zeta_1} \ar[u]^{\tau} & }
\end{equation*}
коммутативна. (Здесь $\zeta_1 \colon H \to \End_\mathbbm{k}(A)$~"--- гомоморфизм, заданный равенством
$\zeta_1(h)(a)=\psi_1(h\otimes a)$ для всех $h\in H$, $a\in A$.)
Отсюда $\tau$~"--- единственный гомоморфизм алгебр,
 такой,
что диаграмма
\begin{equation*}\xymatrix{ H\otimes A \ar[r]^(0.6){\psi} & 
A  \\
H_1 \otimes A \ar[ru]_{\psi_1} \ar[u]^{\tau \otimes \id_A} & }
\end{equation*}
коммутативна.

Остаётся доказать, что $\tau$~"--- гомоморфизм алгебр Хопфа.

Прежде всего заметим, что, в силу того, что $A$~"--- $H$-модульная
алгебра с $1$, из теоремы~\ref{TheoremHopfActionAlgebraUnital}
следует, что $A$ также  $H_1$--модульная
алгебра с $1$.
Поэтому
 $\varepsilon(h)
1_A = h1_A=\tau(h)1_A=\varepsilon(\tau(h))1_A$
для всех $h\in H_1$, откуда $\varepsilon(h)= \varepsilon(\tau(h))$.

Выберем базис $(h_\alpha)_\alpha$ в алгебре $H$.
Докажем, что для всякого набора элементов поля $\lambda_{\alpha\beta}\in \mathbbm{k}$,
такого, что только конечное число этих элементов ненулевые,
условие \begin{equation}\label{EqLambdaAlphaBetaHopfGaloisModUniv}\sum_{\alpha,\beta} \lambda_{\alpha\beta}(h_\alpha a)
(h_\beta b) = 0\text{ для всех }a,b\in A\end{equation}
влечёт $\lambda_{\alpha\beta}=0$ для всех $\alpha,\beta$.

Действительно, пусть выполнено равенство~\eqref{EqLambdaAlphaBetaHopfGaloisModUniv}.
Пусть $\Lambda$~"--- конечное множество индексов, такое, что $\lambda_{\alpha\beta}=0$,
если $\alpha \notin \Lambda$ или $\beta \notin \Lambda$.
 В силу плотности образа отображения $\mathsf{can}$ для любого $\gamma$ существуют $a_{\gamma i},b_{\gamma i} \in A$, такие, что
$$\sum_i a_{\gamma i}(h_\alpha b_{\gamma i})=\left\lbrace
\begin{array}{crl} 1_A & \text{ при } &  \alpha = \gamma,\\
0 & \text{ при } & \alpha \ne \gamma \text{ и } \alpha\in \Lambda.
\end{array} \right.$$
Тогда из~\eqref{EqLambdaAlphaBetaHopfGaloisModUniv}
следует, что для всех
$b\in A$ и всех $\gamma$
выполнено $$ \sum_{\beta\in\Lambda} \lambda_{\gamma\beta} h_\beta b =
\sum_{\substack{\alpha,\beta\in \Lambda, \\ i}} \lambda_{\alpha\beta}a_{\gamma i}(h_\alpha b_{\gamma i})
(h_\beta b) = 0. $$
Используя условие $\ker\zeta=0$, получим $\sum_\beta \lambda_{\gamma\beta} h_\beta = 0$
и $\lambda_{\gamma\beta}=0$ для всех $\beta,\gamma$.

В силу \eqref{EqModCompat} для всех $h\in H_1$ и $a,b \in A$
справедливы равенства $$(\tau(h_{(1)})a)(\tau(h_{(2)})b)=(h_{(1)}a)(h_{(2)}b)=h(ab)=\tau(h)(ab)=
(\tau(h)_{(1)}a)(\tau(h)_{(2)}b).$$
Из~\eqref{EqLambdaAlphaBetaHopfGaloisModUniv}
следует, что $\Delta(\tau(h))=(\tau\otimes \tau)\Delta(h)$
для всех $h\in H_1$ и $\tau$~"--- гомоморфизм алгебр Хопфа.
Поэтому $H$~"--- действительно универсальная алгебра Хопфа модульной структуры $\psi$.
\end{proof}
\begin{remark}
В доказательстве теоремы~\ref{TheoremHModHopfGaloisUnivHopfAlg} использовалась только плотность образа отображения $\mathsf{can}$.
\end{remark}

\section{Универсальные кокоммутативные алгебры Хопфа}
 \label{SectionActionsCocommHopfAlgebras}

В отличие от алгебры Хопфа $\mathcal{H}^\square (A,V)$,
которая задаётся своими порождающими и определяющими соотношениями (см. теорему~\ref{TheoremUnivComeasExistence} и \S\ref{SectionCoactions}), структура алгебры Хопфа ${}_\square \mathcal{H}(A,V)$ не так прозрачна. В связи с этим возникает необходимость разработки альтернативных методов вычисления
$V$-универсальных действующих алгебр Хопфа. Одним из таких способов является рассмотрение
действий кокоммутативных алгебр Хопфа. Оказывается, что для таких
действий существует своя $V$-универсальная кокоммутативная действующая алгебра Хопфа ${}_\square \mathcal{H}(A,V)_\mathrm{coc}$,
причём если удаётся доказать (см., в частности, пример~\ref{ExampleDoubleNumbersUniversal}), что универсальная алгебра Хопфа ${}_\square \mathcal{H}(A,\cosupp \rho)$ для заданного кодействия $\rho$ кокоммутативна, то это означает, что ${}_\square \mathcal{H}(A,\cosupp \rho)$ совпадает с универсальной кокоммутативной алгеброй Хопфа ${}_\square \mathcal{H}(A,\cosupp \rho)_\mathrm{coc}$,
которая вычисляется ниже в теореме~\ref{TheoremUniversalCocommutative}.
 
Для коалгебры $C$ обозначим через $C_{\mathrm{coc}}$ сумму всех её кокоммутативных подкоалгебр.
Легко видеть, что если $H$~"--- алгебра Хопфа, то $H_{\mathrm{coc}}$~"--- её подалгебра Хопфа.
Отсюда, если $A$~"--- 
$\Omega$-алгебра над полем $\mathbbm{k}$, а $V\subseteq \End_\mathbbm{k}(A)$~"--- подалгебра, содержащая тождественный оператор,
то ограничение $\psi^{\mathbf{Hopf}}_{A,V,\mathrm{coc}}$ действия
$\psi^{\mathbf{Hopf}}_{A,V}$
 подалгебры Хопфа ${}_\square \mathcal{H}(A,V)_\mathrm{coc} \subseteq {}_\square \mathcal{H}(A,V)$ удовлетворяет тому же универсальному свойству, что и действие 
алгебры Хопфа ${}_\square \mathcal{H}(A,V)$, но только среди действий кокоммутативных
алгебр Хопфа. В связи с этим назовём алгебру Хопфа ${}_\square \mathcal{H}(A,V)_\mathrm{coc}$
\textit{$V$-универсальной кокоммутативной
 действующей алгеброй Хопфа} на алгебре $A$.\label{NotationsquareHAVcoc}

 Если $\psi \colon H \otimes A \to A$~"--- модульная структура, то 
 алгебру Хопфа ${}_\square \mathcal{H}(A,\cosupp\psi)_\mathrm{coc}$
 будем называть
\textit{универсальной кокоммутативной
 алгеброй Хопфа} модульной структуры $\psi$.

В случае алгебраически замкнутого поля характеристики $0$ теорема
Габриэля--Картье--Костанта--Милнора--Мура (см. теорему~\ref{TheoremCartierGabrielKostant}) позволяет дать конкретное описание алгебры Хопфа ${}_\square \mathcal{H}(A,V)_\mathrm{coc}$. (Определение смэш-произведения см. в \S\ref{SubsectionModComodAlg}.)

\begin{theorem}\label{TheoremUniversalCocommutative} Пусть $A$~"--- 
$\Omega$-алгебра над алгебраически замкнутым полем $\mathbbm{k}$ 
характеристики $0$, а $V\subseteq \End_\mathbbm{k}(A)$~"--- подалгебра, содержащая тождественный оператор.
Пусть $$G := \mathcal U(V) \cap \Aut(A)\text{ и }L :=  V  \cap \Der(A),$$
где $\mathcal U(V)$~"--- группа обратимых элементов алгебры $V$,
а $\Der(A)$~"--- алгебра Ли дифференцирований алгебры $A$, рассматриваемая как подпространство пространства $\End_\mathbbm{k}(A)$.
Тогда ${}_\square \mathcal{H}(A,V)_\mathrm{coc} \cong U(L) \mathbin{\#} \mathbbm{k}G$ (группа $G$ действует на $L$ сопряжениями),
причём действие $\psi$ алгебры Хопфа $U(L) \mathbin{\#} \mathbbm{k}G$ на $A$, которое удовлетворяет универсальному свойству, индуцировано стандартными действиями
группы $G$ и алгебры Ли $L$.
\end{theorem}
\begin{proof} 
Обозначим через $\psi$ действие алгебры Хопфа $U(L) \mathbin{\#} \mathbbm{k}G$ на $A$, индуцированное стандартными действиями группы $G$ и алгебры Ли $L$.
Операторы из $G$ и $L$ являются элементами подпространства $V$,
откуда $\cosupp \psi \subseteq V$.

Пусть $\psi_1 \colon H_1 \otimes A \to A$~"--- действие некоторой кокоммутативной 
алгебры Хопфа $H_1$, такое, что $\cosupp \psi_1 \subseteq V$. В силу теоремы~\ref{TheoremCartierGabrielKostant} существует изоморфизм  $H_1 \cong U(L_1) \mathbin{\#} \mathbbm{k}G_1$, где $G_1$~"---
группа группоподобных элементов алгебры Хопфа $H_1$, а $L_1$~"--- алгебра Ли
примитивных элементов алгебры Хопфа $H_1$.
Требуется показать, что существует единственный гомоморфизм алгебр Хопфа $\tau \colon H_1 \to U(L) \mathbin{\#} \mathbbm{k}G$,
такой, что $\psi(\tau \otimes \id_A) =\psi_1$. Снова определим гомоморфизмы алгебр $\zeta \colon U(L) \mathbin{\#} \mathbbm{k}G \to \End_\mathbbm{k}(A)$
и $\zeta_1 \colon H_1 \to \End_\mathbbm{k}(A)$ через равенства $\zeta(h)(a):=\psi(h\otimes a)$
и $\zeta_1(h_1)(a):=\psi_1(h_1\otimes a)$ при $a\in A$, $h\in U(L) \mathbin{\#} \mathbbm{k}G$ и $h_1 \in H_1$.
Тогда равенство $\psi(\tau \otimes \id_A) =\psi_1$ эквивалентно равенству $\zeta\tau =\zeta_1$.

В силу предложения~\ref{PropositionPrimitiveGrouplikeSmash}
алгебра Ли примитивных элементов алгебры Хопфа $U(L)\mathbin{\#} \mathbbm{k}G$ 
совпадает с $L \otimes 1$ (мы отождествим её с $L$), а группа группоподобных элементов
алгебры Хопфа $H$ совпадает с $1 \otimes G$ (мы отождествим её с $G$).
Поскольку гомоморфизм $\tau$ (если он действительно существует),
обязан отображать 
группоподобные элементы в группоподобные, а примитивные в примитивные,
то должны быть справедливы равенства $\tau(G_1) \subseteq G$ и $\tau(L_1) \subseteq L$.
Заметим, что $\zeta\bigr|_G = \id_G$ и
$\zeta\bigr|_L = \id_L$. Следовательно, $\tau\bigr|_{G_1}$ and $\tau\bigr|_{L_1}$
однозначно определены ограничениями $\zeta_1\bigr|_{G_1} \colon G_1 \to G$ и $\zeta_1\bigr|_{L_1}  \colon L_1 \to L$.
Поскольку $H_1$ как алгебра порождена подмножествами $G_1$ и $L_1$,
существует не более одного гомоморфизма $\tau \colon H_1 \to U(L)\mathbin{\#} \mathbbm{k}G$
алгебр Хопфа, такого, что $\zeta_0\tau =\zeta_1$.
Теперь осталось определить $\tau$ как гомоморфизм алгебр Хопфа индуцированный $\zeta_1\bigr|_{G_1}$ and $\zeta_1\bigr|_{L_1}$.
\end{proof}

Докажем теперь критерий того, что универсальная кокоммутативная и универсальная алгебры Хопфа сопадают:

\begin{theorem}\label{TheoremCocommUnivIsUniv}
Пусть $A$~"--- $\Omega$-алгебра над произвольным полем $\mathbbm{k}$, а $V\subseteq \End_\mathbbm{k}(A)$~"--- подалгебра, содержащая тождественный оператор.
Тогда $ {}_\square \mathcal{H}(A,V) = {}_\square \mathcal{H}(A,V)_\mathrm{coc}$,
если и только если  любое действие $\psi_1 \colon H_1 \otimes A \to A$, где $H_1$~--- алгебра Хопфа,
 а $\cosupp\psi_1 \subseteq V$,
 пропускается через действие $\psi_2$ некоторой кокоммутативной алгебры Хопфа $H_2$,
т.е. существует такой гомоморфизм $\theta$ алгебр Хопфа,
что диаграмма ниже становится коммутативной:

\begin{equation*}\xymatrix{ H_1 \otimes A  \ar[r]^(0.6){\psi_1} \ar@{-->}[d]_{\theta \otimes \id_A} & A  \\
H_2 \otimes A \ar[ru]_{\psi_2} &    
}
\end{equation*}
\end{theorem}
\begin{proof} Необходимость следует из определения универсальной действующей алгебры Хопфа.

Докажем достаточность.
Предположим, что всякое действие алгебры Хопфа с коносителем,
являющимся подалгеброй алгебры $V$, пропускается через  действие
некоторой кокоммутативной алгебры Хопфа $H_2$.
Без ограничения общности можно считать, что соответствующий гомоморфизм $\theta$
сюръективен и, следовательно, коноситель действия алгебры Хопфа $H_2$
является подалгеброй алгебры $V$.
Поскольку всякая такая алгебра Хопфа $H_2$ кокоммутативна,
это означает, что всякое действие алгебры Хопфа с коносителем,
являющимся подалгеброй алгебры $V$, пропускается через $\psi^{\mathbf{Hopf}}_{A,V,\mathrm{coc}}$.

Рассмотрим теперь действие $\psi^{\mathbf{Hopf}}_{A,V} \colon {}_\square \mathcal{H}(A,V) \otimes A \to A$ 
\quad $V$-универсальной действующей алгебры Хопфа $\mathcal{H}(A,V)$.
В силу универсального свойства действия $\psi^{\mathbf{Hopf}}_{A,V}$ и того факта,
что $\psi^{\mathbf{Hopf}}_{A,V}$ пропускается
через $\psi^{\mathbf{Hopf}}_{A,V,\mathrm{coc}}$, следует существование таких
гомоморфизмов $\theta_1, \theta_2$ алгебр Хопфа, что диаграмма ниже становится коммутативной:
\begin{equation*}\xymatrix{ \mathcal{H}(A,V)_\mathrm{coc}\otimes A  \ar[rr]^(0.7){\psi^{\mathbf{Hopf}}_{A,V,\mathrm{coc}}} \ar@<0.5ex>@{-->}[d]^{\theta_1 \otimes A} && A  \\
\mathcal{H}(A,V)\otimes A  \ar[rru]_{\psi^{\mathbf{Hopf}}_{A,V}} \ar@<0.5ex>@{-->}[u]^{\theta_2 \otimes A} &&    
}
\end{equation*}

Из единственности соответствующих отображений в определениях
действий  $\psi^{\mathbf{Hopf}}_{A,V}$ и $\psi^{\mathbf{Hopf}}_{A,V,\mathrm{coc}}$
следует, что $\theta_1 \theta_2 = \id_{\mathcal{H}(A,V)}$
и $\theta_2 \theta_1 = \id_{\mathcal{H}(A,V)_\mathrm{coc}}$.
В частности, алгебра Хопфа $\mathcal{H}(A,V)$
кокоммутативна и $\mathcal{H}(A,V) = \mathcal{H}(A,V)_\mathrm{coc}$.
\end{proof}

В примере~\ref{ExampleCocommUnivDifferent} мы приводим
пример действия $H \otimes A \to A$,
которое хотя и эквивалентно действию $\psi$ некоторой кокоммутативной алгебры Хопфа, однако
для некоторых $h\in H$ и $a,b \in A$ справедливо неравенство
  $(h_{(1)}a)(h_{(2)}b)\ne (h_{(2)}a)(h_{(1)}b)$. Ясно, что такое $H$-действие
не пропускается через действие никакой кокоммутативной алгебры Хопфа,
поэтому в данном случае ${}_\square \mathcal{H}(A,\cosupp\psi) \ncong {}_\square \mathcal{H}(A,\cosupp\psi)_\mathrm{coc}$.

\begin{example}\label{ExampleCocommUnivDifferent} Пусть $\mathbbm{k}$~"--- поле, а $A=\mathbbm{k}1_A\oplus \mathbbm{k}a\oplus \mathbbm{k}b \oplus \mathbbm{k}ab$~"--- ассоциативная $\mathbbm{k}$-алгебра, где
$a^2=b^2=ba=0$. Пусть $S_3$~"--- симметрическая группа, действующая на элементах $\lbrace 1,2,3 \rbrace$. 
Рассмотрим эквивалентные $S_3$-  и $\mathbb Z/4\mathbb Z$-градуировки на $A$,
заданные равенствами
$$A^{(\mathrm{id})}:=A^{(\bar 0)}:= \mathbbm{k} 1_A,\
A^{\bigl((12) \bigr)}:=A^{(\bar 1)}:= \mathbbm{k} a,$$
$$A^{\bigl((23) \bigr)}:=A^{(\bar 2)}:= \mathbbm{k} b,\
A^{\bigl((123) \bigr)}:=A^{(\bar 3)}:= \mathbbm{k} ab.$$
Поскольку градуировки эквивалентны, согласно теореме~\ref{TheoremGradEquivCriterion} соответствующие  $(\mathbbm{k}S_3)^*$- и $(\mathbbm{k}(\mathbb Z/4\mathbb Z))^*$-действия также эквивалентны,
в то время как $(\mathbbm{k}(\mathbb Z/4\mathbb Z))^*$~"--- коммутативная 
кокоммутативная алгебра Хопфа, а $(\mathbbm{k}S_3)^*$~"--- коммутативная 
некокоммутативная алгебра Хопфа, и
существует элемент $h\in (\mathbbm{k}S_3)^*$, такой, что
$(h_{(1)}a)(h_{(2)}b) \ne (h_{(2)}a)(h_{(1)}b)$.
\end{example}
\begin{proof} Для конечной группы $G$ обозначим через $(h_g)_{g\in G}$ базис алгебры $(\mathbbm{k}G)^*$,
двойственный к базису $(g)_{g\in G}$. Тогда $\Delta h_g = \sum\limits_{st=g} h_s \otimes h_t$.
В частности,
\begin{equation*}\begin{split}((h_{(123)})_{(1)}a)((h_{(123)})_{(2)}b)
=\sum_{\sigma\rho=(123)} (h_\sigma a)(h_\rho b) = ab  \ne 
\\ ((h_{(123)})_{(2)}a)((h_{(123)})_{(1)}b) = \sum_{\sigma\rho=(123)} (h_\rho a)(h_\sigma b)=0.\end{split}
\end{equation*}
\end{proof}

Напомним, что в теореме~\ref{TheoremUniversalHopfOfAGradingIsJustUniversalGroup}
было доказано, что универсальная алгебра Хопфа кодействия, соответствующего
групповой градуировке, является групповой алгеброй универсальной группы градуировки.
В примере~\ref{ExampleDoubleNumbersUniversal} ниже мы доказываем, что
для групповых действий аналогичный результат неверен, т.е. что универсальная алгебра Хопфа
действия группы не обязана быть групповой алгеброй.
Кроме того, из примера~\ref{ExampleDoubleNumbersUniversal} видно, как понятие универсальной кокоммутативной алгебры Хопфа
может быть использовано для вычисления универсальной алгебры Хопфа.

\begin{example}\label{ExampleDoubleNumbersUniversal}
Пусть $A := \mathbbm{k}[x]/(x^2)$,
где $\mathbbm{k}$~"--- алгебраически замкнутое поле характеристики $0$. 
Пусть $C_2$~"--- циклическая группа порядка $2$ с порождающим $c$. Определим $C_2$-действие на $A$
по формуле $c \bar x = -\bar x$ и обозначим через $\psi$ соответствующее $\mathbbm{k}C_2$-действие
$\mathbbm{k}C_2 \otimes A \to A$.
Обозначим через $H$ алгебру Хопфа,
которая как векторное пространство
совпадает с $\mathbbm{k}[y]\otimes \mathbbm{k}\mathbbm{k}^\times$, а
структура алгебры и коалгебры на $H$ наследуется при помощи тензорного произведения с соответствующих
структур на $\mathbbm{k}[y]$ и $\mathbbm{k}\mathbbm{k}^\times$,
причём структура коалгебры на $\mathbbm{k}[y]$ задаётся формулами
$\Delta(y)=1 \otimes y + y \otimes 1$ и $\varepsilon(y)=0$.
Антипод $S$ и действие $H$ на $A$ определяются формулами
 $$S(y^k\otimes \lambda)=(-1)^k y^k\otimes \lambda^{-1}\text{ и }(y^k\otimes \lambda) \bar x= \lambda$$ при $k\in\mathbb Z_+$ и $\lambda \in \mathbbm{k}^\times$. 
 Тогда $H$ является универсальной алгеброй Хопфа для $\psi$.
\end{example}
\begin{proof} Обозначим через $\zeta$ гомоморфизм $\mathbbm{k}C_2 \to \End_\mathbbm{k}(A)$,
заданный равенством $\zeta(g)a:=ga$ для всех $g\in C_2$, $a\in A$.
Очевидно, что $\cosupp\psi = \zeta(\mathbbm{k}C_2)$.

Отождествим алгебру $\End_\mathbbm{k}(A)$ с алгеброй
$M_2(\mathbbm{k})$ всех матриц $2\times 2$,
фиксировав базис 
$\bar 1, \bar x$ алгебры $A$. 
Тогда $\zeta(1) = \left(\begin{smallmatrix} 1 & 0 \\
0 & 1\end{smallmatrix}\right)$ и $\zeta(c)= \left(\begin{smallmatrix} 1 & 0 \\
0 & -1\end{smallmatrix}\right)$. Поскольку $\zeta(\mathbbm{k}C_2)$ 
является линейной оболочкой
элементов $\zeta(1)$
и $\zeta(c)$, алгебра $\zeta(\mathbbm{k}C_2)$ является подалгеброй всех диагональных матриц алгебры $M_2(\mathbbm{k})$.
Более того, $\Aut(A) = \left\lbrace\left(\begin{smallmatrix} 1 & 0 \\
0 & \lambda \end{smallmatrix}\right) \mathbin{\bigl|} \lambda \in \mathbbm{k}^\times \right\rbrace$,
а $\Der(A) = \left\lbrace\left(\begin{smallmatrix} 0 & 0 \\
0 & \lambda \end{smallmatrix}\right) \mathbin{\bigl|} \lambda \in \mathbbm{k} \right\rbrace$.
Действительно, всякий автоморфизм обязан сохранять единицу и переводить максимальный нильпотентный идеал $\mathbbm{k}\bar x$ в себя. Для всякого дифференцирования $\delta$ равенство $\delta(\bar 1^2)=\delta(\bar 1)$
влечёт $\delta(\bar 1)=0$, а из равенства $\delta(\bar x^2)=0$ при $\chr \mathbbm{k} 
\ne 2$ следует, что $\bar x \delta(\bar x)=0$, откуда $\delta(\bar x) \in \mathbbm{k}\bar x$.

Следовательно, $\mathcal U(\zeta(\mathbbm{k}C_2)) \cap \Aut(A)=\Aut(A) \cong \mathbbm{k}^\times$,
а $\zeta(\mathbbm{k}C_2) \cap \Der(A)$~"--- это одномерная алгебра Ли,
чья универсальная обёртывающая изоморфна алгебре $\mathbbm{k}[y]$.
Из теоремы~\ref{TheoremUniversalCocommutative}
теперь следует, что $H$~"--- действительно универсальная кокоммутативная алгебра Хопфа
действия $\psi$.

Для того, чтобы показать, что алгебра Хопфа $H$ универсальна и как необязательно кокоммутативная
алгебра Хопфа, в силу теоремы~\ref{TheoremCocommUnivIsUniv} достаточно показать, что для любого
 действия любой алгебры Хопфа $H_1$, эквивалентного
действию $\psi$, соответствующий гомоморфизм $\zeta_1 \colon H_1 \to \End_\mathbbm{k}(A)$ пропускается через 
некоторую кокоммутативную алгебру Хопфа.

Пусть гомоморфизм $\zeta_1 \colon H_1 \to \End_\mathbbm{k}(A)$ 
отвечает действию алгебры Хопфа $H_1$, эквивалентному действию $\psi$.
Тогда $\zeta_1(H_1)$ является алгеброй всех диагональных матриц $2\times 2$.
В частности, $\bar x$~"--- это общий собственный вектор для всех операторов из $H_1$.
Определим $\varphi \in H_1^*$ через $h\bar x = \varphi(h) \bar x$ для всех $h\in H_1$. 
Тогда $\varphi \colon H_1 \to \mathbbm{k}$~"--- гомоморфизм алгебр с единицей
или, что в данном случае то же самое, группоподобный элемент алгебры Хопфа $H^\circ$. 
Обозначим через $\Delta$ и $\varepsilon$, соответственно, коумножение
и коединицу в $H_1$.
В силу предложения~\ref{PropositionHopfUnivModUnitality} имеем
 $h\bar 1 = \varepsilon(h) \bar 1$  для всех $h\in H_1$.
 
 Степени $\varphi^k$, где $k\in \mathbb Z$, элемента $\varphi$ группы $G(H^\circ)$
группоподобных элементов алгебры Хопфа $H^\circ$ могут быть определены по формуле
  $$\varphi^k(h)=\left\lbrace\begin{array}{ccc}\varphi(h_{(1)})\dots \varphi(h_{(k)})
& \text{ при } & k \geqslant 1, \\
\varepsilon(h)  & \text{ при } & k = 0,\\
 \varphi(Sh_{(1)})\dots \varphi(Sh_{(k)})
& \text{ при } & k \leqslant -1,\end{array}\right.$$
где $h\in H_1$.
Пусть $$I := \bigcap\limits_{k\in\mathbb Z} \ker \bigl(\varphi^k\bigr).$$
Поскольку все $\varphi^k \colon H_1 \to \mathbbm{k}$~"--- гомоморфизмы алгебр с единицей, $I$ является идеалом.
В силу того, что в определении идеала $I$ 
используются и отрицательные степени элемента $\varphi$,
справедливо также включение $SI\subseteq I$.

Докажем, что $I$~"--- коидеал,
т. е., что $\Delta(I) \subseteq H_1 \otimes I + I \otimes H_1$.
Для начала покажем, что \begin{equation}\label{EqH_1IIH_1}H_1 \otimes I + I \otimes H_1 = \bigcap\limits_{k,\ell\in\mathbb Z} \ker \bigl(\varphi^k \otimes \varphi^\ell\bigr).\end{equation} 
Включение 
 $$H_1 \otimes I + I \otimes H_1 \subseteq \bigcap\limits_{k,\ell\in\mathbb Z} \ker \bigl(\varphi^k \otimes \varphi^\ell\bigr)$$ очевидно.
Для того, чтобы доказать обратное включение,
выберем в $H_1$
базис $(a_\alpha)_{\alpha \in \Lambda_1 \cup \Lambda_2}$
 таким образом, чтобы $(a_\alpha)_{\alpha \in \Lambda_1}$ являлся базисом идеала $I$.
Пусть $$w:=\sum_{\alpha,\beta \in \Lambda_1 \cup \Lambda_2} \gamma_{\alpha\beta}\, a_{\alpha}\otimes a_{\beta}
\in \bigcap\limits_{k,\ell\in\mathbb Z} \ker \bigl(\varphi^k \otimes \varphi^\ell\bigr)$$
 для некоторых $\gamma_{\alpha\beta} \in \mathbbm{k}$, из которых лишь конечное число коэффициентов ненулевые.
Тогда в силу определения идеала $I$ для любого $k\in\mathbb Z$
справедливо включение $$\sum_{\alpha,\beta \in \Lambda_1 \cup \Lambda_2} \gamma_{\alpha\beta} \varphi^k (a_{\alpha}) a_{\beta} \in I,$$ откуда в силу определения элементов $a_\alpha$ для всех $\beta \in \Lambda_2$ справедливо равенство $$\sum_{\alpha \in \Lambda_1 \cup \Lambda_2} \gamma_{\alpha\beta} \varphi^k (a_{\alpha})=0,$$
т.е. $$\sum_{\alpha \in \Lambda_1 \cup \Lambda_2} \gamma_{\alpha\beta} a_{\alpha}\in I.$$
Поэтому $\gamma_{\alpha\beta} = 0$ при $\alpha,\beta \in \Lambda_2$.
Следовательно, $w\in H_1 \otimes I + I \otimes H_1$
и равенство \eqref{EqH_1IIH_1} доказано.
Поскольку для всех $k,\ell \in\mathbb Z$ и $h\in I$ справедливо равенство
  $$(\varphi^k \otimes \varphi^\ell)(h_{(1)}\otimes h_{(2)})
  =\varphi^{k+\ell}(h)=0,$$  $S$-инвариантный идеал $I$ является
  коидеалом, а значит, и идеалом Хопфа.
  
  Напомним, что в силу определения идеала $I$ справедливы включения $I\subseteq\ker\varepsilon$ и $I\subseteq \ker\varphi$. Следовательно, гомоморфизм $\zeta_1$ пропускается через $H_1/I$, т.е.
  диаграмма ниже коммутативна:
\begin{equation*}\xymatrix{ H _1 \ar[r]^(0.45){\pi} \ar[rd]_{\zeta_1}& H_1/I \ar@{-->}[d] \\
& \End_\mathbbm{k}(A)
}
\end{equation*}
(Здесь $\pi \colon H_1 \twoheadrightarrow H_1/I$~"--- естественный сюръективный
гомоморфизм.) 
 
  Для того, что бы показать, что алгебра Хопфа $H_1/I$ кокоммутативна,
достаточно доказать, что $$h_{(1)}\otimes h_{(2)}-
 h_{(2)}\otimes h_{(1)} \in H_1 \otimes I + I \otimes H_1
\text{ для всех }h\in H_1.$$ 
В силу~\eqref{EqH_1IIH_1} достаточно проверить, что
$$(\varphi^k \otimes \varphi^\ell)(h_{(1)}\otimes h_{(2)}-h_{(2)}\otimes h_{(1)})=0.$$
Действительно, $$(\varphi^k \otimes \varphi^\ell)(h_{(1)}\otimes h_{(2)}-h_{(2)}\otimes h_{(1)})=
\varphi^{k+\ell}(h) - \varphi^{k+\ell}(h) = 0.$$
  Следовательно, алгебра Хопфа $H_1/I$ кокоммутативна, откуда и следует доказываемое утверждение.
 \end{proof}

     \section{Действия аффинных алгебраических групп}
     
     Пусть $G$~"--- аффинная алгебраическая группа над алгебраически замкнутым полем $\mathbbm{k}$
и пусть $\mathcal O(G)$~"--- алгебра регулярных
функций на $G$. Предположим, что $G$ действует рационально автоморфизмами на конечномерной алгебре $A$.
  
  Как было отмечено в примере~\ref{ExampleRegularAffActAutAlgebra}, на $A$ действуют три алгебры Хопфа:
  $\mathcal O(G)^\circ$, $\mathbbm{k}G$ и $U(\mathfrak g)$ (универсальная обёртывающая алгебра алгебры Ли~$\mathfrak g$). Ниже в теореме~\ref{TheoremAffAlgGrAllEquiv} доказывается, что все эти три действия
эквивалентны в случае, когда группа $G$ связна. Для того, чтобы это доказать, нам потребуется вспомогательная лемма:
  
  \begin{lemma}\label{LemmaDensityFinDimImage} Пусть $V$~"--- конечномерный комодуль над коалгеброй $C$, $\rho \colon V \to V \otimes C$~"---
  соответствующее  комодульное отображение, а $\zeta \colon C^* \to \End_\mathbbm{k}(V)$~"---
  соответствующее действие алгебры $C^*$
  на $V$, заданное при помощи равенства $c^* v := c^*(v_{(1)})v_{(0)}$ для $c^* \in C^*$ и $v\in V$.
  Пусть $A \subseteq C^*$~"--- \textbf{плотная}\footnote{Именно такое определение плотности даётся в монографии~\cite{Montgomery}, на результаты из которой мы ссылаемся ниже, однако из доказательства леммы видно, что это определение эквивалентно обычному определению плотности в конечной топологии.} подалгебра, т.е. $$A^\perp := \lbrace c\in C \mid a(c)=0 
  \text{ для всех } a\in A\rbrace = 0.$$ Тогда $\zeta(A)=\zeta(C^*)$.
  \end{lemma}
  \begin{proof}
  В силу того, что пространство $V$ конечномерно, 
существует конечномерная подкоалгебра $D \subseteq C$, такая, что
  $\rho(V)\subseteq V \otimes D$.
Следовательно, достаточно показать, что ограничение линейных функций из
обоих множеств $C^*$ и $A$ на подпространство на $D$ 
   совпадает с $D^*$. Для множества $C^*$ это очевидно, так как любая линейная функция, заданная на подпространстве  $D$, продолжается на всё пространство $C$. Для множества $A$ равенство доказывается следующим образом.
   Линейные функции на $D$, соответствующие элементам алгебры $A$, образуют подпространство $W$ в пространстве $D^*$. Если $W \ne D^*$,
   то в силу конечномерности пространства $D$ существует элемент $d\in D$, $d\ne 0$, такой, что $w(d)=0$
  для всех $w\in W$. Следовательно, $d \in A^\perp$, и мы получаем
  противоречие с тем, что $A$ плотно в $C^*$. Следовательно, ограничение линейных функций из $A$ на $D$ совпадает с $D^*$. Из равенства $\zeta(c^*)=\zeta(a)$ для всех $a\in A$, $c^*\in C^*$,
  таких, что $c^*\bigr|_{D} = a\bigr|_{D}$, следует, что $\zeta(A)=\zeta(C^*)$.
  \end{proof}
  
Теперь мы готовы доказать теорему:

\begin{theorem}\label{TheoremAffAlgGrAllEquiv}
Пусть $G$~"--- связная аффинная алгебраическая 
группа над алгебраически замкнутым полем $\mathbbm{k}$ характеристики $0$,
действующая рационально автоморфизмами на конечномерной алгебре $A$. Пусть $\mathfrak g$~"--- алгебра Ли
 группы $G$. Тогда соответствующие действия алгебр Хопфа
 $\mathbbm{k}G$, $U(\mathfrak g)$ и $\mathcal O(G)^\circ$ на $A$ эквивалентны.
\end{theorem}
\begin{proof}
В силу леммы~\ref{LemmaDensityFinDimImage}
достаточно доказать, что образы алгебр $\mathbbm{k}G$,
$U(\mathfrak g)$ и $\mathcal O(G)^\circ$ плотны в $\mathcal O(G)^*$.
 Для $\mathbbm{k}G$ это следует из определения алгебры $\mathcal O(G)$. Если мы докажем, что
 алгебра  $U(\mathfrak g)$ плотна в $\mathcal O(G)^*$, отсюда сразу будет следовать, что $\mathcal O(G)^\circ$ также плотна в $\mathcal O(G)^*$, так как $U(\mathfrak g) \subseteq \mathcal O(G)^\circ$.

Пусть $I := \ker \varepsilon$, множество всех полиномиальных функций из алгебры $\mathcal O(G)$,
которые принимают нулевое значение в точке $1_G$. В силу предложения~9.2.5 из~\cite{Montgomery} $$U(\mathfrak g)=\lbrace a \in
\mathcal O(G)^\circ \mid a(I^n)=0 \text{ для некоторого } n\in\mathbb N \rbrace=\mathcal O(G)',$$
где $\mathcal O(G)'$ неприводимая компонента для $\varepsilon$ в $\mathcal O(G)^\circ$ (см. определение~5.6.1 в~\cite{Montgomery}). В силу~\cite[глава~II, \S~7.3]{HumphreysAlgGr} 
из связности группы $G$ следует, что $G$ неприводимо как многообразие. 
Отсюда, применив следствие из теоремы Крулля~\cite[следствие~10.18]{AtiyahMacdonald},
получаем, что $\bigcap_{n\geqslant 1} I^n = 0$. Из предложения~9.2.10 монографии~\cite{Montgomery}
следует, что $U(\mathfrak g)=\mathcal O(G)'$ плотно в $\mathcal O(G)^*$. Отсюда $\mathcal O(G)^\circ$
также плотно в $\mathcal O(G)^*$ и действия алгебр Хопфа $\mathbbm{k}G$, $U(\mathfrak g)$ и $\mathcal O(G)^\circ$
на $A$ эквивалентны.\end{proof}
\begin{remark}
Доказательство теоремы~\ref{TheoremAffAlgGrAllEquiv} основано на идее, предложенной М.\,В.~Кочетовым
в ходе совместной работы с автором над статьёй~\cite{ASGordienko6Kochetov}.
\end{remark}
       
\section{Действия алгебр Хопфа на алгебре двойных чисел}\label{SectionDoubleNumbers}

В этом параграфе мы покажем на примере алгебры $\mathbbm{k}[x]/(x^2)$ двойных чисел,
как можно использовать понятие эквивалентности модульных структур для классификации последних.
Для краткости будем обозначать алгебру Хопфа Свидлера $H_4(-1)$ (см. пример~\ref{ExampleTaftAlgebra}) просто через $H_{4}$. 
Напомним, что как алгебра $H_{4}$ 
порождена двумя элементами $c$ и $v$, удовлетворяющими соотношениям $c^{2} = 1$, $v^{2} = 0$ и $vc = -cv$.
Структура коалгебры и антипод задаются равенствами
$$
\Delta(c)=c\otimes c,\quad \Delta(v)=c\otimes v + v\otimes 1, \quad S(c)=c,\quad S(v)=-cv.
$$

\begin{theorem}\label{TheoremDoubleNumbersClassify}
Пусть $\psi \colon H \otimes A \to A$~"--- структура $H$-модульной алгебры с $1$ на $A = \mathbbm{k}[x]/(x^2)$,
где $H$~"--- некоторая алгебра Хопфа, $\chr \mathbbm{k} \ne 2$.
Тогда $\psi$ эквивалентно одной из следующих модульных структур на $A$:
\begin{enumerate}
\item действие поля $\mathbbm{k}$ на алгебре $A$ умножением на скаляры;

\item действие групповой алгебры $\mathbbm{k}G$, где $G=\langle c \rangle_2$, заданное равенством  $c \bar x = -\bar x$;

\item $H_4$-действие, заданное равенствами
$c \bar 1 = \bar 1$, $c \bar x = -\bar x$, $v \bar 1 = 0$,
$v\bar x = \bar 1$. 
\end{enumerate}
\end{theorem}
\begin{proof}
Как и прежде, определим отображение $\zeta \colon H \to \End_\mathbbm{k}(A)$ через $\zeta(h)(a)=\psi(h\otimes a)$,
где $a\in A$, $h\in H$. Кроме того, фиксируем базис $\bar 1$, $\bar x$ в алгебре $A$ и отождествим $\End_\mathbbm{k}(A)$ с алгеброй $M_2(\mathbbm{k})$
квадратных матриц $2 \times 2$.
В силу того, что $A$~"--- $H$-модульная алгебра с $1$, существуют $\alpha,\beta \in H^*$,
такие, что $\zeta(h)=\left(\begin{smallmatrix} \varepsilon(h) & \beta(h) \\
0 & \alpha(h) \end{smallmatrix}\right)$ для всех $h\in H$.

Из~(\ref{EqModCompat}) следует, что
\begin{equation*}\begin{split}
0 = h(\bar x^2)=(h_{(1)} \bar x)(h_{(2)} \bar x)
=(\beta(h_{(1)})\bar 1 + \alpha(h_{(1)}) \bar x)(\beta(h_{(2)})\bar 1 + \alpha(h_{(2)}) \bar x) =
\\ =\beta(h_{(1)})\beta(h_{(2)})\bar 1 + (\alpha(h_{(1)})\beta(h_{(2)}) + \beta(h_{(1)})\alpha(h_{(2)}))\bar x.
\end{split}\end{equation*}
Отсюда
\begin{equation}\label{EqBetaBeta}\beta(h_{(1)})\beta(h_{(2)})=0,\end{equation}
\begin{equation}\label{EqAlphaBeta}\alpha(h_{(1)})\beta(h_{(2)}) + \beta(h_{(1)})\alpha(h_{(2)})=0\end{equation}
для всех $h\in H$.

Рассмотрев ранг матрицы линейного
оператора $\zeta$, заключаем, что $\dim \zeta(H)=\dim \langle \alpha, \beta, \varepsilon \rangle_\mathbbm{k}$.

Если $\dim \zeta(H) = 3$, то $\zeta(H)$~"--- подалгебра
всех верхнетреугольных матриц, и модульная структура $\psi$ эквивалентна структуре 3.

Пусть $\dim \zeta(H) \leqslant 2$.
Если $\alpha$ и $\varepsilon$ линейно зависимы,
то из $\varepsilon\ne 0$ следует, что $\alpha = \gamma \varepsilon$ для некоторого $\gamma \in \mathbbm{k}$.
Отсюда $1=\alpha(1_H) = \gamma \varepsilon(1_H)=\gamma$ и $\alpha = \varepsilon$.
В силу~(\ref{EqAlphaBeta}) имеем $2\beta(h)=0$ и $\beta = 0$. Следовательно, $\zeta(H)$~"--- это алгебра всех скалярных матриц, и модульная структура $\psi$ эквивалентна структуре 1.

Предположим, что $\dim \zeta(H) \leqslant 2$, однако $\alpha$ и $\varepsilon$ линейно независимы.
Тогда $\beta = \lambda \varepsilon + \gamma \alpha$ и
из равенства $\zeta(1_H)=\left(\begin{smallmatrix} 1 & 0 \\
0 & 1 \end{smallmatrix}\right)$ следует, что
$0=\lambda+\gamma$ и $\beta= \lambda(\varepsilon-\alpha)$.
Из~(\ref{EqBetaBeta}) и~(\ref{EqAlphaBeta})
следует, что
\begin{equation}\label{EqDoubleNumbersLambdaVarepsilon}\lambda^2\bigl(\varepsilon(h)-2
\alpha(h)+ \alpha(h_{(1)})\alpha(h_{(2)})\bigr)=0,\end{equation}
\begin{equation}\label{EqDoubleNumbersLambdaAlpha}2\lambda\left(\alpha(h) - \alpha(h_{(1)})\alpha(h_{(2)})\right)=0\end{equation}
для всех $h\in H$. Предположим, что $\lambda \ne 0$.
В силу того, что $\chr \mathbbm{k} \ne 2$, из (\ref{EqDoubleNumbersLambdaAlpha}) следует, что
$$\alpha(h_{(1)})\alpha(h_{(2)})=\alpha(h).$$
Поэтому из~(\ref{EqDoubleNumbersLambdaVarepsilon}) получаем, что 
$$ \alpha(h)=\varepsilon(h),$$
что противоречит предположению о линейной независимости
 $\alpha$ и $\varepsilon$.
Следовательно, $\lambda = 0$ и $\beta = 0$.
Отсюда $\zeta(H)$ совпадает с алгеброй диагональных матриц
и модульная структура $\psi$ эквивалентна структуре 2.
\end{proof}

\newpage

   \chapter{Алгебры с обобщённым $H$-действием}\label{ChapterGenHActions}
   
   К сожалению, не все важные дополнительные структуры на алгебрах можно представить в виде (ко)модульных структур над алгебрами Хопфа. С другой стороны, например, при изучении полиномиальных тождеств в алгебрах для многих построений можно не требовать, чтобы алгебра была модульной над некоторой алгеброй Хопфа, достаточно выполнения более слабых условий, сформулированных ниже. Соответствующие действия мы будем называть обобщёнными $H$-действиями, где $H$~"--- это ассоциативная алгебра с единицей, действующая на заданной алгебре. В случае, когда алгебра над полем $\mathbbm{k}$ градуирована некоторой бесконечной группой $G$ и эта градуировка имеет конечный носитель, при работе с градуированными тождествами удобно бывает заменить градуировку на действие алгебры $(\mathbbm{k}G)^*$, которое также является обобщённым. (См. пример~\ref{ExampleGr} ниже.)
   
     \section{Обобщённые $H$-действия}

Пусть $H$~"--- произвольная ассоциативная алгебра с $1$ над полем $\mathbbm{k}$.
Будем говорить, что (необязательно ассоциативная) алгебра $A$ является
алгеброй с \textit{обобщённым $H$-действием}, если $A$ является левым $H$-модулем
и для любого $h \in H$ существует такое $k\in \mathbb N$ и такие $h'_i, h''_i, h'''_i, h''''_i \in H$, где $1\leqslant i \leqslant k$,
что
\begin{equation}\label{EqGeneralizedHopf}
h(ab)=\sum_{i=1}^k\bigl((h'_i a)(h''_i b) + (h'''_i b)(h''''_i a)\bigr) \text{ для всех } a,b \in A.
\end{equation}

Эквивалентным условием является существование линейных отображений $\Delta, \Theta \colon H \to H\otimes H$ (необязательно коассоциативных), таких, что
$$ h(ab)=\sum\bigl((h_{(1)} a)(h_{(2)} b) + (h_{[1]} b)(h_{[2]} a)\bigr) \text{ для всех } a,b \in A.$$ (Здесь мы используем обозначения $\Delta(h)= \sum h_{(1)} \otimes h_{(2)}$ и $\Theta(h)= \sum  h_{[1]} \otimes h_{[2]}$.)

\begin{remark} Одна и та же алгебра $H$ может действовать на разных алгебрах $A$. При этом отображения $\Delta$ и $\Theta$ также могут быть разными для разных алгебр $A$.
\end{remark}

\begin{example}\label{ExampleHmodule} 
Если $A$~"--- $H$-модульная алгебра, где $H$~"--- некоторая биалгебра,
то $A$~"--- алгебра с обобщённым $H$-действием.
\end{example}

Напомним, что линейное отображение $\varphi\colon A \to A$
называется \textit{антиэндоморфизмом} алгебры $A$, если $\varphi(ab)=\varphi(b)\varphi(a)$
для всех $a,b\in A$. Если моноид $T$ действует на $A$ только эндоморфизмами,
т.е. $t(ab)=(ta)(tb)$ для всех $a,b\in A$, то алгебра $A$ является $\mathbbm{k}T$-модульной, где $\mathbbm{k}T$~"--- полугрупповая биалгебра моноида $T$ (см. примеры~\ref{ExampleFTBialgebra} и~\ref{ExampleFTmodule}).
Однако если $T$ действует ещё и антиэндоморфизмами, то $A$ уже не будет $\mathbbm{k}T$-модульной алгеброй в классическом смысле:

\begin{example}\label{ExampleFTEndAntiEnd}
Пусть $A$~"--- алгебра над полем $\mathbbm{k}$. 
Если моноид $T$ действует на алгебре $A$ эндоморфизмами и антиэндоморфизмами,
то $A$ является алгеброй с обобщённым $\mathbbm{k}T$-действием.
\end{example}

Другим классом алгебр, для которых понятие обобщённого $H$-действия оказывается особенно полезным,
являются градуированные алгебры с конечным носителем градуировки:

\begin{example}\label{ExampleGr}
Пусть $A=\bigoplus_{t\in T} A^{(t)}$~"--- алгебра над полем $\mathbbm{k}$,
градуированная множеством $T$, т.е. для всех $s,t \in T$ существует такое $r\in T$,
что $A^{(s)}A^{(t)}\subseteq A^{(r)}$. Обозначим эту градуировку через $\Gamma$. 
Рассмотрим алгебру $\mathbbm{k}^T$ всех функций из $T$ в $\mathbbm{k}$ с поточечными операциями.
Тогда $\mathbbm{k}^T$ естественным образом действует на $A$: $qa = q(t)a$ для всех $a\in A^{(t)}$
и $q\in \mathbbm{k}^T$,
причём при таком действии $T$-градуированные подпространства в $A$
являются в точности $\mathbbm{k}^T$-подмодулями.
Обозначим через $\zeta \colon \mathbbm{k}^T \to \End_\mathbbm{k}(A)$
соответствующий гомоморфизм.
Заметим, что градуировка $\Gamma$ 
при помощи равенства $s\star t := r$
задаёт на множестве $T$
частичную операцию $\star$ с областью определения $T_0:=\lbrace (s,t) \mid A^{(s)}A^{(t)} \ne 0 \rbrace$.
Пусть $q_t(s):=\left\lbrace\begin{smallmatrix} 1 & \text{при} & s=t,\\ 0 & \text{при} & s\ne t.\end{smallmatrix} \right.$ \label{DefHT}
Если носитель $$\supp \Gamma := \lbrace t\in T \mid A^{(t)}\ne 0\rbrace$$ градуировки $\Gamma$
конечен, множество $T_0$ также конечно и справедливы равенства
\begin{equation}\label{EqIdentityHFiniteSupp}q_r(ab)=\sum_{\substack{(s,t)\in T_0,\\ r=s\star t}}
(q_s a)(q_t b). \end{equation}
(Поскольку обе части равенства линейны по $a$ и $b$, достаточно проверить это равенство только для однородных элементов  $a,b$.) В силу того, что 
алгебра $\mathbbm{k}^T$ является  по модулю $\ker \zeta$ линейной оболочкой множества $(q_t)_{t\in \supp \Gamma}$,
из~\eqref{EqIdentityHFiniteSupp} следует~(\ref{EqGeneralizedHopf}) для всех $h\in \mathbbm{k}^T$.
Отсюда в случае, когда носитель градуировки конечен, $A$
является алгеброй с обобщённым $\mathbbm{k}^T$-действием. 
\end{example}
\begin{example}
Если в предыдущем примере $T$~"--- группа, то $\mathbbm{k}^T \cong (\mathbbm{k}T)^*$. Если группа $T$ конечна, то
$(\mathbbm{k}T)^*$ является алгеброй Хопфа и структура алгебры с обобщённым $\mathbbm{k}^T$-действием на $A$, введённая выше, совпадает со структурой $(\mathbbm{k}T)^*$-модульной алгебры на $A$.
\end{example}
\begin{example}
Если $A$~"--- конечномерная $B$-комодульная алгебра, где $B$~"--- некоторая биалгебра,
то на $A$ также можно ввести структуру алгебры с обобщённым $B^*$-действием.
Достаточно применить лемму~\ref{LemmaSubstituteComult}, взяв в качестве $H_1 \subseteq B$ любое такое конечномерное подпространство, что $\rho(A)\subseteq A \otimes H_1$. (Здесь, как обычно, $\rho \colon A \to A \otimes B$~"--- линейное отображение, задающее на $A$ структуру $B$-комодуля.)
\end{example}

\begin{example} \label{ExampleGenIdGenHAction}
Пусть $A$~"--- ассоциативная алгебра, а $A^+:= A +  \mathbbm{k}\cdot 1$.
Тогда $A$ является алгеброй с обобщённым $A^{+} \otimes \left(A^{+}\right)^{\mathrm{op}}$-действием,
где $$(b\otimes c)a:= bac\text{ для всех }a \in A,\ b \in A^{+}\text{ и }c \in \left(A^{+}\right)^{\mathrm{op}}.$$
Действительно,
$$(b\otimes c)(a_1 a_2):= \bigl((b \otimes 1)a_1 \bigr) \bigl((1 \otimes c)a_2\bigr)\text{ для всех }a_1,a_2 \in A,\ b \in A^{+}\text{ и }c \in \left(A^{+}\right)^{\mathrm{op}}.$$
\end{example}
\medskip

Пусть $A$~"--- алгебра с обобщённым $H$-действием для некоторой ассоциативной алгебры $H$ с единицей $1$
над полем $\mathbbm{k}$. Как и в случае модульных алгебр, подпространство $V\subseteq A$ называется  \textit{инвариантным относительно $H$-действия}, если $HV=V$, т.е. если $V$ является $H$-подмодулем.
Если $A^2\ne 0$ и алгебра $A$ не содержит нетривиальных $H$-инвариантных двусторонних идеалов, то $A$ называется \textit{$H$-простой} алгеброй.

\section{Обобщённые действия, согласованные с градуировками}\label{SectionGradedActions}

В ряде случаев алгебра $A$ над полем $\mathbbm{k}$ бывает наделена, с одной стороны, градуировкой некоторым множеством $T$,
а с другой стороны, обобщённым $H$-действием некоторой ассоциативной алгебры $H$ с $1$,
причём $H$-действие сохраняет компоненты градуировки.  Будем называть такое обобщённое
$H$-действие \textit{согласованным} с $T$-градуировкой.

\begin{example}\label{ExampleSuperInvolution} Пусть $A=A^{(0)}\oplus A^{(1)}$~"--- некоторая $\mathbb Z/2\mathbb Z$-градуированная алгебра. Линейное отображение $\star \colon A \to A$ называется \textit{суперинволюцией} (см. \cite{RacineSuperInvolution}),
если $\left(A^{(k)}\right)^\star = A^{(k)}$ при $k=0,1$, $a^{\star\star}=a$ для всех $a\in A$,
$(ab)^\star=(-1)^{k\ell}b^\star a^\star$ для всех $a\in A^{(k)}$, $b\in A^{(\ell)}$, где $k,\ell \in \lbrace 0,1\rbrace$. На всякой алгебре с суперинволюцией
задано естественное действие групповой алгебры циклической группы второго порядка,
согласованное с $\mathbb Z/2\mathbb Z$-градуировкой.
\end{example}
\begin{example}\label{ExamplePseudoInvolution} Пусть $A=A^{(0)}\oplus A^{(1)}$~"--- некоторая $\mathbb Z/2\mathbb Z$-градуированная алгебра. Линейное отображение $\star \colon A \to A$ называется \textit{псевдоинволюцией} (см. \cite{MartinezZelmanov}),
если $\left(A^{(k)}\right)^\star = A^{(k)}$ при $k=0,1$, $a^{\star\star}=(-1)^k a$,
$(ab)^\star=(-1)^{k\ell}b^\star a^\star$ для всех $a\in A^{(k)}$, $b\in A^{(\ell)}$, где $k,\ell \in \lbrace 0,1\rbrace$. 
На всякой алгебре с псевдоинволюцией
задано естественное действие групповой алгебры циклической группы четвертого порядка,
согласованное с $\mathbb Z/2\mathbb Z$-градуировкой.
\end{example}

Как было впервые отмечено Р.\,Б. дос Сантосом~\cite{dosSantos}, всякая алгебра с суперинволюцией является
алгеброй с обобщённым $H$-действием. Подобное утверждение справедливо
и в общем случае.  Напомним, что через $\mathbbm{k}^T$ обозначается алгебра всех функций из $T$ в $\mathbbm{k}$ с поточечными операциями.
\begin{theorem}\label{TheoremGradGenActionReplace} Пусть $\Gamma \colon A=\bigoplus_{t\in T} A^{(t)}$~"--- $T$-градуировка на алгебре $A$ над полем $\mathbbm{k}$, где $T$~"--- некоторое множество, причём на $A$ также заданное такое обобщённое $H$-действие, 
что $H$~"--- ассоциативная алгебра с $1$, $HA^{(t)}\subseteq A^{(t)}$ для всех $t\in T$, а носитель $\supp \Gamma$ градуировки $\Gamma$ конечен.
Тогда формула $$(q\otimes h)a:=q(t)(ha)\text{ при }a\in A^{(t)},\ t\in T,\ q\in \mathbbm{k}^T,\ h\in H$$
задаёт на $A$ обобщённое $\mathbbm{k}^T \otimes H$-действие, причём $T$-градуированные $H$-подмодули в $A$
являются при таком действии в точности $\mathbbm{k}^T \otimes H$-подмодулями.
\end{theorem}
\begin{proof}
Как и в примере~\ref{ExampleGr}, определим элементы $q_t \in \mathbbm{k}^T$ при $t\in T$
по формуле $q_t(s):=\left\lbrace\begin{smallmatrix} 1 & \text{при} & s=t,\\ 0 & \text{при} & s\ne t\end{smallmatrix} \right.$ и зададим на множестве $T$
частичную операцию $\star$ с областью определения $T_0:=\lbrace (s,t) \mid A^{(s)}A^{(t)} \ne 0 \rbrace$
при помощи равенства $s\star t := r$,
где $A^{(s)}A^{(t)}\subseteq A^{(r)}$.
Тогда из условия~\eqref{EqGeneralizedHopf} для обобщённого $H$-действия следует, что
для любого $h \in H$ существует такое $k\in \mathbb N$ и такие $h'_i, h''_i, h'''_i, h''''_i \in H$, где $1\leqslant i \leqslant k$,
что
\begin{equation}\label{EqGradedActionGenHAction}
(q_r \otimes h)(ab)=\sum_{i=1}^k \sum_{\substack{(s,t)\in T_0,\\ r=s\star t}} \Bigl(\bigl((q_s\otimes h'_i) a\bigr)\bigl((q_t\otimes h''_i) b\bigr) + \bigl((q_t\otimes h'''_i) b\bigr)\bigl((q_s\otimes h''''_i) a\bigr)\Bigr)
\end{equation}
 для всех $r\in T$ и $a,b\in A$.
(В силу линейности обеих частей равенства~\eqref{EqGradedActionGenHAction}
по $a$ и $b$ его достаточно проверить для однородных элементов $a$ и $b$.)
Теперь~\eqref{EqGeneralizedHopf} для $\mathbbm{k}^T \otimes H$-действия следует из того, что $\mathbbm{k}^T \otimes H$ является по модулю ядра действия $\mathbbm{k}^T \otimes H
\to \End_\mathbbm{k}(A)$ линейной оболочкой
элементов $q_r \otimes h$, где $r\in \supp \Gamma$, а $h\in H$.
\end{proof}

Понятно, что когда носитель градуировки состоит из одного элемента, такое обобщённое $\mathbbm{k}^T \otimes H$-действие
сводится просто к обобщённому $H$-действию для той же самой алгебры $H$. Отсюда примеры обобщённых $H$-действий на конечномерных ассоциативных
алгебрах с не $H$-инвариантными радикалами Джекобсона (см. примеры~\ref{ExampleQ1}--\ref{ExampleQ3} вместе с примером~\ref{ExampleGr})
 показывают, что в алгебрах, в которых $T$-градуировка согласована с обобщённым $H$-действием,
 радикал Джекобсона не обязан быть $\mathbbm{k}^T \otimes H$-подмодулем.

Сформулируем теперь достаточные условия инвариантности радикала Джекобсона.
Для этого, отталкиваясь от примеров~\ref{ExampleSuperInvolution}
и~\ref{ExamplePseudoInvolution}, рассмотрим следующий достаточно общий случай.

\begin{definition}\label{DefGradedAction}
 Пусть $G$ и $T$~"--- группы, а $A$~"--- алгебра над полем $\mathbbm{k}$.
Будем говорить, что $A$ наделена \textit{$T$-градуированным $G$-действием}, если на $A$ задана групповая градуировка $A=\bigoplus_{t\in T} A^{(t)}$, гомоморфизм групп $G \to \GL(A)$ и функции
$\alpha, \beta \colon G \times T \times T \to \mathbbm{k}$, такие, что $gA^{(t)}\subseteq A^{(t)}$ и
\begin{equation}\label{EqGradedAction}
g(ab)=\alpha(g,s,t)(ga)(gb)+\beta(g,s,t)(gb)(ga)
\end{equation}
  для всех $g\in G$, $s,t\in T$, $a\in A^{(s)}$ и $b\in A^{(t)}$.
  \end{definition}
  
  \begin{remark}
  Из~\eqref{EqGradedAction} следует, что $\beta(g,s,t)=0$ для всех $g\in G$ и $s,t\in T$,
таких, что $A^{(t)}A^{(s)}\ne 0$ и $st\ne ts$.
  \end{remark}
  
  В силу теоремы~\ref{TheoremGradGenActionReplace} всякое $T$-градуированное $G$-действие
  сводится к обобщённому $\mathbbm{k}^T\otimes \mathbbm{k}G$-действию.

Докажем теперь, что радикал Джекобсона конечномерной ассоциативной алгебры инвариантен относительно
такого действия.

Для начала покажем, что всякий двусторонний градуированный идеал переходит под градуированным действием группы
в некий двусторонний градуированный идеал, причём если этот идеал был нильпотентным, его образ также будет 
нильпотентным идеалом:

\begin{lemma}\label{LemmaGradedActionIdeal} Пусть $\Gamma \colon A=\bigoplus_{t\in T} A^{(t)}$~"--- $T$-градуировка на
алгебре $A$ над полем $\mathbbm{k}$, причём на $A$ также задано $T$-градуированное $G$-действие, где $T$ и $G$~"--- произвольные группы. Тогда для всякого $g\in G$ и двустороннего градуированного идеала $I\subseteq A$
пространство $gI$ также является двусторонним градуированным идеалом, причём если идеал $I$ нильпотентен,
идеал $gI$ также нильпотентен.
\end{lemma}
\begin{proof} Всякий элемент $a \in I$ можно представить в виде суммы его однородных компонент,
причём в силу градуированности идеала $I$ все эти компоненты принадлежат идеалу $I$. Применяя к $a$ произвольный элемент $g\in G$ и учитывая, что группа $G$ сохраняет компоненты градуировки, получаем,
что $ga$ является суммой однородных элементов, каждый из которых принадлежит пространству $gI$. Отсюда
пространство $gI$ градуированное.

В силу градуированности пространства $gI$ для того, чтобы доказать, что $gI$~"--- двусторонний идеал,
достаточно показать, что $ab, ba \in gI$ для всех $a \in A^{(s)}$, $b\in (gI)\cap A^{(t)}$, где $s,t \in T$.
Из~\eqref{EqGradedAction} следует, что
$$g^{-1}(ab)= \alpha(g^{-1},s,t)(g^{-1}a)(g^{-1}b)+\beta(g^{-1},s,t)(g^{-1}b)(g^{-1}a) \in I$$
и
$$g^{-1}(ba)= \alpha(g^{-1},t,s)(g^{-1}b)(g^{-1}a)+\beta(g^{-1},t,s)(g^{-1}a)(g^{-1}b) \in I.$$
Отсюда $gI$~"--- градуированный двухсторонний идеал.

Пусть теперь любое произведение из $n$ элементов идеала $I$ с любой расстановкой скобок равно $0$.
Докажем, что идеал $gI$ обладает тем же свойством. Действительно, пусть $a_1, \ldots, a_n \in gI$~"---
однородные элементы. Рассматривая элемент $g^{-1}(a_1 \ldots a_n)$, где на произведении
$a_1 \ldots a_n$ задана произвольная расстановка скобок, и применяя $(n-1)$ раз формулу~\eqref{EqGradedAction}, получаем, что  $g^{-1}(a_1 \ldots a_n)$
является линейной комбинацией произведений элементов $g^{-1}a_i$ взятых в произвольном порядке.
Поскольку $g^{-1}a_i \in I$ для всех $i=1,\ldots, n$, а число множителей равно $n$,
все такие произведения равны нулю и $g^{-1}(a_1 \ldots a_n)=0$,
откуда $a_1 \ldots a_n=0$. Следовательно, в этом случае идеал $gI$ также нильпотентен.
\end{proof}

Теперь выведем отсюда инвариантность радикала Джекобсона:

\begin{theorem}\label{TheoremTGradedGActionInvRad} Пусть $\Gamma \colon A=\bigoplus_{t\in T} A^{(t)}$~"--- $T$-градуировка на
конечномерной ассоциативной алгебре $A$ над полем $\mathbbm{k}$, причём на $A$ также задано $T$-градуированное $G$-действие, где $T$ и $G$~"--- произвольные группы, и либо $\chr \mathbbm{k} = 0$, либо $\chr \mathbbm{k} > \dim_\mathbbm{k} A$.
Тогда $$\bigl(\mathbbm{k}^T\otimes \mathbbm{k}G\bigr)J(A)\subseteq J(A).$$
\end{theorem}
\begin{proof}
В силу следствия~\ref{CorollaryRadicalGraded} идеал $J(A)$ является градуированным и, следовательно,
$\mathbbm{k}^T$-инвариантным идеалом. Инвариантность идеала $J(A)$ относительно действия группы $G$
 теперь следует из леммы~\ref{LemmaGradedActionIdeal}.
\end{proof}

Для того, чтобы доказать инвариантный аналог теоремы Веддербёрна~"--- Артина, нам
потребуется следующая лемма:

  \begin{lemma}\label{LemmaAnnIdealTGradedGAction}
  Пусть $\Gamma \colon A=\bigoplus_{t\in T} A^{(t)}$~"--- $T$-градуировка на
ассоциативной алгебре $A$ над полем $\mathbbm{k}$, причём на $A$ также задано $T$-градуированное $G$-действие, где $T$ и $G$~"--- произвольные группы.
Пусть $M\subseteq A$~"--- $G$-инвариантное $T$-градуированное подпространство.
Тогда $$\Ann_{\mathrm{lr}}(M) := \lbrace a \in A \mid ab=ba=0 \text{ для всех } b \in M \rbrace$$ также
является $G$-инвариантным $T$-градуированным подпространством.
  \end{lemma}
  \begin{proof} В силу леммы~\ref{LemmaAnnIdealHcomod}
  подпространство $\Ann_{\mathrm{lr}}(M)$ является $T$-градуированным.
  Докажем теперь, что оно $G$-инвариантно.
  
  Пусть $a \in \Ann_{\mathrm{lr}}(M) \cap A^{(s)}$, $b\in M \cap A^{(t)}$, $s,t \in T$, $g\in G$. Тогда
$$ g^{-1}\bigl((ga)b\bigr)=
 \alpha(g^{-1},s,t)a(g^{-1}b)+\beta(g^{-1},s,t)(g^{-1}b)a
=0,$$
откуда $(ga)b=0$. 

Аналогично,
$$g^{-1}(b(ga))= \alpha(g^{-1},t,s)(g^{-1}b)a+\beta(g^{-1},t,s)a(g^{-1}b)=0,$$
откуда $b(ga)=0$.

В силу линейности равенств \begin{equation}\label{EqGabbga}(ga)b=b(ga)=0\end{equation} по $a$ и $b$
и градуированности подпространств $M$ и $\Ann_{\mathrm{lr}}(M)$
равенства~\eqref{EqGabbga} справедливы для всех $a \in \Ann_{\mathrm{lr}}(M)$, $b\in M$.
    Таким образом, $ga \in \Ann_{\mathrm{lr}}(M)$, откуда $\Ann_{\mathrm{lr}}(M)$
является $G$-инвариантным подпространством.
\end{proof}

\begin{theorem}\label{TheoremTGradedGActionInvWedderburnArtin} Пусть $\Gamma \colon B=\bigoplus_{t\in T} B^{(t)}$~"--- $T$-градуировка на
конечномерной полупростой ассоциативной алгебре $B$ над полем $\mathbbm{k}$, причём на $B$ также задано $T$-градуированное $G$-действие, где $T$ и $G$~"--- произвольные группы.
Тогда $$B = B_1 \oplus B_2 \oplus \ldots \oplus B_s$$
  (прямая сумма $T$-градуированных $G$-инвариантных идеалов) для некоторых алгебр $B_i$
  с $T$-градуированным $G$-действием, каждая из которых проста в соответствующем смысле, т.е. не содержит нетривиальных $T$-градуированных $G$-инвариантных идеалов.
\end{theorem}
\begin{proof}
Достаточно повторить доказательство теоремы~\ref{TheoremWedderburnHmod},
использовав вместо леммы~\ref{LemmaAnnIdealHmod}
лемму~\ref{LemmaAnnIdealTGradedGAction}.
\end{proof}

     \section{Слабое разложение Веддербёрна~"--- Мальцева}
     
Аналогично тому, как это было сделано в случае $H$-модульных алгебр (см. \S\ref{SectionH(co)invWedArtMalcev}), определим \textit{$H$-радикал}
$J^H(A)$ конечномерной ассоциативной алгебры $A$ с обобщённым $H$-действием
как максимальный нильпотентный $H$-инвариантный идеал алгебры $A$.

$H$-инвариантный аналог теоремы Веддербёрна~"--- Мальцева можно было бы сформулировать следующим образом:
если $A$~"--- конечномерная ассоциативная алгебра с обобщённым $H$-действием над полем характеристики $0$,
то существует такой гомоморфизм $\varkappa \colon A/J^H(A) \hookrightarrow A$ алгебр и $H$-модулей,
что $\pi\varkappa = \id_{A/J^H(A)}$, где $\pi \colon A \twoheadrightarrow A/J^H(A)$~"--- естественный сюръективный гомоморфизм.
К сожалению, в такой формулировке теорема неверна даже для $H$-модульных алгебр,
см. примеры~\ref{ExampleGnoninvWedderburnMalcev} и~\ref{ExampleU(L)noninvWedderburnMalcev}.

Докажем слабый вариант теоремы Веддербёрна~"--- Мальцева, а именно, что
существует такое $\mathbbm{k}$-линейное отображение $\varkappa \colon A/J^H(A) \hookrightarrow A$,
что $\pi\varkappa = \id_{A/J^H(A)}$ и $\varkappa$ является гомоморфизмом $(B,B)$-бимодулей
для некоторой максимальной полупростой подалгебры 
$B \subseteq A/J^H(A)$. Этот вариант теоремы Веддербёрна~"--- Мальцева будет затем использоваться в доказательстве аналога гипотезы Амицура
для алгебр с обобщённым $H$-действием (см. главу~\ref{ChapterGenHAssocCodim}).

\begin{theorem}[{\cite[лемма 2.6]{ASGordienko15}}]\label{TheoremWeakWedderburnMalcevHRad}
Пусть $A$~"--- конечномерная ассоциативная алгебра с обобщённым $H$-действием
для некоторой ассоциативной алгебры $H$ с $1$ над алгебраически замкнутым полем $\mathbbm{k}$.
Обозначим через $\pi$ естественный сюръективный гомоморфизм $A \twoheadrightarrow A/J^H(A)$.
Тогда существует такое $\mathbbm{k}$-линейное вложение $\varkappa \colon A/J^H(A) \hookrightarrow A$
и такая полупростая (в обычном смысле)
подалгебра $B \subseteq A/J^H(A)$, что \begin{enumerate}\item $\pi\varkappa=\id_{A/J^H(A)}$;
\item $A/J^H(A) = B \oplus J(A/J^H(A))$ (прямая сумма подпространств);
\item
справедливы равенства
 $\varkappa(ba)=\varkappa(b)\varkappa(a)$ и $\varkappa(ab)=\varkappa(a)\varkappa(b)$ для всех $a\in A/J^H(A)$ и $b\in B$.
 \end{enumerate}
\end{theorem}
\begin{proof}
В силу обычной теоремы Веддербёрна~"--- Артина
существует такая максимальная полупростая подалгебра
 $B_0 \subseteq A$, что $A = B_0 \oplus J(A)$ (прямая сумма подпространств).
 Будем рассматривать $A$ как $(B_0, B_0)$-бимодуль.
Докажем, что $A$ является прямой суммой
неприводимых $(B_0, B_0)$-подбимодулей.  

Действительно, поскольку $B_0$~"--- конечномерная полупростая
алгебра над алгебраически замкнутым полем $\mathbbm{k}$,
она изоморфна прямой сумме полных матричных алгебр над $\mathbbm{k}$,
т.е. алгебра $B_0 \otimes B_0^\mathrm{op}$ также является полупростой,
а $A$~"--- вполне приводимый левый $B_0 \otimes B_0^\mathrm{op}$-модуль,
где алгебра $B_0^\mathrm{op}$ антиизоморфна алгебре $B_0$, а $(b_1\otimes b_2)a :=b_1 a b_2$ для всех $b_1 \otimes b_2 \in B_0 \otimes B_0^\mathrm{op}$ и $a\in A$.
Поскольку мы не требуем от алгебры $A$ наличия единицы,
отсюда пока только следует полная приводимость $(B_0, B_0)$-бимодуля
$1_{B_0} A\ 1_{B_0}$.

Рассмотрим разложение Пирса $$A=(1-1_{B_0})A(1-1_{B_0})\oplus 1_{B_0}A(1-1_{B_0})\oplus
(1-1_{B_0})A\ 1_{B_0}\oplus 1_{B_0} A\ 1_{B_0}$$ (прямая сумма $(B_0, B_0)$-подбимодулей),
где $1$~"--- формальная (или присоединённая) единица.
Здесь $1_{B_0}A(1-1_{B_0})$~"--- вполне приводимый левый $B_0$-модуль,
$(1-1_{B_0})A\ 1_{B_0}$~"--- вполне приводимый правый $B_0$-модуль,
а $(1-1_{B_0})A(1-1_{B_0})$~"--- векторное пространство с нулевым $(B_0, B_0)$-действием.
Отсюда $A$~"--- прямая сумма неприводимых $(B_0, B_0)$-бимодулей.
Следовательно, существует такой $(B_0, B_0)$-подбимодуль $N\subseteq A$,
что $J(A)=N\oplus J^H(A)$. 
Заметим, что отображение $$\pi\bigr|_{(B_0\oplus N)} \colon (B_0\oplus N) \mathrel{\widetilde\to}
A/J^H(A)$$ является $\mathbbm{k}$-линейной биекцией.
Введём обозначения $\varkappa := \left(\pi\bigr|_{(B_0\oplus N)}\right)^{-1}$ и $B := \pi(B_0)$. 
Тогда $$A/J^H(A) = \pi(B_0) \oplus \pi(J(A))=B \oplus J(A/J^H(A)).$$

Рассмотрим произвольные элементы $a \in A/J^H(A)$ и $b\in B$. Тогда
$$\pi(\varkappa(ab))=ab=\pi\varkappa(a)\pi\varkappa(b)=\pi(\varkappa(a)\varkappa(b)).$$
Поскольку $\varkappa(ab), \varkappa(a)\varkappa(b) \in B_0\oplus N$,
отсюда следует, что $\varkappa(ab)=\varkappa(a)\varkappa(b)$.
Аналогично, $\varkappa(ba)=\varkappa(b)\varkappa(a)$.
\end{proof}

\newpage

\chapter{Свободные алгебры, полиномиальные тождества и их коразмерности}
\label{ChapterFreePICodim}

В данной главе вводится понятие полиномиального $H$-тождества
и градуированного тождества (в самой общей формулировке), рассматриваются пары сопряжённых функторов, отвечающие соответствующим свободным алгебрам, доказываются оценки для коразмерностей тождеств, а также существование $H$-PI-экспоненты у любой конечномерной $H$-простой алгебры и градуированной PI-экспоненты у любой конечномерной градуированно простой алгебры.

В определениях из \S\ref{SectionHPI}--\ref{SectionGradedPI} мы отталкиваемся от работ А.~Берела~\cite{BereleHopf},
Ю.\,А.~Бахтурина и В.\,В.~Линченко~\cite{BahturinLinchenko} и монографии А.~Джамбруно
и М.\,В.~Зайцева~\cite{ZaiGia}.

  Результаты главы были опубликованы в работах~\cite{ASGordienko13,ASGordienko16,ASGordienko20ALAgoreJVercruysse}.
 
     \section{Полиномиальные $H$-тождества}\label{SectionHPI}
          
Обозначим через $\mathbbm{k}\lbrace X \rbrace$ \textit{(абсолютно) свободную неассоциативную алгебру}
на множестве $X$,  т.е. алгебру всевозможных неассоциативных многочленов от переменных из множества $X$ с коэффициентами из поля $\mathbbm{k}$.\label{DefFXbrace}
Тогда $\mathbbm{k}\lbrace X \rbrace = \bigoplus_{n=1}^\infty \mathbbm{k}\lbrace X \rbrace^{(n)}$,
где $\mathbbm{k}\lbrace X \rbrace^{(n)}$~"--- линейная оболочка одночленов, имеющих суммарную степень $n$  по всем буквам.

Пусть $H$~"--- ассоциативная алгебра с единицей. Обозначим через $\mathbbm{k}\lbrace X|H \rbrace$
алгебру, которая как векторное пространство совпадает с $$\bigoplus_{n=1}^\infty \mathbbm{k}\lbrace X \rbrace^{(n)} \otimes H^{{}\otimes n},$$ 
а умножение задаётся формулой
$(v_1 \otimes w_1)(v_2 \otimes w_2):=(v_1) (v_2) \otimes w_1 \otimes w_2$
для $v_1\in \mathbbm{k}\lbrace X \rbrace^{(k)}$, $v_2\in \mathbbm{k}\lbrace X \rbrace^{(\ell)}$,
$w_1 \in H^{{}\otimes k}$, $w_2 \in H^{{}\otimes \ell}$, $k,\ell\in\mathbb N$.

Введём обозначение $x_1^{h_1}\ldots x_n^{h_n} := x_1\ldots x_n 
\otimes h_1 \otimes h_2 \otimes \ldots \otimes h_n$
для всех $x_1, \ldots, x_n\in X$ и $h_1,\ldots, h_n \in H$, где
расстановка скобок на одночленах $x_1^{h_1}\ldots x_n^{h_n}$ и $x_1\ldots x_n$ совпадает.

Пусть $(h_{\alpha})_{\alpha\in \Lambda}$~"--- базис алгебры $H$. Тогда $\mathbbm{k}\lbrace X|H \rbrace$
 как алгебра изоморфна свободной неассоциативной алгебре на множестве 
$\lbrace x^{h_\alpha} \mid x\in X,\  \alpha\in \Lambda\rbrace$.

\begin{remark}
В случае, когда $H$~"--- это просто ассоциативная алгебра с единицей, мы не рассматриваем на алгебре
$\mathbbm{k}\lbrace X|H \rbrace$ никакого действия. 
\end{remark}

Отождествим $X$ с подмножеством $$\lbrace x^1 \mid x\in X\rbrace \subseteq \mathbbm{k}\lbrace X|H \rbrace.$$ Обозначим через $\iota \colon X \to \mathbbm{k}\lbrace X|H \rbrace$ соответствующее вложение.
Тогда $\mathbbm{k}\lbrace X|H \rbrace$ удовлетворяет следующему универсальному свойству:
для любого отображения $\varphi \colon X \to A$, где $A$~"--- алгебра с обобщённым $H$-действием,
существует единственный гомоморфизм алгебр $\bar \varphi \colon \mathbbm{k}\lbrace X | H \rbrace \to A$,
такой, что для всех $h\in H$ и $x\in X $ справедливо равенство $\bar\varphi(x^h)=h\bar\varphi(x)$ и 
 диаграмма ниже коммутативна:

\begin{equation*}
\xymatrix{ X  \ar[r]^(0.4){\iota} \ar[rd]_{\varphi}& \mathbbm{k}\lbrace X|H \rbrace \ar@{-->}[d]^{\bar\varphi} \\
& B}
\end{equation*}
Для всякого $f\in \mathbbm{k}\lbrace X|H \rbrace$ назовём элемент $\bar\varphi(f)\in A$
\textit{результатом подстановки} элементов $\varphi(x)\in A$ в $f$ вместо переменных $x\in X$.

Будем называть алгебру $\mathbbm{k}\lbrace X|H \rbrace$ \textit{свободной неассоциативной алгеброй на множестве $X$ с символами операторов из алгебры $H$}.\label{DefFXHbrace}
 Элементы алгебры $\mathbbm{k}\lbrace X|H \rbrace$ называются \textit{неассоциативными $H$-многочленами}
от переменных из множества $X$.

Если $A$~"--- алгебра с обобщённым $H$-действием, то пересечение $\Id^H(A)$ ядер всевозможных гомоморфизмов 
алгебр $\varphi \colon \mathbbm{k}\lbrace X|H \rbrace \to A$, удовлетворяющих условию $\varphi(x^h)=h\varphi(x)$
для всех $x\in X$ и $h\in H$,
 называется множеством \textit{полиномиальных $H$-тождеств алгебры $A$}.\label{DefIdHA}
Учитывая универсальное свойство алгебры $\mathbbm{k}\lbrace X|H \rbrace$, легко видеть, что множество $\Id^H(A)$
состоит из всех $H$-многочленов, которые обращаются в нуль при подстановке
элементов алгебры $A$ вместо своих переменных.

Является ли конкретный неассоциативный $H$-многочлен $H$-тождеством или нет, не зависит от того, какими буквами мы обозначаем переменные этого многочлена. Так как любой многочлен зависит лишь от конечного числа переменных, для описания всех полиномиальных $H$-тождеств достаточно рассматривать алгебру 
$\mathbbm{k}\lbrace X|H \rbrace$, где $X=\lbrace x_1, x_2, x_3, \ldots\rbrace$~"--- счётный набор букв, что мы и делаем ниже.
Этим соглашением устраняется неопределённость, которая заключалась в том, что идеал $\Id^H(A)$, заданный выше, зависел от множества $X$.

\begin{example}\label{ExampleIdH}
Рассмотрим $\mathbb Z/2\mathbb Z$-градуировку на алгебре $M_2(\mathbbm{k})$, заданную равенствами $M_2(\mathbbm{k})^{(\bar 0)} = \left\lbrace\left(\begin{smallmatrix}
* & 0 \\
0 & *
\end{smallmatrix}
 \right)\right\rbrace$
и 
$M_2(\mathbbm{k})^{(\bar 1)} = \left\lbrace\left(\begin{smallmatrix}
0 & * \\
* & 0
\end{smallmatrix}
 \right)\right\rbrace$, и соответствующее $(\mathbbm{k}(\mathbb Z/2\mathbb Z))^*$-действие.
 Тогда $$x^{h_0}y^{h_0}-y^{h_0}x^{h_0}\in \Id^{(\mathbbm{k}(\mathbb Z/2\mathbb Z))^*}(M_2(\mathbbm{k})),$$
 где $h_0\in (\mathbbm{k}(\mathbb Z/2\mathbb Z))^*$, $h_0(\bar 0) = 1$ и $h_0(\bar 1) = 0$.
\end{example}

Пусть $Q \subseteq \mathbbm{k}\lbrace X|H \rbrace$~"--- некоторое множество $H$-многочленов.
\textit{Многообразием $\mathbf{Var}(Q)$ алгебр с обобщённым $H$-действием, заданным множеством $Q$,} называется класс всех алгебр с обобщённым $H$-действием, в которых выполняются все тождества из множества $Q$. \textit{Идеалом $H$-тождеств $\Id^H(\mathcal V)$ многообразия $\mathcal V$} называется пересечение всех идеалов $\Id^H(A)$ для всех алгебр $A\in \mathcal V$.
Элементы множества $\Id^H(\mathbf{Var}(Q))$ называются \textit{следствиями} из $H$-многочленов множества $Q$. Если для двух подмножеств $Q_1, Q_2 \subseteq \mathbbm{k}\lbrace X|H \rbrace$
выполнено условие $\Id^H(\mathbf{Var}(Q_1))=\Id^H(\mathbf{Var}(Q_2))$,
то множества $Q_1$ и $Q_2$ называются \textit{эквивалентными}.
Легко видеть, что для любых $Q \subseteq \mathbbm{k}\lbrace X|H \rbrace$ и $f \in Q$
идеал $\Id^H(\mathbf{Var}(Q))$ включает в себя все многочлены, полученные из
$f$ линейными заменами $x_i \mapsto \alpha_{i1} x_{j_{i1}}^{h_{i1}} + \ldots \alpha_{ik} x_{j_{ik}}^{h_{ik}}$,
где $\alpha_{i\ell} \in \mathbbm{k}$, $h_{i\ell} \in H$, $j_{i\ell} \in\mathbb N$, $1\leqslant \ell \leqslant k$, $i,k \in\mathbb N$. Если $H$~"--- это просто ассоциативная алгебра с единицей, то про нелинейные замены ничего сказать нельзя, так как их результат в $\mathbbm{k}\lbrace X|H \rbrace$ не определён.
Разумеется, используя условие~\eqref{EqGeneralizedHopf}, результат любой подстановки
$x_i \mapsto f_i$, где $f_i \in \mathbbm{k}\lbrace X|H \rbrace$, в переменные любого $H$-многочлена снова переписывается в виде некоторого $H$-многочлена $g$ (это будет многократно использоваться в последующих главах), однако $H$-многочлен $g$, вообще говоря, зависит от алгебры $A$.

Даже в случае обычных полиномиальных тождеств (когда $H=\mathbbm{k}$) классификация всевозможных многообразий алгебр в зависимости от идеалов их полиномиальных тождеств представляется неподъёмной задачей. Поэтому многообразия зачастую классифицируются в зависимости от роста их числовых характеристик. Одной из важнейших числовых последовательностей, связанных с полиномиальными тождествами, является последовательность их коразмерностей.

Пусть $$W_n^H := \langle x_{\sigma(1)}^{h_1} x_{\sigma(2)}^{h_2}\ldots x_{\sigma(n)}^{h_n} \mid \sigma \in S_n,\ h_i\in H\rangle_\mathbbm{k} \subset \mathbbm{k}\lbrace X|H \rbrace,$$\label{DefWnH}
где $S_n$~"--- $n$-я группа подстановок,
$n\in\mathbb N$ и на одночленах рассматриваются всевозможные расстановки скобок. Элементы пространств $W_n^H$ называются \textit{полилинейными неассоциативными $H$-многочленами},
а элементы пространств $W_n^H \cap \Id^H(A)$ называются \textit{полилинейными $H$-тождествами} алгебры $A$.

Используя процесс линеаризации~\cite[\S~1.3]{ZaiGia}, 
нетрудно доказать, что над полем характеристики $0$ всякое полиномиальное $H$-тождество
алгебры $A$ с обобщённым $H$-действием
 эквивалентно
конечному набору полилинейных $H$-тождеств алгебры $A$.
Следовательно, в этом случае пространства
$W_n^H \cap \Id^H(A)$, где $n\in\mathbb N$,
содержат всю информацию о полиномиальных $H$-тождествах алгебры $A$. 
Число $c_n^H(A):=\dim\left(\frac{W_n^H}{W_n^H\cap\, \Id^H(A)}\right)$, где $n\in\mathbb N$, называется \textit{$n$-й коразмерностью полиномиальных $H$-тождеств}
или \textit{$n$-й $H$-коразмерностью}
алгебры $A$. \label{DefHCodim}

Предел $\PIexp^H(A):=\lim\limits_{n\rightarrow\infty} \sqrt[n]{c^H_n(A)}$\label{DefHPIexp} (если он существует)
называется \textit{экспонентой роста полиномиальных $H$-тождеств} или \textit{$H$-PI-экспонентой} алгебры $A$.

   Симметрическая группа $S_n$ действует на пространстве
   $\frac {W^H_n}{W^H_n
   \cap \Id^H(A)}$ перестановками переменных.
   Характер $\chi^H_n(A)$ этого представления
   называется $n$-м
  \textit{кохарактером} полиномиальных $H$-тождеств алгебры $A$.\label{DefHCochar}
  В случае, когда $\chr \mathbbm{k} = 0$, $n$-й кохарактер можно представить в виде суммы
  $$\chi^H_n(A)=\sum_{\lambda \vdash n}
   m(A, H, \lambda)\chi(\lambda)$$ неприводимых характеров $\chi(\lambda)$.\label{DefHMul}
   Иными словами, $m(A, H, \lambda)$~"--- это кратности неприводимых подмодулей $M(\lambda)$
   в разложении $\mathbbm{k}S_n$-модуля $\frac {W^H_n}{W^H_n
   \cap \Id^H(A)}$.
  Число $\ell_n^H(A):=\sum_{\lambda \vdash n}
   m(A, H, \lambda)$ называется $n$-й
  \textit{кодлиной} полиномиальных $H$-тождеств алгебры $A$.\label{DefHColength}

 \begin{remark}\label{RemarkOrdinaryCodim} 
 Если $A$~"--- обычная алгебра над полем $\mathbbm{k}$,
то её \textit{обычные полиномиальные тождества}, \textit{коразмерности} и \textit{кохарактеры} можно
определить как, соответственно, $H$-тождества, $H$-коразмерности и $H$-кохарактеры
алгебры $A$ с тривиальным действием алгебры $H=\mathbbm{k}$. Иными словами,
$W_n := W^\mathbbm{k}_n$,
$\Id(A):=\Id^\mathbbm{k}(A)$,
$c_n(A):=c_n^\mathbbm{k}(A)$,
$\PIexp(A):=\PIexp^\mathbbm{k}(A)$,
$\chi_n(A) := \chi^\mathbbm{k}_n(A)$,
$m(A,\lambda):= m(A,\mathbbm{k},\lambda)$,
$\ell_n(A):= \ell_n^\mathbbm{k}(A)$.
\end{remark}

 \begin{remark}\label{RemarkGCodim} 
 Если $A$~"--- алгебра над полем $\mathbbm{k}$ с действием некоторой группы $G$ автоморфизмами
 и антиавтоморфизмами, то из примера~\ref{ExampleFTEndAntiEnd}
 следует, что $A$~"--- алгебра с обобщённым $\mathbbm{k}G$-действием.
 Её \textit{полиномиальные $G$-тождества}, \textit{$G$-коразмерности} и \textit{$G$-кохарактеры}
определяются как, соответственно, $\mathbbm{k}G$-тождества, $\mathbbm{k}G$-коразмерности и $\mathbbm{k}G$-кохарактеры. Иными словами,
$W^G_n := W^{\mathbbm{k}G}_n$,
$\Id^G(A):=\Id^{\mathbbm{k}G}(A)$,
$c^G_n(A):=c_n^{\mathbbm{k}G}(A)$,
$\PIexp^G(A):=\PIexp^{\mathbbm{k}G}(A)$,
$\chi^G_n(A) := \chi^{\mathbbm{k}G}_n(A)$,
$m(A,G,\lambda):= m(A,\mathbbm{k}G,\lambda)$,
$\ell^G_n(A):= \ell_n^{\mathbbm{k}G}(A)$.
\end{remark}

\begin{remark}\label{RemarkDiffCodim} 
 Если $A$~"--- алгебра над полем $\mathbbm{k}$ с действием некоторой алгебры Ли~$\mathfrak{g}$ дифференцированиями, то из примера~\ref{ExampleUgModule}
 следует, что $A$~"--- $U(\mathfrak g)$-модульная алгебра.
 Её \textit{дифференциальные тождества}, \textit{дифференциальные коразмерности} и \textit{дифференциальные кохарактеры}
определяются как, соответственно, её $U(\mathfrak g)$-тождества, $U(\mathfrak g)$-коразмерности и $U(\mathfrak g)$-кохарактеры.
\end{remark}

 \begin{remark}\label{RemarkHActionOnWHn} Истинная важность условия~\eqref{EqGeneralizedHopf}
кроется в том факте, что при помощи~\eqref{EqGeneralizedHopf}
можно естественным образом определить структуру левого $H$-модуля на векторном пространстве $\frac {W^H_n}{W^H_{n}  \cap \Id^H(A)}$. Это делается следующим образом.
Всякий полилинейный $H$-многочлен
может рассматриваться как полилинейная функция на алгебре $A$,
принимающая значения в $A$. При этом подпространство $W^H_{n}  \cap \Id^H(A)$
является ядром соответствующего гомоморфизма $\mathbbm{k}S_n$-модулей $W^H_n \to \Hom_\mathbbm{k}(A^{{}\otimes n}; A)$, где группа $S_n$ действует перестановками аргументов.
Отсюда получаем вложение $\frac {W^H_n}{W^H_{n}  \cap \Id^H(A)} \subseteq \Hom_\mathbbm{k}(A^{{}\otimes n}; A)$.
Пространство $\Hom_\mathbbm{k}(A^{{}\otimes n}; A)$ является левым $H$-модулем, причём
операторы из алгебры $H$ коммутируют с операторами из группы $S_n$: $(hg)(a_1, \ldots, a_n):= hg(a_1, \ldots, a_n)$ при $h\in H$,
$g\in \Hom_\mathbbm{k}(A^{{}\otimes n}; A)$ и $a_1, \ldots, a_n \in A$.
Если $f\in W^H_{n}$ и $h\in H$,
мы можем, применив несколько раз равенство~\eqref{EqGeneralizedHopf},
переписать функцию $h\bar f$, где $\bar f$~"--- образ $f$ в $\frac {W^H_n}{W^H_{n}  \cap \Id^H(A)}$,
в виде линейной комбинации таких произведений переменных $x_i$,
в которых операторы из $H$
применяются только к самим переменным $x_i$, но не к их произведениям.
Другими словами, функция $h\bar f$ на $A$
может быть представлена (необязательно однозначно)
как функция, соответствующая $H$-многочлену из $W^H_n$.  
Следовательно, $\frac {W^H_n}{W^H_{n}  \cap \Id^H(A)}$ является $H$-подмодулем
$H$-модуля $\Hom_\mathbbm{k}(A^{{}\otimes n}; A)$.
Обозначим через $f^h$ любой из таких $H$-многочленов, что $\overline{f^h}=h\bar f$.
 Если $\chr \mathbbm{k} = 0$,  $f^h \notin \Id^H(A)$ и $\mathbbm{k}S_n \bar f \cong M(\lambda)$ для некоторого $\lambda \vdash n$, то $\mathbbm{k}S_n h\bar f$ является ненулевым гомоморфным образом неприводимого $\mathbbm{k}S_n$-модуля $M(\lambda)$.
Отсюда и $\mathbbm{k}S_n h\bar f \cong M(\lambda)$. Это свойство будет использовано ниже в \S\ref{SectionHPIexpExistHSimple} и последующих главах.
 \end{remark}

Приведём пример такой бесконечномерной $H$-модульной ассоциативной алгебры $A$ для бесконечномерной
алгебры Хопфа $H$, что $H$-коразмерности $c_n^H(A)$ бесконечны:

\begin{example}\label{ExampleInfCodim} Пусть $\mathbbm{k}$~"--- некоторое поле.
Обозначим через $C=\langle c \rangle$ бесконечную циклическую группу,
а через $(G, \cdot)$~"--- группу $(\mathbb Q, +)$, записанную в мультипликативной форме.
При этом через $g^\alpha$, где $\alpha \in \mathbb Q$,
условимся обозначать элемент группы $G$, соответствующий числу $\alpha \in \mathbb Q$.
 Зафиксируем некоторое число $m \in\mathbb N$,
где $m \geqslant 2$. Определим действие группы~$C$ на
групповой алгебре $A: = \mathbbm{k}G$ автоморфизмами при помощи формулы $cg^\alpha:=g^{m\alpha}$
для всех $\alpha\in \mathbb Q$. 
 Покажем, что, несмотря на то, что алгебра $A$ коммутативна,
справедливо равенство 
$c_n^{\mathbbm{k}C}(A) = +\infty$ для всех $n\in\mathbb N$.
Для произвольного $n\in\mathbb N$
рассмотрим полилинейные $\mathbbm{k}C$-многочлены $$f_k(x_1, \ldots, x_n):=x_1^{c^k} x_2 \ldots x_n\text{, где }k\in\mathbb Z.$$
(Поскольку алгебра $A$ ассоциативна, расстановка скобок на одночленах неважна. Можно для определённости считать одночлены левонормированными.)
 Тогда $$f_k(g,g,\ldots, g)=(c^k g)g\ldots g = g^{m^k+(n-1)},$$
т.е.  $\mathbbm{k}C$-многочлены $f_k$ линейно независимы по модулю $\Id^{\mathbbm{k}C}(A)$.
\end{example} 

Пример~\ref{ExampleInfCodim} показывает, что при изучении $H$-коразмерностей имеет смысл ограничиться
случаем, когда одна из двух алгебр $A$ и $H$ конечномерна.

Докажем теперь оценку сверху для коразмерностей полиномиальных $H$-тождеств
конечномерных алгебр:

 \begin{proposition}\label{PropositionCodimDim}
Пусть $A$~"--- конечномерная алгебра с обобщённым
$H$-действием для некоторой ассоциативной алгебры $H$ c единицей
над произвольным полем~$\mathbbm{k}$.
Тогда $c_n^H(A) \leqslant (\dim A)^{n+1}$ для всех  $n \in \mathbb N$.
\end{proposition}
\begin{proof}
Снова рассмотрим $H$-многочлены как $n$-линейные отображения из $A$ в $A$.
Из вложения $\frac {W^H_n}{W^H_{n}  \cap \Id^H(A)} \subseteq \Hom_\mathbbm{k}(A^{{}\otimes n}; A)$ следует неравенство
 $$c^H_n(A)=\dim \left(\frac{W^H_n}{W^H_{n}  \cap \Id^H(A)}\right)
\leqslant \dim \Hom_{\mathbbm{k}}(A^{{}\otimes n}; A)=(\dim A)^{n+1}.$$
\end{proof}
\begin{corollary}
Пусть $A$~"--- конечномерная алгебра над произвольным полем $\mathbbm{k}$ с действием некоторой группы $G$ автоморфизмами
 и антиавтоморфизмами. Тогда
$c_n^G(A) \leqslant (\dim A)^{n+1}$ для всех $n \in \mathbb N$.
\end{corollary}

Предложение ниже устанавливает связь между обычными и $H$-коразмерностями:

\begin{proposition}\label{PropositionOrdinaryAndHopf}
Пусть $A$~"--- алгебра с обобщённым
$H$-действием для некоторой ассоциативной алгебры $H$ c единицей
над произвольным полем~$\mathbbm{k}$, а $\zeta \colon H \to \End_\mathbbm{k}(A)$~"--- гомоморфизм, задающий
это $H$-действие. Тогда
$$c_n(A) \leqslant c^{H}_n(A)
  \leqslant (\dim \zeta(H))^n c_n(A) \text{ для всех } n \in \mathbb N.$$
\end{proposition}
\begin{proof}
Снова рассмотрим многочлены как $n$-линейные отображения из $A$ в $A$
и отождествим  $\frac{W_n}{W_n \cap \Id(A)}$ и $\frac{W^{H}_n}{W^{H}_n \cap \Id^H(A)}$
с соответствующими подпространствами в $\Hom_{\mathbbm{k}}(A^{{}\otimes n}; A)$.
Тогда
 $$\frac{W_n}{W_n \cap \Id(A)} \subseteq \frac{W^{H}_n}{W^{H}_n \cap \Id^H(A)}
\subseteq \Hom_{\mathbbm{k}}(A^{{}\otimes n}; A),$$
откуда и следует оценка снизу.

Выберем такие $f_1, \ldots, f_t \in W_n$,
что их образы являются базисом пространства 
$\frac{W_n}{W_n \cap \Id(A)}$.  Тогда для любого одночлена $x_{\sigma(1)}x_{\sigma(2)}
\ldots x_{\sigma(n)}$ (с некоторой фиксированной расстановкой скобок $\Xi$),  где $\sigma \in S_n$, существуют такие коэффициенты $\alpha_{i,\sigma, \Xi} \in \mathbbm{k}$, что
\begin{equation}\label{EqWnModuloId}x_{\sigma(1)}x_{\sigma(2)}
\ldots x_{\sigma(n)} - \sum_{i=1}^t \alpha_{i,\sigma, \Xi} f_i(x_1, \ldots, x_n) \in \Id(A).
\end{equation}

В случае, когда $\dim \zeta(H)=+\infty$, оценка сверху становится тривиальной, поэтому без ограничения общности можно считать, что пространство $\zeta(H)$ конечномерно.
Пусть $\bigl(\zeta(\gamma_j)\bigr)_{j=1}^m$, где $\gamma_j \in H$,"--- базис пространства $\zeta(H)$.
Тогда для любого $h \in H$ существуют такие $\alpha_j \in \mathbbm{k}$,
что $\zeta(h) = \sum_{j=1}^m \alpha_j \zeta(\gamma_j)$ и 
\begin{equation}\label{Eqhtogammaj}
x^h - \sum_{j=1}^m \alpha_j x^{\gamma_j} \in \Id^H(A).
\end{equation}

Отсюда линейная оболочка
$H$-многочленов $x^{\gamma_{i_1}}_{\sigma(1)}x^{\gamma_{i_2}}_{\sigma(2)}
\ldots x^{\gamma_{i_n}}_{\sigma(n)}$, где $\sigma \in S_n$, $1 \leqslant i_j \leqslant m$, со всевозможными расстановками скобок, совпадает по модулю $\Id^H(A)$ с  $W^H_n$. При этом из (\ref{EqWnModuloId}) следует,
что в случае, когда на одночлене $x^{\gamma_{i_1}}_{\sigma(1)}x^{\gamma_{i_2}}_{\sigma(2)}
\ldots x^{\gamma_{i_n}}_{\sigma(n)}$ фиксирована расстановка
скобок~$\Xi$, выполняется условие
 $$x^{\gamma_{i_1}}_{\sigma(1)}x^{\gamma_{i_2}}_{\sigma(2)}
\ldots x^{\gamma_{i_n}}_{\sigma(n)} - \sum_{i=1}^t \alpha_{i,\sigma, \Xi} f_i(x^{\gamma_{i_1}}_1,
 \ldots, x^{\gamma_{i_n}}_n) \in \Id^H(A).$$
Отсюда любой $H$-многочлен из $W^H_n$ может быть представлен по модулю $\Id^H(A)$
в виде линейной комбинации $H$-многочленов $f_i(x^{\gamma_{i_1}}_1,
 \ldots, x^{\gamma_{i_n}}_n)$. Число таких $H$-многочленов равно $m^n t = (\dim \zeta(H))^n c_n(A)$,
 что и завершает доказательство оценки сверху.
\end{proof}
 
 Выше понятие полиномиального $H$-тождества было введено для необязательно ассоциативной алгебры с 
обобщённым $H$-действием. Однако при работе с ассоциативными алгебрами вместо неассоциативных свободных алгебр, как правило, используются ассоциативные свободные алгебры. Как мы увидим в конце параграфа, числовые характеристики полиномиальных $H$-тождеств при этом подходе остаются прежними.

Пусть $X$~"--- множество, а $\mathbbm{k}$~"--- поле. Тогда свободная ассоциативная алгебра $\mathbbm{k}\langle X \rangle$
без единицы раскладывается в прямую сумму подпространств $ \bigoplus_{n=1}^\infty \mathbbm{k}\langle X \rangle^{(n)}$,
где $\mathbbm{k}\langle X \rangle^{(n)}$~"--- линейная оболочка одночленов, имеющих суммарную степень $n$  по всем буквам.

Пусть $H$~"--- ассоциативная алгебра с единицей. Обозначим через $\mathbbm{k}\langle X|H \rangle$
алгебру, которая как векторное пространство совпадает с $$\bigoplus_{n=1}^\infty \mathbbm{k}\langle X \rangle^{(n)} \otimes H^{{}\otimes n},$$ 
а умножение задаётся формулой
$(v_1 \otimes w_1)(v_2 \otimes w_2):=v_1 v_2 \otimes w_1 \otimes w_2$
для $v_1\in \mathbbm{k}\langle X \rangle^{(k)}$, $v_2\in \mathbbm{k}\langle X \rangle^{(\ell)}$,
$w_1 \in H^{{}\otimes k}$, $w_2 \in H^{{}\otimes \ell}$, $k,\ell\in\mathbb N$.

Как и в случае неассоциативных многочленов, будем использовать обозначение $x_1^{h_1}\ldots x_n^{h_n} := x_1\ldots x_n 
\otimes h_1 \otimes h_2 \otimes \ldots \otimes h_n$
для всех $x_1, \ldots, x_n\in X$ и $h_1,\ldots, h_n \in H$.

Пусть $(h_{\alpha})_{\alpha\in \Lambda}$~"--- базис алгебры $H$. Тогда $\mathbbm{k}\langle X|H \rangle$
 как алгебра изоморфна свободной ассоциативной алгебре без единицы на множестве 
$\lbrace x^{h_\alpha} \mid x\in X,\  \alpha\in \Lambda\rbrace$.

Отождествим $X$ с подмножеством $$\lbrace x^1 \mid x\in X\rbrace \subseteq \mathbbm{k}\langle X|H \rangle.$$ Обозначим через $\iota_0 \colon X \to \mathbbm{k}\langle X|H \rangle$ соответствующее вложение.
Тогда $\mathbbm{k}\langle X|H \rangle$ удовлетворяет следующему универсальному свойству:
для любого отображения $\varphi \colon X \to A$, где $A$~"---  ассоциативная алгебра с обобщённым $H$-действием, существует единственный гомоморфизм алгебр $\bar \varphi \colon \mathbbm{k}\langle X | H \rangle \to A$,
такой, что для всех $h\in H$ и $x\in X $ справедливо равенство $\bar\varphi(x^h)=h\bar\varphi(x)$ и 
 диаграмма ниже коммутативна:

\begin{equation*}
\xymatrix{ X  \ar[r]^(0.4){\iota_0} \ar[rd]_{\varphi}& \mathbbm{k}\langle X|H \rangle \ar@{-->}[d]^{\bar\varphi} \\
& B}
\end{equation*}

Будем называть алгебру $\mathbbm{k}\langle X|H \rangle$ \textit{свободной ассоциативной алгеброй без единицы на множестве $X$ с символами операторов из алгебры $H$}.\label{DefFXHangle}

Элементы алгебры $\mathbbm{k}\langle X|H \rangle$ называются \textit{ассоциативными $H$-многочленами}
от переменных из множества $X$.
Если $A$~"--- ассоциативная алгебра с обобщённым $H$-действием, то пересечение $\Id_{\mathrm{assoc}}^H(A)$ ядер всевозможных гомоморфизмов 
алгебр $\varphi \colon \mathbbm{k}\langle X|H \rangle \to A$, удовлетворяющих условию $\varphi(x^h)=h\varphi(x)$
для всех $x\in X$ и $h\in H$,
 называется множеством \textit{ассоциативных полиномиальных $H$-тождеств алгебры $A$}.
Учитывая универсальное свойство алгебры $\mathbbm{k}\langle X|H \rangle$, легко видеть, что множество $\Id_{\mathrm{assoc}}^H(A)$
состоит из всех ассоциативных $H$-многочленов, которые обращаются в нуль при подстановке
элементов алгебры $A$ вместо своих переменных.

Как и в случае неассоциативных тождеств, можно считать, что 
$X=\lbrace x_1, x_2, x_3, \ldots\rbrace$.

Пусть $$P_n^H := \langle x_{\sigma(1)}^{h_1} x_{\sigma(2)}^{h_2}\ldots x_{\sigma(n)}^{h_n} \mid \sigma \in S_n,\ h_i\in H\rangle_\mathbbm{k} \subset \mathbbm{k}\langle X|H \rangle.$$\label{DefPnH}
 Элементы пространств $P_n^H$ называются \textit{полилинейными ассоциативными $H$-многочленами},
а элементы пространств $P_n^H \cap \Id_{\mathrm{assoc}}^H(A)$ называются \textit{полилинейными ассоциативными $H$-тождествами} алгебры $A$.

Как мы видим, определение полиномиального $H$-тождества в ассоциативной алгебре $A$ существенно зависит от того, рассматриваем ли мы алгебру $A$ как ассоциативную или как необязательно ассоциативную.
Однако как всякому неассоциативному, так и всякому ассоциативному $H$-многочлену можно поставить в соответствие функцию,
которая сопоставляет набору элементов алгебры $A$ значение $H$-многочлена
на этом наборе элементов. При этом нулевые функции будут соответствовать полиномиальным $H$-тождествам. Отсюда соответствие $x_i^h \mapsto x_i^h$, $i\in \mathbb N$, $h\in H$,
индуцирует изоморфизм
$\mathbbm{k} \lbrace X | H\rbrace/\Id^H(A) \cong \mathbbm{k}\langle X | H \rangle/ \Id^H_{\mathrm{assoc}}(A)$ алгебр
и изоморфизм $H$- и $\mathbbm{k}S_n$-модулей $\frac {W^H_n}{W^H_{n}
  \cap \Id^H(A)} \cong \frac {P^H_n}{P^H_{n}
  \cap \Id_{\mathrm{assoc}}^H(A)}$. Это означает, что можно отождествить ассоциативные и неассоциативные <<нетождества>>.

Следовательно, при всех $n\in\mathbb N$ справедливо равенство $c_n^H(A)=\dim\left(\frac{P_n^H}{P_n^H\cap\, \Id^H_{\mathrm{assoc}}(A)}\right)$, и определение $n$-й коразмерности, $n$-го кохарактера и $n$-й кодлины полиномиальных $H$-тождеств
ассоциативной алгебры $A$ не зависит от того, используем ли мы алгебру $\mathbbm{k} \lbrace X | H\rbrace$ или
алгебру $\mathbbm{k}\langle X | H \rangle$.

Ниже при работе с ассоциативными алгебрами идеал $\Id^H_{\mathrm{assoc}}(A)$ для краткости также обозначается через $\Id^H(A)$, а полиномиальные $H$-многочлены берутся из алгебры $\mathbbm{k}\langle X | H \rangle$. 

\section{$H$-тождества $H$-модульных алгебр}\label{SectionHIdHMod}

Рассмотрим конструкции предыдущего параграфа в случае модульных алгебр над алгебрами Хопфа.

Итак, пусть $H$~"--- алгебра Хопфа, а $X$~"--- некоторое множество. Тогда алгебра $\mathbbm{k}\lbrace X|H \rbrace$ является левой $H$-модульной алгеброй, где действие алгебры $H$ задаётся формулой $$h \cdot \left(v\otimes h_1 \otimes \ldots \otimes h_n\right) := v
\otimes h_{(1)}h_1 \otimes \ldots \otimes h_{(n)}h_n \text{ для } v\in \mathbbm{k}\lbrace X \rbrace^{(n)}
\text{ и } h_i\in H.$$

Для $H$-модульных алгебр универсальное свойство алгебры $H$ записывается следующим образом:
для любого отображения $\varphi \colon X \to A$, где $A$~"--- необязательно ассоциативная $H$-модульная алгебра,
существует единственный гомоморфизм $H$-модульных алгебр $\bar \varphi \colon \mathbbm{k}\lbrace X | H \rbrace \to A$,
такой, что диаграмма ниже коммутативна:

\begin{equation*}
\xymatrix{ X  \ar[r]^(0.4){\iota} \ar[rd]_{\varphi}& \mathbbm{k}\lbrace X|H \rbrace \ar@{-->}[d]^{\bar\varphi} \\
& B}
\end{equation*}

Иными словами,
алгебра $\mathbbm{k}\lbrace X|H \rbrace$ является в этом случае \textit{свободной неассоциативной $H$-модульной алгеброй},
а $\mathbbm{k}\lbrace - | H \rbrace$ 
оказывается левым сопряжённым функтором к забывающему функтору из
категории необязательно ассоциативных $H$-модульных алгебр в категорию множеств.
Действительно, пусть $\mathbf{NAAlg}({}_H\mathcal M)$~"--- категория, объектами которой являются
необязательно ассоциативные $H$-модульные алгебры, а морфизмами~"--- всевозможные гомоморфизмы $H$-модульных алгебр. Если обозначить через $U \colon \mathbf{NAAlg}({}_H\mathcal M) \to \mathbf{Sets}$ забывающий функтор, то мы получаем биекцию
$$\mathbf{NAAlg}({}_H\mathcal M)(\mathbbm{k}\lbrace X|H \rbrace, A)\cong \mathbf{Sets}(X, UA),$$
естественную по множеству $X$ и алгебре $A$,
существование которой и говорит о том, что функторы $\mathbbm{k}\lbrace -|H \rbrace$ и $U$ являются сопряжёнными. 

Если $A$~"--- $H$-модульная алгебра, то идеал её $H$-тождеств
 $\Id^H(A)$ является таким $H$-инвариантным идеалом алгебры $\mathbbm{k}\lbrace X|H \rbrace$, что $\psi(\Id^H(A))\subseteq \Id^H(A)$
 для любого эндоморфизма $\psi \colon \mathbbm{k}\lbrace X|H \rbrace \to \mathbbm{k}\lbrace X|H \rbrace$,
 сохраняющего на $\mathbbm{k}\lbrace X|H \rbrace$ структуру 
 алгебры и  $H$-модуля.
Обратно, если $I$~"--- идеал алгебры $\mathbbm{k}\lbrace X|H \rbrace$, инвариантный относительно
действия алгебры $H$ и эндоморфизмов $H$-модульной алгебры $\mathbbm{k}\lbrace X|H \rbrace$,
то $\Id^H(\mathbbm{k}\lbrace X|H \rbrace/I)=I$.

В случае, когда $H$~"--- алгебра Хопфа, мы будем рассматривать многообразия, содержащие лишь $H$-модульные
алгебры. В частности, если $Q \subseteq \mathbbm{k}\lbrace X|H \rbrace$, символом $\mathbf{Var}(Q)$ мы будем обозначать \textit{многообразие $H$-модульных алгебр, заданное множеством $Q$,} т.е. класс $H$-модульных алгебр $A$, таких, что $Q\subseteq \Id^H(A)$.

\textit{Идеалом $H$-тождеств $\Id^H(\mathcal V)$ многообразия $\mathcal V$} называется пересечение всех идеалов $\Id^H(A)$ для всех алгебр $A\in \mathcal V$.
Элементы множества $\Id^H(\mathbf{Var}(Q))$ называются \textit{следствиями} из $H$-многочленов множества $Q$. Если для двух подмножеств $Q_1, Q_2 \subseteq \mathbbm{k}\lbrace X|H \rbrace$
выполнено условие $\Id^H(\mathbf{Var}(Q_1))=\Id^H(\mathbf{Var}(Q_2))$,
то множества $Q_1$ и $Q_2$ называются \textit{эквивалентными}.

Определения следствия из множества $H$-тождеств и эквивалентных множеств $H$-тождеств, данные выше, являются дословным повторением  аналогичных определений, сформулированных в предыдущем параграфе для алгебр с обобщённым $H$-действием. Однако благодаря тому, что теперь мы рассматриваем многообразия лишь $H$-модульных алгебр, эти определения, вообще говоря, неэквивалентны определениям, данным в предыдущем параграфе. В силу того, что теперь алгебра $\mathbbm{k}\lbrace X|H \rbrace$ сама является $H$-модульной, вместо переменных любого $H$-многочлена допускается подставлять любые 
$H$-многочлены. Отсюда следствиями конкретного $H$-многочлена $f$ являются в точности $H$-многочлены,
принадлежащие $H$-инвариантному идеалу алгебры $\mathbbm{k}\lbrace X|H \rbrace$, порождённому результатами всевозможных подстановок $H$-многочленов вместо переменных $H$-многочлена $f$.

Если $\mathcal V$~"--- некоторое многообразие необязательно ассоциативных $H$-модульных алгебр, то $\mathbbm{k}\lbrace X|H \rbrace / \Id^H(\mathcal V)$ называется \textit{относительно свободной $H$-модульной алгеброй многообразия $\mathcal V$}.
Легко видеть, что для любой $H$-модульной алгебры $A\in \mathcal V$ и для любого отображения $\varphi \colon X \to A$
существует единственный гомоморфизм $H$-модульных алгебр $\bar \varphi \colon \mathbbm{k}\lbrace X | H \rbrace/ \Id^H(\mathcal V) \to A$,
такой, что диаграмма ниже коммутативна:

\begin{equation*}
\xymatrix{ X  \ar[r]^(0.25){\bar\iota} \ar[rd]_{\varphi}& \mathbbm{k}\lbrace X|H \rbrace/ \Id^H(\mathcal V) \ar@{-->}[d]^{\bar\varphi} \\
& A}
\end{equation*}
(Здесь через $\bar\iota$ обозначено отображение, являющееся композицией отображения
$\iota \colon X \to \mathbbm{k}\lbrace X|H \rbrace$ и гомоморфизма $\mathbbm{k}\lbrace X|H \rbrace \twoheadrightarrow \mathbbm{k}\lbrace X|H \rbrace / \Id^H(\mathcal V)$.)

Если рассматривать $\mathcal V$ как полную подкатегорию категории $\mathbf{NAAlg}({}_H\mathcal M)$, объектами которой являются элементы класса $\mathcal V$ и обозначить через $U \colon \mathcal V \to \mathbf{Sets}$ забывающий функтор,
то мы получаем биекцию
$$\mathcal V(\mathbbm{k}\lbrace X|H \rbrace / \Id^H(\mathcal V), A)\cong \mathbf{Sets}(X, UA),$$
естественную по множеству $X$ и алгебре $A\in\mathcal V$,
существование которой говорит о том, что функторы $\mathbbm{k}\lbrace -|H \rbrace / \Id^H(\mathcal V)$ и $U$ являются сопряжёнными. 

Как и в случае ассоциативных алгебр с обобщённым $H$-действием, в случае ассоциативных $H$-модульных алгебр вместо алгебры $\mathbbm{k}\lbrace X | H \rbrace$ можно использовать алгебру $\mathbbm{k}\langle X | H \rangle$,
которая оказывается \textit{свободной ассоциативной $H$-модульной алгеброй без единицы}.

Обратимся теперь к случаю алгебр Ли. Фиксируем множество $X$ и обозначим через $I$ идеал алгебры $\mathbbm{k}\lbrace X \rbrace$ порождённый элементами $(ab)c+(bc)a+(ca)b$ и $aa$ для всех $a,b,c \in \mathbbm{k}\lbrace X \rbrace$.
Очевидно, алгебра $L(X):=\mathbbm{k}\lbrace X \rbrace/I$ является алгеброй Ли.\label{DefLX} Алгебра
$L(X)$ называется \textit{свободной алгеброй Ли} на множестве $X$, поскольку
для любой алгебры $L$ и для любого отображения $\varphi \colon X \to L$
существует единственный гомоморфизм алгебр Ли $\bar \varphi \colon L(X) \to L$,
такой, что диаграмма ниже коммутативна:

\begin{equation*}
\xymatrix{ X  \ar[r]^(0.4){\iota_1} \ar[rd]_{\varphi}& L(X) \ar@{-->}[d]^{\bar\varphi} \\
& L}
\end{equation*}
(Здесь через $\iota_1$ обозначено отображение, являющееся композицией вложения
$X \hookrightarrow \mathbbm{k}\lbrace X \rbrace$ и гомоморфизма $\mathbbm{k}\lbrace X \rbrace \twoheadrightarrow L(X)$.)

Пусть $H$~"--- алгебра Хопфа. Обозначим через $I_0$ пересечение $H$-инвариантных идеалов
алгебры $\mathbbm{k}\lbrace X|H \rbrace$, содержащих элементы
$(ab)c+(bc)a+(ca)b$ и $aa$ для всех $a,b,c \in \mathbbm{k}\lbrace X|H \rbrace$.
Тогда алгебра $L(X|H):=\mathbbm{k}\lbrace X|H \rbrace/I_0$\label{DefLXH} называется \textit{свободной $H$-модульной алгеброй Ли} на множестве $X$, поскольку
для любой $H$-модульной алгебры $L$ и для любого отображения $\varphi \colon X \to L$
существует единственный гомоморфизм $H$-модульных алгебр Ли $\bar \varphi \colon L(X|H) \to L$,
такой, что диаграмма ниже коммутативна:

\begin{equation*}
\xymatrix{ X  \ar[r]^(0.4){\iota_2} \ar[rd]_{\varphi}& L(X|H) \ar@{-->}[d]^{\bar\varphi} \\
& L}
\end{equation*}
(Здесь через $\iota_2$ обозначено отображение, являющееся композицией вложения
$X \hookrightarrow \mathbbm{k}\lbrace X|H \rbrace$ и гомоморфизма $\mathbbm{k}\lbrace X|H \rbrace \twoheadrightarrow L(X|H)$.)

Для обозначения умножения в алгебре $L(X|H)$ также, как и в других алгебрах Ли, используется символ коммутатора:
$[f,g]:=fg$.

\begin{remark} Обозначим через $I_1$ идеал алгебры $\mathbbm{k}\lbrace X|H \rbrace$,
порождённый элементами $(ab)c+(bc)a+(ca)b$ и $ab+ba$ для всех $a,b,c \in \mathbbm{k}\lbrace X|H \rbrace$. 
 Если алгебра Хопфа $H$ кокоммутативна, т.е. $h_{(1)}\otimes h_{(2)}=h_{(2)}\otimes h_{(1)}$
 для всех $h\in H$, то идеал $I_1$ является $H$-подмодулем. Если $\chr \mathbbm{k}\ne 2$,
 то тождества $[x,x]\equiv 0$ и $[x,y]+[y,x]\equiv 0$ эквивалентны,
 т.е. $I_1 =I_0$. В этом случае  $L(X|H)$
 как алгебра Ли изоморфна свободной алгебре Ли на множестве 
$\lbrace x^{h_\alpha} \mid x\in X,\  \alpha\in \Lambda\rbrace$, где $(h_{\alpha})_{\alpha\in \Lambda}$~"--- базис алгебры $H$.
Однако, если $h_{(1)} \otimes h_{(2)} \ne h_{(2)} \otimes h_{(1)}$ для некоторого $h \in H$,
мы по-прежнему имеем $$[x^{h_{(1)}}, y^{h_{(2)}}]=h[x, y]=-h[y, x]=-[y^{h_{(1)}}, x^{h_{(2)}}]
= [x^{h_{(2)}}, y^{h_{(1)}}]$$ в $L(X | H)$ для всех $x,y \in X$,
т.е. в случае, когда множество $X$ состоит более чем из одного элемента и $h_{(1)} \otimes h_{(2)} \ne h_{(2)} \otimes h_{(1)}$, алгебра Ли $L(X | H)$ 
не является свободной в качестве обычной алгебры Ли.
\end{remark}

Элементы алгебры $L(X | H)$ называются \textit{лиевскими $H$-многочленами}
от переменных из множества $X$.
Если $L$~"--- $H$-модульная алгебра, то пересечение $\Id_{\mathrm{Lie}}^H(L)$ ядер всевозможных гомоморфизмов $H$-модульных алгебр Ли $\varphi \colon L(X | H) \to L$
 называется множеством \textit{лиевских полиномиальных $H$-тождеств алгебры $L$}.
Учитывая универсальное свойство алгебры $L(X | H)$, легко видеть, что множество $\Id_{\mathrm{Lie}}^H(L)$
состоит из всех лиевских $H$-многочленов, которые обращаются в нуль при подстановке
элементов алгебры $L$ вместо своих переменных. 
 Множество $\Id_{\mathrm{Lie}}^H(L)$ является таким $H$-инвариантным идеалом алгебры $L(X | H)$, что $\psi(\Id_{\mathrm{Lie}}^H(L))\subseteq \Id_{\mathrm{Lie}}^H(L)$
 для любого эндоморфизма $\psi \colon L(X | H) \to L(X | H)$,
 сохраняющего на $L(X | H)$ структуру 
 алгебры Ли и  $H$-модуля.
Обратно, если $I$~"--- идеал алгебры Ли $L(X | H)$, инвариантный относительно
действия алгебры $H$ и эндоморфизмов $H$-модульной алгебры Ли $L(X | H)$,
то $\Id_{\mathrm{Lie}}^H(L(X | H)/I)=I$.

Как и выше, для описания всех лиевских полиномиальных $H$-тождеств достаточно положить
$X=\lbrace x_1, x_2, x_3, \ldots\rbrace$.

Пусть $$V_n^H := \langle [x_{\sigma(1)}^{h_1}, x_{\sigma(2)}^{h_2},\ldots, x_{\sigma(n)}^{h_n}] \mid \sigma \in S_n,\ h_i\in H\rangle_\mathbbm{k} \subset L(X | H).$$\label{DefVnH}
 Элементы пространств $V_n^H$ называются \textit{полилинейными лиевскими $H$-многочленами},
а элементы пространств $V_n^H \cap \Id_{\mathrm{Lie}}^H(L)$ называются \textit{полилинейными лиевскими $H$-тождествами} алгебры $L$.

Если рассмотреть функции на алгебре $L$ со значениями в алгебре $L$,
соответствующие неассоциативным и лиевским $H$-многочленам,
то соответствие $x_i^h \mapsto x_i^h$, $i\in \mathbb N$, $h\in H$,
индуцирует изоморфизм
$\mathbbm{k} \lbrace X | H\rbrace/ \Id^H(L) \cong L(X | H)/\Id^H_{\mathrm{Lie}}(L)$ алгебр Ли
и изоморфизм $H$- и $\mathbbm{k}S_n$-модулей $\frac {W^H_n}{W^H_{n}
  \cap \Id^H(L)} \cong \frac {V^H_n}{V^H_{n}
  \cap \Id_{\mathrm{Lie}}^H(L)}$. Это означает, что для алгебры Ли можно отождествить неассоциативные и лиевские <<нетождества>>.

Следовательно, для всех $n\in\mathbb N$ справедливо равенство $c_n^H(L)=\dim\left(\frac{V_n^H}{V_n^H\cap\, \Id^H_{\mathrm{Lie}}(L)}\right)$ и определение $n$-й коразмерности, $n$-го кохарактера и $n$-й кодлины полиномиальных $H$-тождеств
алгебры Ли $L$ не зависит от того, используем ли мы алгебру $\mathbbm{k} \lbrace X | H\rbrace$ или
алгебру $L(X | H)$.

Ниже при работе с алгебрами Ли идеал $\Id^H_{\mathrm{Lie}}(L)$ для краткости также обозначается через $\Id^H(L)$, а $H$-многочлены берутся из алгебры $L(X | H)$.

\section{Градуированные полиномиальные тождества}\label{SectionGradedPI}

Пусть $T$~"--- множество, а $\mathbbm{k}$~"--- поле.

Рассмотрим (абсолютно) свободную неассоциативную алгебру $\mathbbm{k}\lbrace X^{T\text{-}\mathrm{gr}} \rbrace$
на объединении $$X^{T\text{-}\mathrm{gr}}:=\bigsqcup_{t \in T}X^{(t)}$$
непересекающихся множеств $X^{(t)} = \{ x^{(t)}_1,
x^{(t)}_2, \ldots \}$ \label{DefXTgr}

Говорят, что $f=f(x^{(t_1)}_{i_1}, \ldots, x^{(t_s)}_{i_s})$~"--- \textit{градуированное полиномиальное тождество} для $T$-градуированной алгебры $A=\bigoplus_{t\in T}
A^{(t)}$ и пишут $f\equiv 0$,
если $f(a^{(t_1)}_1, \ldots, a^{(t_s)}_s)=0$
для всех $a^{(t_j)}_j \in A^{(t_j)}$, $1 \leqslant j \leqslant s$.
  Множество $\Id^{T\text{-}\mathrm{gr}}(A)$ градуированных полиномиальных
  тождеств алгебры $A$ является
идеалом алгебры $\mathbbm{k}\lbrace
 X^{T\text{-}\mathrm{gr}}\rbrace$.
 
 \begin{remark}
 Если множество $T$ не наделено структурой (полу)группы, мы не рассматриваем никакой градуировки на самой
 алгебре $\mathbbm{k}\lbrace
 X^{T\text{-}\mathrm{gr}}\rbrace$. Если $T$~"--- полугруппа, то  алгебра $\mathbbm{k}\lbrace
 X^{T\text{-}\mathrm{gr}}\rbrace$ также является $T$-градуированной,
 где $T$-степень всякого одночлена является произведением $T$-степеней множителей.
 В этом случае всякая подстановка однородных элементов алгебры $A$ вместо переменных соответствующих $T$-степеней индуцирует (строго) градуированный гомоморфизм $\mathbbm{k}\lbrace
 X^{T\text{-}\mathrm{gr}}\rbrace \to A$.
 \end{remark}

\begin{example}\label{ExampleIdGr}
Рассмотрим полугруппу $T=\mathbb Z_2 = \lbrace \bar 0, \bar 1 \rbrace$
по умножению и $T$-градуировку $\UT_2(\mathbbm{k})=\UT_2(\mathbbm{k})^{(\bar 0)}\oplus \UT_2(\mathbbm{k})^{(\bar 1)}$
на алгебре $\UT_2(\mathbbm{k})$ верхнетреугольных матриц $2\times 2$ над полем $\mathbbm{k}$,
заданную равенствами
 $\UT_2(\mathbbm{k})^{(\bar 1)}=\left(
\begin{array}{cc}
\mathbbm{k} & 0 \\
0 & \mathbbm{k}
\end{array}
 \right)$ и $\UT_2(\mathbbm{k})^{(\bar 0)}=\left(
\begin{array}{cc}
0 & \mathbbm{k} \\
0 & 0
\end{array}
 \right)$. Тогда $$[x^{(\bar 1)}, y^{(\bar 1)}]:=x^{(\bar 1)} y^{(\bar 1)} - y^{(\bar 1)} x^{(\bar 1)}
\in \Id^{T\text{-}\mathrm{gr}}(\UT_2(\mathbbm{k}))$$
и $x^{(\bar 0)} y^{(\bar 0)} 
\in \Id^{T\text{-}\mathrm{gr}}(\UT_2(\mathbbm{k}))$.
\end{example}

Пусть
$$W^{T\text{-}\mathrm{gr}}_n := \langle x^{(t_1)}_{\sigma(1)}
x^{(t_2)}_{\sigma(2)}\ldots x^{(t_n)}_{\sigma(n)}
\mid t_i \in T, \sigma\in S_n \rangle_\mathbbm{k} \subset \mathbbm{k} \lbrace X^{T\text{-}\mathrm{gr}} \rbrace$$
\label{DefWTgrn}
(на одночленах рассматриваются всевозможные расстановки скобок), $n \in \mathbb N$.
Элементы пространства $W^{T\text{-}\mathrm{gr}}_n$ называются \textit{полилинейными $T$-градуированными многочленами степени $n$}.\footnote{Элементы пространств $W^{T\text{-}\mathrm{gr}}_n$,
названные нами полилинейными многочленами, конечно же, таковыми не являются, если их рассматривать как
многочлены от переменных $x_i^{(t)}$. Однако они являются полилинейными функциями от переменных
$x_i$, если рассматривать $x_i^{(t)}$ как $t$-компоненту переменной $x_i$ (см. ниже).}
Число $$c^{T\text{-}\mathrm{gr}}_n(A):=\dim\left(\frac{W^{T\text{-}\mathrm{gr}}_n}{W^{T\text{-}\mathrm{gr}}_n \cap \Id^{T\text{-}\mathrm{gr}}(A)}\right)$$
называется $n$-й \textit{коразмерностью градуированных полиномиальных тождеств}
или $n$-й \textit{градуированной коразмерностью} алгебры $A$.\label{DefCodimGr}

Предел $\PIexp^{T\text{-}\mathrm{gr}}(A):=\lim\limits_{n\rightarrow\infty} \sqrt[n]{c^{T\text{-}\mathrm{gr}}_n(A)}$ (если он существует) называется \textit{экспонентой роста $T$-градуированных тождеств} или \textit{$T$-градуированной PI-экспонентой} алгебры $A$.
\label{DefTgrPIexp}

Симметрическая группа $S_n$ действует на пространстве $\frac {W^{T\text{-}\mathrm{gr}}_n}{W^{T\text{-}\mathrm{gr}}_{n}
  \cap \Id^{T\text{-}\mathrm{gr}}(A)}$
  перестановками переменных внутри каждого множества $X^{(t)}$: $$\sigma x^{(t_1)}_{i_1}\ldots x^{(t_n)}_{i_n}
  := x^{(t_1)}_{\sigma(i_1)}\ldots x^{(t_n)}_{\sigma(i_n)}$$
  при $n\in\mathbb N$, $\sigma \in S_n$, $1\leqslant i_k\leqslant n$, $1\leqslant k \leqslant n$.
   Характер $\chi^{T\text{-}\mathrm{gr}}_n(A)$ представления группы $S_n$
   на пространстве $\frac {W^{T\text{-}\mathrm{gr}}_n}{W^{T\text{-}\mathrm{gr}}_n
   \cap \Id^{T\text{-}\mathrm{gr}}(A)}$\label{DefIdTgrA} называется
   $n$-м
  \textit{кохарактером} градуированных полиномиальных тождеств алгебры $A$.\label{DefCocharGr}
  Если $\chr \mathbbm{k} = 0$, $n$-й кохарактер представляется
  в виде суммы $$\chi^{T\text{-}\mathrm{gr}}_n(A)=\sum_{\lambda \vdash n}
   m(A, T\text{-}\mathrm{gr}, \lambda)\chi(\lambda)$$ \label{DefMulGr}
  неприводимых характеров $\chi(\lambda)$.
  Число $\ell_n^{T\text{-}\mathrm{gr}}(A):=\sum_{\lambda \vdash n}
   m(A, T\text{-}\mathrm{gr}, \lambda)$ называется $n$-й\label{DefColengthGr}
  \textit{кодлиной} градуированных полиномиальных тождеств алгебры $A$.
  
При изучении тождеств в ассоциативных алгебрах и алгебрах Ли вместо свободных неассоциативных алгебр можно использовать, соответственно, свободные ассоциативные алгебры
 и свободные алгебры Ли.
 
 В случае $T$-градуированных ассоциативных алгебр
 градуированные многочлены являются элементами свободной ассоциативной алгебры $\mathbbm{k}\langle X^{T\text{-}\mathrm{gr}}\rangle$ на множестве $X^{T\text{-}\mathrm{gr}}$, причём в случае,
 когда $T$~"--- полугруппа, алгебра $\mathbbm{k}\langle X^{T\text{-}\mathrm{gr}}\rangle$ сама оказывается
 $T$~"-градуированной. Вводятся пространства $P^{T\text{-}\mathrm{gr}}_n \subset \mathbbm{k}\langle X^{T\text{-}\mathrm{gr}}\rangle$ полилинейных ассоциативных $T$-градуированных многочленов
 степени $n$. Множество  $\Id^{T\text{-}\mathrm{gr}}_\mathrm{assoc}(A)$ ассоциативных градуированных полиномиальных тождеств $T$-градуированной ассоциативной алгебры $A$ является идеалом алгебры $\mathbbm{k}\langle X^{T\text{-}\mathrm{gr}}\rangle$.
 
 В случае градуировок алгебр Ли можно, с одной стороны, использовать
 свободную алгебру Ли $L(X^{T\text{-}\mathrm{gr}})$ на множестве $X^{T\text{-}\mathrm{gr}}$.
 C другой стороны, в случае, когда $T$~"--- полугруппа, можно, кроме этого,
использовать алгебру $L(X)^{T\text{-}\mathrm{gr}}$, которая по определению является факторалгеброй
алгебры $\mathbbm{k}\lbrace X^{T\text{-}\mathrm{gr}} \rbrace$ по градуированному идеалу, порождённому 
тождеством Якоби и тождеством антикоммутативности. Если $\chr \mathbbm{k} \ne 2$, то из тождества антикоммутативности
следует, что $L(X)^{T\text{-}\mathrm{gr}} \cong L(X^{T\text{-}\mathrm{gr}})$, если и только если (полу)группа $T$
коммутативна. Обозначим через $V^{T\text{-}\mathrm{gr}}_n$ подпространство полилинейных лиевских $T$-градуированных многочленов степени $n$ одной из этих двух алгебр. Пусть $\Id^{T\text{-}\mathrm{gr}}_\mathrm{Lie}(L)$~"--- идеал лиевских градуированных полиномиальных тождеств $T$-градуированной алгебры Ли $L$.

Точно так же, как и в случае $H$-тождеств (см. предыдущие два параграфа), коразмерности и кохарактеры не зависят от того, какую свободную алгебру мы используем. Отождествление происходит следующим образом: всякий $T$-градуированный многочлен $f$ (вне зависимости от свободной алгебры, элементом которой он является)
можно рассматривать как функцию на $T$-градуированной алгебре $A$ со значениями в $A$.
Пусть задан счётный набор $(a_n)_{n=1}^\infty$ элементов алгебры $A$.
Обозначим через $a^{(t)}_n$ однородные компоненты этих элементов: $a_n = \sum_{t\in T} a^{(t)}_n$,
где $a_n^{(t)}\in A^{(t)}$ и для всякого $n$ лишь конечное число элементов $a_n^{(t)}$ ненулевые.
Тогда значение многочлена $f$ на наборе $(a_n)_{n=1}^\infty$ определяется
как результат подстановки $x_n^{(t)}\mapsto a_n^{(t)}$.
При этом оказывается, что функции с тождественно нулевыми значениями соответствуют градуированным полиномиальным тождествам.
Для ассоциативной $T$-градуированной алгебры $A$ соответствие $x_i^{(t)} \mapsto x_i^{(t)}$, $i\in \mathbb N$, $t\in T$,
индуцирует изоморфизм
$\mathbbm{k} \lbrace X^{T\text{-}\mathrm{gr}}\rbrace/\Id^{T\text{-}\mathrm{gr}}(A) \cong \mathbbm{k}\langle X^{T\text{-}\mathrm{gr}} \rangle/ \Id^{T\text{-}\mathrm{gr}}_{\mathrm{assoc}}(A)$ алгебр
и изоморфизм $\mathbbm{k}S_n$-модулей $\frac {W^{T\text{-}\mathrm{gr}}_n}{W^{T\text{-}\mathrm{gr}}_{n}
  \cap \Id^{T\text{-}\mathrm{gr}}(A)} \cong \frac {P^{T\text{-}\mathrm{gr}}_n}{P^{T\text{-}\mathrm{gr}}_{n}
  \cap \Id_{\mathrm{assoc}}^{T\text{-}\mathrm{gr}}(A)}$.
    Для $T$-градуированной алгебры Ли $L$ это же соответствие
индуцирует изоморфизм
$\mathbbm{k} \lbrace X^{T\text{-}\mathrm{gr}}\rbrace/ \Id^{T\text{-}\mathrm{gr}}(L) \cong L(X^{T\text{-}\mathrm{gr}})/\Id^{T\text{-}\mathrm{gr}}_{\mathrm{Lie}}(L)$ алгебр Ли
и изоморфизм $\mathbbm{k}S_n$-модулей $\frac {W^{T\text{-}\mathrm{gr}}_n}{W^{T\text{-}\mathrm{gr}}_{n}
  \cap \Id^{T\text{-}\mathrm{gr}}(L)} \cong \frac {V^{T\text{-}\mathrm{gr}}_n}{V^{T\text{-}\mathrm{gr}}_{n}
  \cap \Id_{\mathrm{Lie}}^{T\text{-}\mathrm{gr}}(L)}$,
  причём аналогичные изоморфизмы имеют место и в случае, когда вместо алгебры Ли $L(X^{T\text{-}\mathrm{gr}})$ рассматривается алгебра Ли
  $L(X)^{T\text{-}\mathrm{gr}}$.
  
  В силу выше сделанных замечаний при изучении градуированных <<нетождеств>>, коразмерностей и кохарактеров оказывается несущественным,
  какие свободные алгебры использовать. В связи с этим далее идеалы $\Id^{T\text{-}\mathrm{gr}}_{\mathrm{assoc}}(A)$ и $\Id_{\mathrm{Lie}}^{T\text{-}\mathrm{gr}}(L)$ также будут обозначаться, соответственно, через $\Id^{T\text{-}\mathrm{gr}}(A)$ и $\Id^{T\text{-}\mathrm{gr}}(L)$.

Предложение ниже устанавливает связь между обычными и градуированными коразмерностями:
\begin{proposition}\label{PropositionOrdinaryAndGradedCodim}
Пусть $A$~"--- алгебра над полем $\mathbbm{k}$, градуированная
некоторым множеством $T$ (необязательно конечным).
Тогда $c_n(A) \leqslant c^{T\text{-}\mathrm{gr}}_n(A)$.
Кроме того, если множество $T$ конечно, то $c^{T\text{-}\mathrm{gr}}_n(A)
\leqslant |T|^n c_n(A)$ для всех $n\in\mathbb N$.
\end{proposition}
\begin{proof}
Пусть $t_1, \ldots, t_n \in T$.
Обозначим через $W_{t_1, \ldots, t_n}$
пространство полилинейных
неассоциативных многочленов от переменных $x_1^{(t_1)}, \ldots, x_n^{(t_n)}$
с коэффициентами из поля $\mathbbm{k}$.
Тогда $W_n^{T\text{-}\mathrm{gr}}=\bigoplus_{t_1, \ldots, t_n \in T}
W_{t_1, \ldots, t_n}$. Более того, \begin{equation}\label{EqDecompWnGradt1t2tn}
\frac{W^{T\text{-}\mathrm{gr}}_n}{W^{T\text{-}\mathrm{gr}}_n \cap \Id^{T\text{-}\mathrm{gr}}(A)} \cong 
\bigoplus_{t_1, \ldots, t_n \in T}
\frac{W_{t_1, \ldots, t_n}}{W_{t_1, \ldots, t_n} \cap \Id^{T\text{-}\mathrm{gr}}(A)},
\end{equation} поскольку компонента любого градуированного тождества,
отвечающая подпространству $W_{t_1, \ldots, t_n}$ выделяется подстановкой
 $x^{(t)}_i = 0$ для всех $1\leqslant i \leqslant n$ и
 $t \ne t_i$. Отметим, что изоморфизм~\eqref{EqDecompWnGradt1t2tn}
существует для любой характеристики поля $\mathbbm{k}$.

 Пусть $f_1,\ldots, f_{c_n(A)} \in W_n$ такие многочлены, что $\bar f_1, \ldots, \bar f_{c_n(A)}$~"--- базис
пространства $\frac{W_n}{W_n \cap \Id(A)}$.
Тогда для любого одночлена $w=x_{\sigma(1)}\ldots x_{\sigma(n)}$ (с любой расстановкой скобок), $\sigma\in S_n$, существуют $\alpha_{w,i}\in \mathbbm{k}$, такие, что $$x_{\sigma(1)}\ldots x_{\sigma(n)} - \sum_{i=1}^{c_n(A)}\alpha_{w,i} f_i(x_1,\ldots, x_n) \in \Id(A).$$
Заменяя для каждого набора элементов $t_1, \ldots, t_n \in T$ переменные $x_i$ на $x^{\left(t_i\right)}_i$,
получим $$x^{\left(t_{\sigma(1)}\right)}_{\sigma(1)}\ldots x^{\left(t_{\sigma(n)}\right)}_{\sigma(n)} - \sum_{i=1}^{c_n(A)}\alpha_{w,i} f_i\left(x^{(t_1)}_1,\ldots, x^{(t_n)}_n\right) \in \Id^{T\text{-}\mathrm{gr}}(A)$$
и $$\frac{W^{T\text{-}\mathrm{gr}}_n}{W^{T\text{-}\mathrm{gr}}_n \cap \Id^{T\text{-}\mathrm{gr}}(A)}
=\left\langle \bar f_i\left(x^{(t_1)}_1,\ldots, x^{(t_n)}_n\right)\mathrel{\Bigl|}
1\leqslant i \leqslant c_n(A),\  t_1, \ldots, t_n \in T\right\rangle_\mathbbm{k}.$$
Отсюда получаем оценку сверху для $c_n^{T\text{-}\mathrm{gr}}(A)$ в случае, когда множество $T$ конечно.

Для того, чтобы получить нижнюю оценку для $c_n^{T\text{-}\mathrm{gr}}(A)$, для заданного набора $(t_1, \ldots, t_n) \in T^n$ рассмотрим отображение $\varphi_{t_1, \ldots, t_n} \colon W_n \rightarrow \frac{W^{T\text{-}\mathrm{gr}}_n}{W^{T\text{-}\mathrm{gr}}_n \cap \Id^{T\text{-}\mathrm{gr}}(A)}$,
где $\varphi_{t_1, \ldots, t_n}(f)$~"--- образ многочлена $f\left(x^{(t_1)}_1, \ldots, x^{(t_n)}_n\right)$
в $ \frac{W^{T\text{-}\mathrm{gr}}_n}{W^{T\text{-}\mathrm{gr}}_n \cap \Id^{T\text{-}\mathrm{gr}}(A)}$ при $f=f(x_1, \ldots, x_n) \in W_n$.
Из полилинейности многочлена $f$ следует, что $f(x_1, \ldots, x_n) \equiv 0$ 
является обычным тождеством,
если и только если $$f\left(x^{(t_1)}_1, \ldots, x^{(t_n)}_n\right) \equiv 0 $$
является градуированным полиномиальным тождеством для всех $t_1, \ldots, t_n \in T$.
Другими словами, $W_n \cap \Id(A)
= \bigcap\limits_{(t_1, \ldots, t_n) \in T^n} \ker \varphi_{t_1, \ldots, t_n}$.
Поскольку пространство $W_n$ конечномерно, существует конечное подмножество $\Lambda \subseteq T^n$,
такое, что $W_n \cap \Id(A)
= \bigcap\limits_{(t_1, \ldots, t_n) \in \Lambda} \ker \varphi_{t_1, \ldots, t_n}$.

Рассмотрим вложение $$W_n \hookrightarrow
 W_n^{T\text{-}\mathrm{gr}}=\bigoplus_{t_1, \ldots, t_n \in T}
W_{t_1, \ldots, t_n},$$ где образ многочлена $f(x_1, \ldots, x_n)\in W_n$ равен $\sum_{(t_1, \ldots, t_n) \in \Lambda} f\left(x^{(t_1)}_1, \ldots, x^{(t_n)}_n\right)$. 
Тогда из~\eqref{EqDecompWnGradt1t2tn} и нашего выбора множества $\Lambda$
следует, что индуцированное отображение $\frac{W_n}{W_n \cap \Id(A)} \hookrightarrow \frac{W^{T\text{-}\mathrm{gr}}_n}{W^{T\text{-}\mathrm{gr}}_n \cap \Id^{T\text{-}\mathrm{gr}}(A)}$
является вложением, откуда и следует оценка снизу.
\end{proof}

Выше в примере~\ref{ExampleGr}
мы показали, что любая $T$-градуированная алгебра $A$ с конечным носителем градуировки
является алгеброй с обобщённым $\mathbbm{k}^T$-действием.
Покажем, что  изучение градуированных коразмерностей
и кохарактеров такой алгебры~$A$
можно свести к изучению коразмерностей и кохарактеров
её полиномиальных $\mathbbm{k}^T$-тождеств.

\begin{proposition}\label{PropositionCnGrCnGenH}
Пусть $\Gamma \colon A=\bigoplus_{t\in T} A^{(t)}$~"--- градуировка
алгебры $A$ над полем $\mathbbm{k}$ множеством $T$ с конечным носителем $\supp \Gamma$. Тогда $c_n^{T\text{-}\mathrm{gr}}(A)=c_n^{\mathbbm{k}^T}(A)$ и
$\chi_n^{T\text{-}\mathrm{gr}}(A)=\chi_n^{\mathbbm{k}^T}(A)$ для всех $n\in \mathbb N$.
 Если, кроме того, $\chr \mathbbm{k} = 0$, то $\ell_n^{T\text{-}\mathrm{gr}}(A)=\ell_n^{\mathbbm{k}^T}(A)$.
\end{proposition}
\begin{proof} 
Пусть $\xi \colon \mathbbm{k}\lbrace X | \mathbbm{k}^T \rbrace \to \mathbbm{k}\lbrace X^{T\text{-}\mathrm{gr}} \rbrace$~"--- гомоморфизм алгебр, заданный равенствами $\xi(x_i^q) = \sum\limits_{t\in\supp \Gamma} q(t)x^{(t)}_i$ для всех $i\in\mathbb N$ и $q\in \mathbbm{k}^T$. Пусть $f\in \Id^{\mathbbm{k}^T}(A)$. Рассмотрим произвольный гомоморфизм $\psi \colon  
\mathbbm{k}\lbrace X^{T\text{-}\mathrm{gr}} \rbrace \to A$, такой, что
$\psi(x^{(t)}_i)\in A^{(t)}$
для всех $t\in T$ и $i\in\mathbb N$. Тогда гомоморфизм алгебр $\psi\xi \colon \mathbbm{k}\lbrace X | \mathbbm{k}^T \rbrace
\to A$ удовлетворяет условию $$\psi\xi(x_i^q)=\sum\limits_{t\in\supp \Gamma} q(t)\psi\left(x^{(t)}_i\right)=
q\left(\sum\limits_{t\in\supp \Gamma} \psi\left(x^{(t)}_i\right)\right)=q\,\psi\xi(x_i).$$ Следовательно,
для всякого  такого гомоморфизма $\psi$ справедливо равенство $\psi\xi(f) =0$. Отсюда
 $\xi(f)\in \Id^{T\text{-}\mathrm{gr}}(A)$ и $\xi\left(\Id^{\mathbbm{k}^T}(A)\right)\subseteq \Id^{T\text{-}\mathrm{gr}}(A)$.
Обозначим через $$\tilde \xi \colon \mathbbm{k}\lbrace X | \mathbbm{k}^T \rbrace/\Id^{\mathbbm{k}^T}(A) \to \mathbbm{k}\lbrace X^{T\text{-}\mathrm{gr}} \rbrace/\Id^{T\text{-}\mathrm{gr}}(A)$$ гомоморфизм, индуцированный $\xi$.

Пусть $\eta \colon \mathbbm{k}\lbrace X^{T\text{-}\mathrm{gr}} \rbrace \to \mathbbm{k}\lbrace X | \mathbbm{k}^T \rbrace$~"--- гомоморфизм
алгебр, заданный равенствами $\eta\left(x^{(t)}_i\right) = x^{q_t}_i$ при $i\in \mathbb N$
и $t\in T$, где $q_t(s):=\left\lbrace\begin{smallmatrix} 1 & \text{при} & s=t,\\ 0 & \text{при} & s\ne t.\end{smallmatrix} \right.$ Рассмотрим произвольное градуированное полиномиальное тождество $f\in \Id^{T\text{-}\mathrm{gr}}(A)$.
Пусть $\psi \colon  \mathbbm{k}\lbrace X | \mathbbm{k}^T \rbrace \to A$~"--- гомоморфизм,
удовлетворяющий условию
$\psi(x_i^q)=q\psi(x_i)$ для всех $i\in\mathbb N$ и $q\in \mathbbm{k}^T$.
Тогда для всех $i\in\mathbb N$ и $g, t \in T$ справедливы равенства
$$q_g \psi\eta\left(x^{(t)}_i\right) = q_g\psi(x^{q_t}_i)=q_g q_t \psi(x_i)
=\left\lbrace \begin{array}{lll} 0 & \text{ при } & g\ne t,\\
                              \psi\eta\left(x^{(t)}_i\right) & \text{ при } & g=t. \end{array}\right.$$
 Следовательно, $\psi\eta\left(x^{(t)}_i\right) \in A^{(t)}$.
 Отсюда $\psi\eta(f)=0$ и $\eta(\Id^{T\text{-}\mathrm{gr}}(A)) \subseteq \Id^{\mathbbm{k}^T}(A)$.
Обозначим через $\tilde\eta \colon  \mathbbm{k}\lbrace X^{T\text{-}\mathrm{gr}} \rbrace/\Id^{T\text{-}\mathrm{gr}}(A) \to
\mathbbm{k}\lbrace X | \mathbbm{k}^T \rbrace/\Id^{\mathbbm{k}^T}(A)$ индуцированный гомоморфизм.

Ниже используются обозначения $\bar f = f + \Id^{\mathbbm{k}^T}(A) \in \mathbbm{k}\lbrace X | \mathbbm{k}^T \rbrace/\Id^{\mathbbm{k}^T}(A)$ для $f\in
\mathbbm{k}\lbrace X | \mathbbm{k}^T \rbrace$ и  $\bar f = f + \Id^{T\text{-}\mathrm{gr}}(A) \in \mathbbm{k}\lbrace X^{T\text{-}\mathrm{gr}} \rbrace/\Id^{T\text{-}\mathrm{gr}}(A)$ для $f\in \mathbbm{k}\lbrace X^{T\text{-}\mathrm{gr}} \rbrace$.
Заметим, что $$x^q_i - \sum\limits_{t\in\supp \Gamma} q(t) x^{q_t}_i\in \Id^{\mathbbm{k}^T}(A)$$ для всех $q\in \mathbbm{k}^T$ и $i\in\mathbb N$.
Следовательно, $$\tilde\eta\tilde\xi\left(\bar x^q_i\right)=\tilde\eta\left(
\sum\limits_{t\in\supp \Gamma} q(t) \bar x^{(t)}_i\right)
=\sum\limits_{t\in\supp \Gamma} q(t) \bar x^{q_t}_i = \bar x^q_i$$
для всех $q\in \mathbbm{k}^T$ и $i\in\mathbb N$. 
Отсюда получаем, что $\tilde\eta\tilde\xi=\id_{\mathbbm{k}\lbrace X | \mathbbm{k}^T \rbrace/\Id^{\mathbbm{k}^T}(A)}$,
поскольку алгебра $\mathbbm{k}\lbrace X | \mathbbm{k}^T \rbrace/\Id^{\mathbbm{k}^T}(A)$ порождена элементами $\bar x^q_i$, где $q\in \mathbbm{k}^T$ и $i\in\mathbb N$.
Более того, $\tilde\xi\tilde\eta\left(\bar x^{(t)}_i\right)=
\tilde\xi\left(\bar x^{q_t}_i\right)=\bar x^{(t)}_i$ для всех $t\in \supp \Gamma$ и $i\in \mathbb N$.
Следовательно, $\tilde\xi\tilde\eta=\id_{\mathbbm{k}\lbrace X^{T\text{-}\mathrm{gr}} \rbrace/\Id^{T\text{-}\mathrm{gr}}(A)}$
и $\mathbbm{k}\lbrace X^{T\text{-}\mathrm{gr}} \rbrace/\Id^{T\text{-}\mathrm{gr}}(A) \cong \mathbbm{k}\lbrace X | \mathbbm{k}^T \rbrace/\Id^{\mathbbm{k}^T}(A)$
как алгебры. Ограничение изоморфизма $\tilde\xi$ является изоморфизмом $\mathbbm{k}S_n$-модулей
$\frac{W^{\mathbbm{k}^T}_n}{W^{\mathbbm{k}^T}_n \cap \Id^{\mathbbm{k}^T}(A)}$ и $\frac{W^{T\text{-}\mathrm{gr}}_n}{W^{T\text{-}\mathrm{gr}}_n\cap \Id^{T\text{-}\mathrm{gr}}(A)}$.  Следовательно, $$c^{\mathbbm{k}^T}_n(A)=\dim \frac{W^{\mathbbm{k}^T}_n}{W^{\mathbbm{k}^T}_n \cap \Id^{\mathbbm{k}^T}(A)}
= \dim\frac{W^{T\text{-}\mathrm{gr}}_n}{W^{T\text{-}\mathrm{gr}}_n\cap \Id^{T\text{-}\mathrm{gr}}(A)}=c^{T\text{-}\mathrm{gr}}_n(A)$$
и  $\chi_n^{T\text{-}\mathrm{gr}}(A)=\chi_n^{\mathbbm{k}^T}(A)$
для всех $n\in \mathbb N$.
 Если, кроме того, $\chr \mathbbm{k} = 0$, тогда выполнены и равенства $\ell_n^{T\text{-}\mathrm{gr}}(A)=\ell_n^{\mathbbm{k}^T}(A)$.
\end{proof}
\begin{corollary}
Пусть $A$~"--- конечномерная алгебра над произвольным полем,
градуированная произвольным множеством $T$.
Тогда $c^{T\text{-}\mathrm{gr}}_n(A) \leqslant (\dim A)^{n+1}$ для всех  $n \in \mathbb N$.
\end{corollary}
\begin{proof}
В силу конечномерности алгебры $A$ носитель градуировки конечен,
 поэтому достаточно применить предложения~\ref{PropositionCodimDim} и~\ref{PropositionCnGrCnGenH}. \end{proof}
 
 При работе с градуированными тождествами в алгебрах Ли нам будет удобно использовать
 также действие группы линейных характеров градуирующей группы.
 
 Результат леммы~\ref{LemmaAbelianDual} не является новым, однако для удобства читателя мы приведём
его доказательство:

\begin{lemma}\label{LemmaAbelianDual}
Пусть $\mathbbm{k}$~"--- алгебраически замкнутое поле характеристики $0$, $G$~"--- конечнопорождённая
абелева группа, а $\hat G := \Hom(G, \mathbbm{k}^\times)$~"--- группа гомоморфизмов
из $G$ в мультипликативную группу $\mathbbm{k}^\times$ поля $\mathbbm{k}$. Рассмотрим элементы групповой алгебры $\mathbbm{k}\hat G$
как функции на $G$. Тогда для всех попарно различных элементов $\gamma_1, \ldots, \gamma_m \in G$ 
существуют такие элементы $h_1, \ldots, h_m \in \mathbbm{k}\hat G$,
что $h_i(\gamma_j)=\left\lbrace
  \begin{array}{lll}
  1 & \text{при} & i = j, \\
  0 & \text{при} & i \ne j.
  \end{array} \right.$
\end{lemma}
\begin{proof} В случае, когда группа $G$ конечна, из соотношений ортогональности для характеров
следует, что пространство всех функций на $G$ со значениями из $\mathbbm{k}$ является линейной оболочкой
группы $\hat G$. Отсюда такие функции $h_i$ всегда можно найти.

Теперь рассмотрим случай, когда $G = \langle g\rangle$~"--- бесконечная циклическая группа,
достаточно определить элемент $\chi \in \hat G$ через $\chi(g^k)=\lambda^k$, где $\lambda \in \mathbbm{k}^\times$~"--- фиксированный элемент бесконечного порядка.
Тогда из невырожденности матрицы Вандермонда следует, что элементы
 $1$, $\chi$, $\chi^2$, \ldots, $\chi^{m-1}$ линейно независимы как
 функции от $\gamma_1, \ldots, \gamma_m$ и требуемые $h_1, \ldots, h_m$ можно найти и в этом случае. 

В общем случае в силу теоремы о строении конечнопорождённых абелевых групп $G = G_1 \times \mathbb Z^s$, где $G_1$~"--- конечная группа.
Следовательно, $\hat G = \hat G_1 \times (\mathbbm{k}^\times)^s$, где $\hat G_1$~"--- группа линейных характеров
группы $G_1$. Выберем теперь свои элементы для каждой компоненты прямого произведения, рассмотрим их произведения
и получим соответствующие элементы  $h_1, \ldots, h_m \in \mathbbm{k}\hat G$.
\end{proof}

Для любой группы $G$ на любом $G$-градуированном
пространстве $V=\bigoplus_{g\in G} V^{(g)}$
можно задать структуру $\mathbbm{k}\hat G$-модуля:
$$\chi v = \chi(g)v\text{ для всех }\chi \in \hat G,\ g\in G\text{ и }v \in V^{(g)}.$$
При этом всякая $G$-градуированная алгебра $A$ оказывается алгеброй с действием группы $\hat G$
автоморфизмами, причём в силу леммы~\ref{LemmaAbelianDual} градуированные подпространства
являются в точности $\mathbbm{k}\hat G$-подмодулями.

%

\begin{remark}
Можно доказать, что для всякой конечнопорождённой абелевой группы $G$ над алгебраически замкнутым полем $\mathbbm{k}$
характеристики $0$ алгебра Хопфа $\mathbbm{k}G$ не содержит ненулевых нильпотентных элементов,
откуда будет следовать, что $\mathbbm{k}G=\mathcal O(\hat G)$, где $\hat G := \Hom(G, \mathbbm{k}^\times)=\mathbf{Alg}(\mathbbm{k}G, \mathbbm{k})$~"--- множество максимальных идеалов алгебры Хопфа $\mathbbm{k}G$, которое наделено естественной структурой аффинной алгебраической группы,
причём групповая операция на $\hat G$, единица и взятие обратного являются двойственными
операциями к, соответственно, коумножению, коединице и антиподу алгебры Хопфа $\mathbbm{k}G$.
(См.~\cite[\S 4.2.1]{Abe}.) Конечномерные $\mathbbm{k}G$-комодули, т.е. конечномерные $G$-градуированные пространства, будут являться в точности рациональными представлениями группы $\hat G$.
Поскольку группа $\hat G$ является прямым произведением конечной абелевой группы на тор, любое такое конечномерное представление раскладывается в прямую сумму одномерных инвариантных подпространств, т.е. является вполне приводимым, откуда группа $\hat G$ редуктивна.
\end{remark}

Оказывается, что $\hat G$-коразмерности совпадают с $G$-градуированными коразмерностями:

\begin{proposition}\label{PropositionGrToHatG}
Пусть $\Gamma\colon A=\bigoplus_{g\in G} A^{(g)}$~"--- необязательно ассоциативная алгебра 
над алгебраически замкнутым полем $\mathbbm{k}$ характеристики $0$,
градуированная конечнопорождённой 
абелевой группой $G$, причём носитель градуировки $\supp \Gamma$ конечен.
Рассмотрим на $A$  действие группы $\mathbbm{k}\hat G$ автоморфизмами, заданное выше. Тогда $c^{G\text{-}\mathrm{gr}}_n(A)=c^{\mathbbm{k}\hat G}(A)$,
$\chi_n^{G\text{-}\mathrm{gr}}(A)=\chi_n^{\mathbbm{k}\hat G}(A)$ и
 $\ell_n^{G\text{-}\mathrm{gr}}(A)=\ell_n^{\mathbbm{k}\hat G}(A)$  для всех $n\in \mathbb N$.
\end{proposition}
\begin{proof}
Пусть $\lbrace \gamma_1, \ldots, \gamma_m \rbrace := \supp \Gamma$.
Теперь достаточно выбрать элементы $h_1, \ldots, h_m \in \mathbbm{k}\hat G$
в соответствии с леммой~\ref{LemmaAbelianDual}
и дословно повторить доказательство предложения~\ref{PropositionCnGrCnGenH},
используя вместо $q_t$ элемент $h_i$ при $t=\gamma_i$ и элемент $0$ при $t\notin \supp \Gamma$.
\end{proof}

\section{Градуированные $H$-тождества}

Для алгебры, на которой задано обобщённое $H$-действие, согласованное с некоторой $T$-градуировкой, 
можно ввести понятие градуированного $H$-тождества

Пусть $H$~"--- ассоциативная алгебра с единицей над полем $\mathbbm{k}$, а $T$~"--- некоторое множество. Рассмотрим
свободную неассоциативную алгебру $\mathbbm{k}\lbrace X^{T\text{-}\mathrm{gr}}|H \rbrace$ на множестве $X^{T\text{-}\mathrm{gr}}$ с символами операторов из алгебры $H$ (см. \S\ref{SectionHPI}),
где $X^{T\text{-}\mathrm{gr}}:=\bigsqcup_{t \in T}X^{(t)}$~"--- объединение
непересекающихся множеств $X^{(t)} = \{ x^{(t)}_1,
x^{(t)}_2, \ldots \}$.

Будем говорить, что $f=f(x^{(t_1)}_{i_1}, \ldots, x^{(t_s)}_{i_s}) \in \mathbbm{k}\lbrace X^{T\text{-}\mathrm{gr}}|H \rbrace$~"--- \textit{градуированное полиномиальное $H$-тождество} для $T$-градуированной алгебры $A=\bigoplus_{t\in T}
A^{(t)}$ с согласованным обобщённым $H$-действием и писать $f\equiv 0$,
если $f(a^{(t_1)}_1, \ldots, a^{(t_s)}_s)=0$
для всех $a^{(t_j)}_j \in A^{(t_j)}$, $1 \leqslant j \leqslant s$.
  Множество $\Id^{T\text{-}\mathrm{gr}, H}(A)$ градуированных полиномиальных
  $H$-тождеств алгебры $A$ является
идеалом алгебры $\mathbbm{k}\lbrace
 X^{T\text{-}\mathrm{gr}}|H\rbrace$.

Пусть $$W_n^{T\text{-}\mathrm{gr}, H} := \left\langle\left. \left(x_{\sigma(1)}^{(t_1)}\right)^{h_1} \left(x_{\sigma(2)}^{(t_2)}\right)^{h_2}\ldots \left(x_{\sigma(n)}^{(t_n)}\right)^{h_n} \,\right|\, \sigma \in S_n,\ t_i\in T,\ h_i\in H\right\rangle_\mathbbm{k} \subset \mathbbm{k}\lbrace X|H \rbrace,$$
где $S_n$~"--- $n$-я группа подстановок,
$n\in\mathbb N$ и на одночленах рассматриваются всевозможные расстановки скобок. Элементы пространств $W_n^{T\text{-}\mathrm{gr}, H}$ называются \textit{полилинейными неассоциативными $T$-градуированными $H$-многочленами},
а элементы пространств $W_n^{T\text{-}\mathrm{gr}, H} \cap \Id^{T\text{-}\mathrm{gr}, H}(A)$ называются \textit{полилинейными градуированными $H$-тождествами} алгебры $A$.

Число $c_n^{T\text{-}\mathrm{gr}, H}(A):=\dim\left(\frac{W_n^{T\text{-}\mathrm{gr}, H}}{W_n^{T\text{-}\mathrm{gr}, H}\cap\, \Id^{T\text{-}\mathrm{gr}, H}(A)}\right)$, где $n\in\mathbb N$, называется \textit{$n$-й коразмерностью $T$-градуированных полиномиальных $H$-тождеств}
или \textit{$n$-й $T$-градуированной $H$-коразмерностью}
алгебры $A$.

Предел $\PIexp^{T\text{-}\mathrm{gr}, H}(A):=\lim\limits_{n\rightarrow\infty} \sqrt[n]{c^{T\text{-}\mathrm{gr}, H}_n(A)}$ (если он существует)
называется \textit{экспонентой роста $T$-градуированных  полиномиальных $H$-тождеств} или \textit{$T$-градуированной $H$-PI-экспонентой} алгебры $A$.

Симметрическая группа $S_n$ действует на пространстве $\frac {W^{T\text{-}\mathrm{gr}, H}_n}{W^{T\text{-}\mathrm{gr}, H}_{n}
  \cap \Id^{T\text{-}\mathrm{gr}, H}(A)}$
  перестановками переменных внутри каждого множества $X^{(t)}$: $$\sigma \left(x^{(t_1)}_{i_1}\right)^{h_1}\ldots \left(x^{(t_n)}_{i_n}\right)^{h_n}
  := \left(x^{(t_1)}_{\sigma(i_1)}\right)^{h_1}\ldots \left(x^{(t_n)}_{\sigma(i_n)}\right)^{h_n}$$
  при $n\in\mathbb N$, $\sigma \in S_n$, $1\leqslant i_k\leqslant n$, $1\leqslant k \leqslant n$.
   Характер $\chi^{T\text{-}\mathrm{gr}, H}_n(A)$ представления группы $S_n$
   на пространстве $\frac {W^{T\text{-}\mathrm{gr}, H}_n}{W^{T\text{-}\mathrm{gr}, H}_n
   \cap \Id^{T\text{-}\mathrm{gr}, H}(A)}$ называется
   $n$-м
  \textit{кохарактером} градуированных полиномиальных $H$-тождеств алгебры $A$.
  Если $\chr \mathbbm{k} = 0$, $n$-й кохарактер представляется
  в виде суммы $$\chi^{T\text{-}\mathrm{gr}, H}_n(A)=\sum_{\lambda \vdash n}
   m(A, T\text{-}\mathrm{gr}, H, \lambda)\chi(\lambda)$$ 
  неприводимых характеров $\chi(\lambda)$.
  Число $\ell_n^{T\text{-}\mathrm{gr}, H}(A):=\sum_{\lambda \vdash n}
   m(A, T\text{-}\mathrm{gr}, H, \lambda)$ называется $n$-й
  \textit{кодлиной} градуированных полиномиальных $H$-тождеств алгебры $A$.
  
  По аналогии с градуированными тождествами и $H$-тождествами можно ввести лиевские и ассоциативные
  свободные алгебры $L(
 X^{T\text{-}\mathrm{gr}}|H)$ и $\mathbbm{k}\langle
 X^{T\text{-}\mathrm{gr}}|H\rangle$, определить
 соответствующие
  градуированные полиномиальные $H$-тождества и показать,
  что коразмерности, кодлины и кохарактеры не зависят от того, из какой свободной
  алгебры берутся градуированные $H$-многочлены.

В теореме~\ref{TheoremGradGenActionReplace} было доказано, что любое обобщённое $H$-действие,
согласованное с некоторой $T$-градуировкой, сводится к обобщённому $\mathbbm{k}^T\otimes H$-действию.
По аналогии с предложением~\ref{PropositionCnGrCnGenH} покажем, что коразмерности
соответствующих тождеств также совпадают.

\begin{proposition}\label{PropositionCnTGrHCnFTotimesH}
Пусть $\Gamma \colon A=\bigoplus_{t\in T} A^{(t)}$~"--- $T$-градуировка на алгебре $A$ над полем $\mathbbm{k}$, где $T$~"--- некоторое множество, причём на $A$ также заданное такое обобщённое $H$-действие, 
что $H$~"--- ассоциативная алгебра с $1$, $HA^{(t)}\subseteq A^{(t)}$ для всех $t\in T$, а носитель $\supp \Gamma$ градуировки $\Gamma$ конечен.
 Тогда $$c_n^{T\text{-}\mathrm{gr}, H}(A)=c_n^{\mathbbm{k}^T\otimes H}(A)\text{\quad и\quad}\chi_n^{T\text{-}\mathrm{gr}, H}(A)=\chi_n^{\mathbbm{k}^T\otimes H}(A)\text{ для всех }n\in \mathbb N.$$
 Если, кроме того, $\chr \mathbbm{k} = 0$, то $\ell_n^{T\text{-}\mathrm{gr}, H}(A)=\ell_n^{\mathbbm{k}^T\otimes H}(A)$.
\end{proposition}
\begin{proof} 
Пусть $\xi \colon \mathbbm{k}\lbrace X | \mathbbm{k}^T \otimes H \rbrace \to \mathbbm{k}\lbrace X^{T\text{-}\mathrm{gr}} | H \rbrace$~"--- гомоморфизм алгебр, заданный равенствами $\xi(x_i^{q\otimes h}) = \sum\limits_{t\in\supp \Gamma} q(t)\left(x^{(t)}_i\right)^h$ для всех $i\in\mathbb N$, $q\in \mathbbm{k}^T$, $h\in H$. Пусть $f\in \Id^{\mathbbm{k}^T\otimes H}(A)$. Рассмотрим произвольный гомоморфизм $\psi \colon  
\mathbbm{k}\lbrace X^{T\text{-}\mathrm{gr}} | H \rbrace \to A$, такой, что
$\psi(x^{(t)}_i)\in A^{(t)}$ и $\psi\left( \left(x^{(t)}_i\right)^h\right)=h\psi(x^{(t)}_i)$
для всех $t\in T$, $i\in\mathbb N$ и $h\in H$. Тогда гомоморфизм алгебр $\psi\xi \colon \mathbbm{k}\lbrace X | \mathbbm{k}^T \otimes H \rbrace
\to A$ удовлетворяет условию $$\psi\xi(x_i^{q\otimes h})=\sum\limits_{t\in\supp \Gamma} q(t)\psi\left(\left(x^{(t)}_i\right)^h\right)=
(q\otimes h)\sum\limits_{t\in\supp \Gamma} \psi\left(x^{(t)}_i\right)=(q\otimes h)\,\psi\xi(x_i).$$ Следовательно,
для всякого  такого гомоморфизма $\psi$ справедливо равенство $\psi\xi(f) =0$. Отсюда
 $\xi(f)\in \Id^{T\text{-}\mathrm{gr}, H}(A)$ и $\xi\left(\Id^{\mathbbm{k}^T\otimes H}(A)\right)\subseteq \Id^{T\text{-}\mathrm{gr}, H}(A)$.
Обозначим через $$\tilde \xi \colon \mathbbm{k}\lbrace X | \mathbbm{k}^T \otimes H \rbrace/\Id^{\mathbbm{k}^T\otimes H}(A) \to \mathbbm{k}\lbrace X^{T\text{-}\mathrm{gr}} | H \rbrace/\Id^{T\text{-}\mathrm{gr}, H}(A)$$ гомоморфизм, индуцированный $\xi$.

Пусть $\eta \colon \mathbbm{k}\lbrace X^{T\text{-}\mathrm{gr}} | H \rbrace \to \mathbbm{k}\lbrace X | \mathbbm{k}^T \otimes H \rbrace$~"--- гомоморфизм
алгебр, заданный равенствами $\eta\left(\left(x^{(t)}_i\right)^h\right) = x_i^{q_t \otimes h}$ при $i\in \mathbb N$
и $t\in T$. Рассмотрим произвольное градуированное полиномиальное $H$-тождество $f\in \Id^{T\text{-}\mathrm{gr}, H}(A)$.
Пусть $\psi \colon  \mathbbm{k}\lbrace X | \mathbbm{k}^T \otimes H \rbrace \to A$~"--- гомоморфизм,
удовлетворяющий условию
$\psi(x_i^{q\otimes h})=(q\otimes h)\psi(x_i)$ для всех $i\in\mathbb N$, $q\in \mathbbm{k}^T$ и $h\in H$.
Тогда для всех $i\in\mathbb N$ и $g, t \in T$ справедливы равенства
$$q_g \psi\eta\left(x^{(t)}_i\right) = q_g\psi\left(x^{q_t\otimes 1}_i\right)=q_g q_t \psi(x_i)
=\left\lbrace \begin{array}{lll} 0 & \text{ при } & g\ne t,\\
                              \psi\eta\left(x^{(t)}_i\right) & \text{ при } & g=t. \end{array}\right.$$
 Следовательно, $\psi\eta\left(x^{(t)}_i\right) \in A^{(t)}$.
 Кроме того, $$ \psi\eta\left(\left(x^{(t)}_i\right)^h\right) = 
 \psi\left(x_i^{q_t\otimes h}\right)=h\,\psi\left(x_i^{q_t\otimes 1}\right)= h\, \psi\eta\left(x^{(t)}_i\right).$$
 Отсюда $\psi\eta(f)=0$ и $\eta(\Id^{T\text{-}\mathrm{gr}, H}(A)) \subseteq \Id^{\mathbbm{k}^T\otimes H}(A)$.
Обозначим через $$\tilde\eta \colon  \mathbbm{k}\lbrace X^{T\text{-}\mathrm{gr}} | H \rbrace/\Id^{T\text{-}\mathrm{gr}, H}(A) \to
\mathbbm{k}\lbrace X | \mathbbm{k}^T \otimes H \rbrace/\Id^{\mathbbm{k}^T\otimes H}(A)$$ индуцированный гомоморфизм.

Ниже используются обозначения $\bar f = f + \Id^{\mathbbm{k}^T\otimes H}(A) \in \mathbbm{k}\lbrace X | \mathbbm{k}^T \otimes H \rbrace/\Id^{\mathbbm{k}^T\otimes H}(A)$ для $f\in
\mathbbm{k}\lbrace X | \mathbbm{k}^T \otimes H \rbrace$ и  $\bar f = f + \Id^{T\text{-}\mathrm{gr}, H}(A) \in \mathbbm{k}\lbrace X^{T\text{-}\mathrm{gr}} | H \rbrace/\Id^{T\text{-}\mathrm{gr}, H}(A)$ для $f\in \mathbbm{k}\lbrace X^{T\text{-}\mathrm{gr}} | H \rbrace$.
Заметим, что $$x^{q\otimes h}_i - \sum\limits_{t\in\supp \Gamma} q(t) x^{q_t\otimes h}_i\in \Id^{\mathbbm{k}^T\otimes H}(A)$$ для всех $q\in \mathbbm{k}^T$, $h\in H$ и $i\in\mathbb N$.
Следовательно, $$\tilde\eta\tilde\xi\left(\overline{x^{q\otimes h}_i}\right)=\tilde\eta\left(
\sum\limits_{t\in\supp \Gamma} q(t) \overline{\left(x^{(t)}_i\right)^h}\right)
=\sum\limits_{t\in\supp \Gamma} q(t) \overline{x^{q_t\otimes h}_i} = \overline{x^{q\otimes h}_i}$$
для всех $q\in \mathbbm{k}^T$, $h\in H$ и $i\in\mathbb N$. 
Отсюда получаем, что $$\tilde\eta\tilde\xi=\id_{\mathbbm{k}\lbrace X | \mathbbm{k}^T \otimes H \rbrace/\Id^{\mathbbm{k}^T\otimes H}(A)},$$
поскольку алгебра $\mathbbm{k}\lbrace X | \mathbbm{k}^T \otimes H \rbrace/\Id^{\mathbbm{k}^T\otimes H}(A)$ порождена элементами $\overline{x^{q\otimes h}_i}$, где $q\in \mathbbm{k}^T$, $h\in H$ и $i\in\mathbb N$.
Более того, $\tilde\xi\tilde\eta\Biggl(\overline{\left(x^{(t)}_i\right)^h}\Biggr)=
\tilde\xi\left(\overline{x^{q_t\otimes h}_i}\right)=\overline{\left(x^{(t)}_i\right)^h}$ для всех $t\in \supp \Gamma$, $h\in H$ и $i\in \mathbb N$.
Следовательно, $\tilde\xi\tilde\eta=\id_{\mathbbm{k}\lbrace X^{T\text{-}\mathrm{gr}} | H \rbrace/\Id^{T\text{-}\mathrm{gr}, H}(A)}$
и $$\mathbbm{k}\lbrace X^{T\text{-}\mathrm{gr}} | H \rbrace/\Id^{T\text{-}\mathrm{gr}, H}(A) \cong \mathbbm{k}\lbrace X | \mathbbm{k}^T \otimes H \rbrace/\Id^{\mathbbm{k}^T\otimes H}(A)$$
как алгебры. Ограничение изоморфизма $\tilde\xi$ является изоморфизмом $\mathbbm{k}S_n$-модулей
$\frac{W^{\mathbbm{k}^T\otimes H}_n}{W^{\mathbbm{k}^T\otimes H}_n \cap \Id^{\mathbbm{k}^T\otimes H}(A)}$ и $\frac{W^{T\text{-}\mathrm{gr}, H}_n}{W^{T\text{-}\mathrm{gr}, H}_n\cap \Id^{T\text{-}\mathrm{gr}, H}(A)}$.  Следовательно, $$c^{\mathbbm{k}^T\otimes H}_n(A)=\dim \frac{W^{\mathbbm{k}^T\otimes H}_n}{W^{\mathbbm{k}^T\otimes H}_n \cap \Id^{\mathbbm{k}^T\otimes H}(A)}
= \dim\frac{W^{T\text{-}\mathrm{gr}, H}_n}{W^{T\text{-}\mathrm{gr}, H}_n\cap \Id^{T\text{-}\mathrm{gr}, H}(A)}=c^{T\text{-}\mathrm{gr}, H}_n(A)$$
и  $\chi_n^{T\text{-}\mathrm{gr}, H}(A)=\chi_n^{\mathbbm{k}^T\otimes H}(A)$
для всех $n\in \mathbb N$.
 Если, кроме того, $\chr \mathbbm{k} = 0$, тогда выполнены и равенства $\ell_n^{T\text{-}\mathrm{gr}, H}(A)=\ell_n^{\mathbbm{k}^T\otimes H}(A)$.
\end{proof}
 
 \section{Свободно-забывающие сопряжения, соответствующие градуировкам и обобщённым $H$-действиям}\label{SectionGrHAdjunction}

В~\S\ref{SectionHPI} была построена алгебра $\mathbbm{k}\lbrace X|H \rbrace$, причём
в~\S\ref{SectionHIdHMod} было отмечено, что если $H$~"--- алгебра Хопфа, то $\mathbbm{k}\lbrace X|H \rbrace$
является свободной $H$-модульной алгеброй, т.е. функцией объектов левого сопряжённого функтора к забывающему функтору из категории $\mathbf{NAAlg}({}_H\mathcal M)$ в категорию $\mathbf{Sets}$. При этом в случае, 
когда $H$ не является алгеброй Хопфа, алгебра $\mathbbm{k}\lbrace X|H \rbrace$
сама уже не является алгеброй с обобщённым $H$-действием, поскольку никакое $H$-действие
на $\mathbbm{k}\lbrace X|H \rbrace$ не определено.

В~\S\ref{SectionGrAdjoint} было показано, что для всякой группы $G$ свободная ассоциативная $G$-градуированная
алгебра (там она обозначалась через $\mathbbm{k}\langle X\backslash \lbrace 0 \rbrace  \rangle$) является функцией объектов для левого сопряжённого функтора к забывающему функтору из категории
$G$-градуированных алгебр в категорию  $G$-градуированных множеств. Аналог этого сопряжения очевидным образом строится для всякой полугруппы $T$ и между категориями необязательно ассоциативных $T$-градуированных алгебр и категорией $T$-градуированных множеств.

В данном параграфе мы покажем, что если $H$~"--- произвольная ассоциативная алгебра с единицей, а $T$~"--- множество, то можно так расширить категории
алгебр с обобщённым $H$-действием и $T$-градуированных алгебр, что алгебры $\mathbbm{k}\lbrace X|H \rbrace$ и $\mathbbm{k}\lbrace X^{T\text{-}\mathrm{gr}} \rbrace$ также окажутся значениями левых сопряжённых функторов к забывающим функторам из соответствующих категорий.

Напомним, что значения на объектах левых сопряженных функторов к забывающим  функторам из категорий алгебр в категорию множеств называются \textit{свободными} алгебрами, а значения на объектах левых сопряженных функторов к забывающим  функторам из категорий алгебр в категорию векторных пространств~"--- \textit{тензорными} алгебрами. Между свободными и тензорными алгебрами существуют естественные изоморфизмы. Разница между ними только в том, какие объекты используются для построения таких алгебр, множества или векторные пространства. В двух примерах, приведённых выше, рассматриваемые алгебры оказывались свободными. В данном параграфе нам будет удобно использовать язык тензорных алгебр.

Для определённости в параграфе рассматриваются категории необязательно ассоциативных алгебр, хотя, конечно, аналоги конструкций, приводимых ниже, существуют и для соответствующих категорий ассоциативных алгебр и алгебр Ли.

%
%
%

\subsection{Градуировки}\label{SectionGrAdjunction}

Пусть $T$~"--- множество, а $\mathbbm{k}$~"--- поле.
Обозначим через $\mathbf{Vect}^{T\text{-}\mathrm{gr}}_\mathbbm{k}$
категорию, объектами которой являются всевозможные \textit{$T$-градуированные векторные пространства}
над полем $\mathbbm{k}$, т.е. 
векторные пространства
$V$ с фиксированным разложением $V= \bigoplus_{t\in T} V^{(t)}$,
а множество $\mathbf{Vect}^{T\text{-}\mathrm{gr}}_\mathbbm{k}(V,W)$ морфизмов
между $V=\bigoplus_{t\in T} V^{(t)}$ и $W=\bigoplus_{t\in T} W^{(t)}$
состоит из всех линейных отображений $\varphi \colon V \to W$,
таких, что $\varphi\left(V^{(t)}\right)\subseteq W^{(t)}$
для всех $t\in T$.

Обозначим через $\mathbf{NAAlg}^{T\text{-}\mathrm{pgr}}_\mathbbm{k}$
(от англ. ``not necessarily associative partially $T$-graded algebras'')
категорию, объектами которой являются всевозможные необязательно ассоциативные
алгебры $A$ над полем $\mathbbm{k}$ с фиксированными $T$-градуированными подпространствами
$\bigoplus_{t\in T} A^{(t)} \subseteq A$
(вложение может быть строгим),
причём если $A \supseteq \bigoplus_{t\in T} A^{(t)}$ и $B\supseteq \bigoplus_{t\in T} B^{(t)}$~"--- два таких
объекта, то, по определению, множество $\mathbf{NAAlg}^{T\text{-}\mathrm{pgr}}_\mathbbm{k}(A,B)$ морфизмов $A\to B$ состоит из всех гомоморфизмов алгебр $\varphi\colon A\to B$, таких, что
$\varphi(A^{(t)})\subseteq B^{(t)}$ для всех $t\in T$.

Обозначим через $U \colon \mathbf{NAAlg}^{T\text{-}\mathrm{pgr}}_\mathbbm{k} 
\to \mathbf{Vect}^{T\text{-}\mathrm{gr}}_\mathbbm{k}$
забывающий функтор, который сопоставляет каждому объекту $A \supseteq \bigoplus_{t\in T} A^{(t)}$  подпространство $\bigoplus_{t\in T} A^{(t)}$, градуированное множеством $T$,
и ограничивает гомоморфизмы на такие подпространства.

Пусть $V= \bigoplus_{t\in T} V^{(t)}$~"--- $T$-градуированное пространство.
Пусть $Y^{(t)}$~"--- некоторые базисы подпространств $V^{(t)}$. Обозначим через $KV$
абсолютно свободную неассоциативную алгебру $\mathbbm{k}\lbrace Y\rbrace$
с множеством свободных порождающих $Y=\bigsqcup_{t\in T} Y^{(t)}$.
В инвариантной форме $$KV = \bigoplus_{n=1}^\infty \bigoplus_{\substack{\text{всевозможные} \\ \text{расстановки} \\
\text{скобок}}} \underbrace{V \otimes \ldots \otimes V}_n,$$
причём умножение задаётся формулой $vw=v\otimes w$ 
(расстановки скобок в обеих частях совпадают).
Отождествим $V$ с соответствующим подпространством в $KV$
и будем рассматривать $KV \supseteq V= \bigoplus_{t\in T} V^{(t)}$
как объект категории $\mathbf{NAAlg}^{T\text{-}\mathrm{pgr}}_\mathbbm{k}$.

Для любого $\varphi \in \mathbf{Vect}^{T\text{-}\mathrm{gr}}_\mathbbm{k}(V,W)$
существует единственный гомоморфизм алгебр $K\varphi \colon KV \to KW$,
такой, что $(K\varphi)\bigl|_{V}=\varphi$.

\begin{proposition}
Функтор $K \colon \mathbf{Vect}^{T\text{-}\mathrm{gr}}_\mathbbm{k} 
\to \mathbf{NAAlg}^{T\text{-}\mathrm{pgr}}_\mathbbm{k}$ является левым сопряжённым к функтору $U \colon \mathbf{NAAlg}^{T\text{-}\mathrm{pgr}}_\mathbbm{k} 
\to \mathbf{Vect}^{T\text{-}\mathrm{gr}}_\mathbbm{k}$.
\end{proposition}
\begin{proof}
Если $V$~"--- объект категории $\mathbf{Vect}^{T\text{-}\mathrm{gr}}_\mathbbm{k}$,
а $A$~"--- объект категории $\mathbf{NAAlg}^{T\text{-}\mathrm{pgr}}_\mathbbm{k}$,
тогда всякий морфизм $KV \to A$ однозначно определяется своим ограничением на $V$.
Отсюда получаем естественную биекцию $\mathbf{NAAlg}^{T\text{-}\mathrm{pgr}}_\mathbbm{k}(KV,A)\mathrel{\widetilde\to} \mathbf{Vect}^{T\text{-}\mathrm{gr}}_\mathbbm{k}(V, UA)$.
\end{proof}

Пусть теперь $V=\bigoplus_{t\in T} V^{(t)}$,
где $V^{(t)}$~"--- векторные пространства с формальными базисами $\left(x_i^{(t)}\right)_{i\in\mathbb N}$.
Тогда алгебра $KV$ может быть отождествлена с алгеброй $\mathbbm{k}\lbrace X^{T\text{-}\mathrm{gr}}\rbrace$
из~\S\ref{SectionGradedPI}.
Всякая $T$-градуированная алгебра $A$ может рассматриваться как объект
категории $\mathbf{NAAlg}^{T\text{-}\mathrm{pgr}}_\mathbbm{k}$,
где подпространство $\bigoplus_{t\in T} A^{(t)}$ совпадает с $A$.
В этом случае биекция $\mathbf{NAAlg}^{T\text{-}\mathrm{pgr}}_\mathbbm{k}(KV,A)\mathrel{\widetilde\to} \mathbf{Vect}^{T\text{-}\mathrm{gr}}_\mathbbm{k}(V, UA)$
означает, что всякое отображение $\psi \colon X^{T\text{-}\mathrm{gr}} \to A$, такое, что
$\psi\left(X^{(t)}\right) \subseteq A^{(t)}$ для всех $t\in T$,
может быть однозначно продолжено до гомоморфизма алгебр $\bar \psi \colon KV \to A$,
такого, что $\bar\psi\left(X^{(t)}\right) \subseteq A^{(t)}$.

\subsection{Обобщённые $H$-действия}\label{SectionHAdjunction}

Пусть $H$~"--- ассоциативная алгебра с единицей над полем $\mathbbm{k}$.
Обозначим через ${}_H \mathcal M$ категорию левых $H$-модулей,
а через ${}_H \mathbf{NAAlgSubMod}$ 
(от англ. ``not necessarily associative algebras with subspaces that are $H$-modules'')
категорию, объектами которой являются все необязательно ассоциативные алгебры $A$ над $\mathbbm{k}$
с фиксированными подпространствами $A_0 \subseteq A$ (включение может быть строгим),
являющимися левыми $H$-модулями, причём для объектов $A \supseteq A_0$
и $B \supseteq B_0$ множество ${}_H \mathbf{NAAlgSubMod}(A,B)$
морфизмов $A\to B$ состоит из всевозможных гомоморфизмов алгебр $\varphi \colon A \to B$,
где $\varphi(A_0)\subseteq B_0$ и $\varphi\bigl|_{A_0}$
является гомоморфизмом $H$-модулей. 

Снова существует очевидный забывающий функтор $U \colon 
{}_H \mathbf{NAAlgSubMod} \to {}_H \mathcal M$, где $UA:=A_0$ и $U\varphi := \varphi\bigl|_{A_0}$.

Пусть $K\colon {}_H \mathcal M \to {}_H \mathbf{NAAlgSubMod}$~"--- функтор,
который сопоставляет каждому левому $H$-модулю $V$ абсолютно свободную неассоциативную алгебру $KV:=\mathbbm{k}\lbrace Z \rbrace$, где $Z$~"--- базис пространства $V$. 
 Другими словами, 
 $$KV = \bigoplus_{n=1}^\infty \bigoplus_{\substack{\text{всевозможные} \\ \text{расстановки} \\
\text{скобок}}} \underbrace{V \otimes \ldots \otimes V}_n,$$
причём умножение задаётся равенством $vw=v\otimes w$ (расстановки скобок в обеих частях совпадают).
Отождествим $V$ с соответствующим подпространством в $KV$
и будем рассматривать $KV \supseteq V$ как объект категории ${}_H \mathbf{NAAlgSubMod}$.
Для любого $\varphi \in {}_H \mathcal M(V,W)$
существует единственный гомоморфизм алгебр $K\varphi \colon KV \to KW$,
такой, что $(K\varphi)\bigl|_{V}=\varphi$.

\begin{proposition}
Функтор $K \colon {}_H \mathcal M \to {}_H \mathbf{NAAlgSubMod}$ является левым сопряжённым
к функтору $U \colon 
{}_H \mathbf{NAAlgSubMod} \to {}_H \mathcal M$.
\end{proposition}
\begin{proof}
Если $V$~"--- объект категории ${}_H \mathcal M$,
а $A$~"--- объект категории ${}_H \mathbf{NAAlgSubMod}$,
то любой морфизм $KV \to A$
однозначно определяется своим ограничением на $V$.
Следовательно, получаем естественную биекцию ${}_H \mathbf{NAAlgSubMod}(KV,A)\mathrel{\widetilde\to} {}_H \mathcal M(V, UA)$.
\end{proof}

Пусть теперь $X$~"--- произвольное множество, а $V$~"--- свободный левый $H$-модуль со свободным $H$-базисом $X$, т.е. $V=\bigoplus_{x\in X} Hx$.
Тогда алгебру $KV$ можно отождествить с алгеброй $\mathbbm{k}\lbrace X | H \rbrace$
из~\S\ref{SectionHPI}.
Любую алгебру $A$ с обобщённым $H$-действием можно рассматривать как объект категории 
${}_H \mathbf{NAAlgSubMod}$, где $H$-модуль $A_0$ совпадает с $A$.
В этом случае биекция ${}_H \mathbf{NAAlgSubMod}(KV,A)\mathrel{\widetilde\to} {}_H \mathcal M(V, UA)$
означает, что любое отображение $\psi \colon X \to A$
однозначно продолжается до гомоморфизма алгебр $\bar \psi \colon KV \to A$,
такого, что $\bar\psi\left(h x\right) = h \bar\psi\left( x\right)$
для всех $x\in X$.

\section{Гипотеза Амицура и её аналоги}\label{SectionAmitsurDefinition}

Напомним, что ассоциативная алгебра называется \textit{PI-алгеброй}, если в ней выполняется хотя бы одно нетривиальное тождество, т.е. такое тождество, которое справедливо не во всех ассоциативных алгебрах.

В 1980-х годах Ш.~Амицур
 выдвинул следующую гипотезу:

\begin{conjecture}[Ш.~Амицур]
Пусть $A$~"--- ассоциативная PI-алгебра над полем характеристики~$0$,
а
$c_n(A)$~"--- последовательность коразмерностей ее полиномиальных
тождеств. Тогда существует
$\PIexp(A):=\lim\limits_{n \to \infty} \sqrt[n]{c_n(A)} \in \mathbb Z_+$.
\end{conjecture}

 Гипотеза Амицура была
доказана М.\,В.~Зайцевым и А.~Джамбруно~\cite[теорема~6.5.2]{ZaiGia}
 в~1999 году для всех ассоциативных алгебр.
Кроме того, в~2002 году М.\,В.~Зайцев~\cite{ZaicevLie}
доказал аналог гипотезы Амицура для коразмерностей полиномиальных тождеств
 конечномерных алгебр Ли. 
 
 Для полиномиальных $H$-тождеств $H$-модульных ассоциативных алгебр аналог гипотезы Амицура
может быть сформулирован в следующей форме, которая принадлежит Ю.\,А.~Бахтурину:

\begin{conjecture}[Ш.~Амицур~"--- Ю.\,А.~Бахтурин]
 Пусть $A$~"--- конечномерная
ассоциативная $H$-модульная алгебра, где $H$~"--- алгебра Хопфа на полем характеристики $0$.
Тогда существует предел $\PIexp^H(A):=\lim\limits_{n\to\infty}
 \sqrt[n]{c^H_n(A)}$, который является целым числом.
\end{conjecture}

В теоремах~\ref{TheoremHOreAmitsur} и~\ref{TheoremHmoduleAssoc} гипотеза Амицура~"--- Бахтурина
доказывается для всех конечномерных $H$-модульных ассоциативных
алгебр $A$, где $H$~"--- алгебра Хопфа над полем характеристики $0$ и либо радикал Джекобсона 
$J(A)$ является $H$-подмодулем, либо $H$ получена при помощи (возможно, многократного)
расширения Оре конечномерной полупростой алгебры Хопфа косопримитивными элементами.
В теореме~\ref{TheoremDoubleNumbersAmitsurPIexpH} гипотеза Амицура~"--- Бахтурина
доказывается для всех $H$-модульных структур с единицей на алгебре двойных чисел~$\mathbbm{k}[x]/(x^2)$.
В общем же случае вопрос о справедливости гипотезы Амицура~"--- Бахтурина остаётся открытым.
 
 Вообще, под справедливостью аналога гипотезы Амицура для полиномиальных тождеств с дополнительной
 структурой ($H$- и $G$-тождеств, дифференциальных, градуированных, \ldots)
 будем понимать существование целой экспоненты роста коразмерностей
 соответствующих типов тождеств. Главы~\ref{ChapterGenHAssocCodim} и~\ref{ChapterHModLieCodim},
 а также часть главы~\ref{ChapterSGGrAssocCodim} посвящены доказательству аналогов гипотезы
 Амицура для широкого класса ассоциативных алгебр и алгебр Ли с дополнительной структурой.
 Кроме этого, в главе~\ref{ChapterSGGrAssocCodim} приводятся примеры конечномерных
 ассоциативных алгебр, градуированных полугруппами, с нецелой градуированной PI-экспонентой.

 \section{Совпадение $H$-коразмерностей для эквивалентных $H$-модульных структур}\label{SectionEquivCodimCoincide}

Понятие эквивалентности модульных структур существенно уменьшает число случаев,
которые необходимо рассматривать при доказательстве аналога гипотезы Амицура,
поскольку $H$-коразмерности эквивалентных $H$-модульных структур совпадают.

Другим приложением понятия эквивалентности модульных структур является
возможность замены рационального действия связной аффинной алгебраической группы автоморфизмами на действие её алгебры Ли дифференцированиями и наоборот, основанная на теореме~\ref{TheoremAffAlgGrAllEquiv}.

\begin{lemma}\label{LemmaHEquivCodimTheSame}
Пусть $\zeta_1 \colon H_1 \to \End_\mathbbm{k}(A_1)$ и $\zeta_2 \colon H_2 \to \End_\mathbbm{k}(A_2)$~"---
гомоморфизмы алгебр, отвечающие эквивалентным модульным структурам на (необязательно ассоциативных) алгебрах $A_1$ и $A_2$.
Тогда существует изоморфизм алгебр $$\mathbbm{k}\lbrace X | H_1\rbrace/\Id^{H_1}(A_1) 
\mathrel{\widetilde\to} \mathbbm{k}\lbrace X | H_2\rbrace/\Id^{H_2}(A_2),$$ который
при всяком $n\in \mathbb N$ сужается до изоморфизма $\mathbbm{k}S_n$-модулей $$\frac{W_n^{H_1}}{W_n^{H_1}\cap\,\Id^{H_1}(A_1)}\mathrel{\widetilde\to} \frac{W_n^{H_2}}{W_n^{H_2}\cap\,\Id^{H_2}(A_2)}.$$
В частности, $c_n^{H_1}(A_1) = c_n^{H_2}(A_2)$ и
$\chi_n^{H_1}(A)=\chi_n^{H_2}(A)$ для всех $n\in \mathbb N$.
 Если, кроме того, $\chr \mathbbm{k} = 0$, то $\ell_n^{H_1}(A)=\ell_n^{H_2}(A)$.
\end{lemma}
\begin{proof} Пусть $\varphi \colon A_1 \mathrel{\widetilde\to} A_2$~"--- эквивалентность
$H_1$- и $H_2$-модульных структур, отвечающих гомоморфизмам $\zeta_1$ и $\zeta_2$, и пусть $\tilde\varphi \colon \End_\mathbbm{k}(A_1) \mathrel{\widetilde\to} \End_\mathbbm{k}(A_2)$~"--- соответствующий изоморфизм алгебр линейных операторов.
Согласно нашим предположениям, $\tilde\varphi(\zeta_1(H_1))=\zeta_2(H_2)$. 
Следовательно, существуют $\mathbbm{k}$-линейные отображения $\xi \colon H_1 \to H_2$ и $\theta \colon H_2 \to H_1$ (которые необязательно являются гомоморфизмами), такие, что диаграмма ниже коммутирует в обоих направлениях:
\begin{equation*}
\xymatrix{ H_1  \ar@<0.5ex>[r]^\xi \ar[d]_{\zeta_1}& H_2 \ar@<0.5ex>[l]^\theta \ar[d]^{\zeta_2} \\
\End_\mathbbm{k}(A_1) \ar[r]^{\tilde\varphi} & \End_\mathbbm{k}(A_2)}
\end{equation*}

Тогда \begin{equation}\label{EqXiActingCodimTheSame}\xi(h)\varphi(a)=\varphi(ha)\text{ для всех }h\in H_1\text{ и }a\in A_1\end{equation}
и  \begin{equation}\label{EqZetaActingCodimTheSame}\varphi(\theta(h)a)=h\varphi(a)\text{ для всех }h\in H_2\text{ и }a\in A_1.\end{equation}

Определим гомоморфизмы алгебр $\tilde\xi \colon \mathbbm{k}\lbrace X | H_1\rbrace
\to \mathbbm{k}\lbrace X | H_2\rbrace$ и $\tilde\theta \colon \mathbbm{k}\lbrace X | H_2\rbrace
\to \mathbbm{k}\lbrace X | H_1\rbrace$
при помощи равенств $\tilde \xi\left(x_k^{h}\right):=x_k^{\xi(h)}$
для всех $h\in H_1$, $k\in\mathbb N$
и $\tilde \theta\left(x_k^h\right):=x_k^{\theta(h)}$
для всех $h\in H_2$, $k\in\mathbb N$. Из равенств~\eqref{EqXiActingCodimTheSame}
и~\eqref{EqZetaActingCodimTheSame} следует, что 
$\tilde\xi\left( \Id^{H_1}(A_1) \right) \subseteq \Id^{H_2}(A_2)$
и 
$\tilde\theta\left( \Id^{H_2}(A_2) \right) \subseteq \Id^{H_1}(A_1)$.
Отсюда $\tilde\xi$ и $\tilde\theta$ индуцируют гомоморфизмы
$\bar\xi \colon \mathbbm{k}\lbrace X | H_1\rbrace/\Id^{H_1}(A_1)
\to \mathbbm{k}\lbrace X | H_2\rbrace/\Id^{H_2}(A_2)$ и $\bar\theta \colon \mathbbm{k}\lbrace X | H_2\rbrace/\Id^{H_2}(A_2)
\to \mathbbm{k}\lbrace X | H_1\rbrace/\Id^{H_1}(A_1)$.
Заметим, что из~\eqref{EqXiActingCodimTheSame}
и~\eqref{EqZetaActingCodimTheSame} также следует, что $\theta(\xi(h))a=ha$ для всех $a\in A_1$
и $h\in H_1$ и $\xi(\theta(h))a=ha$ для всех $a\in A_2$
и $h\in H_2$. Следовательно $x^h - x^{\theta(\xi(h))}  \in \Id^{H_1}(A_1)$ и $x^h - x^{\xi(\theta(h))} \in \Id^{H_2}(A_2)$. Отсюда $\bar\theta\bar\xi = \id_{\mathbbm{k}\lbrace X | H_1\rbrace/\Id^{H_1}(A_1)}$
и $\bar\xi\bar\theta = \id_{\mathbbm{k}\lbrace X | H_2\rbrace/\Id^{H_2}(A_2)}$, т.е. мы получаем требуемый изоморфизм.

Сравнивая степени $H$-многочленов, получаем, что 
$\bar\xi \left(\frac{W_n^{H_1}}{W_n^{H_1}\cap\,\Id^{H_1}(A_1)}\right) \subseteq \frac{W_n^{H_2}}{W_n^{H_2}\cap\,\Id^{H_2}(A_2)}$
и $\bar\theta \left(\frac{W_n^{H_2}}{W_n^{H_2}\cap\,\Id^{H_2}(A_2)}\right) \subseteq \frac{W_n^{H_1}}{W_n^{H_1}\cap\,\Id^{H_1}(A_1)}$ для всех $n\in\mathbb N$, откуда и следует изоморфизм факторпространств и равенство коразмерностей.
\end{proof}

Лемма~\ref{LemmaHEquivCodimTheSame} будет использована ниже в \S\ref{SectionEquivApplToPolyIden}.

Докажем теперь вариант этой леммы для эквивалентных градуировок:

\begin{lemma}\label{LemmaGrEquivCodimTheSame}
Пусть $\Gamma_1 \colon A_1 = \bigoplus_{t\in T_1} A_1^{(t)}$ и 
$\Gamma_2 \colon A_2 = \bigoplus_{t\in T_2} A_2^{(t)}$~"--- групповые градуировки,
а $\varphi \colon A_1 \mathrel{\widetilde\to} A_2$~"--- изоморфизм алгебр над полем, являющийся
эквивалентностью градуировок $\Gamma_1$ и $\Gamma_2$.
Тогда существует изоморфизм алгебр $$\mathbbm{k}\lbrace X^{T_1\text{-}\mathrm{gr}}\rbrace/\Id^{T_1\text{-}\mathrm{gr}}(A_1) 
\mathrel{\widetilde\to} \mathbbm{k}\lbrace X^{T_2\text{-}\mathrm{gr}}\rbrace/\Id^{T_2\text{-}\mathrm{gr}}(A_2),$$ который
при всяком $n\in \mathbb N$ сужается до изоморфизма $\mathbbm{k}S_n$-модулей $$\frac{W_n^{T_1\text{-}\mathrm{gr}}}{W_n^{T_1\text{-}\mathrm{gr}}\cap\,\Id^{T_1\text{-}\mathrm{gr}}(A_1)}\mathrel{\widetilde\to} \frac{W_n^{T_2\text{-}\mathrm{gr}}}{W_n^{T_2\text{-}\mathrm{gr}}\cap\,\Id^{T_2\text{-}\mathrm{gr}}(A_2)}.$$
В частности, $c_n^{T_1\text{-}\mathrm{gr}}(A_1) = c_n^{T_2\text{-}\mathrm{gr}}(A_2)$ и
$\chi_n^{T_1\text{-}\mathrm{gr}}(A)=\chi_n^{T_2\text{-}\mathrm{gr}}(A)$ для всех $n\in \mathbb N$.
 Если, кроме того, $\chr \mathbbm{k} = 0$, то $\ell_n^{T_1\text{-}\mathrm{gr}}(A)=\ell_n^{T_2\text{-}\mathrm{gr}}(A)$.
\end{lemma}
\begin{proof} Эквивалентность $\varphi$ индуцирует биекцию $\psi \colon \supp \Gamma_1 \mathrel{\widetilde\to} \supp \Gamma_2$,
  где $\varphi\left(A_1^{(t)}\right)=A_2^{\bigl(\psi(t)\bigr)}$ для всех $t\in \supp \Gamma_1$.
  
Определим гомоморфизмы алгебр $\xi \colon \mathbbm{k}\lbrace X^{T_1\text{-}\mathrm{gr}}\rbrace
\to \mathbbm{k}\lbrace X^{T_2\text{-}\mathrm{gr}}\rbrace$ и $\theta \colon \mathbbm{k}\lbrace X^{T_2\text{-}\mathrm{gr}}\rbrace
\to \mathbbm{k}\lbrace X^{T_1\text{-}\mathrm{gr}}\rbrace$
при помощи равенств $$\xi\left(x_k^{(t)}\right):=\left\lbrace\begin{array}{lll}x_k^{\bigl( \psi(t) \bigr)} & \text{при} & t\in \supp \Gamma_1, \\
0 & \text{при} & t\notin \supp \Gamma_1
\end{array}\right.$$ 
и
$$\theta\left(x_k^{(t)}\right):=\left\lbrace\begin{array}{lll}x_k^{\bigl( \psi^{-1}(t) \bigr)} & \text{при} & t\in \supp \Gamma_2, \\
0 & \text{при} & t\notin \supp \Gamma_2,
\end{array}\right.$$ 
$k\in\mathbb N$. Из определения биекции $\psi$ следует, что 
$\xi\left( \Id^{T_1\text{-}\mathrm{gr}}(A_1) \right) \subseteq \Id^{T_2\text{-}\mathrm{gr}}(A_2)$
и 
$\theta\left( \Id^{T_2\text{-}\mathrm{gr}}(A_2) \right) \subseteq \Id^{T_1\text{-}\mathrm{gr}}(A_1)$.
Отсюда $\xi$ и $\theta$ индуцируют гомоморфизмы
$$\bar\xi \colon \mathbbm{k}\lbrace X^{T_1\text{-}\mathrm{gr}}\rbrace/\Id^{T_1\text{-}\mathrm{gr}}(A_1)
\to \mathbbm{k}\lbrace X^{T_2\text{-}\mathrm{gr}}\rbrace/\Id^{T_2\text{-}\mathrm{gr}}(A_2)$$ и $$\bar\theta \colon \mathbbm{k}\lbrace X^{T_2\text{-}\mathrm{gr}}\rbrace/\Id^{T_2\text{-}\mathrm{gr}}(A_2)
\to \mathbbm{k}\lbrace X^{T_1\text{-}\mathrm{gr}}\rbrace/\Id^{T_1\text{-}\mathrm{gr}}(A_1).$$

Заметим, что $\theta\left(\xi\left( x^{(t)}_k \right)\right) = x^{(t)}_k$ для всех $t\in \supp \Gamma_1$,
а $\xi\left(\theta\left( x^{(t)}_k \right)\right) = x^{(t)}_k$ для всех $t\in \supp \Gamma_2$
и $k\in\mathbb N$.
 Отсюда $\bar\theta\bar\xi = \id_{\mathbbm{k}\lbrace X^{T_1\text{-}\mathrm{gr}}\rbrace/\Id^{T_1\text{-}\mathrm{gr}}(A_1)}$
и $\bar\xi\bar\theta = \id_{\mathbbm{k}\lbrace X^{T_2\text{-}\mathrm{gr}}\rbrace/\Id^{T_2\text{-}\mathrm{gr}}(A_2)}$, т.е. мы получаем требуемый изоморфизм.

Сравнивая степени градуированных многочленов, получаем, что 
$\bar\xi \left(\frac{W_n^{T_1\text{-}\mathrm{gr}}}{W_n^{T_1\text{-}\mathrm{gr}}\cap\,\Id^{T_1\text{-}\mathrm{gr}}(A_1)}\right) \subseteq \frac{W_n^{T_2\text{-}\mathrm{gr}}}{W_n^{T_2\text{-}\mathrm{gr}}\cap\,\Id^{T_2\text{-}\mathrm{gr}}(A_2)}$
и $\bar\theta \left(\frac{W_n^{T_2\text{-}\mathrm{gr}}}{W_n^{T_2\text{-}\mathrm{gr}}\cap\,\Id^{T_2\text{-}\mathrm{gr}}(A_2)}\right) \subseteq \frac{W_n^{T_1\text{-}\mathrm{gr}}}{W_n^{T_1\text{-}\mathrm{gr}}\cap\,\Id^{T_1\text{-}\mathrm{gr}}(A_1)}$ для всех $n\in\mathbb N$, откуда и следует изоморфизм факторпространств и равенство коразмерностей и кохарактеров.
\end{proof}

  \section{Оценка сверху для $H$-кодлин}\label{SectionUpperHSimple}


В теореме~1 из работы~\cite{ZaiMishchGiaIntermediate} 
А.~Джамбруно, М.\,В.~Зайцев и С.\,П.~Мищенко
доказали, что \begin{equation*}
\ell_n(A)=\sum_{\lambda \vdash n} m(A,\lambda)
\leqslant (\dim A) (n+1)^{(\dim A)^2+\dim A}
\end{equation*}
 для всех $n\in \mathbb N$.
 
 В теореме~\ref{TheoremUpperBoundColengthHNAssoc}
 ниже мы по аналогии доказываем такую же оценку сверху для $H$-кодлин конечномерных алгебр с обобщённым $H$-действием.

Пусть $A$~"--- некоторая (необязательно ассоциативная) конечномерная
алгебра с обобщённым $H$-действием для некоторой ассоциативной алгебры $H$ с $1$
над полем $\mathbbm{k}$ характеристики $0$.

\begin{lemma}\label{LemmaHTensorCommutative} Пусть $C$~"--- коммутативная
ассоциативная $\mathbbm{k}$-алгебра с единицей. Определим на $A\otimes C$
структуру алгебры с обобщённым $H$-действием при помощи равенств $h(a \otimes c)
:= ha \otimes c$, где $a\in A$ и $c\in C$. Тогда $\Id^H(A\otimes C)=\Id^H(A)$.
\end{lemma}
\begin{proof}
Поскольку алгебра $C$ с единицей,
алгебра
 $A\otimes C$ содержит $H$-инвариантную подалгебру, изоморфную $A$,
 откуда $\Id^H(A\otimes C) \subseteq \Id^H(A)$.
Доказательство обратного включения аналогично доказательству
соответствующего включения для ассоциативных алгебр, не наделённых никаким действием,
см. лемму~1.4.2 в~\cite{ZaiGia}.
\end{proof}

Пусть $a_1, \ldots, a_s$~"--- базис в алгебре $A$. Выберем некототрое число $k\in\mathbb N$.
Обозначим через $\mathbbm{k}[\xi_{ij} \mid 1\leqslant i \leqslant s,\ 1\leqslant j \leqslant k ]$
алгебру коммутативных ассоциативных многочленов
от переменных $\xi_{ij}$ со свободными членами и
с коэффициентами из поля $\mathbbm{k}$.
 Алгебра $$A \otimes \mathbbm{k}[\xi_{ij} \mid 1\leqslant i \leqslant s,\ 1\leqslant j \leqslant k ]$$
 является алгеброй с обобщённым $H$-действием, заданным равенством $$h(a \otimes f):=ha \otimes f\text{ при }a\in A\text{ и }f\in  \mathbbm{k}[\xi_{ij} \mid 1\leqslant i \leqslant s,\ 1\leqslant j \leqslant k ].$$
  Обозначим через $\tilde A_k$ пересечение всех $H$-инвариантных подалгебр алгебры $A \otimes \mathbbm{k}[\xi_{ij} \mid 1\leqslant i \leqslant s,\ 1\leqslant j \leqslant k ]$, содержащих элементы $$\xi_j := \sum_{i=1}^s a_i \otimes \xi_{ij}\text{, где }1\leqslant j \leqslant k.$$ (Буквы $\xi_j$ символизируют 
  <<элементы общего вида>> алгебры $A$.)
  
  \begin{lemma}\label{LemmaHRelativeK}
  Пусть $f=f(x_1, \ldots, x_k) \in \mathbbm{k}\lbrace X | H \rbrace$.
  Тогда $f\in \Id^H(A)$, если и только если $f(\xi_1, \ldots, \xi_k)=0$
  (как элемент алгебры $\tilde A_k$).
  \end{lemma}
  \begin{proof} Из леммы~\ref{LemmaHTensorCommutative} следует, что $$\Id^H(A)=\Id^H(A\otimes \mathbbm{k}[\xi_{ij} \mid 1\leqslant i \leqslant s,\ 1\leqslant j \leqslant k ])\subseteq \Id^H(\tilde A_k).$$
  В частности, если $f\in \Id^H(A)$, то $f(\xi_1, \ldots, \xi_k)=0$.
  
  Обратно, предположим, что $f(\xi_1, \ldots, \xi_k)=0$.
  Докажем, что $f(b_1, \ldots, b_k)=0$ для всех $b_j \in A$. Действительно, $b_j=\sum_{i=1}^s \alpha_{ij} a_i$ для некоторых $\alpha_{ij} \in \mathbbm{k}$.
  Рассмотрим гомоморфизм $$\varphi \colon A\otimes \mathbbm{k}[\xi_{ij} \mid 1\leqslant i \leqslant s,\ 1\leqslant j \leqslant k ] \to A$$
  $H$-модульных алгебр, заданный условием $a \otimes \xi_{ij}\mapsto \alpha_{ij} a$ для всех $a\in A$. Тогда 
  $$f(b_1, \ldots, b_k)=f(\varphi(\xi_1),\ldots, \varphi(\xi_k))=\varphi(f(\xi_1, \ldots, \xi_k))=0$$
  и $f\in \Id^H(A)$.   
  \end{proof}
  
    \begin{lemma}\label{LemmaHRelativeKdim}
  Обозначим через $R_{kn}$ линейную оболочку в $\tilde A_k$ всевозможных произведений $(h_1 \xi_{i_1})\ldots (h_n \xi_{i_n})$, где $h_j \in H$ и $1\leqslant i_j \leqslant k$ при $1\leqslant j\leqslant n$.
  Тогда $$\dim R_{kn}\leqslant (\dim A) (n+1)^{k\dim A}\text{ для всех }n\in\mathbb N.$$
  \end{lemma}
\begin{proof}
Пространство $R_{kn} \subseteq A\otimes \mathbbm{k}[\xi_{ij} \mid 1\leqslant i \leqslant s,\ 1\leqslant j \leqslant k ]$ является подпространством линейной оболочки элементов $a_\ell \otimes \prod\limits_{\substack{
1\leqslant i \leqslant s,\\ 1\leqslant j \leqslant k}} \xi_{ij}^{s_{ij}}$,
где $1\leqslant \ell\leqslant s =\dim A$, $s_{ij}\in\mathbb Z_+$, $\sum\limits_{\substack{
1\leqslant i \leqslant s,\\ 1\leqslant j \leqslant k}} s_{ij} = n$.
Число таких элементов не превосходит $(\dim A) (n+1)^{k\dim A}$,
откуда и следует требуемая оценка сверху для $\dim R_{kn}$.
\end{proof}

Теперь докажем, что все неприводимые $\mathbbm{k}S_n$-подмодули,
которые возникают в разложении модуля $\frac{W^H_n}{W^H_n \cap\, \Id^H(A)}$
с ненулевыми кратностями, отвечают диаграммам Юнга высоты равной или меньшей, чем $\dim A$.
  
  \begin{lemma}\label{LemmaHStripeDimATheorem}
    Пусть $\lambda \vdash n$, где $n\in\mathbb N$, причём $\lambda_{(\dim A)+1} > 0$. Тогда $m(A, H,\lambda)= 0$.
  \end{lemma}
  \begin{proof} В силу теоремы~\ref{ThKratnost}
  достаточно доказать, что $e^{*}_{T_\lambda} f \in \Id^H(A)$
    для всех $f\in W_n^H$.
Выберем в $A$ некоторый базис.
Поскольку многочлены полилинейны,
для проверки того, что они являются тождествами, достаточно подставлять
в них только базисные элементы.
 Заметим, что
$e^{*}_{T_\lambda} = b_{T_\lambda} a_{T_\lambda}$,
где оператор $b_{T_\lambda}$ делает многочлены кососимметричными по переменным каждого столбца таблицы Юнга $T_\lambda$.
Следовательно, если после некоторой подстановки элементов алгебры $A$ многочлен $
e^{*}_{T_\lambda} f$ не обращается в нуль,
то это означает, что вместо переменных каждого столбца таблицы Юнга $T_\lambda$
подставляются разные элементы. Однако если $\lambda_{(\dim A)+1} > 0$,
то высота первого столбца больше, чем $\dim A$. Отсюда
 $e^{*}_{T_\lambda} f \in \Id^H(A)$.
  \end{proof}
  
Докажем теперь главный результат данного параграфа:
  
  \begin{theorem}\label{TheoremUpperBoundColengthHNAssoc}
  Пусть $A$~"--- некоторая (необязательно ассоциативная) конечномерная
алгебра с обобщённым $H$-действием для некоторой ассоциативной алгебры $H$ с $1$
над полем $\mathbbm{k}$ характеристики $0$.
Тогда $$\ell_n^H(A) \leqslant (\dim A) (n+1)^{(\dim A)^2+\dim A}$$ для всех $n\in\mathbb N$.
  \end{theorem}
  \begin{proof}
  Для всякого разбиения $\lambda \vdash n$ фиксируем таблицу Юнга $T_\lambda$
  формы $\lambda$.
  Тогда в силу теоремы~\ref{ThKratnost}
кратность $m(A,H,\lambda)$ неприводимого подмодуля $M(\lambda)=\mathbbm{k}S_n e_{T_\lambda}$
в разложении $\mathbbm{k}S_n$-модуля $\frac{W^H_n}{W^H_n \cap\, \Id^H(A)}$ равна $\dim e_{T_\lambda} \frac{W^H_n}{W^H_n \cap\, \Id^H(A)}$.
Другими словами, число $m(A,H,\lambda)$
равно максимальному числу $m$ таких $H$-многочленов $f_1,\ldots, f_m \in W_n^H$,
что из равенства $$g=\alpha_1 e_{T_\lambda}f_1 + \ldots +\alpha_m e_{T_\lambda}f_m \in \Id^H(A)$$
для некоторых $\alpha_\ell\in \mathbbm{k}$ всегда следует, что $\alpha_1=\ldots = \alpha_m = 0$.
Обозначим через $k_{ij}$ число в клетке $(i,j)$ таблицы $T_\lambda$.
Тогда для любого фиксированного $i$ каждый многочлен $e_{T_\lambda}f_\ell$
симметричен по переменным $ x_{k_{i1}}, \ldots, x_{k_{i\lambda_i}} $.
Используя процесс линеаризации (см., например, \cite[\S 1.3]{ZaiGia}), получаем,
что $H$-многочлен $g$ является полиномиальным $H$-тождеством, если и только если 
$\tilde g$ является полиномиальным $H$-тождеством, где $H$-многочлен $\tilde g$
получен из $g$ подстановкой $x_{k_{ij}} \mapsto x_i$ для всех $i$ и $j$.
Обозначим число строк в таблице Юнга $T_\lambda$ через $k$. 
В силу леммы~\ref{LemmaHStripeDimATheorem} можно считать, что $k\leqslant \dim A$.
 Тогда $H$-многочлен $\tilde g$ зависит от переменных $x_1,\ldots, x_k$ и
 из леммы~\ref{LemmaHRelativeK} следует, что $\tilde g \in \Id^H(A)$
 если и только $\tilde g(\xi_1, \ldots, \xi_k) = 0$ как элемент алгебры $\tilde A_k$.
 Заметим, что $\tilde g(\xi_1, \ldots, \xi_k) = \alpha_1  u_1 + \ldots + \alpha_m u_m$, где
 $u_\ell$~"--- значение $H$-многочлена $e_{T_\lambda}f_\ell$ при подстановке
 $x_{k_{ij}} \mapsto \xi_i$, где $1\leqslant i \leqslant k$ и $1\leqslant j \leqslant \lambda_i$. Следовательно, все $u_i \in R_{kn}$, и если $m >  (\dim A) (n+1)^{k\dim A}$,
 то согласно лемме~\ref{LemmaHRelativeKdim} для любого выбора $H$-многочленов $f_i$ элементы $u_i$
линейно зависимы и $\tilde g(\xi_1, \ldots, \xi_k) = \alpha_1  u_1 + \ldots + \alpha_m u_m = 0$
 для некоторых нетривиальных коэффициентов $\alpha_i$. В частности, $\alpha_1 e_{T_\lambda}f_1 + \ldots \alpha_m e_{T_\lambda}f_m \in \Id^H(A)$ и $m(A,H,\lambda) < m$. Следовательно, для всех $\lambda \vdash n$
 справедливо неравенство $$m(A,H,\lambda)\leqslant (\dim A) (n+1)^{k\dim A} \leqslant (\dim A) (n+1)^{(\dim A)^2}.$$ Поскольку число всех разбиений $\lambda \vdash n$,
 высота которых не больше $\dim A$, не превосходит $(n+1)^{\dim A}$,
 мы получаем требуемую верхнюю оценку для $\ell_n^H(A)$.
\end{proof}

 Если конечномерная алгебра $A$
градуирована некоторым множеством $T$, то
в силу предложения~\ref{PropositionCnGrCnGenH}
 кодлины $\ell_n^{T\text{-}\mathrm{gr}, H}(A)$
градуированных полиномиальных тождеств алгебры $A$ равны её $\mathbbm{k}^T$-кодлинам
  $\ell_n^{\mathbbm{k}^T}(A)$. 
Отсюда получаем следующее следствие из теоремы~\ref{TheoremUpperBoundColengthHNAssoc}:
\begin{corollary}\label{CorollaryUpperBoundColengthGraded}
Пусть $A$~"--- конечномерная алгебра над полем характеристики $0$,
градуированная некоторым множеством $T$. Тогда 
  $$\ell_n^{T\text{-}\mathrm{gr}}(A) \leqslant (\dim A) (n+1)^{(\dim A)^2+\dim A}$$ для всех $n\in\mathbb N$.
\end{corollary}

\section{Разбиения, ограниченные выпуклыми многогранниками}

В данном параграфе мы применим идеи из~\cite{VerZaiMishch} 
для того, чтобы доказать лемму~\ref{LemmaRegionUpperFd}, которая будет затем использована при исследовании
алгебр с необязательно целой PI-экспонентой.

Для начала нам понадобится следующая лемма о выпуклых многогранниках:
\begin{lemma}\label{LemmaDistanceToPolyhedron}
Пусть $P$~"--- непустой выпуклый многогранник в $\mathbb R^n$, заданный системой линейных неравенств  $f_i \leqslant 0$, где $1\leqslant i \leqslant m$.
Тогда для любого $\varepsilon > 0$ существует такое $\delta > 0$, что
расстояние  от любой точки многогранника $$P_\delta = \lbrace  x \mid f_i(x) \leqslant \delta \text{ для всех } 1\leqslant i \leqslant m\rbrace$$
до $P$
 меньше либо равно, чем $\varepsilon$.
\end{lemma}
\begin{proof} Для произвольного подмножества $I\subseteq \lbrace 1, 2, \ldots, m \rbrace$
введём обозначение $\pi_I = \lbrace x\in\mathbb R^n \mid f_i(x)=0 \text{ для всех } i\in I \rbrace$.

Если для какого-то $1\leqslant i \leqslant m$ пересечение $P\cap \pi_{\lbrace i \rbrace}$
пусто, это означает, что неравенство $f_i(x)\leqslant 0$ в определении многоугольника $P$
лишнее и его можно исключить. Поэтому без ограничения общности можно считать,
что $P\cap \pi_{\lbrace i \rbrace} \ne \varnothing$ для всех $1\leqslant i \leqslant m$.

В силу линейности функций $f_i$ можно выбрать такую
общую константу $M_1 > 0$, что $$|f_i(x) - f_i(y)|\leqslant M_1 \rho(x,y)$$
для всех $1\leqslant i \leqslant m$ и $x,y\in\mathbb R^n$,
где $\rho(x,y)=\sqrt{(x^1-y^1)^2+\ldots+(x^n-y^n)^2}$~"--- евклидова метрика в $\mathbb R^n$.

Используя формулу расстояния от точки до гиперплоскости,
выберем такое $M_2 > 0$, что для всех $I\subseteq \lbrace 1, 2, \ldots, m \rbrace$ 
и всех $1\leqslant j \leqslant m$, таких, что функция $f_j$ непостоянна на $\pi_I$,
расстояние от произвольной точки $x \in \pi_I$ до плоскости $\pi_{I\cup \lbrace j \rbrace}$
меньше либо равно $M_2|f_j(x)|$.

Положим $\delta = \frac{\varepsilon}{M_2 \sum_{j=1}^m  (1+M_1M_2)^{j-1}}$ и рассмотрим произвольную точку
$x_0 \in P_\delta$. Если для всех $1\leqslant i \leqslant m$ справедливо
неравенство $f_i(x_0) \leqslant 0$, то $x_0 \in P$ и доказывать
нечего. Предположим, что $f_i(x_0) > 0$ для некоторого $1\leqslant i \leqslant m$.
Пусть $x_1$~"--- ближайшая к $x_0$ точка гиперплоскости~$\pi_{\lbrace i \rbrace}$.
Тогда $\rho(x_0, x_1) \leqslant M_2 |f_i(x_0)| \leqslant M_2\delta$.
При этом $x_1 \in P_{(1 + M_1 M_2) \delta}$.
Теперь рассмотрим вместо $x_0$ точку $x_1$, вместо $\mathbb R^n$~"--- гиперплоскость
$\pi_{\lbrace i \rbrace}$, а вместо $P$~"--- выпуклый многогранник $P\cap \pi_{\lbrace i \rbrace}$
и повторим рассуждения, построив точку $x_2 \in P_{(1 + M_1 M_2)^2 \delta}$,
и т.д. На каждом шаге размерность пространства и число неравенств будут убывать, и в концев
концов мы прийдём к случаю, когда очередная точка $x_k$, где $0\leqslant k \leqslant m$, окажется внутри многогранника. Тогда $$\rho(x_0,x_k) \leqslant \sum_{j=1}^{k} \rho(x_{j-1}, x_j)
\leqslant  M_2 \sum_{j=1}^k  (1+M_1M_2)^{j-1}  \delta \leqslant  \varepsilon.$$
\end{proof}

Покажем теперь, что если для алгебры $A$ с обобщённым $H$-действием все разбиения $\lambda\vdash n$, которые отвечают неприводимым $\mathbbm{k}S_n$-модулям,
встречающимся с ненулевыми кратностями $m(A,H,\lambda)$ в разложении
$\mathbbm{k}S_n$-модуля $\frac{W_n^H}{W_n^H\cap\, \Id^H(A)}$, принадлежат некоторому выпуклому многограннику $\Omega_n$, то число $\mathop{\overline\lim}_{n\to\infty}\sqrt[n]{c_n^{H}(A)}$ 
ограниченно сверху максимумом некоторой функции $\Phi$ на <<непрерывной>> версии $\Omega$ этого многогранника. (По сути данный результат носит сугубо комбинаторный характер, поскольку, как это будет видно из доказательства, он справедлив для размерностей любой последовательности $\mathbbm{k}S_n$-модулей, у которых кратности неприводимых подмодулей ограничены сверху полиномиальной функцией от $n$.)

Пусть $q\in\mathbb N$. Определим функцию
$\Phi(x_1, \ldots, x_q):=\frac{1}{x_1^{x_1} \ldots x_s^{x_q}}$ при $x_1, \ldots, x_q > 0$.
Поскольку $\lim\limits_{x\to +0} x^x = 1$, без ограничения общности можно считать, что $\Phi$
является непрерывной функцией на множестве 
$\lbrace (x_1, \ldots, x_q ) \mid x_i \geqslant 0\rbrace$.

Для фиксированных чисел $\gamma_{ij} \in \mathbb R$, где $1\leqslant i \leqslant m$,
$m \in \mathbb Z_+$, $0\leqslant j \leqslant q$,
определим множество
$$\Omega := \left\lbrace
(\alpha_1, \ldots, \alpha_q)\in \mathbb R^q \mathrel{\Biggl|}\sum_{i=1}^q \alpha_i=1,\ \alpha_1\geqslant\alpha_2\geqslant \ldots \geqslant \alpha_q\geqslant 0,
\ \sum_{j=1}^q \gamma_{ij}\alpha_j \geqslant 0\text{ при } 1\leqslant i \leqslant m\right\rbrace.$$

Пусть также заданы числа $\theta_k \in \mathbb N$, где $q< k \leqslant r$, а $r\in \mathbb Z_+$.

Для всякого $n\in \mathbb N$ определим множество
$$\Omega_n := \left\lbrace
\lambda \vdash n \mathrel{\Biggl|} \sum_{j=1}^q \gamma_{ij}\lambda_j+\gamma_{i0} \geqslant 0 \text{ при } 1\leqslant i \leqslant m,\ \lambda_i \leqslant \theta_i \text{ при } q < i \leqslant r,\ 
\lambda_{r+1}=0 \right\rbrace.$$

В дальнейшем будем относиться к $\Omega$ и $\Omega_n$
как к, соответственно, <<непрерывной>> и <<дискретной>> версиям одного и того же многогранника.

Обозначим через $d$ максимум функции $\Phi$ на компактном множестве $\Omega$, которое предполагаем непустым.

\begin{lemma}\label{LemmaRegionUpperFd}
Пусть $A$~"--- некоторая (необязательно ассоциативная) конечномерная
алгебра с обобщённым $H$-действием для некоторой ассоциативной алгебры $H$ с $1$
над полем~$\mathbbm{k}$ характеристики $0$, причём $m(A, H, \lambda)=0$ для всех $\lambda\vdash n$, где $\lambda 
\notin \Omega_n$, $n\in\mathbb N$. Тогда 
$\mathop{\overline\lim}_{n\to\infty}\sqrt[n]{c_n^{H}(A)}
\leqslant d$.
\end{lemma}
\begin{proof}  Рассмотрим такое разбиение $\lambda \vdash n$, что $m(A,H,\lambda)\ne 0$.
В силу формулы крюков $\dim M(\lambda)=\frac{n!}{\prod_{i,j} h_{ij}}$,
где $h_{ij}$~"--- длина крюка с вершиной в клетке $(i,j)$
диаграммы Юнга $D_\lambda$. Следовательно,
$\dim M(\lambda) \leqslant \frac{n!}{\lambda_1! \ldots \lambda_r!}$.
Теперь заметим, что $\left(x^x\right)'=\left(e^{x\ln x}\right)'=
(\ln x+1)e^{x\ln x}$ и функция $x^x$ убывает при $x\leqslant \frac{1}{e}$.
В силу формулы Стирлинга при некоторых $C_1, C_2 > 0$ и $r_1, r_2 \in\mathbb R$, которые не зависят от $\lambda_i$, и всех достаточно больших $n$ справедливы неравенства \begin{equation}\begin{split}\label{EqMlambdaUpperFdOmega}\dim M(\lambda) \leqslant \frac{C_1 n^{r_1}
\left(\frac{n}{e}\right)^n}{\left(\frac{\lambda_1}{e}\right)^{\lambda_1}\ldots
\left(\frac{\lambda_r}{e}\right)^{\lambda_r}}=C_1 n^{r_1}\left(\frac{1}
{\left(\frac{\lambda_1}{n}\right)^{\frac{\lambda_1}{n}}\ldots
\left(\frac{\lambda_r}{n}\right)^{\frac{\lambda_r}{n}}}\right)^n
\leqslant \\ \leqslant C_1 n^{r_1} \left(\Phi\left(\frac{\lambda_1}{n}, \ldots, \frac{\lambda_q}{n}\right)\right)^n
\frac{n^{\theta_{q+1}+\ldots +\theta_r}}{\theta_{q+1}^{\theta_{q+1}} \ldots \theta_r^{\theta_r}}=C_2 n^{r_2} \left(\Phi\left(\frac{\lambda_1}{n}, \ldots, \frac{\lambda_q}{n}\right)\right)^n.\end{split}\end{equation}

Пусть $\varepsilon > 0$. Поскольку функция $\Phi$ непрерывна на компакте $[0; 1]^q \supseteq \Omega$, она равномерно непрырывна на нём и существует такое $\delta > 0$,
что для всех $x$ из $[0; 1]^q$, таких, что расстояние между
$x$ и $\Omega$ меньше, чем $\delta$, справедливо неравенство $\Phi(x) < d + \varepsilon$.

Следовательно, в силу~(\ref{EqMlambdaUpperFdOmega}) и леммы~\ref{LemmaDistanceToPolyhedron} существует такое $n_0\in\mathbb N$,
что для всех $n \geqslant n_0$ и $\lambda\vdash n$, таких, что $m(A,H,\lambda)\ne 0$,
справедливо неравенство
 $ \dim M(\lambda) \leqslant C_2 n^{r_2} (d+\varepsilon)^n$.

Согласно теореме~\ref{TheoremUpperBoundColengthHNAssoc}
существуют такие $C_3 > 0$  и $r_3\in\mathbb Z_+$, что
 $$\sum_{\lambda \vdash n} m(A,H,\lambda)
\leqslant C_3 n^{r_3}\text{ для всех }n \in \mathbb N.$$

Отсюда $$ c^{H}_n(A) = \sum_{\lambda \vdash n} m(A,H,\lambda) \dim M(\lambda)
 \leqslant  C_2 C_3 n^{r_2+r_3} (d+\varepsilon)^n$$
 и $\mathop{\overline\lim}_{n\to\infty}\sqrt[n]{c_n^H(A)}
\leqslant d+\varepsilon$. Поскольку число $\varepsilon > 0$ было выбрано произвольным, утверждение леммы доказано.
\end{proof}
  
     \section{Существование $H$-PI-экспоненты у $H$-простых алгебр}\label{SectionHPIexpExistHSimple}

Оказывается,
что у всякой конечномерной $H$-простой алгебры существует
 $H$-PI-экспонента.
  
   \begin{theorem}\label{TheoremHSimpleHPIexpHNAssoc} Пусть $A$~"--- конечномерная $H$-простая
(необязательно ассоциативная) алгебра с обобщённым $H$-действием для некоторой ассоциативной алгебры $H$ с $1$ над полем $\mathbbm{k}$ характеристики $0$, $\dim A = s$. Введём
обозначение $$d(A) := \mathop{\overline\lim}_{n\to\infty}\max\limits_{\substack{\lambda \vdash n, \\ m(A,H,\lambda)\ne 0}}
  \Phi\left(\frac{\lambda_1}{n},\ldots, \frac{\lambda_s}{n}\right).$$
  Тогда существует (необязательно целая) $$\PIexp^{H}(A) := \lim_{n\to \infty} \sqrt[n]{c_n^H(A)} = d(A).$$
  \end{theorem}
  
  Теорема~\ref{TheoremHSimpleHPIexpHNAssoc} будет доказана ниже.
  
  Снова используя предложение~\ref{PropositionCnGrCnGenH},
  получаем из теоремы~\ref{TheoremHSimpleHPIexpHNAssoc} следствие для градуированных тождеств:
  
  \begin{corollary}\label{CorollaryGradedExistsExponent}
  Пусть $A$~"--- конечномерная алгебра над полем характеристики $0$,
  градуированная некоторым множеством $T$ таким образом, что $A$
является $T$-градуированно простой, т.е. не содержит нетривиальных градуированных идеалов. Тогда существует $\PIexp^{T\text{-}\mathrm{gr}}(A)=\lim_{n\to \infty} \sqrt[n]{c_n^{T\text{-}\mathrm{gr}}(A)}$.
  \end{corollary}
   
   Докажем сперва, что последовательность $H$-коразмерностей
   не убывает для любой $H$-простой алгебры.
  
   \begin{lemma}\label{LemmaCodimNondecrHsimpleHNAssoc} Пусть $A$~"--- $H$-простая 
алгебра для некоторой ассоциативной алгебры $H$ с $1$ над произвольным полем $\mathbbm{k}$.
  Тогда $c_n^H(A) \leqslant c_{n+1}^H(A)$ для всех $n\in\mathbb N$.
  \end{lemma}
  \begin{proof} Пусть $n\in\mathbb N$,
  а $f_1(x_1, \ldots, x_n)$, \ldots, $f_{c_n^H(A)}(x_1, \ldots, x_n)$~"--- такие $H$-многочлены,
  что их образы являются базисом пространства $\frac{W_n^H}{W_n^H \cap \Id^H(A)}$.
  Предположим, что $H$-многочлены \begin{equation}\label{EqPolynomialsMonotoneCodim}f_1(x_1, \ldots, x_n x_{n+1}), \ldots, f_{c_n^H(A)}(x_1, \ldots, x_n x_{n+1})\end{equation} линейно зависимы по модулю $\Id^H(A)$.
  Тогда существуют такие $\alpha_1, \ldots, \alpha_{c_n^H(A)} \in \mathbbm{k}$, что $$\alpha_1 f_1(a_1, \ldots, a_n a_{n+1})+ \ldots + \alpha_{c_n^H(A)} f_{c_n^H(A)}(a_1, \ldots, a_n a_{n+1}) = 0$$
  для всех $a_i \in A$. Поскольку алгебра $A$ является $H$-простой, справедливо равенство $AA=A$ и 
  $$\alpha_1 f_1(a_1, \ldots, a_n)+ \ldots + \alpha_{c_n^H(A)} f_{c_n^H(A)}(a_1, \ldots, a_n) = 0$$
  для всех $a_i \in A$. Однако $$f_1(x_1, \ldots, x_n), \ldots, f_{c_n^H(A)}(x_1, \ldots, x_n)$$
  линейно независимы по модулю $\Id^H(A)$. Отсюда $$\alpha_1= \ldots= \alpha_{c_n^H(A)}=0,$$
  $H$-многочлены~\eqref{EqPolynomialsMonotoneCodim} линейно независимы по модулю $\Id^H(A)$ и $c_n^H(A) \leqslant c_{n+1}^H(A)$.
  \end{proof} 
  
 Теперь докажем для $c_n^H(A)$ оценку сверху:
  
  \begin{theorem}\label{TheoremUpperBoundCodimPhiHNAssoc} 
  Пусть $A$~"--- конечномерная $H$-простая
(необязательно ассоциативная) алгебра для некоторой ассоциативной алгебры $H$ с $1$ над полем $\mathbbm{k}$ характеристики $0$, $\dim A = s$.
  Тогда существуют такие $C > 0$ и $r\in \mathbb R$, что
  $$c_n^H(A) \leqslant Cn^r \left(\max\limits_{\substack{\lambda \vdash n, \\ m(A,H,\lambda)\ne 0}}
  \Phi\left(\frac{\lambda_1}{n},\ldots, \frac{\lambda_s}{n}\right)\right)^n \text{ для всех } n\in\mathbb N.$$
  \end{theorem} 
\begin{proof}
Пусть $\lambda \vdash n$~"--- такое разбиение, что $m(A,H,\lambda)\ne 0$.
Согласно формуле крюков $\dim M(\lambda)=\frac{n!}{\prod_{i,j} h_{ij}}$,
где $h_{ij}$~"--- длина крюка с вершиной в клетке $(i,j)$
диаграммы Юнга $D_\lambda$. Отсюда
$\dim M(\lambda) \leqslant \frac{n!}{\lambda_1! \ldots \lambda_s!}$.
В силу формулы Стирлинга для всех достаточно больших $n$ справедливы неравенства \begin{equation}\begin{split}\label{EqMlambdaUpperFd}\dim M(\lambda) \leqslant \frac{C_1 n^{r_1}
\left(\frac{n}{e}\right)^n}{\left(\frac{\lambda_1}{e}\right)^{\lambda_1}\ldots
\left(\frac{\lambda_s}{e}\right)^{\lambda_s}}=C_1 n^{r_1}\left(\frac{1}
{\left(\frac{\lambda_1}{n}\right)^{\frac{\lambda_1}{n}}\ldots
\left(\frac{\lambda_s}{n}\right)^{\frac{\lambda_s}{n}}}\right)^n
\leqslant \\ \leqslant C_1 n^{r_1} \left(\Phi\left(\frac{\lambda_1}{n}, \ldots, \frac{\lambda_s}{n}\right)\right)^n \end{split}\end{equation}
при некоторых $C_1 > 0$ и $r_1 \in\mathbb R$, которые не зависят от $\lambda_i$.
Теперь утверждение теоремы следует из теоремы~\ref{TheoremUpperBoundColengthHNAssoc}.
\end{proof}

До конца параграфа считаем выполненными предположения теоремы~\ref{TheoremHSimpleHPIexpHNAssoc}.

Пусть $\lambda \vdash n$, $\mu \vdash m$, 
 $\mathbbm{k}S_n \bar f_1 \cong M(\lambda)$ и $\mathbbm{k}S_m \bar f_2 \cong M(\mu)$
для некоторых $m,n \in\mathbb N$, $f_1 \in W_n^H$ и $f_2 \in W_m^H$.
Тогда образ $H$-многочлена $f_1(x_1, \ldots, x_n)f_2(x_{n+1}, \ldots, x_{m+n})$
порождает $\mathbbm{k}S_{m+n}$-подмодуль в $\frac{W^H_{m+n}}{W^H_{m+n} \cap\, \Id^H(A)}$,
являющийся гомоморфным образом $\mathbbm{k}S_{m+n}$-модуля $$M(\lambda) \hatotimes M(\mu) := (M(\lambda) \otimes M(\mu))
\uparrow S_{m+n} := \mathbbm{k}S_{m+n} \mathbin{\otimes_{\mathbbm{k}(S_n \times S_m)}} (M(\lambda) \otimes M(\mu)).$$
В силу правила Литтлвуда~"--- Ричардсона
все неприводимые компоненты в разложении модуля $M(\lambda) \hatotimes M(\mu)$
отвечают диаграммам Юнга $D_\nu$, полученным из
$D_{\lambda+\mu}$ перемещением некоторых клеток вниз.
Согласно нашим предположениям высота диаграммы $D_{\nu}$
не может быть больше, чем $s=\dim A$.
Другое замечание, которое здесь нужно сделать, заключается в том,
что при перемещении клеток диаграммы $D_{\nu}$ вниз
значение функции $\Phi\left(\frac{\nu_1}{n}, \frac{\nu_2}{n}, \ldots, \frac{\nu_s}{n}\right)$
не убывает, поскольку функция $\frac{1}{x^x(\xi-x)^{\xi-x}}$ возрастает на интервале $x \in \left(0; \frac{\xi}2 \right)$ при фиксированном $0 < \xi \leqslant 1$.

\begin{lemma}\label{LemmaSequenceLambdaLowerHSimpleNAssoc}
Для любого $\varepsilon > 0$ 
существует такое число $t \in \mathbb N$,
такие натуральные числа $n_1 < n_2 < n_3 < \ldots$, где $n_{i+1} - n_i \leqslant t$,
и такие разбиения $\lambda^{(i)} \vdash n_i$, где $m\left(A, H, \lambda^{(i)} \right)\ne 0$,
что $\Phi\left(\frac{\lambda^{(i)}_1}{n_i}, \ldots, \frac{\lambda^{(i)}_s}{n_i}\right) \geqslant
d(A)-\varepsilon$ для всех $i\in \mathbb N$.
\end{lemma}
\begin{proof}
Обозначим через $B$ подалгебру ассоциативной алгебры $\End_\mathbbm{k}(A)$,
порождённую всеми операторами левого и правого умножения на элементы алгебры $A$
и образами элементов алгебры $H$, которые рассматриваются как операторы на $A$.
Поскольку $A$ является $H$-простой алгеброй,
алгебра $A$~"--- неприводимый левый $B$-модуль.
Обозначим через $\mathcal B$ множество всех таких операторов $u\in B$,
что для некоторых $a_1, \ldots, a_m, \tilde a_1, \ldots, \tilde a_{\tilde m}
 \in A$, $m, \tilde m\in \mathbb Z_+$, $h \in H$
  и некоторой расстановки скобок оператор $u$ представляется в виде
$ua = a_1 \ldots a_m (h a) \tilde a_1 \ldots \tilde a_{\tilde m}$ для всех $a\in A$.
 Используя условие~\eqref{EqGeneralizedHopf},
 можно включить операторы $H$-действия внутрь
операторов левого и правого умножения на элементы алгебры $A$ и
представить любой оператор из $B$ как линейную
комбинацию операторов из $\mathcal B$. 
Поскольку алгебра $\End_\mathbbm{k}(A)$ конечномерна,
можно выбрать некоторый конечный базис $u_1, \ldots, u_{\dim B} \in \mathcal B$ пространства $B$.
 Обозначим через $N$ максимальное из всех чисел
$2(m + \tilde m)$, соответствующих элементам $u_i$.

Поскольку алгебра $A$ является $H$-простой, $A^2\ne 0$ и для всех $a,b\in A$, где $a\ne 0$, $b\ne 0$, 
справедливы равенство $A=Ba=Bb$ и неравенство $(Ba)(Bb)\ne 0$. В силу выбора числа $N$
существуют такие $a_1, \ldots, a_m, \tilde a_1, \ldots, \tilde a_{\tilde m}, b_1,\ldots, b_k,
\tilde b_1, \ldots, \tilde b_{\tilde k}
 \in A$, $k,\tilde k, m, \tilde m\in \mathbb Z_+$, $h_1, h_2 \in H$, что
 для некоторой расстановки скобок
$$(a_1 \ldots a_m (h_1 a) \tilde a_1 \ldots \tilde a_{\tilde m})(b_1 \ldots b_k (h_2 b) 
\tilde b_1 \ldots \tilde b_{\tilde k})\ne 0,$$
а $k+\tilde k + m +\tilde m \leqslant N$.

Будем теперь выбирать такое $q \in\mathbb N$, что $\Phi\left(\frac{\mu_1}{q}, \ldots, \frac{\mu_s}{q}\right) \geqslant d(A)-\varepsilon/2$ и $m(A,H,\mu)\ne 0$ для некоторого $\mu \vdash q$. 
Напомним, что функция $\Phi$ является непрерывной на множестве $[0;1]^s$
и, следовательно, равномерно непрерывна на
$[0;1]^s$, поскольку множество $[0;1]^s$ компактно.
Отсюда значение $\Phi$ мало изменяется при малых изменениях значений аргумента.

В силу того, что число $q$ можно брать сколь угодно большим,
мы можем также считать, что $\frac{\sum_{j=1}^i d_j}{iq}$
будет мало при любом $i\in\mathbb N$ и любом выборе
значений $0 \leqslant d_i \leqslant N$.
Следовательно, существуют такие $q$ и $\mu \vdash q$,
что
 \begin{equation}\begin{split}\label{EquationPhiLambda}\Phi\left(\frac{i\mu_1+\sum_{j=1}^i d_j}{
 iq+\sum_{j=1}^i d_j}, \frac{i\mu_2}{iq+\sum_{j=1}^i d_j}, \ldots, \frac{i\mu_s}{iq+\sum_{j=1}^i d_j}\right)= \\ =
\Phi\left(\frac{\frac{\mu_1}{q}+\frac{\sum_{j=1}^i d_j}{iq}}{
1+\frac{\sum_{j=1}^i d_j}{iq}}, \frac{\left(\frac{\mu_2}{q}\right)}{
1+\frac{\sum_{j=1}^i d_j}{iq}}, \ldots, \frac{\left(\frac{\mu_s}{q}\right)}{1+\frac{\sum_{j=1}^i d_j}{iq}}\right)
\geqslant d(A)-\varepsilon\end{split}\end{equation} для всех $i\in\mathbb N$ и всех $0 \leqslant d_i \leqslant N$.

Положим $t:= N+q$.

Выберем такой $H$-многочлен $\tilde f \in W^H_{q} \backslash \Id^H(A)$, что $\mathbbm{k}S_q \tilde f \cong M(\mu)$.
Тогда для некоторой расстановки скобок, 
некоторых $h_1, h_2 \in H$ и некоторых $k, \tilde k, m, \tilde m \geqslant 0$, таких, что $d_1:=k+ \tilde k + m+ \tilde m
\leqslant N$, выполнено условие $$f_1:=\left(y_1 \ldots y_k \tilde f^{h_1}(x_1, \ldots, x_q)
\tilde y_1 \ldots \tilde y_{\tilde k} \right)\left( z_1 \ldots z_m \tilde f^{h_2}(\tilde x_1, \ldots, \tilde x_q)
\tilde z_1 \ldots \tilde z_{\tilde m}\right) \notin \Id^H(A).$$
(Обозначение $f^h$ было введено в замечании~\ref{RemarkHActionOnWHn} выше.)

Рассмотрим $\mathbbm{k}S_{q+k+\tilde k}$-подмодуль $M$ модуля $\frac{W^H_{q+k+\tilde k}}{
W^H_{q+k+\tilde k} \cap\,\Id^H(A)}$
порождённый образом $H$-многочлена
$$y_1 \ldots y_k \tilde f^{h_1}(x_1, \ldots, x_q)
\tilde y_1 \ldots \tilde y_{\tilde k}.$$
В силу замечания ~\ref{RemarkHActionOnWHn}
образ $H$-многочлена $\tilde f^{h_1}(x_1, \ldots, x_q)$ порождает $\mathbbm{k}S_q$-подмодуль,
изоморфный $M(\mu)$,
откуда $M$ является гомоморфным образом модуля $$M(\mu)\mathrel{\widehat\otimes}
\mathbbm{k}S_{k+\tilde k} := (M(\mu)\otimes
\mathbbm{k}S_{k+\tilde k}) \uparrow \mathbbm{k}S_{q+k+\tilde k}.$$
Поскольку все диаграммы Юнга, отвечающие разбиениям числа $k+\tilde k$, получены
из строки длиной $k+\tilde k$ перемещением некоторых клеток,
в силу правила Литтлвуда~"--- Ричардсона
все диаграммы Юнга, отвечающие неприводимым модулям, встречающимся в разложении модуля $M$,
получены из диаграммы Юнга, отвечающей разбиению $(\mu_1+ k+\tilde k, \mu_2, \ldots, \mu_s)$
перемещением некоторых клеток вниз.
Те же самые аргументы можно применить к  $H$-многочлену
$z_1 \ldots z_m \tilde f^{h_2}(\tilde x_1, \ldots, \tilde x_q)
\tilde z_1 \ldots \tilde z_{\tilde m}$. 

Пусть $n_1 := 2q+d_1$, а
$\lambda^{(1)}$~"--- одно из разбиений, отвечающих
неприводимым компонентам в разложении модуля $\mathbbm{k}S_{n_1}\bar f_1$.
Тогда в силу~(\ref{EquationPhiLambda}) и замечаний, сделанных выше и перед леммой, справедливо неравенство $\Phi\left(\frac{\lambda^{(1)}_1}{n_1}, \ldots, \frac{\lambda^{(1)}_s}{n_s}\right) \geqslant d(A)-\varepsilon$.

Снова получаем, что $$f_2:=\left(y_1 \ldots y_k f_1^{h_1}(x_1, \ldots, x_q)
\tilde y_1 \ldots \tilde y_{\tilde k} \right)\left( z_1 \ldots z_m \tilde f^{h_2}(\tilde x_1, \ldots, \tilde x_q)
\tilde z_1 \ldots \tilde z_{\tilde m}\right) \notin \Id^H(A)$$
для некоторой расстановки скобок, некоторых элементов $h_1, h_2 \in H$ и некоторых чисел $k, \tilde k, m, \tilde m \geqslant 0$ (которые могут отличными от тех, что использовались в $f_1$),  таких, что $d_2:=k+ \tilde k + m+ \tilde m
\leqslant N$. Как и в случае с $n_1$, определим
$n_2$ по формуле $n_2 := 3q+d_1+d_2$.
Обозначим через $\lambda^{(2)}$
одно из разбиений, отвечающих неприводимым компонентам в разложении модуля $\mathbbm{k}S_{n_2}\bar f_2$.
Продолжение данной процедуры до бесконечности завершает доказательство леммы.
\end{proof}

\begin{proof}[Доказательство теоремы~\ref{TheoremHSimpleHPIexpHNAssoc}.]
Фиксируем $\varepsilon > 0$.
Возьмём $n_i \in \mathbb N$ и $\lambda^{(i)} \vdash n_i$ из леммы~\ref{LemmaSequenceLambdaLowerHSimpleNAssoc}.
Тогда
\begin{equation}\label{EqCnHi(A)LowerPhi}\begin{split} 
c_{n_i}^H(A) \geqslant \dim M(\lambda^{(i)}) = \frac{n_i!}{\prod_{j,k} h_{jk}}
  \geqslant \frac{n_i!}{(\lambda^{(i)}_1+s-1)! \ldots (\lambda^{(i)}_s+s-1)!} \geqslant \\ \geqslant 
  \frac{n_i!}{n_i^{s(s-1)}\lambda^{(i)}_1! \ldots \lambda^{(i)}_s!} \geqslant
  \frac{C_1 n_i^{r_1} 
\left(\frac{n_i}{e}\right)^{n_i}}{\left(\frac{\lambda^{(i)}_1}{e}\right)^{\lambda^{(i)}_1}\ldots
\left(\frac{\lambda^{(i)}_s}{e}\right)^{\lambda^{(i)}_s}} \geqslant \\ \geqslant C_1 n_i^{r_1}\left(\frac{1}
{\left(\frac{\lambda^{(i)}_1}{n_i}\right)^{\frac{\lambda^{(i)}_1}{n_i}}\ldots
\left(\frac{\lambda^{(i)}_s}{n_i}\right)^{\frac{\lambda^{(i)}_s}{n_i}}}\right)^{n_i} = C_1 n_i^{r_1}
 \left(\Phi\left(\frac{\lambda^{(i)}_1}{n_i}, \ldots, \frac{\lambda^{(i)}_s}{n_i}\right)\right)^{n_i}\end{split}\end{equation}
 для некоторых $C_1 > 0$ и $r_1 \leqslant 0$, которые не зависят от $i$.
 
 Пусть $n \geqslant n_1$. Тогда $n_i \leqslant n < n_{i+1}$ для некоторого $i\in\mathbb N$.
 Используя~\eqref{EqCnHi(A)LowerPhi}, лемму~\ref{LemmaCodimNondecrHsimpleHNAssoc}
 и тот факт, что $\Phi(x_1, x_2, \ldots, x_s) \geqslant 1$
 при $0\leqslant x_1,\ldots, x_s \leqslant 1$, получаем
   \begin{equation*}\begin{split} \sqrt[n]{c_n^H(A)} \geqslant \sqrt[n]{c_{n_i}^H(A)}\geqslant
  \sqrt[n]{C_1 n_i^{r_1}} \left(\Phi\left(\frac{\lambda^{(i)}_1}{n_i}, \ldots, \frac{\lambda^{(i)}_s}{n_i}\right)\right)^\frac{n_i}{n} \geqslant \\
\geqslant  
     \sqrt[n]{C_1 n^{r_1}}
 \left(\Phi\left(\frac{\lambda^{(i)}_1}{n_i}, \ldots, \frac{\lambda^{(i)}_s}{n_i}\right)\right)^{\frac{n-t}{n}}  \geqslant  \sqrt[n]{C_1 n^{r_1}}
 \left(d(A)-\varepsilon\right)^{\frac{n-t}{n}}.
 \end{split}\end{equation*}
  Следовательно, $$\mathop{\underline\lim}_{n\to\infty} \sqrt[n]{c_n^H(A)} \geqslant d(A)-\varepsilon.$$
  Поскольку число $\varepsilon > 0$ было произвольным, получаем, что $\mathop{\underline\lim}_{n\to\infty}
  \sqrt[n]{c_n^H(A)} \geqslant d(A)$.
  Теперь из теоремы~\ref{TheoremUpperBoundCodimPhiHNAssoc} следует, что
   $\lim_{n\to\infty} \sqrt[n]{c_n^H(A)} = d(A)$.
 \end{proof}

\newpage

\chapter{Рост коразмерностей полиномиальных $H$-тождеств  в ассоциативных алгебрах с (обобщённым) $H$-действием}\label{ChapterGenHAssocCodim}

В данной главе аналог гипотезы Амицура (существование целочисленной экспоненты) доказывается для достаточно широкого класса конечномерных ассоциативных алгебр с обобщённым $H$-действием,
включающего в себя $H$-модульные алгебры с $H$-инвариантным радикалом, алгебры с действием произвольной группы автоморфизмами и антиавтоморфизмами, алгебры с действием произвольной алгебры Ли дифференцированиями, алгебры, градуированные произвольными группами, и алгебры с действием
(возможно, многократного) расширения Оре конечномерной полупростой алгебры Хопфа косопримитивными элементами.

Результаты главы были опубликованы в работах~\cite{ASGordienko3, ASGordienko8, ASGordienko6Kochetov, ASGordienko11, ASGordienko12, ASGordienko9, ASGordienko15, ASGordienko20ALAgoreJVercruysse}.

\section{Кососимметрические многочлены}\label{SectionAssocAlt}

Оценка снизу для коразмерностей $H$-тождеств, которая будет получена в \S\ref{SectionAssocUpperLower}--\ref{SectionAssocFinishingTheProof}, опирается
на существование $H$-многочленов, кососимметричных по большому числу наборов переменных.
В доказательстве существования таких многочленов мы будем
увеличивать число наборов кососимметричных переменных при помощи формы следа.
По всей видимости, впервые эта идея была использована Ю.\,П.~Размысловым~\cite[\S 14]{RazmyslovBook}.

Пусть $B$~"--- конечномерная полупростая (в обычном смысле) $H$-простая ассоциативная алгебра над алгебраически замкнутым полем $\mathbbm{k}$ характеристики $0$ с обобщённым действием некоторой ассоциативной
алгебры $H$ с $1$.
 
Обозначим через $\varphi \colon B \to \End_\mathbbm{k}(B)$
левое регулярное представление алгебры $B$, т.е.
$\varphi(a)b=ab$ для всех $a,b \in B$,
а через $\psi \colon B \to \End_\mathbbm{k}(B)$~"--- 
правое регулярное представление алгебры $B$, т.е.
$\psi(a)b=ba$ для всех $a,b \in B$.
Через $\rho$ обозначим действие $H \to \End_\mathbbm{k}(B)$.

Из~(\ref{EqGeneralizedHopf}) следует, что \begin{equation}\label{EqOpComm1}
\rho(h)\varphi(a)=
 \sum_i \Bigl(\varphi(h'_i a)\rho(h''_i)+\psi(h''''_i a)\rho(h'''_i)\Bigr),
\end{equation}
\begin{equation}\label{EqOpComm2}
\rho(h)\psi(a)=
 \sum_i \Bigl(\psi(h''_i a)\rho(h'_i)+\varphi(h'''_i a)\rho(h''''_i)\Bigr),
\end{equation}
\begin{equation}\label{EqOpComm3}
\varphi(a)\psi(b)=\psi(b)\varphi(a)
\end{equation}
для всех $a,b \in B$.

\begin{lemma}\label{LemmaAssocForm}
Билинейная форма $\tr(\varphi(\cdot)\varphi(\cdot))$
невырождена на $B$.
\end{lemma}
\begin{proof}
Напомним, что алгебра $B$ полупроста.
Следовательно, в силу теоремы Веддербёрна~"--- Артина
  $$B\cong M_{k_1}(\mathbbm{k})\oplus M_{k_2}(\mathbbm{k})
\oplus \ldots \oplus M_{k_s}(\mathbbm{k})$$
для некоторых $k_i \in \mathbb N$, где
 $M_{k_i}(\mathbbm{k})$~"--- алгебры все матриц размера $k_i\times k_i$.

Фиксируем в $B$  базис, состоящий из матричных
единиц
 $e_{\alpha\beta}^{(i)}$ алгебр $M_{k_i}(\mathbbm{k})$.
Заметим, что $$\tr(\varphi(e_{\alpha\beta}^{(i)}))=\left\lbrace
\begin{array}{lll} 0, & \text{ если } & \alpha \ne \beta, \\
 k_i, & \text{ если } & \alpha = \beta.
\end{array} \right.$$ Следовательно, $\tr(\varphi(e_{\alpha\beta}^{(i)})
 \varphi(b)) \ne 0$
 для некоторого базисного элемента
 $b$,  если и только
если $b=e_{\beta\alpha}^{(i)}$.
Следовательно, матрица билинейной формы $\tr(\varphi(\cdot)\varphi(\cdot))$
невырождена.
\end{proof}

\begin{remark}\label{RemarkOtherLetters}
По определению пространство $P^H_n$ состоит из всех полилинейных ассоциативных $H$-многочленов
от переменных $x_1, \ldots, x_n$. Однако для удобства в качестве переменных многочленов $f\in P^H_n$ мы будем использовать и другие буквы, например, $y_i$, $z_i$, $u_i$, $v_i$, $w_i$ и т.д., считая последние буквы просто другими обозначениями для $x_j$.
\end{remark}

Лемма~\ref{LemmaAssocAlternateFirst} является аналогом
леммы~1 из работы~\cite{GiaSheZai}:

\begin{lemma}\label{LemmaAssocAlternateFirst}
Пусть $b_1, \ldots, b_\ell$~"--- базис алгебры $B$.
 Тогда для некоторого $T\in\mathbb N$ 
 существует $H$-многочлен
 $$f=f(x_1, \ldots, x_\ell,
y_1, \ldots, y_\ell, z_1, \ldots, z_T, z) \in P^H_{2\ell+T+1},$$
кососимметричный по переменным
  $\lbrace x_1, \ldots, x_\ell \rbrace$
   и по переменным $\lbrace y_1, \ldots, y_\ell \rbrace$,
     где $\ell=\dim B$, удовлетворяющий следующему условию:
     существуют такие $ \bar z_1, \ldots, \bar z_T \in B$,
     что для любого $\bar z \in B$ справедливо равенство
$f(b_1, \ldots, b_\ell,
b_1, \ldots, b_\ell, \bar z_1, \ldots, \bar z_T, \bar z) = \bar z$.
\end{lemma}
\begin{proof}
Поскольку алгебра $B$ является $H$-простой,
в силу теоремы плотности алгебра
 $\End_\mathbbm{k}(B) \cong M_\ell(\mathbbm{k})$
порождена операторами из подпространств $\rho(H)$,
$\varphi(B)$ и $\psi(B)$.
Из~(\ref{EqOpComm1})--(\ref{EqOpComm3}) следует, что
\begin{equation}\label{EqEndB}
\End_\mathbbm{k}(B) = \langle \varphi(a)\psi(b)\rho(h)
\mid a,b \in B, h \in H \rangle_\mathbbm{k}.
\end{equation}

Рассмотрим многочлен Регева
$$ f_\ell(x_1, \ldots, x_{\ell^2}; y_1, \ldots, y_{\ell^2})
=\sum_{\substack{\sigma \in S_\ell, \\ \tau \in S_\ell}}
 (\sign(\sigma\tau))
x_{\sigma(1)}\ y_{\tau(1)}\ x_{\sigma(2)}x_{\sigma(3)}x_{\sigma(4)}
\ y_{\tau(2)}y_{\tau(3)}y_{\tau(4)}\ldots
$$
$$ x_{\sigma\left(\ell^2-2\ell+2\right)}\ldots
 x_{\sigma\left(\ell^2\right)}
\ y_{\tau\left(\ell^2-2\ell+2\right)}\ldots
 y_{\tau\left(\ell^2\right)}.
$$ Этот многочлен является центральным (см., например, \cite[теорема~5.7.4]{ZaiGia})
для алгебры $M_\ell(\mathbbm{k})$, т.е. $f_\ell$ не является полиномиальным тождеством для алгебры $M_\ell(\mathbbm{k})$
и все значения многочлена $f_\ell$ принадлежат центру алгебры $M_\ell(\mathbbm{k})$.
В силу кососимметричности многочлена $f_\ell$
по $x_1, \ldots, x_{\ell^2}$ и $y_1, \ldots, y_{\ell^2}$ этот многочлен не обращается
в нуль, если только если в переменные каждого набора подставляются линейно независимые элементы. 

Заметим, что $\ker \varphi =0 $, поскольку алгебра $B$ полупроста и, следовательно,
с единицей. 
Отсюда отображение $\varphi$ является инъективным.
В силу~(\ref{EqEndB}) в $\End_\mathbbm{k}(B)$ можно выбрать базис, состоящий из элементов 
$$\varphi(b_1), \ldots, \varphi(b_\ell),
 \quad \varphi(b_{i_1})
\psi(b_{k_1})\rho(h_1),\quad
\ldots, \quad  \varphi(b_{i_s})
\psi(b_{k_s})\rho(h_s)$$ для подходящих $i_t,k_t \in \lbrace 1,2, \ldots, \ell \rbrace$,
 $h_t\in H$.
 Заменим в функции $f_\ell$
 переменные $x_i$ на $\varphi(x_i)$, а переменные $y_i$~"--- на $\varphi(y_i)$ для всех $1 \leqslant i \leqslant \ell$. Кроме того, заменим
 $x_{\ell+j}$
 на $\varphi(z_j)
\psi(u_j)\rho(h_j)$, а
$y_{\ell+j}$~"--- на $\varphi(v_j)
\psi(w_{j})\rho(h_j)$ для всех $1 \leqslant j \leqslant s$.
Теперь $x_i$, $y_i$, $z_i$, $u_i$, $v_i$, $w_i$~"--- переменные, принимающие значения
в $B$. Обозначим получившуюся функцию
через $\tilde f_\ell$.
Если подставить $z_t=v_t=b_{i_t}$,
$u_t=w_t=b_{k_t}$ при
$1 \leqslant t \leqslant s$ и
$x_i=y_i = b_i$ при $1 \leqslant i \leqslant \ell$,
тогда $\tilde f_\ell$ скалярным оператором $(\mu \id_{B})$ на $B$,
где $\mu\in \mathbbm{k}$, $\mu \ne 0$.
Введём теперь новую переменную $z$.
Пользуясь~(\ref{EqOpComm1})--(\ref{EqOpComm3}), переместим в выражении
$\tilde f_\ell \cdot z$ все $\rho(h_i)$
вправо и перепишем $\varphi(\ldots)$
и $\psi(\ldots)$ как операторы левого и правого умножения.
Тогда $f := \mu^{-1} \tilde f_\ell \cdot z$ станет $H$-многочленом, принадлежащим пространству $P^H_n$
для подходящего $n \in \mathbb N$. Пусть $T := 4s$. Переименуем переменные $u_i$, $v_i$, $w_i$
в $z_j$ для $s+1 \leqslant j \leqslant T$.
 Тогда $f$ удовлетворяет всем условиям леммы.
\end{proof}

Пусть $k\ell \leqslant n$ для некоторых $k,n \in \mathbb N$.
 Обозначим через $Q^H_{\ell,k,n} \subseteq P^H_n$
подпространство, состоящее из всех многочленов, кососимметричных
по каждому из $k$ попарно непересекающихся наборов переменных $\{x^i_1, \ldots, x^i_\ell \}
\subseteq \lbrace x_1, x_2, \ldots, x_n\rbrace$, где $1 \leqslant i \leqslant k$.

Теорема~\ref{TheoremAssocAlternateFinal}
является аналогом теоремы~1 из~\cite{GiaSheZai}.

\begin{theorem}\label{TheoremAssocAlternateFinal}
Пусть $B$~"--- конечномерная полупростая (в обычном смысле) $H$-простая ассоциативная алгебра над алгебраически
замкнутым полем $\mathbbm{k}$ характеристики $0$ с обобщённым $H$-действием для некоторой ассоциативной
алгебры $H$ с $1$, а $b_1, \ldots, b_\ell$~"--- базис алгебры $B$.
Тогда существуют такие $T \in \mathbb Z_+$ и  $\bar z_1, \ldots, \bar z_T \in B$,
что для любого $k \in \mathbb N$
существует такой $H$-многочлен $$f=f(x_1^1, \ldots, x_\ell^1; \ldots;
x^{2k}_1, \ldots,  x^{2k}_\ell;\ z_1, \ldots, z_T;\ z) \in Q^H_{\ell, 2k, 2k\ell+T+1},$$
что для любого $\bar z \in B$ справедливо равенство
$f(b_1, \ldots, b_\ell; \ldots;
b_1, \ldots, b_\ell; \bar z_1, \ldots, \bar z_T; \bar z) = \bar z$.
\end{theorem}
\begin{proof}
Пусть $f_1=f_1(x_1,\ldots, x_\ell,\ y_1,\ldots, y_\ell,
z_1, \ldots, z_T, z)$~"--- $H$-многочлен из леммы~\ref{LemmaAssocAlternateFirst},
кососимметричный по $x_1,\ldots, x_\ell$ и по $y_1,\ldots, y_\ell$.
В силу того, что $f_1$ удовлетворяет всем условиям теоремы при $k=1$,
можно считать, что $k > 1$. Заметим, что \begin{equation*}\begin{split}
f^{(1)}_1(u_1, v_1, x_1, \ldots, x_\ell,\ y_1,\ldots, y_\ell,
z_1, \ldots, z_T, z) :=\\=
\sum^\ell_{i=1} f_1(x_1, \ldots, u_1 v_1 x_i,  \ldots, x_\ell,\ y_1,\ldots, y_\ell,
z_1, \ldots, z_T, z)\end{split}\end{equation*}
также кососимметричен по $x_1,\ldots, x_\ell$ и по $y_1,\ldots, y_\ell$.
Кроме того, \begin{equation*}\begin{split}
f^{(1)}_1(\bar u_1, \bar v_1, \bar x_1, \ldots, \bar x_\ell,\
\bar y_1,\ldots, \bar y_\ell,
\bar z_1, \ldots, \bar z_T, \bar z) =\\=
 \tr(\varphi( \bar u_1) \varphi(\bar v_1))
f_1(\bar x_1, \bar x_2, \ldots, \bar x_\ell,\ \bar y_1,\ldots, \bar y_\ell,
\bar z_1, \ldots, \bar z_T, \bar z)
\end{split}\end{equation*}
 для любой подстановки элементов из $B$
 вместо своих переменных, поскольку без ограничения общности
 мы можем считать, что $\bar x_1, \ldots, \bar x_\ell$~"--- различные элементы базиса.

Положим \begin{equation*}\begin{split}
f^{(j)}_1(u_1, \ldots, u_j, v_1, \ldots, v_j, x_1, \ldots, x_\ell,\ y_1,\ldots, y_\ell,
z_1, \ldots, z_T, z) :=\\=
\sum^\ell_{i=1} f^{(j-1)}_1(u_1, \ldots,  u_{j-1}, v_1, \ldots, v_{j-1},
 x_1, \ldots, u_j v_j x_i,  \ldots, x_\ell,\ y_1,\ldots, y_\ell,
z_1, \ldots, z_T, z)\end{split}\end{equation*}
для всех $2 \leqslant j \leqslant \ell$.
Снова
\begin{equation}\begin{split}\label{EqKillingAssoc}
f^{(j)}_1(\bar u_1, \ldots, \bar u_j, \bar v_1, \ldots, \bar v_j, \bar x_1, \ldots, \bar x_\ell,\ \bar y_1,\ldots, \bar y_\ell, \bar z_1, \ldots, \bar z_T, \bar z) =\\=
\tr(\varphi( \bar u_1) \varphi(\bar v_1))
 \tr(\varphi( \bar u_2) \varphi(\bar v_2))
 \ldots
 \tr(\varphi( \bar u_j) \varphi(\bar v_j))
 \cdot
\\
\cdot
f_1(\bar x_1, \bar x_2, \ldots, \bar x_\ell,\ \bar y_1,\ldots, \bar y_\ell,
\bar z_1, \ldots, \bar z_T, \bar z).
\end{split}\end{equation}

Заметим, что $\det(\tr(\varphi(b_i) \varphi(b_j)))_{i,j=1}^\ell \ne 0$, поскольку форма $\tr(\varphi(\cdot) \varphi(\cdot))$ в силу леммы~\ref{LemmaAssocForm} невырожденная.
Положим
\begin{equation*}\begin{split} f_2(u_1, \ldots, u_\ell, v_1, \ldots, v_\ell,
x_1, \ldots, x_\ell, y_1, \ldots, y_\ell, z_1, \ldots, z_T, z) :=
\\=\frac{1}{\ell!\det(\tr(\varphi(b_i) \varphi(b_j)))_{i,j=1}^\ell}\sum_{\sigma, \tau \in S_\ell}
\sign(\sigma\tau)
f^{(\ell)}_1(u_{\sigma(1)}, \ldots, u_{\sigma(\ell)}, v_{\tau(1)}, \ldots, v_{\tau(\ell)}, \\
 x_1, \ldots, x_\ell,\ y_1,\ldots, y_\ell,
z_1, \ldots, z_T, z).\end{split}\end{equation*}
Тогда $f_2 \in Q^H_{\ell, 4, 4\ell+T+1}$.
Подставим $x_i=y_i=u_i=v_i=b_i$, $1 \leqslant i \leqslant \ell$.
Выберем также значения $z_j=\bar z_j$, где $1 \leqslant j \leqslant T$, таким
образом, чтобы $$f_1(b_1, \ldots, b_\ell, b_1, \ldots, b_\ell,
\bar z_1, \ldots, \bar z_T, \bar z) = \bar z \text{\quad для всех\quad} \bar z \in B.$$
Докажем, что $$f_2(b_1, \ldots, b_\ell, b_1, \ldots, b_\ell, b_1, \ldots, b_\ell, b_1, \ldots, b_\ell,
\bar z_1, \ldots, \bar z_T, \bar z) = \bar z.$$

Действительно,
\begin{equation*}\begin{split} f_2(b_1, \ldots, b_\ell, b_1, \ldots, b_\ell,
b_1, \ldots, b_\ell, b_1, \ldots, b_\ell, \bar z_1, \ldots, \bar z_T, \bar z) = \\ =\frac{1}{\ell!\det(\tr(\varphi(b_i) \varphi(b_j)))_{i,j=1}^\ell}\sum_{\sigma, \tau \in S_\ell}
\sign(\sigma\tau)
f^{(\ell)}_1(b_{\sigma(1)}, \ldots, b_{\sigma(\ell)}, b_{\tau(1)}, \ldots, b_{\tau(\ell)}, \\
b_1, \ldots, b_\ell,\ b_1,\ldots, b_\ell,
\bar z_1, \ldots, \bar z_T, \bar z).\end{split}\end{equation*} Используя~(\ref{EqKillingAssoc}), получаем
\begin{equation*}\begin{split} f_2(b_1, \ldots, b_\ell, b_1, \ldots, b_\ell,
b_1, \ldots, b_\ell, b_1, \ldots, b_\ell, \bar z_1, \ldots, \bar z_T, \bar z) = \\ = \frac{1}{\ell!\det(\tr(\varphi(b_i) \varphi(b_j)))_{i,j=1}^\ell}
\sum_{\sigma, \tau \in S_\ell}
\sign(\sigma\tau) \tr(\varphi(b_{\sigma(1)}) \varphi(b_{\tau(1)}))
  \ldots \tr(\varphi(b_{\sigma(\ell)}) \varphi(b_{\tau(\ell)}))\cdot \\ \cdot
f_1(b_1, \ldots, b_\ell,\ b_1,\ldots, b_\ell,
\bar z_1, \ldots, \bar z_T, \bar z).
\end{split}\end{equation*}
 Заметим, что
\begin{equation*}\begin{split}\sum_{\sigma, \tau \in S_\ell}
\sign(\sigma\tau) \tr(\varphi(b_{\sigma(1)}) \varphi(b_{\tau(1)}))
  \ldots \tr(\varphi(b_{\sigma(\ell)}) \varphi(b_{\tau(\ell)}))
=\\=\sum_{\sigma, \tau \in S_\ell}
\sign(\sigma\tau) \tr(\varphi(b_{1}) \varphi(b_{\tau\sigma^{-1}(1)}))
  \ldots \tr(\varphi(b_{\ell}) \varphi(b_{\tau\sigma^{-1}(\ell)}))
\mathrel{\stackrel{(\tau'=\tau\sigma^{-1})}{=}}\\=\sum_{\sigma, \tau' \in S_\ell}
\sign(\tau') \tr(\varphi(b_{1}) \varphi(b_{\tau'(1)}))
  \ldots \tr(\varphi(b_{\ell}) \varphi(b_{\tau'(\ell)}))
=\\=\ell!\det(\tr(\varphi(b_i) \varphi(b_j)))_{i,j=1}^\ell.\end{split}\end{equation*}
Следовательно, $$ f_2(b_1, \ldots, b_\ell, b_1, \ldots, b_\ell,
b_1, \ldots, b_\ell, b_1, \ldots, b_\ell, \bar z_1, \ldots, \bar z_T, \bar z) = \bar z. $$
Заметим, что если $H$-многочлен $f_1$ был кососимметричен по каким-то из переменных $z_1,\ldots, z_T$,
то $H$-многочлен $f_2$ также кососимметричен по этим переменным.
Поэтому применяя вышеизложенную процедуру к
$f_2$ вместо $f_1$, получаем $f_3 \in Q^H_{\ell, 6, 6\ell+T+1}$.
Аналогичным образом $H$-многочлен $f_4$ определяется с использованием $H$-многочлена $f_3$, $H$-многочлен $f_5$~"--- с использованием $H$-многочлена $f_4$ и т.д.
В конце концов получается требуемый $H$-многочлен
$f:=f_k \in Q^H_{\ell, 2k, 2k\ell+T+1}$.
\end{proof}

    \section{Свойство (*)} \label{SectionPropertyStar}
    
    Для того, чтобы сформулировать предположения, при которых будет справедлив основной результат главы,
    дадим следующее определение:
 
\begin{enumerate}
\item[(*)] Пусть $B$~"--- конечномерная ассоциативная
алгебра с обобщённым $H$-действием для некоторой ассоциативной алгебры $H$ с единицей
над полем $\mathbbm{k}$. Пусть $b_1, \ldots, b_\ell$~"--- базис алгебры $B$.
Будем говорить, что для алгебры $B$ выполняется свойство (*), если $B$~"--- алгебра
с единицей и существует такое число $n_0 \in \mathbb N$, что для любого $k \in\mathbb N$
существует полилинейный $H$-многочлен $$f=f(x_1^{(1)}, \ldots, x_\ell^{(1)}; \ldots;
x^{(2k)}_1, \ldots,  x^{(2k)}_\ell;\ z_1, \ldots, z_{n_1})$$
и такие элементы $\bar z_i \in B$, где $1\leqslant i \leqslant n_1$, $0\leqslant n_1 \leqslant n_0$,
что $H$-многочлен $f$ кососимметричен по $x_1^{(i)}, \ldots, x_\ell^{(i)}$
для любого $1\leqslant i  \leqslant 2k$
и $f(b_1, \ldots, b_\ell; \ldots;
b_1, \ldots, b_\ell; \bar z_1, \ldots, \bar z_{n_1}) \ne 0$.
\end{enumerate}

\begin{example}\label{ExampleHSimpleGenHSemiSimplePropertyStar}
Если $B$~"--- конечномерная полупростая $H$-простая
алгебра с обобщённым $H$-действием для некоторой ассоциативной алгебры $H$ с единицей
над алгебраически замкнутым полем характеристики $0$,
то в силу теоремы~\ref{TheoremAssocAlternateFinal} для алгебры $B$ выполнено свойство (*). 
\end{example}

\begin{example}\label{ExampleHmodHSemiSimplePropertyStar}
Если $B$~"--- конечномерная $H$-простая $H$-модульная алгебра
для конечномерной полупростой алгебры Хопфа $H$
над алгебраически замкнутым полем характеристики $0$,
то в силу следствия~\ref{CorollaryRadicalHSSSubMod} и примера~\ref{ExampleHSimpleGenHSemiSimplePropertyStar} для алгебры $B$ выполнено свойство (*). 
\end{example}

Покажем теперь, что свойство (*) для конечномерных $H$-простых алгебр наследуется при расширениях Оре алгебры Хопфа $H$ косопримитивными элементами:

\begin{theorem}\label{TheoremHOrePropertyStar}
Пусть $H$~"--- алгебра Хопфа над алгебраически замкнутым полем $\mathbbm{k}$ характеристики $0$,
порождённая как алгебра подалгеброй Хопфа $\tilde H$ и косопримитивным элементом $v\in H$,
где $\Delta v = g_1 \otimes v + v \otimes g_2$ для некоторых $g_1, g_2 \in G(\tilde H)$.
Кроме того, предположим, что существует такой автоморфизм алгебры $\varphi \colon \tilde H \to \tilde H$,
что $vh-\varphi(h)v\in \tilde H$ для всех $h\in\tilde H$.
Допустим, что для всех конечномерных $\tilde H$-простых
алгебр Хопфа выполнено свойство (*).
Тогда для всех конечномерных $H$-простых
алгебр Хопфа также выполнено свойство (*).
\end{theorem}
\begin{corollary}\label{CorollaryOreExtStar}
Если $B$~"--- конечномерная $H$-простая $H$-модульная алгебра
для такой алгебры Хопфа $H$
над алгебраически замкнутым полем характеристики $0$,
что $H$ получена при помощи (возможно, многократного)
расширения Оре конечномерной полупростой алгебры Хопфа косопримитивными элементами,
то для алгебры $B$ выполнено свойство (*).
\end{corollary}
\begin{corollary}\label{CorollaryTaftPropertyStar}
Если $B$~"--- конечномерная $H_{m^2}(\zeta)$-простая $H_{m^2}(\zeta)$-модульная алгебра
над алгебраически замкнутым полем характеристики $0$ для
алгебры Тафта $H_{m^2}(\zeta)$, то для алгебры $B$ выполнено свойство (*).
\end{corollary}

\begin{proof}[Доказательство теоремы~\ref{TheoremHOrePropertyStar}.]
Пусть $A$~"--- конечномерная $H$-простая алгебра. В силу леммы~\ref{LemmaHSemiSimpleIsUnital}
алгебра $A$ обладает единицей.

Если $A$ является $\tilde H$-простой, то согласно предположениям теоремы алгебра $A$ удовлетворяет свойству (*). Поэтому в силу леммы~\ref{LemmaHOreRadical} можно без ограничения общности считать, что $J^{\tilde H}(A) \ne 0$. Из теоремы~\ref{TheoremHSimpleOreExtStructure}
следует, что
существует такое $m\in\mathbb N$, что
$A = \bigoplus_{k=0}^{m-1}v^k \tilde J$, а $J^{\tilde H}(A)=
\bigoplus_{k=0}^{m-2}v^k \tilde J$.
Более того, $A/J^{\tilde H}(A)$~"--- $\tilde H$-простая алгебра.
Согласно предположениям теоремы~\ref{TheoremHOrePropertyStar}
для алгебры $A/J^{\tilde H}(A)$ выполнено свойство (*).
Пусть $a_1, \ldots, a_\ell$~"--- базис идеала $\tilde J$.
Тогда образы $b_1, \ldots, b_\ell$ элементов $v^{m-1}a_1, \ldots, v^{m-1}a_\ell$
в $A/J^{\tilde H}(A)$ являются базисом алгебры $A/J^{\tilde H}(A)$.

В силу свойства (*) существует такое
$n_0 \in \mathbb N$, что для любого $k \in\mathbb N$
существует полилинейный $\tilde H$-многочлен $$\tilde f= \tilde f(x_1^{(1)}, \ldots, x_\ell^{(1)}; \ldots;
x^{(2k)}_1, \ldots,  x^{(2k)}_\ell;\ z_1, \ldots, z_{n_1})$$
и элементы $\bar z_i \in A/J^{\tilde H}(A)$, где $1\leqslant i \leqslant n_1$, $0\leqslant n_1 \leqslant n_0$, что $\tilde f$ кососимметричен по $x_1^{(i)}, \ldots, x_\ell^{(i)}$
для любого $1\leqslant i  \leqslant 2k$,
а $$b:=\tilde f(b_1, \ldots, b_\ell; \ldots;
b_1, \ldots, b_\ell; \bar z_1, \ldots, \bar z_{n_1}) \ne 0.$$
Поскольку $\tilde H$-инвариантный идеал, порождённый
элементом $b$ совпадает с $A/J^{\tilde H}(A)$, а $\left(A/J^{\tilde H}(A)\right)^2=A/J^{\tilde H}(A)$, 
существуют такие $h_1, \ldots, h_m
\in\tilde H$ и $\bar w_1, \ldots, \bar w_{m-1} \in A/J^{\tilde H}(A)$, что $$(h_1 b) \bar w_1 (h_2 b) \bar w_2 \ldots (h_{m-1} b) \bar w_{m-1} (h_m b) \ne 0.$$
Определим следующие полилинейные функции, которые с учётом замечания~\ref{RemarkHActionOnWHn} реализуются
как, соответственно, $\tilde H$- и $H$-многочлены:
 \begin{equation*}\begin{split}
\tilde f_0 := \left(\tilde f(x_{11}^{(1)}, \ldots, x_{1\ell}^{(1)}; \ldots;
x^{(2k)}_{11}, \ldots,  x^{(2k)}_{1\ell};\ z_{11}, \ldots, z_{1n_1})\right)^{h_1}\cdot\\ \cdot \prod_{i=2}^m w_{i-1} \left(\tilde f(x_{i1}^{(1)}, \ldots, x_{i\ell}^{(1)}; \ldots;
x^{(2k)}_{i1}, \ldots,  x^{(2k)}_{i\ell};\ z_{i1}, \ldots, z_{in_1})\right)^{h_i}
\end{split}\end{equation*}
и 
\begin{equation*}\begin{split}
f_0 := \left(\tilde f(x_{11}^{(1)}, \ldots, x_{1\ell}^{(1)}; \ldots;
x^{(2k)}_{11}, \ldots,  x^{(2k)}_{1\ell};\ z_{11}, \ldots, z_{1n_1})\right)^{h_1}\cdot\\ \cdot \prod_{i=2}^m w_{i-1} \left(\tilde f\left(\left(x_{i1}^{(1)}\right)^{v^{i-1}}, \ldots, \left( x_{i\ell}^{(1)}\right)^{v^{i-1}}; \ldots;
\left(x^{(2k)}_{i1} \right)^{v^{i-1}}, \ldots,  \left( x^{(2k)}_{i\ell} \right)^{v^{i-1}};\ z_{i1}, \ldots, z_{in_1}\right)\right)^{h_i}.
\end{split}\end{equation*}

Тогда $\tilde f_0$
не обращается в нуль при следующей подстановке:

\begin{enumerate}
\item $x_{ij}^{(t)}=b_j$ при $1\leqslant i \leqslant m$, $1\leqslant j \leqslant \ell$, $1\leqslant t \leqslant 2k$;

\item $z_{ij}=\bar z_j$ при $1\leqslant i \leqslant m$, $1\leqslant j \leqslant n_1$;

\item $w_i = \bar w_i$ при $1\leqslant i \leqslant m-1$.
\end{enumerate}

Рассмотривая некоторые прообразы $u_i$ элементов $\bar w_i$
и $q_i$ элементов $\bar z_i$ в $A$,
получаем, что значение многочлена $f_0$ при подстановке

\begin{enumerate}
\item $x_{ij}^{(t)}=v^{m-i}a_j$ при $1\leqslant i \leqslant m$, $1\leqslant j \leqslant \ell$, $1\leqslant t \leqslant 2k$;
\item $z_{ij}=q_j$ при $1\leqslant i \leqslant m$, $1\leqslant j \leqslant n_1$;
\item $w_i = u_i$ при $1\leqslant i \leqslant m-1$
\end{enumerate}
не принадлежит идеалу $J^{\tilde H}(A)$.

Обозначим последнюю подстановку через $\Xi$, а значение $H$-многочлена $f_0$ при подстановке $\Xi$~"--- через $a$.

Пусть $f := \Alt_1 \ldots \Alt_{2k} f_0$, где $\Alt_t$~"---
оператор альтернирования по переменным $x^{(t)}_{ij}$, 
$1\leqslant i \leqslant m$, $1\leqslant j \leqslant \ell$.
Рассмотрим образ значения многочлена $f$ при подстановке $\Xi$
в $A/J^{\tilde H}(A)$. Если альтернирование заменяет $x^{(t)}_{ij}$
на $x^{(t)}_{i'j'}$, где $i < i'$, то значение
выражения $(x^{(t)}_{i'j'})^{v^{i-1}}$ при подстановке $\Xi$ равно $v^{m-1+(i-i')}a_{j'}
\in J^{\tilde H}(A)$, т.е. образ всего одночлена в 
$A/J^{\tilde H}(A)$ равен нулю. Следовательно, по модулю $J^{\tilde H}(A)$
можно считать, что альтернирования заменяют $x^{(t)}_{ij}$
на $x^{(t)}_{ij'}$ для тех же самых $i,t$. Учитывая, что
$\tilde f$ кососимметричен по $x_1^{(i)}, \ldots, x_\ell^{(i)}$
для любого $1\leqslant i  \leqslant 2k$,
получаем, что значение многочлена $f$ при подстановке $\Xi$
принадлежит классу $(\ell!)^{2km}a+J^{\tilde H}(A)$, т.е. не равно нулю.

$H$-многочлен $f$ удовлетворяет всем требованиям свойства (*)
для алгебры $H$. 
Следовательно, для всех конечномерных $H$-простых алгебр выполняется свойство (*). 
\end{proof}

 \section{Основная теорема и её следствия}
 \label{SectionGrowthHIdAssoc}
 
Теперь мы готовы к тому, чтобы сформулировать основную теорему данной главы и важнейшие её следствия.
Напомним, что через $J^H(A)$ обозначается максимальный нильпотентный $H$-инвариантный идеал алгебры $A$.

\begin{theorem}\label{TheoremHmodHRadAmitsurPIexpHBdimB}
Пусть $A$~"--- конечномерная ассоциативная алгебра с обобщённым $H$-действием
для некоторой ассоциативной алгебры $H$ с $1$ над полем $\mathbbm{k}$ характеристики $0$. Предположим, что $A/J^H(A) = B_1 \oplus B_2 \oplus \ldots
\oplus B_q$ (прямая сумма $H$-инвариантных идеалов) для некоторых $H$-простых алгебр $B_i$,
удовлетворяющих свойству~(*), причём для некоторого обычного разложения Веддербёрна~"--- Мальцева $$A/J^H(A)=B\oplus J\bigl(A/J^H(A)\bigr)\text{ (прямая сумма подпространств)}$$
существует такое вложение $\varkappa \colon A/J^H(A) \hookrightarrow A$, что
 $\pi\varkappa = \id_{A/J}$, где $\pi \colon A \twoheadrightarrow A/J$~"--- естественный сюръективный гомоморфизм,
 $$\varkappa(a)\varkappa(b)=\varkappa(ab)\text{\quad и\quad}\varkappa(b)\varkappa(a)=\varkappa(ba)\text{ \quad для всех }a\in A/J^H(A)\text{ и }b\in B.$$
Положим
\begin{equation}\begin{split}\label{EqdAssoc}d:= \max\dim\left( B_{i_1}\oplus B_{i_2} \oplus \ldots \oplus B_{i_r}
 \mathbin{\Bigl|}  r \geqslant 1,\right.\\
   \left. (H\varkappa(B_{i_1}))A^+ \,(H\varkappa(B_{i_2})) A^+ \ldots (H\varkappa(B_{i_{r-1}})) A^+\,(H\varkappa(B_{i_r}))\ne 0\right),
  \end{split}\end{equation} где
 $A^+:=A+\mathbbm{k}\cdot 1$.
 Тогда \begin{enumerate}
\item при $d=0$ существует такое $n_0$, что $c_n^H(A)=0$ при всех $n\geqslant n_0$;
\item при $d>0$ существуют такие $C_1,C_2 > 0$ и $r_1,r_2 \in \mathbb R$, что $$C_1 n^{r_1} d^n \leqslant c^H_n(A) \leqslant C_2 n^{r_2} d^n\text{ для всех }n\in\mathbb N.$$
\end{enumerate}
\end{theorem}

Теорема~\ref{TheoremHmodHRadAmitsurPIexpHBdimB} будет доказана в \S\ref{SectionAssocFinishingTheProof}.

\begin{remark}
Из теоремы~\ref{TheoremHmodHRadAmitsurPIexpHBdimB} следует, что существует  $$\PIexp^H(A):=\lim_{n\to\infty}\sqrt[n]{c_n^H(A)}= d \in\mathbb Z_+$$ и
для полиномиальных $H$-тождеств алгебры $A$ справедлив аналог гипотезы Амицура.
\end{remark}
\begin{remark}
В случае, когда основное поле $\mathbbm{k}$ алгебраически замкнуто, существование
вложения $\varkappa$ с требуемыми свойствами следует из теоремы~\ref{TheoremWeakWedderburnMalcevHRad}.
\end{remark}

Докажем теперь основные следствия из теоремы~\ref{TheoremHmodHRadAmitsurPIexpHBdimB}:
\begin{theorem}\label{TheoremHOreAmitsur}
Пусть $H$~"--- алгебра Хопфа над полем $\mathbbm{k}$ характеристики $0$,
причём $H$ либо сама является конечномерной полупростой алгеброй Хопфа,
либо получена при помощи (возможно, многократного)
расширения Оре конечномерной полупростой алгебры Хопфа косопримитивными элементами. (Например,
$H=H_{m^2}(\zeta)$.)
Тогда для любой конечномерной ассоциативной $H$-модульной алгебры $A$
\begin{enumerate}
\item либо существует такое $n_0$, что $c_n^H(A)=0$ при всех $n\geqslant n_0$;
\item либо существуют такие константы $C_1, C_2 > 0$, $r_1, r_2 \in \mathbb R$,
  $d \in \mathbb N$, что $$C_1 n^{r_1} d^n \leqslant c^{H}_n(A)
   \leqslant C_2 n^{r_2} d^n\text{ для всех }n \in \mathbb N.$$
\end{enumerate}
   В частности, для любой конечномерной ассоциативной $H$-модульной алгебры $A$
   существует $\PIexp^H(A)\in\mathbb Z_+$ и, таким образом,
   для $A$ справедлива гипотеза Амицура~"--- Бахтурина.
\end{theorem}
\begin{proof}
Заметим, что $H$-коразмерности
не меняются при расширении основного поля.
Доказательство полностью повторяет соответствующие рассуждения
в случае обычных коразмерностей~\cite[теорема~4.1.9]{ZaiGia}.
Поэтому без ограничения общности можно считать, что основное поле $\mathbbm{k}$
алгебраически замкнуто.
В силу примера~\ref{ExampleHmodHSemiSimplePropertyStar} и следствия~\ref{CorollaryOreExtStar} все конечномерные $H$-простые алгебры над $\mathbbm{k}$
удовлетворяют свойству (*).
Из теоремы~\ref{TheoremSkryabinVanOystaeyen} и леммы~\ref{LemmaHSemiSimpleIsUnital}
следует, что $A/J^H(A) = B_1 \oplus \ldots \oplus B_q$ (прямая сумма $H$-инвариантных идеалов)
для некоторых $H$-простых алгебр $B_i$. Теперь достаточно применить теоремы~\ref{TheoremWeakWedderburnMalcevHRad} и~\ref{TheoremHmodHRadAmitsurPIexpHBdimB}.
\end{proof}

Если потребовать, чтобы радикал Джекобсона алгебры $A$ был её $H$-подмодулем,
то справедливость гипотезы Амицура~"--- Бахтурина получается и в случае произвольной алгебры Хопфа $H$:

\begin{theorem}\label{TheoremHmoduleAssoc}
Пусть $A$~"--- конечномерная 
ассоциативная $H$-модульная алгебра для некоторой алгебры Хопфа $H$
над полем $\mathbbm{k}$ характеристики $0$, причём радикал Джекобсона $J(A)$ является $H$-подмодулем.
Тогда \begin{enumerate}
\item либо существует такое $n_0$, что $c_n^H(A)=0$ при всех $n\geqslant n_0$;
\item либо существуют такие константы $C_1, C_2 > 0$, $r_1, r_2 \in \mathbb R$,
  $d \in \mathbb N$, что $$C_1 n^{r_1} d^n \leqslant c^{H}_n(A)
   \leqslant C_2 n^{r_2} d^n\text{ для всех }n \in \mathbb N.$$
\end{enumerate}
В частности, для любой конечномерной ассоциативной $H$-модульной алгебры $A$
с $H$-инвариантным радикалом Джекобсона
   существует $\PIexp^H(A)\in\mathbb Z_+$ и, таким образом,
   для $A$ справедлива гипотеза Амицура~"--- Бахтурина.
\end{theorem}
\begin{proof}
Пусть $K \supset \mathbbm{k}$~"--- расширение поля $\mathbbm{k}$.
Поскольку алгебра $A/J(A)$ полупроста и $\chr \mathbbm{k}= 0$, $$(A \otimes_\mathbbm{k} K)/(J(A) \otimes_\mathbbm{k} K) \cong (A/J(A)) \otimes_\mathbbm{k} K$$ снова является полупростой алгеброй (см., например, \cite[\S 10.7, следствие b]{PierceAssoc}). В силу того, что идеал $J(A) \otimes_\mathbbm{k} K$ нильпотентен,
справедливо равенство $J(A \otimes_\mathbbm{k} K) = J(A) \otimes_\mathbbm{k} K$.
В частности, $J(A \otimes_\mathbbm{k} K)$ по-прежнему является $H\otimes_\mathbbm{k} K$-инвариантным идеалом.
Как было уже замечено в доказательстве предыдущей теоремы,
 $H$-коразмерности не меняются при расширении основного поля.
 Поэтому без ограничения общности можно считать поле $\mathbbm{k}$ алгебраически
 замкнутым.
 В силу теоремы~\ref{TheoremSkryabinVanOystaeyen}
существует разложение 
  $A/J(A) = B_1 \oplus \ldots \oplus B_q$ (прямая сумма $H$-инвариантных идеалов) для некоторых
   $H$-простых алгебр $B_i$. Поскольку алгебра $A/J(A)$ полупроста,
   алгебры $B_i$ также полупросты. В силу теоремы~\ref{TheoremAssocAlternateFinal}
   для всех алгебр $B_i$ выполняется свойство~(*).
   Теперь для завершения доказательства достаточно применить теоремы~\ref{TheoremWeakWedderburnMalcevHRad} и~\ref{TheoremHmodHRadAmitsurPIexpHBdimB}.
\end{proof}

В случае, когда существует $H$-инвариантный аналог разложения Веддербёрна~"--- Мальцева,
формула~\eqref{EqdAssoc} принимает особенно простой вид, напоминающий
формулу для обычной PI-экспоненты, полученную А.~Джамбруно и М.\,В.~Зайцевым~\cite[\S 6.2]{ZaiGia}:

\begin{theorem}\label{TheoremGenHAmitsurHWedederburn} 
Пусть $A$~"--- конечномерная ассоциативная алгебра с обобщённым $H$-действием
для некоторой ассоциативной алгебры $H$ с $1$ над полем $\mathbbm{k}$ характеристики $0$,
причём существует разложение 
$A=B_0\oplus J^H(A)$ (прямая сумма $H$-подмодулей), где алгебра
$B_0$ представляется в виде
$B_1 \oplus B_2 \oplus \ldots
\oplus B_q$ (прямая сумма $H$-инвариантных идеалов) для некоторых $H$-простых алгебр $B_i$,
удовлетворяющих свойству (*).
Положим 
\begin{equation}\begin{split}\label{EqPIexpH(A)InvWedMalcev}d: = 
 \max\dim\left( B_{i_1}\oplus B_{i_2} \oplus \ldots \oplus B_{i_r}
 \mathbin{\Bigl|}  \right. \\ \left.
    B_{i_1}J^H(A)B_{i_2} J^H(A) \ldots B_{i_{r-1}} J^H(A) B_{i_r} \ne 0,\ r \geqslant 1\right).
  \end{split}\end{equation}
 Тогда \begin{enumerate}
\item при $d=0$ существует такое $n_0$, что $c_n^H(A)=0$ при всех $n\geqslant n_0$;
\item при $d>0$ существуют такие $C_1,C_2 > 0$ и $r_1,r_2 \in \mathbb R$, что $$C_1 n^{r_1} d^n \leqslant c^H_n(A) \leqslant C_2 n^{r_2} d^n\text{ для всех }n\in\mathbb N.$$
\end{enumerate}
\end{theorem}
\begin{proof} Достаточно положить $\varkappa = \left(\pi\bigr|_{B_0} \right)^{-1}$
и заметить, что в этом случае $\varkappa$ является гомоморфизмом алгебр и $H$-модулей,
т.е. $(H B_k) A^+ (H B_\ell)\ne 0$ для некоторых $1\leqslant k,\ell \leqslant q$,
только если либо $k = \ell$, либо $B_k J^H(A) B_\ell \ne 0$. 
Отсюда числа $d$ в теоремах~\ref{TheoremHmodHRadAmitsurPIexpHBdimB} и~\ref{TheoremGenHAmitsurHWedederburn} совпадают, и достаточно применить теорему~\ref{TheoremHmodHRadAmitsurPIexpHBdimB}.
\end{proof}

Теперь выведем из теоремы~\ref{TheoremHmoduleAssoc} справедливость аналога гипотезы Амицура для
дифференциальных тождеств.

\begin{theorem}\label{TheoremDiffAssoc}
Пусть $A$~"--- конечномерная
ассоциативная алгебра над полем $\mathbbm{k}$ характеристики $0$ с действием алгебры Ли $\mathfrak g$ дифференцированиями.
Тогда \begin{enumerate}
\item либо существует такое $n_0$, что $c_n^{U(\mathfrak g)}(A)=0$ при всех $n\geqslant n_0$;
\item либо существуют такие константы $C_1, C_2 > 0$, $r_1, r_2 \in \mathbb R$,
  $d \in \mathbb N$, что $$C_1 n^{r_1} d^n \leqslant c^{U(\mathfrak g)}_n(A)
   \leqslant C_2 n^{r_2} d^n\text{ для всех }n \in \mathbb N.$$
\end{enumerate}
В частности, для любой такой алгебры $A$
   существует $\PIexp^{U(\mathfrak g)}(A)\in\mathbb Z_+$ и, таким образом,
   справедлив аналог гипотезы Амицура.
\end{theorem}
\begin{proof} Как было показано в примере~\ref{ExampleUgModule},
алгебра $A$ является $U(\mathfrak g)$-модульной алгеброй.
В силу теоремы~\ref{TheoremRadicalHSubMod}
радикал Джекобсона $J(A)$
является $U(\mathfrak g)$-подмодулем алгебры $A$.
Теперь неравенство для коразмерностей следует из теоремы~\ref{TheoremHmoduleAssoc}.
\end{proof}

\begin{remark}
В случае, когда $\mathfrak g$~"--- конечномерная полупростая алгебра Ли, справедливо равенство $\PIexp^{U(\mathfrak g)}(A)=\PIexp(A)$ (см. теорему~\ref{TheoremAssDerPIexpEqual} ниже).
 \end{remark}

\begin{remark}
Из теоремы~\ref{TheoremDiffAssoc} следует справедливость аналога гипотезы Амицура для
дифференциальных коразмерностей даже для тех алгебр, для которых не существует
 инвариантного разложения Веддербёрна~"--- Мальцева, см. пример~\ref{ExampleU(L)noninvWedderburnMalcev}.
 \end{remark}

Докажем теперь справедливость аналогов гипотезы Амицура для
градуированных тождеств, $G$-тождеств и градуированных $G$-тождеств:

\begin{theorem}\label{TheoremMainGrAssoc}
Пусть $A$~"--- конечномерная 
ассоциативная алгебра над полем $\mathbbm{k}$ характеристики $0$, градуированная произвольной
группой $G$. 
Тогда \begin{enumerate}
\item либо существует такое $n_0$, что $c_n^{G\text{-}\mathrm{gr}}(A)=0$ при всех $n\geqslant n_0$;
\item либо существуют такие константы $C_1, C_2 > 0$, $r_1, r_2 \in \mathbb R$,
  $d \in \mathbb N$, что $$C_1 n^{r_1} d^n \leqslant c^{G\text{-}\mathrm{gr}}_n(A)
   \leqslant C_2 n^{r_2} d^n\text{ для всех }n \in \mathbb N.$$
\end{enumerate}
В частности, для любой такой алгебры $A$
   существует $\PIexp^{G\text{-}\mathrm{gr}}(A)\in\mathbb Z_+$ и, таким образом,
   справедлив аналог гипотезы Амицура.
\end{theorem}

\begin{theorem}\label{TheoremGAssoc}
Пусть $A$~"--- конечномерная 
ассоциативная алгебра над полем $\mathbbm{k}$ характеристики $0$ с действием произвольной
группы $G$ автоморфизмами и антиавтоморфизмами. 
Тогда \begin{enumerate}
\item либо существует такое $n_0$, что $c_n^{G}(A)=0$ при всех $n\geqslant n_0$;
\item либо существуют такие константы $C_1, C_2 > 0$, $r_1, r_2 \in \mathbb R$,
  $d \in \mathbb N$, что $$C_1 n^{r_1} d^n \leqslant c^{G}_n(A)
   \leqslant C_2 n^{r_2} d^n\text{ для всех }n \in \mathbb N.$$
\end{enumerate}
В частности, для любой такой алгебры $A$
   существует $\PIexp^{G}(A)\in\mathbb Z_+$ и, таким образом,
   справедлив аналог гипотезы Амицура.
\end{theorem}

\begin{remark}
В случае, когда $G$~"--- связная редуктивная аффинная алгебраическая группа, которая рационально действует на $G$, справедливо равенство $\PIexp^{G}(A)=\PIexp(A)$ (см. теорему~\ref{TheoremAssGPIexpEqual} ниже).
 \end{remark}

Теоремы~\ref{TheoremMainGrAssoc} и~\ref{TheoremGAssoc}
являются частными случаями (соответственно, при $G=\lbrace e \rbrace$
и при $T=\lbrace e \rbrace$) следующей теоремы, где мы по аналогии с замечанием~\ref{RemarkGCodim} полагаем $c_n^{T\text{-}\mathrm{gr}, G}(A):= c_n^{T\text{-}\mathrm{gr}, \mathbbm{k}G}(A)$:

\begin{theorem}\label{TheoremTGrGActionAssoc}
Пусть $A$~"--- конечномерная
ассоциативная алгебра над полем $\mathbbm{k}$ характеристики $0$, градуированная произвольной
группой $T$, причём на $A$ задано градуированное действие некоторой группы $G$ (см. определение~\ref{DefGradedAction}).
Тогда \begin{enumerate}
\item либо существует такое $n_0$, что $c_n^{T\text{-}\mathrm{gr}, G}(A)=0$ при всех $n\geqslant n_0$;
\item либо существуют такие константы $C_1, C_2 > 0$, $r_1, r_2 \in \mathbb R$,
  $d \in \mathbb N$, что $$C_1 n^{r_1} d^n \leqslant c^{T\text{-}\mathrm{gr}, G}_n(A)
   \leqslant C_2 n^{r_2} d^n\text{ для всех }n \in \mathbb N.$$
\end{enumerate}
В частности, для любой такой алгебры $A$
   существует $$\PIexp^{T\text{-}\mathrm{gr},G}(A):=\lim\limits_{n\to\infty} \sqrt[n]{c^{T\text{-}\mathrm{gr},G}_n(A)}\in\mathbb Z_+$$ и, таким образом,
   справедлив аналог гипотезы Амицура.
\end{theorem}
\begin{proof}
Как было доказано в теореме~\ref{TheoremGradGenActionReplace},
алгебра $A$ является алгеброй с обобщённым $\mathbbm{k}^T\otimes \mathbbm{k}G$-действием,
причём в силу предложения~\ref{PropositionCnTGrHCnFTotimesH}
для всех $n\in\mathbb N$ справедливы 
равенства $c_n^{T\text{-}\mathrm{gr}, G}(A)=c_n^{\mathbbm{k}^T\otimes \mathbbm{k}G}(A)$.
Как и выше, коразмерности $c_n^{\mathbbm{k}^T\otimes \mathbbm{k}G}(A)$ не меняются
при расширении основного поля $\mathbbm{k}$, поэтому без ограничения общности
можно считать поле $\mathbbm{k}$ алгебраически замкнутым.
В силу теоремы~\ref{TheoremTGradedGActionInvRad}
радикал $J(A)$ является $\mathbbm{k}^T\otimes \mathbbm{k}G$-подмодулем,
а в силу теоремы~\ref{TheoremTGradedGActionInvWedderburnArtin}
алгебра $A/J(A)$ является прямой суммой $\mathbbm{k}^T\otimes \mathbbm{k}G$-инвариантных
идеалов $B_i$, каждый из которых является $\mathbbm{k}^T\otimes \mathbbm{k}G$-простой алгеброй.
Поскольку алгебры $B_i$ полупросты, в силу теоремы~\ref{TheoremAssocAlternateFinal}
для каждой из них выполняется свойство (*).
Теперь достаточно применить теоремы~\ref{TheoremWeakWedderburnMalcevHRad} и~\ref{TheoremHmodHRadAmitsurPIexpHBdimB}.
\end{proof}
\begin{remark}
Из теоремы~\ref{TheoremGAssoc} следует справедливость аналога гипотезы Амицура для
$G$-коразмерностей даже для тех алгебр, для которых не существует
 $G$-инвариантного разложения Веддербёрна~"--- Мальцева, см. пример~\ref{ExampleGnoninvWedderburnMalcev}.
 \end{remark}
 
\section{Оценки сверху и снизу}\label{SectionAssocUpperLower}

В данном параграфе доказываются оценки сверху и снизу из теоремы~\ref{TheoremHmodHRadAmitsurPIexpHBdimB}.
Оказывается, что оценку сверху можно получить при более слабых предположениях, чем предположения
теоремы~\ref{TheoremHmodHRadAmitsurPIexpHBdimB}. Оценка снизу доказывается в предположении, что
существует $H$-многочлен, не являющийся полиномиальным $H$-тождеством, который кососимметричен
по достаточному числу наборов из $d$ переменных.

Пусть $A$~"--- конечномерная ненильпотентная ассоциативная алгебра с обобщённым $H$-действием
для некоторой ассоциативной алгебры $H$ с $1$ над полем $\mathbbm{k}$ характеристики $0$,
а $J$~"--- некоторый такой $H$-инвариантный идеал алгебры $A$, что $J^p=0$ для некоторого $p\in\mathbb N$. Фиксируем разложение $A/J = B_1 \oplus B_2 \oplus \ldots
\oplus B_q$, где $B_i$~"--- некоторые подпространства.
Пусть $\varkappa \colon A/J \hookrightarrow A$~"--- некоторое такое $\mathbbm{k}$-линейное
отображение, что $\pi\varkappa = \id_{A/J}$, где $\pi \colon A \twoheadrightarrow A/J$~"--- естественный сюръективный гомоморфизм. Зададим теперь число $d$ формулой~(\ref{EqdAssoc}).

\begin{lemma}\label{LemmaAssocUpperCochar}
 Пусть $n\in\mathbb N$, а $\lambda = (\lambda_1, \ldots, \lambda_s) \vdash n$, причём $\sum_{k=d+1}^s \lambda_k \geqslant p$. Тогда $m(A, H, \lambda)=0$. 
\end{lemma}
\begin{proof}  В силу теоремы~\ref{ThKratnost}
  достаточно доказать, что $e^{*}_{T_\lambda}f \in \Id^H(A)$ для всех $f \in P^H_n$ и для всех таблиц Юнга $T_\lambda$, отвечающих $\lambda$.

Выберем в алгебре $A$ базис, который является объединением базисов пространств~$\varkappa(B_1),\ldots, \varkappa(B_q)$ и~$J$. Поскольку $H$-многочлен $e^{*}_{T_\lambda}f$ полилинеен, 
достаточно показать, что $e^{*}_{T_\lambda}f$ обращается в нуль при подстановке базисных элементов вместо своих неизвестных. Фиксируем некоторую подстановку базисных элементов
и выберем такие $1 \leqslant i_1,\ldots,i_r \leqslant q$,
что все подставляемые элементы принадлежат подпространству $\varkappa(B_{i_1})\oplus \ldots \oplus \varkappa(B_{i_r}) \oplus J$ и для любого $k$ подставляется хотя бы один элемент из $\varkappa(B_{i_k})$.
Тогда можно считать, что $\dim(B_{i_1}\oplus \ldots \oplus B_{i_r}) \leqslant d$,
поскольку в противном случае значение $H$-многочлена $e^{*}_{T_\lambda}f$
равно нулю в силу определения числа $d$.
 Напомним, что
$e^{*}_{T_\lambda} = b_{T_\lambda} a_{T_\lambda}$, а оператор $b_{T_\lambda}$ делает многочлены кососимметричными
по переменным, отвечающим каждому столбцу таблицы Юнга $T_\lambda$. Отсюда для того, чтобы $H$-многочлен $e^{*}_{T_\lambda} f$ не обратился в нуль, 
в переменные каждого столбца должны поставляться различные элементы базиса.
Следовательно, в переменные должно подставляться по крайней мере $\sum_{k=d+1}^s \lambda_k \geqslant p$ элементов из идеала $J$.
В силу того, что $J$ является $H$-подмодулем и $J^p = 0$, значение $H$-многочлена $e^{*}_{T_\lambda}f$
равно нулю, откуда $e^{*}_{T_\lambda} f \in \Id^H(A)$ и $m(A, H, \lambda)=0$.
\end{proof}

\begin{theorem}\label{TheoremAssocUpper} 
Если $d > 0$, то существуют такие константы $C_2 > 0$ и $r_2 \in \mathbb R$,
что $c^H_n(A) \leqslant C_2 n^{r_2} d^n$
для всех $n \in \mathbb N$. В случае $d=0$ алгебра $A$ совпадает с нильпотентным идеалом $J$
и $c^H_n(A)=0$ при $n\geqslant p$.
\end{theorem}
\begin{proof} Рассмотрим такое $\lambda \vdash n$,
что $m(A,H,\lambda) \ne 0$.
Согласно формуле крюков $$\dim M(\lambda) = \frac{n!}{\prod_{i,j} h_{ij}},$$
где $h_{ij}$~"--- длина крюка в диаграмме $D_\lambda$ с вершиной в $(i, j)$.
Применив лемму~\ref{LemmaAssocUpperCochar} и
оценивая полиномиальный коэффициент $\frac{(\lambda_1 +\ldots +\lambda_d)!}{(\lambda_1)!\ldots (\lambda_d)!}$ сверху по формуле возведения суммы $\underbrace{1+\ldots+1}_d$ в степень
$\lambda_1 +\ldots +\lambda_d$,
получаем, что
\begin{equation*}\begin{split}\dim M(\lambda) \leqslant 
\frac{n!}{(\lambda_1)!\ldots (\lambda_d)!}
\leqslant \frac{n!}{(\lambda_1 +\ldots +\lambda_d)!}
\frac{(\lambda_1 +\ldots +\lambda_d)!}{(\lambda_1)!\ldots (\lambda_d)!}
\leqslant \frac{n!}{(n-p)!}
d^{\lambda_1 +\ldots +\lambda_d} \leqslant C_3 n^{r_3} d^n.
\end{split}\end{equation*}
для некоторых констант $C_3, r_3 > 0$.
Теперь оценка сверху следует из теоремы~\ref{TheoremUpperBoundColengthHNAssoc}.
\end{proof}

\begin{lemma}\label{LemmaAssocCochar} Предположим, что существует такое число $n_0 \in \mathbb N$, что для любого $n\geqslant n_0$
существуют попарно непересекающиеся наборы переменных $X_1$, \ldots, $X_{2k} \subseteq \lbrace x_1, \ldots, x_n\rbrace$, где $k := \left[\frac{n-n_0}{2d}\right]$,
$|X_1| = \ldots = |X_{2k}|=d$, и $H$-многочлен $f \in P^H_n \backslash
\Id^H(A)$, кососимметричный по переменным каждого множества $X_j$. Тогда для любого $n \geqslant n_0$
существует такое разбиение $\lambda = (\lambda_1, \ldots, \lambda_s) \vdash n$,
где $\lambda_i > 2k-p$ для всех $1 \leqslant i \leqslant d$,
что $m(A, H, \lambda) \ne 0$.
\end{lemma}
\begin{proof}
Достаточно доказать, что $e^*_{T_\lambda} f \notin \Id^H(A)$
для некоторой таблицы Юнга $T_\lambda$ требуемой формы $\lambda$.
Известно (см., например, теорему 3.2.7 из~\cite{Bahturin}), что $$\mathbbm{k}S_n = \bigoplus_{\lambda,T_\lambda} \mathbbm{k}S_n e^{*}_{T_\lambda},$$ где суммирование ведётся по множеству стандартных таблиц Юнга $T_\lambda$
всевозможных форм $\lambda \vdash n$. Отсюда $$\mathbbm{k}S_n f = \sum_{\lambda,T_\lambda} \mathbbm{k}S_n e^{*}_{T_\lambda}f
\not\subseteq \Id^H(A)$$ и $e^{*}_{T_\lambda} f \notin \Id^H(A)$ для некоторого $\lambda \vdash n$.
Докажем, что разбиение $\lambda$ имеет требуемый вид.
Достаточно доказать, что
$\lambda_d > 2k-p$, так как
$\lambda_i \geqslant \lambda_d$ для всех $1 \leqslant i \leqslant d$.
В любой строчке таблицы $T_\lambda$ содержится не более одного 
номера переменной из одного и того же множества $X_i$,
поскольку $e^{*}_{T_\lambda} = b_{T_\lambda} a_{T_\lambda}$,
а $a_{T_\lambda}$ симметризует по переменным, отвечающим каждой строчке таблицы $T_\lambda$.
Отсюда $$\sum_{i=1}^{d-1} \lambda_i \leqslant 2k(d-1) + (n-2kd) = n-2k.$$
В силу леммы~\ref{LemmaAssocUpperCochar} справедливо неравенство
$\sum_{i=1}^d \lambda_i > n-p$. Следовательно,
$\lambda_d > 2k-p$.
\end{proof}

Теперь дополним оценку сверху оценкой снизу:
\begin{theorem}\label{TheoremAssocBounds}
Пусть $A$~"--- конечномерная ассоциативная алгебра с обобщённым $H$-действием
для некоторой ассоциативной алгебры $H$ с $1$ над полем $\mathbbm{k}$ характеристики $0$,
а $J$~"--- некоторый такой $H$-инвариантный идеал алгебры $A$, что $J^p=0$ для некоторого $p\in\mathbb N$. Предположим, что задано разложение $A/J = B_1 \oplus B_2 \oplus \ldots
\oplus B_q$, где $B_i$~"--- некоторые подпространства.
Пусть $\varkappa \colon A/J \hookrightarrow A$~"--- некоторое такое $\mathbbm{k}$-линейное
отображение, что $\pi\varkappa = \id_{A/J}$, где $\pi \colon A \twoheadrightarrow A/J$~"--- естественный сюръективный гомоморфизм.
Предположим, что число $d$, заданное формулой~(\ref{EqdAssoc}), больше $0$ и существует такое число $n_0 \in \mathbb N$, что для любого $n\geqslant n_0$
существуют попарно непересекающиеся наборы переменных $X_1$, \ldots, $X_{2k} \subseteq \lbrace x_1, \ldots, x_n\rbrace$, где $k := \left[\frac{n-n_0}{2d}\right]$,
$|X_1| = \ldots = |X_{2k}|=d$, и $H$-многочлен $f \in P^H_n \backslash
\Id^H(A)$, кососимметричный по переменным каждого множества $X_j$.
 Тогда существуют такие $C_1,C_2 > 0$ и $r_1,r_2 \in \mathbb R$, что $$C_1 n^{r_1} d^n \leqslant c^H_n(A) \leqslant C_2 n^{r_2} d^n\text{ для всех }n\in\mathbb N.$$
\end{theorem}
\begin{proof}
Диаграмма Юнга~$D_\lambda$ из леммы~\ref{LemmaAssocCochar} содержит квадратную
поддиаграмму~$D_\mu$, где $\mu=(\underbrace{2k-p, \ldots, 2k-p}_d)$.
Из правила ветвления для группы $S_n$ следует, что
если рассмотреть сужение $S_n$-действия на $M(\lambda)$ до $S_{n-1}$-действия,
то $\mathbbm{k}S_n$-модуль $M(\lambda)$
оказывается прямой суммой всех неизоморфных
$\mathbbm{k}S_{n-1}$-модулей $M(\nu)$, где $\nu \vdash (n-1)$ и всякая таблица $D_\nu$
получена из $D_\lambda$ удалением одной клетки. В частности,
$\dim M(\nu) \leqslant \dim M(\lambda)$.
Применяя правило ветвления $(n-d(2k-p))$ раз, получаем, что $\dim M(\mu) \leqslant \dim M(\lambda)$.
В силу формулы крюков $$\dim M(\mu) = \frac{(d(2k-p))!}{\prod_{i,j} h_{ij}},$$
где $h_{ij}$~"--- длина крюка с вершиной в $(i, j)$.
По формуле Стирлинга
\begin{equation*}\begin{split}c_n^H(A)\geqslant \dim M(\lambda) \geqslant \dim M(\mu) \geqslant \frac{(d(2k-p))!}{((2k-p+d)!)^d}
\sim \\ \sim \frac{
\sqrt{2\pi d(2k-p)} \left(\frac{d(2k-p)}{e}\right)^{d(2k-p)}
}
{
\left(\sqrt{2\pi (2k-p+d)}
\left(\frac{2k-p+d}{e}\right)^{2k-p+d}\right)^d
} \sim C_4 k^{r_4} d^{2kd}\end{split}\end{equation*}
для некоторых констант $C_4 > 0$, $r_4 \in \mathbb Q$
при $k \to \infty$.
Поскольку $k = \left[\frac{n-n_0}{2d}\right]$,
это доказывает оценку снизу.

Оценка сверху была доказана в теореме~\ref{TheoremAssocUpper}.
\end{proof}

\section{Завершение доказательства}\label{SectionAssocFinishingTheProof}

Для того, чтобы завершить доказательство теоремы~\ref{TheoremHmodHRadAmitsurPIexpHBdimB}, осталось построить для алгебры $A$ полилинейный $H$-многочлен, который кососимметричен
по достаточно большому числу наборов из $d$ переменных и при этом не является $H$-тождеством.

\begin{lemma}\label{LemmaAssocLowerPolynomial}
Пусть $A$, $\varkappa$, $B_i$ и $d$ те же, что и в теореме~\ref{TheoremHmodHRadAmitsurPIexpHBdimB}.
Тогда если $d > 0$, то существует такое число $n_0 \in \mathbb N$, что для любого $n\geqslant n_0$
существуют попарно непересекающиеся наборы переменных $X_1$, \ldots, $X_{2k} \subseteq \lbrace x_1, \ldots, x_n\rbrace$, где $k := \left[\frac{n-n_0}{2d}\right]$,
$|X_1| = \ldots = |X_{2k}|=d$, и $H$-многочлен $f \in P^H_n \backslash
\Id^H(A)$, кососимметричный по переменным каждого множества $X_j$.
\end{lemma}
\begin{proof} Пусть $J := J^H(A)$.
Без ограничения общности можно считать, что $$d = \dim(B_1 \oplus B_2 \oplus \ldots \oplus B_r),$$
где
$(H\varkappa(B_1))A^+ (H\varkappa(B_2))A^+ \ldots (H\varkappa(B_{r-1}))A^+ (H\varkappa(B_r))\ne 0$. 

Поскольку идеал $J$ нильпотентен, можно найти такое максимальное число $\sum_{i=1}^r q_i$, где $q_i \in \mathbb Z_+$,
 что
$$\left(a_1 \prod_{i=1}^{q_1} j_{1i}\right) ({\gamma_1}\varkappa(b_1))
\left(a_2 \prod_{i=1}^{q_2} j_{2i}\right)
 ({\gamma_2}\varkappa(b_2)) \ldots \left(a_r \prod_{i=1}^{q_r} j_{ri}\right)
 ({\gamma_r}\varkappa(b_r))  \left(a_{r+1} \prod_{i=1}^{q_{r+1}} j_{r+1,i}\right) \ne 0$$ для
 некоторых $j_{ki}\in J$, $a_k \in A^+$, $b_k \in B_i$, $\gamma_k \in H$.
 Введём обозначение $j_k := a_k\prod_{i=1}^{q_k} j_{ki}$.
 
 Тогда \begin{equation}\label{EqAssocNonZero}j_1 ({\gamma_1}\varkappa(b_1))
j_2  ({\gamma_2}\varkappa(b_2)) \ldots j_r  ({\gamma_r}\varkappa(b_r))j_{r+1} \ne 0\end{equation}
 для некоторых $b_i \in B_i$, $\gamma_i \in H$,
 однако
 \begin{equation}\label{EqAssocbazero}j_1
  \tilde b_1 j_2 \tilde b_2
  \ldots j_r \tilde b_r j_{r+1} = 0 \end{equation} для всех таких $\tilde b_i\in A^+(H\varkappa(B_i))A^+$,
  что  $\tilde b_k\in J(H\varkappa(B_k))A^+ + A^+(H\varkappa(B_k))J$ хотя бы для одного $k$.

Пусть $a^{(i)}_{k}$, где $1 \leqslant k \leqslant d_i := \dim B_i$,
 "---  базисные элементы алгебры $B_{i}$ для $1 \leqslant i \leqslant r$.
 
 В силу свойства~(*) существуют такие константы
$\tilde m_i \in \mathbb Z_+$, что для любого
$k$ существуют полилинейные многочлены $$f_i=f_i(x^{(i, 1)}_1,
 \ldots, x^{(i, 1)}_{d_i};
 \ldots;  x^{(i, 2k)}_1,
 \ldots, x^{(i, 2k)}_{d_i}; z^{(i)}_1, \ldots, z^{(i)}_{m_i}) \in P^H_{2k d_i+m_i} \backslash \Id^H(B_i),$$
   где $0 \leqslant m_i \leqslant \tilde m_i$,
   кососимметричные по переменным каждого из попарно непересекающихся множеств
$X^{(i)}_{\ell}=\lbrace x^{(i, \ell)}_1, x^{(i, \ell)}_2,
\ldots, x^{(i, \ell)}_{d_i} \rbrace$, $1 \leqslant \ell \leqslant 2k$.
В частности, существуют такие $\bar z^{(i)}_\alpha \in B_i$, где $1 \leqslant \alpha \leqslant m_i$,
что
$$\hat b_i := f_i(a^{(i)}_1,
 \ldots, a^{(i)}_{d_i};
 \ldots;  a^{(i)}_1,
 \ldots, a^{(i)}_{d_i}; \bar z^{(i)}_1, \ldots, \bar z^{(i)}_{m_i})\ne 0.$$
 
 Введём обозначения $n_0 := 3r-1+\sum_{i=1}^r \tilde m_i$, $k := \left[\frac{n-n_0}{2d}\right]$, $\tilde k := \left[\frac{n-2kd}{2d_1}\right]+1$. 
 Выберем многочлены $f_i$, где $1 \leqslant i \leqslant r$, для $k = \left[\frac{n-n_0}{2d}\right]$.
 Кроме этого, снова пользуясь свойством (*), выберем $\tilde f_1=\tilde f_1(x^{(1)}_1,
 \ldots, x^{(1)}_{d_i};
 \ldots;  x^{(2\tilde k)}_1,
 \ldots, x^{(2\tilde k)}_{d_i}; z_1, \ldots, z_{\hat m_1}) \in P^H_{2\tilde k d_1+\hat m_1} \backslash \Id^H(B_1),$ где $0\leqslant \hat m_1 \leqslant \tilde m_1$
 и 
$$\hat b := \tilde f_1(a^{(1)}_1,
 \ldots, a^{(1)}_{d_1};
 \ldots;  a^{(1)}_1,
 \ldots, a^{(1)}_{d_1}; \bar z_1, \ldots, \bar z_{\hat m_1})\ne 0$$
 для некоторых $\bar z_1, \ldots, \bar z_{\hat m_1} \in B_1$.

 Поскольку алгебры $B_i$ являются $H$-простыми,
существуют такие элементы
 $h_{i\ell} \in H$, $b_{i\ell}, \tilde b_{i\ell} \in B_i$, $\tilde b_\ell \in B_1$,
 что $\sum_\ell b_{i\ell} (h_{i\ell}\hat b_i) \tilde b_{i\ell} = b_i$
 для всех $2\leqslant i \leqslant r$ и $\sum_\ell \tilde b_\ell (h_{0\ell}\hat b)  b_{1\ell} (h_{1\ell}\hat b_1) \tilde b_{1\ell} = b_1$.
 
 Тогда выражение \begin{equation*}\begin{split}
 j_1 \Biggl(\gamma_1\varkappa \biggl(\sum_{s_1} \tilde b_{s_1} \left( h_{0 s_1}
\tilde f_1\bigl( a^{(1)}_1,
 \ldots, a^{(1)}_{d_1};
 \ldots;  a^{(1)}_1,
 \ldots, a^{(1)}_{d_1}; \bar z_1, \ldots, \bar z_{\hat m_1}\bigr)\right)  b_{1 s_1}  \cdot \\ \cdot
   \left( h_{1 s_1} f_1\bigl(a^{(1)}_1,
 \ldots, a^{(1)}_{d_1};
 \ldots;  a^{(1)}_1,
 \ldots, a^{(1)}_{d_1}; \bar z^{(1)}_1, \ldots, \bar z^{(1)}_{m_1} \bigr)\right) \tilde b_{1 s_1} \biggr) 
 \Biggr)j_2
 \cdot \\ \cdot
    \prod_{i=2}^{r}\Biggl(\gamma_i\varkappa\biggl(\sum_{s_i} b_{is_i} \left(h_{i s_i}  f_i\bigl(a^{(i)}_1,
 \ldots, a^{(i)}_{d_i};
 \ldots;  a^{(i)}_1,
 \ldots, a^{(i)}_{d_i}; \bar z^{(i)}_1, \ldots, \bar z^{(i)}_{m_i}\bigr)\right)\tilde b_{is_i} \biggr) \Biggr)j_{i+1}
 \end{split}
 \end{equation*} 
  равно левой части неравенства~(\ref{EqAssocNonZero}), которая, как следует из~(\ref{EqAssocNonZero}), не равна нулю.
  Отсюда можно выбрать такие индексы $s_1, \ldots, s_r$, что
\begin{equation*}\begin{split}
 a:= j_1 \Biggl(\gamma_1\varkappa \biggl(\tilde b_{s_1} \left( h_{0 s_1}
\tilde f_1\bigl( a^{(1)}_1,
 \ldots, a^{(1)}_{d_1};
 \ldots;  a^{(1)}_1,
 \ldots, a^{(1)}_{d_1}; \bar z_1, \ldots, \bar z_{\hat m_1}\bigr)\right) b_{1 s_1}  \cdot \\ \cdot
   \left( h_{1 s_1} f_1\bigl(a^{(1)}_1,
 \ldots, a^{(1)}_{d_1};
 \ldots;  a^{(1)}_1,
 \ldots, a^{(1)}_{d_1}; \bar z^{(1)}_1, \ldots, \bar z^{(1)}_{m_1} \bigr)\right) \tilde b_{1 s_1} \biggr)
 \Biggr) 
 j_2 \cdot \\ \cdot
   \prod_{i=2}^{r}\Biggl(\gamma_i\varkappa\biggl(b_{is_i} \left(h_{i s_i}  f_i\bigl(a^{(i)}_1,
 \ldots, a^{(i)}_{d_i};
 \ldots;  a^{(i)}_1,
 \ldots, a^{(i)}_{d_i}; \bar z^{(i)}_1, \ldots, \bar z^{(i)}_{m_i}\bigr)\right)\tilde b_{is_i} \biggr) \Biggr)j_{i+1}\ne 0.
 \end{split}
 \end{equation*}   

Пусть $B$~"--- максимальная полупростая подалгебра алгебры $A/J^H(A)$,
для которой $\varkappa$ удовлетворяет свойствам $\varkappa(a)\varkappa(b)=\varkappa(ab)$
и $\varkappa(b)\varkappa(a)=\varkappa(ba)$ для всех $a\in A/J^H(A)$
и $b\in B$. Поскольку в силу свойства~(*)
все алгебры $B_i$ обладает единицей, алгебра $A/J^H(A)=B_1 \oplus B_2 \oplus \ldots
\oplus B_q$ (прямая сумма $H$-инвариантных идеалов) также обладает единицей, причём $1_B = 1_{A/J^H(A)}$.
Введём обозначение $\tilde B_i := 1_{B_i} B$. 
Поскольку $\tilde B_i$ являются гомоморфными образами полупростой алгебры $B$, они также полупросты.
Теперь из $B \subseteq \tilde B_1 \oplus \tilde B_2 \oplus \ldots
\oplus \tilde B_q$ и максимальности подалгебры $B$ следует, что
 $B= \tilde B_1 \oplus \tilde B_2 \oplus \ldots
\oplus \tilde B_q$ (прямая сумма идеалов), а $\tilde B_i$~"--- максимальные
полупростые подалгебры алгебр $B_i$. Отсюда $1_{\tilde B_i} = 1_{B_i}$ и $1_{B_i} \in B$.
 
 Напомним, что, в частности, $\varkappa(b)=\varkappa(b)\varkappa(1_{B_i})$
для всех $b \in B_i$. Отсюда
\begin{equation*}\begin{split}
 a= j_1 \Biggl(\gamma_1\biggl(\varkappa \Bigl(\tilde b_{s_1} \bigl( h_{0 s_1}
\tilde f_1 ( a^{(1)}_1,
 \ldots, a^{(1)}_{d_1};
 \ldots;  a^{(1)}_1,
 \ldots, a^{(1)}_{d_1}; \bar z_1, \ldots, \bar z_{\hat m_1} ) \bigr) b_{1 s_1}   \cdot \\ \cdot
   \bigl( h_{1 s_1} f_1 (a^{(1)}_1,
 \ldots, a^{(1)}_{d_1};
 \ldots;  a^{(1)}_1,
 \ldots, a^{(1)}_{d_1}; \bar z^{(1)}_1, \ldots, \bar z^{(1)}_{m_1} ) \bigr) \tilde b_{1 s_1} \Bigr)
 \varkappa(1_{B_1}) \biggr) \Biggr) j_2
 \cdot \\ \cdot
    \prod_{i=2}^{r}\Biggl(\gamma_i\biggl(\varkappa\Bigl(b_{is_i} \bigl(h_{i s_i}  f_i (a^{(i)}_1,
 \ldots, a^{(i)}_{d_i};
 \ldots;  a^{(i)}_1,
 \ldots, a^{(i)}_{d_i}; \bar z^{(i)}_1, \ldots, \bar z^{(i)}_{m_i} ) \bigr) \tilde b_{is_i} \Bigr)
  \varkappa(1_{B_i})\biggr)\Biggr) j_{i+1}\ne 0.
 \end{split}
 \end{equation*}   
Более того, из $\pi(h\varkappa(a)-\varkappa(ha))=0$ и $\pi(\varkappa(a)\varkappa(b)-\varkappa(ab))=0$ 
следует, что $h\varkappa(a)-\varkappa(ha) \in J$
и $\varkappa(a)\varkappa(b)-\varkappa(ab) \in J$ для всех $a,b\in A$ и $h\in H$.
Отсюда в силу~(\ref{EqAssocbazero}) в множителях слева от $\varkappa(1_{B_i})$
отображение $\varkappa$ ведёт себя как гомоморфизм $H$-модулей и
\begin{equation*}\begin{split}
 a= j_1 \Biggl(\gamma_1\biggl(\varkappa (\tilde b_{s_1}) \Bigl( h_{0 s_1}
\tilde f_1\bigl( \varkappa (a^{(1)}_1),
 \ldots, \varkappa(a^{(1)}_{d_1});
 \ldots;  \varkappa(a^{(1)}_1),
 \ldots, \varkappa(a^{(1)}_{d_1});
\\ 
  \varkappa(\bar z_1), \ldots, \varkappa(\bar z_{\hat m_1})\bigr) \Bigr)\varkappa(b_{1 s_1}) \cdot \\ \cdot
   \Bigl( h_{1 s_1} f_1\bigl(\varkappa(a^{(1)}_1),
 \ldots, \varkappa (a^{(1)}_{d_1});
 \ldots;  \varkappa (a^{(1)}_1),
 \ldots, \varkappa (a^{(1)}_{d_1}); 
\\  
 \varkappa(\bar z^{(1)}_1), \ldots,  \varkappa(\bar z^{(1)}_{m_1})\bigr)\Bigr)  \varkappa(\tilde b_{1 s_1} )
 \varkappa(1_{B_1}) \biggr) \Biggr) j_2
 \cdot \\ \cdot
    \prod_{i=2}^{r}\Biggl(\gamma_i\biggl(\varkappa(b_{is_i}) \Bigl(h_{i s_i}  f_i\bigl(\varkappa(a^{(i)}_1),
 \ldots, \varkappa(a^{(i)}_{d_i});
 \ldots;  \varkappa(a^{(i)}_1),
 \ldots, \varkappa(a^{(i)}_{d_i});
 \\
  \varkappa(\bar z^{(i)}_1 ), \ldots, \varkappa(\bar z^{(i)}_{m_i})\bigr)\Bigr) \varkappa(\tilde b_{is_i}) 
  \varkappa(1_{B_i})\biggr) \Biggr) j_{i+1}\ne 0.
 \end{split}
 \end{equation*}

  Определим полилинейную функцию \begin{equation*}\begin{split}
 f_0 := v_1 \Biggl(\gamma_1\biggl(y_0 \Bigl( h_{0 s_1}
\tilde f_1\bigl( x^{(1)}_1,
 \ldots, x^{(1)}_{d_1};
 \ldots;  x^{(2\tilde k)}_1,
 \ldots, x^{(2\tilde k)}_{d_1};
  z_1, \ldots, z_{\hat m_1} \bigr) \Bigr) y_1 \cdot \\ \cdot
   \Bigl( h_{1 s_1} f_1\bigl(x^{(1,1)}_1,
 \ldots, x^{(1,1)}_{d_1};
 \ldots;  x^{(1,2k)}_1,
 \ldots, x^{(1,2k)}_{d_1}; 
 z^{(1)}_1, \ldots,  z^{(1)}_{m_1} \bigr)\Bigr) w_1  \biggr)\Biggr) v_2 
 \cdot \\ \cdot
   \prod_{i=2}^{r}\Biggl(\gamma_i\biggl(y_i \Bigl(h_{i s_i}  f_i\bigl(x^{(i,1)}_1,
 \ldots, x^{(i, 1)}_{d_i};
 \ldots;  x^{(i, 2k)}_1,
 \ldots, x^{(i, 2k)}_{d_i};
  z^{(i)}_1, \ldots, z^{(i)}_{m_i}\bigr)\Bigr) w_i\biggr) \Biggr) v_{i+1}.
 \end{split}
 \end{equation*}   
 (если $q_i=0$ для некоторых $i$ и элементы $j_i$ отсутствовали в~(\ref{EqAssocNonZero}),
 то переменная $v_i$ также отсутствует в $f_0$.)
Значение $f_0$ при подстановке
$x^{(\alpha)}_{\beta}=\varkappa(a^{(1)}_\beta)$,
    $x^{(i, \alpha)}_{\beta}=\varkappa(a^{(i)}_\beta)$,
    $z_i=\varkappa(\bar z_i)$,
 $z^{(i)}_{\beta}=\varkappa(\bar z^{(i)}_\beta)$, $v_i=j_i$,
 $y_0 = \varkappa(\tilde b_{s_1})$,
 $y_i = \varkappa(b_{is_i})$,
 $w_i = \varkappa(\tilde b_{is_i})\varkappa(1_{B_i})$
  равно $a\ne 0$. 
    Обозначим эту подстановку через $\Xi$.

Пусть $X_\ell = \bigcup_{i=1}^r X^{(i,\ell)}$, где $X^{(i,\ell)}=\lbrace x^{(i,\ell)}_\alpha \mid 1\leqslant \alpha \leqslant d_i \rbrace$. Обозначим через $\Alt_\ell$
оператор альтернирования по множеству $X_\ell$.
   Введём обозначение $\hat f := \Alt_1 \Alt_2 \ldots \Alt_{2k} f_0$.
   Заметим, что альтернирования не меняют $z_i, 
 z^{(i)}_{\beta}, v_i, y_i, w_i$, а
 $H$-многочлен $f_i$ кососимметричен по каждому из множеств $X^{(i)}_\ell$.
   Следовательно, значение выражения $\hat f$ при подстановке $\Xi$
   равно $\left((d_1)! (d_2)! \ldots (d_{r})!\right)^{2k} a \ne 0$,
   так как $B_1 \oplus \ldots \oplus B_r$ является
   прямой суммой $H$-инвариантных идеалов
    и если
альтернирование перемещает переменную из множества
   $X^{(i)}_\ell$ на место переменной из множества $X^{(i')}_\ell$
   при $i \ne i'$, то соответствующие элементы $h\varkappa(a^{(i)}_\beta)$, где $h\in H$,
обращаются в нуль при умножение на элементы пространства $\varkappa(B_{i'})$.
Здесь мы снова использовали тот факт, что в силу~(\ref{EqAssocbazero})
 в множителях слева от $\varkappa(1_{B_i})$
отображение $\varkappa$ ведёт себя как гомоморфизм $H$-модулей и алгебр.
  
  Заметим, что без дополнительных преобразований выражение $\hat f$
является полилинейной функцией, а не $H$-многочленом.
Однако в силу замечания~\ref{RemarkHActionOnWHn}
функцию $\hat f$ можно представить в виде $H$-многочлена
\begin{equation*}\begin{split}
 \tilde f := \Alt_1 \Alt_2 \ldots \Alt_{2k} v_1 y_0^{\tilde h_0} 
\tilde f'_1\left( x^{(1)}_1,
 \ldots, x^{(1)}_{d_1};
 \ldots;  x^{(2\tilde k)}_1,
 \ldots, x^{(2\tilde k)}_{d_1};
  z_1, \ldots, z_{\hat m_1} \right)  y_1^{\tilde h_1} \cdot \\ \cdot
f'_1\left(x^{(1,1)}_1,
 \ldots, x^{(1,1)}_{d_1};
 \ldots;  x^{(1,2k)}_1,
 \ldots, x^{(1,2k)}_{d_1}; 
 z^{(1)}_1, \ldots,  z^{(1)}_{m_1} \right) w_1^{\hat h_1}   v_2 
 \cdot \\ \cdot
   \prod_{i=2}^{r} y_i^{\tilde h_i} f'_i\left(x^{(i,1)}_1,
 \ldots, x^{(i, 1)}_{d_i};
 \ldots;  x^{(i, 2k)}_1,
 \ldots, x^{(i, 2k)}_{d_i};
  z^{(i)}_1, \ldots, z^{(i)}_{m_i}\right) w_i^{\hat h_i}  v_{i+1},
 \end{split}
 \end{equation*}     
где $f_i'$ и $\tilde f_1'$~"--- некоторые $H$-многочлены,
элементы $\tilde h_i, \hat h_i \in H$ получены из $h_{0s_1}$, $h_{is_i}$ и $\gamma_i$
применением равенства~(\ref{EqGeneralizedHopf}), а значение $H$-многочлена $\tilde f$ при
подстановке $\Xi$ снова равно $\left((d_1)! (d_2)! \ldots (d_{r})!\right)^{2k} a \ne 0$.

Представим теперь $\tilde f'_1$ в виде суммы одночленов
и заметим, что $\tilde f$
является линейной комбинацией полилинейных $H$-многочленов \begin{equation*}\begin{split}
 \tilde f_0 := \Alt_1 \Alt_2 \ldots \Alt_{2k} u_1^{\tau_1} \ldots u_s^{\tau_s}  \cdot \\ \cdot
f'_1\left(x^{(1,1)}_1,
 \ldots, x^{(1,1)}_{d_1};
 \ldots;  x^{(1,2k)}_1,
 \ldots, x^{(1,2k)}_{d_1}; 
 z^{(1)}_1, \ldots,  z^{(1)}_{m_1} \right) w_1^{\hat h_1}   v_2 
 \cdot \\ \cdot
   \prod_{i=2}^{r} y_i^{\tilde h_i} f'_i\left(x^{(i,1)}_1,
 \ldots, x^{(i, 1)}_{d_i};
 \ldots;  x^{(i, 2k)}_1,
 \ldots, x^{(i, 2k)}_{d_i};
  z^{(i)}_1, \ldots, z^{(i)}_{m_i}\right) w_i^{\hat h_i}  v_{i+1},
 \end{split}
 \end{equation*}    
 где $u_1, \ldots, u_s$~"--- переменные $x^{(\alpha)}_\beta$, 
 $y_0$, $y_1$, $z_i$ и, возможно, $v_1$, а $\tau_i\in H$~"--- некоторые элементы. Здесь $s=2\tilde k d_1 + \hat m_1+3$, если переменная $v_1$ присутствовала в $f_0$ и $s=2\tilde k d_1 + \hat m_1+2$, если отсутствовала. По крайней мере один из $H$-многочленов $\tilde f_0$ не является
 полиномиальным $H$-тождеством. Снова обозначим его через $\tilde f_0$.
 Заметим, что
 $$\deg \tilde f_0 \geqslant 2\tilde k d_1 + \hat m_1+1
 + \sum_{i=1}^r (2kd_i+m_i+2)>  2\tilde k d_1 + 2 k d > n.$$
 С другой стороны, $\sum_{i=1}^r (2kd_i+m_i+3)-1 \leqslant n$.
 Пусть \begin{equation*}\begin{split}
 f := \Alt_1 \Alt_2 \ldots \Alt_{2k} u_{(\deg \tilde f_0)-n+1 }^{\tau_{(\deg \tilde f_0)-n+1}} \ldots u_s^{\tau_s}  \cdot \\ \cdot
f'_1\left(x^{(1,1)}_1,
 \ldots, x^{(1,1)}_{d_1};
 \ldots;  x^{(1,2k)}_1,
 \ldots, x^{(1,2k)}_{d_1}; 
 z^{(1)}_1, \ldots,  z^{(1)}_{m_1} \right) w_1^{\hat h_1}   v_2 
 \cdot \\ \cdot
   \prod_{i=2}^{r} y_i^{\tilde h_i} f'_i\left(x^{(i,1)}_1,
 \ldots, x^{(i, 1)}_{d_i};
 \ldots;  x^{(i, 2k)}_1,
 \ldots, x^{(i, 2k)}_{d_i};
  z^{(i)}_1, \ldots, z^{(i)}_{m_i}\right) w_i^{\hat h_i}  v_{i+1}.
 \end{split}
 \end{equation*}    
 Тогда $f$ не обращается в нуль при подстановке $\Xi$. С другой стороны,
 $H$-многочлен $f$ кососимметричен по переменным каждого из множеств $X_\ell$, где $1 \leqslant \ell \leqslant 2k$. Наконец, $\deg f = n$. Теперь осталось переименовать переменные $H$-многочлена $f$ в $x_1, \ldots, x_n$  и заметить, что $f$ удовлетворяет всем условиям леммы.
 \end{proof}
\begin{proof}[Доказательство теоремы~\ref{TheoremHmodHRadAmitsurPIexpHBdimB}]
При $d>0$ достаточно применить лемму~\ref{LemmaAssocLowerPolynomial}
и теорему~\ref{TheoremAssocBounds}. 

 В случае $d=0$ алгебра $A$ совпадает с нильпотентным идеалом $J^H(A)$,
 откуда $c^H_n(A)=0$ при $n\geqslant p$.
\end{proof}

  \section{Применение понятия эквивалентности действий и случаи совпадения PI-экспонент}\label{SectionEquivApplToPolyIden}

В этом параграфе мы, в частности, покажем, как понятие эквивалентности модульных структур может быть
использовано для исследования коразмерностей полиномиальных $H$-тождеств.

Во-первых, из теоремы~\ref{TheoremHmoduleAssoc} можно вывести справедливость гипотезы Амицура~"--- Бахтурина для любой структуры $H$-модульной алгебры с $1$ на алгебре $\mathbbm{k}[x]/(x^2)$, где $H$~"--- некоторая алгебра Хопфа:

\begin{theorem}\label{TheoremDoubleNumbersAmitsurPIexpH}
Пусть на алгебре $\mathbbm{k}[x]/(x^2)$ 
определена структура $H$-модульной алгебры с $1$ 
для некоторой алгебры Хопфа $H$ над полем $\mathbbm{k}$ характеристики $0$. Обозначим через $d$
размерность идеала $J^H\bigl( \mathbbm{k}[x]/(x^2) \bigr)$.
Тогда существуют $C_1,C_2 > 0$ и $r_1,r_2 \in \mathbb R$ такие, что \begin{equation}\label{EquationHAmitsur}C_1 n^{r_1} (2-d)^n \leqslant c^H_n(\mathbbm{k}[x]/(x^2)) \leqslant C_2 n^{r_2} (2-d)^n\end{equation}
для всех $n\in\mathbb N$.

В частности, для $H$-тождеств алгебры $\mathbbm{k}[x]/(x^2)$ справедлива гипотеза Амицура~"--- Бахтурина. 
\end{theorem}
\begin{proof} Из леммы~\ref{LemmaHEquivCodimTheSame} 
следует, что~\eqref{EquationHAmitsur} достаточно доказать
для модульных структур, описанных в теореме~\ref{TheoremDoubleNumbersClassify}.
В первых двух случаях радикал Джекобсона алгебры $\mathbbm{k}[x]/(x^2)$, который совпадает с подпространством $\mathbbm{k}\bar x$, $H$-инвариантен, т.е. $d=1$. В последнем случае $\mathbbm{k}[x]/(x^2)$~"--- $H_4$-простая алгебра,
т.е. $d=0$. Как было уже отмечено в предыдущем параграфе, коразмерности не меняются при 
расширении основного поля. Принимая во внимание структуру алгебр $\mathbbm{k}[x]/(x^2)$
и $H_4$, заключаем, что при расширении основного поля не меняется и число $d = \dim_\mathbbm{k} \Bigl(J^H\bigl( \mathbbm{k}[x]/(x^2) \bigr)\Bigr)$. Отсюда без ограничения общности можно предполагать основное поле $\mathbbm{k}$
алгебраически замкнутым. Теперь в случае $H$-инвариантного радикала Джекобсона неравенства~\eqref{EquationHAmitsur} являются следствием
теоремы~\ref{TheoremHmoduleAssoc} и формулы~\eqref{EqdAssoc}, а в случае, когда алгебра $\mathbbm{k}[x]/(x^2)$ является $H_4$-простой,
неравенства~\eqref{EquationHAmitsur} следуют из теоремы~\ref{TheoremHOreAmitsur}
и формулы~\eqref{EqPIexpH(A)InvWedMalcev}.
\end{proof}

Наша следующая цель~"--- доказать, что в случае рационального действия связной редуктивной аффинной алгебраической группы автоморфизмами (в силу предложения~\ref{PropositionConnectedAffAlgAutNoAnti}, если связная группа действует автоморфизмами и антиавтоморфизмами, то она действует автоморфизмами) или действия конечномерной полупростой алгебры Ли дифференцированиями экспоненты соответствующих типов тождеств совпадают с обычными PI-экспонентами.

Начнём с действий алгебраических групп:

 \begin{theorem}\label{TheoremAssGPIexpEqual}
Пусть $A$~"--- конечномерная ассоциативная алгебра над алгебраически замкнутым полем $\mathbbm{k}$ характеристики $0$
с рациональным действием связной редуктивной аффинной алгебраической группы $G$ автоморфизмами.
Тогда $\PIexp^G(A)=\PIexp(A)$.
  \end{theorem}
  \begin{proof} В силу следствия~\ref{CorollaryGReductWedderburnMalcev} для алгебры $A$ существует $G$-инвариантное разложение Веддербёрна~"--- Мальцева $A=B\oplus J(A)$ (прямая сумма $G$-инвариантных подпространств). Согласно теореме~\ref{TheoremWedderburnHmod} максимальная полупростая подалгебра
  $B$ раскладывается в прямую сумму $G$-инвариантных идеалов $B_i$, являющихся $G$-простыми подалгебрами.
  В силу теоремы~\ref{TheoremGSimpleAssoc} все алгебры $B_i$ просты,
  откуда формула~(\ref{EqPIexpH(A)InvWedMalcev}) даёт для $\PIexp^G(A)$ и $\PIexp(A)$
  одно и то же значение.
  \end{proof}

Теперь выведем отсюда аналогичный результат для действий конечномерных полупростых алгебр Ли:
 
  \begin{theorem}\label{TheoremAssDerPIexpEqual}
  Пусть $A$~"--- конечномерная ассоциативная алгебра над полем $\mathbbm{k}$ характеристики $0$
с действием конечномерной полупростой алгебры Ли $\mathfrak g$ дифференцированиями.
Тогда $\PIexp^{U(\mathfrak g)}(A)=\PIexp(A)$.
  \end{theorem}
  \begin{proof} Как и в случае обычных тожеств, $\mathfrak g$-коразмерности не меняются при расширении
  основного поля, поэтому без ограничения общности можно считать, что поле $\mathbbm{k}$ алгебраически замкнуто.
   В силу теоремы~\ref{TheoremLieDiffActionReplacement} действию
  алгебры Ли $\mathfrak g$ соответствует рациональное действие некоторой связной редуктивной аффинной алгебраической группы $G$ автоморфизмами, причём $\mathfrak g$ является алгеброй Ли группы $G$.
  Согласно теореме~\ref{TheoremAffAlgGrAllEquiv} и лемме~\ref{LemmaHEquivCodimTheSame}
  справедливы равенства $c^{U(\mathfrak g)}_n(A)=c^G_n(A)$ для всех $n\in\mathbb N$.
  Теперь достаточно применить теорему~\ref{TheoremAssGPIexpEqual}.
  \end{proof}
  
  \begin{remark} В теореме~\ref{TheoremAssDerPIexpEqual}
доказана справедливость равенства $\PIexp^{U(\mathfrak g)}(A)=\PIexp(A)$, однако сами коразмерности
могут различаться.  
  Если на алгебре $M_2(\mathbbm{k})$ определить присоединённое представление алгебры Ли $\mathfrak{sl}_2(\mathbbm{k})$, то  $c^{U(\mathfrak{sl}_2(\mathbbm{k}))}_1(M_2(\mathbbm{k}))>1$, хотя $c_1(M_2(\mathbbm{k}))=1$.
  \end{remark}
  
Приведём ещё один пример совпадения PI-экспонент.    
 
Работы автора~\cite{MZ2009, MSB2010} были посвящены изучению асимптотического поведения
коразмерностей обобщённых полиномиальных тождеств. Понятие обобщённого $H$-действия позволяет
дать новое доказательство следующего результата:
\begin{theorem}[\cite{MSB2010}]\label{TheoremGIdAssoc}
Пусть $A$~"--- конечномерная
ассоциативная алгебра над полем $\mathbbm{k}$ характеристики $0$, $d:=\PIexp(A)$~"--- обычная PI-экспонента
алгебры $A$, а $gc_n(A)$~"--- последовательность коразмерностей её обобщённых полиномиальных
тождеств (см. определение в~\cite{MZ2009, MSB2010}).
Тогда \begin{enumerate}
\item при $d=0$ существует такое $n_0$, что $gc_n(A)=0$ при всех $n\geqslant n_0$;
\item при $d>0$ существуют такие константы $C_1, C_2 > 0$, $r_1, r_2 \in \mathbb R$, что $$C_1 n^{r_1} d^n \leqslant gc_n(A)
   \leqslant C_2 n^{r_2} d^n\text{ для всех }n \in \mathbb N.$$
\end{enumerate}
В частности, для любой такой алгебры $A$
   существует $\lim\limits_{n\to\infty} \sqrt[n]{gc_n(A)}=d\in\mathbb Z_+$ и, таким образом,
   справедлив аналог гипотезы Амицура.
\end{theorem}
\begin{proof}
Как было отмечено выше в примере~\ref{ExampleGenIdGenHAction},
всякая ассоциативная алгебра $A$ является алгеброй с обобщённым $H$-действием, где $H=A^{+} \otimes \left(A^{+}\right)^{\mathrm{op}}$. В силу того, что алгебра $H$ действует на $A$ линейными комбинациями композиций умножений на элементы алгебры $A$ справа и слева, справедливо равенство $gc_n(A)=c_n^H(A)$ для всех $n\in\mathbb N$.
Как и остальные типы коразмерностей, коразмерности обобщённых тождеств
и $H$-тождеств не меняются при расширении основного поля. Отсюда можно предполагать основное поле $\mathbbm{k}$
алгебраически замкнутым. Заметим, что $J(A)$ является $H$-подмодулем. Пусть $\varkappa \colon A/J(A) \hookrightarrow A$~"--- вложение, отвечающее некоторому обычному разложению Веддербёрна~"--- Мальцева. Вложение $\varkappa$ является гомоморфизмом алгебр. Пусть $A/J(A)=B_1 \oplus \ldots \oplus B_q$~"---
разложение полупростой алгебры $A/J(A)$ в прямую сумму идеалов, являющихся простыми
алгебрами. Тогда все $B_i$ являются $H$-подмодулями,
где структура $H$-модуля на $A/J(A)$
наследуется с соответствующей структуры на $A$.

Заметим, что вложение $\varkappa$ удовлетворяет всем условиям теоремы~\ref{TheoremHmodHRadAmitsurPIexpHBdimB}. Более того,
\begin{equation*}\begin{split}\max\dim\left( B_{i_1}\oplus B_{i_2} \oplus \ldots \oplus B_{i_r}
 \mathbin{\Bigl|}  r \geqslant 1,\right.\\
   \left. (H\varkappa(B_{i_1}))A^+ \,(H\varkappa(B_{i_2})) A^+ \ldots (H\varkappa(B_{i_{r-1}})) A^+\,(H\varkappa(B_{i_r}))\ne 0\right) = \\ =
   \max\dim\left( B_{i_1}\oplus B_{i_2} \oplus \ldots \oplus B_{i_r}
 \mathbin{\Bigl|}  r \geqslant 1,\right.\\
   \left. A^+ \varkappa(B_{i_1})A^+ \varkappa(B_{i_2}) A^+ \ldots \varkappa(B_{i_{r-1}}) A^+\,\varkappa(B_{i_r})\ne 0\right)=\\ =\max\dim\left( B_{i_1}\oplus B_{i_2} \oplus \ldots \oplus B_{i_r}
 \mathbin{\Bigl|}  r \geqslant 1,\right.\\
   \left.  \varkappa(B_{i_1}) J(A) \varkappa(B_{i_2}) J(A) \ldots \varkappa(B_{i_{r-1}}) J(A)\,\varkappa(B_{i_r})\ne 0\right) = \PIexp(A).
  \end{split}\end{equation*}
  
Теперь утверждение теоремы следует из теоремы~\ref{TheoremHmodHRadAmitsurPIexpHBdimB}.
\end{proof}

\section{Примеры и приложения}\label{SectionAssocExamples}

 Воспользуемся теперь формулой~\eqref{EqPIexpH(A)InvWedMalcev}
для того, чтобы вычислить $H$-PI-экспоненты для некоторых важных примеров алгебр 
с обобщённым $H$-действием. 
 
 \subsection{Суммы $H$-простых алгебр и критерии $H$-простоты}
 
   \begin{example}\label{ExampleHStarSum}
 Пусть $A=B_1 \oplus B_2 \oplus \ldots \oplus B_q$ (прямая сумма $H$-инвариантных идеалов)~"---
 ассоциативная алгебра с обобщённым $H$-действием, где $B_i$~"--- конечномерные $H$-простые
 алгебры, удовлетворяющие условию (*), а $H$~"---  ассоциативная алгебра с~$1$
 над полем $\mathbbm{k}$ характеристики $0$. Введём обозначение $d := \max_{1 \leqslant k
  \leqslant q} \dim B_k$. Тогда существуют такие $C_1, C_2 > 0$, $r_1, r_2 \in \mathbb R$,
  что $$C_1 n^{r_1} d^n \leqslant c_n^{H}(A)
  \leqslant C_2 n^{r_2} d^n \text{ для всех } n\in\mathbb N.$$
 \end{example}
 \begin{proof}  Докажем сперва, что все $H$-инвариантные идеалы алгебры $A$ являются прямыми суммами некоторых из идеалов $B_j$. Действительно, если $I \subseteq A$~"--- некоторый $H$-инвариантный идеал,
 то $I B_j, B_j I \subseteq I \cap B_j$, где $I \cap B_j$~"--- также двусторонний $H$-инвариантный идеал, т.е. в силу $H$-простоты алгебр $B_j$ для всякого $j$ либо $B_j \subseteq I$, либо $ I B_j = B_j I = 0$.
 В последнем случае $I \subseteq \bigoplus_{i\ne j} B_i$, так как в силу условия (*)
 алгебра $B_j$ с единицей. Отсюда $I$ является прямой суммой всех таких $B_j$,
 что $B_j \subseteq I$.
  В частности, $J^H(A)=0$ и $\varkappa = \id_A$.
 Теперь осталось применить теорему~\ref{TheoremGenHAmitsurHWedederburn}.
 \end{proof}
 
 С учётом примера~\ref{ExampleHmodHSemiSimplePropertyStar} и следствия~\ref{CorollaryOreExtStar}
 получаем следующий пример:
 
    \begin{example}\label{ExampleHOreExtSum}
 Пусть $A=B_1 \oplus B_2 \oplus \ldots \oplus B_q$ (прямая сумма $H$-инвариантных идеалов)~"---
 ассоциативная $H$-модульная алгебра 
 над алгебраически замкнутым полем $\mathbbm{k}$ характеристики $0$, где $B_i$~"--- конечномерные $H$-простые
 алгебры, а алгебра Хопфа $H$ либо сама является конечномерной полупростой алгеброй Хопфа,
либо получена при помощи (возможно, многократного)
расширения Оре конечномерной полупростой алгебры Хопфа косопримитивными элементами. (Например,
$H=H_{m^2}(\zeta)$.) Введём обозначение $d := \max_{1 \leqslant k
  \leqslant q} \dim B_k$. Тогда существуют такие $C_1, C_2 > 0$, $r_1, r_2 \in \mathbb R$,
  что $$C_1 n^{r_1} d^n \leqslant c_n^{H}(A)
  \leqslant C_2 n^{r_2} d^n \text{ для всех } n\in\mathbb N.$$
 \end{example}
 
 Отсюда вытекает следующий критерий $H$-простоты:
 
     \begin{theorem}\label{TheoremHOreHSimplicityCriterion}
 Пусть $A$~"---конечномерная ассоциативная $H$-модульная алгебра
 над алгебраически замкнутым полем $\mathbbm{k}$ характеристики $0$, где алгебра Хопфа $H$ либо сама является конечномерной полупростой алгеброй Хопфа, либо получена при помощи (возможно, многократного)
расширения Оре конечномерной полупростой алгебры Хопфа косопримитивными элементами. 
Тогда алгебра $A$ является $H$-простой, если и только если $\PIexp^H(A)=\dim A$.
 \end{theorem}
 \begin{proof} Если $J^H(A)\ne 0$, то из~\eqref{EqdAssoc}
 следует, что $$\PIexp^H(A) \leqslant \dim A - \dim J^H(A) < \dim A.$$
 Если же $J^H(A) = 0$, то в силу теоремы~\ref{TheoremSkryabinVanOystaeyen}
 и леммы~\ref{LemmaHSemiSimpleIsUnital} алгебра $A$ является прямой суммой $H$-инвариантных
 идеалов, каждый из которых является $H$-простой алгеброй, т.е. можно воспользоваться
примером~\ref{ExampleHOreExtSum}.
 \end{proof}
 
Используя предложение~\ref{PropositionCnGrCnGenH}, следствие~\ref{CorollaryRadicalGraded}, теоремы~\ref{TheoremTGradedGActionInvWedderburnArtin} и \ref{TheoremAssocAlternateFinal}, получаем аналогичные примеры и результаты для градуированных алгебр
и алгебр с действием некоторой группы автоморфизмами и антиавтоморфизмами:
 
    \begin{example}\label{ExampleGradedSum}
 Пусть $A=B_1 \oplus B_2 \oplus \ldots \oplus B_q$ (прямая сумма градуированных идеалов)~"---
 ассоциативная алгебра 
 над алгебраически замкнутым полем $\mathbbm{k}$ характеристики $0$, градуированная некоторой группой $G$, где $B_i$~"--- конечномерные 
 градуированно простые
 алгебры. Введём обозначение $d := \max_{1 \leqslant k
  \leqslant q} \dim B_k$. Тогда существуют такие $C_1, C_2 > 0$, $r_1, r_2 \in \mathbb R$,
  что $$C_1 n^{r_1} d^n \leqslant c_n^{G\text{-}\mathrm{gr}}(A)
  \leqslant C_2 n^{r_2} d^n \text{ для всех } n\in\mathbb N.$$
 \end{example}
 
      \begin{theorem}\label{TheoremGradedSimplicityCriterion}
 Пусть $A$~"---конечномерная ассоциативная алгебра над алгебраически замкнутым полем $\mathbbm{k}$ характеристики $0$, градуированная некоторой группой $G$. 
Тогда алгебра $A$ является градуированно простой, если и только если $\PIexp^{G\text{-}\mathrm{gr}}(A)=\dim A$.
 \end{theorem}

   \begin{example}\label{ExampleGActionSum}
 Пусть $A=B_1 \oplus B_2 \oplus \ldots \oplus B_q$ (прямая сумма $G$-инвариантных идеалов)~"---
 ассоциативная алгебра  
 над алгебраически замкнутым полем $\mathbbm{k}$ характеристики $0$ с действием некоторой группы $G$ автоморфизмами и антиавтоморфизмами, где $B_i$~"--- конечномерные $G$-простые
 алгебры. Введём обозначение $d := \max_{1 \leqslant k
  \leqslant q} \dim B_k$. Тогда существуют такие $C_1, C_2 > 0$, $r_1, r_2 \in \mathbb R$,
  что $$C_1 n^{r_1} d^n \leqslant c_n^{G}(A)
  \leqslant C_2 n^{r_2} d^n \text{ для всех } n\in\mathbb N.$$
 \end{example}

      \begin{theorem}\label{TheoremGActionSimplicityCriterion}
 Пусть $A$~"---конечномерная ассоциативная алгебра 
 над алгебраически замкнутым полем $\mathbbm{k}$ характеристики $0$ с действием некоторой группы $G$ автоморфизмами и антиавтоморфизмами. 
Тогда алгебра $A$ является $G$-простой, если и только если $\PIexp^{G}(A)=\dim A$.
 \end{theorem}

\subsection{Примеры алгебр, градуированных группами}

В примерах~\ref{Example2M2S3}--\ref{ExampleUT} основное поле $\mathbbm{k}$ является произвольным полем характеристики~$0$.

Напомним, что коразмерности тождеств не меняются при расширении основного поля.
В примерах, приведённых ниже, оценки для коразмерностей также не зависят от основного
поля. Поэтому во всех доказательствах можно без ограничения
общности считать поле $\mathbbm{k}$ алгебраически замкнутым.

\begin{example}\label{Example2M2S3}
Пусть $G=S_3$, а $A=M_2(\mathbbm{k})\oplus M_2(\mathbbm{k})$. Рассмотрим следующую $G$-градуировку
на~$A$: $$A^{(e)} = \left\lbrace\left(\begin{array}{rr}
\alpha & 0 \\
 0 & \beta 
\end{array} \right)\right\rbrace \oplus  \left\lbrace\left(\begin{array}{rr}
\gamma & 0 \\
 0 & \mu 
\end{array} \right)\right\rbrace,$$ $$A^{\bigl((12)\bigr)} = \left\lbrace\left(\begin{array}{rr}
0 & \alpha \\
 \beta & 0  
\end{array} \right)\right\rbrace \oplus  0,\qquad A^{\bigl((23)\bigr)} = 0 \oplus \left\lbrace\left(\begin{array}{rr}
0 & \alpha \\
 \beta & 0  
\end{array} \right)\right\rbrace,$$ а остальные компоненты равны $0$.
Тогда существуют такие $C_1, C_2 > 0$, $r_1, r_2 \in \mathbb R$,
  что $$C_1 n^{r_1} 4^n \leqslant c_n^{G\text{-}\mathrm{gr}}(A)
  \leqslant C_2 n^{r_2} 4^n \text{ для всех } n \in\mathbb N.$$
\end{example}
\begin{proof} Достаточно заметить, что обе копии алгебры $M_2(\mathbbm{k})$ являются градуированными идеалами алгебры $A$, и воспользоваться примером~\ref{ExampleGradedSum}.
\end{proof}

\begin{example}\label{ExampleGroupAlgebra}
Пусть $A=\mathbbm{k}G$, где $G$~"--- конечная группа. Рассмотрим естественную $G$-градуировку $A=\bigoplus_{g\in G} A^{(g)}$, где $A^{(g)}=\mathbbm{k}g$. Тогда существуют такие $C_1, C_2 > 0$, $r_1, r_2 \in \mathbb R$,
что $$C_1 n^{r_1} |G|^n \leqslant c_n^{G\text{-}\mathrm{gr}}(A)
  \leqslant C_2 n^{r_2} |G|^n \text{ для всех } n \in\mathbb N.$$
\end{example}
\begin{proof}
Достаточно заметить, что алгебра $\mathbbm{k}G$
является градуированно простой, и воспользоваться примером~\ref{ExampleGradedSum}.
\end{proof}

\begin{remark}
В следствии 3.4 из работы~\cite{AljadeffKanelBelov} Э.~Альхадефф и А.\,Я.~Канель-Белов
показали, что $$ |G|^n \leqslant c_n^{G\text{-}\mathrm{gr}}(\mathbbm{k}G) \leqslant |G'|\cdot |G|^n \text{ для всех } n\in\mathbb N,$$
где $G'$~"--- коммутант группы $G$. Кроме того, они доказали, что $\lim\limits_{n\to\infty}\left(\frac{c_n^{G\text{-}\mathrm{gr}}(\mathbbm{k}G)}{|G'|\cdot |G|^n}\right) = 1$. 
\end{remark}

\subsection{Примеры алгебр с действием групп и алгебр Ли}

\begin{example}
Пусть $\mathbbm{k}$~"--- поле характеристики $0$,
а
$$A = \left\lbrace\left(\begin{array}{cc} C & D \\
0 & 0
  \end{array}\right) \mathrel{\biggl|} C, D\in M_m(\mathbbm{k})\right\rbrace
  \subseteq M_{2m}(\mathbbm{k}),$$ $m \geqslant 2$.
       Определим $\varphi \in \Aut(A)$ по формуле
  $$\varphi\left(\begin{array}{cc} C & D \\
0 & 0
  \end{array}\right)=\left(\begin{array}{cc} C & C+D \\
0 & 0
  \end{array}\right).$$
  Тогда группа $G=\langle \varphi \rangle
  \cong \mathbb Z$ действует на алгебре $A$ автоморфизмами, т.е. $A$ является $\mathbbm{k}G$-модульной
  алгеброй. Как было показано в примере~\ref{ExampleGnoninvWedderburnMalcev},
  для алгебры $A$ не существует $G$-инвариантного разложения Веддербёрна~"--- Мальцева.
  В то же время существуют такие константны $C_1, C_2 > 0$, $r_1, r_2 \in \mathbb R$, что $$C_1 n^{r_1} m^{2n} \leqslant c^G_n(A) \leqslant C_2 n^{r_2} m^{2n}\text{ для всех }n \in \mathbb N.$$
\end{example}
\begin{proof} Как уже было отмечено, коразмерности не меняются при расширении основного поля.
Более того, при расширении основного поля алгебра $A$ остаётся алгеброй того же типа.
Следовательно, без ограничения общности можно считать, что основное поле алгебраически замкнуто.
Теперь достаточно заметить, что $A/J(A) \cong M_m(\mathbbm{k})$,
откуда в силу теоремы~\ref{TheoremHmodHRadAmitsurPIexpHBdimB} справедливо равенство  $\PIexp^G(A)=\dim M_m(\mathbbm{k})= m^2$.
\end{proof}

Проводя аналогичные рассуждения, получаем следующий пример:
\begin{example}
Пусть $A$~"---  та же ассоциативная алгебра, что и в предыдущем примере.
Определим на векторном пространстве $A$ структуру алгебры Ли при помощи коммутатора  $[x,y]=xy-yx$
и обозначим соответствующую алгебру Ли через $\mathfrak g$.
Рассмотрим стандартное представление алгебры Ли $\mathfrak g$ на $A$ дифференцированиями.
Тогда $A$ оказывается $U(\mathfrak g)$-модульной ассоциативной алгеброй.
Как было показано в примере~\ref{ExampleU(L)noninvWedderburnMalcev}, в $A$ не существует $U(\mathfrak g)$-инвариантного разложения Веддербёрна~"--- Мальцева.  В то же время существуют такие константны $C_1, C_2 > 0$, $r_1, r_2 \in \mathbb R$, что $$C_1 n^{r_1} m^{2n} \leqslant c^{U(\mathfrak g)}_n(A) \leqslant C_2 n^{r_2} m^{2n}\text{ для всех }n \in \mathbb N.$$
\end{example}

  \begin{example}\label{ExampleFields}
  Пусть $A = \mathbbm{k} e_1 \oplus \ldots \oplus \mathbbm{k}{e_m}$ (прямая сумма идеалов),
  где $e_i^2=e_i$, $m \in \mathbb N$. Предположим, что группа $G \subseteq S_m$ действует
  на $A$ по формуле $\sigma e_i := e_{\sigma(i)}$, $\sigma \in G$.
  Рассмотрим разложение $\left\lbrace 1,2, \ldots, m \right\rbrace = \coprod_{i=1}^q O_i$,
  где $O_i$~"--- орбиты $G$-действия на $\left\lbrace 1,2, \ldots, m \right\rbrace$.
  Введём обозначение $d := \max_{1 \leqslant i
  \leqslant q} |O_i|$. Тогда существуют такие $C_1, C_2 > 0$, $r_1, r_2 \in \mathbb R$,
  что $C_1 n^{r_1} d^n \leqslant c_n^G(A)
  \leqslant C_2 n^{r_2} d^n$ для всех $n \in\mathbb N$.
 \end{example}
 \begin{proof}
 Заметим, что $A=B_1 \oplus \ldots \oplus B_q$,
 где $B_i := \langle e_j \mid j \in O_i \rangle_\mathbbm{k}$~"--- $G$-инвариантные идеалы.
 Докажем, что $B_i$ являются $G$-простыми алгебрами для всех $1 \leqslant i \leqslant q$.
 Действительно, если $I$~"--- нетривиальный $G$-инвариантный идеал алгебры $B_i$,
 то существует такое $a = \sum_{j \in O_i} \alpha_j e_j \in I$, что $\alpha_j \in \mathbbm{k}$
 и $\alpha_k \ne 0$ для некоторого $k\in O_i$.
 Следовательно, $e_k = \frac{1}{\alpha_k} e_k a \in I$. Более того,
 для любого $j\in O_i$ существует такое $\sigma \in G$, что  $e_j = \sigma e_{k}$.
  Отсюда $I=B_i$ и алгебра $B_i$ является $G$-простой.

 Теперь из примера~\ref{ExampleGActionSum} следует, что $\PIexp^G(A)=\max_{1 \leqslant i
  \leqslant q} \dim B_i = \max_{1 \leqslant i
  \leqslant q} |O_i|$.
 \end{proof}

Ниже в примере~\ref{ExampleMatrix} группа $G$ может действовать не только автоморфизмами, но и антиавтоморфизмами:

  \begin{example}\label{ExampleMatrix}
  Пусть $A = A_1 \oplus \ldots \oplus A_m$ (прямая сумма идеалов),
   $A_i \cong M_k(\mathbbm{k})$, $1\leqslant i \leqslant m$, и $k,m \in \mathbb N$.
   Обозначим через $\Aut^*(M_k(\mathbbm{k}))$ группу автоморфизмов и антиавтоморфизмов
   алгебры  $M_k(\mathbbm{k})$.
   Тогда группа $\Aut^*(M_k(\mathbbm{k})) \times S_m$
   действует на $A$ следующим образом: если $(\varphi,\sigma)\in \Aut^*(M_k(\mathbbm{k})) \times S_m$
   и $(a_1, \ldots, a_m) \in A$, то
   $$(\varphi, \sigma) \cdot (a_1, \ldots, a_m)
   := \bigl(\varphi (a_{\sigma^{-1}(1)}), \ldots,\varphi(a_{\sigma^{-1}(m)})\bigr).$$
   Пусть $G \subseteq \Aut^*(M_k(\mathbbm{k})) \times S_m$~"--- некоторая  подгруппа.
   Обозначим через $$\pi \colon \Aut^*(M_k(\mathbbm{k})) \times S_m \to S_m$$ естественную проекцию на вторую компоненту.
  Рассмотрим разложение $\left\lbrace 1,2, \ldots, m \right\rbrace = \coprod_{i=1}^q O_i$,
  где $O_i$~"--- орбиты $\pi(G)$-действия на $\left\lbrace 1,2, \ldots, m \right\rbrace$.
  Введём обозначение $d := k^2 \max_{1 \leqslant i
  \leqslant q} |O_i|$. Тогда существуют такие $C_1, C_2 > 0$, $r_1, r_2 \in \mathbb R$,
  что $C_1 n^{r_1} d^n \leqslant c_n^G(A)
  \leqslant C_2 n^{r_2} d^n$ для всех $n \in\mathbb N$.
 \end{example}
 \begin{proof}
 Заметим, что $A=B_1 \oplus \ldots \oplus B_q$,
 где $B_i := \bigoplus_{j \in O_i} A_j$~"--- $G$-инвариантные идеалы.
 Докажем, что $B_i$ являются $G$-простыми алгебрами для всех $1 \leqslant i \leqslant q$.
 Действительно, если $I$~"--- нетривиальный $G$-инвариантный идеал алгебры $B_i$,
 то существует такое $a = \sum_{j \in O_i} a_j \in I$, что $a_j \in A_j$
 и $a_\ell \ne 0$ для некоторого $\ell \in O_i$. Обозначим через $e_\ell$
 единичную матрицу компоненты $A_\ell$, которая, напомним,
 изоморфна алгебре $M_k(\mathbbm{k})$. Тогда $e_\ell a = a_\ell \in I$
 и $I \cap A_\ell \ne 0$. Поскольку алгебра $A_\ell$ проста,
 $I \cap A_\ell = A_\ell$. Заметим теперь, что для любого $j\in O_i$
 существует такое $g\in G$, что $A_j = g A_\ell$.
  Отсюда $I=B_i$, и алгебра $B_i$ является $G$-простой.

 Теперь из примера~\ref{ExampleGActionSum} следует, что $\PIexp^G(A)=\max_{1 \leqslant i
  \leqslant q} \dim B_i = k^2 \max_{1 \leqslant i
  \leqslant q} |O_i|$.
 \end{proof}

Ниже в примере~\ref{ExampleUT} алгебра не является полупростой.

  \begin{example}\label{ExampleUT}
  Пусть $A = A_1 \oplus \ldots \oplus A_m$ (прямая сумма идеалов), где
   $A_i \cong \UT_k(\mathbbm{k})$, $1\leqslant i \leqslant m$,  $k,m \in \mathbb N$,
   а $\UT_k(\mathbbm{k})$~"--- ассоциативная алгебра верхнетреугольных матриц размера $k\times k$.
   Пусть группа $G \subseteq S_m$ действует на $A$ следующим образом: если $\sigma\in G$
   и $(a_1, \ldots, a_m) \in A$, то
   $$\sigma \cdot (a_1, \ldots, a_m)
   := (a_{\sigma^{-1}(1)}, \ldots, a_{\sigma^{-1}(m)}).$$
   Рассмотрим разложение $\left\lbrace 1,2, \ldots, m \right\rbrace = \coprod_{i=1}^s O_i$,
   где $O_i$~"--- орбиты $G$-действия на $\left\lbrace 1,2, \ldots, m \right\rbrace$.
  Введём обозначение $d := k \cdot\max_{1 \leqslant i
  \leqslant s} |O_i|$. Тогда существуют такие $C_1, C_2 > 0$, $r_1, r_2 \in \mathbb R$,
  что $$C_1 n^{r_1} d^n \leqslant c_n^G(A)
  \leqslant C_2 n^{r_2} d^n
   \text{ для всех } n \in\mathbb N.$$
 \end{example}
 \begin{proof} Обозначим через $e_{ij}^{(t)}$, где $1\leqslant i \leqslant j \leqslant k$,
 матричные единицы алгебры $A_t$, где $1 \leqslant t \leqslant m$.
 Тогда $\sigma e_{ij}^{(t)} = e_{ij}^{(\sigma(t))}$
 Заметим, что
\begin{equation}\label{EqWedderburnUT}
 A=\left( \bigoplus_{i=1}^s \bigoplus_{j=1}^k B_{ij}\right) \oplus J,
 \end{equation}
 где $B_{ij} := \langle e_{jj}^{(t)}  \mid t \in O_i \rangle_\mathbbm{k}$ являются $G$-инвариантными
 подалгебрами, а $$J:=\langle e_{ij}^{(t)}  \mid 1\leqslant i < j \leqslant k,
 1 \leqslant t \leqslant m \rangle_\mathbbm{k}$$~"--- $G$-инвариантный нильпотентный идеал.
 Докажем, что $B_{ij}$ являются $G$-простыми алгебрами для всех $1 \leqslant i \leqslant s$
 и $1 \leqslant j \leqslant k$.
 Действительно, если $I$~"--- нетривиальный $G$-инвариантный идеал алгебры $B_{ij}$,
 то существует такое $a = \sum_{t \in O_i} \alpha_t e_{jj}^{(t)} \in I$, что $\alpha_t \in \mathbbm{k}$
 и $\alpha_\ell \ne 0$ для некоторого $\ell \in O_i$.
  Тогда $e^{(\ell)}_{jj} = \frac{1}{\alpha_\ell} e^{(\ell)}_{jj} a  \in I$.
   Теперь заметим, что для любого $t\in O_i$
 существует такое $\sigma \in G$, что $e_{jj}^{(t)} = \sigma\left(e^{(\ell)}_{jj}
 \right)$. Отсюда $I=B_i$, и алгебра $B_i$ является $G$-простой. Следовательно, (\ref{EqWedderburnUT})~"--- $G$-инвариантное разложение Веддербёрна~"--- Мальцева алгебры $A$.

  Пусть  $1 \leqslant i \leqslant s$, $t \in O_i$. Тогда
  $$e_{11}^{(t)} e_{12}^{(t)} e_{22}^{(t)} e_{23}^{(t)} \ldots e_{kk}^{(t)} = e_{1k}^{(t)} \ne 0,$$
  откуда $$B_{i1}JB_{i2}J\ldots J B_{ik} \ne 0$$
  и в силу~\eqref{EqPIexpH(A)InvWedMalcev} справедливо неравенство
\begin{equation}\label{EqExUTLower}
  \PIexp^G(A)\geqslant \dim(B_{i1} \oplus \ldots \oplus B_{ik})
  = k |O_i|.\end{equation}

  Предположим, что $$B_{i_1 j_1}JB_{i_2 j_2}J\ldots J B_{i_r j_r} \ne 0$$
  для некоторых $1 \leqslant i_\ell \leqslant s$, $1 \leqslant j_\ell \leqslant k$.
  Тогда можно выбрать такие $e^{(q_\ell)}_{j_\ell j_\ell} \in B_{i_\ell j_\ell}$,
  $q_\ell \in O_{i_\ell}$, где $1 \leqslant j_\ell \leqslant k$,
  и такие $e^{(q'_\ell)}_{i'_\ell j'_\ell} \in J$, где $1 \leqslant i'_\ell < j'_\ell
  \leqslant k$, что
  $$e^{(q_1)}_{j_1 j_1} e^{(q'_1)}_{i'_1 j'_1} e^{(q_2)}_{j_2 j_2} e^{(q'_2)}_{i'_2 j'_2}
  \ldots e^{(q'_{r-1})}_{i'_{r-1} j'_{r-1}} e^{(q_r)}_{j_r j_r} \ne 0.$$
  В этом случае $q_1 = q'_1 = q_2 = q'_2=\ldots =q'_{r-1}= q_r$ и
  $j_\ell = i'_\ell$, $j'_{\ell-1} = j_\ell$.
  Следовательно, $i_1=\ldots = i_r$, $r \leqslant k$ и
  $ \dim(B_{i_1 j_1} \oplus \ldots \oplus B_{i_r j_r})
  \leqslant k |O_{i_1}|$.
  Отсюда $\PIexp^G(A) \leqslant k\cdot \max_{1 \leqslant i
  \leqslant s} |O_i|$. Оценка снизу была получена в~(\ref{EqExUTLower}).
 \end{proof}

\newpage

\chapter{Рост коразмерностей  полиномиальных  $H$-тождеств в $H$-модульных алгебрах Ли}
\label{ChapterHModLieCodim}

В данной главе аналог гипотезы Амицура (существование целочисленной экспоненты) доказывается для достаточно широкого класса конечномерных алгебр Ли с дополнительной структурой,
включающего в себя $H$-модульные алгебры Ли для конечномерных полупростых алгебр Хопфа $H$ (теорема~\ref{TheoremMainLieHSS}), алгебры Ли с рациональным действием аффинной алгебраической группы автоморфизмами и антиавтоморфизмами (теорема~\ref{TheoremMainLieGAffAlg}), алгебры Ли с действием конечномерной полупростой алгебры Ли дифференцированиями (теорема~\ref{TheoremMainDiffLie}), алгебры Ли, градуированные произвольными группами (теорема~\ref{TheoremMainLieGr}),
$H$-модульные алгебры Ли,  разрешимый радикал которых является нильпотентным $H$-инвариантным идеалом
(теорема~\ref{TheoremMainLieNRSame}), и  алгебры Ли, простые  по отношению к действию алгебры Тафта
(теорема~\ref{TheoremTaftSimpleLieHPIexpExists}).

Результаты главы были опубликованы в работах~\cite{ASGordienko2, ASGordienko7, ASGordienko6Kochetov,
ASGordienko5, ASGordienko17}.

   \section{$H$-хорошие алгебры Ли}\label{SectionHnice}
   
Напомним, что если $M$~"--- модуль над некоторой алгеброй Хопфа $H$,
то алгебра $\End_\mathbbm{k}(M)$ наделена структурой ассоциативной $H$-модульной алгебры
(см. пример~\ref{ExampleHModEnd}), однако соответствующая
алгебра Ли $\mathfrak{gl}(M)$, хотя и наследует с $\End_\mathbbm{k}(M)$ структуру
$H$-модуля, вообще говоря, не обязана являться $H$-модульной алгеброй Ли.

   Для того, чтобы сформулировать требования, предъявляемые к $H$-действию в основной теореме главы (см. теорему~\ref{TheoremMainLieH}),
   дадим следующее определение.

   Пусть $L$~"--- конечномерная $H$-модульная алгебра Ли,
   где $H$~"--- алгебра Хопфа
   над алгебраически замкнутым полем $\mathbbm{k}$ характеристики $0$.
   Будем говорить, что алгебра Ли $L$ является \textit{$H$-хорошей\footnote{В работе автора~\cite{ASGordienko5}, в которой впервые было введено это понятие, использовался
   английский термин \textit{$H$-nice}.}},
   если
   выполнены следующие условия:
\begin{enumerate}
\item \label{ConditionRNinv}
нильпотентный радикал $N$ и разрешимый радикал $R$ алгебры Ли $L$ являются $H$-подмодулями;
\item \label{ConditionLevi} \textit{(существование разложения Леви)}
существует такая $H$-инвариантная максимальная полупростая подалгебра $B \subseteq L$,
что $L=B\oplus R$ (прямая сумма $H$-модулей);
\item \label{ConditionWedderburn} \textit{(существование разложений Веддербёрна~"--- Мальцева)}
для любого $H$-подмодуля $W \subseteq L$
и любой ассоциативной $H$-модульной подалгебры $A_1 \subseteq \End_\mathbbm{k}(W)$
радикал Джекобсона $J(A_1)$ является $H$-подмодулем
и существует такая  $H$-инвариантная максимальная полупростая ассоциативная подалгебра $\tilde A_1 \subseteq A_1$, что $A_1 = \tilde A_1 \oplus J(A_1)$ (прямая сумма $H$-модулей);
\item \label{ConditionLComplHred}
для любой $H$-инвариантной подалгебры Ли $L_0 \subseteq \mathfrak{gl}(L)$,
такой, что $L_0$ является $H$-модульной алгеброй Ли, а пространство $L$ является вполне приводимым $L_0$-модулем без учёта $H$-действия, пространство $L$ является вполне приводимым $(H,L_0)$-модулем.
\end{enumerate}

Приведём основные примеры $H$-хороших алгебр Ли:

\begin{example}\label{ExampleHniceHSS}
Пусть $L$~"--- конечномерная $H$-модульная алгебра Ли,
   где $H$~"--- конечномерная полупростая алгебра Хопфа
   над алгебраически замкнутым полем $\mathbbm{k}$ характеристики $0$.
Тогда алгебра Ли $L$ является $H$-хорошей.
\end{example}
\begin{proof}
В силу следствия~\ref{CorollaryHModRadicalsLie} нильпотентный радикал $N$ и разрешимый радикал $R$ алгебры Ли $L$ являются $H$-подмодулями.
Согласно теореме~\ref{TheoremHLevi}
существует такая $H$-инвариантная максимальная полупростая подалгебра $B \subseteq L$,
что $L=B\oplus R$ (прямая сумма $H$-модулей).
В силу следствий~\ref{CorollaryRadicalHSSSubMod} и~\ref{CorollaryHSSWedderburnMalcev}
существуют $H$-инвариантные разложения Веддербёрна~"--- Мальцева.
Из теоремы~\ref{TheoremHWeyl}
вытекает выполнение условия~\ref{ConditionLComplHred}.

 Следовательно, алгебра Ли $L$ является $H$-хорошей.
\end{proof}

\begin{example}\label{ExampleHniceAffAlg}
Пусть $L$~"--- конечномерная алгебра Ли
   над алгебраически замкнутым полем $\mathbbm{k}$ характеристики $0$,
   на которой рационально действует автоморфизмами
   некоторая редуктивная аффинная алгебраическая группа $G$.
  Тогда алгебра Ли $L$ является $\mathbbm{k}G$-хорошей,
где действие алгебры Хопфа $\mathbbm{k}G$ на $L$ является продолжением
$G$-действия по линейности.
\end{example}
\begin{proof} Во-первых, разрешимый
и нильпотентный радикалы являются $\mathbbm{k}G$-подмодулями,
поскольку радикалы переходят в себя при любых автоморфизмах.
Отсюда условие~\ref{ConditionRNinv} выполнено. 
В силу теоремы~\ref{TheoremAffAlgGrLevi} существует $G$-инвариантное 
разложение Леви.
 
 Рассмотрим произвольный $\mathbbm{k}G$-подмодуль $W \subseteq L$
 и $\mathbbm{k}G$-действие, заданное на ассоциативной алгебре $\End_\mathbbm{k}(W)$ формулой~(\ref{EqHActionOnEnd}).
Это действие соответствует естественному рациональному
$G$-действию на $\End_\mathbbm{k}(W)$: \begin{equation}\label{EquationGEndAction}(g\psi)(w) = g\bigl(\psi({g^{-1}}w)\bigr)\text{ для всех }w \in W,\ \psi \in \End_\mathbbm{k}(W)\text{ и }g \in G.\end{equation}
  Отсюда в силу следствия~\ref{CorollaryGReductWedderburnMalcev}
  для любой
  $\mathbbm{k}G$-модульной подалгебры $A_1 \subseteq \End_\mathbbm{k}(W)$
  существует $\mathbbm{k}G$-инвариантное разложение Веддербёрна~"--- Мальцева.

   Заметим, что в силу кокоммутативности алгебры Хопфа $\mathbbm{k}G$
   алгебра Ли $\mathfrak{gl}(L)$ является $\mathbbm{k}G$-модульной,
   откуда её $\mathbbm{k}G$-модульные подалгебры Ли~"--- это в точности её 
   $G$-инвариантные подалгебры Ли.
  Для любой $G$-инвариантной подалгебры Ли $L_0 \subseteq \mathfrak{gl}(L)$
  пространство $L$ является $(G, L_0)$-модулем.
  Если пространство $L$
является вполне приводимым  $L_0$-модулем
без учёта $G$-действия, то в силу теоремы~\ref{TheoremAffAlgGrWeyl} пространство $L$ является вполне приводимым $(G,L_0)$-модулем.
  
  Следовательно, $L$ является $\mathbbm{k}G$-хорошей алгеброй Ли.
\end{proof}

\begin{example}\label{ExampleHniceAffAlgDiff}
Пусть $L$~"--- конечномерная алгебра Ли
   над алгебраически замкнутым полем $\mathbbm{k}$ характеристики $0$,
   на которой рационально действует автоморфизмами
   некоторая связная редуктивная аффинная алгебраическая группа $G$.
  Тогда алгебра Ли $L$ является $U(\mathfrak g)$-хорошей, где $\mathfrak{g}$~"--- алгебра Ли группы $G$.
\end{example}
\begin{proof} Из~(\ref{EqHActionOnEnd}) следует, что для всякого $U(\mathfrak{g})$-подмодуля $W \subseteq L$ действие группы $G$ на пространстве  $\End_\mathbbm{k}(W)$
задаётся формулой~\eqref{EquationGEndAction}, а $\mathfrak{g}$-действие~"--- формулой
$$(\delta \psi)(w) = \delta \psi(w) - \psi(\delta w)\text{ для всех }w \in W,\ \psi \in \End_\mathbbm{k}(W)\text{ и }\delta \in W.$$
Отсюда следует, что $\mathfrak{g}$-действие на $\End_\mathbbm{k}(W)$ также является
дифференциалом $G$-действия на $\End_\mathbbm{k}(W)$.
Поскольку в силу предложения~\ref{PropositionDerAutConnection}
алгебра Ли $\mathfrak g$ действует дифференцированиями и
оба действия имеют одинаковые инвариантные подпространства,
достаточно использовать пример~\ref{ExampleHniceAffAlg}.
\end{proof}

\begin{example}\label{ExampleHniceDiff}
Пусть $L$~"--- конечномерная алгебра Ли
   над алгебраически замкнутым полем $\mathbbm{k}$ характеристики $0$,
   на которой действует дифференцированиями конечномерная полупростая алгебра Ли $\mathfrak{g}$.
  Тогда алгебра Ли $L$ является $U(\mathfrak g)$-хорошей.
\end{example}
\begin{proof} В силу теоремы~\ref{TheoremLieDiffActionReplacement} 
действие алгебры Ли $\mathfrak{g}$ на $L$ дифференцированиями является
дифференциалом действия некоторой такой односвязной полупростой аффинной алгебраической
группы $G$ на $L$ автоморфизмами, что $\mathfrak{g}$~"--- алгебра Ли группы $G$.
Теперь
достаточно использовать пример~\ref{ExampleHniceAffAlgDiff}.
\end{proof}

Другим важным примером $H$-хорошей алгебры является пример~\ref{ExampleHniceGr},
который будет приведён ниже и окажется частным случаем примера~\ref{ExampleHniceAffAlg}.

\section{Основная теорема и её следствия}\label{SectionMainLie}

Следующая теорема является основным результатом данной главы:

\begin{theorem}\label{TheoremMainLieH}
Пусть $L$~"--- конечномерная $H$-хорошая алгебра Ли
над алгебраически замкнутым полем $\mathbbm{k}$ характеристики $0$,
где $H$~"--- некоторая алгебра Хопфа.
Тогда 
\begin{enumerate}
\item либо существует такое $n_0$, что $c_n^H(L)=0$ при всех $n\geqslant n_0$;
\item либо существуют такие константы $C_1, C_2 > 0$, $r_1, r_2 \in \mathbb R$,
  $d \in \mathbb N$, что $$C_1 n^{r_1} d^n \leqslant c^{H}_n(L)
   \leqslant C_2 n^{r_2} d^n\text{ для всех }n \in \mathbb N.$$
\end{enumerate}
   В частности, существует $\PIexp^H(L)\in\mathbb Z_+$ и, таким образом,
   для $c^H_n(L)$ справедлив аналог гипотезы Амицура.
\end{theorem}

\begin{theorem}\label{TheoremMainLieHSum}
Пусть $L=L_1 \oplus \ldots \oplus L_s$ (прямая сумма $H$-инвариантных идеалов)~"---
$H$-модульная алгебра Ли
над алгебраически замкнутым полем $\mathbbm{k}$ характеристики $0$,
где $H$~"--- некоторая алгебра Хопфа, а $L_i$~"--- $H$-хорошие алгебры.
Тогда существует $\PIexp^H(L)=\max_{1 \leqslant i \leqslant s}
\PIexp^H(L_i)$.
\end{theorem}

Теоремы~\ref{TheoremMainLieH} и~\ref{TheoremMainLieHSum} будут доказаны в конце~\S\ref{SectionLowerLie}.
Формула для $d=d(L,H)=\PIexp^H(L)$ будет дана в \S\ref{SectionHPIexpLie}.

Сформулируем теперь важнейшие следствия из теорем~\ref{TheoremMainLieH} и~\ref{TheoremMainLieHSum}:

\begin{theorem}\label{TheoremMainLieHSS}
Пусть $L$~"--- конечномерная $H$-модульная алгебра Ли
над произвольным полем $\mathbbm{k}$ характеристики $0$,
где $H$~"--- конечномерная полупростая алгебра Хопфа.
Тогда 
\begin{enumerate}
\item либо существует такое $n_0$, что $c_n^H(L)=0$ при всех $n\geqslant n_0$;
\item либо существуют такие константы $C_1, C_2 > 0$, $r_1, r_2 \in \mathbb R$,
  $d \in \mathbb N$, что $$C_1 n^{r_1} d^n \leqslant c^{H}_n(L)
   \leqslant C_2 n^{r_2} d^n\text{ для всех }n \in \mathbb N.$$
\end{enumerate}
   В частности, существует $\PIexp^H(L)\in\mathbb Z_+$ и, таким образом,
   для $c^H_n(L)$ справедлив аналог гипотезы Амицура.
\end{theorem}
\begin{theorem}\label{TheoremMainLieHSSSum}
Пусть $L=L_1 \oplus \ldots \oplus L_s$ (прямая сумма $H$-инвариантных идеалов)~"---
конечномерная $H$-модульная алгебра Ли
над произвольным полем $\mathbbm{k}$ характеристики $0$,
где $H$~"--- конечномерная полупростая алгебра Хопфа.
Тогда $\PIexp^H(L)=\max_{1 \leqslant i \leqslant s}
\PIexp^H(L_i)$.
\end{theorem}
\begin{proof}[Доказательство теорем~\ref{TheoremMainLieHSS} и~\ref{TheoremMainLieHSSSum}]
$H$-коразмерности не меняются при расширении основного поля.
Доказательство полностью повторяет соответствующее доказательство 
для коразмерностей обычных тождеств ассоциативных алгебр~\cite[теорема~4.1.9]{ZaiGia} и
алгебр Ли~\cite[\S 2]{ZaicevLie}.
Отсюда без ограничения общности можно считать основное поле $\mathbbm{k}$ алгебраически замкнутым.
Теперь теоремы~\ref{TheoremMainLieHSS} и~\ref{TheoremMainLieHSSSum}
получаются из теорем~\ref{TheoremMainLieH} и~\ref{TheoremMainLieHSum}
при помощи примера~\ref{ExampleHniceHSS}.
 \end{proof}

\begin{theorem}\label{TheoremMainLieGr}
Пусть $L$~"--- конечномерная алгебра Ли
над полем $\mathbbm{k}$ характеристики $0$,
градуированная произвольной группой $G$.
Тогда 
\begin{enumerate}
\item либо существует такое $n_0$, что $c_n^{G\text{-}\mathrm{gr}}(L)=0$ при всех $n\geqslant n_0$;
\item либо существуют такие константы $C_1, C_2 > 0$, $r_1, r_2 \in \mathbb R$,
  $d \in \mathbb N$, что $$C_1 n^{r_1} d^n \leqslant c^{G\text{-}\mathrm{gr}}_n(L)
   \leqslant C_2 n^{r_2} d^n\text{ для всех }n \in \mathbb N.$$
\end{enumerate}
   В частности, существует $\PIexp^{G\text{-}\mathrm{gr}}(L)\in\mathbb Z_+$ и, таким образом,
   для градуированных тождеств справедлив аналог гипотезы Амицура.
\end{theorem}
\begin{theorem}\label{TheoremMainLieGrSum}
Пусть $L=L_1 \oplus \ldots \oplus L_s$ (прямая сумма градуированных идеалов)~"---
конечномерная алгебра Ли
над полем $\mathbbm{k}$ характеристики $0$,
градуированная произвольной группой $G$.
Тогда $\PIexp^{G\text{-}\mathrm{gr}}(L)=\max_{1 \leqslant i \leqslant s}
\PIexp^{G\text{-}\mathrm{gr}}(L_i)$.
\end{theorem}
Теоремы~\ref{TheoremMainLieGr} и~\ref{TheoremMainLieGrSum} будут выведены
из теорем~\ref{TheoremMainLieH} и~\ref{TheoremMainLieHSum}
в \S\ref{SectionGrLie}.  

Алгебры, на которых некоторая группа действует не только автоморфизмами, но и антиавтоморфизмами
формально не являются модульными алгебрами над алгебрами Хопфа, однако в случае алгебр Ли
следующий приём позволяет свести такой случай к случаю, когда группа действует только автоморфизмами:

\begin{proposition}\label{PropositionGtoAut} 
Пусть $L$~"--- алгебра Ли, на которой некоторая группа $G$ действует автоморфизмами и антиавтоморфизмами.
Обозначим через~$G_0$ подгруппу группы~$G$,
состоящую из элементов, действующих на~$L$ автоморфизмами.
Определим на $L$ действие группы $\tilde G$, изоморфной группе $G$,
при помощи равенства $${\tilde g}a=\left\lbrace\begin{array}{rll}
 ga & \text{при} & g\in G_0, \\
 -ga & \text{при} & g\in G\backslash G_0, 
 \end{array}\right. $$
 где  $a\in L$.
(Через $\tilde g \in \tilde G$ обозначается элемент, соответствующий элементу $g\in G$
при изоморфизме $\tilde G \cong G$.)
 Тогда группа $\tilde G$ действует
 на $\tilde G$ автоморфизмами.
 Обратно, на всякой алгебре Ли $L$,
 на которой группа $\tilde G$
действует автоморфизмами, в случае, когда в группе~$G$ фиксирована подгруппа $G_0$ индекса ${}\leqslant 2$,
 можно определить такое действие группы $G$
автоморфизмами и антиавтоморфизмами, что
элементы подгруппы $G_0$ действуют автоморфизмами,
а элементы подмножества $G\backslash G_0$ действуют антиавтоморфизмами:
$$ga=\left\lbrace\begin{array}{rll}
 {\tilde g}a & \text{при} & g\in G_0, \\
 -{\tilde g}a & \text{при} & g\in G\backslash G_0
 \end{array}\right.$$
  для всех  $a\in L$.
  Кроме того, $c_n^{\tilde G}(L)=c_n^G(L)$ для всех $n \in \mathbb N$. Если $G$~"---
  аффинная алгебраическая группа, рационально действующая
  на $L$, а алгебра Ли $L$ конечномерна, то группа $\tilde G$ также действует на $L$ рационально.
  (Считаем, что на $\tilde G$ задана та же структура аффинного алгебраического
  многообразия, что и на~$G$.)
\end{proposition}
\begin{proof} 
Заметим, что для всех $\tilde g \in \tilde G$
и $a,b\in L$ справедливо равенство
$$\tilde g[a,b]=[\tilde g a, \tilde gb].$$
Это в точности означает, что группа $\tilde G$
действует на $L$ автоморфизмами.
Переход от $\tilde G$- к $G$-действию рассматривается аналогично.

Рассмотрим теперь изоморфизм свободных алгебр $\psi \colon L(X | G) \to L(X|\tilde G)$, заданный равенствами $$\psi\bigl(x^g\bigr) = \left\lbrace\begin{array}{rll}
 x^{\tilde g} & \text{при} & g\in G_0, \\
 -x^{\tilde g} & \text{при} & g\in G\backslash G_0, 
 \end{array}\right.
$$ где $x\in X$.
Тогда $\psi(\Id^G(L))=\Id^{\tilde G}(L)$, $\psi(V^G_n)=V^{\tilde G}_n$
и, следовательно, $c_n^{\tilde G}(L)=c_n^G(L)$ для всех $n \in \mathbb N$.

Если $G$~"--- аффинная алгебраическая группа, то $G_0$~"--- замкнутая подгруппа
индекса ${}\leqslant 2$, а $G$ является объединением неперескающихся замкнутых подмножеств $G_0$ и $G\backslash G_0$.
Обозначим через $\mathcal O(G)$ алгебру регулярных функций на $G$.
Пусть $$I_1 = \lbrace f \in \mathcal O(G) \mid f(g)=0 \text{ для всех } g \in G_0 \rbrace,$$
а $$I_2 = \lbrace f \in \mathcal O(G) \mid f(g)=0 \text{ для всех } g \in G\backslash G_0 \rbrace.$$
Поскольку $G_0 \cap (G\backslash G_0) = \varnothing$, в силу теоремы Гильберта о нулях
справедливо равенство $I_1 + I_2 = \mathcal O(G)$.
Отсюда для некоторых $f_1 \in I_1$ и $f_2 \in I_2$ справедливо равенство $1 = f_1 + f_2$. Следовательно, $f_2(g) - f_1(g) = 1$ для всех $g\in G_0$
и $f_2(g) - f_1(g)= -1$ для всех $g\in G\backslash G_0$.
Отсюда, для того, чтобы домножить каждый оператор из $G\backslash G_0$ на $(-1)$, 
достаточно домножить представление на $\bigl(f_2(g) - f_1(g)\bigr)$,
и оно по-прежнему останется рациональным.
\end{proof}

С учётом предложения~\ref{PropositionGtoAut}
и примера~\ref{ExampleHniceAffAlg}
теоремы~\ref{TheoremMainLieGAffAlg}
и~\ref{TheoremMainLieGAffAlgSum}, сформулированные ниже, оказываются немедленными следствиями теорем~\ref{TheoremMainLieH}
и~\ref{TheoremMainLieHSum}:

\begin{theorem}\label{TheoremMainLieGAffAlg}
Пусть $L$~"--- конечномерная алгебра Ли
над полем $\mathbbm{k}$ характеристики~$0$
с рациональным действием некоторой редуктивной аффинной алгебраической
группы $G$ автоморфизмами и антиавтоморфизмами.
Тогда 
\begin{enumerate}
\item либо существует такое $n_0$, что $c_n^{G}(L)=0$ при всех $n\geqslant n_0$;
\item либо существуют такие константы $C_1, C_2 > 0$, $r_1, r_2 \in \mathbb R$,
  $d \in \mathbb N$, что $$C_1 n^{r_1} d^n \leqslant c^{G}_n(L)
   \leqslant C_2 n^{r_2} d^n\text{ для всех }n \in \mathbb N.$$
\end{enumerate}
   В частности, существует $\PIexp^G(L)\in\mathbb Z_+$ и, таким образом,
   для $G$-тождеств справедлив аналог гипотезы Амицура.
\end{theorem}
\begin{theorem}\label{TheoremMainLieGAffAlgSum}
Пусть $L=L_1 \oplus \ldots \oplus L_s$ (прямая сумма $G$-инвариантных идеалов)~"---
конечномерная алгебра Ли
над полем $\mathbbm{k}$ характеристики $0$
с рациональным действием некоторой редуктивной аффинной алгебраической
группы $G$ автоморфизмами и антиавтоморфизмами.
Тогда $\PIexp^{G}(L)=\max_{1 \leqslant i \leqslant s}
\PIexp^{G}(L_i)$.
\end{theorem}

При этом, как и в случае $G$-тождеств 
ассоциативных алгебр (см.  теорему~\ref{TheoremAssGPIexpEqual}),
если группа $G$ связная, $G$-PI-экспонента совпадает с обычной
(см. теорему~\ref{TheoremLieGConnPIexpEqual} ниже).

Следующие две теоремы можно было бы вывести из
теорем~\ref{TheoremMainLieGAffAlg} и~\ref{TheoremMainLieGAffAlgSum},
пользуясь тем, что всякая конечная группа является редуктивной аффинной алгебраической группой
и всякое её конечномерное представление рационально, однако мы выведем их
из теорем~\ref{TheoremMainLieHSS} и~\ref{TheoremMainLieHSSSum}:

\begin{theorem}\label{TheoremMainLieGFin}
Пусть $L$~"--- конечномерная алгебра Ли
над полем $\mathbbm{k}$ характеристики $0$
с действием некоторой конечной
группы $G$ автоморфизмами и антиавтоморфизмами.
Тогда 
\begin{enumerate}
\item либо существует такое $n_0$, что $c_n^{G}(L)=0$ при всех $n\geqslant n_0$;
\item либо существуют такие константы $C_1, C_2 > 0$, $r_1, r_2 \in \mathbb R$,
  $d \in \mathbb N$, что $$C_1 n^{r_1} d^n \leqslant c^{G}_n(L)
   \leqslant C_2 n^{r_2} d^n\text{ для всех }n \in \mathbb N.$$
\end{enumerate}
   В частности, существует $\PIexp^G(L)\in\mathbb Z_+$ и, таким образом,
   для $G$-тождеств справедлив аналог гипотезы Амицура.
\end{theorem}
\begin{theorem}\label{TheoremMainLieGFinSum}
Пусть $L=L_1 \oplus \ldots \oplus L_s$ (прямая сумма $G$-инвариантных идеалов)~"---
конечномерная алгебра Ли
над полем $\mathbbm{k}$ характеристики $0$
с действием некоторой конечной
группы $G$ автоморфизмами и антиавтоморфизмами.
Тогда $\PIexp^{G}(L)=\max_{1 \leqslant i \leqslant s}
\PIexp^{G}(L_i)$.
\end{theorem}
\begin{proof}[Доказательство теорем~\ref{TheoremMainLieGFin} и~\ref{TheoremMainLieGFinSum}]
В силу теоремы Машке алгебра Хопфа $\mathbbm{k}G$ полупроста.
Теперь достаточно применить теоремы~\ref{TheoremMainLieHSS} и~\ref{TheoremMainLieHSSSum}
для $H=\mathbbm{k}G$.
\end{proof}

 Как и в случае дифференциальных тождеств ассоциативных алгебр (см.  теорему~\ref{TheoremAssDerPIexpEqual}),
 при действии конечномерной полупростой алгебры Ли
 дифференцированиями дифференциальная PI-экспонента совпадает с обычной:

\begin{theorem}\label{TheoremMainDiffLie}
Пусть $L$~"--- конечномерная алгебра Ли
над полем $\mathbbm{k}$ характеристики $0$,
на которой действует дифференцированиями
некоторая конечномерная полупростая алгебра Ли~$\mathfrak g$.
Обозначим через $d:=\PIexp(L)$ экспоненту роста коразмерностей обычных тождеств алгебры Ли $L$.
Тогда 
\begin{enumerate}
\item при $d=0$ существует такое $n_0$, что $c_n^{U(\mathfrak g)}(L)=0$ для всех $n\geqslant n_0$;
\item при $d > 0$ существуют такие константы $C_1, C_2 > 0$, $r_1, r_2 \in \mathbb R$, что $$C_1 n^{r_1} d^n \leqslant c^{U(\mathfrak g)}_n(L)
   \leqslant C_2 n^{r_2} d^n\text{ для всех }n \in \mathbb N.$$
\end{enumerate}
   В частности, существует $\PIexp^{U(\mathfrak g)}(L)\in\mathbb Z_+$ и, таким образом,
   для дифференциальных тождеств справедлив аналог гипотезы Амицура.
\end{theorem}

Теорема~\ref{TheoremMainDiffLie} будет доказана в \S\ref{SectionDiffLie}.

\section{Формулы для $H$-PI-экспоненты}\label{SectionHPIexpLie}

Предположим, что $L$~"--- $H$-хорошая алгебра Ли.

Для всевозможных $r \in \mathbb Z_+$ рассмотрим такие семейства $H$-инвариантных идеалов $I_1, I_2, \ldots, I_r$, $J_1, J_2, \ldots, J_r$ алгебры Ли $L$, что $J_k \subseteq I_k$
и выполнены следующие условия:
\begin{enumerate}
\item всякий фактормодуль $I_k/J_k$ является неприводимым $(H,L)$-модулем;
\item для любых таких $H$-инвариантных $B$-подмодулей $T_k$,
что $I_k = J_k\oplus T_k$, существуют такие числа
$q_i \geqslant 0$, что $$\bigl[[T_1, \underbrace{L, \ldots, L}_{q_1}], [T_2, \underbrace{L, \ldots, L}_{q_2}], \ldots, [T_r,
 \underbrace{L, \ldots, L}_{q_r}]\bigr] \ne 0.$$
\end{enumerate}

Напомним, что через $\Ann(M)$ обозначается аннулятор модуля $M$.
Положим $$d(L, H) := \max \left(\dim \frac{L}{\Ann(I_1/J_1) \cap \ldots \cap \Ann(I_r/J_r)}
\right),$$ где максимум берётся по всем $r \in \mathbb Z_+$ и всем $H$-инвариантным
идеалам $I_1, \ldots, I_r$, $J_1, \ldots, J_r$,
для которых выполнены условия 1--2, приведённые выше. Утверждается, что $\PIexp^H(L)=d(L, H)$.
Соответственно, теорема~\ref{TheoremMainLieH} будет доказываться для $d=d(L,H)$.

Для доказательства равенства дифференциальной PI-экспоненты с обычной PI-экспонентой
в случае действия конечномерной полупростой алгебры Ли дифференцированиями, а также аналогичного результата для $G$-PI-экспоненты в случае рационального действия связной редуктивной аффинной алгебраической
группы автоморфизмами (напомним, что в силу предложений~\ref{PropositionConnectedAffAlgAutNoAnti}
и~\ref{PropositionGtoAut} случай антиавтоморфизмов можно не рассматривать),
нам потребуется ещё одна формула для $H$-PI-экспоненты.

Докажем сперва следующую лемму:

\begin{lemma}\label{LemmaLBQN} Пусть $L$~"--- $H$-хорошая
алгебра Ли. Рассмотрим $H$-инвариантное разложение Леви
$L=B\oplus R$  из условия~\ref{ConditionLevi}
\S\ref{SectionHnice} и ограничение присоединённого
представления $\ad \colon L \to \mathfrak{gl}(L)$
на полупростую алгебру Ли $B$.
Тогда алгебра Ли $L$ является вполне приводимым $(H,B)$-модулем.
Более того, существует такой $(H,B)$-подмодуль $Q \subseteq R$,
что $L=B\oplus Q\oplus N$ (прямая сумма $(H,B)$-подмодулей)
и $[B,Q]=0$.
\end{lemma}
\begin{proof}
Поскольку $\ad\colon L \to \mathfrak{gl}(L)$
является гомоморфизмом алгебр Ли и $H$-модулей,
алгебра Ли $(\ad B)$, будучи гомоморфным образом полупростой $H$-модульной алгебры Ли,
сама является полупростой $H$-модульной алгеброй Ли.
Отсюда в силу теоремы Вейля алгебра Ли~$L$
является вполне приводимым $(\ad B)$-модулем.
Следовательно, в силу условия~\ref{ConditionLComplHred} из \S\ref{SectionHnice}
алгебра Ли~$L$ является вполне приводимым $(H, \ad B)$-модулем
и существует такой $(H,B)$-подмодуль $Q \subseteq R$, что
$R=Q\oplus N$, т.е. $L=B\oplus Q \oplus N$.
При этом в силу предложения 2.1.7 из~\cite{GotoGrosshans}
справедливо включение $[L, R] \subseteq N$.
Отсюда $[B, Q] \subseteq [L, R] \subseteq N$,
в то время как $[B,Q]\subseteq Q$. Следовательно, $[B, Q]=0$.
\end{proof}

Рассмотрим теперь ассоциативную подалгебру $A_0$ алгебры $\End_\mathbbm{k}(L)$,
порождённую подпространством $\ad Q$. Заметим, что $A_0$ является $H$-модульной
алгеброй, поскольку $Q$ является $H$-подмодулем. В силу условия~\ref{ConditionWedderburn} из \S\ref{SectionHnice} существует разложение
$A_0 = \tilde A_0 \oplus J(A_0)$ (прямая сумма $H$-подмодулей),
где $\tilde A_0$~"--- максимальная полупростая подалгебра алгебры $A_0$.

Докажем следующее утверждение:
\begin{lemma}\label{LemmaA0DirectSumOfFields}
Существует такое число $q\in\mathbb N$, что
$$\tilde A_0=\mathbbm{k}e_1 \oplus \ldots \oplus \mathbbm{k}e_q \text{ (прямая сумма идеалов)} $$
    для некоторых идемпотентов $e_i \in A_0$.
\end{lemma}
\begin{proof}
Поскольку разрешимый радикал $R$ является разрешимой алгеброй Ли,
в силу теоремы Ли в $L$  существует такой базис,
что в этом базисе матрицы всех операторов $\ad a$, где $a\in R$, являются
верхнетреугольными.
Обозначим соответствующий изоморфизм алгебр $\End_\mathbbm{k}(L) \to M_s(\mathbbm{k})$ через $\psi$, где $s := \dim L$.
Поскольку $\psi(\ad R) \subseteq \UT_s(\mathbbm{k})$,
справедливо включение $\psi(A_0) \subseteq \UT_s(\mathbbm{k})$, где $\UT_s(\mathbbm{k})$~"--- ассоциативная алгебра верхнетреугольных матриц размера $s\times s$. Однако $$\UT_s(\mathbbm{k}) = \mathbbm{k}e_{11}\oplus \mathbbm{k}e_{22}\oplus
 \ldots\oplus \mathbbm{k}e_{ss}\oplus \tilde N,$$
 где $$\tilde N := \langle e_{ij} \mid 1 \leqslant i < j \leqslant s \rangle_\mathbbm{k}$$
 является нильпотентным идеалом. Поскольку $\psi$~"--- изоморфизм, алгебра $A_0$
 не содержит подалгебр, изоморфных $M_2(\mathbbm{k})$, откуда
  $\tilde A_0=\mathbbm{k}e_1 \oplus \ldots \oplus \mathbbm{k}e_q$ (прямая сумма идеалов)
    для некоторых идемпотентов $e_i \in A_0$.
\end{proof}

В силу того, что $[B,Q]=0$, а элементы $e_i$ являются многочленами от элементов $\ad a$, где $a \in Q$, 
справедливо равенство $[\ad B, \tilde A_0]=0$.
Из полупростоты алгебры Ли $B$ следует, что $(\ad B) \cap \tilde A_0 = 0$.
Будем рассматривать $(\ad B)\oplus \tilde A_0$ как $H$-модульную алгебру Ли.

\begin{lemma}\label{LemmaLHBA0ComplReducible}
Алгебра Ли $L$ является вполне приводимым $(\ad B)\oplus \tilde A_0$- и $(H, (\ad B)\oplus \tilde A_0)$-модулем.
\end{lemma}
\begin{proof}
Заметим, что $e_i$
являются коммутирующими диагонализуемыми на $L$ операторами.
Отсюда для них существует общий базис из собственных векторов и
 $L=\bigoplus_j W_j$, где $W_j$~"--- пересечения собственных подпространств операторов $e_i$.
 Каждый из операторов $e_i$ коммутирует с операторами из $\ad B$. Следовательно, $W_j$ являются $(\ad B)$-подмодулями. В силу того, что алгебра Ли $B$ полупроста, каждое из подпространств $W_j$
является прямой суммой неприводимых $(\ad B)$-подмодулей.
Поскольку операторы $e_i$ действуют на каждом из $W_j$
как скалярные операторы,  алгебра Ли $L$
является прямой суммой неприводимых $(\ad B)\oplus \tilde A_0$-подмодулей. 
Теперь утверждение леммы следует из условия~\ref{ConditionLComplHred} \S\ref{SectionHnice}.
\end{proof}

Теперь мы готовы дать альтернативную формулу для PI-экспоненты, в которой условие~2 из формулы для $d(L,H)$
заменено на более слабое условие~2'.

Пусть $L$~"--- $H$-хорошая алгебра Ли. Для всевозможных $r \in \mathbb Z_+$ рассмотрим такие семейства $H$-инвариантных идеалов $I_1, I_2, \ldots, I_r$, $J_1, J_2, \ldots, J_r$ алгебры Ли $L$, что $J_k \subseteq I_k$ и выполнены следующие условия:
\begin{enumerate}
\item всякий фактормодуль $I_k/J_k$ является неприводимым $(H,L)$-модулем;
\item[2'.] существуют такие $H$-инвариантные $(\ad B)\oplus \tilde A_0$-подмодули $\tilde T_k$, где $I_k = J_k\oplus \tilde T_k$,
 и числа
$q_i \geqslant 0$, что $$\bigl[[\tilde T_1, \underbrace{L, \ldots, L}_{q_1}], [\tilde T_2, \underbrace{L, \ldots, L}_{q_2}], \ldots, [\tilde T_r,
 \underbrace{L, \ldots, L}_{q_r}]\bigr] \ne 0.$$
\end{enumerate}

Положим $$d'(L, H) := \max \left(\dim \frac{L}{\Ann(I_1/J_1) \cap \ldots \cap \Ann(I_r/J_r)}
\right),$$
где максимум берётся по всем $r \in \mathbb Z_+$ и всем наборам $I_1, \ldots, I_r$, $J_1, \ldots, J_r$,
удовлетворяющим условиям~1 и~2'.

Как было только что доказано в лемме~\ref{LemmaLHBA0ComplReducible},
алгебра Ли $L$ является вполне приводимым $(H,(\ad B)\oplus \tilde A_0)$-модулем,
откуда всегда можно выбрать такие $H$-инвариантные $(\ad B)\oplus \tilde A_0$-подмодули $\tilde T_k$,
что $I_k = J_k \oplus \tilde T_k$. Следовательно,
$d'(L, H) \geqslant d(L, H)$.

Отсюда для доказательства равенства $\PIexp^H(L)=d'(L, H)=d(L,H)$
достаточно доказать оценку сверху для коразмерностей с использованием числа $d(L,H)$,
а оценку снизу~"--- с использованием числа $d'(L,H)$.
Именно это и будет сделано в следующих параграфах.

\section{Полиномиальные $H$-тождества представлений и кососимметрические $H$-многочлены}
\label{SectionAltLie}

В данном параграфе доказываются вспомогательные утверждения, которые будут затем использованы для получения оценки снизу.

Пусть $L$~"--- $H$-модульная алгебра Ли для некоторой алгебры Хопфа $H$, $M$~"--- некоторый $(H,L)$-модуль, а $\varphi \colon  L \to \mathfrak{gl}(M)$~"--- соответствующее представление.
Ассоциативный $H$-многочлен $f(x_1, \ldots, x_n)\in \mathbbm{k}\langle X | H \rangle$ называется \textit{ полиномиальным $H$-тождеством} представления $\varphi$, если $f(\varphi(a_1), \ldots, \varphi(a_n))=0$ 
для всех $a_i \in L$.
 Другими словами, $f$~"--- $H$-тождество представления $\varphi$,
 если для любого гомоморфизма $H$-модульных алгебр $\psi \colon 
 \mathbbm{k}\langle X | H \rangle \to \End_\mathbbm{k}(M)$, такого, что $\psi(X)\subseteq \varphi(L)$,
 справедливо равенство $\psi(f)=0$.  
  Множество $\Id^{H}(\varphi)$ всех $H$-тождеств представления $\varphi$ является $H$-инвариантным
  двусторонним идеалом алгебры $\mathbbm{k}\langle X | H \rangle$.

Лемма~\ref{LemmaLieAlternateFirst} является аналогом леммы 1 из~\cite{GiaSheZai}
и леммы~\ref{LemmaAssocAlternateFirst}. 

\begin{lemma}\label{LemmaLieAlternateFirst}
Пусть $M$~"--- некоторый конечномерный неприводимый $(H,L)$-модуль, точный как $L$-модуль,
где $L$~"--- $H$-модульная алгебра Ли, а $H$~"--- алгебра Хопфа
над алгебраически замкнутым полем $\mathbbm{k}$ характеристики $0$. 
Обозначим через $\varphi \colon L \to \mathfrak{gl}(M)$
и $\zeta \colon H \to \End_\mathbbm{k}(M)$ соответствующие гомоморфизмы.
Пусть $(\zeta(\gamma_j))_{j=1}^m$~"--- базис алгебры $\zeta(H)$.
Тогда для некоторого $n\in\mathbb N$ существуют такие $H$-многочлены
$f_j \in P^H_n$, где $1 \leqslant j \leqslant m$,
  кососимметричные по переменным каждого из множеств $\lbrace x_1, \ldots, x_\ell \rbrace$  и $\lbrace y_1, \ldots, y_\ell \rbrace \subseteq \lbrace x_{\ell+1}, \ldots, x_n \rbrace$, где $\ell:=\dim L$,
  что для некоторых $\bar x_i \in L$ справедливо равенство $$\sum_{j=1}^m{f_j(\varphi(\bar x_1), \ldots,  \varphi(\bar x_n))\, \zeta(\gamma_j)}=\id_M.$$
В частности, существует $H$-многочлен $f \in P^H_n \backslash \Id^H(\varphi)$,
кососимметричный по $\lbrace x_1, \ldots, x_\ell \rbrace$  и по $\lbrace y_1, \ldots, y_\ell \rbrace \subseteq \lbrace x_{\ell+1}, \ldots, x_n \rbrace$.
\end{lemma}
\begin{proof}
Поскольку $(H,L)$-модуль $M$ неприводим, в силу теоремы плотности
 алгебра $\End_\mathbbm{k}(M) \cong M_q(\mathbbm{k})$ порождена операторами из $\zeta(H)$
 и $\varphi(L)$. Здесь $q := \dim M$.
 Рассмотрим многочлен Регева
\begin{equation*}\begin{split}
 \hat f(x_1, \ldots, x_{q^2}; y_1, \ldots, y_{q^2})
:=\sum_{\substack{\sigma \in S_q, \\ \tau \in S_q}} (\sign(\sigma\tau))
x_{\sigma(1)}\ y_{\tau(1)}\ x_{\sigma(2)}x_{\sigma(3)}x_{\sigma(4)}
\ y_{\tau(2)}y_{\tau(3)}y_{\tau(4)}\ldots
\cdot
\\ \cdot
 x_{\sigma\left(q^2-2q+2\right)}\ldots x_{\sigma\left(q^2\right)}
\ y_{\tau\left(q^2-2q+2\right)}\ldots y_{\tau\left(q^2\right)}.
\end{split}\end{equation*}
Как уже было отмечено в главе~\ref{ChapterGenHAssocCodim},
этот многочлен является центральным многочленом 
для $M_q(\mathbbm{k})$, т.е. $\hat f$ не является полиномиальным тождеством для $M_q(\mathbbm{k})$,
а его значения принадлежат центру алгебры $M_q(\mathbbm{k})$. (См., например, \cite[теорема~5.7.4]{ZaiGia}.)
Поскольку многочлен $\hat f$ кососимметричный,
его значение является ненулевым скалярным оператором при подстановке всех элементов любого
базиса вместо переменных каждого из множеств  
$\lbrace x_1, \ldots, x_\ell \rbrace$ и $\lbrace y_1, \ldots, y_\ell \rbrace$.

Пусть $a_1, \ldots, a_\ell$~"--- базис алгебры Ли $L$.
Напомним, что в произведении элементов пространств
$\varphi(L)$ и $\zeta(H)$ всегда можно переместить
элементы пространства $\zeta(H)$
вправо, используя равенство~(\ref{EqHLModule2}).
Следовательно,
$\varphi(a_1), \ldots, \varphi(a_\ell)$, $\left(\varphi\left(a_{i_{11}}\right)
\ldots \varphi\left(a_{i_{1,m_1}}\right)\right)\zeta(h_1)$,
\ldots, $\left(\varphi\left(a_{i_{r,1}}\right)
\ldots \varphi\left(a_{i_{r,m_r}}
\right)\right)\zeta(h_r)$ является базисом пространства $\End_\mathbbm{k}(M)$ для подходящих $i_{jk} \in \lbrace 1,2, \ldots, \ell \rbrace$, $h_j \in H$, поскольку алгебра $\End_\mathbbm{k}(M)$ порождена операторами из $\zeta(H)$
и $\varphi(L)$. Заменим в $\hat f$ переменные $x_{\ell+j}$ на $z_{j1}
z_{j2} \ldots z_{j,m_j} \zeta(h_j)$, а переменные $y_{\ell+j}$~"--- на $z'_{j1}
z'_{j2} \ldots z'_{j,m_j} \zeta(h_j)$
и обозначим получившееся выражение через $\tilde f$. Используя~(\ref{EqHLModule2}) ещё раз,
переместим в $\tilde f$ все $\zeta(h)$, где $h \in H$, вправо и представим
 $\tilde f$ в виде суммы $\sum_{j=1}^m{f_j\, \zeta(\gamma_j)}$,
 где $f_j \in P^H_{2\ell+2\sum_{i=1}^r m_i}$ являются
 $H$-многочленами, кососимметричными по $x_1, \ldots, x_\ell$  и по $y_1, \ldots, y_\ell$. 
Заметим, что $\tilde f$ становится ненулевым скалярным оператором на~$M$
при подстановке $x_i=y_i=\varphi(a_i)$ при $ 1 \leqslant i \leqslant \ell$
и $z_{jk}=z_{jk}'=\varphi(a_{i_{jk}})$ при $1 \leqslant j \leqslant r$, $1 \leqslant k \leqslant m_j$. Поделив все $f_j$ на соответствующий скаляр, мы получим требуемые многочлены.
 В частности, по крайней мере для одного из $1 \leqslant j \leqslant m$ 
выполнено условие 
 $f_j \notin \Id^H(\varphi)$, т.е. можно положить $f:=f_j$.
\end{proof}

\begin{lemma}\label{LemmaTwoColumns}
Пусть $\alpha_1, \alpha_2, \ldots, \alpha_q$, $\beta_1, \ldots, \beta_q \in \mathbbm{k}$,
где $\mathbbm{k}$~"--- бесконечное поле, $
1 \leqslant k \leqslant q$,
 $\alpha_i\ne 0$ для всех $1 \leqslant i < k$,
$\alpha_k=0$, а $\beta_k\ne 0$. Тогда существует
такое $\gamma \in \mathbbm{k}$, что $\alpha_i + \gamma \beta_i \ne 0$
для всех $1 \leqslant i \leqslant k$.
\end{lemma}
\begin{proof} Достаточно выбрать
$\gamma \notin \left\lbrace -\frac{\alpha_1}{\beta_1},
\ldots, -\frac{\alpha_{k-1}}{\beta_{k-1}}, 0\right\rbrace$. 
Это возможно сделать, поскольку поле $\mathbbm{k}$ бесконечно.
\end{proof}

\begin{lemma}\label{LemmaS}
Пусть $L=B \oplus Z(L)$~"--- редуктивная $H$-модульная
алгебра Ли, где $H$~"--- алгебра Хопфа
над алгебраически замкнутым полем $\mathbbm{k}$ характеристики $0$, $B$~"---
$H$-инвариантная максимальная полупростая подалгебра, а
$Z(L)$~"--- центр алгебры Ли $L$ с базисом $r_1$, $r_2$, \ldots, $r_t$.
Пусть $M$~"--- конечномерный неприводимый $(H,L)$-модуль, точный как $L$-модуль, 
а $\varphi \colon L \to \mathfrak{gl}(M)$
и $\zeta \colon H \to \End_\mathbbm{k}(M)$~"--- соответствующие гомоморфизмы.
Предположим, что для всех $H$-инвариантных ассоциативных
подалгебр в $\End_\mathbbm{k}(M)$ существует $H$-инвариантное разложение Веддербёрна~"--- Мальцева.
Тогда существует такой $H$-многочлен $f \in P^H_t$,
кососимметричный по переменным $x_1, x_2, \ldots, x_t$,
что $f(\varphi(r_1), \ldots, \varphi(r_t))$ является на~$M$ невырожденным оператором.
\end{lemma}
\begin{proof} В силу леммы~\ref{LemmaAnnHLmoduleHMod} центр $Z(L)$ является $H$-подмодулем.
Согласно лемме~\ref{LemmaRedIrr} существует разложение
$M=M_1\oplus\ldots \oplus M_q$, где $M_j$~"--- $L$-подмодули,
а элементы $r_i$ действуют на каждом из $M_j$ как скалярные операторы.
Докажем, что для любого $j$ существует
такой
$H$-многочлен $f_j \in P^H_t$,  кососимметричный по $x_1, x_2, \ldots, x_t$, что
$f_j(\varphi(r_1), \ldots, \varphi(r_t))$
умножает каждый элемент пространства~$M_j$ на ненулевой элемент поля.

В силу леммы~\ref{LemmaLieAlternateFirst} существует такое $n\in\mathbb N$
и $H$-многочлены $\hat f_k \in P^H_n$, кососимметричные по переменным из множества
$\lbrace x_1, x_2, \ldots, x_\ell\rbrace$, где $\ell := \dim L$,
что
\begin{equation}\label{EqSIdM}
\sum_{k=1}^m{\hat f_k(\varphi(r_1), \varphi(r_2), \ldots, \varphi(r_t), \varphi(\bar x_{t+1}), \varphi(\bar x_{t+2}), \ldots,  \varphi(\bar x_n))\, \zeta(\gamma_k)}=\id_M
\end{equation}
  для некоторых $\bar x_i \in L$. (Вместо переменных множества $\lbrace x_1, x_2, \ldots, x_\ell\rbrace$
  можно подставлять элементы любого базиса пространства $L$,
поскольку $H$-многочлен $\hat f_k$ является кососимметричным по первым $\ell$ переменным.)
 Заметим, что $[\varphi(r_i), \varphi(a)]=0$ для всех $a \in L$ и $1 \leqslant i
\leqslant t$. Следовательно, можно переместить все элементы $\varphi(r_i)$
влево и переписать~(\ref{EqSIdM}) как
$$\sum_k \tilde f_k(\varphi(r_1), \varphi(r_2), \ldots, \varphi(r_t)) b_k = \id_M,$$
где $\tilde f_k \in P^H_t$ являются $H$-многочленами,
кососимметричными по переменным множества $\lbrace x_1, x_2, \ldots, x_t\rbrace$, а $b_k \in \End_\mathbbm{k}(M)$.
Следовательно, для любого $1 \leqslant j \leqslant q$ существует такое $k$, что  $\tilde f_k(\varphi(r_1), \varphi(r_2), \ldots, \varphi(r_t))\bigr|_{M_j}\ne 0$.
Поскольку $\tilde f_j(\varphi(r_1), \varphi(r_2), \ldots, \varphi(r_t))$ является
на $M_j$ скалярным оператором, достаточно взять $f_j := \tilde f_k$.

Теперь докажем по индукции, что для любого $1\leqslant k \leqslant q$
существует такой
$H$-многочлен $g_k =\lambda_1 f_1 +\ldots + \lambda_k f_k$, где $\lambda_i \in \mathbbm{k}$,
что $g_k(\varphi(r_1), \ldots, \varphi(r_t))$ действует на каждом из $M_i$, где $1\leqslant i \leqslant k$, как ненулевой скалярный оператор.

База индукции $k=1$ очевидна, поскольку достаточно положить $\lambda_1 := 1$.
Предположим теперь существование $H$-многочлена $g_{k-1}$. Выведем отсюда
существование $H$-многочлена $g_k$.
Обозначим через $\alpha_j$ скаляры, умножением на которые действует
оператор $g_{k-1}(\varphi(r_1), \ldots, \varphi(r_t))$ на $M_j$, а через $\beta_j$~"--- скаляры, умножением на которые действует оператор $f_k(\varphi(r_1), \ldots, \varphi(r_t))$ на $M_j$, где
$1 \leqslant j \leqslant k$. Согласно предположению индукции
$\alpha_j \ne 0$ при $1 \leqslant j < k$, а в силу свойств $H$-многочлена $f_k$ справедливо равенство  $\beta_k\ne 0$. Если $\alpha_k\ne 0$, то достаточно
положить $g_k:=g_{k-1}$, а если $\alpha_k= 0$, то достаточно
положить $g_k:=g_{k-1} + \gamma f_k$, где $\gamma \in \mathbbm{k}$
берётся из леммы~\ref{LemmaTwoColumns}.

Таким образом, по индукции следует существование $H$-многочленов $g_1, \ldots, g_q$
с заданными свойствами.
Теперь достаточно положить $f:=g_k$.
\end{proof}

Пусть $k\ell \leqslant n$, где $k,\ell,n \in \mathbb N$~"--- некоторые числа.
Как и в \S\ref{SectionAssocAlt}, обозначим через $Q^H_{\ell,k,n} \subseteq P^H_n$
подпространство, состоящее из всех многочленов, кососимметричных
по каждому из $k$ попарно непересекающихся наборов переменных $\{x^i_1, \ldots, x^i_\ell \}
\subseteq \lbrace x_1, x_2, \ldots, x_n\rbrace$, где $1 \leqslant i \leqslant k$.

Теорема~\ref{TheoremLieAlternateFinal}
 является аналогом теоремы 1 из~\cite{GiaSheZai} и теоремы~\ref{TheoremAssocAlternateFinal}.

\begin{theorem}\label{TheoremLieAlternateFinal}
Пусть $L=B \oplus Z(L)$~"--- конечномерная редуктивная $H$-модульная
алгебра Ли, где $H$~"--- алгебра Хопфа
над алгебраически замкнутым полем $\mathbbm{k}$ характеристики~$0$, $B$~"---
$H$-инвариантная максимальная полупростая подалгебра, $\ell := \dim L$.
Пусть $M$~"--- конечномерный неприводимый $(H,L)$-модуль, точный как $L$-модуль,
а $\varphi \colon L \to \mathfrak{gl}(M)$
и $\zeta \colon H \to \End_\mathbbm{k}(M)$~"--- соответствующие гомоморфизмы.
Предположим, что либо $Z(L)=0$, либо для всех $H$-инвариантных ассоциативных
подалгебр в $\End_\mathbbm{k}(M)$ существует $H$-инвариантное разложение Веддербёрна~"--- Мальцева.
Тогда существует такое $T \in \mathbb Z_+$, что для любого $k \in \mathbb N$
существует $f \in Q^H_{\ell, 2k, 2k\ell+T} \backslash \Id^H(\varphi)$.
\end{theorem}
\begin{proof} 
Пусть $f_1=f_1(x_1,\ldots, x_\ell,\ y_1,\ldots, y_\ell,
z_1, \ldots, z_T)$~"--- $H$-многочлен $f$ из леммы~\ref{LemmaLieAlternateFirst},
кососимметричный по $x_1,\ldots, x_\ell$ и по $y_1,\ldots, y_\ell$.
Поскольку $f_1 \in Q^H_{\ell, 2, 2\ell+T} \backslash \Id^H(\varphi)$,
 можно считать, что $k > 1$. Заметим, что
\begin{equation*}\begin{split}
f^{(1)}_1(u_1, v_1, x_1, \ldots, x_\ell,\ y_1,\ldots, y_\ell,
z_1, \ldots, z_T) := \\ =
\sum^\ell_{i=1} f_1(x_1, \ldots, [u_1, [v_1, x_i]],  \ldots, x_\ell,\ y_1,\ldots, y_\ell,
z_1, \ldots, z_T)\end{split}\end{equation*}
кососимметричен по $x_1,\ldots, x_\ell$ и по $y_1,\ldots, y_\ell$,
причём \begin{equation*}\begin{split}
f^{(1)}_1(\bar u_1, \bar v_1, \bar x_1, \ldots, \bar x_\ell,\
\bar y_1,\ldots, \bar y_\ell,
\bar z_1, \ldots, \bar z_T) =\\=
 \tr(\ad_{\varphi(L)} \bar u_1 \ad_{\varphi(L)} \bar v_1)
f_1(\bar x_1, \bar x_2, \ldots, \bar x_\ell,\ \bar y_1,\ldots, \bar y_\ell,
\bar z_1, \ldots, \bar z_T)
\end{split}\end{equation*}
для любой подстановки элементов из пространства $\varphi(L)$,
поскольку мы можем считать, что $\bar x_1, \ldots, \bar x_\ell$~"--- различные элементы базиса.

Пусть \begin{equation*}\begin{split}
f^{(j)}_1(u_1, \ldots, u_j, v_1, \ldots, v_j, x_1, \ldots, x_\ell,\ y_1,\ldots, y_\ell,
z_1, \ldots, z_T) :=\\=
\sum^\ell_{i=1} f^{(j-1)}_1(u_1, \ldots,  u_{j-1}, v_1, \ldots, v_{j-1},
 x_1, \ldots, [u_j, [v_j, x_i]],  \ldots, x_\ell,\ y_1,\ldots, y_\ell,
z_1, \ldots, z_T),\end{split}\end{equation*}
$2 \leqslant j \leqslant s$, $s := \dim B$. Заметим, что
если вместо $u_i$ или $v_i$
подставляются элементы из $\varphi(Z(L))$,
то $H$-многочлены $f^{(j)}_1$ обращаются в нуль, поскольку $Z(L)$ является центром алгебры Ли $L$.
Снова получаем, что
\begin{equation}\begin{split}\label{EqKilling}
f^{(j)}_1(\bar u_1, \ldots, \bar u_j, \bar v_1, \ldots, \bar v_j, \bar x_1, \ldots, \bar x_\ell,\ \bar y_1,\ldots, \bar y_\ell, \bar z_1, \ldots, \bar z_T) =\\= \tr(\ad_{\varphi(L)} \bar u_1 \ad_{\varphi(L)} \bar v_1)
 \tr(\ad_{\varphi(L)} \bar u_2 \ad_{\varphi(L)} \bar v_2)
 \ldots
 \tr(\ad_{\varphi(L)} \bar u_j \ad_{\varphi(L)} \bar v_j)\cdot
\\
\cdot
f_1(\bar x_1, \bar x_2, \ldots, \bar x_\ell,\ \bar y_1,\ldots, \bar y_\ell,
\bar z_1, \ldots, \bar z_T).
\end{split}\end{equation}

Пусть $\eta$~"--- $H$-многочлен $f$ из леммы~\ref{LemmaS}.
Положим
\begin{equation*}\begin{split} f_2(u_1, \ldots, u_\ell, v_1, \ldots, v_\ell,
x_1, \ldots, x_\ell, y_1, \ldots, y_\ell, z_1, \ldots, z_T) :=
\\=\sum_{\sigma, \tau \in S_\ell}
\sign(\sigma\tau)
f^{(s)}_1(u_{\sigma(1)}, \ldots, u_{\sigma(s)}, v_{\tau(1)}, \ldots, v_{\tau(s)}, x_1, \ldots, x_\ell,\ y_1,\ldots, y_\ell,
z_1, \ldots, z_T)\cdot \\ \cdot \eta(u_{\sigma(s+1)}, \ldots, u_{\sigma(\ell)})
\eta(v_{\tau(s+1)}, \ldots, v_{\tau(\ell)}).\end{split}\end{equation*}
Тогда $f_2 \in Q^H_{\ell, 4, 4\ell+T}$. Предположим, что элементы $a_1, \ldots, a_s \in \varphi(B)$
и $a_{s+1}, \ldots, a_\ell \in \varphi(Z(L))$ образуют базис пространства $\varphi(L)$.
Рассмотрим подстановку $x_i=y_i=u_i=v_i=a_i$, где $1 \leqslant i \leqslant \ell$.
Предположим, что значения $z_j=\bar z_j$, где $1 \leqslant j \leqslant T$, выбраны
таким образом, что $f_1(a_1, \ldots, a_\ell, a_1, \ldots, a_\ell,
\bar z_1, \ldots, \bar z_T)\ne 0$. Докажем, что $f_2$ также не обращается в нуль.
Действительно,
\begin{equation*}\begin{split} f_2(a_1, \ldots, a_\ell, a_1, \ldots, a_\ell,
a_1, \ldots, a_\ell, a_1, \ldots, a_\ell, \bar z_1, \ldots, \bar z_T) = \\= \sum_{\sigma, \tau \in S_\ell}
\sign(\sigma\tau)
f^{(s)}_1(a_{\sigma(1)}, \ldots, a_{\sigma(s)}, a_{\tau(1)}, \ldots, a_{\tau(s)}, a_1, \ldots, a_\ell,\ a_1,\ldots, a_\ell,
\bar z_1, \ldots, \bar z_T)\cdot \\ \cdot \eta(a_{\sigma(s+1)}, \ldots, a_{\sigma(\ell)})
\eta(a_{\tau(s+1)}, \ldots, a_{\tau(\ell)})=\\=\left(\sum_{\sigma, \tau \in S_s}
\sign(\sigma\tau)
f^{(s)}_1(a_{\sigma(1)}, \ldots, a_{\sigma(s)}, a_{\tau(1)}, \ldots, a_{\tau(s)},
 a_1, \ldots, a_\ell,\ a_1,\ldots, a_\ell,
\bar z_1, \ldots, \bar z_T)\right)\cdot\\ \cdot \left(
\sum_{ \pi, \omega \in S\lbrace s+1, \ldots,
\ell \rbrace} \sign(\pi\omega)
 \eta(a_{\pi(s+1)}, \ldots, a_{\pi(\ell)})
\eta(a_{\omega(s+1)}, \ldots, a_{\omega(\ell)})\right). \end{split}\end{equation*}
 поскольку элементы $a_j$ при $s < j \leqslant \ell$
 принадлежат центру алгебры Ли $\varphi(L)$,
 а $H$-многочлен $f^{(s)}_j$ 
обращается в нуль, если подставить такой элемент $a_i$ вместо $u_i$ или $v_i$.
Здесь $S\lbrace s+1, \ldots,
\ell \rbrace$~"--- группа подстановок на множестве $\lbrace s+1, \ldots,
\ell \rbrace$. Заметим, что многочлен $\eta$ кососимметрический. Используя~(\ref{EqKilling}), 
получаем, что
\begin{equation*}\begin{split} f_2(a_1, \ldots, a_\ell, a_1, \ldots, a_\ell,
a_1, \ldots, a_\ell, a_1, \ldots, a_\ell, \bar z_1, \ldots, \bar z_T) = \\= \left(
\sum_{\sigma, \tau \in S_s}
\sign(\sigma\tau) \tr(\ad_{\varphi(L)} a_{\sigma(1)}
\ad_{\varphi(L)} a_{\tau(1)})  \ldots \tr(\ad_{\varphi(L)} a_{\sigma(s)}
\ad_{\varphi(L)} a_{\tau(s)}) \right)\cdot \\ \cdot
f_1(a_1, \ldots, a_\ell,\ a_1,\ldots, a_\ell,
\bar z_1, \ldots, \bar z_T)  ((\ell-s)!)^2
\left(\eta(a_{s+1}, \ldots, a_\ell)\right)^2.
\end{split}\end{equation*}
 Заметим, что
\begin{equation*}\begin{split}\sum_{\sigma, \tau \in S_s}
\sign(\sigma\tau) \tr(\ad_{\varphi(L)} a_{\sigma(1)}
\ad_{\varphi(L)} a_{\tau(1)})  \ldots \tr(\ad_{\varphi(L)} a_{\sigma(s)}
\ad_{\varphi(L)} a_{\tau(s)})
=\\=\sum_{\sigma, \tau \in S_s}
\sign(\sigma\tau) \tr(\ad_{\varphi(L)} a_{1}
\ad_{\varphi(L)} a_{\tau\sigma^{-1}(1)})  \ldots
 \tr(\ad_{\varphi(L)} a_{s}
\ad_{\varphi(L)} a_{\tau\sigma^{-1}(s)})\mathrel{\stackrel{(\tau'=\tau\sigma^{-1})}{=}}\\=\sum_{\sigma, \tau' \in S_s}
\sign(\tau') \tr(\ad_{\varphi(L)} a_{1}
\ad_{\varphi(L)} a_{\tau'(1)})  \ldots
 \tr(\ad_{\varphi(L)} a_{s}
\ad_{\varphi(L)} a_{\tau'(s)})=\\=s!\det(\tr(\ad_{\varphi(L)} a_i \ad_{\varphi(L)} a_j))_{i,j=1}^s=
s!\det(\tr(\ad_{\varphi(B)} a_i \ad_{\varphi(B)} a_j))_{i,j=1}^s \ne 0,\end{split}\end{equation*}
поскольку форма Киллинга $\tr(\ad x \ad y)$ полупростой алгебры Ли
 $\varphi(B)$ невырождена.
Отсюда $$ f_2(a_1, \ldots, a_\ell, a_1, \ldots, a_\ell,
a_1, \ldots, a_\ell, a_1, \ldots, a_\ell, \bar z_1, \ldots, \bar z_T) \ne 0. $$
Заметим, что если $H$-многочлен $f_1$ был кососимметричен по некоторым из переменных $z_1,\ldots, z_T$,
$H$-многочлен $f_2$ также будет кососимметричен по
этим переменным.
Следовательно, если применить ту же процедуру не к $f_1$,
а уже к $f_2$, мы получим $f_3 \in Q^H_{\ell, 6, 6\ell+T}$.
Аналогично, определим $f_4$, используя $f_3$, $f_5$ используя $f_4$ и т.д.
В конце концов мы получим
$f:=f_k \in Q^H_{\ell, 2k, 2k\ell+T} \backslash \Id^H(\varphi)$.
\end{proof}

\section{Оценка сверху}
\label{SectionUpperLie}

В \S\ref{SectionUpperLie}--\ref{SectionLowerLie}
рассматривается некоторая $H$-хорошая алгебра Ли $L$ (см. \S\ref{SectionHnice}).

Лемма~\ref{LemmaIrrAnnBQ} является $H$-инвариантным аналогом леммы 4 из~\cite{ZaicevLie}.

\begin{lemma} \label{LemmaIrrAnnBQ}
Пусть $J \subseteq I \subseteq L$~"--- такие $H$-инвариантные идеалы,
что  $I/J$ является неприводимым $(H,L)$-модулем.
Тогда \begin{enumerate}
\item $\Ann (I/J)\cap B$ и $\Ann (I/J)\cap Q$ являются $H$-подмодулями
алгебры Ли $L$; \label{BQInv}
\item $\Ann (I/J)=(\Ann (I/J)\cap B)\oplus (\Ann (I/J)\cap Q)
\oplus N$.  (См. лемму~\ref{LemmaLBQN}.) \label{BQDecomp}
\end{enumerate}
\end{lemma}
\begin{proof}
В силу леммы~\ref{LemmaAnnHLmoduleHMod} пространства
$\Ann (I/J)$,
$\Ann (I/J)\cap B$ и $\Ann (I/J)\cap Q$ являются $H$-подмодулями.

В силу того, что нильпотентный радикал $N$ является нильпотентным идеалом,
для некоторого  $p\in\mathbb N$ 
справедливо равенство $[\underbrace{N, [N, \ldots,[N}_p, I]\ldots]=0$. Отсюда $N^p(I/J)=0$. 
Однако $I/J$ является неприводимым $(H,L)$-модулем,
откуда либо $N(I/J)=0$, либо $N(I/J)=I/J$. 
Следовательно, $N(I/J)=0$ и $N \subseteq \Ann(I/J)$. 
Для того, чтобы доказать лемму,
достаточно показать, что если $b+s \in \Ann (I/J)$, где $b \in B$, $s \in Q$, то
$b,s \in \Ann (I/J)$.
Обозначим через
$\varphi \colon L \to \mathfrak{gl}(I/J)$
гомоморфизм, отвечающий структуре $L$-модуля
на пространстве $I/J$. Тогда $\varphi(b)+\varphi(s)=0$
и $$[\varphi(b), \varphi(B)]=[-\varphi(s), \varphi(B)]=0,$$
поскольку $[B,Q]=0$.
Следовательно, элемент $\varphi(b)$ принадлежит центру алгебры Ли $\varphi(B)$ и $\varphi(b)=\varphi(s)=0$,
так как алгебра Ли $\varphi(B)$ полупроста.
Отсюда $b,s \in \Ann (I/J)$, и лемма доказана.
\end{proof}

Зафиксируем в $L$ композиционный ряд из $H$-инвариантных идеалов
$$L=L_0 \supsetneqq L_1 \supsetneqq L_2 \supsetneqq \ldots \supsetneqq
N\supsetneqq \ldots \supsetneqq L_{\theta-1}
 \supsetneqq L_\theta = 0.$$
 Введём обозначение $\height a := \max_{a \in L_k} k$ для всех $a \in L$.

\begin{remark}
Если $d=d(L,H)=0$, то $L = \Ann(L_{i-1}/L_i)$
для всех $1 \leqslant i \leqslant \theta$ и
 $[a_1, a_2, \ldots, a_n] =0$ для всех $a_i \in L$
 и $n \geqslant \theta +1$. В этом случае $c^H_n(L)=0$
 для всех $n \geqslant \theta +1$. Следовательно, ниже мы можем без ограничения общности считать, что $d > 0$.
\end{remark}

 Пусть $Y:=\lbrace y_{11}, y_{12}, \ldots, y_{1j_1};\,
 y_{21}, y_{22}, \ldots, y_{2j_2}; \ldots;\,
 y_{m1}, y_{m2}, \ldots, y_{mj_m}\rbrace$,
 $Y_1$, \ldots, $Y_q$ и $\lbrace z_1, \ldots, z_m\rbrace$
 являются такими подмножествами множества $\lbrace x_1, x_2, \ldots, x_n\rbrace$,
 что $Y_i \subseteq Y$, $|Y_i|=d+1$, $ Y_i \cap Y_j = \varnothing$
 при $i \ne j$,
 $Y \cap \lbrace z_1, \ldots, z_m\rbrace = \varnothing$,
  $j_i \geqslant 0$.
  Введём обозначение \begin{equation*}\begin{split}f_{m,q}:=\Alt_{1} \ldots \Alt_{q} \Bigl[[z_1^{h_1}, y_{11}^{h_{11}}, y_{12}^{h_{12}},
  \ldots, y_{1j_1}^{h_{1j_1}}],
 [z_2^{h_2},y_{21}^{h_{21}},y_{22}^{h_{22}},\ldots, y_{2j_2}^{h_{2j_2}}], \ldots,
 \\ [ z_m^{h_m}, y_{m1}^{h_{m1}}, y_{m2}^{h_{m2}}, \ldots, y_{mj_m}^{h_{mj_m}}]\Bigr],\end{split}\end{equation*}
 где $\Alt_i$~"--- оператор альтернирования 
по переменным из множества $Y_i$, а
 $h_i, h_{ij}\in H$~"--- произвольные фиксированные элементы.

Пусть $\xi \colon L(X | H) \to L$~"--- гомоморфизм
 $H$-модульных алгебр, индуцированный некоторой подстановкой $\lbrace x_1, x_2, \ldots, x_n, \ldots \rbrace \to L$. Будем называть гомоморфизм $\xi$ \textit{собственным} для  $f_{m,q}$, если
 $\xi(z_1) \in B \cup Q \cup N$,
 $\xi(z_i) \in N$ при $2\leqslant i \leqslant m$
 и
  $\xi(y_{ik})\in B \cup Q$ при $1\leqslant i \leqslant m$,
   $1 \leqslant k \leqslant j_i$.
   
   Воспользуемся идеей из леммы 7 из работы~\cite{ZaicevLie}:

\begin{lemma}\label{LemmaReduct}
Пусть $\xi$~"--- гомоморфизм, собственный для $H$-многочлена $f_{m,q}$.
Тогда выражение $\xi(f_{m,q})$ 
является линейной комбинацией выражений $\psi(f_{m+1,q'})$, где $\psi$~"---
собственный гомоморфизм для $f_{m+1,q'}$, а $q' \geqslant q - (\dim L)m - 2$.
  (Соответствующие множества $Y'$, $Y'_i$, переменные $z'_1, \ldots, z'_{m+1}$
  и элементы $h'_i$, $h'_{ij}$ могут быть различными для различных слагаемых $f_{m+1,q'}$.)
\end{lemma}
\begin{proof}
Пусть $\alpha_i := \height \xi(z_i)$.
Будем доказывать утверждение леммы индукцией по $\sum_{i=1}^m \alpha_i$,
причём сумма будет расти.
Заметим, что $ \alpha_i \leqslant \theta \leqslant \dim L$, и
индукция рано или поздно прекратится.
 Введём обозначения $I_i := L_{\alpha_i}$, $J_i := L_{\alpha_{i}+1}$.

Во-первых, рассмотрим случай, когда для $I_1, \ldots, I_m$,
$J_1, \ldots, J_m$ нарушены условия 1 и 2 из определения числа $d(L,H)$.
В этом случае можно выбрать такие $H$-инвариантные $B$-подмодули
$T_i$, что $I_i = T_i \oplus J_i$ и
\begin{equation}\label{EqTqUpperZero}
 \bigl[[T_1, \underbrace{L, \ldots, L}_{q_1}], [T_2, \underbrace{L, \ldots, L}_{q_2}], \ldots, [T_m,
 \underbrace{L, \ldots, L}_{q_m}]\bigr] = 0
 \end{equation}
 для всех $q_i \geqslant 0$.
Определим элементы $a'_i \in T_i$ и $a''_i \in J_i$ при помощи равенств $\xi(z_i)=a'_i+a''_i$.
Заметим, что $\height a''_i > \height \xi(z_i)$.
Поскольку $H$-многочлен $f_{m,q}$ полилинеен,  выражение
$\xi(f_{m,q})$ можно представить в виде суммы аналогичных слагаемых $\tilde\xi(f_{m,q})$,
 где $\tilde\xi(z_i)$ равно либо
$a'_i$, либо $a''_i$. В силу~(\ref{EqTqUpperZero})
слагаемое, в котором все $\tilde\xi(z_i)=a'_i \in T_i$, равно $0$.
Для всех остальных слагаемых $\tilde\xi(f_{m,q})$ выполняется
условие $\sum_{i=1}^m \height \tilde\xi(z_i) > \sum_{i=1}^m \height \xi(z_i)$.

Следовательно, без ограничения общности можно считать,
что $H$-инвариантные идеалы  $I_1, \ldots, I_m$,
$J_1, \ldots, J_m$ удовлетворяют условиям 1 и 2  из определения числа $d(L,H)$.
В этом случае $\dim(\Ann(I_1/J_1) \cap \ldots \cap \Ann(I_m/J_m))
\geqslant \dim(L)-d$.
В силу леммы~\ref{LemmaIrrAnnBQ}
\begin{equation*}\begin{split}\Ann(I_1/J_1) \cap \ldots \cap \Ann(I_m/J_m)
= (B \cap \Ann(I_1/J_1) \cap \ldots \cap \Ann(I_m/J_m)) \oplus \\ \oplus
(Q \cap \Ann(I_1/J_1) \cap \ldots \cap \Ann(I_m/J_m))\ \oplus\ N.\end{split}\end{equation*}
Выберем в $B$ базис, который включает в себя базис пространства
$B \cap \Ann(I_1/J_1) \cap \ldots \cap \Ann(I_m/J_m)$,
а в $Q$~"--- базис, который включает в себя базис пространства $Q \cap \Ann(I_1/J_1) \cap \ldots \cap \Ann(I_m/J_m)$.
Поскольку $H$-многочлен $f_{m,q}$ полилинеен,
без ограничения общности можно считать, что вместо $y_{k\ell}$ подставляются
только элементы базиса. Напомним, что $H$-многочлен  $f_{m,q}$ 
кососимметричен по переменным из $Y_i$. Следовательно, в случае, когда $\xi(f_{m,q})\ne 0$,
для любого $1 \leqslant i \leqslant q$
существуют такие $y_{kj} \in Y_i$, что либо $$\xi(y_{kj}) \in B \cap \Ann(I_1/J_1) \cap \ldots \cap \Ann(I_m/J_m),$$ либо $$\xi(y_{kj}) \in Q \cap \Ann(I_1/J_1) \cap \ldots \cap \Ann(I_m/J_m).$$

Рассмотрим случай, когда $\xi(y_{kj}) \in B \cap \Ann(I_1/J_1) \cap \ldots \cap \Ann(I_m/J_m)$
для некоторого $y_{kj}$.
Поскольку $L$ является вполне приводимым $(H,B)$-модулем (см. лемму~\ref{LemmaLBQN}),
можно выбрать такие $H$-инвариантные $B$-подмодули $T_k$, что $I_k = T_k \oplus J_k$.
Можно считать, что $\xi(z_k) \in T_k$, поскольку элементы идеалов $J_k$
имееют б\'ольшие высоты.
Следовательно, $[\xi(z_k^{h_k}), a] \in T_k \cap J_k$
 для всех $a \in B \cap \Ann(I_1/J_1) \cap \ldots \cap \Ann(I_m/J_m)$.
Отсюда $[\xi(z_k^{h_k}),a]=0$. Более того, $B \cap \Ann(I_1/J_1) \cap \ldots \cap \Ann(I_m/J_m)$ является $H$-инвариантным идеалом алгебры Ли $B$, a $[B,Q]=0$. 
Поэтому, применяя необходимое число раз тождество Якоби,
приходим к тому, что $$\xi([z_k^{h_{k}},y_{k1}^{h_{k1}},\ldots, y_{kj_k}^{h_{kj_k}}]) = 0.$$
Расписывая альтернирования, получаем, что $\xi(f_{m,q})=0$.

Рассмотрим теперь случай, когда $\xi(y_{k\ell}) \in Q \cap \Ann(I_1/J_1) \cap \ldots \cap \Ann(I_m/J_m)$
для некоторого $y_{k\ell} \in Y_q$. Распишем  в $f_{m,q}$ альтернирование $\Alt_q$
и представим $f_{m,q}$ в виде суммы слагаемых
\begin{equation*}\begin{split}\tilde f_{m,q-1} :=\Alt_{1} \ldots \Alt_{q-1} \bigl[[z_1^{h_1}, y_{11}^{h_{11}}, y_{12}^{h_{12}},
  \ldots, y_{1j_1}^{h_{1j_1}}],
 [z_2^{h_2},y_{21}^{h_{21}},y_{22}^{h_{22}},\ldots, y_{2j_2}^{h_{2j_2}}], \ldots,
 \\ [ z_m^{h_m}, y_{m1}^{h_{m1}}, y_{m2}^{h_{m2}}, \ldots, y_{mj_m}^{h_{mj_m}}]\bigr].
 \end{split}\end{equation*}
 Оператор $\Alt_q$ может менять индексы переменных, однако мы сохраним
 обозначение $y_{k\ell}$ за переменной, обладающей 
 свойством
 $\xi(y_{k\ell}) \in Q \cap \Ann(I_1/J_1) \cap \ldots \cap \Ann(I_m/J_m)$.
Теперь альтернирования оставляют $y_{k\ell}$ на месте.
Заметим, что \begin{equation*}\begin{split}[z_k^{h_{k}},y_{k1}^{h_{k1}},\ldots, y_{k\ell}^{h_{k\ell}}, \ldots, y_{kj_k}^{h_{kj_k}}] = [z_k^{h_{k}},y_{k\ell}^{h_{k\ell}},y_{k1}^{h_{k1}},\ldots, y_{kj_k}^{h_{kj_k}}] +\\ +
\sum\limits_{\beta=1}^{\ell-1}
[z_k^{h_{k}},y_{k1}^{h_{k1}},\ldots, y_{k,\beta-1}^{h_{k,\beta-1}}, [y_{k\beta}^{h_{k\beta}}, y_{k\ell}^{h_{k\ell}}],
 y_{k,\beta+1}^{h_{k,\beta+1}},\ldots, y_{k,{\ell-1}}^{h_{k,{\ell-1}}}, y_{k,{\ell+1}}^{h_{k,{\ell+1}}}, \ldots, y_{kj_k}^{h_{kj_k}}].\end{split}\end{equation*}

Заменим в первом слагаемом $[z_k^{h_{k}},y_{k\ell}^{h_{k\ell}}]$
на $z'_k$ и положим $\xi'(z'_k)
:= \xi([z_k^{h_{k}},y_{k\ell}^{h_{k\ell}}])$ и $\xi'(x) := \xi(x)$ для всех
остальных переменных~$x$. Тогда $\height\xi'(z'_k) >
\height\xi(z_k)$, и мы можем воспользоваться предположением
индукции.
Если для некоторого $j$ выполнено условие $y_{k\beta}\in Y_j$,
то мы расписываем в соответствующем
слагаемом в $\tilde f_{m,q-1}$
альтернирование $\Alt_j$.
Если $\xi(y_{k\beta}) \in B$, то соответствующее слагаемое равно нулю.
Если $\xi(y_{k\beta}) \in Q$, то $\xi([y_{k\beta}^{h_{k\beta}}, y_{k\ell}^{h_{k\ell}}]) \in N$.
Заменим $[y_{k\beta}^{h_{k\beta}}, y_{k\ell}^{h_{k\ell}}]$
 на дополнительную переменную $z'_{m+1}$
и положим $\psi(z'_{m+1}):=\xi([y_{k\beta}^{h_{k\beta}},
 y_{k\ell}^{h_{k\ell}}])$ и $\psi(x):=\xi(x)$
 для остальных переменных $x$.
 Теперь для того, чтобы получить $H$-многочлен требуемого вида, осталось применить необходимое число раз тождество Якоби.
На каждом шаге индукции число $q$
уменьшается не более чем на
 $1$, а максимальное число индукционных шагов равно $(\dim L)m$.
\end{proof}

Поскольку $N$ является нильпотентным идеалом,
для некоторого $p\in \mathbb N$ справедливо равенство
 $N^{p} = 0$.

\begin{lemma}\label{LemmaUpperLie}
Если $\lambda = (\lambda_1, \ldots, \lambda_s) \vdash n$,
причём $\lambda_{d+1} \geqslant p((\dim L)p+3)$ или $\lambda_{\dim L+1} > 0$, то
$m(L, H, \lambda) = 0$.
\end{lemma}

\begin{proof} В силу теоремы~\ref{ThKratnost}
  достаточно доказать, что $e^{*}_{T_\lambda}f \in \Id^H(A)$ для всех $f \in V^H_n$ и для всех таблиц Юнга $T_\lambda$, где $\lambda \vdash n$, $\lambda_{d+1} \geqslant p((\dim L)p+3)$ или $\lambda_{\dim L+1} > 0$.

Выберем в алгебре Ли $L$ базис,
который является объединением базисов пространств $B$, $Q$ и $N$.
Поскольку $H$-многочлены $e^{*}_{T_\lambda}f$ полилинейны,
достаточно подставлять вместо их переменных
только базисные переменные.
Напомним, что
$e^{*}_{T_\lambda} = b_{T_\lambda} a_{T_\lambda}$, а оператор $b_{T_\lambda}$ делает многочлены кососимметричными
по переменным, отвечающим каждому столбцу таблицы Юнга $T_\lambda$. Отсюда для того, чтобы $H$-многочлен $e^{*}_{T_\lambda} f$ не обратился в нуль, 
в переменные каждого столбца должны поставляться различные элементы базиса.
Однако при $\lambda_{\dim L+1} > 0$ длина первого
столбца больше, чем $(\dim L)$. Следовательно, $e^{*}_{T_\lambda} f \in \Id^H(L)$.

Рассмотрим теперь случай $\lambda_{d+1} \geqslant p((\dim L)p+3)$.
Пусть $\xi$~"--- гомоморфизм, отвечающий некоторой подстановке базисных элементов
вместо переменных $x_1, \ldots, x_n$.
Тогда $e^{*}_{T_\lambda}f$ можно представить в виде линейной комбинации $H$-многочленов $f_{m,q}$,
где $1 \leqslant m \leqslant p$, $q \geqslant p((\dim L)p+2)$, а
переменные $z_i$, где $2\leqslant i \leqslant m$, заменяются
на элементы идеала $N$.  (Для различных слагаемых $f_{m,q}$
числа $m$ и $q$,
 переменные $z_i$, $y_{ij}$, множества $Y_i$ и элементы $h_i$, $h_{ij}$ могут быть, вообще говоря, различными.)
Действительно, распишем симметризацию по всех переменным и альтернирование
тем переменным, в которые подставляются элементы из~$N$.
Если ни в одну из переменных элементы из $N$
не подставляются, положим $m=1$, представим
$H$-многочлен $f$
в виде суммы длинных коммутаторов
и в каждом слагаемом обозначим через $z_1$ одну из переменных
во внутреннем коммутаторе, расписав альтернирование
по множеству, включающему $z_1$. 
Если же в некоторые переменные подставляются
элементы из $N$, то обозначим все такие переменные через $z_k$, где $1\leqslant k \leqslant m$
для некоторого $m\in\mathbb N$.
Затем, используя тождество Якоби, занесём одну из таких переменных во внутренний
коммутатор и 
сгруппируем
все переменные, в которые подставляются элементы из $B \cup Q$,
вокруг $z_k$ так, что каждая переменная $z_k$ окажется в самом внутреннем коммутаторе некоторого
длинного коммутатора.

 Применяя необходимое число раз лемму~\ref{LemmaReduct}, будем увеличивать число $m$.
 В силу того, что идеал $N$ нильпотентен, а $\xi(f_{p+1,q})=0$
 для любого $q$ и собственного гомоморфизма~$\xi$,
 уменьшая число  $q$ не более чем на $p((\dim L)p+2)$,
 получим, что $\xi(e^{*}_{T_\lambda}f)=0$.
\end{proof}

\begin{theorem}\label{TheoremUpperLie} Пусть $L$~"--- $H$-хорошая алгебра Ли, причём
$d := d(L,H) > 0$.
Тогда существуют такие константы $C_2 > 0$, $r_2 \in \mathbb R$,
что $c^H_n(L) \leqslant C_2 n^{r_2} d^n$ для всех $n \in \mathbb N$. В случае, когда $d=0$, алгебра Ли $L$ является нильпотентной.
\end{theorem}
\begin{proof} Рассмотрим такое $\lambda \vdash n$,
что $m(A,H,\lambda) \ne 0$.
Согласно формуле крюков $$\dim M(\lambda) = \frac{n!}{\prod_{i,j} h_{ij}},$$
где $h_{ij}$~"--- длина крюка в диаграмме $D_\lambda$ с вершиной в $(i, j)$.
Оценивая полиномиальный коэффициент $\frac{(\lambda_1 +\ldots +\lambda_d)!}{(\lambda_1)!\ldots (\lambda_d)!}$ сверху по формуле возведения суммы $\underbrace{1+\ldots+1}_d$ в степень
$\lambda_1 +\ldots +\lambda_d$ и применяя лемму~\ref{LemmaUpperLie},
получаем, что
\begin{equation*}\begin{split}\dim M(\lambda) \leqslant 
\frac{n!}{(\lambda_1)!\ldots (\lambda_d)!}
\leqslant \frac{n!}{(\lambda_1 +\ldots +\lambda_d)!}
\frac{(\lambda_1 +\ldots +\lambda_d)!}{(\lambda_1)!\ldots (\lambda_d)!}
\leqslant \\ \leqslant \frac{n!}{\Bigl(n-p\bigl((\dim L)p+3\bigr)\bigl((\dim L)-d\bigr)\Bigr)!}
d^{\lambda_1 +\ldots +\lambda_d} \leqslant C_3 n^{r_3} d^n.
\end{split}\end{equation*}
для некоторых констант $C_3, r_3 > 0$.
Теперь оценка сверху следует из теоремы~\ref{TheoremUpperBoundColengthHNAssoc}.
\end{proof}
\begin{remark}
Как видно из доказательства, теорема~\ref{TheoremUpperLie}
справедлива для всех $H$-модульных
алгебр Ли $L$, для которых справедливо утверждение леммы~\ref{LemmaLBQN}.
\end{remark}

\section{Оценка снизу}
\label{SectionLowerLie}

Напомним, что оценка снизу будет доказываться для числа $d'=d'(L,H)$ (см. \S\ref{SectionHPIexpLie}).
Согласно определению числа $d'$
существуют такие $H$-инвариантные идеалы $I_1$, $I_2$, \ldots, $I_r$,
$J_1$, $J_2$, \ldots, $J_r$, где $r \in \mathbb Z_+$, алгебры $L$,
удовлетворяющие условия 1--2, что $J_k \subseteq I_k$ и
$$d' = \dim \frac{L}{\Ann(I_1/J_1) \cap \ldots \cap \Ann(I_r/J_r)}.$$
В доказательстве оценки снизу достаточно рассматривать только случай $d' > 0$.

Без ограничения общности можно считать, что для всех $1 \leqslant \ell \leqslant r$
справедливо неравенство
$$ \bigcap\limits_{k=1}^r \Ann(I_k/J_k) \ne
\bigcap\limits_{\substack{\phantom{,}k=1,\\ k\ne\ell}}^r \Ann(I_k/J_k).$$
В частности, действие алгебры Ли $L$ на всяком $I_k/J_k$ ненулевое.

Наша цель~"--- предъявить такое разбиение $\lambda \vdash n$,
что $m(L, H, \lambda)\ne 0$, а $\dim M(\lambda)$
имеет требуемое асимптотическое поведение.
Мы будем перемножать кососимметрические
$H$-многочлены, построенные в теореме~\ref{TheoremLieAlternateFinal}
для точных неприводимых модулей над редуктивными алгебрами Ли. 
Редуктивные 
алгебры Ли выбираются в лемме ниже:

\begin{lemma}\label{LemmaChooseReduct}
Существуют такие $H$-инвариантные идеалы $B_1, \ldots, B_r$
алгебры Ли $B$ и $H$-подмодули
 $\tilde R_1,\ldots,\tilde R_r \subseteq Q$
(некоторые из $\tilde R_i$ и $B_j$ могут оказаться нулевыми),
что
\begin{enumerate}
\item $B_1+ \ldots + B_r=B_1\oplus \ldots \oplus B_r$;
\item $\tilde R_1+ \ldots + \tilde R_r=\tilde R_1\oplus \ldots \oplus \tilde R_r$;
\item $\sum\limits_{k=1}^r \dim (B_k\oplus  \tilde R_k) = d$;
\item $I_k/J_k$ является точным
$(B_k\oplus\tilde R_k\oplus N)/N$-модулем;
\item $I_k/J_k$ является неприводимым
$\left(H, \left(\sum_{i=1}^r (B_i\oplus \tilde R_i)\oplus N
\right)/N\right)$-модулем;
\item $B_i I_k/J_k = \tilde R_i I_k/J_k = 0$ при $i > k$.
\end{enumerate}
\end{lemma}
\begin{proof}
Положим $N_\ell := \bigcap\limits_{k=1}^\ell \Ann (I_k/J_k)$,
$1 \leqslant \ell \leqslant r$, $N_0 := L$.
Заметим, что идеалы $N_{\ell}$ являются $H$-подмодулями.
Поскольку алгебра Ли $B$ полупростая,
в силу теоремы~\ref{TheoremHLieSemiSimple}
можно выбрать такие $H$-инвариантные идеалы $B_\ell$,
что $N_{\ell-1} \cap B =
 B_\ell \oplus (N_\ell \cap B)$.
 Кроме того, применяя условие~\ref{ConditionLComplHred} из \S\ref{SectionHnice}
к нулевой алгебре, получаем, что $L$ является вполне приводимым $H$-модулем.
Отсюда можно выбрать такие $H$-подмодули
 $\tilde R_\ell$, что
 $N_{\ell-1} \cap Q = \tilde R_\ell
  \oplus (N_\ell \cap Q)$.
Таким образом, свойства 1, 2, 6 доказаны.

 В силу леммы~\ref{LemmaIrrAnnBQ} справедливо равенство
 $N_k = (N_k \cap B) \oplus  (N_k \cap Q) \oplus  N$, из которого
 следует свойство~4. Более того,
 $$N_{\ell-1} = B_\ell \oplus  (N_\ell \cap B) \oplus
  \tilde R_\ell \oplus
   (N_\ell \cap Q) \oplus  N
   = (B_\ell\oplus \tilde R_\ell)\oplus N_\ell$$
   (прямая сумма подпространств).
   Отсюда $L = \left(\bigoplus_{i=1}^r (B_i\oplus \tilde R_i)\right)
   \oplus  N_r$,  и свойства~3 и~5 также доказаны.
\end{proof}

Обозначим через $\zeta \colon H \to \End_\mathbbm{k}(L)$ гомоморфизм, отвечающий $H$-действию,
а через $A$~"--- ассоциативную подалгебру алгебры $\End_\mathbbm{k} (L)$,
порождённую операторами из $(\ad L)$
и $\zeta(H)$.  Пусть $A_0$  и $\tilde A_0$~"--- ассоциативные алгебры из \S\ref{SectionHPIexpLie}.
Выберем в каждом из пространств $\tilde R_\ell$ некоторый базис $a_{\ell 1}, \ldots, a_{\ell, k_\ell}$.

\begin{lemma}\label{LemmaGeneralizedJordanAd}
Для всех $1 \leqslant i \leqslant r$ и $1 \leqslant j \leqslant k_i$ существуют разложения $$\ad a_{ij} = c_{ij} + d_{ij},$$ где операторы $c_{ij} \in \tilde A_0$ являются диагонализуемыми на $L$, $d_{ij} \in J(A)$,
 элементы $c_{ij}$ коммутируют друг с другом,
 а $c_{ij}$ и $d_{ij}$ представляются в виде многочленов от элементов $\ad a$, где $a\in Q$, без свободного члена.
 Более того, пространства $R_\ell := \langle c_{\ell1}, \ldots, c_{\ell, k_\ell} \rangle_\mathbbm{k}$
 являются $H$-подмодулями алгебры $A$.
\end{lemma}
\begin{proof} Поскольку $\ad a_{ij} \in A_0$, существует
разложение $\ad a_{ij} = c_{ij} + d_{ij}$, где $c_{ij}\in \tilde A_0$,
а $d_{ij} \in J(A_0)$.

В силу леммы~\ref{LemmaA0DirectSumOfFields} операторы $c_{ij}$ являются коммутирующими диагонализуемыми
на $L$ операторами. Более того,
$$hc_{\ell j}+hd_{\ell j}=h\psi(a_{\ell j}) \in
 \langle \psi(a_{\ell i})
\mid 1 \leqslant i \leqslant k_\ell\rangle_\mathbbm{k}
\subseteq \langle c_{\ell i}
\mid 1 \leqslant i \leqslant k_\ell \rangle_\mathbbm{k}
\oplus \langle d_{\ell i}
\mid  1 \leqslant i \leqslant k_\ell \rangle_\mathbbm{k}
\subseteq \tilde A_0 \oplus J(A_0)$$
для всех $h \in H$. Однако $\tilde A_0$ и $J(A_0)$ являются $H$-подмодулями, откуда $hc_{\ell j} \in \tilde A_0$, а $hd_{\ell j} \in J(A_0)$.
Следовательно, $hc_{\ell j}\in R_\ell := \langle c_{\ell1}, \ldots, c_{\ell, k_\ell} \rangle_\mathbbm{k}$
и $hd_j\in\langle d_{\ell1}, \ldots, d_{\ell, k_\ell} \rangle_\mathbbm{k}$. Отсюда $R_\ell := \langle c_{\ell1}, \ldots, c_{\ell, k_\ell} \rangle_\mathbbm{k}$
и $\langle d_{\ell1}, \ldots, d_{\ell, k_\ell} \rangle_\mathbbm{k}$ являются $H$-подмодулями алгебры $\End_\mathbbm{k}(V)$.

Теперь для завершения доказательства леммы достаточно применить лемму~\ref{LemmaHSolvableElementsJacobsonRadical}.
\end{proof}

Пусть $$\tilde B := \left(\bigoplus_{i=1}^r \ad B_i\right)\oplus \langle c_{ij} \mid 1\leqslant i \leqslant r,  1 \leqslant j \leqslant k_i
  \rangle_\mathbbm{k},$$
  $$\tilde B_0 := (\ad B)\oplus \langle c_{ij} \mid 1\leqslant i \leqslant r,  1 \leqslant j \leqslant k_i
  \rangle_\mathbbm{k} \subseteq A.$$
  В силу определения пространства $Q$ (см. лемму~\ref{LemmaLBQN}), леммы~\ref{LemmaGeneralizedJordanAd}
  и тождества Якоби, справедливы равенства \begin{equation}\label{EqcijCentral}
 [c_{ij}, \ad B]=0.\end{equation} Отсюда для $\tilde B$ и $\tilde B_0$
 справедливо условие~(\ref{EqHmoduleLieAlgebra}), и $\tilde B$ и $\tilde B_0$ являются $H$-модульными
 алгебрами Ли.
  
  \begin{lemma}\label{LemmaBcReducible}
Алгебра Ли $L$ является вполне приводимым $(H,\tilde B_0)$-модулем.
Более того, для любого $1 \leqslant k \leqslant r$ алгебра Ли $L$ является вполне приводимым $(H,(\ad B_k)\oplus R_k)$-модулем.
\end{lemma}
\begin{proof} В силу условия~\ref{ConditionLComplHred} из \S\ref{SectionHnice}
достаточно показать, что $L$ 
является вполне приводимым $\tilde B_0$-
и $(\ad B_k)\oplus R_k$-модулем без учёта $H$-действия.
Элементы $c_{ij}$ являются диагонализуемыми коммутирующими на
$L$ операторами, откуда всякое собственное подпространство
одного из операторов $c_{ij}$ инвариантно относительно
действия остальных операторов $c_{k\ell}$. По индукции
получаем разложение $L = \bigoplus_{i=1}^\alpha W_i$,
где $W_i$ являются пересечениями 
собственных подпространств операторов $c_{k\ell}$
и все элементы $c_{k\ell}$ действуют на $W_i$ как скалярные операторы.
В силу~(\ref{EqcijCentral}) подпространства $W_i$ являются $B$-подмодулями
и $L$~"--- вполне приводимый $\tilde B_0$- и $(\ad B_k)\oplus R_k$-модуль,
так как алгебры $B$ и $B_k$ полупросты.
\end{proof}

Пусть $\tilde T_k$~"--- $(\ad B)\oplus \tilde A_0$-подмодули из условия 2' \S\ref{SectionHPIexpLie}.
Напомним, что $I_k=\tilde T_k \oplus J_k$.

\begin{lemma}\label{LemmaSiProperties}
  Модули $\tilde T_k$ удовлетворяют следующим свойствам:
\begin{enumerate}
\item каждый $\tilde T_k$ является $H$-инвариантным $B$-подмодулем и неприводимым $(H,\tilde B)$-подмодулем;
\item каждый $\tilde T_k$ является вполне приводимым $(H,(\ad B_k)\oplus R_k)$-модулем,
точным как $(\ad B_k)\oplus R_k$-модуль;
\item $\sum\limits_{k=1}^r\dim ((\ad B_k)\oplus R_k) = d$;
\item $B_i \tilde T_k = R_i \tilde T_k = 0$ для всех $i > k$.
\end{enumerate}
\end{lemma}
\begin{proof} В силу того, что все $c_{ij} \in \tilde A_0$, подпространства $\tilde T_k$ являются $(H,B)$- и $(H,\tilde B)$-подмодулями.

Заметим, что $(\ad a_{ij})w=c_{ij}w$ для всех $w \in I_k/J_k$,
так как $I_k/J_k$ являются неприводимыми $A$-модулями и $J(A)\,I_k/J_k = 0$.
Следовательно, в силу леммы~\ref{LemmaChooseReduct}
модули $I_k/J_k$ являются точными $(\ad B_k)\oplus R_k$-модулями,
  $R_i\, I_k/J_k = 0$ при $i > k$,
  и элементы $c_{ij}$ линейно независимы.
  Более того, в силу свойства~5 из леммы~\ref{LemmaChooseReduct}
  модули
 $I_k/J_k$ являются неприводимыми $\left(H,\left(\sum_{i=1}^r
 (B_i\oplus \tilde R_i)\oplus N
\right)/N\right)$-модулями.
Однако алгебра Ли
$\left(\sum_{i=1}^r (B_i\oplus \tilde R_i)\oplus N
\right)/N$ действует на $I_k/J_k$
теми же самыми операторами, что и $\tilde B$. Отсюда
 $\tilde T_k \cong I_k/J_k$  являются неприводимыми $(H,\tilde B)$-модулями
 и свойство 1 доказано. Согласно лемме~\ref{LemmaBcReducible} алгебра Ли $L$ для любого
 $1 \leqslant k \leqslant r$ является вполне приводимым
  $(H,(\ad B_k)\oplus R_k)$-модулем. Используя изоморфизм $\tilde T_k \cong I_k/J_k$, 
из замечаний, сделанных выше, получаем свойства 2 и 4.
Свойство 3 является следствием свойства 3 из леммы~\ref{LemmaChooseReduct}.
\end{proof}

\begin{lemma}\label{LemmaLieChooseSubmodule} Для любого $1 \leqslant k \leqslant r$
существует разложение
 $$\tilde T_k = T_{k1} \oplus T_{k2} \oplus \ldots
\oplus T_{km},$$ где $T_{kj}$~"--- неприводимые $(H, (\ad B_k)\oplus R_k)$-подмодули,
являющиеся точными $(\ad B_k)\oplus R_k$-модулями, $m \in \mathbb N$,
$1 \leqslant j \leqslant m$.
\end{lemma}
\begin{proof}
В силу свойства 2 из леммы~\ref{LemmaSiProperties}
$$\tilde T_k = T_{k1} \oplus T_{k2} \oplus \ldots
\oplus T_{km}$$ для некоторых
$(H,(\ad B_k)\oplus R_k)$-подмодулей.
Предположим, что для некоторого $1 \leqslant j \leqslant m$
модуль $T_{kj}$ не является точным $(\ad B_k)\oplus R_k$-модулем.
Тогда $b T_{kj}=0$ для некоторого $b \in (\ad B_k)\oplus R_k$,
$b \ne 0$.
Заметим, что $\tilde B = ((\ad B_k)\oplus R_k) \oplus \tilde B_k$,
где $$\tilde B_k := \bigoplus_{i\ne k} (\ad B_i )\oplus
\bigoplus_{i\ne k} R_i$$ и $[(\ad B_k)\oplus R_k, \tilde B_k]=0$.
Обозначим через $\widehat B_k$ ассоциативную подалгебру с $1$
алгебры $\End_\mathbbm{k}(\tilde T_k)$, порождённую операторами из
алгебры Ли $\tilde B_k$. Эта подалгебра является $H$-инвариантной,
причём $$[(\ad B_k)\oplus R_k, \widehat B_k]=0$$ и $\sum_{a \in \widehat B_k}
a T_{kj} \supseteq T_{kj}$ является $H$-инвариантным $\tilde B$-подмодулем
пространства $\tilde T_k$, поскольку $$h\left(\sum_{a \in \widehat B_k}
a T_{kj}\right) = \sum_{a \in \widehat B_k}
(h_{(1)}a) (h_{(2)}T_{kj}) \subseteq \sum_{a \in \widehat B_k}
a T_{kj}$$ для всех $h\in H$. Следовательно,
$ \tilde T_k = \sum_{a \in \widehat B_k}
a T_{kj}$ и $$b \tilde T_k
= \sum_{a \in \widehat B_k}
ba T_{kj} = \sum_{a \in \widehat B_k}
a (bT_{kj})=0.$$ Получаем противоречие с точностью $(\ad B_k)\oplus R_k$-модуля $\tilde T_{k}$.
\end{proof}

 В силу условия~2' из определения числа $d'$ (см.~\S\ref{SectionHPIexpLie})
 существуют такие числа $q_1, \ldots, q_{r} \in \mathbb Z_+$,
 что
$$[[\tilde T_1, \underbrace{L, \ldots, L}_{q_1}], [\tilde T_2, \underbrace{L, \ldots, L}_{q_2}] \ldots, [\tilde T_r,
 \underbrace{L, \ldots, L}_{q_r}]] \ne 0$$
 Выберем такие числа $n_i \in \mathbb Z_+$ с максимальной суммой $\sum\limits_{i=1}^r n_i$, что
$$[[\left(\prod_{k=1}^{n_1} j_{1k}\right)\tilde T_1, \underbrace{L, \ldots, L}_{q_1}],
 [\left(\prod_{k=1}^{n_2} j_{2k}\right) \tilde T_2, \underbrace{L, \ldots, L}_{q_2}] \ldots, [\left(\prod_{k=1}^{n_r} j_{rk}\right) \tilde T_r,
 \underbrace{L, \ldots, L}_{q_r}]] \ne 0
$$ для некоторых $j_{ik}\in J(A)$.
Пусть $j_i := \prod_{k=1}^{n_i} j_{ik}$.
Тогда $j_i  \in J(A) \cup \{1\}$ и
$$[[j_1 \tilde T_1, \underbrace{L, \ldots, L}_{q_1}], [j_2 \tilde T_2, \underbrace{L, \ldots, L}_{q_2}], \ldots, [j_r \tilde T_r,
 \underbrace{L, \ldots, L}_{q_r}]] \ne 0,
 $$ но
\begin{equation}\label{EquationJZero}
[[j_1 \tilde T_1, \underbrace{L, \ldots, L}_{q_1}],
\ldots, [j_k (j \tilde T_k), \underbrace{L, \ldots, L}_{q_k}], \ldots, [j_r \tilde T_r,
 \underbrace{L, \ldots, L}_{q_r}]] = 0
\end{equation}
для всех $j \in J(A)$ и $1 \leqslant k \leqslant r$.

   В силу леммы~\ref{LemmaLieChooseSubmodule} 
можно выбрать такие неприводимые
 $(H,(\ad B_k)\oplus R_k)$-модули
$T_k \subseteq \tilde T_k$, точные как $(\ad B_k)\oplus R_k$-модули,
что \begin{equation}\label{EquationqNonZero}
[[j_1 T_1, \underbrace{L, \ldots, L}_{q_1}], [j_2 T_2, \underbrace{L, \ldots, L}_{q_2}] \ldots, [j_r T_r, \underbrace{L, \ldots, L}_{q_r}]] \ne 0.
\end{equation}

\begin{lemma}\label{LemmaChange}
Пусть $\psi \colon \bigoplus_{i=1}^r(B_i \oplus \tilde R_i) \to
\bigoplus_{i=1}^r((\ad B_i)\oplus R_i) $~"--- линейная биекция,
заданная формулами $\psi(b)= \ad b$ при $b \in B_i$ и $\psi(a_{i\ell})=c_{i\ell}$, $1 \leqslant \ell
\leqslant k_\ell$.
Пусть $f_i$~"--- полилинейные ассоциативные $H$-многочлены, а
$\bar x^{(i)}_1, \ldots, \bar x^{(i)}_{n_i}
\in \bigoplus_{i=1}^r B_i \oplus \tilde R_i$, $\bar t_i \in \tilde T_i$, $\bar u_{ik}\in L$~"--- некоторые элементы.
Тогда
\begin{equation*}\begin{split}[[j_1 f_1(\ad \bar x^{(1)}_1, \ldots, \ad \bar x^{(1)}_{n_1})
 \bar t_1, \bar u_{11}, \ldots, \bar u_{1q_1}],  \ldots, [j_r
 f_r(\ad \bar x^{(r)}_1, \ldots, \ad \bar x^{(r)}_{n_r}) \bar t_r,
 \bar u_{r1}, \ldots, \bar u_{rq_r}]]=\\ =
 [[j_1 f_1(\psi (\bar x^{(1)}_1), \ldots, \psi (\bar x^{(1)}_{n_1}))
 \bar t_1, \bar u_{11}, \ldots, \bar u_{1q_1}],  \ldots,\\ [j_r
 f_r(\psi(\bar x^{(r)}_1), \ldots, \psi (\bar x^{(r)}_{n_r})) \bar t_r,
 \bar u_{r1}, \ldots, \bar u_{rq_r}]].\end{split}
\end{equation*} 
  Другими словами, после замены $\ad a_{i\ell}$ на $c_{i\ell}$ результат не меняется.
\end{lemma}
\begin{proof}
Достаточно воспользоваться равенствами $$\ad a_{i\ell}=c_{i\ell}+d_{i\ell}=\psi(a_{i\ell})+d_{i\ell}$$
и полилинейностью $H$-многочленов $f_i$. В силу~(\ref{EquationJZero}) слагаемые с $d_{i\ell}$ обратятся в нуль.
\end{proof}

\begin{lemma}\label{LemmaLieAlt} Если $d' > 0$, то существует такое число $n_0 \in \mathbb N$,
что для любого $n\geqslant n_0$ существуют попарно непересекающиеся множества $X_1$, \ldots, $X_{2k} \subseteq \lbrace x_1, \ldots, x_n
\rbrace$, где $k := \left[\frac{n-n_0}{2d'}\right]$,
$|X_1| = \ldots = |X_{2k}|=d'$, и $H$-многочлен $f \in V^H_n \backslash
\Id^H(L)$, кососимметричный по переменным каждого из множеств $X_j$.
\end{lemma}
\begin{proof}
Обозначим через $\varphi_i \colon (\ad B_i)\oplus R_i \to
\mathfrak{gl}(T_i)$ представление, отвечающее $(\ad B_i)\oplus R_i$-действию на $T_i$.
В силу теоремы~\ref{TheoremLieAlternateFinal}
существуют такие константы $m_i \in \mathbb Z_+$,
что для любого $k$ существуют полилинейные
ассоциативные $H$-многочлены $f_i \in Q^H_{d_i, 2k, 2k d_i+m_i}
  \backslash \Id^H(\varphi_i)$, где
$d_i := \dim ((\ad B_i)\oplus R_i)$,
кососимметричные по переменным из непересекающихся множеств
$X^{(i)}_{\ell}$, где $1 \leqslant \ell \leqslant 2k$, $|X^{(i)}_{\ell}|=d_i$.

В силу~(\ref{EquationqNonZero}) для некоторых $\bar u_{i\ell} \in L$ и $\bar t_i \in T_i$
справедливо неравенство
$$[[j_1 \bar t_1, \bar u_{11}, \ldots, \bar u_{1,q_1}], [j_2 \bar t_2, \bar u_{21}, \ldots, \bar u_{2,q_2}],
 \ldots, [j_r \bar t_r, \bar u_{r1}, \ldots, \bar u_{r,q_r}]] \ne 0.
 $$
 Все $j_i \in J(A)\cup \{1\}$ являются многочленами от элементов из $\zeta(H)$ и $\ad L$. Обозначим через $\tilde m$ наибольшую среди степеней этих многочленов.

Напомним, что каждое из пространств $T_i$
является неприводимым $(H,(\ad B_i)\oplus R_i)$-модулем и
точным $(\ad B_i)\oplus R_i$-модулем.
Следовательно, в силу теоремы плотности
всякая алгебра
 $\End_\mathbbm{k}(T_i)$ порождена операторами из $\zeta(H)$
и $(\ad B_i)\oplus R_i$.
Отождествим $\End_\mathbbm{k}(T_i)$ с $M_{\dim T_i}(\mathbbm{k})$.
Тогда каждая матричная единица $e^{(i)}_{j\ell} \in M_{\dim T_i}(\mathbbm{k})$
может быть представлена как многочлен от элементов пространств $\zeta(H)$
и $(\ad B_i)\oplus R_i$. Выберем такие многочлены для всех $i$
и всех матричных единиц. Обозначим через $m_0$
наибольшую среди степеней этих многочленов.

Пусть $n_0 := r(2m_0+\tilde m+1)+ \sum_{i=1}^r (m_i+q_i)$.
Выберем $H$-многочлены $f_i$ в соответствии с теоремой~\ref{TheoremLieAlternateFinal} для $k = \left[\frac{n-n_0}{2d'}\right]$ и представлений $\varphi_i$.
Кроме того, выберем $H$-многочлен $\tilde f_1$ для $\tilde k = \left[\frac{n-2kd'-m_1}{2d_1}\right]+1$
и представления $\varphi_1$. При помощи $H$-многочленов
 $f_i$ мы получим необходимое количество альтернирований.
 Однако общая степень такого произведения может оказаться меньше $n$. Мы будем использовать
 $H$-многочлен $\tilde f_1$ для того, чтобы увеличить число переменных
 и получить $H$-многочлен степени $n$.

Поскольку $f_i \notin \Id^H(\varphi_i)$ и $\tilde f_1 \notin \Id^H(\varphi_1)$,
существуют такие $$\bar x_{i1}, \ldots, \bar x_{i, 2k d_i+m_i} \in (\ad B_i)\oplus R_i,$$
что $f_i(\bar x_{i1}, \ldots, \bar x_{i, 2k d_i+m_i})\ne 0$,
и такие $$\bar x_1, \ldots, \bar x_{2\tilde k d_1+m_1} \in (\ad B_1)\oplus R_1,$$
что $\tilde f_1(\bar x_1, \ldots, \bar x_{2\tilde k d_1+m_1}) \ne 0$.
Следовательно, для некоторых
матричных единиц $e^{(i)}_{\ell_i \ell_i},
e^{(i)}_{s_i s_i} \in \End_\mathbbm{k}(T_i)$, $1 \leqslant \ell_i, s_i \leqslant \dim {T_i}$,
$e^{(1)}_{\tilde\ell \tilde\ell}, e^{(1)}_{\tilde s \tilde s} \in \End_\mathbbm{k}(T_1)$, $1 \leqslant \tilde \ell,
\tilde s \leqslant \dim T_1$ справедливы неравенства
 $$e^{(i)}_{\ell_i \ell_i} f_i(\bar x_{i1}, \ldots, \bar x_{i, 2k d_i+m_i})
e^{(i)}_{s_i s_i} \ne 0$$
и $$e^{(1)}_{\tilde\ell \tilde\ell}\tilde f_1(\bar x_1, \ldots, \bar x_{2\tilde k d_1+m_1})
e^{(1)}_{\tilde s \tilde s} \ne 0.$$
 Отсюда $$\sum_{\ell=1}^{\dim_{T_i}}
e^{(i)}_{\ell \ell_i} f_i(\bar x_{i1}, \ldots, \bar x_{i, 2k d_i+m_i})
 e^{(i)}_{s_i \ell}$$ является ненулевым скалярным оператором, принадлежащим алгебре $\End_\mathbbm{k}(T_i)$.

Следовательно,
\begin{equation*}\begin{split} [[j_1\left(\sum_{\ell=1}^{\dim {T_1}}
e^{(1)}_{\ell \ell_1} f_1(\bar x_{11}, \ldots, \bar x_{1,2k d_1+m_1})
e^{(1)}_{s_1 \tilde \ell} \tilde f_1(\bar x_1, \ldots, \bar x_{2\tilde k d_1+m_1})
 e^{(1)}_{\tilde s \ell}\right)\bar t_1, \bar u_{11}, \ldots, \bar u_{1q_1}], \\ [j_2\left(\sum_{\ell=1}^{\dim {T_2}}
e^{(2)}_{\ell \ell_2} f_2(\bar x_{21}, \ldots, \bar x_{2,2k d_2+m_2})
 e^{(2)}_{s_2 \ell}\right)\bar t_2, \bar u_{21}, \ldots, \bar u_{2q_2}],
 \ldots, \\
 [j_r\left(\sum_{\ell=1}^{\dim {T_r}}
e^{(r)}_{\ell \ell_r} f_r(\bar x_{r1}, \ldots, \bar x_{r, 2k d_r+m_r})
 e^{(r)}_{s_r \ell}\right)\bar t_r, \bar u_{r1}, \ldots, \bar u_{rq_r}]]\ne 0.\end{split}
 \end{equation*}
 (Считаем, что всякий $f_i$ является $H$-многочленом от переменных $x_{i1}, \ldots,
x_{i,2k d_i+m_i}$, а всякий $\tilde f_1$ является $H$-многочленом от переменных $x_1, \ldots, x_{2\tilde k d_1 + m_1}$.)
Введём обозначение $X_\ell := \bigcup_{i=1}^{r} X^{(i)}_{\ell}$,
где $X^{(i)}_{\ell}$~"--- множества, по переменным которых
был кососимметричен $H$-многочлен $f_i$.
Обозначим через $\Alt_\ell$
оператор альтернирования по переменным из множества $X_\ell$. 
Рассмотрим выражение
\begin{equation*}\begin{split}\tilde f(x_1, \ldots, x_{2\tilde k d_1 + m_1}; x_{11}, \ldots, x_{1, 2k d_1+m_1};
\ldots; \ x_{r1}, \ldots, x_{r, 2k d_r+m_r}) := \end{split}
 \end{equation*}
\begin{equation*}\begin{split} =
 \Alt_1 \Alt_2 \ldots \Alt_{2k} [[j_1\left(\sum_{\ell=1}^{\dim {T_1}}
e^{(1)}_{\ell \ell_1} f_1(x_{11}, \ldots, x_{1, 2k d_1+m_1})
e^{(1)}_{s_1 \tilde \ell}\ \cdot \right. \\ 
 \left. \cdot \tilde f_1(x_1, \ldots, x_{2\tilde k d_1+m_1})
 e^{(1)}_{\tilde s \ell}\right)\bar t_1, \bar u_{11}, \ldots, \bar u_{1q_1}], \\
 [j_2\left(\sum_{\ell=1}^{\dim {T_2}}
e^{(2)}_{\ell \ell_2} f_2(x_{21}, \ldots, x_{2, 2k d_2+m_2})
 e^{(2)}_{s_2 \ell}\right)\bar t_2, \bar u_{21}, \ldots, \bar u_{2q_2}],
 \ldots, \\
 [j_r\left(\sum_{\ell=1}^{\dim {T_r}}
e^{(r)}_{\ell \ell_r} f_r(x_{r1}, \ldots, x_{r,2k d_r+m_r})
 e^{(r)}_{s_r \ell}\right)\bar t_r, \bar u_{r1}, \ldots, \bar u_{rq_r}]].\end{split}
 \end{equation*}
Тогда
\begin{equation*}\begin{split}\tilde f(\bar x_1, \ldots, \bar x_{2\tilde k d_1 + m_1}; \bar x_{11}, \ldots, \bar x_{1, 2k d_1+m_1};
\ldots; \ \bar x_{r1}, \ldots, \bar x_{r, 2k d_r+m_r})
= \\ = (d_1!)^{2k} \ldots (d_r!)^{2k} [[j_1\left(\sum_{\ell=1}^{\dim {T_1}}
e^{(1)}_{\ell \ell_1} f_1(\bar x_{11}, \ldots, \bar x_{1,2k d_1+m_1})
e^{(1)}_{s_1 \tilde \ell}\ \cdot\right. \\ \left. \cdot \tilde f_1(\bar x_1, \ldots, \bar x_{2\tilde k d_1+m_1})
 e^{(1)}_{\tilde s \ell}\right)\bar t_1, \bar u_{11}, \ldots, \bar u_{1q_1}],
 \ldots, \\
 [j_r\left(\sum_{\ell=1}^{\dim {T_r}}
e^{(r)}_{\ell \ell_r} f_r(\bar x_{r1}, \ldots, \bar x_{r, 2k d_r+m_r})
 e^{(r)}_{s_r \ell}\right)\bar t_r, \bar u_{r1},
  \ldots, \bar u_{rq_r}]]\ne 0,\end{split}
 \end{equation*}
поскольку $H$-многочлены $f_i$ кососимметричны по переменным каждого множества $X^{(i)}_{\ell}$
и в силу леммы~\ref{LemmaSiProperties} справедливо равенство $((\ad B_i)\oplus R_i)\tilde T_\ell = 0$
при $i > \ell$. Теперь представим
$e^{(i)}_{\ell j}$ в виде многочленов от переменных из $(\ad B_i)\oplus R_i$
и $\zeta(H)$.
В силу линейности $H$-многочленов $\tilde f$ по $e^{(i)}_{\ell j}$
можно заменить матричные единицы $e^{(i)}_{\ell j}$ произведениями элементов
из $(\ad B_i)\oplus R_i$ и $\zeta(H)$, и для некоторого выбора таких произведений
выражение  не обратится в нуль. Используя~(\ref{EqHLModule2}),
можно переместить все $\zeta(h)$ вправо.
 В силу~\ref{LemmaChange} все элементы из $(\ad B_i)\oplus R_i$
 можно заменить на элементы из $B_i\oplus \tilde R_i$,
 и выражение по-прежнему останется ненулевым.
Обозначим через $\psi \colon \bigoplus_{i=1}^r (B_i \oplus \tilde R_i) \to
\bigoplus_{i=1}^r ((\ad B_i)\oplus R_i) $ соответствующую линейную биекцию.
Теперь представим $j_i$ как многочлены от элементов из $\ad L$ и $\zeta(H)$.
Поскольку выражение $\tilde f$ линейно по $j_i$,
можно заменить каждый из элементов $j_i$ на некоторый одночлен,
т.е. произведение элементов пространств $\ad L$ и $\zeta(H)$.
Используя~(\ref{EqHLModule2}),
снова перемещаем элементы $\zeta(h)$ вправо. 
Теперь заменим элементы из $\ad L$ на новые переменные.
Тогда
\begin{equation*}\begin{split}\hat f :=
 \Alt_1 \Alt_2 \ldots \Alt_{2k} \biggl[\Bigl[\Bigl[y_{11}, [y_{12}, \ldots
 [y_{1 \alpha_1}, \Bigl[z_{11}, [z_{12},
 \ldots, [z_{1 \beta_1}, 
 \\
  \bigl(h_1 f_1(\ad x_{11}, \ldots, \ad x_{1, 2k d_1+m_1})\bigr)
 [w_{11}, [w_{12}, \ldots, [w_{1 \theta_1}, \\
 \bigl(\tilde h \tilde f_1(\ad x_1, \ldots, \ad x_{2\tilde k d_1+m_1})\bigr)
 [w_{1}, [w_{2}, \ldots, [w_{\tilde \theta},
  t_1^{h'_1}]\ldots \Bigr],
  u_{11}, \ldots, u_{1q_1}\Bigr], \\ 
  \Bigl[\Bigl[y_{21}, [y_{22}, \ldots
 [y_{2 \alpha_2}, \Bigl[z_{21}, [z_{22},
 \ldots, [z_{2 \beta_2},
\\
  \bigl(h_2 f_2(\ad x_{21}, \ldots, \ad x_{2, 2k d_2+m_2})\bigr)
 [w_{21}, [w_{22}, \ldots, [w_{2 \theta_2},
  t_2^{h'_2}]\ldots \Bigr],
  u_{21}, \ldots, u_{2q_2}\Bigr],
 \ldots, \\
 \Bigl[\Bigl[y_{r1}, [y_{r2}, \ldots,
 [y_{r \alpha_r}, \Bigr[z_{r1},
  [z_{r2},
 \ldots, [z_{r \beta_r},
 \\
 \bigl(h_r f_r(\ad x_{r1}, \ldots, \ad x_{r, 2k d_r+m_r})\bigr)
 [w_{r1}, [w_{r2}, \ldots, [w_{r \theta_r}, t_r^{h'_r}]\ldots \Bigr],
  u_{r1}, \ldots, u_{rq_r}\Bigr]\biggr]\end{split}
 \end{equation*}
  для некоторых  $0 \leqslant \alpha_i \leqslant \tilde m$,
  \quad
  $0 \leqslant \beta_i, \theta_i, \tilde \theta \leqslant m_0$,
  \quad $h_i, h'_i, \tilde h \in H$,\quad
  $\bar y_{i\ell}, \bar z_{i\ell},
  \bar w_{i\ell}, \bar w_i \in L$
 не обращается в нуль при подстановке
 $$t_i=\bar t_i,\ u_{i\ell}=\bar u_{i\ell},\ x_{i\ell}=\psi^{-1}(\bar x_{i\ell}),\ x_i = \psi^{-1}(\bar x_i),\ y_{i\ell}=\bar y_{i\ell},\ z_{i\ell}=\bar z_{i\ell},\ w_{i\ell}=\bar w_{i\ell},\ w_i = \bar w_i.$$
 
 Следовательно, \begin{equation*}\begin{split} f_0 :=
 \Alt_1 \Alt_2 \ldots \Alt_{2k} \biggl[\Bigl[\Bigl[y_{11}, [y_{12}, \ldots
 [y_{1 \alpha_1}, \Bigl[z_{11}, [z_{12},
 \ldots, [z_{1 \beta_1},
\\
  \bigl(h_1 f_1(\ad x_{11}, \ldots, \ad x_{1, 2k d_1+m_1})\bigr)
 [w_{11}, [w_{12}, \ldots, [w_{1 \theta_1}, t_1]\ldots \Bigr],
  u_{11}, \ldots, u_{1q_1}\Bigr],\\
  \Bigl[\Bigl[y_{21}, [y_{22}, \ldots
 [y_{2 \alpha_2}, \Bigl[z_{21}, [z_{22},
 \ldots, [z_{2 \beta_2},
\\
  \bigl(h_2 f_2(\ad x_{21}, \ldots, \ad x_{2, 2k d_2+m_2})\bigr)
 [w_{21}, [w_{22}, \ldots, [w_{2 \theta_2},
  t_2^{h'_2}]\ldots \Bigr],
  u_{21}, \ldots, u_{2q_2}\Bigr],
 \ldots, \\
 \Bigl[\Bigl[y_{r1}, [y_{r2}, \ldots,
 [y_{r \alpha_r}, \Bigr[z_{r1},
  [z_{r2},
 \ldots, [z_{r \beta_r},
 \\
 \bigl(h_r f_r(\ad x_{r1}, \ldots, \ad x_{r, 2k d_r+m_r})\bigr)
 [w_{r1}, [w_{r2}, \ldots, [w_{r \theta_r}, t_r^{h'_r}]\ldots \Bigr],
  u_{r1}, \ldots, u_{rq_r}\Bigr]\biggr]\end{split}
 \end{equation*}
   не обращается в нуль при подстановке
 $$t_1 = \bigl(\tilde h \tilde f_1(\ad \bar x_1, \ldots, \ad \bar x_{2\tilde k d_1+m_1})\bigr)
 [\bar w_{1}, [\bar w_{2}, \ldots, [\bar w_{\tilde \theta},
 h'_1 \bar t_1]\ldots],$$
 $t_i=\bar t_i$ при $2 \leqslant i \leqslant r$; $u_{i\ell}=\bar u_{i\ell}$,
 $x_{i\ell}=\psi^{-1}(\bar x_{i\ell})$, $y_{i\ell}=\bar y_{i\ell}$,
 $z_{i\ell}=\bar z_{i\ell}$, $w_{i\ell}=\bar w_{i\ell}$.

Заметим, что $f_0 \in V_{\tilde n}^H$,
  $$\tilde n: = 2kd +r+ \sum_{i=1}^r (m_i + q_i + \alpha_i+\beta_i+\theta_i)
  \leqslant n.$$ Если $n=\tilde n$, положим $f:=f_0$.
  Предположим, что $n > \tilde n$.
Выражение $$\bigl(\tilde h \tilde f_1(\ad \bar x_1, \ldots, \ad \bar x_{2\tilde k d_1+m_1})\bigr)
 [\bar w_{1}, [\bar w_{2}, \ldots, [\bar w_{\tilde \theta},
  h'_1 \bar t_1]\ldots]$$ является линейной комбинацией длинных коммутаторов,
  каждый из которых содержит как минимум $2\tilde k d_1+m_1+1 > n-\tilde n+1$
  элементов алгебры Ли $L$.
       Следовательно, $ f_0$ не обращается в нуль при подстановке
 $t_1 = [\bar v_1, [\bar v_2, [\ldots, [\bar v_q,  h'_1\bar t_1]\ldots]$
 для некоторых $q \geqslant n-\tilde n$, $\bar v_i \in L$;
  $t_i=\bar t_i$ при $2 \leqslant i \leqslant r$; $u_{i\ell}=\bar u_{i\ell}$,
 $x_{i\ell}=\psi^{-1}(\bar x_{i\ell})$, $y_{i\ell}=\bar y_{i\ell}$,
 $z_{i\ell}=\bar z_{i\ell}$, $w_{i\ell}=\bar w_{i\ell}$.
Отсюда \begin{equation*}\begin{split} f :=
 \Alt_1 \Alt_2 \ldots \Alt_{2k} \biggl[\Bigl[\Bigl[y_{11}, [y_{12}, \ldots
 [y_{1 \alpha_1}, \Bigl[z_{11}, [z_{12},
 \ldots, [z_{1 \beta_1},
\\
  \bigl(h_1 f_1(\ad x_{11}, \ldots, \ad x_{1, 2k d_1+m_1})\bigr)
 [w_{11}, [w_{12}, \ldots, [w_{1 \theta_1},
  \\
  \bigl[v_1, [v_2, [\ldots, [v_{n-\tilde n}, t_1]\ldots\bigr]\ldots \Bigr],
  u_{11}, \ldots, u_{1q_1}\Bigr],\\
  \Bigl[\Bigl[y_{21}, [y_{22}, \ldots
 [y_{2 \alpha_2}, \Bigl[z_{21}, [z_{22},
 \ldots, [z_{2 \beta_2},
\\
  \bigl(h_2 f_2(\ad x_{21}, \ldots, \ad x_{2, 2k d_2+m_2})\bigr)
 [w_{21}, [w_{22}, \ldots, [w_{2 \theta_2},
  t_2^{h'_2}]\ldots \Bigr],
  u_{21}, \ldots, u_{2q_2}\Bigr],  \\
 \ldots, \Bigl[\Bigl[y_{r1}, [y_{r2},\ldots,
 [y_{r \alpha_r}, \Bigr[z_{r1},
  [z_{r2},
 \ldots, [z_{r \beta_r},
 \\
 \bigl(h_r f_r(\ad x_{r1}, \ldots, \ad x_{r, 2k d_r+m_r})\bigr)
 [w_{r1}, [w_{r2}, \ldots, [w_{r \theta_r}, t_r^{h'_r}]\ldots \Bigr],
  u_{r1}, \ldots, u_{rq_r}\Bigr]\biggr]\end{split}
 \end{equation*}
  не обращается в нуль при подстановке
  $v_\ell = \bar v_\ell$ при $1 \leqslant \ell \leqslant n-\tilde n$,
  $$t_1 = [\bar v_{n-\tilde n +1}, [\bar v_{n-\tilde n +2}, [\ldots, [\bar v_q,  h'_1\bar t_1]\ldots];$$
  $t_i=\bar t_i$ при $2 \leqslant i \leqslant r$; $u_{i\ell}=\bar u_{i\ell}$,
 $x_{i\ell}=\psi^{-1}(\bar x_{i\ell})$, $y_{i\ell}=\bar y_{i\ell}$,
 $z_{i\ell}=\bar z_{i\ell}$, $w_{i\ell}=\bar w_{i\ell}$.
 Теперь осталось заметить, что $H$-многочлен $f$ принадлежит пространству $V_n^H$ и удовлетворяет всем условиям леммы.
\end{proof}

\begin{lemma}\label{LemmaLieCochar} Пусть
 $k, n_0$~"--- числа из леммы~\ref{LemmaLieAlt}.   Тогда для всех $n \geqslant n_0$
 существует  такое разбиение $\lambda = (\lambda_1, \ldots, \lambda_s) \vdash n$,
 что $m(L, H, \lambda) \ne 0$ и $\lambda_i > 2k-C$ для всех $1 \leqslant i \leqslant d'$.
Здесь $C := p((\dim L)p + 3)((\dim L)-d')$, где $p \in \mathbb N$~"--- такое число, что $N^p=0$.
\end{lemma}
\begin{proof} Достаточно доказать, что $e^*_{T_\lambda} f \notin \Id^H(L)$
для некоторой таблицы Юнга $T_\lambda$ требуемой формы $\lambda$.
Известно (см., например, теорему 3.2.7 из~\cite{Bahturin}), что $$\mathbbm{k}S_n = \bigoplus_{\lambda,T_\lambda} \mathbbm{k}S_n e^{*}_{T_\lambda},$$ где суммирование ведётся по множеству стандартных таблиц Юнга $T_\lambda$
всевозможных форм $\lambda \vdash n$. Отсюда $\mathbbm{k}S_n f = \sum_{\lambda,T_\lambda} \mathbbm{k}S_n e^{*}_{T_\lambda}f
\not\subseteq \Id^H(L)$ и $e^{*}_{T_\lambda} f \notin \Id^H(L)$ для некоторого $\lambda \vdash n$.
Докажем, что разбиение $\lambda$ имеет требуемый вид.
Достаточно доказать, что
$\lambda_{d'} > 2k-C$, так как
$\lambda_i \geqslant \lambda_{d'}$ для всех $1 \leqslant i \leqslant d'$.
В любой строчке таблицы $T_\lambda$ содержится не более одного 
номера переменной из одного и того же множества $X_i$,
поскольку $e^{*}_{T_\lambda} = b_{T_\lambda} a_{T_\lambda}$,
а $a_{T_\lambda}$ симметризует по переменным, отвечающим каждой строчке таблицы $T_\lambda$.
Отсюда $$\sum_{i=1}^{d'-1} \lambda_i \leqslant 2k(d'-1) + (n-2kd') = n-2k.$$
Из леммы~\ref{LemmaUpperLie} и неравенства $d' \geqslant d$ следует, что
$$\lambda_{d'+1} \leqslant \lambda_{d+1} < p((\dim L)p + 3),$$ откуда
$\sum_{i=1}^{d'} \lambda_i > n-C$. Следовательно,
$\lambda_{d'} > 2k-C$.
\end{proof}
\begin{proof}[Доказательство теоремы~\ref{TheoremMainLieH}]
Диаграмма Юнга~$D_\lambda$ из леммы~\ref{LemmaLieCochar} содержит квадратную
поддиаграмму~$D_\mu$, где $\mu=(\underbrace{2k-C, \ldots, 2k-C}_{d'})$.
Из правила ветвления для группы $S_n$ следует, что
если рассмотреть сужение $S_n$-действия на $M(\lambda)$ до $S_{n-1}$-действия,
то $\mathbbm{k}S_n$-модуль $M(\lambda)$
оказывается прямой суммой всех неизоморфных
$\mathbbm{k}S_{n-1}$-модулей $M(\nu)$, где $\nu \vdash (n-1)$ и всякая таблица $D_\nu$
получена из $D_\lambda$ удалением одной клетки. В частности,
$\dim M(\nu) \leqslant \dim M(\lambda)$.
Применяя правило ветвления $(n-d'(2k-C))$ раз, получаем, что $\dim M(\mu) \leqslant \dim M(\lambda)$.
В силу формулы крюков $$\dim M(\mu) = \frac{(d'(2k-C))!}{\prod_{i,j} h_{ij}},$$
где $h_{ij}$~"--- длина крюка с вершиной в $(i, j)$.
По формуле Стирлинга
\begin{equation*}\begin{split}
c_n^H(L)\geqslant \dim M(\lambda) \geqslant \dim M(\mu) \geqslant \frac{(d'(2k-C))!}{((2k-C+d')!)^{d'}}
\sim \\ \sim \frac{
\sqrt{2\pi {d'}(2k-C)} \left(\frac{d'(2k-C)}{e}\right)^{d'(2k-C)}
}
{
\left(\sqrt{2\pi (2k-C+d')}
\left(\frac{2k-C+d'}{e}\right)^{2k-C+d'}\right)^{d'}
} \sim C_4 k^{r_4} \left( d'\right)^{2kd'}
\end{split}\end{equation*}
для некоторых констант $C_4 > 0$, $r_4 \in \mathbb Q$
при $k \to \infty$.
Поскольку $k = \left[\frac{n-n_0}{2d'}\right]$,
это доказывает оценку снизу.
Оценка сверху была доказана в теоремы~\ref{TheoremUpperLie}.
Учитывая, что $d'\geqslant d$, получаем отсюда, что $d=d(L,H)=d'=d'(L,H)$.
Теперь утверждение теоремы~\ref{TheoremMainLieH} следует из полученных оценок.
\end{proof}

\begin{proof}[Доказательство теоремы~\ref{TheoremMainLieHSum}]
Пусть $L = L_1 \oplus \ldots \oplus L_q$, где $L_i$~"--- $H$-хорошие идеалы.
Во-первых, $c^H_n(L) \geqslant c^H_n(L_i)$ для всех $n\in\mathbb N$ и $1 \leqslant i \leqslant q$,
поскольку $L_i$ являются $H$-инвариантными подалгебрами алгебры Ли $L$. 
Отсюда $$\max_{1 \leqslant i \leqslant q} \PIexp^H(L_i) \leqslant \mathop{\underline{\lim}}_{n\to \infty}\sqrt[n]{c^H_n(L)}.$$

Пусть $f_0 \in V_n^H$, где $n\in\mathbb N$. Для того, чтобы проверить, что $f_0 \in \Id^H(L)$, 
достаточно подставлять только элементы базиса.
Выберем в $L$ базис, который является объединением базисов идеалов $L_i$.
Тогда если вместо переменных $H$-многочлена $f_0$
подставляются элементы из разных идеалов $L_i$,
данный $H$-многочлен обращается в нуль.
Отсюда для того, чтобы доказать, что $f_0 \in \Id^H(L)$,
достаточно показать, что $f_0 \in \Id^H(L_i)$
для всех $1\leqslant i \leqslant s$. 
Пусть $d := \max_{1 \leqslant i \leqslant q} d(L_i, H)=\max_{1 \leqslant i \leqslant q} \PIexp^H(L_i)$.
Тогда из леммы~\ref{LemmaUpperLie} следует, что
для всех $\lambda = (\lambda_1, \ldots, \lambda_s) \vdash n$, где
 $\lambda_{d+1} \geqslant p((\dim L)p+3)$ или $\lambda_{\dim L+1} > 0$, справедливо равенство
$m(L, H, \lambda) = 0$. Повторяя рассуждения теоремы~\ref{TheoremUpperLie},
применим теорему~\ref{TheoremUpperBoundColengthHNAssoc} и получим отсюда, что существуют такие $C_2 > 0$, $r_2 \in \mathbb R$,
что $c^H_n(L) \leqslant C_2 n^{r_2} d^n$ для всех $n \in \mathbb N$.
Следовательно, $$\mathop{\overline{\lim}}_{n\to \infty} \sqrt[n]{c^H_n(L)} \leqslant d = \max_{1 \leqslant i \leqslant q} \PIexp^H(L_i).$$
Оценка снизу была получена выше.
\end{proof}

\section{Рост градуированных тождеств}\label{SectionGrLie}

Данный параграф посвящён доказательству теорем~\ref{TheoremMainLieGr} и~\ref{TheoremMainLieGrSum},
сформулированных в \S\ref{SectionMainLie}.

В случае, когда градуирующая группа $G$ конечна, достаточно рассмотреть действие конечномерной полупростой алгебры Хопфа $(\mathbbm{k}G)^*$, которое является двойственным к $\mathbbm{k}G$-кодействию, соответствующему $G$-градуировке.
(См. \S\ref{Section(Co)modules}.) В этом случае теоремы~\ref{TheoremMainLieGr} и~\ref{TheoremMainLieGrSum} получаются из теорем~\ref{TheoremMainLieHSS} и~\ref{TheoremMainLieHSSSum} при помощи предложения~\ref{PropositionCnGrCnGenH} с учётом изоморфизма $(\mathbbm{k}G)^*\cong \mathbbm{k}^G$. Этот переход позволяет
использовать для градуированной PI-экспоненты в случае алгебраически замкнутого поля $\mathbbm{k}$ формулу из \S\ref{SectionHPIexpLie}.

В случае произвольной градуирующей группы используется тот же приём с использованием действия алгебры Хопфа
$\mathbbm{k}\hat G$, где $\hat G := \Hom(G, \mathbbm{k}^\times)$, $$\chi a = \chi(g) a\text{ для всех }\chi \in \hat G,\ a \in L^{(g)} \text{ и } g\in G,$$
$L$~"--- алгебра Ли, градуированная группой $G$.

\begin{example}\label{ExampleHniceGr}
Пусть $L$~"--- конечномерная алгебра Ли
над алгебраически замкнутым полем $\mathbbm{k}$ характеристики $0$, градуированная конечнопорождённой
абелевой группы $G$.
Тогда $L$~"--- $\mathbbm{k}\hat G$-хорошая алгебра.
\end{example}
\begin{proof}
Сперва заметим, что разрешимый и нильпотентный радикалы инвариантны относительно всех автоморфизмов.
Следовательно, они $\hat G$- и $\mathbbm{k}\hat G$-инвариантны.

В силу теоремы~\ref{TheoremGradLevi} существует $G$-градуированное, а значит и $\mathbbm{k}\hat G$-инвариантное разложение Леви.

Из леммы~\ref{LemmaAbelianDual} следует, что $\mathbbm{k}\hat G$-подмодули
любого $G$-градуированного подпространства являются градуированными подпространствами.
Если $V=\bigoplus_{g\in G} V^{(g)}$~"--- $G$-градуированное пространство,
то алгебра $\End_\mathbbm{k}(V)$ наделена следующе естественной градуировкой: $$\Hom_\mathbbm{k}(V^{(g_1)}, V^{(g_2)})
 \subseteq \End_\mathbbm{k}(V)^{(g_2 g_1^{-1})}\text{ для всех }g_1, g_2 \in G.$$ (Результат 
 применения отображений из пространства
 $\Hom_\mathbbm{k}(V^{(g_1)}, V^{(g_2)})$ к элементам компонент $V^{(g)}$, где $g\ne g_1$, равен нулю.)
 Эта $G$-градуировка соответствует стандартному $\mathbbm{k}\hat G$-действию на $\End_\mathbbm{k}(V)$:
 $$(h\psi)(v)=h_{(1)}\psi((Sh_{(2)})v)\text{ для всех }h\in \mathbbm{k}\hat G,\ \psi \in \End_\mathbbm{k}(V),\ v\in V.$$
Отсюда из градуированной версии теоремы Веддербёрна~"--- Мальцева (следствие~\ref{CorollaryGradedWedderburnMalcev})
следует существование $\mathbbm{k}\hat G$-инвариантных разложений Веддербёрна~"--- Мальцева,
которые требуются в условии~\ref{ConditionWedderburn}.
 Аналогично, условие~\ref{ConditionLComplHred} следует из теоремы~\ref{TheoremGradWeyl}. Отсюда алгебра Ли $L$ действительно  является $\mathbbm{k}\hat G$-хорошей.
\end{proof}

Докажем теперь следующий частный случай теоремы~\ref{TheoremMainLieGr}: 
 
\begin{theorem}\label{TheoremFinGenAbelianGr}
Пусть $L$~"--- конечномерная алгебра Ли
над полем $\mathbbm{k}$ характеристики $0$,
градуированная конечнопорождённой абелевой группой $G$.
Тогда 
\begin{enumerate}
\item либо существует такое $n_0$, что $c_n^{G\text{-}\mathrm{gr}}(L)=0$ при всех $n\geqslant n_0$;
\item либо существуют такие константы $C_1, C_2 > 0$, $r_1, r_2 \in \mathbb R$,
  $d \in \mathbb N$, что $$C_1 n^{r_1} d^n \leqslant c^{G\text{-}\mathrm{gr}}_n(L)
   \leqslant C_2 n^{r_2} d^n\text{ для всех }n \in \mathbb N.$$
\end{enumerate}
   В частности, существует $\PIexp^{G\text{-}\mathrm{gr}}(L)\in\mathbb Z_+$ и, таким образом,
   для градуированных тождеств справедлив аналог гипотезы Амицура.
\end{theorem} 
\begin{proof}
Как и прочие виды коразмерностей, градуированные коразмерности
не меняются при расширении основного поля. Доказательство этого факта для
градуированных коразмерностей полностью аналогично
доказательству для случая обычных коразмерностей (см., например, \cite[теорема~4.1.9]{ZaiGia}
или~\cite[\S 2]{ZaicevLie}.
Отсюда без ограничения общности можно считать основное поле
 $\mathbbm{k}$ алгебраически замкнутым.

В примере~\ref{ExampleHniceGr} было показано, что алгебра Ли $L$ является $\mathbbm{k}\hat G$-хорошей,
в силу чего теорема~\ref{TheoremFinGenAbelianGr} следует из теоремы~\ref{TheoremMainLieH} и предложения~\ref{PropositionGrToHatG}.
\end{proof}

Сведём теперь случай произвольной градуирующей группы $G$ к случаю,
когда $G$ является конечнопорождённой абелевой группой.

Нам потребуется следующее утверждение (см., например, \cite[лемма 2.1]{PaReZai}): 
\begin{lemma}\label{LemmaGradedNonZero}
Пусть $L$~"--- алгебра Ли, градуированная некоторой группой~$G$.
Предположим, что $[L^{(g_1)}, \ldots, L^{(g_k)}]\ne 0
$ для некоторых $g_1, \ldots, g_k \in G$. Тогда
$g_i g_j = g_j g_i$ для всех $1 \leqslant i, j \leqslant k$.
\end{lemma}

Пусть $G_0$~"--- подгруппа группы $G$.
Введём обозначение $L_{G_0} := \bigoplus_{g\in G_0} L^{(g)}$.
Тогда силу леммы~\ref{LemmaGrEquivCodimTheSame} для всех $n\in\mathbb N$
справедливо равенство
$c_n^{G\text{-}\mathrm{gr}}(L_{G_0})=c_n^{G_0\text{-}\mathrm{gr}}(L_{G_0})$.

\begin{lemma}\label{LemmaInclExcl}
Пусть $L$~"--- конечномерная алгебра Ли над полем $\mathbbm{k}$ характеристики $0$,
градуированная произвольной группой $G$.
Тогда существуют такие конечнопорождённые абелевы подгруппы $G_1, \ldots, G_r \subseteq G$,
что 
\begin{equation}\begin{split}\label{EqInclExcl}
c^{G\text{-}\mathrm{gr}}_n(L)=
 \sum_{i=1}^r c^{G_i\text{-}\mathrm{gr}}_n(L_{G_i})
- \sum_{i,j=1}^r c^{G_i\cap G_j\text{-}\mathrm{gr}}_n(L_{G_i \cap G_j})+\\+ \sum_{i,j,k=1}^r c^{G_i\cap G_j \cap G_k\text{-}\mathrm{gr}}_n(L_{G_i \cap G_j \cap G_k}) - \ldots +(-1)^{r-1} 
c^{G_1 \cap G_2 \cap \ldots \cap G_r\text{-}\mathrm{gr}}_n(L_{G_1 \cap G_2 \cap \ldots \cap G_r}).\end{split}\end{equation}
\end{lemma}
\begin{proof}
 Сперва заметим, что
 $V^{G\text{-}\mathrm{gr}}_n = \bigoplus_{g_1, \ldots, g_n \in G} V_{g_1, \ldots, g_n}$
 где 
$$
 V_{g_1, \ldots, g_n} := \left\langle \left[x^{(g_{\sigma(1)})}_{\sigma(1)},
x^{(g_{\sigma(2)})}_{\sigma(2)}, \ldots, x^{(g_{\sigma(n)})}_{\sigma(n)}\right]
\mathrel{\Bigr|} \sigma\in S_n \right\rangle_\mathbbm{k}.
$$
Используя процесс линеаризациии (см., например, \cite[теорема 4.2.3]{Bahturin}),
получаем, что \begin{equation}\label{EqVnGrDecomp}\frac{V^{G\text{-}\mathrm{gr}}_n}{V^{G\text{-}\mathrm{gr}}_n \cap \Id^{G\text{-}\mathrm{gr}}(L)}
\cong \bigoplus_{g_1, \ldots, g_n \in G} \frac{V_{g_1, \ldots, g_n}}{V_{g_1, \ldots, g_n} \cap \Id^{G\text{-}\mathrm{gr}}(L)}.\end{equation}

Пусть $\lbrace \gamma_1, \ldots, \gamma_m \rbrace := \lbrace g\in G \mid L^{(g)}\ne 0\rbrace$.
Это множество конечно в силу конечномерности
алгебры $L$. 
Введём обозначение $V_{n_1, \ldots, n_m} := V_{\bar g}$, где $$\bar g := (\underbrace{\gamma_1, \ldots, \gamma_1}_{n_1}, \underbrace{\gamma_2, \ldots, \gamma_2}_{n_2},
\ldots, \underbrace{\gamma_m, \ldots, \gamma_m}_{n_m}),\ n_i \in \mathbb Z_+.$$
 Тогда из~(\ref{EqVnGrDecomp}) следует, что
\begin{equation}\label{EqCn1nm}
c^{G\text{-}\mathrm{gr}}_n(L)=\sum_{n_1+\ldots+n_m = n} \binom{n}{n_1, \ldots, n_m}
c_{n_1, \ldots, n_m}(L),
\end{equation} где $c_{n_1, \ldots, n_m}(L) := \dim \frac{V_{n_1,\ldots,n_m}}{V_{n_1,\ldots,n_m} \cap \Id^{G\text{-}\mathrm{gr}}(L)}$.

 Пусть $G_0$~"--- некоторая подгруппа группы $G$.
 Если $n_i = 0$ для всех $\gamma_i \notin G_0$,
 то $c_{n_1, \ldots, n_m}(L)=c_{n_1, \ldots, n_m}(L_{G_0})$.
(В силу леммы~\ref{LemmaGrEquivCodimTheSame} неважно,
рассматриваем мы $G_0$-
или $G$-градуировку на $L_{G_0}$.)
 Фиксируем $n\in\mathbb N$ и рассмотрим множество $\Theta(G_0)$ всех таких наборов
$(n_1, \ldots, n_m)$, что $n_i\geqslant 0$, $n_1 + \ldots + n_m = n$ и $n_i = 0$ для всех $\gamma_i \notin G_0$.

Введём теперь на множестве $\Theta(G)$ дискретную меру
$\mu$ по формуле $$\mu\bigl((n_1, \ldots, n_m)\bigr)=\binom{n}{n_1, \ldots, n_m}
c_{n_1, \ldots, n_m}(L).$$ Тогда из~(\ref{EqCn1nm}) следует, что
$c^{G\text{-}\mathrm{gr}}_n(L_{G_0})=\mu(\Theta(G_0))$.
Согласно лемме~\ref{LemmaGradedNonZero}
в случае, когда $n_i,n_j \ne 0$ для некоторых таких $i,j$,
что $\gamma_i \gamma_j \ne \gamma_j \gamma_i$,
справедливо равенство $c_{n_1, \ldots, n_m}(L)=0$.
Обозначим множество всех таких наборов через $\Theta_0$. Тогда $\mu(\Theta_0)=0$.
Следовательно, всякая ненулевая коразмерность
$c_{n_1, \ldots, n_m}(L)$ равна $c_{n_1, \ldots, n_m}(L_{G_0})$
для некоторой конечнопорождённой абелевой подгруппы $G_0$ группы $G$. Предположим, что $G_1, \ldots, G_r$~"--- это все абелевы подгруппы группы $G$, порождённые подмножествами
множества $\lbrace \gamma_1, \ldots, \gamma_m \rbrace$. Тогда $\Theta(G)=\Theta_0 \cup \bigcup_{i=1}^r\Theta(G_i)$.
Используя формулу включений-исключений, получаем, что
 \begin{equation*}\begin{split}c^{G\text{-}\mathrm{gr}}_n(L)=\mu(\Theta(G))= \mu\left(\bigcup_{i=1}^r\Theta(G_i)\right)
= \\ = \sum_{i=1}^r \mu(\Theta(G_i))
- \sum_{i,j=1}^r \mu(\Theta(G_i) \cap \Theta(G_j)) + \ldots +(-1)^{r-1} 
\mu(\Theta(G_1) \cap \Theta(G_2) \cap \ldots \cap \Theta(G_r))= \\ =\sum_{i=1}^r \mu(\Theta(G_i))
- \sum_{i,j=1}^r \mu(\Theta(G_i \cap G_j)) + \ldots +(-1)^{r-1} 
\mu(\Theta(G_1 \cap G_2 \cap \ldots \cap G_r))
 =\\ =
 \sum_{i=1}^r c^{G\text{-}\mathrm{gr}}_n(L_{G_i})
- \sum_{i,j=1}^r c^{G\text{-}\mathrm{gr}}_n(L_{G_i \cap G_j}) + \ldots +(-1)^{r-1} 
c^{G\text{-}\mathrm{gr}}_n(L_{G_1 \cap G_2 \cap \ldots \cap G_r}).\end{split}\end{equation*}
\end{proof}
 
 \begin{proof}[Доказательство теоремы~\ref{TheoremMainLieGr}]
В силу теоремы~\ref{TheoremFinGenAbelianGr} всякое слагаемое в формуле~(\ref{EqInclExcl})
имеет требуемую асимптотику для некоторого $d$. Выберем среди всех таких $d$
наибольшее. Тогда при $d>0$ существуют такие $C_2 > 0$ и $r_2 \in\mathbb R$, что
$c_n^{G\text{-}\mathrm{gr}}(L)\leqslant C_2 n^{r_2} d^n$ для всех $n\in\mathbb N$, а
при $d=0$ существует такое $n_0 \in\mathbb N$, что $c_n^{G\text{-}\mathrm{gr}}(L)=0$
для всех $n\geqslant n_0$.
 
В силу выбора числа $d$ существует такая конечнопорождённая абелева подгруппа
$G_0 \subseteq G$, что $\PIexp^{G_0\text{-}\mathrm{gr}}(L_{G_0})=d$.
Отсюда существуют такие $C_1 > 0$ и $r_1 \in\mathbb R$, что
$$C_1 n^{r_1} d^n \leqslant c_n^{G_0\text{-}\mathrm{gr}}(L_{G_0}) = c_n^{G\text{-}\mathrm{gr}}(L_{G_0})\leqslant  c_n^{G\text{-}\mathrm{gr}}(L)\text{ для всех }n\in\mathbb N.$$
\end{proof}
 \begin{proof}[Доказательство теоремы~\ref{TheoremMainLieGrSum}]
Пусть $$L = L_1 \oplus \ldots \oplus L_q,$$ где $L_i$~"--- градуированные идеалы.
Во-первых,
$\PIexp^{G\text{-}\mathrm{gr}}(L)\geqslant  \PIexp^{G\text{-}\mathrm{gr}}(L_i)$ для всех $1 \leqslant i \leqslant q$,
поскольку $L_i$~"--- градуированные подалгебры алгебры $L$.
Во-вторых, в силу леммы~\ref{LemmaInclExcl} достаточно доказать, что
 $$\PIexp^{G\text{-}\mathrm{gr}}(L_{G_0})\leqslant\max_{1\leqslant i\leqslant q}  \PIexp^{G\text{-}\mathrm{gr}}(L_i)$$
 для всех конечнопорождённых абелевых подгрупп $G_0 \subseteq G$.
 Однако $L_{G_0} = (L_1)_{G_0} \oplus \ldots \oplus (L_q)_{G_0}$.
 Используя предложение~\ref{PropositionGrToHatG}, пример~\ref{ExampleHniceGr}
 и теорему~\ref{TheoremMainLieHSum}, получаем, что
 \begin{equation*}\begin{split}
 \PIexp^{G\text{-}\mathrm{gr}}(L_{G_0})
 =\PIexp^{\mathbbm{k}\hat G_0}(L_{G_0})=\max_{1\leqslant i\leqslant q}  \PIexp^{\mathbbm{k}\hat G_0}((L_i)_{G_0})
=\\ = \max_{1\leqslant i\leqslant q}  \PIexp^{G\text{-}\mathrm{gr}}((L_i)_{G_0}) \leqslant \max_{1\leqslant i\leqslant q}  \PIexp^{G\text{-}\mathrm{gr}}(L_i).\end{split}\end{equation*}
\end{proof}

\section{Рост дифференциальных тождеств}\label{SectionDiffLie}

В данном параграфе доказывается теорема~\ref{TheoremMainDiffLie}
об асимптотическом поведении коразмерностей дифференциальных тождеств
в алгебрах Ли с действием конечномерных полупростых алгебр Ли дифференцированиями,
которая будет получена как следствие следующего более общего результата:

\begin{theorem}\label{TheoremddprimeDiff}
Пусть $L$~"--- конечномерная алгебра Ли
над алгебраически замкнутым полем $\mathbbm{k}$ характеристики $0$. Предположим,
что некоторая алгебра Ли $\mathfrak g$ так действует на $L$ дифференцированиями,
что $L$ является $U(\mathfrak g)$-хорошей алгеброй Ли.
Пусть $d:=\PIexp(L)$~"--- PI-экспонента обычных тождеств алгебры Ли $L$.
Тогда 
\begin{enumerate}
\item при $d=0$ существует такое $n_0$, что $c_n^{U(\mathfrak g)}(L)=0$ при всех $n\geqslant n_0$;
\item при $d > 0$ существуют такие константы $C_1, C_2 > 0$, $r_1, r_2 \in \mathbb R$, что $$C_1 n^{r_1} d^n \leqslant c^{U(\mathfrak g)}_n(L)
   \leqslant C_2 n^{r_2} d^n\text{ для всех }n \in \mathbb N.$$
\end{enumerate}
В частности, $\PIexp(L) = \PIexp^{U(\mathfrak g)}(L)$.
\end{theorem}
\begin{remark}\label{RemarkReductActionLieHNice} Если связная редуктивная аффинная алгебраическая группа $G$
действует на $L$ автоморфизмами, то $L$ является
 $U(\mathfrak g)$-хорошей алгеброй (см. пример~\ref{ExampleHniceAffAlgDiff}). 
Как было показано в примере~\ref{ExampleHniceDiff},
любая конечномерная алгебра Ли
над алгебраически замкнутым полем $\mathbbm{k}$ характеристики $0$  с действием
некоторой конечномерной полупростой алгебры Ли дифференцированиями также является
$U(\mathfrak g)$-хорошей.
\end{remark}
\begin{proof}[Доказательство теоремы~\ref{TheoremddprimeDiff}]
В силу теоремы~\ref{TheoremMainLieH} дифференциальные
коразмерности удовлетворяют требуемым оценкам сверху и снизу
для $d=\PIexp^{U(\mathfrak g)}(L)=d'(L, U(\mathfrak g))$,
причём, кроме этого, ещё существует $\PIexp(L)=d'(L, \mathbbm{k})$.
Отсюда для доказательства теоремы достаточно показать, что $d'(L, U(\mathfrak g))=d'(L, \mathbbm{k})$.

Рассмотрим дифференциальные и обычные лиевские многочлены
как полилинейные функции на $L$ и получим, что $c_n(L) \leqslant c^{U(\mathfrak g)}_n(L)$ для всех $n \in\mathbb N$. Отсюда $\PIexp(L) \leqslant \PIexp^{U(\mathfrak g)}(L)$.

 Предположим, что $\mathfrak g$-инвариантные идеалы $I_1, I_2, \ldots, I_r$ и 
$J_1, J_2, \ldots, J_r$ алгебры Ли $L$, где $r \in \mathbb Z_+$ и $J_k \subseteq I_k$,
удовлетворяют условиям 1 и 2' при $H=U(\mathfrak g)$. 
В силу условия 2' существуют такие $\mathfrak g$-инвариантные $(\ad B)\oplus \tilde A_0$-подмодули $T_k$,
где $I_k = J_k\oplus T_k$, и числа $q_i \geqslant 0$, что $$\bigl[[T_1, \underbrace{L, \ldots, L}_{q_1}], [T_2, \underbrace{L, \ldots, L}_{q_2}], \ldots, [T_r,
 \underbrace{L, \ldots, L}_{q_r}]\bigr] \ne 0.$$
 Из леммы~\ref{LemmaLHBA0ComplReducible} следует, что алгебра Ли $L$ является вполне приводимым $(\ad B)\oplus \tilde A_0$-модулем, откуда $T_k=T_{k1}\oplus T_{k2}\oplus \ldots \oplus T_{kn_k}$
 для некоторых неприводимых $(\ad B)\oplus \tilde A_0$-подмодулей $T_{kj}$.
Отсюда можно выбрать такие $1 \leqslant j_k \leqslant n_k$, что
$$\bigl[[T_{1j_1}, \underbrace{L, \ldots, L}_{q_1}], [T_{2j_2}, \underbrace{L, \ldots, L}_{q_2}], \ldots, [T_{rj_r},
 \underbrace{L, \ldots, L}_{q_r}]\bigr] \ne 0.$$
 Пусть $\tilde I_k = T_{kj_k}\oplus J_k$.
 
 Докажем, что $\tilde I_k$ являются идеалами алгебры Ли $L$, причём $\Ann(\tilde I_k / J_k)=\Ann(I_k/J_k)$ для всех $1 \leqslant k \leqslant r$.
  Обозначим через $L_0$, $B_0$, $R_0$, $\mathfrak g_0$, соответственно,
образы алгебр Ли $L$, $B$, $R$, $\mathfrak g$ в $\mathfrak{gl}(I_k/J_k)$.
 Заметим, что алгебра Ли $B_0$ полупроста, а алгебра Ли $R_0$
 разрешима. Отсюда $L_0=B_0\oplus R_0$
 (прямая сумма $\mathfrak g$-подмодулей), где $\mathfrak g$-действие на $\mathfrak{gl}(I_k/J_k)$ 
индуцировано $\mathfrak g$-действием на $I_k/J_k$ 
и соответствует присоединённому представлению
алгебры Ли $\mathfrak g_0$ на
 $\mathfrak{gl}(I_k/J_k)$. 
В частности, $R_0$ является разрешимым идеалом алгебры Ли $(L_0+\mathfrak g_0)$,
а $B_0$ является идеалом алгебры Ли $(B_0 + \mathfrak g_0)$.
Из того, что $I_k/J_k$~"--- неприводимый $(U(\mathfrak g), L)$-модуль,
следует, что $I_k/J_k$~"--- неприводимый $(L_0+\mathfrak g_0)$-модуль.
В силу теоремы Э.~Картана (см., например, \cite[предложение~1.4.11]{GotoGrosshans}), $L_0+\mathfrak g_0 = B_1
 \oplus R_1$ (прямая сумма идеалов), где $B_1$~"--- полупростая алгебра Ли, а $R_1$ либо равно нулю,
 либо совпадает с центром $Z(\mathfrak{gl}(I_k/J_k))$, состоящим из  скалярных операторов.
  Отсюда идеал $R_0 \subseteq R_1$
состоит из скалярных операторов, а
$B_0$ является идеалом алгебры Ли $(L_0+\mathfrak g_0)$.
Рассматривая сюръективный гомоморфизм $(L_0+\mathfrak g_0) \twoheadrightarrow R_1$ с ядром~$B_1$,
 получаем, что $B_0 \subseteq B_1$.

    Поскольку $\tilde I_k/J_k \cong T_{k j_k}$ является неприводимым $(\ad B)\oplus \tilde A_0$-модулем, алгебра Ли  $\tilde A_0$ действует на $\tilde I_k/J_k$ скалярными операторами и
  $\tilde I_k/J_k$ является неприводимым $B_0$- и $L$-модулем. В частности, $\tilde I_k$ является идеалом.
  
  Если $\Ann(\tilde I_k/J_k)\ne \Ann(I_k/J_k)$, то $a\tilde I_k/J_k=0$
  для некоторого $a \in L_0 \cong L/\Ann(I_k/J_k)$, $a\ne 0$.
 Пусть $\varphi \colon L_0 \to \mathfrak{gl}(\tilde I_{k}/J_k)$~"--- соответствующее представление,
 а $a = b + c$, где $b\in B_0$, $c \in R_0$. 
 Тогда $\varphi(b)=-\varphi(c)$ является скалярным оператором на $\tilde I_{k}/J_k$.
 Следовательно, элемент $\varphi(b)$
принадлежит центру полупростой алгебры Ли $\varphi(B_0)$. Отсюда $\varphi(b)=\varphi(c)=0$, $b\ne 0$.
 Напомним, что алгебра Ли  $B_1$ полупроста.
 Следовательно, $B_1 =  B_0 \oplus B_2$ (прямая сумма идеалов)
 для некоторой полупростой алгебры Ли $B_2$. Поскольку $R_1$ состоит из скалярных операторов, $I_k/J_k$ является неприводимым $B_1$-модулем и $$I_k/J_k =  \sum_{\substack{a_i \in B_2,\\ \alpha\in\mathbb Z_+}} a_1 \ldots a_\alpha\ \tilde I_k/J_k.$$
 Теперь из $[b, B_2]=0$ и $b\tilde I_{k}/J_k=0$ следует, что $bI_k/J_k=0$ и $b=0$. Получаем противоречие,
 откуда $\Ann(\tilde I_k/J_k)=\Ann(I_k/J_k)$.
  
  Теперь заметим, что $\tilde I_1, \tilde I_2, \ldots, \tilde I_r$,
$J_1, J_2, \ldots, J_r$ удовлетворяют условиям 1 и 2' при $H=\mathbbm{k}$, т.е.
для случая обычных полиномиальных тождеств. Более того, $$\dim \frac{L}{\Ann(I_1/J_1) \cap \ldots \cap \Ann(I_r/J_r)}
= \dim \frac{L}{\Ann(\tilde I_1/J_1) \cap \ldots \cap \Ann(\tilde I_r/J_r)}.$$
Отсюда $\PIexp^{U(\mathfrak g)}(L)=\PIexp(L)$.
\end{proof}

\begin{proof}[Доказательство теоремы~\ref{TheoremMainDiffLie}]
Снова воспользуемся тем, что $H$-коразмерности не меняются при расширении основного поля.
(Доказательство полностью повторяет соответствующее доказательство 
для коразмерностей обычных тождеств ассоциативных алгебр~\cite[теорема~4.1.9]{ZaiGia} и
алгебр Ли~\cite[\S 2]{ZaicevLie}.)
Отсюда без ограничения общности можно считать основное поле $\mathbbm{k}$ алгебраически замкнутым.
Теперь достаточно применить теорему~\ref{TheoremddprimeDiff} и замечание~\ref{RemarkReductActionLieHNice}.
\end{proof}

\begin{remark}
Из теоремы~\ref{TheoremMainDiffLie} следует, что
обычные и дифференциальные коразмерности имеют схожее асимптотическое поведение.
Однако сами коразмерности могут не совпадать.
Рассмотрим присоединённое представление алгебры Ли $\mathfrak{sl}_2(\mathbbm{k})$ на себе самой. Тогда $c_1(\mathfrak{sl}_2(\mathbbm{k}))=1
< c_1^{U(\mathfrak{sl}_2(\mathbbm{k}))}(\mathfrak{sl}_2(\mathbbm{k}))$,
поскольку $U(\mathfrak{sl}_2(\mathbbm{k}))$-многочлены $x_1^{e_{11}-e_{22}}$ и $x_1^{e_{12}}$
являются линейно независимыми по модулю $\Id^{U(\mathfrak{sl}_2(\mathbbm{k}))}(\mathfrak{sl}_2(\mathbbm{k}))$.
\end{remark}

\begin{theorem}\label{TheoremLieGConnPIexpEqual}
Пусть $L$~"--- конечномерная алгебра Ли над алгебраически замкнутым полем $\mathbbm{k}$ характеристики $0$.
Предположим, что связная редуктивная аффинная алгебраическая группа $G$
рационально действует на $L$ автоморфизмами.
Тогда $\PIexp^{G}(L)=\PIexp(L)$.
\end{theorem}
\begin{proof}
Алгебра Ли $\mathfrak g$ группы $G$ в силу предложения~\ref{PropositionDerAutConnection}
действует на $L$ дифференцированиями.
В силу теоремы~\ref{TheoremAffAlgGrAllEquiv} и леммы~\ref{LemmaHEquivCodimTheSame} для всех $n\in\mathbb N$ справедливо равенство $c_n^{U(\mathfrak g)}(L)=c_n^{G}(L)$.
Теперь из теоремы~\ref{TheoremddprimeDiff} и замечания~\ref{RemarkReductActionLieHNice} следует, что $\PIexp^{G}(L)=\PIexp(L)$.
\end{proof}

\section{Рост $H$-коразмерностей в алгебрах Ли, в которых нильпотентный и разрешимый радикалы совпадают}
\label{SectionHRTheSame}

Оказывается, что в случае, когда разрешимый радикал $H$-модульной алгебры Ли $L$ является нильпотентным $H$-инвариантным идеалом, для доказательства аналога гипотезы Амицура можно не требовать,
чтобы для алгебры Ли выполнялись свойства~\ref{ConditionLevi}--\ref{ConditionLComplHred} из
определения $H$-хорошей алгебры Ли (см. \S\ref{SectionHnice}). 
Более того, в этом случае формула для $H$-PI-экспоненты оказывается
существенно проще, чем в общем случае (см. \S\ref{SectionHPIexpLie}).

\begin{theorem}\label{TheoremMainLieNRSame}
Пусть $L$~"--- конечномерная $H$-модульная алгебра Ли
над произвольным полем $\mathbbm{k}$ характеристики $0$,
где $H$~"--- некоторая алгебра Хопфа, причём разрешимый радикал алгебры Ли $L$
совпадает с её нильпотентным радикалом $N$ и является $H$-подмодулем.
Тогда 
\begin{enumerate}
\item либо существует такое $n_0$, что $c_n^H(L)=0$ при всех $n\geqslant n_0$;
\item либо существуют такие константы $C_1, C_2 > 0$, $r_1, r_2 \in \mathbb R$,
  $d \in \mathbb N$, что $$C_1 n^{r_1} d^n \leqslant c^{H}_n(L)
   \leqslant C_2 n^{r_2} d^n\text{ для всех }n \in \mathbb N.$$
\end{enumerate}
   В частности, существует $\PIexp^H(L)\in\mathbb Z_+$ и, таким образом,
   для $c^H_n(L)$ справедлив аналог гипотезы Амицура.

В случае, когда поле $\mathbbm{k}$ алгебраически замкнуто, константа  $d$ определяется
следующим образом.
Пусть $$L/N = B_1 \oplus \ldots \oplus B_q \text{ (прямая сумма $H$-инвариантных идеалов)},$$
где $B_i$~"--- $H$-простые алгебры Ли, а $\varkappa \colon L/N \to L$~"---
произвольный такой (необязательно $H$-линейный) гомоморфизм алгебр Ли, что $\pi\varkappa = \id_{L/N}$, где $\pi \colon L \twoheadrightarrow L/N$~"--- естественный сюръективный гомоморфизм. 
 Тогда $$d= \max\left( B_{i_1}\oplus B_{i_2} \oplus \ldots \oplus B_{i_r}
 \mathbin{\Bigl|}  r \geqslant 1,\ \bigl[[ H\varkappa(B_{i_1}), \underbrace{L, \ldots, L}_{q_1}], \right.$$ \begin{equation}\label{EqdRNSame}\left. [  H\varkappa(B_{i_2}), \underbrace{L, \ldots, L}_{q_2}], \ldots, [ H\varkappa(B_{i_r}),
 \underbrace{L, \ldots, L}_{q_r}]\bigr] \ne 0 \text{ для некоторых } q_i \geqslant 0 \right).\end{equation}
\end{theorem}
Теорема~\ref{TheoremMainLieNRSame} будет доказана в конце данного параграфа.
\begin{remark}
Существование разложения $$L/N = B_1 \oplus \ldots \oplus B_q \text{ (прямая сумма $H$-инвариантных идеалов)},$$
где $B_i$~"--- $H$-простые алгебры Ли, следует из теоремы~\ref{TheoremHLieSemiSimple}.
Существование отображения $\varkappa$ следует из обычной теоремы Леви.
\end{remark}
\begin{remark}
Заметим, что в силу теоремы~\ref{TheoremDerSimpleLie}
всякая дифференциально простая алгебра Ли проста.
В силу теоремы~\ref{TheoremGSimpleLie}
в случае, когда алгебра Ли $L$ является $G$-простой по отношению к рациональному действию
некоторой связной аффинной алгебраической группы~$G$ автоморфизмами, алгебра Ли $L$ также является
простой. Отсюда в случае, когда $R=N$, из теоремы~\ref{TheoremMainLieNRSame} получается другое
доказательство теорем~\ref{TheoremMainDiffLie}
и~\ref{TheoremLieGConnPIexpEqual}.
Действительно, из условий теорем~\ref{TheoremMainDiffLie}
и~\ref{TheoremLieGConnPIexpEqual} следует существование $H$-инвариантного разложения Леви, откуда можно выбрать отображение $\varkappa$, являющееся гомоморфизмом $H$-модулей.
\end{remark}

Используя предложение~\ref{PropositionGtoAut}, получаем такое следствие из теоремы~\ref{TheoremMainLieNRSame}:

\begin{corollary}
Пусть $L$~"--- конечномерная алгебра Ли
над произвольным полем $\mathbbm{k}$ характеристики $0$,
на которой действует автоморфизмами и антиавтоморфизмами некоторая
группа $G$, причём разрешимый радикал алгебры Ли $L$
совпадает с её нильпотентным радикалом $N$.
Тогда 
\begin{enumerate}
\item либо существует такое $n_0$, что $c_n^G(L)=0$ при всех $n\geqslant n_0$;
\item либо существуют такие константы $C_1, C_2 > 0$, $r_1, r_2 \in \mathbb R$,
  $d \in \mathbb N$, что $$C_1 n^{r_1} d^n \leqslant c^{G}_n(L)
   \leqslant C_2 n^{r_2} d^n\text{ для всех }n \in \mathbb N.$$
\end{enumerate}
   В частности, существует $\PIexp^G(L)\in\mathbb Z_+$ и, таким образом,
   для $c^G_n(L)$ справедлив аналог гипотезы Амицура.
\end{corollary}
\begin{corollary}
Пусть $L$~"--- конечномерная алгебра Ли
над произвольным полем $\mathbbm{k}$ характеристики $0$,
на которой действует дифференцированиями некоторая
алгебра Ли $\mathfrak g$, причём разрешимый радикал алгебры Ли $L$
совпадает с её нильпотентным радикалом $N$.
Тогда 
\begin{enumerate}
\item либо существует такое $n_0$, что $c_n^{U(\mathfrak g)}(L)=0$ при всех $n\geqslant n_0$;
\item либо существуют такие константы $C_1, C_2 > 0$, $r_1, r_2 \in \mathbb R$,
  $d \in \mathbb N$, что $$C_1 n^{r_1} d^n \leqslant c^{U(\mathfrak g)}_n(L)
   \leqslant C_2 n^{r_2} d^n\text{ для всех }n \in \mathbb N.$$
\end{enumerate}
   В частности, существует $\PIexp^{U(\mathfrak g)}(L)\in\mathbb Z_+$ и, таким образом,
   для $c^{U(\mathfrak g)}_n(L)$ справедлив аналог гипотезы Амицура.
\end{corollary}
\begin{proof} В силу следствия~\ref{CorollaryLieRadicalsDiffInvariant}
радикал инвариантен относительно всех дифференцирований.
Теперь достаточно применить теорему~\ref{TheoremMainLieNRSame}.
\end{proof}

Фиксируем такое число~$p\in\mathbb N$, что $N^p=0$.
Тогда дословно повторив рассуждения из леммы~\ref{LemmaAssocUpperCochar}, получаем следующий результат:

\begin{lemma}\label{LemmaNRSameUpperCochar}
 Пусть  $\lambda = (\lambda_1, \ldots, \lambda_s) \vdash n$
для некоторого $n\in\mathbb N$, причём $\sum_{k=d+1}^s \lambda_k \geqslant p$. Тогда $m(L, H, \lambda)=0$. 
\end{lemma}

Дословно повторив рассуждения из теоремы~\ref{TheoremAssocUpper},
получаем оценку сверху:

\begin{lemma}\label{LemmaNRSameUpper} 
Если $d > 0$, то существуют такие константы $C_2 > 0$, $r_2 \in \mathbb R$,
что $c^H_n(L) \leqslant C_2 n^{r_2} d^n$
для всех $n \in \mathbb N$. В случае, когда $d=0$, алгебра Ли $L$ нильпотентна.
\end{lemma}

Докажем теперь существование $H$-многочлена, кососимметричного по достаточному числу наборов переменных:

\begin{lemma}\label{LemmaHRSameLowerPolynomial}
Предположим, что основное поле $\mathbbm{k}$ алгебраически замкнуто, а $L$, $N$, $\varkappa$, $B_i$ и $d$ те же, что и в теореме~\ref{TheoremMainLieNRSame}.
Тогда если $d > 0$, то существует такое число $n_0 \in \mathbb N$, что для любого $n\geqslant n_0$
существуют попарно непересекающиеся наборы переменных $X_1$, \ldots, $X_{2k} \subseteq \lbrace x_1, \ldots, x_n\rbrace$, где $k := \left[\frac{n-n_0}{2d}\right]$,
$|X_1| = \ldots = |X_{2k}|=d$, и $H$-многочлен $f \in V^H_n \backslash
\Id^H(L)$, кососимметричный по переменным каждого множества $X_j$.
\end{lemma}
\begin{proof} Без ограничения общности можно считать, что $$d = \dim(B_1 \oplus B_2 \oplus \ldots \oplus B_r),$$ где
$$\bigl[[  H\varkappa(B_1), a_{11}, \ldots, a_{1q_1}], [  H\varkappa(B_2), a_{21}, \ldots, a_{2q_2}], \ldots, [ H\varkappa(B_r),
 a_{r1}, \ldots, a_{rq_r}]\bigr] \ne 0$$ для некоторых $q_i\geqslant 0$
 и $a_{kj}\in L$. 
 Поскольку идеал $N$ нильпотентен, мы можем увеличить числа $q_i$, добавляя ко множествам $\lbrace a_{ij} \rbrace$ необходимое количество элементов идеала $N$ так, что
 для некоторых $q_i\geqslant 0$, $b_i \in B_i$, $\gamma_i \in H$
 справедливо неравенство
$$\bigl[[ {\gamma_1}\varkappa(b_1), a_{11}, \ldots, a_{1q_1}], [ {\gamma_2}\varkappa(b_2), a_{21}, \ldots, a_{2q_2}], \ldots, [{\gamma_r}\varkappa(b_r),
 a_{r1}, \ldots, a_{rq_r}]\bigr] \ne 0,$$ 
 однако
 \begin{equation}\label{Eqbazero}\bigl[[ \tilde b_1, a_{11}, \ldots, a_{1q_1}], [ \tilde b_2, a_{21}, \ldots, a_{2q_2}], \ldots, [\tilde b_r,
 a_{r1}, \ldots, a_{rq_r}]\bigr] = 0 \end{equation} для всех таких $t_i \geqslant 0$, $\tilde b_i\in [H\varkappa(B_i), \underbrace{L, \ldots, L}_{t_i}]$,
 что $\tilde b_j\in [H\varkappa(B_j), L, \ldots, L, N, L, \ldots, L]$ хотя бы для одного $j$.
 
 Напомним, что $\varkappa$ является гомоморфизмом алгебр Ли.
 Более того, из $\pi(h\varkappa(a)-\varkappa(ha))=0$ 
 следует, что $$h\varkappa(a)-\varkappa(ha) \in N\text{ для всех }a\in L\text{ и }h\in H.$$
 В силу~(\ref{Eqbazero}) получаем отсюда, что если в $$\bigl[[ \gamma_1\varkappa(b_1), a_{11}, \ldots, a_{1q_1}], [  \gamma_2\varkappa(b_2), a_{21}, \ldots, a_{2q_2}], \ldots, [ \gamma_r\varkappa(b_r),
 a_{r1}, \ldots, a_{rq_r}]\bigr]$$
заменить $\varkappa(b_i)$ на коммутатор элемента $\varkappa(b_i)$
и некоторого выражения с участием $\varkappa$, то отображение $\varkappa$
будет вести себя как гомоморфизм $H$-модулей. Это свойство будет использовано ниже.

Поскольку $B_i$ являются неприводимыми $(H, \ad B_i)$-модулями, в силу теоремы~\ref{TheoremLieAlternateFinal}
существуют такие константы $m_i \in \mathbb Z_+$,
что для любого $k$ существуют полилинейные ассоциативные $H$-многочлены $f_i$ степени $(2kd_i + m_i)$,
где $d_i := \dim B_i$, кососимметричные
по переменным из непересекающихся наборов
$X^{(i)}_{\ell}$, где $1 \leqslant \ell \leqslant 2k$, $|X^{(i)}_{\ell}|=d_i$,
причём для каждого $i$ существуют такие элементы из $\ad B_i$, что $f_i$ не обращается в нуль при подстановке этих элементов вместо своих переменных.

Используя тот факт, что $B_i$ являются неприводимыми $(H, \ad B_i)$-модулями, ещё раз,
получаем из теоремы плотности, что алгебры
 $\End_\mathbbm{k}(B_i)$ порождены операторами из алгебр~$H$ и~$\ad B_i$.
 Произведём отождествление алгебр $\End_\mathbbm{k}(B_i)$ и $M_{d_i}(\mathbbm{k})$
 при помощи изоморфизма $\End_\mathbbm{k}(B_i) \cong M_{d_i}(\mathbbm{k})$.
 Тогда любая матричная единица $e^{(i)}_{j\ell} \in M_{d_i}(\mathbbm{k})$
может быть представлена в виде многочлена от операторов
  из алгебр~$H$ и~$\ad B_i$.
  Выберем такие многочлены для всех $i$ и всех матричных единиц. 
Обозначим через $m_0$ наибольшую среди степеней всех таких многочленов.

Пусть $n_0 := r(2m_0+1)+ \sum_{i=1}^r (m_i+q_i)$.
Выберем теперь $H$-многочлены $f_i$, о которых шла речь выше, для $k = \left[\frac{n-n_0}{2d}\right]$.
Кроме того, снова пользуясь теоремой~\ref{TheoremLieAlternateFinal},
выберем $\tilde f_1$ для $\tilde k = \left[\frac{n-2kd-m_1}{2d_1}\right]+1$ и алгебры Ли $B_{i_1}$.
Тогда $H$-многочлены $f_i$ дадут необходимое число наборов переменных, по которым результирующий многочлен
будет кососимметричен, однако общая степень произведения таких $H$-многочленов
может оказаться меньше $n$. Для того чтобы увеличить число переменных
и получить $H$-многочлен степени $n$, мы воспользуемся $H$-многочленом $\tilde f_1$.

В силу теоремы~\ref{TheoremLieAlternateFinal}
существуют такие элементы $\bar x_{i1}, \ldots, \bar x_{i, 2k d_i+m_i} \in B_i$,
что $$f_i(\ad \bar x_{i1}, \ldots, \ad \bar x_{i, 2k d_i+m_i})\ne 0,$$
и такие элементы $\bar x_1, \ldots, \bar x_{2\tilde k d_1+m_1} \in B_1$,
что $\tilde f_1(\ad \bar x_1, \ldots, \ad \bar x_{2\tilde k d_1+m_1}) \ne 0$.
В частности, $$e^{(i)}_{\ell_i \ell_i} f_i(\ad \bar x_{i1}, \ldots, \ad \bar x_{i, 2k d_i+m_i})
e^{(i)}_{s_i s_i} \ne 0$$
и $$e^{(1)}_{\tilde\ell \tilde\ell}\tilde f_1(\ad \bar x_1, \ldots, \ad \bar x_{2\tilde k d_1+m_1})
e^{(1)}_{\tilde s \tilde s} \ne 0$$
 для некоторых $e^{(i)}_{\ell_i \ell_i},
e^{(i)}_{s_i s_i} \in \End_\mathbbm{k}(B_i)$, $1 \leqslant \ell_i, s_i \leqslant d_i$,
$e^{(1)}_{\tilde\ell \tilde\ell}, e^{(1)}_{\tilde s \tilde s} \in \End_\mathbbm{k}(B_1)$, $1 \leqslant \tilde \ell,
\tilde s \leqslant d_1$.
Отсюда $$\sum_{\ell=1}^{d_i}
e^{(i)}_{\ell \ell_i} f_i(\bar x_{i1}, \ldots, \bar x_{i, 2k d_i+m_i})
 e^{(i)}_{s_i \ell}$$ является ненулевым скалярным оператором, принадлежащим алгебре $\End_\mathbbm{k}(B_i)$.

Следовательно,
\begin{equation*}\begin{split} [[\gamma_1\varkappa\left(\sum_{\ell=1}^{d_1}
e^{(1)}_{\ell \ell_1} f_1(\ad \bar x_{11}, \ldots, \ad \bar x_{1,2k d_1+m_1})
e^{(1)}_{s_1 \tilde \ell} \tilde f_1(\ad \bar x_1, \ldots, \ad \bar x_{2\tilde k d_1+m_1})
 e^{(1)}_{\tilde s \ell}b_1\right), a_{11}, \ldots, a_{1q_1}],\\
  [\gamma_2\varkappa\left(\sum_{\ell=1}^{d_2}
e^{(2)}_{\ell \ell_2} f_2(\ad \bar x_{21}, \ldots, \ad \bar x_{2,2k d_2+m_2})
 e^{(2)}_{s_2 \ell}b_2\right), a_{21}, \ldots, a_{2q_2}],
 \ldots, \\
 [\gamma_r\varkappa\left(\sum_{\ell=1}^{d_r}
e^{(r)}_{\ell \ell_r} f_r(\ad \bar x_{r1}, \ldots, \ad \bar x_{r, 2k d_r+m_r})
 e^{(r)}_{s_r \ell}b_r\right), a_{r1}, \ldots, a_{rq_r}]]\ne 0.\end{split}\end{equation*}

Представим теперь матричные единицы
$e^{(i)}_{\ell j}$ в виде многочленов от элементов $\ad B_i$
и $H$.
Используя линейность выражений по $e^{(i)}_{\ell j}$,
можно заменить  $e^{(i)}_{\ell j}$ такими произведениями
элементов алгебр $\ad B_i$
и $H$, что выражение не обратится в нуль.
По определению $H$-модульной алгебры
для всех $h\in H$ и $a, b \in L$
справедливо равенство
$h(\ad a )b=\ad (h_{(1)}a)(h_{(2)}b)$.
Пользуясь этим равенством, переместим
все элементы из $H$ вправо.
Как было отмечено выше, отображение $\varkappa$ является гомоморфизмом алгебр Ли и в силу~(\ref{Eqbazero})
ведёт себя как гомоморфизм $H$-модулей.
Отсюда получаем, что
\begin{equation*}\begin{split}  a_0 := \biggl[\Bigl[\gamma_1\Bigl[\bar y_{11}, [\bar y_{12}, \ldots
 [\bar y_{1 \alpha_1}, 
\\
  \bigl(h_1 f_1(\ad \varkappa(\bar x_{11}), \ldots, \ad \varkappa(\bar x_{1, 2k d_1+m_1}))\bigr)
 [\bar w_{11}, [\bar w_{12}, \ldots, [\bar w_{1 \theta_1},\\
 \bigl(\tilde h \tilde f_1(\ad \varkappa(\bar x_1), \ldots, \ad \varkappa(\bar x_{2\tilde k d_1+m_1}))\bigr)
 [\bar w_{1}, [\bar w_{2}, \ldots, [\bar w_{\tilde \theta},
  \varkappa({h'_1}b_1)]\ldots \Bigr],
  a_{11}, \ldots, a_{1q_1}\Bigr], \\
  \Bigl[\gamma_2\Bigl[\bar y_{21}, [\bar y_{22}, \ldots
 [\bar y_{2 \alpha_2}, \\
  \bigl(h_2 f_2(\ad \varkappa(\bar x_{21}), \ldots, \ad \varkappa(\bar x_{2, 2k d_2+m_2}))\bigr)
 [\bar w_{21}, [\bar w_{22}, \ldots, [\bar w_{2 \theta_2},
  \varkappa({h'_2}b_2)]\ldots \Bigr],
  a_{21}, \ldots, a_{2q_2}\Bigr],
 \ldots, \\ \Bigl[\gamma_r\Bigl[\bar y_{r1}, [\bar y_{r2}, \ldots,
 [\bar y_{r \alpha_r},  \\
 \bigl(h_r f_r(\ad \varkappa(\bar x_{r1}), \ldots, \ad \varkappa(\bar x_{r, 2k d_r+m_r}))\bigr)
 [\bar w_{r1}, [\bar w_{r2}, \ldots, [\bar w_{r \theta_r}, \varkappa({h'_r}b_r)]\ldots \Bigr],
  a_{r1}, \ldots, a_{rq_r}\Bigr]\biggr]\ne 0\end{split}\end{equation*}
 для некоторых  $0 \leqslant \alpha_i, \theta_i, \tilde \theta \leqslant m_0$,
  \quad $h_i, h'_i, \tilde h \in H$,\quad $\bar y_{ij}, \bar w_{ij} \in \varkappa(B_i)$,
  \quad $\bar w_j \in \varkappa(B_1)$.
 Здесь мы считаем, что каждый $H$-многочлен $f_i$ зависит от переменных $x_{i1}, \ldots,
x_{i,2k d_i+m_i}$, а  $H$-многочлен $\tilde f_1$ зависит от переменных  $x_1, \ldots, x_{2\tilde k d_1 + m_1}$.

Введём обозначение $X_\ell := \bigcup_{i=1}^{r} X^{(i)}_{\ell}$,
где $X^{(i)}_{\ell}$~"--- множества, по переменным которых кососимметричны $H$-многочлены $f_i$.
Обозначим через $\Alt_\ell$
операторы альтернирования по переменным из множества $X_\ell$. 

Рассмотрим выражение
\begin{equation*}\begin{split}\hat f :=
 \Alt_1 \Alt_2 \ldots \Alt_{2k} \biggl[\Bigl[\gamma_1\Bigl[y_{11}, [y_{12}, \ldots
 [y_{1 \alpha_1}, 
\\
  \bigl(h_1 f_1(\ad x_{11}, \ldots, \ad x_{1, 2k d_1+m_1})\bigr)
 [w_{11}, [w_{12}, \ldots, [w_{1 \theta_1},\\
 \bigl(\tilde h \tilde f_1(\ad x_1, \ldots, \ad x_{2\tilde k d_1+m_1})\bigr)
 [w_{1}, [w_{2}, \ldots, [w_{\tilde \theta},
  z_1]\ldots \Bigr],
  u_{11}, \ldots, u_{1q_1}\Bigr],\\ \Bigl[\gamma_2\Bigl[y_{21}, [y_{22}, \ldots
 [y_{2 \alpha_2}, \\
  \bigl(h_2 f_2(\ad x_{21}, \ldots, \ad x_{2, 2k d_2+m_2})\bigr)
 [w_{21}, [w_{22}, \ldots, [w_{2 \theta_2},
  z_2]\ldots \Bigr],
  u_{21}, \ldots, u_{2q_2}\Bigr],
 \ldots, \\ \Bigl[\gamma_r\Bigl[y_{r1}, [y_{r2}, \ldots,
 [y_{r \alpha_r},  \\
 \bigl(h_r f_r(\ad x_{r1}, \ldots, \ad x_{r, 2k d_r+m_r})\bigr)
 [w_{r1}, [w_{r2}, \ldots, [w_{r \theta_r}, z_r]\ldots \Bigr],
  u_{r1}, \ldots, u_{rq_r}\Bigr]\biggr].\end{split}\end{equation*}

Тогда значение выражения $\hat f$
при подстановке
$z_i=\varkappa({h'_i}b_i)$, $u_{i\ell}=a_{i\ell}$,
 $x_{i\ell}=\varkappa(\bar x_{i\ell})$, $x_i = \varkappa(\bar x_i)$, $y_{i\ell}=\bar y_{i\ell}$, $w_{i\ell}=\bar w_{i\ell}$, $w_i = \bar w_i$
 равно $(d_1!)^{2k} \ldots (d_r!)^{2k} a_0 \ne 0$,
 так как $f_i$ кососимметричны по переменным каждого из множеств $X^{(i)}_{\ell}$, $[B_i, B_\ell] = 0$ при $i \ne \ell$, а $\varkappa$ является гомоморфизмом алгебр Ли.
 
 Следовательно, \begin{equation*}\begin{split}f_0 := 
  \Alt_1 \Alt_2 \ldots \Alt_{2k} \biggl[\Bigl[\gamma_1\Bigl[y_{11}, [y_{12}, \ldots
 [y_{1 \alpha_1}, 
\\
  \bigl(h_1 f_1(\ad x_{11}, \ldots, \ad x_{1, 2k d_1+m_1})\bigr)
 [w_{11}, [w_{12}, \ldots, [w_{1 \theta_1},
  z_1]\ldots \Bigr],
  u_{11}, \ldots, u_{1q_1}\Bigr],\\ \Bigl[\gamma_2\Bigl[y_{21}, [y_{22}, \ldots
 [y_{2 \alpha_2}, \\
  \bigl(h_2 f_2(\ad x_{21}, \ldots, \ad x_{2, 2k d_2+m_2})\bigr)
 [w_{21}, [w_{22}, \ldots, [w_{2 \theta_2},
  z_2]\ldots \Bigr],
  u_{21}, \ldots, u_{2q_2}\Bigr],
 \ldots, \\
 \Bigl[\gamma_r\Bigl[y_{r1}, [y_{r2}, \ldots,
 [y_{r \alpha_r},  \\
\bigl(h_r f_r(\ad x_{r1}, \ldots, \ad x_{r, 2k d_r+m_r})\bigr)
 [w_{r1}, [w_{r2}, \ldots, [w_{r \theta_r}, z_r]\ldots \Bigr],
  u_{r1}, \ldots, u_{rq_r}\Bigr]\biggr]\end{split}\end{equation*}
   не обращается в нуль при подстановке
 $$z_1 = \bigl(\tilde h \tilde f_1(\ad \varkappa(\bar x_1), \ldots, \ad \varkappa(\bar x_{2\tilde k d_1+m_1}))\bigr)
 [\bar w_{1}, [\bar w_{2}, \ldots, [\bar w_{\tilde \theta},
 \varkappa(h'_1 b_1)]\ldots],$$
  $z_i=\varkappa(h'_i b_i)$ при $2 \leqslant i \leqslant r$; $u_{i\ell}=a_{i\ell}$,
 $x_{i\ell}=\varkappa(\bar x_{i\ell})$, $y_{i\ell}=\bar y_{i\ell}$, $w_{i\ell}=\bar w_{i\ell}$.
 
Заметим, что $f_0 \in V_{\tilde n}^H$, где
  $$\tilde n: = 2kd +r+ \sum_{i=1}^r (m_i + q_i + \alpha_i+\theta_i)
  \leqslant n.$$ Если $n=\tilde n$, положим $f:=f_0$.
  Предположим, что $n > \tilde n$.
Учитывая, что $$\bigl(\tilde h \tilde f_1(\ad \varkappa(\bar x_1), \ldots, \ad \varkappa(\bar x_{2\tilde k d_1+m_1}))\bigr)
 [\bar w_{1}, [\bar w_{2}, \ldots, [\bar w_{\tilde \theta},
  \varkappa(h'_1 b_1)]\ldots]$$ является линейной комбинацией длинных коммутаторов,
  каждый из которых содержит по крайней мере $2\tilde k d_1+m_1+1 > n-\tilde n+1$
  элементов алгебры Ли $L$, получаем, что $H$-многочлен $ f_0$ не обращается в нуль при подстановке
 $z_1 = [\bar v_1, [\bar v_2, [\ldots, [\bar v_\theta,  \varkappa(h'_1 b_1)]\ldots]$
 для некоторых $\theta \geqslant n-\tilde n$, $\bar v_i \in L$;
  $z_i=\varkappa(h'_i b_i)$ при $2 \leqslant i \leqslant r$; $u_{i\ell}=a_{i\ell}$,
 $x_{i\ell}=\varkappa(\bar x_{i\ell})$, $y_{i\ell}=\bar y_{i\ell}$,
  $w_{i\ell}=\bar w_{i\ell}$.
Следовательно, \begin{equation*}\begin{split}f :=  \Alt_1 \Alt_2 \ldots \Alt_{2k} \biggl[\Bigl[\gamma_1\Bigl[y_{11}, [y_{12}, \ldots
 [y_{1 \alpha_1}, 
\\
  \bigl(h_1 f_1(\ad x_{11}, \ldots, \ad x_{1, 2k d_1+m_1})\bigr)
 [w_{11}, [w_{12}, \ldots, [w_{1 \theta_1},\\ 
 \bigl[v_1, [v_2, [\ldots, [v_{n-\tilde n}, z_1]\ldots\bigr]\ldots \Bigr],
  u_{11}, \ldots, u_{1q_1}\Bigr],
  \\ \Bigl[\gamma_2\Bigl[y_{21}, [y_{22}, \ldots
 [y_{2 \alpha_2}, \\
  \bigl(h_2 f_2(\ad x_{21}, \ldots, \ad x_{2, 2k d_2+m_2})\bigr)
 [w_{21}, [w_{22}, \ldots, [w_{2 \theta_2},
  z_2]\ldots \Bigr],
  u_{21}, \ldots, u_{2q_2}\Bigr],
 \ldots, \\
 \Bigl[\gamma_r\Bigl[y_{r1}, [y_{r2}, \ldots,
 [y_{r \alpha_r},  \\
 \bigl(h_r f_r(\ad x_{r1}, \ldots, \ad x_{r, 2k d_r+m_r}) \bigr)
 [w_{r1}, [w_{r2}, \ldots, [w_{r \theta_r}, z_r]\ldots \Bigr],
  u_{r1}, \ldots, u_{rq_r}\Bigr]\biggr]\end{split}\end{equation*}
  не обращается в нуль при подстановке
  $v_\ell = \bar v_\ell$, $1 \leqslant \ell \leqslant n-\tilde n$,
  $$z_1 = [\bar v_{n-\tilde n +1}, [\bar v_{n-\tilde n +2}, [\ldots, [\bar v_\theta,  \varkappa(h'_1 b_1)]\ldots];$$
  $z_i=\varkappa(h'_i b_i)$ при $2 \leqslant i \leqslant r$; $u_{i\ell}=a_{i\ell}$,
 $x_{i\ell}=\varkappa(\bar x_{i\ell})$, $y_{i\ell}=\bar y_{i\ell}$, $w_{i\ell}=\bar w_{i\ell}$.
 Теперь осталось заметить, что $H$-многочлен $f$ принадлежит множеству $V_n^H$ и удовлетворяет всем условиям леммы.
\end{proof}

Повторяя рассуждения леммы~\ref{LemmaAssocCochar}, получаем из леммы~\ref{LemmaHRSameLowerPolynomial}
следующее утверждение:

\begin{lemma}\label{LemmaHRSameCochar} Пусть $k, n_0$~"--- числа из леммы~\ref{LemmaHRSameLowerPolynomial}. Тогда для любого $n \geqslant n_0$
существует такое разбиение $\lambda = (\lambda_1, \ldots, \lambda_s) \vdash n$,
где $\lambda_i > 2k-p$ для всех $1 \leqslant i \leqslant d$,
что $m(L, H, \lambda) \ne 0$.
Здесь $p \in \mathbb N$~"--- такое число, что $N^p=0$.
\end{lemma}

Теперь мы можем доказать основной результат параграфа:

\begin{proof}[Доказательство теоремы~\ref{TheoremMainLieNRSame}]
Пусть $K \supset \mathbbm{k}$~"--- некоторое расширение поля~$\mathbbm{k}$.
Тогда $$(L \otimes_\mathbbm{k} K)/(N \otimes_\mathbbm{k} K) \cong (L/N) \otimes_\mathbbm{k} K$$
снова является полупростой алгеброй Ли, а идеал  $N \otimes_\mathbbm{k} K$ по-прежнему нильпотентен.
Как уже было неоднократно отмечено,  $H$-коразмерности не меняются при расширении основного поля $\mathbbm{k}$.
Отсюда без ограничения общности можно считать, что поле $\mathbbm{k}$ алгебраически замкнуто.
Теперь достаточно, использовав леммы~\ref{LemmaNRSameUpper} и~\ref{LemmaHRSameCochar}, повторить рассуждения теоремы~\ref{TheoremAssocBounds}.
\end{proof}

\section{Примеры и критерии простоты}
\label{SectionLieExamples}

   \begin{example}\label{ExampleHLieSemiSimplePIexp}
 Пусть $L$~"--- полупростая в обычном смысле $H$-модульная алгебра Ли
 над полем $\mathbbm{k}$ характеристики $0$ для некоторой
 алгебры Хопфа $H$.  Тогда существуют такие $d\in\mathbb Z_+$, $C_1, C_2 > 0$, $r_1, r_2 \in \mathbb R$,
  что $$C_1 n^{r_1} d^n \leqslant c_n^{H}(L)
  \leqslant C_2 n^{r_2} d^n \text{ для всех } n\in\mathbb N,$$
  причём если поле~$\mathbbm{k}$ алгебраически замкнуто, то
  $d = \max_{1 \leqslant k
  \leqslant q} \dim B_k$, где $B_i$~"--- конечномерные $H$-простые
 алгебры Ли из разложения $L=B_1 \oplus B_2 \oplus \ldots \oplus B_q$ (прямая сумма $H$-инвариантных идеалов), см. теорему~\ref{TheoremHLieSemiSimple}.
   \end{example}
\begin{proof}
Как и прежде, коразмерности не меняются при расширении основного поля $\mathbbm{k}$, откуда можно
считать, что поле $\mathbbm{k}$ алгебраически замкнуто. Теперь достаточно применить теорему~\ref{TheoremMainLieNRSame}.
\end{proof}

Отсюда получаем следующий критерий $H$-простоты:

\begin{theorem}\label{TheoremHCrSimpleLie}
Пусть $L$~"--- полупростая или $H$-хорошая алгебра Ли для некоторой алгебры Хопфа $H$
над алгебраически замкнутым полем $\mathbbm{k}$ характеристики $0$.
Тогда $\PIexp^H(L)=\dim L$, если и только если алгебра Ли $L$ полупроста и $H$-проста.
\end{theorem}
\begin{proof} В случае, когда $L$~"--- $H$-простая полупростая алгебра Ли,
равенство $\PIexp^H(L)=\dim L$ следует из примера~\ref{ExampleHLieSemiSimplePIexp}.

Если $L$~"--- полупростая алгебра Ли, не являющаяся $H$-простой, то в силу примера~\ref{ExampleHLieSemiSimplePIexp}
справедливо неравенство $\PIexp^H(L)<\dim L$. Поэтому можно считать, что
алгебра Ли $L$ является $H$-хорошей. 
Пусть $\PIexp^H(L)=\dim L$.
Обозначим через $N$ нильпотентный радикал алгебры Ли $L$. Пусть $I_1, \ldots, I_r$, $J_1, \ldots, J_r$~"---
$H$-инвариантные идеалы алгебры Ли $L$, удовлетворяющие условиям 1--2 из \S\ref{SectionHPIexpLie}.
В силу леммы~\ref{LemmaIrrAnnBQ} справедливо включение $N \subseteq \Ann(I_1/J_1) \cap \ldots \cap \Ann(I_r/J_r)$, откуда $\PIexp^H(L) \leqslant (\dim L) - (\dim N)$. Следовательно, $N=0$
и в силу предложения~2.1.7 из~\cite{GotoGrosshans}
выполняется условие
$[L, R] \subseteq N = 0$, где $R$~"--- разрешимый радикал алгебры Ли $L$. Отсюда $R = Z(L)\subseteq N=0$ 
и алгебра Ли $L$ полупроста. Этот случай был уже разобран выше.
\end{proof}

Обратимся теперь к случаю градуированных алгебр:

   \begin{example}\label{ExampleGrLieSemiSimplePIexp}
 Пусть $L=B_1 \oplus B_2 \oplus \ldots \oplus B_q$ (прямая сумма $G$-градуированных идеалов), 
где $B_i$~"--- конечномерные $G$-градуированно простые алгебры Ли над алгебраически замкнутым полем $\mathbbm{k}$
для некоторой группы $G$. Введём обозначение $d := \max_{1 \leqslant k
  \leqslant q} \dim B_k$.  Тогда существуют такие $C_1, C_2 > 0$, $r_1, r_2 \in \mathbb R$,
  что $$C_1 n^{r_1} d^n \leqslant c_n^{G\text{-}\mathrm{gr}}(L)
  \leqslant C_2 n^{r_2} d^n \text{ для всех } n\in\mathbb N.$$
   \end{example}
   \begin{proof}
   Пусть $B_i^{(g)}\ne 0$ для некоторого $g\in G$ и $1 \leqslant k
  \leqslant q$. Тогда в силу своей градуированной простоты
  всякая алгебра Ли $B_i$ представляется в виде суммы всевозможных ненулевых
  коммутаторов $[B_i^{(g)},B_i^{(g_1)}, \ldots, B_i^{(g_t)}]$,
  где $t\in\mathbb Z_+$, а $g_i \in G$. В силу леммы~\ref{LemmaGradedNonZero}
  все такие элементы $g, g_1, \ldots, g_t$ коммутируют между собой,
  откуда $g$ коммутирует со всеми элементами носителя градировки на $B_i$.
  В силу произвольности элемента $g$ все элементы носителя градировки также коммутируют между собой.
  Следовательно, на любой из градуированно простых алгебр Ли $B_i$
  можно заменить градуирующую группу $G$ на некоторую конечнопорождённую абелеву группу $G_i$,
  причём в силу леммы~\ref{LemmaGrEquivCodimTheSame}
  градуированные коразмерности алгебры Ли $B_i$ при этом не меняются.
  Из леммы~\ref{LemmaAbelianDual} следует, что алгебры Ли $B_i$
  являются $\mathbbm{k}\hat G_i$-простыми. 
   Отсюда из предложения~\ref{PropositionGrToHatG} и примера~\ref{ExampleHLieSemiSimplePIexp} получаем, что $$\PIexp^{G\text{-}\mathrm{gr}}(B_i)=\PIexp^{G_i\text{-}\mathrm{gr}}(B_i)=\PIexp^{\mathbbm{k}\hat G_i}(B_i)=\dim B_i.$$
   Теперь достаточно применить теорему~\ref{TheoremMainLieGrSum}.
   \end{proof}
   
Отсюда получаем следующий критерий градуированной простоты в терминах градуированных PI-экспонент:

\begin{theorem}\label{TheoremGrCrSimpleLie}
Пусть $L$~"--- конечномерная алгебра Ли над алгебраически замкнутым полем $\mathbbm{k}$ характеристики $0$,
градуированная некоторой группой $G$.
Тогда $\PIexp^{G\text{-}\mathrm{gr}}(L)=\dim L$, если и только если алгебра Ли $L$ градуированно проста.
\end{theorem}
\begin{proof} Сперва предположим, что нильпотентный радикал $N$ алгебры Ли $L$ не равен $0$,
а градуирующая группа конечнопорождённая абелева.
Пусть $I_1, \ldots, I_r$, $J_1, \ldots, J_r$~"---
$\mathbbm{k}\hat G$-инвариантные идеалы алгебры Ли $L$, удовлетворяющие условиям 1--2 из \S\ref{SectionHPIexpLie}.
В силу леммы~\ref{LemmaIrrAnnBQ} $$N \subseteq \Ann(I_1/J_1) \cap \ldots \cap \Ann(I_r/J_r),$$ откуда в силу предложения~\ref{PropositionGrToHatG}
справедливо неравенство $\PIexp^{G\text{-}\mathrm{gr}}(L) \leqslant (\dim L) - (\dim N)$.
В силу леммы~\ref{LemmaInclExcl} это неравенство справедливо и в случае произвольной
градуирующей группы.
 Следовательно, в случае, когда $\PIexp^{G\text{-}\mathrm{gr}}(L)=\dim L$, справедливо нильпотентный радикал $N$ равен нулю
и в силу предложения~2.1.7 из~\cite{GotoGrosshans}
выполняется условие
$[L, R] \subseteq N = 0$, где $R$~"--- разрешимый радикал алгебры Ли $L$. Отсюда $R = Z(L)\subseteq N=0$ 
и алгебра Ли $L$ полупроста. Теперь достаточно применить теорему~\ref{TheoremCoHLieSemiSimple} и воспользоваться примером~\ref{ExampleGrLieSemiSimplePIexp}.
\end{proof}

Обратимся теперь к алгебрам с $G$-действиями:

   \begin{example}\label{ExampleGLieSemiSimplePIexp}
 Пусть $L=B_1 \oplus B_2 \oplus \ldots \oplus B_q$ (прямая сумма $G$-инвариантных идеалов), 
где $B_i$~"--- конечномерные $G$-простые алгебры Ли над алгебраически замкнутым полем $\mathbbm{k}$
с действием некоторой группы $G$ автоморфизмами и антиавтоморфизмами. Введём обозначение $d := \max_{1 \leqslant k
  \leqslant q} \dim B_k$.  Тогда существуют такие $C_1, C_2 > 0$, $r_1, r_2 \in \mathbb R$,
  что $$C_1 n^{r_1} d^n \leqslant c_n^{G}(L)
  \leqslant C_2 n^{r_2} d^n \text{ для всех } n\in\mathbb N.$$
   \end{example}
   \begin{proof} В силу предложения~\ref{PropositionGtoAut}
можно считать, что группа $G$ действует только автоморфизмами.
Теперь достаточно заметить, что радикалы инвариантны относительно всех автоморфизмов, и применить
теорему~\ref{TheoremMainLieNRSame}.
\end{proof}

Отсюда получаем следующий критерий $G$-простоты:

\begin{theorem}\label{TheoremGCrSimpleLie}
Пусть $L$~"--- конечномерная алгебра Ли над алгебраически замкнутым полем $\mathbbm{k}$ характеристики $0$
с действием некоторой группы $G$ автоморфизмами и антиавтоморфизмами.
Тогда $\PIexp^{G}(L)=\dim L$, если и только если алгебра Ли $L$ является $G$-простой.
\end{theorem}
\begin{proof}
Повторяя рассуждения леммы~\ref{LemmaAssocUpperCochar} и теоремы~\ref{TheoremAssocUpper},
получаем, что если некоторая конечномерная (необязательно ассоциативная) алгебра $A$ с обобщённым $H$-действием обладает $H$-инвариантным нильпотентным идеалом $J$,
то $\overline{\lim}_{n\to\infty} \sqrt[n]{c_n^H(A)} \leqslant (\dim A)-(\dim J)$.
Учитывая, что нильпотентный радикал $N$ алгебры Ли $L$ инвариантен относительно автоморфизмов и антиавтоморфизмов, получаем, что в случае $\PIexp^{G}(L)=\dim L$
справедливо равенство $N=0$ и в силу предложения~2.1.7 из~\cite{GotoGrosshans}
выполняется условие
$[L, R] \subseteq N = 0$, где $R$~"--- разрешимый радикал алгебры Ли $L$. Отсюда $R = Z(L)\subseteq N=0$ 
и алгебра Ли $L$ полупроста. Теперь достаточно применить теорему~\ref{TheoremHLieSemiSimple} и воспользоваться примером~\ref{ExampleGLieSemiSimplePIexp}.
\end{proof}

Обратимся теперь к случаю действий алгебр Ли дифференцированиями:

   \begin{example}\label{ExampleDiffLieSemiSimplePIexp}
 Пусть $L=B_1 \oplus B_2 \oplus \ldots \oplus B_q$ (прямая сумма $\mathfrak g$-инвариантных идеалов), 
где $B_i$~"--- конечномерные простые алгебры Ли над алгебраически замкнутым полем $\mathbbm{k}$
с действием некоторой алгебры Ли $\mathfrak g$ дифференцированиями. Введём обозначение $d := \max_{1 \leqslant k
  \leqslant q} \dim B_k$.  Тогда существуют такие $C_1, C_2 > 0$, $r_1, r_2 \in \mathbb R$,
  что $$C_1 n^{r_1} d^n \leqslant c_n^{U(\mathfrak g)}(L)
  \leqslant C_2 n^{r_2} d^n \text{ для всех } n\in\mathbb N.$$
   \end{example}
\begin{proof} Достаточно применить следствие~\ref{CorollaryLieRadicalsDiffInvariant} и
теорему~\ref{TheoremMainLieNRSame}.
\end{proof}
   
Выведем отсюда критерий простоты в терминах градуированных коразмерностей.
(Напомним, что в силу теоремы~\ref{TheoremDerSimpleLie} понятия простоты и дифференциальной простоты совпадают.)

\begin{theorem}\label{TheoremDiffCrSimpleLie}
Пусть $L$~"--- конечномерная алгебра Ли над алгебраически замкнутым полем $\mathbbm{k}$ характеристики~$0$
с действием некоторой алгебры Ли $\mathfrak g$ дифференцированиями.
Тогда $\PIexp^{U(\mathfrak g)}(L)=\dim L$, если и только если алгебра Ли $L$ является простой.
\end{theorem}
\begin{proof}
В силу следствия~\ref{CorollaryLieRadicalsDiffInvariant}
нильпотентный и разрешимы радикалы являются $\mathfrak g$-подмодулями.
Теперь достаточно повторить доказательство теоремы~\ref{TheoremGCrSimpleLie},
использовав вместо примера~\ref{ExampleGLieSemiSimplePIexp} пример~\ref{ExampleDiffLieSemiSimplePIexp}.
\end{proof}

Сформулируем теперь достаточные условие разрешимости алгебр Ли в терминах их PI-экспонент:

\begin{theorem}\label{TheoremHGenLieSolvable}
Пусть $L$~"--- конечномерная алгебра Ли
над полем $\mathbbm{k}$ характеристики~$0$
с обобщённым $H$-действием
некоторой ассоциативной алгебры $H$ с единицей,
причём $\overline{\lim}_{n\to\infty}\sqrt[n]{c_n^H(L)} < 3$. Тогда алгебра Ли $L$ разрешима.
\end{theorem}
\begin{proof}
 Если алгебра  Ли $L$ неразрешима,
то она останется таковой и при расширении основного поля, откуда можно предполагать основное
поле $\mathbbm{k}$ алгебраически замкнутым. Рассмотрим обычное разложение Леви алгебры Ли $L$
и предположим, что
алгебра  Ли $L$ неразрешима. Тогда $L$ содержит некоторую простую подалгебру $B_0$,
причём $\dim B_0 \geqslant 3$.
В силу примера~\ref{ExampleGLieSemiSimplePIexp} для $G=\langle 1 \rangle$
и предложения~\ref{PropositionOrdinaryAndHopf} существуют такие $C_1 > 0$ и $r_1 \in\mathbb R$,
что
$$C_1 n^{r_1} 3^n \leqslant C_1 n^{r_1} (\dim B_0)^n
\leqslant c_n(B_0)\leqslant c_n(L)\leqslant c^H_n(L),$$
откуда неравенство
$\overline{\lim}_{n\to\infty}\sqrt[n]{c_n^H(L)} < 3$ никак не может быть выполнено.
\end{proof}

\begin{theorem}\label{TheoremGrLieSolvable}
Пусть $L$~"--- конечномерная алгебра Ли
над полем $\mathbbm{k}$ характеристики $0$,
градуированная произвольным множеством $T$,
причём $\overline{\lim}_{n\to\infty}\sqrt[n]{c_n^{T\text{-}\mathrm{gr}}(L)} < 3$. Тогда алгебра Ли $L$ разрешима.
\end{theorem}
\begin{proof}
Достаточно повторить доказательство теоремы~\ref{TheoremHGenLieSolvable},
использовав вместо предложения~\ref{PropositionOrdinaryAndHopf}
предложение~\ref{PropositionOrdinaryAndGradedCodim}.
\end{proof}

Приведём теперь пример неполупростой алгебры Ли, градуированной неабелевой группой:
\begin{example}\label{Example2gl2S3}
Пусть $\mathbbm{k}$~"--- поле характеристики $0$, $G=S_3$, а $L=\mathfrak{gl}_2(\mathbbm{k})\oplus \mathfrak{gl}_2(\mathbbm{k})$. Рассмотрим на~$L$
следующую $G$-градуировку:
 $$L^{(e)} = \left\lbrace\left(\begin{array}{rr}
\alpha & 0 \\
 0 & \beta 
\end{array} \right)\right\rbrace \oplus  \left\lbrace\left(\begin{array}{rr}
\gamma & 0 \\
 0 & \mu 
\end{array} \right)\right\rbrace,$$ $$L^{\bigl((12)\bigr)} = \left\lbrace\left(\begin{array}{rr}
0 & \alpha \\
 \beta & 0  
\end{array} \right)\right\rbrace \oplus  0,\qquad L^{\bigl((23)\bigr)} = 0 \oplus \left\lbrace\left(\begin{array}{rr}
0 & \alpha \\
 \beta & 0  
\end{array} \right)\right\rbrace,$$ остальные компоненты равны $0$.
Тогда существуют такие $C_1, C_2 > 0$, $r_1, r_2 \in \mathbb R$,
что $$C_1 n^{r_1} 3^n \leqslant c_n^{G\text{-}\mathrm{gr}}(L)
  \leqslant C_2 n^{r_2} 3^n \text{ для всех } n \in\mathbb N.$$
\end{example}
\begin{proof} Как обычно, воспользуемся тем фактом, что коразмерности не меняются
при расширении основного поля. При этом отметим, что при расширении поля алгебра Ли $L$ остаётся алгеброй того же типа.
Отсюда можно без ограничения общности считать основное поле $\mathbbm{k}$ алгебраически замкнутым.

Воспользуемся разложением $$L=\mathfrak{sl}_2(\mathbbm{k})\oplus \mathfrak{sl}_2(\mathbbm{k})\oplus Z(L),$$
где центр $Z(L)$ состоит из скалярных матриц из обеих копий алгебры Ли $\mathfrak{gl}_2(\mathbbm{k})$.
Отсюда $$V^{G\text{-}\mathrm{gr}}_n\cap\Id^{G\text{-}\mathrm{gr}}(L)=V^{G\text{-}\mathrm{gr}}_n\cap\Id^{G\text{-}\mathrm{gr}}\bigl(\mathfrak{sl}_2(\mathbbm{k}) \oplus \mathfrak{sl}_2(\mathbbm{k})\bigr)\text{ для всех }n\in \mathbb N.$$
Теперь достаточно заметить, что обе копии простой алгебры Ли $\mathfrak{sl}_2(\mathbbm{k})$
являются градуированными идеалами алгебры Ли $L$,
и воспользоваться примером~\ref{ExampleGrLieSemiSimplePIexp}.
\end{proof}

Приведём теперь два примера с метабелевыми алгебрами Ли:

\begin{example}\label{ExampleLieMetabelian}
Пусть $m\in \mathbb N$, $G\subseteq S_m$~"--- некоторая подгруппа,
а $O_i$~"--- орбиты $G$-действия на множестве
$$\lbrace1, 2, \ldots, m \rbrace = \coprod_{i=1}^s O_i.$$
 Введём обозначение $$d:=\max_{1\leqslant i \leqslant s} |O_i|.$$
Обозначим через $L$ алгебру Ли над произвольным полем $\mathbbm{k}$ характеристики $0$ с базисом $a_1, \ldots, a_m$, $b_1,
\ldots, b_m$, где $\dim L = 2m$, и коммутатором,
заданным формулами $[a_i, a_j]=[b_i,b_j]=0$ и
 $$[a_i, b_j] = \left\lbrace \begin{array}{rrr} b_j & \text{при} &
  i = j,\\
0 & \text{при} & i \ne j.
\end{array} \right.$$
 Предположим, что группа $G$ действует на $L$
следующим образом:
 $\sigma a_i :=a_{\sigma(i)}$ и
 $\sigma b_j :=b_{\sigma(j)}$ для всех $\sigma \in G$.
Тогда существуют такие $C_1, C_2 > 0$ и $r_1, r_2 \in \mathbb R$, что
$$C_1 n^{r_1} d^n \leqslant c_n^{G}(L) \leqslant C_2 n^{r_2} d^n$$
для всех $n\in\mathbb N$.
В частности, если $$G=\langle\tau\rangle \cong 
 \mathbb Z/(m\mathbb Z)%
,$$
где $\tau = (1\,2\, 3\, \ldots\, m)$ (цикл длины $m$),
то
$$C_1 n^{r_1} m^n \leqslant c_n^{G}(L) \leqslant C_2 n^{r_2} m^n.$$
В то же время $c_n(L)=n-1$ для всех $n\in\mathbb N$.
\end{example}
\begin{proof} Если $K \supseteq \mathbbm{k}$~"--- расширение основного поля, то алгебра Ли $K \mathbin{\otimes_\mathbbm{k}} L$
задаётся теми же формулами, что и алгебра Ли $L$.
Поскольку при расширении основного поля коразмерности не меняются,
можно без ограничения общности считать основное поле $\mathbbm{k}$ алгебраически замкнутым.

Пусть $B_i:=\langle b_j \mid j \in O_i \rangle_\mathbbm{k}$, где
$1 \leqslant i \leqslant s$.
Предположим, что $I$~"--- некоторый $G$-инвариантный идеал алгебры Ли $L$.
 Если $b_i \in I$, то
 $b_{\sigma(i)}=\sigma b_i \in I$ для всех $\sigma \in G$.
 Отсюда если $i \in O_j$, то $b_k \in I$ для всех $k \in O_j$.
 Пусть $c:=\sum\limits_{i=1}^m (\alpha_i a_i + \beta_i b_i) \in I$
 для некоторых $\alpha_i, \beta_i \in \mathbbm{k}$.
 Тогда $\beta_i b_i = [a_i, c] \in I$ для всех $1\leqslant i \leqslant m$.
Следовательно, $I=A_0 \oplus B_{i_1} \oplus \ldots \oplus B_{i_k}$
для некоторых $1 \leqslant i_j \leqslant s$ и $A_0 \subseteq \langle
a_1, \ldots, a_m \rangle_\mathbbm{k}$.

Если $I,J \subseteq L$~"--- $G$-инвариантные идеалы, то
$J \subseteq J+[L,L]\cap I \subseteq I$ также является $G$-инвариантным идеалом.
 Предположим, что $I/J$~"--- неприводимый $(G,L)$-модуль. Тогда либо $[L,L]\cap I \subseteq J$
и $\Ann(I/J)=L$, либо $I \subseteq J + [L,L]$,
где $[L,L]=\langle b_1, \ldots, b_m \rangle_\mathbbm{k}$.
Следовательно, в случае $\Ann(I/J)\ne L$ получаем, что $J=A_0 \oplus B_{i_1} \oplus \ldots \oplus B_{i_k}$
и $I=B_\ell \oplus J$ для некоторого $1 \leqslant \ell \leqslant s$.
В этом случае $$\dim(L/\Ann(I/J))=\dim\bigl(\langle a_j \mid j \in O_\ell \rangle_\mathbbm{k}\bigr)=|O_\ell|.$$

Заметим, что если $I_1 = B_{i_1}\oplus J_1$,
а $I_2 = B_{i_2}\oplus J_2$, то $$[[B_{i_1}, L, \ldots, L], [B_{i_2}, L, \ldots, L]]=0.$$
Следовательно, $G$-инвариантные идеалы $I_1, \ldots, I_r$, $J_1, \ldots, J_r$ могут удовлетворять условиям 1--2 из \S\ref{SectionHPIexpLie} только при $r=1$. Отсюда $$d(L, \mathbbm{k}G)=\max_{1 \leqslant i \leqslant s} |O_\ell|$$ и требуемые оценки следуют из теоремы~\ref{TheoremMainLieGFin}.

Рассмотрим обычные полиномиальные тождества.
Используя тождество Якоби, всякий одночлен из пространства $V_n := V_n^\mathbbm{k}$
можно представить в виде линейной комбинации левонормированных
коммутаторов $[x_1, x_j, x_{i_3},
\ldots, x_{i_n}]$. Поскольку в алгебре Ли $L$ выполняется полиномиальное тождество $$[[x,y],[z,t]]\equiv 0,$$ можно считать, что $i_3 < i_4 < \ldots < i_n$.
Теперь осталось заметить, что одночлены $$f_j=[x_1, x_j, x_{i_3},
\ldots, x_{i_n}]\text{, где }2\leqslant j \leqslant n,$$ линейно независимы по модулю $\Id(L)$.
Действительно, если $\sum_{k=2}^n \alpha_k f_k \equiv 0$ для некоторых $\alpha_k \in \mathbbm{k}$,
то при подстановке $x_j=b_1$ и $x_i=a_1$ при $i\ne j$
единственным слагаемым, которое не обратится в нуль будет слагаемое $f_j$.
Отсюда $\alpha_j=0$ и $c_n(L)=n-1$.
\end{proof}

\begin{example}\label{ExampleMetabelianGr}
Пусть $m\in \mathbb N$, $L=\bigoplus\limits_{\bar k \in \mathbb Z/m\mathbb Z} L^{(\bar k)}$~"---
$\mathbb Z/m\mathbb Z$-градуированная алгебра Ли, где
 $L^{(\bar k)}= \langle c_{\bar k}, d_{\bar k} \rangle_\mathbbm{k}$,
  $\dim L^{(\bar k)} = 2$, коммутатор задан формулами
  $[c_{\bar\imath}, c_{\bar\jmath}]=
  [d_{\bar\imath},d_{\bar\jmath}]=0$ и
  $[c_{\bar\imath}, d_{\bar\jmath}] =
   d_{\bar\imath+\bar\jmath}$, а $\mathbbm{k}$~"--- произвольное поле
характеристики $0$.
Тогда существуют такие $C_1, C_2 > 0$ и $r_1, r_2 \in \mathbb R$, что
$$C_1 n^{r_1} m^n \leqslant c_n^{\mathbb Z/m\mathbb Z\text{-}\mathrm{gr}}(L) \leqslant C_2 n^{r_2} m^n$$
для всех $n\in\mathbb N$.
\end{example}
\begin{proof}
Как и прежде, можно без ограничения общности считать поле $\mathbbm{k}$ алгебраически замкнутым.
Пусть $\zeta \in \mathbbm{k}$~"--- примитивный корень из единицы степени $m$.
Тогда для $G=\mathbb Z/m\mathbb Z$ справедливо равенство $\widehat G
=\lbrace \psi_0, \ldots, \psi_{m-1}\rbrace$, где $\psi_\ell(\bar\jmath):=\zeta^{\ell j}$.
Оказывается, что алгебру Ли $L$ из данного примера можно отождествить с
алгеброй Ли из примера~\ref{ExampleLieMetabelian}
при помощи формул $c_{\bar\jmath} = \sum_{k=1}^m \zeta^{-jk} a_k$ и
$d_{\bar\jmath} = \sum_{k=1}^m \zeta^{-jk} b_k$. 
При этом $\hat G$-действию, отвечающему $\mathbb Z/m\mathbb Z$-градуировке на $L$,
соответствует $\langle\tau\rangle$-действие из примера~\ref{ExampleLieMetabelian}:
 $\tau^\ell c_{\bar\jmath} = \zeta^{\ell j} c_{\bar\jmath}
=\psi_\ell(\bar\jmath) c_{\bar\jmath}$
и $\tau^\ell d_{\bar\jmath} = \zeta^{\ell j} d_{\bar\jmath}=\psi_\ell(\bar\jmath)
 d_{\bar\jmath}$.
В силу предложения~\ref{PropositionGrToHatG} справедливо равенство $c_n^{\mathbb Z/m\mathbb Z\text{-}\mathrm{gr}}(L)=c_n^{\langle\tau\rangle}(L)$, и оценки сверху и снизу получаются из примера~\ref{ExampleLieMetabelian}.
\end{proof}

Приведём теперь примеры существования PI-экспонент в алгебрах Ли, для которых не существует инвариантных разложений Леви:

\begin{example}
Пусть $$L = \left\lbrace\left(\begin{array}{cc} C & D \\
0 & 0
  \end{array}\right) \mathrel{\biggl|} C \in \mathfrak{sl}_m(\mathbbm{k}), D\in M_m(\mathbbm{k})\right\rbrace
  \subseteq \mathfrak{sl}_{2m}(\mathbbm{k}),\ m \geqslant 2.$$
  Тогда идеал $$R=\left\lbrace\left(\begin{array}{cc} 0 & D \\
0 & 0
  \end{array}\right) \mathrel{\biggl|} D\in M_m(\mathbbm{k})\right\rbrace
  $$
  является разрешимым (и нильпотентным) идеалом алгебры Ли $L$.
    Определим $\varphi \in \Aut(L)$ по формуле
  $$\varphi\left(\begin{array}{cc} C & D \\
0 & 0
  \end{array}\right)=\left(\begin{array}{cc} C & C+D \\
0 & 0
  \end{array}\right).$$
  Тогда группа $G=\langle \varphi \rangle
  \cong \mathbb Z$ действует на $L$ автоморфизмами, т.е. $L$ является $\mathbbm{k}G$-модульной
  алгеброй.
    Как было показано в примере~\ref{ExampleGnoninvLevi},
  для алгебры $A$ не существует $G$-инвариантного разложения Леви.
  В то же время существуют такие константны $C_1, C_2 > 0$, $r_1, r_2 \in \mathbb R$, что $$C_1 n^{r_1} (m^2-1)^n \leqslant c^G_n(L) \leqslant C_2 n^{r_2} (m^2-1)^n\text{ для всех }n \in \mathbb N.$$
\end{example}
\begin{proof} Как уже было отмечено, коразмерности не меняются при расширении основного поля.
Более того, при расширении основного поля алгебра Ли $L$ остаётся алгеброй того же типа.
Следовательно, без ограничения общности можно считать, что основное поле алгебраически замкнуто.

Заметим, что
 $$N=\left\lbrace\left(\begin{array}{cc} 0 & D \\
0 & 0
  \end{array}\right) \mathrel{\biggl|} D\in M_m(\mathbbm{k})\right\rbrace
  $$
  является разрешимым и нильпотентным радикалом алгебры Ли $L$, причём $L/N \cong \mathfrak{sl}_m(\mathbbm{k})$~"--- простая алгебра Ли. Отсюда $\PIexp^G(L)=\dim\mathfrak{sl}_m(\mathbbm{k})= m^2-1$
  в силу теоремы~\ref{TheoremMainLieNRSame}.
\end{proof}

Проводя аналогичные рассуждения, получаем следующий пример:

\begin{example} Пусть $L$~"---  та же алгебра Ли, что и в предыдущем примере.
Рассмотрим присоединённое представление алгебры Ли $L$ на себе самой дифференцированиями.
Тогда $L$ оказывается $U(L)$-модульной алгеброй Ли.
Как было показано в примере~\ref{ExampleDiffnoninvLevi}, в $A$ не существует $U(L)$-инвариантного разложения Леви.  В то же время существуют такие константны $C_1, C_2 > 0$, $r_1, r_2 \in \mathbb R$, что $$C_1 n^{r_1} (m^{2}-1)^n \leqslant c^{U(L)}_n(L) \leqslant C_2 n^{r_2} (m^2-1)^{n}\text{ для всех }n \in \mathbb N.$$
\end{example}

     \section{Асимптотика $H_{m^2}(\zeta)$-коразмерностей $H_{m^2}(\zeta)$-простых алгебр Ли}\label{SectionHPI-expTaftSLie}

До этого момента в этой главе изучалось асимптотическое поведение $H$-коразмерностей лишь $H$-модульных алгебр Ли с $H$-инвариантными радикалами. В данном параграфе доказывается, что если $L$~"--- конечномерная $H_{m^2}(\zeta)$-простая
$H_{m^2}(\zeta)$-модульная алгебра Ли над алгебраически замкнутым полем $\mathbbm{k}$ характеристики $0$,
то  $\PIexp^{H_{m^2}(\zeta)}(L)=\dim L$. (Напомним, что через $H_{m^2}(\zeta)$ обозначается алгебра Тафта.) В частности,  $H_{m^2}(\zeta)$-PI-экспонента такой алгебры Ли $L$
является целым числом, и для полиномиальных $H_{m^2}(\zeta)$-тождеств
алгебры Ли $L$ справедлив аналог гипотезы Амицура. При этом разрешимый радикал алгебры Ли $L$ может и не являться $H_{m^2}(\zeta)$-подмодулем. (См. алгебры Ли $L(B,0)$ из~\S\ref{SectionClassTaftSLieSimple}.)
 
   \begin{theorem}\label{TheoremTaftSimpleLieHPIexpExists}
Пусть $L$~"--- конечномерная $H_{m^2}(\zeta)$-простая
$H_{m^2}(\zeta)$-модульная алгебра Ли над алгебраически замкнутым полем $\mathbbm{k}$ характеристики $0$.   
   Тогда существуют такие $C>0$ и $r\in \mathbb R$, что
   $$C n^r (\dim L)^n \leqslant c_n^{H_{m^2}(\zeta)}(L) \leqslant (\dim L)^{n+1}\text{ для всех }n\in \mathbb N.$$
   В частности, $\PIexp^{H_{m^2}(\zeta)}(L) = \dim L$, и для $L$ справедлив аналог гипотезы Амицура.
  \end{theorem} 

Сперва докажем существование полилинейного $H$-многочлена, не являющегося
полиномиальным $H$-тождеством, кососимметричного по достаточному числу наборов переменных:

\begin{lemma}\label{LemmaTaftSimpleLieAlt} Пусть $L$~"--- конечномерная
неполупростая $H_{m^2}(\zeta)$-простая Ли над алгебраически замкнутым полем $\mathbbm{k}$ характеристики $0$.
Введём обозначение $\ell := \dim L^{(0)}$.
Тогда существует такое число $r \in \mathbb N$, что для всех $n\geqslant \ell m r + 1$
существуют попарно непересекающиеся подмножества $X_1$, \ldots, $X_{kr} \subseteq \lbrace x_1, \ldots, x_n
\rbrace$, где $k := \left[\frac{n-1}{\ell m r}\right]$,
$|X_1| = \ldots = |X_{kr}|=\ell m$, и $H$-многочлен $f \in V^{H_{m^2}(\zeta)}_n \backslash
\Id^{H_{m^2}(\zeta)}(L)$, кососимметричный по переменным каждого из множеств $X_j$.
\end{lemma}
  \begin{proof}
  Поскольку алгебра Ли $L$ не является полупростой,
  из теоремы~\ref{TheoremTaftSimpleNonSemiSimpleLieClassify} следует,
  что
  $L \cong L(B, 0)$ для некоторой простой алгебры Ли $B$, где $\dim B = \ell$.
  Алгебра Ли $L(B,0)$ наделена $\mathbb Z/m\mathbb Z$-градуировкой $L=\bigoplus_{k=0}^{m-1} L^{(k)}$
  (см. теорему~\ref{TheoremTaftSimpleLiePresent}),
  где компонента $L^{(0)}$ может быть отождествлена с $B$.
  В силу теоремы Ю.\,П.~Размыслова~\cite[теорема~12.1]{RazmyslovBook}
  существует такое число $r\in\mathbb N$, что для любого $k\in\mathbb N$
  существует такой полилинейный ассоциативный многочлен $$f_0=f_0(x_{11}, \ldots, x_{1\ell};\ldots; x_{kr,1}, \ldots, x_{kr, \ell}),$$
  кососимметричный по переменным каждого из множеств $\lbrace x_{i1}, \ldots, x_{i\ell} \rbrace$, $1\leqslant i \leqslant kr$, что
  $f_0(\ad a_1, \ldots, \ad a_\ell;\ldots; \ad a_1, \ldots, \ad a_\ell)$
  является ненулевым скалярным оператором на $B$ для любого базиса $a_1, \ldots, a_\ell$ алгебры Ли $B$.
  
  Пусть $n\in\mathbb N$. Положим $k := \left[ \frac{n-1}{\ell m r}\right]$. Используя теорему Ю.\,П.~Размыслова, выберем многочлен $f_0$,  кососимметричный по $\ell$ переменным каждого из $kr$ множеств, а также многочлен $\tilde f_0$, кососимметричный по $mr$ множествам из $\ell$ переменных.
  Рассмотрим теперь полилинейный лиевский ${H_{m^2}(\zeta)}$-многочлен $$f_1:=\tilde f_0(\ad y_{11}, \ldots, \ad y_{1\ell};\ldots; \ad y_{mr,1}, \ldots, \ad y_{mr, \ell})f_2,$$ где $$f_2=\left(\prod_{i=1}^m f_0(\ad (x_{11i}^{v^{i-1}}), \ldots, \ad (x_{1\ell i}^{v^{i-1}});\ldots; \ad (x_{kr,1,i}^{v^{i-1}}), \ldots, \ad (x_{kr, \ell,i}^{v^{i-1}}))\right)z.$$
  Пусть $b_1, \ldots, b_\ell$~"--- базис компоненты $L^{(m-1)}$.
  Тогда $v^{m-1} b_1, \ldots, v^{m-1} b_\ell$ является базисом алгебры Ли $L^{(0)}=B$.
  Следовательно, ${H_{m^2}(\zeta)}$-многочлен $f_1$
не обращается в нуль при подстановке $x_{jti}=v^{m-i} b_t$, $y_{jt}=v^{m-1} b_t$ и $z=\bar z$
для любого ненулевого элемента $\bar z \in B$. Фиксируем некоторый такой элемент $\bar z$
и обозначим это подстановку через $\Xi$.
  Пусть $b$~"--- значение ${H_{m^2}(\zeta)}$-многочлена $f_1$ при подстановке $\Xi$.
  Рассмотрим теперь ${H_{m^2}(\zeta)}$-многочлен $f_3 := \Alt_1 \ldots \Alt_{kr} f_1$
  где $\Alt_j$~"--- оператор альтернирования по переменным из множества $X_j = \lbrace x_{jti} \mid 1\leqslant t \leqslant \ell,\ 1\leqslant i \leqslant m \rbrace$.
  Заметим, что из  $v^m=0$ следует что все слагаемые, где $x_{jti}$ заменяется при альтернировании
  на  $x_{jt'i'}$ с $i' < i$, обращаются в нуль.
Следовательно, без ограничения общности можно считать, что при альтернировании
 переменные каждого из множеств 
  $\lbrace x_{jti} \mid 1\leqslant t \leqslant \ell \rbrace$ при фиксированных $j$ и $i$
  перемешиваются между собой.
  Поскольку ${H_{m^2}(\zeta)}$-многочлен $f_1$
кососимметричен по переменным каждого из этих множеств, значение ${H_{m^2}(\zeta)}$-многочлена $f_3$
при подстановке $\Xi$ равно $(\ell!)^{kmr} b \ne 0$.
  
  Заметим, что $ k\ell m r+1 \leqslant n < \deg f_3 = (k+1)\ell mr+1$. Расписывая
   $\tilde f_0$, представим $f_3$  в виде линейной комбинации одночленов 
  \begin{equation*}\begin{split}f_4 := \left[w_1, \left[w_2, \ldots, \left[w_{\ell m r}, \Alt_1 \ldots \Alt_{kr}
  \left(\prod_{i=1}^m f_0(\ad (x_{11i}^{v^{i-1}}), \ldots, \ad (x_{1\ell i}^{v^{i-1}}); \right.\right.\right.\right. \\ \left.\left.\left. \ldots; \ad (x_{kr,1,i}^{v^{i-1}}), \ldots, \ad (x_{kr, \ell,i}^{v^{i-1}}))\right)z\right]\ldots\right],\end{split}\end{equation*}
  где переменные $w_i$~"--- это переменные $y_{jt}$, взятые в некотором порядке в зависимости от слагаемого.
  Поскольку $f_3 \notin \Id^{H_{m^2}(\zeta)}(L)$,
  по крайней мере одно из слагаемых $f_4$ не обращается в нуль при подстановке $\Xi$.
  Тогда
  \begin{equation*}\begin{split}f := \left[w_{(k+1)\ell m r-n+2}, \left[w_{(k+1)\ell m r-n+3}, \ldots, \left[w_{\ell m r}, \Alt_1 \ldots \Alt_{kr}
  \left(\prod_{i=1}^m f_0(\ad (x_{11i}^{v^{i-1}}), \ldots, \ad (x_{1\ell i}^{v^{i-1}});\right.\right.\right.\right. \\ \left.\left.\left. \ldots; \ad (x_{kr,1,i}^{v^{i-1}}), \ldots, \ad (x_{kr, \ell,i}^{v^{i-1}}))\right)z\right]\ldots\right] \notin \Id^{H_{m^2}(\zeta)}(L).
  \end{split}\end{equation*}
  Заметим, что $\deg f = n$ и если переименовать переменные ${H_{m^2}(\zeta)}$-многочлена $f$ в $x_1, x_2, \ldots, x_n$, то $f$ удовлетворяет всем условиям леммы.
  \end{proof}
  \begin{proof}[Доказательство теоремы~\ref{TheoremTaftSimpleLieHPIexpExists}]
Оценка сверху $c_n^{H_{m^2}(\zeta)}(L) \leqslant (\dim L)^{n+1}$
следует из предложения~\ref{PropositionCodimDim}.
  
  Если алгебра Ли $L$ полупроста, то оценка снизу следует из примера~\ref{ExampleHLieSemiSimplePIexp}.
    Рассмотрим случай, когда алгебра Ли $L$ неполупроста.  
  
    Пусть $r$~"--- число из леммы~\ref{LemmaTaftSimpleLieAlt}. Введём обозначения $\ell :=\frac{\dim L}m$ и $k := \left[\frac{n-1}{\ell m r}\right]$.
  Докажем, что для любого $n\in\mathbb N$ существует такое $\lambda \vdash n$, что $m(L, {H_{m^2}(\zeta)}, \lambda)\ne 0$ и $\lambda_i \geqslant kr$ для всех $1\leqslant i \leqslant \ell m$.
    Рассмотрим $H$-многочлен $f$ из леммы~\ref{LemmaTaftSimpleLieAlt}.
Достаточно доказать, что $e^*_{T_\lambda} f \notin \Id^H(L)$
для некоторой таблицы Юнга $T_\lambda$ требуемой формы $\lambda$.
Известно (см., например, теорему 3.2.7 из~\cite{Bahturin}), что $$\mathbbm{k}S_n = \bigoplus_{\lambda,T_\lambda} \mathbbm{k}S_n e^{*}_{T_\lambda},$$ где суммирование ведётся по множеству стандартных таблиц Юнга $T_\lambda$
всевозможных форм $\lambda \vdash n$. Отсюда $$\mathbbm{k}S_n f = \sum_{\lambda,T_\lambda} \mathbbm{k}S_n e^{*}_{T_\lambda}f
\not\subseteq \Id^{H_{m^2}(\zeta)}(L)$$ и $e^{*}_{T_\lambda} f \notin \Id^{H_{m^2}(\zeta)}(L)$ для некоторого $\lambda \vdash n$.
Докажем, что разбиение $\lambda$ имеет требуемый вид.
Достаточно доказать, что
$\lambda_{\ell m} \geqslant kr$, так как
$\lambda_i \geqslant \lambda_{\ell m}$ для всех $1 \leqslant i \leqslant \ell m$.
В любой строчке таблицы $T_\lambda$ содержится не более одного 
номера переменной из одного и того же множества $X_i$,
поскольку $e^{*}_{T_\lambda} = b_{T_\lambda} a_{T_\lambda}$,
а $a_{T_\lambda}$ симметризует по переменным, отвечающим каждой строчке таблицы $T_\lambda$.
Отсюда $\sum_{i=1}^{\ell m-1} \lambda_i \leqslant kr(\ell m-1) + (n-k\ell m r) = n-kr$.
 Из леммы~\ref{LemmaHStripeDimATheorem} следует, что если $\lambda \vdash n$ и $\lambda_{\ell m+1} > 0$, то $m(L, {H_{m^2}(\zeta)}, \lambda)=0$.
Следовательно,
$\lambda_{\ell m} \geqslant kr$.

Диаграмма Юнга~$D_\lambda$ содержит прямоугольную
поддиаграмму~$D_\mu$, где $\mu=(\underbrace{kr, \ldots, kr}_{\ell m})$.
Из правила ветвления для группы подстановок $S_n$
следует, что если ограничить $S_n$-действие на $M(\lambda)$ до $S_{n-1}$,
то
$M(\lambda)$ оказывается прямой суммой всех неизоморфных
$\mathbbm{k}S_{n-1}$-модулей $M(\nu)$,  где $\nu \vdash (n-1)$ и всякая диаграмма Юнга $D_\nu$
получена из $D_\lambda$ удалением одной клетки. В частности,
$\dim M(\nu) \leqslant \dim M(\lambda)$.
Применяя правило ветвления $(n-k\ell m r)$ раз, получаем неравенство $\dim M(\mu) \leqslant \dim M(\lambda)$.
В силу формулы крюков $$\dim M(\mu) = \frac{(k\ell m r)!}{\prod_{i,j} h_{ij}},$$
где $h_{ij}$~"--- длина крюка диаграммы $D_\mu$ с вершиной в $(i, j)$.
Согласно формуле Стирлинга
\begin{equation*}\begin{split}c_n^{H_{m^2}(\zeta)}(L)\geqslant \dim M(\lambda) \geqslant \dim M(\mu) \geqslant \frac{(k\ell m r)!}{((kr+\ell m)!)^{\ell m}}
\sim \\ \sim \frac{
\sqrt{2\pi (k\ell m r)} \left(\frac{k\ell m r}{e}\right)^{k\ell m r}
}
{
\left(\sqrt{2\pi (k r +\ell m)}
\left(\frac{k r +\ell m}{e}\right)^{k r +\ell m}\right)^{\ell m}
} \sim C_1 k^{r_1} (\ell m)^{k\ell m r}\end{split}\end{equation*}
при $k \to \infty$
для некоторых констант $C_1 > 0$, $r_1 \in \mathbb Q$. (Как и прежде, мы пишем $f\sim g$, если $\lim \frac{f}{g} = 1$.)
В силу того, что $k = \left[\frac{n-1}{\ell m r}\right]$,
оценка снизу доказана.
\end{proof}

\newpage

\chapter{Рост градуированных коразмерностей в ассоциативных алгебрах, градуированных полугруппами}\label{ChapterSGGrAssocCodim}

В данной главе вычисляются градуированные PI-экспоненты для
таких ассоциативных градуированно простых алгебр $A$, градуированных лентами правых нулей,
что $A/J(A)\cong M_2(\mathbbm{k})$ (теоремы~\ref{TheoremGrPIexpTriangleFracM2} и~\ref{TheoremGrPIexpNonTriangleFracM2}), а в общем случае получается оценка сверху на градуированные коразмерности (теорема~\ref{TheoremUpperFrac}). При этом оказывается, что далеко не всегда градуированная PI-экспонента
является целым числом. Для каждой из четырёх полугрупп  из двух элементов, не являющихся группой (см. предложение~\ref{PropositionSemigroupTwoElements}),
строится пример конечномерной градуированной ассоциативной алгебры с дробной градуированной
PI-экспонентой (теоремы~\ref{TheoremQ3GradFractPI}, \ref{TheoremQ1GradFractPI}, \ref{TheoremQ2GradFractPI} 
и замечание~\ref{RemarkQ3OpGradFractPI}). Для конечномерных ассоциативных алгебр,
градуированных полугруппами с сокращениями (теорема~\ref{TheoremTCancelAmitsur}), а также для
конечномерных ассоциативных алгебр с $1$, градуированных лентами левых и правых нулей (теорема~\ref{TheoremTIdemAmitsur}), доказывается существование целой градуированной PI-экспоненты,
причём оказывается, что в случае конечномерных ассоциативных алгебр с $1$, градуированных лентами левых и правых нулей,
градуированная PI-экспонента совпадает с обычной PI-экспонентой.

Результаты главы были опубликованы
в работах~\cite{ASGordienko13, ASGordienko14JanssensJespers}.

\section{Оценка сверху для $\mathbbm{k}^T$-коразмерностей $T$-градуированно простых алгебр}
\label{SectionUpperSGGr}

Пусть $A=\bigoplus_{t\in T} A^{(t)}$~"--- такая конечномерная ассоциативная $T$-градуированная алгебра над полем $\mathbbm{k}$
характеристики $0$ для некоторого множества $T$, что $A/J(A) \cong M_k(\mathbbm{k})$
для некоторого $k\in\mathbb N$ и для
всех $t\in T$ справедливо равенство $A^{(t)}\cap J(A) = 0$.
(Последнее условие справедливо, например, в случае, когда алгебра $A$ является $T$-градуированно простой,
поскольку тогда ни один идеал не содержит ненулевых однородных элементов, так как любой однородный элемент порождает градуированный идеал.)

Напомним, что $T$-градуированные коразмерности алгебры $A$  в силу предложения~\ref{PropositionCnGrCnGenH} совпадают с её $\mathbbm{k}^T$-коразмерностями, где 
алгебра $\mathbbm{k}^T$ действует на $A$ в соответствии с примером~\ref{ExampleGr}.

В данном параграфе доказывается оценка сверху для $T$-градуированных
коразмерностей алгебры $A$.

Для всякого $t\in T$ выберем в пространстве $A^{(t)}$ некоторый базис $\mathcal B^{(t)}$.
Тогда $\mathcal B = \bigcup_{t\in T} \mathcal B^{(t)}$~"--- базис в алгебре $A$.
Фиксируем также изоморфизм $\psi \colon A/J(A) \mathrel{\widetilde\rightarrow} M_k(\mathbbm{k})$,
а через $\pi \colon A \twoheadrightarrow A/J(A)$ обозначим естественный сюръективный гомоморфизм.
Обозначим через $\theta \colon \mathcal B \rightarrow \mathbb Z$ функцию, заданную равенством $$\theta(a)= \min\left\lbrace i-j 
\mid \alpha_{ij}\ne 0,\ 1\leqslant i,j \leqslant k\right\rbrace,$$ где $\alpha_{ij} \in \mathbbm{k}$
определяются при помощи равенства $\psi\pi(a)=\sum\limits_{1\leqslant i,j \leqslant k} \alpha_{ij} e_{ij}$.

Следующее утверждение играет центральную роль  в данном параграфе:

\begin{lemma}\label{LemmaThetaCondition}
Если для некоторого $\mathbbm{k}^T$-многочлена $f\in P_n^{\mathbbm{k}^T}$ и некоторых $a_i \in \mathcal B$,
где $1\leqslant i \leqslant n$, а $n\in\mathbb N$, справедливо неравенство
 $$f(a_1, \ldots, a_n) \ne 0,$$
 то $$\sum_{i=1}^n \theta(a_i) \leqslant k-1.$$
\end{lemma}
\begin{proof} Без ограничения общности можно считать, что $f$
является линейной комбинацией одночленов $x^{q_{t_1}}_{\sigma(1)} x^{q_{t_2}}_{\sigma(2)}
\ldots x^{q_{t_n}}_{\sigma(n)}$, где, как и прежде, $q_t\in \mathbbm{k}^T$
заданы формулами
$q_t(s):=\left\lbrace\begin{smallmatrix} 1 & \text{при} & s=t,\\ 0 & \text{при} & s\ne t,\end{smallmatrix} \right.$ $t_i\in \supp \Gamma$, $\sigma \in S_n$. 
Поскольку все элементы $a_i$ однородны, значение одночлена $x^{q_{t_1}}_{\sigma(1)} x^{q_{t_2}}_{\sigma(2)}
\ldots x^{q_{t_n}}_{\sigma(n)}$ равно $a_{\sigma(1)}\ldots a_{\sigma(n)}$,
если $a_{\sigma(i)} \in A^{(t_i)}$ для всех  $1\leqslant i \leqslant n$, и $0$ в противном случае.
Однако $a_{\sigma(1)}\ldots a_{\sigma(n)}$ снова является однородным элементом. Напомним, что $J(A)\cap A^{(t)}=0$ для всех $t\in T$.
 Следовательно, $a_{\sigma(1)}\ldots a_{\sigma(n)} \ne 0$, если и только если $$\psi\pi (a_{\sigma(1)}\ldots a_{\sigma(n)})  = \psi\pi (a_{\sigma(1)})\psi\pi (a_{\sigma(2)})\ldots \psi\pi (a_{\sigma(n)}) \ne 0.$$
 Заметим теперь, что $e_{i_1 j_1} e_{i_2 j_2} \ldots e_{i_n j_n} \ne 0$ для некоторых $1\leqslant i_\ell, j_\ell
 \leqslant k$, только если $j_1 = i_2$, $j_2 = i_3$, \ldots, $j_{n-1}=i_n$ и, в частности,
 $1-k \leqslant \sum_{\ell = 1}^n (i_\ell - j_\ell)=i_1 - j_n \leqslant k-1$.
 Следовательно, $a_{\sigma(1)}\ldots a_{\sigma(n)} \ne 0$, только если
 $\sum_{\ell = 1}^n \theta(a_i) \leqslant k-1$.
 \end{proof}

Пусть $r:= \dim A$. Определим $$\beta_\ell := \min \left\lbrace \sum\limits_{i=1}^\ell \theta(a_i)
\mathrel{\biggl|} a_i \in \mathcal B,\ a_i\ne a_j \text{ при } i\ne j \right\rbrace,$$
$$\gamma_\ell := \beta_\ell - \beta_{\ell-1},\quad 1\leqslant \ell \leqslant r,\quad \beta_0 := 0.$$
Без ограничения общности можно считать, что $\mathcal B = (a_1, \ldots, a_r)$,
где $$\theta(a_1) \leqslant \theta(a_2) \leqslant \ldots \leqslant \theta(a_r).$$
Тогда $\beta_\ell = \sum\limits_{i=1}^\ell \theta(a_i)$,
где $\gamma_\ell = \theta(a_\ell)$.
В частности, 
\begin{equation*}1-k=\gamma_1 \leqslant \gamma_2 \leqslant \ldots \leqslant \gamma_r.
\end{equation*}

Равенство $\gamma_1 = 1-k$ следует из того, что матричная единица $e_{1k}$,
которая  должна входить с ненулевым коэффициентом в разложение элемента $\varphi\pi(a)$
для некоторого $a\in\mathcal B$,
имеет наименьшее значение $(i-j)$ среди всех матричных единиц $e_{ij}$.

Докажем теперь основное неравенство для $\mathbbm{k}^T$-кохарактеров алгебры $A$:

\begin{lemma}\label{LemmaInequalityLambdaUpperFrac}
Пусть $f\in P_n^{\mathbbm{k}^T}$, а $\lambda \vdash n$ для некоторого $n\in\mathbb N$, причём $\sum_{i=1}^r \gamma_i \lambda_i \geqslant k$ или $\lambda_{r+1} > 0$. Тогда $e^*_{T_\lambda}f \in \Id^{\mathbbm{k}^T}(A)$ для любой таблицы Юнга $T_\lambda$ формы $\lambda$.
\end{lemma}
\begin{proof}
Заметим, что $\mathbbm{k}^T$-многочлен $e^*_{T_\lambda}f = b_{T_\lambda}a_{T_\lambda}f$
кососимметричен по переменным, отвечающим любому столбцу таблицы $T_\lambda$.
Для того, чтобы определить, является ли $\mathbbm{k}^T$-многочлен $e^*_{T_\lambda}f$
полиномиальным $\mathbbm{k}^T$-тождеством для алгебры $A$, достаточно
подставлять в $e^*_{T_\lambda}f$ лишь элементы из базиса
 $\mathcal B$. Если вместо двух переменных, которые отвечают одному и тому же столбцу
 таблицы $T_\lambda$, подставить одинаковые значения, 
то в силу кососимметричности $\mathbbm{k}^T$-многочлена $e^*_{T_\lambda}f$ по таким переменным
результат подстановки будет равен нулю.
Отсюда при $\lambda_{r+1} > 0$ всегда выполняется условие $e^*_{T_\lambda}f \in \Id^{\mathbbm{k}^T}(A)$,
так как в этом случае высота первого столбца больше, чем $\dim A = r$,
и в его переменные подставляется как минимум два одинаковых элемента.
 
Предположим теперь, что $\sum_{i=1}^r \gamma_i \lambda_i \geqslant k$.
Это неравенство можно переписать в виде
\begin{equation}\label{EquationBetaLambdaUpperFrac}\sum_{i=1}^r (\beta_i-\beta_{i-1}) \lambda_i=\sum_{i=1}^r \beta_i (\lambda_i-\lambda_{i+1}) \geqslant k.\end{equation}
(В силу выше сказанного можно считать, что $\lambda_{r+1}=0$.)
 Заметим, что число $(\lambda_i-\lambda_{i+1})$ равно числу столбцов
 высоты $i$ в таблице $T_\lambda$.
 Предположим, что вместо $x_1, \ldots, x_n$ подставляются элементы $b_1, \ldots, b_n \in\mathcal B$.
В силу замечания, сделанного выше, мы можем считать, что вместо переменных, индексы которых
стоят в одном и том же столбце, подставляются разные элементы базиса.
Из определения чисел $\beta_i$ следует, что $\sum_{i=1}^n \theta(b_i) \geqslant \sum_{i=1}^r \beta_i (\lambda_i-\lambda_{i+1})$.
С учётом~(\ref{EquationBetaLambdaUpperFrac})
получаем, что $\sum_{i=1}^n \theta(b_i) \geqslant k$. Теперь из леммы~\ref{LemmaThetaCondition}
следует, что $(e^*_{T_\lambda}f)(b_1, \ldots, b_n)=0$ и $e^*_{T_\lambda}f \in \Id^{\mathbbm{k}^T}(A)$.
\end{proof}

Из предложения~\ref{PropositionCnGrCnGenH} и лемм~\ref{LemmaRegionUpperFd}
и~\ref{LemmaInequalityLambdaUpperFrac} получаем следующую оценку сверху на рост $T$-градуированных
коразмерностей алгебры $A$:

\begin{lemma}\label{LemmaUpperBoundFUpperFrac}
Введём обозначение $$
\Omega:=\lbrace(\alpha_1, \alpha_2, \ldots, \alpha_r)\in \mathbb R^r \mid \alpha_1 \geqslant \ldots\geqslant \alpha_r
\geqslant 0,\ \alpha_1+\alpha_2+\ldots+\alpha_r=1,\ \gamma_1\alpha_1 + \ldots + \gamma_r\alpha_r \leqslant 0 \rbrace.$$
Тогда $$\mathop{\overline\lim}\limits_{n\rightarrow\infty}
 \sqrt[n]{c_n^{T\text{-}\mathrm{gr}}(A)}\leqslant \max_{x\in\Omega} \Phi(x),$$
 где, как и прежде, $\Phi(x_1, \ldots, x_q):=\frac{1}{x_1^{x_1} \ldots x_s^{x_q}}$.
\end{lemma}

Оставшаяся часть параграфа посвящена вычислению максимума функции $\Phi$.
Начнём с самого простого многогранника ограничений:

\begin{lemma}\label{LemmaMaximumFPlainUpperFrac}
Пусть $r\in\mathbb N$, а $$
\Omega_0:=\left\lbrace(\alpha_1, \alpha_2, \ldots, \alpha_r)\in \mathbb R^r \mathrel{\bigl|} \alpha_1, \ldots, \alpha_r
\geqslant 0,\ \alpha_1+\alpha_2+\ldots+\alpha_r=1 \right\rbrace.$$
Тогда $\max_{x\in\Omega_0} \Phi(x) = r$,
а $$\mathop{\mathrm{argmax}}_{x\in\Omega_0} \Phi(x) = \left(\frac{1}{r}, \frac{1}{r}, \ldots,  \frac{1}{r}\right).$$
\end{lemma}
\begin{proof}
Докажем лемму индукцией по $r$. Случай $r=1$ тривиален. Предположим, что $r\geqslant 2$.
Прежде всего выразим $\alpha_1 = 1-\sum_{i=2}^r \alpha_i$ через $\alpha_2, \alpha_3, \ldots, \alpha_r$
и изучим поведение функции
$\Phi_1(\alpha_2, \ldots, \alpha_r) = \frac{1}{\left(1-\sum_{i=2}^r \alpha_i
\right)^{\left(1-\sum_{i=2}^r \alpha_i\right)} \alpha_2^{\alpha_2} \ldots 
\alpha_r^{\alpha_r}}$ на
$$
\tilde \Omega_0:=\left\lbrace(\alpha_2, \ldots, \alpha_r)\in \mathbb R^{r-1} \mathrel{\biggl|} \alpha_2, \ldots, \alpha_r
\geqslant 0,\ 1-\sum_{i=2}^r \alpha_i\geqslant 0 \right\rbrace.$$
Заметим, что функция $\Phi_1$ непрерывна на компактном множестве $\tilde \Omega_0$
и дифференцируема во всех внутренних точках множества $\tilde \Omega_0$.
Отсюда функция $\Phi_1$
может достигать своего экстремума только в своих критических точках,
являющихся внутренними для $\tilde \Omega_0$, и на границе $\partial\tilde\Omega_0$.
В силу предположения индукции $\Phi_1(x) \leqslant r-1$ для всех  $x\in \partial \tilde\Omega_0$.
Рассмотрим $$\frac{\partial \Phi_1}{\partial \alpha_\ell}(\alpha_2, \ldots, \alpha_r)=
\left(\ln\left(1-\sum_{i=2}^r \alpha_i\right)-\ln \alpha_\ell\right)\Phi_1(\alpha_2, \ldots, \alpha_r).$$
Тогда $\frac{\partial \Phi}{\partial \alpha_\ell}(\alpha_2, \ldots, \alpha_r) = 0$ для всех  $2\leqslant \ell \leqslant r$ только при $\alpha_2=\ldots = \alpha_r = 1-\sum_{i=2}^r \alpha_i = \frac{1}{r}$.
Поскольку $\Phi\left(\frac{1}{r}, \frac{1}{r}, \ldots,  \frac{1}{r}\right)=r > r-1$, отсюда следует утверждение леммы.
\end{proof}

Положительный корень $\zeta$ многочлена $P$, определённого ниже,
будет затем использован
для получения оценки сверху для градуированных коразмерностей:

\begin{lemma}\label{LemmaEquationZetaUpperFrac}
Пусть $r\in\mathbb N$, а $\gamma_i\in\mathbb Z$ при $1\leqslant i \leqslant r$,
причём 
\begin{equation}\label{EquationGammaUpperFrac}
\gamma_1 \leqslant \gamma_2 \leqslant \ldots \leqslant \gamma_r,
\end{equation} $\gamma_1 < 0$. Рассмотрим уравнение 
\begin{equation}\label{EquationZetaUpperFrac}P(\zeta):=\sum_{i=1}^r \gamma_i \zeta^{\gamma_i-\gamma_1} = 0,\end{equation}
где $\zeta$~"--- неизвестное.
Если $\sum_{i=1}^r
\gamma_i \geqslant 0$,
то уравнение~(\ref{EquationZetaUpperFrac}) имеет только один положительный корень $\zeta$,
причём оказывается, что $\zeta \in (0; 1]$.
Если $\sum_{i=1}^r
\gamma_i < 0$, то $P(y) < 0$ для всех $y\in[0;1]$.
\end{lemma}
\begin{proof} 
Если $\gamma_r \leqslant 0$, то
все ненулевые коэффициенты многочлена $P$ отрицательны и $P(y) < 0$ для всех  $y\in[0;1]$. 

Пусть $\gamma_r > 0$. Из неравенства~(\ref{EquationGammaUpperFrac})
следует, что в последовательности коэффициентов многочлена $P$ есть только одна перемена
знака. Отсюда в силу теоремы Декарта уравнение (\ref{EquationZetaUpperFrac}) 
имеет только одно положительное решение~$\zeta$.
Определим число $m\in\mathbb N$ при помощи условия $$\gamma_1=\ldots=\gamma_m<\gamma_{m+1}.$$
Заметим, что $P(0) = m \gamma_1 < 0$, а $P(1)=\sum_{i=1}^r \gamma_i$.
Следовательно, если $\sum_{i=1}^r \gamma_i \geqslant 0$, то $\zeta \in (0; 1]$.  
Если же $P(1)=\sum_{i=1}^r \gamma_i < 0$, то $\zeta > 1$ и $P(y) < 0$ для всех  $y\in[0;1]$.
\end{proof}

Оказывается, что корень $\zeta$ многочлена $P$ является точкой экстремума функции $\Psi$,
которая вводится ниже.
Это свойство будет затем использовано при вычислении максимума функции $\Phi$
на многограннике $\Omega$. 

\begin{lemma}\label{LemmaMinimumGammaSumUpperFrac}
Пусть $r\in\mathbb N$, а $\gamma_i\in\mathbb Z$ при $1\leqslant i \leqslant r$,
причём $\gamma_1 \leqslant \gamma_2 \leqslant \ldots \leqslant \gamma_r$, $\gamma_1 < 0$.
Введём обозначение $\Psi(y):=\sum_{i=1}^r y^{\gamma_i}$.
Тогда \begin{enumerate}
\item
если $\sum_{i=1}^r
\gamma_i \geqslant 0$, то $\min_{y\in(0;1]} \Psi(y) =\Psi(\zeta)$,
где $\zeta \in (0; 1]$~"--- положительный корень многочлена~(\ref{EquationZetaUpperFrac});
\item если $\sum_{i=1}^r
\gamma_i \leqslant 0$, то $\min_{y\in(0;1]} \Psi(y) = r$.
\end{enumerate}
\end{lemma}
\begin{proof}
Заметим, что $\Psi'(y)= \sum_{i=1}^r \gamma_i y^{\gamma_i-1}$.
Отсюда во всех точках полуинтервала $(0;1]$ у функции $\Psi'(y)$ 
тот же знак, что и у многочлена $P(y)=\sum_{i=1}^r \gamma_i y^{\gamma_i-\gamma_1}$.  Кроме того, $\lim_{y\rightarrow +0} \Psi(y)= +\infty$. Из леммы~\ref{LemmaEquationZetaUpperFrac} следует, что если $\sum_{i=1}^r
\gamma_i \geqslant 0$, то $\min_{y\in(0;1]} \Psi(y) =\Psi(\zeta)$, а если $\sum_{i=1}^r
\gamma_i \leqslant 0$, то $\min_{y\in(0;1]} \Psi(y) = \Psi(1)=r$. (В случае $\sum_{i=1}^r
\gamma_i = 0$ справедливы равенства $\zeta = 1$ и $\Psi(\zeta)=r$.)
\end{proof}

Теперь мы готовы вычислить максимум функции $\Phi$. Для удобства мы заменим $\Omega$
на б\'ольшее множество $\tilde\Omega$,
а затем покажем, что максимум на обоих множествах один и тот же.

\begin{lemma}\label{LemmaMaximumFUpperFrac}
Пусть $r\in\mathbb N$, а $\gamma_i\in\mathbb Z$, где $1\leqslant i \leqslant r$,
причём $\gamma_1 \leqslant \gamma_2 \leqslant \ldots \leqslant \gamma_r$,
$\gamma_1 < 0$, $\sum_{i=1}^r
\gamma_i \geqslant 0$.
Обозначим $$
\Omega:=\lbrace(\alpha_1, \alpha_2, \ldots, \alpha_r)\in \mathbb R^r \mid \alpha_1 \geqslant \ldots\geqslant \alpha_r
\geqslant 0,\ \alpha_1+\alpha_2+\ldots+\alpha_r=1,\ \gamma_1\alpha_1 + \ldots + \gamma_r\alpha_r \leqslant 0 \rbrace$$
и $$
\tilde\Omega:=\lbrace(\alpha_1, \alpha_2, \ldots, \alpha_r)\in \mathbb R^r \mid \alpha_1,\ldots,\alpha_r \geqslant 0,\ \alpha_1+\alpha_2+\ldots+\alpha_r=1,\ \gamma_1\alpha_1 + \ldots + \gamma_r\alpha_r \leqslant 0 \rbrace.$$
Тогда $\max_{x\in\Omega} \Phi(x) = \max_{x\in\tilde\Omega} \Phi(x) =\sum_{i=1}^r \zeta^{\gamma_i}$,
где $\zeta \in (0; 1]$~"--- положительный корень уравнения~(\ref{EquationZetaUpperFrac}).
\end{lemma}
\begin{proof}
Как и в лемме~\ref{LemmaMaximumFPlainUpperFrac}, проведём доказательство индукцией по $r$.
Из условий $\gamma_1 < 0$ и $\sum_{i=1}^r
\gamma_i \geqslant 0$ следует, что $r \geqslant 2$.

База индукции не будет доказываться отдельно. Тем не менее из приводимых ниже рассуждений
будет следовать справедливость
утверждения леммы и при $r=2$,
поскольку в случае $r=2$ предположение индукции не будет использовано.
 
Снова выразим $\alpha_1 = 1-\sum_{i=2}^r \alpha_i$ через $\alpha_2, \alpha_3, \ldots, \alpha_r$
и исследуем поведение функции
$\Phi_1(\alpha_2, \ldots, \alpha_r) = \frac{1}{\left(1-\sum_{i=2}^r \alpha_i
\right)^{\left(1-\sum_{i=2}^r \alpha_i\right)} \alpha_2^{\alpha_2} \ldots 
\alpha_r^{\alpha_r}}$ на множестве
$$
\Omega_1:=\left\lbrace(\alpha_2, \ldots, \alpha_r)\in \mathbb R^{r-1} \mathrel{\biggl|} \alpha_2, \ldots, \alpha_r
\geqslant 0,\ 1-\sum_{i=2}^r \alpha_i\geqslant 0,\ \gamma_1+\sum_{i=2}^r(\gamma_i-\gamma_1)\alpha_i \leqslant 0 \right\rbrace.$$
Теперь из доказательства леммы~\ref{LemmaMaximumFPlainUpperFrac} 
следует, что единственной критической точкой функции $\Phi_1$
является точка $\left(\frac{1}{r},  \ldots,  \frac{1}{r}\right)$. Эта точка принадлежит множеству $\Omega_1$, если и только если
$\sum_{i=1}^r \gamma_i \leqslant 0$. Если $\sum_{i=1}^r
\gamma_i = 0$, то в силу леммы~\ref{LemmaMaximumFPlainUpperFrac} справедливо равенство $\max_{x\in\Omega} \Phi(x)=r$.
Поскольку в этом случае $\zeta = 1$, то лемма доказана.

Предположим теперь, что $\sum_{i=1}^r \gamma_i > 0$. Тогда непрерывная функция $\Phi_1$
достигает своего максимума на границе $\partial \Omega_1$. Заметим, что $\partial \Omega_1 = \Omega_2 \cup \bigcup_{i=1}^r \Omega_{1i}$,
где 
$$\Omega_{11}=\left\lbrace(\alpha_2, \ldots, \alpha_r)\in \mathbb R^{r-1} \mathrel{\biggl|} \alpha_2, \ldots, \alpha_r
\geqslant 0,\ 1-\sum_{i=2}^r \alpha_i= 0,\ \gamma_1+\sum_{i=2}^r(\gamma_i-\gamma_1)\alpha_i \leqslant 0 \right\rbrace,$$
\begin{equation*}\begin{split}\Omega_{1\ell}=\left\lbrace(\alpha_2, \ldots, \alpha_r)\in \mathbb R^{r-1} \mathrel{\biggl|} \alpha_2, \ldots, \alpha_r 
\geqslant 0,\ \alpha_\ell=0,\right. \\ \left. 1-\sum_{i=2}^r \alpha_i\geqslant 0,\ \gamma_1+\sum_{i=2}^r(\gamma_i-\gamma_1)\alpha_i \leqslant 0 \right\rbrace
\text{ при } 2\leqslant \ell \leqslant r,\end{split}\end{equation*}
$$\Omega_2=\left\lbrace(\alpha_2, \ldots, \alpha_r)\in \mathbb R^{r-1} \mathrel{\biggl|} \alpha_2, \ldots, \alpha_r
\geqslant 0,\ 1-\sum_{i=2}^r \alpha_i\geqslant 0,\ \gamma_1+\sum_{i=2}^r(\gamma_i-\gamma_1)\alpha_i = 0 \right\rbrace.$$

Докажем, что \begin{equation}\label{EquationOmega1iUpperFrac}\max\limits_{x\in \bigcup_{i=1}^r \Omega_{1i}} \Phi_1(x) < \sum_{j=1}^r \zeta^{\gamma_j},\end{equation} где $\zeta \in (0; 1]$~"---
положительный корень уравнения~(\ref{EquationZetaUpperFrac}).
Переходя к переменным $\alpha_1, \ldots, \hat \alpha_i, \ldots, \alpha_r$, получаем, что $\max_{x\in \Omega_{1i}} \Phi_1(x) = \max_{x\in \Omega'_i} \Phi(x)$ при $1\leqslant i \leqslant r$,
где \begin{equation*}\begin{split}\Omega_i'=\left\lbrace(\alpha_1, \ldots, \hat \alpha_i, \ldots, \alpha_r)\in \mathbb R^{r-1} \mathrel{\biggl|} \alpha_1, \ldots, \hat \alpha_i,\ldots, \alpha_r 
\geqslant 0, \right. \\ \left.  \alpha_1+\ldots+\hat\alpha_i+\ldots+\alpha_r=1,\ \gamma_1\alpha_1 + \ldots+ \widehat{\gamma_i\alpha_i}+\ldots + \gamma_r\alpha_r \leqslant 0 \right\rbrace.\end{split}\end{equation*}
(Для удобства будем обозначать функцию $\Phi(\theta_1, \ldots, \theta_m)=\frac{1}{\theta_1^{\theta_1}
\ldots \theta_m^{\theta_m}}$ одной и той же буквой $\Phi$ для всех  $m$.)

Заметим, что множество $\Omega_i'$ имеет тот же вид, что и множество $\tilde\Omega$
из формулировки леммы. Сперва рассмотрим случай, когда можно применить предположение
индукции, т.е. когда $\sum\limits_{\substack{\ell=1,\\ \ell\ne i}}^r
\gamma_\ell \geqslant 0$ и при этом либо $i\geqslant 2$,
либо $\gamma_2 < 0$. (Ясно, что в этом случае $r \geqslant 3$.)
Применяя предположение
индукции для $(r-1)$, получаем, что $\max_{x\in \Omega'_i} \Phi(x) = \sum\limits_{\substack{\ell=1,\\ \ell\ne i}}^r \left(\zeta'\right)^{\gamma_\ell}$, где $\sum\limits_{\substack{\ell=1,\\ \ell\ne i}}^r \gamma_\ell \left(\zeta'\right)^{\gamma_\ell-\gamma_1}=0$ при $i > 1$
и $\sum_{\substack{\ell=2}}^r \gamma_\ell \left(\zeta'\right)^{\gamma_\ell-\gamma_2}=0$
при $i=1$. В силу леммы~\ref{LemmaMinimumGammaSumUpperFrac}
справедливо равенство
\begin{equation}\label{EquationOmegaMaxPsiUpperFrac} \max_{x\in \Omega_{1i}} \Phi_1(x)=\max_{x\in \Omega'_i} \Phi(x) = \min_{y\in(0;1]} \sum_{\substack{\ell=1,\\ \ell\ne i}}^r y^{\gamma_\ell} < \min_{y\in(0;1]} \sum_{\ell=1}^r y^{\gamma_\ell}=\sum_{j=1}^r \zeta^{\gamma_j}.\end{equation}

Если $i=1$, а $\gamma_2 \geqslant 0$, то в силу того, что
для всех оставшихся $2 \leqslant j \leqslant r$ выполнено условие $\gamma_j, \alpha_j\geqslant 0$, неравенство
$$\gamma_1\alpha_1 + \ldots+ \widehat{\gamma_i\alpha_i}+\ldots + \gamma_r\alpha_r \leqslant 0$$
равносильно равенству нулю всех $\alpha_j$, для которых $\gamma_j\ne 0$.
Отсюда
 \begin{equation*}\begin{split}\Omega_i'=\Omega_1'=\left\lbrace(\alpha_2, \ldots, \alpha_r)\in \mathbb R^{r-1} \mathrel{\biggl|} \alpha_2, \ldots,  \alpha_r 
\geqslant 0,\right. \\ \left. \alpha_2+\ldots+\ldots+\alpha_r=1,\ \alpha_\ell = \alpha_{\ell+1}=\ldots=\alpha_r=0 \right\rbrace,\end{split}\end{equation*}
где число $2\leqslant \ell\leqslant r$ задано условиями $\gamma_2=\ldots=\gamma_{\ell-1}=0$
и $\gamma_\ell > 0$. Тогда из леммы~\ref{LemmaMaximumFPlainUpperFrac}
следует, что $\max_{x\in \Omega'_1} \Phi(x)= \ell-2$. Поскольку $\sum_{j=1}^r \zeta^{\gamma_j} > \sum_{j=2}^{\ell-1} \zeta^{\gamma_j} = \ell-2$, получаем, что $\max_{x\in \Omega_{11}} \Phi_1(x) < \sum_{j=1}^r \zeta^{\gamma_j}$.

Рассмотрим теперь последний случай, когда
$$\sum\limits_{\substack{\ell=1,\\ \ell\ne i}}^r
\gamma_\ell < 0.$$
 Тогда $\left(\frac{1}{r-1}, \frac{1}{r-1}, \ldots,  \frac{1}{r-1}\right)\in \Omega'_i $, 
и в силу леммы~\ref{LemmaMaximumFPlainUpperFrac} справедливо равенство $$\max_{x\in \Omega'_i} \Phi(x) = r-1.$$
Снова применяя лемму~\ref{LemmaMinimumGammaSumUpperFrac}, получаем неравенство~(\ref{EquationOmegaMaxPsiUpperFrac}).
Таким образом, неравенство~(\ref{EquationOmega1iUpperFrac}) доказано.

Докажем теперь, что $\max_{x\in \Omega_2} \Phi_1(x) = \sum_{i=1}^r \zeta^{\gamma_i}$,
где $\zeta \in (0; 1]$~"--- положительный корень уравнения~(\ref{EquationZetaUpperFrac}).
 
  Если $r=2$, то $\Omega_2=\left\lbrace -\frac{\gamma_1}{\gamma_2-\gamma_1} \right\rbrace$,
   $\zeta=\left(-\frac{\gamma_1}{\gamma_2}\right)^{\frac{1}{\gamma_2-\gamma_1}}$,
  \begin{equation*}\begin{split}
  \Phi_1\left(-\frac{\gamma_1}{\gamma_2-\gamma_1}\right)
   =\left(\frac{\gamma_2}{\gamma_2-\gamma_1}\right)^{-\frac{\gamma_2}{\gamma_2-\gamma_1}}\left(-\frac{\gamma_1}{\gamma_2-\gamma_1}\right)^{\frac{\gamma_1}{\gamma_2-\gamma_1}}=\\
 =(\gamma_2-\gamma_1)
   \gamma_2^{-\frac{\gamma_2}{\gamma_2-\gamma_1}}(-\gamma_1)^{\frac{\gamma_1}{\gamma_2-\gamma_1}}
    = \left(-\frac{\gamma_1}{\gamma_2}\right)^{\frac{\gamma_1}{\gamma_2-\gamma_1}}
   +\left(-\frac{\gamma_1}{\gamma_2}\right)^{\frac{\gamma_2}{\gamma_2-\gamma_1}}
   =\zeta^{\gamma_1}+\zeta^{\gamma_2}
  \end{split}\end{equation*}

   Следовательно, без ограничения общности мы можем считать, что $r\geqslant 3$.
   Определим $1\leqslant m < r$ при помощи условия $\gamma_1=\ldots=\gamma_m<\gamma_{m+1}.$
   Тогда для всех  $(\alpha_2, \ldots, \alpha_r)\in \Omega_2$
   получаем, что $\gamma_1+\sum_{i=m+1}^r(\gamma_i-\gamma_1)\alpha_i = 0$
   и \begin{equation}\label{EquationAlpha(m+1)UpperFrac}\alpha_{m+1}= -\frac{1}{\gamma_{m+1}-\gamma_1}\left(\gamma_1+\sum_{i=m+2}^r(\gamma_i-\gamma_1)\alpha_i
   \right).\end{equation}
  Подставим $\alpha_{m+1}$ и заметим, что
  $\max_{x\in \Omega_2} \Phi_1(x) = \max_{x\in \Omega_3} \Phi_2(x)$,
  где 
$$\Phi_2(\alpha_2, \ldots, \alpha_m, \alpha_{m+2},\ldots, \alpha_r) = 
\Phi(\alpha_1, \ldots, \alpha_r),$$
$$
\Omega_3:=\left\lbrace(\alpha_2, \ldots, \alpha_m, \alpha_{m+2},\ldots, \alpha_r)\in \mathbb R^{r-2} \mathrel{\biggl|} \alpha_1, \ldots, \alpha_r
\geqslant 0 \right\rbrace,$$
  $$\alpha_1 = 1-\sum\limits_{\substack{i=2, \\ i\ne m+1}}^r \alpha_i + \frac{1}{\gamma_{m+1}-\gamma_1}\left(\gamma_1+\sum_{i=m+2}^r(\gamma_i-\gamma_1)\alpha_i
\right),$$  
  а $\alpha_{m+1}$ задано равенством~(\ref{EquationAlpha(m+1)UpperFrac}).

Рассмотрим частные производные
\begin{equation*}\begin{split} \frac{\partial \Phi_2}{\partial \alpha_i}(\alpha_2, \ldots, \alpha_m, \alpha_{m+2},\ldots, \alpha_r)=
\left( (-\ln \alpha_1-1)\frac{\partial \alpha_1}{\partial \alpha_i} -\ln \alpha_i-1+ \right.\\
\left.+
(-\ln \alpha_{m+1}-1)\frac{\partial \alpha_{m+1}}{\partial \alpha_i}
 \right) \Phi_2(\alpha_2, \ldots, \alpha_m, \alpha_{m+2},\ldots, \alpha_r).
\end{split}\end{equation*}

При $2\leqslant i \leqslant m$ справедливы равенства $\frac{\partial \alpha_1}{\partial \alpha_i} = -1$,
$\frac{\partial \alpha_{m+1}}{\partial \alpha_i} = 0$ и
$$\frac{\partial \Phi_2}{\partial \alpha_i}(\alpha_2, \ldots, \alpha_m, \alpha_{m+2},\ldots, \alpha_r)=
(\ln \alpha_1 -\ln \alpha_i) \Phi_2(\alpha_2, \ldots, \alpha_m, \alpha_{m+2},\ldots, \alpha_r).$$

При $m+2\leqslant i \leqslant r$ справедливы равенства $\frac{\partial \alpha_1}{\partial \alpha_i} = 
\frac{\gamma_i-\gamma_{m+1}}{\gamma_{m+1}-\gamma_1}$,
$\frac{\partial \alpha_{m+1}}{\partial \alpha_i} = \frac{\gamma_1-\gamma_i}{\gamma_{m+1}-\gamma_1}$ и
\begin{equation*}\begin{split}
\frac{\partial \Phi_2}{\partial \alpha_i}(\alpha_2, \ldots, \alpha_m, \alpha_{m+2},\ldots, \alpha_r)=\\
=\left(-\frac{\gamma_i-\gamma_{m+1}}{\gamma_{m+1}-\gamma_1}\ln \alpha_1 -\ln \alpha_i-
\frac{\gamma_1-\gamma_i}{\gamma_{m+1}-\gamma_1}\ln \alpha_{m+1}\right) \Phi_2(\alpha_2, \ldots, \alpha_m, \alpha_{m+2},\ldots, \alpha_r)
\end{split}\end{equation*}

Если $(\alpha_2, \ldots, \alpha_m, \alpha_{m+2},\ldots, \alpha_r)\in \Omega_3$ 
является критической точкой для функции $\Phi_2$, то $$\frac{\partial \Phi_2}{\partial \alpha_i}(\alpha_2, \ldots, \alpha_m, \alpha_{m+2},\ldots, \alpha_r)=0$$ для всех  $2\leqslant i \leqslant m$ и $m+2\leqslant i \leqslant r$, что эквивалентно равенствам
$$\left\lbrace \begin{array}{lllll} \alpha_i & = & \alpha_1 & \text{ при } & 2\leqslant i \leqslant m,  \\
\alpha_i & = & \alpha_1^{\left(\frac{\gamma_{m+1}-\gamma_i}{\gamma_{m+1}-\gamma_1}\right)}
\alpha_{m+1}^{\left(\frac{\gamma_i-\gamma_1}{\gamma_{m+1}-\gamma_1}\right)}
& \text{ при } & m+2\leqslant i \leqslant r, \\
\alpha_1 & = & 1-\sum\limits_{\substack{i=2, \\ i\ne m+1}}^r \alpha_i + \frac{1}{\gamma_{m+1}-\gamma_1}\left(\gamma_1+\sum_{i=m+2}^r(\gamma_i-\gamma_1)\alpha_i
\right), \\
\alpha_{m+1} & = & -\frac{1}{\gamma_{m+1}-\gamma_1}\left(\gamma_1+\sum_{i=m+2}^r(\gamma_i-\gamma_1)\alpha_i
   \right),
 \end{array}\right.$$
  т.е.
$$\left\lbrace \begin{array}{lllll} \alpha_i & = & \alpha_1 & \text{ при } & 2\leqslant i \leqslant m,  \\
\alpha_i & = & \alpha_1
\left(\frac{\alpha_{m+1}}{\alpha_1}\right)^{\frac{\gamma_i-\gamma_1}{\gamma_{m+1}-\gamma_1}}
& \text{ при } & m+2\leqslant i \leqslant r, \\
\alpha_1 & = & 1-\sum\limits_{\substack{i=2, \\ i\ne m+1}}^r \alpha_i + \frac{1}{\gamma_{m+1}-\gamma_1}\left(\gamma_1+\sum_{i=m+2}^r(\gamma_i-\gamma_1)\alpha_i
\right), \\
\alpha_{m+1} & = & -\frac{1}{\gamma_{m+1}-\gamma_1}\left(\gamma_1+\sum_{i=m+2}^r(\gamma_i-\gamma_1)\alpha_i
   \right).
 \end{array}\right.$$

Заметим, что поскольку нас интересуют только внутренние критические точки функции~$\Phi_2$ на множестве $\Omega_3 \subset \mathbb R^{r-2}$, можно считать, что все $\alpha_i > 0$.

Совершая равносильные преобразования, получаем, что

\begin{equation}\label{EquationAlphaRelationUpperFrac1}\left\lbrace \begin{array}{lllll} \alpha_i & = & \alpha_1
\left(\frac{\alpha_{m+1}}{\alpha_1}\right)^{\frac{\gamma_i-\gamma_1}{\gamma_{m+1}-\gamma_1}}
& \text{ для всех } 1\leqslant i \leqslant r, \\
\sum\limits_{i=1}^r \alpha_i & = & 1, \\
\sum_{i=1}^r \gamma_i \alpha_i & = & 0.
 \end{array}\right.\end{equation}

Введём теперь дополнительную переменную $\xi := \left(\frac{\alpha_{m+1}}{\alpha_1}\right)^{\frac{1}{\gamma_{m+1}-\gamma_1}}$ и получим, что

\begin{equation}\label{EquationAlphaRelationUpperFrac2}\left\lbrace \begin{array}{lllll}
\xi &= & \left(\frac{\alpha_{m+1}}{\alpha_1}\right)^{\frac{1}{\gamma_{m+1}-\gamma_1}}, \\
 \alpha_i & = & \alpha_1
\xi^{\gamma_i-\gamma_1}
& \text{ для всех } 1\leqslant i \leqslant r, \\
\alpha_1 \sum\limits_{i=1}^r \xi^{\gamma_i-\gamma_1} & = & 1, \\
 \sum_{i=1}^r \gamma_i \xi^{\gamma_i-\gamma_1} & = & 0.
 \end{array}\right.\end{equation}
 Заметим теперь, что первое уравнение является следствием второго уравнения
 при $i=m+1$. Следовательно, первоначальная система эквивалентна
 равенствам
\begin{equation}\label{EquationAlphaSolutionUpperFrac}\left\lbrace \begin{array}{lllll}
 \alpha_i & = & \frac{\xi^{\gamma_i-\gamma_1}}{\sum_{i=1}^r \xi^{\gamma_i-\gamma_1}}
& \text{ при } 1\leqslant i \leqslant r, \\
 \sum_{i=1}^r \gamma_i \xi^{\gamma_i-\gamma_1} & = & 0.
 \end{array}\right. \end{equation}
В силу леммы~\ref{LemmaEquationZetaUpperFrac}
последнее уравнение имеет единственное решение $\zeta \in (0;1]$.
Следовательно, точка $$(\alpha_2,\ldots,\alpha_m, \alpha_{m+2}, \ldots, \alpha_r),$$
заданная равенствами~(\ref{EquationAlphaSolutionUpperFrac}),
является единственной внутренней критической точкой функции $\Phi_2$.
Используя~(\ref{EquationAlphaRelationUpperFrac1}) и~(\ref{EquationAlphaRelationUpperFrac2}), получаем, что
  \begin{equation*}\begin{split}\Phi_2(\alpha_2,\ldots,\alpha_m, \alpha_{m+2}, \ldots, \alpha_r)
  = \Phi(\alpha_1,\ldots, \alpha_r)=\frac{1}{\alpha_1^{\alpha_1} \alpha_2^{\alpha_2} \ldots \alpha_r^{\alpha_r}}=
   \\ = \frac{1}{\alpha_1^{\alpha_1+\ldots+\alpha_r} \zeta^{\alpha_1(\gamma_1-\gamma_1)} \ldots \zeta^{\alpha_r(\gamma_r-\gamma_1)}}=
  \frac{1}{\alpha_1 \zeta^{-\gamma_1}}=\sum_{i=1}^r \zeta^{\gamma_i}.
  \end{split}\end{equation*}
 
 Заметим теперь, что значения функции $\Phi_2$ на $\partial \Omega_3$ равны значениям функции $\Phi_1$
 в соответствующих точках множества $\bigcup_{i=1}^r \Omega_{1i}$.
 Следовательно, $$\max_{x\in \tilde\Omega} \Phi(x)=\max_{x\in \Omega_1} \Phi_1(x)=\max_{x\in \Omega_3} \Phi_2(x) = \sum_{i=1}^r \zeta^{\gamma_i}.$$
 Поскольку из~(\ref{EquationAlphaSolutionUpperFrac}) следует, что $\alpha_1 \geqslant \alpha_2\geqslant \ldots
 \geqslant \alpha_r$, получаем, что $$\max_{x\in \Omega} \Phi(x)=\max_{x\in \tilde\Omega} \Phi(x) = 
 \sum_{i=1}^r \zeta^{\gamma_i}.$$
 \end{proof}

Отсюда сразу же получается требуемая оценка сверху:

\begin{theorem}\label{TheoremUpperFrac}
Пусть $A=\bigoplus_{t\in T} A^{(t)}$~"--- такая конечномерная ассоциативная $T$-градуированная алгебра над полем $\mathbbm{k}$
характеристики $0$ для некоторого множества $T$, что $A/J(A) \cong M_k(\mathbbm{k})$
для некоторого $k\in\mathbb N$ и $A^{(t)}\cap J(A) = 0$ для всех  $t\in T$.
Предположим, что числа
 $\gamma_\ell$, где $1\leqslant \ell \leqslant r$, $r:=\dim A$, введённые перед
 леммой~\ref{LemmaInequalityLambdaUpperFrac}, удовлетворяют условию
 $\sum_{i=1}^r
\gamma_i \geqslant 0$.  Тогда $$\mathop{\overline\lim}\limits_{n\rightarrow\infty}
 \sqrt[n]{c_n^{T\text{-}\mathrm{gr}}(A)}\leqslant \sum_{i=1}^r \zeta^{\gamma_i},$$
где $\zeta$~"--- положительный корень
уравнения~(\ref{EquationZetaUpperFrac}).
\end{theorem}
\begin{proof}
Достаточно применить леммы~\ref{LemmaUpperBoundFUpperFrac} и~\ref{LemmaMaximumFUpperFrac}.
\end{proof}
\begin{remark}
Если $A$~"--- конечномерная $T$-градуированно 
простая алгебра над алгебраически замкнутым полем $\mathbbm{k}$
характеристики $0$ для некоторой конечной $0$-простой полугруппы $T$ с тривиальными максимальными подгруппами (например, для некоторой ленты правых нулей), то в силу лемм~\ref{LemmaRadicalSemigroupGradedSimple} и \ref{LemmaReesGrSimpleOtherProperties} для алгебры $A$
существует такое $k\in\mathbb N$, что $A/J(A) \cong M_k(\mathbbm{k})$, и для всех $t\in T$
выполнены условия $A^{(t)}\cap J(A) = 0$. 
 \end{remark}

\section{Случай $A/J(A)\cong M_2(\mathbbm{k})$ и $\PIexp^{T\text{-}\mathrm{gr}}(A)=\dim A$}\label{SectionFracM2Equal}

Пусть $A=\bigoplus_{t\in T} A^{(t)}$~"--- конечномерная
$T$-градуированно простая 
алгебра над полем $\mathbbm{k}$ характеристики $0$
для некоторой ленты правых нулей $T$.
В этом и следующем параграфах вычисляется $\PIexp^{T\text{-}\mathrm{gr}}(A)$ в случае, когда $A/J(A) \cong M_2(\mathbbm{k})$.
 
 Для $t\in T$ введём обозначение $I_t := \pi(A^{(t)})$,
 где $\pi \colon A \twoheadrightarrow A/J(A)$~"--- естественный сюръективный
гомоморфизм. 

 Заметим, что, поскольку $\dim A < +\infty$, только конечное число левых идеалов $I_t$ ненулевые.
  Введём обозначения $T_0 := \lbrace t\in T \mid \dim I_t=2 \rbrace$ и $T_1 = \lbrace t\in T \mid I_t = A/J(A) \rbrace$. При этом $I_t = 0$ для всех  $t\notin T_0 \sqcup T_1$. Более того, из $A^{(t)} \cap \ker \pi =0$ для всех $t\in T$ следует, что $r:= \dim A = 2|T_0|+4|T_1|$.

  Введём на множестве $T_0$ следующее отношение эквивалентности:  $t_1 \sim t_2$, если $I_{t_1}=I_{t_2}$. Поскольку $A/J(A) \cong M_2(\mathbbm{k})$, все неприводимые $A/J(A)$-модули двумерны
  и изоморфны друг другу. Отсюда $A/J(A)=I_{t_1}\oplus I_{t_2}$
  для всех  $I_{t_1}\ne I_{t_2}$, $t_1, t_2 \in T_0$.
 
 Покажем теперь, что если мощности классов эквивалентности
 удовлетворяют некоторой разновидности неравенства треугольника,
 то все элементы множества $T_0$ можно разбить по парам и, возможно, одной тройке так, что внутри каждой пары или тройки
 все элементы попарно неэквивалентны.
   
   \begin{lemma}\label{LemmaTrianglePairsAltFracM2}
Пусть $T_0$~"--- конечное непустое множество, на котором задано отношение эквивалентности $\sim$. 
    Предпололожим, что \begin{equation}\label{EquationEquivTriangleFracM2}|\bar t_0| \leqslant \sum\limits_{\substack{\bar t \in {T_0}/{\sim}, \\
   \bar t \ne \bar t_0}} |\bar t| \text{ для всех  }\bar t_0 \in {T_0}/{\sim}
   .\end{equation}
   Тогда можно выбрать такие элементы $t_i$, что $\lbrace t_1, \ldots, t_{|T_0|} \rbrace = T_0$ и \begin{enumerate}
\item при $2 \mid |T_0|$ для всех $1\leqslant i \leqslant \frac{|T_0|}{2}$ выполнено условие $t_{2i-1} \nsim t_{2i}$;
\item при $2 \nmid |T_0|$ для всех $1\leqslant i \leqslant \frac{|T_0|-1}{2}$ выполнено условие $t_{2i-1} \nsim t_{2i}$, а элементы $t_{|T_0|-2}$, $t_{|T_0|-1}$, $t_{|T_0|}$ попарно неэквивалентны.
\end{enumerate}
   \end{lemma}
  \begin{remark}
  Заметим, что условие~(\ref{EquationEquivTriangleFracM2}) не выполняется,
если и только если существует такой класс эквивалентности
  $\bar t_0 \in {T_0}/{\sim}$, что $|\bar t_0| > \frac{|T_0|}{2}$.
  \end{remark}
  \begin{proof}[Доказательство леммы~\ref{LemmaTrianglePairsAltFracM2}.]
  Проведём доказательство индукцией по $|T_0|$.
  
  Прежде всего заметим, что из~(\ref{EquationEquivTriangleFracM2}) следует, что $|{T_0}/{\sim}| \geqslant 2$. Предположим, что $|{T_0}/{\sim}| = 2$. Тогда из~(\ref{EquationEquivTriangleFracM2})
  следует, что $|\bar t_1| = |\bar t_2|$, где  ${T_0}/{\sim} = \lbrace \bar t_1, \bar t_2 \rbrace$.
  Следовательно, можно определить $\lbrace t_1, \ldots, t_{|T_0|}\rbrace = T_0$
  при помощи равенств $\lbrace t_1,  t_3, \ldots, t_{|T_0|-1}\rbrace := \bar t_1$
  и $\lbrace t_2,  t_4, \ldots, t_{|T_0|}\rbrace := \bar t_2$, и лемма доказана.
  
  Если $|T_0| = 3$, то все элементы множества $T_0$
попарно неэквивалентны и лемма снова доказана.
  
  Предположим теперь, что $|T_0| > 3$. Выберем классы $\bar t_1$ и $\bar t_2$
  с  наибольшим числом элементов. Фиксируем произвольные элементы $t_1 \in \bar t_1$
  и $t_2 \in \bar t_2$. Заметим, что для множества $T_0\backslash \lbrace t_1, t_2 \rbrace$
  и прежнего отношения эквивалентности $\sim$ по-прежнему выполняется условие~(\ref{EquationEquivTriangleFracM2}).
 В силу предположения индукции можно выбрать такие $\lbrace t_3, \ldots, t_{|T|} \rbrace = T_0\backslash \lbrace t_1, t_2 \rbrace$, что элементы $t_1, \ldots, t_{|T_0|}$ удовлетворяют всем условиям леммы.
  \end{proof}
  
Неравенство~(\ref{EquationEquivTriangleFracM2}) будет использовано ниже для того, чтобы различать случаи $\PIexp^{T\text{-}\mathrm{gr}}(A) = \dim A$ и $\PIexp^{T\text{-}\mathrm{gr}}(A) < \dim A$.

Докажем сперва следующий вспомогательный результат, касающийся матричных единиц алгебры $M_2(\mathbbm{k})$:
\begin{lemma}\label{LemmaAlphaBetaDeterminant}
Пусть $\alpha, \beta, \tilde \alpha, \tilde \beta \in \mathbbm{k}$.
Тогда
\begin{enumerate}
\item $$[\alpha e_{11} + \beta e_{12}, \alpha e_{21} + \beta e_{22}]
= \left(\begin{array}{rr}  \alpha\beta & \beta^2 \\  -\alpha^2 & -\alpha\beta \end{array}\right)
= \left(\begin{array}{rr}  \beta & 0 \\ -\alpha  & 0 \end{array}\right)
\left(\begin{array}{rr}  \alpha & \beta \\  0 & 0 \end{array}\right);$$
\item
\begin{equation*}\begin{split}
[\alpha e_{11} + \beta e_{12}, \alpha e_{21} + \beta e_{22}][\tilde\alpha e_{11} + \tilde\beta e_{12}, \tilde\alpha e_{21} + \tilde\beta e_{22}]+\\
 +[\tilde\alpha e_{11} + \tilde\beta e_{12}, \tilde\alpha e_{21} + \tilde\beta e_{22}][\alpha e_{11} + \beta e_{12}, \alpha e_{21} + \beta e_{22}]=\\
= -\left|\begin{array}{rr}  \alpha & \beta \\ \tilde\alpha  & \tilde \beta \end{array}\right|^2(e_{11}+e_{22})
 \end{split}\end{equation*}
\end{enumerate}
\end{lemma}
\begin{proof}
Первое равенство проверяется непосредственно.

Для доказательства второго равенство достаточно заметить, что
\begin{equation*}\begin{split}[\alpha e_{11} + \beta e_{12}, \alpha e_{21} + \beta e_{22}][\tilde\alpha e_{11} + \tilde\beta e_{12}, \tilde\alpha e_{21} + \tilde\beta e_{22}] =\\=
\left(\begin{array}{rr}  \beta & 0 \\ -\alpha  & 0 \end{array}\right)
\left(\begin{array}{rr}  \alpha & \beta \\  0 & 0 \end{array}\right)
\left(\begin{array}{rr}  \tilde\beta & 0 \\ -\tilde\alpha  & 0 \end{array}\right)
\left(\begin{array}{rr}  \tilde\alpha & \tilde\beta \\  0 & 0 \end{array}\right)
=\\=\left(\begin{array}{rr}  \beta & 0 \\ -\alpha  & 0 \end{array}\right)
\left(\begin{array}{rr}  \alpha\tilde\beta - \beta\tilde\alpha  & 0 \\  0 & 0 \end{array}\right)
\left(\begin{array}{rr}  \tilde\alpha & \tilde\beta \\  0 & 0 \end{array}\right)
=\\
=\left|\begin{array}{rr}  \alpha & \beta \\ \tilde\alpha  & \tilde \beta \end{array}\right|
\left(\begin{array}{rr}  \beta & 0 \\ -\alpha  & 0 \end{array}\right)
\left(\begin{array}{rr}  \tilde\alpha & \tilde\beta \\  0 & 0 \end{array}\right)
=
\left|\begin{array}{rr}  \alpha & \beta \\ \tilde\alpha  & \tilde \beta \end{array}\right|
\left(\begin{array}{rr}  \tilde\alpha \beta & \beta \tilde\beta \\ -\alpha\tilde\alpha & -\alpha\tilde\beta \end{array}\right)
.\end{split}\end{equation*}
\end{proof} 
  
  Теперь мы можем доказать существование полилинейного $\mathbbm{k}^T$-многочлена,
не являющегося  $\mathbbm{k}^T$-тождеством, с достаточным числом наборов переменных,
по которым он кососимметричен:
  \begin{lemma}\label{LemmaGrPIexpTriangleAltFracM2}
  Пусть $T_0, T_1 \subseteq T$ и $\sim$~"--- соответственно, подмножества и отношение эквивалентности,
  заданные в начале \S\ref{SectionFracM2Equal}.
  Предположим, что либо справедливо неравенство~(\ref{EquationEquivTriangleFracM2}),
  либо $T_0=\varnothing$. Тогда можно выбрать такое число $n_0 \in \mathbb N$, что для любого $n\geqslant n_0$ существуют такие $\mathbbm{k}^T$-многочлен $f \in P^{\mathbbm{k}^T}_n \backslash
\Id^{\mathbbm{k}^T}(A)$ и попарно непересекающиеся подмножества $X_1$, \ldots, $X_{2k} \subseteq \lbrace x_1, \ldots, x_n\rbrace$, что $k = \left[\frac{n-n_0}{2\dim A}\right]$, 
$|X_1| = \ldots = |X_{2k}|=\dim A$, а $\mathbbm{k}^T$-многочлен $f$ кососимметричен по переменным каждого множества $X_j$.
  \end{lemma}
  \begin{proof}
    Заметим, что многочлен $$f_0(x_1,\ldots, x_4, y_1, \ldots, y_4)=\sum_{\sigma,\rho \in S_4} \sign(\sigma\rho) x_{\sigma(1)}\ y_{\rho(1)}\ x_{\sigma(2)}x_{\sigma(3)}x_{\sigma(4)}\ y_{\rho(2)}y_{\rho(3)}y_{\rho(4)}$$
    не является полиномиальным тождеством для алгебры $M_2(\mathbbm{k})$,
а его значения пропорциональны единичной матрице.
  (См., например, теорему~5.7.4 в~\cite{ZaiGia}.)
  
  Введём обозначения \begin{equation*}\begin{split}f_{t_1, t_2}(x_{t_1,1}, x_{t_1,2}, x_{t_2,1}, x_{t_2,2}) :=\left[x^{q_{t_1}}_{t_1,1}, x^{q_{t_1}}_{t_1,2}\right]\left[x^{q_{t_2}}_{t_2,1}, x^{q_{t_2}}_{t_2,2}\right]+
\left[x^{q_{t_2}}_{t_2,1}, x^{q_{t_2}}_{t_2,2}\right]  \left[x^{q_{t_1}}_{t_1,1}, x^{q_{t_1}}_{t_1,2}\right]\end{split}\end{equation*}
  и \begin{equation*}\begin{split}f_{t_1, t_2, t_3}(x_{t_1,1}, x_{t_1,2}, x_{t_2,1}, x_{t_2,2}, x_{t_3,1}, x_{t_3,2}) := \left[x^{q_{t_1}}_{t_1,1}, x^{q_{t_1}}_{t_1,2}\right] \left[x^{q_{t_3}}_{t_3,1}, x^{q_{t_3}}_{t_3,2}\right] \left[x^{q_{t_2}}_{t_2,1}, x^{q_{t_2}}_{t_2,2}\right] -\\ -
   \left[x^{q_{t_2}}_{t_2,1}, x^{q_{t_2}}_{t_2,2}\right] \left[x^{q_{t_3}}_{t_3,1}, x^{q_{t_3}}_{t_3,2} \right]
\left[x^{q_{t_1}}_{t_1,1}, x^{q_{t_1}}_{t_1,2}\right].\end{split}\end{equation*}
  
   Пусть $\lbrace t_1, \ldots, t_{|T_0|}\rbrace = T_0$~"--- элементы из
   леммы~\ref{LemmaTrianglePairsAltFracM2}.
  
  В случае, когда $2 \mid |T_0|$, положим \begin{equation*}\begin{split}
  f_1 := z_1 \ldots z_{n-(\dim A)2k}\prod_{i=1}^k
   \left(
   \prod_{t\in T_1}
    f_0(x^{q_t}_{i,t,1},\ldots, x^{q_t}_{i,t,4}, y^{q_t}_{i,t,1}, \ldots, y^{q_t}_{i,t,4})
     \right) \cdot \\ \cdot
     \left(
     \prod_{\ell=1}^{\frac{|T_0|}2} f_{t_{2\ell-1}, t_{2\ell}}(x_{i,t_{2\ell-1},1}, x_{i,t_{2\ell-1},2}, x_{i,t_{2\ell},1}, x_{i,t_{2\ell},2}) \cdot\right. \\ \cdot
    f_{t_{2\ell-1}, t_{2\ell}}(y_{i,t_{2\ell-1},1}, y_{i,t_{2\ell-1},2}, y_{i,t_{2\ell},1}, y_{i,t_{2\ell},2})
      \Biggr).\end{split}\end{equation*}
  В случае, когда $2 \nmid |T_0|$, положим \begin{equation*}\begin{split}
  f_1 := z_1 \ldots z_{n-(\dim A)2k}\prod_{i=1}^k
   \left(
   \prod_{t\in T_1} f_0(x^{q_t}_{i,t,1},\ldots, x^{q_t}_{i,t,4}, y^{q_t}_{i,t,1}, \ldots, y^{q_t}_{i,t,4})
    \right) \cdot \\ \cdot
         \left(
    \prod_{\ell=1}^{\frac{|T_0|-3}2} f_{t_{2\ell-1}, t_{2\ell}}(x_{i,t_{2\ell-1},1}, x_{i,t_{2\ell-1},2}, x_{i,t_{2\ell},1}, x_{i,t_{2\ell},2}) \cdot \right. \\ \cdot
     f_{t_{2\ell-1}, t_{2\ell}}(y_{i,t_{2\ell-1},1}, y_{i,t_{2\ell-1},2}, y_{i,t_{2\ell},1}, y_{i,t_{2\ell},2})
    \Biggr)\cdot \\ \cdot
    f_{t_{|T_0|-2}, t_{|T_0|-1}, t_{|T_0|}}
    (x_{i,t_{|T_0|-2},1},
     x_{i,t_{|T_0|-2},2};\ 
      x_{i,t_{|T_0|-1},1},
       x_{i,t_{|T_0|-1},2};\ 
        x_{i,t_{|T_0|},1}, 
        x_{i,t_{|T_0|},2})\cdot \\ \cdot
    f_{t_{|T_0|-2}, t_{|T_0|-1}, t_{|T_0|}}
    (y_{i,t_{|T_0|-2},1},
     y_{i,t_{|T_0|-2},2};\ 
      y_{i,t_{|T_0|-1},1},
       y_{i,t_{|T_0|-1},2};\ 
        y_{i,t_{|T_0|},1}, 
        y_{i,t_{|T_0|},2}).\end{split} \end{equation*}

Докажем, что $f_1 \notin \Id^{\mathbbm{k}^T}(A)$. 
В случае, когда $2 \mid |T_0|$, в качестве $\psi$ возьмём произвольный изоморфизм  $A/J(A) \mathrel{\widetilde\rightarrow} M_2(\mathbbm{k})$.
В случае, когда $2 \nmid |T_0|$, определим $\psi$ следующим образом. Во-первых, заметим, что из $t_{|T_0|-2} \nsim t_{|T_0|-1}$ следует, что $A/J(A)= I_{t_{|T_0|-2}}\oplus I_{t_{|T_0|-1}}$. В силу теоремы~\ref{TheoremSumLeftIdealsMatrix} существует такой изоморфизм $\psi \colon A/J(A) \mathrel{\widetilde\rightarrow} M_2(\mathbbm{k})$,
что $\psi(I_{t_{|T_0|-2}})=\langle e_{11}, e_{21} \rangle_\mathbbm{k}$
и $\psi(I_{t_{|T_0|-1}})=\langle e_{12}, e_{22} \rangle_\mathbbm{k}$.
В случае, когда $2 \nmid |T_0|$, под $\psi$ будем понимать изоморфизм, выбранный именно таким образом.

Из леммы~\ref{LemmaLeftIdealMatrix} следует, что для любого $t \in T_0$ существуют такие
$\alpha_t,\beta_t \in \mathbbm{k}$, что  $$\psi(I_t)=\langle\alpha_t e_{i1} + \beta_t e_{i2} \mid i=1,2 \rangle_\mathbbm{k}.$$
В случае, когда $2 \nmid |T_0|$, в силу выбора изоморфизма $\psi$ можно считать, что $$(\alpha_{t_{|T_0|-2}}, \beta_{t_{|T_0|-2}})=(1, 0)\text{ и }(\alpha_{t_{|T_0|-1}}, \beta_{t_{|T_0|-1}})=
(0, 1).$$ Заметим при этом, что $I_{t_1} = I_{t_2}$, если и только если строки $(\alpha_{t_1}, \beta_{t_1})$ и  $(\alpha_{t_2}, \beta_{t_2})$ пропорциональны.

Выберем некоторый такой элемент $e\in A$, что $\psi\pi(e)$~"--- единичная матрица.
Подставим $z_1 = \ldots  = z_{n-(\dim A)2k} = e$, $$x_{i,t,1}=y_{i,t,1}=\left(\pi{\Bigr|}_{A^{(t)}}\right)^{-1}\psi^{-1}(e_{11}),\ 
x_{i,t,2}=y_{i,t,2}=\left(\pi{\Bigr|}_{A^{(t)}}\right)^{-1}\psi^{-1}(e_{12}),$$
$$x_{i,t,3}=y_{i,t,3}=\left(\pi{\Bigr|}_{A^{(t)}}\right)^{-1}\psi^{-1}(e_{21}),\ 
x_{i,t,4}=y_{i,t,4}=\left(\pi{\Bigr|}_{A^{(t)}}\right)^{-1}\psi^{-1}(e_{22})$$
для всех $t\in T_1$ и $1\leqslant i \leqslant k$
и $$x_{i,t,1}=y_{i,t,1}=\left(\pi{\Bigr|}_{A^{(t)}}\right)^{-1}\psi^{-1}(\alpha_t e_{11} + \beta_t e_{12}),$$ 
$$x_{i,t,2}=y_{i,t,2}=\left(\pi{\Bigr|}_{A^{(t)}}\right)^{-1}\psi^{-1}(\alpha_t e_{21} + \beta_t e_{22})$$
для всех $t\in T_0$ и $1\leqslant i \leqslant k$.

Для того, чтобы показать, что $\mathbbm{k}^T$-многочлен $f_1$ 
не обращается в нуль при данной подстановке,
применим к результату подстановки гомоморфизм $\psi\pi$.
Значение $\mathbbm{k}^T$-многочлена $f_{t_{2\ell-1},t_{2\ell}}$ ненулевое,
так как в силу леммы~\ref{LemmaAlphaBetaDeterminant} \begin{equation*}\begin{split}[\alpha_{t_{2\ell-1}} e_{11} + \beta_{t_{2\ell-1}} e_{12},\ \alpha_{t_{2\ell-1}} e_{21} + \beta_{t_{2\ell-1}} e_{22}][\alpha_{t_{2\ell}} e_{11} + \beta_{t_{2\ell}} e_{12},\ \alpha_{t_{2\ell}} e_{21} + \beta_{t_{2\ell}} e_{22}]+\\
+[\alpha_{t_{2\ell}} e_{11} + \beta_{t_{2\ell}} e_{12},\ \alpha_{t_{2\ell}} e_{21} + \beta_{t_{2\ell}} e_{22}][\alpha_{t_{2\ell-1}} e_{11} + \beta_{t_{2\ell-1}} e_{12},\ \alpha_{t_{2\ell-1}} e_{21} + \beta_{t_{2\ell-1}} e_{22}]=\\=-\left|\begin{array}{cc}\alpha_{t_{2\ell-1}} & \beta_{t_{2\ell-1}} \\
\alpha_{t_{2\ell}} & \beta_{t_{2\ell}} \end{array}\right|^2(e_{11}+e_{22})\ne 0.\end{split} \end{equation*}
$\mathbbm{k}^T$-многочлен
$f_{t_{|T_0|-2}, t_{|T_0|-1}, t_{|T_0|}}$ не обращается в нуль при данной подстановке,
так как в силу леммы~\ref{LemmaAlphaBetaDeterminant} \begin{equation*}\begin{split}[e_{11}, e_{21}][\alpha_{t_{|T_0|}} e_{11} + \beta_{t_{|T_0|}} e_{12},
\alpha_{t_{|T_0|}}e_{21} + \beta_{t_{|T_0|}} e_{22}]
[e_{12}, e_{22}] -\\- [e_{12}, e_{22}]
[\alpha_{t_{|T_0|}} e_{11} + \beta_{t_{|T_0|}} e_{12},
\alpha_{t_{|T_0|}}e_{21} + \beta_{t_{|T_0|}} e_{22}]
[e_{11}, e_{21}] = -\alpha_{t_{|T_0|}}\beta_{t_{|T_0|}}(e_{11}+e_{22})\ne 0.\end{split} \end{equation*}
Таким образом, $f_1 \notin \Id^{\mathbbm{k}^T}(A)$.
Определим теперь $f := 
\Alt_1 \ldots \Alt_{2k} f_1$, где $\Alt_i$~"--- оператор альтернирования
по множеству
$X_i$, где $$X_{2i-1} := \lbrace x_{i,t,j} \mid t\in T,\
1\leqslant j \leqslant 2\text{ при }t\in T_0, \
1\leqslant j \leqslant 4\text{ при }t\in T_1 \rbrace,$$
а
$$ X_{2i} := \lbrace y_{i,t,j} \mid t\in T,\
1\leqslant j \leqslant 2\text{ при }t\in T_0,\
1\leqslant j \leqslant 4\text{ при }t\in T_1 \rbrace,$$ $1\leqslant i \leqslant 2k$.
Заметим, что $f$ не обращается в нуль при той же самой подстановке,
которая использовалась для $\mathbbm{k}^T$-многочлена $f_1$,
поскольку в случае, когда альтернирование заменяет $x_{i,t_1,j_1}$ на $x_{i,t_2,j_2}$ при $t_1 \ne t_2$,
значение выражения $x^{q_{t_1}}_{i,t_2,j_2}$ становится равным нулю
и соответствующее слагаемое также обращается в нуль.

Переименуем теперь переменные $\mathbbm{k}^T$-многочлена $f$ в $x_1, \ldots, x_n$.
Тогда $\mathbbm{k}^T$-многочлен удовлетворяет всем условиям леммы.
  \end{proof}

Докажем, что в случае, когда
справедливо неравенство~(\ref{EquationEquivTriangleFracM2}),
$T$-градуированная PI-экспонента алгебры $A$ равна
$\dim A$:
  
  \begin{theorem}\label{TheoremGrPIexpTriangleFracM2}
Пусть $A=\bigoplus_{t\in T} A^{(t)}$~"--- конечномерная
$T$-градуированно простая 
алгебра над полем $\mathbbm{k}$ характеристики $0$
для некоторой ленты правых нулей $T$, причём $A/J(A) \cong M_2(\mathbbm{k})$.
Пусть $T_0, T_1 \subseteq T$ и $\sim$~"--- соответственно, подмножества и отношение эквивалентности,
  заданные в начале \S\ref{SectionFracM2Equal}.  
  Предположим, что либо справедливо неравенство~(\ref{EquationEquivTriangleFracM2}),
  либо $T_0=\varnothing$.
  Тогда существуют такие $C > 0$ и $r\in \mathbb R$, что
  $$C n^r (\dim A)^n \leqslant c^{T\text{-}\mathrm{gr}}_n(A) \leqslant (\dim A)^{n+1}
  \text{ для всех }n\in \mathbb N.$$
В частности, $\PIexp^{T\text{-}\mathrm{gr}}(A) =\dim A$.
  \end{theorem}
\begin{proof}
В силу предложения~\ref{PropositionCnGrCnGenH}  для всех  $n\in \mathbb N$
справедливо равенство $c_n^{T\text{-}\mathrm{gr}}(A)=c_n^{\mathbbm{k}^T}(A)$.
Теперь достаточно применить предложение~\ref{PropositionCodimDim}, лемму~\ref{LemmaGrPIexpTriangleAltFracM2}
и теорему~\ref{TheoremAssocBounds}.
\end{proof}

\section{Случай $A/J(A)\cong M_2(\mathbbm{k})$ и $\PIexp^{T\text{-}\mathrm{gr}}(A)<\dim A$}\label{SectionFracM2Less}

Пусть $A=\bigoplus_{t\in T} A^{(t)}$~"--- конечномерная
$T$-градуированно простая 
алгебра над полем $\mathbbm{k}$ характеристики $0$
для некоторой ленты правых нулей $T$.
Пусть $T_0, T_1 \subseteq T$ и $\sim$~"--- соответственно, подмножества и отношение эквивалентности,
  заданные в начале \S\ref{SectionFracM2Equal}.

Предположим, что $T_0 \ne \varnothing$ и неравенство~(\ref{EquationEquivTriangleFracM2}) для алгебры $A$ не выполнено. Это эквивалентно существованию такого $t_0 \in T_0$, что 
$|\bar t_0| > \frac{|T_0|}{2}$. Используя теорему~\ref{TheoremSumLeftIdealsMatrix},
фиксируем такой изоморфизм $\psi \colon A/J(A) \rightarrow M_2(\mathbbm{k})$,
что $\psi(I_{t_0})=\langle e_{11}, e_{21} \rangle_\mathbbm{k}$.
В силу леммы~\ref{LemmaLeftIdealMatrix}
для любого $t\in T_0$ можно выбрать такие $\alpha_t, \beta_t\in \mathbbm{k}$,
что набор $( \alpha_t e_{i1}+\beta_t e_{i2} \mathrel{|} 1\leqslant i \leqslant 2)$ является базисом
пространства $\psi(I_t)$. Без ограничения общности можно считать, что $(\alpha_t, \beta_t)=(1,0)$ при $t\sim t_0$. Заметим, что если $I_t\ne I_{t_0}$, то $\beta_t \ne 0$.

Для каждого $t\in T_0$ выберем в компоненте $A^{(t)}$ базис
 $$\left(\left(\pi\bigr|_{A^{(t)}}\right)^{-1}\psi^{-1}(\alpha_t e_{i1}+\beta_t e_{i2}) \mathrel{\Bigl|} 1\leqslant i \leqslant 2 \right)$$, а для каждого $t\in T_1$ выберем в компоненте $A^{(t)}$ базис
 $$\left(\left(\pi\bigr|_{A^{(t)}}\right)^{-1}\psi^{-1}(e_{ij}) \mathrel{\Bigl|} 1\leqslant i,j\leqslant 2 \right).$$ 
 Определим базис $\mathcal B$ в алгебре $A$ как объединение базисов компонент $A^{(t)}$,
 выбранных выше.

Вычислим теперь числа $\gamma_i$, введённые в начале \S\ref{SectionUpperSGGr}.
Заметим, что $$\theta\left(\left(\pi\bigr|_{A^{(t)}}\right)^{-1}\psi^{-1}(\alpha_t e_{i1}+\beta_t e_{i2})\right) 
=i-2$$ при $1\leqslant i \leqslant 2$ и $t\in T_0$, $t\nsim t_0$.
Отсюда \begin{equation}\label{EquationGammaFracM2}(\gamma_1, \ldots, \gamma_r)=(\underbrace{-1,\ldots, -1}_{|T_0|+|T_1|-|\bar t_0|},
 \underbrace{0,\ldots,0}_{|T_0|+2|T_1|},
\underbrace{1,\ldots,1}_{|T_1|+|\bar t_0|}),\end{equation}
и разбиения $\lambda \vdash n$, отвечающие неприводимым $S_n$-характерам, встречающимся
в разложении $T$-градуированного кохарактера алгебры $A$, удовлетворяют неравенству из леммы~\ref{LemmaInequalityLambdaUpperFrac}.

Ниже доказываются три леммы, которые позволят выбрать элементы $b_1,\ldots, b_m$,
которые затем будут подставлены вместо переменных с номерами из заданного столбца
некоторой диаграммы Юнга. При выборе элементов $b_1,\ldots, b_m$
важно контролировать значение суммы $\sum_{j=1}^m \theta(b_i)$.

\begin{lemma}\label{LemmaMSummaGammaFracM2}
Пусть $1 \leqslant m \leqslant r$.
Тогда
$$m-\sum_{j=1}^m \gamma_j \leqslant  3|T_0|+4|T_1|-2|\bar t_0|.$$
 \end{lemma}
 \begin{proof} Если $m\leqslant |T_0|+|T_1|-|\bar t_0|$,
 то $\gamma_j=-1$ при $1\leqslant j \leqslant m$ и $\sum_{j=1}^m \gamma_j = -m$. Следовательно, $$m-\sum_{j=1}^m \gamma_j =  2m \leqslant 2|T_0|+2|T_1|-2|\bar t_0| \leqslant  3|T_0|+4|T_1|-2|\bar t_0|.$$
 
 Если $|T_0|+|T_1|-|\bar t_0| \leqslant m \leqslant 2|T_0|+3|T_1|-|\bar t_0|$,
 то $\gamma_j= 0$ при $|T_0|+|T_1|-|\bar t_0| < j \leqslant m$ и
 $\sum_{j=1}^m \gamma_j = -(|T_0|+|T_1|-|\bar t_0|)$.
 Следовательно, $$m-\sum_{j=1}^m \gamma_j =  m+(|T_0|+|T_1|-|\bar t_0|)\leqslant 3|T_0|+4|T_1|-2|\bar t_0|.$$
 
 Если $m \geqslant 2|T_0|+3|T_1|-|\bar t_0|$,
 то из равенства $\gamma_j=1$ при $j> 2|T_0|+3|T_1|-|\bar t_0|$ следует, что $$m-\sum_{j=1}^m \gamma_j=3|T_0|+4|T_1|-2|\bar t_0|.$$
 \end{proof}

\begin{lemma}\label{LemmaChooseForAColumnPositiveFracM2}
Пусть $\sum_{j=1}^m \gamma_j > 0$ для некоторого $m\in\mathbb N$.
Тогда существуют такие $b_1, \ldots, b_m \in \mathcal B$,  где $b_i\ne b_j$ при $i\ne j$,
что $\sum_{j=1}^m \theta(b_j) = \sum_{j=1}^m \gamma_j$
и 
\begin{itemize}
\item если $\lbrace b_1, \ldots, b_m \rbrace \cap A^{(t)}=\lbrace b_i \rbrace$
для некоторых $1\leqslant i \leqslant m$ и $t\in T_0 \sqcup T_1$,
то $t\in T_0$, $t\sim t_0$ и $b_i=\left(\pi\bigr|_{A^{(t)}}\right)^{-1}\psi^{-1}(e_{11})$;
\item если $\lbrace b_1, \ldots, b_m \rbrace \cap A^{(t)}=\lbrace b_i, b_j \rbrace$
для некоторых $1\leqslant i,j \leqslant m$ и $t\in T_0 \sqcup T_1$,
то либо $\theta(b_i)\ne 0$, либо $\theta(b_j)\ne 0$.
\end{itemize}
\end{lemma}
\begin{proof}
В силу определения чисел $\gamma_i$ (см. начало \S\ref{SectionUpperSGGr}),
существуют такие $b_1, \ldots, b_m \in \mathcal B$, где $b_i\ne b_j$ при $i\ne j$,
что $\sum_{j=1}^m \theta(b_j) = \sum_{j=1}^m \gamma_j$
и $\sum_{j=1}^m \theta(a_j) \geqslant \sum_{j=1}^m \gamma_j$
для всех $a_1, \ldots, a_m \in \mathcal B$, где $a_i\ne a_j$ при $i\ne j$.
Поскольку $\sum_{j=1}^m \gamma_j > 0$, из минимальности $\sum_{j=1}^m \theta(b_j)$
следует, что множество $\lbrace b_1, \ldots, b_m \rbrace$
содержит все элементы $b\in \mathcal B$ с $\theta(b) \leqslant 0$.
Теперь утверждение леммы следует из выбора множества $\mathcal B$.
\end{proof}

\begin{lemma}\label{LemmaChooseForAColumnNegativeFracM2}
Пусть $\sum_{j=1}^m \gamma_j \leqslant q \leqslant 0$ для некоторых $m\in\mathbb N$, $q\in\mathbb Z$.
Тогда существуют такие $b_1, \ldots, b_m \in \mathcal B$, где $b_i\ne b_j$ при $i\ne j$,
что $\sum_{j=1}^m \theta(b_j) = q$
и \begin{itemize}
\item если $\lbrace b_1, \ldots, b_m \rbrace \cap A^{(t)}=\lbrace b_i \rbrace$,
$\theta(b_i)=0$
для некоторых $1\leqslant i \leqslant m$ и $t\in T_0 \sqcup T_1$,
то $t\in T_0$, $t\sim t_0$ и $b_i=\left(\pi\bigr|_{A^{(t)}}\right)^{-1}\psi^{-1}(e_{11})$;
\item если $\lbrace b_1, \ldots, b_m \rbrace \cap A^{(t)}=\lbrace b_i, b_j \rbrace$
для некоторых $1\leqslant i,j \leqslant m$ и $t\in T_0 \sqcup T_1$,
то либо $\theta(b_i)\ne 0$, либо $\theta(b_j)\ne 0$.
\end{itemize}
\end{lemma}
\begin{proof}
Напомним, что $$|\lbrace b\in\mathcal B \mid \theta(b)=-1  \rbrace|=|T_0|+|T_1|-|\bar t_0|,$$
$$|\lbrace b\in\mathcal B \mid \theta(b)=0  \rbrace|=|T_0|+2|T_1|,$$
$$|\lbrace b\in\mathcal B \mid \theta(b)=1  \rbrace|=|T_1|+|\bar t_0|.$$
Пусть $\lbrace t_1, \ldots, t_{|T_1|} \rbrace := T_1$,
$\lbrace \tilde t_1, \ldots, \tilde t_{|T_0|-|\bar t_0|} \rbrace := T_0 \backslash \bar t_0$,
$\lbrace \hat t_1, \ldots, \hat t_{|\bar t_0|} \rbrace := \bar t_0$.

Выделим следующие два случая:

1. Пусть $m < 2(|T_0|+|T_1|-|\bar t_0|)+q$. Положим $\ell = \left[
\frac{m+q}{2}\right]$. Заметим, что тогда $$\ell - q \leqslant \frac{m-q}{2} < |T_0|+|T_1|-|\bar t_0|$$
и $\ell \leqslant \ell - q < |T_0|+|T_1|-|\bar t_0| < |T_1|+|\bar t_0|$.
Из этих неравенств будет следовать, что ниже у нас окажется достаточно
элементов из $(T_0\backslash \bar t_0) \sqcup T_1$
и $\bar t_0 \sqcup T_1$, соответственно.

 Предположим сперва, что $m=2\ell-q$. 
Если $\ell \leqslant |T_1|+q$, то возьмём
\begin{equation*}\begin{split}\lbrace b_1, \ldots, b_m \rbrace = 
\left\lbrace\left(\pi\bigr|_{A^{(t_i)}}\right)^{-1}\psi^{-1}(e_{12})
\mathrel{\bigl|} 1\leqslant i \leqslant \ell-q
\right\rbrace  \cup \\ \cup \left\lbrace\left(\pi\bigr|_{A^{(t_i)}}\right)^{-1}\psi^{-1}(e_{21})
\mathrel{\bigl|} 1\leqslant i \leqslant \ell
\right\rbrace.\end{split}\end{equation*}
Если же $|T_1|+q < \ell \leqslant |T_1|$, то возьмём
\begin{equation*}\begin{split}\lbrace b_1, \ldots, b_m \rbrace = 
\left\lbrace\left(\pi\bigr|_{A^{(t_i)}}\right)^{-1}\psi^{-1}(e_{12})
\mathrel{\bigl|} 1\leqslant i \leqslant |T_1|
\right\rbrace \cup \\ \cup \left\lbrace\left(\pi\bigr|_{A^{(\tilde t_i)}}\right)^{-1}\psi^{-1}(\alpha_{\tilde t_i} e_{11}+\beta_{\tilde t_i} e_{12})
\mathrel{\bigl|} 1\leqslant i \leqslant \ell-q-|T_1|
\right\rbrace \cup \\ \cup \left\lbrace\left(\pi\bigr|_{A^{(t_i)}}\right)^{-1}\psi^{-1}(e_{21})
\mathrel{\bigl|} 1\leqslant i \leqslant \ell
\right\rbrace.\end{split}\end{equation*}
При $\ell > |T_1|$ возьмём
\begin{equation*}\begin{split}\lbrace b_1, \ldots, b_m \rbrace = 
\left\lbrace\left(\pi\bigr|_{A^{(t_i)}}\right)^{-1}\psi^{-1}(e_{12})
\mathrel{\bigl|} 1\leqslant i \leqslant |T_1|
\right\rbrace \cup \\ \cup \left\lbrace\left(\pi\bigr|_{A^{(\tilde t_i)}}\right)^{-1}\psi^{-1}(\alpha_{\tilde t_i} e_{11}+\beta_{\tilde t_i} e_{12})
\mathrel{\bigl|} 1\leqslant i \leqslant \ell-q-|T_1|
\right\rbrace \cup \\ \cup \left\lbrace\left(\pi\bigr|_{A^{(t_i)}}\right)^{-1}\psi^{-1}(e_{21})
\mathrel{\bigl|} 1\leqslant i \leqslant |T_1|
\right\rbrace  \cup \\ \cup \left\lbrace\left(\pi\bigr|_{A^{(\hat t_i)}}\right)^{-1}\psi^{-1}(e_{21})
\mathrel{\bigl|} 1\leqslant i \leqslant \ell-|T_1|
\right\rbrace.\end{split}\end{equation*}
При $m=2\ell-q+1$ в каждом из трёх случаев, рассмотренных выше,
добавим элемент $\left(\pi\bigr|_{A^{(\hat t_1)}}\right)^{-1}\psi^{-1}(e_{11})$.

2. Пусть $m \geqslant 2(|T_0|+|T_1|-|\bar t_0|)+q$.

В силу леммы~\ref{LemmaMSummaGammaFracM2}
$$m-\sum_{j=1}^m \gamma_j \leqslant  3|T_0|+4|T_1|-2|\bar t_0|.$$
Следовательно,
$$m-q - 2(|T_0|+|T_1|-|\bar t_0|) \leqslant  |T_0|+2|T_1|
$$ 
и можно выбрать такие
$$0\leqslant k, \ell \leqslant |T_1|,\qquad 0\leqslant s \leqslant |\bar t_0|,\qquad 0\leqslant u \leqslant |T_0|-|\bar t_0|,$$ что $2(|T_0|+|T_1|-|\bar t_0|)+q + k + \ell + s + u = m$.

Заметим, что $$(|T_0|+|T_1|-|\bar t_0|)+q \geqslant (|T_0|+|T_1|-|\bar t_0|) + \sum_{j=1}^m \gamma_j
\geqslant 0.$$
Если $(|T_0|+|T_1|-|\bar t_0|)+q \leqslant |T_1|$, определим
\begin{equation*}\begin{split}\lbrace b_1, \ldots, b_m \rbrace = 
\left\lbrace\left(\pi\bigr|_{A^{(t_i)}}\right)^{-1}\psi^{-1}(e_{12})
\mathrel{\bigl|} 1\leqslant i \leqslant |T_1|
\right\rbrace \cup \\ \cup \left\lbrace\left(\pi\bigr|_{A^{(t_i)}}\right)^{-1}\psi^{-1}(e_{11})
\mathrel{\bigl|} 1\leqslant i \leqslant k
\right\rbrace \cup \\ \cup \left\lbrace\left(\pi\bigr|_{A^{(t_i)}}\right)^{-1}\psi^{-1}(e_{21})
\mathrel{\bigl|} 1\leqslant i \leqslant (|T_0|+|T_1|-|\bar t_0|)+q
\right\rbrace  \cup \\ \cup \left\lbrace\left(\pi\bigr|_{A^{(t_i)}}\right)^{-1}\psi^{-1}(e_{22})
\mathrel{\bigl|} 1\leqslant i \leqslant \ell
\right\rbrace \cup \\ \cup \left\lbrace\left(\pi\bigr|_{A^{(\tilde t_i)}}\right)^{-1}\psi^{-1}(\alpha_{\tilde t_i} e_{11}+\beta_{\tilde t_i} e_{12})
\mathrel{\bigl|} 1\leqslant i \leqslant |T_0|-|\bar t_0|
\right\rbrace  \cup \\ \cup \left\lbrace\left(\pi\bigr|_{A^{(\tilde t_i)}}\right)^{-1}\psi^{-1}(\alpha_{\tilde t_i} e_{21}+\beta_{\tilde t_i} e_{22})
\mathrel{\bigl|} 1\leqslant i \leqslant u
\right\rbrace  \cup \\ \cup \left\lbrace\left(\pi\bigr|_{A^{(\hat t_i)}}\right)^{-1}\psi^{-1}(e_{11})
\mathrel{\bigl|} 1\leqslant i \leqslant s
\right\rbrace.\end{split}\end{equation*}
При $(|T_0|+|T_1|-|\bar t_0|)+q > |T_1|$ определим
\begin{equation*}\begin{split}\lbrace b_1, \ldots, b_m \rbrace = 
\left\lbrace\left(\pi\bigr|_{A^{(t_i)}}\right)^{-1}\psi^{-1}(e_{12})
\mathrel{\bigl|} 1\leqslant i \leqslant |T_1|
\right\rbrace   \cup \\ \cup \left\lbrace\left(\pi\bigr|_{A^{(t_i)}}\right)^{-1}\psi^{-1}(e_{11})
\mathrel{\bigl|} 1\leqslant i \leqslant k
\right\rbrace   \cup \\ \cup \left\lbrace\left(\pi\bigr|_{A^{(t_i)}}\right)^{-1}\psi^{-1}(e_{21})
\mathrel{\bigl|} 1\leqslant i \leqslant |T_1|
\right\rbrace   \cup \\ \cup \left\lbrace\left(\pi\bigr|_{A^{(t_i)}}\right)^{-1}\psi^{-1}(e_{22})
\mathrel{\bigl|} 1\leqslant i \leqslant \ell
\right\rbrace  \cup \\ \cup \left\lbrace\left(\pi\bigr|_{A^{(\tilde t_i)}}\right)^{-1}\psi^{-1}(\alpha_{\tilde t_i} e_{11}+\beta_{\tilde t_i} e_{12})
\mathrel{\bigl|} 1\leqslant i \leqslant |T_0|-|\bar t_0|
\right\rbrace  \cup \\ \cup \left\lbrace\left(\pi\bigr|_{A^{(\tilde t_i)}}\right)^{-1}\psi^{-1}(\alpha_{\tilde t_i} e_{21}+\beta_{\tilde t_i} e_{22})
\mathrel{\bigl|} 1\leqslant i \leqslant u
\right\rbrace  \cup \\ \cup \left\lbrace\left(\pi\bigr|_{A^{(\hat t_i)}}\right)^{-1}\psi^{-1}(e_{11})
\mathrel{\bigl|} 1\leqslant i \leqslant s
\right\rbrace  \cup \\ \cup \left\lbrace\left(\pi\bigr|_{A^{(\hat t_i)}}\right)^{-1}\psi^{-1}(e_{21})
\mathrel{\bigl|} 1\leqslant i \leqslant (|T_0|-|\bar t_0|)+q
\right\rbrace.\end{split}
\end{equation*}
\end{proof}

В лемме ниже доказывается существование $\mathbbm{k}^T$-многочлена, который не является $\mathbbm{k}^T$-тождеством, но в
то же время кососимметричен по достаточному числу наборов переменных.
Этот $\mathbbm{k}^T$-многочлен и будет порождать $\mathbbm{k}S_n$-подмодуль с размерностью, достаточной
для доказательства оценки снизу.

\begin{lemma}\label{LemmaAltNonTriangleFracM2}
Пусть $\lambda=(\lambda_1, \ldots, \lambda_r) \vdash n$ для некоторого $n\in\mathbb N$. Если $\sum_{i=1}^r \gamma_i \lambda_i \leqslant 0$ и $2 \mid \lambda_i$ при $i \geqslant 2$, 
то существует диаграмма Юнга $T_\lambda$ формы $\lambda$ и $\mathbbm{k}^T$-многочлен $f\in P_n^{\mathbbm{k}^T}$,
такой, что  $b_{T_\lambda}f =\sum_{\sigma\in C_{T_\lambda}}(\sign \sigma)\sigma f\notin \Id^{\mathbbm{k}^T}(A)$.
\end{lemma}
\begin{proof}
Обозначим через $\mu_i$ 
число клеток в $i$-м столбце таблицы $T_\lambda$.
Пусть $m_i := \sum_{j=1}^{\mu_i} \gamma_j$ для $1\leqslant i \leqslant \lambda_1$.

Заметим, что в силу того, что число $(\lambda_i-\lambda_{i+1})$
равно числу столбцов высоты $i$ (полагаем $\lambda_{r+1}:=0$), а $\lambda_1$ 
равно числу всех столбцов, неравенство $\sum_{i=1}^r \gamma_i \lambda_i \leqslant 0$
может быть переписано как
\begin{equation}\label{EquationGammaLambdaFracM2}\sum_{i=1}^r \left(\sum_{j=1}^i \gamma_j -
\sum_{j=1}^{i-1} \gamma_j\right) \lambda_i=\sum_{i=1}^r \left(\sum_{j=1}^i \gamma_j\right) (\lambda_i-\lambda_{i+1}) = \sum_{i=1}^{\lambda_1} m_i \leqslant 0. \end{equation} 
В силу леммы~\ref{LemmaThetaCondition},
если сумма значений функции $\theta$ на элементах базиса,
подставленных вместо переменных некоторого полилинейного $\mathbbm{k}^T$-многочлена
больше, чем $1$,
такой $\mathbbm{k}^T$-многочлен обращается в нуль.
Ниже будут выбраны такие элементы $b_{itj}\in A^{(t)} \cap \mathcal B$,
что сумма значений на них функции $\theta$ равна $0$. В случае, когда $m_i > 0$
для некоторого $i$, нам потребуется сделать сумму значений функции
 $\theta$ на элементах, подставляемых в некоторые другие столбцы,
 отрицательной.

Определим число $1\leqslant \ell \leqslant \lambda_1$ при помощи условий $m_\ell > 0$ и $m_{\ell+1} \leqslant 0$.
Поскольку $2 \mid \lambda_i$ для всех  $i \geqslant 2$, для всех $1 \leqslant j \leqslant \frac{\lambda_2}2$ выполняется равенство $m_{2j-1}=m_{2j}$.
Напомним, что $m_i=-1$ при $\lambda_2+1\leqslant i \leqslant \lambda_1$.
 Отсюда $\ell \leqslant \lambda_2$, $2 \mid \ell$ и $2 \mid \sum_{i=1}^{\ell} m_i$.
 В силу~(\ref{EquationGammaLambdaFracM2})
$$2\sum_{i=1}^{\lambda_2/2} m_{2i} -(\lambda_1-\lambda_2)=\sum_{i=1}^{\lambda_1} m_i \leqslant 0$$
и
$$\sum_{i=1}^{\lambda_2/2} m_{2i} - \left[\frac{\lambda_1-\lambda_2}2\right] \leqslant 0.$$

Отсюда можно выбрать такие числа $N$ и $k_{2i}$, где $\ell < 2i \leqslant \lambda_2$,
что $0\leqslant N \leqslant \left[\frac{\lambda_1-\lambda_2}2\right]$, $m_{2i}\leqslant k_{2i} \leqslant 0$
и $\sum_{i=1}^{\ell/2} m_{2i}+\sum_{i=\ell/2+1}^{\lambda_2} k_{2i} - N = 0$.
Пусть $k_{2i-1}:=k_{2i}$ при $\frac{\ell}2 < i \leqslant \frac{\lambda_2}2$,
$k_i:=-1$ при $\lambda_2+1 \leqslant i \leqslant \lambda_2+2N$,
$k_i:=0$ при $\lambda_2+2N+1 \leqslant i \leqslant \lambda_1$.
Тогда \begin{equation}\label{EquationSumMiQiZeroFracM2}
\sum_{i=1}^{\ell} m_i+\sum_{i=\ell+1}^{\lambda_1} k_i=0.
\end{equation}

Теперь для всякого $1\leqslant i \leqslant \ell$ воспользуемся леммой~\ref{LemmaChooseForAColumnPositiveFracM2}, положив $m=\mu_i$,
а для всякого $\ell+1\leqslant i \leqslant \lambda_1$
воспользуемся  леммой~\ref{LemmaChooseForAColumnNegativeFracM2},
положив $m=\mu_i$ и $q=k_i$.
Отсортируем получившиеся элементы $b_j$
в соответствии с теми однородными компонентами $A^{(t)}$,
которым эти элементы принадлежат.
Для всякого $1\leqslant i \leqslant \lambda_1$
получим элементы $b_{itj}\in A^{(t)} \cap \mathcal B$, где $t\in T_0 \sqcup T_1$, $1\leqslant j \leqslant n_{it}$, общее число элементов $b_{itj}$ при фиксированных $i$ и $t$ равно $n_{it}\geqslant 0$, $$\sum_{t\in T_0 \sqcup T_1} n_{it}= \mu_i,$$
 $$b_{it_1j_1} \ne b_{it_2j_2} \text{ при } (t_1,j_1)\ne (t_2,j_2),$$
$$\sum_{t\in T_0 \sqcup T_1} \sum_{j=1}^{n_{it}} \theta(b_{itj}) =m_i \text{ при } m_i > 0$$ и
$$m_i \leqslant \sum_{t\in T_0 \sqcup T_1} \sum_{j=1}^{n_{it}} \theta(b_{itj})=k_i \leqslant 0 \text{ при } m_i \leqslant 0.$$

В силу~(\ref{EquationSumMiQiZeroFracM2})
справедливо равенство
\begin{equation}\label{EquationSumThetaBitj0FracM2}
\sum_{i=1}^{\lambda_1} 
\sum_{t\in T_0 \sqcup T_1} \sum_{j=1}^{n_{it}} \theta(b_{itj}) = 0.
\end{equation}

Поскольку для всех $2i \leqslant \lambda_2$ справедливы равенства $k_{2i-1}=k_{2i}$ и $\mu_{2i-1}=\mu_{2i}$,
можно считать, что $n_{{2i-1},t}=n_{{2i},t}$,
$b_{{2i-1},t,j}=b_{{2i},t,j}$
для всех $1 < 2i \leqslant \lambda_2$, $t\in T_0 \sqcup T_1$, $1\leqslant j \leqslant
n_{{2i},t}$.

Элементы $b_{itj}$
будут впоследствии подставляться вместо тех переменных, номера которых
находятся в $i$-м столбце таблицы $T_\lambda$. 

Обозначим через $W_{-1}$ 
множество всех таких пар $(i, t)$, где $1\leqslant i \leqslant \lambda_1$, $t\in T_0 \sqcup T_1$, что 
$\theta(b_{itj})\leqslant 0$ для всех  $1\leqslant j \leqslant n_{ti}$,
причём для некоторого $1\leqslant j \leqslant n_{ti}$ справедливо равенство $\theta(b_{itj})=-1$.
Через $W_1$ обозначим  множество всех таких пар $(i, t)$, где $1\leqslant i \leqslant \lambda_1$, $t\in T_0 \sqcup T_1$, что $\theta(b_{itj})\geqslant 0$ для всех  $1\leqslant j \leqslant n_{ti}$,
причём для некоторого $1\leqslant j \leqslant n_{ti}$ справедливо равенство 
$\theta(b_{itj})=1$.
Для каждого $1\leqslant i \leqslant \lambda_1$
определим $W_1^{(i)} := \lbrace t\in T_0 \sqcup T_1 \mid
(i,t)\in W_1 \rbrace$ и
и $W_0^{(i)} := \lbrace t\in T_0 \sqcup T_1 \mid
(i,t)\notin W_{-1} \sqcup W_1, n_{it} > 0 \rbrace$.

В силу~(\ref{EquationSumThetaBitj0FracM2}) справедливо равенство $|W_{-1}| = |W_1|$.
Следовательно, существуют такие отображения $\varkappa \colon W_1 \rightarrow \lbrace 1,\ldots, \lambda_1\rbrace$ и $\rho \colon W_1 \rightarrow T_0 \sqcup T_1$, что $(i,t) \mapsto (\varkappa(i,t),\rho(i,t))$ является биекцией $W_1 \rightarrow W_{-1}$.

  Определим теперь многочлены $f_{it}$ и $\tilde f_{it}$, где $1\leqslant i \leqslant \lambda_2$, а $t\in T_0 \sqcup T_1$, следующим образом.

При $n_{it}=1$ положим $f_{it}(x_1):=x_1$.

При $n_{it}=2$ положим $f_{it}(x_1, x_2):=x_1 x_2 - x_2 x_1$.

При $n_{it}=3$ положим $f_{it}(x_1,x_2, x_3):=\sum_{\sigma \in S_3} \sign(\sigma) x_{\sigma(1)}x_{\sigma(2)}x_{\sigma(3)}$.

При $1\leqslant n_{it}\leqslant 3$ положим $$\tilde f_{it}(x_1,\ldots, x_{n_{it}}; y_1,\ldots, y_{n_{it}})
:=f_{it}(x_1,\ldots, x_{n_{it}})f_{it}(y_1,\ldots, y_{n_{it}}).$$

При $n_{it}=4$ положим $$\tilde f_{it}(x_1,\ldots, x_4, y_1, \ldots, y_4):=\sum_{\sigma,\tau \in S_4} \sign(\sigma\tau) x_{\sigma(1)}\ y_{\tau(1)}\ x_{\sigma(2)}x_{\sigma(3)}x_{\sigma(4)}\ y_{\tau(2)}y_{\tau(3)}y_{\tau(4)}.$$
Напомним, что такой многочлен не является для алгебры $M_2(\mathbbm{k})$
полиномиальным тождеством, а его значения пропорциональны единичной матрице.
  (См., например, \cite[теорема 5.7.4]{ZaiGia}.)

Пусть $X_i:=\lbrace x_{itj} \mid 1\leqslant j \leqslant n_{it},\ t\in T_0 \sqcup T_1 \rbrace$, $1\leqslant i \leqslant \lambda_1$.
Обозначим через $\Alt_i$ оператор альтернирования по переменным из множества $X_i$.

Положим \begin{equation*}\begin{split}f := \Alt_1 \Alt_2 \ldots \Alt_{\lambda_1}
\prod_{i=1}^{\lambda_2/2} \left(\prod_{t\in W_0^{(2i-1)}}
 \tilde f_{2i-1,t}(x_{2i-1,t,1}^{q_t},\ldots, x_{2i-1,t,n_{2i-1,t}}^{q_t};\ x_{2i,t,1}^{q_t},\ldots, x_{2i,t,n_{2i,t}}^{q_t}) \right. \cdot \\  \cdot \left.
 \left(\prod_{t\in W_1^{(2i-1)}}
 f_{\varkappa(2i-1,t)\rho(2i-1,t)}\left(x_{\varkappa(2i-1,t)\rho(2i-1,t)1}^{q_{\rho(2i-1,t)}},\ldots, x_{\varkappa(2i-1,t)\rho(2i-1,t)n_{\varkappa(2i-1,t)\rho(2i-1,t)}}^{q_{\rho(2i-1,t)}}\right)
\right.\right. \cdot \\ \cdot \left.\left.
 f_{2i-1,t}(x_{2i-1,t,1}^{q_t},\ldots, x_{2i-1,t,n_{2i-1,t}}^{q_t})
  \right.\right. \cdot \\  \cdot 
  f_{\varkappa(2i,t)\rho(2i,t)}\left(x_{\varkappa(2i,t)\rho(2i,t)1}^{q_{\rho(2i,t)}},\ldots, x_{\varkappa(2i,t)\rho(2i,t)n_{\varkappa(2i,t)\rho(2i,t)}}^{q_{\rho(2i,t)}}\right)
   \cdot \\  \cdot  \left.\left.
 f_{2i,t}(x_{2i,t,1}^{q_t},\ldots, x_{2i,t,n_{2i,t}}^{q_t}) 
 \right) \right)
 \prod_{\substack{i=\lambda_2+1, \\ t\in W^{(i)}_0}}^{\lambda_1} x_{it1}^{q_t}.\end{split}\end{equation*}
 
Заметим, что в силу лемм~\ref{LemmaChooseForAColumnPositiveFracM2} и \ref{LemmaChooseForAColumnNegativeFracM2} в случае, когда $(i,t)\in W_{-1}$, множество $\lbrace\psi\pi(b_{it1}),\ldots, \psi\pi(b_{itn_{it}})\rbrace$
совпадает с одним из следующих множеств: $\lbrace e_{12}\rbrace$, $\lbrace \alpha_t e_{11}+\beta_t e_{12} \rbrace$, $\lbrace e_{12}, e_{22}\rbrace$, $\lbrace e_{11}, e_{12}\rbrace$,
$\lbrace \alpha_t e_{11}+\beta_t e_{12}, \alpha_t e_{21}+\beta_t e_{22} \rbrace$,
 $\lbrace e_{11}, e_{12}, e_{22}\rbrace$.

Если $t\in W_0^{(i)}$, то $\lbrace\psi\pi(b_{it1}),\ldots, \psi\pi(b_{itn_{it}})\rbrace$
совпадает с одним из следующих множеств: $\lbrace e_{11}\rbrace$, $\lbrace e_{12}, e_{21}\rbrace$,
$\lbrace e_{12}, e_{22}, e_{21}\rbrace$, $\lbrace e_{11}, e_{12}, e_{21}\rbrace$,
$\lbrace e_{11}, e_{12}, e_{22}, e_{21}\rbrace$.

Если $t\in W_1^{(i)}$, то $\lbrace\psi\pi(b_{it1}),\ldots, \psi\pi(b_{itn_{it}})\rbrace$
совпадает с одним из следующих множеств: $\lbrace e_{21}\rbrace$, $\lbrace e_{22}, e_{21}\rbrace$,
 $\lbrace e_{21}, e_{11}\rbrace$,
$\lbrace e_{22}, e_{21}, e_{11}\rbrace$.

В силу леммы~\ref{LemmaAlphaBetaDeterminant} и всех замечаний, сделанных выше,
значение $\mathbbm{k}^T$-многочлена $f$ при подстановке $x_{itj}=b_{itj}$, где $1\leqslant i \leqslant \lambda_1$,
$t\in T_0 \sqcup T_1$, $1\leqslant j \leqslant n_{it}$,
не обращается в нуль под действием гомоморфизма $\psi\pi$, поскольку
если для каких-то $t_1 \ne t_2$ альтернирование заменяет $x_{i,t_1,j_1}$ на $x_{i,t_2,j_2}$,
то значение выражения $x^{q_{t_1}}_{i,t_2,j_2}$ 
нулевое и всё соответствующее слагаемое обращается в нуль.
Переименуем теперь для удобства переменные $\mathbbm{k}^T$-многочлена $f$ в $x_1, \ldots, x_n$.
Тогда $f$ удовлетворяет всем требования данной леммы.
\end{proof}

Теперь мы готовы вычислить $\PIexp^{T\text{-}\mathrm{gr}}(A)$ 
в случае, когда неравенство~(\ref{EquationEquivTriangleFracM2}) не выполнено:

\begin{theorem}\label{TheoremGrPIexpNonTriangleFracM2}
Пусть $A$~"--- конечномерная
$T$-градуированно простая 
алгебра над полем $\mathbbm{k}$ характеристики $0$
для некоторой ленты правых нулей $T$, причём $A/J(A) \cong M_2(\mathbbm{k})$.
Пусть $T_0, T_1 \subseteq T$ и $\sim$~"--- соответственно, подмножества и отношение эквивалентности,
  заданные в начале \S\ref{SectionFracM2Equal}.  
  Предположим, что $|\bar t_0| > \frac{|T_0|}{2}$ для некоторого $\bar t_0 \in T_0/\sim$.
  Тогда $$\PIexp^{T\text{-}\mathrm{gr}}(A) = |T_0|+2|T_1|
  + 2\sqrt{(|T_1|+|\bar t_0|)(|T_0|+|T_1|-|\bar t_0|)}<2|T_0|+4|T_1|=\dim A.$$
\end{theorem}
\begin{proof}
Пусть $\mathcal B$~"--- базис алгебры $A$,
выбранный в начале \S\ref{SectionFracM2Less}
Из равенства (\ref{EquationGammaFracM2}) следует, что
  $0<\zeta=\sqrt{\frac{|T_0|+|T_1|-|\bar t_0|}{|T_1|+|\bar t_0|}}<1$ является корнем уравнения~(\ref{EquationZetaUpperFrac}).

Пусть $$\Omega = \left\lbrace (\alpha_1, \ldots, \alpha_r)\in \mathbb R^r \mathrel{\biggl|} \sum_{i=1}^r \alpha_i = 1,\ 
\alpha_1 \geqslant \alpha_2 \geqslant \ldots \geqslant \alpha_r\geqslant 0,\ \sum_{i=1}^r \gamma_i \alpha_i \leqslant 0\right\rbrace.$$
В силу леммы~\ref{LemmaMaximumFUpperFrac}
\begin{equation}\label{EquationMaximumGammaOmegaFracM2} d:=\max_{x\in \Omega} \Phi(x) =\sum_{i=1}^r \zeta^{\gamma_i}= |T_0|+2|T_1|
  + 2\sqrt{(|T_1|+|\bar t_0|)(|T_0|+|T_1|-|\bar t_0|)}.\end{equation}
Обозначим через $(\alpha_1, \ldots, \alpha_r) \in \Omega$ такую точку, что
$\Phi(\alpha_1, \ldots, \alpha_r)=d$.

Для любого $n\in\mathbb N$ определим $\mu\vdash n$ при помощи равенств
$\mu_i := 2\left[\frac{\alpha_i n}{2}\right]$ при $2\leqslant i \leqslant r$
и $\mu_1 :=  n-\sum_{i=2}^r \mu_i$.

В силу~(\ref{EquationAlphaRelationUpperFrac1}) справедливо равенство $\sum_{i=1}^r \gamma_i \alpha_i = 0$.
Поскольку $$\gamma_1 = \ldots = \gamma_{|T_0|+|T_1|-|\bar t_0|}=-1,$$
$$\gamma_{|T_0|+|T_1|-|\bar t_0|+1} = \ldots = \gamma_{2|T_0|+3|T_1|-|\bar t_0|}=0,$$
и $$\gamma_{2|T_0|+3|T_1|-|\bar t_0|+1} = \ldots = \gamma_r=1,$$
из~(\ref{EquationAlphaRelationUpperFrac2}) следует, что
$$\alpha_1 = \ldots = \alpha_{|T_0|+|T_1|-|\bar t_0|},$$
$$\alpha_{|T_0|+|T_1|-|\bar t_0|+1} = \ldots = \alpha_{2|T_0|+3|T_1|-|\bar t_0|},$$
$$\alpha_{2|T_0|+3|T_1|-|\bar t_0|+1} = \ldots = \alpha_r,$$
откуда
$$\alpha_1 n-2 \leqslant \mu_2 = \ldots=\mu_{|T_0|+|T_1|-|\bar t_0|} \leqslant \alpha_1 n,$$
$$\alpha_{|T_0|+|T_1|-|\bar t_0|+1} n-2 \leqslant \mu_{|T_0|+|T_1|-|\bar t_0|+1} = \ldots=\mu_{2|T_0|+3|T_1|-|\bar t_0|} \leqslant \alpha_{|T_0|+|T_1|-|\bar t_0|+1} n,$$
$$\alpha_r n-2 \leqslant \mu_{2|T_0|+3|T_1|-|\bar t_0|+1} = \ldots=\mu_r \leqslant \alpha_r n.$$
Теперь из $\sum_{i=1}^r \alpha_i = 1$ следует, что
$\alpha_1 n \leqslant \mu_1 \leqslant \alpha_1 n+2r$.

Заметим, что
\begin{equation}
\label{EquationSumMuGammaFracM2}\begin{split} \sum_{i=1}^r \gamma_i \mu_i=
-\left(n-\sum_{i=2}^r \mu_i\right)-\sum_{i=2}^{|T_0|+|T_1|-|\bar t_0|} \mu_i +  \sum_{i=2|T_0|+3|T_1|-|\bar t_0|+1}^r \mu_i =\\= \left(\sum_{i=|T_0|+|T_1|-|\bar t_0|+1}^{2|T_0|+3|T_1|-|\bar t_0|} \mu_i\right)+ 2\left(\sum_{i=2|T_0|+3|T_1|-|\bar t_0|+1}^r \mu_i\right) - n
\leqslant \\
\leqslant n\left(\sum_{i=|T_0|+|T_1|-|\bar t_0|+1}^{2|T_0|+3|T_1|-|\bar t_0|} \alpha_i\right)
+2n \left(\sum_{i=2|T_0|+3|T_1|-|\bar t_0|+1}^r \alpha_i \right)
 -n
\sum_{i=1}^r \alpha_i = n \sum_{i=1}^r \gamma_i \alpha_i = 0.\end{split}\end{equation}
Аналогичным образом получаем, что
\begin{equation}
\label{EquationSumMuGamma2FracM2}\begin{split} \sum_{i=1}^r \gamma_i \mu_i=
 \left(\sum_{i=|T_0|+|T_1|-|\bar t_0|+1}^{2|T_0|+3|T_1|-|\bar t_0|} \mu_i\right)+ 2\left(\sum_{i=2|T_0|+3|T_1|-|\bar t_0|+1}^r \mu_i\right) - n
 \geqslant \\\\
\geqslant n\left(\sum_{i=|T_0|+|T_1|-|\bar t_0|+1}^{2|T_0|+3|T_1|-|\bar t_0|} \alpha_i\right)
+2n \left(\sum_{i=2|T_0|+3|T_1|-|\bar t_0|+1}^r \alpha_i \right)
 -n
\sum_{i=1}^r \alpha_i-4r = \\\\= n \sum_{i=1}^r \gamma_i \alpha_i -4r = -4r.\end{split}\end{equation}

 Из~(\ref{EquationSumMuGammaFracM2}) и леммы~\ref{LemmaAltNonTriangleFracM2} следует, что $b_{T_\mu}f \notin \Id^{\mathbbm{k}^T}(A)$
 для некоторого $f\in P^{\mathbbm{k}^T}_n$. Обозначим через $\bar f$ образ $\mathbbm{k}^T$-многочлена $f$ в факторпространстве $\frac{P^{\mathbbm{k}^T}_n}{P^{\mathbbm{k}^T}_n
 \cap\ \Id^{\mathbbm{k}^T}(A)}$. Рассмотрим подмодуль $\mathbbm{k}S_n b_{T_\mu}\bar f \subseteq \frac{P^{\mathbbm{k}^T}_n}{P^{\mathbbm{k}^T}_n
 \cap\ \Id^{\mathbbm{k}^T}(A)}$. Поскольку все $S_n$-представления над полями характеристики  $0$
 вполне приводимы, $$\mathbbm{k}S_n b_{T_\mu}\bar f \cong \mathbbm{k}S_n e_{T_{\lambda^{(1)}}}
 \oplus \ldots \oplus \mathbbm{k}S_n e_{T_{\lambda^{(s)}}}$$
 для некоторых $\lambda^{(i)}\vdash n$ и некоторых таблиц Юнга $T_{\lambda^{(i)}}$
 формы $\lambda^{(i)}$, $1\leqslant i \leqslant s$, $s\in\mathbb N$.
В частности, $e^*_{T_{\lambda^{(1)}}} \mathbbm{k}S_n b_{T_\mu} \ne 0$.

Теперь заметим, что из $e^*_{T_{\lambda^{(1)}}} \mathbbm{k}S_n b_{T_\mu} \ne 0$ следует, что
$a_{T_{\lambda^{(1)}}} \sigma b_{T_\mu} = \sigma a_{\sigma^{-1}T_{\lambda^{(1)}}}  b_{T_\mu} \ne 0$
для некоторой подстановки
 $\sigma \in S_n$. В силу того, что $a_{\sigma^{-1}T_{\lambda^{(1)}}}$
является оператором симметризации
по числам из строк таблицы Юнга $\sigma^{-1}T_{\lambda^{(1)}}$,
а $b_{T_\mu}$ является оператором альтернирования 
по числам из столбцов таблицы Юнга $T_\mu$, все числа из первой строчки таблицы $\sigma^{-1}T_{\lambda^{(1)}}$
должны находиться в разных столбцах таблицы $T_\mu$. Отсюда $\left(\lambda^{(1)}\right)_1 \leqslant \mu_1$.
Более того, все числа из каждого из $\mu_r$ первых столбцов таблицы $T_\mu$
должны находиться в разных строчках таблицы $\sigma^{-1}T_{\lambda^{(1)}}$. Поскольку в силу леммы~\ref{LemmaInequalityLambdaUpperFrac}
справедливо равенство $\left(\lambda^{(1)}\right)_{r+1}=0$,
получаем отсюда, что $\left(\lambda^{(1)}\right)_r \geqslant \mu_r$.

Теперь из~(\ref{EquationSumMuGamma2FracM2})
и того, что $\mu_1-\mu_i \leqslant 2(r+1)$
 для всех $2\leqslant i \leqslant |T_0|+|T_1|-|\bar t_0|$,
следует, что
\begin{equation*}\begin{split}\sum_{i=1}^r \gamma_i \left(\lambda^{(1)}\right)_i \geqslant \sum_{i=1}^r \gamma_i \left(\lambda^{(1)}\right)_i
- \sum_{i=1}^r \gamma_i \mu_i - 4r =\\
 =
\sum_{i=1}^{|T_0|+|T_1|-|\bar t_0|} \left(\mu_i-\left(\lambda^{(1)}\right)_i\right)+\sum_{i=2|T_0|+3|T_1|-|\bar t_0|+1}^r \left(\left(\lambda^{(1)}\right)_i-\mu_i\right)
-4r \geqslant \\ \geqslant\sum_{i=1}^{|T_0|+|T_1|-|\bar t_0|} \left(\mu_1-\left(\lambda^{(1)}\right)_i\right)+\sum_{i=2|T_0|+3|T_1|-|\bar t_0|+1}^r \left(\left(\lambda^{(1)}\right)_i-\mu_r\right)
-6r-2r^2\end{split}
\end{equation*}  

Поскольку в силу леммы~\ref{LemmaInequalityLambdaUpperFrac} справедливо неравенство $\sum_{i=1}^r \gamma_i \left(\lambda^{(1)}\right)_i \leqslant 1$ и оба числа $\left(\mu_1-\left(\lambda^{(1)}\right)_i\right)$
и $\left(\left(\lambda^{(1)}\right)_i-\mu_r\right)$ неотрицательные, получаем, что
 $\mu_1 - (2r^2+6r+1) \leqslant \left(\lambda^{(1)}\right)_i \leqslant \mu_1$
для всех $1\leqslant i \leqslant |T_0|+|T_1|-|\bar t_0|$
и
$\mu_r \leqslant \left(\lambda^{(1)}\right)_i \leqslant \mu_r+(2r^2+6r+1)$
для всех $2|T_0|+3|T_1|-|\bar t_0|+1\leqslant i \leqslant r$.

Напомним, что из $a_{\sigma^{-1}T_{\lambda^{(1)}}}  b_{T_\mu} \ne 0$ следует,
что все числа из любого столбца таблицы $T_\mu$ находятся в разных строчках
таблицы $\sigma^{-1}T_{\lambda^{(1)}}$.
Применяя это свойство к первым $\mu_{2|T_0|+3|T_1|-|\bar t_0|}$ столбцам,
получаем, что в последних  $|\bar t_0|+|T_1|+1$
строчках таблицы $T_{\lambda^{(1)}}$ находятся как минимум $\sum_{i=2|T_0|+3|T_1|-|\bar t_0|}^r\mu_i$
клеток и $$\sum_{i=2|T_0|+3|T_1|-|\bar t_0|}^r \left(\lambda^{(1)}\right)_i
\geqslant \sum_{i=2|T_0|+3|T_1|-|\bar t_0|}^r\mu_i=(|\bar t_0|+|T_1|)\mu_r
+ \mu_{2|T_0|+3|T_1|-|\bar t_0|}.$$
Следовательно, \begin{equation*}\begin{split}
\lambda^{(1)}_{2|T_0|+3|T_1|-|\bar t_0|} \geqslant(|\bar t_0|+|T_1|)\mu_r + \mu_{2|T_0|+3|T_1|-|\bar t_0|} - \sum_{i=2|T_0|+3|T_1|-|\bar t_0|+1}^r \left(\lambda^{(1)}\right)_i \geqslant \\
\geqslant \mu_{2|T_0|+3|T_1|-|\bar t_0|} - (2r^3+6r^2+r)
\end{split}
\end{equation*}

Пусть $\lambda=(\lambda_1, \lambda_2, \ldots, \lambda_r)$,
где $$ \lambda_i=\left\lbrace\begin{array}{rrr}
\mu_1 - (2r^2+6r+1) & \text{ при } & 1\leqslant i \leqslant |T_0|+|T_1|-|\bar t_0|,\\
\mu_{2|T_0|+3|T_1|-|\bar t_0|} - (2r^3+6r^2+r) & \text{ при } & |T_0|+|T_1|-|\bar t_0|+1 \leqslant i \leqslant 2|T_0|+3|T_1|-|\bar t_0|, \\
\mu_r & \text{ при } & 2|T_0|+3|T_1|-|\bar t_0|+1\leqslant i \leqslant r.
\end{array} \right.$$ Введем обозначение $n_1 :=\sum_{i=1}^r \lambda_i$.
Тогда $n - (2r^4+6r^3+r^2) \leqslant n_1 \leqslant n$.

Для любого $\varepsilon > 0$ существует такое $n_0\in\mathbb N$,
что для любого $n\geqslant n_0$ справедливы неравенства $\lambda_1\geqslant \ldots \geqslant \lambda_r$ 
и $\Phi\left(\frac{\lambda_1}{n_1},\ldots,\frac{\lambda_r}{n_1}\right) > d-\varepsilon$.
Поскольку $D_\lambda$ является поддиаграммой диаграммы Юнга $D_{\lambda^{(1)}}$,
справедливо неравенство $c^{\mathbbm{k}^T}_n(A)\geqslant \dim M_\lambda$
и в силу формулы крюков и формулы Стирлинга существуют
такие $C_1 > 0$ и $r_1\in\mathbb R$, что
 \begin{equation}\begin{split} c^{\mathbbm{k}^T}_n(A) \geqslant \dim M(\lambda) = \frac{n_1!}{\prod_{i,j} h_{ij}}
  \geqslant \frac{n_1!}{(\lambda_1+r-1)! \ldots (\lambda_r+r-1)!} \geqslant \\
  \geqslant \frac{n_1!}{n_1^{r(r-1)}\lambda_1! \ldots \lambda_r!} \geqslant
  \frac{C_1 n_1^{r_1} 
\left(\frac{n_1}{e}\right)^{n_1}}{\left(\frac{\lambda_1}{e}\right)^{\lambda_1}\ldots
\left(\frac{\lambda_r}{e}\right)^{\lambda_r}} \geqslant \\ \geqslant C_1 n_1^{r_1}\left(\frac{1}
{\left(\frac{\lambda_1}{n_1}\right)^{\frac{\lambda_1}{n_1}}\ldots
\left(\frac{\lambda_r}{n_1}\right)^{\frac{\lambda_r}{n_1}}}\right)^{n_1} \geqslant C_1 n_1^{r_1}
 (d-\varepsilon)^{n - (2r^4+6r^3+r^2)}.\end{split}\end{equation}
 Следовательно, $\mathop{\underline\lim}_{n\rightarrow\infty}\sqrt[n]{c_n^{\mathbbm{k}^T}(A)}
\geqslant d-\varepsilon$. Из произвольности числа $\varepsilon > 0$ получаем, что $\mathop{\underline\lim}_{n\rightarrow\infty}\sqrt[n]{c_n^{\mathbbm{k}^T}(A)}
\geqslant d$. В силу предложения~\ref{PropositionCnGrCnGenH} и~(\ref{EquationMaximumGammaOmegaFracM2}) 
$$\mathop{\underline\lim}\limits_{n\rightarrow\infty}
 \sqrt[n]{c_n^{T\text{-}\mathrm{gr}}(A)}\geqslant  |T_0|+2|T_1|
  + 2\sqrt{(|T_1|+|\bar t_0|)(|T_0|+|T_1|-|\bar t_0|)}.$$

Применяя теорему~\ref{TheoremUpperFrac} и~(\ref{EquationMaximumGammaOmegaFracM2}), получаем, что
$$\mathop{\overline\lim}\limits_{n\rightarrow\infty}
 \sqrt[n]{c_n^{T\text{-}\mathrm{gr}}(A)}\leqslant  |T_0|+2|T_1|
  + 2\sqrt{(|T_1|+|\bar t_0|)(|T_0|+|T_1|-|\bar t_0|)}.$$
  
Из условия $|\bar t_0| > \frac{|T_0|}2$
следует, что
 $|T_1|+|\bar t_0| > |T_0|+|T_1|-|\bar t_0|$,
$$2\sqrt{(|T_1|+|\bar t_0|)(|T_0|+|T_1|-|\bar t_0|)} < (|T_1|+|\bar t_0|)+(|T_0|+|T_1|-|\bar t_0|)$$
  и $$\lim\limits_{n\rightarrow\infty}
 \sqrt[n]{c_n^{T\text{-}\mathrm{gr}}(A)} = |T_0|+2|T_1|
  + 2\sqrt{(|T_1|+|\bar t_0|)(|T_0|+|T_1|-|\bar t_0|)} < 2|T_0|+4|T_1|=\dim A.$$
\end{proof}
\begin{remark}\label{RemarkGrPIexpNonTriangleFracM2}
Конечномерная $T$-градуированно простая алгебра $A$
с заданными $|T_0|$, $|T_1|$ и $\frac{|T_0|}{2} < |\bar t_0| \leqslant |T_0|$ существует в силу предложения~\ref{TheoremImagesOfGradedComponentsReesExistence}.
\end{remark}
\begin{remark}
Алгебра $A$ из теоремы~\ref{TheoremGrPIexpNonTriangleFracM2}
не содержит единицу, поскольку иначе в силу предложения~\ref{PropositionTIdemGradedIdeals}
существовал бы изоморфизм $A \cong M_2(\mathbbm{k})$,
откуда в силу предложения~\ref{PropositionOrdinaryAndGradedCodim} и теоремы~\ref{TheoremGradedSimplicityCriterion} (применённой для тривиальной группы $G$) было бы справедливо равенство $$\PIexp^{T\text{-}\mathrm{gr}}(A)=\
\PIexp(A)=\dim A=4.$$ Вообще, как мы увидим чуть ниже в теореме~\ref{TheoremTIdemAmitsur}, если алгебра, градуированная лентой правых нулей является алгеброй с единицей, её градуированная PI-экспонента всегда целая и совпадает с её обычной PI-экспонентой.
\end{remark}

Напомним, что в \S\ref{SectionSemigroupTwoElements}
через $Q_3$ была обозначена лента правых нулей, состоящая из двух элементов.
Посльзуясь теоремой~\ref{TheoremGrPIexpNonTriangleFracM2}, мы можем построить пример такой конечномерной $Q_3$-простой алгебры $A$, что $\PIexp^{Q_3\text{-}\mathrm{gr}}(A) < \dim A$ и не является целым числом:

\begin{theorem}\label{TheoremQ3GradFractPI}
Пусть $\mathbbm{k}$~"--- поле характеристики $0$.
Обозначим через $I$ неприводимый левый $M_2(\mathbbm{k})$-модуль,
изоморфный минимальному левому идеалу $\langle e_{12}, e_{22}\rangle_\mathbbm{k} \subset M_2(\mathbbm{k})$.
Пусть $A = M_2(\mathbbm{k})\oplus I$ (прямая сумма левых идеалов), где $IM_2(\mathbbm{k}) := 0$, а $I^2:=0$.
Определим на $A$ градуировку полугруппой $Q_3$ при помощи равенств $A^{(e_1)}=(M_2(\mathbbm{k}),0)$, $A^{(e_2)}=\lbrace (\varphi(a), a) \mid a \in I \rbrace$, где $\varphi \colon I \hookrightarrow M_2(\mathbbm{k})$~"--- естественное вложение, которое является
гомоморфизмом $ M_2(\mathbbm{k})$-модулей. Тогда алгебра $A$ является $Q_3$-градуированно простой и существует  $$\PIexp^{Q_3\text{-}\mathrm{gr}}(A)=\lim\limits_{n\to \infty} \sqrt[n]{c_n^{Q_3\text{-}\mathrm{gr}}(A)} =3+2\sqrt 2= 5{,}8284\ldots$$
\end{theorem}

\begin{remark}\label{RemarkQ3OpGradFractPI}
Аналогично, $\lim\limits_{n\to \infty} \sqrt[n]{c_n^{Q_3^{\,\mathrm{op}}\text{-}\mathrm{gr}}(A^{\,\mathrm{op}})} =3+2\sqrt 2= 5{,}8284\ldots$
\end{remark}
\begin{proof}[Доказательство теоремы~\ref{TheoremQ3GradFractPI}]
Сперва заметим, что $J(A)=I$.
Предположим, что $W\ne 0$~"--- некоторый градуированный идеал алгебры $A$.
Тогда существует ненулевой однородный элемент $(a_1, b_1) \in W$, где $a_1 \in M_2(\mathbbm{k})$, а $b_1\in I$.
Поскольку $I$ не содержит ненулевых однородных элементов, выполнено условие $a_1 \ne 0$.
Тогда $(a_1,b_1)(E,0) = (a_1,0) \in (M_2(\mathbbm{k}),0) \cap W$.
Поскольку $W$ является двусторонним идеалом, а алгебра $M_2(\mathbbm{k})$~"--- проста, справедливо включение $(M_2(\mathbbm{k}),0) \subseteq W$. Более того, $(0,I) = (M_k(\mathbbm{k}),0) (0,I) \subseteq W$. Отсюда $W = A$, и алгебра $A$
действительно $Q_3$-градуированно проста.

Определяя множества $T_0$ и $T_1$ так, как было указано в начале \S\ref{SectionFracM2Equal},
получаем, что $T_0=\lbrace e_2 \rbrace$, $T_1=\lbrace e_1 \rbrace$.
Отсюда согласно теореме~\ref{TheoremGrPIexpNonTriangleFracM2} справедливо равенство
$$\PIexp^{Q_3\text{-}\mathrm{gr}}(A)=\lim\limits_{n\to \infty} \sqrt[n]{c_n^{Q_3\text{-}\mathrm{gr}}(A)} =3+2\sqrt 2= 5{,}8284\ldots$$
\end{proof}

\section{$Q_1$- и $Q_2$-градуированные алгебры с нецелой градуированной PI-экспонентой}

Напомним, что в \S\ref{SectionSemigroupTwoElements} были введены полугруппы~$Q_1=(\lbrace 0, 1 \rbrace, \cdot)$ и~$Q_2=\lbrace 0, v \mid v^2 = 0\rbrace$.
В данном параграфе строятся примеры конечномерных $Q_1$- и $Q_2$-градуированно простых алгебр с дробной
градуированной PI-экспонентой.

%

\begin{theorem}\label{TheoremQ1GradFractPI}
Пусть $A = M_2(\mathbbm{k})\oplus \UT_2(\mathbbm{k})$ (прямая сумма идеалов), где $\mathbbm{k}$~"--- поле характеристики $0$. 
Зададим на $A$ градуировку полугруппой $Q_1$ равенствами $A^{(0)}=(M_2(\mathbbm{k}),0)$ и $A^{(1)}=\lbrace (\varphi(a), a) \mid a \in \UT_2(\mathbbm{k}) \rbrace$, где $\varphi \colon \UT_2(\mathbbm{k}) \hookrightarrow M_2(\mathbbm{k})$~"---
естественное вложение. Другими словами, $A$ является алгеброй из примера~\ref{ExampleQ1} при $k=2$. Тогда существует  $$\lim\limits_{n\to \infty} \sqrt[n]{c_n^{Q_1\text{-}\mathrm{gr}}(A)} =4+2\sqrt 2= 6{,}8284\ldots$$
\end{theorem}

Для того, чтобы доказать теорему~\ref{TheoremQ1GradFractPI},
нам потребуется следующая лемма.  Будем для краткости опускать $\varphi$
и писать $(e_{ij}, e_{ij})$ вместо $(e_{ij}, \varphi^{-1}(e_{ij}))$.
\begin{lemma}\label{LemmaAltQ1}
Пусть $\lambda \vdash n$, $n\in\mathbb N$, $\lambda_8 = 0$, а $\lambda_6+\lambda_7 \leqslant \lambda_1$.
Тогда $m(A,\mathbbm{k}^{Q_1},\lambda) \ne 0$.
\end{lemma}
\begin{proof}
Достаточно доказать, что для некоторых $f\in P^{\mathbbm{k}^{Q_1}}_n$ и $T_\lambda$
выполнено условие $e_{T_\lambda}f \not\equiv 0$ на $A$.

Заметим, что $e_{T_\lambda}f$ кососимметричен по $\lambda_7$ попарно непересекающимся множествам переменных,
каждое из которых состоит из $7$ переменных. Из леммы~\ref{LemmaThetaCondition}
следует, что каждый такой столбец будет вносить вклад как минимум $1$ в сумму значений функции $\theta$
на элементах, которые подставляются вместо переменных множества $e_{T_\lambda}f$.
Следовательно, это вклад нужно будет компенсировать.

Пусть $\beta_2=\lambda_6-\lambda_7$. Выберем произвольные числа $\beta_3,\ldots,\beta_{12} \geqslant 0$,
такие, что $$\beta_3+\beta_5+\beta_7+\beta_9+\beta_{11} = \lambda_7,\quad
\beta_3+\beta_4=\lambda_5-\lambda_6,\quad \beta_5+\beta_6=\lambda_4-\lambda_5,$$
$$\beta_7+\beta_8=\lambda_3-\lambda_4,\quad \beta_9+\beta_{10}=\lambda_2-\lambda_3\text{\quad и\quad}\beta_{11}+\beta_{12}=\lambda_1-\lambda_2.$$ Другими словами,
$$D_\lambda=\begin{array}{|c|c|c|c|c|c|c|c|c|c|c|c|}
\multicolumn{1}{c}{\lambda_7} & \multicolumn{1}{c}{\beta_2} & \multicolumn{1}{c}{\beta_3} & \multicolumn{1}{c}{\beta_4} & \multicolumn{1}{c}{\beta_5} & \multicolumn{1}{c}{\beta_6} & \multicolumn{1}{c}{\beta_7} & \multicolumn{1}{c}{\beta_8} & \multicolumn{1}{c}{\beta_9} &
\multicolumn{1}{c}{\beta_{10}} & \multicolumn{1}{c}{\beta_{11}} & \multicolumn{1}{c}{\beta_{12}} \\
\hline
 \ldots & \ldots & \ldots & \ldots & \ldots & \ldots & \ldots & \ldots & \ldots & \ldots & \ldots & \ldots \\
 \cline{1-12}
 \ldots & \ldots & \ldots & \ldots & \ldots & \ldots & \ldots & \ldots & \ldots & \ldots \\
 \cline{1-10}
 \ldots & \ldots & \ldots & \ldots & \ldots & \ldots & \ldots & \ldots \\
 \cline{1-8}
 \ldots & \ldots & \ldots & \ldots & \ldots & \ldots \\
 \cline{1-6}
 \ldots & \ldots & \ldots & \ldots \\
 \cline{1-4}
 \ldots & \ldots \\
 \cline{1-2}
 \ldots  \\
 \cline{1-1}
\end{array}.$$
(Здесь числа $\beta_i$ означают число столбцов в каждом блоке.)

Каждый из $\lambda_7$ первых столбцов будет вносить вклад $1$ в $\theta$,
который будет компенсироваться столбцами, помеченными числами $\beta_3$, $\beta_5$, $\beta_7$, $\beta_9$
и $\beta_{11}$. Столбцы, помеченные $\beta_2$, $\beta_3$, $\beta_4$, $\beta_6$, $\beta_8$, $\beta_{10}$ и
 $\beta_{12}$ будут давать в $\theta$ нулевой вклад.

Выберем некоторую таблицу Юнга $T_\lambda$ формы $\lambda$,
заполненную числами от $1$ до $n$.
Для каждого столбца таблицы $T_\lambda$
определим полилинейный кососимметричный $\mathbbm{k}^{Q_1}$-многочлен, зависящий от переменных с индексами из этого столбца.
Для краткости будем обозначать $\mathbbm{k}^{Q_1}$-многочлены,
отвечающие столбцам блока с номером $i$
одним и тем же символом $f_i$.
Через $(i_1, \ldots, i_\ell)$ 
обозначим набор, состоящий из чисел, записанных в клетках соответствующего столбца, взятых сверху вниз.
Обозначим через $S\lbrace i_1, \ldots, i_\ell\rbrace$ группу подстановок на множестве $i_1, \ldots, i_\ell$.
Положим $$f_1 := \sum_{\sigma\in S\lbrace i_1, \ldots, i_7\rbrace} (\sign \sigma)
x^{q_0}_{\sigma(i_3)}
x^{q_1}_{\sigma(i_2)}
x^{q_0}_{\sigma(i_6)}
x^{q_1}_{\sigma(i_4)}
x^{q_0}_{\sigma(i_5)}
x^{q_0}_{\sigma(i_1)}
x^{q_1}_{\sigma(i_7)},
$$ 
$$f_2 := \sum_{\sigma\in S\lbrace i_1, \ldots, i_6\rbrace} (\sign \sigma)
x^{q_0}_{\sigma(i_3)}
x^{q_1}_{\sigma(i_2)}
x^{q_0}_{\sigma(i_6)}
x^{q_1}_{\sigma(i_4)}
x^{q_0}_{\sigma(i_5)}
x^{q_0}_{\sigma(i_1)}
,$$ 
$$f_3 := \sum_{\sigma\in S\lbrace i_1, \ldots, i_5\rbrace} (\sign \sigma)
x^{q_0}_{\sigma(i_5)}
x^{q_1}_{\sigma(i_4)}
x^{q_0}_{\sigma(i_1)}
x^{q_1}_{\sigma(i_2)}
x^{q_0}_{\sigma(i_3)}
,$$
$$f_4 := \sum_{\sigma\in S\lbrace i_1, \ldots, i_5\rbrace} (\sign \sigma)
x^{q_1}_{\sigma(i_2)}
x^{q_0}_{\sigma(i_3)}
x^{q_1}_{\sigma(i_5)}
x^{q_1}_{\sigma(i_4)}
x^{q_0}_{\sigma(i_1)}
,$$
$$f_5 := \sum_{\sigma\in S\lbrace i_1, \ldots, i_4\rbrace} (\sign \sigma)
x^{q_1}_{\sigma(i_4)}
x^{q_0}_{\sigma(i_1)}
x^{q_1}_{\sigma(i_2)}
x^{q_0}_{\sigma(i_3)}
,$$ $$f_6 := \sum_{\sigma\in S\lbrace i_1, \ldots, i_4\rbrace} (\sign \sigma)
x^{q_1}_{\sigma(i_2)}
x^{q_1}_{\sigma(i_3)}
x^{q_1}_{\sigma(i_4)}
x^{q_0}_{\sigma(i_1)}
,$$ $$f_7 := \sum_{\sigma\in S\lbrace i_1, i_2, i_3\rbrace} (\sign \sigma)
x^{q_0}_{\sigma(i_1)}
x^{q_1}_{\sigma(i_2)}
x^{q_0}_{\sigma(i_3)}
,\qquad f_8 := \sum_{\sigma\in S\lbrace i_1, i_2, i_3\rbrace} (\sign \sigma)
x^{q_1}_{\sigma(i_2)}
x^{q_1}_{\sigma(i_3)}
x^{q_0}_{\sigma(i_1)}
,$$ $$f_9 := \sum_{\sigma\in S\lbrace i_1, i_2\rbrace} (\sign \sigma)
x^{q_0}_{\sigma(i_1)}
x^{q_1}_{\sigma(i_2)}
,\qquad f_{10} := \sum_{\sigma\in S\lbrace i_1, i_2\rbrace} (\sign \sigma)
x^{q_0}_{\sigma(i_2)}
x^{q_0}_{\sigma(i_1)}
,$$ $$f_{11} := x^{q_0}_{i_1},\qquad f_{12} := x^{q_1}_{i_1}.$$

Теперь определим следующий $\mathbbm{k}^{Q_1}$-многочлен:
$$f=(f_1 f_3)^{\beta_3}(f_1 f_5)^{\beta_5}(f_1 f_7)^{\beta_7}(f_1 f_9)^{\beta_9}(f_1 f_{11})^{\beta_{11}}
f_2^{\beta_2} f_4^{\beta_4}f_6^{\beta_6}f_8^{\beta_8}f_{10}^{\beta_{10}}f_{12}^{\beta_{12}} \in P_n.$$
Как было уже отмечено, различные копии $\mathbbm{k}^{Q_1}$-многочлена $f_i$ зависят от разных переменных.

Копии $\mathbbm{k}^{Q_1}$-многочлена $f_1$ являются кососимметричными $\mathbbm{k}^{Q_1}$-многочленами 
степени $7$, отвечающими первым $\lambda_7$ столбцам высоты $7$.

Копии $\mathbbm{k}^{Q_1}$-многочлена  $f_2$ являются кососимметричными $\mathbbm{k}^{Q_1}$-многочленами 
степени $6$, отвечающими следующим $\beta_2$ столбцам высоты $6$.

\ldots

Копии $\mathbbm{k}^{Q_1}$-многочлена $f_{12}$ являются $\mathbbm{k}^{Q_1}$-многочленами 
степени $1$, индексы переменных которых принадлежат последним $\beta_{12}$ столбцам высоты $1$.

Докажем, что $e_{T_\lambda}f \not\equiv 0$. Для того, чтобы это проверить, заполним диаграмму $D_\lambda$
определёнными однородными элементами и обозначим получившуюся таблицу через $\tau$. (См. рис.~\ref{FigureTauQ1}.)
\begin{landscape}

\begin{figure}\caption{Подстановка элементов алгебры $A$ вместо переменных многочлена $e_{T_\lambda}f$.}\label{FigureTauQ1}
$$\tau=\begin{array}{|c|c|c|c|c|c|c|c|c|c|c|c|}
\multicolumn{1}{c}{\lambda_7} & \multicolumn{1}{c}{\beta_2} & \multicolumn{1}{c}{\beta_3} & \multicolumn{1}{c}{\beta_4} & \multicolumn{1}{c}{\beta_5} & \multicolumn{1}{c}{\beta_6} & \multicolumn{1}{c}{\beta_7} & \multicolumn{1}{c}{\beta_8} & \multicolumn{1}{c}{\beta_9} &
\multicolumn{1}{c}{\beta_{10}} & \multicolumn{1}{c}{\beta_{11}} & \multicolumn{1}{c}{\beta_{12}} \\
\hline
 (e_{21},0) & (e_{21},0) & (e_{21},0) & (e_{21},0) & (e_{21},0) & (e_{21},0) & (e_{21},0) & (e_{21},0) & (e_{21},0) & (e_{21},0) & (e_{21},0) & (e_{11}, e_{11}) \\
 \cline{1-12}
 (e_{11}, e_{11}) & (e_{11}, e_{11}) & (e_{11}, e_{11}) & (e_{11}, e_{11}) & (e_{11}, e_{11}) & (e_{11}, e_{11}) & (e_{11}, e_{11}) & (e_{11}, e_{11}) & (e_{11}, e_{11}) & (e_{12}, 0) \\
 \cline{1-10}
 (e_{11}, 0) & (e_{11}, 0) & (e_{11}, 0) & (e_{11}, 0) & (e_{11}, 0) & (e_{12}, e_{12}) & (e_{11}, 0) & (e_{12}, e_{12}) \\
 \cline{1-8}
 (e_{22}, e_{22}) & (e_{22}, e_{22}) & (e_{22}, e_{22}) & (e_{22}, e_{22}) & (e_{22}, e_{22}) & (e_{22}, e_{22}) \\
 \cline{1-6}
 (e_{22}, 0) & (e_{22}, 0) & (e_{22}, 0) & (e_{12}, e_{12}) \\
 \cline{1-4}
 (e_{12}, 0) & (e_{12}, 0) \\
 \cline{1-2}
 (e_{12}, e_{12})  \\
 \cline{1-1}
\end{array}$$ 
\end{figure}

 (Здесь в $i$-м блоке имеется $\beta_i$ столбцов с одними и теми же значениями.
Для краткости запишем всякое значение для каждого блока один раз.
При этом таблица $\tau$ остаётся таблицей формы $\lambda$.)

\end{landscape}

Теперь подставим в каждую переменную значение из соответствующей клетки
таблицы~$\tau$. При такой подстановке $\mathbbm{k}^{Q_1}$-многочлен $f$ не обращается 
в нуль.

Напомним, что $e_{T_\lambda} = a_{T_\lambda}b_{T_\lambda}$,
где $a_{T_\lambda}$~"--- оператор симметризации по переменным каждой строчки,
а $b_{T_\lambda}$~"--- оператор альтернирования по переменным каждого столбца.
Поскольку $\mathbbm{k}^{Q_1}$-многочлены $f_i$ кососимметричные, $\mathbbm{k}^{Q_1}$-многочлен $b_{T_\lambda} f$
равен $\mathbbm{k}^{Q_1}$-многочлену $f$, умноженному на ненулевой коэффициент.

Второй строчке таблицы $T_\lambda$ отвечает два подмножества переменных.
В переменные первого подмножества подставляются переменные $(e_{11},e_{11})  \in A^{(1)}$,
а в переменные второго подмножества~"--- переменные $(e_{12},0) \in A^{(0)}$.
Отсюда, если некоторое слагаемое элемента $a_{T_\lambda}$
перемешивает переменные из этих двух подмножеств,
по крайней мере одна переменная из второй группы (она входит в $f_{10}$)
заменяется на переменную из первой группы. Однако $\mathbbm{k}^{Q_1}$-многочлен $f_{10}$
обращается в нуль, если вместо хотя бы одной из его переменных подставляется
элемент из $A^{(1)}$, поскольку к обеим переменным $\mathbbm{k}^{Q_1}$-многочлена 
применяется оператор $q_0$.
Отсюда все слагаемые $\mathbbm{k}^{Q_1}$-многочлена $a_{T_\lambda}b_{T_\lambda}f$,
в которых переменные из этих двух групп перемешиваются друг с другом, обращаются при подстановке, заданной таблицей $\tau$, в нуль.

Отсюда, если некоторое слагаемое элемента $a_{T_\lambda}$ 
переводит одну из переменных первых двух столбцов в переменную с другим значением из таблицы $\tau$,
вместо переменных $\mathbbm{k}^{Q_1}$-многочленов $f_1$ и $f_2$ подставляется
слишком много элементов из
 $A^{(1)}$ и результат такой подстановки равен нулю из-за действия оператора $q_0$.
 Следовательно, все слагаемые $\mathbbm{k}^{Q_1}$-многочлена  $a_{T_\lambda}b_{T_\lambda}f$,
 в которых переменные с разными значениями меняются местами, обращаются в нуль. 
Продолжая эту процедуру, в конце концов получаем,
что если слагаемое элемента $a_{T_\lambda}$
не переводит множества переменных с одними и теми же значениями из таблицы $\tau$ в себя,
соответствующее слагаемое в $a_{T_\lambda}b_{T_\lambda}f$ обращается в нуль.
Отсюда значение $\mathbbm{k}^{Q_1}$-многочлена $a_{T_\lambda}b_{T_\lambda}f$ 
пропорционально значению $\mathbbm{k}^{Q_1}$-многочлена $b_{T_\lambda}f$ с ненулевым коэффициентом, т.е. также не равно нулю. Лемма доказана.
\end{proof}

\begin{proof}[Доказательство теоремы~\ref{TheoremQ1GradFractPI}.]
В обозначениях \S\ref{SectionUpperSGGr} $$\gamma_1 = -1,\ \gamma_2=\gamma_3=\gamma_4=\gamma_5=0,\ \gamma_6=\gamma_7=1,\ r=7,$$ откуда 
в силу леммы~\ref{LemmaInequalityLambdaUpperFrac}
для всех $\lambda\vdash n$, где $n\in\mathbb N$, таких, что $m(A, \mathbbm{k}^{Q_1}, \lambda)\ne 0$,
выполнены условия $\lambda_8 = 0$ и $\lambda_6+\lambda_7 \leqslant \lambda_1 + 1$.
Кроме того, $\zeta=\frac{\sqrt 2}{2}$ и для
$$\Omega := \left\lbrace (\alpha_1, \ldots, \alpha_7)\in \mathbb R^7 \mathrel{\biggl|} \sum_{i=1}^7 \alpha_i = 1,\ 
\alpha_1 \geqslant \alpha_2 \geqslant \ldots \geqslant \alpha_7\geqslant 0,\ \alpha_7
+\alpha_6 \leqslant \alpha_1\right\rbrace$$
в силу леммы~\ref{LemmaMaximumFUpperFrac} 
справедливо равенство
$ d:=\max_{x\in \Omega} \Phi(x) = 4+2\sqrt 2= 6{,}8284\ldots$
Обозначим через $(\alpha_1, \ldots, \alpha_7) \in \Omega$ такую точку, что
$\Phi(\alpha_1, \ldots, \alpha_7)=d$.
Для любого $n\in\mathbb N$ зададим $\mu\vdash n$ при помощи равенств
$\mu_i = [\alpha_i n]$ при $2\leqslant i \leqslant 7$
и $\mu_1 =  n-\sum_{i=1}^7 \mu_i$.
Для любого $\varepsilon > 0$ существует такое $n_0\in\mathbb N$,
что для любого $n\geqslant n_0$
справедливо неравенство $\Phi\left(\frac{\mu_1}{n},\ldots,\frac{\mu_7}{n}\right) > d-\varepsilon$.
 В силу леммы~\ref{LemmaAltQ1} выполнено условие $m(A,\mathbbm{k}^{Q_1},\mu) \ne 0$.
 Применяя формулу крюков и формулу Стирлинга, получаем, что существуют
 такие $C_1 > 0$ и $r_1\in\mathbb R$, что
 \begin{equation}\begin{split} c^{\mathbbm{k}^{Q_1}}_n(A) \geqslant \dim M(\mu) = \frac{n!}{\prod_{i,j} h_{ij}}
  \geqslant \frac{n!}{(\mu_1+6)! \ldots (\mu_7+6)!} \geqslant \\ \geqslant
  \frac{n!}{n^{42}\mu_1! \ldots \mu_7!} \geqslant
  \frac{C_1 n^{r_1} 
\left(\frac{n}{e}\right)^n}{\left(\frac{\mu_1}{e}\right)^{\mu_1}\ldots
\left(\frac{\mu_7}{e}\right)^{\mu_7}}\geqslant \\ \geqslant C_1 n^{r_1}\left(\frac{1}
{\left(\frac{\mu_1}{n}\right)^{\frac{\mu_1}{n}}\ldots
\left(\frac{\mu_7}{n}\right)^{\frac{\mu_7}{n}}}\right)^n \geqslant C_1 n^{r_1}
 (d-\varepsilon)^n.\end{split}\end{equation}
 Отсюда $\mathop{\underline\lim}_{n\to\infty}\sqrt[n]{c_n^{\mathbbm{k}^{Q_1}}(A)}
\geqslant d-\varepsilon$. Поскольку число $\varepsilon > 0$ было выбрано произвольным,
справедливо неравенство $$\mathop{\underline\lim}_{n\to\infty}\sqrt[n]{c_n^{\mathbbm{k}^{Q_1}}(A)}
\geqslant d,$$
которое в силу предложения~\ref{PropositionCnGrCnGenH}
эквивалентно неравенству
$$\mathop{\underline\lim}_{n\to\infty}\sqrt[n]{c_n^{Q_1\text{-}\mathrm{gr}}(A)}
\geqslant d.$$
 Из теоремы~\ref{TheoremUpperFrac} следует неравенство $$\mathop{\overline\lim}_{n\to\infty}\sqrt[n]{c_n^{Q_1\text{-}\mathrm{gr}}(A)}
\leqslant d,$$
откуда $$\lim\limits_{n\to \infty} \sqrt[n]{c_n^{Q_1\text{-}\mathrm{gr}}(A)} =4+2\sqrt 2= 6{,}8284\ldots$$
\end{proof}

Приведём теперь пример $Q_2$-градуированной алгебры с нецелой градуированной
 PI-экспонентой:

\begin{theorem}\label{TheoremQ2GradFractPI}
Пусть $A_2 = M_2(\mathbbm{k})\oplus \mathbbm{k} j_{11} \oplus \mathbbm{k} j_{12} \oplus \mathbbm{k} j_{22}$ (прямая сумма идеалов), где
$\langle j_{11}, j_{12}, j_{22}\rangle_\mathbbm{k}^2=0$, а
 $\mathbbm{k}$~"--- поле характеристики $0$. 
Зададим на $A$ градуировку полугруппой $Q_2$ равенствами  $A_2^{(0)}=(M_2(\mathbbm{k}),0)$, $A_2^{(v)}=\langle (e_{11}, j_{11}), (e_{12}, j_{12}), (e_{22}, j_{22})\rangle$. Тогда существует  $$\lim\limits_{n\to \infty} \sqrt[n]{c_n^{Q_2\text{-}\mathrm{gr}}(A_2)} =
4+2\sqrt 2= 6{,}8284\ldots$$
\end{theorem}
\begin{proof}
Обозначим через $A$ алгебру из теоремы~\ref{TheoremQ1GradFractPI}. 
Зададим линейную биекцию
$\psi \colon A \mathrel{\widetilde{\rightarrow}} A_2$ формулами $\psi(e_{ij}, e_{k\ell})=(e_{ij}, j_{k\ell})$
и $\psi(e_{ij}, 0)=(e_{ij}, 0)$.
Тогда $\psi(A^{(0)})=A_2^{(0)}$ и $\psi(A^{(1)})=A_2^{(v)}$.
 Определим теперь изоморфизм алгебр $\Theta \colon \mathbbm{k}\langle X^{Q_1\text{-}\mathrm{gr}} \rangle
\to \mathbbm{k}\langle X^{Q_2\text{-}\mathrm{gr}} \rangle$ при помощи равенств
 $\Theta(x^{(0)}_i)=x^{(0)}_i$ и $\Theta(x^{(1)}_i)=x^{(v)}_i$.
Докажем, что \begin{equation}\label{EqIdQ1IdQ2}\Theta(\Id^{Q_1\text{-}\mathrm{gr}}(A))=\Id^{Q_2\text{-}\mathrm{gr}}(A_2).\end{equation}
Поскольку характеристика поля $\mathbbm{k}$ равна $0$, оба идеала $\Id^{Q_1\text{-}\mathrm{gr}}(A)$ и $\Id^{Q_2\text{-}\mathrm{gr}}(A_2)$
порождены как идеалы градуированных тождеств полилинейными многочленами.
(Доказательство этого факта повторяет доказательство теоремы~1.3.8 из~\cite{ZaiGia}.)
Другими словами, для того, чтобы доказать равенство~(\ref{EqIdQ1IdQ2}),
достаточно показать, что
если $f\in \mathbbm{k}\langle X^{Q_1\text{-}\mathrm{gr}} \rangle$ 
полилинеен как обычный многочлен от переменных
$x_1^{(t_1)}, x_2^{(t_2)},\ldots, x_n^{(t_n)}$, где $t_i\in Q_1$, то
$f \in \Id^{Q_1\text{-}\mathrm{gr}}(A)$, если и только если $\Theta(f)\in
\Id^{Q_2\text{-}\mathrm{gr}}(A_2)$.
Напомним, что для проверки градуированного тождества подставляются только однородные элементы.
Заметим, что
$\pi \psi(a) = \pi(a)$, где $\pi$~"--- проекция на первую компоненту: $\pi(a,b)=a$ для всех  $(a,b)\in A$ и $(a,b)\in A_2$. Более того, если $a \in A^{(0)} \cup A^{(1)}$, то $\pi(a)=0$, если и только если $a=0$. Поскольку обе полугруппы $Q_1$ и $Q_2$ коммутативны, значение градуированного многочлена $f$ 
при подстановке однородных элементов снова является однородным элементом.
Применяя $\pi$, получаем, что $f\in P_n^{Q_1\text{-}\mathrm{gr}}$
обращается в нуль при подстановке однородных элементов $a^{(t_i)}_i\in A^{(t_i)}$,  где $t_i\in Q_1$,
если и только если $\Theta(f)$ обращается в нуль при подстановке элементов $\psi\left(a^{(\tilde t_i)}_i\right)\in A_2^{(\tilde t_i)}$. (Здесь $\tilde 0 = 0$,
а $\tilde 1 = v$.) Следовательно, справедливо равенство~(\ref{EqIdQ1IdQ2}),
откуда
$$c_n^{Q_1\text{-}\mathrm{gr}}(A)=c_n^{Q_2\text{-}\mathrm{gr}}(A_2)\text{ для всех  }n\in\mathbb N.$$
Теперь достаточно применить теорему~\ref{TheoremQ1GradFractPI}.
\end{proof}

\section{Достаточные условия справедливости аналога гипотезы Амицура для полиномиальных $T$-градуированных тождеств}

Напомним, что под справедливостью аналога гипотезы Амицура понимается существование целой градуированной PI-экспоненты. Одно семейство алгебр, градуированных лентами правых нулей, с целой  градуированной PI-экспонентой уже было построено в теореме~\ref{TheoremGrPIexpTriangleFracM2}. Приведём теперь другие достаточные условия справедливости аналога гипотезы Амицура.

\begin{theorem}\label{TheoremTCancelAmitsur}
Пусть $A$~"--- конечномерная 
ассоциативная алгебра над полем $\mathbbm{k}$ характеристики $0$, градуированная некоторой полугруппой $T$
с сокращениями. 
Тогда \begin{enumerate}
\item либо существует такое $n_0$, что $c_n^{T\text{-}\mathrm{gr}}(A)=0$ при всех $n\geqslant n_0$;
\item либо существуют такие константы $C_1, C_2 > 0$, $r_1, r_2 \in \mathbb R$,
  $d \in \mathbb N$, что $$C_1 n^{r_1} d^n \leqslant c^{T\text{-}\mathrm{gr}}_n(A)
   \leqslant C_2 n^{r_2} d^n\text{ для всех }n \in \mathbb N.$$
\end{enumerate}
В частности, для любой такой алгебры $A$
   существует $\PIexp^{T\text{-}\mathrm{gr}}(A)\in\mathbb Z_+$ и, таким образом,
   справедлив аналог гипотезы Амицура.
   \end{theorem}
\begin{proof}
Как обычно (см., например, \cite[теорема~4.1.9]{ZaiGia}),
коразмерности не меняются при расширении основного поля,
поэтому без ограничения общности можно считать, что основное поле
$\mathbbm{k}$ алгебраически замкнуто.

Согласно следствию 4.1 из~\cite{KelarevBook} идеал $J(A)$ является градуированным.
В силу предложения~\ref{PropositionTCancelWedderburn}
факторалгебра 
 $A/J(A)$ раскладывается в прямую сумму $T$-градуированных идеалов, являющихся $T$-градуированно простыми алгебрами. Теперь достаточно применить теоремы~\ref{TheoremWeakWedderburnMalcevHRad}, \ref{TheoremAssocAlternateFinal} и~\ref{TheoremHmodHRadAmitsurPIexpHBdimB} и  предложение~\ref{PropositionCnGrCnGenH}.
 \end{proof}

Если конечномерная алгебра с единицей градуирована лентой левых или правых нулей, то её градуированная PI-экспонента совпадает с обычной:

\begin{theorem}\label{TheoremTIdemAmitsur}
Пусть $A$~"--- конечномерная 
ассоциативная алгебра с единицей над полем $\mathbbm{k}$ характеристики $0$, градуированная некоторой лентой $T$
левых или правых нулей. 
Тогда \begin{enumerate}
\item либо существует такое $n_0$, что $c_n^{T\text{-}\mathrm{gr}}(A)=0$ при всех $n\geqslant n_0$;
\item либо существуют такие константы $C_1, C_2 > 0$, $r_1, r_2 \in \mathbb R$, что $$C_1 n^{r_1} d^n \leqslant c^{T\text{-}\mathrm{gr}}_n(A)
   \leqslant C_2 n^{r_2} d^n\text{ для всех }n \in \mathbb N,$$
\end{enumerate}
где число $d=\PIexp(A)$ является обычной PI-экспонентой алгебры $A$.
В частности, для любой такой алгебры $A$
   справедлив аналог гипотезы Амицура.
\end{theorem}
\begin{proof}
Поскольку алгебра $A$ конечномерна, можно без ограничения общности
считать, что полугруппа $T$ конечна. Кроме того, можно снова предполагать
основное поле $\mathbbm{k}$ алгебраически замкнутым.

В силу предложения~\ref{PropositionTIdemGradedIdeals}
радикал Джекобсона $J(A)$ является градуированным идеалом.
 В силу теоремы~\ref{TheoremTIdemGradedWeddMalcev}
 можно выбрать такую максимальную градуированную подалгебру $B$, что $A=B\oplus J(A)$ (прямая сумма градуированных подпространств). Пусть $B=B_1 \oplus \ldots \oplus B_s$ (прямая сумма идеалов) для некоторых простых алгебр $B_i$.
В силу предложения~\ref{PropositionTIdemGradedIdeals} идеалы~$B_i$ градуированные, а из
 теорем~\ref{TheoremAssocAlternateFinal}, \ref{TheoremGenHAmitsurHWedederburn} и  предложения~\ref{PropositionCnGrCnGenH} следует, что
для
\begin{equation*}\begin{split}d:= \max\dim\left( B_{i_1}\oplus B_{i_2} \oplus \ldots \oplus B_{i_r}
 \mathbin{\Bigl|}  r \geqslant 1,\ B_{i_1} J(A) \,B_{i_2} J(A) \ldots B_{i_{r-1}} J(A)\,B_{i_r}\ne 0\right)\end{split}\end{equation*}
 справедливы неравенства из формулировки теоремы.
 При этом $d=\PIexp^{T\text{-}\mathrm{gr}}(A)= \PIexp(A)$,
 поскольку алгебры $B_i$ просты как обычные алгебры.
 \end{proof}
 \begin{remark}
 Равенство $\PIexp^{T\text{-}\mathrm{gr}}(A)=\PIexp(A)$ не означает равенства коразмерностей.
  Действительно, пусть  $k\in\mathbb N$, $k\geqslant 2$, $A=M_k(\mathbbm{k})$,
 $T=\lbrace t_1, \ldots, t_k\rbrace$, где $t_i t_\ell = t_\ell$
 для всех  $1\leqslant i,\ell \leqslant k$, а $A^{(t_i)}=\langle e_{1 i},\ldots, e_{k i}\rangle_\mathbbm{k}$.
  Тогда элементы $x_1^{(t_i)}$
  линейно независимы по модулю $\Id^{T\text{-}\mathrm{gr}}(A)$.
  Для того, чтобы в этом убедиться, достаточно подставить $x_1^{(t_i)} = e_{ii}$.
 Отсюда $c_1(A)=1 < c_1^{T\text{-}\mathrm{gr}}(A)=k$.
 \end{remark}

\newpage

\chapter*{Заключение}
\addcontentsline {toc}{chapter}{Заключение}

Перечислим основные результаты работы:

\begin{enumerate}

\item построен пример элементарной градуировки на матричной алгебре,
которая не реализуется конечной группой (теорема~\ref{TheoremRegradeElementaryOmega});

\item введено понятие эквивалентности (ко)действий алгебр Хопфа
на $\Omega$-алгебрах, получены достаточные условия существования $V$-универсальных
(ко)действующих биалгебр и алгебр Хопфа, а также установлена двойственность между ними (глава~\ref{ChapterOmegaAlg});

\item найдены приложения понятия эквивалентности действий алгебр Хопфа к полиномиальным тождествам
(\S\ref{SectionEquivApplToPolyIden});

\item получены достаточные условия $H$-коинвариантности радикала Джекобсона в $H$-комодульных ассоциативных алгебрах (теорема~\ref{TheoremRadicalHSubComod});

\item получены достаточные условия $H$-(ко)инвариантности радикалов в $H$-(ко)модульных алгебрах Ли
(теоремы~\ref{TheoremHModRadicalsLie} и~\ref{TheoremLieRadicalHSubComod});

\item доказан $H$-(ко)инвариантный аналог теоремы Леви (теоремы~\ref{TheoremHcoLevi} и~\ref{TheoremHLevi});

\item классифицированы конечномерные ассоциативные алгебры и алгебры Ли, простые по отношению к действию алгебр Тафта (теоремы~\ref{TheoremTaftSimpleNonSemiSimplePresent}, \ref{TheoremTaftSimpleNonSemiSimpleClassify}, \ref{TheoremTaftSimpleSemisimple} и \S\ref{SectionClassTaftSLieSS}--\ref{SectionClassTaftSLieNSS});

\item классифицированы конечномерные ассоциативные градуированно простые алгебры, градуированные конечными полугруппами с тривиальными максимальными подгруппами (\S\ref{SectionReesSemigroupSimpleDescription}--\ref{SectionTGradedReesExistence});

\item доказаны аналоги гипотезы Амицура для следующих классов конечномерных алгебр:

\begin{itemize}
\item  ассоциативных алгебр и алгебр Ли, градуированных произвольными группами (теоремы~\ref{TheoremMainGrAssoc} и~\ref{TheoremMainLieGr});

\item  ассоциативных алгебр c действием произвольной группы автоморфизмами и антиавтоморфизмами (теорема~\ref{TheoremGAssoc});

\item  ассоциативных алгебр c градуированным действием произвольной группы (теорема~\ref{TheoremTGrGActionAssoc});

\item  ассоциативных алгебр c действием произвольной алгебры Ли дифференцированиями (теорема~\ref{TheoremDiffAssoc});

\item  алгебр Ли c рациональным действием редуктивной аффинной алгебраической группы автоморфизмами и антиавтоморфизмами (теорема~\ref{TheoremMainLieGAffAlg});

\item  алгебр Ли c действием конечномерной полупростой алгебры Ли дифференцированиями (теорема~\ref{TheoremMainDiffLie});

\item  ассоциативных $H$-модульных алгебр в случае, когда алгебра Хопфа $H$ либо полупроста и конечномерна, либо получена при помощи (возможно, многократного) расширения Оре конечномерной полупростой алгебры Хопфа косопримитивными элементами (теорема~\ref{TheoremHOreAmitsur});

\item  ассоциативных $H$-модульных алгебр $A$ в случае, когда $H$~"--- произвольная алгебра Хопфа,
а радикал Джекобсона $J(A)$ является $H$-подмодулем (теорема~\ref{TheoremHmoduleAssoc});

\item $H$-модульных алгебр Ли в случае, когда алгебра Хопфа $H$ полупроста и конечномерна (теорема~\ref{TheoremMainLieHSS});

\item $H$-модульных алгебр Ли в случае, когда  $H$~"--- произвольная алгебра Хопфа,
а разрешимый радикал является $H$-инвариантным нильпотентным идеалом (теорема~\ref{TheoremMainLieNRSame});
\end{itemize}

\item построена серия конечномерных градуированно простых ассоциативных алгебр, градуированных лентами правых нулей, с дробной градуированной PI-экспонентой (теорема~\ref{TheoremGrPIexpNonTriangleFracM2} и замечание~\ref{RemarkGrPIexpNonTriangleFracM2});

\item доказано существование градуированной PI-экспоненты у любой конечномерной градуированно простой
(необязательно ассоциативной) алгебры, градуированной произвольным множеством (следствие~\ref{CorollaryGradedExistsExponent}), и $H$-PI-экспоненты у любой конечномерной $H$-простой
(необязательно ассоциативной) алгебры с обобщённым $H$-действием (теорема~\ref{TheoremHSimpleHPIexpHNAssoc}).

\end{enumerate}

Данные исследования могут получить следующее развитие.

Во-первых, не лишено интереса обобщение понятия эквивалентности для (ко)действий
бимоноидов и моноидов Хопфа в произвольных заплетённых моноидальных категориях. В данном случае
особенный интерес представляют условия, которые необходимо наложить на соответствующую
категорию, чтобы в ней существовали универсальные бимоноиды и моноиды Хопфа.

Во-вторых, остаётся недоказанным существание разложения конечномерной $H$-модульной
алгебры Ли, не содержащей ненулевых $H$-инвариантных нильпотентных идеалов, в прямую
сумму $H$-простых идеалов в случае, когда алгебра Хопфа $H$ является неполупростой, например, когда $H$~"--- алгебра Тафта.

Во-третьих, гипотеза Амицура~"--- Бахтурина для конечномерных
$H$-модульных ассоциативных алгебр по-прежнему остаётся открытой в случае произвольных
алгебр Хопфа $H$. 
В случае произвольных
алгебр Хопфа $H$ остаётся открытым и аналог гипотезы Амицура для конечномерных
$H$-модульных алгебр Ли.

В-четвертых, представляет интерес вопрос о существовании PI-экспоненты у конечномерных ассоциативных
алгебр и алгебр Ли с произвольной дополнительной структурой. 

В-пятых, до сих пор неизвестны примеры конечномерных ассоциативных алгебр и алгебр Ли с действием некоторой конечномерной алгебры Ли $\mathfrak g$ дифференцированиями (в силу теорем~\ref{TheoremAssDerPIexpEqual} и~\ref{TheoremMainDiffLie}  такая алгебра Ли $\mathfrak g$ не может быть полупростой), у которых экспонента роста коразмерностей дифференциальных тождеств не совпадает
с экспонентой роста коразмерностей обычных тождеств.

\end{document}